\documentclass[11pt]{amsbook}

\newtheorem{thm}{Theorem}[section]
\newtheorem*{thm*}{Theorem}

\newtheorem{claim}[thm]{Claim}

\newtheorem{cor}[thm]{Corollary}

\newtheorem{lem}[thm]{Lemma}
\newtheorem*{lem*}{Lemma}
\newtheorem{mainthm}{Theorem}
\newtheorem*{mainthm*}{Theorem}
\newtheorem{maincor}[mainthm]{Corollary}
\newtheorem{prop}[thm]{Proposition}

\theoremstyle{definition}

\newtheorem*{case*}{Case}
\newtheorem{conj}[thm]{Conjecture}

\newtheorem{defn}[thm]{Definition}
\newtheorem*{defn*}{Definition}
\newtheorem{exmp}[thm]{Example}
\newtheorem*{exmp*}{Example}

\newtheorem{hyp}[thm]{Hypothesis}

\newtheorem{step}{Step}\renewcommand{\thestep}{}

\theoremstyle{remark}

\newtheorem{case}{Case}\renewcommand{\thecase}{}

\newtheorem{rmk}[thm]{Remark}
\newtheorem*{rmk*}{Remark}

\makeatletter
\def\alphenumi{
  \def\theenumi{\alph{enumi}}
  \def\p@enumi{\theenumi}
  \def\labelenumi{(\@alph\c@enumi)}}
\makeatother

\makeatletter
\def\thecase{\@arabic\c@case}
\makeatother

\makeatletter
\def\thestep{\@arabic\c@step}
\makeatother

\newcount\hh
\newcount\mm
\mm=\time
\hh=\time
\divide\hh by 60
\divide\mm by 60
\multiply\mm by 60
\mm=-\mm
\advance\mm by \time
\def\hhmm{\number\hh:\ifnum\mm<10{}0\fi\number\mm}

\let\oldmarginpar\marginpar
\renewcommand\marginpar[1]{\-\oldmarginpar[\raggedleft\footnotesize #1]%
{\raggedright\footnotesize #1}}

\renewcommand\emptyset{\varnothing}

\newcommand\too{\longrightarrow}

\newcommand\CC{\mathbb{C}}

\newcommand\HH{\mathbb{H}}

\newcommand\KK{\mathbb{K}}

\newcommand\NN{\mathbb{N}}

\newcommand\PP{\mathbb{P}}

\newcommand\RR{\mathbb{R}}

\newcommand\VV{\mathbb{V}}

\newcommand\ZZ{\mathbb{Z}}

\newcommand\cA{{\mathcal{A}}}
\newcommand\cB{{\mathcal{B}}}

\newcommand\cF{{\mathcal{F}}}
\newcommand\cG{{\mathcal{G}}}
\newcommand\cH{{\mathcal{H}}}

\newcommand\cT{{\mathcal{T}}}

\newcommand\calV{{\mathcal{V}}}
\newcommand\cW{{\mathcal{W}}}

\newcommand\fA{{\mathfrak{A}}}

\newcommand\fB{{\mathfrak{B}}}

\newcommand\fF{{\mathfrak{F}}}
\newcommand\fg{{\mathfrak{g}}}

\newcommand\fm{{\mathfrak{m}}}

\newcommand\fN{{\mathfrak{N}}}

\newcommand\fV{{\mathfrak{V}}}

\newcommand\fW{{\mathfrak{W}}}

\newcommand\sA{{\mathscr{A}}}
\newcommand\sB{{\mathscr{B}}}
\newcommand\sC{{\mathscr{C}}}
\newcommand\sD{{\mathscr{D}}}
\newcommand\sE{{\mathscr{E}}}
\newcommand\sF{{\mathscr{F}}}
\newcommand\sG{{\mathscr{G}}}
\newcommand\sH{{\mathscr{H}}}
\newcommand\sI{{\mathscr{I}}}

\newcommand\sK{{\mathscr{K}}}
\newcommand\sL{{\mathscr{L}}}
\newcommand\sM{{\mathscr{M}}}
\newcommand\sN{{\mathscr{N}}}
\newcommand\sO{{\mathscr{O}}}
\newcommand\sP{{\mathscr{P}}}

\newcommand\sR{{\mathscr{R}}}
\newcommand\sS{{\mathscr{S}}}
\newcommand\sT{{\mathscr{T}}}
\newcommand\sU{{\mathscr{U}}}
\newcommand\sV{{\mathscr{V}}}
\newcommand\sW{{\mathscr{W}}}
\newcommand\sX{{\mathscr{X}}}
\newcommand\sY{{\mathscr{Y}}}
\newcommand\sZ{{\mathscr{Z}}}

\newcommand\bkappa{{\boldsymbol{\kappa}}}

\newcommand\bsigma{{\boldsymbol{\sigma}}}

\newcommand\bR{{\mathbf{R}}}

\newcommand\bx{{\mathbf{x}}}

\newcommand{\cov}{\nabla}

\newcommand\eps{\varepsilon}

\newcommand\la{\lambda}
\newcommand\La{\Lambda}

\newcommand\Om{\Omega}

\newcommand\gl{{\mathfrak{g}\mathfrak{l}}}

\newcommand\GL{\operatorname{GL}}
\newcommand\Or{\operatorname{O}}

\newcommand\SL{\operatorname{SL}}
\newcommand\SO{\operatorname{SO}}

\newcommand\SU{\operatorname{SU}}
\newcommand\U{\operatorname{U}}

\newcommand\less{\setminus}

\newcommand\ad{{\operatorname{ad}}}
\newcommand\Ad{{\operatorname{Ad}}}

\newcommand\Aut{\operatorname{Aut}}

\DeclareMathOperator{\BMO}{BMO}

\newcommand\Coker{\operatorname{Coker}}

\newcommand\CS{\operatorname{CS}}

\newcommand\diag{\operatorname{diag}}

\newcommand\Diff{\operatorname{Diff}}
\newcommand\diam{\operatorname{diam}}
\DeclareMathOperator{\mydirac}{\slashed{\partial}}
\newcommand\dist{\operatorname{dist}}
\newcommand\divg{\operatorname{div}}

\newcommand\End{\operatorname{End}}

\newcommand{\esssup}{\operatornamewithlimits{ess\ sup}}

\newcommand\grad{\operatorname{grad}}

\newcommand\Hess{\operatorname{Hess}}

\newcommand\Ind{\operatorname{Index}}
\DeclareMathOperator{\Inj}{Inj}

\newcommand\Imag{\operatorname{Im}}

\newcommand\Ker{\operatorname{Ker}}

\DeclareMathOperator{\Lip}{Lip}

\newcommand\Met{\operatorname{Met}}

\newcommand\Real{\operatorname{Re}}
\newcommand\Ric{\operatorname{Ric}}

\newcommand\Ran{\operatorname{Ran}}
\newcommand\rank{\operatorname{rank}}

\newcommand\red{\operatorname{red}}

\newcommand\Riem{\operatorname{Riem}}
\DeclareMathOperator{\Scal}{Scal}

\newcommand\supp{\operatorname{supp}}

\newcommand\Sym{\operatorname{Sym}}

\newcommand\tr{\operatorname{tr}}

\DeclareMathOperator{\VMO}{VMO}
\newcommand\vol{\operatorname{vol}}
\newcommand\Vol{\operatorname{Vol}}

\DeclareMathOperator{\YM}{\mathscr{Y}\!\!\mathscr{M}}

\newcommand\afortiori{{\emph{a fortiori }}}

\newcommand\apriori{{\emph{a priori }}}
\newcommand\Apriori{{\emph{A priori }}}

\newcommand\asd{{\mathrm{asd}}}

\newcommand\deRham{{\mathrm{de Rham}}}

\newcommand\id{{\mathrm{id}}}

\newcommand\loc{{\mathrm{loc}}}

\newcommand\mutatis{{\emph{mutatis mutandis }}}

\newcommand\myref{{\textsc{ref}}}
\newcommand\RF{{\mathrm{RF}}}

\newcommand\ym{\textsc{ym}}

\newcommand\tg{{\tilde g}}

\numberwithin{equation}{section}
\newcounter{sectioncontinuation}

\usepackage{amssymb}
\usepackage{chngcntr}
\usepackage{graphicx}
\usepackage{hyperref}
\usepackage{mathrsfs}
\usepackage[usenames]{color}
\usepackage{paralist}
\usepackage{slashed}
\usepackage{url}
\usepackage{verbatim}
\usepackage{xr}

\counterwithout{section}{chapter}

\hypersetup{pdftitle={Global Existence and Convergence of Solutions to Gradient Systems and Applications to Yang-Mills Gradient Flow}}
\hypersetup{pdfauthor={Paul M. N. Feehan}}

\makeindex

\begin{document}

\frontmatter

\title[Global Existence and Convergence for Gradient Systems]{Global Existence and Convergence of Solutions to Gradient Systems and Applications to Yang-Mills Gradient Flow}

\author[Paul M. N. Feehan]{Paul M. N. Feehan}
\address{Department of Mathematics, Rutgers, The State University of New Jersey, 110 Frelinghuysen Road, Piscataway, NJ 08854-8019, United States of America}
\email{feehan@math.rutgers.edu}

\curraddr{School of Mathematics, Institute for Advanced Study, Princeton, NJ 08540}
\email{feehan@math.ias.edu}


\subjclass[2010]{Primary 35A01, 47J35, 58D25, 58E15, 58J35; secondary 35K55, 37B35, 37L15}


\maketitle



\setcounter{page}{5}

\tableofcontents



\chapter*{Preface}
Our primary goal in this monograph is to develop new results on global existence and convergence of Yang-Mills gradient flow over closed, smooth Riemannian manifolds of arbitrary dimension. Our secondary goal is to develop general methods for analyzing global existence and convergence questions for gradient-like flows on Banach spaces, so they can apply equally well to the many different gradient (or gradient-like) flows arising in geometric analysis, mathematical physics, and applied mathematics --- including Chern-Simons flow, harmonic map flow, knot energy flow, mean curvature flow, Ricci flow, and Yamabe scalar curvature flow in geometric analysis and many other examples in mathematical physics and applied mathematics.

Given $T > 0$, a Banach space $\sX$, a smooth function $\sE:\sX\to\RR$, and a point $u_0 \in \sX$, a smooth map, $u:[0, T) \to \sX$, is called a \emph{gradient flow} for $\sE$ if it is a solution to the Cauchy problem for the \emph{gradient system},
\[
\dot u(t) = -\sE'(u(t)), \quad\text{for all } t\in [0,T), \quad u(0) = u_0,
\]
as an identity in $\sX^*$ (the continuous dual space of $\sX$), where we abbreviate $\dot u = dt/dt$. The best known examples of gradient flows occurring in geometric analysis include pure and coupled Yang-Mills flows, harmonic map flow, and Yamabe scalar curvature flow. The relatively recent discovery that Ricci flow can be viewed (after some manipulation) as a gradient flow is due to Perelman \cite{Perelman_2002}. More generally, $u:[0, T) \to \sX$, is called a \emph{gradient-like flow} for $\sE$ if it is a solution to the Cauchy problem for the \emph{gradient-like system},
\[
\dot u(t) = -\sE'(u(t)) + R(t), \quad\text{for all } t\in [0,T), \quad u(0) = u_0,
\]
where $R:[0,\infty) \to \sX^*$ is a smooth map such that $R(t)$ converges to zero as $t\to\infty$.

Our monograph is partly inspired by Leon Simon's landmark article \cite{Simon_1983} on asymptotic analysis for a certain class of nonlinear evolution equations and thus is written with applications to geometric analysis in mind. However, the methods discussed in this book are not limited to geometric analysis. At the time of writing, over 250 articles cite Simon's ideas \cite{Simon_1983} and collectively they address a very broad spectrum of global existence and convergence problems arising in all areas of geometric analysis, mathematical physics, and applied mathematics.

Simon's approach relies on his celebrated generalization to infinite dimensions of the {\L}ojasiewicz gradient inequality to a specific class of analytic energy functionals on $C^{2,\alpha}$ H\"older spaces of sections of a vector bundle over a closed, smooth Riemannian manifold. Over the intervening years, this {\L}ojasiewicz-Simon gradient inequality has since been generalized by many authors --- see Feehan and Maridakis \cite{Feehan_Maridakis_Lojasiewicz-Simon_harmonic_maps} and references cited therein. For example, if $\sX$ is continuously embedded in a Hilbert space $\sH$ and $\sE$ is analytic and $x_\infty \in \sX$ is a critical point such that the Hessian operator, $\sE''(x_\infty):\sX\to\sX^*$, is Fredholm with index zero, then there exist constants $c \in [1,\infty)$, $\sigma \in (0,1]$, and $\theta \in [1/2, 1)$ such that gradient map obeys \cite[Theorem 1]{Feehan_Maridakis_Lojasiewicz-Simon_harmonic_maps}
\[
\|\sE'(x)\|_{\sX^*} \geq c|\sE(x) - \sE(x_\infty)|^\theta,
\quad\text{for all } x \in \sX \text{ such that } \|x-x_\infty\|_\sX < \sigma.
\]
However, in applications of the concepts in \cite{Simon_1983}, many authors proceed by adapting Simon's proofs or results on convergence to their particular setting without developing a more general theory in an abstract setting with universal applicability. The monograph of Huang \cite{Huang_2006} is a notable exception, but we found his abstract hypotheses on gradient-like flow difficult to verify in practice and, in particular, were unable to verify them in the case of Yang-Mills gradient flow. Fortunately, Johan R\r{a}de discovered in his beautiful analysis \cite{Rade_1992} of Yang-Mills gradient flow over manifolds of dimension two or three that the following \apriori interior length estimate (which we abstract here from the specific setting of \cite[Lemma 7.3]{Rade_1992} in Yang-Mills gauge theory) is the key property required to apply the {\L}ojasiewicz-Simon gradient inequality to establish convergence, convergence rates, global existence, and stability properties for gradient flows,
\[
\int_{\delta}^{T} \|\dot{u}(t)\|_\sX \,dt \leq C\int_0^{T} \|\dot{u}(t)\|_\sH,
\]
where $C = C(\delta) \in [1,\infty)$. While \apriori estimates for solutions to quasilinear, second-order, parabolic partial differential equations are certainly employed in \cite{Simon_1983}, this length estimate itself is not explicitly stated. Thus, one goal of our monograph is to reformulate the analysis by Huang \cite{Huang_2006} for abstract gradient-like flows by replacing his hypotheses by a hypothesis that R\r{a}de's length estimate holds along with hypotheses on short-time behavior of solutions that are typically known or easily checked in applications.

These alternative hypotheses are very convenient because it is relatively easy to prove that R\r{a}de's length estimate holds for a gradient flow if the induced evolution equation for $v(t) = \dot u(t)$,
\[
\ddot u(t) = -\sE''(u(t))\dot u(t), \quad\text{for all } t\in [0,T),
\]
has the form of a linear evolution equation \cite{Sell_You_2002},
\[
\dot v(t) + \cA v(t) = f(t), \quad\text{for all } t\in [0,T),
\]
where $\cA$ is a (positive) sectorial operator on a Banach space $\calV$ and hence the generator of an analytic semigroup on $\calV$. Typically, the function $f:[0,T)\to \cW$, which depends on $u$, is a Lipschitz map into a Banach space $\cW$. It is interesting to note that $u(t)$ itself does \emph{not} need to obey a nonlinear evolution equation of the above form with $f(t)$ replaced by $\cF(t,u(t))$ and, in the case of Yang-Mills gradient flow, one sees a nice example of this phenomenon. Indeed, it is well-known that certain flows --- such as mean curvature flow, Ricci flow, or Yang-Mills gradient flow --- only become solutions to quasilinear, second-order, parabolic partial differential equations after applying a version of the \emph{DeTurck trick} \cite{DeTurck_1983} with a suitably-chosen time-varying family of diffeomorphisms or gauge transformations. However, it is not always desirable to apply the DeTurck trick and R\r{a}de does not require it for his analysis in \cite{Rade_1992}. In the case of Ricci flow, short-time existence of solutions can be obtained either by application of the DeTurck trick \cite{DeTurck_1983} or the \emph{Nash-Moser Implicit Function Theorem} \cite{Hamilton_1982, Hamilton_1982bams, Moser_1966a, Moser_1966b, Nash_1956} and it is likely that the same should hold for Yang-Mills gradient flow.

A gradient system for $u:[0,T) \to \sX$ often has the form,
\[
\dot u(t) + \cA u(t) = \cF(t,u(t)), \quad\text{for all } t\in [0,T),
\]
where $\cF:[0,\infty)\times\calV^{2\beta}\to \cW$ is a Lipschitz map into a Banach space $\cW$ and $\calV^{2\beta}$, for $\beta \in [0,1)$, is a fractional power of the domain $\calV^2 = \sD(\cA)$ of an unbounded linear operator $\cA:\sD(\cA) \subset \cW\to\cW$. For example, the equation for harmonic map gradient flow has this form and, after an application of the DeTurck trick, so do the equations for Ricci flow and Yang-Mills gradient flow.

The development of \apriori estimates, short-time existence, and regularity for mild solutions to such nonlinear evolution equations has been established in considerable generality --- see, for example, the wonderful exposition by Sell and You \cite{Sell_You_2002} --- and this can frequently eliminate the need for such analysis to be repeated in specific applications. In our monograph, we lean heavily on the development in \cite{Sell_You_2002} and extend their analysis in directions suitable for application to gradient systems with properties frequently encountered in differential geometry. The general approach to nonlinear evolution equations in \cite{Sell_You_2002} and elsewhere using analytic semigroups effectively reduces the analysis to one of showing that the relevant elliptic differential or pseudo-differential operators are sectorial on a desired Banach space, generalizing the important \emph{resolvent estimates} due to Agmon \cite{AgmonLecturesEllipticBVP}. This analysis requires a careful treatment of existence, uniqueness, regularity, and \apriori estimates for solutions to elliptic systems, whether for sections of vector bundles over closed, smooth manifolds or domains in Euclidean space. In the case of scalar elliptic or parabolic partial differential equations, Gilbarg and Trudinger \cite{GilbargTrudinger}, Krylov \cite{Krylov_LecturesHolder, Krylov_LecturesSobolev}, and Lieberman \cite{Lieberman} provide excellent, self-contained references. By contrast, well-known references such as Agmon \cite{AgmonLecturesEllipticBVP}, Lady{\v{z}}enskaja, Solonnikov, and Ural$'$ceva \cite{LadyzenskajaSolonnikovUralceva} or Morrey \cite{Morrey} notwithstanding, a sufficiently general analysis for elliptic or parabolic systems required in applications is very difficult to find in the literature and results, when they can be found, are widely scattered among many references. Therefore, we give a largely self-contained development of the results we need for elliptic systems, including the key resolvent (Agmon) estimates required to show that elliptic systems on Sobolev spaces define sectorial operators in a wide variety of settings. While we restrict our attention to partial differential operators rather than pseudo-differential operators more generally, recent applications (see Blatt \cite{Blatt_2016arxiv} and references cited therein) to gradient flows for knot energy functionals may point to a growth in the role of elliptic pseudo-differential operators in geometric analysis.

The {\L}ojasiewicz-Simon gradient inequality can be used to analyze convergence, convergence rates, global existence, and stability of gradient flows started near a critical point of an analytic potential function $\sE$ and thus provides most information when that critical point is a local minimum. When a geometric gradient flow is started at a point that is not close to a critical point, then such flows often acquire singularities, at finite or infinite time. In the case of Yang-Mills gradient flow over four-dimensional manifolds, we give a self-contained, careful treatment of the energy bubbling phenomenon and the formation of singularities at finitely distinct points in the base four-dimensional manifold. As observed by Struwe \cite{Struwe_1985, Struwe_1996}, the bubbling phenomenon for harmonic map gradient flow is similar in many (though not all) respects.

As we noted at the beginning of our preface, our primary goal in this monograph is to develop new results on global existence and convergence of solutions to the gradient flow equation for the Yang-Mills energy functional on a principal $G$-bundle, $P$, over a closed, Riemannian, smooth manifold, $X$, where $G$ is a compact Lie group. Much of our analysis allows $X$ to have arbitrary dimension, but ultimately we focus on the case where $X$ has dimension four. The study of this gradient flow had been initiated by Atiyah and Bott \cite{Atiyah_Bott_1983} in the case where $X$ is a closed Riemann surface. The first major results in (real) dimension four are due to Donaldson \cite{DonASD}, who proved global existence for Yang-Mills gradient flow on a Hermitian vector bundle, $E$, over a compact, K\"ahler surface, $X$, when the initial connection, $A_0$, on $E$ is unitary with curvature of type $(1,1)$ and in addition subsequential convergence, modulo gauge transformations, to a Hermitian Yang-Mills connection on $E$ if the holomorphic bundle, $(E,\bar\partial_{A_0})$, is also stable.

The challenge ever since Donaldson's celebrated article \cite{DonASD} has been to understand to what extent his global existence and convergence results continue to hold when $X$ is allowed to be any closed, four-dimensional, Riemannian, smooth manifold and there are no assumptions of stability for $P$. In this monograph, we prove the following:
\begin{enumerate}
\item
\label{item:Global_existence_and_convergence_near_minimum}
When the initial connection, $A_0$, on $P$ is close enough
to a local minimum of the Yang-Mills energy functional, the Yang-Mills gradient flow exists for all time and converges to a Yang-Mills connection on $P$ as time tends to infinity.

\item
\label{item:Global_existence_arbitrary_initial_energy}
If the initial connection, $A_0$, is allowed to have arbitrary energy but we restrict to the setting of a Hermitian vector bundle over a compact, complex, Hermitian (but not necessarily K\"ahler) surface and $A_0$ has curvature of type $(1,1)$, then the Yang-Mills gradient flow exists for all time, though bubble singularities may (and in certain cases must) occur in the limit as time tends to infinity.
\end{enumerate}

In order to better understand some the shared properties and the differences in as broad a class of nonlinear evolutionary equations as possible, our development in this monograph proceeds from the very general to the specific: nonlinear evolutionary equations in Banach spaces, abstract gradient and gradient-like systems in Banach spaces, Yang-Mills gradient flow over closed, Riemannian manifolds of arbitrary dimension $d \geq 2$, Yang-Mills gradient flow over closed, four-dimensional Riemannian manifolds, and finally Yang-Mills gradient flow over complex surfaces. We hope that our monograph will be of interest to researchers in all areas of geometric analysis, mathematical physics, and applied mathematics who are concerned with questions of global existence and convergence of solutions to nonlinear evolutionary equations that can be analyzed as gradient or gradient-like flows.

\chapter*{Acknowledgments}
We became interested in the questions of global existence and convergence of Yang-Mills gradient flow while participating in the workshop \emph{Geometry and topology of smooth 4-manifolds}, at the Max Planck Institute for Mathematics, Bonn, June 3--7, 2013. Our approach in this monograph to these questions was described in one half of our research proposal submitted in October 2013 to the Analysis program of the Division of Mathematical Sciences at the National Science Foundation.

We are very grateful to Professor Peter Teichner and the staff of the Max Planck Institute for Mathematics, Bonn, for their generous hospitality during our visit in Summer 2013. Much our research was conducted while a visiting scholar in the Department of Mathematics at Columbia University, New York, and we warmly thank Professor Ioannis Karatzas and his staff for their kind hospitality. We would like to thank Professor Ari Laptev and the staff of the Institut Mittag-Leffler, Stockholm, for a very useful visit in January 2014. We are most grateful to Professors Michael Struwe and Tristan Rivi\`ere and the Forschungsinstitut f\"ur Mathematik at ETH Z\"urich and Professor Nigel Hitchin at the Mathematical Institute, University of Oxford, for much appreciated visits in May and June 2014, respectively, and early opportunities to present the main results of our monograph in their seminars. Part of our research was completed while visiting the Isaac Newton Institute for Mathematical Sciences, Cambridge, during June and July 2014, for their program, \emph{Free Boundary Problems and Related Topics}; we are extremely grateful to Professor Henrik Shahgholian for his invitation and the staff of the Newton Institute for their support. We warmly thank Professor Herbert Koch and staff of the Hausdorff Research Institute for Mathematics, Bonn, for an opportunity to participate in July and August 2014 in their trimester program, \emph{Harmonic Analysis and Partial Differential Equations}, and present our research. We thank Professor J\o rgen Andersen and the Center for Quantum Geometry of Moduli Spaces at \r{A}rhus University, Denmark, for their generous hospitality during our visit in August 2014. Lastly, we thank Professors Helmut Hofer and Thomas Spencer at the Institute for Advanced Study, Princeton for facilitating our visit to Institute for Advanced Study during the academic year 2015-16 and Professor Simon Thomas for facilitating our competitive fellowship leave from Rutgers University.

Since we commenced work on this monograph, many mathematicians kindly responded to our questions about their articles or related research. In that regard, we are most grateful to Simon Brendle, Piermarco Cannarsa, Huai-Dong Cao, David Groisser, Robert Haslhofer, Adam Jacob, Sergiu Klainerman, John Lott, Andre Neves, Johan R\r{a}de, Lorenzo Sadun, Andreas Schlatter, Peter Takac, Peter Topping, and Michael Weinstein. We are especially grateful to our colleagues at Rutgers, Natasa Sesum and Shadi Tahvildar-Zadeh, for explanations of their articles and related research. We are indebted to Richard Wentworth for his invaluable help in clarifying several important points and his patient explanations.

We are very grateful to the mathematicians who provided us with additional opportunities to present our results in their seminars commencing in Spring 2014, including Huai-Dong Cao at Lehigh University, Adam Levine and Rafael Montezuma at Princeton University, Tian-Jun Li and Camelia Pop at the University of Minnesota, Cl\'ement Mouhot at the University of Cambridge, Thomas Parker at Michigan State University, Brendan Owens at the University of Glasgow, Christian Saemann at Heriot-Watt University and Antony Maciocia at the University of Edinburgh, Rafe Mazzeo at Stanford University, Benjamin Sharpe and Richard Thomas at Imperial College London, Daniela De Silva at Columbia University, and Jeffrey Streets at the University of California at Irvine.

Many of the participants in our seminar lectures provided invaluable assistance by way of thoughtful questions and feedback; we would especially like to thank Sir Michael Atiyah, Sir Simon Donaldson, Richard Hamilton, Nigel Hitchin, Dominic Joyce, Duong Phong, and Michael Struwe. We are very grateful to Manousos Maridakis for a very helpful collaboration to understand the abstract {\L}ojasiewicz-Simon gradient inequality and its applications to coupled Yang-Mills energy functionals. We thank Peter Tak{\'a}{\v{c}} for many helpful conversations regarding the {\L}ojasiewicz-Simon gradient inequality, for explaining his proof of \cite[Proposition 6.1]{Feireisl_Takac_2001} and how it can be generalized, and for his kindness when hosting his visit to the Universit{\"a}t R{\"o}stock in July 2015.

Work on our monograph consumed an ever increasing proportion of our time since June 2013 and we gratefully acknowledge the patience and tolerance of our collaborators on unrelated projects during this period, including Ruoting Gong, Camelia Pop, Jian V. Song and especially Tom Leness.

We are most grateful for the longstanding encouragement of Ronald Fintushel, Peter Kronheimer, Tom Mrowka, Duong Phong, and Cliff Taubes for our research in gauge theory and to Karen Uhlenbeck for her interest in our work and engaging conversations on Yang-Mills theory.

We were partially supported by National Science Foundation grant DMS-1510064 and the Oswald Veblen Fund and Fund for Mathematics (Institute for Advanced Study, Princeton) during the preparation of this monograph.
\bigskip
\bigskip

\leftline{October 14, 2016}
\leftline{New Brunswick, New Jersey}

\mainmatter

\chapter{Introduction}
\label{chapter:Introduction}

\section{Introduction to Yang-Mills gradient flow}
\label{sec:Introduction_Yang-Mills_gradient_flow}
Let $G$ be a compact Lie group and $P$ a principal $G$-bundle over a closed, connected, smooth manifold, $X$, with Riemannian metric, $g$. The quotient $\sB(P) = \sA(P)/\sG(P)$ of the affine space, $\sA(P)$, of connections on $P$, modulo the action of the group $\sG(P) = \Aut P$ of gauge transformations of $P$, has fundamental significance in algebraic geometry, low-dimensional topology, classification of four-dimensional, smooth  manifolds, and high-energy physics.

Given a connection $A$ on $P$, we denote its curvature by $F_A \in \Omega^2(X; \ad P) \equiv C^\infty(X,\Lambda^2\otimes\ad P)$ and consider the associated \emph{Yang-Mills energy functional},
\index{Yang-Mills energy functional}
\begin{equation}
\label{eq:Yang-Mills_energy_functional}
\sE_g(A)  := \frac{1}{2}\int_X |F_A|^2\,d\vol_g.
\end{equation}
Atiyah and Bott \cite{Atiyah_Bott_1983} had proposed that the gradient flow of this energy functional with respect to the natural $L^2$ Riemannian metric on $\sB(P)$ would prove to be of vital importance in understanding the topology of $\sB(P)$
via an infinite-dimensional Morse theory and explored this approach in the case where $X$ is a Riemann surface.

The following conjecture essentially goes back to Atiyah and Bott \cite{Atiyah_Bott_1983}, Sedlacek \cite{Sedlacek}, Taubes \cite{TauPath, TauFrame, TauStable}, and Uhlenbeck \cite{UhlLp, UhlRem}; it appears explicitly in an article by Schlatter, Struwe, and Tahvildar-Zadeh \cite[p. 118]{Schlatter_Struwe_Tahvildar-Zadeh_1998} and elsewhere.

\begin{conj}[Global existence and convergence of Yang-Mills gradient flow over closed four-dimensional manifolds]
\label{conj:Yang-Mills_gradient_flow_global_existence}
Let $G$ be a compact Lie group and $P$ a principal $G$-bundle over a closed, connected, four-dimensional, smooth manifold, $X$, with Riemannian metric, $g$.
If $A_0$ is a smooth connection on $P$, then there is a smooth solution, $A(t)$ for $t\in [0,\infty)$, to the gradient flow,
\begin{align}
\label{eq:Yang-Mills_gradient_flow}
\frac{\partial A}{\partial t} &= -d_{A(t)}^{*_g}F_{A(t)},
\\
\label{eq:Yang-Mills_gradient_flow_initial_condition}
A(0) &= A_0,
\end{align}
for the Yang-Mills energy functional with respect to the $L^2$ Riemannian metric on the affine space of connections on $P$. Moreover, as $t\to\infty$, the flow, $A(t)$, converges to a smooth Yang-Mills connection, $A_\infty$, on $P$.
\end{conj}

The importance of the question embodied in Conjecture \ref{conj:Yang-Mills_gradient_flow_global_existence} had begun to gain further recognition in the early mathematical approaches to Yang-Mills gauge theory in \cite{AHS, Bourguignon_Lawson_1981, Bourguignon_Lawson_Simons_1979, UhlRem, UhlLp} (during the period 1978 to 1982) and more specifically in \cite{Atiyah_Bott_1983, DonASD, TauPath, TauFrame, TauStable} (during the period 1982 to 1989).

As we shall explain, four is the critical dimension for the base manifold, $X$, and the focus of this monograph. When $X$ has dimension two or three, much more is known. The initial energy $\sE(A_0)$, Lie group $G$, topology of $P$ and $X$, Riemannian metric $g$ on $X$ and possible local or global symmetries in the flow all bear on the questions of long-time existence and convergence, so we must allow refinements to Conjecture \ref{conj:Yang-Mills_gradient_flow_global_existence}. A connection, $A$, is a \emph{critical point} of $\sE$ if and only if it obeys the Yang-Mills equation,
\begin{equation}
\label{eq:Yang-Mills_equation}
d_A^{*,g}F_A = 0 \quad\hbox{on }X,
\end{equation}
since $d_A^{*,g}F_A = \sE_g'(A)$ when the gradient of $\sE$ is defined by the $L^2$ metric.

The first major advance towards Conjecture \ref{conj:Yang-Mills_gradient_flow_global_existence} was due to Donaldson \cite{DonASD} in the case of a Hermitian vector bundle, $E$, over a complex projective algebraic surface, $X$, later generalized in \cite{DK} to the case of a K\"ahler surface; see Theorem \ref{mainthm:Donaldson_Kronheimer_6-2-7_and_6-2-14_plus} for a statement of (a special case of) his results on global existence and subsequential convergence modulo gauge transformations.

Donaldson's results in \cite{DonASD} established the existence of anti-self-dual connections on stable, holomorphic vector bundles over complex algebraic surfaces. Those results marked the first significant extensions of earlier results due to Taubes \cite{TauSelfDual, TauIndef} on the existence of anti-self-dual connections over closed Riemannian four-manifolds; taken together, they led to enormous advances in our understanding of the smooth classification of four-manifolds \cite{DK, Friedman_Morgan_1998, KMStructure}.

Daskalopoulos and Wentworth \cite{Daskalopoulos_Wentworth_2004, Daskalopoulos_Wentworth_2007} proved certain extensions of Donaldson's results for Yang-Mills gradient flow on a Hermitian vector bundle, $E$, over a K\"ahler surface, $X$, when the stability condition assumed by Donaldson is relaxed. In particular, it follows from results of Daskalopoulos and Wentworth in \cite{Daskalopoulos_Wentworth_2004, Daskalopoulos_Wentworth_2007} that one can construct examples of unstable holomorphic vector bundles, $E$, and initial connections, $A_0$, such that the Yang-Mills gradient flow necessarily develops bubble singularities at $T = \infty$; see Section \ref{subsec:Daskalopoulos_Wentworth} for a discussion of one such example. On the other hand, even when the stability condition employed by Donaldson is relaxed, they show that the Yang-Mills gradient flow cannot develop bubble singularities in finite time, $T < \infty$. Hong and Tian \cite{Hong_Tian_2004} independently established related results on the asymptotic behavior of Yang-Mills gradient flow, relying more on analytical methods.

When $X$ instead has dimension two or three, R\r{a}de \cite[Theorems 1, $1'$, and 2]{Rade_1992} has shown that Conjecture \ref{conj:Yang-Mills_gradient_flow_global_existence} is true. Daskalopoulos \cite{Daskalopoulos_1992} proved a similar result when $X$ has dimension two, namely global existence and convergence modulo gauge transformations.

Struwe \cite[Theorem 2.3]{Struwe_1994} established existence and uniqueness of solutions to the Yang-Mills gradient flow \eqref{eq:Yang-Mills_gradient_flow}, up to a finite time $T_1>0$ characterized by the first occurrence of bubble singularities and conjectured the long-time existence and uniqueness of the gradient flow \cite[Theorem 2.4]{Struwe_1994}, modulo blow-ups and bundle topology changes at finitely many times, $0<T_1<\cdots<T_K\leq \infty$, a result later proved by his Ph.D. student, Schlatter
\cite[Theorems 1.2 and 1.3]{Schlatter_1997}. Similar results were proved independently by
Kozono, Maeda, and Naito \cite{Kozono_Maeda_Naito_1995}; see also results of Chen and Shen
\cite{Chen_Shen_1993, Chen_Shen_1994, Chen_Shen_1995, Chen_Shen_Zhou_2002}.

Using methods of stochastic analysis, Pulemotov \cite{Pulemotov_2008} has established long-time existence for the Yang-Mills gradient flow over the unit ball, $B\subset\RR^4$, when $F_A(t)$ has suitable boundary conditions and the initial energy, $\sE(A_0)$, is sufficiently small.

Kozono, Maeda, and Naito \cite[Corollary 5.7]{Kozono_Maeda_Naito_1995} established the global existence in Conjecture \ref{conj:Yang-Mills_gradient_flow_global_existence} when $P$ is a product bundle, $X\times G$, and the initial energy, $\sE(A_0)$, is sufficiently small. Schlatter \cite[Theorem 1.6]{Schlatter_1996} claimed to extend the preceding results for $\sE(A_0)$ small and $P = X\times G$ to the case of a principal $\SU(2)$-bundle $P$ over $S^4$, with $c_2(E) = 1$ and $\|F_{A_0}^+\|_{L^2(S^4)}$ small and $S^4$ having its standard, round metric of radius one. Schlatter's proof-by-contradiction argument tries to adapt the proof of a related result for harmonic maps due to Struwe \cite{Struwe_1996}.
His argument hinges on there being a positive lower bound for the energy of a connection on a principal $\SU(2)$-bundle $\widetilde P \to S^4$ that he builds in his construction \cite[Section 5]{Schlatter_1996}. However, a close examination of his proof reveals that he does not use the hypothesis that $A(t)$ is Yang-Mills gradient flow and, after examining his construction of $\widetilde P$, one discovers that $\widetilde P \cong S^4 \times \SU(2)$.

Schlatter, Struwe, and Tahvildar-Zadeh \cite[Theorem 3.1]{Schlatter_Struwe_Tahvildar-Zadeh_1998} have established long-time existence for the $\SO(4)$-equivariant Yang-Mills gradient flow on the product $\SU(2)$-bundle over the unit ball, $B\subset\RR^4$, for any finite initial energy, $\sE(A_0)$, and boundary condition on $\partial B$. In contrast, Grotowski \cite{Grotowski_2001} has shown that for the product bundle, $\RR^d\times \SO(d)$ with $d \geq 5$, Yang-Mills gradient flow develops singularities in finite time for a class of $\SO(d)$-equivariant initial connections, $A_0$.

When $d=4$, the Yang-Mills energy functional is invariant with respect to conformal changes of the Riemannian metric on $X$. Using the fact that the Yang-Mills energy functional is not conformally invariant when $d \geq 5$, one can employ rescaling to obtain connections with arbitrarily small energy on nontrivial principal bundles over the $d$-dimensional sphere, $S^d$, with its standard, round metric of radius one. Naito \cite[Theorem 1.3]{Naito_1994} exploited this property to show that if $G \subset \SO(d)$, for $d \geq 5$, and $P$ is a nontrivial principal $G$-bundle over $S^d$ (with its standard, round metric of radius one), there is a positive constant, $\eps_1$, such that if $\|F_{A_0}\|_{L^2(S^4)} < \eps_1$, then a solution to \eqref{eq:Yang-Mills_gradient_flow} and \eqref{eq:Yang-Mills_gradient_flow_initial_condition} blows up in finite time.

Related finite-time blow-up results for $S^1$-equivariant harmonic map gradient flow, from $B\subset\RR^2$ to $S^2\subset\RR^3$, were established by Chang, Ding, and Ye \cite{Chang_Ding_Ye_1992}.

Waldron \cite{Waldron_2014arxiv, WaldronThesis} has recently established several significant results concerning Conjecture \ref{conj:Yang-Mills_gradient_flow_global_existence} by methods that are entirely different from those employed in our monograph.

Chen and Zhang \cite{Chen_Zhang_2015} and Kelleher and Streets have an ongoing programs to study the important problem of singularity formulation in Yang-Mills gradient flow over base manifolds of dimension four and higher \cite{Kelleher_Streets_2014arxiv, Kelleher_Streets_2016preprint}.

\section{Main results}
In this section, we highlight some of our main results in this monograph for abstract gradient flow on Banach spaces and Yang-Mills gradient flow in particular.

\subsection{Convergence, global existence, and stability for solutions to abstract gradient systems on Banach spaces}
\label{subsec:Main_results_gradient_system_Banach_space}
For the convenience of the reader, we shall summarize in this section some of the main results for abstract gradient systems on Banach spaces that are stated in more generality in Sections \ref{sec:Huang_3_and_5_gradient-system} and \ref{sec:Huang_3_and_5_gradientlike_system} for \emph{gradient} and \emph{gradient-like systems}, respectively. (See Palais \cite{Palais_1966LS} or Mawhin and Willem \cite[Section 4]{Mawhin_Willem_2010} for the related concept of \emph{pseudo-gradient system} due to Palais.) Our results could easily be stated for gradient systems on Banach manifolds but, as the results are `local' in nature, we may confine our attention to gradient flow on Banach spaces without loss of generality. The conclusions of our results resemble those of
\begin{inparaenum}[\itshape a\upshape)]
\item Simon in \cite{Simon_1983}, for a certain class of parabolic, quasi-linear, second-order partial differential equations or systems on a closed, finite-dimensional, smooth Riemannian manifold, and R\r{a}de \cite{Rade_1992}, for Yang-Mills gradient flow over a closed, smooth Riemannian manifold of dimension two or three, but hold in far greater generality, and
\item Huang in \cite{Huang_2006}, but replace hypotheses that are exceedingly difficult to verify by ones that are relatively straightforward to confirm.
\end{inparaenum}
The collection of results in this section comprise a `toolkit' that may be directly and easily applied to analyze a wide range of gradient systems in geometric analysis --- including harmonic map gradient flow, knot energy flow, mean curvature flow, Ricci flow, Yamabe flow, and (coupled and pure) Yang-Mills flow, as well as numerous other gradient systems in applied mathematics and mathematical physics --- \emph{without} the need to reprove analogues for those applications of the original results due to Simon or R\r{a}de.

\begin{hyp}[\Apriori interior estimate for a trajectory]
\label{hyp:Abstract_apriori_interior_estimate_trajectory_main_introduction}
(Compare Hypothesis \ref{hyp:Abstract_apriori_interior_estimate_trajectory}.)
Let $\sX$ be a Banach space that is continuously embedded in a Hilbert space $\sH$. If $\delta \in (0,\infty)$ is a constant, then there is a constant $C_1 = C_1(\delta) \in [1,\infty)$ with the following significance. If $S, T \in \RR$ are constants obeying $S+\delta \leq T$ and $u \in C^\infty([S,T); \sX)$, we say that $\dot u \in C^\infty([S,T); \sX)$ obeys an \apriori \emph{interior estimate on $(0, T]$} if
\begin{equation}
\label{eq:Abstract_apriori_interior_estimate_trajectory_main_introduction}
\int_{S+\delta}^T \|\dot u(t)\|_\sX\,dt \leq C_1\int_S^T \|\dot u(t)\|_\sH\,dt.
\end{equation}
\end{hyp}

In applications, $u \in C^\infty([S,T); \sX)$ in Hypothesis \ref{hyp:Abstract_apriori_interior_estimate_trajectory_main_introduction} will often be a solution to a quasi-linear parabolic partial differential system, from which an \apriori estimate \eqref{eq:Abstract_apriori_interior_estimate_trajectory_main_introduction} may be deduced. For example, Hypothesis \ref{hyp:Abstract_apriori_interior_estimate_trajectory_introduction} can be verified by Lemma \ref{lem:Rade_7-3_abstract_interior_L1_in_time_V2beta_space_time_derivative_interior} for a nonlinear evolution equation on a Banach space $\calV$ of the form (see Caps \cite{Caps_evolution_equations_scales_banach_spaces}, Henry \cite{Henry_geometric_theory_semilinear_parabolic_equations}, Pazy \cite{Pazy_1983}, Sell and You \cite{Sell_You_2002}, Tanabe \cite{Tanabe_1979, Tanabe_1997} or Yagi \cite{Yagi_abstract_parabolic_evolution_equations_applications})
\begin{equation}
\label{eq:Nonlinear_evolution_equation_Banach_space}
\frac{du}{dt} + \cA u = \cF(t,u(t)), \quad t\geq 0, \quad u(0) = u_0,
\end{equation}
where $\cA$ is a positive, sectorial, unbounded operator on a Banach space, $\cW$, with domain $\calV^2 \subset \cW$ and the nonlinearity, $\cF$, has suitable properties.

For example, given a compatible connection, $\nabla$, on a finite-rank, Riemannian vector bundle over a closed, smooth, Riemannian manifold $(X,g)$, and a choice of Banach space $\cW = L^p(X;V)$ for $p \in (1,\infty)$, one can choose $\calV^{2\alpha}$, for $\alpha\geq 0$, to be the domain $W^{2\alpha,p}(X;\RR)$ of the fractional powers, $(\nabla^{*_g}\nabla + 1)^\alpha$, of the unbounded operator, $\cA = \nabla^{*_g}\nabla + 1: L^p(X;V) \to L^p(X;V)$. When $p=2$ and $\cH = L^2(X;V)$, one can choose $\calV^{2\alpha}$, for $\alpha\in\RR$, to be the domain $H^{2\alpha}(X;V)$ of the fractional powers, $(\nabla^{*_g}\nabla + 1)^\alpha$, of the unbounded operator, $\nabla^{*_g}\nabla + 1: L^2(M;\RR) \to L^2(M;\RR)$. (See Kre{\u\i}n, Petun{\={\i}}n, and Sem{\"e}nov \cite{Krein_Petunin_Semenov_interpolation_linear_operators} for a discussion of such \emph{scales} of Banach spaces.)

Harmonic map gradient flow (as a quasi-linear parabolic partial differential system) may be placed in the form \eqref{eq:Nonlinear_evolution_equation_Banach_space} directly (for suitable choices of Banach spaces) whereas the non-parabolic Ricci or Yang-Mills gradient flows acquire this form after application of a version of the DeTurck trick \cite{DeTurck_1983}, by choosing a suitable path of diffeomorphisms of the base manifold or the automorphisms of the principal $G$-bundle, respectively,

\begin{rmk}[Regularity and \apriori interior estimates for solutions to Yang-Mills gradient flow]
\label{rmk:Abstract_apriori_interior_estimate_trajectory_introduction}
The Hypothesis \ref{hyp:Abstract_apriori_interior_estimate_trajectory_introduction} is verified for nonlinear evolution equations on Banach spaces by Lemma \ref{lem:Rade_7-3_abstract_interior_L1_in_time_V2beta_space_time_derivative_interior} and for  solutions to Yang-Mills gradient flow in Lemmata \ref{lem:Rade_7-3},  \ref{lem:Rade_7-3_L1_in_time_H2beta_in_space_apriori_estimate_by_L1_in_time_L2_in_space}, \ref{lem:Rade_7-3_arbitrary_dimension}, and Corollary \ref{cor:Rade_7-3_arbitrary_dimension_L1_time_W1p_space}.
\end{rmk}

We recall our

\begin{thm}[{\L}ojasiewicz-Simon gradient inequality for analytic functions on Banach spaces]
\label{thm:Huang_2-4-5_introduction}
\cite[Theorem 1]{Feehan_Maridakis_Lojasiewicz-Simon_harmonic_maps}, \cite[Theorem 2.4.5]{Huang_2006}
Let $\sX$ be a Banach space that is continuously embedded in a Hilbert space $\sH$. Let $\sU \subset \sX$ be an open subset, $\sE:\sU\to\RR$ be an analytic function, and $\varphi\in\sU$ be a critical point of $\sE$, that is, $\sE'(\varphi) = 0$. Assume that $\sE''(\varphi):\sX\to \sX'$ is a Fredholm operator with index zero. Then there are constants, $c \in [1,\infty)$, and $\sigma \in (0,1]$,
and $\theta \in [1/2,1)$ such that
\begin{equation}
\label{eq:Simon_2-2_dualspacenorm_introduction}
\|\sE'(u)\|_{\sX'} \geq c|\sE(u) - \sE(\varphi)|^\theta, \quad \forall\, u \in \sU \hbox{ such that } \|u-\varphi\|_\sX < \sigma.
\end{equation}
\end{thm}

We have the following analogue of Huang \cite[Theorems 3.3.3 and 3.3.6]{Huang_2006} and abstract analogue of Simon \cite[Corollary 2]{Simon_1983}.

\begin{mainthm}[Convergence of a subsequence implies convergence for a smooth solution to a gradient system]
\label{mainthm:Simon_corollary_2_introduction}
(Compare Theorem \ref{thm:Simon_corollary_2}.)
Let $\sU$ be an open subset of a real Banach space, $\sX$, that is continuously embedded and dense in a Hilbert space, $\sH$. Let $\sE:\sU\subset \sX\to\RR$ be an analytic function with gradient map $\sE':\sU\subset \sX \to \sH$. Assume that $\varphi \in \sU$ is a critical point of $\sE$, that is $\sE'(\varphi)=0$. If $u \in C^\infty([0,\infty); \sX)$ is a solution to the gradient system,
\begin{equation}
\label{eq:gradient_system}
\dot u(t) = -\sE'(u(t)), \quad t \in (0,\infty),
\end{equation}
and the orbit $O(u) = \{u(t): t\geq 0\} \subset \sX$ is
precompact\footnote{Recall that \emph{precompact} (or \emph{relatively compact}) subspace $Y$ of a topological space $X$ is a subset whose closure is compact. If the topology on $X$ is metrizable, then a subspace $Z \subset X$ is compact if and only if $Z$ is \emph{sequentially compact} \cite[Theorem 28.2]{Munkres2}, that is, every infinite sequence in $Z$ has a convergent subsequence in $Z$ \cite[Definition, p. 179]{Munkres2}.}, and $\varphi$ is a cluster point of $O(u)$, then $u(t)$ converges to $\varphi$ as $t\to\infty$ in the sense that
\[
\lim_{t\to\infty}\|u(t)-\varphi\|_\sX = 0
\quad\hbox{and}\quad
\int_0^\infty \|\dot u\|_\sH\,dt < \infty.
\]
Furthermore, if $u$ satisfies Hypothesis \ref{hyp:Abstract_apriori_interior_estimate_trajectory_main_introduction} on $(0, \infty)$, then
\[
\int_1^\infty \|\dot u\|_\sX\,dt < \infty.
\]
\end{mainthm}

We next have the following abstract analogue of R\r{a}de's \cite[Proposition 7.4]{Rade_1992}, in turn a variant the \emph{Simon Alternative}, namely \cite[Theorem 2]{Simon_1983}.

\begin{mainthm}[Simon Alternative for convergence for a smooth solution to a gradient system]
\label{mainthm:Huang_3-3-6_introduction}
(Compare Theorem \ref{thm:Huang_3-3-6}.)
Let $\sU$ be an open subset of a real Banach space, $\sX$, that is continuously embedded and dense in a Hilbert space, $\sH$. Let $\sE:\sU\subset \sX\to\RR$ be an analytic function with gradient map $\sE':\sU\subset \sX \to \sH$. Assume that
\begin{enumerate}
\item $\varphi \in \sU$ is a critical point of $\sE$, that is $\sE'(\varphi)=0$; and

\item Given positive constants $b$, $\eta$, and $\tau$, there is a constant $\delta = \delta(\eta, \tau, b) \in (0, \tau]$ such that if $v$ is a smooth solution to \eqref{eq:gradient_system} on $[t_0, t_0 + \tau)$ with $t_0 \in \RR$ and $\|v(t_0)\|_\sX \leq b$, then
\begin{equation}
\label{eq:Gradient_solution_near_initial_data_at_t0_for_short_enough_time_introduction}
\sup_{t\in [t_0, t_0+\delta]}\|v(t) - v(t_0)\|_\sX < \eta.
\end{equation}
\end{enumerate}
If $(c,\sigma,\theta)$ are the {\L}ojasiewicz-Simon constants for $(\sE,\varphi)$, then there is a constant
\[
\eps = \eps(c, C_1, \delta, \theta, \rho, \sigma, \tau, \varphi) \in (0, \sigma/4)
\]
with the following significance.  If $u:[0, \infty)\to \sU$ is a smooth solution to \eqref{eq:gradient_system} that satisfies Hypothesis \ref{hyp:Abstract_apriori_interior_estimate_trajectory_main_introduction} on $(0, \infty)$ and there is a constant $T \geq 0$ such that
\begin{equation}
\label{eq:Rade_7-2_banach_introduction}
\|u(T) - \varphi\|_\sX < \eps,
\end{equation}
then either
\begin{enumerate}
\item
\label{item:Theorem_3-3-6_energy_u_at_time_t_below_energy_critical_point_introduction}
$\sE(u(t)) < \sE(\varphi)$ for some $t>T$, or
\item
\label{item:Theorem_3-3-6_u_converges_to_limit_u_at_infty_introduction}
$u(t)$ converges in $\sX$ to a limit $u_\infty \in \sX$ as $t\to\infty$ in the sense that
$$
\lim_{t\to\infty}\|u(t)-u_\infty\|_\sX =0
\quad\hbox{and}\quad
\int_1^\infty \|\dot u\|_\sX\,dt < \infty.
$$
If $\varphi$ is a cluster point of the orbit $O(u) = \{u(t): t\geq 0\}$, then $u_\infty = \varphi$.
\end{enumerate}
\end{mainthm}

In applications, the short-time estimate \eqref{eq:Gradient_solution_near_initial_data_at_t0_for_short_enough_time_introduction} for $v \in C^\infty([t_0,t_0+\tau); \sX)$ will usually follow from the fact that $v$ is a solution to a quasi-linear parabolic partial differential system, from which \eqref{eq:Gradient_solution_near_initial_data_at_t0_for_short_enough_time_introduction} may be deduced. We have the following enhancement of Huang \cite[Theorem 3.4.8]{Huang_2006}.

\begin{mainthm}[Convergence rate under the validity of a {\L}ojasiewicz-Simon gradient inequality]
\label{mainthm:Huang_3-4-8_introduction}
(Compare Theorem \ref{thm:Huang_3-4-8}.)
Let $\sU$ be an open subset of a real Banach space, $\sX$, that is continuously embedded and dense in a Hilbert space, $\sH$. Let $\sE:\sU\subset \sX\to\RR$ be an analytic function with gradient map $\sE':\sU\subset \sX \to \sH$. Let $u:[0,\infty) \to \sX$ be a smooth solution to the gradient system \eqref{eq:gradient_system} and assume that $O(u) \subset \sU_\sigma \subset \sU$, where $(c,\sigma,\theta)$ are the {\L}ojasiewicz-Simon constants for $(\sE,\varphi)$ and $\sU_\sigma := \{x\in \sX: \|x-\varphi\|_\sX < \sigma\}$. Then there exists $u_\infty \in \sH$ such that
\begin{equation}
\label{eq:Huang_3-45_H_introduction}
\|u(t) - u_\infty\|_\sH \leq \Psi(t), \quad t\geq 0,
\end{equation}
where
\begin{equation}
\label{eq:Huang_3-45_growth_rate_introduction}
\Psi(t)
:=
\begin{cases}
\displaystyle
\frac{1}{c(1-\theta)}\left(c^2(2\theta-1)t + (\gamma-a)^{1-2\theta}\right)^{-(1-\theta)/(2\theta-1)},
& 1/2 < \theta < 1,
\\
\displaystyle
\frac{2}{c}\sqrt{\gamma-a}\exp(-c^2t/2),
&\theta = 1/2,
\end{cases}
\end{equation}
and $a, \gamma$ are constants such that $\gamma > a$ and
\[
a \leq \sE(v) \leq \gamma, \quad\forall\, v \in \sU.
\]
If in addition $u$ obeys Hypothesis \ref{hyp:Abstract_apriori_interior_estimate_trajectory_main_introduction}, then $u_\infty \in \sX$ and
\begin{equation}
\label{eq:Huang_3-45_X_introduction}
\|u(t+1) - u_\infty\|_\sX \leq 2C_1\Psi(t), \quad t\geq 0,
\end{equation}
where $C_1 \in [1,\infty)$ is the constant in Hypothesis \ref{hyp:Abstract_apriori_interior_estimate_trajectory_main_introduction} for $\delta=1$.
\end{mainthm}

One calls a critical point $\varphi \in \sU$ of $\sE$ a \emph{ground state} if $\sE$ attains its minimum value on $\sU$ at this point, that is,
\[
\sE(\varphi) = \inf_{u\in \sU}\sE(u).
\]
We have the following analogue of Huang \cite[Theorem 5.1.1]{Huang_2006}.

\begin{mainthm}[Existence and convergence of a global solution to a gradient system near a local minimum]
\label{mainthm:Huang_5-1-1_introduction}
(Compare Theorem \ref{thm:Huang_5-1-1}.)
Let $\sU$ be an open subset of a real Banach space, $\sX$, that is continuously embedded and dense in a Hilbert space, $\sH$. Let $\sE:\sU\subset \sX\to\RR$ be an analytic function with gradient map $\sE':\sU\subset \sX \to \sH$. Let $\varphi \in \sU$ be a ground state of $\sE$ on $\sU$ and suppose that $(c,\sigma,\theta)$ are the {\L}ojasiewicz-Simon constants for $(\sE,\varphi)$. Assume that
\begin{enumerate}
\item For each $u_0 \in \sU$, there exists a unique smooth solution to the Cauchy problem \eqref{eq:gradient_system} with $u(0)=u_0$, on a time interval $[0, \tau)$ for some positive constant, $\tau$;

\item Hypothesis \ref{hyp:Abstract_apriori_interior_estimate_trajectory_main_introduction} holds for smooth solutions to the gradient system \eqref{eq:gradient_system}; and

\item Given positive constants $b$ and $\eta$, there is a constant $\delta = \delta(\eta, \tau, b) \in (0, \tau]$ such that if $v$ is a smooth solution to the gradient system \eqref{eq:gradient_system} on $[0, \tau)$ with $\|v(0)\|_\sX \leq b$, then
\begin{equation}
\label{eq:Gradient_solution_near_initial_data_at_time_zero_for_short_enough_time_introduction}
\sup_{t\in [0, \delta]}\|v(t) - v(0)|_\sX < \eta.
\end{equation}
\end{enumerate}
Then there is a constant $\eps = \eps(c,C_1,\delta, \theta, \rho, \sigma, \tau, \varphi) \in (0, \sigma/4)$ with the following significance. For each $u_0 \in \sU_\eps$, the Cauchy problem \eqref{eq:gradient_system} with $u(0)=u_0$ admits a global smooth solution, $u:[0,\infty) \to \sU_{\sigma/2}$, that converges to a limit $u_\infty \in \sX$ as $t\to\infty$ with respect to the $\sX$ norm in the sense that
$$
\lim_{t \to \infty} \|u(t) - u_\infty\|_\sX = 0 \quad\hbox{and}\quad \int_1^\infty\|\dot u(t)\|_\sX\,dt < \infty.
$$
\end{mainthm}

Finally, we have the following analogue of Huang \cite[Theorem 5.1.2]{Huang_2006}.

\begin{mainthm}[Convergence to a critical point and stability of a ground state]
\label{mainthm:Huang_5-1-2_introduction}
(Compare Theorem \ref{thm:Huang_5-1-2}.)
Assume the hypotheses of Theorem \ref{mainthm:Huang_5-1-1_introduction}. Then, for each $u_0 \in \sU_\eps$, the Cauchy problem \eqref{eq:gradient_system} with $u(0)=u_0$ admits a global smooth solution $u:[0,\infty) \to \sU_{\sigma/2}$ that converges in $\sX$ as $t\to\infty$ to some critical point $u_\infty \in \sU_\sigma$. The critical point, $u_\infty$, satisfies $\sE(u_\infty) = \sE(\varphi)$. As an equilibrium of \eqref{eq:gradient_system} , the point $\varphi$ is Lyapunov stable (see Definition \ref{defn:Sell_You_page_32_Lyapunov_and_uniform_asymptotic_stability}). If $\varphi$ is isolated or a cluster point of the orbit $O(u)$, then $\varphi$ is uniformly asymptotically stable (see Definition \ref{defn:Sell_You_page_32_Lyapunov_and_uniform_asymptotic_stability}).
\end{mainthm}

It is appropriate at this point to recall an important existence result for critical points of functions on Banach spaces. We first review the well-known

\begin{defn}[Palais-Smale Condition]
\label{defn:Palais-Smale_Condition}
\cite[Section 4]{Mawhin_Willem_2010}
Let $\sE$ be a real-valued $C^1$ function on a Banach space, $\sX$. If $\{x_n\}_{n\in\NN} \subset \sX$ is a sequence such that $\{\sE(x_n)\}_{n\in\NN}$ is bounded and $\sE'(x_n) \to 0$ as $n \to \infty$, then $\{x_n\}_{n\in\NN}$ contains a convergent subsequence (whose limit is a critical point of $\sE$).
\end{defn}

Here, $\sE'(x) \in \sX'$ denotes the Fr\'echet derivative of $\sE$ at $x \in \sX$. The Palais-Smale Condition is closely related to the Condition (C) of Palais and Smale:

\begin{defn}[Condition C of Palais and Smale]
\label{defn:Condition_C_Palais_and_Smale}
\cite[Section 4]{Mawhin_Willem_2010}
Let $\sE$ be a real-valued $C^1$ function on a Banach space, $\sX$. If $S \subset X$ is a non-empty subset on which $\sE$ is bounded but $\|\sE'(\cdot)\|_{\sX'}$ is not bounded away from zero, then the closure of $S$ contains a critical point of $\sE$.
\end{defn}

The results of Section \ref{subsec:Main_results_gradient_system_Banach_space} apply regardless of whether the function $\sE:\sU\subset\sX\to\RR$ obeys the Palais-Smale Condition. For example, the Yang-Mills $L^2$-energy functional (when the base manifold has dimension four or greater) and harmonic map $L^2$-energy functional (when the source manifold has dimension two or greater) famously do \emph{not} obey the Palais-Smale Condition. However, when the Palais-Smale Condition \emph{does} hold, one has the celebrated \emph{Mountain Pass Theorem} due to Ambrosetti and Rabinowitz.

\begin{thm}[Mountain Pass Theorem]
\label{thm:Mountain_pass}
\cite{Ambrosetti_Rabinowitz_1973},
\cite[Theorem 3]{Mawhin_Willem_2010}
Let $\sE$ be a real-valued $C^1$ function on a Banach space, $\sX$, that obeys the Palais-Smale Condition and assume there exist $x,y \in \sX$ and $R > 0$ and $b \in \RR$ such that
\begin{gather*}
0 < R < \|x-y\|_\sX,
\\
\sE(w) \geq b > \max\{\sE(x), \sE(y)\}, \quad\forall\, w \in \sX \text{ such that } \|w - y\|_\sX = R.
\end{gather*}
If
\[
\Gamma := \{g \in C([0, 1]; \sX) : g(0) = x,\  g(1) = y\},
\]
then
\[
c = \inf_{g \in \Gamma} \sup_{t\in [0,1]} \sE(g(t)) \geq b
\]
is a critical value of $\sE$.
\end{thm}

\subsection{Preliminaries required for statements of main results on Yang-Mills gradient flow}
\label{subsec:Preliminaries_statements_main_results_Yang-Mills_gradient_flow}
Let $\ad P$ denote the real vector bundle associated to $P$ by the adjoint representation of $G$. Given a fixed $C^\infty$ reference connection, $A_1$, on a principal $G$-bundle, $P$, over a closed, connected, smooth manifold, $X$, of dimension $d \geq 2$ and Riemannian metric, $g$, we have covariant derivative operators, for any integer $i\geq 0$,
$$
\nabla_{A_1}:\Omega^i(X;\ad P) \to \Omega^{i+1}(X;\ad P),
$$
defined by the connection, $A_1$, on $P$ and Levi-Civita connection on $TX$ for the Riemannian metric $g$, where $\Omega^i(X;\ad P) := C^\infty(X;\Lambda^i\otimes\ad P)$ and we abbreviate the wedge product, $\Lambda^i(T^*X)$, by $\Lambda^i$ for any $i \geq 1$, so $\Lambda^1(T^*X) = T^*X$. For each $p \in [1,\infty]$ and integer $k\geq 0$, we define Sobolev spaces, $W_{A_1}^{k,p}(X;\Lambda^i\otimes\ad P)$, as completions of $\Omega^i(X;\ad P)$ with respect to the norms,
$$
\|a\|_{W_{A_1}^{k,p}(X)}
:=
\left(\sum_{j=0}^k \int_X |\nabla_{A_1}^j a|^p\,d\vol_g \right)^{1/p},
\quad \hbox{for } 1\leq p < \infty,
$$
and
$$
\|a\|_{W_{A_1}^{k,\infty}(X)}
:=
\sum_{j=0}^k \esssup_X|\nabla_{A_1}^j a|,
\quad\hbox{for } p = \infty,
$$
while for $p \in (1,\infty)$ and integer $k<0$, one uses Banach-space duality to define
$$
W_{A_1}^{-k,p}(X;\Lambda^i\otimes\ad P)
:=
\left(W_{A_1}^{k,p}(X;\Lambda^i\otimes\ad P)\right)'.
$$
More generally, one defines the Sobolev spaces, $W_{A_1}^{s,p}(X;\Lambda^i\otimes\ad P)$
for $s \in \RR$ and $p \in (1,\infty)$ via the fractional powers, $(\nabla_{A_1}^*\nabla_{A_1} + 1)^s$ when $s\geq 0$ and duality when $s<0$. When $p = 2$ and $s\in\RR$, we denote
$$
H_{A_1}^s(X;\Lambda^i\otimes\ad P) := W_{A_1}^{s,2}(X;\Lambda^i\otimes\ad P).
$$
Because the Yang-Mills gradient flow equation \eqref{eq:Yang-Mills_gradient_flow} is not parabolic and has a nonlinearity which is critical in dimension four, questions surrounding its solution are highly sensitive to the regularity of the initial data. Thus, for the sake of clarity and simplicity, we state our main results for Yang-Mills gradient flow in our Introduction in the case of an initial connection, $A_0$, of class $C^\infty$ and defer the corresponding statements for an initial connection in a Sobolev class to the indicated sections in our monograph, such as Section \ref{sec:Application_abstract_gradient_system_results_Yang-Mills_energy_functional}. As we shall see in Section \ref{subsec:Yang-Mills_gradient_flow_near_local_minimum}, we recover the previous results on Yang-Mills gradient flow due to Daskalopoulos \cite{Daskalopoulos_1992} when $d=2$ and to R\r{a}de \cite{Rade_1992} when $d=2,3$, by different methods.

\subsection{Yang-Mills gradient flow near critical points: base manifolds of arbitrary dimension}
\label{subsec:Yang-Mills_gradient_flow_near_local_minimum}
We describe the key results for Yang-Mills gradient flow started near a local minimum, $A_{\min}$, for the Yang-Mills energy functional and then describe convergence properties for Yang-Mills gradient flow near an arbitrary critical point, $A_\ym$, that is known to be a cluster point of a Yang-Mills gradient flow line, $\{A(t): t \geq 0\}$, whose existence is already established. When $X$ has dimension four, absolute minima for the Yang-Mills energy functional are provided by $g$-anti-self-dual connections when $P$ has non-positive Pontrjagin numbers and $g$-self-dual connections when $P$ has non-negative Pontrjagin numbers (see Section \ref{sec:Taubes_1982_Appendix} for a discussion of the classification of principal $G$-bundles and absolute minima of the Yang-Mills energy functional in dimension four). Because the signs of the Pontrjagin numbers are simply reversed when one reverses the orientation of $X$, there is no loss in generality if we restrict, as we shall henceforth when $X$ has dimension four, our attention to principal $G$-bundles with non-positive Pontrjagin numbers and absolute minima of the Yang-Mills energy functional correspond to $g$-anti-self-dual connections. The results in this subsection slightly specialize those in Section \ref{sec:Application_abstract_gradient_system_results_Yang-Mills_energy_functional}. If $A$ is a connection on $P$, we let $[A]$ denote its gauge-equivalence class in the quotient space, $\sB(P)$. Recall that $\Aut P$, the group of $C^\infty$ automorphisms (or gauge transformations) of $P$ may be equivalently viewed as the group of $C^\infty$ sections of $\Ad P$ under fiberwise multiplication \cite[Section 3.1]{FrM}.

\begin{mainthm}[Global existence and convergence of Yang-Mills gradient flow near a local minimum]
\label{mainthm:Yang-Mills_gradient_flow_global_existence_and_convergence_started_near_minimum}
Let $G$ be a compact Lie group and $P$ a principal $G$-bundle over a closed, connected, oriented, smooth manifold, $X$, of dimension $d\geq 2$ and with Riemannian metric, $g$. Let $A_1$ and $A_{\min}$ be $C^\infty$ connections on $P$, with $A_{\min}$ being a local minimum, and let
\begin{inparaenum}[\itshape a\upshape)]
\item $k=1$ and $p = 2$ if $2\leq d\leq 4$, or
\item $k=2$ and $p > d/2$ if $d \geq 5$.
\end{inparaenum}
Then there are constants\footnote{These are the {\L}ojasiewicz-Simon constants for the Yang-Mills energy functional --- see Theorem \ref{thm:Huang_3-3-6_Yang-Mills}.}
$c \in [1,\infty)$, and $\sigma \in (0,1]$, and $\theta \in [1/2,1)$,
depending on $(A_1, A_{\min},g,p)$, with the following significance.
\begin{enumerate}
\item \emph{Global existence:} There is a constant $\eps \in (0,\sigma/4)$, depending on $(A_1, A_{\min},g,p)$, with the following significance. If $A_0$ is a $C^\infty$ connection on $P$ such that
$$
\|A_0 - A_{\min}\|_{W_{A_1}^{k,p}(X)} < \eps,
$$
then there exists a solution, $A(t) = A_0 + a(t)$ for $t\in [0,\infty)$, with
$$
a \in C^\infty([0,\infty)\times X; \Lambda^1\otimes\ad P),
$$
to the Yang-Mills gradient flow \eqref{eq:Yang-Mills_gradient_flow} with initial data, $A(0) = A_0$, and
$$
\|A(t) - A_{\min}\|_{W_{A_1}^{k,p}(X)} < \sigma/2, \quad\forall\, t \in [0,\infty).
$$

\item \emph{Dependence on initial data:} The solution, $A(t)$ for $t \in [0,\infty)$, varies continuously with respect to $A_0$ in the $C_{\loc}([0,\infty); W_{A_1}^{k,p}(X;\Lambda^1\otimes\ad P))$ topology and, more generally, smoothly for all non-negative integers, $l,m$, in the $C_{\loc}^l([0,\infty); H_{A_1}^m(X;\Lambda^1\otimes\ad P))$ topology.

\item \emph{Convergence:} As $t\to\infty$, the flow, $A(t)$, converges strongly with respect to the norm on $W_{A_1}^{k,p}(X;\Lambda^1\otimes\ad P)$ to a Yang-Mills connection, $A_\infty$, of class $C^\infty$ on $P$, and the gradient-flow line has finite length in the sense that
$$
\int_0^\infty \left\|\frac{\partial A}{\partial t}\right\|_{W_{A_1}^{k,p}(X)}\,dt < \infty.
$$
If $A_{\min}$ is a cluster point of the orbit, $O(A) = \{A(t):t\geq 0\}$, then $A_\infty = A_{\min}$.

\item \emph{Convergence rate:} For all $t \geq 1$,
\begin{multline}
\label{eq:Rade_Proposition_7-4_convergence_rate_introduction}
\|A(t) - A_\infty\|_{W_{A_1}^{k,p}(X)}
\\
\leq
\begin{cases}
\displaystyle
\frac{1}{c(1-\theta)}\left(c^2(2\theta-1)(t-1) + (\sE(A_0)-\sE(A_\infty))^{1-2\theta}\right)^{-(1-\theta)/(2\theta-1)},
& 1/2 < \theta < 1,
\\
\displaystyle
\frac{2}{c}\sqrt{\sE(A_0)-\sE(A_\infty)}\exp(-c^2(t-1)/2),
&\theta = 1/2.
\end{cases}
\end{multline}

\item \emph{Stability:} As an equilibrium of the Yang-Mills gradient flow \eqref{eq:Yang-Mills_gradient_flow}, the point $A_\infty$ is Lyapunov stable; if $A_\infty$ is isolated or a cluster point of the orbit $O(A)$, then $A_\infty$ is uniformly asymptotically stable.

\item \emph{Uniqueness:} Any two solutions are equivalent modulo a path of gauge transformations,
$$
u \in C^\infty([0,\infty)\times X; \Ad P), \quad u(0) = \id_P.
$$
\end{enumerate}
\end{mainthm}

\begin{rmk}[On the role of the $C^\infty$ reference connection in Theorem \ref{mainthm:Yang-Mills_gradient_flow_global_existence_and_convergence_started_near_minimum}]
\label{rmk:Role_smooth_reference_connection}
While the $C^\infty$ connections $A_1$ (used to define Sobolev norms) and $A_{\min}$ (the critical point) play very distinct roles in the statement of Theorem \ref{mainthm:Yang-Mills_gradient_flow_global_existence_and_convergence_started_near_minimum}, one could simplify the statement slightly by choosing $A_1 = A_{\min}$ with little loss of generality. While we have assumed that the connection, $A_0$ on $P$, defining the initial data for the Yang-Mills gradient flow is $C^\infty$, one could allow the initial data, $A_0$, more generally to belong to a Sobolev class, as we do in the body of our monograph.
\end{rmk}

It is worth noting that Theorem \ref{mainthm:Yang-Mills_gradient_flow_global_existence_and_convergence_started_near_minimum} does not contradict an important result due to Naito:

\begin{thm}[Finite-time blow-up for Yang-Mills gradient flow over a sphere of dimension greater than or equal to five]
\label{thm:Naito_1-3}
\cite[Theorem 1.3]{Naito_1994}
Let $G \subset\SO(n)$ be a compact Lie group and $P$ a non-trivial principal $G$-bundle over the sphere, $S^d$, of dimension $d\geq 5$ and with its standard round Riemannian metric of radius one. Then there is a constant $\eps_1 \in (0,1]$ with the following significance. If $A_0$ is a $C^\infty$ connection on $P$ such that $\|F_{A_0}\|_{L^2(S^d)} < \eps_1$, then the solution $A(t)$ to Yang-Mills gradient flow \eqref{eq:Yang-Mills_gradient_flow} with initial data, $A(0) = A_0$, blows up in finite time.
\end{thm}

Since $P$ in Theorem \ref{thm:Naito_1-3} is a \emph{non-trivial} bundle over $S^d$, it does not support a product connection (the only flat connection on $P$ since the base manifold is simply connected). On the other hand, by an energy gap result due to Nakajima \cite[Corollary 1.2]{Nakajima_1987} (compare \cite[Theorem 1]{Feehan_yangmillsenergygapflat}), for small enough $\eps_0 \in (0,1]$, if $A_\ym$ is a Yang-Mills connection on $P$ with $\|F_{A_\ym}\|_{L^2(S^d)} < \eps_0$, then $A_\ym$ is necessarily flat and therefore is the product connection, which is impossible when $P$ is non-trivial. Hence, Theorem \ref{mainthm:Yang-Mills_gradient_flow_global_existence_and_convergence_started_near_minimum} does not apply in the setting of Theorem \ref{thm:Naito_1-3} since it is not possible to choose a Yang-Mills connection, $A_{\min}$, on a non-trivial $G$-bundle $P$ over $S^d$ and initial data, $A_0$, such that \emph{both} $\|A_0 - A_{\min}\|_{W_{A_1}^{k,p}(S^d)}$ \emph{and} $\|F_{A_0}\|_{L^2(S^d)}$ are small because this would force $\|F_{A_{\min}}\|_{L^2(S^d)}$ to also be small, and thus force $A_{\min}$ to be the product connection.

In view of the crucial role played by local minima of the Yang-Mills energy functional in Theorem 1, it is important to review some results regarding their \emph{non-existence} due to Bourguignon, Lawson, and Simons \cite{Bourguignon_Lawson_1981, Bourguignon_Lawson_Simons_1979}; in Section \ref{sec:Critical_points_Yang-Mills_energy_functional}, we survey some of the results regarding existence of non-minimal Yang-Mills connections. A Yang-Mills connection $A_\ym$ is \emph{weakly stable} if the Hessian $\sE''(A_\ym)$ of the Yang-Mills $L^2$-energy functional is a non-negative operator \cite[Equation (1.1)]{Bourguignon_Lawson_1981}.

\begin{thm}
\label{thm:Bourguignon_Lawson_A}
\cite[Theorem A]{Bourguignon_Lawson_1981}
Let $G$ be a compact Lie group, $S^d$ be the sphere of dimension $d$ with its standard round metric of radius one, and $P$ a principal $G$-bundle over $S^d$. If $d \geq 5$, then there are no weakly stable Yang-Mills connections on $P$.
\end{thm}

\begin{thm}
\label{thm:Bourguignon_Lawson_B}
\cite[Theorem B]{Bourguignon_Lawson_1981}
Let $G=\SU(2)$, $\SU(3)$, or $\U(2)$, and $S^4$ be the four-dimensional sphere with its standard round metric of radius one, and $P$ a principal $G$-bundle over $S^4$. If $A$ is a weakly stable Yang-Mills connection on $P$, then $A$ is anti-self-dual, self-dual, or Abelian.
\end{thm}

\begin{thm}
\label{thm:Bourguignon_Lawson_Bprime}
\cite[Theorem B${}'$]{Bourguignon_Lawson_1981}
Let $G=\SU(2)$, and $X$ be a closed, four-dimensional, homogenous, orientable, Riemannian manifold, and $P$ a principal $G$-bundle over $X$. If $A$ is a weakly stable Yang-Mills connection on $P$, then $A$ is anti-self-dual, self-dual, or Abelian.
\end{thm}

A Yang-Mills connection $A_\ym$ is \emph{strictly stable} if the Hessian $\sE''(A_\ym)$ of the Yang-Mills $L^2$-energy functional is a strictly positive operator on restriction to a Coulomb-gauge slice through $A_\ym$ \cite[p. 190]{Bourguignon_Lawson_1981}. If $S^d/\Gamma$, for $d \geq 4$, is a non-trivial quotient (for example, real projective space), then there exist strictly stable Yang-Mills connections on principal $G$-bundles over $S^d/\Gamma$ \cite[p. 190 and Section 9]{Bourguignon_Lawson_1981}.

Further information regarding the behavior of Yang-Mills gradient flow near an arbitrary critical point can sometimes be deduced from the following analogue of Theorem \ref{mainthm:Huang_3-3-6_introduction}.

\begin{mainthm}[Convergence of a subsequence implies convergence for a solution to Yang-Mills gradient flow]
\label{mainthm:Yang-Mills_gradient_flow_global_with_critical_point_in_orbit_closure}
Let $G$ be a compact Lie group and $P$ a principal $G$-bundle over a closed, connected, oriented, smooth manifold, $X$, of dimension $d\geq 2$ and with Riemannian metric, $g$. Let $A_1$ and $A_\ym$ be $C^\infty$ connections on $P$, with $A_\ym$ being a Yang-Mills connection, and let
\begin{inparaenum}[\itshape a\upshape)]
\item $k=1$ and $p = 2$ if $2\leq d\leq 4$, or
\item $k=2$ and $p > d/2$ if $d \geq 5$.
\end{inparaenum}
If $A(t) = A_1+a(t)$, for $t \in [0,\infty)$, with $a\in C^\infty([0,\infty); W_{A_1}^{k,p}(X;\Lambda^1\otimes\ad P))$, is a solution to Yang-Mills gradient flow \eqref{eq:Yang-Mills_gradient_flow} and $[A_\ym]$ is a cluster point of the orbit $O([A]) = \{[A(t)]: t\geq 0\}$, in the sense that there exist sequences of times, $\{t_m\}_{m\in\NN} \subset [0,\infty)$ with $t_m\to\infty$ as $m\to\infty$, and $W^{k+1,p}$ gauge transformations, $\{\Phi_m\}_{m\in\NN} \subset \Aut(P)$, such that
\[
\Phi_m^*A(t_m) \to A_\ym \quad\text{in } W_{A_1}^{k,p}(X;\Lambda^1\otimes\ad P) \text{ as } m \to \infty,
\]
then there exists a $W^{k+1,p}$ gauge transformation, $\Phi \in \Aut(P)$, such that $A(t)$ converges to $\tilde A_\ym = \Phi^*A_\ym$ as $t\to\infty$ in the sense that
\[
\lim_{t\to\infty}\|A(t)-\tilde A_\ym\|_{W_{A_1}^{k,p}(X)} = 0
\quad\hbox{and}\quad
\int_0^\infty \|\dot A\|_{W_{A_1}^{k,p}(X)}\,dt < \infty.
\]
\end{mainthm}

\subsection{Yang-Mills gradient flow for an initial connection with almost minimal energy: four-dimensional base manifolds}
\label{subsec:Yang-Mills_gradient_flow_initial_connection_almost_minimal_energy}
It is more useful in applications if one can recast Theorem \ref{mainthm:Yang-Mills_gradient_flow_global_existence_and_convergence_started_near_minimum} in the context of an initial connection, $A_0$, which is close to a minimum in an energy rather than a norm sense. Recall that, when $X$ has dimension four, one may define (see \cite[Sections 1.1.5 and 2.1.3]{DK}) the $g$-self-dual and $g$-anti-self-dual components,
$$
F_A^{\pm,g} = \frac{1}{2}(1\pm*_g)F_A \in \Omega^\pm(X; \ad P),
$$
of the curvature, $F_A$, of a connection, $A$, on $P$,
$$
F_A \in \Omega^2(X; \ad P) = \Omega^+(X; \ad P)\oplus \Omega^-(X; \ad P).
$$
Let $d_A^+$ denote the composition of the covariant exterior derivative, $d_A:\Omega^1(X; \ad P) \to \Omega^2(X; \ad P)$, with the projection $\Omega^2(X; \ad P) \to \Omega^+(X; \ad P)$, so $d_A^+ = \frac{1}{2}(1\pm*)d_A$. The second-order partial differential operator, $d_A^+d_A^{+,*}$ on $\Omega^+(X; \ad P)$, is self-adjoint and is well-known to be elliptic (see, for example, \cite{AHS, DK, FU, TauSelfDual}), with a discrete spectrum of non-negative, real eigenvalues. Let $\mu(A)$ denote the least eigenvalue of the Laplace operator, $d_A^+d_A^{+,*}$ on $\Omega^+(X; \ad P)$.

If $\cG(\cdot,\cdot)$ denotes the Green kernel of the Laplace operator, $d^*d$, on $\Omega^2(X)$, we define (see Section \ref{subsec:Feehan_slice_4_and_5} for further details),
\begin{align*}
\|v\|_{L^\sharp(X)} &:= \sup_{x\in X}\int_X \cG(x,y) |v|(y)\,d\vol_g(y),
\\
\|v\|_{L^{\sharp,2}(X)} &:= \|v\|_{L^\sharp(X)} + \|v\|_{L^2(X)},
\quad\forall\, v \in \Omega^2(X;\ad P).
\end{align*}
We recall that $\cG(x,y)$ has a singularity comparable with $\dist_g(x,y)^{-2}$, when $x,y\in X$ are close \cite{Chavel}. The norm $\|v\|_{L^2(X)}$ is conformally invariant and $\|v\|_{L^\sharp(X)}$ is scale invariant. One can show that $\|v\|_{L^\sharp(X)} \leq c_p\|v\|_{L^p(X)}$ for every $p>2$, where $c_p$ depends at most on $p$ and the Riemannian metric, $g$, on $X$. See \cite{FeehanSlice} or \cite{TauGluing} and references cited therein for further discussion of this `critical-exponent' norm.

\begin{maincor}[Global existence and convergence of Yang-Mills gradient flow for an initial connection with almost minimal energy and positive eigenvalue]
\label{maincor:Yang-Mills_gradient_flow_global_existence_and_convergence_started_near_connection with almost minimal energy_and_vanishing_H2+}
Let $G$ be a compact Lie group, $P$ a principal $G$-bundle over a closed, connected, four-dimensional, oriented, smooth manifold, $X$, with Riemannian metric, $g$, and $E_0, \mu_0 \in (0,\infty)$ constants. Then there are constants, $C_0 = C_0(E_0,g,\mu_0)\in (0,\infty)$ and $\eta = \eta(E_0,g,\mu_0) \in (0,1]$, with the following significance. If $A_0$ is a $C^\infty$ connection on $P$ such that
\begin{subequations}
\begin{align}
\mu(A_0) &\geq \mu_0,
\\
\|F_{A_0}^+\|_{L^{\sharp,2}(X)} &\leq \eta,
\\
\|F_{A_0}\|_{L^2(X)} &\leq E_0,
\end{align}
\end{subequations}
then there is a $g$-anti-self-dual connection, $A_\asd$ on $P$, of class $C^\infty$ such that
$$
\|A_\asd - A_0\|_{H_{A_0}^1(X)} \leq C_0\|F_{A_0}^+\|_{L^{\sharp,2}(X)},
$$
and the following holds. Let $c>0$, $\sigma>0$, and $\theta\in [1/2,1)$ be the constants in Theorem \ref{mainthm:Yang-Mills_gradient_flow_global_existence_and_convergence_started_near_minimum} defined by $([A_1],[A_\asd],g)$. For any $C^\infty$ connection, $A_0'$ on $P$, there is a constant $\eps \in (0,\sigma/4)$, depending on
$([A_1],[A_\asd],\|A_0' - A_1\|_{H_{A_1}^1(X)}, g)$, such that the conclusions of Theorem \ref{mainthm:Yang-Mills_gradient_flow_global_existence_and_convergence_started_near_minimum} hold for the Yang-Mills gradient flow \eqref{eq:Yang-Mills_gradient_flow} with initial data, $A(0) = A_0'$, and where $A_{\min} = A_\asd$.
\end{maincor}

\begin{rmk}[Selection of the $g$-anti-self-dual connection in Corollary \ref{maincor:Yang-Mills_gradient_flow_global_existence_and_convergence_started_near_connection with almost minimal energy_and_vanishing_H2+}]
\label{rmk:Selection_g-ASD_connection_near_A0}
The $g$-anti-self-dual connection, $A_\asd$ on $P$ in Corollary \ref{maincor:Yang-Mills_gradient_flow_global_existence_and_convergence_started_near_connection with almost minimal energy_and_vanishing_H2+}, may be constructed as the unique solution,
$A_\asd = A_0+d_{A_0}^{+,*}v$, to the $g$-anti-self-dual equation, $F^{+,g}(A_0+d_{A_0}^{+,*}v)=0$ on $X$, with $v \in \Omega^+(X;\ad P)$ obeying
$$
\|v\|_{H_{A_0}^2(X)} + \|v\|_{C(X)} \leq C_0\|F_{A_0}^+\|_{L^{\sharp,2}(X)}.
$$
See Theorem \ref{thm:Proposition_Feehan_Leness_7-6} in the sequel.
\end{rmk}

Corollary \ref{maincor:Yang-Mills_gradient_flow_global_existence_and_convergence_started_near_connection with almost minimal energy_and_vanishing_H2+} is a consequence of Theorem \ref{mainthm:Yang-Mills_gradient_flow_global_existence_and_convergence_started_near_minimum} and existence of solutions to the anti-self-dual equation (Theorem \ref{thm:Proposition_Feehan_Leness_7-6}). The hypothesis\footnote{This is equivalent to $\Coker d_{A_0}^{+,g} = 0$, where $\Coker d_{A_0}^{+,g} := \Omega^+(X;\ad P)/\Ran d_{A_0}^{+,g}$.}
$\mu(A_0)>0$ in Corollary \ref{maincor:Yang-Mills_gradient_flow_global_existence_and_convergence_started_near_connection with almost minimal energy_and_vanishing_H2+} is easily achieved in many situations of practical interest. To describe those, we recall some facts concerning `positive' or `generic' metrics, respectively.

\begin{defn}[Good Riemannian metric]
\label{defn:Good_Riemannian_metric_introduction}
Let $G$ be a compact Lie group, $X$ a closed, connected, four-dimensional, oriented, smooth manifold, and $\eta \in H^2(X;\pi_1(G))$ an obstruction class. We say that a Riemannian metric, $g$, on $X$ is \emph{good} if for every principal $G$-bundle, $P$, over $X$ with $\eta(P) = \eta$ and non-positive Pontrjagin vector, $\bkappa(P)$, and every connection, $A$, on $P$ of class $H^1$ with $F_A^{+,g}=0$, then $\Coker d_A^{+,g} = 0$.
\end{defn}

We refer the reader to Section \ref{sec:Taubes_1982_Appendix} for the classification of principal $G$-bundles, in terms of Pontrjagin vectors and obstruction classes, over closed, connected, four-dimensional manifolds.

A Riemannian metric, $g$, on $X$ is good in the sense of Definition \ref{defn:Good_Riemannian_metric_introduction} if
\begin{inparaenum}[\itshape a\upshape)]
\item $G$ is an arbitrary compact Lie group and $g$ is \emph{positive} in the sense of \eqref{eq:Freed_Uhlenbeck_page_174_positive_metric} (see Lemma \ref{lem:Positive_metric_implies_positive_lower_bound_small_eigenvalues} in the sequel) or,
\item $G$ is $\SU(2)$ or $\SO(3)$ and $g$ is \emph{generic} in the sense of Freed and Uhlenbeck (see, for example, Theorem \ref{thm:Donaldson-Kronheimer_Corollary_4-3-15_and_18_and_Proposition_4-3-20} in the sequel).
\end{inparaenum}
In particular, Corollary \ref{maincor:Yang-Mills_gradient_flow_global_existence_and_convergence_started_near_connection with almost minimal energy_and_vanishing_H2+} yields the following

\begin{maincor}[Global existence and convergence of Yang-Mills gradient flow for an initial connection with almost minimal energy and a good Riemannian metric]
\label{maincor:Yang-Mills_gradient_flow_global_existence_and_convergence_started_near_connection_with_almost_minimal_energy_and_g_is_good_Lp_version}
Let $G$ be a compact Lie group, $P$ a principal $G$-bundle over a closed, connected, four-dimensional, oriented, smooth manifold, $X$, with \emph{good} Riemannian metric, $g$, in the sense of Definition \ref{defn:Good_Riemannian_metric_introduction}, and $K \in (0,\infty)$ and $p \in (2,4)$ constants, and $A_1$ a $C^\infty$ reference connection on $P$. Then there are constants, $\delta = \delta(g,p) \in (0,1]$ and $\eta = \eta([A_1],g,K,P,p) \in (0,1]$, with the following significance. If $A_0$ is a $C^\infty$ connection on $P$ such that
\begin{subequations}
\label{eq:FA0+_Lp_small_and_FA0_Lp_bounded}
\begin{align}
\label{eq:FA0+_Lp_small}
\|F_{A_0}^+\|_{L^p(X)} &\leq \eta,
\\
\label{eq:FA0_Lp_bounded}
\|F_{A_0}\|_{L^p(X)} &\leq K,
\end{align}
\end{subequations}
and
\begin{equation}
\label{eq:A0_minus_A1_W14_small}
\|A_0 - A_1\|_{W_{A_0}^{1,4}(X)} \leq \delta,
\end{equation}
then the conclusions of Theorem \ref{mainthm:Yang-Mills_gradient_flow_global_existence_and_convergence_started_near_minimum} hold for the Yang-Mills gradient flow \eqref{eq:Yang-Mills_gradient_flow} with initial data, $A(0) = A_0$.
\end{maincor}

\begin{rmk}[Discreteness of critical values implied by the {\L}ojasiewicz-Simon gradient inequality and anti-self-duality of the limiting Yang-Mills connection]
\label{rmk:LS_gradient_inequality_and_ASD_YM_limit}
The {\L}ojasiewicz-Simon gradient inequality for the Yang-Mills energy functional (Theorem \ref{thm:Rade_proposition_7-2}) implies that there are constants $c \in (0,\infty)$, $\sigma \in (0,1]$, and $\theta \in [1/2,1)$, depending on the triple $([A_\infty],g,P)$, with the following significance: If $A$ is an $H^1$ connection on $P$ such that $\|A - A_\infty\|_{H^1_{A_\infty}(X)} < \sigma$, then
$$
\|\sE'(A)\|_{L^2(X)} \geq c|\sE(A) - \sE(A_\infty)|^\theta.
$$
In particular, if there is a $g$-anti-self-dual connection $A_{\asd}$ on $P$ such that $\|A_{\asd} - A_\infty\|_{H^1_{A_\infty}(X)} < \sigma$, then we can apply the preceding gradient inequality with $A = A_\asd$ to conclude that $\sE(A_\infty) = \sE(A_{\asd})$, so the energy $\sE(A_\infty)$ is minimal and $A_\infty$ is also $g$-anti-self-dual.
\end{rmk}

\begin{rmk}[Anti-self-duality of the limiting Yang-Mills connection for a positive Riemannian metric]
\label{rmk:Positive_Riemannian_metric_and_ASD_YM_limit}
Our hypothesis \eqref{eq:FA0+_Lp_small} that $\|F_{A_0}^+\|_{L^p(X)} \leq \eta$ and Lemma \ref{lem:Nonincreasing_ASD_curvature_Yang-Mills_gradient_flow}, which implies
$$
\|F_{A(t)}^+\|_{L^2(X)} \leq \sqrt{2}\|F_{A_0}^+\|_{L^2(X)}, \quad\forall\, t \geq 0,
$$
ensures that $\|F_{A_\infty}^+\|_{L^2(X)} < \sqrt{2}c_p\eta$. Hence, the $L^2$-isolation theorem of Min-Oo \cite[Theorem 2]{Min-Oo_1982} also yields the conclusion that the limit, $A_\infty$, is $g$-anti-self-dual when $g$ is positive.
\end{rmk}

It is possible to relax the hypothesis, $\|F_{A_0}^+\|_{L^p(X)} \leq \eta$, in Corollary \ref{maincor:Yang-Mills_gradient_flow_global_existence_and_convergence_started_near_connection_with_almost_minimal_energy_and_g_is_good_Lp_version}, to the scale invariant condition, $\|F_{A_0}^+\|_{L^{\sharp,2}(X)} \leq \eta$, together with a cumbersome auxiliary hypothesis.

\begin{maincor}[Global existence and convergence of Yang-Mills gradient flow for an initial connection with almost minimal energy and a good Riemannian metric]
\label{maincor:Yang-Mills_gradient_flow_global_existence_and_convergence_started_near_connection_with_almost_minimal_energy_and_g_is_good_L2sharp_version}
Assume the hypotheses of Corollary \ref{maincor:Yang-Mills_gradient_flow_global_existence_and_convergence_started_near_connection_with_almost_minimal_energy_and_g_is_good_Lp_version}, except replace the hypothesis \eqref{eq:FA0+_Lp_small} by
\begin{subequations}
\begin{align}
\label{eq:FA0+_Lsharp2_small}
\|F_{A_0}^+\|_{L^{\sharp,2}(X)} &\leq \eta,
\\
\label{eq:FAasd_Lp_bounded}
\|F(A_\asd)\|_{L^p(X)} &\leq K,
\end{align}
\end{subequations}
where $A_\asd = A_0+d_{A_0}^{+,*}v$ is the unique solution to the $g$-anti-self-dual equation, $F^{+,g}(A_0+d_{A_0}^{+,*}v)=0$, with $v \in \Omega^+(X;\ad P)$ obeying
$$
\|v\|_{H_{A_0}^2(X)} + \|v\|_{C(X)} \leq C_0\|F_{A_0}^+\|_{L^{\sharp,2}(X)}.
$$
Then the conclusions of Corollary \ref{maincor:Yang-Mills_gradient_flow_global_existence_and_convergence_started_near_connection_with_almost_minimal_energy_and_g_is_good_Lp_version} continue to hold.
\end{maincor}

\begin{rmk}[Dependence of {\L}ojasiewicz-Simon constants on the critical point]
\label{rmk:Dependence_LS_constants_precompact_subset}
Corollary \ref{maincor:Yang-Mills_gradient_flow_global_existence_and_convergence_started_near_connection_with_almost_minimal_energy_and_g_is_good_L2sharp_version} should hold without the auxiliary hypothesis \eqref{eq:FAasd_Lp_bounded}. The inclusion of the condition \eqref{eq:FAasd_Lp_bounded} originates from the facts that
\begin{inparaenum}[\itshape a\upshape)]
\item the moduli space, $M(P,g)$, of $g$-anti-self-dual connections is (except in rare cases) non-compact (as a topological space equipped with its usual Uhlenbeck topology); and
\item the triple of {\L}ojasiewicz-Simon constants (namely, the \emph{proportionality constant}, $c>0$, \emph{exponent}, $\theta \in (1/2,1]$, and \emph{radius}, $\sigma > 0$) depend on the critical point, $[A_\asd] \in M(P,g)$.
\end{inparaenum}
The hypotheses \eqref{eq:FA0+_Lp_small} in Corollary \ref{maincor:Yang-Mills_gradient_flow_global_existence_and_convergence_started_near_connection_with_almost_minimal_energy_and_g_is_good_Lp_version} or \eqref{eq:FAasd_Lp_bounded} in Corollary \ref{maincor:Yang-Mills_gradient_flow_global_existence_and_convergence_started_near_connection_with_almost_minimal_energy_and_g_is_good_L2sharp_version} allow us to select a precompact open subset of $M(P,g)$, depending only on $(g,K,P,p)$.

In addition to R\r{a}de in \cite{Rade_1992, RadeThesis}, the second early application of the {\L}ojasiewicz-Simon gradient inequality in a gauge-theory setting is due to Morgan, Mrowka, and Ruberman \cite{MMR}, in the context of the Chern-Simons functional on the affine space of connections on a principal $G$-bundle, $Q$, over a closed, smooth three-dimensional manifold, $Y$, and the question of dependence of the {\L}ojasiewicz-Simon constants on the critical points also arises in their application. However, the critical points of the Chern-Simons functional are \emph{flat} connections and the moduli space, $M(Q)$, of flat connections on $Q$ is compact. Thus, it suffices to cover $M(Q)$ by finitely many `Simon coordinate patches' (see \cite[Definition 4.3.2]{MMR}) defined in turn by finitely many gauge-equivalence classes of flat connections, with their associated triples of {\L}ojasiewicz-Simon constants.

The proof of the {\L}ojasiewicz-Simon gradient inequality (Theorem \ref{thm:Rade_proposition_7-2}) for the Yang-Mills energy functional, $\sE(A)$, requires a splitting of $\sE$ into a finite and infinite-dimensional part using $L^2$-orthogonal projection onto the kernel of the Hodge Laplace operator, $\Delta_A = d_A^*d_A + d_Ad_A^*$ on $\Omega^1(X;\ad P)$, and its orthogonal complement, with the infinite-dimensional gradient inequality then being deduced from the finite-dimensional version due to {\L}ojasiewicz \cite{Lojasiewicz_1963, Lojasiewicz_1965, Lojasiewicz_1984, Simon_1996}). (For example, see the direct proofs due to R\r{a}de of his \cite[Proposition 7.2]{Rade_1992} when $X$ has dimension $2$ or $3$.) This splitting is reminiscent of the method employed by Kuranishi \cite{Kuranishi} in his analysis of complex structures and adapted by Taubes in \cite{TauIndef, TauFrame}, Donaldson \cite{DonConn}, and Donaldson and Kronheimer \cite{DK} in their construction of open neighborhoods of points in the Uhlenbeck boundary, $\partial M(P,g) = \bar M(P,g) \less M(P,g)$. The behavior of the eigenvalues of $\Delta_A$ as $[A]$ converges to a point in $\partial M(P,g)$ is analyzed by Taubes in \cite{TauFrame}.

It is not unreasonable to expect that $M(P,g)$ has a finite cover by suitably defined Simon coordinate patches, as in \cite{MMR}. Some evidence for this is provided by Corollary \ref{cor:Rade_proposition_7-2_dimension_four}, which specializes Theorem \ref{thm:Rade_proposition_7-2} to a manifold, $X$, of dimension four and where we observe that the {\L}ojasiewicz-Simon radius, $\sigma$, and exponent, $\theta$, depend only on the conformal equivalence class, $[g]$, of the Riemannian metric on $X$. In Corollary \ref{cor:Rade_proposition_7-2_S4}, we further specialize Theorem \ref{thm:Rade_proposition_7-2} to the case of $S^4$ with its standard round metric of radius one and note that if $c$ is the {\L}ojasiewicz-Simon proportionality constant for a Yang-Mills connection, $A_\ym$, with scale one in the sense of \cite[Equation (3.10)]{FeehanGeometry} or \cite[p. 343]{TauFrame}, then $c\lambda$ is the {\L}ojasiewicz-Simon proportionality constant for the Yang-Mills connection, $f_\lambda^*A_\ym$, with scale $\lambda > 0$, where $f_\lambda$ is the conformal diffeomorphism of $S^4$ defined by the rescaling map $\RR^4 \ni x \mapsto x/\lambda \in \RR^4$.
\end{rmk}


\subsection{Yang-Mills gradient flow for an initial connection with arbitrary energy: complex surfaces}
\label{subsec:Yang-Mills_gradient_flow_initial_connection_arbitrary_energy}
Finally, when we relax the condition that the initial data, $A_0$, have suitably small $F_{A_0}^+$, we retain global existence but may no longer have convergence without energy loss at $t=\infty$. To set our own results and ensuing remarks in context, we first describe a special case of the global existence and convergence results due to Donaldson \cite{DonASD, DonInfDet, DK}. We recall \cite[Section 6.1.4]{DK} that a unitary connection, $A$, on a Hermitian vector bundle, $E$, over a compact, complex, Hermitian surface, $(X,h)$, is \emph{Hermitian Yang-Mills} if
\begin{equation}
\label{eq:Donaldson_Kronheimer_6-1-4_HYM_connection_equation}
i\widehat F_A = \lambda\id_E,
\end{equation}
where $\widehat F_A := \langle F_A, \omega\rangle$ (pointwise inner product), and $\deg E := \langle c_1(E)\smile\omega,[X]\rangle$, and $\omega$ denotes the K\"ahler form on $X$, and
\begin{equation}
\label{eq:Donaldson_Kronheimer_6-1-4_HYM_connection_lambda}
\lambda := \frac{2\pi}{\Vol_h(X)}\frac{\deg E}{\rank_\CC E}.
\end{equation}
See, for example, Kobayashi \cite{Kobayashi} or L\"ubke and Teleman \cite{Lubke_Teleman_1995} for additional information on Hermitian Yang-Mills connections. That understood, for the sake of consistency with the remainder of our monograph, which focuses on Yang-Mills rather than Hermitian Yang-Mills connections more generally, we shall only describe Donaldson's results for the case where $E$ is a complex rank two, Hermitian vector bundle of degree zero (so $G=\SU(2)$) and $X$ is a compact, K\"ahler surface; we refer to Section \ref{subsec:Preliminaries_curvature_1-1_Kaehler_identities_Bochner} for definitions and notation in this context, which closely follows those of \cite{DK}. We recall that \cite[Definition 6.1.3]{DK} a holomorphic $\SL(2,\CC)$ bundle $\sE$ over a compact, Hermitian surface, $X$, is \emph{stable} if for each holomorphic line bundle $\sL$ over $X$ for which there is a non-trivial holomorphic map $\sE \to \sL$ we have $\deg \sL > 0$; see also \cite[Section 6.1.4]{DK}. Recall that in \cite{DonASD, DK}, the group of $C^\infty$ bundle automorphisms of $E$ preserving the Hermitian metric on $E$ is denoted by $\sG_E$ and that $\fg_E$ is the real vector bundle of skew-adjoint endomorphisms of $E$. Throughout our monograph, $\NN$ denotes the set of non-negative integers.

\begin{mainthm}[Global existence and subsequential convergence modulo gauge transformations of Yang-Mills gradient flow for an initial connection of type $(1,1)$ over a K\"ahler surface]
\cite{DonASD}, \cite[Propositions 6.1.10, 6.2.7 and 6.2.14, Sections 6.2.5 and 6.2.6]{DK}
\label{mainthm:Donaldson_Kronheimer_6-2-7_and_6-2-14_plus}
Let $E$ be a complex rank two, Hermitian vector bundle with $c_1(E)=0$ over a compact, complex, Hermitian surface, $X$, and $A_0$ a $C^\infty$ unitary connection on $E$ with curvature, $F_{A_0}$, of type $(1,1)$. Then the following hold.
\begin{enumerate}
\item \emph{Global existence:} If the Hermitian metric on $X$ is \emph{K\"ahler}, then there exists a solution, $A(t) = A_0 + a(t)$ for $t\in [0,\infty)$, with
$$
a \in C^\infty([0,\infty)\times X; \Lambda^1\otimes\fg_E),
$$
to the Yang-Mills gradient flow \eqref{eq:Yang-Mills_gradient_flow} with initial data, $A(0) = A_0$, and $F_A(t)$ is of type $(1,1)$ for all $t\geq 0$;

\item \emph{Subsequential convergence:} If in addition $(E,\bar\partial_{A_0})$ is a \emph{stable} holomorphic vector bundle, then there are
\begin{enumerate}
\item a sequence of times, $\{t_m\}_{m\in\NN} \subset (0, \infty)$ with $t_m \to \infty$ as $m \to \infty$;

\item an irreducible $C^\infty$ anti-self-dual connection, $A_\infty$, on $E$; and

\item a sequence of gauge transformations, $\{\Phi_m\}_{m\in\NN} \subset \sG_E$,
\end{enumerate}
such that the sequence, $\{\Phi_m^*A(t_m)\}_{m\in\NN}$, converges strongly to $A_\infty$ over $X$ as $m \to \infty$ in the sense of $H_{A_0}^2(X;\Lambda^1\otimes\fg_E)$.

\item \emph{Uniqueness of the limit:} The limit, $A_\infty$, is unique up to the action of an element of $\sG_E$.
\end{enumerate}
\end{mainthm}

It is natural to ask whether the convergence statement in Theorem \ref{mainthm:Donaldson_Kronheimer_6-2-7_and_6-2-14_plus} can be strengthened to full convergence, $\Phi(t)^*A(t) \to \infty$ as $t\to \infty$, for a $C^\infty$ path of $C^\infty$ unitary gauge transformations, $\Phi(t) \in \sG_E$ for $t\in [0,\infty)$, or even if $A(t) \to A_\infty$ as $t\to \infty$, where convergence is again in the sense of $H_{A_0}^2(X;\Lambda^1\otimes\fg_E)$. One approach is to follow the example of R\r{a}de \cite[Section 7]{Rade_1992}, for Yang-Mills gradient flow over a closed Riemannian manifold, $X$, of dimension $d=2$ or $3$, and appeal to our results on the {\L}ojasiewicz-Simon gradient inequality for the Yang-Mills energy functional. Indeed, by exact analogy with R\r{a}de's proof of his \cite[Theorem 2]{Rade_1992} via \cite[Proposition 7.4]{Rade_1992}, we can combine Theorem \ref{mainthm:Donaldson_Kronheimer_6-2-7_and_6-2-14_plus} with our Theorem \ref{thm:Huang_3-3-6_Yang-Mills} to give

\begin{maincor}[Convergence of Yang-Mills gradient flow for an initial connection of type $(1,1)$ over a K\"ahler surface]
\label{cor:Donaldson_Kronheimer_6-2-7_and_6-2-14_convergence}
Assume the hypotheses of Theorem \ref{mainthm:Donaldson_Kronheimer_6-2-7_and_6-2-14_plus}, with $(E,\bar\partial_{A_0})$ a stable holomorphic vector bundle over a compact, K\"ahler surface, $X$. Then $A(t) \to A_\infty$ as $t\to \infty$ in the sense of $H_{A_0}^1(X;\Lambda^1\otimes\fg_E)$.
\end{maincor}

\begin{rmk}[Convergence modulo gauge transformations of Yang-Mills gradient flow for an initial connection of type $(1,1)$ over a K\"ahler surface]
\label{rmk:Donaldson_Kronheimer_6-2-7_and_6-2-14_convergence_mod_gauge}
It is possible to prove a weaker version of Corollary \ref{cor:Donaldson_Kronheimer_6-2-7_and_6-2-14_convergence} without appealing to the {\L}ojasiewicz-Simon gradient inequality, namely that $\Phi(t)^*A(t) \to A_\infty$ as $t \to \infty$ in the sense of $H_{A_0}^2(X;\Lambda^1\otimes\fg_E)$, for a $C^\infty$ path of $C^\infty$ unitary gauge transformations, $\Phi(t) \in \sG_E$ for $t\in [0,\infty)$. To see this, one can exploit the fact that the limit, $A_\infty$, is unique up to the action of $\sG_E$ and the concept is explained more fully in Section \ref{subsec:Proof_corollary_Donaldson_Kronheimer_6-2-7_and_6-2-14_convergence}; we are grateful to Richard Wentworth for this idea and helpful discussions.
\end{rmk}

\begin{mainthm}[Global existence and weak convergence of Yang-Mills gradient flow for an initial connection of type $(1,1)$ over a complex, Hermitian surface]
\label{mainthm:Yang-Mills_gradient_flow_global_existence_and_convergence_started_with arbitrary_initial_energy}
Let $E$ be a Hermitian vector bundle over a compact, complex, Hermitian surface, $X$, and $A_0$ a $C^\infty$ unitary connection on $E$ with curvature, $F_{A_0}$, of type $(1,1)$. Then the following hold.
\begin{enumerate}
\item \emph{Global existence:} There exists a solution, $A(t) = A_0 + a(t)$ for $t\in [0,\infty)$, with
$$
a \in C^\infty([0,\infty)\times X; \Lambda^1\otimes\fg_E),
$$
to the Yang-Mills gradient flow \eqref{eq:Yang-Mills_gradient_flow} with initial data, $A(0) = A_0$.

\item \emph{Dependence on initial data:} The solution, $A(t)$ for $t \in [0,\infty)$, varies continuously with respect to $A_0$ in the $C_{\loc}([0,\infty); H_{A_0}^1(X;\Lambda^1\otimes\fg_E))$ topology and, more generally, smoothly for all non-negative integers, $k,l$, in the $C_{\loc}^l([0,\infty); H_{A_0}^k(X;\Lambda^1\otimes\fg_E))$ topology.

\item \emph{Uniqueness:} Any two solutions are equivalent modulo a $C^\infty$ path of unitary gauge transformations, $\{u(t)\}_{t\geq 0} \subset \sG_E$, with $u(0) = \id_E$.

\item \emph{Weak convergence:} There are
\begin{enumerate}
\item a sequence of times, $\{t_m\}_{m\in\NN} \subset (0, \infty)$ with $t_m \to \infty$ as $m \to \infty$;

\item a finite (possibly empty) subset of points, $\Sigma := \{x_1\ldots,x_L\} \subset X$;

\item a $C^\infty$ Yang-Mills connection, $A_\infty$, on a Hermitian vector bundle, $E_\infty$, over $X$; and

\item a sequence of $H_{\loc}^3$ isomorphisms of Hermitian vector bundles, $\Phi_m: E_\infty \restriction \Sigma \to E \restriction \Sigma$,
\end{enumerate}
such that the sequence, $\{\Phi_m^*A(t_m)\}_{m\in\NN}$, converges to $A_\infty$ over $X \less \Sigma$ as $m \to \infty$ weakly in the sense of $H_{A_0,\loc}^2(X\less\Sigma;\Lambda^1\otimes\fg_E)$ and strongly in the sense of $W_{A_0,\loc}^{1,p}(X\less\Sigma;\Lambda^1\otimes\fg_E)$ for any $p \in [2,4)$.
\end{enumerate}
In addition, the set, $\Sigma$, is uniquely determined by the initial data.
\end{mainthm}

Theorem 1 in \cite{Daskalopoulos_Wentworth_2007}, where it is assumed in addition that $X$ is K\"ahler, suggests that $A_\infty$ in Theorem \ref{mainthm:Yang-Mills_gradient_flow_global_existence_and_convergence_started_with arbitrary_initial_energy} may be unique up to the action of $\sG_E$ and that certain integer multiplicities associated with the points $x_i \in \Sigma$ are uniquely determined by the flow, $A(t)$. As explained in \cite[Section 2]{Daskalopoulos_Wentworth_2007}, these integers have both an algebraic, sheaf-theoretic interpretation and an analytic interpretation in terms of weights associated with delta measures arising in limits of curvature densities, $|F_A(t_m)|^2$, as $m \to \infty$. The integer multiplicities have a third interpretation as characteristic numbers of the bundles arising in bubble-tree limits of $A(t)$ associated with each singular point, as explained by our Theorem \ref{thm:Kozono_Maeda_Naito_5-4} and a partial analogue of the Daskalopoulos-Wentworth characterization of those multiplicities is developed in Lemma \ref{lem:Daskalopoulos_Wentworth_2007_lemma_5_analogue} in the sequel.

The fact that $A(t)$ will generally have \emph{bubble-tree limits} --- requiring a finite rather than a single level of rescaling to converge (modulo gauge transformations and after passing to subsequences) without further energy loss --- appears to have been overlooked in \cite{Kozono_Maeda_Naito_1995, Schlatter_1997}.

We expect that Theorem \ref{mainthm:Yang-Mills_gradient_flow_global_existence_and_convergence_started_with arbitrary_initial_energy} should admit generalizations. In the following remarks, we describe some of the possible extensions that we will explore elsewhere, together with known limitations due to Daskalopoulos and Wentworth (see Section \ref{subsec:Daskalopoulos_Wentworth}) which prevent the convergence given in Theorem \ref{mainthm:Donaldson_Kronheimer_6-2-7_and_6-2-14_plus} or Corollary \ref{cor:Donaldson_Kronheimer_6-2-7_and_6-2-14_convergence} when $(E,\bar\partial_{A_0})$ is an unstable holomorphic vector bundle.

\begin{rmk}[More general Lie groups]
\label{rmk:Compact_Lie_groups}
In Theorem \ref{mainthm:Yang-Mills_gradient_flow_global_existence_and_convergence_started_with arbitrary_initial_energy}, we restricted our attention to $G = \U(n)$ for convenience, but the proof should extend, with minor modifications, to the case of any compact Lie structure group using the framework of Ramanathan and Subramanian \cite{Ramanathan_Subramanian_1988}.
\end{rmk}

\begin{rmk}[Uniqueness of $A_\infty$ and uniqueness of $\Sigma$ with multiplicity]
\label{rmk:Uniqueness_bubble_points_and_multiplicities}
With the additional hypothesis that the metric, $h$, on $X$ is K\"ahler, Daskalopoulos and Wentworth showed (see their \cite[Theorem 1]{Daskalopoulos_Wentworth_2007}) that the set of bubble points, $\Sigma = \{x_1,\ldots,x_L\} \subset X$, and certain associated integers (the \emph{analytic multiplicities}) are uniquely determined by the initial data, $A_0$. Moreover, they showed that the limit, $A_\infty$, is also unique (up to the action of $\sG_E$) and may be precisely identified in terms of $A_0$ and $E$ \cite[Theorem 1]{Daskalopoulos_Wentworth_2004}. It seems likely that their results would extend to the case where $h$ is merely Hermitian. Analogues of these uniqueness results for harmonic heat flow have been established by Irwin \cite{IrwinThesis}, Kwon \cite[Theorem 1.16]{KwonThesis}, and Topping \cite{ToppingThesis, Topping_2004am}, albeit under restrictive hypotheses. In the case of harmonic heat flow, versions of the {\L}ojasiewicz-Simon gradient inequality are employed by Irwin, Kwon, and Topping to establish their main results, so it is possible that our version of the {\L}ojasiewicz-Simon gradient inequality for the Yang-Mills energy functional (Theorem \ref{thm:Huang_3-3-6_Yang-Mills}) may play a similar role, just as it does in the simpler setting of our proof of Corollary \ref{cor:Donaldson_Kronheimer_6-2-7_and_6-2-14_convergence}. Moreover, the {\L}ojasiewicz-Simon gradient inequality may point to generalizations of the main theorems of Daskalopoulos and Wentworth in \cite{Daskalopoulos_Wentworth_2004, Daskalopoulos_Wentworth_2007} from the setting of K\"ahler surfaces to Riemannian four-manifolds. Our proof of Theorem \ref{mainthm:Yang-Mills_gradient_flow_global_existence_and_convergence_started_with arbitrary_initial_energy} suggests that the main theorems of \cite{Daskalopoulos_Wentworth_2004, Daskalopoulos_Wentworth_2007} may extend to the case of compact, complex, Hermitian surfaces, without appealing to the {\L}ojasiewicz-Simon gradient inequality.
\end{rmk}

\begin{rmk}[Convergence of $A(t)$ as $t\to\infty$ to an anti-self-dual connection $A_\infty$ on $E$]
\label{rmk:Non_Kaehler_convergence}
Given global existence, it is natural to ask whether our approach using Yang-Mills gradient-like flow can also give convergence, with $T=\infty$ in \eqref{eq:Yang-Mills_gradient-like_flow_Kaehler_surface}, when $(E,\bar\partial_{A_0})$ is further assumed to be a stable holomorphic vector bundle. When $X$ is K\"ahler and $A(t)$ is pure Yang-Mills gradient flow, the proof of convergence is presented in \cite[Sections 6.1.3, 6.1,4, and 6.2.3--6]{DK}. The proof of \cite[Proposition 6.2.14]{DK} relies on \cite[Corollary 6.2.12]{DK} for pure Yang-Mills gradient flow (our replacement, Corollary \ref{cor:Donaldson_Kronheimer_6-2-12}, for Yang-Mills gradient-like flow is not immediately applicable for that purpose). However, the role of \cite[Proposition 6.2.14, page 222]{DK} is to show that the limiting connection, $A_\infty$, if it exists is Yang-Mills on a bundle $E_\infty$ over $X$ and this we know for an arbitrary closed, Riemannian, smooth four-dimensional manifold by Theorem \ref{thm:Kozono_Maeda_Naito_5-3}. The result \cite[Corollary 6.2.12]{DK} is again used in the proof of \cite[Proposition 6.2.14, page 225]{DK} to show that either
\begin{inparaenum}[\itshape a\upshape)]
\item $\widehat F_{A_\infty} = 0$ and $A_\infty$ is anti-self-dual by
\cite[Proposition 2.1.59]{DK}, or
\item $\widehat F_{A_\infty}$ is non-zero and $A_\infty$ is a reducible connection on a holomorphic bundle, $\sL\oplus\sL^{-1}$ (when $n=2$), induced from a constant curvature connection on $\sL$ with $\deg\sL > 0$.
\end{inparaenum}
It may be possible to localize the proof of \cite[Proposition 6.2.14, page 225]{DK} for our application of Yang-Mills gradient-like flow.
Lastly, for our application of Yang-Mills gradient-like flow, it may be possible to localize the argument in \cite[Sections 6.2.5 and 6.2.6]{DK}, where it is shown that $A_\infty$ is anti-self-dual on $E_\infty$ and $E_\infty \cong E$, so no bubbling occurs as $t\to\infty$. Again, it is not obvious that all of the required arguments in \cite[Section 6.2]{DK} will extend, though there are reasons for optimism.

Jacob \cite{Jacob_2015conm} has shown that Donaldson's heat flow \cite{DonASD} exhibits both global existence and convergence when $X$ is a compact, complex surface with a Gauduchon metric (which always exists) and $(E,\bar\partial_{A_0})$ is stable, but Donaldson's heat flow does not necessarily coincide with Yang-Mills gradient flow when $X$ is non-K\"ahler. McNamara and Zhao \cite{McNamara_Zhao_2014arxiv} extend Jacob's results in \cite{Jacob_2015conm} to the case where $(E,\bar\partial_{A_0})$ is unstable, by analogy with the articles of Daskalopoulos and Wentworth \cite{Daskalopoulos_Wentworth_2004, Daskalopoulos_Wentworth_2007} when $X$ is K\"ahler.
\end{rmk}

\begin{rmk}[Global existence for solutions to pure Yang-Mills gradient flow in the case of initial data with small $F_{A_0}^{0,2}$]
\label{rmk:Small_FA02_global_existence}
Rather than assume that $A_0$ has curvature $F_{A_0}$ of type $(1,1)$ over $X$, it may suffice to assume that the component $F_{A_0}^{0,2}$ is suitably small over $X$, keeping in mind the identification $F_{A_0}^{0,2} = \bar\partial_{A_0}^2$ in \eqref{eq:FA02_is_barpartialA_squared} and the fact that $F_{A_0}^{2,0} = -(F_{A_0}^{0,2})^*$ when $A$ is unitary. We recall from \cite[Lemma 2.1.57]{DK} that
$$
F_A^+ = F_A^{2,0} + \langle F_A,\omega \rangle \omega + F_A^{0,2}
\quad\hbox{and}\quad
F_A^- = F_A^{1,1} - \langle F_A,\omega \rangle \omega,
$$
where $|\omega| = 1$ and $\omega$ is a $(1,1)$ form defined by the almost complex structure and Riemannian metric on $X$, and $\Omega^{p,q}(X)$ denotes the usual decomposition of $\Omega^r(X)$, for $p+q=r$, on a complex manifold, and $F_A^{p,q}$ the corresponding components of  the curvature, $F_A \in \Omega^2(X;\fg_E)$, of a unitary connection, $A$, on $E$.
\end{rmk}

\begin{rmk}[Global existence for solutions to pure Yang-Mills gradient flow in the case of an almost complex, four-dimensional, Riemannian manifold]
\label{rmk:Almost_complex_manifold_global_existence}
Because we use a perturbation argument it is possible, though certainly not obvious, that our proof of global existence for solutions to pure Yang-Mills gradient flow when $X$ is a compact, complex, possibly non-K\"ahler surface may extend to the case when $X$ is a symplectic or even an almost complex, four-dimensional, Riemannian smooth manifold, at least if $\bar\partial^2$ is suitably small over $X$, even if not identically zero. It is worth noting that by the decompositions \eqref{eq:Donaldson_Kronheimer_lemma_2-1-57} of $\Omega^2(X;\fg_E)$, we can define a connection $A$ on $E$ to have curvature of type $(1,1)$ in this setting by requiring that
$$
F_A^+ \in \Omega^0(X;\fg_E)\omega,
$$
where $\omega$ is a $(1,1)$ form defined by the almost complex structure and Riemannian metric on $X$.
\end{rmk}

\begin{rmk}[Previous results on existence of Yang-Mills connections on Hermitian vector bundles over complex surfaces]
When $E$ has complex rank $n$ and $c_1(E)=0$ (so the concepts of Yang-Mills and Hermitian Yang-Mills connections coincide
\cite[Section 6.1.4]{DK}) and $E$ is stable, then existence of an irreducible Yang-Mills connection on $E$ when $X$ is K\"ahler follows from more general results (proved without using Yang-Mills gradient flow) by Uhlenbeck and Yau \cite{Uhlenbeck_Yau_1986, Uhlenbeck_Yau_1989} and Ramanathan and Subramanian \cite{Ramanathan_Subramanian_1988}. When $X$ is more generally a compact, complex Hermitian, possibly non-K\"ahler surface, then existence of an irreducible Yang-Mills connection on $E$ follows from more general results of Buchdahl \cite{Buchdahl_1988} by adapting Donaldson's approach in \cite{DonNS} for vector bundles over Riemann surfaces (and again proved without using Yang-Mills gradient flow).
\end{rmk}

\subsection{Counterexamples to convergence of Yang-Mills gradient flow without bubbling}
\label{subsec:Daskalopoulos_Wentworth}
That bubble singularities can --- and in certain examples will --- occur at $T=\infty$ in Theorem \ref{mainthm:Yang-Mills_gradient_flow_global_existence_and_convergence_started_with arbitrary_initial_energy}, even when $X$ is K\"ahler, is illustrated by results of Daskalopoulos and Wentworth \cite{Daskalopoulos_Wentworth_2004, Daskalopoulos_Wentworth_2004correction, Daskalopoulos_Wentworth_2007} when $(E,\bar\partial_{A_0})$ is an \emph{unstable} holomorphic vector bundle over a compact, K\"ahler surface, $X$. We are indebted to Richard Wentworth for his explanations of their results.

Theorem 1 in \cite{Daskalopoulos_Wentworth_2007} extends Theorem \ref{mainthm:Donaldson_Kronheimer_6-2-7_and_6-2-14_plus} to the case where $(E,\bar\partial_{A_0})$ is unstable. The weak limit, $A_\infty$, is a (uniquely identified) Yang-Mills connection on a Hermitian vector bundle $E_\infty$, with $E_\infty \restriction \Sigma \cong E \restriction \Sigma$ and possibly non-empty singular set, $\Sigma$, just as in Theorem \ref{mainthm:Yang-Mills_gradient_flow_global_existence_and_convergence_started_with arbitrary_initial_energy}.

That $\Sigma$ can be non-empty is illustrated by the following example kindly provided to us by Wentworth. Consider nontrivial extensions of coherent sheaves on $\CC\PP^2$,
\begin{equation}
\label{eqn:extension}
0 \too \sO_{\CC\PP^2} \too \sE \too \sI_\Sigma \too 0,
\end{equation}
where $\sI_\Sigma$ is the ideal sheaf of a zero-dimensional subscheme.  It suffices to take $\Sigma$ to be a collection of distinct points.  By Friedman \cite[p. 38, \emph{Example}]{FriedmanBundleBook}, there are choices $\Sigma\neq \emptyset$ so that $\sE$ is locally free. Since $\sI_\Sigma^{\ast\ast}=\sO_{\CC\PP^2}$, we have $c_1(\sE)=0$ and $c_2(\sE)=\ell(\sO_{\CC\PP^2}/\sI_\Sigma)$, where $\ell$ denotes the length of the torsion sheaf. In the simple case where $\Sigma$ is a collection of distinct points each with multiplicity one, then $c_2(\sE)$ is just the cardinality of $\Sigma$.

Now it is easy to show that $\sE$ is strictly semistable (with respect the natural polarization on $\CC\PP^2$, say). In fact, the Seshadri filtration of $\sE$ is exactly $\sO_{\CC\PP^2}\subset \sE$, since $\sO_{\CC\PP^2}$ is the maximal destabilizing saturated subsheaf. The associated graded of the filtration  is  ${\rm Gr}(\sE)=\sO_{\CC\PP^2}\oplus \sI_\Sigma$, ${\rm Gr}(\sE)^{\ast\ast}=\sO_{\CC\PP^2}^{\oplus 2}$.

We choose the Fubini-Study K\"ahler metric $\omega $ on $\CC\PP^2$, a Hermitian metric $H$ on $\sE$, and let $A_0$ be the \emph{Chern connection}, that is, the unique unitary connection on $\sE$ that is compatible with the holomorphic structure on $\sE$ (see \cite[Proposition 4.2.14]{Huybrechts_2005}.  By the results of \cite{Daskalopoulos_Wentworth_2004}, any Uhlenbeck limit of Yang-Mills gradient flow, $A(t)$, with initial condition $A_0$ (and with respect to the metrics $\omega$ and $H$) must coincide with the Yang-Mills connection on ${\rm Gr}(\sE)^{\ast\ast}$; in this case, a trivial connection.  Consequently, there must be bubbling, since $c_2(\sE)\neq 0$.  In fact, \cite[Theorem 1]{Daskalopoulos_Wentworth_2007} says the bubbling occurs precisely at $\Sigma$ (with multiplicities).

Examples like \eqref{eqn:extension} were the original motivation for \cite[Theorem 1]{Daskalopoulos_Wentworth_2007}.  Notice that away from $\Sigma$ the metric gives a smooth splitting of \eqref{eqn:extension} with respect to which the $\bar\partial$-operator has the form
$$
\bar\partial_\sE=\left(\begin{matrix} \bar\partial_{\sO_{\CC\PP^2}} & \beta \\ 0 &  \bar\partial_{\sI_\Sigma}
\end{matrix}\right),
$$
and the connection looks like
$$
d_\sE=\left(\begin{matrix} d_{\sO_{\CC\PP^2}} & \beta \\ -\beta^\ast &  d_{\sI_\Sigma}
\end{matrix}\right).
$$
The philosophy is that along the flow, $A(t)$, the second fundamental form, $\beta(t)$, converges to zero. This is indeed the case but near $\Sigma$, the form $\beta(t)$ acquires singularities and estimates are difficult to obtain. The splitting above no longer makes sense, and one expects this to be reflected in the analysis.

\subsection{R\r{a}de's results on global existence and convergence of Yang-Mills gradient flow over base manifolds of dimensions two or three}
\label{subsec:Yang-Mills_gradient_flow_global_existence_and_convergence_rade}
When the base manifold, $X$, has dimension two or three, we obtain the following improvement of Theorem \ref{mainthm:Yang-Mills_gradient_flow_global_existence_and_convergence_started_near_minimum} due to R\r{a}de \cite{Rade_1992}. In place of R\r{a}de's rather complicated proof of global existence of the flow (see the proof of \cite[Theorems 1 and $1'$]{Rade_1992} in \cite[Sections 4, 5, and 6]{Rade_1992}), we shall instead observe that global existence can be deduced almost immediately
from a simple adaption of Struwe's analysis of the Yang-Mills gradient flow over four-dimensional base manifold \cite{Struwe_1994}.

\begin{mainthm}[Global existence and convergence of Yang-Mills gradient flow over base manifolds of dimensions two or three]
\label{mainthm:Rade_1_and_2}
(See \cite[Theorems 1 and 2]{Rade_1992}.)
Let $G$ be a compact Lie group and $P$ a principal $G$-bundle over a closed, connected, oriented, smooth manifold, $X$, of dimension $2$ or $3$ and with Riemannian metric, $g$. Let $A_1$ be a $C^\infty$ connection on $P$. Then there are constants $c \in [1,\infty)$, and $\sigma \in (0,1]$, and $\theta \in [1/2,1)$, depending on $A_1, G, g$, with the following significance.
\begin{enumerate}
\item
\label{item:Rade_theorem_1_and_2_global_existence}
\emph{Global existence:} There is a constant $\eps \in (0,\sigma/4)$, depending on $(A_1, G, g)$, with the following significance. If $A_0$ is a $C^\infty$ connection on $P$, then there exists a solution, $A(t) = A_0 + a(t)$ for $t\in [0,\infty)$, with
$$
a \in C^\infty([0,\infty)\times X; \Lambda^1\otimes\ad P),
$$
to the Yang-Mills gradient flow \eqref{eq:Yang-Mills_gradient_flow} with initial data, $A(0) = A_0$.

\item
\label{item:Rade_theorem_1_and_2_dependence_initial_data}
\emph{Dependence on initial data:} The solution, $A(t)$ for $t \in [0,\infty)$, varies continuously with respect to $A_0$ in the $C_{\loc}([0,\infty); H_{A_1}^1(X;\Lambda^1\otimes\ad P))$ topology and, more generally, smoothly for all non-negative integers, $l,m$, in the $C_{\loc}^l([0,\infty); H_{A_1}^m(X;\Lambda^1\otimes\ad P))$ topology.

\item
\label{item:Rade_theorem_1_and_2_convergence}
\emph{Convergence:} As $t\to\infty$, the flow, $A(t)$, converges strongly with respect to the norm on $H_{A_1}^1(X;\Lambda^1\otimes\ad P)$ to a Yang-Mills connection, $A_\infty$, of class $C^\infty$ on $P$, and the gradient-flow line has finite length in the sense that
$$
\int_0^\infty \left\|\frac{\partial A}{\partial t}\right\|_{H_{A_1}^1(X)}\,dt < \infty.
$$
If $A_\ym$ is a cluster point of the orbit, $O(A) = \{A(t):t\geq 0\}$, then $A_\infty = A_\ym$.

\item
\label{item:Rade_theorem_1_and_2_convergence_rate}
\emph{Convergence rate:} For all $t \geq 1$,
\begin{multline}
\label{eq:Rade_Proposition_7-4_convergence_rate_dimension_2_or_3}
\|A(t) - A_\infty\|_{H_{A_1}^1(X)}
\\
\leq
\begin{cases}
\displaystyle
\frac{1}{c(1-\theta)}\left(c^2(2\theta-1)(t-1) + (\sE(A_0)-\sE(A_\infty))^{1-2\theta}\right)^{-(1-\theta)/(2\theta-1)},
& 1/2 < \theta < 1,
\\
\displaystyle
\frac{2}{c}\sqrt{\sE(A_0)-\sE(A_\infty)}\exp(-c^2(t-1)/2),
&\theta = 1/2.
\end{cases}
\end{multline}

\item
\label{item:Rade_theorem_1_and_2_stability}
\emph{Stability:} If the critical point, $A_\infty$, is a local minimum then, as an equilibrium of the Yang-Mills gradient flow \eqref{eq:Yang-Mills_gradient_flow}, the point $A_\infty$ is Lyapunov stable; if $A_\infty$ is isolated or a cluster point of the orbit $O(A)$, then $A_\infty$ is uniformly asymptotically stable.

\item
\label{item:Rade_theorem_1_and_2_uniqueness}
\emph{Uniqueness:} Any two solutions are equivalent modulo a path of gauge transformations,
$$
u \in C^\infty([0,\infty)\times X; \Ad P), \quad u(0) = \id_P.
$$
\end{enumerate}
\end{mainthm}

\begin{rmk}[Comparison of Theorem \ref{mainthm:Rade_1_and_2} with results of R\r{a}de]
\label{rmk:Rade_theorem_1_and_2_comparison}
Items \eqref{item:Rade_theorem_1_and_2_stability} is not explicitly proved by R\r{a}de in \cite{Rade_1992}, but this assertion can be derived with the aid of the {\L}ojasiewicz-Simon gradient inequality, just as in the proof of \cite[Theorem 5.1.2]{Huang_2006}. R\r{a}de establishes a more general version of Theorem \ref{mainthm:Rade_1_and_2}, in that he allows initial data, $A_0$, of class $H^1$ and a weaker concept of solution to Yang-Mills gradient flow \eqref{eq:Yang-Mills_gradient_flow} \cite[Definition, p. 127]{Rade_1992}. Finally, if $d = 2$ and $G = \U(n)$ and $A_\infty$ is irreducible, then R\r{a}de obtains $\theta = 1/2$, yielding exponential convergence in \eqref{eq:Rade_Proposition_7-4_convergence_rate_dimension_2_or_3}.
\end{rmk}

\subsection{Harmonic map gradient flow near critical points}
\label{subsec:Harmonic_map_gradient_flow_near_critical points}
The results of this section are essentially due to Simon \cite{Simon_1983, Simon_1985} (compare Theorem \ref{thm:Simon_in_Kwon_thesis_1-15}), but the specific statements we give here and their justifications are difficult to find in the literature. However, as in the more difficult case of Yang-Mills gradient flow, the results of this section may be obtained as direct consequences of the general methods that we develop in this monograph and thus have self-contained proofs. For an introduction to harmonic maps and harmonic map gradient flow, we refer to Eells and Sampson \cite{Eells_Sampson_1964ajm}, Eells and Lemaire \cite{Eells_Lemaire_1978}, Hamilton \cite{Hamilton_1975}, H{\'e}lein \cite{Helein_harmonic_maps}, Jost \cite{Jost_two_dim_geom_var_probs, Jost_riemannian_geometry_geometric_analysis}, Lin and Wang \cite{Lin_Wang_2008}, Simon \cite{Simon_1996}, Struwe \cite{Struwe_1985, Struwe_1996, Struwe_variational_methods}, and references cited therein. For a selection of more recent articles on harmonic map gradient flow, see Li and Zhu \cite{Li_Zhu_2012}, Luo \cite{Luo_2012}, and Topping \cite{Topping_2004am, Topping_2004mz}. Biernat \cite{Biernat_2015}, Bizo{\'n} and Wasserman \cite{Bizon_Wasserman_2015}, Boling, Kelleher, and Streets \cite{Boling_Kelleher_Streets_2015arxiv}, Fan \cite{Fan_1999} and Zhang \cite{Zhang_2012} have results on singularity formulation in harmonic map gradient flow for source manifolds of dimension two and higher.

We begin by recalling a consequence of Theorem \ref{thm:Huang_2-4-5_introduction} for the harmonic map $L^2$-energy functional.

\begin{defn}[Harmonic map energy functional]
\label{defn:Harmonic_map_energy_functional}
Let $(M,g)$ and $(N,h)$ be a pair of closed, Riemannian, smooth manifolds. One defines the \emph{harmonic map $L^2$-energy functional} by
\begin{equation}
\label{eq:Harmonic_map_energy_functional}
\sE_{g,h}(f)
:=
\frac{1}{2} \int_M |df|_{g,h}^2 \,d\vol_g,
\end{equation}
for smooth maps, $f:M\to N$, where $df:TM \to TN$ is the differential map.
\end{defn}

When clear from the context, we omit explicit mention of the Riemannian metrics $g$ on $M$ and $h$ on $N$ and write $\sE = \sE_{g,h}$. Although initially defined for smooth maps, the energy functional $\sE$ in Definition \ref{defn:Harmonic_map_energy_functional}, extends to the case of Sobolev maps of class $W^{1,2}$. To define the gradient of the energy functional $\sE$ in \eqref{eq:Harmonic_map_energy_functional} with respect to the $L^2$ metric on $C^\infty(M;N)$, we first choose an isometric embedding, $(N,h) \hookrightarrow \RR^n$ for a sufficiently large $n$ (courtesy of the isometric embedding theorem due to Nash \cite{Nash_1956}), and recall that by \cite[Equations (8.1.10) and (8.1.13)]{Jost_riemannian_geometry_geometric_analysis} we have
\begin{align*}
\left(\sE'(f),u\right)_{L^2(M,g)}
&:=
\left.\frac{d}{dt}\sE(\exp_{f}(tu))\right|_{t=0}
\\
&\,= \left(\Delta_g f,u\right)_{L^2(M,g)}
\\
&\,= \left(\Pi_h(f)\Delta_g f,u\right)_{L^2(M,g)},
\end{align*}
for all $u \in C^\infty(M;f^*TN)$, where $\Pi_h(y):\RR^n \to T_yN$ is orthogonal projection and $\exp_y:T_yN \to N$ is the exponential map, so $\exp_y(0) = y \in N$, for all $y \in N$. (Note that one could alternatively define
\[
\left(\sE'(f),u\right)_{L^2(M,g)}
=
\left.\frac{d}{dt}\sE(\pi(f + tu))\right|_{t=0}
\]
as implied by \cite[Equations (2.2)(i) and (ii)]{Simon_1996}, where $\pi$ is the nearest
point projection onto $N$ from a normal tubular neighborhood.) Thus, viewing the gradient as an operator and applying \cite[Lemma 1.2.4]{Helein_harmonic_maps},
\begin{equation}
\label{eq:Gradient_harmonic_map_operator}
\sE'(f) = \Pi_h(f)\Delta_g = \Delta_g f - A_h(df,df),
\end{equation}
as in \cite[Equations (2.2)(iii) and (iv)]{Simon_1996}. Here, $A_h$ denotes the second fundamental form of the isometric embedding, $(N,h) \subset \RR^n$ and
\begin{equation}
\label{eq:Laplace-Beltrami_operator}
\Delta_g
:=
-\divg_g \grad_g
=
d^{*,g}d
=
-\frac{1}{\sqrt{\det g}} \frac{\partial}{\partial x^\beta}
\left(\sqrt{\det g}\, \frac{\partial f}{\partial x^\alpha} \right)
\end{equation}
denotes the Laplace-Beltrami operator for $(M,g)$ (with the opposite sign convention to that of \cite[Equations (1.14) and (1.33)]{Chavel}) acting on the scalar components $f^i$ of $f = (f^1,\ldots,f^n)$ and $\{x^\alpha\}$ denote local coordinates on $M$. As usual, the gradient vector field, $\grad_g f^i \in C^\infty(TM)$, is defined by $\langle\grad_g f^i,\xi\rangle_g := df^i(\xi)$ for all $\xi \in C^\infty(TM)$ and $1
\leq i \leq n$ and the divergence function, $\divg_g\xi \in C^\infty(M;\RR)$, by the pointwise trace,
$\divg_g\xi := \tr(\eta \mapsto \nabla_\xi^g\eta)$, for all $\eta \in C^\infty(TM)$.

One says that a smooth map $f:M\to N$ is \emph{harmonic} if it is a critical point of the $L^2$ energy functional \eqref{eq:Harmonic_map_energy_functional}, that is
\[
\sE'(f) = \Delta_g f - A_h(df,df) = 0.
\]
Given a smooth map $f:M\to N$, an isometric embedding $(N,h) \hookrightarrow \RR^n$, a non-negative integer $k$, and $p \in [1,\infty)$, we define the Sobolev norms,
\[
\|f\|_{W^{k,p}(M)} := \left(\sum_{i=1}^n \|f^i\|_{W^{k,p}(M)}^p\right)^{1/p},
\]
with
\[
\|f^i\|_{W^{k,p}(M)} := \left(\sum_{j=0}^k \int_M |(\nabla^g)^j f^i|^p \,d\vol_g\right)^{1/p},
\]
where $\nabla^g$ denotes the Levi-Civita connection on $TM$ and all associated bundles (that is, $T^*M$ and their tensor products), and if $p = \infty$, we define
\[
\|f\|_{W^{k,\infty}(M)} = \|f\|_{C^k(M)}
:=
\sum_{i=1}^{n}\sum_{j=0}^k \esssup_M |(\nabla^g)^j f^i|.
\]
If $k=0$, then we denote $\|f\|_{W^{0,p}(M)} = \|f\|_{L^p(M)}$. For $p \in [1,\infty)$ and nonnegative integers $k$, we use \cite[Theorem 3.12]{AdamsFournier} (applied to $W^{k,p}(M;\RR^n)$ and noting that $M$ is a closed manifold) and Banach space duality to define
\[
W^{-k,p'}(M;\RR^n) := \left(W^{k,p}(M;\RR^n)\right)^*,
\]
where $p'\in (1,\infty]$ is the dual exponent defined by $1/p+1/p'=1$. Elements of the Banach space dual $(W^{k,p}(M;\RR^n))^*$ may be characterized via \cite[Section 3.10]{AdamsFournier} as distributions in the Schwartz space $\sD'(M;\RR^n)$ \cite[Section 1.57]{AdamsFournier}.

In particular, when $p = 1$ and $p' = \infty$ and $k$ is a non-negative integer, we have
\[
W^{-k,\infty}(M;\RR^n) := \left(W^{k,1}(M;\RR^n)\right)^*.
\]
Lastly, we note that if $(N,h)$ is real analytic, then the isometric embedding $(N,h) \hookrightarrow \RR^n$ may also be chosen to be analytic by the analytic isometric embedding theorem due to Nash \cite{Nash_1966}, with a simplified proof due to Greene and Jacobowitz \cite{Greene_Jacobowitz_1971}).

The statement of the forthcoming Theorem \ref{thm:Lojasiewicz-Simon_gradient_inequality_energy_functional_Riemannian_manifolds} includes the most delicate dimension for the source Riemannian manifold, $(M,g)$, namely the case where $M$ has dimension $d=2$ and allows a Sobolev norm for the definition of the {\L}ojasiewicz-Simon neighborhood of a harmonic map that appears to be optimal for that case, namely, $W^{2,1}(M;N)$, as well as the suboptimal $W^{1,p}(M;N)$ with $p > 2$. Following the landmark articles by Sacks and Uhlenbeck \cite{Sacks_Uhlenbeck_1981, Sacks_Uhlenbeck_1982}, the case where the domain manifold $M$ has dimension two is well-known to be critical.

\begin{thm}[{\L}ojasiewicz-Simon gradient inequality for the energy functional for maps between pairs of Riemannian manifolds]
\label{thm:Lojasiewicz-Simon_gradient_inequality_energy_functional_Riemannian_manifolds}
\cite[Theorem 4]{Feehan_Maridakis_Lojasiewicz-Simon_harmonic_maps}
Let $d\geq 2$ and $k \geq 1$ be integers and $p\in [1,\infty)$ be such that
\[
kp > d \quad\text{or}\quad k=d \text{ and } p=1.
\]
Let $(M,g)$ and $(N,h)$ be closed, Riemannian, smooth manifolds, with $M$ of dimension $d$. If $(N,h)$ is real analytic (respectively, $C^\infty$) and $f\in W^{k,p}(M;N)$, then the gradient map\footnote{Thus $T^*_fW^{k,p}(M;N)$ is the dual of the tangent space $T_fW^{k,p}(M;N)$ of the Banach manifold $W^{k,p}(M;N)$ at the point $f$.}
\[
\sE'(f): T_fW^{k,p}(M; N) \to T_f^*W^{k,p}(M; N),
\]
is a real analytic (respectively, $C^\infty$) map of Banach spaces. If $f_\infty \in W^{k,p}(M; N)$ is a harmonic map, then there are positive constants $c \in [1, \infty)$, $\sigma \in (0,1]$, and $\theta \in [1/2,1)$, depending on $f_\infty$, $g$, $h$, $k$, $p$, $M$, and $N$ with the following significance. If $f\in W^{k,p}(M;N)$ obeys the $W^{k,p}$ \emph{{\L}ojasiewicz-Simon neighborhood} condition,
\begin{equation}
\label{eq:Lojasiewicz-Simon_gradient_inequality_harmonic_map_neighborhood_Riemannian_manifold}
\|f - f_\infty\|_{W^{k,p}(M)} < \sigma,
\end{equation}
then the harmonic map energy functional
\eqref{eq:Harmonic_map_energy_functional} obeys the
\emph{{\L}ojasiewicz-Simon gradient inequality},
\begin{equation}
\label{eq:Lojasiewicz-Simon_gradient_inequality_harmonic_map_energy_functional_Riemannian_manifold}
\|\sE'(f)\|_{W^{-k,p'}(M)}
\geq
c|\sE(f) - \sE(f_\infty)|^\theta.
\end{equation}
\end{thm}

\begin{rmk}[Previous versions of the {\L}ojasiewicz-Simon gradient inequality for the harmonic map energy functional]
Topping \cite[Lemma 1]{Topping_1997} proved a {\L}ojasiewicz-type gradient inequality for maps, $f:S^2 \to S^2$, with small $L^2$ energy, with the latter criterion replacing the usual small $C^{2,\alpha}(M;\RR^n)$ norm criterion of Simon for the difference between a map and a critical point \cite[Theorem 3]{Simon_1983}. Simon uses a $C^2(M;\RR^n)$ norm to measure distance between maps, $f:M \to N$,  in \cite[Equation (4.27)]{Simon_1985}. Topping's result is generalized by Liu and Yang in \cite[Lemma 3.3]{Liu_Yang_2010}. Kwon \cite[Theorem 4.2]{KwonThesis} obtains a {\L}ojasiewicz-type gradient inequality for maps, $f:S^2 \to N$, that are $W^{2,p}(S^2;\RR^n)$-close to a harmonic map,  with $1 < p \leq 2$. However, her proof explicitly uses the fact that $p > 1$.
\end{rmk}

We now turn to questions of convergence, global existence, and stability of solutions to harmonic map gradient flow. We begin by recalling the following result due to Eells and Sampson on the existence of minima for the harmonic map energy functional.

\begin{thm}[Existence of harmonic maps whose energy is absolutely minimizing in a homotopy class]
\label{thm:Eells_Sampson_page_158}
\cite[Corollary, p. 158]{Eells_Sampson_1964ajm}
Let $(M,g)$ and $(N,h)$ be a pair of closed, Riemannian, smooth manifolds and assume that $(N,h)$ has non-positive sectional curvature. If $f_0 \in C(M;N)$ is a continuous map, then there exists a $C^\infty$ harmonic map $f_{\min}$ that is homotopic to $f_0$ and which is absolutely energy minimizing in the homotopy class $[f_0]$.
\end{thm}

The curvature hypothesis in Theorem \ref{thm:Eells_Sampson_page_158} was relaxed by Sacks and Uhlenbeck in their celebrated

\begin{thm}[Existence of harmonic maps whose energy is absolutely minimizing in a homotopy class]
\label{thm:Sacks_Uhlenbeck_5-1}
\cite[Theorem 5.1]{Sacks_Uhlenbeck_1981}
Let $(M,g)$ and $(N,h)$ be a pair of closed, Riemannian, smooth manifolds and assume that $\pi_2(N)=0$. If $f_0 \in C(M;N)$ is a continuous map, then there exists a $C^\infty$ harmonic map $f_{\min}$ that is homotopic to $f_0$ and which is absolutely energy minimizing in the homotopy class $[f_0]$.
\end{thm}

The proof of Theorem \ref{thm:Eells_Sampson_page_158} relies on fundamental results for harmonic map gradient flow due to Eells and Sampson \cite[Theorem, p. 156]{Eells_Sampson_1964ajm}; the version of their results that we state below appears as \cite[Theorem 5.3.1]{Lin_Wang_2008}.

\begin{thm}[Global existence and convergence of harmonic map gradient flow into a closed, Riemannian, smooth manifold with non-positive sectional curvature]
\label{thm:Eells_Sampson}
\cite[Theorem, p. 156]{Eells_Sampson_1964ajm},
\cite[Theorem 5.3.1]{Lin_Wang_2008}
Let $(M,g)$ and $(N,h)$ be a pair of closed, Riemannian, smooth manifolds and assume that $(N,h)$ has non-positive sectional curvature. If $f_0 \in C^\infty(M;N)$, then there exists a unique solution, $f \in C^\infty([0,\infty)\times M; N)$, to the harmonic map gradient flow equation,
\begin{equation}
\label{eq:Harmonic_map_gradient_flow}
\dot f(t) = - \sE'(f), \quad f(0) = f_0,
\end{equation}
a harmonic map, $f_\infty \in C^\infty(M;N)$, and an unbounded sequence $\{t_m\}_{m=0}^\infty \subset [0,\infty)$ such that $\sE(f_\infty) \leq \sE(f_0)$ and
\[
\|f(t_m) - f_\infty\|_{C^2(M;\RR^n)} \to 0, \quad\text{as } m \to \infty.
\]
\end{thm}

If in addition to the hypotheses of Theorem \ref{thm:Eells_Sampson}, we assume that $(N,h)$ is real analytic, it then follows from Theorem \ref{thm:Eells_Sampson} and \cite[Corollary 2]{Simon_1983} that
\[
\|f(t) - f_\infty\|_{C^2(M;\RR^n)} \to 0, \quad\text{as } t \to \infty.
\]
For a recent exposition of the proof of Theorem \ref{thm:Eells_Sampson} and certain generalizations, we refer the reader to Lin and Wang \cite[Section 5.3]{Lin_Wang_2008}. Short-time existence, uniqueness, and regularity of solutions to the harmonic map gradient flow equation \eqref{eq:Harmonic_map_gradient_flow} can be established either by appealing to general results for nonlinear parabolic systems developed in this monograph or by appealing to previous expositions or original results due to Eells and Sampson \cite{Eells_Sampson_1964ajm}, Lin and Wang \cite{Lin_Wang_2008}, Simon \cite{Simon_1996}, or Struwe \cite{Struwe_1985}.
Lastly, in order to apply our results for abstract gradient flows, we need to verify Hypothesis \ref{hyp:Abstract_apriori_interior_estimate_trajectory_main_introduction} with the following analogue of Corollary \ref{cor:Rade_7-3_arbitrary_dimension_L1_time_W1p_space} for the Yang-Mills $L^2$-energy functional.

\begin{lem}[\Apriori $L^1$-in-time-$W^{1,p}$-in-space interior estimate for a solution to harmonic map gradient flow]
\label{lem:Harmonic_map_L1_time_W1p_space_apriori_interior_estimate_trajectory}
Let $(M,g)$ and $(N,h)$ be closed, Riemannian, smooth manifolds, with $M$ of dimension $d \geq 2$, and $f_\infty \in C^\infty(M; N)$, and $p \in (2\vee d/2,\infty)$. Then there are positive constants, $C = C(d, f_\infty, g, h, p) \in [1,\infty)$ and $\eps_1 = \eps_1(d, f_\infty, g, h, p) \in (0, 1]$, such that if $f \in C^\infty((S, T)\times M; N)$ is a solution to harmonic map gradient flow \eqref{eq:Harmonic_map_gradient_flow} on an interval $(S, T)$, where $S \in \RR$ and $\delta > 0$ and $T$ obey $S + 2\delta \leq T \leq \infty$, and
\begin{equation}
\label{eq:Harmonic_map_Linfinity_in_time_W1p_in_space_small_norm_ft_minus_f1_condition}
\|f(t) - f_\infty\|_{W^{1,p}(M)} \leq \eps_1, \quad\forall\, t \in (S, T),
\end{equation}
then there is an integer $n = n(d,p) \geq 1$ such that
\begin{equation}
\label{eq:Harmonic_map_apriori_interior_estimate}
\int_{S+\delta}^T \|\dot f(t)\|_{W^{1,p}(M)}\,dt
\leq
C\left(1 + \delta^{-n}\right)\int_S^T \|\dot f(t)\|_{L^2(M)}\,dt.
\end{equation}
\end{lem}

In applications of Lemma \ref{lem:Harmonic_map_L1_time_W1p_space_apriori_interior_estimate_trajectory}, the map $f_\infty \in C^\infty(M; N)$ will be harmonic, but that is not a hypothesis of Lemma \ref{lem:Harmonic_map_L1_time_W1p_space_apriori_interior_estimate_trajectory}. We omit the proof of Lemma \ref{lem:Harmonic_map_L1_time_W1p_space_apriori_interior_estimate_trajectory} as the ideas are very similar to those employed in the proof of Corollary \ref{cor:Rade_7-3_arbitrary_dimension_L1_time_W1p_space} and it again relies on the abstract Lemma \ref{lem:Rade_7-3_abstract_interior_L1_in_time_V2beta_space_time_derivative_interior}.
Consequently, by virtue of Theorem \ref{thm:Lojasiewicz-Simon_gradient_inequality_energy_functional_Riemannian_manifolds} together with Theorems \ref{mainthm:Huang_3-4-8_introduction}, \ref{mainthm:Huang_5-1-1_introduction}, and \ref{mainthm:Huang_5-1-2_introduction}, we obtain the following analogue of Theorem \ref{mainthm:Yang-Mills_gradient_flow_global_existence_and_convergence_started_near_minimum}.

\begin{mainthm}[Global existence and convergence of harmonic map flow near a local minimum]
\label{mainthm:Harmonic_map_gradient_flow_global_existence_and_convergence_started_near_minimum}
Let $(M,g)$ and $(N,h)$ be closed, Riemannian, smooth manifolds, with $M$ of dimension $d\geq 2$ and $(N,h)$ real analytic. If $f_{\min} \in C^\infty(M;N)$ is a local minimum for the harmonic map energy functional $\sE$ and $c \in [1,\infty)$, and $\sigma \in (0,1]$, and $\theta \in [1/2,1)$ are the {\L}ojasiewicz-Simon constants for $(\sE,f_{\min})$ given by Theorem \ref{thm:Lojasiewicz-Simon_gradient_inequality_energy_functional_Riemannian_manifolds}, and $p \in (2\vee d/2,\infty)$, then there is a constant $\eps \in (0,\sigma/4)$ such that the following hold:
\begin{enumerate}
\item \emph{Global existence and uniqueness:} If $f_0 \in C^\infty(M;N)$ obeys
$$
\|f_0 - f_{\min}\|_{W^{1,p}(M)} < \eps,
$$
then there exists a unique solution, $f \in C^\infty([0,\infty)\times M; N)$, to the harmonic map gradient flow equation \eqref{eq:Harmonic_map_gradient_flow} with initial data, $f(0) = f_0$, and
$$
\|f(t) - f_{\min}\|_{W^{1,p}(M)} < \sigma/2, \quad\forall\, t \in [0,\infty).
$$

\item \emph{Dependence on initial data:} The solution, $f(t)$ for $t \in [0,\infty)$, varies continuously with respect to $f_0$ in the $C_{\loc}([0,\infty); W^{1,p}(M;N))$ topology and, more generally, smoothly for all integers $k\geq 1$ and $l\geq 0$ in the $C_{\loc}^l([0,\infty); W^{k,p}(M;N))$ topology.

\item \emph{Convergence:} As $t\to\infty$, the flow, $f(t)$, converges strongly with respect to the norm on $W^{1,p}(M;N)$ to a harmonic map, $f_\infty \in C^\infty(M;N)$, and the gradient-flow line has finite length in the sense that
$$
\int_0^\infty \|\dot f(t)\|_{W^{1,p}(M)}\,dt < \infty.
$$
If $f_{\min}$ is a cluster point of the orbit, $O(f) = \{f(t):t\geq 0\}$, then $f_\infty = f_{\min}$.

\item \emph{Convergence rate:} For all $t \geq 1$,
\begin{multline*}
\|f(t) - f_\infty\|_{W^{1,p}(M)}
\\
\leq
\begin{cases}
\displaystyle
\frac{1}{c(1-\theta)}\left(c^2(2\theta-1)(t-1) + (\sE(f_0)-\sE(f_\infty))^{1-2\theta}\right)^{-(1-\theta)/(2\theta-1)},
& 1/2 < \theta < 1,
\\
\displaystyle
\frac{2}{c}\sqrt{\sE(f_0)-\sE(f_\infty)}\exp(-c^2(t-1)/2),
&\theta = 1/2.
\end{cases}
\end{multline*}

\item \emph{Stability:} As an equilibrium of the harmonic map gradient flow \eqref{eq:Harmonic_map_gradient_flow}, the point $f_\infty$ is Lyapunov stable; if $f_\infty$ is isolated or a cluster point of the orbit $O(f)$, then $f_\infty$ is uniformly asymptotically stable.
\end{enumerate}
\end{mainthm}

When $M = N = S^2$ (with its standard round metric of radius one) and $\sE_\infty := \lim_{t\to\infty} \sE(u(t))$, Topping \cite[Theorem 1.7]{Topping_2004am} has shown that
$$
|\sE(u(t)) - \sE_\infty| \leq C_0\exp(-t/C_0), \quad\forall\, t\geq 0,
$$
for a positive constant, $C_0$; he allows that $u(t)$ may bubble at finitely many points in $S^2$ as $t\to\infty$.

It is worth noting that Theorem \ref{mainthm:Harmonic_map_gradient_flow_global_existence_and_convergence_started_near_minimum} does not contradict the example, due to Chang, Ding, and Ye \cite{Chang_Ding_Ye_1992}, of finite-time blow-up for harmonic map gradient flow for maps from $S^2$ to $S^2$ and initial data $f_0$ of sufficiently high energy.

Further information regarding the behavior of harmonic map flow near an arbitrary critical point can sometimes be deduced from the following consequence of Theorem \ref{mainthm:Huang_3-3-6_introduction} and analogue of Theorem \ref{mainthm:Yang-Mills_gradient_flow_global_with_critical_point_in_orbit_closure} for Yang-Mills gradient flow.

\begin{mainthm}[Convergence of a subsequence implies convergence for a solution to harmonic map gradient flow]
\label{mainthm:Harmonic_map_gradient_flow_global_with_critical_point_in_orbit_closure}
Let $(M,g)$ and $(N,h)$ be closed, Riemannian, smooth manifolds, with $M$ of dimension $d\geq 2$ and $(N,h)$ real analytic, $f_\infty \in C^\infty(M;N)$ be a harmonic map, and $p \in (2\vee d/2,\infty)$.
If $f \in C^\infty([0,\infty)\times M; N)$ is a solution to the harmonic map gradient flow equation \eqref{eq:Harmonic_map_gradient_flow} and $f_\infty$ is a cluster point of the orbit $O(f) = \{f(t): t\geq 0\}$, in the sense that there exists a sequence of times, $\{t_m\}_{m\in\NN} \subset [0,\infty)$ with $t_m\to\infty$ as $m\to\infty$, such that
\[
f(t_m) \to f_\infty \quad\text{in } W^{1,p}(M;\RR^n) \text{ as } m \to \infty,
\]
then $f(t)$ converges to $f_\infty$ as $t\to\infty$ in the sense that
\[
\lim_{t\to\infty}\|f(t)-f_\infty\|_{W^{1,p}(M)} = 0
\quad\hbox{and}\quad
\int_0^\infty \|\dot f\|_{W^{1,p}(M)}\,dt < \infty.
\]
\end{mainthm}

\section{Summary}
For the benefit of the reader, we outline the remainder of our monograph.

In order to highlight methods which may have application to Yang-Mills gradient flow (and vice versa), Chapter \ref{chapter:Comparison_global_existence_convergence_results_gradient_and_heat_flows_geometric_analysis} surveys other nonlinear evolution equations in geometric analysis and mathematical physics, including
\begin{inparaenum}[\itshape a\upshape)]
\item harmonic map gradient flow from Riemann surfaces into a target Riemannian manifold (see Section \ref{sec:Gradient_flow_harmonic_map_energy_functional}), and
\item conformal Yamabe scalar curvature heat flow (see Section \ref{sec:Yamabe_scalar_curvature_heat_flow}).
\end{inparaenum}
We also briefly mention (see Section \ref{sec:Other_gradient_flows}) Chern-Simons gradient flow, Donaldson's heat flow, Ginzburg-Landau energy gradient flow, knot energy gradient flow, Lagrangian mean curvature flow, mean curvature flow, Ricci curvature flow, Taubes' anti-self-dual curvature flow, and other gradient or gradient-like flows in applied mathematics and fluid dynamics, pointing the reader to more detailed references in each case.

In Chapter \ref{chapter:Preliminaries}, we gather a number of preliminary definitions and facts that we shall need throughout the course of our monograph. In Section \ref{sec:Yang-Mills_gradient_flow_solution_concepts}, we define different concepts of \emph{classical solution} different concepts of \emph{weak solution} to Yang-Mills gradient flow. Section \ref{sec:Heat_equation_method} surveys the origins and philosophy of the heat equation method from the perspective of solving an elliptic equation. Section \ref{sec:Taubes_1982_Appendix} reviews the classification of principal $G$-bundles, the Chern-Weil formula, and characterization of minima of the Yang-Mills energy functional in dimension four. In Section \ref{sec:Critical_points_Yang-Mills_energy_functional}, we briefly review what is known about minimal and non-minimal critical points of the Yang-Mills energy functional over a base manifold of dimension four, together with consequences for Yang-Mills gradient flow.

Chapter \ref{chapter:Sell_You_4} provides a development of linear and nonlinear evolutionary equations in Banach spaces from the perspective of analytic semigroups defined by positive, sectorial operators on Banach spaces. Much of this chapter closely follows the presentation due to Sell and You \cite{Sell_You_2002}, but we extend their local well-posedness results for nonlinear evolutionary equations by
\begin{inparaenum}[\itshape a\upshape)]
\item considering polynomial nonlinearities;
\item allowing initial data of lower regularity than that assumed \cite{Sell_You_2002} using Banach spaces with temporal weights; and
\item deriving explicit lower bounds for the minimal lifetimes of solutions to nonlinear evolutionary equations in Banach spaces.
\end{inparaenum}
Section \ref{sec:Sell_You_4-2} summarizes the linear theory required from \cite{Sell_You_2002}.  Section \ref{sec:Local_existence_nonlinear_evolution_equation_Banach_space} summarizes the nonlinear theory required from \cite{Sell_You_2002} and develops the aforementioned extensions.

In Chapter \ref{chapter:Elliptic_and_parabolic_partial_differential_systems}, we consider elliptic differential systems on domains in $\RR^d$ and on vector bundles over closed manifolds, develop \apriori estimates, existence and uniqueness results, resolvent (or Agmon) estimates, and consequently show that these elliptic differential systems define analytic semigroups on standard $L^p$, $C^0$, and $L^1$ Banach spaces in Section \ref{sec:Sell_You_3-8-2_standard_Sobolev_spaces} and certain critical-exponent Banach spaces in Section \ref{sec:Sell_You_3-8-2_critical_exponent_elliptic_Sobolev_spaces}. While there are many references for elliptic and parabolic scalar differential operators (such as Gilbarg and Trudinger \cite{GilbargTrudinger} and Krylov \cite{Krylov_LecturesSobolev}) the available references which discuss elliptic and parabolic differential systems (such as Agmon \cite{AgmonLecturesEllipticBVP}, Chen and Wu \cite{Chen_Wu_1998}, Lady{\v{z}}enskaja, Solonnikov and Ural$'$ceva \cite{Ladyzhenskaya_Uraltseva_1968, LadyzenskajaSolonnikovUralceva}, and Morrey \cite{Morrey}) do not develop the results we need for our application, so Chapter \ref{chapter:Elliptic_and_parabolic_partial_differential_systems} provides an essentially self-contained treatment.

The Yang-Mills heat equation --- a quasilinear parabolic equation defined by a choice of reference connection, as distinct from the non-parabolic Yang-Mills gradient flow equation --- takes the form,
\begin{equation}
\label{eq:Yang-Mills_heat_equation_introduction}
\frac{\partial a}{\partial t} + d_A^*F_A + d_Ad_A^*a = 0,
\end{equation}
where $A(t) = A_1 + a(t)$, and $A_1$ is a fixed $C^\infty$ reference connection on $P$, and $a(t) = A(t) - A_1 \in \Omega^1(X;\ad P)$ for $t\in [0, \tau)$. The questions of local existence, uniqueness, and regularity of solutions to the Yang-Mills heat equation are developed from many different points of view in Chapter \ref{chapter:Existence_uniqueness_regularity_Yang-Mills_heat_equation} and for initial data of different regularities. While the question of local well-posedness for quasilinear parabolic equations is often dismissed as `standard', it is difficult to find a standard reference, so we establish all of the results one might need in detail and in as much generality as possible. The elegant approach due to Struwe \cite{Struwe_1994}, particular to a base manifold of dimension less than or equal to four, requires a package of existence, uniqueness, and regularity results for a solution $a(t) \in \Omega^1(X;\ad P)$, for $t\in [0,\infty)$, to the linear heat equation,
$$
\frac{\partial a}{\partial t} + \left(d_{A_1}^*d_{A_1} + d_{A_1}d_{A_1}^*\right)a = f, \quad a(0) = a_0,
$$
or more simply,
\begin{equation}
\label{eq:Linear_heat_equation_rough_Laplacian}
\frac{\partial a}{\partial t} + \nabla_{A_1}^*\nabla_{A_1} a = f, \quad a(0) = a_0,
\end{equation}
and those are developed in Section \ref{sec:Sell_You_4-2-3_heat_equation_vector_bundle}.

Section \ref{sec:Local_well-posedness_yang_mills_heat_equation} develops the required local posedness results for the Yang-Mills heat equation when the initial data is assumed to be of sufficient regularity that more standard versions of the theory of analytic semigroups and nonlinear evolution equations in Banach spaces \cite{Sell_You_2002} are applicable. In addition, we derive lower bounds for the minimal lifetime of a solution in terms of suitable norms of the initial data and an \apriori estimate for the length of the trajectory defined by a solution.

In Section \ref{sec:Critical-exponent_parabolic_Sobolev_spaces_linear_parabolic_operator_vector_bundle_manifold}, we construct a family of
critical-exponent parabolic Sobolev spaces, over a base manifold $X$ of arbitrary dimension $d\geq 2$, by analogy with the families of critical-exponent elliptic Sobolev spaces frequently employed by Taubes \cite{TauPath, TauFrame, TauStable, TauConf, TauGluing} and developed further by the author in \cite{FeehanSlice} when $X$ has dimension $d=4$. We then develop an \apriori estimate for a linear parabolic operator on sections of a vector bundle over a closed, Riemannian, smooth manifold, $X$, of dimension $d\geq 2$.

Section \ref{sec:Struwe_3and4} develops local well-posedness for the Yang-Mills heat equation with initial data of minimal regularity, based partly on ideas of Kozono, Maeda, and Naito \cite{Kozono_Maeda_Naito_1995}, using the theory of analytic semigroups and nonlinear evolution equations in Banach spaces \cite{Sell_You_2002}, when $X$ has dimension $d\geq 2$, and partly on ideas of Struwe \cite{Struwe_1994}, when $X$ has dimension $d\leq 4$. Our version of Struwe's approach \cite{Struwe_1994} employs the Contraction Mapping Principle, based on existence of strong solutions,
$$
a \in L^2(0,T;H_{A_1}^2(X;\Lambda^1\otimes\ad P)) \cap H^1(0,T; L^2(X;\Lambda^1\otimes\ad P)),
$$
to the linear heat equation \eqref{eq:Linear_heat_equation_rough_Laplacian}.

Finally, in Section \ref{sec:Local_existence_yang_mills_gradient_coulomb_gauge_flow}, we establish local existence, uniqueness, and regularity for solutions to Yang-Mills gradient flow, employing the Donaldson-DeTurck trick \cite{DonASD, DK} to pass from a solution to the Yang-Mills heat equation to the Yang-Mills gradient flow equation. The question of uniqueness of a solution to the Yang-Mills gradient flow equation can be approached either via the method of Kozono, Maeda, and Naito \cite{Kozono_Maeda_Naito_1995} or that of Struwe \cite{Struwe_1994}, who derives a complementary, but essentially stronger result by a different method. We describe both approaches and add detail to selected calculations.

All of these approaches rely on Donaldson's version \cite{DonASD}, \cite[Equation (6.3.3)]{DK} of the DeTurck Trick for Ricci flow \cite{DeTurck_1983} to convert the Yang-Mills \emph{gradient flow} equation in $\Omega^1(X; \Lambda^1\otimes\ad P)$,
\begin{equation}
\label{eq:Yang-Mills_gradient_flow_equation_introduction}
\frac{\partial A}{\partial t} + d_A^*F_A = 0,
\end{equation}
with initial data $A(0) = A_0 \in \Omega^1(X; \ad P)$, to the Yang-Mills \emph{heat} equation \eqref{eq:Yang-Mills_heat_equation_introduction}. For $a(t)$, with $t \in (0,\tau)$, solving \eqref{eq:Yang-Mills_heat_equation_introduction}, one defines a family of gauge transformations, $u(t) \in \Aut P$, by solving
\begin{equation}
\label{eq:Struwe_18_introduction}
u(t)^{-1}\circ \frac{\partial u(t)}{\partial t} = - d_{A(t)}^*a(t), \quad\forall\, t \in (0,\tau), \quad u(0) = \id_P,
\end{equation}
and discovers that $\widetilde A(t) = u(t)^*A(t)$ solves the Yang-Mills gradient flow equation \eqref{eq:Yang-Mills_gradient_flow_equation_introduction}.

Chapter \ref{chapter:Lojasiewicz-Simon_gradient_inequality_and_stability_and_convergence} comprises our development of the {\L}ojasiewicz-Simon gradient inequality, convergence, and stability for (pure) abstract gradient systems in Banach spaces, together with a derivation of two versions of the {\L}ojasiewicz-Simon gradient inequality for the Yang-Mills energy functional over a base manifold, $X$, of dimension in the range $d \geq 2$. While our approach to convergence and stability for abstract gradient (and later gradient-like) systems is inspired by that of Huang \cite{Huang_2006}, we employ a fundamentally different paradigm (a generalization of R\r{a}de's \cite[Lemma 7.3]{Rade_1992}) that, we hope, lends itself to a much wider range of applications than permitted by the approach of Huang \cite{Huang_2006}.

By way of introduction to the infinite-dimensional {\L}ojasiewicz-Simon gradient inequality and its application to gradient systems, we recall the

\begin{thm}[Long-time existence and convergence of a solution to a gradient system in $\RR^n$]
\label{thm:Lojasiewicz_1984_introduction}
\cite{Lojasiewicz_1963}, \cite[Theorem 1]{Lojasiewicz_1984}
Let $f$ be an analytic, non-negative function on a neighborhood of the origin in $\RR^n$ such that $f(0) = 0$. Then there exists a neighborhood, $U$, of the origin such that each trajectory, $y_{x_0}(t)$, with $y_{x_0}(0) = x_0$, of the system,
$$
\frac{dy}{dt} = -(\grad f)(y(t)),
$$
is defined on $[0, \infty)$, has finite length, and converges uniformly to a point in $Z := \{x \in U: \nabla f(x) = 0\}$ as $t \to \infty$.
\end{thm}

To prove his result, {\L}ojasiewicz uses the following version of his gradient inequality \cite{Lojasiewicz_1965}: \emph{If $F$ is an analytic function on a neighborhood of the origin in $\RR^n$ such that $F(0) = 0$, then there exists a constant, $\theta \in (0, 1)$, and a neighborhood, $U$, of the origin such that}
$$
|\grad F(x)| \geq |F(x)|^\theta, \quad \forall\, x \in U.
$$
In \cite{Simon_1983}, Simon developed an infinite-dimensional version of the finite-dimensional gradient inequality due to {\L}ojasiewicz. Simon used his gradient inequality to establish global existence and convergence of smooth solutions to a certain class of gradient-like (or pseudo-gradient) flow equations arising in geometric analysis, including harmonic map gradient flow when the target Riemannian manifold is real analytic. However, in the form developed in \cite{Simon_1983}, Simon's results are not applicable to the gradient flow for the Yang-Mills energy functional.

Since the publication of Simon's seminal article \cite{Simon_1983}, there have been many attempts to generalize his gradient inequality to abstract gradient inequalities for functionals defined on open subsets of Hilbert or Banach spaces, including work of Chill, Huang, and Takac and many others cited in \cite{Huang_2006}. One variant of the {\L}ojasiewicz-Simon gradient inequality, developed specifically for the Yang-Mills energy functional in his Ph.D. thesis \cite{RadeThesis} for a base manifold of dimension two or three, is due to R\r{a}de \cite{Rade_1992}, though his version appears not well known. Section \ref{sec:Huang_2} reviews one of the most general, abstract versions of the {\L}ojasiewicz-Simon gradient inequality (\cite[Theorem 2.4.2 (i)]{Huang_2006}) following the presentation due to Huang \cite[Chapter 2]{Huang_2006}. In Section \ref{sec:Huang_2_yangmills}, we then use this to generalize R\r{a}de's gradient inequality for the Yang-Mills energy functional when $X$ has dimensions $d=2, 3$ and also extend it to the case of $d\geq 4$.

Suppose that $\sE:\sU\subset\sX \to \RR$ is a functional on an open subset, $\sU$, of a real, reflexive Banach space, $\sX$, that is in turn continuously embedded and dense in a real Hilbert space, $\sH$, with analytic gradient map, $\sE':\sU\subset\sX \to \sH$. Since $\sX \hookrightarrow \sH$ is a continuous embedding, we also have $\sH \cong \sH' \hookrightarrow \sX'$, where $\sX'$ denotes the Banach-space dual of $\sX$ and the isometric isomorphism between $\sH$ and $\sH'$ is given by the canonical identification. One requires, in addition, that $\sE$ obeys certain technical properties listed in Hypothesis \ref{hyp:Huang_2-4_H1_H2_H3}. If $\varphi \in \sU$ is a critical point of $\sE$, that is, $\sE'(\varphi) = 0$, then there are positive constants, $c$, $\sigma$, and $\theta \in [1/2, 1)$ such that
\begin{equation}
\label{eq:Simon_2-2_introduction}
\|\sE'(u)\|_\sH \geq c|\sE(u) - \sE(\varphi)|^\theta, \quad \forall\, u \in \sU \hbox{ such that } \|u-\varphi\|_\sX < \sigma.
\end{equation}
We refer the reader to Theorem \ref{thm:Huang_2-4-2} for the precise statement. Our Theorems \ref{thm:Rade_proposition_7-2_L2} and \ref{thm:Rade_proposition_7-2} specialize Theorem \ref{thm:Huang_2-4-2} to the case of the Yang-Mills energy functional with
$$
\sH = L^2(X; \Lambda^1\otimes\ad P),
$$
and different choices of Sobolev space for the Banach space, $\sX$. When $\sX = H_{A_1}^1(X; \Lambda^1\otimes\ad P)$ and $d=4$, where $A_1$ a $C^\infty$ reference connection on $P$, then our Theorem \ref{thm:Rade_proposition_7-2} provides the desired result: If $A_\infty$ is a $C^\infty$ Yang-Mills connection on $P$, then there are positive constants $c$, $\sigma$, and $\theta \in [1/2,1)$, depending on $A_1$, $g$, and $P$ such that
\begin{equation}
\label{eq:Rade_7-1_introduction}
\|\sE'(A)\|_{L^2(X)} \geq c|\sE(A) - \sE(A_\infty)|^\theta,
\end{equation}
for all connections, $A$, of class $H^2$ on $P$ such that
\begin{equation}
\label{eq:Rade_7-1_neighborhood_introduction}
\|A - A_\ym\|_{H_{A_1}^1(X)} < \sigma,
\end{equation}
recalling that $\sE(A) = \frac{1}{2}\|F_A\|_{L^2(X)}^2$ and $\sE'(A) = d_A^*F_A$.

Section \ref{sec:Huang_3_and_5_gradient-system} returns to the abstract setting of Section \ref{sec:Huang_2} and develops convergence, global existence, and stability results for strong solutions, $u(t)\in \sX$ for $t \geq 0$ with $\dot u(t)\in \sH$ for a.e. $t > 0$, to the Cauchy problem for a gradient system,
\begin{equation}
\label{eq:Gradient_system_Cauchy_problem}
\dot u(t) = -\sE'(u(t)) \quad\hbox{(in $\sH$, for a.e. $t >0$)}, \quad u(0) = u_0.
\end{equation}
While our development formally mirrors that of Huang \cite{Huang_2006}, we replace his hypotheses (either \cite[Equation (3.10)]{Huang_2006} for gradient-like flow or \cite[Equation (3.10$'$)]{Huang_2006} for gradient flow) with the following abstract hypothesis modeled on R\r{a}de's \apriori estimate \cite[Lemma 7.3]{Rade_1992} for the length of a Yang-Mills gradient flow trajectory.

\begin{hyp}[\Apriori interior estimate for a trajectory]
\label{hyp:Abstract_apriori_interior_estimate_trajectory_introduction}
Let $C$ and $\mu$ be positive constants and let $T \in (0, \infty]$. Given an open subset, $\sU \subset \sX$, let $u:[0,T)\to \sU$ and $\dot u:[0,T) \to \sH$ be continuous. We say that $\dot u$ obeys an \apriori \emph{interior estimate on $(0, T]$} if, for every $S \geq 0$ and $\delta > 0$ obeying $S+\delta \leq T$, the map $\dot u:[S+\delta,T) \to \sX$ is Bochner integrable and there holds
\begin{equation}
\label{eq:Abstract_apriori_interior_estimate_trajectory_introduction}
\int_{S+\delta}^T \|\dot u(t)\|_\sX\,dt \leq C(1+\delta^{-\mu})\int_S^T \|\dot u(t)\|_\sH\,dt.
\end{equation}
\end{hyp}

See Hypothesis \ref{hyp:Abstract_apriori_interior_estimate_trajectory} for a slightly more general statement, with weaker regularity requirements on the path, $u(t)$.
The abstract \apriori estimate \eqref{eq:Abstract_apriori_interior_estimate_trajectory_introduction} is realized for a solution to Yang-Mills gradient flow by our Lemmata \ref{lem:Rade_7-3_abstract_interior_L1_in_time_V2beta_space_time_derivative_interior}, \ref{lem:Rade_7-3}, and \ref{lem:Rade_7-3_L1_in_time_H2beta_in_space_apriori_estimate_by_L1_in_time_L2_in_space}. For example, Lemma \ref{lem:Rade_7-3} yields the estimate \eqref{eq:Abstract_apriori_interior_estimate_trajectory_introduction} with $\sX = H_{A_1}^1(X;\Lambda^1\otimes\ad P)$ and $\sH = L^2(X;\Lambda^1\otimes\ad P)$ and $\mu=1/2$, when $X$ has dimension $d$ in the range $2\leq d\leq 5$ and $A(t)$ is a strong solution to Yang-Mills gradient flow that obeys
$$
\|A(t) - A_1\|_{H^1_{A_1}(X)} < \eps_1 \quad\forall\, t \in (S, T),
$$
for a small enough positive constant, $\eps_1$, depending on the $C^\infty$ reference connection, $A_1$, $d$, and the Riemannian metric, $g$, on $X$, so
$$
\int_{S+\delta}^T \|\dot A(t)\|_{H_{A_1}^1(X)}\,dt
\leq
C(1 + \delta^{-1/2})\int_S^T \|\dot A(t)\|_{L^2(X)}\,dt,
$$
where $C$ depends at most on $A_1$, $d$, and $g$.

Given Hypothesis \ref{hyp:Abstract_apriori_interior_estimate_trajectory_introduction}, we then obtain an abstract version (Lemma \ref{lem:Huang_3-3-4}) of Simon's \cite[Lemma 1]{Simon_1983} and a replacement of Huang's \cite[Lemma 3.3.4]{Huang_2006}, giving an estimate for the integral,
$$
\int_0^T \|\dot u(t)\|_\sH\,dt
\leq
\int_{\sE(u(T))}^{\sE(u(0))} \frac{1}{c|s - \sE(\varphi)|^\theta} \,ds,
$$
and hence an estimate for the integral,
$$
\int_\delta^T \|\dot u(t)\|_\sX\,dt,
$$
when $u$ is a solution to a gradient system obeying a {\L}ojasiewicz-Simon gradient inequality along the trajectory, $u(t)$ for $t\in [0,T)$, and $c>0$, and $\theta \in [1/2,1)$ and $\sigma>0$ are the {\L}ojasiewicz-Simon constants. These estimates are the key ingredient which allow us to obtain Theorem \ref{thm:Huang_3-3-6}, an abstract version of R\r{a}de's \cite[Proposition 7.4]{Rade_1992}, Simon's \cite[Theorem 2]{Simon_1983} and a replacement for Huang's \cite[Theorem 3.3.6]{Huang_2006}, giving a convergence alternative for a global solution, $u:[0,\infty)\to\sX$, with initial data, $u(0) = u_0$, sufficiently close to a critical point, $\varphi$, of $\sE$. Furthermore, the rate of convergence of the solution, $u(t)$, to a limit, $u_\infty \in \sX$, as $t\to\infty$ can be estimated using our enhancement (Theorem \ref{thm:Huang_3-4-8}) of Huang's \cite[Theorem 3.4.8]{Huang_2006}.

One calls a critical point $\varphi \in \sU$ of $\sE$ a \emph{ground state} if $\sE$ attains its minimum on $\sU$ at this point, that is,
$$
\sE(\varphi) = \inf_{u\in \sU}\sE(u).
$$
When $u_0$ is sufficiently close to a ground state, $\varphi$, our version, Theorem \ref{thm:Huang_5-1-1}, of Huang's \cite[Theorem 5.1.1]{Huang_2006} yields global existence of a solution, $u:[0,\infty)\to\sX$, to the Cauchy problem \eqref{eq:Gradient_system_Cauchy_problem} for an abstract gradient system obeying a {\L}ojasiewicz-Simon gradient inequality, together with convergence of $u(t)$ to a limit, $u_\infty\in\sX$, as $t\to\infty$ in the sense that
$$
\lim_{t \to \infty} \|u(t) - u_\infty\|_\sX = 0
\quad\hbox{and}\quad
\int_1^\infty\|\dot u(t)\|_\sX\,dt < \infty.
$$
Moreover, our analogue, Theorem \ref{thm:Huang_5-1-2}, of \cite[Theorem 5.1.2]{Huang_2006} yields Lyapunov stability for a ground state, $\varphi$; if $\varphi$ is isolated or a cluster point of the path, $\{u(t):t\geq 0\}$, then $\varphi$ is uniformly asymptotically stable.

Chapter \ref{chapter:Huang_5-1_application_pure_Yang-Mills_gradient_flow} establishes global existence and convergence for Yang-Mills gradient flow near a local minimum. Section \ref{sec:Simon_4_linear_theory} establishes certain key \apriori estimates for a variational solution to a linear heat equation generalizing those of R\r{a}de \cite{Rade_1992} and Simon \cite{Simon_1983}. In Section \ref{sec:Rade_Lemma_7-3_generalization}, we establish the required \apriori estimates for lengths of Yang-Mills gradient-like flow lines and justify Hypothesis \ref{hyp:Abstract_apriori_interior_estimate_trajectory_introduction} in the setting of Yang-Mills gradient and gradient-like flow, whose applications we discuss further in the sequel. Finally, in Section \ref{sec:Application_abstract_gradient_system_results_Yang-Mills_energy_functional} we complete the proof of Theorem \ref{mainthm:Yang-Mills_gradient_flow_global_existence_and_convergence_started_near_minimum}.

In Chapter \ref{chapter:Gradient_flow_arbitrary_initial_energy}, we commence our examination of the problems of global existence and convergence of smooth solutions to Yang-Mills gradient flow when the energy, $\sE(A_0)$, of the initial connection is arbitrary. For this purpose, we revert once more to the abstract setting of Section \ref{sec:Huang_3_and_5_gradient-system} and in Section \ref{sec:Huang_3_and_5_gradientlike_system} we generalize those results from the case of gradient to gradient-like systems. In particular, we introduce a modification of the pure gradient flow for $\sE$ modeled on that of Simon's \cite[Equations (0.1) or (3.1)]{Simon_1983}.

As in Definition \ref{defn:Strong_solution_to_gradientlike_system}, let $\sE':\sU\subset \sX \to \sH$ be a gradient map associated with a $C^1$ functional, $\sE:\sU \to \RR$, where $\sU$ is an open subset of a real, reflexive Banach space, $\sX$, that is continuously embedded and dense in a Hilbert space, $\sH$, and let $R:[0, T) \to \sH$ be a continuous map for some $0 < T \leq \infty$. We call a trajectory, $u:[0,T)\to \sU\subset \sX$ a \emph{strong solution} of a \emph{pseudogradient} or \emph{gradient-like system} for $\sE$ if
\begin{equation}
\label{eq:Simon_0-1_introduction}
\dot u(t) = -\sE'(u(t)) + R(t), \quad\hbox{a.e. } t \in (0,T), \quad u(0) = u_0 \in \sU,
\end{equation}
as an equation in $\sH$. As one might expect from the more concrete situation considered \cite{Simon_1983}, the perturbation, $R(t)$, must obey certain estimates in order to allow us to again draw any conclusions regarding convergence or global existence and these are described in Hypothesis \ref{hyp:Huang_3-10}, which requires $R$ to be such that
\begin{subequations}
\label{eq:Huang_3-10a_introduction}
\begin{align}
\label{eq:Huang_3-10a_gradE_innerproduct_dotu_introduction}
(-\sE'(u(t)), \dot u(t))_\sH &\geq \|\sE'(u(t))\|_\sH \|\dot u(t)\|_\sH + F'(t),
\\
\label{eq:Huang_3-10a_gradE_norm_H_introduction}
\|\sE'(u(t))\|_\sH &\geq \|\dot u(t)\|_\sH + G'(t),  \quad\hbox{for a.e. } t \geq 0,
\end{align}
\end{subequations}
where the functions $F, G: [0, \infty) \to [0, \infty)$ are absolutely continuous, non-increasing (non-negative) functions satisfying
\begin{subequations}
\label{eq:Huang_3-10b_introduction}
\begin{align}
\label{eq:Huang_3-10b_limit_F+G_zero_introduction}
\lim_{t\to \infty} F(t) + G(t) &= 0,
\\
\label{eq:Huang_3-10b_integral_zero_to_infinity_of_phiF_finite_introduction}
\int_0^\infty \phi(F(t))\,dt &< \infty,
\end{align}
\end{subequations}
where $\phi(x) := c|x - a|^\theta$, $x \in \RR$, for constants $c > 0$ and $a \in \RR$ and $\theta \in [1/2, 1)$. For the gradient-like system \eqref{eq:Simon_0-1_introduction}, one can choose
\begin{subequations}
\label{eq:Huang_F_and_G_in_terms_of_Simon_R_introduction}
\begin{align}
\label{eq:Huang_F_in_terms_of_Simon_R_introduction}
F(t) &= \frac{1}{2}\int_t^\infty\|R(s)\|_\sH^2\,ds,
\\
\label{eq:Huang_G_in_terms_of_Simon_R_introduction}
G(t) &= \int_t^\infty\|R(s)\|_\sH\,ds, \quad\forall\, t \geq 0.
\end{align}
\end{subequations}
We then generalize our results from gradient to gradient-like systems, albeit with more complicated statements and proofs, including:
\begin{inparaenum}[\itshape a\upshape)]
\item a growth estimate for a solution to a gradient-like system;
\item convergence of gradient-like flow near a local minimum;
\item convergence rate near a critical point; and
\item global existence and convergence of a solution to a gradient-like system started near a local minimum.
\end{inparaenum}

In order to try to apply our results for abstract gradient-like flow to Yang-Mills gradient flow we must carefully examine the nature of the singularities that can potentially occur in Yang-Mills gradient flow, whether in finite or infinite time. Consequently, in Section \ref{sec:Local_apriori_estimates_Yang-Mills_heat_and_gradient_flows} we review and in some cases extend the local \apriori estimates for Yang-Mills heat and gradient flows developed by Kozono, Maeda, and Naito \cite{Kozono_Maeda_Naito_1995}, Schlatter \cite{Schlatter_1997}, and Struwe \cite{Struwe_1994}. A local version, Lemma \ref{lem:Schlatter_2-4_and_Struwe_3-6_local}, of the continuous extension Lemma \ref{lem:Schlatter_2-4_and_Struwe_3-6} of Struwe plays a particularly important role in our analysis in the sequel and this is proved in Section \ref{subsec:Local_version_continuous_extension_lemma_Schlatter_2-4}.

As we noted earlier, the possibility of bubble-tree limits appears to have been overlooked in previous treatments \cite{Kozono_Maeda_Naito_1995, Schlatter_1997} of the limiting behavior of Yang-Mills gradient flow near bubble points, so we devote Section \ref{sec:Uhlenbeck_and_bubbletree_limits_Yang-Mills_gradient_flow_four-manifold} to a detailed analysis of
Uhlenbeck and bubble-tree limits for solutions to Yang-Mills gradient flow over four-dimensional manifolds, building on the local \apriori estimates reviewed in Section \ref{sec:Local_apriori_estimates_Yang-Mills_heat_and_gradient_flows}.

Given a Riemannian metric, $g$, on $X$, and a solution, $A(t)$, to Yang-Mills gradient flow on a principal $G$-bundle, $P$, over $X$ with initial data, $A_0$, it is natural to compare $A(t)$ with a solution, $\bar A(t)$, to Yang-Mills gradient flow on $P$ for a nearby Riemannian metric, $\bar g$, and the same initial data, $A_0$. Thus, in Section \ref{sec:Continuity_solution_Yang-Mills_heat_equation_wrt_Riemannian_metric}, we examine continuity and stability of solutions to the Yang-Mills heat equation with respect to variations in the Riemannian metric.

In Chapter \ref{chapter:Yang-Mills_gradient-like_flow_four_manifolds_applications}, we consider Yang-Mills gradient-like flow over four-dimensional manifolds and applications, motivated by the results of our analysis of Uhlenbeck limits for Yang-Mills gradient flow in Section \ref{sec:Uhlenbeck_and_bubbletree_limits_Yang-Mills_gradient_flow_four-manifold}. Section \ref{sec:Yang-Mills_gradient-like_flow_over_four_sphere} develops a selection of results for Yang-Mills gradient-like flow over the four-dimensional sphere, $S^4$, with its standard round metric of radius one. In Section \ref{sec:Yang-Mills_gradient-like_flow_Kahler_surface}, we consider Yang-Mills gradient-like flow over compact K{\"a}hler surfaces and apply our analysis to complete the proof of Theorem \ref{mainthm:Yang-Mills_gradient_flow_global_existence_and_convergence_started_with arbitrary_initial_energy}. We conclude Section \ref{sec:Yang-Mills_gradient-like_flow_Kahler_surface} with a proof Corollary \ref{cor:Donaldson_Kronheimer_6-2-7_and_6-2-14_convergence}.

Finally, in Chapter \ref{chapter:Solution_to_anti-self-dual_equation_applications}, we develop a series of results related to the problem of solving the anti-self-dual equation,
$$
F^{+,g}(A + a) = 0,
$$
for a perturbation, $a \in \Omega^1(X;\ad P)$, such that the connection $A+a$ is exactly $g$-anti-self-dual, given an approximately $g$-anti-self-dual connection, $A$, on a principal $G$-bundle, $P$, over a closed, four-dimensional, oriented, smooth manifold, $X$, with Riemannian metric, $g$. Our analysis in Section \ref{sec:Solution_to_anti-self-dual_equation} builds on previous results of the author and Leness in \cite{FLKM1}, Sedlacek \cite{Sedlacek}, and especially Taubes
\cite{TauSelfDual, TauPath, TauIndef, TauFrame, TauStable}. In Section \ref{sec:Yang-Mills_gradient_flow_initial_connection_almost_minimal_energy_proofs}, we combine the results of Section \ref{sec:Solution_to_anti-self-dual_equation} and Theorem \ref{mainthm:Yang-Mills_gradient_flow_global_existence_and_convergence_started_near_minimum} to establish Corollaries \ref{maincor:Yang-Mills_gradient_flow_global_existence_and_convergence_started_near_connection with almost minimal energy_and_vanishing_H2+}, \ref{maincor:Yang-Mills_gradient_flow_global_existence_and_convergence_started_near_connection_with_almost_minimal_energy_and_g_is_good_Lp_version}, and \ref{maincor:Yang-Mills_gradient_flow_global_existence_and_convergence_started_near_connection_with_almost_minimal_energy_and_g_is_good_L2sharp_version}.

\chapter[Comparison of global existence and convergence results]{Comparison of global existence and convergence results for gradient and heat flows in geometric analysis}
\label{chapter:Comparison_global_existence_convergence_results_gradient_and_heat_flows_geometric_analysis}
Difficulties in the analysis of Yang-Mills gradient flow also arise in other nonlinear evolution equations in geometric analysis. We cannot hope in this short survey to provide anything remotely approaching a comprehensive survey of the main results on global existence and convergence for nonlinear evolution equations in geometric analysis and mathematical physics. Instead, we focus on a few well-known results for some of the gradient flows whose behavior appears to most closely emulate that of Yang-Mills gradient flow in dimension four, including harmonic map gradient flow for Riemann surfaces and Yamabe scalar curvature heat flow. We shall mainly concentrate on harmonic map gradient flow for Riemann surfaces, indicate some parallel results for Yamabe scalar curvature flow, and briefly mention results for a few other geometric flows. Our purpose is solely to try to help us understand why one gradient flow might exhibit global existence or convergence and another while an apparently similar gradient flow might not possess those properties.

\section{Gradient flow for the harmonic map energy functional}
\label{sec:Gradient_flow_harmonic_map_energy_functional}
The harmonic map energy functional is conformally invariant when the source manifold has dimension two, like the Yang-Mills energy functional in dimension four, in the case of maps from a Riemann surface, $M$, into a Riemannian manifold of dimension $d \geq 2$. Therefore it is especially interesting to try to understand similarities and differences between their respective gradient flows. Indeed, a comparison between harmonic map gradient flow in dimension two and the Yang-Mills gradient flow in dimension four was undertaken by Grotowski and Shatah in \cite{Grotowski_Shatah_2007} for precisely this reason.

\subsection{Comparison between the harmonic map and Yang-Mills $L^2$-energy functionals}
\label{subsec:Comparison_between_harmonic_map_and_Yang-Mills_energy_functionals}
Modulo an application of the Donaldson-DeTurck trick in the case of Yang-Mills gradient flow, the harmonic map and Yang-Mills gradient flow equations are both quasilinear parabolic but the nonlinearity in the case of the harmonic map gradient flow appears slightly less well behaved and that may account for the sharp distinction between the finite-time blow up behaviors in the presence of rotational symmetry exhibited by Chang, Ding, and Ye in \cite{Chang_Ding_Ye_1992} and Schlatter, Struwe, and Tahvildar-Zadeh \cite{Schlatter_Struwe_Tahvildar-Zadeh_1998}.

To try to better understand the reasons for the differences between the harmonic map and Yang-Mills gradient flows, let us examine some differences between the nonlinearities appearing in the respective $L^2$-energy functionals. By comparing the gradient of the harmonic map energy functional \eqref{eq:Gradient_harmonic_map_operator} with that of the Yang-Mills energy functional, it is interesting to note that, in the latter case, the nonlinearity is less severe, as one encounters terms of the form $a\times a\times a$ and $a\times\nabla_{A_1}a$, but not $\nabla_{A_1}a\times\nabla_{A_1}a$ or even $a\times\nabla_{A_1}a\times\nabla_{A_1}a$. In calculations involving the gradient flows of these functionals, the $L^2$ norms of the gradients, $\|\sE'(f)\|_{L^2(M)}$ and $\|\sE'(A)\|_{L^2(X)}$, play a significant role and so these differences are important.

\subsection{The {\L}ojasiewicz-Simon gradient inequality and Simon's results on harmonic map gradient flow}
\label{subsec:Lojasiewicz-Simon_gradient_inequality_harmonic_map_energy_functional}
Simon applied his {\L}ojasiewicz-Simon gradient inequality \cite{Simon_1983} to obtain the following global existence and convergence result\footnote{Our summary of Simon's results is taken from \cite[Theorem 1.15]{KwonThesis}.} for harmonic map gradient flow in \cite{Simon_1985}.

\begin{thm}[Global existence and convergence for harmonic map gradient flow]
\cite{Simon_1983, Simon_1985}
\label{thm:Simon_in_Kwon_thesis_1-15}
Let $M$ and $N$ be closed, Riemannian, smooth manifolds. Suppose $f: M \times [0,\infty) \to N$ is a smooth solution to harmonic map gradient flow such that there exist a sequence of times, $\{t_m\}_{m\in\NN} \subset [0,\infty)$ with $t_m\to \infty$ as $m\to\infty$, and a smooth harmonic map, $f_\infty: M \to N$, such that $\lim_{m\to\infty} \|f(t_m) - f_\infty\|_{C^k(M)} = 0$, for any $k \in \NN$. If either the target $N$ is \emph{real analytic} or $f_\infty$ is \emph{integrable} (in the sense of \cite{Simon_1985}), then $\lim_{t\to\infty} \|f(t) - f_\infty\|_{C^k(M)} = 0$, for any $k \in \NN$.  Furthermore, if the initial map, $f_0$, is sufficiently close to a locally energy minimizing map (in the $C^k(M)$ topology for suitably large $k$), then a solution, $f(t)$, exists for all time and asymptotically converges to a smooth, locally energy minimizing map, $f_\infty$.
\end{thm}

Versions of the {\L}ojasiewicz-Simon gradient inequality were subsequently used by Irwin \cite{IrwinThesis} and Kwon \cite{KwonThesis}, both former Ph.D. students of Simon, as well as by Liu and Yang \cite{Liu_Yang_2010} and Topping \cite{ToppingThesis, Topping_1997}.

\subsection{Finite-time blow up for $S^1$-equivariant harmonic map gradient flow from a two-dimensional disk into $S^2$}
\label{subsec:Chang_Ding_Ye}
In the context of harmonic maps from a Riemann surface, there are the following well-known results of Chang, Ding, and Ye \cite{Chang_Ding_Ye_1992}. Let
$$
D := \{x = (x_1, x_2 , 0) \in \RR^3: |x|^2 < 1\} \quad\hbox{and}\quad S^2 := \{x \in \RR^3: |x|^2 = 1\}.
$$
Given $u_0 \in C^1(\bar D , S^2)$, they consider the initial-boundary value problem for the harmonic map gradient flow equation,
\begin{equation}
\label{eq:Chang_Ding_Ye_1}
\frac{\partial u}{\partial t} = \Delta u + |\nabla u|^2 u,
\quad
\begin{cases}
u(0, x) = u_0(x), & x \in \bar D,
\\
u(t,x) = u_0(x), & x \in \partial D.
\end{cases}
\end{equation}
It is well known that \eqref{eq:Chang_Ding_Ye_1} admits a unique classical solution, $u$, that solves the problem on $[0, \sT) \times \bar D$, where $\sT \equiv \sT(u_0) \in (0, \infty]$ is the maximal existence time for $u$. If $\sT < \infty$, one says that the solution, $u$, blows up in finite time. In \cite{Chang_Ding_1991}, it is shown that if the initial map, $u_0$, has the following
symmetric form,
\begin{equation}
\label{eq:Chang_Ding_Ye_2}
u_0(x) = \left( \frac{x}{|x|} \sin h_0(|x|), \cos h_0(|x|) \right), \quad x \in \bar D,
\end{equation}
then the solution of \eqref{eq:Chang_Ding_Ye_2} has the form,
$$
u(t,x) = \left( \frac{x}{|x|} \sin h(t,|x|), \cos h(t,|x|) \right), \quad (t,x) \in [0,\sT) \times \bar D.
$$
In particular, from \cite{Chang_Ding_Ye_1992} the study of $S^1$-equivariant harmonic map gradient flow for maps $u:\bar D \to S^2$ reduces to the study of the following initial-boundary value problem for a nonlinear, singular ordinary differential equation \cite[Equation (3)]{Chang_Ding_Ye_1992},
\begin{subequations}
\label{eq:Chang_Ding_Ye_3}
\begin{gather}
\label{eq:Chang_Ding_Ye_3a}
\frac{\partial h}{\partial t} = \frac{\partial^2 h}{\partial r^2} + \frac{1}{r}\frac{\partial h}{\partial r}
- \frac{1}{r^2}\cos h \sin h,
\\
\label{eq:Chang_Ding_Ye_3b}
\begin{cases}
h(0,r) = h_0(r), & 0 < r < 1,
\\
h(t,0) = h_0(0) = 0, & t \geq 0,
\\
h(t,1) = h_0(1) = b, & t \geq 0,
\end{cases}
\end{gather}
\end{subequations}
where $b \in \RR$ and $h_0 \in C^1([0,1];\RR)$. If $|h_0| \leq \pi$ on $[0,1$, Chang and Ding \cite{Chang_Ding_1991} showed that the solution, $h$ of \eqref{eq:Chang_Ding_Ye_3} exists for all $t \geq 0$ and consequently the solution $u$ of \eqref{eq:Chang_Ding_Ye_1} with initial map, $u_0$, given by \eqref{eq:Chang_Ding_Ye_2} is a global solution. However, if $|h_0(1)| > \pi$, Chang, Ding, and Ye \cite{Chang_Ding_Ye_1992} showed that the solution, $h$ of \eqref{eq:Chang_Ding_Ye_3} blows up in finite time and consequently the same is true for the solution $u \in C^1(\bar D, S^2)$ of \eqref{eq:Chang_Ding_Ye_1} with initial map, $u_0$, given by \eqref{eq:Chang_Ding_Ye_2}.

In \cite[pp. 514--515]{Chang_Ding_Ye_1992}, the authors show how to extend their finite-time blow up result to the case of maps $u \in C^1(S^2, S^2)$ with sufficiently large initial energy.

\subsection{No finite-time blow up for $\SO(4)$-equivariant Yang-Mills gradient flow}
\label{subsec:Schlatter_Struwe_Tahvildar-Zadeh}
In \cite{Schlatter_Struwe_Tahvildar-Zadeh_1998}, motivated by the results of Chang, Ding, and Ye \cite{Chang_Ding_Ye_1992}, Schlatter, Struwe, and Tahvildar-Zadeh show that the study of $\SO(4)$-equivariant Yang-Mills gradient flow over the closed unit ball, $\bar B \subset \RR^4$, reduces to the study of the following initial-boundary value problem for a nonlinear, singular ordinary differential equation \cite[Equations (4), (5), and (6)]{Schlatter_Struwe_Tahvildar-Zadeh_1998},
\begin{subequations}
\label{eq:Schlatter_Struwe_Tahvildar-Zadeh-4_and_5_and_6}
\begin{gather}
\label{eq:Schlatter_Struwe_Tahvildar-Zadeh-4_and_5}
\frac{\partial f}{\partial t} = \frac{\partial^2 f}{\partial r^2} + \frac{1}{r}\frac{\partial f}{\partial r}
- \frac{2}{r^2}f(1-f)(2-f),
\\
\label{eq:Schlatter_Struwe_Tahvildar-Zadeh-6}
\begin{cases}
f(0, r) = f_0(r), & 0 < r < 1,
\\
f(t,0) = f_0(0) = 0, & t \geq 0,
\\
f(t,1) = f_0(1) = b, & t \geq 0,
\end{cases}
\end{gather}
\end{subequations}
where $b \in \RR$ is arbitrary and $f_0 \in H^1_{\loc}(0,1;\RR)$. The condition that $A_0$ is a finite-energy, $\SO(4)$-equivariant Yang-Mills connection on $B\times G$ is equivalent to $f_0$ obeying $f_0(0) = 0$ and
$$
\sE(f) := \int_0^1 e(f_0)r\,dr < \infty,
$$
where
$$
e(f_0) := \left|\frac{\partial f}{\partial r}\right|^2 + \frac{f^2(2-f)^2}{2r^2}.
$$
Their principal result, \cite[Theorem 1]{Schlatter_Struwe_Tahvildar-Zadeh_1998}, is that \eqref{eq:Schlatter_Struwe_Tahvildar-Zadeh-4_and_5_and_6} admits a unique, global smooth solution. Consequently, there is no finite-time blow up for $\SO(4)$-equivariant Yang-Mills gradient flow on $\bar B \subset \RR^4$.

\section{Conformal Yamabe scalar curvature heat flow}
\label{sec:Yamabe_scalar_curvature_heat_flow}
There is an extensive literature on questions related to Yamabe scalar curvature heat flow. We shall only focus our attention on results of Brendle \cite{Brendle_2005, Brendle_2007invent, Brendle_2008sdg, Brendle_2011jjm}, Carlotto, Chodosh, and Rubinstein \cite{Carlotto_Chodosh_Rubinstein_2015}, Schwetlick and Struwe \cite{Schwetlick_Struwe_2003}, and Ye \cite{Ye_1994jdg}. In our brief account here we rely on the survey articles by Brendle \cite{Brendle_2008sdg, Brendle_2011jjm}; we refer the reader to those articles for additional background omitted here. We begin by recalling the celebrated

\begin{conj}[Yamabe Conjecture]
\label{conj:Yamabe_conjecture}
(See Yamabe \cite{Yamabe_1960}.)
Let $X$ be a closed manifold of dimension $d\geq 3$, and let $g_0$ be a Riemannian metric on $X$. Then there exists a metric $g$ on $X$ which is conformally equivalent to $g_0$ and has constant scalar curvature.
\end{conj}

Conjecture \ref{conj:Yamabe_conjecture} was proved by Trudinger \cite{Trudinger_1968}, Aubin \cite{Aubin_1976}, and Schoen \cite{Schoen_1984}. Using earlier analysis developed by Trudinger, Aubin proved Conjecture \ref{conj:Yamabe_conjecture} when $d \geq 6$ and $(X,g_0)$ is not locally conformally flat. Schoen completed the proof of Conjecture \ref{conj:Yamabe_conjecture} using the Positive Mass Theorem to settle the rest of cases. Bahri \cite{Bahri_1993na, Bahri_1993lnpam} provided a different proof of Conjecture \ref{conj:Yamabe_conjecture} when $g_0$ is locally conformally flat. Bahri's approach does not rely on the Positive Mass Theorem.

Hamilton pioneered the gradient flow approach to the Yamabe problem, which we now outline. Let $g(t)$, for $t\geq 0$, be a $C^\infty$ path of $C^\infty$ Riemannian metrics on $X$. One says that $g(t)$ solves the \emph{unnormalized Yamabe flow} if
\begin{equation}
\label{Brendle_2011_4}
\frac{\partial g}{\partial t}
=
-R_{g(t)} g(t).
\end{equation}
We use $R_{g(t)}$ to denote the scalar curvature of the metric $g(t)$. One says that $g(t)$ solves the \emph{normalized Yamabe flow} if
\begin{equation}
\label{Brendle_2011_5}
\frac{\partial g}{\partial t}
=
-\left(R_{g(t)} - r_{g(t)}\right) g(t).
\end{equation}
In the preceding equation, $r_{g(t)}$ is the mean value of the scalar curvature of the metric, $g(t)$ and is thus defined by
$$
r_{g(t)}
:=
\frac{\displaystyle\int_X R_{g(t)}\,d\vol_{g(t)}}{\Vol_{g(t)}(X)}.
$$
Equations \eqref{Brendle_2011_4} and \eqref{Brendle_2011_5} are regarded as equivalent because any solution of the equation \eqref{Brendle_2011_4} can be transformed into a solution of \eqref{Brendle_2011_5} using rescaling. Following Brendle in this summary, we restrict our attention to the normalized Yamabe flow \eqref{Brendle_2011_5}. We observe that
\begin{equation}
\label{Brendle_2011_6}
\frac{\partial R_g}{\partial t}
=
(d-1)\Delta_{g(t)} R_{g(t)} - R_{g(t)} \left(R_{g(t)} - r_{g(t)}\right),
\end{equation}
when $g(t)$ solves \eqref{Brendle_2011_5}.

The Yamabe flow preserves the conformal class of a Riemannian metric and this observation allows one to write
$$
g(t) = u(t)^{4/(d-2)}g_0,
$$
with respect to one fixed metric, $g_0$. The scalar curvature of the metric $g(t)$ in the preceding equation is then found to be
$$
R_{g(t)}
=
u(t)^{(d+2)/(d-2)} \left(-\frac{4(n-1)}{n-2}\Delta_{g_0}u(t) + R_{g_0}u(t)\right).
$$
Therefore, Equation \eqref{Brendle_2011_5} for the metric $g(t)$ reduces to an equation for the scalar conformal factor, namely
\begin{equation}
\label{Brendle_2011_7}
\frac{\partial}{\partial t} \left( u(t)^{(d+2)/(d-2)} \right)
=
\frac{d+2}{4}\left(\frac{4(d-1)}{d-2}\Delta_{g_0}u(t) - R_{g_0}u(t)
+ r_{g(t)}u(t)^{(d+2)/(d-2)} \right),
\end{equation}
which is a nonlinear parabolic partial differential equation. One can interpret the Yamabe flow as the \emph{gradient flow} for the \emph{Yamabe energy function} given by
$$
\sE_{g_0}(u)
=
\frac{\displaystyle\int_X\left(\frac{4(d-1)}{d-2}|du|_{g_0}^2 + R_{g_0}u^2\right)\,d\vol_{g_0}}
{\displaystyle\left(\int_X u^{2d/(d-2)}\,d\vol_{g_0}\right)^{(d-2)/d}}.
$$
For any positive $u \in C^\infty(X)$, one can show that
$$
\sE_{g_0}(u)
=
\frac{\displaystyle\int_X R_g\,d\vol_g}
{\displaystyle\left(\Vol_g(X)\right)^{(d-2)/d}},
$$
where we write $g = u^{4/(d-2)}g_0$. Hence the Yamabe function, $\sE_{g_0}$, arises as the restriction of the \emph{normalized Einstein-Hilbert action} to the conformal class of the Riemannian $g_0$.

As usual, one is interested in the longtime behavior of the Yamabe flow, so we recall the

\begin{thm}[Hamilton (circa 1989, unpublished)]
\label{thm:Brendle_2011_2-1}
Let $X$ be a closed manifold of dimension $d\geq 3$. Given any initial metric $g_0$ on $X$, the Yamabe flow \eqref{Brendle_2011_5} admits a solution which is defined for all $t \geq 0$.
\end{thm}

Given Theorem \ref{thm:Brendle_2011_2-1}, which asserts global existence of $g(t)$ for $0 \leq t < \infty$, it is natural to consider the asymptotic behavior of $g(t)$ as $t\to\infty$ and in this context, one has the

\begin{conj}[Hamilton]
\label{conj:Brendle_2011_2-2}
Let $(X, g_0)$ be a closed Riemannian manifold of dimension $d \geq 3$, and let $g(t)$ be the unique solution of the Yamabe flow with initial metric, $g_0$. Then $g(t)$ converges to a metric of constant scalar curvature as $t \to \infty$.
\end{conj}

Chow \cite{Chow_1992} proved a special case of Conjecture \ref{conj:Brendle_2011_2-2}: If $(X,g_0)$ is locally conformally flat and has positive Ricci curvature, then the Yamabe flow, $g(t)$, with initial data, $g_0$, converges to a metric of constant scalar curvature. Ye \cite{Ye_1994jdg} obtained the following improvement of Chow's result.

\begin{thm}[Ye, \cite{Ye_1994jdg}]
\label{thm:Brendle_2011_2-3}
Let $(X,g_0)$ be a closed Riemannian manifold of dimension $d \geq 3$, and let $g(t)$ be the unique solution of the Yamabe flow with initial metric $g_0$. If $(X, g_0)$ is locally conformally flat, then $g(t)$ converges to a metric of constant scalar curvature as $t \to \infty$.
\end{thm}

Schwetlick and Struwe \cite{Schwetlick_Struwe_2003} in turn refined Theorem \ref{thm:Brendle_2011_2-3} by replacing the assumption that $(X, g_0)$ is locally conformally flat with a weaker assumption on the energy of the initial data.

\begin{thm}[Schwetlick and Struwe, \cite{Schwetlick_Struwe_2003}]
\label{thm:Brendle_2011_2-4}
Let $(X, g_0)$ be a closed Riemannian manifold of dimension $d$, where $3 \leq d \leq 5$. Moreover, let $g(t)$ be the unique solution to the Yamabe flow with initial metric $g_0$. If the Yamabe energy of $g_0$ is less than
$$
\left[Y(X,g_0)^{d/2} + Y(S^d)^{d/2}\right]^{2/d},
$$
then $g(t)$ converges to a metric of constant scalar curvature as $t \to \infty$.
\end{thm}

In Theorem \ref{thm:Brendle_2011_2-4}, we use $Y(S^d)$ to denote the infimum of the Yamabe energy function on the sphere, $S^d$, with its standard round metric of radius one and $Y(X,g_0)$ is the infimum of the Yamabe energy function on $(X,g_0)$. Brendle further improved Theorem \ref{thm:Brendle_2011_2-4} using an argument that allowed him to omit the hypothesis on the Yamabe energy of $g_0$, and obtain convergence of the Yamabe flow for any initial data, $g_0$:

\begin{thm}[Brendle, \cite{Brendle_2005}]
\label{thm:Brendle_2011_2-5}
Let $(X,g_0)$ be a closed Riemannian manifold of dimension $d$. We assume that either $3 \leq d \leq 5$ or $d\geq 3$ and $(X, g_0)$ is locally conformally flat. Moreover, let $g(t)$, for $t \geq 0$, be the unique solution to the Yamabe flow with initial metric $g_0$. Then $g(t)$ converges to a metric of constant scalar curvature as $t \to \infty$.
\end{thm}

It is interesting to recall Brendle's summary of his main ideas involved in the proofs of Theorems \ref{thm:Brendle_2011_2-4} and \ref{thm:Brendle_2011_2-5}. Brendle writes $g(t) = u(t)^{4/(d-2)}g_0$, for $u(t) \in C^\infty(X)$ obeying $u(t) > 0$ on $X$ and $t\geq 0$. If $u(t)$ is uniformly bounded from above, convergence of $g(t)$ follows from the asymptotic analysis of gradient flows developed by Simon \cite{Simon_1983}. Therefore, it is enough to prove that $u(t)$ is uniformly bounded from above. To achieve this goal, Brendle employs a blow-up analysis common in energy bubbling problems. Hence, he assumes that $\sup_X u(t_k) \to \infty$ for a sequence $\{t_k\}_{k\geq 1} \subset (0,\infty)$ with $t_k \to \infty$. Brendle then applies a theorem due to Struwe \cite{Struwe_1984} to prove that $g(t)$ develops bubble singularities for at most finite number of times.

In Theorem \ref{thm:Brendle_2011_2-4}, the hypothesis on the energy of the initial metric $g_0$ ensures that $g(t)$ develops at most one bubble. Schwetlick and Struwe excluded the latter possibility by appealing to the Positive Mass Theorem due to Schoen and Yau \cite{Schoen_Yau_1981}. When $g_0$ has higher initial energy, the problem is much more difficult, as one might expect from our experience with harmonic map or Yang-Mills gradient flows. The difficulty is due to the fact that $g(t)$ could form multiple bubbles. Brendle excludes the possibility of multiple bubbles by again applying the Positive Mass Theorem.

Finally, we describe another convergence result due to Brendle for the Yamabe flow in dimension $d \geq 6$. Suppose that $(X,g_0)$ is a closed Riemannian manifold of dimension $d \geq 6$. Define 
$$
\sZ := \left\{p\in X \left|\limsup_{x\to p} \dist_{g_0}(p,x)^{[(d-6)/2]}\left|\sW_{g_0}(x)\right| = 0 \right.\right\},
$$
where $\sW_{g_0}$ is the Weyl curvature of $g_0$ and $\dist_{g_0}(\cdot,\cdot)$ is the Riemannian distance function. One can show that $\sZ$ depends at most on the conformal class of the metric $g_0$.

\begin{thm}[Brendle, \cite{Brendle_2007invent}]
\label{thm:Brendle_2011_2-6}
Let $(X, g_0)$ be a closed Riemannian manifold of dimension $d\geq 6$. Assume that either $X$ is spin or $\sZ = \emptyset$. Moreover, let $g(t)$, for $t \geq 0$, be the unique solution to the Yamabe flow with initial metric $g_0$. Then $g(t)$ converges to a metric of constant scalar curvature as $t\to\infty$.
\end{thm}

The limiting metric provided by Theorem \ref{thm:Brendle_2011_2-6} need not be a global minimum of the Yamabe energy function. 

The previously cited articles do not explicitly use the {\L}ojasiewicz-Simon gradient inequality, but in his proof of Theorem \ref{thm:Brendle_2011_2-5}, Brendle appeals to Simon's convergence result in \cite{Simon_1983}, as he notes in \cite{Brendle_2011jjm}. Moreover, Brendle's \cite[Lemma 6.5]{Brendle_2005} (see \cite[Equation (100)]{Brendle_2005}) relies on a (finite-dimensional) version of the {\L}ojasiewicz gradient inequality. More recently, Carlotto, Chodosh, and Rubinstein \cite{Carlotto_Chodosh_Rubinstein_2015} have used the (infinite-dimensional) {\L}ojasiewicz-Simon gradient inequality to characterize the rate of convergence of a converging volume-normalized Yamabe flow in terms of Morse theoretic properties of the limiting metric. The authors appeal to the abstract formulation of the {\L}ojasiewicz-Simon gradient inequality due to Chill \cite{Chill_2003}.

\section{Ricci flow}
\label{sec:Ricci_flow}
In this section, we comment on how our results in Section \ref{subsec:Main_results_gradient_system_Banach_space} for abstract gradient flow might be applied to study properties of Ricci flow.

\subsection{Ricci flow as a gradient flow}
\label{subsec:Ricci_flow_gradient_flow}
Let $M$ be a closed, smooth Riemannian manifold. Following Perelman \cite{Perelman_2002}, we recall that Hamilton's \emph{Ricci flow},
\begin{align}
\label{eq:Ricci_flow_equation}
\frac{\partial g}{\partial t} &= -2\Ric(g(t)),
\\
\label{eq:Ricci_flow_initial_data}
g(0) &= g_0,
\end{align}
is equivalent, modulo a $C^\infty$ family of time-varying diffeomorphisms of $M$, to the gradient flow of the functional
\begin{equation}
\label{eq:Perelman_lambda_functional}
\lambda(g) = \inf_{\begin{subarray}{c} f \in C^\infty(M) \\ \int_M e^{-f}\,d\vol_g=1 \end{subarray}}
\int_M \left(\Scal(g) + |df|_g^2\right)e^{-f}\,d\vol_g.
\end{equation}
The functional, $\lambda$, is nondecreasing under the Ricci flow and the critical points of $\lambda$ are precisely the Ricci-flat metrics \cite{Kleiner_Lott_2008, Perelman_2002}.
Recall that Ricci-flat metrics were already known to be critical points (but not extrema) of the Einstein-Hilbert functional and fixed points of the Ricci flow \cite{Hamilton_1982}.

For further discussion of Ricci curvature flow and its applications, we refer the reader to the original articles by Hamilton \cite{Hamilton_1982, Hamilton_1995} and Perelman \cite{Perelman_2002, Perelman_2003a, Perelman_2003b}, together with the numerous expositions for detailed discussions of the relevant analysis, including those of Cao and Zhu \cite{Cao_Zhu_2006, Cao_Zhu_2006erratum, Cao_Zhu_2006arxiv}, Chow et al. \cite{Chow_Knopf_2004, Chow_Lu_Ni_2006, Chow_etal_Ricci_Flow_I, Chow_etal_Ricci_Flow_II}, Kleiner and Lott \cite{Kleiner_Lott_2008}, Morgan and Fong \cite{Morgan_Fong_2010}, Morgan and Tian \cite{Morgan_Tian_2007}, M\"uller \cite{Muller_2006}, notes by Tao \cite{Tao_2008lecture1, Tao_2008lecture8}, and Topping \cite{Topping_Lectures_ricci_flow}.

\subsection{{\L}ojasiewicz-Simon gradient inequality for Perelman's $\lambda$ functional}
\label{subsec:Lojasiewicz–Simon gradient inequality_Perelman_lambda}
A {\L}ojasiewicz-Simon gradient inequality for $\lambda$ has been established, by directly adapting the proof due to Simon of his \cite[Theorem 3]{Simon_1983}, by several researchers, including Ache \cite{Ache_2011arxiv}, Haslhofer \cite{Haslhofer_2012cvpde}, Haslhofer and M\"uller \cite{Haslhofer_Muller_2014}, and Sun and Wang \cite{Sun_Wang_2015}; see also Kr{\"o}ncke \cite{Kroncke_2015cvpde, Kroncke_2013arxiv}. The statement below is taken from Haslhofer \cite{Haslhofer_2012cvpde} and Haslhofer and M\"uller \cite{Haslhofer_Muller_2014}.

\begin{thm}[{\L}ojasiewicz-Simon gradient inequality for Perelman's $\lambda$ functional]
\label{thm:Lojasiewicz–Simon gradient inequality_Perelman_lambda}
\cite[Appendix A]{Ache_2011arxiv}
\cite[Theorem 3]{Haslhofer_2012cvpde},
\cite[Theorem 3]{Haslhofer_Muller_2014},
\cite[Lemma 3.1]{Sun_Wang_2015}
Let $g_\RF$ be a Ricci-flat Riemannian metric on a closed, smooth manifold, $M$. Then there exist constants  $c \in [1,\infty)$, and $\sigma \in (0,1]$, and $\theta \in [1/2,1)$ such that, for all $g$ in the space of Riemannian metrics on $M$ obeying
\[
\|g - g_\RF\|_{C_{g_\RF}^{2,\alpha}(M)} < \sigma,
\]
one has
\begin{equation}
\label{eq:Lojasiewicz–Simon gradient inequality_Perelman_lambda}
\|\Ric(g) + \Hess_g(f_g)\|_{L^2(M,e^{-f_g}\, d\vol_g)} \geq c|\lambda(g)|^\theta,
\end{equation}
where $f_g$ is the minimizer in \eqref{eq:Perelman_lambda_functional} realizing $\lambda(g)$. If $g_\RF$ is \emph{integrable}, then $\theta=1/2$.
\end{thm}

Recall that a Ricci-flat Riemannian metric, $g_\RF$, is \emph{integrable} if for every symmetric two-tensor $h$ in the kernel of the
linearization of $\Ric$, one can find a curve of Ricci-flat metrics with initial velocity $h$; see
\cite[Section 12]{Besse_1987} and \cite[Sections 1 and 3]{Haslhofer_2012cvpde} for further details.

The overall strategy of the proof of Theorem \ref{thm:Lojasiewicz–Simon gradient inequality_Perelman_lambda} is originally due to R\r{a}de \cite{Rade_1992} in the context of the Yang-Mills $L^2$-energy functional on a principal $G$-bundle $P$ over a closed manifold of dimension two or three. Namely, one must first restrict the functional to a slice for the action of group of diffeomorphisms of $M$ (automorphisms of the $G$-bundle $P$ in \cite{Rade_1992}) and, in this setting, one appeals to the slice theorem due to Ebin \cite{Ebin_1968, Ebin_1970} and Palais \cite{Palais_1961} for the quotient space of Riemannian metrics on $M$ modulo diffeomorphisms of $M$. The group of $C^{2,\alpha}$ diffeomorphisms, $\Diff(M)$, acts on the open subspace of $C^{2,\alpha}$ Riemannian metrics, $\Met(M)$, of the Banach space, $C^{2,\alpha}(\Sym^2(T^*M))$, by pullback. For a given Riemannian metric, $g \in \Met(M)$, we have an $L^2$-orthogonal decomposition,
\[
C_g^{2,\alpha}(\Sym^2(T^*M)) = \Ker \divg_g \oplus \Ran \divg_g^*.
\]
One then has the

\begin{thm}[Ebin-Palais slice theorem]
\label{thm:Ebin_Palais_slice_theorem}
\cite{Ebin_1968, Ebin_1970}, \cite{Palais_1961}
Let $M$ be a closed, smooth manifold and $\alpha \in (0,1)$. If $g_0 \in \Met(M)$, then there exists a constant $\eps = \eps(g_0,\alpha_0) \in (0,1]$ with the following significance. If
\[
\|g - g_0\|_{C_{g_0}^{2,\alpha}(M)} < \eps,
\]
then there exists a $C^{2,\alpha}$ diffeomorphism, $\varphi \in \Diff(M)$, such that $g = \varphi^*\hat g$, where
\[
\|\hat g - g_0\|_{C_{g_0}^{2,\alpha}(M)} < \eps
\quad\text{and}\quad
\divg_{g_0}(\hat g - g_0) = 0.
\]
\end{thm}

Ebin employs H\"older spaces in his version of the slice theorem (as customary in the study of Ricci flow), whereas Palais employs Sobolev spaces; both are valid choices. Rather than establish Theorem \ref{thm:Lojasiewicz–Simon gradient inequality_Perelman_lambda} by arguing that it can be justified by adapting Simon's proof of his \cite[Theorem 3]{Simon_1983} (as done in the cited references), we show here that it can instead be easily derived as a corollary of our abstract {\L}ojasiewicz-Simon gradient inequality Theorem \ref{thm:Huang_2-4-5_introduction}.

\begin{proof}[Proof of Theorem \ref{thm:Lojasiewicz–Simon gradient inequality_Perelman_lambda} using Theorem \ref{thm:Huang_2-4-5_introduction}]
The gradient of $\lambda$ at $g \in \Met(M)$ is given by the \emph{Perelman first variation formula} (for example, \cite[Equation (1.5)]{Haslhofer_2012cvpde}),
\begin{equation}
\label{eq:Perelman_first_variation_formula}
\lambda'(g)[h] = -\int_M \langle \Ric(g) + \Hess_g(f_g), h \rangle_g e^{-f_g}\,d\vol_g,
\quad\forall\, h \in C_g^{2,\alpha}(\Sym^2(T^*M)).
\end{equation}
The Hessian of $\lambda$ at a Ricci-flat metric, $g_\RF \in \Met(M)$ (thus $\lambda(g_\RF)=0$ and $\lambda'(g_\RF)=0$), is given by the \emph{Cao-Hamilton-Ilmanen second variation formula} \cite[Theorem 1.1]{Cao_Hamilton_Ilmanen_2004arxiv}, \cite[Equations (1.3) and (2.13)]{Haslhofer_2012cvpde},
\begin{equation}
\label{eq:Cao-Hamilton-Ilmanen_second_variation_formula}
\lambda''(g_\RF)[h,h] = \frac{1}{2\vol_{g_\RF}(M)}\int_M \langle L_{g_\RF}^L h, h \rangle_{g_\RF} e^{-f_{g_\RF}}\,d\vol_{g_\RF},
\quad\forall\, h \in \Ker \divg_{g_\RF},
\end{equation}
where $L_g$ is the \emph{Lichnerowicz Laplacian},
\[
L_g h_{\alpha\beta} := \Delta_g h_{\alpha\beta} + 2R_{\alpha\mu\beta\nu}(g)h_{\mu\nu}, \quad\forall\, h \in \Sym^2(T^*X),
\]
and $\lambda''(g_\RF)[h,h] = 0$ for $h \in \Ran \divg_{g_\RF}^*$. The general expression for $\lambda''(g_\RF)[h,k]$, for possibly distinct $h,k \in C_{g_\RF}^{2,\alpha}(\Sym^2(T^*M))$, can be obtained from the expression for $\lambda''(g_\RF)[h,h]$ by polarization, noting that the Hessian operator is symmetric.

Given a Ricci-flat metric, $g_\RF \in \Met(M)$, we choose $\sX = \sX(g_\RF)$ and $\sH = \sH(g_\RF)$ in the statement of Theorem \ref{thm:Huang_2-4-5_introduction} by setting
\begin{align*}
\sX &:= \Ker\left\{\divg_{g_\RF}: C_{g_\RF}^{2,\alpha}(\Sym^2(T^*M)) \to C_{g_\RF}^{1,\alpha}(\Sym^2(T^*M))\right\},
\\
\sH &:= L^2(\Sym^2(T^*M), e^{-f_{g_\RF}}\,d\vol_{g_\RF}),
\end{align*}
and observe that $\sX \subset \sH$ is a continuous embedding (in fact, dense). The functional $\lambda:\Met(M)\to \RR$ is analytic \cite{Haslhofer_2012cvpde}, \cite{Haslhofer_Muller_2014}, \cite{Sun_Wang_2015} and thus its restriction, $\lambda:\sX\to\RR$, to a slice is also analytic. Moreover, the elliptic differential operator,
\[
L_{g_\RF}: C_{g_\RF}^{2,\alpha}(\Sym^2(T^*M)) \to C_{g_\RF}^\alpha(\Sym^2(T^*M)),
\]
is Fredholm by Lemma \ref{lem:Gilkey_1-4-5_Holder} and has index zero since $L_g$ is a compact perturbation of the Laplace operator, $\Delta_g$.
It follows that the Hessian, $\lambda''(g_\RF):\sX\to\sX^*$, is Fredholm with index zero.
Because $\sX \subset \sH$ is a continuous, dense embedding, then $\sH\subset\sX^*$ is a continuous embedding with norm $\kappa \in [1,\infty)$ and Theorem \ref{thm:Huang_2-4-5_introduction} provides constants $c \in [1,\infty)$, $\sigma \in (0,1]$, and $\theta \in [1/2,1)$ such that
\[
\kappa\|\lambda'(g)\|_\sH \geq \|\lambda'(g)\|_{\sX*} \geq c|\lambda(g) - \lambda(g_\RF)|^\theta = c|\lambda(g)|^\theta.
\]
The conclusion is now immediate from the expression for $\lambda'(g)$ and by relabelling $c$.
\end{proof}

\subsection{Ricci flow as a gradient system}
\label{subsec:Ricci_flow_as_gradient_flow}
Given $g \in \Met(M) \subset C^\infty(\Sym^2(T^*M))$, one defines an inner product on the tangent space $T_g C^\infty(\Sym^2(T^*M)) = C^\infty(\Sym^2(T^*M))$ by setting
\begin{equation}
\label{eq:Perelman_inner_product}
(h,k)_{L^2(M,e^{-f_g}\,d\vol_g)} := \int_M \langle h, k \rangle_g e^{-f_g}\,d\vol_g, \quad\forall\, h, k \in C^\infty(\Sym^2(T^*M)),
\end{equation}
where $f_g$ is the corresponding minimizer in \eqref{eq:Perelman_lambda_functional}. We denote $\sH_g := L^2(\Sym^2(T^*M),e^{-f_g}\,d\vol_g)$.

Suppose that $\tilde g(t) \in C^\infty(\Sym^2(T^*M))$ is a $C^\infty$ function of $t \in (0,T)$ and that $f_{\tilde g}(t)$ is the corresponding minimizer in \eqref{eq:Perelman_lambda_functional} for each $t \in (0,T)$. It is known (for example, \cite{Haslhofer_2012cvpde}, \cite{Haslhofer_Muller_2014}) that $f_{\tilde g}(t) \equiv f_{\tilde g(t)}\in C^\infty(M)$ and is a $C^\infty$ function of $t\in (0,T)$. The gradient system for the functional $-2\lambda:\Met(M) \to \RR$ defined by the \emph{time-varying} inner product \eqref{eq:Perelman_inner_product} determined by $\tilde g(t)$ is
\[
\left(\frac{\partial \tilde g}{\partial t}, h\right)_{\sH_{\tilde g(t)}}  =  2\lambda'(\tilde g(t))[h], \quad\forall\, h \in C^\infty(\Sym^2(T^*M)).
\]
From \eqref{eq:Perelman_first_variation_formula}, this gradient system is
\begin{multline*}
\int_M \left\langle \frac{\partial\tilde g}{\partial t}, h\right\rangle_{\tilde g(t)}\,e^{-f_{\tilde g}(t)}\,d\vol_{\tilde g(t)}
\\
=
-2\int_M \left\langle\Ric(\tilde g(t) + \Hess_{\tilde g(t)}(f_{\tilde g(t)}(t)), h\right\rangle_{g(t)}\,e^{-f_{\tilde g}(t)}\,d\vol_{\tilde g(t)}, \quad\forall\, h \in C^\infty(\Sym^2(T^*M)),
\end{multline*}
or equivalently,
\begin{multline*}
\int_M \frac{\partial\tilde g}{\partial t}(\eta)\,e^{-f_{\tilde g}(t)}\,d\vol_{\tilde g(t)}
\\
=
-2\int_M \left(\Ric(\tilde g(t)) + \Hess_{\tilde g(t)}(f_{\tilde g}(t))\right)(\eta)\,e^{-f_{\tilde g}(t)}\,d\vol_{\tilde g(t)}, \quad\forall\, \eta \in C^\infty(\Sym^2(TM)).
\end{multline*}
Therefore, the gradient flow for the functional $-2\lambda:\Met(M) \to \RR$ with respect to the time-varying bilinear pairing,
\[
C^\infty(\Sym^2(T^*M))\times C^\infty(\Sym^2(TM)) \ni (h,\eta) \mapsto \int_M h(\eta)\, e^{-f_{\tilde g}(t)}\,d\vol_{\tilde g(t)},
\]
is given
\begin{equation}
\label{eq:Modified_Ricci_flow}
\frac{\partial\tilde g}{\partial t}(t) = -2\Ric(\tilde g(t)) - 2\Hess_{\tilde g(t)}(f_{\tilde g}(t)), \quad\forall\, t \in (0,T).
\end{equation}
This is called \emph{modified Ricci flow} and it is known to be equivalent to (pure) Ricci flow \eqref{eq:Ricci_flow_equation} modulo pullback by a $C^\infty$ family of diffeomorphisms of $M$. While it is likely that our results for abstract gradient flow summarized in Section \ref{subsec:Main_results_gradient_system_Banach_space} can be applied to \eqref{eq:Modified_Ricci_flow}, the application is not immediate because the inner product \eqref{eq:Perelman_inner_product} on $\sH_{\tilde g(t)}$ is time-varying. We refer the reader to \cite{Haslhofer_2012cvpde}, \cite{Haslhofer_Muller_2014} for results --- obtained via the {\L}ojasiewicz-Simon gradient inequality for the Perelman $\lambda$-functional --- on convergence and stability of Ricci flow near a Ricci-flat connection that is a local maximum for $\lambda$.

\subsection{Alternative approach to interpreting Ricci flow as a gradient system}
\label{subsec:Ricci_flow_as_gradient_flow_alternative}
The approach taken in Section \ref{subsec:Ricci_flow_as_gradient_flow} to interpreting modified Ricci flow as a gradient system has the unattractive feature that one must allow a family of Hilbert spaces with time-varying inner products. There is an alternative approach to interpreting a (possibly different) modified Ricci flow as a gradient flow, but this comes at the cost of now defining $f(t) \in C^\infty(M)$ as the solution to a backward heat equation, so we cannot prescribe initial data at time $t=0$ but only terminal data at some time $T > 0$, given a solution $g(t)$ to Ricci flow \eqref{eq:Ricci_flow_equation} on $(0,T)$. This alternative approach (and the attendant difficulties) are described by Andrews and Hopper \cite[Chapter 10]{Andrews_Hopper_2011}, Kleiner and Lott \cite[Sections 6 and 10]{Kleiner_Lott_2008}, M\"uller \cite{Muller_2006}, Tao \cite{Tao_2008lecture8}, and Topping \cite[Section 6.4]{Topping_Lectures_ricci_flow}.

Following Topping \cite[Section 6.4]{Topping_Lectures_ricci_flow}, we suppose that we are given a solution $g(t)$, for $t\in [0,T)$, to Ricci flow \eqref{eq:Ricci_flow_equation} with initial data $g(0) = g_0$. For some choice $f_T \in C^\infty(M)$ of terminal data, one solves the backward heat equation\footnote{We adhere to the geometer's sign convention for the Laplace operator, $\Delta_g = d^{*_g}d$ on $C^\infty(M)$, opposite to the analyst's sign convention employed by Topping in \cite[Equation (6.4.6)]{Topping_Lectures_ricci_flow}.}
\begin{equation}
\label{eq:Topping_6-4-6}
\frac{\partial f}{\partial t} = \Delta_{g(t)} f(t) - |df|_{g(t)}^2 - \Scal(g(t)), \quad t \in (0, T),
\end{equation}
with terminal data $f(T) = f_T$. Modulo pullback by a $C^\infty$ family of diffeomorphisms of $M$, the system \eqref{eq:Ricci_flow_equation}, \eqref{eq:Topping_6-4-6} is equivalent to the system,
\begin{align}
\label{eq:Topping_6-4-9}
\frac{\partial\hat g}{\partial t}(t) &= -2\Ric(\hat g(t)) - 2\Hess_{\hat g(t)}(\hat f(t)),
\\
\label{eq:Topping_6-4-10}
\frac{\partial \hat f}{\partial t} &= \Delta_{\hat g(t)} \hat f(t)- \Scal(\hat g(t)), \quad t \in (0, T).
\end{align}
Naturally, this process of demonstrating equivalence works in reverse and in the cited references, one starts with system \eqref{eq:Topping_6-4-9}, \eqref{eq:Topping_6-4-10}, following the observation that given a $C^\infty$ function $\hat g(t) \in \Met(M)$ of $t\in(0,T)$, the equation \eqref{eq:Topping_6-4-10} is equivalent to the constraint that the measure (for example, see Topping \cite[Equation (6.4.4)]{Topping_Lectures_ricci_flow}),
\begin{equation}
\label{eq:Topping_6-4-4}
dm := e^{-\hat f}\,d\vol_{\hat g},
\end{equation}
is constant with respect to time. Given a static measure $dm$ on $M$, one defines
\begin{equation}
\label{eq:Perelman_F_functional}
\sF_m(\hat g,\hat f) := \int_M \left(\Scal(\hat g) + |d\hat f|_{\hat g}^2\right)\,dm, \quad\forall\, (\hat g,\hat f) \in \Met(M)\times C^\infty(M).
\end{equation}
The gradient system for the functional $\Met(M) \ni \hat g \mapsto -2\sF_m(\hat g, \hat f) \in \RR$ (with $\hat f$ determined by $\hat g$ and $\omega$ via \eqref{eq:Topping_6-4-4}) is
\[
\int_M \frac{\partial \hat g}{\partial t}(\eta)\,dm  =  2\sF_m'(\hat g(t))[\eta], \quad\forall\, \eta \in C^\infty(\Sym^2(TM)).
\]
More explicitly, using \eqref{eq:Perelman_first_variation_formula}, this is
\[
\int_M \frac{\partial \hat g}{\partial t}(\eta)\,dm  =  -2\int_M\left(\Ric(\hat g(t)) + \Hess_{\hat g(t)}(\hat f(t))\right)(\eta)\,dm, \quad\forall\, \eta \in C^\infty(\Sym^2(TM)),
\]
and this in turn yields \eqref{eq:Topping_6-4-9}.

\subsection{Applications of the {\L}ojasiewicz-Simon gradient inequality to convergence, global existence, and stability for Ricci flow}
\label{subsec:Applications_Lojasiewicz–Simon gradient inequality_Ricci_flow}
Given Theorem \ref{thm:Lojasiewicz–Simon gradient inequality_Perelman_lambda}, one should be able to apply our abstract results in Section \ref{subsec:Main_results_gradient_system_Banach_space} for convergence, global existence, and stability for gradient flows to Ricci flow. We briefly outline how this program could be implemented in order to recover and perhaps extend results due to Ache \cite{Ache_2011arxiv}, Haslhofer \cite{Haslhofer_2012cvpde}, Haslhofer and M\"uller \cite{Haslhofer_Muller_2014}, Sesum \cite{Sesum_2006}, and others.

While the Ricci flow equation \eqref{eq:Ricci_flow_equation} is not parabolic, it can be modified in a standard way to yield a parabolic, quasi-linear, second-order partial differential system (the \emph{Ricci-DeTurck flow}) via the DeTurck trick (see \cite[Section 3]{DeTurck_1983} for $M$ of dimension $d=3$ and \cite[Section 4]{DeTurck_1983} for $M$ of dimension $d\geq 4$) via the pull back action on a path of Riemannian metrics, $g(t)$, by a suitable $C^\infty$ family of diffeomorphisms, $\varphi(t) \in \Diff(M)$, just as we describe elsewhere in this monograph for Yang-Mills gradient flow. The abstract Hypothesis \ref{hyp:Abstract_apriori_interior_estimate_trajectory_introduction} can then be verified for the resulting Ricci-DeTurck flow, $\tilde g(t)$, by Lemma \ref{lem:Rade_7-3_abstract_interior_L1_in_time_V2beta_space_time_derivative_interior} for a solution to a nonlinear evolution equation on a Banach space $\calV$ of the form
\[
\frac{du}{dt} + \cA u = \cF(t,u(t)), \quad t\geq 0, \quad u(0) = u_0,
\]
where $\cA$ is a positive, sectorial, unbounded operator on $\cW$ with domain $\calV^2 \subset \cW$ and the nonlinearity, $\cF$, has suitable properties.

The main technical difficulty in making this program precise is that if one interprets Ricci flow as a gradient system using the approach of Section \ref{subsec:Ricci_flow_as_gradient_flow}, one must extend the abstract results in Section \ref{subsec:Main_results_gradient_system_Banach_space} to the setting of a time-varying family of Hilbert spaces $\sH_t$ depending on the solution $u(t)$. On the other hand, if one interprets Ricci flow as a gradient system using the approach of Section \ref{subsec:Ricci_flow_as_gradient_flow_alternative}, one must address the fact that one (generally) cannot solve \eqref{eq:Topping_6-4-10} forward in time for $\hat f(t)$, given an arbitrary initial condition $\hat f(0) = \hat f_0 \in C^\infty(M)$, unlike the pure Ricci flow equation \eqref{eq:Ricci_flow_equation} for $g(t)$, given an arbitrary initial condition $g(0) = g_0 \in \Met(M)$.

\section{Other gradient flows}
\label{sec:Other_gradient_flows}
We briefly mention here some of the many related geometric flows to which the {\L}ojasiewicz-Simon gradient inequality and techniques developed in this monograph are applicable, together with selected previous results and further references for those flows.

\subsection{Anti-self-dual curvature flow over four-dimensional manifolds}
\label{subsec:Anti-self-dual_curvature_flow_over_4-manifolds}
In \cite{TauPath, TauStable}, Taubes established global existence and convergence for \emph{anti-self-dual curvature flow}. Because we review some of his results in Section \ref{subsec:Perturbation_and_solution_anti-self-dual_equation}, we defer a discussion of that flow to our review in that section. This flow may be viewed as a solution to a gradient system for an energy functional, but Taubes does not use the {\L}ojasiewicz-Simon gradient inequality to obtain his results -- although there appears to be no impediment to such an analysis.

\subsection{Chern-Simons gradient flow over three-dimensional manifolds}
\label{subsec:Chern-Simons_gradient_flow_over_3-manifolds}
Suppose $X$ is a closed, four-dimensional, Riemannian, smooth manifold with an end isometric to $Y\times [0,\infty)$, where $Y$ is a closed, three-dimensional, Riemannian, smooth manifold. A finite-energy solution to the anti-self-dual equation for a connection on a principal $G$-bundle over $X$ may be interpreted, over the infinite tube, as the solution to \emph{Chern-Simons gradient flow} on the product $G$-bundle over $Y$. As discussed by Morgan, Mrowka, and Ruberman in \cite{MMR}, the {\L}ojasiewicz-Simon gradient inequality for the Chern-Simons functional plays an essential role in analyzing the asymptotic behavior of a finite-energy anti-self-dual connection, $A(\cdot,t)$, as $t\to \infty$, where $(y,t)$ denote coordinates on $Y\times [0,\infty)$.

\subsection{Donaldson heat flow}
\label{subsec:Donaldson_heat_flow}
For expositions of the \emph{Donaldson heat flow} of Hermitian metrics\footnote{This flow is only equivalent to the Yang-Mills gradient flow of connections of type $(1,1)$ on $E$ when the Hermitian metric on $X$ is K\"ahler.}
on a complex vector bundle, $E$, over a complex, Hermitian manifold, $X$, we refer the reader to Donaldson's article \cite{DonASD} and the discussions by Jacob \cite{Jacob_2015conm} and Jost \cite{Jost_1991}.
Donaldson does not use the {\L}ojasiewicz-Simon gradient inequality to obtain his results, although there appears to be no impediment to such an analysis, but rather cleverly exploits special features of the flow that are particular to K\"ahler surfaces.

\subsection{Fluid dynamics}
\label{subsec:Fluid dynamics}
For applications of the {\L}ojasiewicz-Simon gradient inequality to fluid dynamics, see the articles by Feireisl, Lauren{\c{c}}ot, and Petzeltov{\'a} \cite{Feireisl_Laurencot_Petzeltova_2007}, Frigeri, Grasselli, and Krej{\v{c}}{\'{\i}} \cite{Frigeri_Grasselli_Krejcic_2013}, Grasselli and Wu \cite{Grasselli_Wu_2013}, and Wu and Xu \cite{Wu_Xu_2013}.

\subsection{Knot energy gradient flows}
\label{subsec:Knot_energy_gradient_flow}
O’Hara \cite{OHara_1991, OHara_1992, OHara_1994} invented a family of knot energies (compare also the knot energy functional proposed by Buck and Orloff \cite{Buck_Orloff_1995}). Blatt and his collaborators consider the gradient flow of the O'Hara and related energy functionals in \cite{Blatt_2012jntr, Blatt_2012cvpde, Blatt_2016arxiv}. In \cite{Blatt_2016arxiv}, Blatt studies the gradient flow of the sum of one of the O'Hara energies and a positive multiple of the length. He shows that the gradients of these knot energies can be written as the normal part of a quasilinear operator, derived short-time existence results for these flows, and applies the {\L}ojasiewicz-Simon gradient inequality to prove long-time existence and convergence to critical points.

\subsection{Lagrangian mean curvature flow}
\label{subsec:Lagrangian_mean_curvature_flow}
For a survey of results on global existence and convergence results and characterizations of first-time singularities in \emph{Lagrangian mean curvature flow}, we refer to the articles by Neves \cite{Neves_2013}, Wang \cite{Wang_2008sdg}, and references cited therein for further details and applications of the {\L}ojasiewicz-Simon gradient inequality.

\subsection{Mean curvature flow}
\label{subsec:Mean_curvature_flow}
For a discussion of applications of a version of the {\L}ojasiewicz-Simon gradient inequality to mean curvature flow, we refer to the articles by Colding and Minicozzi \cite{Colding_Minicozzi_2014sdg}, Schulze \cite{Schulze_2014}, and references cited therein for further details and applications of the {\L}ojasiewicz-Simon gradient inequality.

For interpretations of mean curvature flow as a gradient system, we refer the reader to Bellettini \cite{Bellettini_lecture_notes_on_mean_curvature_flow}, Ilmanen \cite{Ilmanen_1994}, Mantegazza \cite{Mantegazza_2011_lectures_mcf}, Smoczyk \cite{Smoczyk_2012}, and Zaal \cite{Zaal_2015}
and references contained therein. For applications of the DeTurck trick \cite{DeTurck_1983} to convert mean curvature flow to a nonlinear parabolic partial differential equation and establish short-time existence, we refer to Andrews and Baker \cite{Andrews_Baker_2010}, Baker \cite{BakerThesis}, Leng, Zhao, and Zhao \cite{Leng_Zhao_Zhao_2014}, and references contained therein. As in the case of Ricci flow, the interpretation of mean curvature flow as a gradient system leads to the introduction a time-varying family of Hilbert spaces --- a family of $L^2$ spaces defined by a measure that depends on the time-varying family of immersions \cite[Section 1.2, page 7]{Mantegazza_2011_lectures_mcf}.

\subsection{Seiberg-Witten gradient flow}
The \emph{Seiberg-Witten gradient flow} is the gradient flow for the Seiberg-Witten action functional, and thus a type of coupled Yang-Mills gradient flow --- see Feehan and Maridakis \cite{Feehan_Maridakis_Lojasiewicz-Simon_coupled_Yang-Mills} and Hong and Schabrun \cite{Hong_Schabrun_2010} and references cited therein for further details and applications of the {\L}ojasiewicz-Simon gradient inequality.

\subsection{Yang-Mills gradient flow over cylindrical-end manifolds}
\label{subsec:Yang-Mills_gradient_flow_cylindrical-end_manifolds}
Duncan \cite{Duncan_2015preprint} has studied the over four-dimensional manifolds with finitely many cylindrical ends, paving the way for extensions of the ideas in this monograph to complete, non-compact manifolds. While Duncan does not use the {\L}ojasiewicz-Simon gradient inequality in \cite{Duncan_2015preprint}, there does not appear to be any impediment to its application to extend the results of \cite{Duncan_2015preprint} to higher energies, given the now well-established analytical framework for Yang-Mills gauge theory on manifolds with cylindrical ends --- see Morgan, Mrowka, and Ruberman \cite{MMR}, Mrowka \cite{MrowkaThesis}, and Taubes \cite{TauCasson, TauL2}, based on the framework established by Lockhart and McOwen \cite{Lockhart_McOwen_1985}.

\subsection{Yang-Mills-Higgs gradient flow over K\"ahler surfaces}
\label{subsec:Yang-Mills-Higgs_gradient_flow_Kaehler_surfaces}
The \emph{Yang-Mills-Higgs gradient flow} over K\"ahler surfaces has been explored by Li and Zhang \cite{Li_Zhang_2011} and Wang and Zhang \cite{Wang_Zhang_2008}, extending earlier similar results for Yang-Mills gradient flow over K\"ahler surfaces due to Daskalopoulos and Wentworth \cite{Daskalopoulos_Wentworth_2004, Daskalopoulos_Wentworth_2007}. In \cite{Wilkin_2008}, Wilkin has extended the results of Daskalopoulos \cite{Daskalopoulos_1992} and R\r{a}de \cite{Rade_1992} for Yang-Mills gradient flow to the case of Yang-Mills-Higgs gradient flow on a Hermitian vector bundle over a closed Riemann surface; those results are extended to the case of a principal $G$-bundle over a closed Riemann surface by Biswas and Wilkin in \cite{Biswas_Wilkin_2010}. The Yang-Mills-Higgs gradient flow is a type of coupled Yang-Mills gradient flow --- see Feehan and Maridakis \cite{Feehan_Maridakis_Lojasiewicz-Simon_coupled_Yang-Mills} and references cited therein for further details for further details and applications of the {\L}ojasiewicz-Simon gradient inequality.

\subsection{Additional gradient flows in mathematical physics and applied mathematics: Cahn-Hilliard, Ginzburg-Landau, Kirchoff-Carrier and related energy functionals}
\label{subsec:Ginzburg-Landau_energy_gradient_flow}
For applications of the {\L}ojasiewicz-Simon gradient inequality to proofs of global existence, convergence, convergence rate, and stability of nonlinear evolution equations arising in other areas of mathematical physics (including the Cahn-Hilliard, Ginzburg-Landau, Kirchoff-Carrier, porous medium, reaction-diffusion, and semi-linear heat and wave equations), we refer to the monograph by Huang \cite{Huang_2006} for a comprehensive introduction and to the articles by Chill \cite{Chill_2003, Chill_2006}, Chill and Fiorenza \cite{Chill_Fiorenza_2006}, Chill, Haraux, and Jendoubi \cite{Chill_Haraux_Jendoubi_2009}, Chill and Jendoubi \cite{Chill_Jendoubi_2003, Chill_Jendoubi_2007}, Feireisl and Simondon \cite{Feireisl_Simondon_2000}, Feireisl and Tak{\'a}{\v{c}} \cite{Feireisl_Takac_2001}, Grasselli, Wu, and Zheng \cite{Grasselli_Wu_Zheng_2009}, Haraux \cite{Haraux_2012}, Haraux and Jendoubi \cite{Haraux_Jendoubi_1998, Haraux_Jendoubi_2007, Haraux_Jendoubi_2011}, Haraux, Jendoubi, and Kavian \cite{Haraux_Jendoubi_Kavian_2003}, Huang and Tak{\'a}{\v{c}} \cite{Huang_Takac_2001}, Jendoubi \cite{Jendoubi_1998jfa}, Rybka and Hoffmann \cite{Rybka_Hoffmann_1998, Rybka_Hoffmann_1999}, Simon \cite{Simon_1983}, and Tak{\'a}{\v{c}} \cite{Takac_2000}.

\chapter{Preliminaries}
\label{chapter:Preliminaries}

\section{Classical and weak solutions to the Yang-Mills gradient flow equation}
\label{sec:Yang-Mills_gradient_flow_solution_concepts}
We assume throughout this section that $G$ be a compact Lie group and $P$ a principal $G$-bundle over a closed, compact, connected, Riemannian, smooth dimensional manifold, $X$, of dimension $d\geq 2$. We begin with the

\begin{defn}[Classical solution to the Cauchy problem for Yang-Mills gradient flow]
\label{defn:Classical_solution_Yang-Mills_gradient_flow}
Let $G$ be a compact Lie group and $P$ be a principal $G$-bundle over a closed, Riemannian, smooth manifold, $X$ of dimension $d\geq 2$ and $T \in (0,\infty]$. If $A_0$ is a connection on $P$ of class $C^\infty$, we call $A(t)$, for $t \in [0,T)$, a \emph{classical} solution, with \emph{smoothly attained initial data}, to the Cauchy problem for Yang-Mills gradient flow if
$$
A - A_0 \in C^\infty([0,T)\times X;\Lambda^1\otimes\ad P),
$$
and
\begin{align}
\label{eq:Struwe_3_classical}
\frac{\partial A}{\partial t} + d_A^*F_A &= 0 \quad\hbox{on } (0,T)\times X,
\\
\label{eq:Struwe_4_classical}
A(0) &= A_0.
\end{align}
If $A_0$ is a connection on $P$ of class $W^{s,p}$, for $s\geq 1$ and $p\geq 2$, we call $A(t)$, for $t \in [0,T)$, a \emph{classical} solution, with \emph{continuously attained initial data}, to the Cauchy problem for Yang-Mills gradient flow if
$$
A - A_0
\in C([0,T); W_{A_1}^{s,p}(X;\Lambda^1\otimes\ad P))
\cap C^\infty((0,T)\times X;\Lambda^1\otimes\ad P),
$$
for some fixed $C^\infty$ reference connection $A_1$ on $P$
and \eqref{eq:Struwe_3_classical} and \eqref{eq:Struwe_4_classical} hold.
\end{defn}

By a wide margin, when $A_0$ is $C^\infty$, Definition \ref{defn:Classical_solution_Yang-Mills_gradient_flow} provides the simplest solution concept for the Cauchy problem for Yang-Mills gradient flow.

There are several variants of the concept of `weak' solution for an evolution equation as well as different terminologies and that can create some confusion. For Yang-Mills gradient flow, we provide four equivalent formulations in Definition \ref{defn:Weak_solution_Yang-Mills_gradient_flow} based on the corresponding definitions for a first-order linear evolution equation in a Banach space $V$ \cite[Proposition III.2.1]{Showalter}, \cite[Section 30.1]{Zeidler_nfaa_v3}, where $H$ is a Hilbert space with continuous and dense embedding, $V \hookrightarrow H$, so $V \hookrightarrow H \hookrightarrow V'$ forms an evolution triple in the sense of \cite[Section 30.1]{Zeidler_nfaa_v3} and $V'$ denotes the dual of $V$.

Given a $C^\infty$ connection, $A_1$, on $P$, we define $H_{A_1}^1(X;\Lambda^1\otimes\ad P) = W_{A_1}^{1,2}(X;\Lambda^1\otimes\ad P)$ in the usual way via the covariant derivative associated with $A_1$. In order that the energy, $\sE(A(t)) = \frac{1}{2}\|F_A(t)\|_{L^2(X)}^2$, of a solution to Yang-Mills gradient flow be well-defined for a.e. $t \in [0,T)$, it is natural to choose
$$
V := H_{A_1}^1(X;\Lambda^1\otimes\ad P) \quad\hbox{and}\quad H := L^2(X;\Lambda^1\otimes\ad P),
$$
where $A_1$ is a fixed $C^\infty$ reference connection on $P$, and thus
$$
V' = H_{A_1}^{-1}(X;\Lambda^1\otimes\ad P).
$$
Concepts of solution to Yang-Mills gradient flow with regularities intermediate between that of `weak' and classical solution will be defined as needed throughout our monograph and we shall work with solutions of greatest regularity whenever possible. The concept of weak solution, however, provide a useful framework in which to analyze uniqueness of solutions to Yang-Mills gradient flow. Prior to introducing this concept, we review some useful preliminaries.

By analogy with the definition of $W_2(0,T;V,H)$ in \cite[p. 105]{Showalter}, we define the Hilbert space,
\begin{multline}
\label{eq:Showalter_W_2_space}
W_{A_1;2}(0,T; \Lambda^1\otimes \ad P)
\\
:= \{a \in L^2(0,T; H_{A_1}^1(X;\Lambda^1\otimes\ad P)): \dot a \in L^2(0,T; H_{A_1}^{-1}(X;\Lambda^1\otimes\ad P))\},
\end{multline}
with norm
$$
\|a\|_{L^2(0,T;H_{A_1}^1(X))} + \|\dot a\|_{L^2(0,T;H_{A_1}^{-1}(X))}.
$$
We have the following convenient analogues of \cite[Proposition III.1.2 and Corollary III.1.1]{Showalter}.

\begin{prop}[Continuous embedding a Hilbert space into a Banach space of continuous maps]
\label{prop:Showalter_III-1-2}
There is a continuous embedding,
$$
W_{A_1;2}(0,T; \Lambda^1\otimes \ad P) \hookrightarrow C([0,T]:L^2(X;\Lambda^1\otimes\ad P)),
$$
and, for $a \in W_{A_1;2}(0,T; \Lambda^1\otimes \ad P)$, the function $\|a(\cdot)\|_{L^2(X)}^2$ is absolutely continuous on $[0,T]$, with
$$
\frac{d}{dt}\|a(t)\|_{L^2(X)}^2 = (\dot a(t), a(t))_{L^2(X)},
\quad\hbox{a.e. } t \in (0, T).
$$
\end{prop}

\begin{cor}[Absolute continuity of pairings of time-dependent maps]
\label{cor:Showalter_III-1-1}
If $a, b \in W_{A_1;2}(0,T; \Lambda^1\otimes \ad P)$, then the function $(a(\cdot), b(\cdot))_{L^2(X)}$ is absolutely continuous on $[0,T]$, with
$$
\frac{d}{dt}(a(t), b(t))_{L^2(X)} = (\dot a(t), b(t))_{L^2(X)} + (a(t), \dot b(t))_{L^2(X)}, \quad \quad\hbox{a.e. } t \in (0, T).
$$
\end{cor}

We can now provide the

\begin{defn}[Weak solution to the Cauchy problem for Yang-Mills gradient flow]
\label{defn:Weak_solution_Yang-Mills_gradient_flow}
Let $G$ be a compact Lie group, $P$ be a principal $G$-bundle over a closed, Riemannian, smooth manifold, $X$, of dimension $d\geq 2$ and $T \in (0,\infty]$ and $A_0$ a connection on $P$ of class $H^1$ and $A_1$ a fixed $C^\infty$ reference connection on $P$. Denote $X_T = (0,T)\times X$ and $a_0 := A_0-A_1 \in H_{A_1}^1(X;\Lambda^1\otimes\ad P)$. We call a family of connections, $A(t) = A_1 + a(t)$ for $t \in [0, T)$, a \emph{weak} solution to the Cauchy problem for Yang-Mills gradient flow if it satisfies one of the following formulations:
\begin{enumerate}
\item\label{item:Weak_type_1_solution_Yang-Mills_gradient_flow}
\emph{Weak type 1}: One has $a \in W_{A_1;2}(0,T; \Lambda^1\otimes\ad P)$ with $a(0) = a_0$ and
\begin{equation}
\label{eq:YMGF_Zeidler_30-1}
\dot a(t) + d_A^*F_A(t) = 0 \quad\hbox{in } H_{A_1}^{-1}(X;\Lambda^1\otimes\ad P),
\quad\hbox{a.e. } t \in (0,T).
\end{equation}
The identity in \eqref{eq:YMGF_Zeidler_30-1} means
\begin{equation}
\label{eq:YMGF_Zeidler_30-2}
(\dot a(t),b)_{L^2(X)} + (d_A^*F_A(t),b)_{L^2(X)} = 0,
\quad\hbox{a.e. } t \in (0,T), \quad\forall\, b \in H_{A_1}^1(X;\Lambda^1\otimes\ad P).
\end{equation}
This formulation follows \cite[Equations (30.1) and (30.2)]{Zeidler_nfaa_v3}.

\item\label{item:Weak_type_2_solution_Yang-Mills_gradient_flow}
\emph{Weak type 2}: One has $a \in L^2(0,T; H_{A_1}^1(X;\Lambda^1\otimes\ad P))$ and $\dot a \in (L^2(0,T; H_{A_1}^1(X;\Lambda^1\otimes\ad P)))'$ with $a(0) = a_0$ and
\begin{equation}
\label{eq:YMGF_Zeidler_30-7a}
\dot a + d_A^*F_A = 0 \quad\hbox{in } (L^2(0,T; H_{A_1}^1(X;\Lambda^1\otimes\ad P)))'.
\end{equation}
The identity in \eqref{eq:YMGF_Zeidler_30-7a} means
\begin{equation}
\label{eq:YMGF_Showalter_proposition_III-2-1-a}
(\dot a, b)_{L^2(X_T)} + (d_A^*F_A, b)_{L^2(X_T)} = 0,
\quad\forall\, b \in L^2(0,T; H_{A_1}^1(X;\Lambda^1\otimes\ad P)).
\end{equation}
This formulation follows \cite[Proposition III.2.1 (a)]{Showalter},  \cite[Equation (30.7)]{Zeidler_nfaa_v3}.

\item\label{item:Weak_type_3_solution_Yang-Mills_gradient_flow}
\emph{Weak type 3}: One has $a \in W_{A_1;2}(0,T; \Lambda^1\otimes\ad P)$ and
\begin{multline}
\label{eq:YMGF_Showalter_proposition_III-2-1-b}
-(a, \dot b)_{L^2(X_T)} + (d_A^*F_A, b)_{L^2(X_T)} = (a_0,b(0))_{L^2(X)},
\\
\forall\, b \in W_{A_1;2}(0,T; \Lambda^1\otimes\ad P) \hbox{ with } b(T) = 0.
\end{multline}
This formulation follows \cite[Proposition III.2.1 (b)]{Showalter}.

\item\label{item:Weak_type_4_solution_Yang-Mills_gradient_flow}
\emph{Weak type 4}:
One has $a \in W_{A_1;2}(0,T; \Lambda^1\otimes\ad P)$ with $a(0) = a_0$ and
\begin{equation}
\label{eq:YMGF_Showalter_proposition_III-2-1-c_distributions}
(\dot a, b)_{L^2(X)} + (d_A^*F_A, b)_{L^2(X)} = 0 \quad\hbox{in }\sD^*(0,T),
\quad\forall\, b\in H_{A_1}^1(X;\Lambda^1\otimes\ad P).
\end{equation}
The identity in \eqref{eq:YMGF_Showalter_proposition_III-2-1-c_distributions} means
\begin{multline}
\label{eq:YMGF_Showalter_proposition_III-2-1-c}
-\int_0^T (a(t), b)_{L^2(X)}\dot\varphi(t)\,dt + \int_0^T (d_A^*F_A(t), b)_{L^2(X)}\varphi(t)\,dt = 0,
\\
\quad\forall\, b\in H_{A_1}^1(X;\Lambda^1\otimes\ad P),\ \varphi \in C_0^\infty(0,T).
\end{multline}
This formulation follows \cite[Proposition III.2.1 (c)]{Showalter}, \cite[Equation (30.4)]{Zeidler_nfaa_v3}.
\end{enumerate}
\end{defn}

\begin{rmk}[On the interpretation of $d_A^*F_A$ in Definition \ref{defn:Weak_solution_Yang-Mills_gradient_flow}]
We observe that, for any $b \in L^2(0,T; H_{A_1}^1(X;\Lambda^1\otimes\ad P))$,
\begin{align*}
(d_A^*F_A, b)_{L^2(X_T)}
&= \int_0^T (d_A^*F_A(t), b(t))_{L^2(X)}\,dt
\\
&= \int_0^T (F_A(t), d_A^*b(t))_{L^2(X)}\,dt
\\
&= (F_A, d_Ab)_{L^2(X_T)}
\\
&= (F_{A_1} + d_{A_1}a + [a,a], d_Ab)_{L^2(X_T)},
\end{align*}
and so this inner product is well-defined in Definition \ref{defn:Weak_solution_Yang-Mills_gradient_flow} by the (implicit) interpretation of $d_A^*F_A$ as an element of $L^2(0,T; H_{A_1}^{-1}(X;\Lambda^1\otimes\ad P))$. A similar comment applies to $(d_A^*F_A(t), b)_{L^2(X)}$ when $b \in H_{A_1}^1(X;\Lambda^1\otimes\ad P)$.
\end{rmk}

The formulation \eqref{item:Weak_type_1_solution_Yang-Mills_gradient_flow} in Definition \ref{defn:Weak_solution_Yang-Mills_gradient_flow} agrees with the \cite[Definitions 1.1 and 1.3 (ii)]{Huang_2006} due to Huang for a weak solution to the Cauchy problem for a gradient system. The formulation \eqref{item:Weak_type_3_solution_Yang-Mills_gradient_flow} in Definition \ref{defn:Weak_solution_Yang-Mills_gradient_flow} agrees with the \cite[Definition, p. 94 (ii)]{Kozono_Maeda_Naito_1995} due to Kozono, Maeda, and Naito for a weak solution to the Cauchy problem for Yang-Mills gradient flow.

Unfortunately, it is unclear how to prove existence, even for sufficiently short time intervals, of weak solutions to the Cauchy problem for Yang-Mills gradient flow in the sense of Definition \ref{defn:Weak_solution_Yang-Mills_gradient_flow}. As the discussion in Struwe \cite[Section 4.4]{Struwe_1996} illustrates, the difficulty is due to a combination of the minimal regularity assumption on the initial data and the fact that the Yang-Mills gradient flow equation is not parabolic, thus lacking the smoothing behavior of solutions to homogeneous parabolic equations. A similar point is made by R\r{a}de in \cite[p. 125]{Rade_1992}. We shall return to this fundamental issue in our discussion of the `Donaldson-DeTurck trick' in Section \ref{subsec:Struwe_4-4_Donaldson-DeTurck trick}.

The solution whose existence is established by Kozono, Maeda, and Naito in their \cite[Theorem A]{Kozono_Maeda_Naito_1995} is actually a strong solution to their \cite[Equation (1.3)]{Kozono_Maeda_Naito_1995}, the parabolic initial value problem obtained by applying the Donaldson-DeTurck trick, and not a solution (weak or strong) to their \cite[Equations (1.1) and (1.2)]{Kozono_Maeda_Naito_1995}, the Cauchy problem for Yang-Mills gradient flow. When the initial data, $A_0$, is of class $W^{s+1,p}$ with $s\geq 1$ and $p\in[2,\infty]$ obeying $sp>d$ (for a base manifold of dimension $d\geq 2$), one can readily show (as we review in Section \ref{subsec:Struwe_4-4_Donaldson-DeTurck trick}) that the Donaldson-DeTurck trick yields a `strong' and thus also a weak solution in the sense of Definition \ref{defn:Weak_solution_Yang-Mills_gradient_flow}. However, that argument encounters fundamental difficulties when $sp\leq d$.

For these reasons, Struwe employs a different concept of weak solution. Our formulation in Definition \ref{defn:Struwe_weak_solution} below is stronger than that in his \cite[Definition 2.1]{Struwe_1994}, but agrees with the formulation in his \cite[Theorem 2.3 (i) and Section 4.4]{Struwe_1994}.

\begin{defn}[Struwe weak solution to the Cauchy problem for Yang-Mills gradient flow]
\label{defn:Struwe_weak_solution}
\cite[Theorem 2.3 (i)]{Struwe_1994}
Let $G$ be a compact Lie group, $P$ be a principal $G$-bundle over a closed, Riemannian, smooth manifold, $X$, of dimension $d\geq 2$ and $T \in (0,\infty]$ and $A_0$ a connection on $P$ of class $H^1$ and $A_1$ a fixed $C^\infty$ reference connection on $P$. Denote $a_0 := A_0-A_1 \in H_{A_1}^1(X;\Lambda^1\otimes\ad P)$. We call $A(t) = A_1+a(t)$, for $t\in [0,T)$, a \emph{Struwe weak} solution to the Cauchy problem for Yang-Mills gradient flow if $a(0) = a_0$ and
\begin{subequations}
\label{eq:Struwe_Definition_2-1}
\begin{gather}
\label{eq:Struwe_Definition_2-1_A-A1_regularity}
a \in C([0,T); L^2(X;\Lambda^1\otimes\ad P))
\cap H^1(0,T; L^2(X;\Lambda^1\otimes\ad P)),
\\
\label{eq:Struwe_Definition_2-1_F_A_regularity}
F_A \in C([0,T); L^2(X;\Lambda^1\otimes\ad P)),
\\
\label{eq:Struwe_Definition_2-1_variational_equation}
-(a, \dot b)_{L^2(X_T)} + (F_A, d_Ab)_{L^2(X_T)} =  0,
\quad
\forall\, b \in C_0^\infty((0,T)\times X; \Lambda^1\otimes\ad P).
\end{gather}
\end{subequations}
\end{defn}

As comparison with \cite[Definition 2.1]{Struwe_1994} shows, other variants of Definition \ref{defn:Struwe_weak_solution} are possible. It is important to note that Definition \ref{defn:Struwe_weak_solution} does \emph{not} require
$$
a \in L^2(0,T; H_{A_1}^1(X;\Lambda^1\otimes\ad P)),
$$
and this is the key difference with Definition \ref{defn:Weak_solution_Yang-Mills_gradient_flow}.
Ambiguity in the concept of `weak' solution, at least in the sense of Definition \ref{defn:Weak_solution_Yang-Mills_gradient_flow}, is removed by the

\begin{prop}[Equivalence of formulations of weak solution]
\label{prop:Showalter_III-2-1}
The formulations in Definition \ref{defn:Weak_solution_Yang-Mills_gradient_flow} of weak solution to the Cauchy problem for Yang-Mills gradient flow are equivalent.
\end{prop}

\begin{proof}
The corresponding equivalence for any first-order \emph{linear} evolution equation, defined with the aid of an evolution triple $V \hookrightarrow H \hookrightarrow V'$, is given in \cite[Proposition III.2.1]{Showalter}, \cite[Section 30.1]{Zeidler_nfaa_v3}. For the sake of completeness, we provide the proof of the desired equivalence for Yang-Mills gradient flow.

Proposition \ref{prop:Showalter_III-1-2} implies that any $a \in W_{A_1;2}(0,T ;\Lambda^1\otimes\ad P))$ satisfies
\begin{equation}
\label{eq:A_C0_path_in_L2}
a \in C([0,T]; L^2(X;\Lambda^1\otimes\ad P)),
\end{equation}
and so an initial condition $a(0) = a_0$ is well-defined (as an identity in $L^2(X;\Lambda^1\otimes\ad P)$). Moreover, Corollary \ref{cor:Showalter_III-1-1} implies, given $a, b \in W_{A_1;2}(0,T ;\Lambda^1\otimes\ad P))$, that the function,
\begin{equation}
\label{eq:AC_TEMPLABEL}
[0, T] \ni t \mapsto (a(t), b(t))_{L^2(X)} \in \RR,
\end{equation}
is absolutely continuous.

\emph{Weak type \ref{item:Weak_type_2_solution_Yang-Mills_gradient_flow} $\implies$ \ref{item:Weak_type_3_solution_Yang-Mills_gradient_flow}.}
Suppose that $A=A_1+a$ is a weak solution in the sense of Item \eqref{item:Weak_type_2_solution_Yang-Mills_gradient_flow} in Definition \ref{defn:Weak_solution_Yang-Mills_gradient_flow} and $b \in W_{A_1;2}(0,T; \Lambda^1\otimes\ad P)$ with $b(T)=0$. Since
$$
\dot a \in (L^2(0,T; H_{A_1}^1(X;\Lambda^1\otimes\ad P)))'
$$
and
$$
(L^2(0,T; H_{A_1}^1(X;\Lambda^1\otimes\ad P)))'
\cong
L^2(0,T; H_{A_1}^{-1}(X;\Lambda^1\otimes\ad P))
\quad\hbox{(by \cite[Theorem III.1.5]{Showalter})},
$$
and $a\in L^2(0,T; H_{A_1}^1(X;\Lambda^1\otimes\ad P))$, then the regularity hypothesis on $a(t)$ is equivalent to $a \in W_{A_1;2}(0,T ;\Lambda^1\otimes\ad P))$ by \eqref{eq:Showalter_W_2_space}. By the absolute continuity of \eqref{eq:AC_TEMPLABEL}, we have
$$
\frac{d}{dt}(a(t), b(t))_{L^2(X)} = (\dot a(t), b(t))_{L^2(X)} + (a(t), \dot b(t))_{L^2(X)}, \quad\hbox{a.e. } t \in (0, T),
$$
and, integrating over $[0,T]$,
\begin{align*}
{}&(a(T), b(T))_{L^2(X)} - (a(0), b(0))_{L^2(X)}
\\
&\quad = \int_0^T \frac{d}{dt}(a(t), b(t))_{L^2(X)}\,dt
\\
&\quad = \int_0^T (\dot a(t), b(t))_{L^2(X)}\,dt + \int_0^T (a(t), \dot b(t))_{L^2(X)}\,dt \quad\hbox{(by Corollary \ref{cor:Showalter_III-1-1}),}
\end{align*}
and thus, substituting \eqref{eq:YMGF_Showalter_proposition_III-2-1-a}, namely,
$$
(\dot a(t),b(t))_{L^2(X)} + (F_A(t),d_{A(t)}b(t))_{L^2(X)} = 0,
\quad\hbox{a.e. } t \in (0,T),
$$
into the first integral in the last equality, we obtain
\begin{equation}
\label{eq:Integral_0T_Yang-Mills_gradient_flow_strong_variational_solution_equation}
\begin{aligned}
{}&(a(T), b(T))_{L^2(X)} - (a(0), b(0))_{L^2(X)}
\\
&= -\int_0^T (F_{A(t)}, d_{A(t)}b(t))_{L^2(X)}\,dt + \int_0^T (a(t), \dot b(t))_{L^2(X)}\,dt \quad\hbox{(by \eqref{eq:YMGF_Showalter_proposition_III-2-1-a})},
\end{aligned}
\end{equation}
Consequently, $A=A_1+a$ obeys \eqref{eq:YMGF_Showalter_proposition_III-2-1-b} since $a(0) = a_0$ and $b(T)=0$. Hence, $A$ is a weak solution in the sense of Item \eqref{item:Weak_type_3_solution_Yang-Mills_gradient_flow} in Definition \ref{defn:Weak_solution_Yang-Mills_gradient_flow}.

\emph{Weak type \ref{item:Weak_type_3_solution_Yang-Mills_gradient_flow} $\implies$ \ref{item:Weak_type_4_solution_Yang-Mills_gradient_flow}.}
Assume $A=A_1+a$ is a weak solution in the sense of Item \eqref{item:Weak_type_3_solution_Yang-Mills_gradient_flow} in Definition \ref{defn:Weak_solution_Yang-Mills_gradient_flow}. Substituting $b(t) = b\varphi(t)$ into \eqref{eq:YMGF_Showalter_proposition_III-2-1-b}, with $b \in H_{A_1}^1(X;\Lambda^1\otimes\ad P)$ and $\varphi \in C^\infty_0(0,T;\RR)$, yields
$$
-\int_0^T (a(t), b)_{L^2(X)}\dot\varphi(t)\,dt + \int_0^T (F_{A(t)}, d_{A(t)}b)_{L^2(X)}\varphi(t)\,dt = 0,
$$
for all such $b$ and $\varphi$. Thus, \eqref{eq:YMGF_Showalter_proposition_III-2-1-c} holds. On the other hand, substituting $b(t) = b\varphi(t)$ into \eqref{eq:YMGF_Showalter_proposition_III-2-1-b}, with $b \in H_{A_1}^1(X;\Lambda^1\otimes\ad P)$ and $\varphi \in C^\infty_0([0,T;\RR)$, yields
$$
-\int_0^T (a(t), b)_{L^2(X)}\dot\varphi(t)\,dt + \int_0^T (F_{A(t)}, d_{A(t)}b)_{L^2(X)}\varphi(t)\,dt = (a_0, b)_{L^2(X)},
$$
and, integrating by parts with respect to time,
$$
\int_0^T (\dot a(t), b)_{L^2(X)}\varphi(t)\,dt + \int_0^T (F_{A(t)}, d_{A(t)}b)_{L^2(X)}\varphi(t)\,dt + (a(0), b)_{L^2(X)} = (a_0, b)_{L^2(X)}.
$$
Consequently, $(a(0), b)_{L^2(X)} = (a_0, b)_{L^2(X)}$ for all $b \in H_{A_1}^1(X;\Lambda^1\otimes\ad P)$ and so $a(0)=a_0$. Thus, $A=A_1+a$ is a weak solution in the sense of Item \eqref{item:Weak_type_4_solution_Yang-Mills_gradient_flow} in Definition \ref{defn:Weak_solution_Yang-Mills_gradient_flow}.

\emph{Weak type \ref{item:Weak_type_4_solution_Yang-Mills_gradient_flow} $\implies$ \ref{item:Weak_type_1_solution_Yang-Mills_gradient_flow}.}
Assume $A=A_1+a$ is a weak solution in the sense of Item \eqref{item:Weak_type_4_solution_Yang-Mills_gradient_flow} in Definition \ref{defn:Weak_solution_Yang-Mills_gradient_flow}. Integration by parts with respect to time in \eqref{eq:YMGF_Showalter_proposition_III-2-1-c} yields
\begin{multline*}
\int_0^T (\dot a(t), b)_{L^2(X)}\varphi(t)\,dt + \int_0^T (F_{A(t)}, d_{A(t)}b)_{L^2(X)}\varphi(t)\,dt = 0,
\\
\quad\forall\, b \in H_{A_1}^1(X;\Lambda^1\otimes\ad P),\ \varphi \in C_0^\infty(0,T),
\end{multline*}
and so
$$
(\dot a(t), b)_{L^2(X)} + (F_{A(t)}, d_{A(t)}b)_{L^2(X)} = 0, \quad\hbox{a.e. } t \in (0, T),
\quad\forall, b \in H_{A_1}^1(X;\Lambda^1\otimes\ad P),
$$
that is, $A$ obeys \eqref{eq:YMGF_Zeidler_30-2}. Hence, $A=A_1+a$ is a weak solution in the sense of Item \eqref{item:Weak_type_1_solution_Yang-Mills_gradient_flow} in Definition \ref{defn:Weak_solution_Yang-Mills_gradient_flow}.

\emph{Weak type \ref{item:Weak_type_1_solution_Yang-Mills_gradient_flow} $\implies$ \ref{item:Weak_type_2_solution_Yang-Mills_gradient_flow}.} Assume that $A=A_1+a$ is a weak solution in the sense of Item \eqref{item:Weak_type_1_solution_Yang-Mills_gradient_flow} in Definition \ref{defn:Weak_solution_Yang-Mills_gradient_flow}. The regularity hypothesis on $a(t)$ is equivalent to $a \in L^2(0,T;H_{A_1}^1(X;\Lambda^1\otimes\ad P))$ and $\dot a \in (L^2(0,T;H_{A_1}^1(X;\Lambda^1\otimes\ad P)))'$. For any $b \in L^2(0,T;H_{A_1}^1(X;\Lambda^1\otimes\ad P))$ we have $b(t) \in H_{A_1}^1(X;\Lambda^1\otimes\ad P)$ for a.e. $t \in (0, T)$ and so integrating \eqref{eq:YMGF_Zeidler_30-2} over $[0,T]$ yields
\begin{multline*}
\int_0^T (\dot a(t), b(t))_{L^2(X_T)}\,dt + \int_0^T (d_A^*F_A(t), b(t))_{L^2(X_T)}\,dt = 0,
\\
\forall\, b \in L^2(0,T; H_{A_1}^1(X;\Lambda^1\otimes\ad P)),
\end{multline*}
that is, \eqref{eq:YMGF_Showalter_proposition_III-2-1-a} holds and $A=A_1+a$ is a weak solution in the sense of Item \eqref{item:Weak_type_2_solution_Yang-Mills_gradient_flow} in Definition \ref{defn:Weak_solution_Yang-Mills_gradient_flow}.
\end{proof}

\section{On the heat equation method}
\label{sec:Heat_equation_method}
According to R\r{a}de \cite[page 124]{Rade_1992}, the proof of the Hodge Theorem for de Rham cohomology by Milgram and Rosenbloom was one of the first applications of the heat equation method in geometric analysis \cite{Milgram_Rosenbloom_1951a, Milgram_Rosenbloom_1951a}. Conceptual expositions of the heat equation method have appeared elsewhere (see Atiyah and Bott \cite{Atiyah_Bott_1983} and Donaldson and Kronheimer \cite[Section 6.1.2]{DK}, but a few additional, informal comments may be worthwhile.

Suppose first that we wish to solve the nonhomogeneous Hodge Laplace equation on $\Omega^1(X;\ad P)$, namely,
$$
\Delta_A a = f,
$$
where $\Delta_A = d_A^*d_A + d_Ad_A^*$ and $a, f \in \Omega^p(X;\ad P)$ and $X$ is a closed, Riemannian, smooth manifold. There will in general be a finite-dimensional cokernel obstruction to solving such an equation and, because $\Delta_A$ is (formally) self-adjoint, that is equivalent to $\Delta_A$ having a non-zero kernel. Suppose we instead consider the corresponding nonhomogeneous heat equation, namely,
$$
\frac{\partial a}{\partial t} + \Delta_A a = f, \quad a(0) = a_0,
$$
for some choice of initial data, $a_0 \in \Omega^p(X;\ad P)$. Now choose a large enough positive constant, $\mu$, that ensures that the operator, $\Delta_A + \mu$, is positive and consider the related nonhomogeneous heat equation,
$$
\frac{\partial a}{\partial t} + (\Delta_A + \mu)a = f + \mu a, \quad a(0) = a_0.
$$
To avoid having the term $\mu a$ appear on the right-hand side, we may apply the familiar exponential-shift trick, set $b(t) = e^{-\mu t}a(t)$, and observe that $b_0 = b(0) = a(0) = a_0$ and
\begin{align*}
\frac{\partial b}{\partial t} &= e^{-\mu t}\frac{\partial a}{\partial t} - \mu e^{-\mu t}a
\\
&= e^{-\mu t}(-\Delta_A a + f) - \mu e^{-\mu t}a
\\
&= -(\Delta_A + \mu)b + g,
\end{align*}
where we set $g(t) = e^{-\mu t}f(t)$. Thus,
$$
\frac{\partial b}{\partial t} + (\Delta_A + \mu)b = g, \quad b(0) = b_0.
$$
Proceeding formally, we may write the solution as
$$
b(t) = e^{-(\Delta_A + \mu)t}b_0 + \int_0^t e^{-(\Delta_A + \mu)(t-s)}g(s)\,ds, \quad t \geq 0.
$$

\section[Classification of principal $G$-bundles and the Chern-Weil formula]{Classification of principal $G$-bundles, the Chern-Weil formula, and absolute minima of the Yang-Mills energy functional}
\label{sec:Taubes_1982_Appendix}

In this section, we extend the discussion in \cite[Sections 2.1.3 and 2.1.4]{DK} to the case of compact Lie groups and also establish certain conventions for calculations involving the curvature of a connection.

Let $G$ be a Lie group and $P$ be a principal $G$-bundle over a closed, connected, four-dimensional, smooth manifold, $X$. We shall need to appeal to the classification of principal $G$-bundles, $P$, for compact Lie groups, so we recall the facts we shall require from \cite[Appendix]{Sedlacek},  \cite[Appendix]{TauSelfDual}. If $G$ is a compact, connected Lie group and $\widetilde G$ is its universal covering group, then
$$
\widetilde G \cong \RR^m \times G_1\times\cdots \times G_l,
$$
where the $G_i$ are compact, simple, simply-connected Lie groups and consequently the real vector bundle,
\begin{equation}
\label{eq:Taubes_1982_page_139_adjoint_bundle}
\ad P := P \times_{\ad} \fg
\end{equation}
associated to $G$ by the adjoint representation, $\Ad:G \ni u \to \Ad_u \in \Aut\fg$, splits
$$
\ad P = \underline\RR^m\times \ad_i P\oplus \cdots \oplus \ad_l P,
$$
where $\ad_i P := P \times_{\ad}\fg_i$, for $i = 1\,\ldots, l$ and $\underline\RR^m := X \times \RR^m$ (see \cite[Appendix]{Sedlacek}). The simple Lie algebras, $\fg_i$, are the Lie algebras of the Lie groups, $G_i$. We shall henceforth further restrict our attention to \emph{semisimple} Lie groups, $G$, and so we may omit the factor, $\RR^m$, which arises from a torus in $G$ (and corresponding Abelian Lie algebra with trivial Killing form). The $\fg_i$ are nontrivial ideals comprising the semisimple Lie algebra, $\fg$, of $G$.

Assume now that $X$ is oriented. Given a connection, $A$, on $P$, Chern-Weil theory provides representatives for the first Pontrjagin classes of the vector bundles, $\ad P$ and $\ad_i P$, namely \cite[Equation (A.7)]{TauSelfDual}
\begin{equation}
\label{eq:Taubes_1982_A7_deRham_cohomology_classes}
\begin{aligned}
p_1(P) \equiv p_1(\ad P) &= -\frac{1}{4\pi^2}\tr_\fg(F_A\wedge F_A) \in H^4_{\deRham}(X),
\\
p_1^i(P) \equiv p_1(\ad_i P) &= -\frac{1}{4\pi^2}\tr_{\fg_i}(F_A\wedge F_A) \in H^4_{\deRham}(X), \quad 1 \leq i \leq l,
\end{aligned}
\end{equation}
and hence scalar and vector Pontrjagin numbers \cite[Equation (A.7)]{TauSelfDual} (compare \cite[page 121]{Kozono_Maeda_Naito_1995}),
\begin{equation}
\label{eq:Taubes_1982_A7_scalar_and_vector_Pontrjagin_numbers}
\begin{aligned}
p_1(P)[X] \equiv p_1(\ad P)[X] &= -\frac{1}{4\pi^2}\int_X\tr_\fg(F_A\wedge F_A),
\\
p_1^i(P)[X] \equiv p_1(\ad_i P)[X] &= -\frac{1}{4\pi^2}\int_X\tr_{\fg_i}(F_A\wedge F_A), \quad 1 \leq i \leq l.
\end{aligned}
\end{equation}
Principal $G$-bundles, $P$, are classified \cite[Appendix]{Sedlacek}, \cite[Propositions A.1 and A.2]{TauSelfDual} by a cohomology class $\eta(P) \in H^2(X;\pi_1(G))$ and the vector Pontrjagin degree \cite[Equation (A.6)]{TauSelfDual},
\begin{equation}
\label{eq:Taubes_1982_A6}
\bkappa(P) \equiv (\kappa^1(P),\ldots,\kappa^l(P)) := \left(\frac{1}{r_{\fg_1}} p_1^1(\ad P)[X], \ldots, \frac{1}{r_{\fg_l}} p_1^l(\ad P)[X]\right) \in \ZZ^l,
\end{equation}
where the positive integers, $r_{\fg_i}$, depend on the Lie groups, $G_i$ \cite[Equation (A.5)]{TauSelfDual}; for example, if $G = \SU(n)$, then $r_\fg = 4n$. The Pontrjagin degree of $P$ is defined by
\begin{equation}
\label{eq:Taubes_1982_A7_Pontrjagin_degree}
\kappa(P) := \sum_{i=1}^l \kappa^i(P) = -\frac{1}{4\pi^2}\sum_{i=1}^l \frac{1}{r_{\fg_i}}\int_X\tr_{\fg_i}(F_A\wedge F_A).
\end{equation}
For $G = \Or(n)$ or $\SO(n)$, then $\eta(P) = w_2(P) \in H^2(X; \ZZ/2\ZZ)$, where $w_2(P) \equiv w_2(V)$ and $V = P\times_{\Or(n)}\RR^n$ or $P\times_{\SO(n)}\RR^n$ is the real vector bundle associated to $P$ via the standard representation, $\Or(n) \hookrightarrow \GL(n;\RR)$ or $\SO(n) \hookrightarrow \GL(n;\RR)$; for $G = \U(n)$, then $\eta(P) = c_1(P) \in H^2(X; \ZZ)$, where $c_1(P) \equiv c_1(E)$ and $E = P\times_{\U(n)}\CC^n$ is the complex vector bundle associated to $P$ via the standard representation, $\U(n) \hookrightarrow \GL(n;\CC)$ \cite[Theorem 2.4]{Sedlacek}. The topological invariant, $\eta \in H^2(X; \pi_1(G))$, is the obstruction to the existence of a principal $G$-bundle, $P$ over $X$, with a specified vector Pontrjagin number.

Assume in addition that $X$ has a Riemannian metric. To relate the Chern-Weil formulae \eqref{eq:Taubes_1982_A7_scalar_and_vector_Pontrjagin_numbers} to the $L^2(X)$-norms of $F_A^\pm$, we need to recall some facts concerning the Killing form \cite{Knapp_2002}. Every element $\xi$ of a Lie algebra $\fg$ over a field $\KK$ defines an adjoint endomorphism, $\ad\,\xi \in \End_\KK\fg$, with the help of the Lie bracket via $(\ad\,\xi)(\zeta) := [\xi, \zeta]$, for all $\zeta \in \fg$. For a finite-dimensional Lie algebra, $\fg$, its Killing form is the symmetric bilinear form,
\begin{equation}
B(\xi, \zeta) := \tr(\ad\,\xi \circ \ad\,\zeta), \quad\forall\xi, \zeta \in \fg,
\end{equation}
with values in $\KK$. Since we restrict to compact Lie groups, their Lie algebras are real. The Lie algebra, $\fg$, is semisimple (according to the Cartan criterion) if and only if its Killing form is non-degenerate. If a Lie algebra, $\fg$, is a direct sum of its ideals, $\fg_1 \ldots, \fg_l$, then the Killing form of $\fg$ is the direct sum of the Killing forms of the individual summands, $\fg_i$. The Killing form of a semisimple Lie algebra is negative definite. For example, if $G = \SU(n)$, then $B(M, N) = 2n\tr(MN)$ for matrices $M, N \in \CC^{n\times n}$, while if $G = \SO(n)$, then $B(M, N) = (n-2)\tr(MN)$ for matrices $M, N \in \RR^{n\times n}$. In particular, if $B_{\fg_i}$ is the Killing form on $\fg_i$, then it defines an inner product on $\fg_i$ via $\langle \cdot, \cdot \rangle_{\fg_i} = -B_{\fg_i}(\cdot,\cdot)$ and thus a norm $|\cdot|_{\fg_i}$ on $\fg_i$. But
$$
F_A = F_A^+ \oplus F_A^- \in \Omega^2(X;\ad P) = \Omega^+(X;\ad P) \oplus \Omega^-(X;\ad P),
$$
corresponding to the positive and negative eigenspaces of $*:\Lambda^2 \to \Lambda^2$, so
$$
F_A^\pm = \frac{1}{2}(1 \pm *)F_A \in \Omega^\pm(X; \ad P),
$$
and
$$
F_A\wedge F_A = (F_A^+\oplus F_A^-)\wedge (F_A^+\oplus F_A^-) = F_A^+ \wedge *F_A^+ - F_A^- \wedge *F_A^-.
$$
Hence, for $i = 1, \ldots, l$,
$$
\tr_{\fg_i}(F_A\wedge F_A) = \tr_{\fg_i}(F_A^+ \wedge *F_A^+) - \tr_{\fg_i}(F_A^- \wedge *F_A^-) = \left(|F_A^+|_{\fg_i}^2 - |F_A^-|_{\fg_i}^2\right)\,d\vol.
$$
The components of the vector Pontrjagin degree may thus be computed by
$$
\kappa^i(P) = \frac{1}{4\pi^2 r_{\fg_i}} \int_X \left(|F_A^+|_{\fg_i}^2 - |F_A^-|_{\fg_i}^2\right)\,d\vol, \quad 1 \leq i \leq l.
$$
If $A$ is self-dual, then $F_A^- \equiv 0$ over $X$ and
$$
\kappa^i(P) = \frac{1}{4\pi^2 r_{\fg_i}} \int_X |F_A|_{\fg_i}^2\,d\vol \geq 0, \quad 1 \leq i \leq l,
$$
while if $A$ is anti-self-dual, then $F_A^+ \equiv 0$ over $X$ and
$$
\kappa^i(P) = -\frac{1}{4\pi^2 r_{\fg_i}} \int_X |F_A|_{\fg_i}^2\,d\vol \leq 0, \quad 1 \leq i \leq l.
$$
Consequently, if $P$ admits a self-dual connection, then $\bkappa(P) \leq 0$ while if $P$ admits an anti-self-dual connection, then $\bkappa(P) \geq 0$, where we use $\bkappa(P) \leq 0$ ($\geq 0$) as an abbreviation for $\kappa^i(P) \leq 0$ ($\geq 0$) for $1 \leq i \leq l$.

On the other hand, as the pointwise norm $|F_A|_\fg$ over $X$ of $F_A \in \Omega^2(X; \ad P)$ is defined by
$$
|F_A|_\fg^2 = -\sum_{i=1}^l\frac{1}{r_{\fg_i}}\tr_{\fg_i}(F_A\wedge F_A) = \sum_{i=1}^l \frac{1}{r_{\fg_i}}|F_A|_{\fg_i}^2,
$$
then
\begin{align*}
\int_X |F_A|_\fg^2 \,d\vol &= \sum_{i=1}^l \frac{1}{r_{\fg_i}}\int_X \left(|F_A^+|_{\fg_i}^2 + |F_A^-|_{\fg_i}^2\right)\,d\vol
\\
&\geq \left|\sum_{i=1}^l \frac{1}{r_{\fg_i}}\int_X \left(|F_A^+|_{\fg_i}^2 - |F_A^-|_{\fg_i}^2\right)\,d\vol\right|
\\
&= 4\pi^2|\kappa(P)|.
\end{align*}
Hence, $4\pi^2|\kappa(P)|$ represents a topological lower bound for the Yang-Mills energy functional,
$$
2\sE(A) = \int_X |F_A|_\fg^2 \,d\vol.
$$
If $\bkappa(P) \geq 0$, then $2\sE(A)$ achieves its lower bound,
$$
\sum_{i=1}^l \frac{1}{r_{\fg_i}}\int_X \left(|F_A^+|_{\fg_i}^2 - |F_A^-|_{\fg_i}^2\right)\,d\vol = 4\pi^2|\kappa(P)|,
$$
if and only if
$$
\sum_{i=1}^l \frac{1}{r_{\fg_i}}\int_X \left(|F_A^+|_{\fg_i}^2 + |F_A^-|_{\fg_i}^2\right)\,d\vol
=
\sum_{i=1}^l \frac{1}{r_{\fg_i}}\int_X \left(|F_A^+|_{\fg_i}^2 - |F_A^-|_{\fg_i}^2\right)\,d\vol,
$$
that is, if and only if
$$
\int_X |F_A^-|_\fg^2 \,d\vol = 0,
$$
in other words, if and only if $F_A^- \equiv 0$ over $X$ and $A$ is self-dual. Similarly, if $\bkappa(P) \leq 0$, then $2\sE(A)$ achieves its lower bound $-4\pi^2|\kappa(P)|$ if and only if $F_A^+ \equiv 0$ over $X$ and $A$ is anti-self-dual.

\section[Critical points of the Yang-Mills energy functional]{Critical points of the Yang-Mills energy functional and asymptotic limits of Yang-Mills gradient flow}
\label{sec:Critical_points_Yang-Mills_energy_functional}
It is well-known that there are topological obstructions to the existence of absolute minima, or anti-self-dual connections, on a principal $G$-bundle, $P$, over a closed, oriented, Riemannian, smooth four-dimensional manifold, $X$, for reasons evident from Donaldson and Kronheimer \cite{DK}, Freed and Uhlenbeck \cite{FU}, and Taubes \cite{TauSelfDual, TauIndef}. (The question of existence of self-dual connections is symmetric, by reversing the orientation of $X$, and so we may confine our discussion, without loss of generality, to the question of existence of anti-self-dual connections.) The nature of those topological obstructions can be inferred from Taubes' proofs of existence of anti-self-dual connections in \cite{TauSelfDual, TauIndef}, the generic metrics theorem of Freed and Uhlenbeck \cite{FU}, and discussions concerning the dimension of the moduli space of anti-self-dual connections in \cite{DK}.

Min-Oo's $L^2$-isolation result \cite{Min-Oo_1982} for Yang-Mills connections asserts that, if the curvature of the Riemannian metric on $X$ is pointwise positive over $X$ in the sense of \eqref{eq:Freed_Uhlenbeck_page_174_positive_metric}, then there is a positive constant, $\eps$, such that if $A$ is a Yang-Mills connection on $P$ with the property that
\begin{equation}
\label{eq:Minoo_theorem_2}
\|F_A^{+,g}\|_{L^2(X,g)} < \eps,
\end{equation}
then $A$ is necessarily anti-self-dual. (The positivity condition \eqref{eq:Freed_Uhlenbeck_page_174_positive_metric} on the Riemannian metric is the same as that exploited by Atiyah, Hitchin, and Singer \cite{AHS} and by Taubes in \cite{TauSelfDual}.) Therefore, when this curvature condition is obeyed, one can start the Yang-Mills gradient flow at an approximately anti-self-dual connection, $A_0$, produced by splicing onto $X$ an exact anti-self-dual connection from $S^4$ in the manner of \cite{TauSelfDual}, and such that $F_{A_0}^{+,g}$ obeys the $L^2$ bound \eqref{eq:Minoo_theorem_2}. It can be shown that the $L^2$-norm of the self-dual curvature is non-increasing with time (see Lemma \ref{lem:Nonincreasing_ASD_curvature_Yang-Mills_gradient_flow} in the sequel). Thus, if the gradient flow converges to a Yang-Mills connection, $A_\infty$, on $P$, this limit must necessarily be anti-self-dual. However, the existence of anti-self-dual connections in this setting is already known from work of Taubes \cite{TauSelfDual}, whose results establish existence not only when the geometric conditions of Min-Oo or Atiyah, Hitchin, and Singer are obeyed but also under the much weaker, purely topological condition that $b^+(X) = 0$, whenever existence is also known on a principal $G$-bundle over $S^4$ with the same Pontrjagin number.

It is more interesting to ask what happens when the Riemannian metric $g$ on $X$ does not obey the positivity condition of \cite{AHS, Min-Oo_1982} (which is always satisfied by the standard round metric on $S^4$) or even when $X$ does not obey the topological condition $b^+(X) = 0$ exploited by Taubes in \cite{TauSelfDual}. When $b^+(X) > 0$, Taubes still gives existence of anti-self-dual connections \cite{TauIndef} when $G$ is $\SU(2)$ or $\SO(3)$, but only when the Pontrjagin number is sufficiently large. Taubes' results were extended by Graham \cite{Taylor_2002} to allow for $G = \SO(n)$, for all integers $n \geq 4$. Such results, at least when $G$ is $\SU(2)$ or $\SO(3)$, are consistent with the generic metrics theorem of Freed and Uhlenbeck \cite{FU}, which asserts that the moduli space of anti-self-dual connections is a smooth (open) manifold, away from finitely points representing gauge-equivalence classes of reducible connections when $b^+(X) \geq 0$, and a smooth (open) manifold when $b^+(X) > 0$, provided $g$ is suitably generic. Thus, when the Pontrjagin number is such that the moduli space of anti-self-dual connections has positive expected dimension by the moduli space dimension formulas given, for example, in \cite{AHS, DK, FU, FrM, Lawson}, it is possible that the limiting Yang-Mills connection, $A_\infty$, produced by Yang-Mills gradient flow starting at an initial connection, $A_0$, with sufficiently small $\|F_{A_0}^{+,g}\|_{L^2(X,g)}$ is also an absolute minimum of the Yang-Mills energy functional, or an anti-self-dual connection, in other words. However, there is no guarantee that $A_\infty$ is anything other than a non-minimal Yang-Mills connection unless the Riemannian metric, $g$, obeys the positivity condition described in \cite{AHS, Min-Oo_1982, TauSelfDual}. As we know from work of Sibner, Sibner, and Uhlenbeck \cite{SibnerSibnerUhlenbeck}, with many further extensions developed by Bor and Montgomery \cite{Bor_1992, Bor_Montgomery_1990} and Sadun and Segert \cite{Sadun_Segert_1991, Sadun_Segert_1992cmp, Sadun_Segert_1992cpam, Sadun_Segert_1993, Sadun_1994}, non-minimal Yang-Mills connections exist on $S^4$ and thus one can not rule out their existence on arbitrary four-dimensional manifolds, $X$. (We remark that geometric questions concerning non-minimal Yang-Mills connections have also been explored by Stern \cite{Stern_2010}.) We can conclude that, when the Riemannian metric, $g$, on $X$ does \emph{not} obey the positivity condition positivity condition \eqref{eq:Freed_Uhlenbeck_page_174_positive_metric} but is generic in the Freed-Uhlenbeck sense, the group $G$ is $\SU(2)$ or $\SO(3)$, and the Pontrjagin number of $P$ is such that the moduli space of anti-self-dual connections has \emph{negative} expected dimension, then no matter how small the value of $\|F_{A_0}^{+,g}\|_{L^2(X,g)}$ for the initial connection, $A_0$, on $P$, the limiting connection, $A_\infty$, produced by Yang-Mills gradient flow must necessarily be a \emph{non-minimal} Yang-Mills connection.

\chapter{Linear and nonlinear evolutionary equations in Banach spaces}
\label{chapter:Sell_You_4}

\section{Linear evolutionary equations in Banach spaces}
\label{sec:Sell_You_4-2}
In this section, we shall review the \apriori estimates and existence, uniqueness, and regularity results which we shall need for linear evolutionary equations in Banach and Hilbert spaces \cite[Section 4.2]{Sell_You_2002}. In a later section, we apply this framework to deduce the corresponding results for the heat equation defined by the connection Laplace operator \eqref{eq:Connection_Laplacian}
on Sobolev spaces of sections of a vector bundle over a closed, Riemannian, smooth manifold, generalizing results for linear scalar parabolic equations described by Evans in \cite[Section 7.1]{Evans2} or for linear parabolic systems described by  Lady{\v{z}}enskaja, Solonnikov, and Ural$'$ceva in \cite[Chapter 8]{LadyzenskajaSolonnikovUralceva}.

\subsection{Analytic semigroups and sectorial operators}
\label{subsec:Sell_You_3-6}
We review some of the key concepts concerning analytic semigroups and sectorial operators, closely following \cite[Section 3.6]{Sell_You_2002}. A special class of $C^0$-semigroups (see \cite[Section 3.1]{Sell_You_2002}), namely, the \emph{analytic semigroups}, plays a fundamental role in the study of the dynamics of infinite-dimensional systems. There are two principal reasons why analytic semigroups are important in the study of systems of nonlinear parabolic partial differential equations.

The first reason is owing to the good information one has on the behavior of solutions as time $t \searrow 0$. It is known that, under reasonable conditions, a $C^0$-semigroup, $(e^{\cA t}, \cA)$, on a Banach space, $\cW$, is an analytic semigroup if and only if there are constants, $M_0 \geq 1$ and $M_1 > 0$, and $a \in \RR$, such that
$$
\|e^{\cA t}w\|_\cW \leq M_0 e^{at}\|w\|_\cW
\quad\hbox{and}\quad
\|\cA e^{\cA t}w\|_\cW \leq M_1 t^{-1} e^{at}\|w\|_\cW
$$
for all $w \in \cW$ and $t > 0$. Since the nonlinear evolutionary equations we consider are perturbations of a linear equation, by having good information about the behavior of solutions of the linear problem, as $t \searrow 0$, one can use this to study the initial value problem for related nonlinear problems.

The second reason for the importance of analytic semigroups is given in Theorem \ref{thm:Sell_You_37-5} below. A vector-valued function, $f : D \to \cW$, where $D$ is an open subset of the complex plane, $\CC$, and $\cW$ is a Banach space, is said to be \emph{analytic} if, for any $z_0 \in D$, the strong limit
$$
\lim_{z\to 0} \frac{1}{z}(f(z_0 + z) - f(z_0))
$$
exists in $\cW$. Recall the

\begin{defn}[Resolvent set]
\label{defn:Resolvent_set}
\cite[Section 13.26]{Rudin}, \cite[Section 8.1]{Yosida}
Let $\rho(\cA) \subset \CC$ denote the \emph{resolvent set} for $\cA$, that is, the set of all $\lambda \in \CC$ such that $\lambda  - \cA:\sD(\cA) \subset \cW \to \cW$ is a one-to-one map with dense range $\Ran(\lambda - \cA) \subset \cW$ and bounded inverse, $R(\lambda, \cA) := (\lambda - \cA)^{-1}$, the \emph{resolvent operator} on $\cW$.
\end{defn}

Given $a \in \RR$ and $\delta, \sigma \in (0, \pi)$, one defines sectors in the complex plane, $\CC$, by \cite[p. 77]{Sell_You_2002}
\begin{align}
\label{eq:Sell_You_page_77_complex_plane_sector_definition_Delta_of_a}
\Delta_\delta(a) := \{z \in \CC: |\arg(z-a)| < \delta \hbox{ and } z \neq a\},
\\
\label{eq:Sell_You_page_77_complex_plane_sector_definition_Sigma_of_a}
\Sigma_\sigma(a) := \{z \in \CC: |\arg(z-a)| > \sigma \hbox{ and } z \neq a\}.
\end{align}
When $a=0$, we abbreviate $\Delta_\delta(0) = \Delta_\delta$ and $\Sigma_\sigma(0) = \Sigma_\sigma$. Note that $\Delta_\delta(-a) = -\Sigma_\sigma(a)$, provided that $\delta = \pi -\sigma$ (see \cite[Figure 3.1]{Sell_You_2002}). As a result one has,
\begin{equation}
\label{eq:Sell_You_36-1}
\rho(\cA) \supset \Sigma_\sigma(a) \iff \rho(-\cA) \supset \Delta_\delta(-a), \quad\hbox{where } \delta = \pi - \sigma,
\end{equation}
for any linear operator $\cA : \sD(\cA) \subset \cW$, where $\sD(\cA)$ is dense in $\cW$.

\begin{defn}[Analytic semigroup]
\label{defn:Analytic_semigroup}
\cite[p. 77]{Sell_You_2002}
Let $(\cT(t), \cA)$ be a $C^0$-semigroup on $\cW$. One says say that $(\cT(t), \cA)$ is an
\emph{analytic semigroup} if there exists an extension of $\cT(t)$ to a mapping $\cT(z)$ defined for $z$ in some sector $\Delta_\delta \cup \{0\}$ and satisfying the following conditions:
\begin{enumerate}
\item The mapping $z \mapsto \cT(z)$ is a mapping of $\Delta_\delta \cup \{0\}$ into $\sL(\cW)$, the Banach space of bounded linear operators on $\cW$.

\item $\cT(z_1 + z_2) = \cT(z_1)\cT(z_2)$ for all $z_1$ and $z_2$ in $\Delta_\delta \cup \{0\}$.

\item For each $w \in \cW$, one has $\cT(z)w \to w$, as $z \to 0$ in $\Delta_\delta \cup \{0\}$.

\item For each $w \in \cW$, the function $z \mapsto \cT(z)w$ is an analytic mapping from $\Delta_\delta \cup \{0\}$ into $\cW$.
\end{enumerate}
\end{defn}

There remains the key question of determining which properties on the infinitesimal generator, $\cA$, will guarantee that a given $C^0$-semigroup is an analytic semigroup. This prompts the following

\begin{defn}[Sectorial operator]
\label{defn:Sell_You_page_78_definition_of_sectorial_operator}
\cite[p. 78]{Sell_You_2002}
A linear operator, $\cA : \sD(\cA) \subset \cW \to \cW$, is \emph{sectorial} if it obeys the following two conditions:
\begin{enumerate}
\item $\cA$ is densely defined and closed;
\item There exist real numbers $a \in \RR$, $\sigma \in (0, \pi/2)$, and $M \geq 1$, such that one has $\Sigma_\sigma(a) \subset \rho(\cA)$, and
\begin{equation}
\label{eq:Sell_You_36-2}
\|R(\lambda, \cA)\| \leq \frac{M}{|\lambda - a|}, \quad \forall\, \lambda \in \Sigma_\sigma(a),
\end{equation}
or equivalently by \eqref{eq:Sell_You_36-1}, that $\Delta_\delta(-a) \subset \rho(-\cA)$ where $\delta = \pi-\sigma$, and
\begin{equation}
\label{eq:Sell_You_36-3}
\|R(\lambda, -\cA)\| \leq \frac{M}{|\lambda + a|}, \quad \forall\, \lambda \in \Delta_\delta(-a),
\end{equation}
\end{enumerate}
A sectorial operator, $\cA$, is said to be \emph{positive} if it satisfies \eqref{eq:Sell_You_36-2} for some $a > 0$.
\end{defn}

See Haase \cite{Haase_2006} for a comprehensive treatment of sectorial operators.

\subsection{Fractional powers and interpolation spaces}
\label{subsec:Sell_You_3-7}
We review several key concepts and results from \cite[Section 3.7]{Sell_You_2002}.

Let $I$ be an interval in $\RR$. Following \cite[p. 93]{Sell_You_2002}, a family of Banach spaces $\calV^\alpha$ with norms $\|\cdot\|_\alpha$, defined for $\alpha \in I$, is said to be a family of \emph{interpolation spaces} on $I$ if the following property holds: For $\alpha, \beta \in I$ with $\alpha \geq \beta$, one has $\calV^\alpha \hookrightarrow \calV^\beta$, and $\calV^\alpha$ is dense in $\calV^\beta$. Moreover, for every $\theta \in [0,1]$, there is a constant, $C = C(\alpha, \beta, \theta)$, such that
$$
\|v\|_\gamma \leq C\|v\|_\alpha^\theta \|v\|_\beta^{1-\theta}, \quad\forall\, v \in V^\alpha,
$$
where $\gamma = \theta\alpha + (1-\theta)\beta$.

Let $\cW$ denote a given Banach space and let $\cA$ be a positive, sectorial operator on $\cW$. For any $\alpha > 0$, one defines \cite[Equation (37.8)]{Sell_You_2002}
\begin{equation}
\label{eq:Sell_You_37-8}
\cA^{-\alpha} := \frac{1}{\Gamma(\alpha)} \int_0^\infty t^{\alpha-1} e^{-\cA t}\, dt,
\end{equation}
where $e^{-\cA t}$ is the analytic semigroup generated by $-\cA$, and $\Gamma(\alpha)$ is the Gamma
function.

The \emph{fractional power} $\cA^\alpha$ of the operator $\cA$ is defined to be
$$
\cA^{-\alpha} := \left(\cA^{-\alpha}\right)^{-1}, \quad\forall\, \alpha > 0,
$$
and the domain is given by $\sD(\cA^\alpha) := \Ran(\cA^{-\alpha})$. Also one defines $\cA^0 := I$, the
identity on $\cW$. For $\alpha > 0$, the domain $\sD(\cA^\alpha)$ becomes a Banach space with respect to the graph norm,
\begin{equation}
\label{eq:Sell_You_37-9}
\|u\|_\alpha := \|\cA^\alpha u\|_\cW, \quad\forall\, u \in \sD(\cA^\alpha).
\end{equation}
Note that $\|u\|_0 = \|u\|_\cW$, for all $u \in \cW$. The interpolation space of order $2\alpha$ is defined to be $\calV^{2\alpha} := \sD(\cA^\alpha)$, where $\sD(\cA^\alpha)$ has the graph norm.

\begin{lem}[Basic properties of fractional power spaces defined by positive, sectorial operators]
\label{lem:Sell_You_37-4}
\cite[Lemma 37.4]{Sell_You_2002}
Let $\cA$ be a positive, sectorial operator on a Banach space $\cW$, and let $e^{-\cA t}$ denote the analytic semigroup on $\cW$ generated by $-\cA$. Then for $\alpha, \beta \geq 0$, the following properties are valid:
\begin{enumerate}
\item The operator $\cA^\alpha$ is a densely defined, closed linear operator.

\item For $\alpha \geq \beta$, one has $\sD(\cA^\alpha) \hookrightarrow \sD(\cA^\beta)$,
in terms of the graph norms on these spaces, and $\sD(\cA^\alpha)$ is dense in $\sD(\cA^\beta)$.

\item If in addition, $\cA$ has compact resolvent, then one has $\sD(\cA^\alpha) \hookrightarrow \sD(\cA^\beta)$, whenever $\alpha > \beta$.

\item One has $\cA^\alpha \cA^\beta = \cA^\beta \cA^\alpha = \cA^{\alpha+\beta}$ on $\sD(\cA^\gamma)$, for any $\alpha, \beta \in \RR$, where $\gamma := \max\{\alpha, \beta, \alpha+\beta\}$.

\item One has $\cA^\alpha e^{-\cA t}u = e^{-\cA t} \cA^\alpha u$, for all $u \in \sD(\cA^\alpha)$ and $t \geq 0$. Furthermore,
$e^{-\cA t}$ is an analytic semigroup on $\calV^\alpha$, for each $\alpha \in \RR$.

\item The mapping,
$$
[0, \infty) \times \calV^{2\alpha} \ni (t, u) \mapsto e^{-\cA t}u \in \calV^{2\alpha},
$$
is continuous, for every $\alpha \in \RR$.
\end{enumerate}
\end{lem}

\begin{thm}[Fundamental theorem on sectorial operators]
\label{thm:Sell_You_37-5}
\cite[Theorem 37.5]{Sell_You_2002}
Let $\cA$ be a positive, sectorial operator on a Banach space $\cW$, and let $\cT(t) = e^{-\cA t}$ be the analytic semigroup generated by $-\cA$. Then the following statements hold:
\begin{enumerate}
\item For any $r \geq 0$ and $t > 0$, the semigroup $e^{-\cA t}$ maps $\cW$ into $\sD(\cA^r)$, and it is strongly continuous in $t > 0$.

\item For any $r \geq 0$, there is a constant $M_r > 0$ such that
\begin{equation}
\label{eq:Sell_You_37-11}
\|e^{-\cA t}\|_{\sL(\cW, \sD(\cA^r))} = \|\cA^r e^{-\cA t}\|_{\sL(\cW)} \leq M_r t^{-r} e^{-at}, \quad\forall\, t > 0,
\end{equation}
where $a > 0$ is given by \eqref{eq:Sell_You_36-2}.

\item For $0 < \alpha \leq 1$, there is a constant $K_\alpha > 0$ such that
\begin{equation}
\label{eq:Sell_You_37-12}
\|e^{-\cA t}w - w\|_\cW \leq K_\alpha t^\alpha \|\cA^\alpha w\|_\cW, \quad\forall\, t \geq 0 \hbox{ and } w \in \sD(\cA^\alpha).
\end{equation}

\item The functions $e^{-\cA t}w$ are Lipschitz continuous in $t$, for $t > 0$ and $w \in \cW$. More precisely, for every $r \geq 0$, there is a constant $C_r > 0$ such that
\begin{equation}
\label{eq:Sell_You_37-13}
\left\|\cA^r \left(e^{-\cA (t+h)}w - e^{-\cA t}w\right)\right\|_\cW \leq C_r |h| t^{-(1+r)} \|w\|_\cW,
\end{equation}
for all $t > 0$ and $w \in \cW$.

\item For $w \in \cW$, one has
$$
e^{-\cA t}w \in C([0, \infty); \cW) \cap C((0, \infty; \calV^{2r}),
$$
for all $r \geq 0$. If in addition, one has $w \in \calV^{2\rho}$, for some $\rho > 0$, then
\begin{equation}
\label{eq:Sell_You_37-14}
e^{-\cA t}w \in C([0, \infty); \calV^{2\rho}) \cap C_{\loc}^{0,\rho-\alpha}([0, \infty); \calV^{2\alpha})
\cap C_{\loc}^{0,1}((0, \infty); \calV^{2r}),
\end{equation}
for all $r \geq 0$, and all $\alpha$ with $0 \leq \alpha < \rho$.

\item If in addition, the sectorial operator $\cA$ has compact resolvent, then the analytic semigroup $e^{-\cA t}w$ is compact, for $t > 0$.
\end{enumerate}
\end{thm}

\subsection{Solution concepts and the variation of constants formula}
\label{subsec:Sell_You_4-2-preamble}
We recall the

\begin{hyp}[Standing Hypothesis A]
\label{hyp:Sell_You_4_standing_hypothesis_A}
\cite[Standing Hypothesis A, p. 141]{Sell_You_2002}
Let $\cA$ be a positive, sectorial operator on a Banach space $\cW$ with associated analytic semigroup $e^{-\cA t}$. Let $\calV^{2\alpha}$ be the family of interpolation spaces generated by the fractional powers of $\cA$, where $\calV^{2\alpha} = \sD(\cA^\alpha)$, for $\alpha \geq 0$. Let $\|\cA^\alpha u\|_\cW = \|u\|_{\calV^{2\alpha}}$ denote the norm on $\calV^{2\alpha}$. See Lemma \ref{lem:Sell_You_37-4} for more information.
\end{hyp}

\begin{hyp}[Standing Hypothesis B]
\label{hyp:Sell_You_4_standing_hypothesis_B}
\cite[Standing Hypothesis B, p. 142]{Sell_You_2002}
The operator $\cA$ is a positive, self-adjoint, linear operator, with compact resolvent, on a Hilbert space $\cH$. Consequently, $\cA$ satisfies the Hypothesis \ref{hyp:Sell_You_4_standing_hypothesis_A}. Moreover, the fractional power spaces $\calV^\alpha$ are defined for all $\alpha \in \RR$ and
\cite[Equation (37.2)]{Sell_You_2002} defines the Hilbert space structure on each $\calV^\alpha$. Also the semigroup $e^{-\cA t}$ is compact, for $t > 0$. See \cite[Theorem 37.2]{Sell_You_2002} for more information.
\end{hyp}


We will focus first on the \emph{linear inhomogeneous evolution equation} on a Banach space $\cW$ and initial condition,
\begin{align}
\label{eq:Sell_You_42-1}
\frac{\partial u}{\partial t} + (\cA u)(t) &= f(t) \quad\hbox{in } \cW, \quad \hbox{for } t > 0,
\\
\label{eq:Sell_You_42-2}
u(0) & = u_0 \in \cW.
\end{align}
We assume that $-\cA$ is the infinitesimal generator of an analytic semigroup, $e^{-\cA t}$, and that $f \in L_{\loc}([0, T); \cW)$, where $0 < T \leq \infty$. We call the expression
\begin{equation}
\label{eq:Sell_You_42-3}
u(t) := e^{-\cA t}u_0 + \int_0^t e^{-\cA (t - s)}f(s)\,ds, \quad t \in [0, T).
\end{equation}
the \emph{variation of constants formula}, where the integral is in the Bochner sense and represents a point in $\cW$ for
each $t > 0$ \cite[Appendix C]{Sell_You_2002}.

\subsection{Solutions via analytic semigroup theory}
\label{subsec:Sell_You_4-2-2}
The \emph{Newton-Leibnitz} formula in the space $\cW$ takes the form,
\begin{equation}
\label{eq:Sell_You_42-7}
\frac{d}{dt} \int_0^t e^{-\cA (t-s)}f(s)\,ds = -\cA \int_0^t e^{-\cA (t-s)}f(s)\,ds + f(t),
\quad\hbox{a.e } t > 0.
\end{equation}
We recall the following result, which yields further information on the regularity of mild solutions, when $p > 1$.

\begin{lem}[Regularity of mild solutions to a linear evolution equation]
\cite[Lemma 42.7]{Sell_You_2002}
Let $f \in L^p_{\loc}([0, T); \cW)$ and let the Hypothesis \ref{hyp:Sell_You_4_standing_hypothesis_A} hold, where $1 \leq p \leq \infty$ and $0 < T \leq \infty$. Then the following are valid:
\begin{enumerate}
\item The implications
\begin{equation}
\label{eq:Sell_You_42-12}
\begin{aligned}
\hbox{classical solution in } \calV^{2\alpha} &\implies \hbox{strong solution in } \calV^{2\alpha}
\\
&\implies \hbox{mild solution in } \calV^{2\alpha}
\end{aligned}
\end{equation}
for a solution $u$, are valid for any $\alpha$ with $\alpha \geq 0$.

\item Assume that $p > 1$. Let $u_0 \in \cW$ be fixed and let $u$ satisfy \eqref{eq:Sell_You_42-3} in $\cW$ on $[0, T)$. Then $u$ is a mild solution of equation \eqref{eq:Sell_You_42-1}
in $\cW$ and
\begin{equation}
\label{eq:Sell_You_42-13}
u \in C([0, T); \cW) \cap C^{0, \theta}((0, T); \calV^{2r}),
\end{equation}
for every $r$ with $0 \leq r < 1 - 1/p$, where $0 < \theta_0 < 1 - 1/p - r$.

\item Assume that $p > 1$ and $u_0 \in \calV^{2\rho}$, where $0 < \rho < 1 - 1/p$. Then $u$ is a mild solution of equation \eqref{eq:Sell_You_42-1} in $\calV^{2\alpha}$ and
\begin{equation}
\label{eq:Sell_You_42-14}
u \in C([0, T); \calV^{2\alpha}) \cap C^{0, \theta_1}((0, T); \calV^{2r}) \cap C^{0, \theta_2}((0, T); \calV^{2\sigma}),
\end{equation}
for every $\alpha$, $r$, and $\sigma$ with $0 \leq \alpha \leq \rho$, and $0 \leq r < \rho$, and $\rho \leq \sigma < 1 - 1/p$, where $\theta_1 = \theta_1(r)$ and $\theta_2 = \theta_2(r)$ are positive.
\end{enumerate}
\end{lem}

\begin{thm}[Existence of mild and strong solutions to a linear evolution equation]
\label{thm:Sell_You_42-9}
\cite[Theorem 42.9]{Sell_You_2002}
Assume that Hypothesis \ref{hyp:Sell_You_4_standing_hypothesis_A} is satisfied and
\begin{equation}
\label{eq:Sell_You_42-16}
f \in L^p_{\loc}([0, T); \cW) \cap C^{0,\theta}((0, T); \cW),
\end{equation}
where $1 \leq p \leq \infty$ and $\theta \in (0, 1)$ and $0 < T \leq \infty$. For any $u_0 \in \cW$, define $u$ by the variation of constants formula \eqref{eq:Sell_You_42-3}. Then $u$ is mild solution of equation \eqref{eq:Sell_You_42-1} in $\cW$ that satisfies
\begin{equation}
\label{eq:Sell_You_42-18}
u \in C([0, T); \cW) \cap C^{0, 1-r}((0, T); \calV^{2r}) \cap C((0, T); \sD(\cA)).
\end{equation}
for every $r \in [0, 1)$. In addition, the Newton-Leibnitz Formula \eqref{eq:Sell_You_42-7} is valid in $\cW$, and $u$ is a strong solution in $\cW$, for $0 \leq t < T$. Moreover, if $u_0 \in \calV^{2\rho}$, for some $\rho \in (0, 1]$, then $u$ is a strong solution in $\calV^{2\rho}$ on $[0, T)$.
\end{thm}

In the next result we show that by assuming additional spatial regularity for the forcing function $f$, one obtains a significant improvement in the spatial regularity of the mild solutions.

\begin{thm}[Higher spatial regularity of mild solutions]
\label{thm:Sell_You_42-10}
\cite[Theorem 42.10]{Sell_You_2002}
In addition to Hypothesis \ref{hyp:Sell_You_4_standing_hypothesis_A}, let
\begin{equation}
\label{eq:Sell_You_42-20}
f \in L^p_{\loc}([0, T); \calV^{2\nu}), \quad\hbox{where } p > 1 \hbox{ and } 1 < \nu p \leq \infty.
\end{equation}
Let $\beta$ satisfy $0 \leq \beta < \nu - 1/p$. Assume that $u_0 \in \calV^{2\mu}$, where $0 < \mu \leq 1 + \beta$, and let $u$ be the mild solution of equation \eqref{eq:Sell_You_42-1} in $\cW$ on $[0, T)$. Then $u$ is a strong solution to equation \eqref{eq:Sell_You_42-1} in $\calV^{2\mu}$ on $[0, T)$, and $u$ satisfies
\begin{equation}
\label{eq:Sell_You_42-21}
u \in C([0, T); \calV^{2\mu}) \cap C^{0,\theta_0}_{\loc}([0, T); \calV^{2\rho}) \cap C_{\loc}^{0,\theta_1}((0, T); \calV^{2+2\beta}),
\end{equation}
where $0 \leq \rho < \mu$ and $0 \leq \beta < \nu - 1/p$ and $\theta_0$ and $\theta_1$ are positive. Assume instead that \eqref{eq:Sell_You_42-20} is replaced by the condition
\begin{equation}
\label{eq:Sell_You_42-22}
f \in L^1_{\loc}([0, T); \cW) \cap L^p_{\loc}((0, T); \calV^{2\nu}), \quad\hbox{where } p > 1 \hbox{ and } 1 < \nu p \leq \infty.
\end{equation}
and that $u_0 \in \cW$. Then for every $\tau \in (0,T)$, the \emph{translate} $u_\tau := u(\cdot + \tau)$ \cite[Section 2.1.3]{Sell_You_2002} satisfies \eqref{eq:Sell_You_42-21}, where the interval $[0, T)$ is replaced by $[0, T - \tau)$. Furthermore, $u_\tau$ is a strong solution to equation \eqref{eq:Sell_You_42-21} in $\calV^{2\nu}$ on $[0, T-\tau)$, for every $\mu$ with $0 \leq \mu \leq 1 + \beta$.
\end{thm}

\subsection{Weak solutions to evolutionary linear equations in Hilbert spaces}
\label{subsec:Sell_You_4-2-3}
We assume that Hypothesis \ref{hyp:Sell_You_4_standing_hypothesis_B} (`Standing Hypothesis B') is satisfied. We review a fourth solution concept, that of a weak solution, and describe its relationship to the other solution concepts given in Section \ref{subsec:Sell_You_4-7-1} (for evolutionary nonlinear equations on Banach spaces --- see \cite[Section 4.2]{Sell_You_2002} for linear equations), namely mild, strong, and classical. In particular, given $\alpha \in \RR$, one says that a function $u = u(t)$ is a \emph{weak solution} of the linear evolution equation \eqref{eq:Sell_You_42-1} \emph{in the space} $\calV^\alpha$, on the interval $[0, T)$, where $0 < T \leq \infty$, provided that $u(0) = u_0 \in \calV^\alpha$ and $f \in L_{\loc}^2([0, T); \calV^{\alpha-1})$, and that the following properties are satisfied:
\begin{enumerate}
\item one has
\begin{equation}
\label{eq:Sell_You_42-25}
u \in L^\infty_{\loc}([0,T); \calV^\alpha) \cap L^2_{\loc}([0,T); \calV^{\alpha+1});
\end{equation}

\item the function $u$ has a time derivative $\partial_t u \in L_{\loc}^p([0, T); \calV^{\alpha-1})$, for some $p$ with $1 \leq p < \infty$, such that the equation
\begin{equation}
\label{eq:Sell_You_42-26}
u(t) - u(t_0) = \int_{t_0}^t \frac{\partial u}{\partial t}(s) \,ds \quad\hbox{holds in } \calV^{\alpha-1},
\quad\forall\, t_0, t \in [0, T);
\end{equation}

\item one has
\begin{equation}
\label{eq:Sell_You_42-27}
\|u(t)\|_{\calV^\alpha} + \int_{t_0}^t \|\cA^{\frac{1}{2}} u(s)\|_{\calV^\alpha}^2\,ds
\leq
\|u(t_0)\|_{\calV^\alpha} + \int_{t_0}^t \|\cA^{-\frac{1}{2}} f(s)\|_{\calV^\alpha}^2\,ds,
\end{equation}
for almost all $t_0$ and $t$ with $0 \leq t_0 \leq t < T$; and

\item the function $u$ satisfies the weak-integrated form of Equation \eqref{eq:Sell_You_42-1},
\begin{equation}
\label{eq:Sell_You_42-28}
\begin{split}
(\cA^{-\frac{1}{2}}(u(t) - u(t_0)), \cA^{\frac{1}{2}}v)_{\calV^\alpha} + \int_{t_0}^t (\cA^{\frac{1}{2}} u(s), \cA^{\frac{1}{2}} v(s))_{\calV^\alpha}\,ds
\\
= \int_{t_0}^t (\cA^{-\frac{1}{2}} f(s), \cA^{\frac{1}{2}} v(s))_{\calV^\alpha}\,ds,
\end{split}
\end{equation}
for all $v \in \calV^{\alpha+1}$ and all $t_0$ and $t$ with $0 \leq t_0 \leq t < T$.
\end{enumerate}

There are several properties inherited by every weak solution of Equation \eqref{eq:Sell_You_42-1}. In particular, one has the following result.

\begin{lem}[Properties of weak solutions]
\cite[Lemma 42.11]{Sell_You_2002}
\label{lem:Sell_You_42-11}
In addition to the Hypothesis \ref{hyp:Sell_You_4_standing_hypothesis_B}, let $u = u(t)$ be a weak solution of Equation \eqref{eq:Sell_You_42-1} in $\calV^\alpha$, on the interval $[0, T)$, where $0< T \leq \infty$, and
$$
u_0 \in \calV^\alpha \quad\hbox{and}\quad f \in L_{\loc}^p([0, T); \calV^{\alpha-1}), \quad\hbox{for } 2 \leq p \leq \infty.
$$
Then the following properties are valid:
\begin{enumerate}
\item The time derivative satisfies $\partial_t u \in L_{\loc}^2([0, T); \calV^{\alpha-1})$ and $u$ satisfies
$$
u \in C([0, T); \calV^\alpha) \cap C^{0,\frac{1}{2}}_{\loc}([0, T); \calV^{\alpha-1}).
$$
\item One has $\partial u_t + \cA u = f$ a.e. on $(0, T)$ in the space $\calV^{\alpha-1}$, and $u$ is a mild solution
of Equation \eqref{eq:Sell_You_42-1} in $\calV^\beta$, for each $\beta < \alpha + (p - 2)p^{-1}$.
\end{enumerate}
\end{lem}

The next result provides existence of weak solutions and the connection between these solutions and mild solutions. One uses the
inequality $\lambda_1\|v\|_{\calV^\alpha}^2 \leq \|\cA^{\frac{1}{2}}v\|_{\calV^\alpha}^2$, where $\lambda_1 > 0$ is the first eigenvalue of the positive, self-adjoint operator $\cA$.

\begin{thm}[Existence, uniqueness, and regularity of weak solutions to a linear evolutionary equation in a Hilbert space]
\cite[Theorem 42.12]{Sell_You_2002}
\label{thm:Sell_You_42-12}
In addition to the Hypothesis \ref{hyp:Sell_You_4_standing_hypothesis_B}, let
$$
u_0 \in \calV^\alpha, \quad\hbox{and}\quad f \in L_{\loc}^p([0, T); \calV^{\alpha-1}), \quad\hbox{for } 2 \leq p \leq \infty.
$$
Then the following properties hold:
\begin{enumerate}
\item
\label{item:Theorem_Sell_You_42-12_1}
There is a unique weak solution $u = u(t)$ of Equation \eqref{eq:Sell_You_42-1}, in the space $\calV^\alpha$, on the interval $[0, T)$, with $u(0) = u_0$,
\begin{equation}
\label{eq:Sell_You_42-29}
\begin{gathered}
\frac{\partial u}{\partial t} \in L_{\loc}^2([0, T); \calV^{\alpha-1}), \quad\hbox{and}
\\
u \in C([0, T); \calV^\beta) \cap C^{0,\theta_1}_{\loc}([0, T); \calV^\sigma) \cap L_{\loc}^2([0, T); \calV^{\alpha+1}),
\end{gathered}
\end{equation}
for each $\beta$ and $\sigma$, where $\beta \leq \alpha$, and $\sigma < \alpha$, and $\theta_1 = \theta_1(\sigma) > 0$. Furthermore, $u$ is a mild solution of Equation \eqref{eq:Sell_You_42-1} in $\calV^\beta$, for each $\beta < \alpha + (p - 2)p^{-1}$.

\item
\label{item:Theorem_Sell_You_42-12_2}
For $0 \leq t < T$, the following inequalities are valid in the space $\calV^\alpha$:
\begin{subequations}
\label{eq:Sell_You_42-30}
\begin{align}
\label{eq:Sell_You_42-30_Linfty_time_Valpha_space_u}
\|u(t)\|_{\calV^\alpha}^2
&\leq
e^{-\lambda_1 t}\|u_0\|_{\calV^\alpha}^2 + \int_0^t e^{-\lambda_1(t-s)}\|\cA^{-\frac{1}{2}}f(s)\|_{\calV^\alpha}^2 \,ds,
\\
\label{eq:Sell_You_42-30_L2_time_Valpha+1_space_u}
\int_0^t\|\cA^{\frac{1}{2}}u(s)\|_{\calV^\alpha}^2 \,ds
&\leq
\|u_0\|_{\calV^\alpha}^2 + \int_0^t \|\cA^{-\frac{1}{2}}f(s)\|_{\calV^\alpha}^2 \,ds,
\\
\label{eq:Sell_You_42-30_L2_time_Valpha-1_space_partial_t_u}
\int_0^t\|\cA^{-\frac{1}{2}}\partial_tu(s)\|_{\calV^\alpha}^2 \,ds
&\leq
2\|u_0\|_{\calV^\alpha}^2 + 4\int_0^t \|\cA^{-\frac{1}{2}}f(s)\|_{\calV^\alpha}^2 \,ds.
\end{align}
\end{subequations}
\item
\label{item:Theorem_Sell_You_42-12_3}
If $p > 2$, then $u$ is a mild solution of Equation \eqref{eq:Sell_You_42-1} in $\calV^\alpha$, and in addition to \eqref{eq:Sell_You_42-29}, one has
\begin{equation}
\label{eq:Sell_You_42-31}
C^{0,\theta_2}_{\loc}((0, T); \calV^\sigma),
\end{equation}
for each $\sigma$, where $\alpha \leq \sigma < \alpha + 1 - 2/p$ and $\theta_2 = \theta_2(\sigma) > 0$.
\end{enumerate}
\end{thm}

Theorem \ref{thm:Sell_You_42-12} can used as a part of a bootstrap argument, where one has $u_0 \in \calV^{\alpha+1}$ and $f \in L_{\loc}^p([0, T); \calV^\alpha)$, for some $p$ with $2 \leq p \leq \infty$. Since one has the embeddings,
$$
u_0 \in \calV^{\alpha+1} \hookrightarrow \calV^\alpha \quad\hbox{and}\quad
f \in L_{\loc}^p([0, T); \calV^{\alpha+1}) \hookrightarrow L_{\loc}^p([0, T); \calV^\alpha),
$$
then Theorem \ref{thm:Sell_You_42-12} is applicable for both $\alpha$ and $\alpha + 1$. Furthermore, the uniqueness of the weak solutions implies that the weak solution $u = u(t)$ in $\calV^\alpha$ is, in fact, a weak solution in $\calV^{\alpha+1}$. Hence, it is a mild solution in $\calV^\beta$, for $\beta < \alpha + 1 + (p - 2)p^{-1}$ by Theorem \ref{thm:Sell_You_42-12}. As a result, inequalities \eqref{eq:Sell_You_42-30} are valid, where $\alpha$ is replaced by $\alpha + 1$; and from \eqref{eq:Sell_You_42-29}, one has
\begin{equation}
\label{eq:Sell_You_42-34}
u \in C([0, T); \calV^\beta) \cap C^{0,\theta_1}_{\loc}([0, T); \calV^\sigma) \cap L_{\loc}^2([0, T); \calV^{\alpha+2}),
\end{equation}
for $\beta \leq \alpha + 1$ and $\sigma < \alpha + 1$, where $\theta_1 = \theta_1(\sigma) > 0$. If in addition, $p > 2$,
then \eqref{eq:Sell_You_42-31} implies that for some $\theta_2 = \theta_2(\sigma) > 0$, one has
\begin{equation}
\label{eq:Sell_You_42-35}
u \in C^{0,\theta_1}_{\loc}((0, T); \calV^\sigma) \quad\hbox{for } \alpha + 1 \leq \sigma < \alpha + 2 - \frac{2}{p}.
\end{equation}
Moreover, in this setting, one can also apply \cite[Theorem 42.10]{Sell_You_2002}, which leads to valuable information about the strong solutions of Equation \eqref{eq:Sell_You_42-1}.

\begin{cor}[Regularity of weak solutions]
\cite[Corollary 42.13]{Sell_You_2002}
\label{cor:Sell_You_42-13}
In addition to the Hypothesis \ref{hyp:Sell_You_4_standing_hypothesis_B}, let $u_0 \in \calV^\alpha$ and let $f \in L_{\loc}^p([0, T); \calV^\alpha)$, for some $p$ with $2 \leq p \leq \infty$. Then the following statements are valid:
\begin{enumerate}
\item The weak solution $u$ of Equation \eqref{eq:Sell_You_42-1} in $\calV^\alpha$ on $[0, T)$ is a strong solution to equation  \eqref{eq:Sell_You_42-1} in $\calV^\alpha$ on $[0, T)$, and $u$ satisfies
\begin{equation}
\label{eq:Sell_You_42-36}
u \in C([0, T); \calV^\alpha) \cap C^{0,\theta_0}_{\loc}([0, T); \calV^\sigma)
\cap C^{0,\theta_1}_{\loc}((0, T); \calV^{\alpha+1+2\beta}),
\end{equation}
for every $\sigma$ and $\beta$ with $\sigma < \alpha$ and $0 \leq \beta < 1 - 1/p$, where $\theta_0 = \theta_0(\sigma)$ and $\theta_1 = \theta_1(\beta)$ are positive.

\item For every $\tau \in (0, T)$, one has $u(\tau) \in \calV^{\alpha+1}$, and the conclusions in
Items \eqref{item:Theorem_Sell_You_42-12_1}--\eqref{item:Theorem_Sell_You_42-12_3} in Theorem \ref{thm:Sell_You_42-12} are valid on the interval $\tau \leq t < T$, with $\alpha$ replaced by $\alpha+1$ and $u_0$ replaced by $u(\tau)$. In addition, \eqref{eq:Sell_You_42-34} and \eqref{eq:Sell_You_42-35} are valid for the \emph{translate} $u_\tau := u(\cdot + \tau)$ \cite[Section 2.1.3]{Sell_You_2002}, where the interval $[0, T)$ is replaced by $[0, T - \tau)$. Furthermore, the following inequalities are valid, for $0<\tau<t<T$:
\begin{subequations}
\label{eq:Sell_You_42-37}
\begin{align}
\label{eq:Sell_You_42-37_Linfty_time_Valpha+1_space_u}
\|u(t)\|_{\calV^{\alpha+1}}^2
&\leq
e^{-\lambda_1 (t-\tau)}\|u(\tau)\|_{\calV^{\alpha+1}}^2 + \int_\tau^t e^{-\lambda_1(t-s)}\|f(s)\|_{\calV^\alpha}^2 \,ds,
\\
\label{eq:Sell_You_42-37_L2_time_Valpha+2_space_u}
\int_\tau^t\|\cA^{\frac{1}{2}}u(s)\|_{\calV^{\alpha+1}}^2 \,ds
&\leq
\|u(\tau)\|_{\calV^{\alpha+1}}^2 + \int_\tau^t \|f(s)\|_{\calV^\alpha}^2 \,ds,
\\
\label{eq:Sell_You_42-37_L2_time_Valpha_space_partial_t_u}
\int_\tau^t\|\cA^{-\frac{1}{2}}\partial_tu(s)\|_{\calV^{\alpha+1}}^2 \,ds
&\leq
2\|u(\tau)\|_{\calV^{\alpha+1}}^2 + 4\int_{\tau}^t \|f(s)\|_{\calV^\alpha}^2 \,ds.
\end{align}
\end{subequations}
\end{enumerate}
\end{cor}

If the source function $f$ has additional \emph{temporal} regularity, the mild solution $u$ to \eqref{eq:Sell_You_42-1} acquires
additional spatial regularity.

\begin{thm}[Spatial regularity implied by temporal regularity of the source function and spatial regularity of the initial data]
\cite[Theorem 42.14]{Sell_You_2002}
\label{thm:Sell_You_42-14}
In addition to the Hypothesis \ref{hyp:Sell_You_4_standing_hypothesis_B}, let
$$
f \in C([0, T); \calV^\alpha) \cap W^{1,p}([0, T); \calV^{\alpha-1}),
$$
for some $p$ with $2 \leq p \leq \infty$, and let $g = \partial_tf$. Let $u_0 \in \calV^{\alpha+2}$ and define $v_0 := f(0) - \cA u_0 \in \calV^\alpha$. Let $u = u(t)$ be the weak solution to Equation \eqref{eq:Sell_You_42-1} in the space $\calV^\alpha$, with $u(0) = u_0$. Then the following properties hold:
\begin{enumerate}
\item $u$ is a mild solution of Equation \eqref{eq:Sell_You_42-1} in $\calV^\sigma$, for each $\sigma < \alpha + 2$, and $u$ is a strong solution in $\calV^\beta$, for each $\beta \leq \alpha$. Moreover, one has
\begin{equation}
\label{eq:Sell_You_42-38}
u \in C^1([0, T); \calV^\beta) \cap C([0, T); \calV^\nu) \cap C^{0,\theta_1}_{\loc}([0, T); \calV^\sigma),
\end{equation}
for all $\beta$, $\nu$, and $\sigma$ with $\beta \leq \alpha$, and $\nu \leq \alpha + 2$, and $\sigma < \alpha + 2$, where
$\theta_1 = \theta_1(\sigma) > 0$.

\item The function $v := \partial_t u$ satisfies
\begin{equation}
\label{eq:Sell_You_42-39}
v(t) = e^{-\cA t}v_0 + \int_0^t e^{-\cA(t-r)}g(r)\,dr, \quad\hbox{for } t \geq 0,
\end{equation}
in any space $\calV^\beta$, with $\beta < \alpha$, and $v$ is a weak solution of $\partial_tv + \cA v = g(t)$ in the space $\calV^\alpha$ with
\begin{equation}
\label{eq:Sell_You_42-40}
v \in C([0, T); \calV^\beta) \cap L_{\loc}^2([0, T); \calV^{\alpha+1}) \cap C^{0,\theta_2}([0, T); \calV^\sigma),
\end{equation}
for all $\beta$ and $\sigma$ with $\beta \leq \alpha$ and $\sigma < \alpha$, where $\theta_2 = \theta_2(\sigma) > O$.
\end{enumerate}
\end{thm}

Finally, we turn to an issue raised in Theorem \ref{thm:Sell_You_42-14}, which shows that if $u_0 \in \calV^{\alpha+2}$, then $u$ satisfies property \eqref{eq:Sell_You_42-38}. One can ask what happens with other mild solutions, where $u_0 \in \calV^\alpha$, but $u_0 \notin \calV^{\alpha+2}$? An answer to this question is provided by

\begin{thm}[Spatial regularity implied by temporal regularity of the source function]
\cite[Theorem 42.15]{Sell_You_2002}
\label{thm:Sell_You_42-15}
In addition to the Hypothesis \ref{hyp:Sell_You_4_standing_hypothesis_B}, let
$$
f \in C([0, T); \calV^\alpha) \cap W^{1,2}([0, T); \calV^{\alpha-1}).
$$
For any $u_0 \in \calV^\alpha$, let $u = u(t)$ denote the weak solution of Equation \eqref{eq:Sell_You_42-1} in $\calV^\alpha$ on $[0, T)$. Then the following hold:
\begin{enumerate}
\item $u$ is a strong solution in $\calV^\alpha$, and it satisfies \eqref{eq:Sell_You_42-29}, \eqref{eq:Sell_You_42-31}, and \eqref{eq:Sell_You_42-36}, with $p = \infty$.

\item $u$ is a mild solution in $\calV^\beta$, for each $\beta < \alpha + 1$.

\item $u$ satisfies
\begin{equation}
\label{eq:Sell_You_42-43}
u \in C^1((0, T); \calV^\beta) \cap C((0, T); \calV^\nu) \cap C^{0,\theta_4}_{\loc}((0, T); \calV^\sigma),
\end{equation}
for every $\beta$, $\nu$, and $\sigma$ with $\beta \leq \alpha$, $\nu \leq \alpha + 2$, and $\sigma < \alpha + 2$, where
$\theta_4 = \theta_4(\sigma) > O$.

\item For every $\tau \in (0, T)$, the translate $u_\tau$ is a strong solution to \eqref{eq:Sell_You_42-1}, where $f$ is replaced by $f_\tau \in \calV^{\alpha+1}$ on $[0,T - \tau)$.
\end{enumerate}
\end{thm}

\section{Local existence for a nonlinear evolution equation in a Banach space}
\label{sec:Local_existence_nonlinear_evolution_equation_Banach_space}
In this section, we develop an approach based on semigroup theory to the question of local existence and uniqueness of a solutions to a nonlinear evolution equation in a Banach space. The results we discuss both review and further develop those described by Sell and You in \cite[Chapter 4]{Sell_You_2002}.

\subsection{Local existence and uniqueness results for mild solutions}
\label{subsec:Sell_You_4-7-1}
We assume that we are in the setting of Hypothesis \ref{hyp:Sell_You_4_standing_hypothesis_A} (the `Standing Hypothesis A' of Sell and You).

\begin{defn}[Maps in $C_{\Lip}([0,\infty)\times \calV; \cW)$]
\label{defn:Sell_You_page_221_CLip_0_to_infinity_times_V_into_W}
\cite[p. 221]{Sell_You_2002}
For Banach spaces, $\calV$ and $\cW$, let $C_{\Lip}([0,\infty)\times \calV; \cW)$ denote the set of continuous, locally bounded, Lipschitz maps,
$$
\cF:[0,\infty)\times \calV \to \cW,
$$
such that for any bounded set $B \subset \calV$, there are positive constants $K_0=K_0(B)$ and $K_1=K_1(B)$ with the properties that
\begin{align}
\label{eq:Sell_You_46-1}
\|\cF(t,v)\|_\cW &\leq K_0, \quad\forall\, t \geq 0 \quad\hbox{and}\quad v \in B,
\\
\label{eq:Sell_You_46-2}
\|\cF(t,v_1) - \cF(t,v_2)\|_\cW &\leq K_1\|v_1-v_2\|_\calV, \quad\forall\, t \geq 0 \quad\hbox{and}\quad v_1, v_2 \in B.
\end{align}
\end{defn}

Let $\calV$ and $\cW$ be Banach spaces and, for $\beta \in [0,1)$, require that
$$
\cF:[0,\infty)\times \calV^{2\beta} \to \cW
$$
belong to $C_{\Lip}([0,\infty)\times \calV^{2\beta}; \cW)$, the set of continuous, locally bounded, Lipschitz maps, such that for each bounded set $B \subset \calV^{2\beta}$, there are positive constants $K_0=K_0(B)$ and $K_1=K_1(B)$ such that,
\begin{align}
\label{eq:Sell_You_46-7}
\|\cF(t,v)\|_\cW &\leq K_0, \quad\forall\, t \geq 0 \quad\hbox{and}\quad v \in B,
\\
\label{eq:Sell_You_46-8}
\|\cF(t,v_1) - \cF(t,v_2)\|_\cW &\leq K_1\|v_1-v_2\|_{\calV^{2\beta}}, \quad\forall\, t \geq 0 \quad\hbox{and}\quad v_1, v_2 \in B.
\end{align}
Let $I = [t_0, t_0 + \tau)$ be an interval in $[0,\infty)$, where $\tau > 0$, and let $\rho\geq 0$. A pair $(u , I)$ is called a \emph{mild solution} of
\begin{equation}
\label{eq:Sell_You_47-1}
\dot u(t) + \cA u(t) = \cF(t,u(t)), \quad\hbox{for } u(t_0) = u_0 \in \cW \hbox{ and all } t \geq t_0 \geq 0,
\end{equation}
\emph{in the space} $\calV^{2\rho}$ on $I$ if $u \in C(I; \calV^{2\rho})$ and is a solution of the integral equation in $\calV^{2\rho}$,
\begin{equation}
\label{eq:Sell_You_47-2}
u(t) = e^{-(t-t_0)\cA}u_0 + \int_{t_0}^t e^{-(t-s)\cA} \cF(s,u(s))\,ds, \quad\hbox{for all } t \in I.
\end{equation}
Note that the initial condition obeys $u(t_0) = u_0 \in \calV^{2\rho}$.

A pair $(u,I)$, where $u \in C(I; \calV^{2\rho})$ and satisfies $u(t_0) = u_0$, is (strongly) differentiable almost everywhere with $\dot u$ and $\cA u$ in $L^1_{\loc}(I, \cW)$,
and satisfies the differential equation in $\cW$,
\begin{equation}
\label{eq:Sell_You_47-3}
\dot u(t) + \cA u(t) = \cF(t,u(t)), \quad\hbox{for a.e. } t \in (t_0, t_0 + \tau),
\end{equation}
is called a \emph{strong solution} of \eqref{eq:Sell_You_47-1} \emph{in the space} $\calV^{2\rho}$ on $I$. If in addition, one has $\dot u \in C(I, \calV^{2\rho})$ and the differential equation in \eqref{eq:Sell_You_47-3} is satisfied for all $t \in (t_0, t_0 + \tau)$ , then $(u,I)$ is called a \emph{classical solution} of \eqref{eq:Sell_You_47-1} \emph{in the space $\calV^{2\rho}$} on $I$.

Notice that $(u,I)$ is a mild solution of \eqref{eq:Sell_You_47-1} if and only if $v(t) := u(t)$ is a mild solution of the \emph{linear} inhomogeneous problem,
$$
\dot v(t) + \cA v(t) = \cF(t,u(t)), \quad\hbox{for } v(t_0) = u_0 \in \calV^{2\rho} \hbox{ and all } t \geq t_0 \geq 0.
$$
Consequently, \cite[Lemma 42.1]{Sell_You_2002} implies that a classical solution, or a strong solution, if it exists, must be a mild solution. We begin by recalling the following local existence and uniqueness result for mild solutions of \eqref{eq:Sell_You_47-1} (compare \cite[Theorem 46.1]{Sell_You_2002}).

\begin{thm}[Existence and uniqueness of mild solutions to a nonlinear evolution equation in a Banach space]
\label{thm:Sell_You_lemma_47-1}
\cite[Lemma 47.1]{Sell_You_2002}
Assume the setup of the preceding paragraphs and that, for some $\beta \in [0, 1)$,
\begin{equation}
\label{eq:Sell_You_47-4}
\cF \in C_{\Lip}([0,\infty)\times \calV^{2\beta}; \cW).
\end{equation}
Given $b > 0$, there exists a positive constant,
$$
\tau = \tau\left(b, K_0, K_1, M_0, M_\beta, \beta\right),
$$
with the following significance.  For every $u_0 \in \calV^{2\beta}$ obeying $\|u_0\|_{\calV^{2\beta}} \leq b$ and every $t_0 \geq 0$, the initial value problem \eqref{eq:Sell_You_47-1} has a unique, mild solution in $\calV^{2\beta}$ on an interval $[t_0, t_0 + \tau)$, and
\begin{equation}
\label{eq:Sell_You_47-5}
u \in C([t_0, t_0+\tau); \calV^{2\beta}) \cap C_{\loc}^{0,\theta_1}([t_0, t_0+\tau); \calV^{2\alpha}) \cap C_{\loc}^{0,\theta}((t_0, t_0+\tau); \calV^{2r}),
\end{equation}
for all $\alpha$ and $r$ with $0 \leq \alpha < r$ and $0 \leq r < 1$, where $\theta_1 > 0$ and $\theta > 0$.
\end{thm}

Lemma \ref{lem:Sell_You_lemma_47-1_tau_formula} below provides an explicit formula for $\tau$ and an \apriori estimate for $u$ in Theorem \ref{thm:Sell_You_lemma_47-1} in terms of known constants.

By imposing more explicit polynomial growth and Lipschitz conditions on the nonlinearity, $\sF$, we obtain more a precise lower bound, $\tau$, on the lifetime of the mild solution, $u$, to \eqref{eq:Sell_You_47-1} as well as a more precise \apriori estimate for $u$ on the interval $[t_0, t_0 + \tau]$. In particular, we replace \eqref{eq:Sell_You_46-7} and \eqref{eq:Sell_You_46-8} with following more precise growth and Lipschitz conditions, for some $n \geq 1$ and positive constants $\kappa_0, \kappa_1$:
\begin{align}
\label{eq:Sell_You_46-7_polynomial_nonlinearity}
\|\cF(t, v)\|_\cW  &\leq \kappa_0\left(1 + \|v\|_{\calV^{2\beta}}^n\right),
\quad \forall\, t\geq 0 \hbox{ and } v \in \calV^{2\beta},
\\
\label{eq:Sell_You_46-8_polynomial_nonlinearity}
\|\cF(t, v_1) - \cF(t, v_2)\|_\cW  &\leq \kappa_1\left(1 + \|v_1\|_{\calV^{2\beta}}^{n-1} + \|v_2\|_{\calV^{2\beta}}^{n-1}\right)
\|v_1 - v_2\|_{\calV^{2\beta}},
\\
&\notag\qquad \forall\, t\geq 0 \hbox{ and } v_1, v_2 \in \calV^{2\beta}.
\end{align}
In the case of the nonlinearity defined by the Yang-Mills heat equation
\eqref{eq:Yang-Mills_heat_equation_as_perturbation_rough_Laplacian_plus_one_heat_equation},
we have $n=3$ but in general $n$ need not be an integer.

This yields the following improvement on Theorem \ref{thm:Sell_You_lemma_47-1}; we include a detailed proof, since the proof that $u$ in Theorem \ref{thm:Sell_You_lemma_47-1} belongs to $C([t_0, t_0+\tau); \calV^{2\beta})$ is omitted in \cite{Sell_You_2002}.

\begin{thm}[Existence and uniqueness of mild solutions to a nonlinear evolution equation in a Banach space]
\label{thm:Sell_You_lemma_47-1_polynomial_nonlinearity}
Assume the setup
of the paragraphs preceding Theorem \ref{thm:Sell_You_lemma_47-1} and that, for some $\beta \in [0, 1)$ and $n \geq 1$, the nonlinearity,
\begin{equation}
\label{eq:Sell_You_47-4_polynomial_nonlinearity}
\cF \in C_{\Lip}([0,\infty)\times \calV^{2\beta}; \cW),
\end{equation}
obeys \eqref{eq:Sell_You_46-7_polynomial_nonlinearity} and \eqref{eq:Sell_You_46-8_polynomial_nonlinearity}.
Given $b > 0$, there exists a positive
constant\footnote{This constant $\tau$ is given explicitly via \eqref{eq:Sell_You_47-6_polynomial_nonlinearity}.},
$$
\tau = \tau\left(b, M_0, M_\beta, n, \beta, \kappa_0, \kappa_1\right),
$$
with the following significance.  For every $u_0 \in \calV^{2\beta}$ obeying $\|u_0\|_{\calV^{2\beta}} \leq b$ and every $t_0 \geq 0$, the initial value problem \eqref{eq:Sell_You_47-1} has a unique, mild solution in $\calV^{2\beta}$ on an interval $[t_0, t_0 + \tau)$, which obeys
\begin{equation}
\label{eq:Sell_You_47-5_polynomial_nonlinearity}
u \in C([t_0, t_0+\tau]; \calV^{2\beta}) \cap
C_{\loc}^{0,\theta_1}([t_0, t_0+\tau); \calV^{2\alpha}) \cap C_{\loc}^{0,\theta}((t_0, t_0+\tau); \calV^{2r}),
\end{equation}
for all $\alpha$ and $r$ with $0 \leq \alpha < r$ and $0 \leq r < 1$, where $\theta_1 > 0$ and $\theta > 0$.
Moreover, the solution $u$ obeys the \apriori estimate,
\begin{equation}
\label{eq:Sell_You_lemma_47-1_polynomial_nonlinearity_apriori_estimate}
\|u(t)\|_{\calV^{2\beta}}
\leq
M_0\|u_0\|_{\calV^{2\beta}}
+
\frac{2M_\beta \kappa_0}{1-\beta} \left(1 + M_0 b\right)^n (t-t_0)^{1-\beta},
\quad \forall\, t\in[t_0, t_0+\tau].
\end{equation}
If $v_0 \in \calV^{2\beta}$ obeys $\|v_0\|_{\calV^{2\beta}} \leq b$ and $v$ is the unique, mild solution to \eqref{eq:Sell_You_47-1} in $\calV^{2\beta}$ on $[t_0, t_0 + \tau]$ with $v(0) = v_0$ and satisfying \eqref{eq:Sell_You_47-5_polynomial_nonlinearity}, then
\begin{equation}
\label{eq:Sell_You_lemma_47-1_polynomial_nonlinearity_continuity_with_respect_to_initial_data}
\sup_{t\in [t_0, t_0+\tau]}\|u(t)-v(t)\|_{\calV^{2\beta}} \leq 2M_0\|u_0 - v_0\|_{\calV^{2\beta}}.
\end{equation}
\end{thm}

\begin{proof}
We employ a contraction mapping argument based on the Banach space,
\begin{equation}
\label{eq:Sell_You_page_234_Banach_space_polynomial_nonlinearity}
\begin{gathered}
\fV := C([t_0, t_0+\tau]; \calV^{2\beta}) \quad\hbox{with norm}
\\
\|v\|_\fV \equiv \|v\|_{C_\gamma([t_0, t_0+\tau]; \calV^{2\beta})}
:= \sup_{t\in[t_0, t_0+\tau]}\|v(t)\|_{\calV^{2\beta}},
\end{gathered}
\end{equation}
just like Sell and You in their proof of \cite[Lemma 47.1]{Sell_You_2002} (recalled in our monograph as Theorem \ref{thm:Sell_You_lemma_47-1}). Remember that $\|v\|_{\calV^{2\beta}} \equiv \|\cA^\beta v\|_\cW$. We set
\begin{equation}
\label{eq:Sell_You_47-6_polynomial_nonlinearity_spatial_radius}
R :=  M_0b + 1,
\end{equation}
and, noting that $\beta < 1$ by hypothesis, choose $\tau > 0$ to be the largest constant such that
\begin{subequations}
\label{eq:Sell_You_47-6_polynomial_nonlinearity}
\begin{gather}
\label{eq:Sell_You_47-6_polynomial_nonlinearity_bound}
M_\beta \kappa_0 \left(1 + R^n \right) \frac{\tau^{1-\beta}}{1 - \beta} \leq 1,
\\
\label{eq:Sell_You_47-6_polynomial_nonlinearity_Lipschitz}
M_\beta \kappa_1 \left(1 + 2R^{n-1}\right) \frac{\tau^{1-\beta}}{1 - \beta} \leq\frac{1}{2}.
\end{gather}
\end{subequations}
We define $\fF := \{v \in \fV: \|v\|_\fV \leq R\}$, a closed ball in $\fV$. We shall seek a fixed point of a map $\sT:\fF \to \fF$ defined by
\begin{equation}
\label{eq:Sell_You_page_233_contraction_map_definition}
\hat u(t) \equiv \sT u(t) := e^{-\cA t}u_0 + \int_{t_0}^t e^{-\cA (t - s)} \cF(s, u(s))\, ds,
\quad \forall\, t \in [t_0, t_0 + \tau],
\end{equation}
regarding which we need to show that
\begin{inparaenum}[\itshape a\upshape)]
\item $\|\hat u(t)\|_{\calV^{2\beta}} \leq R$, for all $t \in [t_0, t_0 + \tau]$, and
\item $\hat u(t)$ is continuous with respect to $t \in [t_0, t_0 + \tau]$, so $\hat u \in \fF$, and
\item $\sT$ is a contraction map.
\end{inparaenum}
It is convenient to divide the proof of existence and uniqueness of $u$ into three corresponding steps.

\setcounter{step}{0}
\begin{step}[Boundedness of $\hat u(t)$ in $\calV^{2\beta}$ for $t_0 \leq t \leq t_0+\tau$]
\label{step:Sell_You_lemma_47-1_polynomial_nonlinearity_proof_of_boundedness}
Let $u \in \fF$ and define $\hat u$ by \eqref{eq:Sell_You_page_233_contraction_map_definition}. We calculate that, for $t_0 \leq t \leq t_0+\tau$,
\begin{align*}
\|\cA^\beta \hat u(t)\|_\cW
&\leq \|e^{-\cA (t-t_0)}\cA^\beta u_0\|_\cW
+ \int_{t_0}^t \|e^{-\cA(t-s)} \cA^\beta \cF(s, u(s))\|_\cW \, ds
\\
& = M_0 e^{-a(t-t_0)}\|u_0\|_{\calV^{2\beta}}
+ M_\beta \int_{t_0}^t (t-s)^{-\beta} e^{-a(t-s)}\|\cF(s, u(s))\|_\cW \, ds
\quad\hbox{(by \eqref{eq:Sell_You_37-11})}
\\
& \leq M_0\|u_0\|_{\calV^{2\beta}}
+ M_\beta \kappa_0 \int_{t_0}^t (t-s)^{-\beta} \left(1 + \|u(s)\|_{\calV^{2\beta}}^n \right)\, ds
\quad\hbox{(by \eqref{eq:Sell_You_46-7_polynomial_nonlinearity} and $a\geq 0$)}
\\
& \leq M_0 \|u_0\|_{\calV^{2\beta}}
+ M_\beta \kappa_0 \left(1 + \|u\|_\fV^n \right)
\int_{t_0}^t (t-s)^{-\beta}\,ds
\quad\hbox{(by \eqref{eq:Sell_You_page_234_Banach_space_polynomial_nonlinearity})}
\\
& = M_0 \|u_0\|_{\calV^{2\beta}}
+ M_\beta \kappa_0 \left(1 + \|u\|_\fV^n \right) \frac{(t-t_0)^{1-\beta}}{1-\beta}
\quad\hbox{(as $\beta<1$).}
\end{align*}
Therefore, because $\|u\|_\fV \leq R$ for $u \in \fF$,
\begin{equation}
\label{eq:Sell_You_lemma_47-1_polynomial_nonlinearity_apriori_estimate_preliminary}
\|\cA^\beta \hat u(t)\|_\cW
\leq M_0 \|u_0\|_{\calV^{2\beta}} + \frac{M_\beta \kappa_0}{1-\beta}(1 + R^n) (t-t_0)^{1-\beta},
\quad\forall\, t \in [t_0, t_0 + \tau].
\end{equation}
By definition of $R$ in \eqref{eq:Sell_You_47-6_polynomial_nonlinearity_spatial_radius} and as $\tau$ obeys \eqref{eq:Sell_You_47-6_polynomial_nonlinearity_bound}, we obtain
$$
\|\cA^\beta \hat u(t)\|_\cW \leq R, \quad\forall\, t \in [t_0, t_0 + \tau],
$$
and hence $\hat u(t) \in \bar B := \{v \in \calV^{2\beta}: \|v\|_{\calV^{2\beta}} \leq R\}$ for all $t \in [t_0, t_0 + \tau]$. This completes Step \ref{step:Sell_You_lemma_47-1_polynomial_nonlinearity_proof_of_boundedness}.
\end{step}

\begin{step}[Continuity of $\hat u(t)$ for $t_0 \leq t \leq t_0 + \tau$]
\label{step:Sell_You_lemma_47-1_polynomial_nonlinearity_proof_of_continuity}
The proof that $[t_0, t_0 + \tau] \ni t \mapsto \hat u(t) \in \calV^{2\beta}$ is a continuous mapping is similar to the proof of the corresponding fact in \cite[Theorem 46.1]{Sell_You_2002}. We note that
$$
\cA^\beta \hat u(t + h) - \cA^\beta \hat u(t) = E_1(t,h) + E_2(t,h) + E_3(t,h),
$$
where
\begin{align*}
E_1(t,h) &:= e^{-\cA(t+h-t_0)} \cA^\beta u_0 - e^{-\cA(t-t_0)} \cA^\beta u_0,
\\
E_2(t,h) &:=
\begin{cases}
\displaystyle
\int_{t_0}^t \left(\cA^\beta e^{-\cA(t+h-s)} \cF(s, u(s)) - \cA^\beta e^{-\cA(t-s)} \cF(s, u(s))\right)\, ds,
&\hbox{for } h\geq 0,
\\
\displaystyle
\int_{t_0}^{t+h} \left(\cA^\beta e^{-\cA(t+h-s)} \cF(s, u(s)) - \cA^\beta e^{-\cA(t-s)} \cF(s, u(s))\right)\, ds,
&\hbox{for } h < 0,
\end{cases}
\\
E_3(t,h) &:=
\begin{cases}
\displaystyle\int_t^{t+h} \cA^\beta e^{-\cA(t+h-s)} \cF(s, u(s))\, ds,
&\hbox{for } h\geq 0,
\\
-\displaystyle\int_{t+h}^t \cA^\beta e^{-\cA(t-s)} \cF(s, u(s))\, ds,
&\hbox{for } h < 0.
\end{cases}
\end{align*}
Without loss of generality, we assume that $|h| \leq 1$. If $t = t_0$, we further restrict $h$ to satisfy $0 \leq h \leq 1$. For $t > t_0$, we have
$$
\|E_1(t,h)\|_\cW \leq C_\beta |h| (t-t_0)^{-1}\|u_0\|_{\calV^{2\beta}}
\quad\hbox{(by \cite[Equation (37.13)]{Sell_You_2002})}.
$$
Thus, $\|E_1(t,h)\|_\cW \to 0$, as $h \to 0$, when $t>t_0$. For the case $t = t_0$, we have $h \geq 0$ and observe that $\cA^\beta u_0 \in \cW$ and
$$
\|E_1(t,h)\|_\cW = \|e^{-\cA h}\cA^\beta u_0 - \cA^\beta u_0\|_\cW \to 0 \quad\hbox{as } h\downarrow 0
\quad\hbox{by \cite[Lemma 31.2]{Sell_You_2002}.}
$$
For the term $E_2$, our calculations in Step \ref{step:Sell_You_lemma_47-1_polynomial_nonlinearity_proof_of_boundedness} for the estimate of $\|\cA^\beta \hat u(t)\|_\cW$ show that
$$
\|E_2(t,h)\|_\cW
\leq
\begin{cases}
\displaystyle
M_\beta\int_{t_0}^t (t-s)^{-\beta} \|e^{-\cA h} \cF(s, u(s)) - \cF(s, u(s))\|_\cW \, ds,
&\hbox{for } h \geq 0,
\\
\displaystyle
M_\beta\int_{t_0}^t (t+h-s)^{-\beta} \|\cF(s, u(s)) - e^{\cA h} \cF(s, u(s))\|_\cW \, ds,
&\hbox{for } h < 0.
\end{cases}
$$
Again, \cite[Lemma 31.2]{Sell_You_2002} implies that
$$
\|e^{-\cA |h|} \cF(s, u(s)) - \cF(s, u(s))\|_\cW \to 0 \quad\hbox{as } h\to 0
\quad\hbox{by \cite[Lemma 31.2]{Sell_You_2002}.}
$$
Therefore, $\|E_2(t,h)\|_\cW \to 0$, as $h \to 0$, for each $t \in [t_0, t_0 + \tau]$ and the Lebesgue Dominated Convergence Theorem for Bochner integrals \cite[Theorem C.4 (4)]{Sell_You_2002}. By \eqref{eq:Sell_You_46-7_polynomial_nonlinearity} and our calculations in Step \ref{step:Sell_You_lemma_47-1_polynomial_nonlinearity_proof_of_boundedness} for the estimate of $\|\cA^\beta \hat u(t)\|_\cW$, the term $E_3$ satisfies, for all $t \in [t_0, t_0 + \tau]$,
$$
\|E_3(t,h)\|_\cW
\leq
\begin{cases}
\displaystyle
M_0M_\beta\kappa_0 \left(1 + \sup_{s\in[t,t+h]}\|u(s)\|_{\calV^{2\beta}}^n \right) \frac{|h|^{1-\beta}}{1-\beta},
&\hbox{for } h \geq 0,
\\
\displaystyle
M_\beta\kappa_0 \left(1 + \sup_{s\in[t,t+h]}\|u(s)\|_{\calV^{2\beta}}^n \right) \frac{|h|^{1-\beta}}{1-\beta},
&\hbox{for } h < 0,
\end{cases}
$$
so $\|E_3(t,h)\|_\cW \to 0$, as $h \to 0$, for all $t \in [t_0, t_0 + \tau]$, recalling that $\beta < 1$. (In the case $h\geq 0$, we also use the fact that $\|e^{-\cA h}\|_{\sL(\cW)} \leq M_0$ by \cite[Equation (37.11)]{Sell_You_2002} with $r=0$, where $M_0\geq 1$ .) Thus, $\hat u \in C([t_0, t_0 + \tau]; \calV^{2\beta}) \equiv \fV$ and this completes Step \ref{step:Sell_You_lemma_47-1_polynomial_nonlinearity_proof_of_continuity}.
\end{step}

\begin{step}[Contraction mapping property of $\sT$]
\label{step:Sell_You_lemma_47-1_polynomial_nonlinearity_proof_of_contraction property}
Next we show that for $\tau$ obeying \eqref{eq:Sell_You_47-6_polynomial_nonlinearity_Lipschitz}, the mapping $\sT$ is a contraction on $\fF$ with contraction coefficient less than or equal to $1/2$. Indeed, let $u_1, u_2 \in \fF$ and set $\hat u_i := \sT u_i$ for $i = 1,2$. Then, for $t_0 \leq t \leq t_0 + \tau$, we have
\begin{align*}
{}& \|\cA^\beta (\hat u_1(t) - \hat u_2(t))\|_\cW
\\
&\quad \leq \int_{t_0}^t \|e^{-\cA(t-s)} \cA^\beta (\cF(s, u_1(s)) - \cF(s, u_2(s)))\|_\cW \, ds
\\
&\quad\leq M_\beta \int_{t_0}^t (t-s)^{-\beta} e^{-a(t-s)}\|\cF(s, u_1(s)) - \cF(s, u_2(s))\|_\cW \, ds
\quad\hbox{(by \eqref{eq:Sell_You_37-11})}
\\
&\quad \leq M_\beta \kappa_1 \int_{t_0}^t (t-s)^{-\beta}
\left(1 + \|u_1(s)\|_{\calV^{2\beta}}^{n-1} + \|u_2(s)\|_{\calV^{2\beta}}^{n-1} \right)\|u_1(s) - u_2(s)\|_{\calV^{2\beta}} \,ds
\quad\hbox{(by \eqref{eq:Sell_You_46-8_polynomial_nonlinearity})}
\\
&\quad \leq M_\beta \kappa_1 \left(1 + \|u_1\|_\fV^{n-1} + \|u_2\|_\fV^{n-1} \right)\|u_1 - u_2\|_\fV
\int_{t_0}^t (t-s)^{-\beta}\,ds
\quad\hbox{(by \eqref{eq:Sell_You_page_234_Banach_space_polynomial_nonlinearity})}
\\
&\quad \leq M_\beta \kappa_1 \left(1 + 2R^{n-1}\right)\|u_1 - u_2\|_\fV
\frac{\tau^{1-\beta}}{1-\beta}
\quad\hbox{(as $\beta<1$ and $\|u_i\|_\fV \leq R$ for $u_i \in \fF$).}
\end{align*}
Therefore, by \eqref{eq:Sell_You_page_234_Banach_space_polynomial_nonlinearity} and \eqref{eq:Sell_You_47-6_polynomial_nonlinearity_Lipschitz} we obtain
\begin{equation}
\label{eq:Sell_You_lemma_47-1_polynomial_nonlinearity_contraction_mapping_property}
\|\hat u_1 - \hat u_2\|_\fV \leq \frac{1}{2}\|u_1 - u_2\|_\fV,
\end{equation}
and as a result, the map $\sT$ has a unique fixed point $u \in \fF$. This fixed point is the mild solution of \eqref{eq:Sell_You_47-1} on $[t_0, t_0+\tau]$, and because of the contraction property, this solution is uniquely determined. This completes Step \ref{step:Sell_You_lemma_47-1_polynomial_nonlinearity_proof_of_contraction property}.
\end{step}

\begin{step}[\Apriori and continuity estimates]
\label{step:Sell_You_lemma_47-1_polynomial_nonlinearity_apriori_and_continuity_estimates}
From \eqref{eq:Sell_You_lemma_47-1_polynomial_nonlinearity_apriori_estimate_preliminary}, the definition \eqref{eq:Sell_You_47-6_polynomial_nonlinearity_spatial_radius} of $R$, and the fact that $\hat u = u$, we obtain \eqref{eq:Sell_You_lemma_47-1_polynomial_nonlinearity_apriori_estimate}. If $u, v$ are the unique mild solutions in $\calV^{2\beta}$ on $[t_0, t_0 + \tau]$ to \eqref{eq:Sell_You_47-1} with initial data $u_0, v_0 \in \calV^{2\beta}$, respectively, then a slight modification of the derivation of the contraction mapping property \eqref{eq:Sell_You_lemma_47-1_polynomial_nonlinearity_contraction_mapping_property} (simply observe that $u_0-v_0$ may be nonzero) for $\sT:\fF \to \fF$ now yields, noting that $\hat u = \sT u = u$ and $\hat v = \sT v = v$ and recalling the definition of $\fV$ in \eqref{eq:Sell_You_page_234_Banach_space_polynomial_nonlinearity},
$$
\|u - v\|_\fV \leq M_0\|u_0 - v_0\|_{\calV^{2\beta}} + \frac{1}{2}\|u - v\|_\fV,
$$
and rearrangement gives the continuity estimate \eqref{eq:Sell_You_lemma_47-1_polynomial_nonlinearity_continuity_with_respect_to_initial_data}.
\end{step}

\begin{step}[Regularity]
\label{step:Sell_You_lemma_47-1_polynomial_nonlinearity_regularity}
Then \cite[Lemma 42.7]{Sell_You_2002},
with $p = \infty$, implies that \eqref{eq:Sell_You_47-5_polynomial_nonlinearity} holds.
\end{step}

This completes the proof of Theorem \ref{thm:Sell_You_lemma_47-1_polynomial_nonlinearity}.
\end{proof}

In the proof of \cite[Lemma 47.1]{Sell_You_2002} provided by Sell and You, they give an explicit formula for the positive constant, $\tau$, in the statement of Theorem \ref{thm:Sell_You_lemma_47-1} and an \apriori estimate for the solution, $u$. We record this formula and \apriori estimate in Lemma \ref{lem:Sell_You_lemma_47-1_tau_formula}, but in order to introduce them and define all the relevant constants, we require a digression on semigroups and, for convenience, we shall closely follow
\cite[Section 36.2]{Sell_You_2002}. A special class of $C^0$-semigroups, namely the \emph{analytic semigroups}, plays a fundamental role in the study of the dynamics of infinite-dimensional systems. There are two principal reasons why analytic semigroups are important in the study of systems of nonlinear parabolic partial differential equations. The first reason is owing to the good information one has on the behavior of solutions as time $t \downarrow 0$. In particular, under reasonable conditions, a $C^0$-semigroup, $(e^{\cA t}, \cA)$, on a Banach space, $\cW$, is an analytic semigroup if and only if there are constants $M_0 \geq 1$ and $M_1 > 0$ and a constant $a \in \RR$ such that
$$
\|e^{\cA t}\| \leq M_0 e^{-a t}\|w\| \quad\hbox{and}\quad \|\cA e^{\cA t}\| \leq M_1 t^{-1} e^{-a t}\|w\|,
$$
for all $w \in \cW$ and $t > 0$. Since the nonlinear evolutionary equations under discussion are perturbations of a linear equation, by having good information about the behavior of solutions of the linear problem, as $t \downarrow 0$, one can
use this to study the initial value problem for related nonlinear problems. The second reason for the importance of analytic semigroups is given in the Fundamental Theorem on Sectorial Operators (Theorem \ref{thm:Sell_You_37-5}).

We now recall a fundamental characterization of the generator of an analytic semigroup; for additional details on analytic semigroups, see Hille and Phillips \cite[Chapter 17]{Hille_Phillips}, Kato \cite[Section 9.1.6]{Kato}, or Yosida \cite[Section 9.9]{Yosida}.
In \cite[Theorem 1.7.7]{Pazy_1983}, Pazy provides a sufficient condition for an operator to be the generator of a uniformly bounded $C^0$ semigroup, while \cite[Theorem 2.5.2]{Pazy_1983} provides several equivalent characterizations for the generator of an analytic semigroup. The characterization in Theorem \ref{thm:Renardy_Rogers_12-31} below given in Renardy and Rogers \cite{Renardy_Rogers_2004} will be most convenient for our application and is a stronger result than its nearest equivalent \cite[Theorem 2.5.2 (c)]{Pazy_1983} provided by Pazy or \cite[Theorem 36.2 (3) and Equation (36.3)]{Sell_You_2002} provided by Sell and You.

\begin{thm}[Generator of an analytic semigroup]
\label{thm:Renardy_Rogers_12-31}
\cite[Theorem 12.31]{Renardy_Rogers_2004}
A closed, densely defined operator $\cA$ on a Banach space $\cW$ is the generator of an analytic semigroup if and only if there exists $\lambda_0 \in \RR$ such that the resolvent set contains a half-plane,
$$
\rho(\cA) \supset \{z \in \CC: \Real z > \lambda_0\} \equiv \Delta_{\pi/2}(\lambda_0),
$$
and, moreover, there is a positive constant $C$ such that
\begin{equation}
\label{eq:Renardy_Rogers_2004_12-65}
\|\sR(\lambda, \cA)\|_{\sL(\cW)} \leq \frac{C}{|\lambda - \lambda_0|},
\quad \forall\, \lambda \in \Delta_{\pi/2}(\lambda_0).
\end{equation}
If this is the case, then actually the resolvent set contains a sector,
$$
\rho(\cA) \supset \Delta_{\pi/2 + \eps}(\lambda_0),
$$
for some $\eps \in (0, \pi/2)$, and a resolvent estimate analogous to \eqref{eq:Renardy_Rogers_2004_12-65} holds in this sector. Moreover, the semigroup is represented by
\begin{equation}
\label{eq:Renardy_Rogers_2004_12-66}
e^{\cA t} = \frac{1}{2\pi i}\oint_\Gamma e^{\lambda t}\sR(\lambda, -\cA) \, d\lambda, \quad t \geq 0,
\end{equation}
where $\Gamma \subset \CC$ is any curve from $e^{-i\phi}\infty$ to $e^{i\phi}\infty$ such that $\Gamma$ lies entirely in the closed sector $\bar\Delta_\phi(\lambda_0) = \{\lambda \in \CC: \arg(\lambda - \lambda_0) \leq \phi\}$, where $\phi$ is any angle such that $\pi/2 < \phi < \pi/2 + \eps$.
\end{thm}

\begin{rmk}[Sectorial operators and generators of analytic semigroups]
\label{rmk:Renardy_Rogers_2004_12-31}
If $\cA$ is a sectorial operator on a Banach space $\cW$ with sector $\Delta_{\pi/2}(\lambda_0)$ for some $\lambda_0 \in \RR$ in the sense of Definition \ref{defn:Sell_You_page_78_definition_of_sectorial_operator}, then $\cA$ satisfies the hypotheses of Theorem \ref{thm:Renardy_Rogers_12-31}.
\end{rmk}

Recall from the Fundamental Theorem on Sectorial Operators (Theorem \ref{thm:Sell_You_37-5}) that for any $r \geq 0$, there is a positive constant, $M_r$, such that the inequality \eqref{eq:Sell_You_37-11} holds, where $e^{-\cA t}$ is the analytic semigroup generated  by $-\cA$ and $a$ is given by \eqref{eq:Sell_You_36-2}. This completes our digression on semigroups and we can now review the promised \apriori estimate and formula for the existence time of the mild solution, $u$, provided in the proof (see \cite[pp. 234--235]{Sell_You_2002}) of Theorem \ref{thm:Sell_You_lemma_47-1}.

\begin{lem}[An \apriori estimate and a positive lower bound for the existence time of a mild solution]
\label{lem:Sell_You_lemma_47-1_tau_formula}
Assume the hypotheses of Theorem \ref{thm:Sell_You_lemma_47-1}. Then the solution, $u$, obeys
\begin{equation}
\label{eq:Sell_You_page_234_contraction_mapping_apriori_bound_solution}
\|u(t)\|_{\calV^{2\beta}} \leq M_0b + \frac{M_\beta K_0}{(1-\beta)}(t-t_0)^{1-\beta},
\quad\forall\, t \in [t_0, t_0+\tau],
\end{equation}
where $M_0$ is given by \eqref{eq:Sell_You_37-11}, and $\tau$ is given by
\begin{equation}
\label{eq:Sell_You_47-6}
\tau^{1-\beta} = \min\left\{\frac{M_0(1-\beta)}{2M_\beta K_0}, \frac{1-\beta}{2M_\beta K_1}\right\},
\end{equation}
where $K_0=K_0(B)$ is given by \eqref{eq:Sell_You_46-7} with
\begin{equation}
\label{eq:Sell_You_page_234_ball}
B := \{v \in \calV^{2\beta}: \|v\|_{\calV^{2\beta}} < M_0(b+1)\},
\end{equation}
and $K_1=K_1(B)$ is given by \eqref{eq:Sell_You_46-8}, and $M_\beta$ is given by \eqref{eq:Sell_You_37-11}.
\end{lem}

\begin{rmk}[On the significance of the formula \eqref{eq:Sell_You_47-6} for $\tau$]
\label{rmk:Sell_You_lemma_47-1_tau_formula}
One can gain further intuition for the meaning of \eqref{eq:Sell_You_47-6} by noting that the positive constants, $K_0, K_1$, pertain to properties of the nonlinearity, $\cF(t,v)$, whereas the positive constants, $M_0, M_\beta$, pertain to properties of the linear operator, $\cA$. If in fact $\cF \equiv 0$, then one could take $K_0 = K_1 = 0$ in \eqref{eq:Sell_You_47-6}, giving $\tau = \infty$, just as we would expect from linear theory. It is important to also note that $\tau$ depends on the upper bound $b$ on the norm of initial data, $u_0 \in \calV^{2\beta}$, through the dependence of the constants, $K_0(B)$ and $K_1(B)$, on the radius, $M_0(b + 1)$, of the ball, $B$.
\end{rmk}

\subsection{Strong solutions}
\label{subsec:Sell_You_4-7-2}
Suppose that $\calV$ and $\cW$ are two Banach spaces and $\theta \in (0, 1]$. Recall from \cite[Equation (46.6) or p. 658]{Sell_You_2002} that $\cF \in C_{\Lip;\theta}([0,\infty)\times \calV; \cW)$ if $\cF \in C_{\Lip}([0,\infty)\times \calV; \cW)$ (in the sense of \eqref{eq:Sell_You_46-7} and \eqref{eq:Sell_You_46-8}) and, for every bounded set $B\subset \calV$ and compact set $J\subset [0,\infty)$, there is a nonnegative constant, $K_2 = K_2(B,J)$, such that
\begin{equation}
\label{eq:Sell_You_46-6}
\|\cF(t_1,v_1) - \cF(t_2,v_2)\|_\cW \leq K_2\left(\|v_1-v_2\|_\calV + |t_1-t_2|^\theta\right), \quad\forall\, t_1, t_2 \in J \hbox{ and } v_1, v_2 \in B.
\end{equation}
By imposing an additional regularity condition of this kind on the nonlinearity $\cF(t,u)$, one can show that the mild solution is a strong solution (compare \cite[Theorem 46.2]{Sell_You_2002}). Recall that $\calV^2 = \sD(\cA)$.

\begin{thm}
\label{thm:Sell_You_lemma_47-2}
\cite[Lemma 47.2]{Sell_You_2002}
Assume the setup in Section \ref{subsec:Sell_You_4-7-1} and that, for some $\beta \in [0, 1)$ and $\theta \in (0, 1]$,
\begin{equation}
\label{eq:Sell_You_lemma_47-2_F_C_Lipschitz_theta}
\cF \in C_{\Lip;\theta}([0,\infty)\times \calV^{2\beta}; \cW).
\end{equation}
If $u_0 \in \calV^{2\beta}$ and $u$ is a mild solution of Equation \eqref{eq:Sell_You_47-1} in $\calV^{2\beta}$ on an interval $[0, T)$ for some $T > 0$, then $u$ is a strong solution in $\calV^{2\beta}$ on the interval $[0, T)$, and it satisfies
\begin{equation}
\label{eq:Sell_You_47-7}
u \in C([0, T); \calV^{2\alpha}) \cap C_{\loc}^{0,1-r}((0, T); \calV^{2r}) \cap C((0, T); \sD(\cA)),
\end{equation}
for all $\alpha$ and $r$ with $0 \leq \alpha \leq \beta$ and $0 \leq r < 1$.
\end{thm}

\begin{rmk}[Regularity and the lower bound for the lifetime of the solution]
The original statement of \cite[Lemma 47.2]{Sell_You_2002} only asserts that $T>0$ exists but it is clear from its proof in \cite{Sell_You_2002} that the stronger conclusion follows. This observation explains our phrasing of Theorem \ref{thm:Sell_You_lemma_47-2}.
\end{rmk}

\subsection{Maximally defined solutions}
\label{subsec:Sell_You_4-7-3}
Let $(u_1, I_1)$ and $(u_2,I_2)$ be two mild solutions of \eqref{eq:Sell_You_47-1}, where $I_i = [t_0, t_0 + \tau_i)$, for $i = 1,2$, and $\tau_1 \leq \tau_2$. Owing to the uniqueness of solutions, one must have $u_1(t) = u_2(t)$, for $t \in I_1$. Hence, $(u_2, I_2)$ is an \emph{extension} of $(u_1, I_1)$. When $\tau_1 < \tau_2$, calls $(u_2, I_2)$ a \emph{proper extension} of $(u_1, I_1)$. A solution $(u, I)$ of \eqref{eq:Sell_You_47-1} is said to be a \emph{maximally defined solution} if $(u, I)$ has no proper extension. We recall that (compare \cite[Theorem 46.3]{Sell_You_2002})

\begin{thm}
\label{thm:Sell_You_lemma_47-4}
\cite[Lemma 47.4]{Sell_You_2002}
Assume the setup in Section \ref{subsec:Sell_You_4-7-1} and let $\cF$ obey \eqref{eq:Sell_You_47-4} for some $\beta \in [0,1)$. Then for every $u_0 \in \calV^{2\beta}$ and $t_0 \geq 0$, there is a unique, maximally defined, mild solution $(u, I)$ of \eqref{eq:Sell_You_47-1} in $\calV^{2\beta}$, where $I = [t_0, t_0 + T)$. Furthermore, either $T = \infty$, or
\begin{equation}
\label{eq:Sell_You_47-8}
\lim_{t\uparrow T} \|\cA^\beta u(t)\|_\cW = \infty.
\end{equation}
\end{thm}


Recall that $\|\cA^\beta v\|_\cW = \|v\|_{\calV^{2\beta}}$.

\subsection{Continuous dependence of solutions}
\label{subsec:Sell_You_4-7-4}
We continue the setup in Section \ref{subsec:Sell_You_4-7-1}. For $\beta \in [0,1)$ and any $\cF \in C_{\Lip}([0,\infty)\times \calV^{2\beta};\cW)$ and $u_0 \in \calV^{2\beta}$, we let $\phi(u_0, \cF, t)$ denote the maximally defined, mild solution of \eqref{eq:Sell_You_47-1} in $\calV^{2\beta}$, with $t_0 = 0$, that satisfies $\phi(u_0, \cF, 0) = u_0$, and let $I = [0, T)$ denote the interval of definition of $\phi(u_0, \cF, t)$, where $T = T(u_0, \cF)$ satisfies $0 < T \leq \infty$. Next define
\begin{equation}
\label{eq:Sell_You_47-9}
\Xi := \{(u_0, \cF, t) \in \calV^{2\beta} \times C_{\Lip}([0,\infty)\times \calV^{2\beta};\cW) \times [0,\infty) : 0 \leq t < T(u_0, \cF)\}.
\end{equation}
We equip $\Xi$ with the topologies $\sT_A^0$ or $\sT_{bo} = \sT_{bo}^0$ \cite[p. 222]{Sell_You_2002}. We recall the (compare \cite[Theorem 46.4]{Sell_You_2002})

\begin{thm}
\label{thm:Sell_You_47-5}
\cite[Theorem 47.5]{Sell_You_2002}
Assume the setup in Section \ref{subsec:Sell_You_4-7-1}, let $\cF$ obey \eqref{eq:Sell_You_47-4} for some $\beta \in [0,1)$, and let $\Xi$ be given by \eqref{eq:Sell_You_47-9}. Then the following hold:
\begin{enumerate}
\item The mapping $(u_0, \cF, t) \mapsto \phi(u_0, \cF, t)$ of $(\Xi, \sT_A^0)$ or $(\Xi, \sT_{bo}^0)$ into $\calV^{2\beta}$ is continuous, and $\phi(u_0, \cF, t)$ is locally Lipschitz continuous in $\cF$ and $u_0$;
\item The set $\Xi$ is open in $\calV^{2\beta} \times C_{\Lip}([0,\infty)\times \calV^{2\beta};\cW) \times [0,\infty)$;
\item If $\tau \in [0,T(u_0,\cF))$ and $t \in [0,T(\phi(u_0, \cF, \tau)), \cF_\tau))$, then $\tau + t \in [0,T(u_0,\cF))$ and one has
\begin{equation}
\label{eq:Sell_You_47-10}
\phi(\phi(u_0, \cF, \tau), \cF_\tau, t) = \phi(u_0, \cF, \tau + t),
\end{equation}
where $\cF_\tau(u, t) := \cF(u, \tau + t)$. In particular, if $\cF \in C_{\Lip}( \calV^{2\beta};\cW)$ is time-independent, then $\cF_\tau = \cF$, for all $t\geq 0$, and
\begin{equation}
\label{eq:Sell_You_46-24}
\phi(\phi(u_0, \cF, \tau), \cF, t) = \phi(u_0, \cF, \tau + t).
\end{equation}
\end{enumerate}
\end{thm}

\subsection{A standard result for long time existence}
\label{subsec:Sell_You_4-7-6}
In order to highlight the \emph{nonstandard} nature of the question of long-time existence of a solution to the Yang-Mills heat equation
\eqref{eq:Yang-Mills_heat_equation_as_perturbation_rough_Laplacian_plus_one_heat_equation},
it is useful to recall a standard sufficient criterion for long time existence to \eqref{eq:Sell_You_47-1}. We shall later explain why simple criteria of this kind fail for the Yang-Mills heat equation \eqref{eq:Yang-Mills_heat_equation_as_perturbation_rough_Laplacian_plus_one_heat_equation}.

\begin{thm}
\label{thm:Sell_You_47-7}
\cite[Theorem 47.7]{Sell_You_2002}
Assume the setup in Section \ref{subsec:Sell_You_4-7-1}
and let $\cF$ obey \eqref{eq:Sell_You_47-4} for some $\beta \in [0,1)$. Then a sufficient condition for $T(u_0, \cF) = \infty$, for all $u_0 \in \calV^{2\beta}$, is that there are nonnegative constants $C_0$ and $C_1$ such that
\begin{equation}
\label{eq:Sell_You_47-15}
\|\cF(t,v)\|_\cW \leq C_0 + C_1\|\cA^\beta v\|_\cW, \quad\forall\, (t, v) \in [0,\infty) \times \calV^{2\beta}.
\end{equation}
\end{thm}

\subsection{Regularity in space and time}
\label{subsec:Sell_You_4-8-2}
For $\theta \in (0, 1]$ and $T > 0$ and, temporarily assuming that $\calV, \cW$ are Banach spaces, we recall that $C^{0,\theta}_{\loc}([0,T);\calV)$ is the vector subspace of functions $v \in C([0,T); \calV)$ such that, for each compact subset $J \subset [0, T)$, there is a nonnegative constant $K(J)$ with
$$
\|v(t_1) - v(t_2)\|_\calV \leq K|t_1 - t_2|, \quad \forall \, t_1, t_2 \in J,
$$
and the vector space $C^{0,\theta}_{\loc}((0,T);\calV)$ is similarly defined \cite[Appendix B.1 and p. 658]{Sell_You_2002}. Moreover, we recall that $C^1_F([0,T)\times \calV; \cW)$ is the vector space of continuously Fr\'echet-differentiable functions \cite[Appendix C.3 and p. 655]{Sell_You_2002}.

\begin{thm}
\label{thm:Sell_You_48-5}
\cite[Theorem 48.5]{Sell_You_2002}
Assume that Hypothesis \ref{hyp:Sell_You_4_standing_hypothesis_B} holds for a Hilbert space $\cH$ and, for some $\beta \in [0,1)$, that $\cF$ obeys
\begin{equation}
\label{eq:theorem_Sell_You_48-5_F_in_CLip_R_times_V2beta_to_H_and_C1Frechet_R_times_V2beta_to_H}
\cF \in C_{\Lip}([0,\infty)\times \calV^{2\beta}; \cH) \cap C^1_F([0,\infty)\times \calV^{2\beta}; \cH).
\end{equation}
Assume further that there is a $p \in [2,\infty]$ such that if $v : [0, T) \to \calV^{2\beta}$ is continuous and strongly differentiable for $0 < t < T$, then $f(t) := \cF(t, v(t))$ satisfies
\begin{equation}
\label{eq:theorem_Sell_You_48-5_f_in_C_0_T_to_H_and_W1_p_loc_0_T_to_Vminus_one}
f \in C([0, T); \cH) \cap W^{1,p}_{\loc}([0, T; \calV^{-1}).
\end{equation}
Let $\dot f = \partial_u \cF \dot u + \partial_t \cF$, where $\partial_t \cF$ and $\partial_u \cF$ are the (partial) Fr\'echet derivatives of $\cF$. Let $u_0 \in \sD(\cA) = \calV^2$ and define $v_0 := \cF(0, u_0) - \cA u_0 \in \cH$. Let $u = u(t) = \phi(u_0, \cF, t)$ denote the maximally defined mild solution of \eqref{eq:Sell_You_47-1} in $\calV^{2\beta}$ on the interval $[0, T)$, where $0 < T \leq \infty$. Then $u$ satisfies the following properties:
\begin{enumerate}
\item The function $u$ is a strong solution to \eqref{eq:Sell_You_47-1} in $\calV^{2\beta}$. Moreover,
for each $r \in [0, 1)$, there is a $\theta = \theta(r) > 0$ such that one has
$$
u \in C^1([0, T); \cH) \cap C([0,T); \calV^2) \cap C^{0, \theta}_{\loc}([0,T); \calV^{2r}).
$$

\item Set $g := \partial_u \cF \dot u + \partial_t \cF$. Then, for each $\alpha < 0$, the function $v := \dot u$ satisfies
$$
v(t) = e^{-t\cA}v_0 + \int_0^t e^{-(t-s)\cA} g(s)\,ds, \quad \forall\, t \geq 0.
$$
Also $v$ is a mild solution of $\dot v + \cA v = g$ in $\calV^\alpha$, and it is a weak solution in $\cH$ with
$$
v \in C([0, T); \cH) \cap L^2_{\loc}([0,T); \calV^1) \cap C^{0, \theta_1}_{\loc}([0,T); \calV^{2\alpha}),
$$
where $\theta_1 = \theta_1(\alpha) > 0$.

\item If instead, one has $u_0 \in \calV^{2\beta}$, then $u = u(t) = \phi(u_0, \cF, t)$ satisfies
\begin{equation}
\label{eq:Sell_You_48-3}
u \in C^1([0, T); \calV^{2\beta}) \cap C((0,T); \calV^{2(1+\beta)}) \cap C^{0, \theta_2}_{\loc}([0,T); \calV^{2r}),
\end{equation}
where $\theta_2 = \theta_2(r) > 0$, for $r \in [0, 1+\beta)$. Furthermore, for each such $r$, the function $u$ is a strong solution to \eqref{eq:Sell_You_47-1} in $\calV^{2r}$ on $[t_1, T)$, for any $t_1$ with $0<t_1<T$.
\end{enumerate}
\end{thm}

\subsection{Existence and uniqueness of mild solutions to a nonlinear evolution equation in a Banach space with initial data of minimal regularity I}
\label{subsec:Sell_You_lemma_4-7-1_Kozono_Maeda_Naito}
Theorem \ref{thm:Sell_You_lemma_47-1} and Theorem \ref{thm:Sell_You_lemma_47-1_polynomial_nonlinearity} provide useful existence and uniqueness results for mild solutions, in $\calV^{2\beta}$ on $[t_0, t_0+\tau]$, to the nonlinear evolution equation \eqref{eq:Sell_You_47-1}, but they lack the flexibility we would like to have with $u(t) \in \calV^{2\beta}$ for $t > 0$ but a more relaxed requirement on the regularity of the initial data than $u_0 \in \calV^{2\beta}$. For the linear evolution equation \eqref{eq:Sell_You_42-1}, we have seen a statement of this kind in Theorem \ref{thm:Sell_You_42-15} and in this subsection, we establish an analogue of this statement for \eqref{eq:Sell_You_47-1}.

The results \cite[Theorem 3.1 and Lemma 3.4]{Kozono_Maeda_Naito_1995} of Kozono, Maeda, and Naito, together with their proofs, suggest that the introduction of time-weighted function spaces will allow the greatest flexibility and we employ such ideas in our proof of Theorem \ref{thm:Kozono_Maeda_Naito_lemma_3-4}.

We shall need to replace \eqref{eq:Sell_You_46-7_polynomial_nonlinearity} and \eqref{eq:Sell_You_46-8_polynomial_nonlinearity} with following more refined growth and Lipschitz conditions, for some $n \geq 1$ (for the Yang-Mills heat equation
\eqref{eq:Yang-Mills_heat_equation_as_perturbation_rough_Laplacian_plus_one_heat_equation}, we have $n=3$) and $\alpha, \gamma \in [0,1)$, and positive constants $\kappa_2, \kappa_3$:
\begin{align}
\label{eq:Sell_You_46-7_Kozono_Maeda_Naito}
\|\cA^{-\gamma}\cF(t, v)\|_\cW  &\leq \kappa_2\left(1 + \|v\|_{\calV^{2\alpha}}^n\right),
\quad \forall\, t\geq 0 \hbox{ and } v \in \calV^{2\alpha},
\\
\label{eq:Sell_You_46-8_Kozono_Maeda_Naito}
\|\cA^{-\gamma}\cF(t, v_1) - \cF(t, v_2)\|_\cW
&\leq \kappa_3\left(1 + \|v_1\|_{\calV^{2\alpha}}^{n-1} + \|v_2\|_{\calV^{2\alpha}}^{n-1}\right) \|v_1 - v_2\|_{\calV^{2\alpha}},
\\
&\notag\qquad \forall\, t\geq 0 \hbox{ and } v_1, v_2 \in \calV^{2\alpha}.
\end{align}
Again, in the case of the nonlinearity defined by the Yang-Mills heat equation \eqref{eq:Yang-Mills_heat_equation_as_perturbation_rough_Laplacian_plus_one_heat_equation},
we have $n=3$ but in general $n$ need not be an integer. Our Theorem \ref{thm:Kozono_Maeda_Naito_lemma_3-4} extends \cite[Theorem 3.1 and Lemma 3.4]{Kozono_Maeda_Naito_1995}, who restrict their attention to the Yang-Mills heat equation \eqref{eq:Yang-Mills_heat_equation_as_perturbation_rough_Laplacian_plus_one_heat_equation},
taking $\cW = L^p(X; \Lambda^1\otimes \ad P)$ and $\calV = W^{2,p}_{A_1}(X; \Lambda^1\otimes \ad P)$ with $p=\dim X$. Our proof of Theorem \ref{thm:Kozono_Maeda_Naito_lemma_3-4} employs many of the ideas developed by Kozono, Maeda, and Naito in \cite{Kozono_Maeda_Naito_1995}, but we abstract their result and simplify their argument. We prove a stronger version of their result as Theorem \ref{thm:Kozono_Maeda_Naito_lemma_3-4_plus_uniqueness} in Section \ref{subsec:Sell_You_lemma_4-7-1_Kozono_Maeda_Naito_II}.

To prove existence and uniqueness of a solution $u$ to \eqref{eq:Sell_You_47-1} in the space \eqref{eq:Sell_You_47-5_Kozono_Maeda_Naito_lemma_3-4}, given $u_0 \in \cW$, we shall employ a contraction mapping argument for the map $\sT$ in \eqref{eq:Sell_You_page_233_contraction_map_definition}, based on the choice of Banach subspace, $\fV$, of functions $v \in C([t_0, t_0+\tau]; \cW)\cap C((0,T); \calV^{2\beta})$ with finite norm
\begin{equation}
\label{eq:Kozono_Maeda_Naito_Banach_space}
\|v\|_\fV
:=
\sup_{\begin{subarray}{c}t\in [t_0, t_0+\tau] \\ \delta \in \{0, \alpha, \beta\} \end{subarray}} (t-t_0)^\delta \|v(t)\|_{\calV^{2\delta}}.
\end{equation}
Our choice of $\fV$ is motivated by \cite[Equation (3.9)]{Kozono_Maeda_Naito_1995}.

\begin{thm}[Existence and uniqueness of mild solutions to a nonlinear evolution equation in a Banach space]
\label{thm:Kozono_Maeda_Naito_lemma_3-4}
Assume the setup
of the paragraphs preceding Theorem \ref{thm:Sell_You_lemma_47-1} and that, for some $\beta \in [0, 1)$, the nonlinearity,
\begin{equation}
\label{eq:Sell_You_47-4_Kozono_Maeda_Naito}
\cF \in C_{\Lip}([0,\infty)\times \calV^{2\beta}; \cW),
\end{equation}
obeys \eqref{eq:Sell_You_46-7_Kozono_Maeda_Naito} and \eqref{eq:Sell_You_46-8_Kozono_Maeda_Naito} for some $n \geq 1$, and $\alpha, \gamma \in [0,1)$ obeying $\gamma+n\alpha < 1$ and $\beta+\gamma < 1$, and positive constants $\kappa_2, \kappa_3$. Given $b > 0$, there exists a positive constant\footnote{The constant $\tau$ is given explicitly via \eqref{eq:Sell_You_47-6_Kozono_Maeda_Naito}.}
$$
\tau = \tau\left(b, M_0, M_\alpha, M_\beta, M_\gamma, M_{\alpha+\gamma}, M_{\beta+\gamma}, n, \alpha, \beta, \gamma, \kappa_2, \kappa_3 \right) > 0,
$$
with the following significance.  For every $u_0 \in \cW$ obeying $\|u_0\|_\cW \leq b$ and every $t_0 \geq 0$,
the initial value problem \eqref{eq:Sell_You_47-1} has a unique, mild solution in $\cW$, on the interval $[t_0, t_0 + \tau]$, which obeys
\begin{equation}
\label{eq:Sell_You_47-5_Kozono_Maeda_Naito_lemma_3-4}
u \in C([t_0, t_0+\tau]; \cW) \cap C((t_0, t_0+\tau); \calV^{2\beta}).
\end{equation}
Moreover, for $\delta=0,\alpha,\beta$, the solution $u$ obeys the \apriori estimate,
\begin{multline}
\label{eq:Kozono_Maeda_Naito_lemma_3-4_apriori_estimate_with_rearrangement}
(t - t_0)^\delta\|u(t)\|_{\calV^{2\delta}}
\leq
M_\delta \|u_0\|_\cW + M_{\delta+\gamma} B(1-n\alpha, 1-\delta-\gamma) \kappa_2
\\
\times \left( \tau^{n\alpha} + \left(\bar M b  + 1\right)^n \right)
(t - t_0)^{1-\gamma-n\alpha}, \quad \forall\, t \in [t_0, t_0 + \tau],
\end{multline}
where $\bar M := \max_{\delta=0,\alpha,\beta} M_\delta$ and $B(x,y)$ is Euler's Beta integral \cite[Section 5.12]{Olver_Lozier_Boisvert_Clark}; if $v_0 \in \cW$ obeys $\|v_0\|_\cW \leq b$ and $v$ is the unique, mild solution to \eqref{eq:Sell_You_47-1} in $\cW$ on $[t_0, t_0 + \tau]$ with $v(0) = v_0$ and
and satisfying \eqref{eq:Sell_You_47-5_Kozono_Maeda_Naito_lemma_3-4}, then
\begin{equation}
\label{eq:Kozono_Maeda_Naito_lemma_3-4_continuity_with_respect_to_initial_data}
\|u - v\|_\fV \leq 2\bar M \|u_0 - v_0\|_\cW.
\end{equation}
\end{thm}

\begin{rmk}[Application to the Yang-Mills heat equation]
\label{rmk:Kozono_Maeda_Naito_lemma_3-4}
While we include Theorem \ref{thm:Kozono_Maeda_Naito_lemma_3-4} for the sake of completeness, it appears difficult to usefully apply the result to the Yang-Mills heat equation \eqref{eq:Yang-Mills_gradient_flow_equation} because of the relatively high degree ($n =3$) of the Yang-Mills heat equation polynomial nonlinearity
\eqref{eq:Yang-Mills_heat_equation_nonlinearity_relative_rough_Laplacian_plus_one}.
\end{rmk}

\begin{proof}[Proof of Theorem \ref{thm:Kozono_Maeda_Naito_lemma_3-4}]
We proceed as in the proof of Theorem \ref{thm:Sell_You_lemma_47-1_polynomial_nonlinearity}, but set
\begin{equation}
\label{eq:Kozono_Maeda_Naito_radius}
R :=  \left(\max_{\delta=0,\alpha,\beta} M_\delta\right)b  + 1,
\end{equation}
and, noting that $\gamma + n\alpha < 1$ by hypothesis, choose $\tau > 0$ to be the largest constant such that
\begin{subequations}
\label{eq:Sell_You_47-6_Kozono_Maeda_Naito}
\begin{gather}
\label{eq:Sell_You_47-6_Kozono_Maeda_Naito_bound}
\max_{\delta=0,\alpha,\beta} M_{\delta+\gamma} B(1-n\alpha, 1-\delta-\gamma) \kappa_2
(\tau^{n\alpha} + R^n) \tau^{1-\gamma-n\alpha} \leq 1,
\\
\label{eq:Sell_You_47-6_Kozono_Maeda_Naito_Lipschitz}
\max_{\delta=0,\alpha,\beta} M_{\delta+\gamma} B(1-n\alpha, 1-\delta-\gamma) \kappa_3
(\tau^{(n-1)\alpha} + R^{n-1}) \tau^{1-\gamma-n\alpha} \leq \frac{1}{2}.
\end{gather}
\end{subequations}
It is convenient to divide the proof of existence and uniqueness of $u$ into three steps.

\setcounter{step}{0}
\begin{step}[Boundedness of the map $\sT$]
\label{step:Kozono_Maeda_Naito_lemma_3-4_proof_of_boundedness}
We define $\fF := \{v \in \fV: \|v\|_\fV \leq R\}$, a closed ball in $\fV$, and now show that $\sT$ in maps $\fF$ into itself. Let $u \in \fF$ and set $\hat u = \sT u$. For $t_0 \leq t \leq t_0+\tau$,
\begin{align*}
{}& \|\cA^\beta \hat u(t)\|_\cW
\\
&\quad \leq \|\cA^\beta e^{-\cA (t-t_0)}u_0\|_\cW
+ \int_{t_0}^t \|\cA^{\beta+\gamma} e^{-\cA(t-s)} \cA^{-\gamma}\cF(s, u(s))\|_\cW \, ds
\\
&\quad \leq M_\beta (t-t_0)^{-\beta} e^{-a(t-t_0)}\|u_0\|_\cW
+ M_{\beta+\gamma} \int_{t_0}^t (t-s)^{-\beta-\gamma} e^{-a(t-s)}\|\cA^{-\gamma}\cF(s, u(s))\|_\cW \, ds
\\
&\qquad\hbox{(by \eqref{eq:Sell_You_37-11})}
\\
&\quad \leq M_\beta (t-t_0)^{-\beta} \|u_0\|_\cW
+ M_{\beta+\gamma} \kappa_2 \int_{t_0}^t (t-s)^{-\beta-\gamma}
\left(1 + \|u(s)\|_{\calV^{2\alpha}}^n \right)\, ds
\quad\hbox{(by \eqref{eq:Sell_You_46-7_Kozono_Maeda_Naito} and $a\geq 0$)}
\\
&\quad = M_\beta (t-t_0)^{-\beta} \|u_0\|_\cW
\\
&\qquad + M_{\beta+\gamma} \kappa_2 \int_{t_0}^t (t-s)^{-\beta-\gamma}
(s-t_0)^{-n\alpha}\left((s-t_0)^{n\alpha} + (s-t_0)^{n\alpha}\|u(s)\|_{\calV^{2\alpha}}^n \right) \,ds
\\
&\quad \leq M_\beta (t-t_0)^{-\beta} \|u_0\|_\cW
+ M_{\beta+\gamma} \kappa_2 \left((t-t_0)^{n\alpha}
+ \sup_{s\in[t_0,t]}(s-t_0)^{n\alpha}\|u(s)\|_{\calV^{2\alpha}}^n \right)
\\
&\qquad \times \int_{t_0}^t (t-s)^{-\beta-\gamma} (s-t_0)^{-n\alpha}\,ds.
\end{align*}
The integral on the right-hand side in the preceding inequality may be calculated via Euler's Beta integral \cite[Section 5.12]{Olver_Lozier_Boisvert_Clark}:
\begin{align*}
{}&\int_{t_0}^t (t-s)^{-\beta-\gamma} (s-t_0)^{-n\alpha}\,ds
\\
&\quad =
\int_0^{t-t_0} (t-t_0-r)^{-\beta-\gamma} r^{-n\alpha}\,dr
\\
&\quad = (t-t_0)^{1-\beta-\gamma-n\alpha}\int_0^1 q^{(1-n\alpha) - 1} (1-q)^{(1-\beta-\gamma) - 1} \,dq
\\
&\quad = (t-t_0)^{1-\beta-\gamma-n\alpha} B(1-n\alpha, 1-\beta-\gamma)
\\
&\quad = (t-t_0)^{1-\beta-\gamma-n\alpha} \frac{\Gamma(1-n\alpha)\Gamma(1-\beta-\gamma)}{\Gamma(2-n\alpha-\beta-\gamma)}
\quad\hbox{(by \cite[Equation 5.12.1]{Olver_Lozier_Boisvert_Clark}),}
\end{align*}
provided $n\alpha < 1$ and $\beta+\gamma < 1$, as assured by our hypotheses. Therefore, when $\delta=\beta$,
\begin{multline}
\label{eq:Kozono_Maeda_Naito_lemma_3-4_apriori_estimate_preliminary}
(t-t_0)^\delta \|\cA^\delta \hat u(t)\|_\cW
\leq
M_\delta \|u_0\|_\cW + M_{\delta+\gamma} B(1-n\alpha, 1-\delta-\gamma) \kappa_2 (t-t_0)^{1-\gamma-n\alpha}
\\
\times \left(\tau^{n\alpha} + \sup_{s\in[t_0, t_0+\tau]}\left((s-t_0)^\alpha\|u(s)\|_{\calV^{2\alpha}}\right)^n \right),
\quad \forall\, t \in [t_0, t_0+\tau],
\end{multline}
where we recall that $\gamma+n\alpha < 1$ by hypothesis; naturally, one can repeat the same calculation with $\beta$ replaced by $0$ or $\alpha$ to conclude that \eqref{eq:Kozono_Maeda_Naito_lemma_3-4_apriori_estimate_preliminary} holds for $\delta = 0, \alpha, \beta$ and $\delta+\gamma < 1$. Therefore, by definition of $\fV$ in \eqref{eq:Kozono_Maeda_Naito_Banach_space}, and $R$ in \eqref{eq:Kozono_Maeda_Naito_radius}, the facts that $\|u_0\|_\cW \leq b$ by hypothesis, $\tau$ obeys \eqref{eq:Sell_You_47-6_Kozono_Maeda_Naito_bound}, and $u \in \fF$ so $\|u\|_\fV \leq R$, we obtain
$$
\|(t-t_0)^\delta \cA^\delta \hat u(t)\|_\cW \leq R, \quad\forall\, t \in [t_0, t_0 + \tau], \quad \delta = 0, \alpha, \beta,
$$
and thus,
$$
\|\sT u\|_\fV \equiv \|\hat u\|_\fV \leq R, \quad \forall\, u \in \fF.
$$
This completes Step \ref{step:Kozono_Maeda_Naito_lemma_3-4_proof_of_boundedness}.
\end{step}

\begin{step}[Contraction mapping property of $\sT$]
\label{step:Kozono_Maeda_Naito_lemma_3-4_proof_of_contraction property}
Next we show that for $\tau$ obeying \eqref{eq:Sell_You_47-6_Kozono_Maeda_Naito_Lipschitz}, the mapping $\sT$ is a contraction on $\fF$ with contraction coefficient less than or equal to $1/2$. We defer the proof that $\hat u \in C([t_0, t_0+\tau]; \cW)$, when $u\in \fF$, to Step \ref{step:Kozono_Maeda_Naito_lemma_3-4_proof_of_continuity}. Let $u_1, u_2 \in \fF$ and set $\hat u_i := \sT u_i$ for $i = 1,2$. For $t_0 \leq t \leq t_0 + \tau$, we have
\begin{align*}
{}& \|\cA^\beta (\hat u_1(t) - \hat u_2(t))\|_\cW
\\
&\quad \leq \int_{t_0}^t \|\cA^\beta e^{-\cA(t-s)} (\cF(s, u_1(s)) - \cF(s, u_2(s)))\|_\cW \, ds
\\
&\quad = \int_{t_0}^t \|\cA^{\beta+\gamma} e^{-\cA(t-s)} \cA^{-\gamma} (\cF(s, u_1(s)) - \cF(s, u_2(s)))\|_\cW \, ds
\\
&\quad\leq M_{\beta+\gamma} \int_{t_0}^t (t-s)^{-\beta-\gamma} e^{-a(t-s)}\|\cA^{-\gamma}
(\cF(s, u_1(s)) - \cF(s, u_2(s)))\|_\cW \, ds
\quad\hbox{(by \eqref{eq:Sell_You_37-11})}
\\
&\quad \leq M_{\beta+\gamma} \kappa_3 \int_{t_0}^t (t-s)^{-\beta-\gamma}
\left(1 + \|u_1(s)\|_{\calV^{2\alpha}}^{n-1} + \|u_2(s)\|_{\calV^{2\alpha}}^{n-1} \right)\|u_1(s) - u_2(s)\|_{\calV^{2\alpha}} \,ds
\\
&\qquad\hbox{(by \eqref{eq:Sell_You_46-8_Kozono_Maeda_Naito} and $a\geq 0$)}
\\
&\quad = M_{\beta+\gamma} \kappa_3 \int_{t_0}^t (t-s)^{-\beta-\gamma} (s - t_0)^{-n\alpha}
(s - t_0)^\alpha \|u_1(s) - u_2(s)\|_{\calV^{2\alpha}}
\\
&\qquad \times \left((s - t_0)^{(n-1)\alpha} + (s - t_0)^{(n-1)\alpha}\|u_1(s)\|_{\calV^{2\alpha}}^{n-1}
+ (s - t_0)^{(n-1)\alpha}\|u_2(s)\|_{\calV^{2\alpha}}^{n-1} \right) \,ds
\\
&\quad \leq M_{\beta+\gamma} \kappa_3 \sup_{s\in [t_0, t]}(s - t_0)^\alpha \|u_1(s) - u_2(s)\|_{\calV^{2\alpha}}
\int_{t_0}^t (t-s)^{-\beta-\gamma} (s - t_0)^{-n\alpha}\,ds
\\
&\qquad \times \left( (t - t_0)^{(n-1)\alpha}
+ \sup_{\begin{subarray}{c}s\in [t_0, t] \\ i=1,2 \end{subarray}}
\left((s - t_0)^\alpha \|u_i(s)\|_{\calV^{2\alpha}}\right)^{n-1} \right).
\end{align*}
Thus, applying \cite[Equation 5.12.1]{Olver_Lozier_Boisvert_Clark} as before and recalling the definition \eqref{eq:Kozono_Maeda_Naito_Banach_space} of $\fV$ and that $\|u_i\|_\fV \leq R$ for $i = 1, 2$, we obtain, when $\delta=\beta$,
\begin{multline}
\label{eq:Kozono_Maeda_Naito_lemma_3-4_continuity_with_respect_to_initial_data_preliminary}
(t-t_0)^\delta \|\cA^\delta (\hat u_1(t) - \hat u_2(t))\|_\cW
\leq M_{\delta+\gamma} \kappa_3 \tau^{1-\gamma-n\alpha} B(1-n\alpha, 1-\delta-\gamma)
\\
\times \left( \tau^{(n-1)\alpha} + R^{n-1} \right)
\sup_{s\in [t_0, t_0+\tau]}(s - t_0)^\alpha \|u_1(s) - u_2(s)\|_{\calV^{2\alpha}},
\quad \forall\, t \in [t_0, t_0 + \tau],
\end{multline}
and that by repeating the preceding argument with $\beta$ replaced by $0$ or $\alpha$, the inequality
\eqref{eq:Kozono_Maeda_Naito_lemma_3-4_continuity_with_respect_to_initial_data_preliminary} also holds when $\delta = 0$ or $\alpha$. Therefore, by the fact that $\tau$ obeys \eqref{eq:Sell_You_47-6_Kozono_Maeda_Naito_Lipschitz}, the definition \eqref{eq:Kozono_Maeda_Naito_Banach_space} of $\fV$, and the fact that $u_i \in \fF$ implies $\|u_i\|_\fV \leq R$, we see that
$$
\|\sT u_1 - \sT u_2\|_\fV \equiv \|\hat u_1 - \hat u_2\|_\fV \leq \frac{1}{2}\|u_1 - u_2\|_\fV, \quad \forall\, u_1, u_2 \in \fF.
$$
This completes Step \ref{step:Kozono_Maeda_Naito_lemma_3-4_proof_of_contraction property}.
\end{step}

\begin{step}[Continuity of $\hat u(t)$ for $t_0\leq t\leq t_0+\tau$ and existence and uniqueness of the solution]
\label{step:Kozono_Maeda_Naito_lemma_3-4_proof_of_continuity}
By Theorem \ref{thm:Sell_You_lemma_47-1_polynomial_nonlinearity}, we know already that a mild solution to \eqref{eq:Sell_You_47-1} exists, belongs to \eqref{eq:Sell_You_47-5_Kozono_Maeda_Naito_lemma_3-4} (because, \afortiori, it already has greater regularity \eqref{eq:Sell_You_47-5_polynomial_nonlinearity}) and is unique when $u_0 \in \calV^{2\beta}$, so we may employ an approximation argument (similar to the proof of \cite[Theorem 8.30]{GilbargTrudinger}) when it is only known that $u_0 \in \cW$ to prove that $\hat u$ has the claimed $C([t_0, t_0+\tau]; \cW)$-continuity property.

Recall that $\calV^{2\beta}$ is dense in $\cW$, by the definition $V^{2\beta} := \sD(\cA^\beta)$ in \cite[p. 96]{Sell_You_2002} and the fact that $\sD(\cA^\beta)$ is dense in $\cW$ for any $\beta\geq 0$ by \cite[Lemma 34.7 (1)]{Sell_You_2002} since $\cA$ is a positive, sectorial operator on the Banach space $\cW$ by the Hypothesis \ref{hyp:Sell_You_4_standing_hypothesis_A}. Hence, given $u_0 \in \cW$ with $\|u_0\|_\cW \leq b$, we may choose a sequence $\{u_0^m\}_{m\in\NN} \subset \calV^{2\beta}$ such that $u_0^m \to u_0$ in $\cW$ as $m \to \infty$ and $\|u_0^m\|_\cW \leq b$ for all $m \in \NN$. Let $\{u^m\}_{m\in\NN}$ denote the corresponding sequence of mild solutions in $V^{2\beta}$ to \eqref{eq:Sell_You_47-1} on $[t_0, t_0+\tau]$ and $u^m(0) = u_0^m$ that are produced by Theorem \ref{thm:Sell_You_lemma_47-1_polynomial_nonlinearity}. From Step \ref{step:Kozono_Maeda_Naito_lemma_3-4_proof_of_boundedness}, we have the \apriori estimates,
$$
\|u^m\|_\fV \leq R, \quad \forall\, m \in \NN,
$$
and thus from the derivation of \eqref{eq:Kozono_Maeda_Naito_lemma_3-4_apriori_estimate_preliminary} and \eqref{eq:Kozono_Maeda_Naito_lemma_3-4_continuity_with_respect_to_initial_data_preliminary} and noting that $\hat u^m \equiv \sT u^m = u^m$ for all $m\in\NN$, we obtain the following estimate for all $k, l \in \NN$:
\begin{multline*}
\|u^k - u^l\|_\fV
\leq
\left(\max_{\delta=0,\alpha,\beta} M_\delta\right)\|u_0^k - u_0^l\|_\cW
\\
+ \max_{\delta=0,\alpha,\beta} M_{\delta+\gamma}
\kappa_3 \tau^{1-\gamma-n\alpha} B(1-n\alpha, 1-\delta-\gamma) \left(\tau^{(n-1)\alpha} + R^{n-1} \right)\|u^k - u^l\|_\fV.
\end{multline*}
In particular, by the fact that $\tau$ obeys \eqref{eq:Sell_You_47-6_Kozono_Maeda_Naito_Lipschitz}, we have
$$
\|u^k - u^l\|_\fV \leq \left(\max_{\delta=0,\alpha,\beta} M_\delta\right)\|u_0^k - u_0^l\|_\cW
+ \frac{1}{2}\|u^k - u^l\|_\fV,
$$
and so rearrangement yields,
\begin{equation}
\label{eq:Kozono_Maeda_Naito_lemma_3-4_continuity_with_respect_to_initial_data_sequence}
\|u^k - u^l\|_\fV \leq 2\left(\max_{\delta=0,\alpha,\beta} M_\delta\right)\|u_0^k - u_0^l\|_\cW, \quad\forall\, k, l \in \NN.
\end{equation}
Therefore, $\{u^m\}_{m\in\NN}$ is Cauchy in $\fV$, so $u^m \to \tilde u$ in $\fV$ as $m \to \infty$, for some $\tilde u \in \fV$. Taking the limit of $u^m = \sT u^m$ as $m \to \infty$ gives $\tilde u = \sT \tilde u$, that is, $\tilde u$ solves \eqref{eq:Sell_You_47-2} and obeys $\|\tilde u\|_\fV \leq R$, by taking the limit of our \apriori estimate for the sequence $\{u^m\}_{m\in\NN}$. The estimate \eqref{eq:Kozono_Maeda_Naito_lemma_3-4_continuity_with_respect_to_initial_data_sequence} remains valid when $u^k$ is replaced by $\hat u$ and $u^l$ is replaced by $\tilde u$ and $u_0^k$ and $u_0^l$ are replaced by $u_0$. Thus, we must have $\tilde u = \hat u$. Since $\{u^m\}_{m\in\NN}$ is Cauchy in $C([t_0, t_0+\tau]; \cW)$ (because, \afortiori, it is Cauchy in $\fV$), we necessarily also have $\tilde u = \hat u \in C([t_0, t_0+\tau]; \cW)$.

From Step \ref{step:Kozono_Maeda_Naito_lemma_3-4_proof_of_boundedness}, we have $\|\sT u\|_\fV \leq R$ for all $u \in \fF$ and thus, $\sT u \in \fF$ for all $u \in \fF$ since we now know in addition that $\sT u \equiv \hat u \in C([t_0, t_0+\tau]; \cW)$. Because $\sT : \fF \to \fF$ is a contraction by Step \ref{step:Kozono_Maeda_Naito_lemma_3-4_proof_of_contraction property}, the map $\sT$ has a unique fixed point $u \in \fF$. This fixed point is the mild solution of \eqref{eq:Sell_You_47-1} in $\cW$ on $[t_0, t_0+\tau]$, and because $\sT$ is a contraction, it is unique. This completes Step \ref{step:Kozono_Maeda_Naito_lemma_3-4_proof_of_continuity}.
\end{step}

Our proof of continuity of $\sT u(t)$ with respect to $t \in [t_0, t_0+\tau]$ also yields the regularity result,
$$
u \in C([0, \tau]; \cW) \cap C((0, \tau]; \calV^{2\beta}),
$$
that is, $u$ obeys \eqref{eq:Sell_You_47-5_Kozono_Maeda_Naito_lemma_3-4}.

\begin{step}[\Apriori and continuity estimates]
\label{step:Kozono_Maeda_Naito_lemma_3-4_apriori_and_continuity_estimates}
The \apriori estimate \eqref{eq:Kozono_Maeda_Naito_lemma_3-4_apriori_estimate_with_rearrangement} for $u$ follows from \eqref{eq:Kozono_Maeda_Naito_lemma_3-4_apriori_estimate_preliminary}, the definition \eqref{eq:Kozono_Maeda_Naito_Banach_space} of $\fV$, the definition \eqref{eq:Kozono_Maeda_Naito_radius} of $R$, and the facts that $\hat u = u$ and $\|u\|_\fV \leq R$ since $u \in \fF$. The continuity estimate \eqref{eq:Kozono_Maeda_Naito_lemma_3-4_continuity_with_respect_to_initial_data} follows from \eqref{eq:Kozono_Maeda_Naito_lemma_3-4_continuity_with_respect_to_initial_data_sequence} by replacing $u^k$ by $u$ and $u^l$ by $v$.
\end{step}

This completes the proof of Theorem \ref{thm:Kozono_Maeda_Naito_lemma_3-4}.
\end{proof}

Of course, when $t > 0$, the unique mild solution provided by Theorem \ref{thm:Kozono_Maeda_Naito_lemma_3-4} has the same regularity properties as discovered in our results originating with Theorem \ref{thm:Sell_You_lemma_47-1}.

\subsection{Existence and uniqueness of mild solutions to a nonlinear evolution equation in a Banach space with initial data of minimal regularity II}
\label{subsec:Sell_You_lemma_4-7-1_Kozono_Maeda_Naito_II}
In Section \ref{subsec:Sell_You_lemma_4-7-1_Kozono_Maeda_Naito}, we proved existence and uniqueness of mild solutions to a nonlinear evolution equation in a Banach space $\calV^{2\alpha}$ on an interval $(0, \tau)$, given initial data in $\cW$, when $\alpha$ obeys $\gamma + n\alpha < 1$, for some $\alpha, \gamma \in [0, 1)$ and $n \geq 1$ and the nonlinearity, $\cF(t, v)$, obeys the structural conditions defined by $n, \alpha, \gamma$ and $\cA$ in \eqref{eq:Sell_You_46-7_Kozono_Maeda_Naito} and \eqref{eq:Sell_You_46-8_Kozono_Maeda_Naito}, for all $(t, v) \in [0, \calV^{2\alpha})$. However, for later applications to the Yang-Mills heat equation \eqref{eq:Yang-Mills_heat_equation_as_perturbation_rough_Laplacian_plus_one_heat_equation},
the constraint $\gamma + n\alpha < 1$ is too strong. Indeed, keeping in mind \cite[Lemma 3.3]{Kozono_Maeda_Naito_1995} and the cubic polynomial structure of the (time-homogeneous) Yang-Mills heat equation nonlinearity
\eqref{eq:Yang-Mills_heat_equation_nonlinearity_relative_rough_Laplacian_plus_one}, we would like to choose $n = 3$ and $\alpha = \gamma = \frac{1}{4}$ and have $\cF \in C_{\Lip}(\calV; \cW)$ obey the following bounds, for some positive constants, $\mu_1, \mu_2, \mu_3$,
\begin{subequations}
\label{eq:Kozono_Maeda_Naito_3-7}
\begin{align}
\label{eq:Kozono_Maeda_Naito_3-7_constant}
f_0 &\in \cW,
\\
\label{eq:Kozono_Maeda_Naito_3-7_linear}
\|\cF_1(v)\|_\cW &\leq \mu_1\|v\|_\cW
\\
\label{eq:Kozono_Maeda_Naito_3-7_quadratic}
\|\cA^{-\frac{1}{4}}\cF_2(v)\|_\cW &\leq \mu_2\|\cA^{\frac{1}{2}}v\|_\cW \|\cA^{\frac{1}{4}}v\|_\cW
\\
\notag
&= \mu_2\|v\|_\calV\|v\|_{\calV^{\frac{1}{2}}},
\\
\label{eq:Kozono_Maeda_Naito_3-7_cubic}
\|\cA^{-\frac{1}{4}}\cF_3(v)\|_\cW &\leq \mu_3\|\cA^{\frac{1}{4}}v\|_\cW^3
\\
\notag
&= \mu_3\|v\|_{\calV^{\frac{1}{2}}}^3, \quad\forall\, v \in \calV,
\end{align}
\end{subequations}
and where the nonlinearity has the structure,
\begin{equation}
\label{eq:Kozono_Maeda_Naito_3-2_nonlinearity}
\cF(v) = f_0 + \cF_1(v) + \cF_2(v) + \cF_3(v), \quad\forall\, v \in \calV,
\end{equation}
and $\cA$ obeys Hypothesis \ref{hyp:Sell_You_4_standing_hypothesis_A}.

Kozono, Maeda, and Naito choose $\cW = L^d(X; \Lambda^1\otimes\ad P)$ in \cite[pp. 100-101]{Kozono_Maeda_Naito_1995}, so $\calV^2 = W^{2,d}_{A_1}(X; \Lambda^1\otimes\ad P)$, the domain of $\sA_2: \sD(\sA_2) \subset L^d(X; \Lambda^1\otimes\ad P) \to L^d(X; \Lambda^1\otimes\ad P)$, when $\cA = \sA_2$ is the realization of $\sA = \nabla_{A_1}^*\nabla_{A_1} + 1$, the (augmented) connection Laplace operator \eqref{eq:Connection_Laplacian} on $C^\infty(X; \Lambda^1\otimes\ad P)$. The Sobolev Embedding Theorem  \cite[Theorem 4.12]{AdamsFournier} ensures that $W^{2\delta,d}(X) \subset C(X)$ and that $W^{2\delta,d}(X)$ is a Banach algebra when $\delta > \frac{1}{2}$.

Again keeping in mind \cite[Lemma 3.3]{Kozono_Maeda_Naito_1995}, we shall also require that $\cF$ obey the following Lipschitz continuity conditions, for some positive constants, $\eta_1, \eta_2$, and $v_1, v_2 \in \calV$,
\begin{subequations}
\label{eq:Kozono_Maeda_Naito_3-8}
\begin{align}
\label{eq:Kozono_Maeda_Naito_3-8_quadratic}
\|\cA^{-\frac{1}{4}}(\cF_2(v_1) - \cF_2(v_2))\|_\cW
&\leq
\eta_2\left(\|\cA^{\frac{1}{2}}v_1\|_\cW + \|\cA^{\frac{1}{2}}v_2\|_\cW\right) \|\cA^{\frac{1}{4}}(v_1-v_2)\|_\cW
\\
\notag
&\quad + \eta_2\left(\|\cA^{\frac{1}{4}}v_1\|_\cW + \|\cA^{\frac{1}{4}}v_2\|_\cW\right) \|\cA^{\frac{1}{2}}(v_1-v_2)\|_\cW
\\
\notag
&= \eta_2\left(\|v_1\|_\calV + \|v_2\|_\calV\right)\|v_1 - v_2\|_{\calV^{\frac{1}{2}}}
\\
\notag
&\quad + \eta_2\left(\|v_1\|_{\calV^{\frac{1}{2}}} + \|v_2\|_{\calV^{\frac{1}{2}}}\right)\|v_1 - v_2\|_\calV,
\\
\label{eq:Kozono_Maeda_Naito_3-8_cubic}
\|\cA^{-\frac{1}{4}}(\cF_3(v_1) - \cF_3(v_2))\|_\cW
&\leq \eta_3\left(\|\cA^{\frac{1}{4}}v_1\|_\cW^2 + \|\cA^{\frac{1}{4}}v_2\|_\cW^2\right) \|\cA^{\frac{1}{4}}(v_1-v_2)\|_\cW
\\
\notag
&= \eta_3\left(\|v_1\|_{\calV^{\frac{1}{2}}}^2 + \|v_2\|_{\calV^{\frac{1}{2}}}^2\right)\|v_1 - v_2\|_{\calV^{\frac{1}{2}}}.
\end{align}
\end{subequations}
These inequalities are, of course, motivated by the structure of the Yang-Mills heat equation nonlinearity
\eqref{eq:Yang-Mills_heat_equation_nonlinearity_relative_rough_Laplacian_plus_one}.

One should realize that our Theorem \ref{thm:Kozono_Maeda_Naito_lemma_3-4} above is both stronger than \cite[Lemma 3.4]{Kozono_Maeda_Naito_1995}, in the sense that it also asserts uniqueness of the mild solution, but also weaker in that its hypothesis, $\gamma + n\alpha < 1$, excludes the important case $n = 3$ and $\alpha = \gamma = \frac{1}{4}$. The lack of uniqueness in \cite[Lemma 3.4]{Kozono_Maeda_Naito_1995} can perhaps be interpreted as the result of an appeal to the Schauder Fixed-point Theorem \cite[Theorem 11.1 or Corollary 11.2]{GilbargTrudinger} rather than the Banach Contraction Mapping Theorem.

To prove existence and uniqueness of a solution $u$ to \eqref{eq:Sell_You_47-1} in the space \eqref{eq:Sell_You_47-5_Kozono_Maeda_Naito_lemma_3-4_plus_uniqueness}, given $u_0 \in \cW$ and a positive constant $\tau$, we shall again employ a contraction mapping argument for a certain map $\sT$, based on the choice of Banach subspace, $\fV$, of functions $v \in C([0, \tau]; \cW)\cap C((0,\tau]; \calV)$ with finite norm, by analogy with \eqref{eq:Kozono_Maeda_Naito_Banach_space} but now defined as
\begin{equation}
\label{eq:Kozono_Maeda_Naito_Banach_space_plus_uniqueness}
\|v\|_\fV
:=
\sup_{\begin{subarray}{c}t\in [0, \tau] \\ \delta \in \{0, \frac{1}{4}, \frac{1}{2}\} \end{subarray}} t^\delta \|v(t)\|_{\calV^{2\delta}}
\equiv
\sup_{\begin{subarray}{c}t\in [0, \tau] \\ \delta \in \{0, \frac{1}{4}, \frac{1}{2}\} \end{subarray}}
t^\delta \|\cA^\delta v(t)\|_\cW.
\end{equation}
Basic properties of the fractional power spaces, $\calV^{2\alpha} = \sD(\cA^\alpha)$, are given by Lemma \ref{lem:Sell_You_37-4} and \cite[Theorems 37.6 and 37.7]{Sell_You_2002}, together with a key example in \cite[Lemma 37.8]{Sell_You_2002}.

We can now state the following significant improvement of both our Theorem \ref{thm:Kozono_Maeda_Naito_lemma_3-4} and \cite[Lemma 3.4]{Kozono_Maeda_Naito_1995}, in that we now obtain both existence \emph{and} uniqueness. Moreover, our proof of Theorem \ref{thm:Kozono_Maeda_Naito_lemma_3-4_plus_uniqueness} is considerably simpler than the proof of \cite[Lemma 3.4]{Kozono_Maeda_Naito_1995}.

\begin{thm}[Existence and uniqueness of mild solutions to a nonlinear evolution equation in a Banach space]
\label{thm:Kozono_Maeda_Naito_lemma_3-4_plus_uniqueness}
Assume that Hypothesis \ref{hyp:Sell_You_4_standing_hypothesis_A} holds and
\begin{equation}
\label{eq:Sell_You_47-4_Kozono_Maeda_Naito_plus_uniqueness}
\cF \in C_{\Lip}(\calV; \cW),
\end{equation}
obeys \eqref{eq:Kozono_Maeda_Naito_3-7}, \eqref{eq:Kozono_Maeda_Naito_3-2_nonlinearity}, and \eqref{eq:Kozono_Maeda_Naito_3-8} for $f_0 \in \cW$ and some positive constants, $\mu_1, \mu_2, \mu_3, \eta_2, \eta_3$. Given $\delta_0 \in [\frac{1}{2}, \frac{3}{4})$ and $u_0 \in \cW$, there are positive constants\footnote{The constant $\tau$ is given explicitly via \eqref{eq:Sell_You_47-6_Kozono_Maeda_Naito_plus_uniqueness} and $C_0$ is equal to the solution $R$ to \eqref{eq:Kozono_Maeda_Naito_radius_plus_uniqueness_bound_and_Lipschitz}.}
$$
\tau = \tau\left(u_0, \cA, \|f_0\|_\cW, \delta_0, \mu_1, \mu_2, \mu_3, \eta_2, \eta_3 \right)
\quad\hbox{and}\quad
C_0 = C_0(\cA, \delta_0, \eta_2, \eta_3),
$$
with the following significance.  The initial value problem \eqref{eq:Sell_You_47-1} has a unique, mild solution in $\cW$, on the interval $[0, \tau]$, which obeys
\begin{equation}
\label{eq:Sell_You_47-5_Kozono_Maeda_Naito_lemma_3-4_plus_uniqueness}
u \in C([0, \tau]; \cW) \cap C((0, \tau]; \calV^{2\delta_0}).
\end{equation}
Moreover, for every $\delta \in [0, \delta_0]$, the solution $u$ obeys the \apriori estimate,
\begin{equation}
\label{eq:Kozono_Maeda_Naito_lemma_3-4_plus_uniqueness_apriori_estimate_with_rearrangement}
t^\delta\|u(t)\|_{\calV^{2\delta}}
\leq
C_0,
\quad
\forall\, t \in (0, \tau].
\end{equation}
\end{thm}

\begin{rmk}[Initial data in Theorem \ref{thm:Kozono_Maeda_Naito_lemma_3-4_plus_uniqueness} and well-posedness]
\label{rmk:Kozono_Maeda_Naito_lemma_3-4_plus_uniqueness_initial_data}
It is important to realize that while Theorem \ref{thm:Kozono_Maeda_Naito_lemma_3-4_plus_uniqueness} asserts existence and uniqueness of a solution, $u$, satisfying the regularity property \eqref{eq:Sell_You_47-5_Kozono_Maeda_Naito_lemma_3-4_plus_uniqueness}, it does \emph{not} assert well-posedness in the sense that $u$ also depends continuously on the initial data $u_0 \in \cW$. The reason is that our method of proof (this is also true of the proof of \cite[Lemma 3.4]{Kozono_Maeda_Naito_1995}) requires us to choose an approximation in $\cW$ of $u_0 \in \cW$ by $\tilde u_0 \in \calV^2$ (or at least in $\calV^{2\delta_0}$) and the minimum lifetime, $\tau$, of the solution $u$ depends on the $\calV$-norm of $\tilde u_0$ and not just the $\cW$-norm of $u_0$. Note the contrast with Theorem \ref{thm:Kozono_Maeda_Naito_lemma_3-4}, which asserts existence, uniqueness, and continuity with respect to the initial data.
\end{rmk}

\begin{proof}[Proof of Theorem \ref{thm:Kozono_Maeda_Naito_lemma_3-4_plus_uniqueness}]
We again follow the strategy of the proof of Theorem \ref{thm:Sell_You_lemma_47-1_polynomial_nonlinearity}, but now choose $R$ to be the largest positive constant such that
\begin{equation}
\label{eq:Kozono_Maeda_Naito_radius_plus_uniqueness_bound_and_Lipschitz}
(2\eta_2 R + \eta_3 R^2) \times \max_{\delta\in [0, \delta_0]} M_{\delta+\frac{1}{4}} B\left(\textstyle\frac{3}{4} - \delta, \textstyle\frac{1}{4}\right) = \frac{1}{16},
\end{equation}
where $M_{\delta+\frac{1}{4}}$ is the positive constant associated with $\cA$ via Theorem \ref{thm:Sell_You_37-5}.

Our proof of Theorem \ref{thm:Kozono_Maeda_Naito_lemma_3-4} did not take advantage of the more refined structure of the nonlinearity $\cF$ in \eqref{eq:Kozono_Maeda_Naito_3-2_nonlinearity} but instead relied on the condition $\gamma + n\alpha < 1$ to guarantee that \eqref{eq:Sell_You_47-6_Kozono_Maeda_Naito_plus_uniqueness} holds. In the present case, we shall adapt an idea partially used, but not fully exploited in the proof of \cite[Lemma 3.4]{Kozono_Maeda_Naito_1995} and decompose the variation of constants formula to define the map $\sT:\fF \to \fF$. Recall that $\calV^{2\alpha}$ is dense in $\cW$, by the definition $\calV^{2\alpha} := \sD(\cA^\alpha)$ in \cite[p. 96]{Sell_You_2002} and the fact that $\sD(\cA^\alpha)$ is dense in $\cW$ for any $\alpha\geq 0$ by
\cite[Lemma 34.7 (1)]{Sell_You_2002} since $\cA$ is a positive, sectorial operator on the Banach space $\cW$ by the Hypothesis \ref{hyp:Sell_You_4_standing_hypothesis_A}. Therefore, we may choose $\tilde u_0 \in \calV^2$ (depending on $\eps$) such that
\begin{equation}
\label{eq:Sell_You_47-6_Kozono_Maeda_Naito_plus_uniqueness_choice_of_tilde_u0}
\|\tilde u_0 - u_0\|_\cW \times \max_{\delta\in [0, \delta_0]} M_\delta \leq \frac{R}{2},
\end{equation}
where again $M_\delta$ is the positive constant associated with $\cA$ via Theorem \ref{thm:Sell_You_37-5}.

Next, noting that $M$ is the constant in \eqref{eq:u_1_norm_C_0_to_T_into_V2delta} and depending at most on $\cA$, we choose $\tau$ to be the largest positive constant such that
\begin{subequations}
\label{eq:Sell_You_47-6_Kozono_Maeda_Naito_plus_uniqueness}
\begin{multline}
\label{eq:Sell_You_47-6_Kozono_Maeda_Naito_plus_uniqueness_bound}
\left( \|f_0\|_\cW + \mu_1M\|\tilde u_0\|_\cW \right) \times \max_{\delta \in  [0, \delta_0] }\frac{M_\delta \tau^{1-\delta}}{1-\delta}
+ \left( \tau^{\frac{3}{4}} \mu_2 M^2\|\tilde u_0\|_\calV \|\tilde u_0\|_{\calV^{\frac{1}{2}}} + \tau^{\frac{3}{4}} \mu_3 M^3\|\tilde u_0\|_{\calV^{\frac{1}{2}}}^3 \right.
\\
+ \left. \mu_1\tau^{\frac{3}{4}} R + \eta_2 \tau^{\frac{1}{2}} R M\|\tilde u_0\|_\calV + \eta_2 \tau^{\frac{1}{4}} R M\|\tilde u_0\|_{\calV^{\frac{1}{2}}}
+ \eta_3\tau^{\frac{1}{2}} R M^2\|\tilde u_0\|_{\calV^{\frac{1}{2}}}^2 \right)
\\
\times \max_{\delta\in [0, \delta_0]} M_{\delta+\frac{1}{4}} B\left(\textstyle\frac{3}{4} - \delta, \textstyle\frac{1}{4}\right) \leq \frac{R}{4},
\end{multline}
\begin{multline}
\label{eq:Sell_You_47-6_Kozono_Maeda_Naito_plus_uniqueness_Lipschitz}
\left( \mu_1 \tau^{\frac{3}{4}} + 2\eta_2 \tau^{\frac{1}{2}} M\|\tilde u_0\|_\calV + 2\eta_2 \tau^{\frac{1}{4}} M\|\tilde u_0\|_{\calV^{\frac{1}{2}}}
+ 4\eta_3 \tau^{\frac{1}{2}} M^2\|\tilde u_0\|_{\calV^{\frac{1}{2}}}^2 \right)
\\
\times \max_{\delta \in  [0, \delta_0] }
M_{\delta+\frac{1}{4}} B\left(\textstyle\frac{3}{4} - \delta, \textstyle\frac{1}{4}\right) \leq \frac{1}{4}.
\end{multline}
\end{subequations}
We divide the proof of existence and uniqueness of $u$ into several steps.

We define $\fF := \{v \in \fV: \|v\|_\fV \leq R\}$, a closed ball in $\fV$. Our goal is to apply the Banach Contraction Mapping Theorem to a map $\sT:\fF \to \fF$ suggested by the
the nonlinear integral equation (compare \eqref{eq:Sell_You_47-2}),
\begin{equation}
\label{eq:Sell_You_47-2_time_homogeneous}
u(t) = e^{-\cA t}u_0 + \int_0^t e^{-\cA (t - s)}\cF(u(s))\,ds, \quad\forall\, t \in [0, \tau],
\end{equation}
but our choice will not be the same as \eqref{eq:Sell_You_page_233_contraction_map_definition} because of the lower spatial regularity of $u_0$.

Recalling that $\tilde u_0 \in \calV$ was chosen to satisfy \eqref{eq:Sell_You_47-6_Kozono_Maeda_Naito_plus_uniqueness_choice_of_tilde_u0}, we set
$$
u_1(t) = e^{-\cA t}\tilde u_0, \quad\forall\, t \in [0, \infty),
$$
and observe that $u_1 \in C([0, \infty); \calV)$ by Theorem \ref{thm:Sell_You_37-5}. It is now natural to seek a solution, $u = u_1 + v$, to the nonlinear integral equation \eqref{eq:Sell_You_47-2_time_homogeneous} with the aid of a small perturbation, $v \in C([0, \tau]; \cW) \cap C((0, \tau]; \calV^{2\beta})$ for a suitable $\tau > 0$, that solves
\begin{multline*}
u_1(t) + v(t) =  e^{-\cA t}(u_0 - \tilde u_0) + u_1(t) + \int_0^t e^{-\cA (t - s)} \cF(u_1(s)) \,ds
\\
+ \int_0^t e^{-\cA (t - s)}[ \cF(u_1(s) + v(s)) - \cF(u_1(s)) ] \,ds,
\quad\forall\, t \in [0, \tau],
\end{multline*}
that is, for $v$ solving
\begin{multline*}
v(t) = e^{-\cA t}(u_0 - \tilde u_0) + \int_0^t e^{-\cA (t - s)} \cF(u_1(s)) \,ds
\\
+ \int_0^t e^{-\cA (t - s)}[ \cF(u_1(s) + v(s)) - \cF(u_1(s)) ] \,ds,
\quad\forall\, t \in [0, \tau].
\end{multline*}
A related device is used, for similar reasons, by Kozono, Maeda, and Naito in \cite[pp. 105--106]{Kozono_Maeda_Naito_1995} and by Struwe in \cite[p. 135]{Struwe_1994}. We therefore define
\begin{multline}
\label{eq:Sell_You_page_233_contraction_map_definition_Kozono_Maeda_Naito_plus_uniqueness}
\hat v(t) := (\sT v)(t) \equiv e^{-\cA t}(u_0 - \tilde u_0) + \int_0^t e^{-\cA (t - s)} \cF(u_1(s)) \,ds
\\
+ \int_0^t e^{-\cA (t - s)}[ \cF(u_1(s) + v(s)) - \cF(u_1(s)) ] \,ds, \quad\forall\, t \in [0, \tau], v \in \fV.
\end{multline}
As usual, we aim to show that the preceding formal expression for $\sT$ defines a contraction mapping on the closed ball $\fF \subset \fV$.

\setcounter{step}{0}
\begin{step}[Boundedness of the map $\sT$]
\label{step:Kozono_Maeda_Naito_lemma_3-4_plus_uniqueness_proof_of_boundedness}
Observe that, for all $\delta \in [0, \frac{3}{4})$ and $t \in (0, \tau]$ and $v \in \fV$, the expressions \eqref{eq:Kozono_Maeda_Naito_3-2_nonlinearity} for $\cF$ and \eqref{eq:Sell_You_page_233_contraction_map_definition_Kozono_Maeda_Naito_plus_uniqueness} for $\hat v$ yield
\begin{align*}
\|\cA^\delta \hat v(t)\|_\cW &\leq \|\cA^\delta e^{-\cA t}(u_0 - \tilde u_0)\|_\cW
\\
&\quad + \int_0^t \|\cA^\delta e^{-\cA (t - s)} \cF(u_1(s))\|_\cW \,ds
 + \int_0^t \|\cA^{\delta+\frac{1}{4}} e^{-\cA (t - s)}\cA^{-\frac{1}{4}} \cF_1(v(s))\|_\cW \,ds
\\
&\quad + \int_0^t \|\cA^{\delta+\frac{1}{4}} e^{-\cA (t - s)}\cA^{-\frac{1}{4}} [ \cF_2(u_1(s) + v(s)) - \cF_2(u_1(s)) ] \|_\cW \,ds
\\
&\quad + \int_0^t \|\cA^{\delta+\frac{1}{4}} e^{-\cA (t - s)}\cA^{-\frac{1}{4}} [ \cF_3(u_1(s) + v(s)) - \cF_3(u_1(s)) ] \|_\cW \,ds
\\
&\leq M_\delta t^{-\delta}\|u_0 - \tilde u_0\|_\cW + M_\delta \int_0^t (t-s)^{-\delta} \left( \|f_0\|_\cW + \mu_1\|u_1(s)\|_\cW \right) \,ds
\\
&\quad + M_{\delta+\frac{1}{4}} \int_0^t (t-s)^{-\delta - \frac{1}{4}} \left[\mu_2\|u_1(s)\|_\calV\|u_1(s)\|_{\calV^{\frac{1}{2}}}
 + \mu_3 \|u_1(s)\|_{\calV^{\frac{1}{2}}}^3 \right.
\\
&\quad + \mu_1\|v(s)\|_\cW + \eta_2\left( \|u_1(s)\|_\calV + \|v(s)\|_\calV \right) \|v(s)\|_{\calV^{\frac{1}{2}}}
\\
&\quad + \eta_2\left( \|u_1(s)\|_{\calV^{\frac{1}{2}}} + \|v(s)\|_{\calV^{\frac{1}{2}}} \right) \|v(s)\|_\calV
\\
&\quad + \left. \eta_3\left( \|u_1(s)\|_{\calV^{\frac{1}{2}}}^2 + \|v(s)\|_{\calV^{\frac{1}{2}}}^2) \|v(s)\|_{\calV^{\frac{1}{2}}} \right) \right]\,ds,
\end{align*}
where, in obtaining the second inequality, we again used the decomposition of $\cF$ in \eqref{eq:Kozono_Maeda_Naito_3-2_nonlinearity}, the inequalities in \eqref{eq:Kozono_Maeda_Naito_3-7} and \eqref{eq:Kozono_Maeda_Naito_3-8}, and appealed to
Theorem \ref{thm:Sell_You_37-5} to estimate,
\begin{align*}
\int_0^t \|\cA^\delta e^{-\cA (t - s)} \cF(u_1(s))\|_\cW \,ds
&\leq
\int_0^t \|\cA^\delta e^{-\cA (t - s)} ( f_0 + \cF_1(u_1(s)) ) \|_\cW \,ds
\\
&\quad + \int_0^t \|\cA^{\delta+\frac{1}{4}} e^{-\cA (t - s)}\cA^{-\frac{1}{4}} ( \cF_2(u_1(s)) + \cF_3(u_1(s)) ) \|_\cW \,ds
\\
&\leq
M_\delta \int_0^t (t-s)^{-\delta}\|f_0 + \cF_1(u_1(s))\|_\cW \,ds
\\
&\quad + M_{\delta+\frac{1}{4}} \int_0^t (t-s)^{-\delta - \frac{1}{4}} \|\cA^{-\frac{1}{4}} ( \cF_2(u_1(s)) + \cF_3(u_1(s)) ) \|_\cW \,ds.
\end{align*}
Therefore, for all $\delta \in [0, \frac{3}{4})$ and $t \in (0, \tau]$ and $v \in \fV$,
\begin{align*}
{}&\|\cA^\delta \hat v(t)\|_\cW
\\
&\quad \leq M_\delta t^{-\delta}\|u_0 - \tilde u_0\|_\cW + \left(\|f_0\|_\cW + \mu_1\sup_{s\in(0,\tau)}\|u_1(s)\|_\cW \right) M_\delta \int_0^t (t-r)^{-\delta}\,dr
\\
&\qquad + \sup_{s\in(0,\tau)}\left[\mu_2\|s^{\frac{1}{2}} u_1(s)\|_\calV\|s^{\frac{1}{4}} u_1(s)\|_{\calV^{\frac{1}{2}}}
+ \mu_3 \|s^{\frac{1}{4}} u_1(s)\|_{\calV^{\frac{1}{2}}}^3 \right.
\\
&\qquad\quad + \mu_1\|s^{\frac{3}{4}} v(s)\|_\cW + \eta_2\left(\|s^{\frac{1}{2}} u_1(s)\|_\calV + \|s^{\frac{1}{2}} v(s)\|_\calV\right) \|s^{\frac{1}{4}} v(s)\|_{\calV^{\frac{1}{2}}}
\\
&\qquad\quad + \eta_2\left(\|s^{\frac{1}{4}} u_1(s)\|_{\calV^{\frac{1}{2}}}  + \|s^{\frac{1}{4}} v(s)\|_{\calV^{\frac{1}{2}}} \right) \|s^{\frac{1}{2}} v(s)\|_\calV
\\
&\qquad\quad + \left. \eta_3\left(\|s^{\frac{1}{4}} u_1(s)\|_{\calV^{\frac{1}{2}}}^2 + \|s^{\frac{1}{4}} v(s)\|_{\calV^{\frac{1}{2}}}^2\right) \|s^{\frac{1}{4}} v(s)\|_{\calV^{\frac{1}{2}}} \right]
M_{\delta+\frac{1}{4}} \int_0^t (t-r)^{-\delta-\frac{1}{4}} r^{-\frac{3}{4}}\,dr.
\end{align*}
We now use the facts that \cite[Section 5.12]{Olver_Lozier_Boisvert_Clark}
$$
\int_0^t (t-r)^{-\delta-\frac{1}{4}} r^{-\frac{3}{4}}\,dr = t^{-\delta}B\left(\textstyle\frac{3}{4} - \delta, \textstyle\frac{1}{4}\right), \quad\forall\, t > 0,
$$
and $u_1 \in C([0,\infty; \calV) \cap C([0,\infty; \calV^{\frac{1}{2}}) = C([0,\infty; \calV)$. Hence, keeping in mind the definition \eqref{eq:Kozono_Maeda_Naito_Banach_space_plus_uniqueness} of the Banach space, $\fV$, we see that, for all $\delta \in [0, \frac{3}{4})$ and $t \in (0, \tau]$ and $v \in \fV$,
\begin{align*}
{}&\|\cA^\delta \hat v(t)\|_\cW
\\
&\quad\leq M_\delta t^{-\delta}\|u_0 - \tilde u_0\|_\cW + \frac{M_\delta t^{1-\delta}}{1-\delta}\left( \|f_0\|_\cW + \mu_1\|u_1\|_{C[0,\tau]; \cW)} \right)
\\
&\qquad + t^{-\delta} \tau^{\frac{3}{4}} \left( \mu_2\|u_1\|_{C([0,\tau]; \calV)}\|u_1\|_{C([0,\tau]; \calV^{\frac{1}{2}})}
+ \mu_3\|u_1\|_{C([0,\tau]; \calV^{\frac{1}{2}})}^3 \right) M_{\delta+\frac{1}{4}} B\left(\textstyle\frac{3}{4} - \delta, \textstyle\frac{1}{4}\right)
\\
&\qquad + t^{-\delta}\left( \mu_1\tau^{\frac{3}{4}} \|v\|_\fV + \eta_2\tau^{\frac{1}{2}} \|u_1\|_{C([0,\tau]; \calV)} \|v\|_\fV
 + \eta_2\tau^{\frac{1}{4}} \|u_1\|_{C([0,\tau]; \calV^{\frac{1}{2}})} \|v\|_\fV + 2\eta_2\|v\|_\fV^2 \right.
\\
&\quad\qquad + \left. \eta_3\tau^{\frac{1}{2}} \|u_1\|_{C([0,\tau]; \calV^{\frac{1}{2}})}^2 \|v\|_\fV + \eta_3\|v\|_\fV^3 \right)
M_{\delta+\frac{1}{4}} B\left(\textstyle\frac{3}{4} - \delta, \textstyle\frac{1}{4}\right),
\end{align*}
and thus,
\begin{multline*}
\sup_{t \in (0, \tau)} t^\delta \|\cA^\delta \hat v(t)\|_\cW
\leq
M_\delta \|u_0 - \tilde u_0\|_\cW + \frac{M_\delta t^{1-\delta}}{1-\delta}\left( \|f_0\|_\cW + \mu_1\|u_1\|_{C[0,\tau]; \cW)} \right)
\\
+ \tau^{\frac{3}{4}} \left( \mu_2\|u_1\|_{C([0,\tau]; \calV)}\|u_1\|_{C([0,\tau]; \calV^{\frac{1}{2}})}
+ \mu_3\|u_1\|_{C([0,\tau]; \calV^{\frac{1}{2}})}^3 \right) M_{\delta+\frac{1}{4}} B\left(\textstyle\frac{3}{4} - \delta, \textstyle\frac{1}{4}\right)
\\
+ \left( \mu_1\tau^{\frac{3}{4}} \|v\|_\fV + \eta_2\tau^{\frac{1}{2}} \|u_1\|_{C([0,\tau]; \calV)}\|v\|_\fV + \eta_2\tau^{\frac{1}{4}} \|u_1\|_{C([0,\tau]; \calV^{\frac{1}{2}})}\|v\|_\fV
+ 2\eta_2\|v\|_\fV^2 \right.
\\
+ \left. \eta_3\tau^{\frac{1}{2}} \|u_1\|_{C([0,\tau]; \calV^{\frac{1}{2}})}^2 \|v\|_\fV + \eta_3\|v\|_\fV^3 \right)
M_{\delta+\frac{1}{4}} B\left(\textstyle\frac{3}{4} - \delta, \textstyle\frac{1}{4}\right).
\end{multline*}
In particular, we may choose $\delta \in  \{0, \frac{1}{4}, \frac{1}{2}\}$ and hence
\begin{multline}
\label{eq:Kozono_Maeda_Naito_lemma_3-4_plus_uniqueness_apriori_estimate_preliminary}
\|\sT v\|_\fV
\leq
\max_{\delta \in  [0, \delta_0]} M_\delta \|u_0 - \tilde u_0\|_\cW
+ \max_{\delta \in  [0, \delta_0]} \frac{M_\delta \tau^{1-\delta}}{1-\delta}\left( \|f_0\|_\cW + \mu_1\|u_1\|_{C[0,\tau]; \cW)} \right)
\\
+ \left[ \tau^{\frac{3}{4}} \left( \mu_2\|u_1\|_{C([0,\tau]; \calV)}\|u_1\|_{C([0,\tau]; \calV^{\frac{1}{2}})}
+ \mu_3\|u_1\|_{C([0,\tau]; \calV^{\frac{1}{2}})}^3 \right) \right.
\\
+ \mu_1\tau^{\frac{3}{4}} \|v\|_\fV + \eta_2\tau^{\frac{1}{2}} \|u_1\|_{C([0,\tau]; \calV)}\|v\|_\fV + \eta_2\tau^{\frac{1}{4}} \|u_1\|_{C([0,\tau]; \calV^{\frac{1}{2}})}\|v\|_\fV
\\
+ \left. 2\eta_2\|v\|_\fV^2 + \eta_3\tau^{\frac{1}{2}} \|u_1\|_{C([0,\tau]; \calV^{\frac{1}{2}})}^2 \|v\|_\fV + \eta_3\|v\|_\fV^3 \right]
\max_{\delta \in  [0, \delta_0]} M_{\delta+\frac{1}{4}} B\left(\textstyle\frac{3}{4} - \delta, \textstyle\frac{1}{4}\right).
\end{multline}
Therefore, by definition of $\fV$ in \eqref{eq:Kozono_Maeda_Naito_Banach_space_plus_uniqueness}, and $R$ in \eqref{eq:Kozono_Maeda_Naito_radius_plus_uniqueness_bound_and_Lipschitz}, the fact that $u_0$ and $\tilde u_0$ obey \eqref{eq:Sell_You_47-6_Kozono_Maeda_Naito_plus_uniqueness_choice_of_tilde_u0}, and $\tau$ obeys \eqref{eq:Sell_You_47-6_Kozono_Maeda_Naito_plus_uniqueness_bound}, and $v \in \fF$ if and only if $\|v\|_\fV \leq R$, we obtain
$$
\|\sT v\|_\fV \equiv \|\hat v\|_\fV \leq \frac{R}{2} + \frac{R}{16} + \frac{R}{4} \leq R, \quad \forall\, v \in \fF.
$$
Note that $R$ now depends on $u_0 \in \cW$ and \emph{not} just an upper bound on $\|u_0\|_\cW$, since we needed to choose $\tilde u_0 \in \calV$ satisfying \eqref{eq:Sell_You_47-6_Kozono_Maeda_Naito_plus_uniqueness_choice_of_tilde_u0}, and using
\begin{equation}
\label{eq:u_1_norm_C_0_to_T_into_V2delta}
\|u_1\|_{C([0, T]; \calV^{2\delta})} = \sup_{t\in [0, T]}\|e^{-\cA t}\tilde u_0\|_{\calV^{2\delta}} \leq M\|\tilde u_0\|_{\calV^{2\delta}},
\quad\forall\, \delta \in \{0, \textstyle\frac{1}{4}, \textstyle\frac{1}{2}\} \hbox{ and } T>0,
\end{equation}
where $M \geq 1$ is a constant (depending only on $\cA$) and the inequalities follow from Lemma \ref{lem:Sell_You_37-4} and Theorem \ref{thm:Sell_You_37-5}, since $e^{-\cA t}$ is an analytic semigroup on $\calV^\alpha$ for any $\alpha \geq 0$. Consequently, aside from the term $M_\delta \|u_0 - \tilde u_0\|_\cW$, the remainder of the right-hand side of the inequality \eqref{eq:Kozono_Maeda_Naito_lemma_3-4_plus_uniqueness_apriori_estimate_preliminary} obeys, for $\delta \in  [0, \delta_0]$,
\begin{multline*}
\frac{M_\delta \tau^{1-\delta}}{1-\delta}\left( \|f_0\|_\cW + \mu_1M\|\tilde u_0\|_\cW \right)
+ \left[ \tau^{\frac{3}{4}} \left( \mu_2M^2 \|\tilde u_0\|_\calV \|\tilde u_0\|_{\calV^{\frac{1}{2}}}
+ \mu_3M^3 \|\tilde u_0\|_{\calV^{\frac{1}{2}}}^3 \right) \right.
\\
+ \mu_1\tau^{\frac{3}{4}} R + \eta_2 \tau^{\frac{1}{2}} R M\|\tilde u_0\|_\calV + \eta_2 \tau^{\frac{1}{4}} R M\|\tilde u_0\|_{\calV^{\frac{1}{4}}}
+ \eta_3\tau^{\frac{1}{2}} R M^2\|\tilde u_0\|_{\calV^{\frac{1}{2}}}^2
\\
+ \left. \left( 2\eta_2 R^2 + \eta_3 R^3 \right) \right] M_{\delta+\frac{1}{4}} B\left(\textstyle\frac{3}{4} - \delta, \textstyle\frac{1}{4}\right),
\end{multline*}
and this explains our choice of $R$ in \eqref{eq:Kozono_Maeda_Naito_radius_plus_uniqueness_bound_and_Lipschitz} and $\tau$ in \eqref{eq:Sell_You_47-6_Kozono_Maeda_Naito_plus_uniqueness_bound}. This completes Step \ref{step:Kozono_Maeda_Naito_lemma_3-4_proof_of_boundedness}.
\end{step}

\begin{step}[Contraction mapping property of $\sT$]
\label{step:Kozono_Maeda_Naito_lemma_3-4_plus_uniqueness_proof_of_contraction property}
For any $v_1, v_2 \in \fV$, we set $\hat v_i = \sT v_i$ for $i = 1, 2$ and observe that the definition of $\sT$ in \eqref{eq:Sell_You_page_233_contraction_map_definition_Kozono_Maeda_Naito_plus_uniqueness} yields
$$
\hat v_1(t) - \hat v_2(t) = \int_0^t e^{-\cA (t - s)}[ \cF(u_1(s) + v_1(s)) - \cF(u_1(s) + v_2(s)) ] \,ds,
\quad\forall\, t \in [0, \tau].
$$
The expression \eqref{eq:Kozono_Maeda_Naito_3-2_nonlinearity} for $\cF$ gives
\begin{multline*}
\cF(u_1(s) + v_1(s)) - \cF(u_1(s) + v_2(s))
\\
=
\cF_1(v_1(s) - v_2(s)) + \cF_2(u_1(s) + v_1(s)) - \cF_2(u_1(s) + v_2(s))
\\
+ \cF_3(u_1(s) + v_1(s)) - \cF_3(u_1(s) + v_2(s)).
\end{multline*}
Therefore, applying the Lipschitz inequalities \eqref{eq:Kozono_Maeda_Naito_3-7_linear} and \eqref{eq:Kozono_Maeda_Naito_3-8} yields, for all $\delta \in [0, \frac{3}{4})$ and $t \in (0, \tau]$ and $v_1, v_2 \in \fV$,
\begin{align*}
\|\cA^\delta (\hat v_1(t) - \hat v_2(t))\|_\cW
& \leq M_{\delta+\frac{1}{4}} \int_0^t (t-s)^{-\delta - \frac{1}{4}} \left[ \mu_1\|v_1(s) - v_2(s)\|_\cW \right.
\\
&\quad + \eta_2\left( \|u_1(s) + v_1(s)\|_\calV + \|u_1(s) + v_2(s)\|_\calV \right) \|v_1(s) - v_2(s)\|_{\calV^{\frac{1}{2}}}
\\
&\quad + \eta_2\left( \|u_1(s) + v_1(s)\|_{\calV^{\frac{1}{2}}} + \|u_1(s) + v_2(s)\|_{\calV^{\frac{1}{2}}} \right)
\|v_1(s) - v_2(s)\|_\calV
\\
&\quad + \left. \eta_3\left( \|u_1(s) + v_1(s)\|_{\calV^{\frac{1}{2}}}^2 + \|u_1(s) + v_2(s)\|_{\calV^{\frac{1}{2}}}^2 \right)
\|v_1(s) - v_2(s)\|_{\calV^{\frac{1}{2}}} \right] \,ds.
\end{align*}
Consequently, by again inserting appropriate powers of $s \in [0, t]$ in each of the integral factors, we find that, for all $\delta \in [0, \frac{3}{4})$ and $t \in (0, \tau]$ and $v_1, v_2 \in \fV$,
\begin{align*}
{}&\|\cA^\delta (\hat v_1(t) - \hat v_2(t))\|_\cW
\\
&\quad \leq \sup_{s \in [0,t]} \left[ \mu_1 s^{\frac{3}{4}}\|v_1(s) - v_2(s)\|_\cW \right.
\\
&\qquad + \eta_2\left( 2s^{\frac{1}{2}} \|u_1(s)|_\calV + s^{\frac{1}{2}}\|v_1(s)\|_\calV + s^{\frac{1}{2}} \|v_2(s)\|_\calV \right)
s^{\frac{1}{4}} \|v_1(s) - v_2(s)\|_{\calV^{\frac{1}{2}}}
\\
&\qquad + \eta_2\left( 2s^{\frac{1}{4}} \|u_1(s)\|_{\calV^{\frac{1}{2}}} + s^{\frac{1}{4}}\|v_1(s)\|_{\calV^{\frac{1}{2}}} + s^{\frac{1}{4}} \|v_2(s)\|_{\calV^{\frac{1}{2}}} \right)
s^{\frac{1}{2}} \|v_1(s) - v_2(s)\|_\calV
\\
&\qquad + \left. \eta_3\left( \left(s^{\frac{1}{4}} \|u_1(s) + v_1(s)\|_{\calV^{\frac{1}{2}}}\right)^2
+ \left(s^{\frac{1}{4}} \|u_1(s) + v_2(s)\|_{\calV^{\frac{1}{2}}} \right)^2 \right)
s^{\frac{1}{4}} \|v_1(s) - v_2(s)\|_{\calV^{\frac{1}{2}}} \right]
\\
&\qquad \times t^{-\delta} M_{\delta+\frac{1}{4}} B\left(\textstyle\frac{3}{4} - \delta, \textstyle\frac{1}{4}\right).
\end{align*}
By again appealing to the fact that $u_1 \in C([0,\infty; \calV) \cap C([0,\infty; \calV^{\frac{1}{2}})$ and applying the definition \eqref{eq:Kozono_Maeda_Naito_Banach_space_plus_uniqueness} of the norm on $\fV$, we discover that, for all $\delta \in [0, \frac{3}{4})$ and $t \in (0, \tau]$ and $v_1, v_2 \in \fV$,
\begin{align*}
t^\delta \|\cA^\delta (\hat v_1(t) - \hat v_2(t))\|_\cW
&\leq
\left[ \mu_1 \tau^{\frac{3}{4}} + \eta_2\left( 2\tau^{\frac{1}{2}} \|u_1\|_{C([0,\infty; \calV)} + \|v_1\|_\fV + \|v_2\|_\fV \right) \right.
\\
&\quad + \eta_2\left( 2\tau^{\frac{1}{4}} \|u_1\|_{C([0,\infty; {\calV^{\frac{1}{2}}})} + \|v_1\|_\fV + \|v_2\|_\fV \right)
\\
&\quad + \left. 2\eta_3\left( 2\tau^{\frac{1}{2}} \|u_1\|_{C([0,\infty; {\calV^{\frac{1}{2}}})}^2 + \|v_1\|_\fV^2
+ \|v_2\|_\fV^2 \right) \right]
\\
&\qquad \times M_{\delta+\frac{1}{4}} B\left(\textstyle\frac{3}{4} - \delta, \textstyle\frac{1}{4}\right)\|v_1 - v_2\|_\fV.
\end{align*}
In particular, we may choose $\delta \in  \{0, \frac{1}{4}, \frac{1}{2}\}$ and, recalling that $\hat v_i = \sT v_i$ for $i = 1, 2$ and the definition \eqref{eq:Kozono_Maeda_Naito_Banach_space_plus_uniqueness} of the norm on $\fV$, we obtain
\begin{multline}
\label{eq:Kozono_Maeda_Naito_lemma_3-4_plus_uniqueness_contraction_mapping_preliminary}
\|\sT v_1 - \sT v_2\|_\fV
\leq
\left[ \mu_1 \tau^{\frac{3}{4}} + \eta_2\left( 2\tau^{\frac{1}{2}} \|u_1\|_{C([0,\infty; \calV)} + \|v_1\|_\fV + \|v_2\|_\fV \right) \right.
\\
+ \eta_2\left( 2\tau^{\frac{1}{4}} \|u_1\|_{C([0,\infty; {\calV^{\frac{1}{2}}})} + \|v_1\|_\fV + \|v_2\|_\fV \right)
\\
+ \left. 2\eta_3\left( 2\tau^{\frac{1}{2}} \|u_1\|_{C([0,\infty; {\calV^{\frac{1}{2}}})}^2 + \|v_1\|_\fV^2
+ \|v_2\|_\fV^2 \right) \right]
\\
\times \max_{\delta \in  [0, \delta_0]} M_{\delta+\frac{1}{4}} B\left(\textstyle\frac{3}{4} - \delta, \textstyle\frac{1}{4}\right) \|v_1 - v_2\|_\fV.
\end{multline}
Therefore, by definition of $R$ in \eqref{eq:Kozono_Maeda_Naito_radius_plus_uniqueness_bound_and_Lipschitz}, the bounds on $u_1$ in \eqref{eq:u_1_norm_C_0_to_T_into_V2delta}, the fact that $\tau$ obeys \eqref{eq:Sell_You_47-6_Kozono_Maeda_Naito_plus_uniqueness_Lipschitz} and $v_1, v_2 \in \fF$ if and only if $\|v_i\|_\fV \leq R$ for $i = 1, 2$, we see that
$$
\|\sT v_1 - \sT v_2\|_\fV \leq \left(\frac{1}{4} + \frac{1}{4}\right)\|v_1 - v_2\|_\fV = \frac{1}{2}\|v_1 - v_2\|_\fV, \quad \forall\, v_1, v_2 \in \fF,
$$
as desired for our application of the Banach Contraction Mapping Theorem.
\end{step}

\begin{step}[Continuity of $\hat v(t)$ for $0\leq t\leq \tau$ and existence and uniqueness of the solution]
\label{step:Kozono_Maeda_Naito_lemma_3-4_plus_uniqueness_proof_of_continuity}
The proof of continuity of continuity of $\hat v(t)$ with respect to $t \in [0, \tau]$ follows the pattern of Step \ref{step:Kozono_Maeda_Naito_lemma_3-4_proof_of_continuity} in the proof of Theorem \ref{thm:Kozono_Maeda_Naito_lemma_3-4}, \emph{mutatis mutandis}, except that we now choose a sequence, $\{\tilde u_0^m\}_{m \in \NN} \subset \calV^2$, such that $\tilde u_0^m \to \tilde u_0$ in $\calV^2$ as $m \to \infty$. As in the proof of Theorem \ref{thm:Kozono_Maeda_Naito_lemma_3-4}, we now obtain existence and uniqueness of a solution $v \in \fV$ by the Banach Contraction Mapping Theorem and the proof of continuity of $\hat v(t)$ with respect to $t \in [0, \tau]$ also yields the regularity result,
$$
v \in C([0, \tau]; \cW) \cap C((0, \tau]; \calV^{2\delta_0}).
$$
Since $u = e^{-cA t}\tilde u_0 + v$ and $\tilde u_0 \in \calV^{2\delta_0}$, we see that $u$ has the same regularity as $v$, that is, $u$ obeys \eqref{eq:Sell_You_47-5_Kozono_Maeda_Naito_lemma_3-4_plus_uniqueness}.
\end{step}

\begin{step}[\Apriori estimate]
\label{step:Kozono_Maeda_Naito_lemma_3-4_plus_uniqueness_apriori_and_continuity_estimates}
The \apriori estimate \eqref{eq:Kozono_Maeda_Naito_lemma_3-4_plus_uniqueness_apriori_estimate_with_rearrangement} follows immediately from our derivation of \eqref{eq:Kozono_Maeda_Naito_lemma_3-4_plus_uniqueness_apriori_estimate_preliminary} by taking $C_0 = R$ and noting that $v = \sT v$ because the solution $v$ is a fixed point of the mapping $\sT:\fF\to\fF$ and recalling that $u = e^{-cA t}\tilde u_0 + v$.
\end{step}

This completes the proof of Theorem \ref{thm:Kozono_Maeda_Naito_lemma_3-4_plus_uniqueness}.
\end{proof}

\chapter[Elliptic and parabolic partial differential systems on manifolds]{A priori estimates, existence, uniqueness, and regularity for elliptic and parabolic partial differential systems on manifolds}
\label{chapter:Elliptic_and_parabolic_partial_differential_systems}

\section[Elliptic partial differential systems and analytic semigroups]{Elliptic partial differential systems and analytic semigroups on $L^p$, $C^0$, and $L^1$ Banach spaces}
\label{sec:Sell_You_3-8-2_standard_Sobolev_spaces}
In order to bring the abstract theory of evolution equations in Banach spaces which we have discussed thus far to bear on the Yang-Mills heat equation \eqref{eq:Yang-Mills_heat_equation_as_perturbation_rough_Laplacian_plus_one_heat_equation}
or its linearization, we will need show that the Laplace operator, $\Delta_A$ on $\Omega^1(X; \ad P)$, appearing in those parabolic equations has realizations in various useful (and not necessarily standard) Sobolev spaces which are sectorial in the sense of Definition \ref{defn:Sell_You_page_78_definition_of_sectorial_operator} and hence determine infinitesimal generators, $-\Delta_A$, of analytic semigroups, $e^{-\Delta_A t}$, by Theorem \ref{thm:Renardy_Rogers_12-31}.

While there are several well-known treatments of \apriori $L^p$ estimates and existence, uniqueness, and regularity results for elliptic partial differential systems, none are entirely suited to our applications in this monograph or as comprehensive and well-developed as their counterparts for scalar, second-order elliptic partial differential operators on open subsets of $\RR^d$, such as can be found in the references due to Gilbarg and Trudinger \cite{GilbargTrudinger} or Krylov \cite{Krylov_LecturesSobolev}. When one in addition requires results for elliptic partial differential operators on vector bundles over closed manifolds and resolvent estimates and analytic semigroup generation results on $L^p$, $C^0$, and $L^1$ spaces, then suitable references become even harder to locate. For this reason, we shall provide in this section a largely self-contained review together with our own further development of the results we shall need in this monograph for elliptic partial differential systems, building on relatively recent articles due to Cannarsa, Terreni, and Vespri \cite{Cannarsa_Terreni_Vespri_1985} and Denk, Haller-Dintelmann, Heck, Hieber, Pr\"uss, and Stavrakidis \cite{Denk_Hieber_Pruss_2003, Haller-Dintelmann_Heck_Hieber_2006, Heck_Hieber_Stavrakidis_2010}. In addition, when this can be done easily in a self-contained manner, we extend existing \apriori $L^p$ estimates and existence, uniqueness, and regularity results described by Krylov \cite{Krylov_LecturesSobolev} for scalar elliptic partial differential operators on open subsets of $\RR^d$ to the case of certain elliptic partial differential systems. Our review and development builds on references such as those of Agmon \cite[Section 6]{AgmonLecturesEllipticBVP}, Agmon, Douglis, and Nirenberg \cite{AgmonDouglisNirenberg2} Chen and Wu \cite[Part II]{Chen_Wu_1998}, Ladyzhenskaya and Ural$'$tseva \cite[Chapter 7]{Ladyzhenskaya_Uraltseva_1968}, and Morrey \cite[Chapter 6]{Morrey}.

For resolvent estimates and results on analytic semigroup generation on $L^p$, $C^0$, and $L^1$ spaces determined by elliptic partial differential operators, the theory is by far most well-developed in the case of scalar elliptic partial differential operators of order $m \geq 1$ on open subsets of $\RR^d$, beginning with early work of Browder \cite{Browder_1961} when $m =2$ and $p = 2$, Agmon for even $m \geq 2$ and $1 < p < \infty$ \cite{Agmon_1962}, and Stewart \cite{Stewart_1974, Stewart_1980} even $m \geq 2$ and $p = \infty$ via the Stewart-Masuda method. Expositions of their and related further results are due to Jacob \cite{Jacob_v1, Jacob_v2, Jacob_v3}, Pazy \cite[Section 7.3]{Pazy_1983}, Sell and You \cite[Section 3.8.2]{Sell_You_2002}, and Tanabe \cite{Tanabe_1979, Tanabe_1997}. Grubb \cite{Grubb_1996} and Jacob \cite{Jacob_v1, Jacob_v2, Jacob_v3} provide resolvent estimates and analytic semigroup generation results for pseudo-differential operators.

Although the \apriori $L^p$ estimates and existence, uniqueness, and regularity results for elliptic partial differential operators which we describe in this section provide the analytical foundation for our monograph, our treatment is developed with the goal of proving the existence of analytic semigroups on useful Banach spaces, with infinitesimal generators given by elliptic partial differential operators. Expositions of the abstract treatment of sectorial operators and analytic semigroups can be found in references due to Banasiak and Arlotti \cite{Banasiak_Arlotti_2006}, Engel and Nagel \cite{Engel_Nagel_2000, Engel_Nagel_2006}, Jacob \cite{Jacob_v1, Jacob_v2, Jacob_v3}, Lorenzi and Bertoldi \cite{Lorenzi_Bertoldi_2007}, \cite{Lunardi_1995}, Pazy \cite{Pazy_1983}, Sell and You \cite{Sell_You_2002}, Tanabe \cite{Tanabe_1979, Tanabe_1997}, and Vrabie \cite{Vrabie_2003}.

Our development in this section focuses on elliptic partial differential operators acting on standard Sobolev spaces but, as we shall see in Section \ref{sec:Sell_You_3-8-2_critical_exponent_elliptic_Sobolev_spaces}, we shall be able to extend the theory in this section in a natural way there to provide \apriori estimates, existence and uniqueness regularity results, resolvent estimates, and analytic semigroup generation results for second-order elliptic partial differential operators, with scalar principal symbol, acting on certain Banach spaces --- whose definitions are inspired by ideas of Taubes --- of sections of vector bundles over closed manifolds.


\subsection{Sobolev embedding and multiplication theorems for real derivative exponents}
\label{subsec:Sobolev_embedding_multiplication}
Suppose we take $\cW = L^p(X;\RR)$ and $\cA = \Delta + 1$, where $X$ is a closed, Riemannian, smooth dimensional manifold and $\Delta = d^*d: C^\infty(X) \to C^\infty(X)$ is the Laplace operator (on scalar functions) defined by the Riemannian metric, in which case $\calV = \sD(\cA) = W^{2,p}(X)$ and
$\calV^s = W^{s, p}(X)$, where the Sobolev spaces $W^{s, p}(X)$, for $1<p<\infty$ and $s\in\RR$, may be defined as in \cite[Section 13.6]{Taylor_PDE3}.

We now review several different cases of the `Sobolev multiplication theorems' described by Freed and Uhlenbeck \cite[pp. 95--96]{FU} and Palais \cite[Section 9]{PalaisFoundationGlobal}.
We shall state all Sobolev embedding and multiplication results for the case real or complex-valued functions on a closed, Riemannian, smooth manifold and so, for example, $L^p(X)$ may denote $L^p(X;\RR)$ or $L^p(X;\CC)$, but all of these results extend to the case of sections of real Riemannian or complex Hermitian vector bundles over $X$, with pointwise scalar multiplication, $(u_1, u_2) \mapsto u_1u_2$, replaced by pointwise tensor product,  $(u_1, u_2) \mapsto u_1\otimes u_2$.

\begin{lem}[Sobolev multiplication theorem for nonnegative integer derivative exponents and neither factor in continuous range]
\label{lem:Freed_Uhlenbeck_equation_6-34}
\cite[Equation (6.34)]{FU}, \cite[Section 9]{PalaisFoundationGlobal}
Let $X$ be a closed, Riemannian, smooth manifold of dimension $d \geq 2$. If $k, k_1, k_2$ and $p, p_1, p_2$ satisfy
\begin{equation}
\label{eq:Freed_Uhlenbeck_6-34_hypothesis}
\begin{gathered}
k, k_1, k_2 \in \NN \hbox{ with } k_1, k_2 \geq k, \quad p_1, p_2 > 1, \quad 1 \leq p < \infty, \quad p_1k_1, p_1k_2 < d,
\\
\hbox{and}\quad k_1 - \frac{d}{p_1} + k_2 - \frac{d}{p_2} \geq k - \frac{d}{p},
\end{gathered}
\end{equation}
then the following multiplication map is defined and continuous,
\begin{equation}
\label{eq:Freed_Uhlenbeck_6-34}
W^{k_1,p_1}(X) \times W^{k_2,p_2}(X) \to W^{k,p}(X).
\end{equation}
\end{lem}

Clearly, extensions of Lemma \ref{lem:Freed_Uhlenbeck_equation_6-34} to arbitrarily many multiplicative factors or nonnegative real derivative exponents are possible, just as in \cite[Section 9]{PalaisFoundationGlobal}. However, the following extension will be especially useful and illustrates the method of proof of more general results along these lines.

\begin{lem}[Sobolev multiplication theorem for three factors, nonnegative integer derivative exponents, and no factor in continuous range]
\label{lem:Freed_Uhlenbeck_equation_6-34_nonnegative_cubic}
Let $X$ be a closed, Riemannian, smooth manifold of dimension $d \geq 2$. If $k, k_i \in \ZZ$ and $p, p_i \geq 1$, for $i=1,2,3$, satisfy
\begin{equation}
\label{eq:Freed_Uhlenbeck_6-34_nonnegative_cubic_hypothesis}
\begin{gathered}
k_i \geq k \geq 0, \quad p_i > 1, \quad p < \infty,
\quad p_ik_i < d \quad\hbox{for } i = 1,2,3,
\\
\hbox{and}\quad k_1 - \frac{d}{p_1} + k_2 - \frac{d}{p_2} + k_3 - \frac{d}{p_3} \geq k - \frac{d}{p},
\end{gathered}
\end{equation}
then the following multiplication map is defined and continuous,
\begin{equation}
\label{eq:Freed_Uhlenbeck_6-34_nonnegative_cubic}
W^{k_1,p_1}(X) \times W^{k_2,p_2}(X) \times W^{k_3,p_3}(X) \to W^{k,p}(X).
\end{equation}
\end{lem}

\begin{proof}
Set $m_i := k_i - k$ for $i = 1,2,3$, so the $m_i$ are non-negative integers by our hypotheses on $k$ and the $k_i$. We shall first consider the case $k = 0$. Our hypothesis \eqref{eq:Freed_Uhlenbeck_6-34_nonnegative_cubic_hypothesis} ensures that $m_ip_i < d$, when $k = 0$, and as $p_i > 1$, for $i = 1,2,3$, the Sobolev Embedding Theorem \cite[Theorem 4.12]{AdamsFournier} gives
\begin{equation}
\label{eq:Sobolev_embedding_Wm_ip_i_into_Lq_i}
W^{m_i, p_i}(X) \hookrightarrow L^{q_i}(X), \quad\hbox{for } q_i := \frac{dp_i}{d - m_ip_i} \hbox{ and } i = 1,2,3.
\end{equation}
Assume temporarily that the $q_i$ obey
\begin{equation}
\label{eq:Exponent_condition_for_cubic_generalized_holder}
\frac{1}{q_1} + \frac{1}{q_2} + \frac{1}{q_3} \leq \frac{1}{p}.
\end{equation}
Then the generalized H\"older inequality \cite[Equation (7.11)]{GilbargTrudinger} yields a continuous multiplication map,
$$
L^{q_1}(X) \times L^{q_2}(X) \times L^{q_3}(X) \to L^p(X).
$$
Composing the preceding continuous multiplication map with the Sobolev embeddings \eqref{eq:Sobolev_embedding_Wm_ip_i_into_Lq_i} gives a continuous Sobolev multiplication map,
\begin{equation}
\label{eq:Freed_Uhlenbeck_6-34_nonnegative_cubic_k_is_zero}
W^{m_1, p_1}(X) \times W^{m_2, p_2}(X) \times W^{m_3, p_3}(X) \to L^p(X),
\end{equation}
and this is just the map \eqref{eq:Freed_Uhlenbeck_6-34_nonnegative_cubic} when $k = 0$.

To complete our proof of \eqref{eq:Freed_Uhlenbeck_6-34_nonnegative_cubic} when $k = 0$, it remains to verify that the $q_i$ obey \eqref{eq:Exponent_condition_for_cubic_generalized_holder} for $i = 1,2,3$. Using $d/q_i = d/p_i - m_i$ for $i = 1,2,3$, the main inequality in \eqref{eq:Freed_Uhlenbeck_6-34_nonnegative_cubic_hypothesis} for $k_i$ and $p_i$ can be rewritten as
$$
k_1 - m_1 - \frac{d}{q_1} + k_2 - m_2 - \frac{d}{q_2} + k_3 - m_3 - \frac{d}{q_3} \geq k - \frac{d}{p},
$$
and thus, again because $k_i - m_i = 0$ for $i =  1,2,3$ when $k = 0$, this inequality is equivalent to
$$
- \frac{d}{q_1}- \frac{d}{q_2} - \frac{d}{q_3} \geq - \frac{d}{p},
$$
and this is in turn equivalent to \eqref{eq:Exponent_condition_for_cubic_generalized_holder}.

Continuity of the Sobolev multiplication map \eqref{eq:Freed_Uhlenbeck_6-34_nonnegative_cubic} for arbitrary integer $k \geq 0$ follows by applying \eqref{eq:Freed_Uhlenbeck_6-34_nonnegative_cubic_k_is_zero} to the derivative expressions, $\nabla^j(f_1f_2f_3)$, for $0 \leq j \leq k$, where $f_i \in W^{k+m_i, p_i}(X) = W^{k_i, p_i}(X)$ for $i = 1,2,3$.
\end{proof}

Recall from \cite[Section 3.7]{AdamsFournier} that, given $s \in \RR$ with $s \geq 0$ (Adams and Fournier assume $s \in \NN$) and $p \in (1, \infty)$ and for dual exponent $p' \in (1, \infty)$ given by $1 = 1/p + 1/p'$, one defines
\begin{equation}
W^{-s, p'}(X) := (W^{s, p}(X))',
\end{equation}
where $(W^{s, p}(X))'$ denotes the Banach-space dual of $W^{s, p}(X)$, and one finds that (see \cite[Section 3.14]{AdamsFournier} when $s \in \NN$) the norm on $W^{-s, p'}(X)$ may be characterized by
\begin{equation}
\label{eq:Adams_Fournier_3-14}
\|u\|_{W^{-s, p'}(X)}
=
\sup_{\begin{subarray}{c} v \in W^{s, p}(X) \\ \|v\|_{W^{s, p}(X)} \leq 1 \end{subarray}} (u, v)_{L^2(X)}.
\end{equation}
The relation \eqref{eq:Freed_Uhlenbeck_6-34} extends by duality to give

\begin{lem}[Sobolev multiplication theorem for one negative and one non-negative integer derivative exponent]
\label{lem:Freed_Uhlenbeck_equation_6-34_minus_k_k2}
\cite[p. 96]{FU}, \cite[Section 9]{PalaisFoundationGlobal}
Let $X$ be a closed, Riemannian, smooth manifold of dimension $d \geq 2$. The multiplication map \eqref{eq:Freed_Uhlenbeck_6-34} is defined and continuous for $k, k_1, k_2$ and $p, p_1, p_2$ if and only if following multiplication map is defined and continuous, where $p'$ and $p_2'$ are dual exponents obeying $1/p + 1/p' = 1$ and  $1/p_2 + 1/p_2' = 1$,
\begin{equation}
\label{eq:Freed_Uhlenbeck_6-34_minus_k_k2}
W^{k_1,p_1}(X) \times W^{-k,p}(X) \to W^{-k_2,p_2'}(X).
\end{equation}
\end{lem}

The multiplication theorems simplify in the \emph{continuous range}, $kp>d$.

\begin{lem}[Sobolev multiplication theorem for non-negative integer derivative exponents and two factors in continuous range]
\label{lem:Freed_Uhlenbeck_equation_6-34_algebra}
\cite[Theorem 4.39]{AdamsFournier}, \cite[p. 96]{FU}, \cite[Corollary 9.7]{PalaisFoundationGlobal}
Let $X$ be a closed, Riemannian, smooth manifold of dimension $d \geq 2$. If $k$ and $p$ obey
\begin{equation}
\label{eq:Freed_Uhlenbeck_6-34_algebra_hypothesis}
k \in \NN \quad\hbox{and}\quad kp > d,
\end{equation}
then $W^{k,p}(X)$ is an algebra, that is, the following multiplication map is defined and continuous,
\begin{equation}
\label{eq:Freed_Uhlenbeck_6-34_algebra}
W^{k,p}(X) \times W^{k,p}(X) \to W^{k,p}(X).
\end{equation}
\end{lem}

\begin{lem}[Sobolev multiplication theorem for non-negative integer derivative exponents and one factor in continuous range]
\label{lem:Freed_Uhlenbeck_equation_6-34_module}
\cite[p. 96]{FU}, \cite[Corollary 9.7]{PalaisFoundationGlobal}
Let $X$ be a closed, Riemannian, smooth manifold of dimension $d \geq 2$. If $k, k_1$ and $p$ obey
\begin{equation}
\label{eq:Freed_Uhlenbeck_6-34_module_hypothesis}
k_1,k \in \NN \hbox{ with } k_1\geq k, \quad 1 \leq p < \infty, \quad k_1p > d,
\end{equation}
then $W^{k,p}(X)$ is a $W^{k_1,p}(X)$-module, that is, the following multiplication map is defined and continuous,
\begin{equation}
\label{eq:Freed_Uhlenbeck_6-34_module}
W^{k_1,p}(X) \times W^{k,p}(X) \to W^{k,p}(X).
\end{equation}
\end{lem}

These multiplication results are consequences of the H\"older inequality and the standard Sobolev Embedding Theorem (for integer $k$) \cite[Theorem 4.12]{AdamsFournier}. The latter result may be extended to non-integral $s\in\RR$. Indeed, we recall the

\begin{prop}[Sobolev embedding theorem for real derivative exponents]
\label{prop:Taylor_proposition_13-6-3and4}
\cite{Calderon_1961}, \cite[Theorems 9.1 and 9.2]{PalaisFoundationGlobal}, \cite[Propositions 13.6.3 and 13.6.4]{Taylor_PDE3}
Let $X$ be a closed, Riemannian, smooth manifold of dimension $d \geq 2$. If $p \in (1,\infty)$ and $s \in \RR$, then the following embeddings are continuous,
\begin{equation}
\label{eq:Taylor_proposition_13-6-3and4}
W^{s,p}(X)
\hookrightarrow
\begin{cases}
L^{dp/(d-sp)}(X), &\hbox{for } 0 \leq sp < d,
\\
C(X), &\hbox{for } sp > d.
\end{cases}
\end{equation}
\end{prop}


Palais provides Sobolev multiplication results \cite[Theorems 9.4--9.6]{PalaisFoundationGlobal} for the Sobolev spaces $W^{s,p}(X)$ for non-integral $s\in\RR$, obtaining those from the case of integral $k \in \ZZ$ by interpolation theory. We shall instead derive the few special cases we shall need from the Sobolev embedding \eqref{eq:Taylor_proposition_13-6-3and4} and the previous Sobolev multiplication results for integer $k\in\ZZ$.

\begin{lem}[Sobolev multiplication theorem for nonnegative real derivative exponents and neither factor in continuous range]
\label{lem:Freed_Uhlenbeck_equation_6-34_nonnegative_real}
Let $X$ be a closed, Riemannian, smooth manifold of dimension $d \geq 2$. If $s, s_1, s_2$ and $p, p_1, p_2$ satisfy
\begin{equation}
\label{eq:Freed_Uhlenbeck_6-34_nonnegative_real_hypothesis}
\begin{gathered}
s, s_1, s_2 \in \RR \hbox{ with } s_1, s_2 \geq s \geq 0, \quad p_1, p_2 > 1, \quad 1 \leq p < \infty, \quad p_1s_1, p_1s_2 < d,
\\
\hbox{and}\quad s_1 - \frac{d}{p_1} + s_2 - \frac{d}{p_2} \geq s - \frac{d}{p},
\end{gathered}
\end{equation}
then the following multiplication map is defined and continuous,
\begin{equation}
\label{eq:Freed_Uhlenbeck_6-34_nonnegative_real}
W^{s_1,p_1}(X) \times W^{s_2,p_2}(X) \to W^{s,p}(X).
\end{equation}
\end{lem}

\begin{proof}
First, writing $s = k + \alpha$, where $k \in \NN$ is the largest integer with $k \leq s$ and $\alpha := s - k \in [0,1)$, we find that
\begin{equation}
\label{eq:Taylor_proposition_13-6-4_s_and_2_into_k_and_p}
W^{s, p}(X) = W^{k+\alpha, p}(X) \hookrightarrow W^{k,q}(X),
\end{equation}
by \eqref{eq:Taylor_proposition_13-6-3and4} with $q := dp/(d - \alpha p) \in [p, \infty)$. To see this, recall that for $p \in (1,\infty)$ and $s\in\RR$ (by \cite[Equation (13.6.1) and Proposition 13.6.1]{Taylor_PDE3} and \cite[Section 4.1]{Taylor_PDE1})
\begin{equation}
\label{eq:Taylor_13-6-1}
W^{s,p}(X) := (\Delta + 1)^{-s/2}L^p(X),
\end{equation}
and so the composition,
$$
W^{s, p}(X) = W^{k+\alpha, p}(X) \xrightarrow{(\Delta + 1)^{k/2}} W^{\alpha, p}(X)
\xrightarrow{\text{(by \eqref{eq:Taylor_proposition_13-6-3and4})}} L^q(X) \xrightarrow{(\Delta + 1)^{-k/2}} W^{k,q}(X),
$$
is a continuous map for $q = dp/(d - \alpha p)$. Indeed, by hypothesis, we have $(k+\alpha)=sp < d$ and $k\geq 0$ (since $s \geq 0$ by hypothesis) and hence $\alpha p < d$, so the preceding application of Proposition \ref{prop:Taylor_proposition_13-6-3and4} is valid.

For $i=1,2$ and writing $s_i = k_i + \alpha_i$, where $k_i \in \NN$ are the largest integers with $k_i \leq s_i$ and $\alpha_i := s_i - k_i \in [0,1)$, we see that \eqref{eq:Taylor_proposition_13-6-4_s_and_2_into_k_and_p} yields
\begin{equation}
\label{eq:Taylor_proposition_13-6-4_si_and_2_into_ki_and_pi}
W^{s_i, p_i}(X) = W^{k_i+\alpha_i, p_i}(X) \hookrightarrow W^{k_i,q_i}(X),
\end{equation}
for $q_i := dp_i/(d - \alpha_i p_i) \in [p_i, \infty)$.

Given $s \in \RR$ with $s \geq 0$ and $p \in (1, \infty)$, we define $k \in \NN$ as the smallest integer with $k \geq s$ and $\beta := k - s \in [0,1)$ and define $q \in (1, \infty)$ by $p = dq/(d-\beta q)$, and find that
\begin{equation}
\label{eq:Taylor_proposition_13-6-4_k_and_q_into_s_and_2}
W^{k, q}(X) = W^{s+\beta, q}(X) \hookrightarrow W^{s, p}(X).
\end{equation}
To see this, observe that
$$
W^{k, q}(X) = W^{s+\beta, q}(X) \xrightarrow{(\Delta + 1)^{s/2}} W^{\beta, q}(X)
\xrightarrow{\text{(by \eqref{eq:Taylor_proposition_13-6-3and4})}} L^p(X) \xrightarrow{(\Delta + 1)^{-s/2}} W^{s, p}(X),
$$
is a continuous map, as desired.

We now claim that the $k, k_1, k_2$ and $q, q_1, q_2$ obey the hypotheses in \eqref{eq:Freed_Uhlenbeck_6-34_hypothesis}. Indeed, $k, k_1, k_2 \in \NN$ by construction and the hypothesis $s_1, s_2 \geq s$ implies $k_1, k_2 \geq k$. By construction, the $q, q_1, q_2$ obey $q_i \geq p_i > 1$ for $i=1,2$ and $q < \infty$ and because $s_ip_i = (k_i+\alpha_i)p_i < d$ for $i=1,2$ by the hypothesis \eqref{eq:Freed_Uhlenbeck_6-34_nonnegative_real_hypothesis}, they also obey $k_iq_i < d$ for $i=1,2$. Lastly, using our definitions of $q, q_1, q_2$ via
$$
\frac{1}{q_i} = \frac{1}{p_i} - \frac{\alpha_i}{d} \quad\hbox{for } i = 1, 2 \quad\hbox{and}\quad \frac{1}{p} = \frac{1}{q} - \frac{\beta}{d},
$$
and our definitions of $k, k_1, k_2$ and $\alpha, \alpha_1, \alpha_2$, we see that
\begin{align*}
k_1 - \frac{d}{q_1} + k_2 - \frac{d}{q_2} &= s_1 - \alpha_1 - \frac{d}{p_1} + \alpha_1 + s_2 - \alpha_2 - \frac{d}{p_2} + \alpha_2
\\
&= s_1 - \frac{d}{p_1} + s_2 - \frac{d}{p_2}
\\
&\geq s - \frac{d}{p} \quad\hbox{(by \eqref{eq:Freed_Uhlenbeck_6-34_nonnegative_real_hypothesis})}
\\
&= k - \beta - \frac{d}{q} + \beta
\\
&= k - \frac{d}{q},
\end{align*}
as claimed. Therefore, by composing the continuous maps \eqref{eq:Freed_Uhlenbeck_6-34}, \eqref{eq:Taylor_proposition_13-6-4_si_and_2_into_ki_and_pi}, and \eqref{eq:Taylor_proposition_13-6-4_k_and_q_into_s_and_2}, we have
\begin{align*}
W^{s_1,p_1}(X) \times W^{s_2,p_2}(X) &= W^{k_1+\alpha_1, p_1}(X) \times W^{k_2+\alpha_2, p_2}(X)
\\
&\hookrightarrow W^{k_1,q_1}(X) \times W^{k_1,q_1}(X) \quad\hbox{(by \eqref{eq:Taylor_proposition_13-6-4_si_and_2_into_ki_and_pi})}
\\
&\to W^{k,q}(X) \quad\hbox{(by \eqref{eq:Freed_Uhlenbeck_6-34})}
\\
&= W^{s+\beta, q}(X)
\\
&\hookrightarrow W^{s, p}(X) \quad\hbox{(by \eqref{eq:Taylor_proposition_13-6-4_k_and_q_into_s_and_2})}.
\end{align*}
This yields \eqref{eq:Freed_Uhlenbeck_6-34_nonnegative_real}, as desired.
\end{proof}

Note that the multiplication map \eqref{eq:Freed_Uhlenbeck_6-34_nonnegative_real} is just the map \eqref{eq:Freed_Uhlenbeck_6-34} with $k, k_1, k_2 \in \NN$ replaced by $s, s_1, s_2 \in \RR$ with $s, s_1, s_2 \geq 0$. The following elementary extension of will also be useful.

\begin{lem}[Sobolev multiplication theorem for three factors, nonnegative real derivative exponents, and no factor in continuous range]
\label{lem:Freed_Uhlenbeck_equation_6-34_nonnegative_real_cubic}
Let $X$ be a closed, Riemannian, smooth manifold of dimension $d \geq 2$. If $s, s_i \in \RR$ and $p, p_i \geq 1$, for $i=1,2,3$, satisfy
\begin{equation}
\label{eq:Freed_Uhlenbeck_6-34_nonnegative_real_cubic_hypothesis}
\begin{gathered}
s_i \geq s \geq 0, \quad p_i > 1, \quad p < \infty, \quad p_is_i < d \quad\hbox{for } i=1,2,3,
\\
\hbox{and}\quad s_1 - \frac{d}{p_1} + s_2 - \frac{d}{p_2} + s_3 - \frac{d}{p_3} \geq s - \frac{d}{p},
\end{gathered}
\end{equation}
then the following multiplication map is defined and continuous,
\begin{equation}
\label{eq:Freed_Uhlenbeck_6-34_nonnegative_real_cubic}
W^{s_1,p_1}(X) \times W^{s_2,p_2}(X) \times W^{s_3,p_3}(X) \to W^{s,p}(X).
\end{equation}
\end{lem}

\begin{proof}
The result may be proved by repeating the proof of Lemma \ref{lem:Freed_Uhlenbeck_equation_6-34_nonnegative_real}, replacing role of Lemma \ref{lem:Freed_Uhlenbeck_equation_6-34} (a quadratic multiplication map with nonnegative integer derivative exponents) by Lemma \ref{lem:Freed_Uhlenbeck_equation_6-34_nonnegative_cubic} (a cubic multiplication map with nonnegative integer derivative exponents), \emph{mutatis mutandis}.
\end{proof}

Given Lemma \ref{lem:Freed_Uhlenbeck_equation_6-34_nonnegative_real}, the proof of Lemma \ref{lem:Freed_Uhlenbeck_equation_6-34_minus_k_k2} via duality and the proofs of Lemmata \ref{lem:Freed_Uhlenbeck_equation_6-34_algebra} and \ref{lem:Freed_Uhlenbeck_equation_6-34_module} now extend to give

\begin{lem}[Sobolev multiplication theorem for one negative and one non-negative integer derivative exponent]
\label{lem:Freed_Uhlenbeck_equation_6-34_minus_s_s2_real}
Let $X$ be a closed, Riemannian, smooth manifold of dimension $d \geq 2$. The multiplication map \eqref{eq:Freed_Uhlenbeck_6-34_nonnegative_real} is defined and continuous for $s, s_1, s_2$ and $p, p_1, p_2$ if and only if following multiplication map is defined and continuous, where $p'$ and $p_2'$ are dual exponents obeying $1/p + 1/p' = 1$ and  $1/p_2 + 1/p_2' = 1$,
\begin{equation}
\label{eq:Freed_Uhlenbeck_6-34_minus_s_s2_real}
W^{s_1,p_1}(X) \times W^{-s,p}(X) \to W^{-s_2,p_2'}(X).
\end{equation}
\end{lem}

\begin{lem}[Sobolev multiplication theorem for non-negative real derivative exponents and two factors in continuous range]
\label{lem:Freed_Uhlenbeck_equation_6-34_nonnegative_real_algebra}
\cite[Corollary 9.7]{PalaisFoundationGlobal}
Let $X$ be a closed, Riemannian, smooth manifold of dimension $d \geq 2$. If $s$ and $p$ obey
\begin{equation}
\label{eq:Freed_Uhlenbeck_6-34_algebra_nonnegative_real_hypothesis}
s \in \RR \quad\hbox{and}\quad sp > d,
\end{equation}
then $W^{s,p}(X)$ is an algebra, that is, the following multiplication map is defined and continuous,
\begin{equation}
\label{eq:Freed_Uhlenbeck_6-34_nonnegative_real_algebra}
W^{s,p}(X) \times W^{s,p}(X) \to W^{s,p}(X).
\end{equation}
\end{lem}

\begin{lem}[Sobolev multiplication theorem for non-negative real derivative exponents and one factor in continuous range]
\label{lem:Freed_Uhlenbeck_equation_6-34_module_nonnegative_real}
\cite[Corollary 9.7]{PalaisFoundationGlobal}
Let $X$ be a closed, Riemannian, smooth manifold of dimension $d \geq 2$. If $s, s_1$ and $p$ obey
\begin{equation}
\label{eq:Freed_Uhlenbeck_6-34_module_nonnegative_real_hypothesis}
s_1,s \in \RR \hbox{ with } s_1\geq s \geq 0, \quad 1 \leq p < \infty, \quad s_1p > d,
\end{equation}
then $W^{s, p}(X)$ is a $W^{s_1, p}(X)$-module, that is, the following multiplication map is defined and continuous,
\begin{equation}
\label{eq:Freed_Uhlenbeck_6-34_module_nonnegative_real}
W^{s_1, p}(X) \times W^{s, p}(X) \to W^{s, p}(X).
\end{equation}
\end{lem}

More generally, suppose we are given a $C^\infty$ complex Hermitian (real Riemannian) vector bundle $E$ with $C^\infty$ Hermitian (Riemannian) connection $A$ and a covariant derivative, $\nabla_A: C^\infty(X; E)\to C^\infty(X; \Lambda^1\otimes E)$. For $1 < p < \infty$ and $s \in \RR$, we may define the Sobolev space $W^{s,p}_A(X; E)$
with the aid of the (augmented) connection Laplace operator \eqref{eq:Connection_Laplacian}, that is, $\nabla_A^*\nabla_A+1: C^\infty(X; E) \to C^\infty(X; E)$. When $p = 2$ and $s\in\RR$, we write $H_A^s(X; E) := W_A^{s,2}(X; E)$ for brevity.

All of the preceding embedding and pointwise multiplication results for real or complex valued functions on $X$ then extend, \emph{mutatis mutandis}, to the case of sections of smooth complex Hermitian (real Riemannian) vector bundles $E_i$, for $i=1,2,3$, and pointwise tensor-product multiplication. Indeed, this is the situation considered by Freed and Uhlenbeck in \cite[pp. 95--96]{FU} and Palais in \cite[Chapter 9]{PalaisFoundationGlobal}.

\subsection{$L^p$ theory for scalar elliptic partial differential and pseudo-differential operators on $\RR^d$ and applications to elliptic systems}
\label{subsec:Krylov_LecturesSobolev_13.4.5}
In this subsection, we derive an \apriori $L^p$ estimate and unique solvability result (Theorem \ref{thm:Krylov_LecturesSobolev_13.4.5_diagonal_principal_symbol}) for an elliptic partial differential system with a diagonal principal symbol from corresponding results for a scalar elliptic partial differential equation. Because of our focus in this monograph on elliptic partial differential operators on sections of vector bundles over closed manifolds, it will suffice here to consider elliptic partial differential systems on $\RR^d$.

While our main goal in this subsection is to prove Theorem \ref{thm:Krylov_LecturesSobolev_13.4.5_diagonal_principal_symbol}, it is convenient to introduce the relevant background by first recalling the corresponding results for scalar elliptic partial differential operators. Suppose first that $a = (a^{ij})$ is a constant, real, symmetric matrix obeying
\begin{equation}
\label{eq:Krylov_LecturesSobolev_page_93}
\kappa^{-1} \leq a^{ij}\xi_i\xi_j \leq \kappa, \quad\forall\, \xi \in \RR^d \hbox{ with } |\xi|=1,
\end{equation}
for some constant $\kappa > 0$, and consider the elliptic, homogeneous, scalar, second-order partial differential operator as in \cite[Equation (4.0.1)]{Krylov_LecturesSobolev}, albeit with the opposite sign convention,
\begin{equation}
\label{eq:Krylov_LecturesSobolev_4-0-1}
\sA := -a^{ij}\frac{\partial^2}{\partial x_i\partial x_j}.
\end{equation}
Given $p \in (1, \infty)$, and $f \in L^p(\RR^d)$, and $\lambda > 0$, unique solvability of the scalar equation,
\begin{equation}
\label{eq:Krylov_LecturesSobolev_4-3-6_elliptic}
(\sA + \lambda)u = f \quad\hbox{a.e. on } \RR^d \quad \hbox{for } u \in W^{2,p}(\RR^d),
\end{equation}
is provided by \cite[Theorem 4.3.8 (ii)]{Krylov_LecturesSobolev}. The \apriori estimate,
\begin{multline}
\label{eq:Krylov_LecturesSobolev_4-3-8}
\lambda\|u\|_{L^p(\RR^d)} + \lambda^{1/2}\|Du\|_{L^p(\RR^d)} + \|D^2u\|_{L^p(\RR^d)} \leq C\|(\sA + \lambda)u\|_{L^p(\RR^d)},
\\
\forall\, u \in W^{2,p}(\RR^d),
\end{multline}
with $C = C(d, p, \kappa)$ is also provided by \cite[Theorem 4.3.8 (ii)]{Krylov_LecturesSobolev} with the aid of the interpolation inequality \cite[Theorem 1.5.1]{Krylov_LecturesSobolev},
$$
\|Du\|_{L^p(\RR^d)} \leq \|u\|_{L^p(\RR^d)}^{1/2}\|D^2u\|_{L^p(\RR^d)}^{1/2}, \quad \forall\, u \in W^{2,p}(\RR^d).
$$
When $p=2$, unique solvability of \eqref{eq:Krylov_LecturesSobolev_4-3-6_elliptic} is obtained by the Fourier transform in the shape of \cite[Theorem 1.3.16]{Krylov_LecturesSobolev}. The proof for general case $p \in (1, \infty)$ follows from the case $p=2$, the \apriori estimate \eqref{eq:Krylov_LecturesSobolev_4-3-8}, and the fact that $(\Delta + \lambda)C^\infty_0(\RR^d)$ is dense in $L^p(\RR^d)$ by \cite[Theorem 1.1.6]{Krylov_LecturesSobolev}, for any $p \in [1, \infty)$. Here, we denote $\sA = \Delta$ when $a^{ij} = \delta^{ij}$ in \eqref{eq:Krylov_LecturesSobolev_4-0-1}, again opposite to the sign convention in \cite{Krylov_LecturesSobolev}.

\begin{defn}[Strongly elliptic scalar partial differential operator with complex constant coefficients]
\label{defn:Krylov_LecturesSobolev_12-2-1}
\cite[Definition 12.2.1]{Krylov_LecturesSobolev}
Let $m \geq 1$ and $a_\alpha \in \CC$, for multi-indices $\alpha \in \NN^d$ with $|\alpha| \leq m$, and denote\footnote{Krylov denotes $D = (\partial_{x_1},\ldots,\partial_{x_d})$, but we include the `$-i$' in $D$ and exclude the $i^{|\alpha|}$ factor in the definition of strong ellipticity for the sake of consistency with \cite{Denk_Hieber_Pruss_2003, Haller-Dintelmann_Heck_Hieber_2006, Heck_Hieber_Stavrakidis_2010}.}
$D = -i(\partial_{x_1},\ldots,\partial_{x_d})$, with $i = \sqrt{-1}$, and $D^\alpha = (-i\partial_{x_1})^{\alpha_1}\cdots(-i\partial_{x_d})^{\alpha_d}$. The operator,
$$
\sA := a(D) = \sum_{|\alpha|\leq m}a_\alpha D^\alpha,
$$
is called a \emph{scalar partial differential operator of order $m\geq 1$ with complex constant coefficients} and \emph{strongly elliptic} if
\begin{equation}
\label{eq:Krylov_LecturesSobolev_page_271}
\begin{aligned}
\mathring{a}(\xi) := \sum_{|\alpha| = m}a_\alpha \xi^\alpha &\neq 0, \quad\forall\, \xi \in \RR^d \less \{0\}, \quad\hbox{and}
\\
a(\xi) := \sum_{|\alpha|\leq m}a_\alpha \xi^\alpha &\neq 0, \quad\forall\, \xi \in \RR^d,
\end{aligned}
\end{equation}
where $\mathring{a}(\xi)$ and $a(\xi)$ are the \emph{principal symbol} and \emph{symbol} of $\sA$, respectively.
\end{defn}

For example, when $m=2$, the operator $\sA = \Delta + 1$ is strongly elliptic in the sense of Definition \ref{defn:Krylov_LecturesSobolev_12-2-1} since it has symbol $a(\xi) = |\xi|^2 + 1$.
It will be useful to recall the

\begin{defn}[Scalar elliptic symbol and scalar pseudo-differential operator]
\label{defn:Krylov_LecturesSobolev_12-2-1_and_7}
\cite[Definitions 7.8.1 and 18.1.1]{Hormander_v3},
\cite[Definitions 12.4.1 and 12.4.7]{Krylov_LecturesSobolev}
Given $\mu \in \RR$, a function $a \in C^\infty(\RR^d\times\RR^d; \CC)$ is a \emph{symbol of order $\mu$} if, for each pair of multi-indices $\alpha,\beta \in \NN^d$, there is a positive constant, $C_{\alpha,\beta}$, such that
\begin{equation}
\label{eq:Krylov_LecturesSobolev_12-4-1}
|\partial_\xi^\alpha \partial_x^\beta a(x,\xi)| \leq C_{\alpha,\beta}(1 + |\xi|^2)^{(\mu - |\alpha|)/2}, \quad\forall\, x, \xi \in \RR^d.
\end{equation}
If in addition there is a positive constant, $\kappa$, such that
\begin{equation}
\label{eq:Krylov_LecturesSobolev_12-4-2}
|a(x,\xi)| \geq \kappa(1 + |\xi|^2)^{\mu/2}, \quad\forall\, x, \xi \in \RR^d,
\end{equation}
then the symbol, $a$, is \emph{elliptic}, with \emph{constant of ellipticity} $\kappa$. The operator, $\sP:\sS \to \sS$, on the Schwartz space, $\sS$, of rapidly decaying functions \cite[Definition 12.2.1]{Krylov_LecturesSobolev} defined by the composition,
$$
\sP u := \sF^{-1}(a \sF(u)), \quad\forall\, u \in \sS,
$$
where
$$
\sF(u)(\xi) := \frac{1}{(2\pi)^{d/2}}\int_{\RR^d} e^{-ix\cdot\xi} u(x)\, dx
$$
denotes the Fourier transform, is called a \emph{pseudo-differential operator of order $\mu$ with symbol $a$.}
\end{defn}

The pseudo-differential operator,
$$
(\Delta + 1)^{\mu/2},
$$
obtained in Definition \ref{defn:Krylov_LecturesSobolev_12-2-1_and_7} with symbol
$$
a(x,\xi) := (|\xi|^2 + 1)^{\mu/2}, \quad\forall\, x, \xi \in \RR^d,
$$
is a strongly elliptic pseudo-differential operator of order $\mu \in \RR$ and agrees with the naive definition of $(\Delta + 1)^{\mu/2}$ when $\mu = 2n$ and $n\geq 0$ is an integer \cite[Example 12.4.10]{Krylov_LecturesSobolev}.

For $s\in \RR$ and $p \in (1,\infty)$, we recall that the Banach space of \emph{Bessel potentials} is given by \cite[Definition 13.3.1]{Krylov_LecturesSobolev}
\begin{equation}
\label{eq:Krylov_LecturesSobolev_definition_13-3-1}
H^{s,p}(\RR^d) := (\Delta + 1)^{-s/2}L^p(\RR^d),
\end{equation}
with norm
\begin{equation}
\label{eq:Krylov_LecturesSobolev_13-3-1}
\|u\|_{H^{s,p}(\RR^d)} := \|(\Delta + 1)^{s/2}u\|_{L^p(\RR^d)}.
\end{equation}
By \cite[Theorem 13.3.7 (ii)]{Krylov_LecturesSobolev}, we are assured that $C^\infty_0(\RR^d)$ is dense in $H^{s,p}(\RR^d)$. When $m \in \NN$, then \cite[Theorem 13.3.12]{Krylov_LecturesSobolev} implies that
$$
H^{m,p}(\RR^d) = W^{m,p}(\RR^d).
$$
When $\sP$ is a pseudo-differential operator of order $\mu \in \RR$ in the sense of Definition \ref{defn:Krylov_LecturesSobolev_12-2-1_and_7}, then \cite[Theorem 13.3.10]{Krylov_LecturesSobolev} asserts that $\sP$ defines a bounded map from $H^{s+\mu,p}(\RR^d)$ into $H^{s,p}(\RR^d)$,
\begin{equation}
\label{eq:Krylov_LecturesSobolev_13-3-3}
\|\sP u\|_{H^{s, p}(\RR^d)} \leq C\|u\|_{H^{s+\mu, p}(\RR^d)}, \quad\forall\, u \in H^{s+\mu,p}(\RR^d).
\end{equation}
If $\sP$ is strongly elliptic in the sense of \cite[Definitions 12.4.1 and 7]{Krylov_LecturesSobolev}, then \cite[Theorem 13.3.10]{Krylov_LecturesSobolev} also asserts that $\sP$ defines a bounded, \emph{one-to-one} map from $H^{s+\mu,p}(\RR^d)$ \emph{onto} $H^{s,p}(\RR^d)$, with \apriori estimate,
\begin{equation}
\label{eq:Krylov_LecturesSobolev_13-3-4}
\|u\|_{H^{s+\mu, p}(\RR^d)} \leq C\|\sP u\|_{H^{s, p}(\RR^d)}, \quad\forall\, u \in H^{s+\mu,p}(\RR^d).
\end{equation}
For $s \in \RR$, then \cite[Theorem 13.3.7 (iii)]{Krylov_LecturesSobolev} provides a result that is simpler than
\cite[Theorem 13.3.10]{Krylov_LecturesSobolev}, but whose systems analogue is adequate for our applications in this monograph: The operator $(\Delta + 1)^{\mu/2}$ is an isometric isomorphism from the Banach space $H^{s+\mu,p}(\RR^d)$ onto $H^{s,p}(\RR^d)$. In the case $\mu=1$ (the Cauchy operator), the latter result is given by \cite[Theorem 12.9.3]{Krylov_LecturesSobolev}.

\begin{defn}[Uniformly strongly elliptic scalar partial differential operator]
\label{defn:Krylov_LecturesSobolev_13_4_0_uniformly_strongly_elliptic_partialdo}
\cite[Section 13.4]{Krylov_LecturesSobolev}
If $m \geq 1$ is an integer and $a_\alpha:\RR^d \to \CC$ is a collection of measurable functions, for multi-indices $\alpha \in \NN^d$ with $|\alpha| \leq m$, and $\kappa$ is a positive constant such that the symbol,
$$
a(x, \xi) := \sum_{|\alpha|\leq m}a_\alpha(x) \xi^\alpha, \quad \forall\, \xi \in \RR^d \hbox{ and a.e. } x\in \RR^d,
$$
obeys
$$
\kappa^{-1}(1 + |\xi|^m) \geq |a(x,\xi)| \geq \kappa(1 + |\xi|^m), \quad \forall\, \xi \in \RR^d \hbox{ and a.e. } x\in \RR^d,
$$
then the scalar partial differential operator of order $m\geq 1$,
$$
\sA := a(x, D),
$$
is called \emph{uniformly strongly elliptic}.
\end{defn}

Given $\sA$ as in Definition \ref{defn:Krylov_LecturesSobolev_13_4_0_uniformly_strongly_elliptic_partialdo} and $\lambda \geq 0$, one defines (compare \cite[p. 328]{Krylov_LecturesSobolev})
\begin{equation}
\label{eq:Krylov_LecturesSobolev_page_328_L_lambda}
\sA_\lambda :=  \sum_{|\alpha|\leq m} \lambda^{1-|\alpha|/m} a_\alpha(x) D^\alpha,
\end{equation}
and given a collection of functions, $b^\alpha \in L^\infty(\RR^d)$ for $|\alpha| \leq m-1$ such that
\begin{equation}
\label{eq:Krylov_LecturesSobolev_page_329_b_alpha_bound}
|b^\alpha| \leq K \quad\hbox{a.e. on } \RR^d, \quad \forall\, \alpha \in \NN^d \hbox{ with } |\alpha| \leq m-1,
\end{equation}
one defines \cite[p. 329]{Krylov_LecturesSobolev}
$$
B := \sum_{|\alpha|\leq m-1} b^\alpha(x) D^\alpha.
$$
Suppose in addition that the coefficients, $a_\alpha$, satisfy the following continuity property: There is an increasing function, $\omega:[0,\infty) \to \RR$, such that $\omega(t) \to 0$ as $t\downarrow 0$ and for all $\alpha \in \NN^d$ with $|\alpha| \leq m$,
\begin{equation}
\label{eq:Krylov_LecturesSobolev_assumption_13-4-1}
|a_\alpha(x) - a_\alpha(y)| \leq \omega(|x-y|), \quad\forall\, x, y \in \RR^d.
\end{equation}
For $p \in (1,\infty)$ and $f \in L^p(\RR^d)$ and $\lambda \in [\lambda_0, \infty)$ for a sufficiently large positive constant $\lambda_0$, unique solvability of the equation,
\begin{equation}
\label{eq:Krylov_LecturesSobolev_page_331_L_lambda_plus_B_equation}
(\sA_\lambda + B)u = f \quad\hbox{a.e. on } \RR^d \quad\hbox{for } u \in W^{m,p}(\RR^d),
\end{equation}
is provided by \cite[Theorem 13.4.5]{Krylov_LecturesSobolev}, with \apriori estimate given by \cite[Lemma 13.4.4]{Krylov_LecturesSobolev},
\begin{multline}
\label{eq:Krylov_LecturesSobolev_13-4-2}
\sum_{|\alpha|\leq m} \lambda^{1 - |\alpha|/m}\|D^\alpha u\|_{L^p(\RR^d)} \leq C\|(\sA_\lambda + B)u\|_{L^p(\RR^d)},
\\
\quad\forall\, \lambda \in [\lambda_0,\infty) \hbox{ and } u \in W^{m,p}(\RR^d),
\end{multline}
for $C = C(d, K, m, p, \kappa, \omega)$.

General results for suitably-defined elliptic partial differential systems can be found, for example, in work of Agmon, Douglis, and Nirenberg \cite{AgmonDouglisNirenberg2}.
However, when combined later with the method of continuity \cite[Theorem 5.2]{GilbargTrudinger} and \apriori $L^p$ estimates in \cite{Denk_Hieber_Pruss_2003, Haller-Dintelmann_Heck_Hieber_2006} for elliptic partial differential operators with non-diagonal, matrix-valued principal symbols, the following consequence of \cite[Theorem 13.4.5]{Krylov_LecturesSobolev} will be very useful, notwithstanding the simplicity of its proof.

\begin{thm}[Unique solvability and  \apriori $L^p$ estimate for a system of uniformly strongly elliptic partial differential equations on $\RR^d$]
\label{thm:Krylov_LecturesSobolev_13.4.5_diagonal_principal_symbol}
Let $N\geq 1$ and for $n \in \{1,\ldots,N\}$, let $a_\alpha^n \in C_b(\RR^d;\CC)$, for all $\alpha \in \NN^d$ with $|\alpha|=m$, and assume that the homogeneous partial differential operators, $\sum_{|\alpha| = m} a_\alpha^n(x)D^\alpha$, of order $m$ are uniformly strongly elliptic in the sense of Definition \ref{defn:Krylov_LecturesSobolev_13_4_0_uniformly_strongly_elliptic_partialdo}, for $1\leq n \leq N$. Let $b^\alpha \in L^\infty(\RR^d;\CC^{N\times N})$ obey \eqref{eq:Krylov_LecturesSobolev_page_329_b_alpha_bound}, for all $\alpha \in \NN^d$ with $|\alpha|\leq m-1$. Suppose $p \in (1,\infty)$. Then there are positive constants, $C$ and $\lambda_0$, depending at most on $d, K, N, p, \kappa, \omega$, with the following significance. Let $a_\alpha := \diag(a_\alpha^1, \ldots, a_\alpha^N) \in C_b(\RR^d;\CC^{N\times N})$ and
\begin{equation}
\label{eq:Krylov_LecturesSobolev_13_4_0_diagonal_principal_matrix_lower_order_coefficients}
\sA := \sum_{|\alpha| = m} a_\alpha(x)D^\alpha + \sum_{|\alpha| \leq m-1} b^\alpha(x)D^\alpha.
\end{equation}
Then the following inequality holds,
\begin{multline}
\label{eq:Krylov_LecturesSobolev_13-4-2_diagonal_principal_matrix_lower_order_coefficients}
\sum_{|\alpha|\leq m} \lambda^{1 - |\alpha|/m}\|D^\alpha u\|_{L^p(\RR^d)} \leq C\|(\sA + \lambda)u\|_{L^p(\RR^d)},
\\
\quad\forall\, \lambda \in [\lambda_0, \infty) \hbox{ and } u \in W^{m,p}(\RR^d).
\end{multline}
Moreover, for each $f \in L^p(\RR^d; \CC^N)$ and $\lambda \in [\lambda_0, \infty)$, there exists a unique solution $u \in W^{m,p}(\RR^d; \CC^N)$ to
\begin{equation}
\label{eq:Krylov_LecturesSobolev_page_331_diagonal_principal_matrix_lower_order_coefficients_equation}
(\sA + \lambda)u = f \quad\hbox{a.e. on } \RR^d.
\end{equation}
\end{thm}

\begin{proof}
Writing $\sA^n := \sum_{|\alpha| = m} a_\alpha^n(x)D^\alpha$ for $1\leq n\leq N$ and $u = (u_1,\ldots, u_n)$, the inequality \eqref{eq:Krylov_LecturesSobolev_13-4-2} yields
$$
\sum_{|\alpha|\leq m} \lambda^{1 - |\alpha|/m}\|D^\alpha u_n\|_{L^p(\RR^d)} \leq C\|(\sA^n + \lambda)u_n\|_{L^p(\RR^d)}, \quad 1\leq n\leq N,
\quad\forall\, u \in W^{m,p}(\RR^d),
$$
for all $\lambda \in [\lambda_0, \infty)$. For sufficiently large $\lambda_0$, one then obtains \eqref{eq:Krylov_LecturesSobolev_13-4-2_diagonal_principal_matrix_lower_order_coefficients} (see, for example, the proof of Theorem \ref{thm:Haller-Dintelmann_Heck_Hieber_5-1_full_apriori_estimate} below).

When $b^\alpha \equiv 0$ on $\RR^d$, for all $\alpha \in \NN^d$ with $|\alpha| \leq m-1$, then existence and uniqueness of a solution $u \in W^{m,p}(\RR^d; \CC^N)$ to \eqref{eq:Krylov_LecturesSobolev_page_331_diagonal_principal_matrix_lower_order_coefficients_equation} follows immediately from \cite[Theorem 13.4.5]{Krylov_LecturesSobolev}. For the general case, where the lower-order coefficients $b^\alpha$ may be non-zero, the conclusion follows by the method of continuity exactly as in the proof of \cite[Theorem 13.4.5]{Krylov_LecturesSobolev}.
\end{proof}

For a homogeneous elliptic partial differential operator with a non-diagonal principal symbol which is complex-matrix or even Banach-space operator valued, one may use a Calder\'on-Zygmund theory for vector-valued functions to derive \apriori estimates such as \eqref{eq:Krylov_LecturesSobolev_13-4-2_diagonal_principal_matrix_lower_order_coefficients} --- see, for example, \cite{Denk_Hieber_Pruss_2003, Haller-Dintelmann_Heck_Hieber_2006, Taylor_PDO_and_NLPDE} and the references contained therein --- by analogy with the development for scalar-valued functions as in \cite[Chapter 9]{GilbargTrudinger} or \cite{Krylov_LecturesSobolev}.

\subsection{Analytic semigroups on $L^p(\Omega;\RR)$, $C_0(\bar\Omega;\RR)$, and $L^1(\Omega;\RR)$ defined by scalar, elliptic partial differential operators on bounded open subsets of Euclidean space with a homogeneous Dirichlet boundary condition}
\label{subsec:Pazy_7-3}
In preparation for our subsequent development of existence and uniqueness results, \apriori estimates, resolvent estimates, and analytic semigroup generation results for elliptic systems on open subsets $\Omega \subseteqq \RR^d$, with homogeneous Dirichlet boundary conditions, we first provide an overview of the simpler and well-documented theory of analytic semigroups on $L^p(\Omega;\RR)$, $C_0(\bar\Omega;\RR)$, and $L^1(\Omega;\RR)$ defined by scalar, second-order, elliptic partial differential operators with $C^\infty$ coefficients and Dirichlet boundary conditions on bounded open subsets $\Omega \Subset \RR^d$. We largely follow the outline due to Pazy in \cite[Section 7.3]{Pazy_1983}; the far easier case where $p=2$, is described, for example, by Pazy in \cite[Section 7.2]{Pazy_1983}. In later subsections, we give extensions of these analytic semigroup generation results to the case of elliptic partial differential systems of order $m \geq 1$ on open subsets $\Omega \subseteqq \RR^d$, with homogeneous Dirichlet boundary conditions, together with the supporting existence and uniqueness results, \apriori estimates, resolvent estimates, and their detailed proofs in many cases.

We suppose first that $1 < p < \infty$ and $\Omega \subset \RR^d$ is a bounded domain with $C^\infty$ boundary $\partial\Omega$ and that \cite[Equations (7.2.1) or (7.3.1)]{Pazy_1983}
\begin{equation}
\label{eq:Pazy_7-3-1}
\sA(x, D)u := \sum_{|\alpha| \leq 2m}a_\alpha(x) D^\alpha u, \quad x \in \Omega,
\end{equation}
is a partial differential operator of order $2m$, with $m \geq 1$, where the coefficients $a_\alpha$ belong to $C^\infty(\bar\Omega;\CC)$, and \emph{strongly elliptic} \cite[Definition 7.2.1]{Pazy_1983},
that is, that there is a positive constant $\kappa$ such that
\begin{equation}
\label{eq:Pazy_7-2-2}
\sum_{|\alpha| \leq 2m}a_\alpha(x) \xi^\alpha \geq \kappa|\xi|^{2m}, \quad\forall\, \xi \in \RR^d, \quad x \in \Omega,
\end{equation}
recalling
that we define $D = -i\partial$ and $\partial = (\partial_{x_1},\ldots,\partial_{x_d})$. The following fundamental \apriori estimate is well-known (see, for example, \cite[Theorem 9.13]{GilbargTrudinger} when $m=2$ and $\sA$ has real coefficients).

\begin{thm}[\Apriori estimate]
\label{thm:Pazy_7-3-1}
\cite{AgmonDouglisNirenberg1},
\cite[Theorem 7.3.1]{Pazy_1983}
Let $\sA(x,D)$ as in \eqref{eq:Pazy_7-3-1} be strongly elliptic partial differential operator of order $2m$, for $m \geq 1$, in the sense of \eqref{eq:Pazy_7-2-2} on a bounded open subset $\Omega \subset \RR^d$ with $C^\infty$ boundary $\partial\Omega$ and $1<p<\infty$. Then there is a positive constant, $C$, such that
\begin{equation}
\label{eq:Pazy_7-3-3}
\|u\|_{W^{2m,p}(\Omega)} \leq C\left(\|\sA u\|_{L^p(\Omega)} + \|u\|_{L^p(\Omega)}\right),
\quad\forall\, u \in W^{2m,p}(\Omega;\CC)\cap W^{m,p}_0(\Omega;\CC).
\end{equation}
\end{thm}

Using this \apriori estimate together with an argument of Agmon one proves the following theorem; see Agmon \cite[Equation (2.3)]{Agmon_1962} and his proof of \cite[Theorem 2.1]{Agmon_1962}.

\begin{thm}[Agmon estimate]
\label{thm:Pazy_7-3-2}
\cite[Theorem 2.1 and Equation (2.3)]{Agmon_1962}, \cite{AgmonLecturesEllipticBVP},
\cite[Theorem 7.3.2]{Pazy_1983}
Assume the hypotheses of Theorem \ref{thm:Pazy_7-3-1}. Then there are constants, $C > 0$ and $\lambda_0 \geq 0$ and $\vartheta \in (\pi/2, \pi)$, such that
\begin{multline}
\label{eq:Pazy_7-3-4}
|\lambda|\|u\|_{L^p(\Omega)} \leq C\|(\sA + \lambda)u\|_{L^p(\Omega)},
\\
\forall\, \lambda \in \Delta_\vartheta \hbox{ such that } |\lambda| \geq \lambda_0, \quad
u \in W^{2m,p}(\Omega;\CC)\cap W^{m,p}_0(\Omega;\CC).
\end{multline}
\end{thm}

Given a strongly elliptic operator $\sA(x, D)$ of order $2m$, we associate a linear, unbounded operator $\sA_p$ on $L^p(\Omega;\CC)$ as follows:

\begin{defn}[Domain of a partial differential operator on $L^p(\Omega;\CC)$]
\label{defn:Pazy_7-3-3}
\cite[Definition 7.3.3]{Pazy_1983}
Assume the setup of Theorem \ref{thm:Pazy_7-3-1}. Let
\begin{align}
\label{eq:Pazy_7-3-5}
\sD(\sA_p) &:= W^{2m,p}(\Omega;\CC)\cap W^{m,p}_0(\Omega;\CC),
\\
\label{eq:Pazy_7-3-6}
\sA_p u &:= \sA(x, D)u, \quad \forall\, u \in \sD(\sA_p).
\end{align}
\end{defn}

The domain $\sD(\sA_p)$ contains $C^\infty_0(\Omega;\CC)$ and is therefore dense in $L^p(\Omega;\CC)$. Moreover, from Theorem \ref{thm:Pazy_7-3-1} it follows readily (see Remark \ref{rmk:Domain_elliptic_partial_differential_operator_on_L1}) that $\sA_p$ is a closed operator on $L^p(\Omega;\CC)$. From Theorems \ref{thm:Pazy_7-3-1} and \ref{thm:Pazy_7-3-2} one deduces the

\begin{thm}[Generator of an analytic semigroup on $L^p(\Omega;\CC)$]
\label{thm:Pazy_7-3-5}
\cite[Theorem 7.3.5]{Pazy_1983},
Assume the hypotheses of Theorem \ref{thm:Pazy_7-3-1}. If $\sA_p$ is the operator associated with $\sA(x,D)$ by Definition \ref{defn:Pazy_7-3-3}, then $-\sA_p$ is the infinitesimal generator of an analytic semigroup on $L^p(\Omega;\CC)$.
\end{thm}

See \cite[Theorem 7.3.5]{Pazy_1983} or \cite[Theorem 38.2]{Sell_You_2002} for proofs of Theorems \ref{thm:Pazy_7-3-2} and \ref{thm:Pazy_7-3-5} when $m=2$ and $\sA$ can be expressed in divergence form without lower-order coefficients. We turn now to the cases $p = 1$ and $p = \infty$ and start with $p = \infty$. We have the following analogue of Definition \ref{defn:Pazy_7-3-3}.

\begin{defn}[Domain of a partial differential operator on $L^\infty(\Omega;\CC)$]
\label{defn:Pazy_7-3-3_L_infinity}
\cite[Equations (7.3.20) and (7.3.21)]{Pazy_1983}
Assume the setup of Theorem \ref{thm:Pazy_7-3-1}. Let
\begin{gather}
\label{eq:Pazy_7-3-20}
\begin{split}
\sD(\sA_\infty) := \{u \in L^\infty(\Omega;\CC): u \in W^{2m,p}(\Omega;\CC), \forall\, p > d,
\\
\sA(x, D)u \in L^\infty(\Omega;\CC), u = 0 \hbox{ on } \partial\Omega\},
\end{split}
\\
\label{eq:Pazy_7-3-21}
\sA_\infty u = \sA(x, D)u, \quad\forall\, u \in \sD(\sA_\infty).
\end{gather}
\end{defn}

By the Sobolev Embedding Theorem \cite[Theorem 4.12]{AdamsFournier} it follows that $\sD(\sA_\infty) \subset C(\bar\Omega;\CC)$ and therefore the condition $u = 0$ on $\partial\Omega$ makes sense. Moreover, from the definition \eqref{eq:Pazy_7-3-20} of $\sD(A_\infty)$, it follows that $\sD(\sA_\infty) \subset W^{2m,p}(\Omega;\RR) \cap W^{m,p}_0(\Omega;\RR) = \sD(\sA_p)$ for every $p > d$. The Stewart-Masuda argument \cite{Stewart_1974} leads to the Agmon estimate, for constants $C > 0$ and $\lambda_0 \geq 0$ and $\vartheta \in (\pi/2, \pi)$, such that
\begin{equation}
\label{eq:Pazy_7-3-22}
|\lambda|\|u\|_{L^\infty(\Omega)} \leq C\|(\sA + \lambda)u\|_{L^\infty(\Omega)},
\quad\forall\, \lambda \in \Delta_\vartheta \hbox{ such that } |\lambda| \geq \lambda_0, \quad
u \in \sD(\sA_\infty).
\end{equation}
The estimate \eqref{eq:Pazy_7-3-22} follows immediately from Theorem \ref{thm:Stewart_1974_1_and_2} below,
as do the facts that $\sA + \lambda$ is injective and has closed range for all $\lambda \in \Delta_\vartheta$ such that $|\lambda| \geq \lambda_0$.

However, $\sA_\infty$ cannot be the infinitesimal generator of even a $C^0$ semigroup on $L^\infty(\Omega;\CC)$ because $\sD(\sA_\infty)$ is never dense in $L^\infty(\Omega;\CC)$. Indeed, we have noted earlier that $\sD(\sA_\infty) \subset C(\bar\Omega;\CC)$ and therefore the closure of $\sD(\sA_\infty)$ in $L^\infty(\Omega;\CC)$ is also contained in $C(\bar\Omega;\CC)$. But $C(\bar\Omega;\CC)$ is not dense in $L^\infty(\Omega;\CC)$ and thus $\sD(\sA_\infty)$ cannot be dense in $L^\infty(\Omega;\CC)$. To overcome this difficulty one restricts to spaces of continuous functions on $\Omega$, as in the following analogue of Definition \ref{defn:Pazy_7-3-3}.

\begin{defn}[Domain of a partial differential operator on $C_0(\bar\Omega;\CC)$]
\label{defn:Pazy_7-3-3_C}
\cite[Equations (7.3.23) and (7.3.24)]{Pazy_1983}
Assume the setup of Theorem \ref{thm:Pazy_7-3-1}. Let
\begin{align}
\label{eq:Pazy_7-3-23}
\sD(\sA_c) &:= \{u \in \sD(\sA_\infty): \sA(x, D)u \in C_0(\bar\Omega;\CC)\},
\\
\label{eq:Pazy_7-3-24}
\sA_c u &= \sA(x, D)u, \quad u \in \sD(\sA_c),
\end{align}
where $C_0(\bar\Omega;\CC) := \{u \in C(\bar\Omega;\CC): u = 0 \hbox{ on } \partial\Omega\}$.
\end{defn}

One now has the

\begin{thm}[Generator of an analytic semigroup on $C_0(\bar\Omega;\CC)$]
\label{thm:Pazy_7-3-7}
\cite[Theorem 7.3.7]{Pazy_1983}
Assume the hypotheses of Theorem \ref{thm:Pazy_7-3-1}. If $\sA_c$ is the operator associated with $\sA(x,D)$ by Definition \ref{defn:Pazy_7-3-3_C}, then $-\sA_c$ is the infinitesimal generator of an analytic semigroup on $C_0(\bar\Omega;\CC)$.
\end{thm}

The proof of Theorem \ref{thm:Pazy_7-3-7} is based on the \apriori $L^\infty$ estimate \eqref{eq:Stewart_1974_1-1} in Theorem \ref{thm:Stewart_1974_1_and_2} below.
The \apriori estimate \eqref{eq:Pazy_7-3-3} in Theorem \ref{thm:Pazy_7-3-1} is not known to hold when $p = 1$, but an analytic semigroup generation result for this case can nonetheless be derived by exploiting a duality between continuous functions and $L^1$ functions on $\Omega$. We start with the\footnote{One can also deduce this result from \cite[Lemma 2.7]{AdamsFournier}, noting that the proof of \cite[Lemma 2.7]{AdamsFournier} is valid when $p=1$.}

\begin{lem}[Reverse H\"older inequality for $p=1$]
\label{lem:Pazy_7-3-8}
\cite[Exercise 4.26 (1)]{Brezis}, \cite[Lemma 7.3.8]{Pazy_1983}
Let $\Omega \subseteqq \RR^d$ be an open
subset\footnote{We note that from the proof of \cite[Lemma 7.3.8]{Pazy_1983} or statement of \cite[Exercise 4.26]{Brezis} that $\Omega \subset \RR^d$ may be any open subset, not necessarily a bounded domain as assumed in \cite[Lemma 7.3.8]{Pazy_1983}.}.
If $u \in L^1(\Omega;\CC)$, then
\begin{equation}
\label{eq:Pazy_7-3-26}
\|u\|_{L^1(\Omega)}
=
\sup_{\begin{subarray}{c} \varphi \in C^\infty_0(\Omega; \CC) \\ \|\varphi\|_{L^\infty(\Omega)} \leq 1\end{subarray}}
(u, \varphi)_{L^2(\Omega)}.
\end{equation}
\end{lem}

From \cite[Exercise 4.26 (1)]{Brezis}, one may replace $C^\infty_0(\Omega)$ by $C_0(\Omega)$ in \eqref{eq:Pazy_7-3-26}.

Recall that $(L^1(\Omega;\CC))' = L^\infty(\Omega;\CC)$ by \cite[Theorem 4.14]{Brezis}, but $(L^\infty(\Omega;\CC))' \neq L^1(\Omega;\CC)$. However, if $\sM(\bar\Omega)$ denotes the space of Radon measures on a bounded open subset $\Omega \subset \RR^d$ \cite[pp. 114--115]{Brezis}, then $\sM(\bar\Omega) = (C(\bar\Omega;\CC))'$ by \cite[Theorem 4.31]{Brezis} while $\sM(\Omega) = (C_0(\bar\Omega;\CC))'$.

We turn now to the definition of the operator $\sA_1$, associated with the strongly elliptic operator $\sA(x, D)$ given by \eqref{eq:Pazy_7-3-1}, on the space $L^1(\Omega;\CC)$. We replace\footnote{Pazy defines $\sD(\sA_1)$ to be the set of all $u \in W^{2m-1, 1}(\Omega; \CC) \cap W^{m, 1}_0(\Omega; \CC)$ such that $\sA(x, D) \in L^1(\Omega;\CC)$, but this does not appear justifiable in view of the precise characterizations of $\sD(\sA_1)$ given, for example, by Guidetti in \cite[Theorem 3.3]{Guidetti_1993} and Tanabe in \cite[Theorem 5.8]{Tanabe_1997}.}
\cite[Definition 3.9]{Pazy_1983} (see the discussion in Remark \ref{rmk:Domain_elliptic_partial_differential_operator_on_L1} below) by our

\begin{defn}[Domain of a partial differential operator on $L^1(\Omega;\CC)$]
\label{defn:Pazy_7-3-9_modified}
Assume the setup of Theorem \ref{thm:Pazy_7-3-1}. Let
\begin{equation}
\label{eq:Pazy_7-3-27}
\sD(\sA_1) \subset L^1(\Omega;\CC),
\end{equation}
be the domain of the smallest closed extension of $\sA: C^\infty_0(\Omega;\CC) \subset L^1(\Omega;\CC) \to L^1(\Omega;\CC)$. The operator $\sA_1$ is defined by
\begin{equation}
\label{eq:Pazy_7-3-28}
\sA_1 = \sA(x, D)u, \quad\forall\, u \in \sD(\sA_1).
\end{equation}
\end{defn}

One then has the

\begin{thm}[Generator of an analytic semigroup on $L^1(\Omega;\CC)$]
\label{thm:Pazy_7-3-10}
\cite[Theorem 7.3.10]{Pazy_1983}
Assume the hypotheses of Theorem \ref{thm:Pazy_7-3-1}. If $\sA_1$ is the operator associated with $\sA(x,D)$ by Definition \ref{defn:Pazy_7-3-9_modified}, then $-\sA_1$ is the infinitesimal generator of an analytic semigroup on $L^1(\Omega;\CC)$.
\end{thm}

The proof of Theorem \ref{thm:Pazy_7-3-10} employs Lemma \ref{lem:Pazy_7-3-8} and Theorem \ref{thm:Stewart_1974_1_and_2}.

\begin{rmk}[Alternative proofs of Theorem \ref{thm:Pazy_7-3-10} and further related results]
\label{rmk:Pazy_7-3-10}
See Amann \cite[p. 226]{Amann_1983}, Angiuli, Pallara, and Paronetto \cite[Proposition 2.5]{Angiuli_Pallara_Paronetto_2010}, Cannarsa and Vespri \cite{Cannarsa_Vespri_1988} (for $\Omega = \RR^d$), Di Blasio \cite{DiBlasio_1991}, and Lunardi and Metafune \cite{Lunardi_Metafune_2004} for proofs of Theorem \ref{thm:Pazy_7-3-10} and extensions as well as additional related results, including consideration of boundary conditions other than the homogeneous Dirichlet boundary condition assumed in Theorem \ref{thm:Pazy_7-3-10}.
\end{rmk}

As we noted earlier, the proofs of Theorems \ref{thm:Pazy_7-3-7} and \ref{thm:Pazy_7-3-10} are based on fundamental results of Stewart \cite{Stewart_1974, Stewart_1980}. Because of the important role played by the Stewart-Masuda method in our article and for the sake of comparison with later results, we now review the main result of \cite{Stewart_1974}. For a possibly unbounded open subset $\Omega \subseteqq \RR^d$, one denotes \cite[p. 144]{Stewart_1974}
\begin{equation}
\label{eq:Stewart_page_144_C_sub_0_bar_omega}
C_0(\bar\Omega) := \{u \in C(\bar\Omega): u = 0 \hbox{ on } \partial\Omega \hbox{ and } |u(x)| \to 0 \hbox{ as } |x| \to \infty\}.
\end{equation}
Following Jacob \cite[p. xiv]{Jacob_v1}, we also find it convenient to denote
\begin{equation}
\label{eq:Jacob_page_xiv_C_sub_infinity_omega}
C_\infty(\Omega) := \{u \in C(\Omega): |u(x)| \to 0 \hbox{ as } |x| \to \infty\},
\end{equation}
so that $C_0(\bar\Omega) = C_\infty(\RR^d)$ when $\Omega  = \RR^d$.

\begin{thm}[\Apriori estimate for a uniformly strongly elliptic partial differential operator]
\label{thm:Stewart_1974_1_and_2}
\cite[Theorems 1 and 2]{Stewart_1974}, \cite[Theorem 1]{Stewart_1980}
Let $d \geq 2$ and $m\geq 1$ be integers, $K$, $L$, $q>d$, and $\kappa$ be positive constants, and $\omega$ be a modulus of continuity. Let $\Omega \subseteqq \RR^d$ be an open subset with boundary, $\partial\Omega$, which is uniformly regular of class $C^{2m}$ in the sense of \cite[Section 4.10]{AdamsFournier}, with overlap constant $K$ and coordinate chart and chart inverse derivative bound $L$. Then there are positive constants, $C$ and $r_0$ and $\vartheta \in (\pi/2, \pi)$ and $\lambda_0$, with the following significance. Let $\sA$ be a scalar, uniformly strongly elliptic partial differential operator of order $2m$ in the sense of Definition \ref{defn:Krylov_LecturesSobolev_13_4_0_uniformly_strongly_elliptic_partialdo} with $a_\alpha \in C(\bar\Omega)$ having modulus of continuity $\omega$ when $|\alpha|=2m$, and $a_\alpha \in L^\infty(\Omega)$ when $|\alpha| < 2m$. If
\begin{multline*}
u \in \{v \in W^{2m,q}_{\loc}(\bar\Omega) \cap C_0(\bar\Omega): \sA v \in C_0(\bar\Omega),
\\
D^\alpha v = 0 \hbox{ on } \partial\Omega, \hbox{for all } \alpha \in \NN^d\hbox{ with }|\alpha| < m\}
\end{multline*}
and $\lambda \in \bar\Delta_\vartheta$ with $|\lambda| \geq \lambda_0$, then
\begin{multline}
\label{eq:Stewart_1974_1-2}
\sum_{|\alpha| < 2m} |\lambda|^{1-|\alpha|/2m}\|D^\alpha u\|_{C(\bar\Omega)}
+ \sum_{|\alpha| = 2m} |\lambda|^{d/2mq} \sup_{x_0 \in \Omega} \|D^\alpha u\|_{L^q(\Omega\cap B(x_0, r_0|\lambda|^{-1/2m}))}
\\
\leq C \sup_{x_0 \in \RR^d} \|(\sA + \lambda)u\|_{L^q(\Omega\cap B(x_0, r_0|\lambda|^{-1/2m}))},
\end{multline}
and
\begin{multline}
\label{eq:Stewart_1974_1-1}
\sum_{|\alpha| < 2m} |\lambda|^{1-|\alpha|/2m}\|D^\alpha u\|_{C(\bar\Omega)}
+ \sum_{|\alpha| = 2m} |\lambda|^{d/2mq} \sup_{x_0 \in \RR^d} \|D^\alpha u\|_{L^q(\Omega\cap B(x_0, r_0|\lambda|^{-1/2m}))}
\\
\leq C \|(\sA + \lambda)u\|_{C(\bar\Omega)}.
\end{multline}
\end{thm}

While the \apriori estimate \eqref{eq:Stewart_1974_1-1} follows from the \apriori estimate \eqref{eq:Stewart_1974_1-2}, we include both estimates for the sake of comparison with Theorem \ref{thm:Heck_Hieber_Stavrakidis_lemma_3-2}.


\subsection{Existence, uniqueness, \apriori estimates, and analytic semigroups on $L^p(\Omega;\CC^N)$, and $C_0(\bar\Omega;\CC^N)$, and $L^1(\Omega;\CC^N)$ defined by second-order, elliptic partial differential systems and a homogeneous Dirichlet boundary condition}
\label{subsec:Cannarsa_Terreni_Vespri_1985}
In this subsection, we review the existence and uniqueness results and \apriori estimates for a second-order, elliptic partial differential system on a bounded open subset of $\RR^d$ and a homogeneous Dirichlet boundary condition, due to Cannarsa, Terreni, and Vespri \cite{Cannarsa_Terreni_Vespri_1985}, together with their consequences for resolvent estimates and analytic semigroup generation results on $L^p(\Omega;\CC^N)$, and $C_0(\bar\Omega;\CC^N)$, and $L^1(\Omega;\CC^N)$.

Let $\Omega \subset \RR^d$ be bounded open subset, with boundary of class $C^2$, and consider the second-order differential
operator\footnote{Remember that we denote $\partial := (\partial_{x_1},\ldots,\partial_{x_d})$ and $D := -i\partial$, so $\partial^\alpha := \partial_{x_1}^{\alpha_1}\cdots\partial_{x_d}^{\alpha_d}$ and $D := (-i\partial_{x_1})^{\alpha_1}\cdots (-i\partial_{x_d})^{\alpha_d}$, for all $\alpha = (\alpha_1,\ldots,\alpha_d) \in \NN^d$.},
\begin{equation}
\label{eq:Cannarsa_Terreni_Vespri_1-1}
\sA(x, D) u := \sum_{|\alpha| \leq 2} a_\alpha(x)D^\alpha u, \quad\forall\, u \in C^\infty(\Omega;\CC^N).
\end{equation}
Here, the coefficients $a_\alpha$, for $\alpha \in \NN^d$ with $|\alpha| \leq 2$, are complex matrix-valued functions such that
\begin{align}
\label{eq:Cannarsa_Terreni_Vespri_2-5}
a_\alpha &\in C(\bar\Omega; \CC^{N\times N}), \quad\hbox{for } |\alpha| = 2,
\\
\label{eq:Cannarsa_Terreni_Vespri_2-6}
a_\alpha &\in L^\infty(\Omega; \CC^{N\times N}), \quad\hbox{for } |\alpha| \leq 1.
\end{align}
We then make the

\begin{defn}[Second-order, strictly elliptic partial differential operator]
\label{defn:Cannarsa_Terreni_Vespri_2-7}
One says that $\sA$ as in \eqref{eq:Cannarsa_Terreni_Vespri_1-1} is \emph{(strictly) elliptic} if there is a positive constant, $\kappa$, such that
\begin{equation}
\label{eq:Cannarsa_Terreni_Vespri_2-7}
\sum_{|\alpha| = 2} \Real \langle a_\alpha\eta, \eta \rangle \xi^\alpha > \kappa|\xi|^2|\eta|^2,
\quad\forall\, x \in \Omega, \ \xi \in \RR^d, \ \eta \in \CC^N.
\end{equation}
\end{defn}

The main results we require due to Cannarsa, Terreni, and Vespri \cite{Cannarsa_Terreni_Vespri_1985} are as follows.

\begin{thm}[Existence, uniqueness, and \apriori estimate for a solution to a second-order elliptic system on $W^{2,p}(\Omega; \CC^N)$ when $2\leq p < \infty$]
\label{thm:Cannarsa_Terreni_Vespri_6-7}
\cite[Theorem 6.7]{Cannarsa_Terreni_Vespri_1985}
Let $\Omega \subset \RR^d$ be bounded domain, with boundary of class $C^2$, and $\sA$ as in \eqref{eq:Cannarsa_Terreni_Vespri_1-1} with coefficients obeying \eqref{eq:Cannarsa_Terreni_Vespri_2-5}, \eqref{eq:Cannarsa_Terreni_Vespri_2-6}, and \eqref{eq:Cannarsa_Terreni_Vespri_2-7}. Then there are a constant $\omega_1 \geq 0$ and, given\footnote{As Piermarco Cannarsa points out \cite{Cannarsa_2014_private}, there a misprint in the statement of \cite[Theorem 6.7]{Cannarsa_Terreni_Vespri_1985}: the range allowed for $p$ should read $2 \leq p < + \infty$ and does \emph{not} include the case $p = + \infty$, as is also clear from the proof in \cite[pp. 91--93]{Cannarsa_Terreni_Vespri_1985}.}
$2\leq p < \infty$, a constant $C > 0$ with the following significance. For all $\lambda \in \CC$ obeying $\Real \lambda > \omega_1$ and $f \in L^p(\Omega;\CC^N)$, there is a unique solution $u \in W^{2,p}(\Omega;\CC^N)\cap W^{1,p}_0(\Omega;\CC^N)$ to
\begin{equation}
\label{eq:Cannarsa_Terreni_Vespri_3-1}
(\sA + \lambda)u = f \quad\hbox{a.e. on } \Omega,
\end{equation}
and, moreover, $u$ obeys
\begin{equation}
\label{eq:Cannarsa_Terreni_Vespri_6-15}
(|\lambda|-\omega_1)\|u\|_{L^p(\Omega)} + (|\lambda|-\omega_1)^{1/2}\|u\|_{W^{1,p}(\Omega)} + \|u\|_{W^{2,p}(\Omega)}
\leq C\|f\|_{L^p(\Omega)}.
\end{equation}
\end{thm}

For the case $p = \infty$, Cannarsa, Terreni, and Vespri provide the

\begin{thm}[Existence, uniqueness, and \apriori estimate for a solution to a second-order elliptic system on $W^{1,\infty}(\Omega; \CC^N)$]
\label{thm:Cannarsa_Terreni_Vespri_6-1}
\cite[Theorem 6.1]{Cannarsa_Terreni_Vespri_1985}
Let $\Omega \subset \RR^d$ be bounded domain, with boundary of class $C^2$, and $\sA$ as in \eqref{eq:Cannarsa_Terreni_Vespri_1-1} with coefficients obeying \eqref{eq:Cannarsa_Terreni_Vespri_2-5}, \eqref{eq:Cannarsa_Terreni_Vespri_2-6}, and \eqref{eq:Cannarsa_Terreni_Vespri_2-7}. Then there are constants, $C > 0$ and $\omega_1 \geq 0$, with the following significance. For all $\lambda \in \CC$ obeying $\Real \lambda > \omega_1$ and $f \in L^\infty(\Omega;\CC^N)$, there is a unique solution $u \in W^{2,2}(\Omega;\CC^N)\cap W^{1,\infty}(\Omega;\CC^N)\cap W^{1,2}_0(\Omega;\CC^N)$ to \eqref{eq:Cannarsa_Terreni_Vespri_3-1} and, moreover, $u$ obeys
\begin{equation}
\label{eq:Cannarsa_Terreni_Vespri_6-2}
(|\lambda|-\omega_1)\|u\|_{L^\infty(\Omega)} + (|\lambda|-\omega_1)^{1/2}\|u\|_{W^{1,\infty}(\Omega)}
\leq C\|f\|_{L^\infty(\Omega)}.
\end{equation}
\end{thm}

The \apriori estimate \eqref{eq:Cannarsa_Terreni_Vespri_6-2} in Theorem \ref{thm:Cannarsa_Terreni_Vespri_6-1} can be augmented with an \apriori estimate of the second-order derivatives of $u$, to give an \apriori estimate (albeit using $L^{2,\mu}$ Morrey spaces) analogous to that of Stewart in \cite[Equation (1.1)]{Stewart_1974} (who uses $L^q$ spaces) as described in the following remark.

\begin{rmk}[Supplementary \apriori estimate for the second-order derivatives of the solution $u$ in Theorem \ref{thm:Cannarsa_Terreni_Vespri_6-1}]
\cite[Remark 6.4]{Cannarsa_Terreni_Vespri_1985}
\label{rmk:Cannarsa_Terreni_Vespri_6-4}
Assume the hypotheses of Theorem \ref{thm:Cannarsa_Terreni_Vespri_6-1}. If $\mu \in (d-2, d)$ (see \cite[Lemma 6.2]{Cannarsa_Terreni_Vespri_1985} for this restriction on $\mu$), then
\begin{equation}
\label{eq:Cannarsa_Terreni_Vespri_6-2_second_order_derivatives}
(|\lambda|-\omega_1)^{(d-\mu)/2}\sup_{x \in \bar\Omega}\sum_{|\alpha|=2}\|D^\alpha u\|_{L^{2,\mu}(\Omega\cap B(x,r/2))}
\leq C\|f\|_{L^\infty(\Omega)},
\end{equation}
with $r = K(|\lambda|-\omega_1)^{-1/2}$ as in \cite[Equation (6.11)]{Cannarsa_Terreni_Vespri_1985}, where $K$ is a positive constant. Here, the Morrey space $L^{2,\mu}(\Omega; \CC^N)$, for $\mu \in (0, d)$, is defined via the norm \cite[pp. 59--60]{Cannarsa_Terreni_Vespri_1985}
\begin{equation}
\|u\|_{L^{2,\mu}(\Omega)}^2
:=
\sup_{\begin{subarray}{c}x \in \Omega \\ 0 < \varrho \leq \diam(\Omega)\end{subarray}}
\varrho^{-\mu} \|u\|_{L^2(\Omega\cap B(x,\varrho))}^2.
\end{equation}
See Troianiello \cite[Section 1.4]{Troianiello} for an introduction to the Morrey, John-Nirenberg (BMO), and Campanato spaces, $L^{2,\mu}(\Omega; \CC^N)$, with $0<\mu<d$, and $\mu = d$, and $d < \mu \leq d + 2$, respectively. We recall from \cite[Section 1.4.2]{Troianiello}, that
\begin{inparaenum}[\itshape a\upshape)]
\item $L^{2,\nu}(\Omega; \CC^N) \hookrightarrow L^{2,\mu}(\Omega; \CC^N)$ if $0 \leq \mu < \nu \leq d+2$; and
\item $L^{2,0}(\Omega; \CC^N) \cong L^{2,\nu}(\Omega; \CC^N)$ (as a Banach space); and
\item $L^p(\Omega; \CC^N) \hookrightarrow L^{2,\mu}(\Omega; \CC^N)$ if $p > 2$ and $\mu = d(p-2)p$; and
\item $L^{2,\mu}(\Omega; \CC^N) \cong C^{0,\delta}(\bar\Omega; \CC^N)$ (as a Banach space), when $d < \mu \leq d + 2$ and $\delta = (\mu - d)/2$ \cite[Theorem 1.17]{Troianiello}.
\end{inparaenum}
\end{rmk}

\begin{rmk}[Regularity of the solution $u$ in Theorem \ref{thm:Cannarsa_Terreni_Vespri_6-1}]
\label{rmk:Cannarsa_Terreni_Vespri_theorem_6-1_solution_regularity}
The Sobolev Embedding Theorem \cite[Theorem 4.12]{AdamsFournier} implies that $W^{1,\infty}(\Omega;\CC^N) \hookrightarrow C(\bar\Omega;\CC^N)$ and so we may write more simply in Theorem \ref{thm:Cannarsa_Terreni_Vespri_6-1} that $u \in H^2(\Omega;\CC^N)\cap W^{1,\infty}_0(\Omega;\CC^N)$. Moreover, there is a continuous embedding $L^\infty(\Omega;\CC^N) \hookrightarrow L^{2,\mu}(\bar\Omega;\CC^N)$ for any $\mu \in [0, d)$ by \cite[Remark 2.2 and pp. 87--88]{Cannarsa_Terreni_Vespri_1985}. Therefore, the function $f$ in Theorem \ref{thm:Cannarsa_Terreni_Vespri_6-1} belongs also to $L^{2,\mu}(\bar\Omega;\CC^N)$ and (see the discussion in \cite[p. 87]{Cannarsa_Terreni_Vespri_1985}) therefore $u$ belongs to $W^2_{(\mu)}(\Omega;\CC^N)$ by \cite[Theorem 5.1]{Cannarsa_Terreni_Vespri_1985}. Recall from \cite[Theorem 2.1 and p. 61]{Cannarsa_Terreni_Vespri_1985}) that $W^{2,2}_{(\mu)}(\Omega;\CC^N)$ is the subspace of functions $u \in W^{2,2}(\Omega;\CC^N)$ such that $D^\alpha u \in L^{2,\mu}(\bar\Omega;\CC^N)$ (when $0 \leq \mu < d$) for all $\alpha \in \NN^d$ with $|\alpha|=2$.
\end{rmk}

Following the suggestion in \cite[Remark 6.8]{Cannarsa_Terreni_Vespri_1985}, we now extend Theorem \ref{thm:Cannarsa_Terreni_Vespri_6-7} to allow for $1 < p < 2$ as well.

\begin{thm}[Existence, uniqueness, and \apriori estimate for a solution to a second-order elliptic system on $W^{2,p}(\Omega; \CC^N)$ when $1 < p < 2$]
\label{thm:Cannarsa_Terreni_Vespri_6-7_one_lessthan_p_lessthan_two}
\cite[Remark 6.8]{Cannarsa_Terreni_Vespri_1985}
Assume the hypotheses of Theorem \ref{thm:Cannarsa_Terreni_Vespri_6-7} and in addition\footnote{This additional hypothesis is omitted in \cite[Remark 6.8]{Cannarsa_Terreni_Vespri_1985} but such an assumption, or at least that the $a_\alpha$ belong to $W^{2,\infty}(\Omega)$ for $\alpha \in \NN^d$ with $|\alpha| \leq 2$, appears to be necessary when we rely an argument involving the formal adjoint, $\sA^*$, of $\sA$.}
that the coefficients $a_\alpha$ of $\sA$ in \eqref{eq:Cannarsa_Terreni_Vespri_1-1} belong to $C^2(\bar\Omega; \CC^N)$ for $\alpha \in \NN^d$ with $|\alpha| \leq 2$. Then there are a constant $\omega_1 \geq 0$ and, given $1 < p < 2$, a constant $C > 0$ with the following significance. For all $\lambda \in \CC$ obeying $\Real \lambda > \omega_1$ and $f \in L^p(\Omega;\CC^N)$, there is a unique solution $u \in W^{2,p}(\Omega;\CC^N)\cap W^{1,p}_0(\Omega;\CC^N)$ to \eqref{eq:Cannarsa_Terreni_Vespri_3-1} and, moreover, $u$ obeys the \apriori estimate \eqref{eq:Cannarsa_Terreni_Vespri_6-15}.
\begin{equation}
\label{eq:Cannarsa_Terreni_Vespri_apriori_estimate_one_lessthan_p_lessthan_two}
(|\lambda|-\omega_1)\|u\|_{L^p(\Omega)} + (|\lambda|-\omega_1)^{1/2}\|Du\|_{L^p(\Omega)} + \sum_{|\alpha| = 2} \|D^\alpha u\|_{L^p(\Omega)}
\leq C\|f\|_{L^p(\Omega)}.
\end{equation}
\end{thm}

\begin{proof}
We proceed by analogy with the proof of the Calder\'on-Zygmund Inequality for the case $1 < p < 2$ via duality (see, for example, step (v) in the proof of \cite[Theorem 9.9]{GilbargTrudinger}).

Let $q \in (2, \infty)$ be the Sobolev exponent dual to $p$, so $1/p + 1/q = 1$. Suppose $\alpha \in \NN^d$ with $|\alpha| \leq 2$ and recall that (see, for example, \cite[Lemma 2.7]{AdamsFournier}),
\begin{equation}
\label{eq:reverse_Holder_inequality_one_lessthan_p_lessthan_infinity}
\|D^\alpha u\|_{L^p(\Omega)}
=
\sup_{\begin{subarray}{c} \varphi \in C^\infty_0(\Omega; \CC^N) \\ \|\varphi\|_{L^q(\Omega)} \leq 1\end{subarray}}
(D^\alpha u, \varphi)_{L^2(\Omega)}.
\end{equation}
We begin by verifying the \apriori estimate \eqref{eq:Cannarsa_Terreni_Vespri_6-15} with right-hand side $f = (\sA + \lambda)u$ and $\lambda \in \CC$ with $\Real\lambda > \omega_1$ and note that, by continuity, it suffices to establish \eqref{eq:Cannarsa_Terreni_Vespri_6-15} when $u \in C^\infty_0(\Omega; \CC^N)$, and hence $f \in C^\infty_0(\Omega; \CC^N)$, which we shall now assume. Thus, noting that $C^\infty_0(\Omega; \CC^N) \subset L^p(\Omega; \CC^N)$, we
calculate\footnote{More precisely, if we temporarily use $\sA'$ to denote the formal adjoint of $\sA$ and $\cA^*$ to denote the abstract Banach-space adjoint of an operator $\cA$, then we are implicitly asserting that $(\sA_p)^* = \sA_q'$; see the proof of this equality in \cite[Equation (5.68)]{Tanabe_1997} for the case $N=1$.}
\begin{align*}
(D^\alpha u, \varphi)_{L^2(\Omega)}
&= (D^\alpha (\sA_p + \lambda)^{-1}f, \varphi)_{L^2(\Omega)}
\\
&= (f, (D^\alpha (\sA_p + \lambda)^{-1})^*\varphi)_{L^2(\Omega)}
\\
&= (f, (\sA_q^* + \lambda)^{-1}D^\alpha\varphi)_{L^2(\Omega)}
\\
&= (f, D^\alpha(\sA_q^* + \lambda)^{-1}\varphi)_{L^2(\Omega)}.
\end{align*}
To understand the origin of the last equality, we may suppose that the action of the integral operator, $(\sA_q^* + \lambda)^{-1}$ on $L^q(\Omega; \CC^N)$, is represented by a kernel, $K_\lambda \in C^2(\Omega\times\Omega - \Delta_\Omega; \CC^{N\times N})$, where we denote $\Delta_\Omega := \diag(\Omega\times\Omega)$; see, for example, \cite[Theorem 5.7]{Tanabe_1997} for the case $N=1$. Therefore,
\begin{align*}
(\sA_q^* + \lambda)^{-1}D^\alpha\varphi(x)
&= \int_\Omega K_\lambda(x,y)D^\alpha\varphi(y) \,dy
\\
&= \int_\Omega D^\alpha K_\lambda(x,y)\varphi(y) \,dy \quad \hbox{(integration by parts)}
\\
&= D^\alpha(\sA_q^* + \lambda)^{-1}\varphi(x), \quad\hbox{a.e. } x \in \Omega.
\end{align*}
Consequently,
\begin{align*}
(|\lambda|-\omega_1)^{1-|\alpha|/2}|(D^\alpha u, \varphi)_{L^2(\Omega)}|
&\leq \|f\|_{L^p(\Omega)} (|\lambda|-\omega_1)^{1-|\alpha|/2}\|D^\alpha (\sA_q^* + \lambda)^{-1} \varphi\|_{L^q(\Omega)}
\\
&\leq C\|f\|_{L^p(\Omega)} \|\varphi\|_{L^q(\Omega)}
\quad\hbox{(by Theorem \ref{thm:Cannarsa_Terreni_Vespri_6-7} and \eqref{eq:Cannarsa_Terreni_Vespri_6-15})}
\\
&\leq C\|f\|_{L^p(\Omega)}.
\end{align*}
Therefore, for $\alpha \in \NN^d$ with $|\alpha| \leq 2$, we see that
\begin{equation}
\label{eq:Cannarsa_Terreni_Vespri_3-1_one_leq_p_lessthan_2_prelim}
(|\lambda|-\omega_1)^{1-|\alpha|/2}\|D^\alpha u\|_{L^p(\Omega)}
\leq
C\|f\|_{L^p(\Omega)},
\end{equation}
and this yields the desired \apriori estimate \eqref{eq:Cannarsa_Terreni_Vespri_3-1} when $1 < p < 2$.

For existence (and uniqueness) of a solution $u \in W^{2,p}(\Omega;\CC^N)\cap W^{1,p}_0(\Omega;\CC^N)$ to $(\sA + \lambda)u = f$, observe that if $\tilde f \in C^\infty_0(\Omega)$, then $\tilde f \in L^2(\Omega;\CC^N)$ and so Theorem \ref{thm:Cannarsa_Terreni_Vespri_6-7} implies that there exists a (unique) $\tilde u \in W^{2,2}(\Omega;\CC^N)\cap W^{1,2}_0(\Omega;\CC^N)$ satisfying $(\sA + \lambda)\tilde u = \tilde f$ a.e. on $\Omega$, for $\lambda \in \CC$ with $\Real\lambda > \omega_1$. Clearly, since $\Omega$ is bounded and $1 < p < 2$, we also have $\tilde u \in W^{2,p}(\Omega;\CC^N)\cap W^{1,p}_0(\Omega;\CC^N)$. Because $C^\infty_0(\Omega)$ is dense in $L^p(\Omega;\CC^N)$, we see that range of the operator,
\begin{equation}
\label{eq:A_plus_lambda_from_W2p_cap_W1p0_to_Lp_vector_valued}
\sA + \lambda: W^{2,p}(\Omega;\CC^N)\cap W^{1,p}_0(\Omega;\CC^N) \to L^p(\Omega;\CC^N),
\end{equation}
is dense in $L^p(\Omega;\CC^N)$ for $\lambda \in \CC$ with $\Real\lambda > \omega_1$.

The operator $\sA + \lambda$ in \eqref{eq:A_plus_lambda_from_W2p_cap_W1p0_to_Lp_vector_valued} is clearly continuous but also the \apriori estimate \eqref{eq:Cannarsa_Terreni_Vespri_3-1} means that it is bounded below in the sense of \cite[Definition 2.1]{Abramovich_Aliprantis_2002} and so has closed range by \cite[Theorem 2.5]{Abramovich_Aliprantis_2002}. In particular, the operator \eqref{eq:A_plus_lambda_from_W2p_cap_W1p0_to_Lp_vector_valued} is surjective and one-to-one by \cite[Lemma 2.8]{Abramovich_Aliprantis_2002}. In other words, given $f \in L^p(\Omega;\CC^N)$, there exists a unique solution $u \in W^{2,p}(\Omega;\CC^N)\cap W^{1,p}_0(\Omega;\CC^N)$ to $(\sA + \lambda)u = f$ a.e. on $\Omega$ and this completes the proof.
\end{proof}

Because neither Theorem \ref{thm:Cannarsa_Terreni_Vespri_6-7} nor Theorem \ref{thm:Cannarsa_Terreni_Vespri_6-1} provide a $W^{2,\infty}$ \apriori estimate or solvability of \eqref{eq:Cannarsa_Terreni_Vespri_3-1} for $u \in W^{2,\infty}(\Omega;\CC^N) \cap W^{1,\infty}_0(\Omega;\CC^N)$ given $f \in L^\infty(\Omega;\CC^N)$, the proof of Theorem \ref{thm:Cannarsa_Terreni_Vespri_6-7_one_lessthan_p_lessthan_two} for $1 < p < 2$ does not entirely extend to the case $p = 1$ as well. However, we can still extract a partial analogue of Theorem \ref{thm:Cannarsa_Terreni_Vespri_6-7_one_lessthan_p_lessthan_two} for $p = 1$ which is enough for the purpose of generation of analytic semigroups. First, we note that the proof of Lemma \ref{lem:Pazy_7-3-8} extends with minimal change to give\footnote{The result may also be extracted as a special case of Lemma \ref{lem:Pazy_7-3-8_vector_bundle_manifold_Taubes_weight}, whose proof we do include.}

\begin{lem}[Reverse H\"older inequality for $p=1$]
\label{lem:Pazy_7-3-8_vector_valued}
Let $\Omega \subseteqq \RR^d$ be an open subset and $N\geq 1$ an integer. If $u \in L^1(\Omega;\CC^N)$, then
\begin{equation}
\label{eq:Pazy_7-3-26_vector_valued}
\|u\|_{L^1(\Omega)}
=
\sup_{\begin{subarray}{c} \varphi \in C^\infty_0(\Omega; \CC^N) \\ \|\varphi\|_{L^\infty(\Omega)} \leq 1\end{subarray}}
(u, \varphi)_{L^2(\Omega)}.
\end{equation}
\end{lem}

The following result may be compared with \cite[Theorem 2.4]{Angiuli_Pallara_Paronetto_2010} (scalar, second-order elliptic operator), \cite[Theorem 3.3]{Guidetti_1993} (elliptic system of order $2m$ with $m \geq 1$), and \cite[Theorem 5.8]{Tanabe_1997} (scalar elliptic operator of order $m \geq 1$).

\begin{thm}[Existence, uniqueness, and \apriori $L^1$ estimate and domain for an elliptic system]
\label{thm:Cannarsa_Terreni_Vespri_6-7_p_is_one}
Assume the hypotheses of Theorem \ref{thm:Cannarsa_Terreni_Vespri_6-7} for the case $1 < p < 2$. Then the operator,
$$
\sA: C^2(\bar\Omega; \CC^N)\cap C_0(\bar\Omega; \CC^N) \subset L^1(\Omega; \CC^N) \to L^1(\Omega; \CC^N),
$$
is closable and, if $\sA_1$ denotes the smallest closed extension, one has
\begin{equation}
\label{eq:Cannarsa_Terreni_Vespri_apriori_function_domain_A1_partial_characterization}
\sD(\sA_1) \subset \{u \in L^1(\Omega; \CC^N): \sA u \in L^1(\Omega; \CC^N)\}.
\end{equation}
Moreover, there are constants $C > 0$ and $\omega_1 \geq 0$ with the following significance. If $\lambda \in \CC$ obeys $\Real \lambda > \omega_1$ and $u \in \sD(\sA_1)$, then
\begin{equation}
\label{eq:Cannarsa_Terreni_Vespri_apriori_function_and_gradient_estimate_p_is_one}
(|\lambda|-\omega_1)\|u\|_{L^1(\Omega)} + (|\lambda|-\omega_1)^{1/2}\|Du\|_{L^1(\Omega)}
\leq C\|(\sA_1 + \lambda)\|_{L^1(\Omega)}.
\end{equation}
Finally, for all $\lambda \in \CC$ with $\Real\lambda > \omega_1$ and $f \in L^1(\Omega;\CC^N)$, there exists a unique solution $u \in \sD(\sA_1)$ to \eqref{eq:Cannarsa_Terreni_Vespri_3-1}.
\end{thm}

\begin{proof}
The fact that $\sA: C^2(\bar\Omega; \CC^N)\cap C_0(\bar\Omega; \CC^N) \subset L^1(\Omega; \CC^N) \to L^1(\Omega; \CC^N)$ is closable follows, \emph{mutatis mutandis}, from the proof for the example in \cite[Section 2.6]{Yosida}; compare the proofs of \cite[Lemma 9.1]{Amann_1983} and \cite[Lemma 2.3]{Angiuli_Pallara_Paronetto_2010}. The fact that \eqref{eq:Cannarsa_Terreni_Vespri_apriori_function_domain_A1_partial_characterization} holds is immediate from the definition of the smallest closed extension.

Consider $\lambda \in \CC$ with $\Real \lambda > \omega_1$, where $\omega_1$ is as in Theorem \ref{thm:Cannarsa_Terreni_Vespri_6-1}. Suppose first that $u \in C^\infty_0(\Omega; \CC^N)$ and set $f = (\sA + \lambda)u \in C^\infty_0(\Omega; \CC^N)$. When the multi-index $\alpha \in \NN^d$ obeys $|\alpha|=0$ or $1$, then the proof of the \apriori estimate \eqref{eq:Cannarsa_Terreni_Vespri_apriori_estimate_one_lessthan_p_lessthan_two} provided by Theorem \ref{thm:Cannarsa_Terreni_Vespri_6-7_one_lessthan_p_lessthan_two} extends to give the \apriori estimate \eqref{eq:Cannarsa_Terreni_Vespri_apriori_function_and_gradient_estimate_p_is_one} by making the following changes:
\begin{enumerate}
\item The role of the reverse H\"older inequality \eqref{eq:reverse_Holder_inequality_one_lessthan_p_lessthan_infinity} when $1<p<2$ is replaced by that of \eqref{eq:Pazy_7-3-26_vector_valued} when $p = 1$ in Lemma \ref{lem:Pazy_7-3-8_vector_valued};
\item The role of the \apriori estimate \eqref{eq:Cannarsa_Terreni_Vespri_6-15} in Theorem \ref{thm:Cannarsa_Terreni_Vespri_6-7} is replaced by that of \eqref{eq:Cannarsa_Terreni_Vespri_6-2} in Theorem \ref{thm:Cannarsa_Terreni_Vespri_6-1}.
\end{enumerate}
The fact that the \apriori estimate \eqref{eq:Cannarsa_Terreni_Vespri_apriori_function_and_gradient_estimate_p_is_one} continues to hold for $u \in \sD(\sA_1)$ follows by continuity, the definition of $\sD(\sA_1)$, and the fact that $C^\infty_0(\Omega; \CC^N)$ is dense in $\sD(\sA_1)$.

For existence (and uniqueness) of a solution $u \in \sD(\sA_1)$ to $(\sA + \lambda)u = f$, given $f \in L^1(\Omega; \CC^N)$, observe that if $\tilde f \in C^\infty_0(\Omega)$, then $\tilde f \in L^2(\Omega;\CC^N)$ and so Theorem \ref{thm:Cannarsa_Terreni_Vespri_6-7} implies that there exists a (unique) $\tilde u \in W^{2,2}(\Omega;\CC^N)\cap W^{1,2}_0(\Omega;\CC^N)$ satisfying $(\sA + \lambda)\tilde u = \tilde f$ a.e. on $\Omega$, for $\lambda \in \CC$ with $\Real\lambda > \omega_1$. Clearly, since $\Omega$ is bounded, we also have $\tilde u \in L^1(\Omega;\CC^N)$ and, by definition of $\sD(\sA_1)$, we see that $\tilde u$ belongs to $\sD(\sA_1)$. Because $C^\infty_0(\Omega)$ is dense in $L^1(\Omega;\CC^N)$, we see that range of the operator,
\begin{equation}
\label{eq:A_plus_lambda_from_domain_A1_to_L1_vector_valued}
\sA_1 + \lambda: \sD(\sA_1) \to L^1(\Omega;\CC^N),
\end{equation}
is dense in $L^1(\Omega;\CC^N)$ for $\lambda \in \CC$ with $\Real\lambda > \omega_1$.

By construction, the operator $\sA_1$ and thus also $\sA_1 + \lambda$ have closed range. In particular, the operator $\sA_1 + \lambda$ is surjective since its range is closed and dense in $L^1(\Omega;\CC^N)$ and thus is equal to $L^1(\Omega;\CC^N)$ for $\lambda \in \CC$ with $\Real\lambda > \omega_1$. The \apriori estimate \eqref{eq:Cannarsa_Terreni_Vespri_apriori_function_and_gradient_estimate_p_is_one} ensures that $\sA_1 + \lambda$ is one-to-one for $\lambda \in \CC$ with $\Real\lambda > \omega_1$. In other words, given $f \in L^1(\Omega;\CC^N)$, there exists a unique solution $u \in \sD(\sA_1)$ to $(\sA + \lambda)u = f$ a.e. on $\Omega$ and this completes the proof.
\end{proof}

\begin{rmk}[On the domain of an elliptic partial differential operator on $L^1(\Omega; \CC^N)$]
\label{rmk:Domain_elliptic_partial_differential_operator_on_L1}
A linear operator $\cA: \sD(\cA) \subset \cW \to \calV$, where $\cW$ and $\calV$ are Banach spaces, is said to be \emph{closed} \cite[Section 2.6]{Yosida} if
\begin{multline}
\label{eq:Yosida_2-6-2}
\{w_n\}_{n\in\NN} \subset \sD(\cA) \hbox{ and } w_n \xrightarrow{\cW} w \in \cW \hbox{ and } \cA w_n \xrightarrow{\calV} v \in \calV
\\
\implies w \in \sD(\cW) \hbox{ and } \cA w = v,
\end{multline}
that is, if and only if the graph $G(\cA) := \{(w, \cA w): w \in \cW\}$ is a closed subspace of $\cW\times \calV$. A linear operator $\cA: \sD(\cA) \subset \cW \to \calV$ is called \emph{closable} or \emph{pre-closed} if the closure in $\cW \times \calV$ of the graph $G(\cA)$ is the graph of a linear operator, say $\cB: \sD(\cB) \subset \cW \to \calV$ \cite[Section 2.6]{Yosida}. We recall from \cite[Proposition 2.6.2]{Yosida} and its proof that a linear operator $\cA: \sD(\cA) \subset \cW \to \calV$ closable if and only if the following condition is satisfied:
\begin{equation}
\label{eq:Yosida_2-6-3}
\{w_n\}_{n\in\NN} \subset \sD(\cA) \hbox{ and } w_n \xrightarrow{\cW} 0 \hbox{ and } \cA w_n \xrightarrow{\calV} v \in \calV
\implies v = 0.
\end{equation}
If $\cA$ obeys \eqref{eq:Yosida_2-6-3}, then the \emph{smallest closed extension} of $\cA$ is the linear operator $\cB: \sD(\cB) \subset \cW \to \calV$ defined by
\begin{multline}
\label{eq:Yosida_2-6-4}
w \in \sD(\cB) \iff \exists\, \{w_n\}_{n\in\NN} \subset \sD(\cA) \hbox{ with } w_n \xrightarrow{\cW} w \hbox{ and } \cA w_n \xrightarrow{\calV} v \in \calV,
\\
\hbox{in which case } \cB w := v.
\end{multline}
The condition \eqref{eq:Yosida_2-6-3} ensures that the value $v = \cB w$ in \eqref{eq:Yosida_2-6-4} is defined uniquely by $w \in \sD(\cB)$.

A second-order partial differential operator $\sA$ as in \eqref{eq:Cannarsa_Terreni_Vespri_1-1}, with coefficients in $C^2(\Omega; \CC^{N\times N})$, is closable in $L^2(\Omega; \CC^N)$ by \cite[Section 2.6]{Yosida} and more generally, by an identical argument, in $L^p(\Omega; \CC^N)$ when $1 < p < \infty$. In that situation, the \apriori estimate \eqref{eq:Cannarsa_Terreni_Vespri_6-15} identifies the domain of the smallest closed extension $\sA_p$ of $\sA: C^\infty_0(\Omega; \CC^N) \subset L^p(\Omega; \CC^N) \to L^p(\Omega; \CC^N)$ as $\sD(\sA_p) := W^{2,p}(\Omega; \CC^N) \cap W^{1,p}_0(\Omega; \CC^N)$. The same argument gives closability of partial differential operators with order $m\geq 1$ and coefficients in $C^m(\Omega; \CC^{N\times N})$, just as in \cite[Section 2.6]{Yosida}.

In the absence of an \apriori estimate such as \eqref{eq:Cannarsa_Terreni_Vespri_6-15} when $p=1$, various alternative characterizations of $\sD(\sA_1)$ in the case $N = 1$ (scalar elliptic partial differential operator ) have been provided by Amann \cite[Equation (9.1)]{Amann_1983}, Angiuli, Pallara, and Paronetto \cite[p. 259]{Angiuli_Pallara_Paronetto_2010}, Di Blasio \cite{DiBlasio_1991} (order $m = 2$), Lunardi and Metafune \cite{Lunardi_Metafune_2004}, Pazy \cite[Definition 7.3.9]{Pazy_1983} (order $2m$ with $m \geq 1$), Tanabe \cite[Theorem 5.8]{Tanabe_1997} (order $m \geq 1$), Vespri \cite[p. 103]{Vespri_1991} (order $m = 2$) and by Guidetti \cite[Theorem 3.3]{Guidetti_1993} in the case $N \geq 1$ (elliptic system of order $2m$ with $m \geq 1$); the most refined characterizations of $\sD(\sA_1)$ are provided by Tanabe \cite[Theorem 5.8]{Tanabe_1997} and Guidetti \cite[Theorem 3.3]{Guidetti_1993}.
\end{rmk}

\begin{rmk}[\Apriori $L^1$ estimates]
\label{rmk:L1_apriori_gradient_estimate_for_an_elliptic_partial_differential_operator}
We note that \apriori $L^1$ estimates related to \eqref{eq:Cannarsa_Terreni_Vespri_apriori_function_and_gradient_estimate_p_is_one} are also provided by Amann \cite[Proposition 9.2]{Amann_1983}, Angiuli, Pallara, and Paronetto \cite[Theorem 3.5]{Angiuli_Pallara_Paronetto_2010}(2010), Brezis and Strauss \cite[Theorem 8]{Brezis_Strauss_1973}, Cannarsa and Vespri \cite[Lemma 4.3 and Theorem 4.4]{Cannarsa_Vespri_1988}, Tanabe (implicitly) \cite[Theorem 5.8]{Tanabe_1997}, and Vespri \cite[Remark 1 on page 103]{Vespri_1991}.
\end{rmk}

By analogy with Definition \ref{defn:Pazy_7-3-3} for the case $N=1$, we define, for $1 < p < \infty$,
\begin{align}
\label{eq:Pazy_7-3-5_system}
\sD(\sA_p) &:= W^{2,p}(\Omega;\CC^N)\cap W^{1,p}_0(\Omega;\CC^N) \subset L^p(\Omega;\CC^N),
\\
\label{eq:Pazy_7-3-6_system}
\sA_p u &:= \sA(x, D)u, \quad \forall\, u \in \sD(\sA_p).
\end{align}
For the case $p=\infty$, by analogy with Definition \ref{defn:Pazy_7-3-3_L_infinity} for the case $N=1$, we set
\begin{multline}
\label{eq:Pazy_7-3-20_system}
\sD(\sA_\infty) := \{u \in W^{1,\infty}_0(\Omega;\CC^N): D^\alpha u\in L^{2,\mu}(\Omega;\CC^N),
\\
\forall\, \mu \in [0, d) \hbox{ and } \alpha \in \NN^d \hbox{ with } |\alpha| = 2, \hbox{ and } \sA(x, D)u \in L^\infty(\Omega;\CC^N)\},
\end{multline}
and
\begin{equation}
\label{eq:Pazy_7-3-21_system}
\sA_\infty u := \sA(x, D)u, \quad\forall\, u \in \sD(\sA_\infty).
\end{equation}
See Remark \ref{rmk:Cannarsa_Terreni_Vespri_6-4} for a discussion of the Morrey spaces, $L^{2,\mu}(\Omega;\CC^N)$.

As suggested by Cannarsa, Terreni, and Vespri in \cite{Cannarsa_Terreni_Vespri_1985}, Theorems \ref{thm:Cannarsa_Terreni_Vespri_6-7}, \ref{thm:Cannarsa_Terreni_Vespri_6-1}, and \ref{thm:Cannarsa_Terreni_Vespri_6-7_one_lessthan_p_lessthan_two} quickly lead to the following results on the generation of analytic semigroups. A result similar to Corollary \ref{cor:Cannarsa_Terreni_Vespri_6-7_sectorial_operator_analytic_semigroup_LpOmega} below in the case $p = 1$ has been obtained by Guidetti \cite[Theorem 3.3]{Guidetti_1993}.

\begin{cor}[Resolvent estimate for a second-order elliptic system when $1 \leq p \leq \infty$ and sectorial property and generation of an analytic semigroup on $L^p(\Omega; \CC^N)$ when $1 \leq p < \infty$]
\label{cor:Cannarsa_Terreni_Vespri_6-7_sectorial_operator_analytic_semigroup_LpOmega}
Assume the hypotheses of Theorem \ref{thm:Cannarsa_Terreni_Vespri_6-7} and, in addition when $1 \leq p < 2$, that the coefficients $a_\alpha$ of $\sA$ in \eqref{eq:Cannarsa_Terreni_Vespri_1-1} belong to $C^2(\bar\Omega; \CC^N)$ for $\alpha \in \NN^d$ with $|\alpha| \leq 2$. Suppose $1 \leq p \leq \infty$. Then there are constants, $C>0$ and $\omega_0\geq 0$ and $\vartheta \in (\pi/2, \pi)$, such that $\rho(-\sA_p) \supset \Delta_\vartheta(\omega_0)$ and
\begin{equation}
\label{eq:Cannarsa_Terreni_Vespri_Lp_sectorial_estimate}
\|(\sA_p + \lambda)^{-1}\|_{\sL(L^p(\Omega;\CC^N))} \leq \frac{C}{|\lambda - \omega_0|},
\quad \forall\, \lambda \in \Delta_\vartheta(\omega_0).
\end{equation}
Moreover, restricting to the case $1 \leq p < \infty$, we have that $\sA_p$ is a sectorial operator on $L^p(\Omega; \CC^N)$ in the sense of Definition \ref{defn:Sell_You_page_78_definition_of_sectorial_operator} and $-\sA_p$ generates an analytic semigroup, $e^{-\sA_p t}$, on $L^p(\Omega; \CC^N)$.
\end{cor}

\begin{proof}
By Remarks \ref{rmk:Cannarsa_Terreni_Vespri_6-4} and \ref{rmk:Cannarsa_Terreni_Vespri_theorem_6-1_solution_regularity} and the definition \eqref{eq:Pazy_7-3-20_system} of $\sD(\sA_\infty)$, we see that the \apriori estimate \eqref{eq:Cannarsa_Terreni_Vespri_6-2} holds for all $u \in \sD(\sA_\infty)$.

As a consequence of Theorems \ref{thm:Cannarsa_Terreni_Vespri_6-7}, \ref{thm:Cannarsa_Terreni_Vespri_6-1}, \ref{thm:Cannarsa_Terreni_Vespri_6-7_one_lessthan_p_lessthan_two}, and \ref{thm:Cannarsa_Terreni_Vespri_6-7_p_is_one} one obtains, for $1 \leq p \leq \infty$ and all $\lambda \in \CC$ obeying $\Real \lambda > \omega_1$, that
\begin{equation}
\label{eq:Cannarsa_Terreni_Vespri_6-15_agmon_estimate_preliminary}
\|u\|_{L^p(\Omega)} \leq \frac{C}{|\lambda|-\omega_1}\|(\sA_p + \lambda)u\|_{L^p(\Omega)},
\quad \forall\, u \in \sD(\sA_p).
\end{equation}
Suppose furthermore that $\lambda \in \CC$ obeys $\Real \lambda > 3\omega_1$. Using $|\lambda - 3\omega_1| \leq |\lambda| + 3\omega_1$ and thus $|\lambda| \geq |\lambda - 3\omega_1| - 3\omega_1$, we see that
$$
|\lambda| - 3\omega_1 \geq |\lambda - 3\omega_1| - 6\omega_1 > |\lambda - 3\omega_1| - 2|\lambda|,
$$
where we use $|\lambda| > 3\omega_1$, since $\Real \lambda > 3\omega_1$, to obtain the second inequality. Therefore,
$$
3(|\lambda| - \omega_1) > |\lambda - 3\omega_1|, \quad \Real \lambda > 3\omega_1,
$$
that is,
$$
|\lambda| - \omega_1 > \frac{1}{3}|\lambda - 3\omega_1|, \quad\forall\, \lambda \in \CC \hbox{ such that } \Real \lambda > 3\omega_1.
$$
Hence, by combining the preceding inequality with \eqref{eq:Cannarsa_Terreni_Vespri_6-15_agmon_estimate_preliminary}, we find that, for $1 \leq p \leq \infty$,
\begin{equation}
\label{eq:Cannarsa_Terreni_Vespri_6-15_agmon_estimate}
\|u\|_{L^p(\Omega)} \leq \frac{3C}{|\lambda - 3\omega_1|}\|(\sA_p + \lambda)u\|_{L^p(\Omega)},
\quad \forall\, \lambda \in \Delta_{\pi/2}(3\omega_1) \hbox{ and } u \in \sD(\sA_p),
\end{equation}
where we recall from \eqref{eq:Sell_You_page_77_complex_plane_sector_definition_Delta_of_a} that $\Delta_{\pi/2}(3\omega_1) = \{\lambda \in \CC: \Real \lambda > 3\omega_1\}$. From Theorems \ref{thm:Cannarsa_Terreni_Vespri_6-7}, \ref{thm:Cannarsa_Terreni_Vespri_6-1}, \ref{thm:Cannarsa_Terreni_Vespri_6-7_one_lessthan_p_lessthan_two}, and \ref{thm:Cannarsa_Terreni_Vespri_6-7_p_is_one}, we obtain that, for any $\lambda \in \Delta_{\pi/2}(3\omega_1)$, the operators,
$$
\sA_p + \lambda: W^{2,p}(\Omega;\CC^N)\cap W^{1,p}_0(\Omega;\CC^N) \to L^p(\Omega;\CC^N), \quad 1 \leq p \leq \infty,
$$
are surjective (and necessarily one-to-one by \eqref{eq:Cannarsa_Terreni_Vespri_6-15_agmon_estimate}) and the \apriori estimate \eqref{eq:Cannarsa_Terreni_Vespri_6-15_agmon_estimate} yields the estimate,
\begin{equation}
\label{eq:Cannarsa_Terreni_Vespri_6-15_resolvent_estimate}
\|(\sA_p + \lambda)^{-1}\|_{\sL(L^p(\Omega))} \leq \frac{3C}{|\lambda - 3\omega_1|},
\quad\forall\, \lambda \in \Delta_{\pi/2}(3\omega_1).
\end{equation}
In particular, for $\lambda \in \Delta_{\pi/2}(3\omega_1)$, the operators
$$
(\sA_p + \lambda)^{-1}: L^p(\Omega;\CC^N) \to L^p(\Omega;\CC^N), \quad 1 \leq p \leq \infty,
$$
are bounded and so
$$
\Delta_{\pi/2}(3\omega_1) \subset \rho(-\sA_p), \quad 1 \leq p \leq \infty,
$$
where $\rho(-\sA_p)$ is the resolvent set of $-\sA_p$ (Definition \ref{defn:Resolvent_set}).

Restricting now to the case $1 \leq p < \infty$, we observe that the domain $\sD(\sA_p)$ is dense in $L^p(\Omega;\CC^N)$ when $1 \leq p < \infty$. Therefore, Theorem \ref{thm:Renardy_Rogers_12-31} implies that the resolvent set $\rho(-\sA_p)$ thus actually contains a larger sector, $\Delta_{\pi/2 + \eps}(3\omega_1)$, for some $\eps > 0$, and that $\sA_p$ is a sectorial operator on $L^p(\Omega;\CC^N)$ in the sense of Definition \ref{defn:Sell_You_page_78_definition_of_sectorial_operator}, with constants $a = -3\omega_1$, $M := \max\{1, 3C\}$, and $\delta = \pi/2+\eps$, and that $-\sA_p$ is the generator of an analytic semigroup, $e^{-\sA_p t}$, on $L^p(\Omega;\CC^N)$.
\end{proof}

For the case of an analytic semigroup on $C_0(\bar\Omega;\CC^N) := \{u \in C(\bar\Omega;\CC^N): u = 0 \hbox{ on } \partial\Omega\}$, by analogy with Definition \ref{defn:Pazy_7-3-3_C} for the case $N=1$, we set
\begin{equation}
\label{eq:Pazy_7-3-23_system}
\sD(\sA_c) := \{u \in \sD(\sA_\infty): \sA(x, D)u \in C_0(\bar\Omega;\CC^N)\},
\end{equation}
and observe that $\sD(\sA_c)$ is dense in $C_0(\bar\Omega;\CC^N)$ by the Sobolev Embedding Theorem \cite[Theorem 4.12]{AdamsFournier}.

\begin{cor}[Sectorial property of a second-order elliptic system and generation of an analytic semigroup on $C_0(\bar\Omega;\CC^N)$]
\label{cor:Cannarsa_Terreni_Vespri_6-7_sectorial_operator_analytic_semigroup_C0barOmega}
Assume the hypotheses of Theorem \ref{thm:Cannarsa_Terreni_Vespri_6-7}. Then there are constants, $C>0$ and $\omega_0\geq 0$ and $\vartheta \in (\pi/2, \pi)$, such that $\rho(-\sA_c) \supset \Delta_\vartheta(\omega_0)$ and
\begin{equation}
\label{eq:Cannarsa_Terreni_Vespri_C0_sectorial_estimate}
\|(\sA_c + \lambda)^{-1}\|_{\sL(C_0(\bar\Omega;\CC^N))} \leq \frac{C}{|\lambda - \omega_0|},
\quad \forall\, \lambda \in \Delta_\vartheta(\omega_0).
\end{equation}
Moreover, $\sA_c$ is a sectorial operator on $C_0(\bar\Omega;\CC^N)$ in the sense of Definition \ref{defn:Sell_You_page_78_definition_of_sectorial_operator} and $-\sA_c$ generates an analytic semigroup, $e^{-\sA_c t}$, on $C_0(\bar\Omega;\CC^N)$.
\end{cor}

\begin{proof}
The conclusion follows from Corollary \ref{cor:Cannarsa_Terreni_Vespri_6-7_sectorial_operator_analytic_semigroup_LpOmega} applying to the case $p=\infty$ and Theorem \ref{thm:Renardy_Rogers_12-31}, given the additional fact that $\sD(\sA_c)$ is dense in $C_0(\bar\Omega;\CC^N)$.
\end{proof}


\subsection{Existence, uniqueness, \apriori $L^p$ estimates, and analytic semigroups on $L^p(\RR^d;\CC^N)$ defined by elliptic partial differential systems of order $m \geq 1$}
\label{subsec:Haller-Dintelmann_Heck_Hieber_2006}
In this subsection, we review the existence and uniqueness results and \apriori $L^p$ estimates, when $1<p<\infty$, for a parameter-elliptic partial differential system of order $m\geq 1$ on $\RR^d$, due to Haller-Dintelmann, Heck, and Hieber \cite{Haller-Dintelmann_Heck_Hieber_2006}, together with their consequences for resolvent estimates and analytic semigroup generation results on $L^p(\RR^d;\CC^N)$. The coefficients of the principal symbol of the partial differential system are allowed to belong to $L^\infty(\RR^d;\CC^{N\times N}) \cap \VMO(\RR^d;\CC^{N\times N})$. Closely related results are due to Denk, and Hieber, and Pr\"uss \cite{Denk_Hieber_Pruss_2003}.

We consider systems of differential operators of the form,
\begin{equation}
\label{eq:Haller-Dintelmann_Heck_Hieber_page_720_definition_differential_system}
\sA := \sum_{|\alpha|\leq m} a_\alpha(x) D^\alpha,
\end{equation}
where $D = -i(\partial_{x_1}, \ldots, \partial_{x_d})$ and $\alpha = (\alpha_1,\ldots,\alpha_d)$ is a multi-index and $a_\alpha \in L^\infty(\RR^d;\CC^{N\times N})$.

\begin{defn}[Parameter-elliptic partial differential operator]
\label{defn:Haller-Dintelmann_Heck_Hieber_page_720_parameter_elliptic}
\cite[p. 720]{Haller-Dintelmann_Heck_Hieber_2006}, \cite[p. 301]{Heck_Hieber_Stavrakidis_2010}
The differential operator $\sA$ in \eqref{eq:Haller-Dintelmann_Heck_Hieber_page_720_definition_differential_system} is $(\sM,\theta)$ \emph{parameter-elliptic} on $\RR^d$ if there are constants $\sM>0$ and $\theta \in [0,\pi)$ such that the principal symbol, $\mathring{a}(x,\xi) := \sum_{|\alpha| = m} a_\alpha(x) \xi^\alpha \in \CC^{N\times N}$ for $\xi \in \RR^d$, obeys
\begin{gather}
\label{eq:Heck_Hieber_Stavrakidis_6}
\sigma(\mathring{a}(x,\xi)) \subset \bar\Delta_\theta,
\\
\label{eq:Heck_Hieber_Stavrakidis_7}
\|\mathring{a}(x,\xi)^{-1}\| \leq \sM, \quad \forall\, \xi \in \RR^d \hbox{ such that } |\xi| = 1, \hbox{ and a.e. } x \in \RR^d,
\end{gather}
where $\Delta_\theta$ is as in \eqref{eq:Sell_You_page_77_complex_plane_sector_definition_Delta_of_a} and $\sigma(B)$ denotes the spectrum of a matrix $B \in \CC^{N\times N}$.
\end{defn}

For example, if the principal symbol, $\mathring{a}(x,\xi)$, is real scalar multiple of the $N\times N$ identity matrix for a.e. $x \in \RR^d$ and all $\xi \in \RR^d$ with $|\xi|=1$, then $\sigma(\mathring{a}(x,\xi)) \subset (\kappa, \sK)$, for some positive constants, $\kappa$ and $\sK$, and $\sA$ will be parameter-elliptic with $\sM \leq \kappa^{-1}$ and $\theta = 0$. Similarly, if $m = 2$ and the principal symbol, $\mathring{a}(x,\xi)$, is strictly elliptic in the sense of Definition \ref{defn:Cannarsa_Terreni_Vespri_2-7} with constant of ellipticity, $\kappa$, then $\sA$ will be parameter-elliptic with $\sM \leq \kappa^{-1}$ and $\theta = \pi/2$.

Given $a \in L_{\loc}^1(\RR^d; \CC^N \times \CC^N)$, denote \cite[p. 721]{Haller-Dintelmann_Heck_Hieber_2006}, for any bounded open subset $G \subset \RR^d$,
$$
a_G := \frac{1}{|G|}\int_G a(x)\, dx,
$$
where $|G|$ denotes the Lebesgue measure of $G$, and define
\begin{equation}
\label{eq:Haller-Dintelmann_Heck_Hieber_page_721_definition_star_norm_and_BMO}
\|a\|_* := \sup_{x_0 \in \RR^d, r > 0} \frac{1}{|B_r(x_0)|}\int_{B_r(x_0)} |a(x) - a_{B_r(x_0)}|\,dx,
\end{equation}
where $B_r(x_0) := \{x\in\RR^d: |x-x_0| < r\}$. One says that a function $a \in L_{\loc}^1(\RR^d; \CC^N \times \CC^N)$ has \emph{bounded mean oscillation} (BMO), or $a \in \BMO(\RR^d; \CC^N \times \CC^N)$, if $\|a\|_* < \infty$. For $r > 0$ and $a \in \BMO(\RR^d; \CC^N \times \CC^N)$, define
\begin{equation}
\label{eq:Haller-Dintelmann_Heck_Hieber_page_721_definition_eta_function}
\eta_a(r) := \sup_{x_0 \in \RR^d} \frac{1}{|B_r(x_0)|}\int_{B_r(x_0)} |a(x) - a_{B_r(x_0)}|\,dx.
\end{equation}
One says that a function $a \in \BMO(\RR^d; \CC^N \times \CC^N)$ has \emph{vanishing mean oscillation} (VMO), or $a \in \VMO(\RR^d; \CC^N \times \CC^N)$, if
$$
\lim_{r \downarrow 0} \eta_a(r) = 0.
$$
Suppose now that $a \in C(\bar\RR^d; \CC^N \times \CC^N)$ with modulus of continuity $\omega$, so
\begin{equation}
\label{eq:Krylov_5-2-2}
|a(x) - a(y)| \leq \omega(|x-y|), \quad \forall\, x, y \in \RR^d.
\end{equation}
Then, for any $x_0 \in \RR^d$,
$$
|a_{B_r(x_0)} - a(x_0)| = \left|\frac{1}{|B_r(x_0)|}\int_{B_r(x_0)} (a(x) - a(x_0))\, dx\right| \leq \omega(r), \quad \forall\, r > 0,
$$
and so the definition \eqref{eq:Haller-Dintelmann_Heck_Hieber_page_721_definition_eta_function} of $\eta_a$ yields
\begin{equation}
\label{eq:Krylov_eta_a_bounded_by_modulus_of_continuity}
\eta_a(r) \leq \omega(r), \quad \forall\, r > 0,
\end{equation}
and clearly $a \in \VMO(\RR^d; \CC^N \times \CC^N)$.

We next recall the main results of Haller-Dintelmann, Heck, and Hieber \cite{Haller-Dintelmann_Heck_Hieber_2006} for the generation of an analytic semigroup on $L^p(\RR^d; \CC^N)$ when $1<p<\infty$. We include Theorem \ref{thm:Haller-Dintelmann_Heck_Hieber_3-1} below for the sake of completeness, though it plays no essential role in the sequel.

\begin{thm}[Resolvent estimate for a parametric-elliptic operator on $L^p(\RR^d; \CC^N)$ when $1<p<\infty$]
\label{thm:Haller-Dintelmann_Heck_Hieber_3-1}
\cite[Theorem 3.1 and Remark 3.3]{Haller-Dintelmann_Heck_Hieber_2006}
Let $d \geq 2$ and $N \geq 1$ be integers, $\sM > 0$, $p \in (1, \infty)$, $\theta_0 \in (0, \pi)$, and $\theta \in [0, \theta_0)$. Let $\sA$ in \eqref{eq:Haller-Dintelmann_Heck_Hieber_page_720_definition_differential_system} be an $(\sM, \pi-\theta_0)$-parameter elliptic partial differential operator of order $m \geq 1$ with $a_\alpha \in L^\infty(\RR^d; \CC^{N\times N}) \cap \VMO(\RR^d; \CC^{N\times N})$ when $|\alpha|=m$, and $a_\alpha \in L^\infty(\RR^d; \CC^{N\times N})$ when $|\alpha| < m$. Then there are non-negative constants, $C$ and $\lambda_0$, such that $\rho(-\sA_p) \supset \Delta_\theta(\lambda_0)$ and
\begin{equation}
\label{eq:Haller-Dintelmann_Heck_Hieber_3-1}
\|(\sA_p + \lambda)^{-1}\|_{\sL(L^p(\RR^d; \CC^N))} \leq \frac{C}{|\lambda-\lambda_0|},
\quad \forall\, \lambda \in \Delta_\theta(\lambda_0),
\end{equation}
where $\Delta_\theta(\lambda_0) \subset \CC$ is as in \eqref{eq:Sell_You_page_77_complex_plane_sector_definition_Delta_of_a} and the realization $\sA_p:\sD(\sA_p) \subset L^p(\RR^d; \CC^N) \to L^p(\RR^d; \CC^N)$ is defined by \cite[Equation (2.2)]{Haller-Dintelmann_Heck_Hieber_2006},
$$
\sD(\sA_p) := W^{m,p}(\RR^d; \CC^N) \quad\hbox{and}\quad \sA_p u = \sA(x,D)u,
\quad\forall\, u \in \sD(\sA_p).
$$
\end{thm}

The proof of Theorem \ref{thm:Haller-Dintelmann_Heck_Hieber_3-1} is driven by the following \apriori estimate due to Haller-Dintelmann, Heck, and Hieber.

\begin{thm}[\Apriori estimate when $1<p<\infty$ for a homogeneous parameter-elliptic partial differential operator]
\label{thm:Haller-Dintelmann_Heck_Hieber_5-1}
\cite[Theorem 5.1 and Remark 3.3]{Haller-Dintelmann_Heck_Hieber_2006}, \cite[Theorem 6.6]{Heck_Hieber_2003}
Assume the hypotheses of Theorem \ref{thm:Haller-Dintelmann_Heck_Hieber_3-1} with $\theta_0 \in (0, \pi)$ and, in addition, that $\sA$ is \emph{homogeneous} of degree $m$. Then there are constants $C > 0$, $\lambda_1 \geq 0$, and $\eta > 0$ with the following significance.  If the coefficients $a_\alpha$ obey
\begin{equation}
\label{eq:Haller-Dintelmann_Heck_Hieber_Theorem_5-1_aij_VMO_condition}
\max_{|\alpha|=m}\|a_\alpha\|_* \leq \eta,
\end{equation}
then for all $u \in W^{m,p}_{\loc}(\RR^d; \CC^N)$ and $\lambda \in \Delta_\theta$ with $|\lambda| > \lambda_1$,
\begin{equation}
\label{eq:Haller-Dintelmann_Heck_Hieber_Theorem_5-1_apriori_estimate}
\sum_{|\alpha|\leq m} |\lambda|^{1-|\alpha|/m}\|D^\alpha u\|_{L^p(\RR^d)}
\leq
C\|(\sA + \lambda)u\|_{L^p(\RR^d)}.
\end{equation}
\end{thm}

In our application, we shall need to relax the requirements in Theorem \ref{thm:Haller-Dintelmann_Heck_Hieber_5-1} that
\begin{inparaenum}[\itshape a\upshape)]
\item $\sA$ be homogeneous of order $m$,
\item that the coefficients $a_\alpha$ obey \eqref{eq:Haller-Dintelmann_Heck_Hieber_Theorem_5-1_aij_VMO_condition}, and
\item replace the condition $\lambda \in \Delta_\theta$ with $|\lambda| > \lambda_1$ by $\lambda \in \Delta_\theta(\lambda_3)$
\end{inparaenum}
to give the generalizations in Theorem \ref{thm:Haller-Dintelmann_Heck_Hieber_5-1_full_apriori_estimate} of the \apriori estimate \eqref{eq:Krylov_LecturesSobolev_13-4-2_diagonal_principal_matrix_lower_order_coefficients} in Theorem \ref{thm:Krylov_LecturesSobolev_13.4.5_diagonal_principal_symbol}. One can use a localization procedure (see \cite[p. 735]{Haller-Dintelmann_Heck_Hieber_2006}) to obtain
\eqref{eq:Haller-Dintelmann_Heck_Hieber_Theorem_5-1_full_apriori_estimate} without the restriction \eqref{eq:Haller-Dintelmann_Heck_Hieber_Theorem_5-1_aij_VMO_condition} imposed in Theorem \ref{thm:Haller-Dintelmann_Heck_Hieber_5-1}; we provide a simplified version of that argument under the assumption that the coefficients of the principal symbol are uniformly continuous on $\RR^d$.

\begin{thm}[\Apriori estimate when $1<p<\infty$]
\label{thm:Haller-Dintelmann_Heck_Hieber_5-1_full_apriori_estimate}
Assume the hypotheses of Theorem \ref{thm:Haller-Dintelmann_Heck_Hieber_3-1} and, in addition, that the coefficients $a_\alpha$ of the principal symbol of $\sA$ in \eqref{eq:Haller-Dintelmann_Heck_Hieber_page_720_definition_differential_system} for $\alpha \in \NN^d$ with $|\alpha|=m$ belong to $C(\bar\RR^d; \CC^{N\times N})$, with modulus of continuity $\omega$. Then there are constants $C > 0$, $\lambda_2 \geq 0$, and $\lambda_3 \geq 0$ with the following significance.
\begin{enumerate}
\item
\label{item:Haller-Dintelmann_Heck_Hieber_5-1_full_apriori_estimate}
If $\lambda \in \Delta_\theta$ and $|\lambda| \geq \lambda_2$, then
\begin{equation}
\label{eq:Haller-Dintelmann_Heck_Hieber_Theorem_5-1_full_apriori_estimate}
\sum_{|\alpha|\leq m} |\lambda|^{1-|\alpha|/m}\|D^\alpha u\|_{L^p(\RR^d)}
\leq
C\|(\sA + \lambda)u\|_{L^p(\RR^d)}, \quad \forall\, u \in W^{m,p}(\RR^d; \CC^N).
\end{equation}
\item
\label{item:Haller-Dintelmann_Heck_Hieber_5-1_full_apriori_estimate_sector}
If $\lambda \in \Delta_\theta(\lambda_3)$, then
\begin{equation}
\label{eq:Haller-Dintelmann_Heck_Hieber_Theorem_5-1_full_apriori_estimate_sector}
\sum_{|\alpha|\leq m} |\lambda - \lambda_3|^{1-|\alpha|/m}\|D^\alpha u\|_{L^p(\RR^d)}
\leq
C\|(\sA + \lambda)u\|_{L^p(\RR^d)}, \quad \forall\, u \in W^{m,p}(\RR^d; \CC^N).
\end{equation}
\end{enumerate}
\end{thm}

\begin{proof}
We divide the proof into steps corresponding to the three issues described in the paragraph preceding the statement of
Theorem \ref{thm:Haller-Dintelmann_Heck_Hieber_5-1_full_apriori_estimate}.

\setcounter{step}{0}
\begin{step}[Relaxing the condition that $\sA$ be homogeneous of order $m$]
For $\sA$ as in \eqref{eq:Haller-Dintelmann_Heck_Hieber_page_720_definition_differential_system}, the operator $\sA - \sum_{|\alpha|< m} a_\alpha D^\alpha = \sum_{|\alpha| = m} a_\alpha D^\alpha$ is homogeneous of degree $m$ and so Theorem \ref{thm:Haller-Dintelmann_Heck_Hieber_5-1} yields, for all $\in \Delta_\theta$ with $|\lambda|>\lambda_1$,
\begin{align*}
\sum_{|\alpha|\leq m} |\lambda|^{1-|\alpha|/m}\|D^\alpha u\|_{L^p(\RR^d)}
&\leq
C\left\|\left(\sA + \lambda  - \sum_{|\alpha|< m} a_\alpha D^\alpha\right)u\right\|_{L^p(\RR^d)}
\\
&\leq C\|(\sA + \lambda)u\|_{L^p(\RR^d)} + \sum_{|\alpha|< m}C\|D^\alpha u\|_{L^p(\RR^d)}
\\
&\leq C\|(\sA + \lambda)u\|_{L^p(\RR^d)} + \frac{1}{2}\sum_{|\alpha|< m} |\lambda|^{1 - |\alpha|/m} \|D^\alpha u\|_{L^p(\RR^d)},
\end{align*}
provided $|\lambda|^{1 - k/m} \geq 2C$ for $k = 0, \ldots, m-1$. The \apriori estimate \eqref{eq:Haller-Dintelmann_Heck_Hieber_Theorem_5-1_full_apriori_estimate}, still assuming that the coefficients $a_\alpha$ obey \eqref{eq:Haller-Dintelmann_Heck_Hieber_Theorem_5-1_aij_VMO_condition} when $|\alpha|=m$, now follows by choosing $\lambda_2 \geq \max\{1, \lambda_1\}$ and $\lambda_2 \geq (2C)^{m/(m - k)}$ for  $k = 0, \ldots, m-1$, where $C$ and $\lambda_1$ are as in Theorem \ref{thm:Haller-Dintelmann_Heck_Hieber_5-1}.
\end{step}

\begin{step}[Relaxing the condition that the coefficients $a_\alpha$ obey \eqref{eq:Haller-Dintelmann_Heck_Hieber_Theorem_5-1_aij_VMO_condition} when $|\alpha|=m$]
By the definition \eqref{eq:Haller-Dintelmann_Heck_Hieber_page_721_definition_star_norm_and_BMO} of $\|a_\alpha\|_*$ and the inequality \eqref{eq:Krylov_eta_a_bounded_by_modulus_of_continuity}, we know that $\|a_\alpha\|_* = 0$ if $a_\alpha$ is constant on $\RR^d$, in which case the inequality \eqref{eq:Krylov_eta_a_bounded_by_modulus_of_continuity} is obviously obeyed. Thus, it suffices to use a version of the method of freezing the coefficients (see, for example, the proofs of \cite[Theorem 6.2]{GilbargTrudinger} or \cite[Lemma 1.6.2]{Krylov_LecturesSobolev}).

We fix a cut-off function $\zeta \in C^\infty_0(\RR^d)$ such that $\zeta = 1$ on the hypercube $Q_1$ and $\supp\zeta \subset Q_2$ and $0 \leq \zeta \leq 1$ on $\RR^d$, where we denote
$$
Q_r(x_0) := \{x \in \RR^d: |x-x_0| < r\}, \quad\hbox{for } x_0 \in \RR^d \hbox{ and } r > 0,
$$
and $Q_r(x_0) = Q_r$ if $x_0$ is the origin in $\RR^d$. Furthermore, we may suppose that for each integer $m \geq 1$, there are positive constants $K_m$ such that
$$
|D^\alpha\zeta| \leq K_m, \quad \forall\, \alpha \in \NN^d \hbox{ with } |\alpha| = m.
$$
Note that $Q_r(x_0) := \{x \in \RR^d: |x-x_0| < r\} \subset B_{r\sqrt{d}}(x_0)$, for any $x_0 \in \RR^d$ and $r > 0$. For a positive constant $r > 0$ to be determined further below, we now express $\RR^d$ as a union, over all $n \in \NN$, of closed hypercubes, $\bar Q_r(x^n)$, centered at grid points, $x^n \in r\ZZ^d \subset \RR^d$, such that the open hypercubes $Q_r(x^{n_1})$ and $Q_r(x^{n_2})$ are disjoint for all $n_1 \neq n_2$. We define a corresponding sequence of cut-off functions, $\{\zeta_n\}_{n \in \NN} \subset C^\infty_0(\RR^d)$, by setting
$$
\zeta_n(x) := \zeta((x - x^n)/r), \quad \forall\, x \in \RR^d.
$$
Thus, $\supp \zeta_n \subset Q_{2r}(x^n)$ and $\zeta_n = 1$ on $Q_r(x^n)$ and
\begin{equation}
\label{eq:Haller-Dintelmann_Heck_Hieber_5-1_cut-off_function_derivative_bounds}
|D^\alpha\zeta_n| \leq r^mK_m, \quad \forall\, \alpha \in \NN^d \hbox{ with } |\alpha| \leq m.
\end{equation}
We now define the partial differential operator $\sA_{(n)}$ with constant principal symbol coefficients $a_\alpha(x_n)$ for all $\alpha \in \NN^d$ when $|\alpha| = m$ and the same coefficients $a_\alpha$ as those of $\sA$ in \eqref{eq:Haller-Dintelmann_Heck_Hieber_page_720_definition_differential_system} when $|\alpha| < m$. The  operator $\sA_{(n)}$ is $(\sM, \pi-\theta_0)$ parameter-elliptic in the sense of Definition \ref{defn:Haller-Dintelmann_Heck_Hieber_page_720_parameter_elliptic}.

The \apriori estimate \eqref{eq:Haller-Dintelmann_Heck_Hieber_Theorem_5-1_full_apriori_estimate}, under the assumption that \eqref{eq:Haller-Dintelmann_Heck_Hieber_Theorem_5-1_aij_VMO_condition} holds, can thus be applied with the operator $\sA_{(n)}$, to give, for $\lambda \in \Delta_\theta$ and $|\lambda| \geq \lambda_2$ and $u \in W^{m,p}(\RR^d; \CC^N)$,
\begin{equation}
\label{eq:Haller-Dintelmann_Heck_Hieber_Theorem_5-1_full_apriori_estimate_constant_coefficients}
\sum_{|\alpha|\leq m} |\lambda|^{1-|\alpha|/m}\|D^\alpha (\zeta_n u)\|_{L^p(\RR^d)}
\leq
C\|(\sA_{(n)} + \lambda)(\zeta_n u)\|_{L^p(\RR^d)}.
\end{equation}
Noting that the coefficients $a_\alpha$ have modulus of continuity $\omega$ when $|\alpha|=m$,
\begin{align*}
\|(\sA_{(n)} + \lambda)(\zeta_n u)\|_{L^p(\RR^d)} &\leq \|(\sA + \lambda)(\zeta_n u)\|_{L^p(\RR^d)} + \|(\sA - \sA_{(n)})(\zeta_n u)\|_{L^p(\RR^d)}
\\
&\leq \|(\sA + \lambda)(\zeta_n u)\|_{L^p(\RR^d)} + \sum_{|\alpha|=m}\|(a_\alpha - a_\alpha(x^n))D^\alpha(\zeta_n u)\|_{L^p(\RR^d)}
\\
&\leq \|(\sA + \lambda)(\zeta_n u)\|_{L^p(\RR^d)} + \sum_{|\alpha|=m} \omega(r\sqrt{d})\|D^\alpha(\zeta_n u)\|_{L^p(\RR^d)},
\end{align*}
and so, choosing $r>0$ small enough that $C\omega(r\sqrt{d}) \leq 1/2$ and supposing without loss of generality that $\lambda_2 \geq 1$ we obtain, via rearrangement, that
\begin{equation}
\label{eq:Haller-Dintelmann_Heck_Hieber_Theorem_5-1_full_apriori_estimate_small_cube}
\sum_{|\alpha|\leq m} |\lambda|^{1-|\alpha|/m}\|D^\alpha (\zeta_n u)\|_{L^p(\RR^d)}
\leq
C\|(\sA + \lambda)(\zeta_n u)\|_{L^p(\RR^d)}.
\end{equation}
We now estimate the left-hand side of the inequality \eqref{eq:Haller-Dintelmann_Heck_Hieber_Theorem_5-1_full_apriori_estimate} with
\begin{align*}
\sum_{|\alpha|\leq m} |\lambda|^{1-|\alpha|/m}\|D^\alpha u\|_{L^p(\RR^d)}
&= \sum_{|\alpha|\leq m} |\lambda|^{1-|\alpha|/m} \left(\sum_{n\in\NN}\|D^\alpha u\|_{L^p(Q_r(x^n))}^p\right)^{1/p}
\\
&\leq \sum_{|\alpha|\leq m} |\lambda|^{1-|\alpha|/m} \left(\sum_{n\in\NN}\|D^\alpha (\zeta_n u)\|_{L^p(Q_r(x^n))}^p\right)^{1/p}
\\
&= \sum_{|\alpha|\leq m} \left(\sum_{n\in\NN}|\lambda|^{p(1-|\alpha|/m)} \|D^\alpha (\zeta_n u)\|_{L^p(Q_r(x^n))}^p\right)^{1/p}
\\
&\leq C\left(\sum_{n\in\NN} \|(\sA + \lambda)(\zeta_n u)\|_{L^p(\RR^d)}^p\right)^{1/p}
\quad\hbox{(by \eqref{eq:Haller-Dintelmann_Heck_Hieber_Theorem_5-1_full_apriori_estimate_small_cube})}
\\
&\leq C\left(\sum_{n\in\NN} \left(\|\zeta_n(\sA + \lambda)u\|_{L^p(\RR^d)} + \|[\sA, \zeta_n] u\|_{L^p(\RR^d)}\right)^p\right)^{1/p}
\\
&\leq C\left(\sum_{n\in\NN} \left(\|\zeta_n(\sA + \lambda)u\|_{L^p(\RR^d)}^p  + \|[\sA, \zeta_n] u\|_{L^p(\RR^d)}^p\right)\right)^{1/p}
\\
&\qquad\hbox{(using $(x + y)^p \leq 2^{p-1}(x^p + y^p)$ for all $x, y  \geq 0$ \cite[Lemma 2.2]{AdamsFournier})}
\\
&\leq C\left(\sum_{n\in\NN} \|(\sA + \lambda)u\|_{L^p(Q_{2r}(x^n))}^p
+ \sum_{|\alpha|\leq m-1} \sum_{n\in\NN} \|D^\alpha u\|_{L^p(Q_{2r}(x^n))}^p \right)^{1/p}
\\
&\qquad \hbox{(since $[\sA, \zeta_n]$ has order $m-1$ and applying \eqref{eq:Haller-Dintelmann_Heck_Hieber_5-1_cut-off_function_derivative_bounds})}
\\
&\leq C\left( \|(\sA + \lambda)u\|_{L^p(\RR^d)}^p + \sum_{|\alpha|\leq m-1} \|D^\alpha u\|_{L^p(\RR^d)}^p \right)^{1/p}.
\end{align*}
Therefore, noting that $(\sum_i y_i^p)^{1/p} \leq \sum_i y_i$ for all $y_i \geq 0$ and $p > 1$ by \cite[Inequality (1.4.1) or Theorem 19]{HardyLittlewoodPolya},
$$
\sum_{|\alpha|\leq m} |\lambda|^{1-|\alpha|/m}\|D^\alpha u\|_{L^p(\RR^d)}
\leq
C\|(\sA + \lambda)u\|_{L^p(\RR^d)} + C\sum_{|\alpha|\leq m-1} \|D^\alpha u\|_{L^p(\RR^d)}.
$$
Hence, for large enough $\lambda_2$ (depending on $C$), we can use rearrangement in the preceding inequality to give \eqref{eq:Haller-Dintelmann_Heck_Hieber_Theorem_5-1_full_apriori_estimate} \emph{without} the restriction  \eqref{eq:Haller-Dintelmann_Heck_Hieber_Theorem_5-1_aij_VMO_condition}, as desired.
\end{step}

\begin{step}[Replacing the condition $\lambda \in \Delta_\theta$ with $|\lambda| > \lambda_2$ by $\lambda \in \Delta_\theta(\lambda_3)$]
Suppose $\lambda_3$ obeys $\lambda_3 \geq \lambda_2$, where $\lambda_2$ is as in Item \eqref{item:Haller-Dintelmann_Heck_Hieber_5-1_full_apriori_estimate}, and consider $\lambda \in \CC\less\{\lambda_3\}$ such that $|\arg(\lambda - \lambda_3)| < \theta$, where $\theta$ is as in Theorem \ref{thm:Haller-Dintelmann_Heck_Hieber_3-1}. We write $\lambda - \lambda_3 = Re^{i\varphi}$ with $|\varphi| < \theta$ and examine the two cases,
\begin{inparaenum}[\itshape a\upshape)]
\item $0 \leq |\varphi| < \pi/2$, and
\item $\pi/2 \leq |\varphi| < \theta$.
\end{inparaenum}
We observe that
$$
|\lambda - \lambda_3| \leq |\lambda| + \lambda_3 < 2|\lambda|,
\quad\forall\, \lambda \in \CC \hbox{ with } |\lambda| > \lambda_3.
$$
If $|\varphi| < \pi/2$, then $\Real \lambda > \lambda_3$ and so
$$
|\lambda| > \lambda_3 \quad\hbox{and}\quad |\lambda-\lambda_3| < 2|\lambda|,
\quad \forall\, \lambda \in \Delta_\theta(\lambda_3)
\hbox{ with } |\arg(\lambda - \lambda_3)| \in [0, \pi/2),
$$
and for this case, it suffices to choose $\lambda_3 = \lambda_2$. The inequality \eqref{eq:Haller-Dintelmann_Heck_Hieber_Theorem_5-1_full_apriori_estimate} thus gives
\begin{multline}
\label{eq:Haller-Dintelmann_Heck_Hieber_Theorem_5-1_full_apriori_estimate_acute}
\sum_{|\alpha|\leq m} |\lambda - \lambda_3|^{1-|\alpha|/m}\|D^\alpha u\|_{L^p(\RR^d)}
\leq
C\|(\sA + \lambda)u\|_{L^p(\RR^d)}.
\\
\forall\, u \in W^{m,p}(\RR^d; \CC^N) \quad\hbox{and}\quad
\lambda \in \Delta_\theta(\lambda_3) \hbox{ with } |\arg(\lambda - \lambda_3)| \in [0, \pi/2).
\end{multline}
If $\pi/2 \leq |\varphi| < \theta < \pi$, then $\Real \lambda = \lambda_3 + |\lambda - \lambda_3|\cos \varphi$ with $0 \leq |\cos\varphi| < |\cos\theta| < 1$ and so
$$
\lambda_3 < |\Real \lambda| + |\lambda - \lambda_3||\cos \theta|,
$$
which gives
$$
|\lambda - \lambda_3| \leq |\lambda| + \lambda_3 < 2|\lambda| + |\lambda - \lambda_3||\cos \theta|,
$$
and therefore,
$$
|\lambda - \lambda_3| < \frac{2}{1-|\cos \theta|}|\lambda|,
\quad\forall\, \lambda \in \Delta_\theta(\lambda_3)
\hbox{ with } |\arg(\lambda - \lambda_3)| \in [\pi/2, \theta).
$$
Consequently, we also have
$$
\lambda_3 < |\lambda| + \frac{2|\cos \theta|}{1-|\cos \theta|}|\lambda|
= \frac{1+|\cos \theta|}{1-|\cos \theta|}|\lambda|,
\quad\forall\, \lambda \in \Delta_\theta(\lambda_3)
\hbox{ with } |\arg(\lambda - \lambda_3)| \in [\pi/2, \theta).
$$
and thus,
$$
|\lambda| > \frac{1-|\cos \theta|}{1+|\cos \theta|}\lambda_3,
\quad\forall\, \lambda \in \Delta_\theta(\lambda_3)
\hbox{ with } |\arg(\lambda - \lambda_3)| \in [\pi/2, \theta).
$$
Hence, if we require $\lambda_3 = \lambda_3(\lambda_2, \theta)$ to be large enough that
$$
\frac{1-|\cos \theta|}{1+|\cos \theta|}\lambda_3 \geq \lambda_2,
$$
then we obtain $|\lambda| > \lambda_2$ in this case too. Therefore, for a possibly larger constant $C$, the inequality
\eqref{eq:Haller-Dintelmann_Heck_Hieber_Theorem_5-1_full_apriori_estimate} gives
\begin{multline}
\label{eq:Haller-Dintelmann_Heck_Hieber_Theorem_5-1_full_apriori_estimate_oblique}
\sum_{|\alpha|\leq m} |\lambda - \lambda_3|^{1-|\alpha|/m}\|D^\alpha u\|_{L^p(\RR^d)}
\leq
C\|(\sA + \lambda)u\|_{L^p(\RR^d)}.
\\
\forall\, u \in W^{m,p}(\RR^d; \CC^N) \quad\hbox{and}\quad
\lambda \in \Delta_\theta(\lambda_3) \hbox{ with } |\arg(\lambda - \lambda_3)| \in [\pi/2, \theta),
\end{multline}
and the combination of \eqref{eq:Haller-Dintelmann_Heck_Hieber_Theorem_5-1_full_apriori_estimate_acute} and
\eqref{eq:Haller-Dintelmann_Heck_Hieber_Theorem_5-1_full_apriori_estimate_oblique} yields \eqref{eq:Haller-Dintelmann_Heck_Hieber_Theorem_5-1_full_apriori_estimate_sector}.
\end{step}
This completes the proof of Theorem \ref{thm:Haller-Dintelmann_Heck_Hieber_5-1_full_apriori_estimate}.
\end{proof}

The following result is essentially \cite[Corollary 3.2 and Remark 3.3]{Haller-Dintelmann_Heck_Hieber_2006}, but our hypotheses are slightly stronger and we include a proof which is omitted in \cite{Haller-Dintelmann_Heck_Hieber_2006}.

\begin{thm}[Sectorial property of a parameter-elliptic system and generation of an analytic semigroup on $L^p(\RR^d; \CC^N)$ when $1<p<\infty$]
\label{thm:Corollary_Haller-Dintelmann_Heck_Hieber_3-2}
Assume the hypotheses of Theorem \ref{thm:Haller-Dintelmann_Heck_Hieber_3-1} and, in addition, that $\theta_0 > \pi/2$, the coefficients $a_\alpha$ of the principal symbol of $\sA$ in \eqref{eq:Haller-Dintelmann_Heck_Hieber_page_720_definition_differential_system} for $\alpha \in \NN^d$ with $|\alpha|=m$ belong to $C(\bar\RR^d; \CC^{N\times N})$, and \emph{one} of the following holds:
\begin{enumerate}
\item $m = 2n$ for an integer $n\geq 1$, or
\item The coefficients $a_\alpha$ of $\sA$ in \eqref{eq:Haller-Dintelmann_Heck_Hieber_page_720_definition_differential_system} belong to $C^m_b(\RR^d;\CC^{N\times N})$ for all $\alpha \in \NN^d$ with $|\alpha|\leq m$.
\end{enumerate}
Then $\rho(-\sA_p) \supset \Delta_\theta(\lambda_0)$ and the resolvent estimate \eqref{eq:Haller-Dintelmann_Heck_Hieber_3-1} holds for $\lambda_0$, $\theta$, and $C$ as in Theorem \ref{thm:Haller-Dintelmann_Heck_Hieber_3-1}, $\sA_p$ is a sectorial operator on $L^p(\RR^d; \CC^N)$ in the sense of Definition \ref{defn:Sell_You_page_78_definition_of_sectorial_operator}, and $-\sA_p$ generates an analytic semigroup $e^{-\sA_p t}$ on $L^p(\RR^d; \CC^N)$.
\end{thm}

\begin{proof}
Let us first check that the resolvent estimate \eqref{eq:Haller-Dintelmann_Heck_Hieber_3-1} follows from Theorem \ref{thm:Haller-Dintelmann_Heck_Hieber_5-1_full_apriori_estimate}, since the details were omitted in \cite{Haller-Dintelmann_Heck_Hieber_2006}. We note that the \apriori estimate
\eqref{eq:Haller-Dintelmann_Heck_Hieber_Theorem_5-1_full_apriori_estimate_sector} in Theorem \ref{thm:Haller-Dintelmann_Heck_Hieber_5-1_full_apriori_estimate}, if we choose $\lambda_0 = \lambda_3$, gives
\begin{equation}
\label{eq:Haller-Dintelmann_Heck_Hieber_3-1apriori_estimate}
|\lambda - \lambda_0| \|u\|_{L^p(\RR^d)} \leq C\|(\sA + \lambda)u\|_{L^p(\RR^d)},
\quad \forall\, u \in \sD(\sA_p) \quad\hbox{and}\quad \lambda \in \Delta_\theta(\lambda_0).
\end{equation}
Clearly, $\sD(\sA_p) = W^{m,p}(\RR^d;\CC^N)$ is dense in $L^p(\RR^d;\CC^N)$, since $\sD(\sA_p)$ contains $C_0^\infty(\RR^d; \CC^N)$ and $C_0^\infty(\RR^d; \CC^N)$ is dense in $L^p(\RR^d;\CC^N)$. The inequality \eqref{eq:Haller-Dintelmann_Heck_Hieber_3-1apriori_estimate} implies that $\sA_p + \lambda$ is one-to-one on $\sD(\sA_p)$ for all $\lambda \in \Delta_\theta(\lambda_0)$.

The \apriori estimate \eqref{eq:Haller-Dintelmann_Heck_Hieber_Theorem_5-1_apriori_estimate} (thanks to Corollary
\ref{thm:Haller-Dintelmann_Heck_Hieber_5-1_full_apriori_estimate}) means that it is bounded below in the sense of \cite[Definition 2.1]{Abramovich_Aliprantis_2002} and so has closed range by \cite[Theorem 2.5]{Abramovich_Aliprantis_2002}.

To show that in fact $\Ran(\sA_p + \lambda) = L^p(\RR^d;\CC^N)$, we consider cases corresponding to our two alternative hypotheses.

\begin{case}[$m = 2n$ for an integer $n\geq 1$]
\label{case:Corollary_Haller-Dintelmann_Heck_Hieber_3-1_even_order}
We consider the one-parameter family of partial differential operators,
$$
\sA(t) := (1-t)(\Delta + 1)^{m/2} + t\sA, \quad t \in [0, 1],
$$
where $\Delta := -\sum_{i=1}^d\partial_{x_ix_i}$ is the Laplace operator, $\kappa$ is a positive constant, and we recall that $m=2n$, for an integer $n\geq 1$, by hypothesis. Then, for a possibly larger $\sM$ than in our hypotheses, the partial differential operator $\sA(t)$ is $(\sM, \theta)$ parameter-elliptic in the sense of
Definition \ref{defn:Haller-Dintelmann_Heck_Hieber_page_720_parameter_elliptic} for all $t \in [0, 1]$ and so the \apriori estimate
\eqref{eq:Haller-Dintelmann_Heck_Hieber_3-1apriori_estimate} holds with $\sA$ replaced by $\sA(t)$. But the operator
$(\Delta + 1)^{m/2}: W^{m,p}(\RR^d;\CC^N) \to L^p(\RR^d;\CC^N)$ is onto by Theorem \ref{thm:Krylov_LecturesSobolev_13.4.5_diagonal_principal_symbol} and so $\sA: W^{m,p}(\RR^d;\CC^N) \to L^p(\RR^d;\CC^N)$ is onto by the method of continuity \cite[Theorem 5.3]{GilbargTrudinger}. This concludes consideration of Case \ref{case:Corollary_Haller-Dintelmann_Heck_Hieber_3-1_even_order}.
\end{case}

\begin{case}[$a_\alpha \in C^m_b(\RR^d;\CC^{N\times N})$  for all $\alpha \in \NN^d$ with $|\alpha|\leq m$]
\label{case:Corollary_Haller-Dintelmann_Heck_Hieber_3-1_smooth_coefficients}
The range, $\Ran(\sA_p + \lambda)$, is closed in $L^p(\RR^d;\CC^N)$ since $\sA_p$ is closed and to show that $\sA_p + \lambda$ maps onto $L^p(\RR^d;\CC^N)$ for $\lambda \in \Delta_\theta(\lambda_0)$, it suffices to show that the range, $\Ran(\sA_p + \lambda)$, is dense in $L^p(\RR^d;\CC^N)$. For this purpose, we adapt the proof of \cite[Theorem 38.2]{Sell_You_2002} (in particular, see \cite[pp. 109--110]{Sell_You_2002}). Let $q \in (1,\infty)$ be the dual exponent to $p$, so $1/p+1/q=1$. Suppose there is a $v \in W^{m,q}(\RR^d;\CC^N)$ which belongs to the $L^2$-orthogonal complement of $\Ran(\sA_p + \lambda)$,
that is, such that
\begin{equation}
\label{eq:L2_inner_product_A+lambda_on_u_and_v_is_zero}
((\sA_p + \lambda)u, v)_{L^2(\RR^d)} = 0, \quad\forall\, u \in \sD(\sA_p).
\end{equation}
Because the coefficients $a_\alpha$ of $\sA$ belong to $C^m_b(\RR^d;\CC^{N\times N})$, we can integrate by parts in \eqref{eq:L2_inner_product_A+lambda_on_u_and_v_is_zero} to give
\begin{equation}
\label{eq:L2_inner_product_u_and_Aadjoint+lambda_on_v_is_zero}
(u, (\sA_q^* + \bar\lambda)v)_{L^2(\RR^d)} = 0, \quad\forall\, u \in \sD(\sA_p),
\end{equation}
where $\sA^*$ is the formal adjoint \cite[Equation (6.2)]{AgmonLecturesEllipticBVP} of $\sA$ defined by $(\sA u, w)_{L^2(\RR^d)} = (u, \sA^*w)_{L^2(\RR^d)}$ for all $u, w \in C^\infty_0(\RR^d; \CC^N)$, that is,
\begin{equation}
\label{eq:Haller-Dintelmann_Heck_Hieber_page_720_definition_differential_system_formal_adjoint}
\sA^*w = \sum_{|\alpha|\leq m} D^\alpha(\bar a_\alpha w), \quad \forall\, w \in C^\infty_0(\RR^d; \CC^N),
\end{equation}
where we recall from the definition \eqref{eq:Haller-Dintelmann_Heck_Hieber_page_720_definition_differential_system} of $\sA$ that $D = -i(\partial_{x_1}, \ldots, \partial_{x_d})$, and let $\sA_q^*$ denotes its realization with $\sD(\sA_q^*) = W^{m,q}(\RR^d; \CC^N)$. Since $\sD(\sA_p) = W^{m,p}(\RR^d;\CC^N)$ is dense in $L^p(\RR^d;\CC^N)$, the relation \eqref{eq:L2_inner_product_u_and_Aadjoint+lambda_on_v_is_zero} implies that $v$ belongs to the kernel of the bounded operator $\sA^*+\bar\lambda: W^{m,q}(\RR^d;\CC^N) \to L^q(\RR^d;\CC^N)$. But $\lambda \in \Delta_\theta(\lambda_0)$ if and only if $\bar\lambda \in \Delta_\theta(\lambda_0)$ and our proof that $\sA_p+\lambda: W^{m,p}(\RR^d;\CC^N) \to L^p(\RR^d;\CC^N)$ is one-to-one if $\lambda \in \Delta_\theta(\lambda_0)$ applies equally well to show that $\sA_q^*+\bar\lambda: W^{m,q}(\RR^d;\CC^N) \to L^q(\RR^d;\CC^N)$ is one-to-one if $\bar\lambda \in \Delta_\theta(\lambda_0)$ and therefore $v = 0$. Consequently, the identity \eqref{eq:L2_inner_product_A+lambda_on_u_and_v_is_zero} implies that $\Ran(\sA_p + \lambda)$ is dense in $L^p(\RR^d;\CC^N)$, as desired.  Thus, we again obtain that $\sA: W^{m,p}(\RR^d;\CC^N) \to L^p(\RR^d;\CC^N)$ is onto. This concludes consideration of Case \ref{case:Corollary_Haller-Dintelmann_Heck_Hieber_3-1_smooth_coefficients}.
\end{case}

The \apriori estimate \eqref{eq:Haller-Dintelmann_Heck_Hieber_3-1apriori_estimate} or the Open Mapping Theorem (see, for example, \cite[Section 2.5]{Yosida}) ensures that the inverse,
$$
(\sA_p + \lambda)^{-1}: L^p(\RR^d;\CC^N) \to L^p(\RR^d;\CC^N),
$$
exists as a bounded, linear operator for all $\lambda \in \Delta_\theta(\lambda_0)$, and thus
$$
\rho(-\sA_p) \supset \Delta_\theta(\lambda_0),
$$
and the resolvent estimate \eqref{eq:Haller-Dintelmann_Heck_Hieber_3-1} now follows from the \apriori estimate \eqref{eq:Haller-Dintelmann_Heck_Hieber_3-1apriori_estimate}.

Therefore, $\sA_p$ is a sectorial operator on $L^p(\RR^d;\CC^N)$ in the sense of Definition \ref{defn:Sell_You_page_78_definition_of_sectorial_operator}, with constants $a = -\lambda_0$, $M := \max\{1, C\}$, and $\delta = \theta$, where $\theta_0 > \pi/2$ and $\theta \in (\pi/2, \theta_0)$. Thus, $-\sA_p$ is the generator of an analytic semigroup, $e^{-\sA_p t}$, on $L^p(\RR^d;\CC^N)$ by Theorem \ref{thm:Renardy_Rogers_12-31}.
\end{proof}

Higher-order \apriori estimates and regularity may be obtained by adapting the proofs of the corresponding results in \cite[Chapters 6 and 9]{GilbargTrudinger} to give

\begin{thm}[Higher-order \apriori $L^p$ estimate for a parameter-elliptic system]
\label{thm:Corollary_Haller-Dintelmann_Heck_Hieber_5-1_apriori_estimate_higher_order}
Assume the hypotheses of Theorem \ref{thm:Haller-Dintelmann_Heck_Hieber_5-1_full_apriori_estimate} and, in addition, that $k\geq 0$ is an integer and that the coefficients $a_\alpha$ of $\sA$ in \eqref{eq:Haller-Dintelmann_Heck_Hieber_page_720_definition_differential_system} belong to $C^k(\bar\RR^d; \CC^{N\times N})$ for all $\alpha \in \NN^d$ when $|\alpha|=m$ and $W^{k,\infty}(\RR^d; \CC^{N\times N})$ when $|\alpha|<m$. Then there are constants $C_0>0$, $\lambda_4\geq 0$, and $\lambda_5\geq 0$ with the following significance.
\begin{enumerate}
\item If $\lambda \in \Delta_\theta$ with $|\lambda| \geq \lambda_4$ and $u \in W^{k+m,p}(\RR^d; \CC^N)$, then
\begin{equation}
\label{eq:Haller-Dintelmann_Heck_Hieber_Theorem_5-1_apriori_estimate_higher_order}
\sum_{\begin{subarray}{c}|\alpha|\leq m \\ |\beta|\leq k\end{subarray}} |\lambda|^{1-|\alpha|/m}\|D^{\alpha+\beta} u\|_{L^p(\RR^d)}
\leq
C_0\|(\sA + \lambda)u\|_{W^{k,p}(\RR^d)}.
\end{equation}
\item If $\lambda \in \Delta_\theta(\lambda_5)$ and $u \in W^{k+m,p}(\RR^d; \CC^N)$, then
\begin{equation}
\label{eq:Haller-Dintelmann_Heck_Hieber_Theorem_5-1_apriori_estimate_higher_order_sector}
\sum_{\begin{subarray}{c}|\alpha|\leq m \\ |\beta|\leq k\end{subarray}} |\lambda - \lambda_5|^{1-|\alpha|/m}\|D^{\alpha+\beta} u\|_{L^p(\RR^d)}
\leq
C_0\|(\sA + \lambda)u\|_{W^{k,p}(\RR^d)}.
\end{equation}
\end{enumerate}
\end{thm}

\begin{proof}
Notice that when $k = 0$ the \apriori estimate \eqref{eq:Haller-Dintelmann_Heck_Hieber_Theorem_5-1_apriori_estimate_higher_order} is given by the inequality \eqref{eq:Haller-Dintelmann_Heck_Hieber_Theorem_5-1_full_apriori_estimate} in Theorem \ref{thm:Haller-Dintelmann_Heck_Hieber_5-1_full_apriori_estimate} and so without loss of generality we may assume that $k\geq 1$. Suppose $\beta \in \NN^d$ with $|\beta| \leq k$ and consider $u \in W^{k+m,p}(\RR^d; \CC^N)$. The inequality \eqref{eq:Haller-Dintelmann_Heck_Hieber_Theorem_5-1_full_apriori_estimate} yields
$$
\sum_{|\alpha|\leq m} |\lambda|^{1-|\alpha|/m}\|D^\alpha D^\beta u\|_{L^p(\RR^d)}
\leq
C\|(\sA + \lambda)D^\beta u\|_{L^p(\RR^d)}, \quad\forall\, \lambda \in \Delta_\theta \hbox{ with } |\lambda| \geq \lambda_2.
$$
Observe that
$$
(\sA + \lambda)D^\beta u = D^\beta(\sA + \lambda)u + [\sA, D^\beta]u,
$$
where the commutator $[\sA, D^\beta]$ is a partial differential operator of order $m-1$ with coefficients in $L^\infty(\RR^d; \CC^{N\times N})$. Thus,
\begin{align*}
\|(\sA + \lambda)D^\beta u\|_{L^p(\RR^d)} &\leq \|D^\beta(\sA + \lambda)u\|_{L^p(\RR^d)} + \|[\sA, D^\beta]u\|_{L^p(\RR^d)}
\\
&\leq \|(\sA + \lambda)u\|_{W^{k,p}(\RR^d)} + C\sum_{|\gamma| \leq k-1}\|D^\gamma u\|_{L^p(\RR^d)}.
\end{align*}
Hence, for a large enough positive constant $\lambda_4 \geq \lambda_2$, the terms $\|D^\gamma u\|_{L^p(\RR^d)}$ on the right-hand side of the preceding inequality with $|\gamma| \leq k-1$ may be absorbed into the left-hand side of the inequality \eqref{eq:Haller-Dintelmann_Heck_Hieber_Theorem_5-1_apriori_estimate_higher_order}  and this
completes the proof that \eqref{eq:Haller-Dintelmann_Heck_Hieber_Theorem_5-1_apriori_estimate_higher_order} holds for $k \geq 1$.

The \apriori estimate \eqref{eq:Haller-Dintelmann_Heck_Hieber_Theorem_5-1_apriori_estimate_higher_order_sector} follows from \eqref{eq:Haller-Dintelmann_Heck_Hieber_Theorem_5-1_apriori_estimate_higher_order} just as in the proof of Theorem \ref{thm:Haller-Dintelmann_Heck_Hieber_5-1_full_apriori_estimate}.
\end{proof}

\begin{thm}[Higher-order regularity for a parameter-elliptic system]
\label{thm:Corollary_Haller-Dintelmann_Heck_Hieber_5-1_higher_order_regularity}
Assume the hypotheses of Theorem \ref{thm:Haller-Dintelmann_Heck_Hieber_5-1_full_apriori_estimate} and, in addition, that $k\geq 0$ is an integer and that the coefficients $a_\alpha$ of $\sA$ in \eqref{eq:Haller-Dintelmann_Heck_Hieber_page_720_definition_differential_system} belong to $C^k(\bar\RR^d; \CC^{N\times N})$ for all $\alpha \in \NN^d$ when $|\alpha|=m$ and $W^{k,\infty}(\RR^d; \CC^{N\times N})$ when $|\alpha|<m$. Then there are constants $\lambda_4\geq 0$ and $\lambda_5\geq 0$ with the following significance. If $\lambda \in \Delta_\theta$ with $|\lambda|\geq \lambda_4$ or $\lambda \in \Delta_\theta(\lambda_5)$ and $u \in W^{m,p}(\RR^d; \CC^N)$ obey
$$
(\sA + \lambda)u \in W^{k,p}(\RR^d; \CC^N),
$$
then $u \in W^{k+m,p}(\RR^d; \CC^N)$.
\end{thm}

\begin{proof}
To prove first that $u \in W^{m+1,p}(\RR^d; \CC^N)$, we can adapt the finite-difference quotient method of proof of \cite[Theorems 6.17, 6.19 and 9.19]{GilbargTrudinger}, which concern the higher-order regularity of solutions to an elliptic, linear, scalar, second-order partial differential equation on an open subset $\Omega \subset \RR^d$. The required technical results on finite-difference quotients are given by \cite[Lemmata 7.23 and 7.24]{GilbargTrudinger}. The role played by the \apriori interior Schauder estimate \cite[Corollary 6.3]{GilbargTrudinger} in the proof of \cite[Theorem 6.17]{GilbargTrudinger} is replaced by that of the \apriori global Sobolev estimates in Theorem \ref{thm:Haller-Dintelmann_Heck_Hieber_5-1_full_apriori_estimate}. The situation here is simpler than in the proof of \cite[Theorem 6.17]{GilbargTrudinger} due to the absence of a domain boundary and the proof of \cite[Theorem 6.17]{GilbargTrudinger} now carries over without significant further change to show that $u \in W^{m+1,p}(\RR^d; \CC^N)$. The conclusion $u \in W^{k+m,p}(\RR^d; \CC^N)$ when $k \geq 2$ follows by induction on $k$, just as in the proof of \cite[Theorem 6.17]{GilbargTrudinger}.
\end{proof}

\subsection{Existence, uniqueness, \apriori $L^\infty$ estimates, and analytic semigroups on $C_\infty(\RR^d;\CC^N)$ defined by elliptic partial differential systems of order $m\geq 1$}
\label{subsec:Heck_Hieber_Stavrakidis_2010}
In this subsection, we review the existence and uniqueness results and \apriori $L^\infty$ estimates for a parameter-elliptic partial differential system of order $m\geq 1$ on $\RR^d$, due to Heck, Hieber, and Stavrakidis \cite{Heck_Hieber_Stavrakidis_2010}, together with their consequences for resolvent estimates and analytic semigroup generation results on $C_\infty(\RR^d;\CC^N)$. The coefficients of the principal symbol of the partial differential system are allowed to belong to $L^\infty(\RR^d;\CC^{N\times N}) \cap \VMO(\RR^d;\CC^{N\times N})$. Closely related results are due to Denk and Dreher \cite{Denk_Dreher_2011}.

For the case $p=\infty$, one defines \cite[p. 301]{Heck_Hieber_Stavrakidis_2010}
\begin{equation}
\label{eq:Heck_Hieber_Stavrakidis_page_301_definition_domain_A_infinity}
\sD(\sA_\infty) := \left\{u \in \bigcap_{p\geq 1} W^{m,p}_{\loc}(\RR^d;\CC^N):
u, \, \sA u \in L^\infty(\RR^d;\CC^N) \right\}.
\end{equation}
Note that an asymptotic boundary condition, such as $|u(x)| \to \infty$ as $|x| \to \infty$, is not incorporated by Heck, Hieber, and Stavrakidis in their definition \eqref{eq:Heck_Hieber_Stavrakidis_page_301_definition_domain_A_infinity} of $\sD(\sA_\infty)$. We include Theorem \ref{thm:Heck_Hieber_Stavrakidis_Theorem_1-1} below for the sake of completeness, though it plays no role in the sequel.

\begin{thm}[Existence, uniqueness, and resolvent estimate for an elliptic partial differential system on $W^{m-1,\infty}(\RR^d; \CC^N)$]
\label{thm:Heck_Hieber_Stavrakidis_Theorem_1-1}
\cite[Theorem 1.1]{Heck_Hieber_Stavrakidis_2010}
Assume the hypotheses of Theorem \ref{thm:Haller-Dintelmann_Heck_Hieber_3-1}. Then there are constants, $C>0$ and $\lambda_0\geq 0$, with the following significance. If $f \in L^\infty(\RR^d; \CC^N)$ and $\lambda \in \Delta_\theta(\lambda_0)$, then there is a unique $u \in \sD(\sA_\infty)$ such that $(\sA + \lambda)u = f$ a.e. on $\RR^d$ and
\begin{equation}
\label{eq:Heck_Hieber_Stavrakidis_Theorem_1-1}
\|(\sA_\infty + \lambda)^{-1}\|_{\sL(L^\infty(\RR^d; \CC^N))} \leq \frac{C}{|\lambda|}.
\end{equation}
\end{thm}

The proof of Theorem \ref{thm:Heck_Hieber_Stavrakidis_Theorem_1-1} depends heavily on the following \apriori estimate due to Heck, Hieber, and Stavrakidis and which generalizes Theorem \ref{thm:Stewart_1974_1_and_2}, due to Stewart, for the case $N=1$ and with coefficients of the principal symbol belonging to $C_b(\RR^d)$.

\begin{thm}[\Apriori estimate for a parameter-elliptic partial differential operator]
\label{thm:Heck_Hieber_Stavrakidis_lemma_3-2}
\cite[Lemma 3.2]{Heck_Hieber_Stavrakidis_2010}
Assume the hypotheses of Theorem \ref{thm:Haller-Dintelmann_Heck_Hieber_3-1} and, in addition, that $q>d$. Then there are constants, $K_q \geq 0$ and $\omega_q \geq 0$, with the following significance.  If $u \in C^{m-1}(\RR^d; \CC^N)\cap W^{m,q}_{\loc}(\RR^d; \CC^N)$ and $\lambda \in \CC$ with $\Real \lambda \geq \min\{1, \omega_q\}$,
then\footnote{There appears to be a typographical error in \cite[Equation (9)]{Heck_Hieber_Stavrakidis_2010} and the first displayed inequality in \cite[Lemma 3.2]{Heck_Hieber_Stavrakidis_2010}, where the term $\lambda-\sA$ should be replaced by $\lambda+\sA$: compare with \cite[Proposition 2.1]{Heck_Hieber_Stavrakidis_2010} and the resolvent estimate in \cite[Theorem 1.1]{Heck_Hieber_Stavrakidis_2010} and its derivation via \cite[Equation (9)]{Heck_Hieber_Stavrakidis_2010} in \cite[p. 307]{Heck_Hieber_Stavrakidis_2010}.}
\begin{multline}
\label{eq:Heck_Hieber_Stavrakidis_9prelim}
\sum_{|\alpha| < m} |\lambda|^{1-|\alpha|/m}\|D^\alpha u\|_{L^\infty(\RR^d)}
+ \sum_{|\alpha| = m} |\lambda|^{d/m q} \sup_{x_0 \in \RR^d} \|D^\alpha u\|_{L^q(B(x_0, |\lambda|^{-1/m})}
\\
\leq K_q |\lambda|^{d/m q} \sup_{x_0 \in \RR^d} \|(\sA + \lambda)u\|_{L^q(B(x_0, |\lambda|^{-1/m})}.
\end{multline}
In particular, if $\sA u \in L^\infty(\RR^d; \CC^N)$, then
\begin{multline}
\label{eq:Heck_Hieber_Stavrakidis_9}
\sum_{|\alpha| < m} |\lambda|^{1-|\alpha|/m}\|D^\alpha u\|_{L^\infty(\RR^d)}
+ \sum_{|\alpha| = m} |\lambda|^{d/m q} \sup_{x_0 \in \RR^d} \|D^\alpha u\|_{L^q(B(x_0, |\lambda|^{-1/m})}
\\
\leq K_q \gamma_d^{1/q} \|(\sA + \lambda)u\|_{L^\infty(\RR^d)},
\end{multline}
where $\gamma_d$ is the volume of the unit ball in $\RR^d$.
\end{thm}

To simplify our proof of Theorem \ref{thm:Heck_Hieber_Stavrakidis_Theorem_1-1_corollary_sectorial_operator} --- which we include since this is omitted \cite{Heck_Hieber_Stavrakidis_2010} --- we restrict our statement to the case of partial differential operators $\sA$ in \eqref{eq:Haller-Dintelmann_Heck_Hieber_page_720_definition_differential_system} with smooth coefficients $a_\alpha$. In addition, we modify the definition \eqref{eq:Heck_Hieber_Stavrakidis_page_301_definition_domain_A_infinity} of $\sD(\sA_\infty)$ in order to include an asymptotic boundary condition as $|x| \to \infty$ for $x \in \RR^d$, by analogy with \eqref{eq:Pazy_7-3-5_system} in the case of a bounded open subset, $\Omega \subset \RR^d$,
\begin{equation}
\label{eq:Heck_Hieber_Stavrakidis_page_301_definition_domain_A_infinity_asymptotically_zero}
\tag{\ref*{eq:Heck_Hieber_Stavrakidis_page_301_definition_domain_A_infinity}$'$}
\sD(\sA_\infty) := \left\{u \in \bigcap_{p\geq 1} W^{m,p}_{\loc}(\RR^d;\CC^N):
u \in C_\infty(\RR^d;\CC^N) \hbox{ and } \sA u \in L^\infty(\RR^d;\CC^N) \right\},
\end{equation}
where $C_\infty(\RR^d;\CC^N) := \{u \in C(\RR^d;\CC^N): |u(x)| \to 0 \hbox{ as } |x| \to \infty\}$. By analogy with \eqref{eq:Pazy_7-3-23_system} and \cite[p. 144]{Stewart_1974} (Stewart addresses the scalar case $N=1$), we define
\begin{equation}
\label{eq:Pazy_7-3-23_system_Euclidean_space}
\sD(\sA_c) := \{u \in \sD(\sA_\infty): \sA(x, D)u \in C_\infty(\RR^d;\CC^N)\}.
\end{equation}
We then have the

\begin{thm}[Sectorial property of a parameter-elliptic system and generation of an analytic semigroup on $C_\infty(\RR^d; \CC^N)$]
\label{thm:Heck_Hieber_Stavrakidis_Theorem_1-1_corollary_sectorial_operator}
Let $d \geq 2$, $\sM > 0$, $N \geq 1$, and $\theta_0 \in (0, \pi)$. Let $\sA$ in \eqref{eq:Haller-Dintelmann_Heck_Hieber_page_720_definition_differential_system} be an $(\sM, \pi-\theta_0)$-parameter elliptic partial differential operator of order $m \geq 1$ whose coefficients $a_\alpha$ belong to $C^\infty(\bar\RR^d;\CC^{N\times N})$ for all $\alpha \in \NN^d$ with $|\alpha|\leq m$. If $\theta_0 > \pi/2$, then there are constants, $C>0$ and $\lambda_0\geq 0$ and $\vartheta \in (\pi/2, \theta_0)$, such that $\rho(-\sA_c) \supset \Delta_\vartheta(\lambda_0)$ and
\begin{equation}
\label{eq:Heck_Hieber_Stavrakidis_Theorem_1-1_sectorial_estimate}
\|(\sA_c + \lambda)^{-1}\|_{\sL(C_\infty(\RR^d;\CC^N))} \leq \frac{C}{|\lambda - \lambda_0|},
\quad \forall\, \lambda \in \Delta_\vartheta(\lambda_0).
\end{equation}
Moreover, $\sA_c$ is a sectorial operator on $C_\infty(\RR^d; \CC^N)$ in the sense of Definition \ref{defn:Sell_You_page_78_definition_of_sectorial_operator} and $-\sA_c$ generates an analytic semigroup, $e^{-\sA_c t}$, on $C_\infty(\RR^d; \CC^N)$.
\end{thm}

\begin{proof}
Let us first check that \eqref{eq:Heck_Hieber_Stavrakidis_Theorem_1-1_sectorial_estimate} follows from Theorem \ref{thm:Heck_Hieber_Stavrakidis_lemma_3-2}. Fix $q>d$, choose $\lambda_0 := \min\{1, \omega_q\}$ and $\theta = \pi/2$, and observe that Theorem \ref{thm:Heck_Hieber_Stavrakidis_lemma_3-2} implies
\begin{equation}
\label{eq:Heck_Hieber_Stavrakidis_Theorem_1-1_prelim}
|\lambda| \|u\|_{L^\infty(\RR^d)} \leq C\|(\sA_c + \lambda)u\|_{L^\infty(\RR^d)},
\quad \forall\, u \in \sD(\sA_c) \quad\hbox{and}\quad \lambda \in \CC \hbox{ with } \Real \lambda > \lambda_0.
\end{equation}
Observe that $\sD(A_c)$, defined by \eqref{eq:Heck_Hieber_Stavrakidis_page_301_definition_domain_A_infinity_asymptotically_zero} and \eqref{eq:Pazy_7-3-23_system_Euclidean_space}, is a dense subset of $C_\infty(\RR^d; \CC^N)$, since $\sD(\sA_c)$ contains $C_0^\infty(\RR^d; \CC^N)$ and $C_0^\infty(\RR^d; \CC^N)$ is dense in $C_\infty(\RR^d; \CC^N)$.

The inequality \eqref{eq:Heck_Hieber_Stavrakidis_Theorem_1-1_prelim} implies that $\sA_c + \lambda$ is one-to-one on $\sD(\sA_c)$ for all $\lambda \in \CC$ with $\Real \lambda > \lambda_0$.

Note that $C_\infty(\RR^d;\CC^N)$ is a Banach subspace of $L^\infty(\RR^d;\CC^N)$. If $\{u_n\}_{n=1}^\infty \subset \sD(\sA_c)$ is a sequence such that $\{\sA_c u_n\}_{n=1}^\infty \subset C_\infty(\RR^d;\CC^N)$ converges in $L^\infty(\RR^d;\CC^N)$ to $u \in  L^\infty(\RR^d;\CC^N)$ and $v \in C_\infty(\RR^d;\CC^N)$, respectively, as $n\to \infty$, then the \apriori estimate \eqref{eq:Heck_Hieber_Stavrakidis_9} also implies that,
$$
\|u_l - u_n\|_{W^{m-1,\infty}(\RR^d)} \to 0, \quad\hbox{as } \min\{l,n\} \to \infty,
$$
and, for any $q > d$,
$$
\sum_{|\alpha| = m} \sup_{x_0 \in \RR^d} \|D^\alpha (u_l-u_n)\|_{L^q(B(x_0, |\lambda|^{-1/m})} \to 0
\quad\hbox{as } \min\{l,n\} \to \infty.
$$
Consequently, the sequence $\{u_n\}_{n=1}^\infty$ converges to a limit in $\sD(\sA_\infty)$ by its definition in \eqref{eq:Heck_Hieber_Stavrakidis_page_301_definition_domain_A_infinity_asymptotically_zero}, that is, we must have $u \in \sD(\sA_\infty)$. Moreover, we have $\sA_\infty u_n \to \sA_\infty u$ in $L^q_{\loc}(\RR^d; \CC^N)$ as $n\to \infty$, and so $v = \sA_\infty u$ a.e. on $\RR^d$ and as $v \in C_\infty(\RR^d;\CC^N)$, then we must have $u \in \sD(\sA_c)$ and $v = \sA_c u$ on $\RR^d$. Thus, $\sA_c$ has closed range. We now make the

\begin{claim}
\label{claim:Corollary_Heck_Hieber_Stavrakidis_Theorem_1-1_density}
If $\lambda \in \CC$ with $\Real \lambda > \lambda_0$, then $\Ran(\sA_c + \lambda)$ is dense in $C_\infty(\RR^d; \CC^N)$
\end{claim}

\begin{proof}[Proof of Claim \ref{claim:Corollary_Heck_Hieber_Stavrakidis_Theorem_1-1_density}]
Noting that $C_0^\infty(\RR^d; \CC^N)$ is dense in $C_\infty(\RR^d; \CC^N)$, choose $f \in C_0^\infty(\RR^d; \CC^N)$. Theorem \ref{thm:Corollary_Haller-Dintelmann_Heck_Hieber_3-2} implies that there exists a (unique) $u \in W^{m,q}(\RR^d; \CC^N)$ such that $(\sA + \lambda) u = f$ a.e. on $\RR^d$. Theorem \ref{thm:Corollary_Haller-Dintelmann_Heck_Hieber_5-1_higher_order_regularity} implies that $u$ belongs to $\cap_{k\geq 0}W^{k+m,q}(\RR^d; \CC^N)$ and hence $u \in C_\infty(\RR^d; \CC^N)$ by the Sobolev Embedding Theorem \cite[Theorem 4.12]{AdamsFournier}, noting that $W^{m,q}(\RR^d; \CC^N) = W^{m,q}_0(\RR^d; \CC^N)$ by \cite[Corollary 3.23]{AdamsFournier} and recalling that $W^{m,q}_0(\RR^d; \CC^N)$ is the completion of $C^\infty_0(\RR^d; \CC^N)$ by \cite[Section 3.2]{AdamsFournier}. In particular, $u$ belongs to $\sD(A_\infty)$ in \eqref{eq:Heck_Hieber_Stavrakidis_page_301_definition_domain_A_infinity_asymptotically_zero} and furthermore, since $f \in C_\infty(\RR^d; \CC^N)$ (because, \afortiori, $f \in C_0^\infty(\RR^d; \CC^N)$), we obtain that $u \in \sD(\sA_c)$ by its definition in \eqref{eq:Pazy_7-3-23_system_Euclidean_space}. Hence, $f \in \Ran(\sA_c + \lambda)$ and the conclusion follows.
\end{proof}

Therefore, because $\Ran(\sA_c + \lambda)$ is closed and dense in $C_\infty(\RR^d; \CC^N)$, we see that $\Ran(\sA_c + \lambda) = C_\infty(\RR^d; \CC^N)$. Furthermore, $\sD(\sA_c + \lambda)$ is contained in $C_\infty(\RR^d; \CC^N)$ (as a dense subset) and so the \apriori estimates in Theorem \ref{thm:Heck_Hieber_Stavrakidis_lemma_3-2} or the Open Mapping Theorem \cite[Section 2.5]{Yosida} ensures that the inverse,
$$
(\sA_c + \lambda)^{-1}: C_\infty(\RR^d; \CC^N) \to C_\infty(\RR^d; \CC^N),
$$
exists as a bounded, linear operator for all $\lambda \in \CC$ with $\Real \lambda > \lambda_0$ and so
$$
\rho(-\sA_c) \supset \{\lambda \in \CC: \Real \lambda > \lambda_0\}.
$$
Moreover, if $\lambda \in \CC$ with $\Real \lambda > \lambda_0$, then
$$
|\lambda - \lambda_0| \leq |\lambda| + \lambda_0 < 2|\lambda|,
$$
and \eqref{eq:Heck_Hieber_Stavrakidis_Theorem_1-1_prelim} thus yields, after replacing $2C$ by $C$,
\begin{multline}
\label{eq:Heck_Hieber_Stavrakidis_Theorem_1-1_halfplane}
\|(\sA_c + \lambda)^{-1}\|_{L^\infty(\RR^d; \CC^N)} \leq \frac{C}{|\lambda - \lambda_0|},
\\
\quad \forall\, \lambda \in \CC\less\{\lambda_0\} \quad\hbox{with}\quad 0 \leq |\arg(\lambda-\lambda_0)| < \pi/2.
\end{multline}
It remains to show that $\rho(-\sA_c)$ actually contains a sector, $\Delta_\vartheta(\lambda_0)$, for some $\vartheta \in (\pi/2, \theta_0)$ and that \eqref{eq:Heck_Hieber_Stavrakidis_Theorem_1-1_halfplane} continues to hold for $\lambda \in \Delta_\vartheta(\lambda_0)$.

Consider $\lambda \in \CC\less\{\lambda_0\}$ such that $|\arg(\lambda - \lambda_0)| < \theta_0$. We write $\lambda - \lambda_0 = Re^{i\varphi}$ with $|\varphi| < \theta_0$ and note that it suffices to consider the remaining case, $\pi/2 \leq |\varphi| < \theta_0$. Choose any $\phi \in (-\pi/2, \pi/2)$, write
$$
\sA_c + \lambda = \sA_c + (\lambda_0 + Re^{i\phi}) + R(e^{i\varphi} - e^{i\phi}),
$$
observe that \eqref{eq:Heck_Hieber_Stavrakidis_Theorem_1-1_sectorial_estimate} (for the already established case of $z = \lambda_0 + Re^{i\phi} \in \CC\less\{\lambda_0\}$ with $|\arg(z - \lambda_0)| < \pi/2$) yields
$$
\|(\sA_c + \lambda_0 + Re^{i\phi})^{-1}\|_{\sL(C_\infty(\RR^d;\CC^N))} \leq \frac{C}{R},
$$
and require (just as in the proof of \cite[Theorem 12.31]{Renardy_Rogers_2004}) that $\varphi$ obeys
$$
R(e^{i\varphi} - e^{i\phi}) < \frac{R}{C} \leq 1/\|(\sA_c + \lambda_0 + Re^{i\phi})^{-1}\|_{\sL(C_\infty(\RR^d;\CC^N))},
$$
that is, $\varphi$ obeys\footnote{There is a minor typographical error in the proof of \cite[Theorem 12.31]{Renardy_Rogers_2004}.}
$|e^{i\varphi} - e^{i\phi}| < 1/C$. Therefore,
\begin{align*}
\|(\sA_c + \lambda_0 + Re^{i\phi}) - (\sA_c + \lambda)\|_{\sL(C_\infty(\RR^d;\CC^N))}
&= |R(e^{i\varphi} - e^{i\phi}) |
\\
&< 1/\|(\sA_c + \lambda_0 + Re^{i\phi})^{-1}\|_{\sL(C_\infty(\RR^d;\CC^N))},
\end{align*}
and so the Neumann series \cite[Theorem 6.3]{Abramovich_Aliprantis_2002} for $(\sA_c + \lambda)^{-1}$ converges.
In particular, we may choose $\vartheta \in (\pi/2, \theta_0)$ such that $|e^{i\vartheta} - e^{i\phi}| < 1/C$ when $\phi \in (0, \pi/2)$ and $|e^{-i\vartheta} - e^{i\phi}| < 1/C$ when $\phi \in (-\pi/2, 0)$, provided $|\phi|$ is close enough to $\pi/2$. We thus find that
$$
\rho(-\sA_c) \supset \{\lambda \in \CC\less\{\lambda_0\}: \pi/2 \leq |\arg(\lambda-\lambda_0)| < \vartheta\},
$$
and
\begin{multline}
\label{eq:Heck_Hieber_Stavrakidis_Theorem_1-1_oblique}
\|(\sA_c + \lambda)^{-1}\|_{\sL(C_\infty(\RR^d;\CC^N))} \leq \frac{C}{|\lambda - \lambda_0|},
\\
\quad \forall\, \lambda \in \CC\less\{\lambda_0\} \quad\hbox{with}\quad \pi/2 \leq |\arg(\lambda-\lambda_0)| < \vartheta.
\end{multline}
By combining \eqref{eq:Heck_Hieber_Stavrakidis_Theorem_1-1_halfplane} and \eqref{eq:Heck_Hieber_Stavrakidis_Theorem_1-1_oblique}, we see that $\rho(-\sA_c) \supset \Delta_\vartheta(\lambda_0)$ and the resolvent estimate \eqref{eq:Heck_Hieber_Stavrakidis_Theorem_1-1_sectorial_estimate} holds.

Observe that $\sD(\sA_c)$ is dense in $C_\infty(\RR^d;\CC^N)$ by the Sobolev Embedding Theorem \cite[Theorem 4.12]{AdamsFournier}. Therefore, $\sA_c$ is a sectorial operator on $C_\infty(\RR^d;\CC^N)$ in the sense of Definition \ref{defn:Sell_You_page_78_definition_of_sectorial_operator}, with constants $a = -\lambda_0$, $M := \max\{1, C\}$, and $\delta = \vartheta$, and $-\sA_c$ is the generator of an analytic semigroup, $e^{-\sA_c t}$, on $C_\infty(\RR^d;\CC^N)$ by Theorem \ref{thm:Renardy_Rogers_12-31}.
\end{proof}

\begin{rmk}[Sectorial property of a parameter-elliptic system and generation of an analytic semigroup on $L^1(\RR^d; \CC^N)$]
Theorems \ref{thm:Heck_Hieber_Stavrakidis_lemma_3-2} and \ref{thm:Heck_Hieber_Stavrakidis_Theorem_1-1_corollary_sectorial_operator} do not immediately yield the analogous results on $L^1(\RR^d; \CC^N)$ by duality, by analogy with our proof of Theorem \ref{thm:Cannarsa_Terreni_Vespri_6-7_p_is_one} (where $m = 2$). This is because
our proof of Theorem \ref{thm:Cannarsa_Terreni_Vespri_6-7_p_is_one} does not extend from bounded $\Omega\Subset\Omega$ to the case $\Omega = \RR^d$ since we used the fact that $W^{m,2}(\Omega;\CC^N) \subset W^{m,p}(\Omega;\CC^N)$, for $1 \leq p < 2$, when $\Omega$ is bounded. However, such a consequence should follows using weighted Sobolev spaces and a strategy in that spirit is used by Cannarsa and Vespri \cite{Cannarsa_Vespri_1988} in their proofs of existence, uniqueness, \apriori estimates, sectoriality, and generation of analytic semigroups on $L^1(\RR^d)$.
\end{rmk}

\subsection{Existence, uniqueness, \apriori $L^p$ estimates, and analytic semigroups defined by elliptic partial differential operators of order $m \geq 1$ on the Banach space $L^p(X;E)$ of sections of a complex vector bundle over a closed manifold}
\label{subsec:Haller-Dintelmann_Heck_Hieber_vector_bundle_manifold}
We now adapt the results of Section \ref{subsec:Haller-Dintelmann_Heck_Hieber_2006} to the case of a elliptic partial differential operator of order $m \geq 1$,
$$
\sA: W^{m, p}(X; E) \to L^p(X; E),
$$
where $1<p<\infty$ and $E$ is a $C^\infty$ complex vector bundle of rank $N \geq 1$ over a closed, smooth manifold, $X$, of dimension $d \geq 2$, and $W^{s, p}(X; E)$, for $s \in \RR$, is the Sobolev space of $W^{s,p}$ sections of $E$.

We shall adopt the following definition of partial differential operator on a vector bundle and parameter-ellipticity, by analogy with the usual definition of a pseudo-differential operator on a vector bundle over a manifold and ellipticity \cite[Section 1.3]{Gilkey2}, \cite[Definitions 18.1.20 and 18.1.32]{Hormander_v3}.

\begin{defn}[Partial differential operator on a vector bundle over a manifold and parameter-ellipticity]
\label{defn:Elliptic_partial_differential_operator_vector_bundle_manifold}
Let $E$ and $F$ be $C^\infty$ complex vector bundles of rank $N \geq 1$ over a smooth manifold, $X$, of dimension $d \geq 2$. A continuous linear map, $\sA: C^\infty(X; E) \to C^\infty(X; F)$, is a \emph{partial differential operator of order $m\geq 1$ (with $C^\infty$ coefficients)} if for every coordinate chart, $\varphi: X \supset U \cong \varphi(U) \subset \RR^d$, and pair of vector bundle trivializations, $\varsigma: E\restriction U \cong U\times\CC^N$ and $\tau: F\restriction U \cong U\times\CC^N$, the induced operator,
\begin{equation}
\label{eq:partial_differential_operator_vector_bundle_manifold_localization}
\varphi_*\tau\sA\varsigma^{-1}\varphi^* : C^\infty(\varphi(U); \CC^N) \to C^\infty(\varphi(U); \CC^N),
\end{equation}
is a partial differential operator of order $m\geq 1$ (with $C^\infty$ coefficients) on $\varphi(U) \subset \RR^d$ in the sense of \eqref{eq:Haller-Dintelmann_Heck_Hieber_page_720_definition_differential_system}, with $\RR^d$ replaced by the open subset, $\varphi(U)$. We say that $\sA$ is $(\sM,\theta)$ \emph{parameter-elliptic} if the local differential operator \eqref{eq:partial_differential_operator_vector_bundle_manifold_localization} is $(\sM,\theta)$ parameter-elliptic in the sense of Definition \ref{defn:Haller-Dintelmann_Heck_Hieber_page_720_parameter_elliptic}.
\end{defn}

\begin{rmk}[Other concepts of elliptic operator on a vector bundle over a manifold]
\label{rmk:Elliptic_differential_operator_vector_bundle_manifold_alternative_definitions}
By replacing the role of the local Definition \ref{defn:Haller-Dintelmann_Heck_Hieber_page_720_parameter_elliptic} in Definition \ref{defn:Elliptic_partial_differential_operator_vector_bundle_manifold}, one could define $\sA$, \emph{mutatis mutandis}, to be a
\begin{inparaenum}[\itshape a\upshape)]
\item \emph{elliptic pseudo-differential operator} using the local Definition \ref{defn:Krylov_LecturesSobolev_12-2-1_and_7},
\item \emph{uniformly strongly elliptic scalar partial differential operator} using the local Definition \ref{defn:Krylov_LecturesSobolev_13_4_0_uniformly_strongly_elliptic_partialdo},
\item \emph{second-order, strictly elliptic partial differential operator} using the local Definition \ref{defn:Cannarsa_Terreni_Vespri_2-7} when $m=2$.
\end{inparaenum}
\end{rmk}

In view of our later applications, we shall assume that $X$ is endowed with a $C^\infty$ Riemannian metric, that $E$ is a complex (respectively, real) Hermitian (respectively, Riemannian) vector bundle endowed with Hermitian (respectively, orthogonal) $C^\infty$ connection (see, for example, \cite{Kobayashi, Kobayashi_Nomizu_v1, Kobayashi_Nomizu_v2}), and that $F=E$ in Definition \ref{defn:Elliptic_partial_differential_operator_vector_bundle_manifold}. We denote the associated covariant derivative by $\nabla: C^\infty(X; E) \to C^\infty(X; E\otimes T^*X)$ and let $\nabla^*:C^\infty(X; E\otimes T^*X) \to C^\infty(X; E)$ denote the $L^2$-adjoint of $\nabla$ defined by the metrics on $X$ and $E$ via
\begin{equation}
\label{eq:Covariant_derivative_adjoint}
(\nabla_\xi^* u, v)_{L^2(X)} = (u, \nabla_\xi v)_{L^2(X)}, \quad\forall\, \xi \in C^\infty(TX), \ u, v \in C^\infty(X;E).
\end{equation}
The principal symbol of the \emph{connection Laplacian} \cite[Section 2.8]{LM},
\begin{equation}
\label{eq:Connection_Laplacian}
\nabla^*\nabla: C^\infty(X;E) \to C^\infty(X;E)
\end{equation}
is expressed in terms of the Riemannian metric on $X$, namely $g \in S^2(T^*X)$ (the space of symmetric two-tensors on $X$) by \cite[Equation (2.8.4)]{LM}
\begin{equation}
\label{eq:Connection_Laplacian_principal_symbol}
\hbox{Symbol}(\nabla^*\nabla) = -\id_E\otimes g \in C^\infty(X;\End_\CC(E)\otimes S^2(T^*X)).
\end{equation}
Thus, with respect to a local coordinate chart on $X$ and a local trivialization of $E$, the principal coefficients in \eqref{eq:Cannarsa_Terreni_Vespri_1-1} are given by
$$
a_{ij} = g^{ij}\, \id_{\CC^N}, \quad 1 \leq i,j \leq d,
$$
where $\id_{\CC^N} \in \CC^{N\times N}$ is the identity matrix, and so $\nabla^*\nabla$ is elliptic in the sense of Cannarsa, Terreni, and Vespri in Definition \ref{defn:Cannarsa_Terreni_Vespri_2-7} and thus $(\sM, \theta)$ parameter-elliptic in the sense of Haller-Dintelmann, Heck, and Hieber in Definition \ref{defn:Haller-Dintelmann_Heck_Hieber_page_720_parameter_elliptic} with $\sM = 1/\kappa$ and $\theta = 0$.

Hence, we may define the Sobolev spaces, $W^{s, p}(X; E)$ and their norms, either via
$$
W^{s, p}(X; E) := (\nabla^*\nabla + 1)^{-s/2}L^p(X; E), \quad s \in \RR, \quad 1 < p < \infty,
$$
or by transferring the usual definition\footnote{Henceforth, we drop the distinction between $H^{s, p}(\RR^d; \CC^N)$ and $W^{s, p}(\RR^d; \CC^N)$.}
\eqref{eq:Krylov_LecturesSobolev_definition_13-3-1} (for $1 < p < \infty$) of $W^{s, p}(\RR^d; \CC^N)$ to $W^{s, p}(X; E)$ via the patching constructions in \cite[Section 1.3]{Gilkey2}, \cite[Section 18.1]{Hormander_v3}, including the case $1 \leq p \leq \infty$ when $s = k \in \NN$.

The goal of this section is then to prove the following analogue of Theorems \ref{thm:Cannarsa_Terreni_Vespri_6-7} and \ref{thm:Cannarsa_Terreni_Vespri_6-7_one_lessthan_p_lessthan_two} (for $1<p<\infty$) and Theorem \ref{thm:Corollary_Haller-Dintelmann_Heck_Hieber_3-2}.

\begin{thm}[Resolvent estimate and generation of an analytic semigroup on $L^p(X; E)$ when $1<p<\infty$]
\label{thm:Haller-Dintelmann_Heck_Hieber_Theorem_3-1_vector_bundle_manifold}
Let $\sM > 0$, and $p \in (1, \infty)$, and $\theta_0 \in (0, \pi)$, and $\theta \in [0, \theta_0)$. Let $E$ be a $C^\infty$ Hermitian vector bundle of complex rank $N\geq 1$ with $C^\infty$ Hermitian connection over a closed, Riemannian, $C^\infty$ manifold, $X$, of dimension $d \geq 2$, and let $\sA:C^\infty(X; E) \to C^\infty(X; E)$ be an $(\sM, \pi-\theta_0)$ parameter-elliptic partial differential operator of order $m \geq 1$ with $C^\infty$ coefficients, as in Definition \ref{defn:Elliptic_partial_differential_operator_vector_bundle_manifold}. Then there are constants, $C>0$ and $\lambda_0\geq 0$, such that
\begin{equation}
\label{eq:Haller-Dintelmann_Heck_Hieber_Theorem_3-1_vector_bundle_manifold}
\|(\sA_p + \lambda)^{-1}\|_{\sL(L^p(X; E))} \leq \frac{C}{|\lambda - \lambda_0|},
\quad \forall\, \lambda \in \Delta_\theta(\lambda_0),
\end{equation}
where the realization $\sA_p: \sD(\sA_p) \subset L^p(X; E) \to L^p(X; E)$ is defined by
$$
\sD(\sA_p) := W^{m,p}(X; E).
$$
If $\theta_0 > \pi/2$, then $-\sA_p$ generates an analytic semigroup, $e^{-\sA_pt}$, on $L^p(X; E)$.
\end{thm}

To prove Theorem \ref{thm:Haller-Dintelmann_Heck_Hieber_Theorem_3-1_vector_bundle_manifold}, we shall adapt the proof of \cite[Theorem 8.5.3]{Krylov_LecturesSobolev}, which employs a standard, universal patching argument (credited by Krylov in \cite[p. 172]{Krylov_LecturesSobolev} to F. Browder) valid for elliptic partial differential equations or systems, on domains in Euclidean space, manifolds, or vector bundles over manifolds. Given an elliptic scalar second-order differential operator, $\sA$, and real $\lambda \geq \lambda$, for a large enough constant $\lambda_0\geq 0$, Krylov uses this patching argument to establish unique solvability of $(\sA + \lambda)u = f$ for $u\in W^{2,p}(\Omega)\cap W^{1,p}_0(\Omega)$ on a domain $\Omega \subset \RR^d$ with $C^2$ boundary $\partial\Omega$, given unique solvability and \apriori estimates when $\Omega = \RR^d$ or $\Omega = \RR^{d-1}\times\RR_+$ (see \cite[Theorem 8.5.3]{Krylov_LecturesSobolev}), together with general \apriori estimates for $u$ (see \cite[Theorem 8.3.7, Lemma 8.5.2, and Theorem 8.5.6]{Krylov_LecturesSobolev}). Krylov applies a similar patching argument to establish unique solvability of $(\sA + \lambda)u = f$ for $u\in C^{2,\alpha}_0(\bar\Omega)$ (see \cite[Theorem 6.5.3]{Krylov_LecturesHolder}).

Note that our sign convention for $\sA$ is opposite to that of Krylov for his elliptic scalar second-order differential operator $L$ in \cite[p. 157]{Krylov_LecturesHolder}, so in adapting his equations and estimates in \cite[Chapter 8]{Krylov_LecturesSobolev}, we replace $L$ by $-\sA$.

Before proceeding to the proofs of a series of auxiliary results and ultimately Theorem \ref{thm:Haller-Dintelmann_Heck_Hieber_Theorem_3-1_vector_bundle_manifold} itself, we begin with the preliminary setup we shall need. We fix a covering of $X$ by coordinate neighborhoods, $U_i\subset X$ for $i=1,\ldots, n$, local coordinate charts on $X$, and local trivializations of $E$,
\begin{equation}
\label{eq:Krylov_Sobolev_lectures_page_169_definition_local_coordinate_chart_and_trivialization}
\begin{gathered}
\psi_i: X \supset U_i \cong \psi_i(B) \subset \RR^d,
\\
\tau_i : E\restriction U_i \cong U_i \times \CC^N, \quad \hbox{for } i=1,\ldots, n.
\end{gathered}
\end{equation}
To apply the local results from Section \ref{subsec:Haller-Dintelmann_Heck_Hieber_2006}, we must choose a $C^\infty$ partition of unity, $\{\zeta_i\}_{i=1}^n$, subordinate to the covering, $\{U_i\}_{i=1}^n$, and $C^\infty$ functions, $\{\eta_i\}_{i=1}^n$, such that (compare \cite[pp. 166, 168, and 171--172]{Krylov_LecturesSobolev}),
\begin{subequations}
\label{eq:Krylov_LecturesSobolev_partition_of_unity}
\begin{gather}
\label{eq:Krylov_LecturesSobolev_zeta_and_eta_properties}
0 \leq \zeta_i, \eta_i \leq 1 \hbox{ on } X, \quad \supp \zeta_i, \eta_i \Subset U_i,
\quad \eta_i = 1 \hbox{ on } \supp \zeta_i
\quad\hbox{for } 1 \leq i \leq n,
\\
\label{eq:Krylov_LecturesSobolev_sum_zetai_is_one}
\hbox{and}\quad \sum_{i=1}^n \zeta_i^2 = 1.
\end{gather}
\end{subequations}
Thus, given $u \in W^{k,p}(X; E)$, for $0\leq k\leq m$, then $\psi_{i,*} (\eta_i u) := \eta_i u \circ \psi_i^{-1}  \in W^{k,p}(\RR^d; E)$ and $\tau_i\psi_{i,*} (\eta_i u) \in W^{k,p}(\RR^d; \CC^N)$. Conversely, given $w \in W^{k,p}(\RR^d; \CC^N)$, we have $\eta_i\psi_i^*w := \eta_i w\circ\psi_i \in W^{k,p}(X; \CC^N)$ and $\tau_i^{-1}\eta_i\psi_i^*w \in W^{k,p}(X; E)$.

By analogy with \cite[p. 169 and p. 172]{Krylov_LecturesSobolev}, we define
\begin{align}
\label{eq:Local_elliptic_system_on_Euclidean space_from_global_operator_vector_bundle}
\sA_i w &:= \psi_{i,*}\tau_i\sA \tau_i^{-1}\psi_i^* w,
\quad \forall\, w \in W^{m,p}(\psi_i(U_i); \CC^N),
\\
\label{eq:Krylov_Sobolev_lectures_page_172_definition_global_resolvent_on_Euclidean space}
\sR(\lambda, -\sA_i) &:= (\sA_i + \lambda)^{-1}: L^p(\RR^d; \CC^N) \to W^{m,p}(\RR^d; \CC^N),
\quad\forall\, \lambda \in \Delta_\theta(\lambda_0),
\\
\label{eq:Krylov_Sobolev_lectures_page_169_definition_local_resolvent_on_patch}
\sR_i(\lambda, -\sA)v &:= \tau_i^{-1}\psi_i^*\sR(\lambda, -\sA_i)\psi_{i,*}\tau_i(\eta_i v),
\quad \forall\, v \in W^{m,p}(X; E), \quad \hbox{for } i=1,\ldots, n,
\end{align}
where $\theta \in (\theta_0, \pi)$, $\theta_0 \in (\pi/2, \pi)$, and $\lambda_0 \geq 0$ are the constants in Theorem \ref{thm:Corollary_Haller-Dintelmann_Heck_Hieber_3-2}, for an elliptic partial differential operator $\sA_i$ on $C^\infty(\RR^d; \CC^N)$ of order $m$ as in \eqref{eq:Haller-Dintelmann_Heck_Hieber_page_720_definition_differential_system}, for $i=1,\ldots, n$.

In our definition \eqref{eq:Krylov_Sobolev_lectures_page_172_definition_global_resolvent_on_Euclidean space}, \emph{we assume throughout this subsection} that the $\psi_i$ and $\tau_i$ in \eqref{eq:Krylov_Sobolev_lectures_page_169_definition_local_coordinate_chart_and_trivialization} are defined on slightly larger coordinate neighborhoods, $V_i \Supset U_i$, and that the coefficients of $\sA_i$ are smoothly extended from $\psi_i(U_i)$ to $\RR^d$, by patching over $\psi(V_i \less U_i)$, so that the following holds.

\begin{lem}
\label{lem:Each_sAi_obeys_corollary_Haller-Dintelmann_Heck_Hieber_5-1}
For $i=1,\ldots,n$, the extended differential operator, $\sA_i$ on $C^\infty(\RR^d; \CC^N)$, obtained from \eqref{eq:Local_elliptic_system_on_Euclidean space_from_global_operator_vector_bundle} obeys the hypotheses of Theorem \ref{thm:Corollary_Haller-Dintelmann_Heck_Hieber_3-2}, with constants $C\geq 0$, and $\theta_0 \in (\pi/2, \pi)$, and $\theta \in (\theta_0, \pi)$, and $\lambda_0 \geq 0$.
\end{lem}

We now develop a series of technical results leading to the proof of Theorem \ref{thm:Haller-Dintelmann_Heck_Hieber_Theorem_3-1_vector_bundle_manifold}, beginning with the following analogue of \cite[Lemma 8.4.1]{Krylov_LecturesSobolev}.

\begin{lem}
\label{lem:Krylov_Sobolev_lectures_8-4-1}
Assume the hypotheses of Theorem \ref{thm:Haller-Dintelmann_Heck_Hieber_Theorem_3-1_vector_bundle_manifold}. Let $f \in L^p(X; E)$ and $\lambda \in \Delta_\theta(\lambda_0)$. If $u \in  W^{m,p}(X; E)$ is a solution to
\begin{equation}
\label{eq:Krylov_Sobolev_lectures_8-5-1}
(\sA + \lambda)u = f \quad\hbox{a.e. on } X,
\end{equation}
then $u$ is also a solution to
\begin{equation}
\label{eq:Krylov_Sobolev_lectures_8-4-2}
u = \sum_i \zeta_i \sR_i(\lambda, -\sA)(\zeta_i f + [\sA, \zeta_i]u),
\end{equation}
where $[\sA, \zeta_i]u \equiv \sA (\zeta_i u) - \zeta_i \sA u$.
\end{lem}

\begin{proof}
Calculating, we find that
\begin{align*}
{}&\sum_i \zeta_i \sR_i(\lambda, -\sA)(\zeta_i f + [\sA, \zeta_i]u)
\\
&\quad = \sum_i \zeta_i \sR_i(\lambda, -\sA)(\zeta_i \sA u + [\sA, \zeta_i]u + \lambda \zeta_i u)
\\
&\quad = \sum_i \zeta_i \sR_i(\lambda, -\sA)(\sA \zeta_i u + \lambda \zeta_i u)
\\
&\quad = \sum_i \zeta_i \psi_i^*\tau_i^{-1}\sR(\lambda, -\sA_i)\tau_i\psi_{i,*}(\eta_i(\sA u + \lambda)(\zeta_i u))
\quad\hbox{(by \eqref{eq:Krylov_Sobolev_lectures_page_169_definition_local_resolvent_on_patch})}
\\
&\quad = \sum_i \zeta_i \psi_i^*\tau_i^{-1} (\sA_i + \lambda)^{-1} \tau_i\psi_{i,*}(\sA u + \lambda)(\zeta_i u)
\quad\hbox{(by \eqref{eq:Krylov_LecturesSobolev_zeta_and_eta_properties} and \eqref{eq:Krylov_Sobolev_lectures_page_172_definition_global_resolvent_on_Euclidean space})}
\\
&\quad = \sum_i \zeta_i \psi_i^*\tau_i^{-1} (\sA_i + \lambda)^{-1} (\sA_i + \lambda)\tau_i\psi_{i,*}(\zeta_i u)
\quad\hbox{(by \eqref{eq:Local_elliptic_system_on_Euclidean space_from_global_operator_vector_bundle})}
\\
&\quad = \sum_i \zeta_i \psi_i^*\tau_i^{-1} \tau_i\psi_{i,*}(\zeta_i u)
= \sum_i \zeta_i^2 u = u \quad\hbox{(by \eqref{eq:Krylov_LecturesSobolev_sum_zetai_is_one}),}
\end{align*}
and thus $u$ solves \eqref{eq:Krylov_Sobolev_lectures_8-4-2}. This completes the proof of Lemma \ref{lem:Krylov_Sobolev_lectures_8-4-1}.
\end{proof}

We have the following analogue of \cite[Definition 8.4.2]{Krylov_LecturesSobolev}.

\begin{defn}[Regularizer of $\sA + \lambda$ on $W^{m,p}(X; E)$]
\label{defn:Krylov_8-4-2}
Assume the setup in Theorem \ref{thm:Haller-Dintelmann_Heck_Hieber_Theorem_3-1_vector_bundle_manifold}. The \emph{regularizer} of the operator $\sA + \lambda: W^{m,p}(X; E) \to L^p(X; E)$ is given by
\begin{equation}
\label{eq:Krylov_Sobolev_lectures_8-4-2_definition_regularizer}
\bR(\lambda, -\sA)f := \sum_i \zeta_i \sR_i(\lambda, -\sA)\zeta_i f, \quad \forall\, f \in L^p(X; E).
\end{equation}
\end{defn}

Clearly, the regularizer $\bR(\lambda, -\sA)$ of $\sA + \lambda$ is obtained by omitting the commutator term on the right-hand side of \eqref{eq:Krylov_Sobolev_lectures_8-4-2} and should be viewed as an approximation to the inverse, $\sR(\lambda, -\sA) = (\sA + \lambda)^{-1}: L^p(X; E) \to W^{m,p}(X; E)$, whose existence need to establish, for $\lambda$ in a suitable sector in $\CC$. By contrast with Lemma \ref{lem:Krylov_Sobolev_lectures_8-4-1}, the proof of the following converse and analogue of
\cite[Lemma 8.5.1]{Krylov_LecturesSobolev} is more involved.

\begin{lem}
\label{lem:Krylov_Sobolev_lectures_8-5-1}
Assume the hypotheses of Theorem \ref{thm:Haller-Dintelmann_Heck_Hieber_Theorem_3-1_vector_bundle_manifold}. Then there is a constant $\lambda_1 \geq \lambda_0$ with the following significance. If $f \in L^p(X; E)$ and $\lambda \in \Delta_\theta(\lambda_0)$ and $u \in  W^{m-1,p}(X; E)$ is a solution to \eqref{eq:Krylov_Sobolev_lectures_8-4-2}, then
\begin{enumerate}
  \item $u \in  W^{m,p}(X; E)$, and
  \item If $|\lambda| \geq \lambda_1$, then $u$ is also a solution to \eqref{eq:Krylov_Sobolev_lectures_8-5-1}.
\end{enumerate}
\end{lem}

\begin{proof}
By hypothesis we have $f \in L^p(X; E)$ and because the commutators, $[\sA, \zeta_i]$, are differential operators of order $m-1$, we also have $[\sA, \zeta_i]u \in L^p(X; E)$ for $i=1,\ldots, n$. By the definition of $\sR_i(\lambda, -\sA)$ through \eqref{eq:Krylov_Sobolev_lectures_page_172_definition_global_resolvent_on_Euclidean space} and \eqref{eq:Krylov_Sobolev_lectures_page_169_definition_local_resolvent_on_patch}, we must have (compare \cite[Exercise 8.3.9]{Krylov_LecturesSobolev})
$$
\sum_i \zeta_i \sR_i(\lambda, -\sA)(\zeta_i f + [\sA, \zeta_i]u) \in W^{m,p}(X; E),
$$
and because $u$ solves \eqref{eq:Krylov_Sobolev_lectures_8-4-2}, we see that $u \in W^{m,p}(X; E)$.

Next, denote $g := (\sA + \lambda)u \in L^p(X; E)$. Lemma \ref{lem:Krylov_Sobolev_lectures_8-4-1} implies that the equality \eqref{eq:Krylov_Sobolev_lectures_8-4-2} holds with $g$ in place of $f$, and to finish the proof, we need only show that if $|\lambda| \geq \lambda_1$, for $\lambda_1$ sufficiently large, then $g = f$. This is equivalent to showing that if $h \in L^p(X; E)$ and $\bR(\lambda, -\sA)h = 0$ a.e. on $X$, then $h = 0$ a.e. on $X$. Indeed, because
$$
\sum_i \zeta_i \sR_i(\lambda, -\sA)(\zeta_i f + [\sA, \zeta_i]u) = u = \sum_i \zeta_i \sR_i(\lambda, -\sA)(\zeta_i g + [\sA, \zeta_i]u),
$$
the Definition \ref{defn:Krylov_8-4-2} of $\bR(\lambda, -\sA)$ implies that $\bR(\lambda, -\sA)f = \bR(\lambda, -\sA)g$. Observe that (compare \cite[Exercise 8.3.9]{Krylov_LecturesSobolev})
\begin{align*}
(\sA + \lambda)\sR_i(\lambda, -\sA)(\zeta_i h)
&=
(\sA + \lambda)\tau_i^{-1}\psi_i^*\sR(\lambda, -\sA_i)\psi_{i,*}\tau_i(\eta_i \zeta_i h)
\quad\hbox{(by \eqref{eq:Krylov_Sobolev_lectures_page_169_definition_local_resolvent_on_patch})}
\\
&=
\tau_i^{-1}\psi_i^*(\sA_i + \lambda)\sR(\lambda, -\sA_i)\psi_{i,*}\tau_i(\eta_i \zeta_i h)
\quad\hbox{(by \eqref{eq:Local_elliptic_system_on_Euclidean space_from_global_operator_vector_bundle})}
\\
&= \eta_i \zeta_i h \quad\hbox{(by Definition \ref{defn:Resolvent_set} of a resolvent)}
\end{align*}
and so, because $\eta_i=1$ on $\supp\zeta_i$ by \eqref{eq:Krylov_LecturesSobolev_zeta_and_eta_properties},
\begin{equation}
\label{eq:sA_plus_lambda_applied_to_sRi}
(\sA + \lambda)\sR_i(\lambda, -\sA)(\zeta_i h) = \zeta_i h, \quad \forall\, h \in L^p(X; E).
\end{equation}
Therefore, since $\bR(\lambda, -\sA)h = 0$ a.e. on $X$, we discover that
\begin{align*}
0 = (\sA + \lambda) \bR(\lambda, -\sA)h
&=
\sum_{i=1}^n (\sA + \lambda) \zeta_i \sR_i(\lambda, -\sA)(\zeta_i h)
\quad\hbox{(by \eqref{eq:Krylov_Sobolev_lectures_8-4-2_definition_regularizer})}
\\
&= \sum_{i=1}^n \zeta_i (\sA + \lambda) \sR_i(\lambda, -\sA)(\zeta_i h)
+ \sum_{i=1}^n [\sA, \zeta_i] \sR_i(\lambda, -\sA)(\zeta_i h)
\\
&= \sum_{i=1}^n \zeta_i^2 h + \sum_{i=1}^n [\sA, \zeta_i] \sR_i(\lambda, -\sA)(\zeta_i h)
\quad\hbox{(by \eqref{eq:sA_plus_lambda_applied_to_sRi})}
\\
&= h - \sT_\lambda h,
\end{align*}
recalling that $\sum_{i=1}^n \zeta_i^2 = 1$ on $X$ by \eqref{eq:Krylov_LecturesSobolev_zeta_and_eta_properties} and where we have defined
$$
\sT_\lambda h := -\sum_{i=1}^n [\sA, \zeta_i] \sR_i(\lambda, -\sA)(\zeta_i h).
$$
To finish the proof it suffices to show that for all $\lambda \in \Delta_\theta$ with $|\lambda| > \lambda_1$ and $\lambda_1$ sufficiently large, $\sT_\lambda$ is a contraction operator on $L^p(X; E)$. We calculate that
\begin{align*}
\|\sT_\lambda h\|_{L^p(X)}
&\leq
\sum_{i=1}^n \| [\sA, \zeta_i] \sR_i(\lambda, -\sA)(\zeta_i h)\|_{L^p(X)}
\\
&\leq C\|\sR_i(\lambda, -\sA)(\zeta_i h)\|_{W^{m-1,p}(X)}
\quad \hbox{(because $[\sA, \zeta_i]$ has order $m-1$)}
\\
&=
\sum_{i=1}^n C\|\tau_i^{-1}\psi_i^*\sR(\lambda, -\sA_i)\psi_{i,*}\tau_i(\eta_i \zeta_i h)\|_{W^{m-1,p}(X)}
\quad\hbox{(by \eqref{eq:Krylov_Sobolev_lectures_page_169_definition_local_resolvent_on_patch})}
\\
&\leq
\sum_{i=1}^n\sum_{|\alpha|\leq m-1}
C\|D^\alpha\sR(\lambda, -\sA_i)\psi_{i,*}\tau_i(\zeta_i h)\|_{L^p(\RR^d)}
\quad\hbox{(by \eqref{eq:Krylov_LecturesSobolev_zeta_and_eta_properties})}
\\
&\leq
\sum_{i=1}^n C\left( \sum_{k=0}^{m-1} |\lambda|^{k/m - 1} \right) \|\psi_{i,*}\tau_i(\zeta_i h)\|_{L^p(\RR^d)}
\quad\hbox{(by Theorem \ref{thm:Haller-Dintelmann_Heck_Hieber_5-1_full_apriori_estimate} and Lemma \ref{lem:Each_sAi_obeys_corollary_Haller-Dintelmann_Heck_Hieber_5-1})}
\\
&\leq C\left( \sum_{k=0}^{m-1} |\lambda|^{k/m - 1} \right) \|h\|_{L^p(X)}.
\end{align*}
Since the constant $C$ is independent of $\lambda$, we may choose $\lambda_1 > 0$ large enough that
$$
C\sum_{k=0}^{m-1} \lambda_1^{k/m - 1} \leq \frac{1}{2},
$$
and thus, for any $\lambda \in \Delta_\theta$ obeying $|\lambda| > \lambda_1$, the same inequality holds with $\lambda_1$ replaced by $|\lambda|$. In particular, this yields $\|\sT_\lambda h\|_{L^p(X)} \leq (1/2)\|h\|_{L^p(X)}$ for $|\lambda| > \lambda_1$ and completes the proof of Lemma \ref{lem:Krylov_Sobolev_lectures_8-5-1}.
\end{proof}

We proceed to establish the following analogue of \cite[Lemma 8.5.2]{Krylov_LecturesSobolev}.

\begin{lem}
\label{lem:Krylov_Sobolev_lectures_8-5-2}
Assume the hypotheses of Theorem \ref{thm:Haller-Dintelmann_Heck_Hieber_Theorem_3-1_vector_bundle_manifold}. Then there are constants, $C > 0$ and $\lambda_2 \geq 0$, with the following significance. If $f \in L^p(X; E)$ and $\lambda \in \Delta_\theta(\lambda_0)$ with $|\lambda| \geq \lambda_2$, then there is a unique solution $u \in  W^{m-1,p}(X; E)$ to \eqref{eq:Krylov_Sobolev_lectures_8-4-2} and
\begin{equation}
\label{eq:Krylov_Sobolev_lectures_lemma_8-5-2_estimate}
\sum_{k=0}^{m-1} |\lambda|^{1 - k/m} \|\nabla^k u\|_{L^p(X)} \leq C\|f\|_{L^p(X)}.
\end{equation}
\end{lem}

\begin{proof}
We use the contraction mapping principle in $W^{m-1,p}(X; E)$. Accordingly, to prove Lemma \ref{lem:Krylov_Sobolev_lectures_8-5-2}, it suffices to show that, for sufficiently large $\lambda_2$, and all $f \in L^p(X; E)$ and $v \in W^{m-1,p}(X; E)$ and $\lambda \in \Delta_\theta(\lambda_2)$ and (compare \eqref{eq:Krylov_Sobolev_lectures_8-4-2})
\begin{equation}
\label{eq:Krylov_Sobolev_lectures_8-4-2_rhs_v}
u := \sum_i \zeta_i \sR_i(\lambda, -\sA)(\zeta_i f + [\sA, \zeta_i]v),
\end{equation}
we have $u \in W^{m-1,p}(X; E)$ and
\begin{equation}
\label{eq:Krylov_Sobolev_lectures_8-5-2}
\sum_{k=0}^{m-1} |\lambda|^{1 - k/m} \|\nabla^k u\|_{L^p(X)} \leq C\left(\|f\|_{L^p(X)} + \|v\|_{W^{m-1,p}(X)}\right).
\end{equation}
We know from the proof of Lemma \ref{lem:Krylov_Sobolev_lectures_8-5-1} that $u \in W^{m-1,p}(X; E)$. Next, notice that
\begin{align*}
{}&\sum_{k=0}^{m-1} |\lambda|^{1 - k/m} \|\nabla^k u\|_{L^p(X)}
\\
&\leq \sum_{k=0}^{m-1} |\lambda|^{1 - k/m}
\sum_i \left( \|\nabla^k \zeta_i \sR_i(\lambda, -\sA) (\zeta_i f)\|_{L^p(X)}
+ \|\nabla^k \zeta_i \sR_i(\lambda, -\sA) [\sA, \zeta_i]v\|_{L^p(X)} \right)
\\
&\qquad\hbox{(by \eqref{eq:Krylov_Sobolev_lectures_8-4-2_rhs_v})}
\\
&\leq \sum_{|\alpha|\leq m-1} |\lambda|^{1 - |\alpha|/m}
\sum_i \left( \|D^\alpha\sR(\lambda, -\sA_i) \psi_{i,*}\tau_i (\zeta_i f)\|_{L^p(\RR^d)}
+ \|D^\alpha\sR(\lambda, -\sA_i)\psi_{i,*}\tau_i [\sA, \zeta_i]v\|_{L^p(\RR^d)} \right)
\\
&\qquad\hbox{(by \eqref{eq:Krylov_Sobolev_lectures_page_169_definition_local_resolvent_on_patch} and \eqref{eq:Krylov_LecturesSobolev_zeta_and_eta_properties} and assuming $|\lambda|\geq 1$)}
\\
&\leq
\sum_i C \left( \|\psi_{i,*}\tau_i (\zeta_i f)\|_{L^p(\RR^d)} + \|\psi_{i,*}\tau_i [\sA, \zeta_i]v\|_{L^p(\RR^d)} \right)
\\
&\quad\hbox{(by Theorem \ref{thm:Haller-Dintelmann_Heck_Hieber_5-1_full_apriori_estimate} with $|\lambda| > \lambda_1$ and Lemma \ref{lem:Each_sAi_obeys_corollary_Haller-Dintelmann_Heck_Hieber_5-1})}
\\
&\leq C\left(\|f\|_{L^p(X)} + \|v\|_{W^{m-1,p}(X)}\right)
\quad\hbox{(since $[\sA, \zeta_i]$ has order $m-1$)}.
\end{align*}
Therefore, \eqref{eq:Krylov_Sobolev_lectures_lemma_8-5-2_estimate} holds. Given $v \in W^{m-1,p}(X; E)$, define $u := \sT_\lambda v \in W^{m-1,p}(X; E)$ by the right-hand side of \eqref{eq:Krylov_Sobolev_lectures_8-4-2_rhs_v} and observe that, for $v_1, v_2 \in W^{m-1,p}(X; E)$, the estimate \eqref{eq:Krylov_Sobolev_lectures_lemma_8-5-2_estimate} yields
$$
\sum_{k=0}^{m-1} |\lambda|^{1 - k/m} \|\nabla^k (\sT_\lambda v_1 - \sT_\lambda v_2)\|_{L^p(X)}
\leq
C\|v_1 - v_2\|_{W^{m-1,p}(X)},
$$
and thus, choosing $\lambda_2 \geq \max\{1, \lambda_1\}$ large enough that
$$
\sum_{k=0}^{m-1} \lambda_2^{1 - k/m} \leq \frac{1}{2},
$$
we obtain, for $\lambda \in \Delta_\theta(\lambda_0)$ with $|\lambda| > \lambda_2$,
$$
\|\sT_\lambda v_1 - \sT_\lambda v_2\|_{W^{m-1,p}(X)} \leq \|v_1 - v_2\|_{W^{m-1,p}(X)},
$$
and this completes the proof of Lemma \ref{lem:Krylov_Sobolev_lectures_8-5-2}.
\end{proof}

Lemmata \ref{lem:Krylov_Sobolev_lectures_8-5-1} and \ref{lem:Krylov_Sobolev_lectures_8-5-2} thus give the following analogue of \cite[Theorem 8.5.3]{Krylov_LecturesSobolev}.

\begin{thm}
\label{thm:Krylov_Sobolev_lectures_8-5-3}
Assume the hypotheses of Theorem \ref{thm:Haller-Dintelmann_Heck_Hieber_Theorem_3-1_vector_bundle_manifold}. Then there are constants, $C > 0$ and $\lambda_3 \geq 0$, with the following significance. If $f \in L^p(X; E)$ and $\lambda \in \Delta_\theta(\lambda_0)$ with $|\lambda| \geq \lambda_3$, then there is a unique solution, $u \in  W^{m,p}(X; E)$, to $(\sA + \lambda)u = f$ a.e. on $X$ and
\begin{equation}
\label{eq:Krylov_Sobolev_lectures_theorems_8-5-3_and_6_estimate}
\sum_{k=0}^m |\lambda|^{1 - k/m} \|\nabla^k u\|_{L^p(X)} \leq C\|f\|_{L^p(X)}.
\end{equation}
\end{thm}

\begin{proof}
Choose $\lambda_3 := \max\{\lambda_1, \lambda_2\}$, where $\lambda_1$ and $\lambda_2$ are the constants in Lemmata \ref{lem:Krylov_Sobolev_lectures_8-5-1} and \ref{lem:Krylov_Sobolev_lectures_8-5-2}, respectively. Then Lemma \ref{lem:Krylov_Sobolev_lectures_8-5-2} provides a unique solution, $u \in  W^{m-1,p}(X; E)$, to \eqref{eq:Krylov_Sobolev_lectures_8-4-2} and Lemma \ref{lem:Krylov_Sobolev_lectures_8-5-1} ensures that this $u$ actually belongs to $W^{m,p}(X; E)$ and solves \eqref{eq:Krylov_Sobolev_lectures_8-5-1}.

The bounds on the terms $|\lambda|^{1 - k/m} \|\nabla^k u\|_{L^p(X)}$ in \eqref{eq:Krylov_Sobolev_lectures_theorems_8-5-3_and_6_estimate}, for $0\leq k \leq m-1$, are given by Lemma \ref{lem:Krylov_Sobolev_lectures_8-5-2}. To estimate the term $\|\nabla^m u\|_{L^p(X)}$, we again proceed as in the proof of Lemma \ref{lem:Krylov_Sobolev_lectures_8-5-2} to give
\begin{align*}
\|\nabla^m u\|_{L^p(X)}
&\leq
\sum_i \left( \|\nabla^m \zeta_i \sR_i(\lambda, -\sA) (\zeta_i f)\|_{L^p(X)}
+ \|\nabla^m \zeta_i \sR_i(\lambda, -\sA) [\sA, \zeta_i]u\|_{L^p(X)} \right)
\\
&\leq
\sum_i \left( \|\sR(\lambda, -\sA_i)\psi_{i,*}\tau_i (\zeta_i f)\|_{W^{m,p}(\RR^d)}
+ \|\sR(\lambda, -\sA_i)\psi_{i,*}\tau_i [\sA, \zeta_i]u\|_{W^{m,p}(\RR^d)} \right)
\\
&\leq
\sum_i C\left( \|\psi_{i,*}\tau_i (\zeta_i f)\|_{L^p(\RR^d)} + \|\psi_{i,*}\tau_i [\sA, \zeta_i]u\|_{L^p(\RR^d)} \right)
\\
&\qquad\hbox{(by Theorem \ref{thm:Haller-Dintelmann_Heck_Hieber_5-1_full_apriori_estimate}
with $|\lambda| > \max\{1, \lambda_1\}$ and Lemma \ref{lem:Each_sAi_obeys_corollary_Haller-Dintelmann_Heck_Hieber_5-1})}
\\
&\leq C\left(\|f\|_{L^p(X)} + \|u\|_{W^{m-1,p}(X)}\right)
\quad\hbox{(since $[\sA, \zeta_i]$ has order $m-1$)}
\\
&\leq C\|f\|_{L^p(X)} + \frac{1}{2}\sum_{k=0}^{m-1} |\lambda|^{1 - k/m} \|\nabla^k u\|_{L^p(X)},
\end{align*}
where the last inequality follows if $|\lambda| \geq \lambda_3$ and we choose $\lambda_3$ large enough that $2C \leq \lambda_3^{1 - k/m}$ for $0 \leq k \leq m-1$. Rearrangement yields \eqref{eq:Krylov_Sobolev_lectures_theorems_8-5-3_and_6_estimate} and completes the proof of Theorem \ref{thm:Krylov_Sobolev_lectures_8-5-3}.
\end{proof}

Finally, we can apply Theorem \ref{thm:Krylov_Sobolev_lectures_8-5-3} to complete the

\begin{proof}[Proof of Theorem \ref{thm:Haller-Dintelmann_Heck_Hieber_Theorem_3-1_vector_bundle_manifold}]
It remains only to establish the resolvent estimate. From the inequality \eqref{eq:Krylov_Sobolev_lectures_theorems_8-5-3_and_6_estimate} with $f = (\sA + \lambda)u$, we obtain
$$
|\lambda| \|u\|_{L^p(X)} \leq C\|(\sA + \lambda)u\|_{L^p(X)},
$$
and, just as in the proof of the \apriori estimate \eqref{eq:Haller-Dintelmann_Heck_Hieber_Theorem_5-1_full_apriori_estimate_sector} in Theorem \ref{thm:Haller-Dintelmann_Heck_Hieber_5-1_full_apriori_estimate}, this gives \eqref{eq:Haller-Dintelmann_Heck_Hieber_Theorem_3-1_vector_bundle_manifold} by increasing $\lambda_0$ if necessary.

Observe that $\sD(\sA_p) = W^{m, p}(X; E)$ is dense in $L^p(X; E)$. Consequently, if $\theta > \pi/2$, then $\sA_p$ is a sectorial operator on $L^p(X; E)$, with constants $a = -\lambda_0$, $M := \max\{1, C\}$, and $\delta = \theta$ in Definition \ref{defn:Sell_You_page_78_definition_of_sectorial_operator}, and the generator of an analytic semigroup, $e^{- \sA_p t}$, on $L^p(X; E)$ by Theorem \ref{thm:Renardy_Rogers_12-31}.
\end{proof}

\begin{rmk}[Alternative proofs of Theorems \ref{thm:Krylov_Sobolev_lectures_8-5-3} and \ref{thm:Haller-Dintelmann_Heck_Hieber_Theorem_3-1_vector_bundle_manifold} when $m=2$]
\label{rmk:Krylov_Sobolev_lectures_8-5-3_alternative_proof}
When $m=2$ and $\sA$ is a second-order, strictly elliptic partial differential operator on $W^{2,p}(X; E)$ with $C^\infty$ coefficients, generalizing the local definition on $\Omega \subseteqq \RR^d$ for such an operator in Definition \ref{defn:Cannarsa_Terreni_Vespri_2-7}, then Theorem \ref{thm:Krylov_Sobolev_lectures_8-5-3} could also be deduced by
\begin{inparaenum}[\itshape a\upshape)]
\item replacing the role of Theorem \ref{thm:Haller-Dintelmann_Heck_Hieber_5-1_full_apriori_estimate} with that of Theorem \ref{thm:Cannarsa_Terreni_Vespri_6-7} when $2\leq p<\infty$ and Theorem \ref{thm:Cannarsa_Terreni_Vespri_6-7_one_lessthan_p_lessthan_two} when $1<p<2$,
\item replacing the role of Theorem \ref{thm:Corollary_Haller-Dintelmann_Heck_Hieber_3-2} by that of Corollary \ref{cor:Cannarsa_Terreni_Vespri_6-7_sectorial_operator_analytic_semigroup_LpOmega} (restricted to $1<p<\infty$).
\end{inparaenum}
\end{rmk}

\subsection{Existence, uniqueness, \apriori $L^\infty$ and $L^1$ estimates, and analytic semigroups on $C(X; E)$ and $L^1(X; E)$ defined by elliptic partial differential operators of order $m \geq 1$ on a complex vector bundle over a closed manifold}
\label{subsec:Heck_Hieber_Stavrakidis_vector_bundle_manifold}
Given a complex (or real) $C^\infty$ vector bundle $E$ of rank $N \geq 1$ over a closed, Riemannian, smooth manifold, $X$, of dimension $d \geq 2$, we next proceed to adapt the results of Sections \ref{subsec:Cannarsa_Terreni_Vespri_1985} and \ref{subsec:Heck_Hieber_Stavrakidis_2010} to the case of realizations of a parameter-elliptic partial differential operator $\sA$ of order $m \geq 1$, as in Definition \ref{defn:Elliptic_partial_differential_operator_vector_bundle_manifold}, acting on sections of $E$ over $X$,
\begin{gather*}
\sA_\infty: \sD(\sA_\infty) \subset L^\infty(X; E) \to L^\infty(X; E),
\\
\sA_c: \sD(\sA_c) \subset C(X; E) \to C(X; E),
\\
\sA_1: \sD(\sA_1) \subset L^1(X; E) \to L^1(X; E),
\end{gather*}
where we define, by analogy with \eqref{eq:Pazy_7-3-20}, \eqref{eq:Pazy_7-3-23}, and  \eqref{eq:Pazy_7-3-27}, respectively,
\begin{gather}
\label{eq:Heck_Hieber_Stavrakidis_page_301_definition_domain_A_infinity_vector_bundle_manifold}
\sD(\sA_\infty) := \left\{u \in \bigcap_{p>d} W^{m,p}(X; E): \ u, \sA(x, D)u \in L^\infty(X; E) \right\},
\\
\label{eq:Pazy_7-3-23_vector_bundle_on_manifold}
\sD(\sA_c) := \{u \in \sD(\sA_\infty): \sA u \in C(X; E)\},
\\
\label{eq:Pazy_7-3-27_vector_bundle_on_manifold}
\begin{split}
\sD(\sA_1) \subset W^{m-1,1}(X; E), \hbox{ where $\sA_1$ is the smallest closed extension of }
\\
\sA: C^\infty(X; E) \subset L^1(X; E) \to L^1(X; E).
\end{split}
\end{gather}
We recall from Remark \ref{rmk:Domain_elliptic_partial_differential_operator_on_L1} that the partial differential operator $\sA: C^\infty(X; E) \subset L^1(X; E) \to L^1(X; E)$ is indeed closable in $L^1(X; E)$. We then have the following vector bundle and manifold version of Theorem \ref{thm:Heck_Hieber_Stavrakidis_lemma_3-2}.

\begin{thm}[\Apriori estimate for a parameter-elliptic partial differential operator on $W^{m-1,\infty}(X; E)$]
\label{thm:Heck_Hieber_Stavrakidis_lemma_3-2_vector_bundle_manifold}
Assume the hypotheses of Theorem \ref{thm:Haller-Dintelmann_Heck_Hieber_Theorem_3-1_vector_bundle_manifold} and, in addition, that $q>d$. Then there are constants, $C_q > 0$ and $\omega_q \geq 0$, with the following significance.  If $u \in C^{m-1}(X; E)\cap W^{m,q}(X; E)$ and $\lambda \in \CC$ with $\Real \lambda \geq \min\{1, \omega_q\}$, then
\begin{multline}
\label{eq:Heck_Hieber_Stavrakidis_9prelim_vector_bundle_manifold}
\sum_{k=0}^{m-1} |\lambda|^{1 - k/m} \|\nabla^k u\|_{L^\infty(X)} + |\lambda|^{d/mq} \sup_{x_0 \in X} \|\nabla^m u\|_{L^q(B(x_0, |\lambda|^{-1/m})}
\\
\leq C_q |\lambda|^{d/m q} \sup_{x_0 \in X} \|(\sA + \lambda)u\|_{L^q(B(x_0, |\lambda|^{-1/m})}.
\end{multline}
In particular, if $\sA u \in L^\infty(X; E)$, then
\begin{multline}
\label{eq:Heck_Hieber_Stavrakidis_9_vector_bundle_manifold}
\sum_{k=0}^{m-1} |\lambda|^{1 - k/m} \|\nabla^k u\|_{L^\infty(X)} + |\lambda|^{d/mq} \sup_{x_0 \in X} \|\nabla^m u\|_{L^q(B(x_0, |\lambda|^{-1/m})}
\\
\leq C_q\|(\sA + \lambda)u\|_{L^\infty(X)}.
\end{multline}
\end{thm}

\begin{proof}
The global \apriori estimates \eqref{eq:Heck_Hieber_Stavrakidis_9prelim_vector_bundle_manifold} and \eqref{eq:Heck_Hieber_Stavrakidis_9_vector_bundle_manifold} on $X$ follow from their analogues \eqref{eq:Heck_Hieber_Stavrakidis_9prelim} and \eqref{eq:Heck_Hieber_Stavrakidis_9} on $\RR^d$ in Theorem \ref{thm:Heck_Hieber_Stavrakidis_lemma_3-2} by the method of proof of the \apriori estimate \eqref{eq:Krylov_Sobolev_lectures_lemma_8-5-2_estimate} in Lemma \ref{lem:Krylov_Sobolev_lectures_8-5-2} or the \apriori estimate \eqref{eq:Krylov_Sobolev_lectures_theorems_8-5-3_and_6_estimate} in Theorem \ref{thm:Krylov_Sobolev_lectures_8-5-3}, \emph{mutatis mutandis}, except that the role played by Theorem \ref{thm:Haller-Dintelmann_Heck_Hieber_5-1_full_apriori_estimate} \eqref{item:Haller-Dintelmann_Heck_Hieber_5-1_full_apriori_estimate} is replaced by that of Theorem \ref{thm:Heck_Hieber_Stavrakidis_lemma_3-2}.
\end{proof}

We next have the following $L^\infty$ version of Theorem \ref{thm:Krylov_Sobolev_lectures_8-5-3}.

\begin{thm}[Existence and uniqueness for an elliptic system on $\sD(\sA_\infty) \subset L^\infty(X; E)$]
\label{thm:Krylov_Sobolev_lectures_8-5-3_Linfinity}
Assume the hypotheses of Theorem \ref{thm:Haller-Dintelmann_Heck_Hieber_Theorem_3-1_vector_bundle_manifold}. Then there are constants, $K_q>0$ and $\omega_q \geq 0$, with the following significance. If $f \in L^\infty(X; E)$ and $\lambda \in \Delta_\theta(\lambda_0)$ with $\Real\lambda \geq \min\{1, \omega_q\}$, then there is a unique solution, $u \in  W^{m,\infty}(X; E)$, to $(\sA + \lambda)u = f$ a.e. on $X$ and $u$ obeys the \apriori estimate \eqref{eq:Heck_Hieber_Stavrakidis_9_vector_bundle_manifold}.
\end{thm}

\begin{proof}
The proof follows the same pattern as that of Theorem \ref{thm:Krylov_Sobolev_lectures_8-5-3}, \emph{mutatis mutandis}, with the following observations:
\begin{enumerate}
  \item Lemma \ref{lem:Each_sAi_obeys_corollary_Haller-Dintelmann_Heck_Hieber_5-1} continues to hold with $p=\infty$, with the role of Theorem \ref{thm:Haller-Dintelmann_Heck_Hieber_5-1_full_apriori_estimate} in its hypothesis replaced by that of Theorem \ref{thm:Heck_Hieber_Stavrakidis_lemma_3-2};

  \item Lemma \ref{lem:Krylov_Sobolev_lectures_8-4-1} continues to hold with $p=\infty$;

  \item The Definition \ref{defn:Krylov_8-4-2} of the regularizer is valid with $p=\infty$;

  \item Lemma \ref{lem:Krylov_Sobolev_lectures_8-5-1} and Lemma \ref{lem:Krylov_Sobolev_lectures_8-5-2} continue to hold with $p=\infty$, where the role in their proofs played by Theorem \ref{thm:Haller-Dintelmann_Heck_Hieber_5-1_full_apriori_estimate} \eqref{item:Haller-Dintelmann_Heck_Hieber_5-1_full_apriori_estimate} is replaced by that of Theorem \ref{thm:Heck_Hieber_Stavrakidis_lemma_3-2};
\end{enumerate}
This completes the proof.
\end{proof}

\begin{rmk}[Alternative proof of Theorem \ref{thm:Krylov_Sobolev_lectures_8-5-3_Linfinity} when $m=2$]
\label{thm:Krylov_Sobolev_lectures_8-5-3_Linfinity_alternative_proof}
When $m=2$ and $\sA$ is a second-order, strictly elliptic partial differential operator on $W^{2,p}(X; E)$ with $C^\infty$ coefficients, generalizing the local definition on $\Omega \subseteqq \RR^d$ for such an operator in Definition \ref{defn:Cannarsa_Terreni_Vespri_2-7}, then Theorem \ref{thm:Krylov_Sobolev_lectures_8-5-3_Linfinity} could also be deduced by replacing the roles of the \apriori estimates in Theorem \ref{thm:Haller-Dintelmann_Heck_Hieber_5-1_full_apriori_estimate} by those of Theorem \ref{thm:Cannarsa_Terreni_Vespri_6-1} and Remark \ref{rmk:Cannarsa_Terreni_Vespri_6-4}.
\end{rmk}

We obtain the following $L^\infty$ version of Theorem \ref{thm:Haller-Dintelmann_Heck_Hieber_Theorem_3-1_vector_bundle_manifold} and vector bundle and manifold analogue of Theorem \ref{thm:Heck_Hieber_Stavrakidis_Theorem_1-1_corollary_sectorial_operator}.

\begin{thm}[Sectorial property of a parameter-elliptic operator and generation of an analytic semigroup on $C(X; E)$]
\label{thm:Heck_Hieber_Stavrakidis_Theorem_1-1_corollary_sectorial_operator_vector_bundle_manifold}
Assume the hypotheses of Theorem \ref{thm:Haller-Dintelmann_Heck_Hieber_Theorem_3-1_vector_bundle_manifold}. If $\theta_0 > \pi/2$, then there are constants, $C > 0$ and $\lambda_0 \geq 0$ and $\vartheta \in (\pi/2, \theta_0)$, such that $\rho(-\sA_c) \supset \Delta_\vartheta(\lambda_0)$ and
\begin{equation}
\label{eq:Heck_Hieber_Stavrakidis_Theorem_1-1_sectorial_estimate_vector_bundle_manifold}
\|(\sA_c + \lambda)^{-1}\|_{\sL(C(X; E))} \leq \frac{C}{|\lambda - \lambda_0|},
\quad \forall\, \lambda \in \Delta_\vartheta(\lambda_0).
\end{equation}
Moreover, $\sA_c$ is a sectorial operator on $C(X; E)$ in the sense of Definition \ref{defn:Sell_You_page_78_definition_of_sectorial_operator} and $-\sA_c$ generates an analytic semigroup, $e^{-\sA_c t}$, on $C(X; E)$.
\end{thm}

\begin{proof}
The proof follows the same pattern as that of Theorem \ref{thm:Heck_Hieber_Stavrakidis_Theorem_1-1_corollary_sectorial_operator}, \emph{mutatis mutandis}, with the following observations:
\begin{enumerate}
\item The role of the \apriori estimate \eqref{eq:Heck_Hieber_Stavrakidis_9} in Theorem \ref{thm:Heck_Hieber_Stavrakidis_lemma_3-2} is replaced by that of the \apriori estimate \eqref{eq:Heck_Hieber_Stavrakidis_9_vector_bundle_manifold} in Theorem \ref{thm:Heck_Hieber_Stavrakidis_lemma_3-2_vector_bundle_manifold};
\item The role of Theorem \ref{thm:Corollary_Haller-Dintelmann_Heck_Hieber_3-2} is replaced by that of Theorem \ref{thm:Haller-Dintelmann_Heck_Hieber_Theorem_3-1_vector_bundle_manifold}.
\end{enumerate}
Note finally that $\sD(\sA_c)$ is dense in $C(X; E)$ since $\sD(\sA_c)$ contains $C^\infty(X; E)$ and $C^\infty(X; E)$ is dense in $C_\infty(X; E)$. This completes the proof.
\end{proof}

The proof of Lemma \ref{lem:Pazy_7-3-8} extends without change to give\footnote{The result may also be extracted as a special case of Lemma \ref{lem:Pazy_7-3-8_vector_bundle_manifold}, whose proof we do include.}


\begin{lem}[Reverse H\"older inequality for $p=1$]
\label{lem:Pazy_7-3-8_vector_bundle_manifold}
Let $E$ be a complex Hermitian (or real Riemannian) $C^\infty$ vector bundle over a $C^\infty$, orientable, Riemannian manifold, $X$. If $u \in L^1(X; E)$, then
\begin{equation}
\label{eq:Pazy_7-3-26_vector_bundle_manifold}
\|u\|_{L^1(X)} = \sup_{\begin{subarray}{c} \varphi \in C^\infty(X; E) \\ \|\varphi\|_{L^\infty(X)} \leq 1\end{subarray}}
(u, \varphi)_{L^2(X)}.
\end{equation}
\end{lem}


We have the following analogue of Theorem \ref{thm:Cannarsa_Terreni_Vespri_6-7_p_is_one} for parameter-elliptic partial differential operators of order $m \geq 1$ acting on sections of a Hermitian smooth vector bundle with Hermitian connection over a closed, Riemannian, smooth manifold .

\begin{thm}[Existence, uniqueness, and \apriori $L^1$ estimate and domain for a parameter-elliptic partial differential operator]
\label{thm:Cannarsa_Terreni_Vespri_6-7_p_is_one_m_geq_one_vector_bundle_manifold}
Assume the hypotheses of Theorem \ref{thm:Haller-Dintelmann_Heck_Hieber_Theorem_3-1_vector_bundle_manifold}. Then the operator,
$$
\sA: C^m(X; E) \subset L^1(X; E) \to L^1(X; E),
$$
is closable and, if $\sA_1$ denotes the smallest closed extension, one has
\begin{equation}
\label{eq:Cannarsa_Terreni_Vespri_apriori_function_domain_A1_partial_characterization_m_geq_one_vector_bundle_manifold}
\sD(\sA_1) \subset \{u \in L^1(X; E): \sA u \in L^1(X; E)\}.
\end{equation}
Moreover, there are constants $C > 0$ and $\omega_1 \geq 0$ with the following significance. If $\lambda \in \CC$ obeys $\Real \lambda > \omega_1$ and $u \in \sD(\sA_1)$, then
\begin{equation}
\label{eq:Cannarsa_Terreni_Vespri_apriori_function_and_gradient_estimate_p_is_one_m_geq_one_vector_bundle_manifold}
\sum_{k=0}^{m-1} |\lambda|^{1 - k/m} \|\nabla^k u\|_{L^1(X)}
\leq C\|(\sA_1 + \lambda)u\|_{L^1(X)}.
\end{equation}
Finally, for all $\lambda \in \CC$ with $\Real\lambda > \omega_1$ and $f \in L^1(X; E)$, there exists a unique solution $u \in \sD(\sA_1)$ to
\begin{equation}
\label{eq:Cannarsa_Terreni_Vespri_3-1_p_is_one_m_geq_one_vector_bundle_manifold}
(\sA + \lambda)u = f \quad\hbox{a.e. on } X.
\end{equation}
\end{thm}

\begin{proof}
The proof by duality follows the same pattern as that of Theorem \ref{thm:Cannarsa_Terreni_Vespri_6-7_one_lessthan_p_lessthan_two} (when $1 < p < 2$) and especially Theorem \ref{thm:Cannarsa_Terreni_Vespri_6-7_p_is_one} (when $p = 1$), \emph{mutatis mutandis}, with the following changes:
\begin{enumerate}
\item The role of Theorem \ref{thm:Cannarsa_Terreni_Vespri_6-7} (providing existence of solutions when $f \in L^2(\Omega; \CC^N)$) is replaced by that of Theorem \ref{thm:Krylov_Sobolev_lectures_8-5-3} (providing existence of solutions when $f \in L^2(X; E)$);

\item The role of Theorem \ref{thm:Cannarsa_Terreni_Vespri_6-1} (providing a $W^{1,\infty}(\Omega;\CC^N)$ \apriori estimate \eqref{eq:Cannarsa_Terreni_Vespri_6-2} for a solution $u$ in terms of $f = (\sA + \lambda)u \in L^\infty(\Omega;\CC^N)$) is replaced by that of Theorem \ref{thm:Heck_Hieber_Stavrakidis_lemma_3-2_vector_bundle_manifold} (providing a $W^{m-1,\infty}(X; E)$ \apriori estimate \eqref{eq:Heck_Hieber_Stavrakidis_9_vector_bundle_manifold} for a solution $u$ in terms of $f = (\sA + \lambda)u \in L^\infty(X; E)$);

\item The role of Lemma \ref{lem:Pazy_7-3-8_vector_valued} (the reverse Holder inequality \eqref{eq:Pazy_7-3-26_vector_valued} for $u \in L^1(\Omega;\CC^N)$) is replaced by that of Lemma \ref{lem:Pazy_7-3-8_vector_bundle_manifold} (the reverse Holder inequality \ref{eq:Pazy_7-3-26_vector_bundle_manifold} for $u \in L^1(X; E)$).
\end{enumerate}
This completes the proof.
\end{proof}

In particular, we obtain the

\begin{cor}[Resolvent estimate for a parameter-elliptic partial differential operator, sectorial property, and generation of an analytic semigroup on $L^1(X; E)$]
\label{cor:Cannarsa_Terreni_Vespri_L1_m_geq_one_sectorial_operator_vector_bundle_manifold}
Assume the hypotheses of Theorem \ref{thm:Cannarsa_Terreni_Vespri_6-7_p_is_one_m_geq_one_vector_bundle_manifold}. Then there are constants, $C > 0$ and $\lambda_0 \geq 0$ and $\vartheta \in (\pi/2, \pi)$, such that $\rho(-\sA_1) \supset \Delta_\vartheta(\lambda_0)$ and
\begin{equation}
\label{eq:Cannarsa_Terreni_Vespri_L1_sectorial_estimate_vector_bundle_manifold}
\|(\sA_1 + \lambda)^{-1}\|_{\sL(L^1(X; E))} \leq \frac{C}{|\lambda - \lambda_0|},
\quad \forall\, \lambda \in \Delta_\vartheta(\lambda_0).
\end{equation}
Moreover, $\sA_1$ is a sectorial operator on $L^1(X; E)$ in the sense of Definition \ref{defn:Sell_You_page_78_definition_of_sectorial_operator} and $-\sA_1$ generates an analytic semigroup, $e^{-\sA_1 t}$, on $L^1(X; E)$.
\end{cor}

\begin{proof}
The proof follows the same pattern as that of Corollary \ref{cor:Cannarsa_Terreni_Vespri_6-7_sectorial_operator_analytic_semigroup_LpOmega} (restricted to $p=1$), \emph{mutatis mutandis}, with the role of Theorem \ref{thm:Cannarsa_Terreni_Vespri_6-7_p_is_one} replaced by that of Theorem \ref{thm:Cannarsa_Terreni_Vespri_6-7_p_is_one_m_geq_one_vector_bundle_manifold}. This completes the proof.
\end{proof}

\subsection{\Apriori interior $L^2$ and $L^p$ estimates for a solution to the linear heat equation}
\label{subsec:Krylov_2-4-4_covariant_derivative_vector_bundle}
Throughout this monograph, when we need to localize an \apriori global estimate, we shall almost always do this with the aid of an explicit (and often carefully chosen) cut-off function and then estimate the resulting commutator terms over the support of the gradient of the cut-off function. However, it can still be useful to have at hand a genuine \apriori interior estimate in the customary sense (for example, of Gilbarg and Trudinger \cite{GilbargTrudinger} or Krylov \cite{Krylov_LecturesHolder, Krylov_LecturesSobolev}). In this section, we point out that one can easily adapt, to the case of elliptic or parabolic systems over closed Riemannian manifolds, an elegant construction described by Krylov \cite{Krylov_LecturesHolder, Krylov_LecturesSobolev} for deriving an \apriori interior estimate from a global estimate for a solution to a scalar elliptic or parabolic partial differential equation on an open subset of $\RR^d$. The method described by Krylov cleverly uses an infinite sequence of cut-off functions and is arguably simpler than the weighted function space approach described by Gilbarg and Trudinger \cite[Section 6.1]{GilbargTrudinger}.

Let $X$ be a $C^\infty$ closed manifold. By analogy with \cite[page 60]{Krylov_LecturesSobolev}, given a time $t_0 \in \RR$ and a point $x_0 \in X$, we define the `parabolic cylinder',
$$
Q_r(t_0,x_0) := (t_0, t_0 + r^2) \times B_r(x_0) \subset \RR \times X.
$$
The following \apriori $L^2$ interior estimate is an analogue of \cite[Lemma 2.4.4]{Krylov_LecturesSobolev} for the scalar parabolic equation \cite[Equation (2.3.1)]{Krylov_LecturesSobolev} on open subsets of $\RR^{d+1}$. The interior estimate in \cite[Lemma 2.4.4]{Krylov_LecturesSobolev} follows from the corresponding global estimates in \cite[Theorems 2.3.1, 2.3.2, or Corollary 2.3.3]{Krylov_LecturesSobolev} for a solution over $(-\infty, T)\times\RR^d$, where $T \in (-\infty, \infty]$.

\begin{lem}
\label{lem:Krylov_2-4-4_covariant_derivative_vector_bundle}
Let $X$ be a $C^\infty$ closed manifold of dimension $d \geq 2$ with Riemannian metric $g$, and $E$ a complex Hermitian (real Riemannian) vector bundle over $X$, and $A$ a $C^\infty$ Hermitian (Riemannian) connection on $E$. Then there exists a positive constant, $C$, with the following significance. If $(t_0, x_0) \in \RR\times X$ and $0 < r < R < \infty$ and $\mu \geq 0$ and $u \in C^\infty(X; E)$, then
\begin{multline}
\label{eq:Krylov_2-4-1_covariant_derivative_vector_bundle}
\mu\|u\|_{L^2(Q_r(t_0,x_0))} + \|\nabla_Au\|_{L^2(Q_r(t_0,x_0))} + \|\nabla_A^2u\|_{L^2(Q_r(t_0,x_0))} + \|\partial_t u\|_{L^2(Q_r(t_0,x_0))}
\\
\leq
C\left(\|(L_A+\mu)u\|_{L^2(Q_R(t_0,x_0))} + \left(1 + (R-r)^{-2}\right)\|u\|_{L^2(Q_R(t_0,x_0))} \right),
\end{multline}
where $L_A = \partial_t + \nabla_A^*\nabla_A$.
\end{lem}

The following \apriori $L^p$ interior estimate is an analogue of \cite[Theorem 5.2.5]{Krylov_LecturesSobolev} for the scalar parabolic equation \cite[Equation (5.2.1)]{Krylov_LecturesSobolev} on open subsets of $\RR^{d+1}$ and is derived from the corresponding global \apriori estimate in \cite[Theorem 4.3.8]{Krylov_LecturesSobolev}.

\begin{thm}
\label{lem:Krylov_5-2-5_covariant_derivative_vector_bundle}
Let $k \geq 0$ be an integer, $p \in (1, \infty)$, and $r, R \in (0, \infty)$ with $r < R$, and $X$ be a $C^\infty$ closed manifold of dimension $d \geq 2$ with Riemannian metric $g$, and $E$ a complex Hermitian (real Riemannian) vector bundle over $X$, and $A$ a $C^\infty$ Hermitian (Riemannian) connection on $E$. Then there exists a positive constant, $C$, with the following significance. If $(t_0, x_0) \in \RR\times X$ and $\mu \geq 0$ and $u \in C^\infty(X; E)$, then,
\begin{multline}
\label{eq:Krylov_2-4-1_covariant_derivative_vector_bundle_higher_order}
\sum_{j=0}^k\left(\mu\|\nabla_A^j u\|_{L^p(Q_r(t_0,x_0))} + \|\nabla_A^{j+1} u\|_{L^p(Q_r(t_0,x_0))}  \right.
\\
+ \left. \|\nabla_A^{j+2} u\|_{L^p(Q_r(t_0,x_0))} + \|\nabla_A^j \partial_t u\|_{L^p(Q_r(t_0,x_0))} \right)
\\
\leq
C\left(\sum_{j=0}^k\|\nabla_A^j (L_A+\mu)u\|_{L^p(Q_R(t_0,x_0))} + \|u\|_{L^p(Q_R(t_0,x_0))} \right),
\end{multline}
where $L_A = \partial_t + \nabla_A^*\nabla_A$.
\end{thm}

\begin{rmk}[On the derivation of \apriori interior Schauder estimates in H\"older spaces for a solution to the linear heat equation]
One could also easily adapt the method described by Krylov for deducing an interior \apriori estimate from a global estimate over $\RR^d$ in the elliptic case \cite[Theorem 7.1.1]{Krylov_LecturesHolder} or $(-\infty,T)\times\RR^d$ in the parabolic case \cite[Theorems 8.11.1 and 8.12.1]{Krylov_LecturesHolder}.
\end{rmk}

\section[Analytic semigroups on critical-exponent Sobolev spaces]{Elliptic partial differential systems and analytic semigroups on critical-exponent Sobolev spaces}
\label{sec:Sell_You_3-8-2_critical_exponent_elliptic_Sobolev_spaces}
In Section \ref{sec:Sell_You_3-8-2_standard_Sobolev_spaces}, we established that a parameter-elliptic partial differential operator, $\sA$, of order $m \geq 1$ --- acting on sections of a complex vector bundle, $E$, of rank $N$ over a closed, Riemannian, smooth manifold, $X$, of dimension $d \geq 2$ --- had realizations $\sA_p$, $\sA_c$, and $\sA_1$ as sectorial operators on $L^p(X; E)$ (for $1 < p < \infty$), $C(X; E)$, and $L^1(X; E)$, respectively, and hence that $-\sA_p$, $-\sA_c$, and $-\sA_1$ are generators of analytic semigroups on $L^p(X; E)$, $C(X; E)$, and $L^1(X; E)$, respectively.

In this section, in order to obtain the optimal results later in our monograph
we wish to extend that analysis further and show, when $d=4$, that the realization, $\sA_\sharp$, of the parameter-elliptic partial differential operator $\sA$ on a certain Banach space --- defined first by Taubes and explored further by the author in \cite{FeehanSlice} --- is also a sectorial operator and hence the generator of analytic semigroup on that Banach space. Thus, in Section \ref{subsec:Feehan_slice_4_and_5} we recall the definition of these Banach spaces and some fundamental \apriori estimates from \cite[Sections 4 and 5]{FeehanSlice}, while in Section \ref{subsec:Sectorial_property_and_analytic_semigroup_on_Taubes_spaces} we show that the proof of Corollary \ref{cor:Cannarsa_Terreni_Vespri_L1_m_geq_one_sectorial_operator_vector_bundle_manifold} readily extends to give the desired sectorial operator property of $\sA_\sharp$.

\subsection{Critical exponent elliptic Sobolev spaces and \apriori estimates}
\label{subsec:Feehan_slice_4_and_5}
We now describe the basic properties of the
critical-exponent norms and corresponding
Banach spaces introduced by Taubes in
\cite{TauPath, TauFrame, TauStable, TauConf, TauGluing}.
In particular, we give the basic embedding, multiplication, and composition lemmata which we can expect from our experience with the analogous results for standard Sobolev spaces. We shall make frequent use of the pointwise
Kato Inequality \eqref{eq:FU_6-20_first-order_Kato_inequality},
that is, $|d|v||\le |\nabla_Av|$ for $v\in \Om^0(E)$, so that the norms of the embedding and multiplication maps depend at most on the
Riemannian manifold $(X,g)$. Moreover, for simplicity, we confine our
attention to the case of closed four-manifolds: there are obvious analogues
of the Sobolev lemmata described here for any $d$-manifold, with
$d>2$. Similarly, extensions are possible to the case of complete manifolds
bounded geometry (bounded curvature and injectivity radius uniformly
bounded from below)---see \cite{AdamsFournier},  \cite{Aubin} for further details for
Sobolev embedding results in those situations and for the construction of
Green kernels.  We refer the reader to the monograph of Adams and Fournier
\cite{AdamsFournier} for a comprehensive treatment of Sobolev spaces and to that of
Stein \cite{Stein} for a treatment based on potential functions.

Throughout this section, $A$, $B$ denote $C^\infty$ orthogonal connections on
Riemannian vector bundles $E$, $F$ over $X$ with $C^\infty$ sections $u$, $v$,
respectively.  We first have the following analogues of the $L^2$ and $L^4$
norms \cite[Equation (4.1)]{FeehanSlice},
\begin{equation}
\label{eq:BasicSharpNorms}
\begin{aligned}
\|u\|_{L^\sharp(X)} &= \sup_{x\in X}\|\dist^{-2}(x,\cdot)|u|\|_{L^1(X)}, \\
\|u\|_{L^{2\sharp}(X)} &= \sup_{x\in X}\|\dist^{-1}(x,\cdot)|u|\|_{L^2(X)},
\end{aligned}
\end{equation}
where $\dist(x,y)$ denotes the geodesic distance between points $x$ and $y$
in $X$ defined by the metric $g$; these norms have the same behavior as the
$L^2$ and $L^4$ norms with respect to constant rescalings of
the metric $g$---the $L^\sharp$ norm on two-forms and the $L^{2\sharp}$
norm on one-forms are {\em scale invariant\/}. Indeed, one sees this by
noting that if $g\mapsto \tg = \la^{-2}g$, then
$\dist_{\tg}(x,y)=\la^{-1}\dist_g(x,y)$ and $dV_{\tg}=\la^{-4}dV_g$, while
for any $a\in\Om^1(E)$ and $v\in\Om^2(E)$, we have
$|a|_{\tg}=\la|a|_g$, and $|v|_{\tg}=\la^2|v|_g$.

Next, we define analogues of the $L^2_1$ and $L^2_2$ norms
\begin{align*}
\|u\|_{L^2_{1,A}(X)} &=\|\nabla_Au\|_{L^2(X)} + \|u\|_{L^2(X)}, \\
\|u\|_{L^2_{2,A}(X)} &=\|\nabla_A^2u\|_{L^2(X)} + \|\nabla_Au\|_{L^2(X)} +
\|u\|_{L^2(X)},
\end{align*}
and set \cite[Equation (4.2)]{FeehanSlice}
\begin{subequations}
\label{eq:CompositeLSharp21Norms}
\begin{align}
\label{eq:FeehanSlice_4-2_W1_sharp}
\|u\|_{L^\sharp_{1,A}(X)}
&=\|\nabla_Au\|_{L^\sharp(X)} + \|u\|_{L^{2\sharp}(X)} + \|u\|_{L^\sharp(X)},
\\
\label{eq:FeehanSlice_4-2_W2_sharp}
\|u\|_{L^\sharp_{2,A}(X)}
&=\|\nabla_A^*\nabla_Au\|_{L^\sharp(X)} + \|u\|_{L^\sharp(X)},
\end{align}
\end{subequations}
where $\nabla_A^* = -*\nabla_A*:\Om^1(E)\to\Om^0(E)$ is the $L^2$-adjoint of
the map $\nabla_A:\Om^0(E)\to\Om^1(E)$.

Finally, we define analogues of the $C^0\cap L^2_2$ norm
$$
\|u\|_{C^0\cap L^2_{2,A}(X)} = \|u\|_{C^0(X)} + \|u\|_{L^2_{2,A}(X)},
$$
and set \cite[Equation (4.3)]{FeehanSlice}
\begin{subequations}
\label{eq:WholeFamilyOfSharpNorms}
\begin{align}
\|u\|_{L^{\sharp,2}(X)} &= \|u\|_{L^\sharp\cap L^2(X)} =
\|u\|_{L^\sharp(X)} + \|u\|_{L^2(X)}, \\
\|u\|_{L^{2\sharp,4}(X)} &= \|u\|_{L^{2\sharp}\cap L^4(X)} =
\|u\|_{L^{2\sharp}(X)} + \|u\|_{L^4(X)},   \\
\|u\|_{L^{\sharp,2}_{1,A}(X)} &= \|u\|_{L^\sharp_{1,A}\cap L^2_{1,A}(X)} =
\|u\|_{L^\sharp_{1,A}(X)} + \|u\|_{L^2_{1,A}(X)},   \\
\label{eq:Feehan_slice_4-3_W2+sharp_2}
\|u\|_{L^{\sharp,2}_{2,A}(X)} &= \|u\|_{L^\sharp_{2,A}\cap L^2_{2,A}(X)} =
\|u\|_{L^\sharp_{2,A}(X)} + \|u\|_{L^2_{2,A}(X)} .
\end{align}
\end{subequations}
It might have appeared, at first glance, a little more natural to continue the
obvious pattern and instead define $\|u\|_{L^\sharp_{2,A}(X)}$ using
$\|\nabla_A^2u\|_{L^\sharp(X)}$: as we shall below, though, the given
definition is most useful in practice. For related reasons, if $u\in\Om^1(E) =
\Om^0(\La^1\otimes E)$,
it is convenient to define the norm $\|u\|_{L^\sharp_{1,A}(X)}$ by \cite[Equation (4.4)]{FeehanSlice}
\begin{equation}
\|u\|_{L^\sharp_{1,A}(X)}
=\|\nabla_A^*u\|_{L^\sharp(X)} + \|u\|_{L^{2\sharp}(X)} +
\|u\|_{L^\sharp(X)}.
\label{eq:LSharp1AdjointCov}
\end{equation}
Let $L^\sharp(X)$ be the Banach space completion of $C^\infty(X)$ with respect to
the norm $\|\cdot\|_{L^\sharp}$ and similarly define the remaining Banach
spaces above.

We have the following extensions of the standard Sobolev embedding theorem
\cite{FU}, \cite{PalaisFoundationGlobal}: their proofs are given in \cite[Section 5]{FeehanSlice}. See also
\cite{DonApprox}, \cite{ParkerTaubes}, \cite{TauPath}, \cite[Section 6]{TauFrame},
\cite[Equation (3.4) and Section 5]{TauStable}, and \cite[Lemma 4.7]{TauConf}.

\begin{lem}
\label{lem:Feehan_4-1}
\cite[Lemma 4.1]{FeehanSlice}
The following are continuous embeddings:
\begin{enumerate}
\item $L^p_k(E)\subset L^\sharp_k(E)$, for $k=0,1,2$ and all $p>2$;
\item $L^q(E)\subset L^{2\sharp}(E)$, for all $q>4$;
\item $L^2_1(E)\subset L^{2\sharp}(E)$.
\end{enumerate}
\end{lem}

In the reverse direction we have:

\begin{lem}
\label{lem:Feehan_4-2}
\cite[Lemma 4.2]{FeehanSlice}
The following are continuous embeddings:
\begin{enumerate}
\item $L^\sharp(E)\subset L^1(E)$ and $L^{2\sharp}(E)\subset L^2(E)$;
\item $L^\sharp_2(E)\subset C^0\cap L^2_1(E)$;
\end{enumerate}
\end{lem}

We next consider the extension of the standard Sobolev multiplication lemma
\cite{FU}, \cite{PalaisFoundationGlobal}. While there is no continuous
multiplication map
$L^2_2\times L^2_2 \to L^2_2$, it is worth observing that there is a
continuous bilinear map $C^0\cap L^2_2(E)\times C^0\cap L^2_2(F)
\to C^0\cap L^2_2(E\otimes F)$ given by $(u,v)\mapsto u\otimes v$.
Note that for $u\in \Om^0(E)$ and $v\in \Om^0(F)$ we have \cite[Equation (4.5)]{FeehanSlice}
\begin{equation}
\begin{aligned}
\cov^2_{A\otimes B}(u\otimes v)
&= (\nabla_A^2 u)\otimes v + 2\nabla_A u\otimes\nabla_B v + u\otimes\cov _B^2v,
\label{eq:SecondCovDerivLeibnitz}
\\
\cov^*_{A\otimes B}\nabla_{A\otimes B}(u\otimes v)
&= (\nabla_A^*\nabla_A u)\otimes v + *((*\nabla_Au)\wedge\nabla_B v)
\\
&\qquad - *(\nabla_A u\wedge *\nabla_B v) + u\otimes\nabla_B^*\nabla_B v.
\end{aligned}
\end{equation}
Similarly, for $u\in\Om^0(\La^1\otimes E)$ and $v\in\Om^0(F)$, we have \cite[Equation (4.6)]{FeehanSlice}
\begin{equation}
\cov^*_{A\otimes B}(u\otimes v)
= (\nabla_A^*u)\otimes v + *(*u\wedge \nabla_B v)
\label{eq:AdjointCovLeibnitz}
\end{equation}
In particular, we see that if $u,v\in\Om^0(\gl(E))$, then \cite[Equation (4.7)]{FeehanSlice}
\begin{equation}
\begin{aligned}
\cov^*_A\nabla_A(uv)
&= (\nabla_A^*\nabla_A u)v + *((*\nabla_Au)\wedge\nabla_A v)
\\
&\qquad - *(\nabla_A u\wedge (*\nabla_A v)) + u(\nabla_A^*\nabla_A v),
\label{eq:LaplaceLeibnitz}
\end{aligned}
\end{equation}
an identity we will need in the next section.

\begin{lem}
\label{lem:Feehan_4-3}
\cite[Lemma 4.3]{FeehanSlice}
Let $\Om^0(E)\times\Om^0(F)\to \Om^0(E\otimes F)$ be given by $(u,v)\mapsto
u\otimes v$. Then the following hold.
\begin{enumerate}
\item The map $C^0(E)\otimes L^\sharp(F)\to
L^\sharp(E\otimes F)$ is continuous;
\item The map $L^{2\sharp}(E)\otimes
L^{2\sharp}(F) \to L^\sharp(E\otimes F)$ is continuous;
\item The spaces $L^\sharp_1(F)$, $L^2_1(F)$, and
$L^\sharp_2(F)$ are $L^\sharp_2(E)$-modules;
\item The spaces $L^2_1(F)$, $L^{\sharp,2}_1(F)$, and
$L^{\sharp,2}_2(F)$ are $L^{\sharp,2}_2(E)$-modules;
\end{enumerate}
The conclusions continue to hold for $\Om^1(E)$ in place of $\Om^0(E)$ and
the norms on $L^{\sharp}_1(\La^1\otimes E)$ and
$L^{\sharp,2}_1(\La^1\otimes E)$ defined via
\eqref{eq:LSharp1AdjointCov}.
\end{lem}

We continue the notation and assumptions of the preceding paragraphs.
The estimates which we collect here are due to Taubes and they arise, in a variety of contexts, in the proofs of
\cite[Lemma 5.4]{ParkerTaubes}, \cite[Equation (2.14) and Lemmata 3.5, 3.6, and A.3]{TauPath},
\cite[Equation (3.4b) and Lemma 6.2]{TauFrame}, \cite[Lemma 5.6]{TauStable},
and \cite[Sections 4 (c), (d), and (e)]{TauConf}.

Let $G(x,y)$ be the kernel function for
the Green's operator $(d^*d+1)^{-1}$ of the Laplace operator $d^*d+1$ on
$C^\infty(X)$. The kernel $G(x,y)$ of $(d^*d+1)^{-1}$
behaves like $\dist^{-2}(x,y)$ as $\dist(x,y)\to 0$ (see
\cite[Lemma 4.7]{TauConf} and \cite[Section 5]{TauStable}):

\begin{lem}
\label{lem:Green}
\cite[Lemma 5.1]{FeehanSlice}
The kernel $G(x,y)$ is a positive $C^\infty$ function away from the diagonal
in $X\times X$ and as $\dist(x,y)\to 0$,
$$
G(x,y) = {\frac{1}{4\pi^2\dist^2(x,y)}} + o(\dist^{-2}(x,y)).
$$
\end{lem}

Lemma \ref{lem:Green} implies that there is a constant $c$ depending
at most on $g$ such that for all $x\ne y$ in $X$,
\begin{equation}
c^{-1}\dist^{-2}(x,y) \le G(x,y)
\le c\dist^{-2}(x,y). \label{eq:GreenKerEst}
\end{equation}
Consequently, for all $u\in \Om^0(E)$, we have
\begin{equation}
c^{-1}\|u\|_{L^\sharp(X)} \le \|G|u|\|_{C^0(X)} \le c\|u\|_{L^\sharp(X)}.
\label{eq:LinftyGreenEst}
\end{equation}
An inequality similar to that in Lemma \ref{lem:Lp*Est} below is given by
\cite[Equation (3.4)]{TauStable}; see \cite[Lemma 5.4(a)]{ParkerTaubes} for
a related result on $\RR^3$.

\begin{lem}
\label{lem:Lp*Est}
\cite[Lemma 5.2]{FeehanSlice}
For all $f\in L^2_1(\RR^4)$, where $\RR^4$ has its standard metric,
$$
\sup_{x\in\RR^4}\|\dist^{-1}(x,\cdot)f\|_{L^2(\RR^4)}
\le \frac{1}{2}\|\cov f\|_{L^2(\RR^4)}.
$$
Suppose $X$ be a closed, oriented, Riemannian four-manifold. Then there is a
positive constant $c$ such that for all $f\in L^2_1(X)$,
$$
\sup_{x\in X}\|\dist^{-1}(x,\cdot)f\|_{L^2(X)} \le c\|f\|_{L^2_1(X)}.
$$
\end{lem}

The key estimates (1) and (2)
in Lemma \ref{lem:Feehan_5-3} below and the
estimates (1), (2), and (3) in Lemma \ref{lem:L22Estu} are essentially
those of \cite[Lemma 6.2]{TauFrame}, except that the dependence of the constant on
$\|F_A\|_{L^2}$ is made explicit, but the argument is the same as that in
\cite{TauFrame}.

\begin{lem}
\label{lem:Feehan_5-3}
\cite[Lemma 5.3]{FeehanSlice}
Let $X$ be a closed, oriented four-manifold with metric $g$. Then there is
a constant $c$ with the following significance. Let $E$ be a Riemannian
vector bundle over $X$ and let $A$ be an orthogonal
$L^2_2$ connection on $E$ with curvature $F_A$. Then
$L^\sharp_2(E)\subset C^0\cap L^2_1(E)$ and the following estimates hold:
\begin{align}
\|\nabla_Au\|_{L^{2\sharp}(X)}+ \|u\|_{C^0(X)}
&\le c\|\nabla_A^*\nabla_Au\|_{L^\sharp(X)} + \|u\|_{L^\sharp(X)}, \tag{1}\\
\|\nabla_Au\|_{L^{2\sharp}(X)}+ \|u\|_{C^0(X)}
&\le c\|\nabla_A^*\nabla_Au\|_{L^\sharp(X)} + \|u\|_{L^2(X)}, \tag{2}\\
\|u\|_{L^1(X)} &\le c\|u\|_{L^{\sharp}(X)}, \tag{3} \\
\|u\|_{L^2(X)} &\le c\|u\|_{L^{2\sharp}(X)}, \tag{4} \\
\|\nabla_Au\|_{L^2(X)} &\le c\|\nabla_Au\|_{L^{2\sharp}(X)}. \tag{5}
\end{align}
\end{lem}

\begin{lem}
\label{lem:L22Estu}
\cite[Lemma 5.4]{FeehanSlice}
Continue the hypotheses of Lemma \ref{lem:Feehan_5-3}. Then for any $u\in
(C^0\cap L^2_2)(E)$, we have
\begin{align}
\|\nabla_A^2u\|_{L^2(X)} &\le \|\nabla_A^*\nabla_Au\|_{L^2(X)} +
c\|F_A\|_{L^2(X)}^{1/2}\|\nabla_Au\|_{L^4(X)} \tag{1}\\
&\qquad + \|F_A\|_{L^2(X)}\|u\|_{C^0(X)}, \notag\\
\|\nabla_Au\|_{L^4(X)} &\le \|u\|_{C^0(X)}^{1/2}
\left(\|\nabla_A^*\nabla_Au\|_{L^2(X)}+2\|\nabla_A^2u\|_{L^2(X)}\right)^{1/2},
\tag{2}\\
\|\nabla_A^2u\|_{L^2(X)}
&\le 2\|\nabla_A^*\nabla_Au\|_{L^2(X)} + c\|F_A\|_{L^2(X)}\|u\|_{C^0(X)}. \tag{3}
\end{align}
\end{lem}

\begin{lem}
\label{lem:L22InfinityEstu}
\cite[Lemma 5.5]{FeehanSlice}
Continue the hypotheses of Lemma \ref{lem:Feehan_5-3}. Then for any $u\in
L^{\sharp,2}_2(E)$, we have:
$$
\|u\|_{L^2_{2,A}(X)} + \|u\|_{C^0(X)}
\le c(1+\|F_A\|_{L^2(X)})\left(\|\nabla_A^*\nabla_Au\|_{L^{\sharp,2}(X)} +
\|u\|_{L^2(X)}\right).
$$
\end{lem}

\subsection{Sectorial property of an elliptic operator and generation of an analytic semigroup on Taubes spaces}
\label{subsec:Sectorial_property_and_analytic_semigroup_on_Taubes_spaces}
By analogy with the definition of standard Sobolev spaces \cite[Section 3.2]{AdamsFournier}, we make the

\begin{defn}[Taubes domain and range spaces]
\label{eq:Taubes_domains_and_range_spaces}
Let $E$ be a complex Hermitian (or real Riemannian) $C^\infty$ vector bundle with $C^\infty$ Hermitian (or orthogonal) connection $\nabla$ over a $C^\infty$, closed, four-dimensional, orientable, Riemannian manifold, $X$. The \emph{Taubes domain space}, $W^{2, 2+\sharp}(X; E)$, is the completion of $C^\infty(X; E)$ with respect to the norm,
\begin{equation}
\label{eq:Taubes_domain_space_norm}
\|u\|_{W^{2, 2+\sharp}(X)} := \|u\|_{C(X)} + \|\nabla u\|_{L^2(X)} + \|\nabla^2 u\|_{L^2(X)} + \|\nabla^*\nabla u\|_{L^\sharp(X)}.
\end{equation}
The \emph{Taubes range space}, $L^{2,\sharp}(X; E)$, is the completion of $C^\infty(X; E)$ with respect to the norm,
\begin{equation}
\label{eq:Taubes_range_space_norm}
\|u\|_{L^{2,\sharp}(X)} := \|u\|_{L^2(X)} + \|u\|_{L^{\sharp}(X)}.
\end{equation}
\end{defn}

We know from \cite{FeehanSlice} that $W^{2, 2+\sharp}(X; E)$ and $L^{2,\sharp}(X; E)$ are indeed Banach spaces. We know also from \cite[Lemmata 4.1, 5.3, and 5.4]{FeehanSlice} that the norm in \eqref{eq:Taubes_domain_space_norm} has many equivalent variants (remembering that $X$ is compact), including
\begin{align}
\label{eq:Taubes_domain_space_norm_maximal}
\tag{\ref*{eq:Taubes_domain_space_norm}$'$}
\|u\|'_{W^{2, 2+\sharp}(X)} &:= \|u\|_{C(X)} + \|\nabla u\|_{L^{2,\sharp}(X)} + \|\nabla^2 u\|_{L^2(X)} + \|\nabla^*\nabla u\|_{L^\sharp(X)},
\\
\label{eq:Taubes_domain_space_norm_minimal_L2}
\tag{\ref*{eq:Taubes_domain_space_norm}$''$}
\|u\|''_{W^{2, 2+\sharp}(X)} &:= \|u\|_{L^2(X)} + \|\nabla^*\nabla u\|_{L^{2,\sharp}(X)},
\\
\label{eq:Taubes_domain_space_norm_minimal_Lsharp}
\tag{\ref*{eq:Taubes_domain_space_norm}$'''$}
\|u\|'''_{W^{2, 2+\sharp}(X)} &:= \|u\|_{L^{\sharp}(X)} + \|\nabla^*\nabla u\|_{L^{2,\sharp}(X)},
\end{align}
We select the choice \eqref{eq:Taubes_domain_space_norm}, rather than its equivalent variants, merely as a matter of convenience for a norm defining $W^{2, 2+\sharp}(X; E)$.

Our goal in this section is prove Theorem \ref{thm:Taubes_existence_uniqueness_apriori_estimate_sA_plus_lambda}, which establishes an \apriori estimate, existence, and uniqueness for a solution $u \in W^{2; \sharp, 2}(X; E)$ to $(\sA + \lambda)u = f$ a.e. on $X$, given $f \in $ and $\lambda \in \CC$ with $\Real\lambda > \lambda_0$, for some constant $\lambda_0 \geq 0$, and hence Corollary \ref{cor:Sectorial_property_and_analytic_semigroup_on_Taubes_spaces}, which establishes the sectorial property of the realization, $\sA_{\sharp, 2}$, of a second-order, elliptic partial differential operator, $\sA$, with principal symbol given by the Riemannian metric on $X$, and hence the fact that $-\sA_{\sharp, 2}$ is the infinitesimal generator of an analytic semigroup on $L^{\sharp, 2}(X; E)$.

Our definition \eqref{eq:BasicSharpNorms} of the Banach space $L^\sharp(X; E)$ and of the Banach spaces $L^{2,\sharp}(X; E)$ and $W^{2; 2, \sharp}(X; E)$ in Definition \ref{eq:Taubes_domains_and_range_spaces} suggest that we should model our proof on the corresponding result, Corollary \ref{cor:Cannarsa_Terreni_Vespri_L1_m_geq_one_sectorial_operator_vector_bundle_manifold}, for the Banach spaces $L^1(X; E)$ and $W^{2, 1}(X; E)$, respectively.

For any $x \in X$, let $\mu_x$ be the measure defined on Borel subsets $B \subset X$ by
\begin{equation}
\label{eq:Taubes_inverse_square_distance_weight_Lebesgue_measure}
\mu_x(B) := \int_B\dist_g^{-2}(x, \cdot)\, d\vol_g,
\end{equation}
where we recall that $\dist_g^{-2}(x, \cdot)$ is the weight function, defined by the Riemannian metric $g$ on $X$, appearing in the definition \eqref{eq:BasicSharpNorms} of the Banach space $L^\sharp(X)$. (We assume throughout our monograph, as will normally be the case in all our applications, that $X$ is orientable but, if not, one can simply replace the Riemannian volume form, $d\vol_g$, by the Riemannian volume density $|d\vol_g|$.) The density, $\dist_g^{-2}(x, \cdot)\, d\vol_g$, in \eqref{eq:Taubes_inverse_square_distance_weight_Lebesgue_measure} is determined by a weight, $\dist_g^{-2}(x, \cdot) \in C^\infty(X\less\{x\})$, which is only defined almost everywhere on $X$ and, in particular, is not $C^1$ on $X$ as assumed in standard expositions of the integration theory needed to define $L^p$ spaces on manifolds (for example, \cite[Proposition 1.23]{BerlineGetzlerVergne}, \cite[Section 11.4]{Folland_realanalysis}, \cite[Section 6.3]{Hormander_v1}), but clearly this does not effect the definition of the integral. We now have the following slight extension of Lemma \ref{lem:Pazy_7-3-8_vector_bundle_manifold}, where the Riemannian volume form is replaced by an arbitrary Borel measure, $\mu$, allowing us to include examples such as \eqref{eq:Taubes_inverse_square_distance_weight_Lebesgue_measure}.

\begin{lem}[Reverse H\"older inequality for $p=1$]
\label{lem:Pazy_7-3-8_vector_bundle_manifold_Taubes_weight}
Let $E$ be a complex Hermitian (or real Riemannian) $C^\infty$ vector bundle over a $C^\infty$, orientable, Riemannian manifold, $X$, and let $\mu$ be a Borel measure on $X$. If $u \in L^1(X, \mu; E)$, then
\begin{equation}
\label{eq:Pazy_7-3-26_vector_bundle_manifold_Taubes_weight}
\|u\|_{L^1(X, \mu)}
=
\sup_{\begin{subarray}{c} \varphi \in C^\infty(X; E) \\ \|\varphi\|_{L^\infty(X)} \leq 1\end{subarray}}
(u, \varphi)_{L^2(X, \mu)}.
\end{equation}
\end{lem}

\begin{proof}
We directly adapt, \emph{mutatis mutandis}, the proof of \cite[Lemma 7.3.8]{Pazy_1983} (whose statement is recorded in this monograph as Lemma \ref{lem:Pazy_7-3-8}) but we shall include the details for the sake of completeness and because of its importance. For every $\varphi \in C^\infty(X; E)$ satisfying $\|\varphi\|_{L^\infty(X)} \leq 1$ we have
$$
\left|(u, \varphi)_{L^2(X, \mu)}\right| \leq \|u\|_{L^1(X,\mu)}\|\varphi\|_{L^\infty(X)} \leq \|u\|_{L^1(X,\mu)},
$$
and so the supremum on the right-hand side of \eqref{eq:Pazy_7-3-26_vector_bundle_manifold_Taubes_weight} is \emph{less} than or equal to $\|u\|_{L^1(X, \mu)}$. Since $C^\infty(X; E)$ is dense in $L^1(X, \mu; E)$ it suffices to prove the result for $u \in C^\infty(X; E)$. Let $\{p_n\}_{n \in \NN} \subset C^\infty(X; \End_\CC(E))$ be a sequence with the properties that, for each $n\in \NN$,
\begin{multline*}
p_n(0) = 0 \quad\hbox{and}\quad |p_n(z)| \leq 1 \quad\forall\, z \in \CC, \quad\hbox{and}\quad
\\
p_n(z) = \frac{\bar z}{|z|} \quad\forall\, x \in X \hbox{ and } z \in E_x \cong \CC^N \hbox{ with } |z| \geq \frac{1}{n+1},
\end{multline*}
where the `$0$' in the expression $p_n(0)$ denotes the zero element of $\End_\CC(E_x)$ for any $x \in X$. Then $p_n\circ u \in C^\infty(X; E)$ and $\|p_n\circ u \|_{L^\infty(X)} \leq 1$ for all $n \in \NN$. The Lebesgue Dominated Convergence Theorem yields
$$
\lim_{n\to\infty} (u, p_n\circ u)_{L^2(X, \mu)} = \int_X |u| \,d\mu = \|u\|_{L^1(X, \mu)},
$$
and hence the supremum on the right-hand side of \eqref{eq:Pazy_7-3-26_vector_bundle_manifold_Taubes_weight} is \emph{greater} than or equal to $\|u\|_{L^1(X, \mu)}$.
\end{proof}

For ease of language, we shall restrict our attention to the case of Hermitian connections on complex Hermitian vector bundles, with the case of orthogonal connections on (real) Riemannian vector bundles following \emph{mutatis mutandis}. Recall from \eqref{eq:Connection_Laplacian_principal_symbol} that the connection Laplacian, $\nabla^*\nabla$, has principal symbol $-\id_E\otimes g \in C^\infty(X;\End_\CC(E)\otimes S^2(T^*X))$ and clearly obeys
\begin{equation}
\label{eq:Connection_Laplacian_parameter_ellipticity}
\sigma(\id_E\otimes g(x;\xi)) \subset (\kappa, \sK), \quad \forall\, x\in X \hbox{ and } \xi \in (TX)_x,
\end{equation}
for some positive constants, $\kappa$ and $\sK$, where $\sigma(B)$ denotes the spectrum of a linear operator $B \in \End_\CC(E_x)$. With respect to a local chart, $\{x^1,\ldots,x^4\}$, on a coordinate neighborhood in $X$ we have $g = \sum_{i,j=1}^4 g^{ij}dx^i\otimes dx^j$ and $\xi =  \sum_{i=1}^4 \xi_i \partial/\partial x_i$ and
$$
g(x;\xi) =  \sum_{i,j=1}^4 g^{ij}(x)\xi_i\xi_j.
$$
Therefore, $\nabla^*\nabla$ is $(\sM,\theta)$-parameter elliptic on $X$ in the sense of Definition \ref{defn:Elliptic_partial_differential_operator_vector_bundle_manifold}, with parameters $\sM \leq \kappa^{-1}$ and $\theta = 0$, and also strictly elliptic with ellipticity constant $\kappa$ in the sense of Remark \ref{rmk:Elliptic_differential_operator_vector_bundle_manifold_alternative_definitions} because, \afortiori, it is $(\sM, \theta)$-parameter elliptic with parameters $\sM \leq \kappa^{-1}$ and $\theta = 0$.

Our proof of Theorem \ref{thm:Taubes_existence_uniqueness_apriori_estimate_sA_plus_lambda} with the aid of a duality argument and weighted $L^1$-spaces has a precedent in the work of Cannarsa and Vespri \cite{Cannarsa_Vespri_1988} in their proofs of existence, uniqueness, \apriori estimates, sectoriality, and generation of analytic semigroups on $L^1(\RR^d)$ --- see their \cite[Remark 3.7 and Theorem 4.4]{Cannarsa_Vespri_1988}.

\begin{thm}[Existence, uniqueness, and \apriori estimate for a second-order elliptic partial differential operator on $W^{2, 2+\sharp}(X; E)$ with scalar principal symbol]
\label{thm:Taubes_existence_uniqueness_apriori_estimate_sA_plus_lambda}
Let $E$ be a $C^\infty$, complex Hermitian vector bundle with Hermitian connection $\nabla$ over a $C^\infty$, orientable, four-dimensional manifold, $X$, with Riemannian metric $g$. Let $\sA$ be a second-order partial differential operator on $C^\infty(X; E)$ with $C^\infty$ coefficients in the sense of Definition \ref{defn:Elliptic_partial_differential_operator_vector_bundle_manifold} and has principal symbol given by the negative of the Riemannian metric, namely $-\id_E \otimes g \in C^\infty(X;\End_\CC(E)\otimes S^2(T^*X))$. Then there are constants, $C > 0$ and $\lambda_0 \geq 0$, with the following significance. For all $\lambda \in \CC$ with $\Real \lambda > \lambda_0$ and $u \in W^{2, 2+\sharp}(X; E)$, one has
\begin{multline}
\label{eq:Taubes_sA_plus_lambda_apriori_estimate}
|\lambda|\|u\|_{L^{2,\sharp}(X)} + |\lambda|^{1/2}\|\nabla u\|_{L^2(X)} + \|u\|_{C(X)} + \|\nabla^2 u\|_{L^2(X)} + \|\nabla^*\nabla u\|_{L^\sharp(X)}
\\
\leq C\|(\sA + \lambda)u\|_{L^{2, \sharp}(X)}.
\end{multline}
Moreover, given $f \in L^{2, \sharp}(X; E)$, there is a unique solution $u \in W^{2, 2+\sharp}(X; E)$ to
\begin{equation}
\label{eq:Taubes_sA_plus_lambda_equation}
(\sA + \lambda)u = f \quad\hbox{a.e. on } X.
\end{equation}
\end{thm}

\begin{proof}
In establishing the \apriori estimate \eqref{eq:Taubes_sA_plus_lambda_apriori_estimate}, it suffices by continuity to assume that $u \in C^\infty(X; E)$. Suppose first that $\sA = \nabla^*\nabla$ and recall that $\nabla^*\nabla$ is $(\sM, \theta)$-parameter elliptic in the sense of Definition \ref{defn:Haller-Dintelmann_Heck_Hieber_page_720_parameter_elliptic},  with parameters $\sM = \kappa^{-1}$ and $\theta = 0$, for a positive constant of ellipticity, $\kappa$, depending only on the Riemannian metric, $g$, on $X$. For any $\theta \in (0,\pi)$, Theorem \ref{thm:Krylov_Sobolev_lectures_8-5-3} (taking $p = 2$) then implies that there are constants, $\lambda_0 \geq 0$ and $C > 0$, such that
\begin{multline}
\label{eq:Taubes_sA_plus_lambda_apriori_estimate_pure_L_2_part}
|\lambda|\|u\|_{L^2(X)} + |\lambda|^{1/2}\|\nabla u\|_{L^2(X)} + \|\nabla^2 u\|_{L^2(X)}
\leq C\|(\nabla^*\nabla + \lambda)u\|_{L^2(X)},
\\
\forall\, \lambda \in \Delta_\theta(\lambda_0) \hbox{ and } u \in C^\infty(X; E).
\end{multline}
On the other hand, the \apriori estimate \eqref{eq:Heck_Hieber_Stavrakidis_9_vector_bundle_manifold} in Theorem \ref{thm:Heck_Hieber_Stavrakidis_lemma_3-2_vector_bundle_manifold} (specialized to $m=2$) implies that, for possibly larger constants, $\lambda_0 \geq 0$ and $C > 0$,
$$
|\lambda|\|u\|_{C(X)} \leq C\|(\nabla^*\nabla + \lambda)u\|_{C(X)}, \quad \forall\, \lambda \in \Delta_{\pi/2}(\lambda_0) \hbox{ and } u \in C^\infty(X; E).
$$
Choosing $\mu = \mu_x = \dist_g^{-2}(\cdot, x)$ in the definition \eqref{eq:Taubes_inverse_square_distance_weight_Lebesgue_measure} (for any $x \in X$) of a Borel measure on $X$ and also in Lemma \ref{lem:Pazy_7-3-8_vector_bundle_manifold_Taubes_weight} yields, using duality just as in the derivation of the \apriori estimate \eqref{eq:Cannarsa_Terreni_Vespri_3-1_one_leq_p_lessthan_2_prelim} in the proof of Theorem \ref{thm:Cannarsa_Terreni_Vespri_6-7_one_lessthan_p_lessthan_two},
$$
|\lambda|\|u\|_{L^1(X,\mu_x)} \leq C\|(\nabla^*\nabla + \lambda)u\|_{L^1(X,\mu_x)},
\quad\forall\, \lambda \in \Delta_{\pi/2}(\lambda_0) \hbox{ and } u \in C^\infty(X; E).
$$
Taking the supremum over all $x \in X$ and applying the definition \eqref{eq:BasicSharpNorms} of the norm $\|\cdot\|_{L^\sharp(X)}$ gives
\begin{equation}
\label{eq:Taubes_sA_plus_lambda_apriori_estimate_pure_L_sharp_part}
|\lambda|\|u\|_{L^\sharp(X)} \leq C\|(\nabla^*\nabla + \lambda)u\|_{L^\sharp(X)},
\quad\forall\, \lambda \in \Delta_{\pi/2}(\lambda_0) \hbox{ and } u \in C^\infty(X; E).
\end{equation}
Assembling what we have achieved in \eqref{eq:Taubes_sA_plus_lambda_apriori_estimate_pure_L_2_part} (and now restricting to $\theta \in [\pi/2, \pi)$) and \eqref{eq:Taubes_sA_plus_lambda_apriori_estimate_pure_L_sharp_part} produces the bound,
\begin{multline}
\label{eq:Taubes_sA_plus_lambda_apriori_estimate_prelim}
|\lambda|\|u\|_{L^{2,\sharp}(X)} + |\lambda|^{1/2}\|\nabla u\|_{L^2(X)} + \|\nabla^2 u\|_{L^2(X)}
\leq C\|(\nabla^*\nabla + \lambda)u\|_{L^{2,\sharp}(X)},
\\
\forall\, \lambda \in \Delta_{\pi/2}(\lambda_0) \hbox{ and } u \in C^\infty(X; E).
\end{multline}
To complete the proof of the \apriori estimate \eqref{eq:Taubes_sA_plus_lambda_apriori_estimate} for $\sA = \nabla^*\nabla$, it remains to show that we can include in \eqref{eq:Taubes_sA_plus_lambda_apriori_estimate_prelim} the terms $\|u\|_{C(X)}$ and $\|\nabla^*\nabla u\|_{L^\sharp(X)}$ appearing in the left-hand side of \eqref{eq:Taubes_sA_plus_lambda_apriori_estimate}. But Lemma \ref{lem:Feehan_5-3} provides the inequality,
$$
\|u\|_{C(X)} \leq C\left(\|\nabla^*\nabla u\|_{L^\sharp(X)} + \|u\|_{L^\sharp(X)}\right), \quad\forall\, u \in C^\infty(X; E),
$$
and, for $\lambda_0 \geq C/2$, we can absorb the term $\|u\|_{L^\sharp(X)}$ in the right-hand side of the preceding inequality into the term $|\lambda|\|u\|_{L^\sharp(X)}$ in the left-hand side of the \apriori estimate \eqref{eq:Taubes_sA_plus_lambda_apriori_estimate}. Similarly, we can include $\|\nabla^*\nabla u\|_{L^\sharp(X)}$ in the left-hand side of \eqref{eq:Taubes_sA_plus_lambda_apriori_estimate} using the elementary inequality,
$$
\|\nabla^*\nabla u\|_{L^\sharp(X)} \leq \|(\nabla^*\nabla + \lambda)u\|_{L^\sharp(X)} + |\lambda|\|u\|_{L^\sharp(X)},
\quad\forall\, \lambda \in \CC \hbox{ and } u \in C^\infty(X; E),
$$
and applying the estimate \eqref{eq:Taubes_sA_plus_lambda_apriori_estimate_pure_L_sharp_part} to bound $|\lambda|\|u\|_{L^\sharp(X)}$ in terms of $\|(\nabla^*\nabla + \lambda)u\|_{L^\sharp(X)}$. Combining the preceding two observations with the inequality \eqref{eq:Taubes_sA_plus_lambda_apriori_estimate_prelim} yields the \apriori estimate \eqref{eq:Taubes_sA_plus_lambda_apriori_estimate} when $\sA = \nabla^*\nabla$.

Now suppose that $\sA$ is as allowed more generally in the hypotheses of our theorem. The principal symbols of $\nabla^*\nabla$ and $\sA$ are equal and thus
\begin{multline*}
\|(\nabla^*\nabla + \lambda)u\|_{L^{2,\sharp}(X))}
\leq
\|(\sA + \lambda)u\|_{L^{2,\sharp}(X)} + C\left(\|u\|_{L^{2,\sharp}(X)} + \|\nabla u\|_{L^{2,\sharp}(X)} \right),
\\
\forall\, \lambda \in \CC \hbox{ and } u \in C^\infty(X; E),
\end{multline*}
where $C$ may depend on the Riemannian metric on $X$ and the Hermitian connection and fiber metric on $E$. By choosing the constant $\lambda_0 \geq 0$ sufficiently large (for example, greater than $C/2$), the terms $\|u\|_{L^{2,\sharp}(X)}$ and $\|\nabla u\|_{L^2(X)}$ on the right-hand side of the preceding inequality maybe absorbed into the left-hand side of the \apriori estimate \eqref{eq:Taubes_sA_plus_lambda_apriori_estimate}, noting that the latter estimate is valid whenever $\lambda \in \CC$ obeys $\Real \lambda > \lambda_0$. Finally, we observe that, for any constant $\eps > 0$,
\begin{align*}
\|\nabla u\|_{L^\sharp(X)} &\leq C\|\nabla u\|_{L^{8/3}(X)}  \quad\hbox{(by Lemma \ref{lem:Feehan_4-1})}
\\
&\leq C\left(\eps^{-1}\|\nabla u\|_{L^2(X)} + \eps\|\nabla u\|_{L^4(X)}\right)  \quad\hbox{(by \cite[Equation (7.10)]{GilbargTrudinger})}
\\
&\leq C\left(\eps^{-1}\|\nabla u\|_{L^2(X)} + \eps\|\nabla^2 u\|_{L^2(X)}\right),
\end{align*}
where the final inequality follows from the Sobolev embedding $W^{1,2}(X; E) \hookrightarrow L^4(X; E)$ \cite[Theorem 4.12]{AdamsFournier}. To obtain the second inequality above, we use the interpolation inequality \cite[Equation (7.10)]{GilbargTrudinger} with interpolation parameter $\theta = 1/2$ and exponents $p = 2$, $q = 8/3$, and $r = 4$ obeying $p \leq q \leq r$ and related by $1/q = \theta/p + (1-\theta)/r$, so $1/q  = 1/4 + 1/8 = 3/8$, while $\mu$ in \cite[Equation (7.10)]{GilbargTrudinger} is given by $\mu = (1/p - 1/q)/(1/q - 1/r) = (1/2 - 3/8)/(3/8 - 1/4) = (1/8)/(1/8) = 1$. Hence, after relabeling the constant $\eps$, we obtain the interpolation inequality,
$$
\|\nabla u\|_{L^\sharp(X)} \leq C\eps^{-1}\|\nabla u\|_{L^2(X)} + \eps\|\nabla^2 u\|_{L^2(X)},
\quad\forall\, \eps > 0 \hbox{ and } u \in C^\infty(X; E).
$$
We now choose $\eps > 0$ small enough that the term $C\eps\|\nabla^2 u\|_{L^2(X)}$ can be absorbed into the left-hand side of the inequality \eqref{eq:Taubes_sA_plus_lambda_apriori_estimate}. We then choose $\lambda_0 \geq 0$ large enough that the terms $C\|u\|_{L^{2,\sharp}(X)}$ and $C\eps^{-1}\|\nabla u\|_{L^2(X)}$ can also be absorbed into the left-hand side of \eqref{eq:Taubes_sA_plus_lambda_apriori_estimate}. This completes the proof of the \apriori estimate \eqref{eq:Taubes_sA_plus_lambda_apriori_estimate}.

To establish existence and uniqueness of a solution $u \in W^{2, 2+\sharp}(X; E)$ to \eqref{eq:Taubes_sA_plus_lambda_equation}, given $f \in L^{2,\sharp}(X; E)$, recall first that Lemma \ref{lem:Feehan_4-1} provides an embedding $L^p(X; E) \hookrightarrow L^{2,\sharp}(X; E)$, for any $p > 2$, with dense range. For $\tilde f \in L^p(X; E)$ and $\lambda \in \CC$ with $\Real\lambda > \lambda_0$, Theorem \ref{thm:Krylov_Sobolev_lectures_8-5-3} provides a unique solution $\tilde u \in  W^{2,p}(X; E)$ to \eqref{eq:Taubes_sA_plus_lambda_equation} and because Lemma \ref{lem:Feehan_4-1} gives an embedding $W^{2,p}(X; E) \hookrightarrow W^{2, 2+\sharp}(X; E)$, we see that $\tilde u \in  W^{2, 2+\sharp}(X; E)$. Hence, the operator
\begin{equation}
\label{eq:A_plus_lambda_from_W2-2-sharp_to_L2-sharp}
\sA + \lambda: W^{2, 2+\sharp}(X; E) \to L^{2,\sharp}(X; E),
\end{equation}
has dense range for $\lambda \in \CC$ with $\Real\lambda > \lambda_0$.

The operator $\sA + \lambda$ in \eqref{eq:A_plus_lambda_from_W2-2-sharp_to_L2-sharp} is clearly continuous but also the \apriori estimate \eqref{eq:Taubes_sA_plus_lambda_apriori_estimate} means that it is bounded below in the sense of \cite[Definition 2.1]{Abramovich_Aliprantis_2002} and so has closed range by \cite[Theorem 2.5]{Abramovich_Aliprantis_2002}. In particular, the operator \eqref{eq:A_plus_lambda_from_W2-2-sharp_to_L2-sharp} is surjective and one-to-one by \cite[Lemma 2.8]{Abramovich_Aliprantis_2002}. Hence, given $f \in L^{2,\sharp}(X; E)$, there exists a unique solution $u \in W^{2, 2+\sharp}(X; E)$ to \eqref{eq:Taubes_sA_plus_lambda_equation} and this completes the proof of the theorem.
\end{proof}

Let $\sA$ be as in the hypotheses of Theorem \ref{thm:Taubes_existence_uniqueness_apriori_estimate_sA_plus_lambda} and, by analogy with Definition \ref{defn:Pazy_7-3-3}, set
\begin{align}
\label{eq:Taubes_domain_for_second_order_elliptic_operator}
\sD(\sA_{2,\sharp}) &:= W^{2; 2,\sharp}(X; E),
\\
\label{eq:Taubes_realization_for_second_order_elliptic_operator}
\sA_{2,\sharp} u &:= \sA(x, D)u, \quad \forall\, u \in \sD(\sA_{2,\sharp}).
\end{align}
We then have the

\begin{cor}[Sectorial property of a second-order elliptic partial differential operator with scalar principal symbol and generation of an analytic semigroup on $L^{2, \sharp}(X; E)$]
\label{cor:Sectorial_property_and_analytic_semigroup_on_Taubes_spaces}
Assume the hypotheses of Theorem \ref{thm:Taubes_existence_uniqueness_apriori_estimate_sA_plus_lambda}. Then there are constants, $C > 0$ and $\lambda_0 \geq 0$ and $\vartheta \in (\pi/2, \pi)$, such that $\rho(-\sA_{2,\sharp}) \supset \Delta_\vartheta(\lambda_0)$ and
\begin{equation}
\label{eq:Resolvent_estimate_Taubes_spaces}
\|(\sA_{2,\sharp} + \lambda)^{-1}\|_{\sL(L^{2,\sharp}(X; E))} \leq \frac{C}{|\lambda - \lambda_0|},
\quad \forall\, \lambda \in \Delta_\vartheta(\lambda_0).
\end{equation}
Moreover, $\sA_{2,\sharp}$ is a sectorial operator on $L^{2,\sharp}(X; E)$ in the sense of Definition \ref{defn:Sell_You_page_78_definition_of_sectorial_operator} and $-\sA_{2,\sharp}$ generates an analytic semigroup, $e^{-\sA_{2,\sharp} t}$, on $L^{2,\sharp}(X; E)$.
\end{cor}

\begin{proof}
The proof follows the same pattern as that of Corollary \ref{cor:Cannarsa_Terreni_Vespri_6-7_sectorial_operator_analytic_semigroup_LpOmega}, \emph{mutatis mutandis}, with the following observations:
\begin{enumerate}
\item The roles of Theorem \ref{thm:Cannarsa_Terreni_Vespri_6-7} (or Theorems \ref{thm:Cannarsa_Terreni_Vespri_6-7_one_lessthan_p_lessthan_two} or \ref{thm:Cannarsa_Terreni_Vespri_6-7_p_is_one}) are replaced by those of Theorem \ref{thm:Taubes_existence_uniqueness_apriori_estimate_sA_plus_lambda};

\item The domain $\sD(\sA_{2,\sharp}) = W^{2; 2, \sharp}(X; E)$ is dense in $L^{2, \sharp}(X; E)$ since $W^{2; 2, \sharp}(X; E)$ contains $C^\infty(X; E)$ and because $C^\infty(X; E)$ is dense in $L^{2, \sharp}(X; E)$.
\end{enumerate}
We note that the \apriori estimate \eqref{eq:Taubes_sA_plus_lambda_apriori_estimate} in Theorem \ref{thm:Taubes_existence_uniqueness_apriori_estimate_sA_plus_lambda} yields
$$
|\lambda|\|u\|_{L^{2,\sharp}(X)} \leq C\|(\sA + \lambda)u\|_{L^{2, \sharp}(X)},
\quad\forall\, \lambda \in \Delta_{\pi/2}(\lambda_0) \hbox{ and } u \in W^{2, 2+\sharp}(X; E).
$$
But $|\lambda - \lambda_0| \leq |\lambda| + |\lambda_0| \leq 2|\lambda|$ for all $\lambda \in \CC$ with $\Real \lambda > \lambda_0$ and thus
$$
|\lambda - \lambda_0|\|u\|_{L^{2,\sharp}(X)} \leq 2C\|(\sA + \lambda)u\|_{L^{2, \sharp}(X)},
\quad\forall\, \lambda \in \Delta_{\pi/2}(\lambda_0) \hbox{ and } u \in W^{2, 2+\sharp}(X; E).
$$
The remainder of the proof now follows just as in the proof of Corollary \ref{cor:Cannarsa_Terreni_Vespri_6-7_sectorial_operator_analytic_semigroup_LpOmega}.
\end{proof}

\subsection{Well-posedness for the linear heat equation in critical-exponent elliptic Sobolev spaces}
\label{subsec:Well-posedness_heat_equation_critical_exponent_elliptic_Sobolev_spaces}
Suppose that $A$ is a $C^\infty$ connection on a principal $G$-bundle $P$ over a $C^\infty$ closed, four-dimensional, Riemannian manifold $X$. One can now specialize the abstract development by Sell and You in \cite[Section 4.2]{Sell_You_2002} and summarized in our monograph in Section \ref{sec:Sell_You_4-2} by taking
\begin{subequations}
\label{eq:Standing_hypothesis_operator_A_Taubes}
\begin{align}
\label{eq:Standing_hypothesis_A_operator_cA_is_rough_Laplacian+1_on_Omega_1_adP_Taubes}
\sA &= \nabla_A^*\nabla_A + 1,
\\
\label{eq:Standing_hypothesis_operator_cA_choice_W_Taubes}
\cW &= L^{2+\sharp}(X; \Lambda^1\otimes\ad P),
\\
\label{eq:Standing_hypothesis_operator_cA_choice_V_Taubes}
\calV &\equiv \sD(\cA) = W^{2, 2+\sharp}(X; \Lambda^1\otimes\ad P).
\end{align}
\end{subequations}
According to Corollary \ref{cor:Sectorial_property_and_analytic_semigroup_on_Taubes_spaces}, the operator $\sA$ in \eqref{eq:Standing_hypothesis_A_operator_cA_is_rough_Laplacian+1_on_Omega_1_adP_Taubes} is sectorial on $\cW = L^{2+\sharp}(X; \Lambda^1\otimes\ad P)$ and the generator of analytic semigroup on $L^{2+\sharp}(X; \Lambda^1\otimes\ad P)$. Clearly, the operator $\sA$ in \eqref{eq:Standing_hypothesis_A_operator_cA_is_rough_Laplacian+1_on_Omega_1_adP_Taubes} is also positive and hence fulfills Hypothesis \ref{hyp:Sell_You_4_standing_hypothesis_A} with $\calV \equiv \sD(\cA) = W^{2, 2+\sharp}(X; \Lambda^1\otimes\ad P)$, noting that the domain $\sD(\cA)$ of the smallest closed extension of \eqref{eq:Standing_hypothesis_A_operator_cA_is_rough_Laplacian+1_on_Omega_1_adP_Taubes} of $\sA$ is identified by the \apriori estimate \eqref{eq:Taubes_sA_plus_lambda_apriori_estimate} in Theorem \ref{thm:Taubes_existence_uniqueness_apriori_estimate_sA_plus_lambda}. Consequently, depending on our choice of temporal regularity and any additional choice of spatial regularity for $f$ in equation \eqref{eq:Linear_heat_equation_with_rough_Laplacian_plus_one_on_Omega_1_adP}
and whether or not any additional regularity condition is imposed on the initial data $a_0$, the conclusions of the abstract Theorems \ref{thm:Sell_You_42-9} and \ref{thm:Sell_You_42-10} on existence, uniqueness, and regularity of solutions to the linear heat equation \eqref{eq:Linear_heat_equation_with_rough_Laplacian_plus_one_on_Omega_1_adP}
now hold for the choice of setup in \eqref{eq:Standing_hypothesis_operator_A_Taubes}.

\chapter[Solutions to the Yang-Mills heat equation]{Existence, uniqueness, and regularity of solutions to the Yang-Mills heat equation}
\label{chapter:Existence_uniqueness_regularity_Yang-Mills_heat_equation}

\section[Weak solutions to the heat equation]{Weak solutions to the heat equation on Sobolev spaces of sections of a vector bundle}
\label{sec:Sell_You_4-2-3_heat_equation_vector_bundle}
In this section, we wish to apply the results of Section \ref{subsec:Sell_You_4-2-3} to the heat equation,
\begin{equation}
\label{eq:Sell_You_42-1_heat_equation_connection_Laplacian}
\frac{\partial v}{\partial t} + \nabla_A^*\nabla_Av = h,
\end{equation}
where $A$ is a connection on a principle $G$-bundle, $P$, over a closed, Riemannian, smooth manifold, $X$, of dimension $d \geq 2$ and $\nabla_A$ is the covariant derivative on an associated Riemannian vector bundle, $V = P\times_\rho \VV$, defined by a representation
$\rho:G \to \End\VV$ for a compact Lie group, $G$, and $v \in C^\infty((0,T)\times X; V)$ while $f \in C^\infty((0,T)\times X; V)$ is a given source function.

\subsection{Existence, uniqueness, and \apriori estimates for weak solutions to the heat equation}
\label{sec:Sell_You_4-2-3_heat_equation_vector_bundle_existence_uniqueness_apriori_estimates}
The connection Laplace operator, $\nabla_A^*\nabla_A$, is self-adjoint on $L^2(X;\Lambda^1\otimes\ad P)$ but not positive, so we instead consider the augmented connection Laplace operator, $\sA := \nabla_A^*\nabla_A + 1$, and note that, setting $v(t) =: e^{t}u(t)$ for $t \geq 0$, Equation \eqref{eq:Sell_You_42-1_heat_equation_connection_Laplacian} for $v$ is equivalent to the following equation
\begin{equation}
\label{eq:Sell_You_42-1_heat_equation}
\frac{\partial u}{\partial t} + (\nabla_A^*\nabla_A + 1)u = f,
\end{equation}
for $u$, where $f(t) := e^{-t}h(t)$ for $t \geq 0$. We set $\sH := L^2(X;V)$ and $\calV^2 := H_A^2(X; V) = \sD(\nabla_A^*\nabla_A+1)$, so that $\|v\|_{H^s_A(X)} = \|(\nabla_A^*\nabla_A + 1)^{s/2}v\|_{L^2(X)}$, for all $v \in H^s_A(X;V)$ and $s\in\RR$.
We obtain the following special cases of Lemma \ref{lem:Sell_You_42-11} with $\alpha = 0$ and $1$, respectively.

\begin{lem}[Properties of weak solutions in the space $L^2(X; V)$]
\label{lem:Sell_You_42-11_heat_equation_alpha_is_zero}
Let $u = u(t)$ be a weak solution of equation \eqref{eq:Sell_You_42-1_heat_equation} in $L^2(X; V)$, on the interval $[0, T)$, where $0< T \leq \infty$, and
$$
u_0 \in L^2(X; V) \quad\hbox{and}\quad f \in L_{\loc}^p([0, T); H_A^{-1}(X; V)), \quad\hbox{for } 2 \leq p \leq \infty.
$$
Then the following properties are valid:
\begin{enumerate}
\item The time derivative satisfies $\partial_t u \in L_{\loc}^2([0, T); H_A^{-1}(X; V))$ and $u$ satisfies
$$
u \in C([0, T); L^2(X; V)) \cap C^{0,\frac{1}{2}}_{\loc}([0, T); H_A^{-1}(X; V)).
$$
\item $u$ obeys Equation \eqref{eq:Sell_You_42-1_heat_equation} a.e. on $(0, T)$ in the space $H_A^{-1}(X; V)$, and $u$ is a mild solution of Equation \eqref{eq:Sell_You_42-1_heat_equation} in $H^\beta_A(X; V)$, for each $\beta < (p - 2)p^{-1}$.
\end{enumerate}
\end{lem}

\begin{lem}[Properties of weak solutions in the space $H_A^1(X; V)$]
\label{lem:Sell_You_42-11_heat_equation_alpha_is_one}
Let $u = u(t)$ be a weak solution of equation \eqref{eq:Sell_You_42-1_heat_equation} in $H_A^1(X; V)$, on the interval $[0, T)$, where $0< T \leq \infty$, and
$$
u_0 \in H_A^1(X; V) \quad\hbox{and}\quad f \in L_{\loc}^p([0, T); L^2(X;V)), \quad\hbox{for } 2 \leq p \leq \infty.
$$
Then the following properties are valid:
\begin{enumerate}
\item The time derivative satisfies $\partial_t u \in L_{\loc}^2([0, T); L^2(X;V))$ and $u$ satisfies
$$
u \in C([0, T); H_A^1(X; V)) \cap C^{0,\frac{1}{2}}_{\loc}([0, T); L^2(X;V)).
$$
\item $u$ obeys Equation \eqref{eq:Sell_You_42-1_heat_equation} a.e. on $(0, T)$ in the space $L^2(X;V)$, and $u$ is a mild solution of Equation \eqref{eq:Sell_You_42-1_heat_equation} in $H^\beta_A(X; V)$, for each $\beta < 1 + (p - 2)p^{-1}$.
\end{enumerate}
\end{lem}

Although Lemmata \ref{lem:Sell_You_42-11_heat_equation_alpha_is_zero} and \ref{lem:Sell_You_42-11_heat_equation_alpha_is_one} comprise the cases we shall appeal to most frequently, it is useful for later reference, to record the general case of arbitrary $\alpha \in \RR$.

\begin{lem}[Properties of weak solutions in the space $H_A^\alpha(X; V)$ for $\alpha \in \RR$]
\label{lem:Sell_You_42-11_heat_equation_alpha_is_real}
For $\alpha \in \RR$, let $u = u(t)$ be a weak solution of Equation \eqref{eq:Sell_You_42-1_heat_equation} in $H_A^\alpha(X; V)$ on the interval $[0, T)$, where $0< T \leq \infty$, and
$$
u_0 \in H_A^\alpha(X; V) \quad\hbox{and}\quad f \in L_{\loc}^p([0, T); H_A^{\alpha-1}(X; V)), \quad\hbox{for } 2 \leq p \leq \infty.
$$
Then the following properties are valid:
\begin{enumerate}
\item The time derivative satisfies $\partial_t u \in L_{\loc}^2([0, T); H_A^{\alpha-1}(X; V))$ and $u$ satisfies
$$
u \in C([0, T); H_A^\alpha(X; V)) \cap C^{0,\frac{1}{2}}_{\loc}([0, T); H_A^{\alpha-1}(X; V)).
$$
\item One has $\partial u_t + \cA u = f$ a.e. on $(0, T)$ in the space $H_A^{\alpha-1}(X; V)$, and $u$ is a mild solution
of Equation \eqref{eq:Sell_You_42-1} in $H_A^\beta(X; V)$, for each $\beta < \alpha + (p - 2)p^{-1}$.
\end{enumerate}
\end{lem}

\begin{proof}[Proofs of Lemmata \ref{lem:Sell_You_42-11_heat_equation_alpha_is_zero},\ref{lem:Sell_You_42-11_heat_equation_alpha_is_one}, and \ref{lem:Sell_You_42-11_heat_equation_alpha_is_real}]
We apply Lemma \ref{lem:Sell_You_42-11} to Equation \eqref{eq:Sell_You_42-1_heat_equation} with $\sA := \nabla_A^*\nabla_A + 1$, and $\alpha = 0$ or $1$ or arbitrary $\alpha \in \RR$, respectively.
\end{proof}

For $\alpha = 0$ or $1$ or arbitrary $\alpha \in \RR$, respectively, Theorem \ref{thm:Sell_You_42-12} also specializes to

\begin{thm}[Existence, uniqueness, and regularity of weak solutions to the heat equation in the space $L^2(X; V)$]
\label{thm:Sell_You_42-12_heat_equation_alpha_is_zero}
Let
$$
u_0 \in L^2(X; V) \quad\hbox{and}\quad f \in L_{\loc}^p([0, T); H_A^{-1}(X;V)), \quad\hbox{for } 2 \leq p \leq \infty.
$$
Then the following properties hold:
\begin{enumerate}
\item
\label{item:Theorem_Sell_You_42-12_1_heat_equation_alpha_is_zero}
There is a unique weak solution $u = u(t)$ of Equation \eqref{eq:Sell_You_42-1_heat_equation}, in the space $L^2(X; V)$, on the interval $[0, T)$, with $u(0) = u_0$,
\begin{equation}
\label{eq:Sell_You_42-29_heat_equation}
\begin{gathered}
\frac{\partial u}{\partial t} \in L_{\loc}^2([0, T); H_A^{-1}(X;V)), \quad\hbox{and}
\\
u \in C([0, T); H^\beta_A(X; V)) \cap C^{0,\theta_1}_{\loc}([0, T); H^\sigma_A(X; V)) \cap L_{\loc}^2([0, T); H_A^1(X; V)),
\end{gathered}
\end{equation}
for each $\beta$ and $\sigma$, where $\beta \leq 0$, and $\sigma < 0$, and $\theta_1 = \theta_1(\sigma) > 0$. Furthermore, $u$ is a mild solution of Equation \eqref{eq:Sell_You_42-1} in $H^\beta_A(X; V)$, for each $\beta < (p - 2)p^{-1}$.

\item
\label{item:Theorem_Sell_You_42-12_2_heat_equation_alpha_is_zero}
For $0 \leq t < T$, the following inequalities are valid in the space $L^2(X; V)$:
\begin{subequations}
\label{eq:Sell_You_42-30_heat_equation_alpha_is_zero}
\begin{align}
\label{eq:Sell_You_42-30_heat_equation_Linfty-L2_u}
\|u(t)\|_{L^2(X)}^2
&\leq
\|u_0\|_{L^2(X)}^2 + \int_0^t \|f(s)\|_{H_A^{-1}(X)}^2 \,ds,
\\
\label{eq:Sell_You_42-30_heat_equation_L2-H1_u}
\int_0^t e^{-2s} \|u(s)\|_{H_A^1(X)}^2 \,ds
&\leq
\|u_0\|_{L^2(X)}^2 + \int_0^t e^{-2s} \|f(s)\|_{H_A^{-1}(X)}^2 \,ds,
\\
\label{eq:Sell_You_42-30_heat_equation_L2-Hminus1_dudt}
\int_0^t\|e^{-2s} \partial_tu(s)\|_{H_A^{-1}(X)}^2 \,ds
&\leq
4\|u_0\|_{L^2(X)}^2 + 8\int_0^t e^{-2s} \|f(s)\|_{H_A^{-1}(X)}^2 \,ds
\\
\notag
&\quad + 2\int_0^t e^{-2s}\|u(s)\|_{H_A^{-1}(X)}^2 \,ds.
\end{align}
\end{subequations}
\item
\label{item:Theorem_Sell_You_42-12_3_heat_equation_alpha_is_zero}
If $p > 2$, then $u$ is a mild solution of Equation \eqref{eq:Sell_You_42-1} in $L^2(X;V)$, and in addition to \eqref{eq:Sell_You_42-29_heat_equation}, one has
\begin{equation}
\label{eq:Sell_You_42-31_heat_equation_alpha_is_zero}
C^{0,\theta_2}_{\loc}((0, T); H^\sigma_A(X;V)),
\end{equation}
for each $\sigma$, where $0 \leq \sigma < 1 - 2/p$ and $\theta_2 = \theta_2(\sigma) > 0$.
\end{enumerate}
\end{thm}

\begin{thm}[Existence, uniqueness, and regularity of weak solutions to the heat equation in the space $H_A^1(X; V)$]
\label{thm:Sell_You_42-12_heat_equation_alpha_is_one}
Let
$$
u_0 \in H_A^1(X; V) \quad\hbox{and}\quad f \in L_{\loc}^p([0, T); L^2(X;V)), \quad\hbox{for } 2 \leq p \leq \infty.
$$
Then the following properties hold:
\begin{enumerate}
\item
\label{item:Theorem_Sell_You_42-12_1_heat_equation_alpha_is_one}
There is a unique weak solution $u = u(t)$ of Equation \eqref{eq:Sell_You_42-1_heat_equation}, in the space $H_A^1(X; V)$, on the interval $[0, T)$, with $u(0) = u_0$,
\begin{equation}
\label{eq:Sell_You_42-29_heat_equation_alpha_is_one}
\begin{gathered}
\frac{\partial u}{\partial t} \in L_{\loc}^2([0, T); L^2(X;V)), \quad\hbox{and}
\\
u \in C([0, T); H^\beta_A(X; V)) \cap C^{0,\theta_1}_{\loc}([0, T); H^\sigma_A(X; V)) \cap L_{\loc}^2([0, T); H_A^2(X; V)),
\end{gathered}
\end{equation}
for each $\beta$ and $\sigma$, where $\beta \leq 1$, and $\sigma < 1$, and $\theta_1 = \theta_1(\sigma) > 0$. Furthermore, $u$ is a mild solution of Equation \eqref{eq:Sell_You_42-1} in $H^\beta_A(X; V)$, for each $\beta < 1 + (p - 2)p^{-1}$.

\item
\label{item:Theorem_Sell_You_42-12_2_heat_equation_alpha_is_one}
For $0 \leq t < T$, the following inequalities are valid in the space $H_A^1(X; V)$:
\begin{subequations}
\label{eq:Sell_You_42-30_heat_equation_alpha_is_one}
\begin{align}
\label{eq:Sell_You_42-30_heat_equation_Linfty-H1_u}
\|u(t)\|_{H_A^1(X)}^2
&\leq
\|u_0\|_{H_A^1(X)}^2 + \int_0^t \|f(s)\|_{L^2(X)}^2 \,ds,
\\
\label{eq:Sell_You_42-30_heat_equation_L2-H2_u}
\int_0^t e^{-2s} \|u(s)\|_{H_A^2(X)}^2 \,ds
&\leq
\|u_0\|_{H_A^1(X)}^2 + \int_0^t e^{-2s} \|f(s)\|_{L^2(X)}^2 \,ds,
\\
\label{eq:Sell_You_42-30_heat_equation_Linfty-L2_dudt}
\int_0^t \|e^{-2s} \partial_tu(s)\|_{L^2(X)}^2 \,ds
&\leq
4\|u_0\|_{H_A^1(X)}^2 + 8\int_0^t e^{-2s} \|f(s)\|_{L^2(X)}^2 \,ds
\\
\notag
&\quad + 2\int_0^t e^{-2s} \|u(s)\|_{L^2(X;V)}^2 \,ds.
\end{align}
\end{subequations}
\item
\label{item:Theorem_Sell_You_42-12_3_heat_equation_alpha_is_one}
If $p > 2$, then $u$ is a mild solution of Equation \eqref{eq:Sell_You_42-1} in $H_A^1(X;V)$, and in addition to \eqref{eq:Sell_You_42-29_heat_equation}, one has
\begin{equation}
\label{eq:Sell_You_42-31_heat_equation_alpha_is_one}
C^{0,\theta_2}_{\loc}((0, T); H^\sigma_A(X;V)),
\end{equation}
for each $\sigma$, where $1 \leq \sigma < 2 - 2/p$ and $\theta_2 = \theta_2(\sigma) > 0$.
\end{enumerate}
\end{thm}


\begin{thm}[Existence, uniqueness, and regularity of weak solutions to a linear evolutionary equation in the space $H_A^\alpha(X; V)$ for $\alpha \in \RR$]
\label{thm:Sell_You_42-12_heat_equation_alpha_is_real}
For $\alpha \in \RR$, let
$$
u_0 \in H_A^\alpha(X;V), \quad\hbox{and}\quad f \in L_{\loc}^p([0, T); H_A^{\alpha-1}(X;V)),
\quad\hbox{for } 2 \leq p \leq \infty.
$$
Then the following properties hold:
\begin{enumerate}
\item
\label{item:Theorem_Sell_You_42-12_1_heat_equation_alpha_is_real}
There is a unique weak solution $u = u(t)$ of Equation \eqref{eq:Sell_You_42-1_heat_equation}, in the space $H_A^\alpha(X;V)$, on the interval $[0, T)$, with $u(0) = u_0$,
\begin{equation}
\label{eq:Sell_You_42-29_heat_equation_alpha_is_real}
\begin{gathered}
\frac{\partial u}{\partial t} \in L_{\loc}^2([0, T); H_A^{\alpha-1}(X;V)), \quad\hbox{and}
\\
u \in C([0, T); H_A^\beta(X;V)) \cap C^{0,\theta_1}_{\loc}([0, T); H_A^\sigma(X;V))
\cap L_{\loc}^2([0, T); H_A^{\alpha+1}(X;V)),
\end{gathered}
\end{equation}
for each $\beta$ and $\sigma$, where $\beta \leq \alpha$, and $\sigma < \alpha$, and $\theta_1 = \theta_1(\sigma) > 0$. Furthermore, $u$ is a mild solution of Equation \eqref{eq:Sell_You_42-1_heat_equation} in $H_A^\beta(X;V)$, for each $\beta < \alpha + (p - 2)p^{-1}$.

\item
\label{item:Theorem_Sell_You_42-12_2_heat_equation_alpha_is_real}
For $0 \leq t < T$, the following inequalities are valid in the space $H_A^\alpha(X;V)$:
\begin{subequations}
\label{eq:Sell_You_42-30_heat_equation_alpha_is_real}
\begin{align}
\label{eq:Sell_You_42-30_Linfty_time_Halpha_space_u}
\|u(t)\|_{H_A^\alpha(X)}^2
&\leq
e^{-\lambda_1 t}\|u_0\|_{H_A^\alpha(X)}^2 + \int_0^t e^{-\lambda_1(t-s)}\|f(s)\|_{H_A^{\alpha-1}(X)}^2 \,ds,
\\
\label{eq:Sell_You_42-30_L2_time_Halpha+1_space_u}
\int_0^t\|u(s)\|_{H_A^{\alpha+1}(X)}^2 \,ds
&\leq
\|u_0\|_{H_A^\alpha(X)}^2 + \int_0^t \|f(s)\|_{H_A^{\alpha-1}(X)}^2 \,ds,
\\
\label{eq:Sell_You_42-30_L2_time_Halpha-1_space_partial_t_u}
\int_0^t\|\partial_tu(s)\|_{H_A^{\alpha-1}(X)}^2 \,ds
&\leq
2\|u_0\|_{H_A^\alpha(X)}^2 + 4\int_0^t \|f(s)\|_{H_A^{\alpha-1}(X)}^2 \,ds.
\end{align}
\end{subequations}
\item
\label{item:Theorem_Sell_You_42-12_3_heat_equation_alpha_is_real}
If $p > 2$, then $u$ is a mild solution of Equation \eqref{eq:Sell_You_42-1_heat_equation} in $H_A^\alpha(X)$, and in addition to \eqref{eq:Sell_You_42-29_heat_equation_alpha_is_real}, one has
\begin{equation}
\label{eq:Sell_You_42-31_heat_equation_alpha_is_real}
C^{0,\theta_2}_{\loc}((0, T); H_A^\sigma(X)),
\end{equation}
for each $\sigma$, where $\alpha \leq \sigma < \alpha + 1 - 2/p$ and $\theta_2 = \theta_2(\sigma) > 0$.
\end{enumerate}
\end{thm}

\begin{proof}[Proof of Theorems \ref{thm:Sell_You_42-12_heat_equation_alpha_is_zero}, \ref{thm:Sell_You_42-12_heat_equation_alpha_is_one}, and \ref{thm:Sell_You_42-12_heat_equation_alpha_is_real}]
It suffices to consider the case $\alpha=1$ as the cases $\alpha = 0$ and $\alpha \in \RR$ arbitrary are similar. We apply Theorem \ref{thm:Sell_You_42-12} to Equation \eqref{eq:Sell_You_42-1_heat_equation} with $\sA := \nabla_A^*\nabla_A + 1$, and $\alpha = 1$ and note the fact that $\|w\|_{H^s_A(X)} = \|(\nabla_A^*\nabla_A + 1)^{s/2}w\|_{L^2(X)}$, for all $w \in H^s_A(X;V)$ and $s \in \RR$.
The inequalities \eqref{eq:Sell_You_42-30_Linfty_time_Valpha_space_u} and \eqref{eq:Sell_You_42-30_L2_time_Valpha+1_space_u} (with $\alpha = 1$ and $\lambda_1 = 1$) yield
\begin{align*}
\|v(t)\|_{H_A^1(X;V)}^2
&\leq
e^{-t}\|v_0\|_{H_A^1(X;V)}^2 + \int_0^t e^{-(t-s)}\|f_1(s)\|_{L^2(X;V)}^2 \,ds,
\\
\int_0^t\|v(s)\|_{H_A^2(X;V)}^2 \,ds
&\leq
\|v_0\|_{H_A^1(X;V)}^2 + \int_0^t \|f_1(s)\|_{L^2(X;V)}^2 \,ds.
\end{align*}
We now substitute $v(t) = e^{-t}u(t)$ and $f_1(t) = e^{-t}f(t)$ in the preceding inequalities to give \eqref{eq:Sell_You_42-30_heat_equation_Linfty-H1_u} and \eqref{eq:Sell_You_42-30_heat_equation_L2-H2_u}.

The inequality \eqref{eq:Sell_You_42-30_L2_time_Valpha-1_space_partial_t_u} (with $\alpha = 1$ and $\lambda_1 = 1$) implies that
$$
\int_0^t\|\partial_tv(s)\|_{L^2(X;V)}^2 \,ds
\leq
2\|v_0\|_{H_A^1(X;V)}^2 + 4\int_0^t \|f_1(s)\|_{L^2(X;V)}^2 \,ds.
$$
Using $\partial_t v = e^{-t}(\partial_t u - u) = e^{-t}\partial_t u - v$ and
$$
e^{-t}\|\partial_tu(t)\|_{L^2(X;V)} \leq \|\partial_tv(t)\|_{L^2(X;V)} + \|v(t)\|_{L^2(X;V)},
$$
so that
$$
e^{-2t}\|\partial_tu(t)\|_{L^2(X;V)}^2 \leq 2\|\partial_tv(t)\|_{L^2(X;V)}^2 + 2\|v(t)\|_{L^2(X;V)}^2.
$$
Combining the preceding inequalities gives \eqref{eq:Sell_You_42-30_heat_equation_Linfty-L2_dudt}.
\end{proof}

Theorems \ref{thm:Sell_You_42-12_heat_equation_alpha_is_zero}, \ref{thm:Sell_You_42-12_heat_equation_alpha_is_one}, and \ref{thm:Sell_You_42-12_heat_equation_alpha_is_real} give the following useful \apriori estimates which one can also derive by other methods that we shall describe in the sequel.

\begin{cor}[\Apriori estimates for weak solutions to the heat equation in the space $L^2(X; V)$]
\label{cor:Sell_You_42-12_heat_equation_alpha_is_zero_apriori_estimates}
Let $u_0 \in L^2(X; V)$ and $f \in L^2(0, T; H_A^{-1}(X;V))$ and $u = u(t)$ be a weak solution of Equation \eqref{eq:Sell_You_42-1_heat_equation}, in the space $L^2(X; V)$, on the interval $[0, T)$, with $u(0) = u_0$. Then the following \apriori estimates hold:
\begin{subequations}
\label{eq:Sell_You_42-30_heat_equation_alpha_is_zero_compact_form}
\begin{align}
\label{eq:Sell_You_42-30_heat_equation_Linfty-L2_u_compact_form}
\|u\|_{L^\infty(0,T;L^2(X))} &\leq \|u_0\|_{L^2(X)} + \|f\|_{L^2(0,T;H_A^{-1}(X))},
\\
\label{eq:Sell_You_42-30_heat_equation_L2-H1_u_compact_form}
\|e^{-2t}u\|_{L^2(0,T;H_A^1(X))} &\leq \|u_0\|_{L^2(X)} + \|e^{-2t}f\|_{L^2(0,T;H_A^{-1}(X))},
\\
\label{eq:Sell_You_42-30_heat_equation_L2-Hminus1_dudt_compact_form}
\|e^{-2t}\partial_tu\|_{L^2(0,T;H_A^{-1}(X))} &\leq 2\|u_0\|_{L^2(X)} + 2\sqrt{2}\|e^{-2t}f\|_{L^2(0,T;H_A^{-1}(X))}
\\
\notag
&\quad + \sqrt{2}\|e^{-2t}u\|_{L^2(0,T;H_A^{-1}(X))},
\end{align}
\end{subequations}
and, for $T < \infty$,
\begin{equation}
\label{eq:Sell_You_42-30_heat_equation_L2-Hminus1_dudt_compact_form_finite_T0}
\|\partial_tu\|_{L^2(0,T;H_A^{-1}(X))}
\leq
\left(2 + \sqrt{2T}\right) e^{2T} \|u_0\|_{L^2(X)} + \left(2\sqrt{2} + \sqrt{2T}\right) e^{2T} \|f\|_{L^2(0,T;H_A^{-1}(X))}.
\end{equation}
\end{cor}

\begin{cor}[\Apriori estimates for weak solutions to the heat equation in the space $H_A^1(X; V)$]
\label{cor:Sell_You_42-12_heat_equation_alpha_is_one_apriori_estimates}
Let $u_0 \in H_A^1(X; V)$ and $f \in L^2(0, T; L^2(X;V))$ and $u = u(t)$ be a weak solution of Equation \eqref{eq:Sell_You_42-1_heat_equation}, in the space $H_A^1(X; V)$, on the interval $[0, T)$, with $u(0) = u_0$. Then the following \apriori estimates hold:
\begin{subequations}
\label{eq:Sell_You_42-30_heat_equation_alpha_is_one_compact_form}
\begin{align}
\label{eq:Sell_You_42-30_heat_equation_Linfty-H1_u_compact_form}
\|u\|_{L^\infty(0,T;H_A^1(X))} &\leq \|u_0\|_{H_A^1(X)} + \|f\|_{L^2(0,T;L^2(X))},
\\
\label{eq:Sell_You_42-30_heat_equation_L2-H2_u_compact_form}
\|e^{-2t}u\|_{L^2(0,T;H_A^2(X))} &\leq \|u_0\|_{H_A^1(X)} + \|e^{-2t}f\|_{L^2(0,T;L^2(X))},
\\
\label{eq:Sell_You_42-30_heat_equation_L2-L2_dudt_compact_form}
\|e^{-2t}\partial_tu\|_{L^2(0,T;L^2(X))} &\leq 2\|u_0\|_{H_A^1(X)} + 2\sqrt{2}\|e^{-2t}f\|_{L^2(0,T;L^2(X))}
\\
\notag
&\quad + \sqrt{2}\|e^{-2t}u\|_{L^2(0,T;L^2(X))},
\end{align}
\end{subequations}
and, for $T < \infty$,
\begin{equation}
\label{eq:Sell_You_42-30_heat_equation_L2-L2_dudt_compact_form_finite_T0}
\|\partial_tu\|_{L^2(0,T;L^2(X))}
\leq
\left(2 + \sqrt{2T}\right) e^{2T} \|u_0\|_{H_A^1(X)} + \left(2\sqrt{2} + \sqrt{2T}\right) e^{2T} \|f\|_{L^2(0,T;L^2(X))}.
\end{equation}
\end{cor}

\begin{cor}[\Apriori estimates for weak solutions to the heat equation in the space $H_A^\alpha(X; V)$ for $\alpha \in \RR$]
\label{cor:Sell_You_42-12_heat_equation_alpha_is_real_apriori_estimates}
For $\alpha \in \RR$, let $u_0 \in H_A^\alpha(X; V)$ and $f \in L^2(0, T; H_A^{\alpha-1}(X;V))$ and $u = u(t)$ be a weak solution of Equation \eqref{eq:Sell_You_42-1_heat_equation}, in the space $H_A^\alpha(X; V)$, on the interval $[0, T)$, with $u(0) = u_0$. Then the following \apriori estimates hold:
\begin{subequations}
\label{eq:Sell_You_42-30_heat_equation_alpha_is_real_compact_form}
\begin{align}
\label{eq:Sell_You_42-30_heat_equation_Linfty-Halpha_u_compact_form}
\|u\|_{L^\infty(0,T;H_A^\alpha(X))} &\leq \|u_0\|_{H_A^\alpha(X)} + \|f\|_{L^2(0,T;H_A^{\alpha-1}(X))},
\\
\label{eq:Sell_You_42-30_heat_equation_L2-Halpha+1_u_compact_form}
\|e^{-2t}u\|_{L^2(0,T;H_A^{\alpha+1}(X))} &\leq \|u_0\|_{H_A^\alpha(X)} + \|e^{-2t}f\|_{L^2(0,T;H_A^{\alpha-1}(X))},
\\
\label{eq:Sell_You_42-30_heat_equation_L2-Halpha-1_dudt_compact_form}
\|e^{-2t}\partial_tu\|_{L^2(0,T;H_A^{\alpha-1}(X))} &\leq 2\|u_0\|_{H_A^\alpha(X)}
+ 2\sqrt{2}\|e^{-2t}f\|_{L^2(0,T;H_A^{\alpha-1}(X))}
\\
\notag
&\quad + \sqrt{2}\|e^{-2t}u\|_{L^2(0,T;H_A^{\alpha-1}(X))},
\end{align}
\end{subequations}
and, for $T < \infty$,
\begin{multline}
\label{eq:Sell_You_42-30_heat_equation_L2-Halpha-1_dudt_compact_form_finite_T0}
\|\partial_tu\|_{L^2(0,T;H_A^{\alpha-1}(X))}
\\
\leq \left(2 + \sqrt{2T}\right) e^{2T} \|u_0\|_{H_A^\alpha(X)} + \left(2\sqrt{2} + \sqrt{2T}\right) e^{2T}
\|f\|_{L^2(0,T;H_A^{\alpha-1}(X))}.
\end{multline}
\end{cor}

\begin{proof}[Proofs of Corollaries \ref{cor:Sell_You_42-12_heat_equation_alpha_is_zero_apriori_estimates}, \ref{cor:Sell_You_42-12_heat_equation_alpha_is_one_apriori_estimates}, and \ref{cor:Sell_You_42-12_heat_equation_alpha_is_real_apriori_estimates}]
It suffices to consider the case $\alpha=1$, since the discussion is identical for the cases $\alpha=0$ or arbitrary $\alpha\in\RR$. The inequalities \eqref{eq:Sell_You_42-30_heat_equation_Linfty-H1_u_compact_form}, \eqref{eq:Sell_You_42-30_heat_equation_L2-H2_u_compact_form}, and \eqref{eq:Sell_You_42-30_heat_equation_L2-L2_dudt_compact_form} simply rewrite those in \eqref{eq:Sell_You_42-30_heat_equation_alpha_is_one} in more compact form. As in the proof of \cite[Theorem 37.6]{Sell_You_2002}, the spectral theory for positive self-adjoint operators on a Hilbert space \cite{Rudin} yields, for all $s \in \RR$,
$$
(1 + \nabla_A^*\nabla_A)^s v = \int_0^\infty (1 + \lambda)^s \,dE(\lambda)v,
\quad\forall\, v \in \sD\left((1 + \nabla_A^*\nabla_A)^s\right),
$$
and, for $s_1 \leq s_2$ and $v \in \sD\left((1 + \nabla_A^*\nabla_A)^{s_2}\right)$,
$$
\|v\|_{H_A^{s_1}(X)}^2 = \int_0^\infty (1 + \lambda)^{2s_1} \,\|dE(\lambda)v\|_{L^2(X)}^2
\leq
\int_0^\infty (1 + \lambda)^{2s_2} \,\|dE(\lambda)v\|_{L^2(X)}^2 = \|v\|_{H_A^{s_2}(X)}^2,
$$
that is,
\begin{equation}
\label{eq:Hilbert_space_scale_monotonicity}
\|v\|_{H_A^{s_1}(X)} \leq \|v\|_{H_A^{s_2}(X)},
\quad\forall\, s_1 \leq s_2 \hbox{ and } v \in \sD\left((1 + \nabla_A^*\nabla_A)^{s_2}\right).
\end{equation}
We then observe that
\begin{align*}
\|e^{-2t}u\|_{L^2(0,T;L^2(X))} &\leq \sqrt{T} \|u\|_{L^\infty(0,T;L^2(X))}
\\
&\leq \sqrt{T}\left( \|u_0\|_{L^2(X)} + \|f\|_{L^2(0,T;H_A^{-1}(X))} \right)
\quad\hbox{(by \eqref{eq:Sell_You_42-30_heat_equation_Linfty-L2_u_compact_form})}
\\
&\leq \sqrt{T}\left( \|u_0\|_{L^2(X)} + \|f\|_{L^2(0,T;L^2(X))} \right)
\quad\hbox{(by \eqref{eq:Hilbert_space_scale_monotonicity})}.
\end{align*}
Combining the preceding inequality with \eqref{eq:Sell_You_42-30_heat_equation_L2-L2_dudt_compact_form} and
$$
e^{-2T}\|\partial_tu\|_{L^2(0,T;L^2(X))} \leq \|e^{-2t}\partial_tu\|_{L^2(0,T;L^2(X))}
$$
yields \eqref{eq:Sell_You_42-30_heat_equation_L2-L2_dudt_compact_form_finite_T0}.
\end{proof}

\begin{rmk}[On the constants in the \apriori estimates in Corollaries \ref{cor:Sell_You_42-12_heat_equation_alpha_is_zero_apriori_estimates}, \ref{cor:Sell_You_42-12_heat_equation_alpha_is_one_apriori_estimates}, and \ref{cor:Sell_You_42-12_heat_equation_alpha_is_real_apriori_estimates}]
\label{rmk:Sell_You_42-12_heat_equation_alpha_is_zero_or_one_or_real_apriori_estimate_constant_dependencies}
The fact that the constants appearing in the \apriori estimates in Corollaries \ref{cor:Sell_You_42-12_heat_equation_alpha_is_zero_apriori_estimates}, \ref{cor:Sell_You_42-12_heat_equation_alpha_is_one_apriori_estimates}, and
\ref{cor:Sell_You_42-12_heat_equation_alpha_is_real_apriori_estimates}
have simple, universal values is due to our choice of norms for the Hilbert spaces, $H_A^s(X; V)$, namely $\|v\|_{H_A^s(X; V)} = \|(1 + \nabla_A^*\nabla_A)^s v\|_{L^2(X)}$ for real $s \geq 0$ and the norm on $\|v\|_{H_A^s(X; V)}$ defined by duality, $H_A^s(X; V) = (H_A^{-s}(X; V))'$, for real $s < 0$. If the norm on $H_A^k(X; V)$ is defined in the usual way for integer $k \geq 0$,
$$
\|v\|_{H_A^k(X; V)} = \left(\sum_{j=0}^k\int_X |\nabla_A^j v|^2\,d\vol_g\right)^{1/2},
$$
then the constants appearing in the \apriori estimates in Corollaries \ref{cor:Sell_You_42-12_heat_equation_alpha_is_zero_apriori_estimates} and \ref{cor:Sell_You_42-12_heat_equation_alpha_is_one_apriori_estimates} may acquire dependencies on the reference connection, $A$, and Riemannian metric, $g$.
\end{rmk}

\subsection{Regularity of weak solutions to the heat equation}
\label{sec:Sell_You_4-2-3_heat_equation_vector_bundle_regularity}
We have the following analogue of the regularity result, Corollary \ref{cor:Sell_You_42-13}, for the heat equation \eqref{eq:Sell_You_42-1_heat_equation}.

\begin{cor}[Regularity of weak solutions]
\label{cor:Sell_You_42-13_heat_equation}
For $\alpha\in\RR$, let $u_0 \in H_A^\alpha(X;V)$ and let $f \in L_{\loc}^p([0, T); H_A^\alpha(X;V))$, for some $p$ with $2 \leq p \leq \infty$. Then the following statements are valid:
\begin{enumerate}
\item The weak solution $u$ of Equation \eqref{eq:Sell_You_42-1_heat_equation} in $H_A^\alpha(X;V)$ on $[0, T)$ is a strong solution to Equation  \eqref{eq:Sell_You_42-1_heat_equation} in $H_A^\alpha(X;V)$ on $[0, T)$, and $u$ satisfies
\begin{equation}
\label{eq:Sell_You_42-36_heat_equation}
u \in C([0, T); H_A^\alpha(X;V)) \cap C^{0,\theta_0}_{\loc}([0, T); H_A^\sigma(X;V))
\cap C^{0,\theta_1}_{\loc}((0, T); H_A^{\alpha+1+2\beta}(X;V)),
\end{equation}
for every $\sigma$ and $\beta$ with $\sigma < \alpha$ and $0 \leq \beta < 1 - 1/p$, where $\theta_0 = \theta_0(\sigma)$ and $\theta_1 = \theta_1(\beta)$ are positive.

\item For every $\tau \in (0, T)$, one has $u(\tau) \in H_A^{\alpha+1}(X;V)$, and the conclusions in
Items \eqref{item:Theorem_Sell_You_42-12_1_heat_equation_alpha_is_real} -- \eqref{item:Theorem_Sell_You_42-12_3_heat_equation_alpha_is_real} in Theorem \ref{thm:Sell_You_42-12_heat_equation_alpha_is_real} are valid on the interval $\tau \leq t < T$, with $\alpha$ replaced by $\alpha+1$ and $u_0$ replaced by $u(\tau)$. In addition, the following regularity properties are valid for the \emph{translate} $u_\tau := u(\cdot + \tau)$ \cite[Section 2.1.3]{Sell_You_2002}, where the interval $[0, T)$ is replaced by $[0, T - \tau)$:
\begin{equation}
\label{eq:Sell_You_42-34_heat_equation}
u_\tau \in C([0, T); H_A^\beta(X;V)) \cap C^{0,\theta_1}_{\loc}([0, T); H_A^\sigma(X;V))
\cap L_{\loc}^2([0, T); H_A^{\alpha+2}(X;V)),
\end{equation}
for $\beta \leq \alpha + 1$ and $\sigma < \alpha + 1$, where $\theta_1 = \theta_1(\sigma) > 0$. If in addition, $p > 2$,
then \eqref{eq:Sell_You_42-31_heat_equation_alpha_is_real} implies that for some $\theta_2 = \theta_2(\sigma) > 0$, one has
\begin{equation}
\label{eq:Sell_You_42-35_heat_equation}
u_\tau \in C^{0,\theta_1}_{\loc}((0, T); H_A^\sigma(X;V)) \quad\hbox{for } \alpha + 1 \leq \sigma < \alpha + 2 - \frac{2}{p}.
\end{equation}
Furthermore, the following inequalities are valid, for $0<\tau<t<T$:
\begin{subequations}
\label{eq:Sell_You_42-37_heat_equation}
\begin{align}
\label{eq:Sell_You_42-37_Linfty_time_Halpha+1_space_u_heat_equation}
\|u(t)\|_{H_A^{\alpha+1}(X)}^2
&\leq
e^{-\lambda_1 (t-\tau)}\|u(\tau)\|_{H_A^{\alpha+1}(X)}^2 + \int_\tau^t e^{-\lambda_1(t-s)}\|f(s)\|_{H_A^\alpha(X)}^2 \,ds,
\\
\label{eq:Sell_You_42-37_L2_time_Halpha+2_space_u_heat_equation}
\int_\tau^t\|u(s)\|_{H_A^{\alpha+2}(X)}^2 \,ds
&\leq
\|u(\tau)\|_{H_A^{\alpha+1}(X)}^2 + \int_\tau^t \|f(s)\|_{H_A^\alpha(X)}^2 \,ds,
\\
\label{eq:Sell_You_42-37_L2_time_Halpha_space_partial_t_u_heat_equation}
\int_\tau^t\|\partial_tu(s)\|_{H_A^\alpha(X)}^2 \,ds
&\leq
2\|u(\tau)\|_{H_A^{\alpha+1}(X)}^2 + 4\int_{\tau}^t \|f(s)\|_{H_A^\alpha(X)}^2 \,ds.
\end{align}
\end{subequations}
\end{enumerate}
\end{cor}

Next, we have the following analogue of the regularity result, Theorem \ref{thm:Sell_You_42-14}, for the heat equation \eqref{eq:Sell_You_42-1_heat_equation}.

\begin{thm}[Spatial regularity implied by temporal regularity of the source function and spatial regularity of the initial data]
\label{thm:Sell_You_42-14_heat_equation}
For $\alpha\in\RR$, let
$$
f \in C([0, T); H_A^\alpha(X; V)) \cap W^{1,p}([0, T); H_A^{\alpha-1}(X; V)),
$$
for some $p$ with $2 \leq p \leq \infty$, and let $g = \partial_tf$. Let $u_0 \in H_A^{\alpha+2}(X;V)$ and, for $\sA := \nabla_A^*\nabla_A + 1$, define
$$
v_0 := f(0) - \sA u_0 \in H_A^\alpha(X; V).
$$
Let $u = u(t)$ be the weak solution to Equation \eqref{eq:Sell_You_42-1_heat_equation} in the space $H_A^\alpha(X; V)$, with $u(0) = u_0$. Then the following properties hold:
\begin{enumerate}
\item $u$ is a mild solution of Equation \eqref{eq:Sell_You_42-1_heat_equation} in $H_A^\sigma(X;V)$, for each $\sigma < \alpha + 2$, and $u$ is a strong solution in $H_A^\beta(X; V)$, for each $\beta \leq \alpha$. Moreover, one has
\begin{equation}
\label{eq:Sell_You_42-38_heat_equation}
u \in C^1([0, T); H_A^\beta(X; V)) \cap C([0, T); H_A^\nu(X;V)) \cap C^{0,\theta_1}_{\loc}([0, T); H_A^\sigma(X;V)),
\end{equation}
for all $\beta$, $\nu$, and $\sigma$ with $\beta \leq \alpha$, and $\nu \leq \alpha + 2$, and $\sigma < \alpha + 2$, where
$\theta_1 = \theta_1(\sigma) > 0$.

\item The function $v := \partial_t u$ satisfies
\begin{equation}
\label{eq:Sell_You_42-39_heat_equation}
v(t) = e^{-\sA t}v_0 + \int_0^t e^{-\sA(t-r)}g(r)\,dr, \quad\hbox{for } t \geq 0,
\end{equation}
in any space $H_A^\beta(X; V)$, with $\beta < \alpha$, and $v$ is a weak solution of $\partial_tv + \sA v = g(t)$ in the space $H_A^\alpha(X; V)$ with
\begin{equation}
\label{eq:Sell_You_42-40_heat_equation}
v \in C([0, T); H_A^\beta(X; V)) \cap L_{\loc}^2([0, T); H_A^{\alpha+1}(X;V)) \cap C^{0,\theta_2}([0, T); H_A^\sigma(X;V)),
\end{equation}
for all $\beta$ and $\sigma$ with $\beta \leq \alpha$ and $\sigma < \alpha$, where $\theta_2 = \theta_2(\sigma) > O$.
\end{enumerate}
\end{thm}

Theorem \ref{thm:Sell_You_42-14_heat_equation} leads to the familiar result that $C^\infty$ spatial and temporal smoothness of the source function and $C^\infty$ spatial smoothness of the initial data imply $C^\infty$ smoothness of a solution to the heat equation for $t \geq 0$, namely

\begin{cor}[$C^\infty$ smoothness of a solution implied by $C^\infty$ spatial and temporal smoothness of the source function and $C^\infty$ spatial regularity of the initial data]
\label{cor:Sell_You_42-14_heat_equation_smoothness_for_t_geq_zero}
Let
$$
f \in C^\infty([0, T)\times X; V) \quad\hbox{and}\quad u_0 \in C^\infty(X; V).
$$
If $u = u(t)$ is a the weak solution to Equation \eqref{eq:Sell_You_42-1_heat_equation} in the space $H_A^{\alpha_0}(X; V)$, for some $\alpha_0 \in \RR$, with $u(0) = u_0$, then
$$
u \in C^\infty([0, T)\times X; V)
$$
and $u$ is a classical solution to Equation \eqref{eq:Sell_You_42-1_heat_equation} on $[0,T)\times X$.
\end{cor}

\begin{proof}
The conclusion follows by repeated application of Theorem \ref{thm:Sell_You_42-14_heat_equation} to Equation \eqref{eq:Sell_You_42-1_heat_equation} and its derivatives with respect to time, for all $\alpha\geq \alpha_0$.
\end{proof}

Lastly, we have the following analogue of the regularity result, Theorem \ref{thm:Sell_You_42-15}, for the heat equation \eqref{eq:Sell_You_42-1_heat_equation}.

\begin{thm}[Spatial regularity implied by temporal regularity of the source function]
\label{thm:Sell_You_42-15_heat_equation}
For $\alpha\in\RR$, let
$$
f \in C([0, T); H_A^\alpha(X; V)) \cap W^{1,2}([0, T); H_A^{\alpha-1}(X; V)).
$$
For any $u_0 \in H_A^\alpha(X; V)$, let $u = u(t)$ denote the weak solution of Equation \eqref{eq:Sell_You_42-1_heat_equation} in $H_A^\alpha(X; V)$ on $[0, T)$. Then the following hold:
\begin{enumerate}
\item $u$ is a strong solution in $H_A^\alpha(X; V)$, and it satisfies \eqref{eq:Sell_You_42-29_heat_equation}, \eqref{eq:Sell_You_42-31_heat_equation_alpha_is_real}, and \eqref{eq:Sell_You_42-36_heat_equation}, with $p = \infty$.

\item $u$ is a mild solution in $H_A^\beta(X; V)$, for each $\beta < \alpha + 1$.

\item $u$ satisfies
\begin{equation}
\label{eq:Sell_You_42-43_heat_equation}
u \in C^1((0, T); H_A^\beta(X; V)) \cap C((0, T); H_A^\nu(X;V)) \cap C^{0,\theta_4}_{\loc}((0, T); H_A^\sigma(X;V)),
\end{equation}
for every $\beta$, $\nu$, and $\sigma$ with $\beta \leq \alpha$, $\nu \leq \alpha + 2$, and $\sigma < \alpha + 2$, where
$\theta_4 = \theta_4(\sigma) > O$.

\item For every $\tau \in (0, T)$, the translate $u_\tau$ is a strong solution to \eqref{eq:Sell_You_42-1_heat_equation}, where $f$ is replaced by $f_\tau \in H_A^{\alpha+1}(X;V)$ on $[0,T - \tau)$.
\end{enumerate}
\end{thm}

Theorem \ref{thm:Sell_You_42-15_heat_equation} also leads to the usual result that $C^\infty$ spatial and temporal smoothness of the source function $C^\infty$ smoothness of a solution to the heat equation for $t > 0$, namely

\begin{cor}[$C^\infty$ smoothness of a solution implied by $C^\infty$ spatial and temporal smoothness of the source function]
\label{cor:Sell_You_42-15_heat_equation_smoothness_for_t_greaterthan_zero}
Let $\alpha_0 \in \RR$ and
$$
f \in C^\infty([0, T)\times X; V) \quad\hbox{and}\quad u_0 \in H_A^{\alpha_0}(X; V).
$$
If $u = u(t)$ is a the weak solution to Equation \eqref{eq:Sell_You_42-1_heat_equation} in the space $H_A^{\alpha_0}(X; V)$, then
$$
u \in C^\infty((0, T)\times X; V) \cap C([0, T); H_A^{\alpha_0}(X; V)),
$$
and $u$ is a classical solution to Equation \eqref{eq:Sell_You_42-1_heat_equation} on $(0,T)\times X$.
\end{cor}

\begin{proof}
The conclusion follows by repeated application of Theorem \ref{thm:Sell_You_42-15_heat_equation} to Equation \eqref{eq:Sell_You_42-1_heat_equation} and its derivatives with respect to time, for all $\alpha\geq \alpha_0$.
\end{proof}

\subsection{Application to the heat equation defined by the Hodge Laplace operator}
\label{sec:Sell_You_4-2-3_heat_equation_vector_bundle_Hodge_Laplacian}
One can relate the \apriori estimates described in this section involving the augmented connection Laplacian, $\nabla_A^*\nabla_A+1$ on $C^\infty(X;\Lambda^p\otimes\ad P) = \Omega^p(X;\ad P)$ in \eqref{eq:Connection_Laplacian}, and the \emph{Hodge Laplace operator}, for any integer $p \geq 1$,
\begin{equation}
\label{eq:Lawson_page_93_Hodge_Laplacian}
\Delta_A := d_A^*d_A + d_Ad_A^* \quad\hbox{on } \Omega^p(X;\ad P),
\end{equation}
via the Bochner-Weitzenb\"ock formulae \cite{Bourguignon_1981, Bourguignon_1990}, \cite[Appendix C]{FU}, \cite[Appendix II]{Lawson} and \cite{Wu_1988}. From \cite[Corollaries II.2 and II.3]{Lawson}, respectively, one has
\begin{align}
\label{eq:Lawson_corollary_II-2}
\Delta_A a = \nabla_A^*\nabla_Aa + \{\Ric_g, a\} + \{F_A, a\}, \quad\forall\, a \in \Omega^1(X;\ad P),
\\
\label{eq:Lawson_corollary_II-3}
\Delta_A v = \nabla_A^*\nabla_Av + \{\Ric_g, v\} + \{\Riem_g, v\} + \{F_A, v\}, \quad\forall\, v \in \Omega^2(X;\ad P),
\end{align}
where $\Ric_g$ and $\Riem_g$ denote the Ricci and Riemann curvature tensors of the Riemannian metric $g$ on the manifold $X$ of dimension $d \geq 2$ and $\{, \}$ denotes universal bilinear expressions (independent of the Riemannian metric on $X$) and which we also may write using `$\times$' for brevity. In applications of \eqref{eq:Lawson_corollary_II-3}, we often combine the terms $\Ric_g \times v$ and $\Riem_g \times v$ and simply write $\Riem_g \times v$ instead of $\Ric_g \times v + \Riem_g \times v$ in \eqref{eq:Lawson_corollary_II-3}.

\section{Local well-posedness for the Yang-Mills heat equation}
\label{sec:Local_well-posedness_yang_mills_heat_equation}
We now return to the setting described in the Introduction to our monograph. Let $G$ be a compact Lie group and $P$ a principal $G$-bundle with $C^\infty$ connection, $A_1$, over a closed, connected, oriented, smooth manifold, $X$, of dimension $d \geq 2$ and Riemannian metric $g$. While our ultimate goal is to develop existence, uniqueness, and regularity theory for a family of connections, $A(t) = A_1 + a(t)$ on $P$ for $t \geq 0$, solving the nonlinear \emph{Yang-Mills gradient flow equation},
\begin{equation}
\label{eq:Yang-Mills_gradient_flow_equation}
\frac{\partial a}{\partial t} + d_{A(t)}^*F_{A(t)} = 0 \quad\hbox{in } \Omega^1(X; \Lambda^1\otimes\ad P), \quad \hbox{for } t > 0,
\end{equation}
with initial data,
\begin{equation}
\label{eq:Yang-Mills_heat_or_gradient_flow_equation_initial_condition}
a(0) = a_0 \in \Omega^1(X; \ad P),
\end{equation}
we shall first consider the closely related nonlinear \emph{Yang-Mills heat equation},
\begin{equation}
\label{eq:Yang-Mills_heat_equation}
\frac{\partial a}{\partial t} + d_{A(t)}^*F_{A(t)} + d_{A(t)}d_{A(t)}^*a(t) = 0
\quad\hbox{in } \Omega^1(X; \Lambda^1\otimes\ad P), \quad \hbox{for } t > 0,
\end{equation}
As we explain in more detail in Section \ref{eq:Yang-Mills_heat_equation}, this is a nonlinear \emph{parabolic} equation for
$$
a(t) \in \Omega^1(X; \ad P), \quad\forall\, t > 0,
$$
unlike the nonlinear Yang-Mills gradient flow equation \eqref{eq:Yang-Mills_gradient_flow_equation}. Consequently, we can apply the abstract theory for linear and nonlinear evolution equations developed in Sections \ref{sec:Sell_You_4-2} and \ref{sec:Local_existence_nonlinear_evolution_equation_Banach_space}, respectively, coupled with the existence, uniqueness, and regularity theory for elliptic systems and associated semigroup theory developed in Section \ref{sec:Sell_You_3-8-2_standard_Sobolev_spaces}.

References for the existence, uniqueness, and regularity theory for parabolic systems include Amann \cite{Amann_1995}, Koshelev \cite{Koshelev_1995}, Lady{\v{z}}enskaja, Solonnikov, and Ural$'$ceva \cite{LadyzenskajaSolonnikovUralceva}, while references for scalar parabolic equations include Krylov \cite{Krylov_LecturesHolder, Krylov_LecturesSobolev}, Lady{\v{z}}enskaja, Solonnikov, and Ural$'$ceva \cite{LadyzenskajaSolonnikovUralceva}, and Lieberman \cite{Lieberman}. However, it will be rarely, if ever, possible in our monograph to appeal to any specific result in those references which one can apply directly to the Yang-Mills heat equation or related linear parabolic equations. For that reason, we shall develop the theory we need in this monograph for these equations from basic principles.

\subsection{The Yang-Mills heat equation}
\label{subsec:Yang-Mills_heat_equation}
We continue the notation and setup of the Introduction to this section. Recall that the Laplace operator, $\Delta_{A_1}$ on $\Omega^1(X; \ad P)$, is defined by the exterior covariant derivative,
$$
d_{A_1}:\Omega^l(X; \ad P) \to \Omega^{l+1}(X; \ad P),
$$
and its formal adjoint with respect to the Riemannian metric,
$$
d_{A_1}^*:\Omega^l(X; \ad P) \to \Omega^{l-1}(X; \ad P),
$$
namely \cite{DK, FU, FrM, Lawson, LM},
\begin{equation}
\label{eq:Laplace_operator_on_adP_valued_one_forms}
\Delta_{A_1} = d_{A_1}^*d_{A_1} + d_{A_1}d_{A_1}^*.
\end{equation}
By writing the nonlinear Yang-Mills heat equation \eqref{eq:Yang-Mills_heat_equation} as a perturbation,
\begin{equation}
\label{eq:Yang-Mills_heat_equation_as_perturbation_rough_Laplacian_plus_one_heat_equation}
\frac{\partial a}{\partial t} + (\nabla_{A_1}^*\nabla_{A_1} + 1)a(t) = \sF(a(t)), \quad\forall\, t > 0,
\end{equation}
of the \emph{linear} heat equation defined by the \emph{augmented connection Laplacian}, $\nabla_{A_1}^*\nabla_{A_1} + 1$, that is,
\begin{equation}
\label{eq:Linear_heat_equation_with_rough_Laplacian_plus_one_on_Omega_1_adP}
\frac{\partial a}{\partial t} + (\nabla_{A_1}^*\nabla_{A_1} + 1)a(t) = f(t)
\quad\hbox{in } \Omega^1(X; \Lambda^1\otimes\ad P), \quad \hbox{for } t > 0,
\end{equation}
with $f(t) := \sF(a(t))$, we obtain
the following schematic expression for the \emph{Yang-Mills heat equation nonlinearity},
\begin{multline}
\label{eq:Yang-Mills_heat_equation_nonlinearity_relative_rough_Laplacian_plus_one}
-\sF(a) := d_{A_1}^*F_{A_1} + (F_{A_1} - 1)\times a + \Ric_g\times a + \nabla_{A_1}a\times a + a\times a\times a,
\\
\forall\, a \in \Omega^1(X; \ad P).
\end{multline}
Here, we have applied the Bochner-Weitzenb\"ock formula \eqref{eq:Lawson_corollary_II-2} to express the Laplace operator $\Delta_{A_1}$ in terms of the connection Laplace operator $\nabla_{A_1}^*\nabla_{A_1}$ plus zeroth-order terms and used the definition of the curvature, $F_A$, to write \cite{DK, FU, FrM, Lawson, LM}
\begin{equation}
\label{eq:FA_1+a_expression}
F_{A_1 + a} = F_{A_1} + d_{A_1}a + [a, a] \in \Omega^2(X; \ad P).
\end{equation}
Rather than keep precise track of universal but otherwise unimportant constants or $C^\infty$ coefficients depending at most on the Riemannian metric, $g$, or Lie group, $G$, we shall instead use schematic expressions such as
$$
F_{A_1 + a} = F_{A_1} + d_{A_1}a + a \times a.
$$
This convention explains the origin of the relatively simple form of the schematic expression for the nonlinearity \eqref{eq:Yang-Mills_heat_equation_nonlinearity_relative_rough_Laplacian_plus_one}.

\subsection{`Standing Hypothesis A' and the augmented connection Laplacian}
\label{subsec:Standing_hypothesis_A_and_the_augmented_connection_Laplacian}
We continue the notation and setup of the preceding subsection. For $p \in (1,\infty)$, we choose
\begin{subequations}
\label{eq:Standing_hypothesis_operator_A_rough_Laplacian_plus_one_on_Lp}
\begin{align}
\label{eq:W_is_Lp_X_Lambda1_adP}
\cW &= L^p(X; \Lambda^1\otimes\ad P),
\\
\label{eq:sA_is_rough_Laplacian_plus_one}
\sA &= \nabla_{A_1}^*\nabla_{A_1} + 1,
\\
\label{eq:V_is_W2p_X_Lambda1_adP}
\calV &= \sD(\sA_p) = W^{2,p}_{A_1}(X; \Lambda^1\otimes\ad P).
\end{align}
\end{subequations}
Here, the domain, $\sD(\sA_p)$, is defined with respect to the range, $L^p(X; \Lambda^1\otimes\ad P)$, as the smallest closed extension of the realization, $\sA_p$, of the partial differential operator, $\sA$, on $L^p(X; \Lambda^1\otimes\ad P)$,
$$
\sA_p: \sD(\sA_p) \subset L^p(X; \Lambda^1\otimes\ad P) \to L^p(X; \Lambda^1\otimes\ad P).
$$
We shall employ the fractional powers,
$$
\calV^{2\alpha} = W_{A_1}^{2\alpha, p}(X; \Lambda^1\otimes\ad P),
$$
for $\alpha \in \RR$ with $\alpha \geq 0$, recalling that $\calV^0 = \cW$ \cite[Section 3.7]{Sell_You_2002}.

To simplify notation, we shall not distinguish between the partial differential operator, $\nabla_{A_1}^*\nabla_{A_1} + 1$ in \eqref{eq:sA_is_rough_Laplacian_plus_one} on $C^\infty(X; \Lambda^1\otimes\ad P)$, and its realization on $L^p(X; \Lambda^1\otimes\ad P)$. According to Theorem \ref{thm:Haller-Dintelmann_Heck_Hieber_Theorem_3-1_vector_bundle_manifold}, the realization $\sA_p$ is sectorial on $L^p(X; \Lambda^1\otimes\ad P)$ and $-\sA_p$ is the infinitesimal generator of an analytic semigroup on $L^p(X; \Lambda^1\otimes\ad P)$. Our choice of $\sA$ in \eqref{eq:sA_is_rough_Laplacian_plus_one} also defines a positive realization, $\sA_p$, on $L^p(X; \Lambda^1\otimes\ad P)$ and hence fulfills Hypothesis \ref{hyp:Sell_You_4_standing_hypothesis_A} (`Standing Hypothesis A') with $\calV \equiv \sD(\sA_p) = W^{2, p}(X; \Lambda^1\otimes\ad P)$, noting that the domain, $\sD(\sA_p)$, of the smallest closed extension of $\sA_p$ is identified by the \apriori estimate \eqref{eq:Krylov_Sobolev_lectures_theorems_8-5-3_and_6_estimate} in Theorem \ref{thm:Krylov_Sobolev_lectures_8-5-3}. Consequently, depending on the spatial and temporal regularity for the source function, $f$, in equation \eqref{eq:Linear_heat_equation_with_rough_Laplacian_plus_one_on_Omega_1_adP} and the regularity of the initial data, $a_0$, the conclusions of the abstract Theorems \ref{thm:Sell_You_42-9} and \ref{thm:Sell_You_42-10} on existence, uniqueness, and regularity of solutions to the linear heat equation \eqref{eq:Linear_heat_equation_with_rough_Laplacian_plus_one_on_Omega_1_adP} will hold for the choice of setup in \eqref{eq:Standing_hypothesis_operator_A_rough_Laplacian_plus_one_on_Lp}.

Naturally, the same conclusions in the abstract Theorems \ref{thm:Sell_You_42-9} and \ref{thm:Sell_You_42-10} would also hold if we had instead selected $\sA = \Delta_{A_1} + \lambda_0$, where $\lambda_0 \geq 0$ is a constant chosen large enough (depending on the curvatures of the connection $A_1$ on $P$ and the Riemannian metric $g$ on $X$) to ensure that $\Delta_{A_1} + \lambda_0$ defines a sectorial operator on $L^p(X; \Lambda^1\otimes\ad P)$ and, in place of \eqref{eq:Linear_heat_equation_with_rough_Laplacian_plus_one_on_Omega_1_adP}, consider
\begin{equation}
\label{eq:Sell_You_42-1_linear_heat_equation_with_Laplacian_on_Omega_1_adP}
\frac{\partial a}{\partial t} + (\Delta_{A_1} + \lambda_0)a(t) = f(t)
\quad\hbox{in } W^{2, p}(X; \Lambda^1\otimes\ad P), \quad \hbox{for } t > 0.
\end{equation}
Alternatively, when $f \equiv 0$, one could employ perturbation theory for analytic semigroups as in \cite[Section 4.4]{Sell_You_2002} and write \eqref{eq:Sell_You_42-1_linear_heat_equation_with_Laplacian_on_Omega_1_adP} as a perturbed equation,
\begin{equation}
\label{eq:Sell_You_44-1_linear_heat_equation_with_Laplacian_on_Omega_1_adP_perturbed}
\frac{\partial a}{\partial t} + \sA_p a(t) = \sB_p a(t)
\quad\hbox{in } W^{2, p}(X; \Lambda^1\otimes\ad P), \quad \hbox{for } t > 0,
\end{equation}
where $\sA = \nabla_{A_1}^*\nabla_{A_1} + 1$ and $\sB = \nabla_{A_1}^*\nabla_{A_1} + 1 - \Delta_{A_1}$, a zeroth-order operator thanks to the Bochner-Weitzenb\"ock formula \eqref{eq:Lawson_corollary_II-2}. However, because such a linear perturbation can be viewed as just a special case of the more general nonlinear perturbations allowed in our treatment of the nonlinear Yang-Mills heat equation \eqref{eq:Yang-Mills_heat_equation_as_perturbation_rough_Laplacian_plus_one_heat_equation}, based in part on \cite[Section 4.7]{Sell_You_2002}, so we can restrict our attention to the simpler linear heat equation \eqref{eq:Linear_heat_equation_with_rough_Laplacian_plus_one_on_Omega_1_adP} rather than \eqref{eq:Sell_You_42-1_linear_heat_equation_with_Laplacian_on_Omega_1_adP} or \eqref{eq:Sell_You_44-1_linear_heat_equation_with_Laplacian_on_Omega_1_adP_perturbed}.

\subsection{Sobolev embedding and multiplication theorems for real derivative exponents and the Yang-Mills heat equation nonlinearity}
\label{subsec:Sobolev_embedding_multiplication_Yang-Mills_heat_equation_nonlinearity}
We continue the notation and setup of the preceding subsection. When choosing suitable Sobolev spaces to represent $\cW$, and thus $\calV$, and the exponent, $\beta$, defining the fractional power, $\calV^{2\beta}$, the structure of the Yang-Mills heat equation nonlinearity \eqref{eq:Yang-Mills_heat_equation_nonlinearity_relative_rough_Laplacian_plus_one} suggests that there are two cases to consider for $p \in (1, \infty)$ and $\beta \in (0, 1)$:
\begin{enumerate}
\item \label{item:2beta_p_greaterthan_d}
$p$ and $\beta$ obey
\begin{equation}
\label{eq:2beta_p_greaterthan_d}
\beta \geq \frac{1}{2} \quad\hbox{and}\quad 2\beta p > d,
\end{equation}
so $W^{2\beta,p}_{A_1}(X; \Lambda^1\otimes\ad P) \hookrightarrow W^{1,p}_{A_1}(X; \Lambda^1\otimes\ad P)$ and $W^{2\beta,p}_{A_1}(X; \Lambda^1\otimes\ad P) \hookrightarrow C(X; \Lambda^1\otimes\ad P)$ and $W^{2\beta,p}(X)$ is a Banach algebra by Lemma \ref{lem:Freed_Uhlenbeck_equation_6-34_nonnegative_real_algebra}; or

\item \label{item:2beta_p_lessthan_d}
$p$ and $\beta$ obey
\begin{equation}
\label{eq:2beta_p_lessthan_d}
0 < 2\beta p < d,
\end{equation}
and one must appeal to the Sobolev multiplication result Lemma \ref{lem:Freed_Uhlenbeck_equation_6-34_nonnegative_real}.
\end{enumerate}
Clearly, Case \eqref{item:2beta_p_greaterthan_d} will mean that handling the Yang-Mills heat equation nonlinearity will be straightforward but imposes a stronger requirement on the regularity of the initial data $a_0 \in W^{2\beta,p}_{A_1}(X; \Lambda^1\otimes\ad P)$ when we appeal to Theorem \ref{thm:Sell_You_lemma_47-1}. The constraint $\beta \geq 1/2$ is included in \eqref{eq:2beta_p_greaterthan_d} to ensure that product terms such as $\nabla_{A_1}a \times a$ belong to $W^{2\beta,p}_{A_1}(X; \Lambda^1\otimes\ad P)$ when $a\in W^{2\beta,p}_{A_1}(X; \Lambda^1\otimes\ad P)$.

The constraint \eqref{eq:2beta_p_greaterthan_d} means that
$$
\frac{d}{2p} < \beta < 1,
$$
and thus also requires $p > d/2$ to ensure that some choice of $\beta$ is possible.

We now consider the more complicated Case \eqref{item:2beta_p_lessthan_d}, where $p \in (1, \infty)$ and $\beta \in (0, 1)$ obey $0 < 2\beta p < d$. From Lemma \ref{lem:Freed_Uhlenbeck_equation_6-34_nonnegative_real}, we have continuous Sobolev multiplication maps,
\begin{gather}
\label{eq:Sobolev_multiplication_Ws1_p_times_Ws2_p_into_Lp_when_d_geq_2}
W^{s_1,p}(X) \times W^{s_2,p}(X) \to L^p(X),
\\
\label{eq:Sobolev_multiplication_Ws1_p_times_Ws2_p_times_Ws3_p_into_Lp_when_d_geq_2}
W^{s_1,p}(X) \times W^{s_2,p}(X) \times W^{s_3,p}(X)\to L^p(X),
\end{gather}
where we choose the minimal $s_1 \geq 0$, $s_2=s_1+1$, and $s_3 \geq 0$, respectively, such that
\begin{gather*}
\left(s_1 - \frac{d}{p}\right) + \left(s_2 - \frac{d}{p}\right) \geq -\frac{d}{p}, \quad\hbox{with } s_1p, s_2p < d,
\\
\left(s_1 - \frac{d}{p}\right) + \left(s_2 - \frac{d}{p}\right) + \left(s_3 - \frac{d}{p}\right) \geq -\frac{d}{p}, \quad\hbox{with } s_1p, s_2p, s_3p < d,
\end{gather*}
In the map \eqref{eq:Sobolev_multiplication_Ws1_p_times_Ws2_p_into_Lp_when_d_geq_2}, we choose (keeping in mind the term $\nabla_{A_1}a\times a$ in the Yang-Mills heat equation nonlinearity)
\begin{equation}
\label{eq:Sobolev_multiplication_Ws1_p_times_Ws2_p_into_Lp_when_d_geq_2_minimal_s1_and_s2}
s_1 = \frac{d}{2p} - \frac{1}{2} \quad\hbox{and}\quad s_2 = s_1 + 1 = \frac{d}{2p} + \frac{1}{2}, \quad\hbox{with } p < d,
\end{equation}
noting that $s_1 \equiv d/2p -1/2 \geq 0 \iff p \leq d$ and $s_2p \equiv d/2 + p/2 < d \iff p < d$. In the map \eqref{eq:Sobolev_multiplication_Ws1_p_times_Ws2_p_times_Ws3_p_into_Lp_when_d_geq_2}, we choose (keeping in mind that our optimal selections in \eqref{eq:Sobolev_multiplication_Ws1_p_times_Ws2_p_into_Lp_when_d_geq_2_minimal_s1_and_s2} restrict further freedom of choice),
\begin{equation}
\label{eq:Sobolev_multiplication_Ws1_p_times_Ws2_p_times_Ws3_p_into_Lp_when_d_geq_2_minimal_s1_and_s2_and_s3}
s_1 = s_2 = \frac{d}{2p} + \frac{1}{2} \quad\hbox{and}\quad s_3 = \frac{d}{p} - 1, \quad\hbox{with } p < d.
\end{equation}
Consequently, we obtain the Sobolev multiplication maps
\begin{gather}
\label{eq:Sobolev_multiplication_Wdover2p-half_p_times_Wdover2p+half_p_into_Lp_when_d_geq_2}
W^{\frac{d}{2p} - \frac{1}{2},p}(X) \times W^{\frac{d}{2p} + \frac{1}{2},p}(X) \to L^p(X),
\\
\label{eq:Sobolev_multiplication_Wdover2p+half_p_times_Wdover2p+half_p_times_Wdoverp-one_p_into_Lp_when_d_geq_2}
W^{\frac{d}{2p} + \frac{1}{2},p}(X) \times W^{\frac{d}{2p} + \frac{1}{2},p}(X) \times W^{\frac{d}{p} - 1,p}(X)\to L^p(X).
\end{gather}
Note that
$$
\frac{d}{p} - 1 \leq \frac{d}{2p} + \frac{1}{2} \iff p \geq \frac{d}{3},
$$
and, without loss of generality, we may further restrict $p \in (1, d)$ to $p \in [d/3, d)$ when $d \geq 4$ or $p \in (1, d)$ when $d = 3$ to guarantee the embedding,
\begin{equation}
\label{eq:Sobolev_embedding_Wdover2p+half_p_into_Wdoverp-one_p_when_d_geq_2}
W^{\frac{d}{2p} + \frac{1}{2}, p}(X) \hookrightarrow W^{\frac{d}{p} - 1, p}(X), \quad\hbox{for } d \geq 3 \hbox{ and } p \geq d/3.
\end{equation}
The constraint $0 < 2\beta p < d$ in Case \eqref{item:2beta_p_lessthan_d} already restricts our choice of $\beta$ to $\beta \in (0, d/2p)$. We now aim to select $\beta$ to ensure that the following embedding holds,
$$
\calV^{2\beta} \hookrightarrow W^{\frac{d}{2p} + \frac{1}{2},p}(X; \Lambda^1\otimes\ad P),
$$
so $2\beta \geq d/2p + 1/2$ and thus
\begin{equation}
\label{eq:beta_function_d_and_p_ensuring_bounded_quadratic_and_cubic_maps}
\beta \geq \frac{d}{4p} + \frac{1}{4}, \quad\hbox{for } d \geq 3 \hbox{ and } \frac{d}{3} < p < d.
\end{equation}
A minimal choice of $\beta = 2/4p + 1/4$ satisfies $\beta \in (1/2, 1)$. A choice of $p = d/3$ would force $\beta \geq 1$, so we must restrict to $p > d/3$, while the fact that $p < d$ ensures $d/2p > 1/2$ and guarantees that \emph{some} choice of $\beta \in (1/2, d/2p)$ is possible via \eqref{eq:beta_function_d_and_p_ensuring_bounded_quadratic_and_cubic_maps}.

In developing the higher-order regularity of a strong solution to the Yang-Mills heat equation, we shall need analogues of the quadratic and cubic maps \eqref{eq:Sobolev_multiplication_Wdover2p-half_p_times_Wdover2p+half_p_into_Lp_when_d_geq_2} and
\eqref{eq:Sobolev_multiplication_Wdover2p+half_p_times_Wdover2p+half_p_times_Wdoverp-one_p_into_Lp_when_d_geq_2} with $L^p(X)$ on the right-hand side replaced by $W^{s, p}(X)$ for $s \geq 0$. However, by applying Lemma \ref{lem:Freed_Uhlenbeck_equation_6-34_nonnegative_real} for $s \geq 0$ such that
\begin{equation}
\label{eq:d_and_p_and_s_unstable_range_for_Sobolev_cubic_and_quadratic_multiplication}
\left(s +\frac{d}{2p} + \frac{1}{2}\right)p < d,
\end{equation}
so the hypothesis \eqref{eq:Freed_Uhlenbeck_6-34_nonnegative_real_hypothesis} of Lemma \ref{lem:Freed_Uhlenbeck_equation_6-34_nonnegative_real} is obeyed, we see that the following Sobolev multiplication maps,
\begin{gather}
\label{eq:Sobolev_multiplication_Ws+dover2p-half_p_times_Ws+dover2p+half_p_into_Wsp_when_d_geq_2}
W^{s + \frac{d}{2p} - \frac{1}{2},p}(X) \times W^{s +\frac{d}{2p} + \frac{1}{2},p}(X) \to W^{s, p}(X),
\\
\label{eq:Sobolev_multiplication_Ws+dover2p+half_p_times_Ws+dover2p+half_p_times_Ws+doverp-one_p_into_Wsp_when_d_geq_2}
W^{s + \frac{d}{2p} + \frac{1}{2},p}(X) \times W^{s + \frac{d}{2p} + \frac{1}{2},p}(X) \times W^{s + \frac{d}{p} - 1,p}(X)\to W^{s, p}(X),
\end{gather}
are defined and continuous. Of course, when $s>0$ is large enough that
$$
\left(s +\frac{d}{2p} + \frac{1}{2}\right)p > d,
$$
then $W^{s + \frac{d}{2p} - \frac{1}{2},p}(X)$ is a Banach algebra by Lemma \ref{lem:Freed_Uhlenbeck_equation_6-34_nonnegative_real_algebra}. In this situation, we may choose
$$
\cW = W^{s, p}(X; \Lambda^1\otimes\ad P),
$$
and $\beta \in (0, 1)$ large enough to ensure
$$
\calV^{2\beta} \hookrightarrow W^{s + \frac{d}{2p} + \frac{1}{2},p}(X; \Lambda^1\otimes\ad P),
$$
so $2\beta \geq s + d/2p + 1/2$, that is,
\begin{equation}
\label{eq:beta_function_d_and_p_ensuring_bounded_quadratic_and_cubic_maps_when_cW_is_Wsp}
\frac{s}{2} + \frac{d}{4p} + \frac{1}{4} \leq \beta < 1,
\end{equation}
for $d \geq 2$ and $p \in (1, \infty)$ obeying \eqref{eq:d_and_p_and_s_unstable_range_for_Sobolev_cubic_and_quadratic_multiplication}. We shall apply these Sobolev multiplication results in the next subsection.

\subsection{Local well-posedness for a strong solution to the Yang-Mills heat equation in $W^{2\beta, p}$ given initial data in $W^{2\beta, p}$}
\label{subsec:Sell_You_4-7-1_and_4-7-2_Yang-Mills}
We continue the notation and setup of the preceding subsection. We now apply Theorems \ref{thm:Sell_You_lemma_47-1} and \ref{thm:Sell_You_lemma_47-2} to establish existence and uniqueness of a mild solution in $L^p(X; \Lambda^1\otimes\ad P)$ and a solution in $W_{A_1}^{2\beta, p}p(X; \Lambda^1\otimes\ad P)$, respectively, on an interval $[0, \tau)$ for some $\tau > 0$. We shall give the precise dependencies of $\tau$ in a later subsection.

\begin{thm}[Existence and uniqueness of mild solutions to the Yang-Mills heat equation in $W^{2\beta, p}$ with initial data in $W^{2\beta, p}$ over a manifold of dimension $d$]
\label{thm:Existence_uniqueness_mild_solution_Yang-Mills_heat_equation_in_W_2beta_p_initial_data_in_W_2beta_p}
Let $G$ be a compact Lie group and $P$ a principal $G$-bundle over a closed, connected, orientable, Riemannian, smooth manifold, $X$,  of dimension $d \geq 2$. Let $A_1$ be a reference connection of class $C^\infty$ on $P$. Let $d \geq 2$ and $p \in (1, \infty)$ and $\beta \in (0, 1)$ obey \emph{one} of the following conditions:
\footnote{Taking $\beta = d/4p + 1/4$ in \eqref{eq:d_greaterthan_3_and_dover3_lessthan_p_lessthan_dover2_and_2beta_p_lessthan_d} ensures that $1/2 < \beta < d/2p$, using the fact that $p < d$ for this case.}
\begin{gather}
\label{eq:p_greaterthan_dover2_and_2beta_p_greaterthan_d}
d \geq 3 \quad\hbox{and}\quad p > \frac{d}{2} \quad\hbox{and}\quad \beta \geq \frac{1}{2} \quad\hbox{and}\quad \frac{d}{2p} < \beta < 1, \quad\hbox{or }
\\
\label{eq:d_greaterthan_3_and_dover3_lessthan_p_lessthan_dover2_and_2beta_p_lessthan_d}
d \geq 3 \quad\hbox{and}\quad \frac{d}{3} < p < d \quad\hbox{and}\quad \frac{d}{4p} + \frac{1}{4} \leq \beta < \frac{d}{2p}.
\end{gather}
Then for every $a_0 \in W_{A_1}^{2\beta, p}(X; \Lambda^1\otimes\ad P)$, there is a unique mild solution $A = A_1 + a$ to the Yang-Mills heat equation \eqref{eq:Yang-Mills_heat_equation_as_perturbation_rough_Laplacian_plus_one_heat_equation} in $W_{A_1}^{2\beta, p}(X; \Lambda^1\otimes\ad P)$ with $A(0) = A_1 + a_0$ and
\begin{multline}
\label{eq:Sell_You_47-5_Yang-Mills_mild_solution_in_W_2beta_p_initial_data_in_W_2beta_p}
a \in C([0, \tau]; W_{A_1}^{2\beta, p}(X; \Lambda^1\otimes\ad P)) \cap C_{\loc}^{0,\theta_1}([0, \tau); W_{A_1}^{2\alpha, p}(X; \Lambda^1\otimes\ad P))
\\
\cap C_{\loc}^{0,\theta}((0, \tau); W_{A_1}^{2r, p}(X; \Lambda^1\otimes\ad P)),
\end{multline}
for some $\tau >0$ and all $\alpha$ and $r$ with $0 \leq \alpha < r$ and $0 \leq r < 1$, where $\theta_1 > 0$ and $\theta > 0$.
\end{thm}

\begin{rmk}[Application of Theorem \ref{thm:Existence_uniqueness_mild_solution_Yang-Mills_heat_equation_in_W_2beta_p_initial_data_in_W_2beta_p} when $X$ has dimension four]
\label{rmk:Existence_uniqueness_mild_solution_Yang-Mills_heat_equation_in_W_2beta_p_initial_data_in_W_2beta_p}
When $d = 4$, the choice $p = 2$ and $\beta = 1$ is not possible in Theorem \ref{thm:Existence_uniqueness_mild_solution_Yang-Mills_heat_equation_in_W_2beta_p_initial_data_in_W_2beta_p}, so this result does not yield an analogue, for mild solutions, of the local well-posedness result for strong solutions in $H_{A_1}^2(X; \Lambda^1\otimes\ad P)$, as described in \cite[Section 4.3]{Struwe_1994}, with initial data in $H_{A_1}^1(X; \Lambda^1\otimes\ad P)$. However, when $d = 4$, we can choose
\begin{inparaenum}[\itshape a\upshape)]
\item $p > 2$ and $\beta \in (2/p, 1)$, or
\item $4/3 < p \leq 2$ and $\beta = 1/p + 1/4 \in (3/4, 1)$.
\end{inparaenum}
\end{rmk}

\begin{proof}[Proof of Theorem \ref{thm:Existence_uniqueness_mild_solution_Yang-Mills_heat_equation_in_W_2beta_p_initial_data_in_W_2beta_p}]
We shall apply Theorem \ref{thm:Sell_You_lemma_47-1} with $\cW = L^p(X; \Lambda^1\otimes\ad P)$ and $\calV = W_{A_1}^{2,p}(X; \Lambda^1\otimes\ad P)$, for $p \in (1, \infty)$ and, noting that $\calV^{2\beta} = W_{A_1}^{2\beta,p}(X; \Lambda^1\otimes\ad P)$, for $\beta \in [0, 1)$. We need only consider the more difficult case, where $\beta$ and $p$ obey \eqref{eq:d_greaterthan_3_and_dover3_lessthan_p_lessthan_dover2_and_2beta_p_lessthan_d}. The case where $\beta$ and $p$ obey \eqref{eq:p_greaterthan_dover2_and_2beta_p_greaterthan_d} follows by the same argument, except we can now use the fact that $W^{2\beta,p}(X)$ is a Banach algebra by Lemma \ref{lem:Freed_Uhlenbeck_equation_6-34_nonnegative_real_algebra} in applications of the Sobolev multiplication theorems.

It remains to verify, for $p$ and $\beta$ obeying \eqref{eq:d_greaterthan_3_and_dover3_lessthan_p_lessthan_dover2_and_2beta_p_lessthan_d}, that the Yang-Mills heat equation nonlinearity $\sF$ in \eqref{eq:Yang-Mills_heat_equation_nonlinearity_relative_rough_Laplacian_plus_one}
obeys the hypothesis \eqref{eq:Sell_You_47-4}, namely that $\sF \in C_{\Lip}([0,\infty)\times \calV^{2\beta}; \cW)$ in the sense that $\sF$ obeys \eqref{eq:Sell_You_46-1} and \eqref{eq:Sell_You_46-2} (equivalently, \eqref{eq:Sell_You_46-7} and \eqref{eq:Sell_You_46-8}), for suitable constants $K_0 = K_0(B)$ and $K_1 = K_1(B)$, given any bounded set $B \subset \calV$.

For this purpose, we shall use our Sobolev multiplication results to estimate the different terms in the expression \eqref{eq:Yang-Mills_heat_equation_nonlinearity_relative_rough_Laplacian_plus_one} for the Yang-Mills heat equation nonlinearity, $\sF(t,a)$, and the difference, $\sF(t,a_1) - \sF(t,a_2)$, when
$$
a, a_1, a_2 \in W_{A_1}^{2\beta,p}(X; \Lambda^1\otimes\ad P).
$$
We begin with the

\setcounter{step}{0}
\begin{step}[Boundedness of the Yang-Mills nonlinearity in $L^p(X)$]
The quadratic and cubic terms obey (for simplicity, we suppress denoting the dependency on $A_1$ in our notation for norms on the Sobolev spaces $W_{A_1}^{s,p}(X; \Lambda^1\otimes\ad P)$)
\begin{equation}
\label{eq:L2_estimate_quadratic_and_cubic_nonlinearities_H2_initial_data}
\begin{aligned}
\|\nabla_{A_1}a\times a\|_{L^p(X)} &\leq C\|\nabla_{A_1}a\|_{W^{\frac{d}{2p} - \frac{1}{2}, p}(X)} \|a\|_{W^{\frac{d}{2p} + \frac{1}{2}, p}(X)}
\quad\hbox{(by \eqref{eq:Sobolev_multiplication_Wdover2p-half_p_times_Wdover2p+half_p_into_Lp_when_d_geq_2})}
\\
&\leq C\|a\|_{W^{\frac{d}{2p} + \frac{1}{2}, p}(X)}^2,
\\
\|a\times a\times a\|_{L^p(X)} &\leq C\|a\|_{W^{\frac{d}{p} -1, p}(X)}\|a\|_{W^{\frac{d}{2p} + \frac{1}{2}, p}(X)}^2
\quad\hbox{(by \eqref{eq:Sobolev_multiplication_Wdover2p+half_p_times_Wdover2p+half_p_times_Wdoverp-one_p_into_Lp_when_d_geq_2})}
\\
&\leq C\|a\|_{W^{\frac{d}{2p} + \frac{1}{2}, p}(X)}^3
\quad\hbox{(by \eqref{eq:Sobolev_embedding_Wdover2p+half_p_into_Wdoverp-one_p_when_d_geq_2})},
\end{aligned}
\end{equation}
while affine term obeys
\begin{equation}
\label{eq:L2_estimate_linear_term_FA_btimesb_H2_initial_data}
\begin{aligned}
{}&\|d_{A_1}^*F_{A_1} + \Ric_g\times a + (F_{A_1}-1)\times a\|_{L^p(X)}
\\
&\quad \leq \|d_{A_1}^*F_{A_1}\|_{L^p(X)} + C\left(\|\Ric_g\|_{L^\infty(X)} + \|F_{A_1}\|_{L^\infty(X)} + 1\right)\|a\|_{L^p(X)}
\\
&\quad \leq \|d_{A_1}^*F_{A_1}\|_{L^p(X)} + C\left(\|\Ric_g\|_{L^\infty(X)} + \|F_{A_1}\|_{L^\infty(X)} + 1\right)\|a\|_{W^{\frac{d}{2p} + \frac{1}{2}, p}(X)}.
\end{aligned}
\end{equation}
We now combine \eqref{eq:L2_estimate_quadratic_and_cubic_nonlinearities_H2_initial_data} and \eqref{eq:L2_estimate_linear_term_FA_btimesb_H2_initial_data} with our expression \eqref{eq:Yang-Mills_heat_equation_nonlinearity_relative_rough_Laplacian_plus_one} for $\sF(t,a)$ to give
\begin{multline}
\label{eq:Yang-Mills_heat_equation_nonlinearity_relative_rough_Laplacian_plus_1_Lp_bound}
\|\sF(t, a)\|_{L^p(X)}
\leq
\|d_{A_1}^*F_{A_1}\|_{L^p(X)} + C\left(\|\Ric_g\|_{L^\infty(X)} + \|F_{A_1}\|_{L^\infty(X)} + 1\right)\|a\|_{W^{\frac{d}{2p} + \frac{1}{2}, p}(X)}
\\
+ C\left(\|a\|_{W^{\frac{d}{2p} + \frac{1}{2}, p}(X)}^2 + \|a\|_{W^{\frac{d}{2p} + \frac{1}{2}, p}(X)}^3\right),
\\
\forall\, (t, a) \in [0,\infty) \times W_{A_1}^{2\beta, p}(X; \Lambda^1\otimes\ad P).
\end{multline}
Therefore, $\sF$ obeys \eqref{eq:Sell_You_46-7}, as required for this step.
\end{step}

\begin{step}[Lipschitz property the Yang-Mills nonlinearity in $L^p(X)$]
We first observe that
$$
\nabla_{A_1}a_1\times a_1 - \nabla_{A_1}a_2\times a_2 = \nabla_{A_1}(a_1 - a_2)\times a_1 + \nabla_{A_1}a_2\times (a_1 - a_2),
$$
Thus,
\begin{align*}
{}&\|\nabla_{A_1}a_1\times a_1 - \nabla_{A_1}a_2\times a_2\|_{L^p(X)}
\\
&\quad \leq \|\nabla_{A_1}(a_1 - a_2)\times a_1\|_{L^p(X)} + \|\nabla_{A_1}a_2\times (a_1 - a_2)\|_{L^p(X)}
\\
&\quad \leq c\|\nabla_{A_1}(a_1 - a_2)\|_{W^{\frac{d}{2p} - \frac{1}{2}, p}(X)} \|a_1\|_{W^{\frac{d}{2p} + \frac{1}{2}, p}(X)}
\\
&\qquad + c\|\nabla_{A_1}a_2\|_{W^{\frac{d}{2p} - \frac{1}{2}, p}(X)} \|a_1 - a_2\|_{W^{\frac{d}{2p} + \frac{1}{2}, p}(X)}
\quad\hbox{(by \eqref{eq:Sobolev_multiplication_Wdover2p-half_p_times_Wdover2p+half_p_into_Lp_when_d_geq_2})}
\\
&\quad \leq c\left(\|a_1\|_{W^{\frac{d}{2p} + \frac{1}{2}, p}(X)} + \|a_2\|_{W^{\frac{d}{2p} + \frac{1}{2}, p}(X)} \right)\|a_1 - a_2\|_{W^{\frac{d}{2p} + \frac{1}{2}, p}(X)}
\quad\hbox{(by \eqref{eq:Sobolev_embedding_Wdover2p+half_p_into_Wdoverp-one_p_when_d_geq_2})},
\end{align*}
where $c$ depends at most on $d$, $p$, and the Riemannian metric $g$ on $X$.

For the cubic term, we observe that
\begin{multline*}
a_1\times a_1\times a_1 - a_2\times a_2\times a_2
\\
= (a_1-a_2)\times a_1\times a_1 + a_2\times (a_1-a_2)\times a_1 + a_2\times a_2\times (a_1-a_2).
\end{multline*}
Thus,
\begin{align*}
{}& \|a_1\times a_1\times a_1 - a_2\times a_2\times a_2\|_{L^p(X)}
\\
&\quad \leq \|(a_1-a_2)\times a_1\times a_1\|_{L^p(X)} + \|a_2\times (a_1-a_2)\times a_1\|_{L^p(X)} + \|a_2\times a_2\times (a_1-a_2)\|_{L^p(X)}
\\
&\quad \leq c\|a_1-a_2\|_{W^{\frac{d}{2p} - \frac{1}{2}, p}(X)} \|a_1\|_{W^{\frac{d}{2p} + \frac{1}{2}, p}(X)}^2
\\
&\qquad + c\|a_2\|_{W^{\frac{d}{2p} + \frac{1}{2}, p}(X)} \|a_1-a_2\|_{W^{\frac{d}{2p} - \frac{1}{2}, p}(X)} \|a_1\|_{W^{\frac{d}{2p} + \frac{1}{2}, p}(X)}
\\
&\qquad + c\|a_2\|_{W^{\frac{d}{2p} + \frac{1}{2}, p}(X)}^2 \|a_1-a_2\|_{W^{\frac{d}{2p} - \frac{1}{2}, p}(X)}
\quad\hbox{(by \eqref{eq:Sobolev_multiplication_Wdover2p+half_p_times_Wdover2p+half_p_times_Wdoverp-one_p_into_Lp_when_d_geq_2})}
\\
&\quad\leq c\left(\|a_1\|_{W^{\frac{d}{2p} + \frac{1}{2}, p}(X)}^2 + \|a_2\|_{W^{\frac{d}{2p} + \frac{1}{2}, p}(X)}^2 \right)
\|a_1-a_2\|_{W^{\frac{d}{2p} - \frac{1}{2}, p}(X)}
\quad\hbox{(by \eqref{eq:Sobolev_embedding_Wdover2p+half_p_into_Wdoverp-one_p_when_d_geq_2})},
\end{align*}
using $2xy \leq x^2 + y^2$ to obtain the last inequality. Therefore, combining the preceding inequalities with our expression \eqref{eq:Yang-Mills_heat_equation_nonlinearity_relative_rough_Laplacian_plus_one} for $\sF(t,a_1)$ and $\sF(t,a_2)$ yields
\begin{multline}
\label{eq:Yang-Mills_heat_equation_nonlinearity_relative_rough_Laplacian_plus_one_Lp_Lipschitz}
\|\sF(t, a_1) - \sF(t, a_2)\|_{L^p(X)}
\\
\leq C\left(\|a_1\|_{W^{\frac{d}{2p} + \frac{1}{2}, p}(X)} + \|a_2\|_{W^{\frac{d}{2p} + \frac{1}{2}, p}(X)}
+ \|a_1\|_{W^{\frac{d}{2p} + \frac{1}{2}, p}(X)}^2 + \|a_2\|_{W^{\frac{d}{2p} + \frac{1}{2}, p}(X)}^2\right)
\\
\times\|a_1 - a_2\|_{W^{\frac{d}{2p} + \frac{1}{2}, p}(X)}
\\
+ C\left(\|\Ric_g\|_{L^\infty(X)} + \|F_{A_1}\|_{L^\infty(X)} + 1\right)\|a_1 - a_2\|_{W^{\frac{d}{2p} + \frac{1}{2}, p}(X)},
\\
\forall\, t \in [0,\infty) \hbox{ and } a_1, a_2 \in W_{A_1}^{2\beta, p}(X; \Lambda^1\otimes\ad P).
\end{multline}
Therefore, $\sF$ obeys \eqref{eq:Sell_You_46-8}, as required for this step.
\end{step}


The conclusions now follows immediately from Theorem \ref{thm:Sell_You_lemma_47-1}.
\end{proof}

The proof of Theorem \ref{thm:Existence_uniqueness_mild_solution_Yang-Mills_heat_equation_in_W_2beta_p_initial_data_in_W_2beta_p} actually yields more, as we see next.

\begin{thm}[Strong solutions and regularity of mild solutions to the Yang-Mills heat equation in $W^{2\beta, p}$ with initial data in $W^{2\beta, p}$ over a manifold of dimension $d$]
\label{thm:Existence_uniqueness_strong_solution_Yang-Mills_heat_equation_in_W_2beta_p_initial_data_in_W_2beta_p}
Assume the hypotheses of Theorem \ref{thm:Existence_uniqueness_mild_solution_Yang-Mills_heat_equation_in_W_2beta_p_initial_data_in_W_2beta_p}. If $a_0 \in W_{A_1}^{2\beta, p}(X; \Lambda^1\otimes\ad P)$ and $a(t)$ is a mild solution of the Yang-Mills heat equation
\eqref{eq:Yang-Mills_heat_equation_as_perturbation_rough_Laplacian_plus_one_heat_equation} in the Banach space $W_{A_1}^{2\beta, p}(X; \Lambda^1\otimes\ad P)$ on an interval $[0, T)$ for some $T > 0$, then $a(t)$ is a strong solution in $W_{A_1}^{2\beta, p}(X; \Lambda^1\otimes\ad P)$ on the interval $[0, T)$, and it satisfies
\begin{multline}
\label{eq:Sell_You_47-7_Yang-Mills_strong_solution_in_W_2beta_p_initial_data_in_W_2beta_p}
a \in C([0, T); W_{A_1}^{2\alpha, p}(X; \Lambda^1\otimes\ad P)) \cap C_{\loc}^{0,1-r}((0, T); W_{A_1}^{2r, p}(X; \Lambda^1\otimes\ad P))
\\
\cap C((0, T); W_{A_1}^{2, p}(X; \Lambda^1\otimes\ad P)),
\end{multline}
for all $\alpha$ and $r$ with $0 \leq \alpha \leq \beta$ and $0 \leq r < 1$.
\end{thm}

\begin{proof}
It suffices to observe that Yang-Mills heat equation nonlinearity $\sF$ in \eqref{eq:Yang-Mills_heat_equation_nonlinearity_relative_rough_Laplacian_plus_one} obeys the hypothesis \eqref{eq:Sell_You_lemma_47-2_F_C_Lipschitz_theta} in Theorem \ref{thm:Sell_You_lemma_47-2}. Indeed, for $\cW = L^p(X; \Lambda^1\otimes\ad P)$ and $\calV = W_{A_1}^{1, p}(X; \Lambda^1\otimes\ad P)$, we have
$$
\sF \in C_{\Lip}([0,\infty)\times \calV^{2\beta}; \cW),
$$
by the Definition \ref{defn:Sell_You_page_221_CLip_0_to_infinity_times_V_into_W} of $C_{\Lip}([0,\infty)\times \calV^{2\beta}; \cW)$ and as we showed that $\sF$ obeys
\eqref{eq:Yang-Mills_heat_equation_nonlinearity_relative_rough_Laplacian_plus_1_Lp_bound} and
\eqref{eq:Yang-Mills_heat_equation_nonlinearity_relative_rough_Laplacian_plus_one_Lp_Lipschitz} in the proof of Theorem \ref{thm:Existence_uniqueness_mild_solution_Yang-Mills_heat_equation_in_W_2beta_p_initial_data_in_W_2beta_p}. (These two properties were established in detail for the case $2\beta p < d$, but the proof for the case $2\beta p > d$ is considerably simpler; note also that $T > 0$ did not need to be sufficiently small for these properties to hold.) Because $\sF(t,v) = \sF(v)$ for all $t \in [0,\infty)$ and $v \in W_{A_1}^{2\beta, p}(X; \Lambda^1\otimes\ad P)$, then $\sF$ automatically satisfies \eqref{eq:Sell_You_46-6} with $\theta = 1$. The conclusion now follows immediately from Theorem \ref{thm:Sell_You_lemma_47-2}.
\end{proof}

\subsection{Higher-order spatial and temporal regularity of a strong solution to the Yang-Mills heat equation in $W^{2\beta, p}$}
\label{subsec:Sell_You_4-8-2_but_with_standing_hypothesis_A_Yang-Mills}
As we can see from Theorem \ref{thm:Existence_uniqueness_strong_solution_Yang-Mills_heat_equation_in_W_2beta_p_initial_data_in_W_2beta_p}, the strong solution, $a(t)$, has greater spatial regularity than the initial data, $a_0$, for any $t > 0$. We can thus iterate
Theorem \ref{thm:Existence_uniqueness_strong_solution_Yang-Mills_heat_equation_in_W_2beta_p_initial_data_in_W_2beta_p} to give $C^\infty$ spatial and temporal regularity for $t \in (0, T)$ and appeal to the regularity result in Theorem \ref{thm:Sell_You_42-10} for a linear evolution equation to give $C^\infty$ spatial and temporal regularity for $t \in [0, T)$.

\begin{thm}[Classical solutions and $C^\infty$ regularity of mild solutions for positive time to the Yang-Mills heat equation with initial data in $W^{2\beta, p}$ over a manifold of dimension $d$]
\label{thm:Smoothness_strong_solution_Yang-Mills_heat_equation_initial_data_in_W_2beta_p}
Assume the hypotheses of Theorem \ref{thm:Existence_uniqueness_mild_solution_Yang-Mills_heat_equation_in_W_2beta_p_initial_data_in_W_2beta_p}. If $a_0 \in W_{A_1}^{2\beta, p}(X; \Lambda^1\otimes\ad P)$ and $a(t)$ is a mild solution of the Yang-Mills heat equation
\eqref{eq:Yang-Mills_heat_equation_as_perturbation_rough_Laplacian_plus_one_heat_equation} in $W_{A_1}^{2\beta, p}(X; \Lambda^1\otimes\ad P)$ on an interval $[0, T)$ for some $T > 0$, then $a(t)$ is a classical $C^\infty$ solution on the interval $(0, T)$, and it satisfies
\begin{equation}
\label{eq:Smoothness_strong_solution_Yang-Mills_heat_equation_initial_data_in_W_2beta_p}
a \in C([0, T); W_{A_1}^{2\beta, p}(X; \Lambda^1\otimes\ad P)) \cap C^\infty((0, T) \times X; \Lambda^1\otimes\ad P)).
\end{equation}
If in addition, $a_0 \in W_{A_1}^{2, p}(X; \Lambda^1\otimes\ad P)$, then $a(t)$ also satisfies
\begin{equation}
\label{eq:C1_on_0_to_T_strong_solution_Yang-Mills_heat_equation_initial_data_in_W_2beta_p}
a \in C^1([0, T); L^p(X; \Lambda^1\otimes\ad P)).
\end{equation}
\end{thm}

\begin{proof}
We divide the proof of the regularity property \eqref{eq:Smoothness_strong_solution_Yang-Mills_heat_equation_initial_data_in_W_2beta_p} into the two cases corresponding to the two alternative hypotheses for $\beta$ and $p$ in Theorem \ref{thm:Existence_uniqueness_mild_solution_Yang-Mills_heat_equation_in_W_2beta_p_initial_data_in_W_2beta_p}.

\setcounter{case}{0}
\begin{case}[$2\beta p > d$]
The case where $\beta$ and $p$ obey \eqref{eq:p_greaterthan_dover2_and_2beta_p_greaterthan_d} is easier because we can now use the fact that $W^{2\beta,p}(X)$ is a Banach algebra by Lemma \ref{lem:Freed_Uhlenbeck_equation_6-34_nonnegative_real_algebra} in applications of the Sobolev multiplication theorem. Theorem \ref{thm:Existence_uniqueness_strong_solution_Yang-Mills_heat_equation_in_W_2beta_p_initial_data_in_W_2beta_p} yields
\begin{equation}
\label{eq:a_continuous_up_to_time0_into_W2beta_p_and_greaterthan_time0_into_W2_p}
a \in C([0, T); W_{A_1}^{2\beta, p}(X; \Lambda^1\otimes\ad P)) \cap C((0, T); W_{A_1}^{2, p}(X; \Lambda^1\otimes\ad P)).
\end{equation}
But $W^{2, p}(X)$ is a Banach algebra because, \afortiori, this is true of $W^{2\beta, p}(X)$ for this case, and as the Yang-Mills heat equation nonlinearity $\sF(a(t))$ in \eqref{eq:Yang-Mills_heat_equation_nonlinearity_relative_rough_Laplacian_plus_one} is a cubic polynomial in $a(t)$, we have
\begin{equation}
\label{eq:Fa_continuous_up_to_time0_into_W2beta_p_and_greaterthan_time0_into_W2_p}
\sF(a) \in C([0, T); W_{A_1}^{2\beta, p}(X; \Lambda^1\otimes\ad P)) \cap C((0, T); W_{A_1}^{2, p}(X; \Lambda^1\otimes\ad P)).
\end{equation}
We can now apply Theorem \ref{thm:Sell_You_42-10}, with source function
$$
f(a(t)) := \sF(a(t)), \quad t \in [0, T).
$$
Because \eqref{eq:Fa_continuous_up_to_time0_into_W2beta_p_and_greaterthan_time0_into_W2_p} gives, for any $t_0 \in (0, T)$,
$$
f \in C([t_0, T); W_{A_1}^{2, p}(X; \Lambda^1\otimes\ad P)),
$$
we may take $\nu = 1$ in the hypothesis \eqref{eq:Sell_You_42-20} of Theorem \ref{thm:Sell_You_42-10}. Noting that $a(t_0) \in W_{A_1}^{2, p}(X; \Lambda^1\otimes\ad P)$ by \eqref{eq:a_continuous_up_to_time0_into_W2beta_p_and_greaterthan_time0_into_W2_p} and hence choosing $\mu = 1+\beta$ in the hypotheses of Theorem \ref{thm:Sell_You_42-10} (applied to the interval $[t_0, T)$ rather than $[0, T)$, we obtain from \eqref{eq:Sell_You_42-21} that
$$
a \in C([t_0, T); W_{A_1}^{2 + 2\beta, p}(X; \Lambda^1\otimes\ad P)),
$$
recalling that $\calV = W_{A_1}^{2, p}(X; \Lambda^1\otimes\ad P)$ and so
$$
\calV^{2\mu} = \calV^{2+2\beta} = W_{A_1}^{2 + 2\beta, p}(X; \Lambda^1\otimes\ad P).
$$
Since $t_0 \in (0, T)$ was arbitrary, we obtain
$$
a \in C((0, T); W_{A_1}^{2 + 2\beta, p}(X; \Lambda^1\otimes\ad P)).
$$
We may assume without loss of generality that $\beta \geq 1/2$ and so, in particular, we have
$$
a \in C((0, T); W_{A_1}^{3, p}(X; \Lambda^1\otimes\ad P)).
$$
According to Lemma \ref{lem:Sell_You_37-4}, we have that $e^{-\cA t}$ is an analytic semigroup on $\calV^\alpha$, for each $\alpha \in \RR$, if $\cA$ is a positive, sectorial operator on a Banach space $\cW$ and $\calV = \sD(\cA)$. We can therefore repeat the preceding application of Theorem \ref{thm:Sell_You_42-10} on the interval $[t_0, T)$, this time with $a(t_0) \in W_{A_1}^{3, p}(X; \Lambda^1\otimes\ad P)$ and $\cW = W_{A_1}^{1, p}(X; \Lambda^1\otimes\ad P))$ and $\calV = W_{A_1}^{3, p}(X; \Lambda^1\otimes\ad P))$, to give
$$
a \in C((0, T); W_{A_1}^{3 + 3\beta, p}(X; \Lambda^1\otimes\ad P)).
$$
By iterating in this way, we find that
$$
a \in C((0, T); W_{A_1}^{k+2, p}(X; \Lambda^1\otimes\ad P)), \quad\forall\, k \in \NN,
$$
and applying the Sobolev Embedding Theorem \cite[Theorem 4.12]{AdamsFournier}, we obtain
$$
a \in C((0, T); C^l(X; \Lambda^1\otimes\ad P)), \quad\forall\, l \in \NN.
$$
But $a(t)$ is a strong solution in $W_{A_1}^{k+2, p}(X; \Lambda^1\otimes\ad P)$, for all $k \geq \NN$, to the equation,
$$
\frac{\partial a}{\partial t} + (\nabla_{A_1}^*\nabla_{A_1} + 1)a(t) = \sF(a(t)), \quad\hbox{a.e. } t \in (0, T).
$$
Therefore,
$$
\frac{\partial a}{\partial t} = -(\nabla_{A_1}^*\nabla_{A_1} + 1)a + \sF(a) \in C((0, T); C^l(X; \Lambda^1\otimes\ad P)), \quad\forall\, l \in \NN,
$$
and by iterating the preceding argument, noting that $\sF(a(t))$ is a cubic polynomial in $a(t)$, we obtain
$$
\frac{\partial^m a}{\partial t^m} \in C((0, T); C^l(X; \Lambda^1\otimes\ad P)), \quad\forall\, l, m \in \NN,
$$
that is
$$
a \in C^m((0, T); C^l(X; \Lambda^1\otimes\ad P)), \quad\forall\, l, m \in \NN,
$$
and the first conclusion \eqref{eq:Smoothness_strong_solution_Yang-Mills_heat_equation_initial_data_in_W_2beta_p} follows in this case.
\end{case}


\begin{case}[$2\beta p < d$]
We now consider the more difficult case, where $\beta$ and $p$ obey \eqref{eq:d_greaterthan_3_and_dover3_lessthan_p_lessthan_dover2_and_2beta_p_lessthan_d}. However, in this situation we again note that, for initial data $a_0 \in W_{A_1}^{2\beta, p}(X; \Lambda^1\otimes\ad P)$ and $\beta \in (0, 1)$ as in the hypotheses of Theorem \ref{thm:Existence_uniqueness_mild_solution_Yang-Mills_heat_equation_in_W_2beta_p_initial_data_in_W_2beta_p} for the case $2\beta p < d$, Theorem \ref{thm:Existence_uniqueness_strong_solution_Yang-Mills_heat_equation_in_W_2beta_p_initial_data_in_W_2beta_p} still yields the spatial regularity improvement \eqref{eq:a_continuous_up_to_time0_into_W2beta_p_and_greaterthan_time0_into_W2_p} for $t > 0$.

If $2p > d$, then we could have chosen $\beta \in (0, 1)$ so that $2\beta p > d$. Thus, we may suppose without loss of generality that $2p < d$ and so we may choose $s > 0$ satisfying \eqref{eq:d_and_p_and_s_unstable_range_for_Sobolev_cubic_and_quadratic_multiplication}, that is,
$$
s +\frac{d}{2p} + \frac{1}{2} < 2 \quad\hbox{and}\quad \left(s +\frac{d}{2p} + \frac{1}{2}\right)p < d.
$$
The quadratic and cubic Sobolev multiplication maps \eqref{eq:Sobolev_multiplication_Ws+dover2p-half_p_times_Ws+dover2p+half_p_into_Wsp_when_d_geq_2} and \eqref{eq:Sobolev_multiplication_Ws+dover2p+half_p_times_Ws+dover2p+half_p_times_Ws+doverp-one_p_into_Wsp_when_d_geq_2} are defined and continuous for this $s$ and may be used in place of \eqref{eq:Sobolev_multiplication_Wdover2p-half_p_times_Wdover2p+half_p_into_Lp_when_d_geq_2} and
\eqref{eq:Sobolev_multiplication_Wdover2p+half_p_times_Wdover2p+half_p_times_Wdoverp-one_p_into_Lp_when_d_geq_2} when estimating the $W_{A_1}^{s,p}(X)$ norm of $\sF(a(t))$ since
$$
a(t) \in W_{A_1}^{2, p}(X; \Lambda^1\otimes\ad P) \hookrightarrow  W_{A_1}^{s +\frac{d}{2p} + \frac{1}{2}, p}(X; \Lambda^1\otimes\ad P), \quad\forall\, t \in (0, T).
$$
We can repeatedly apply the regularity result in Theorem \ref{thm:Sell_You_42-10}, just as in the case $2\beta p > d$, with
$$
\cW = W_{A_1}^{s, p}(X; \Lambda^1\otimes\ad P) \quad\hbox{and}\quad \calV = W_{A_1}^{s +\frac{d}{2p} + \frac{1}{2}, p}(X; \Lambda^1\otimes\ad P),
$$
until we reach the continuous spatial range,
$$
a \in C((0, T); C(X; \Lambda^1\otimes\ad P)),
$$
and then the argument in the case $2\beta p > d$ yields the first conclusion \eqref{eq:Smoothness_strong_solution_Yang-Mills_heat_equation_initial_data_in_W_2beta_p} for this case too.
\end{case}

We now verify the $C^1$ property \eqref{eq:C1_on_0_to_T_strong_solution_Yang-Mills_heat_equation_initial_data_in_W_2beta_p} of $a(t)$ for $t \in [0, 1)$. The continuity property of $\sF(a)$ in \eqref{eq:Fa_continuous_up_to_time0_into_W2beta_p_and_greaterthan_time0_into_W2_p} when $2\beta p > d$ and the polynomial estimate
\eqref{eq:Yang-Mills_heat_equation_nonlinearity_relative_rough_Laplacian_plus_1_Lp_bound} when $2\beta p < d$ gives in particular that
$$
\sF(a) \in
\begin{cases}
C([0, T); W_{A_1}^{2\beta, p}(X; \Lambda^1\otimes\ad P)), &\hbox{if } 2\beta p > d,
\\
C([0, T); L^p(X; \Lambda^1\otimes\ad P)), &\hbox{if } 2\beta p < d,
\end{cases}
$$
and so the hypothesis \eqref{eq:Sell_You_42-20} of Theorem \ref{thm:Sell_You_42-10} is obeyed with source function $f(t) = \sF(a(t))$, for $t \in [0, T)$, and the choices $p = \infty$ and $\nu = 1$ and $\calV = W_{A_1}^{2, p}(X; \Lambda^1\otimes\ad P)$. We can take $\mu = 1$ in Theorem \ref{thm:Sell_You_42-10}, since we are now given $a_0 \in W_{A_1}^{2, p}(X; \Lambda^1\otimes\ad P)$, so its conclusion \eqref{eq:Sell_You_42-21} gives
$$
a \in C([0, T); W_{A_1}^{2, p}(X; \Lambda^1\otimes\ad P)).
$$
Combining the preceding observations yields the inclusions
$$
(\nabla_{A_1}^*\nabla_{A_1} + 1)a, \ \sF(a) \in C([0, T); L^p(X; \Lambda^1\otimes\ad P)),
$$
and thus also
$$
\frac{\partial a}{\partial t} = -(\nabla_{A_1}^*\nabla_{A_1} + 1)a + \sF(a) \in C([0, T); L^p(X; \Lambda^1\otimes\ad P)).
$$
Therefore, the second conclusion \eqref{eq:C1_on_0_to_T_strong_solution_Yang-Mills_heat_equation_initial_data_in_W_2beta_p} follows.
\end{proof}

It remains to examine the higher-order regularity of the strong solution $a(t)$ up to time $t = 0$, given initial data $a_0$ with suitable regularity.

\begin{thm}[Classical solutions and $C^\infty$ regularity of mild solutions up to the initial time for the Yang-Mills heat equation with $C^\infty$ initial data over a manifold of dimension $d$]
\label{thm:Smoothness_strong_solution_Yang-Mills_heat_equation_smooth_initial_data}
Assume the hypotheses of
Theorem \ref{thm:Existence_uniqueness_mild_solution_Yang-Mills_heat_equation_in_W_2beta_p_initial_data_in_W_2beta_p}.
If $a_0 \in C^\infty(X; \Lambda^1\otimes\ad P)$ and $a(t)$ is a mild solution of the Yang-Mills heat equation
\eqref{eq:Yang-Mills_heat_equation_as_perturbation_rough_Laplacian_plus_one_heat_equation} in $W_{A_1}^{2\beta, p}(X; \Lambda^1\otimes\ad P)$ on an interval $[0, T)$ for some $T > 0$, then $a(t)$ is a classical solution on $[0, T)$ and
\begin{equation}
\label{eq:Smoothness_on_0_to_T_strong_solution_Yang-Mills_heat_equation_smooth_initial_data}
a \in C^\infty([0, T) \times X; \Lambda^1\otimes\ad P).
\end{equation}
\end{thm}

\begin{proof}
It suffices to observe that instead of taking
$$
\cW = L^p(X; \Lambda^1\otimes\ad P) \quad\hbox{and}\quad \calV = W_{A_1}^{2,p}(X; \Lambda^1\otimes\ad P),
$$
with $\cA = \nabla_{A_1}^*\nabla_{A_1} + 1$, as usual, in the proofs of Theorems
\ref{thm:Existence_uniqueness_mild_solution_Yang-Mills_heat_equation_in_W_2beta_p_initial_data_in_W_2beta_p} and
\ref{thm:Existence_uniqueness_strong_solution_Yang-Mills_heat_equation_in_W_2beta_p_initial_data_in_W_2beta_p}, we may instead choose, for any $k \in \NN$,
$$
\cW = W_{A_1}^{k,p}(X; \Lambda^1\otimes\ad P) \quad\hbox{and}\quad \calV = W_{A_1}^{k+2,p}(X; \Lambda^1\otimes\ad P).
$$
In particular, it suffices to restrict our attention to $k$ large enough that $(k+2)p > d$ and so $W^{k+2,p}(X)$ is a Banach algebra by Lemma \ref{lem:Freed_Uhlenbeck_equation_6-34_algebra}.

For $a_0 \in W_{A_1}^{k+2,p}(X; \Lambda^1\otimes\ad P)$, the proof of Theorem \ref{thm:Smoothness_strong_solution_Yang-Mills_heat_equation_initial_data_in_W_2beta_p} in the case $2\beta p > d$ and $k = 0$ yields, \emph{mutatis mutandis},
\begin{align*}
a &\in C([0, T); W_{A_1}^{k+2, p}(X; \Lambda^1\otimes\ad P))
\\
(\nabla_{A_1}^*\nabla_{A_1} + 1)a &\in C([0, T); W_{A_1}^{k, p}(X; \Lambda^1\otimes\ad P)),
\\
\sF(a) &\in C([0, T); W_{A_1}^{k+2,p}(X; \Lambda^1\otimes\ad P)),
\end{align*}
and thus also
$$
\frac{\partial a}{\partial t} = -(\nabla_{A_1}^*\nabla_{A_1} + 1)a + \sF(a) \in C([0, T); W^{k, p}(X; \Lambda^1\otimes\ad P)).
$$
Therefore,
$$
a \in C^1([0, T); W_{A_1}^{k,p}(X; \Lambda^1\otimes\ad P)).
$$
But $\dot a = \partial a/\partial t$ is a strong solution to the time derivative of the Yang-Mills heat equation \eqref{eq:Yang-Mills_heat_equation_as_perturbation_rough_Laplacian_plus_one_heat_equation},
$$
\frac{\partial \dot a}{\partial t} + (\nabla_{A_1}^*\nabla_{A_1} + 1)\dot a(t) = \frac{\partial }{\partial t}\sF(a(t)), \quad\hbox{a.e. } t \in (0, T),
$$
with initial data
$$
\dot a(0) \in W_{A_1}^{k, p}(X; \Lambda^1\otimes\ad P),
$$
and source function with regularity
$$
\frac{\partial }{\partial t}\sF(a(t)) \in C([0, T); W_{A_1}^{k, p}(X; \Lambda^1\otimes\ad P)).
$$
The regularity result \eqref{eq:Sell_You_42-21} in Theorem \ref{thm:Sell_You_42-10} now ensures that
$$
\dot a \in C([0, T); W_{A_1}^{k+2,p}(X; \Lambda^1\otimes\ad P)),
$$
and thus
$$
\frac{\partial \dot a}{\partial t} = -(\nabla_{A_1}^*\nabla_{A_1} + 1)\dot a + \frac{\partial }{\partial t}\sF(a) \in C([0, T); W_{A_1}^{k,p}(X; \Lambda^1\otimes\ad P)).
$$
Therefore,
$$
a \in C^2([0, T); W_{A_1}^{k,p}(X; \Lambda^1\otimes\ad P)).
$$
The preceding argument may be repeated, for all $m \in \NN$ and $k \in \NN$ with $(k+2)p > d$, to give
$$
a \in C^m([0, T); W_{A_1}^{k,p}(X; \Lambda^1\otimes\ad P)),
$$
and hence obtain the conclusion \eqref{eq:Smoothness_on_0_to_T_strong_solution_Yang-Mills_heat_equation_smooth_initial_data}.
\end{proof}

\begin{rmk}[Regularity of the initial data in Theorem \ref{thm:Existence_uniqueness_strong_solution_Yang-Mills_heat_equation_in_W_2beta_p_initial_data_in_W_2beta_p} and the result of Kozono, Maeda, and Naito]
\label{rmk:Comment_on_Kozono_Maeda_Naito_Theorem_3-1}
Given $a_0 \in L^d(X; \Lambda^1\otimes\ad P)$, Kozono, Maeda, and Naito \cite[Theorem 3.1]{Kozono_Maeda_Naito_1995} establish the local existence of a strong solution $a(t)$ to \eqref{eq:Yang-Mills_heat_equation_as_perturbation_rough_Laplacian_plus_one_heat_equation} with the properties that
\begin{align*}
a \in C([0,T); L^d(X; \Lambda^1\otimes\ad P) \cap C((0,T); W_{A_1}^{2,d}(X; \Lambda^1\otimes\ad P) \cap C^1((0,T); L^d(X; \Lambda^1\otimes\ad P).
\end{align*}
Their application of semigroup theory to achieve this conclusion is different from ours and this explains the discrepancy between their strong solution regularity result and that of our Theorem \ref{thm:Existence_uniqueness_strong_solution_Yang-Mills_heat_equation_in_W_2beta_p_initial_data_in_W_2beta_p}, which requires $a_0 \in W_{A_1}^{2\beta, p}(X; \Lambda^1\otimes\ad P)$.
\end{rmk}

\subsection{\Apriori estimate for the $W^{2\beta, p}$ norm and minimum lifetime of mild or strong solution to the Yang-Mills heat equation in $W^{2\beta, p}$ given initial data in $W^{2\beta, p}$}
We continue the notation and setup of the preceding subsection. We have seen that, with only the hypotheses for the structure of an abstract nonlinearity, $\cF$, as assumed in the abstract Theorem \ref{thm:Sell_You_lemma_47-1}, the abstract Lemma \ref{lem:Sell_You_lemma_47-1_tau_formula} already provides a \emph{lower} bound \eqref{eq:Sell_You_47-6} on the lifetime, $\tau$, and \emph{upper} bound \eqref{eq:Sell_You_page_234_contraction_mapping_apriori_bound_solution} in the $\calV^{2\beta}$-norm (for $\beta \in (0, 1)$) of the mild solution, $u(t)$ for $t \in [0, \tau)$, and hence the lifetime of the strong solution via Theorem \ref{thm:Sell_You_lemma_47-2}, since the constant, $\tau$, does not change in the passage from mild to strong, as is evident
from the proof Theorem \ref{thm:Sell_You_lemma_47-2}. However, because the Yang-Mills heat equation nonlinearity \eqref{eq:Yang-Mills_heat_equation_nonlinearity_relative_rough_Laplacian_plus_one} is a cubic polynomial (with respect to suitable choices of Sobolev spaces), these two bounds can be sharpened, as we have already demonstrated. Indeed, the abstract Theorem \ref{thm:Sell_You_lemma_47-1_polynomial_nonlinearity}, which assumes a $\cW$-norm estimate for $\cF(t,v)$ which is polynomial in the $\calV^{2\beta}$-norm of $v$, yields the lower bound \eqref{eq:Sell_You_47-6_polynomial_nonlinearity} for the lifetime and $\calV^{2\beta}$-norm bound \eqref{eq:Sell_You_lemma_47-1_polynomial_nonlinearity_apriori_estimate} for the solution.

Although we shall state main result of this subsection only for mild solutions to the Yang-Mills heat equation
\eqref{eq:Yang-Mills_heat_equation_as_perturbation_rough_Laplacian_plus_one_heat_equation}, the \apriori estimate and lower bound for the lifetime of course continue to hold for strong or classical solutions.

\begin{thm}[\Apriori estimate and lifetime of a mild solution to the Yang-Mills heat equation in $W^{2\beta, p}$ given initial data in $W^{2\beta, p}$]
\label{thm:Apriori_estimate_and_lifetime_mild_solution_Yang-Mills_heat_equation_in_W_2beta_p_initial_data_in_W_2beta_p}
Assume the hypotheses of
Theorem \ref{thm:Existence_uniqueness_mild_solution_Yang-Mills_heat_equation_in_W_2beta_p_initial_data_in_W_2beta_p}.
Given $b > 0$, there exists a positive constant\footnote{This constant $\tau$ is given explicitly via \eqref{eq:Sell_You_47-6_polynomial_nonlinearity}.},
$$
\tau = \tau\left(b, M_0, M_\beta, n, \beta, \kappa_0, \kappa_1\right),
$$
with the following significance.  For every $a_0 \in W^{2\beta, p}_{A_1}(X; \Lambda^1\otimes\ad P)$ obeying
$$
\|a_0\|_{W_{A_1}^{2\beta, p}(X)} \leq b,
$$
the Yang-Mills heat equation \eqref{eq:Yang-Mills_heat_equation_as_perturbation_rough_Laplacian_plus_one_heat_equation} with initial data $a(0) = a_0$ has a unique, mild solution $a(t)$ in $W^{2\beta, p}_{A_1}(X; \Lambda^1\otimes\ad P)$ on an interval $[0, \tau)$, with the regularity property \eqref{eq:Sell_You_47-5_Yang-Mills_mild_solution_in_W_2beta_p_initial_data_in_W_2beta_p}. Moreover, the solution $a(t)$ obeys the \apriori estimate,
\begin{equation}
\label{eq:Sell_You_lemma_47-1_polynomial_nonlinearity_apriori_estimate_Yang-Mills}
\|a(t)\|_{W_{A_1}^{2\beta, p}(X)}
\leq
M_0\|a_0\|_{W_{A_1}^{2\beta, p}(X)}
+
\frac{2M_\beta \kappa_0}{1-\beta} \left(1 + M_0 b\right)^n t^{1-\beta},
\quad \forall\, t\in[0, \tau].
\end{equation}
If in addition $\hat a_0 \in W^{2\beta, p}_{A_1}(X; \Lambda^1\otimes\ad P)$ obeys
$$
\|\hat a_0\|_{W_{A_1}^{2\beta, p}(X)} \leq b,
$$
and $\hat a$ is the unique, mild solution to \eqref{eq:Sell_You_47-1} in $W^{2\beta, p}_{A_1}(X; \Lambda^1\otimes\ad P)$ on $[0, \tau]$ with initial data $\hat a(0) = \hat a_0$, then
\begin{equation}
\label{eq:Sell_You_lemma_47-1_polynomial_nonlinearity_continuity_with_respect_to_initial_data_Yang-Mills}
\sup_{t\in [0, \tau]}\|a(t)-\hat a(t)\|_{W_{A_1}^{2\beta, p}(X)} \leq 2M_0\|a_0 - \hat a_0\|_{W_{A_1}^{2\beta, p}(X)}.
\end{equation}
\end{thm}

\begin{proof}
We just need to verify that the Yang-Mills heat equation nonlinearity
\eqref{eq:Yang-Mills_heat_equation_nonlinearity_relative_rough_Laplacian_plus_one} obeys the hypotheses of Theorem \ref{thm:Sell_You_lemma_47-1_polynomial_nonlinearity}. We recall from the proof of Theorem \ref{thm:Existence_uniqueness_mild_solution_Yang-Mills_heat_equation_in_W_2beta_p_initial_data_in_W_2beta_p} that, for either the case $2\beta p > d$ or $2\beta d < d$, the estimate
\eqref{eq:Yang-Mills_heat_equation_nonlinearity_relative_rough_Laplacian_plus_1_Lp_bound}
implies that $\sF$ obeys the inequality \eqref{eq:Sell_You_46-7_polynomial_nonlinearity} with $n = 3$. Similarly, for either the case $2\beta p > d$ or $2\beta d < d$, the estimate
\eqref{eq:Yang-Mills_heat_equation_nonlinearity_relative_rough_Laplacian_plus_one_Lp_Lipschitz} ensures that $\sF$ obeys the inequality \eqref{eq:Sell_You_46-8_polynomial_nonlinearity} with $n = 3$. The conclusions now follow from Theorem \ref{thm:Sell_You_lemma_47-1_polynomial_nonlinearity}, keeping track of the
dependencies of the constants.
\end{proof}

\subsection{Maximally defined solutions to the Yang-Mills heat equation}
\label{subsec:Maximally_defined_solutions_Yang-Mills_heat_equation}
Theorem \ref{thm:Sell_You_lemma_47-4} implies that the maximal interval $[0,T)$ of definition for the solution $a(t)\in H^4(X, \Lambda^1\otimes\ad P)$
to \eqref{eq:Struwe_21} is characterized by
\eqref{eq:Sell_You_47-8}, that is,
\begin{equation}
\label{eq:Sell_You_46-21_Yang-Mills}
\lim_{t\uparrow T}\|a(t)\|_{H^1(X)} = \infty.
\end{equation}
We wish to show that this is equivalent to the characterization by Struwe in \cite[Theorem 2.3]{Struwe_1994}, in which case we will have completed our alternative proof of \cite[Theorem 2.3]{Struwe_1994}.

Clearly, if $T$ is characterized by \cite[Theorem 2.3 (iii) and Equation (5)]{Struwe_1994}, then \eqref{eq:Sell_You_46-21_Yang-Mills} holds.
To establish the reverse implication, we argue by contradiction, much as in \cite[Section 7]{Struwe_1994}. Suppose that $a(t)\in H^4(X, \Lambda^1\otimes\ad P)$ and $T$ obey \eqref{eq:Sell_You_46-21_Yang-Mills}, but that there is an $R \in (0,1]$ such that \cite[Equation (15)]{Struwe_1994} holds, namely
$$
\sup_{\begin{subarray}{}x^0 \in X \\ t\in [0,T)\end{subarray}}\|F_{A(t)}\|_{L^2(B(x^0,R))} < \delta,
$$
where $\delta>0$ is as in Lemma \ref{lem:Schlatter_2-3_and_Struwe_3-3}.
Then Lemma \ref{lem:Schlatter_2-4_and_Struwe_3-6}
implies that
$$
A(t) \to A(T) \quad\hbox{in } H^1(X, \Lambda^1\otimes\ad P), \quad\hbox{as } t\uparrow T,
$$
for some $A(T) = A_0 + a(T)$ with $a(T) \in H^1(X, \Lambda^1\otimes\ad P)$, and so
$$
\lim_{t\uparrow T}\|a(t)\|_{H^1(X)} = \|a(T)\|_{H^1(X)} < \infty,
$$
which contradicts \eqref{eq:Sell_You_46-21_Yang-Mills}.

\subsection{Global $L^1$-in-time a priori estimates for a solution to a nonlinear evolution equation}
\label{subsec:Global_L1-in-time_apriori_estimates_solution_nonlinear_evolution_equation}
In this subsection, we develop abstract global \apriori estimates which complement \cite[Lemma 7.3]{Rade_1992} (whose generalization to a closed, smooth Riemannian manifold $X$ of dimension $d$ with $2 \leq d \leq 4$ is stated as the forthcoming Lemma \ref{lem:Rade_7-3}) by allowing far greater flexibility with regard to the choice of spatial regularity.

\begin{lem}[A global $L^1$-in-time \apriori estimate for a mild solution to a nonlinear evolution equation]
\label{lem:Rade_7-3_abstract_L1_in_time_V2beta_space}
Assume that Hypothesis \ref{hyp:Sell_You_4_standing_hypothesis_A} holds. Let $\beta \in [0, 1)$, and let $\eps$ be a positive constant such that
$$
\eps a^{1+\beta} M_\beta \Gamma(1-\beta) \leq \frac{1}{2},
$$
where the constants $a > 0$ and $M_\beta > 0$ are as in Theorem \ref{thm:Sell_You_37-5}, and $K \in [1,\infty)$. Suppose $t_0 \in \RR$ and $T$ is such that $t_0 < T \leq \infty$ and that $\cF \in C_{\Lip}([t_0,T)\times\calV^{2\beta}; \cW)$ obeys
$$
\cF(t,v) = f_0(t) + \cF_1(t,v) + \cF_2(t,v), \quad \forall\, (t,v) \in [t_0,T)\times\calV^{2\beta},
$$
with $f_0 \in C_{\Lip}([t_0,T); \cW)$ and $\cF_1, \cF_2 \in C_{\Lip}([t_0,T)\times\calV^{2\beta}; \cW)$ and
\begin{align*}
\|\cF_1(t,v)\|_\cW &\leq K\|v\|_\cW,
\\
\|\cF_2(t,v)\|_\cW &\leq \eps\|v\|_{\calV^{2\beta}}, \quad\forall\, (t,v) \in [t_0,T)\times\calV^{2\beta}.
\end{align*}
If $u$ is a mild solution to the nonlinear evolution equation \eqref{eq:Sell_You_47-1} on $[t_0, T)$ defined by $\cA$ and $\cF$, with $u \in C([t_0, T); \cW) \cap C((t_0, T); \calV^{2\beta})$, then
\begin{multline}
\label{eq:Rade_7-3_abstract_L1_in_time_V2beta_space_apriori_estimate}
\int_{t_0}^T \|u(t)\|_{\calV^{2\beta}} \,dt
\\
\leq a^{1+\beta} M_\beta \Gamma(1-\beta) \left(\|u_0\|_\cW
+ \int_{t_0}^T \|f_0(t)\|_\cW\,dt + K\int_{t_0}^T \|u(t)\|_\cW \,dt\right).
\end{multline}
\end{lem}

\begin{proof}
From \eqref{eq:Sell_You_47-2}, we have
$$
u(t) = e^{-\cA (t - t_0)}u_0 + \int_{t_0}^t e^{-\cA (t - s)} \cF(s,u(s)) \,ds, \quad\forall\, t \in [t_0, T).
$$
Thus, in the usual way, we see that for $t \in [t_0, T]$,
\begin{align*}
\|\cA^\beta u(t)\|_\cW &\leq \|\cA^\beta e^{-\cA (t - t_0)}u_0\|_\cW
+ \int_{t_0}^t \|\cA^\beta e^{-\cA (t - s)} \cF(s,u(s))\|_\cW \,ds
\\
&\leq M_\beta (t - t_0)^{-\beta} e^{-a(t - t_0)} \|u_0\|_\cW
 + M_\beta \int_{t_0}^t (t - s)^{-\beta} e^{-a(t - s)}\|\cF(s,u(s))\|_\cW \,ds
\\
&\qquad\hbox{(by Theorem \ref{thm:Sell_You_37-5})}
\\
&\leq M_\beta (t - t_0)^{-\beta} e^{-a(t - t_0)} \|u_0\|_{L^2(X)}
\\
&\quad + M_\beta \int_{t_0}^t (t - s)^{-\beta} e^{-a(t - s)}
\left(\|f_0(s)\|_\cW + K\|u(s)\|_\cW + \eps\|u(s)\|_{\calV^{2\beta}}\right) \,ds,
\end{align*}
where we applied the hypotheses on $\cF$ to obtain the last inequality. Recall that the \emph{gamma} and \emph{incomplete gamma functions} are given by \cite[Sections 5.2 (i) and 8.1]{Olver_Lozier_Boisvert_Clark},
$$
\gamma(\alpha, z) = \int_0^z t^{\alpha-1} e^{-t}\,dt \quad\hbox{and}\quad \Gamma(\alpha) = \int_0^\infty t^{\alpha-1} e^{-t}\,dt,
$$
for $\alpha \in \CC$ with $\Real \alpha > 0$ and $z \in \CC$, and so
\begin{align*}
\int_{t_0}^t (t - s)^{-\beta} e^{-a(t - s)}\,ds &=  a^{1+\beta}\int_0^{a(t-t_0)} r^{(1-\beta) - 1} e^{-r}\,dr
\\
&= a^{1+\beta} \gamma(1-\beta, a(t-t_0)),
\end{align*}
and thus,
\begin{equation}
\label{eq:Incomplete_gamma_function_inequality}
\int_{t_0}^t (t - s)^{-\beta} e^{-a(t - s)}\,ds \leq a^{1+\beta}\Gamma(1-\beta), \quad\forall\, t \geq t_0.
\end{equation}
Note also that
\begin{align*}
\int_s^T (t - s)^{-\beta} e^{-a(t - s)}\,dt &= a^{1+\beta}\int_0^{a(T-s)} r^{(1-\beta) - 1} e^{-r}\,dr
\\
&= a^{1+\beta} \gamma(1-\beta, a(T-s))
\\
&\leq a^{1+\beta} \Gamma(1-\beta), \quad \forall\, T \geq s.
\end{align*}
By integrating with respect to $t\in[t_0,T]$, we see that
\begin{align*}
\int_{t_0}^T \|\cA^\beta u(t)\|_\cW \,dt &\leq M_\beta \|u_0\|_\cW \int_{t_0}^T (t - t_0)^{-\beta} e^{-a(t - t_0)}\,dt
\\
&\quad + M_\beta \int_{t_0}^T \int_{t_0}^t (t - s)^{-\beta} e^{-a(t - s)}
\left(\|f_0(s)\|_\cW + K\|u(s)\|_\cW + \eps\|u(s)\|_{\calV^{2\beta}}\right) \,ds\,dt
\\
&= M_\beta a^{1+\beta} \gamma(1-\beta, a(T-t_0)) \|u_0\|_\cW
\\
&\quad + M_\beta \int_{t_0}^T \left(\|f_0(s)\|_\cW + K\|u(s)\|_\cW + \eps\|u(s)\|_{\calV^{2\beta}}\right)
\int_s^T (t - s)^{-\beta} e^{-a(t - s)}\,dt\,ds
\\
&\qquad\text{(by Fubini's Theorem)}
\\
&= a^{1+\beta} M_\beta \gamma(1-\beta, a(T-t_0)) \|u_0\|_\cW
\\
&\quad + M_\beta \int_{t_0}^T \left(\|f_0(s)\|_\cW + K\|u(s)\|_\cW + \eps\|u(s)\|_{\calV^{2\beta}}\right)
a^{1+\beta} \gamma(1-\beta, a(T-s))\,ds,
\end{align*}
and thus
\begin{align*}
\int_{t_0}^T \|\cA^\beta u(t)\|_\cW \,dt &\leq a^{1+\beta} M_\beta \Gamma(1-\beta) \|u_0\|_\cW
\\
&\quad + a^{1+\beta} M_\beta \Gamma(1-\beta) \int_{t_0}^T
\left(\|f_0(s)\|_\cW + K\|u(s)\|_\cW + \eps\|u(s)\|_{\calV^{2\beta}}\right) \,ds.
\end{align*}
The conclusion follows by rearrangement, noting that $\|\cA^\beta u(t)\|_\cW  = \|u(t)\|_{\calV^{2\beta}}$ for all $t \in [t_0, T]$.
\end{proof}

\begin{rmk}[Inequalities]
The inequality \eqref{eq:Incomplete_gamma_function_inequality} is useful when $t \gg t_0$ and we seek an estimate of the integral which is independent of $t-t_0 \gg 1$; when $t-t_0$ is small and positive, then the following simpler inequality
$$
\int_{t_0}^t (t - s)^{-\beta} e^{-a(t - s)}\,ds \leq \int_{t_0}^t (t - s)^{-\beta} \,ds
= \frac{(t-t_0)^{1-\beta}}{1-\beta},
$$
usually suffices.
\end{rmk}

We have the following corollary of Lemma \ref{lem:Rade_7-3_abstract_L1_in_time_V2beta_space}.

\begin{lem}[A global $L^1$-in-time \apriori estimate for the time derivative of a solution to a nonlinear evolution equation]
\label{lem:Rade_7-3_abstract_global_L1_in_time_V2beta_space_time_derivative}
Assume that Hypothesis \ref{hyp:Sell_You_4_standing_hypothesis_A} holds. Let $\beta \in [0, 1)$ and $\eps \in (0,1]$ be such that
$$
\eps a^{1+\beta} M_\beta \Gamma(1-\beta) \leq \frac{1}{2},
$$
where the constants $a > 0$ and $M_\beta > 0$ are as in Theorem \ref{thm:Sell_You_37-5}, and $K \in [1,\infty)$. Let $\cF \in C_{\Lip}([t_0,T)\times\calV^{2\beta}; \cW)$ for $t_0 \in \RR$ and $T$ obeying $t_0 < T \leq \infty$. Let $u \in C([t_0, T); \cW) \cap C((t_0, T); \calV^{2\beta})$ be a mild solution to the nonlinear evolution equation \eqref{eq:Sell_You_47-1} on $[t_0, T)$ defined by $\cA$ and $\cF$,
\[
\frac{du}{dt}(t) + \cA u(t) = \cF(t, u(t)).
\]
Furthermore, suppose that $\dot u \in C([t_0, T); \cW) \cap C((t_0, T); \calV^{2\beta})$ and that $\dot u$ is a mild solution to the nonlinear evolution equation on $[t_0, T)$,
\begin{equation}
\label{eq:Sell_You_47-1_time_derivative}
\frac{d\dot u}{dt}(t) + \cA\dot u(t) = \cG(t,\dot u(t)),
\end{equation}
where
\[
\cG(t,\dot u(t)) := \frac{d}{\partial t}\cF(t,u(t))
\]
can be expressed as
\[
\cG(t,v) = g_0(t) + \cG_1(t,v) + \cG_2(t,v), \quad\forall\, (t,v) \in [t_0,T) \times \calV^{2\beta},
\]
with $g_0 = \dot f_0 \in C_{\Lip}([t_0,T); \cW)$ and $\cG_1, \cG_2 \in C_{\Lip}([t_0,T)\times\calV^{2\beta}; \cW)$ and
\begin{align*}
\|\cG_1(t,v)\|_\cW &\leq K\|v\|_\cW,
\\
\|\cG_2(t,v)\|_\cW &\leq \eps\|v\|_{\calV^{2\beta}}, \quad\forall\, (t,v) \in [t_0,T) \times \calV^{2\beta}.
\end{align*}
Then
\begin{equation}
\label{eq:Rade_7-3_abstract_L1_in_time_V2beta_space_apriori_estimate_time_derivative}
\int_{t_0}^T \|\dot u(t)\|_{\calV^{2\beta}} \,dt \leq a^{1+\beta} M_\beta \Gamma(1-\beta) \left(\|\dot u_0\|_\cW
+ \int_{t_0}^T \|\dot f_0(t)\|_\cW\,dt + K\int_{t_0}^T \|\dot u(t)\|_\cW \,dt\right).
\end{equation}
\end{lem}

\begin{proof}
The conclusion follows by applying Lemma \ref{lem:Rade_7-3_abstract_L1_in_time_V2beta_space} to the nonlinear evolution equation \eqref{eq:Sell_You_47-1_time_derivative} with $\cG(t,v)$ in place of $\cF(t,v)$.
\end{proof}

\subsection{Interior $L^1$-in-time a priori estimates for a solution to a nonlinear evolution equation}
\label{subsec:Interior_L1-in-time_apriori_estimates_solution_nonlinear_evolution_equation}
In this subsection, we establish abstract interior \apriori estimates which complement the global estimates in Section \ref{subsec:Global_L1-in-time_apriori_estimates_solution_nonlinear_evolution_equation}; these estimates are especially useful when we have limited information regarding initial data.

\begin{lem}[An interior $L^1$-in-time \apriori estimate for a strong solution to a nonlinear evolution equation]
\label{lem:Rade_7-3_abstract_L1_in_time_V2beta_space_interior}
Assume the hypotheses of Lemma \ref{lem:Rade_7-3_abstract_L1_in_time_V2beta_space}. If $\delta$ is a constant obeying $t_0 + 2\delta \leq T$ and $u$ is a strong solution to the nonlinear evolution equation \eqref{eq:Sell_You_47-1}, then
\begin{equation}
\label{eq:Rade_7-3_abstract_L1_in_time_V2beta_space_apriori_estimate_interior}
\int_{t_0+\delta}^T \|u(t)\|_{\calV^{2\beta}} \,dt \leq a^{1+\beta} M_\beta \Gamma(1-\beta)
\left(\int_{t_0}^T \|f_0(t)\|_\cW\,dt + K(1+\delta^{-1})\int_{t_0}^T \|u(t)\|_\cW \,dt\right).
\end{equation}
\end{lem}

\begin{proof}
We define a smooth cutoff function, $\zeta\in C^\infty(\RR;[0,1])$, such that
\[
\zeta(t) = 1 \quad\forall\, t \geq t_0 + \delta \quad\text{and}\quad \zeta(t) = 0 \quad\forall\, t \leq t_0.
\]
We may construct $\zeta$ by setting
\[
\zeta(t) := \kappa\left(\delta^{-1}(t-t_0)\right) \quad (m\geq 0),
\]
where $\kappa \in C^\infty(\RR;[0,1])$ obeys $\kappa(t) = 0$ for $t\leq 0$ and $\kappa(t) = 1$ for $t\geq 1$. We observe that
\[
|\dot\zeta(t)| \leq N\delta^{-1}, \quad\forall\, t \in \RR,
\]
where $N \geq 1$ is a universal constant (depending at most on the choice of $\kappa$). Because $u$ is a strong solution to \eqref{eq:Sell_You_47-1} on $(t_0,T)$, we see that $\zeta u$  is a strong solution to the following nonlinear evolution equation on $(t_0,T)$,
\[
\frac{d(\zeta u)}{dt} + \cA(\zeta u) = \zeta\cF(t, u) + \dot\zeta u.
\]
Denote
\[
\cH(t,\zeta u) := \zeta\cF(t,u) + \dot\zeta u,
\]
and write
\[
\cH(t,\zeta u) := h_0(t) + \cH_1(t,\zeta u) + \cH_2(t,\zeta u)
\]
where
\[
h_0(t) := \zeta f_0(t), \quad \cH_1(t,\zeta u) := \zeta\cF_1(t,u) + \dot\zeta u, \quad \cH_1(t,\zeta u) := \zeta\cF_2(t,u).
\]
By the hypotheses in Lemma \ref{lem:Rade_7-3_abstract_L1_in_time_V2beta_space} on $\cF$, we obtain
\begin{align*}
\|\cH_1(t,\zeta v)\|_\cW &\leq K\|\zeta v\|_\cW + |\dot\zeta|\|v\|_\cW
\\
&\leq K\|\zeta v\|_\cW + N\delta^{-1}\|v\|_\cW,
\\
\|\cH_2(t,\zeta v)\|_\cW &\leq \eps\|\zeta v\|_{\calV^{2\beta}},
\quad\forall\, (t,v) \in (t_0,T)\times\calV^{2\beta}.
\end{align*}
Consequently, the method of proof of the \apriori estimate \eqref{eq:Rade_7-3_abstract_L1_in_time_V2beta_space_apriori_estimate} in the Lemma \ref{lem:Rade_7-3_abstract_L1_in_time_V2beta_space} applied \mutatis to the mild solution $\zeta u$ on $(t_0,T)$ now yields
\[
\int_{t_0}^T \|\zeta(t) u(t)\|_{\calV^{2\beta}} \,dt
\leq
a^{1+\beta} M_\beta \Gamma(1-\beta) \left(\int_{t_0}^T \|f_0(t)\|_\cW\,dt
+ \left(K+N\delta^{-1}\right)\int_{t_0}^T \|\zeta(t) u(t)\|_\cW \,dt \right).
\]
The conclusion thus follows from the fact that $\zeta = 1$ on $(t_0+\delta,T)$.
\end{proof}

Similarly, it will be useful to have an \apriori interior estimate analogue of Lemma \ref{lem:Rade_7-3_abstract_global_L1_in_time_V2beta_space_time_derivative}.

\begin{lem}[An interior $L^1$-in-time \apriori estimate for the time derivative of a solution to a nonlinear evolution equation]
\label{lem:Rade_7-3_abstract_interior_L1_in_time_V2beta_space_time_derivative_interior}
Assume the hypotheses of Lemma \ref{lem:Rade_7-3_abstract_global_L1_in_time_V2beta_space_time_derivative}. If $\delta$ is a constant obeying $t_0 + 2\delta \leq T$ and $\dot u$ is a strong solution to the nonlinear evolution equation \eqref{eq:Sell_You_47-1_time_derivative}, then
\begin{equation}
\label{eq:Rade_7-3_abstract_L1_in_time_V2beta_space_apriori_estimate_time_derivative_interior}
\int_{t_0+\delta}^T \|\dot u(t)\|_{\calV^{2\beta}} \,dt \leq a^{1+\beta} M_\beta \Gamma(1-\beta)
\left( \int_{t_0}^T \|\dot f_0(t)\|_\cW\,dt + K(1+\delta^{-1})\int_{t_0}^T \|\dot u(t)\|_\cW \,dt\right).
\end{equation}
\end{lem}

\begin{proof}
The conclusion follows by replacing the role of Lemma \ref{lem:Rade_7-3_abstract_L1_in_time_V2beta_space} by that of Lemma \ref{lem:Rade_7-3_abstract_global_L1_in_time_V2beta_space_time_derivative}  in the proof of Lemma \ref{lem:Rade_7-3_abstract_L1_in_time_V2beta_space_interior}.
\end{proof}

\section[Critical-exponent parabolic Sobolev spaces]{Critical-exponent parabolic Sobolev spaces and linear parabolic operators on sections of vector bundles over compact manifolds}
\label{sec:Critical-exponent_parabolic_Sobolev_spaces_linear_parabolic_operator_vector_bundle_manifold}
In this section, we shall adapt our development in \cite[Sections 4 and 5]{FeehanSlice} of critical-exponent Sobolev spaces for elliptic operators on sections of vector bundles over closed four-dimensional manifolds, which is based on ideas of Taubes \cite{TauPath, TauFrame, TauStable, TauConf, TauGluing}, to the case of parabolic operators.

We let $X$ be a closed, oriented, Riemannian manifold of dimension $d \geq 2$ and let $K(t,x,y)$ denote the heat kernel, that is, the fundamental solution for the (augmented) heat operator \cite{Berger_Gauduchon_Mazet_1971, Chavel, Davies_1989, Grigoryan_2009},
\begin{equation}
L := \frac{\partial}{\partial t} + \Delta + 2  \quad\hbox{on } C^\infty((0, \infty) \times X),
\end{equation}
where $\Delta = d^*d$ is the Laplace operator on $C^\infty(X)$. See \cite[Definition 6.1]{Chavel} for a precise definition of the heat kernel. We recall that the heat kernel for $\partial_t + \Delta$ on $C^\infty((0, \infty) \times \RR^d)$ is given by
$$
\frac{1}{(4\pi t)^{d/2}}e^{-|x-y|^2/4t}, \quad\forall\, (t, x, y) \in (0, \infty) \times \RR^d \times \RR^d.
$$
As a consequence, for example, of the parametrix construction of $K(t, x, y)$ on $(0, \infty) \times X \times X$ when $X$ is a closed, Riemannian manifold \cite[Section 6.4]{Chavel}, one finds that the heat kernel, $K(t, x, y)$, is comparable to
$$
\frac{1}{(4\pi t)^{d/2}}e^{-\dist_g^2(x, y)/4t},
$$
for $t > 0$, when $x, y \in X$ and $\dist_g(x, y)$ is small relative to the injectivity radius of $(X, g)$. By analogy with our definitions \eqref{eq:BasicSharpNorms} of the $L^\sharp$ and $L^{2\sharp}$ norms of $u \in C^\infty(X; E)$, we define, for $u \in C^\infty([0, T]\times X; E)$ and $0 < T \leq \infty$,
\begin{subequations}
\label{eq:Critical-exponent_parabolic_Sobolev_spaces}
\begin{align}
\label{eq:Critical-exponent_parabolic_Sobolev_space_Llozenge_0_to_T_times_X}
\|u\|_{L^\lozenge((0, T)\times X; E)} &:= \sup_{(t,x)\in(0,T)\times X} \int_0^t\int_X K(t-s, x, y)|u|(s, y)\,d\vol_g(y) \,ds,
\\
\label{eq:Critical-exponent_parabolic_Sobolev_space_L2lozenge_0_to_T_times_X}
\|u\|_{L^{2\lozenge}((0, T)\times X; E)} &:= \sup_{(t,x)\in(0,T)\times X} \left(\int_0^t\int_X K(t-s, x, y)|u(s, y)|^2\,d\vol_g(y) \,ds\right)^{1/2}.
\end{align}
\end{subequations}
Of course, the preceding two norms are related by
$$
\|u\|_{L^{2\lozenge}((0, T)\times X; E)} = \left(\| |u|^2 \|_{L^\lozenge((0, T)\times X)}\right)^{1/2}.
$$
We define the (augmented) heat operator on $C^\infty(X; E)$ by
\begin{equation}
\label{eq:Augmented_heat_operator_on_sections_vectorbundle_over_manifold}
L_A := \frac{\partial}{\partial t} + \nabla_A^*\nabla_A + 1 \quad\hbox{on } C^\infty((0, \infty) \times X; E),
\end{equation}
and give a parabolic analogue of Lemma \ref{lem:Feehan_5-3}. For brevity, we denote $K_{t,x} = K(t,x,\cdot)$, for all $(t,x) \in (0,T)\times X$.

The well-known pointwise identity \cite[Equation (6.18)]{FU},
\begin{equation}
\label{eq:Freed_Uhlenbeck_6-18}
|\nabla_Au|^2 + \frac{1}{2} \Delta|u|^2
=
\langle\nabla_A^*\nabla_A u,u\rangle \quad\hbox{a.e. on } X,
\end{equation}
provides the key technical ingredient in the following lemma.

\begin{lem}
\label{lem:Feehan_5-3_heat_operator}
Let $X$ be a $C^\infty$ closed, oriented
manifold of dimension $d \geq 2$ with Riemannian metric $g$, and $E$ a complex Hermitian (or real Riemannian) vector bundle over $X$, and $A$ a $C^\infty$ Hermitian (Riemannian) connection on $E$ with curvature $F_A$. If $u \in C^\infty(X; E)$, then
\begin{multline}
\label{eq:Feehan_5-3_heat_operator}
\frac{1}{2}\|\nabla_Au\|_{L^{2\lozenge}((0, T)\times X; E)} + \frac{1}{4}\|u\|_{C([0, T]\times X; E)}
\\
\leq
2\|L_Au\|_{L^\lozenge((0, T)\times X)}
+ \frac{1}{\sqrt{2}}\sup_{(t, x) \in (0,T)\times X}\|u(0,\cdot)\|_{L^2(X, K_{t,x}; E)}.
\end{multline}
\end{lem}

\begin{proof}
The identity \eqref{eq:Freed_Uhlenbeck_6-18} gives
\begin{align*}
|\nabla_Au|^2 + \frac{1}{2}\left(\frac{\partial}{\partial t} + \Delta + 2\right)|u|^2
&= \frac{1}{2}\frac{\partial}{\partial t}|u|^2 + \langle\nabla_A^*\nabla_A u,u\rangle + |u|^2
\\
&= \left\langle\frac{\partial u}{\partial t} + \nabla_A^*\nabla_A u + u, u \right\rangle \quad\hbox{a.e. on } X,
\end{align*}
that is,
$$
|\nabla_Au|^2 + \frac{1}{2}L(|u|^2) = \langle L_Au, u \rangle \quad\hbox{a.e. on } X.
$$
The variation of constants formula \cite[Equation (6)]{Chavel}, \cite[Equation (42.3)]{Sell_You_2002} is
\begin{align*}
u(t, x) &= e^{-(\Delta + 2)t}u(0, x) + \int_0^t e^{-(\Delta + 2)(t-s)} Lu(s, x) \,ds
\\
&= \int_X K(t, x, y)u(0, y)\,d\vol_g(y) + \int_0^t \int_X K(t-s, x, y) Lu(s, y)\,d\vol_g(y) \,ds, \quad\forall\, t \geq 0, \ x\in X,
\end{align*}
where the action of the semigroup, $e^{-(\Delta + 2)t}$, on $C^\infty(X)$, is represented by the heat kernel, $K(t,x,y)$. By combining the preceding identities, we obtain
\begin{align*}
{}&\int_0^t \int_X K(t-s, x, y)|\nabla_Au(y)|^2\,d\vol_g(y) + \frac{1}{2}\int_0^t  \int_X K(t, x, y) (L|u|^2)(s, y)\,d\vol_g(y) \,ds
\\
&= \int_0^t \int_X K(t-s, x, y)|\nabla_Au(y)|^2\,d\vol_g(y) + \frac{1}{2}|u(t, x)|^2 - \frac{1}{2}\int_X K(t, x, y)|u(0, y)|^2\,d\vol_g(y)
\\
&= \int_0^t\int_X K(t-s, x, y)\langle L_Au, u\rangle(s,y)\,d\vol_g(y) \,ds,  \quad\forall\, t \geq 0, \ x\in X.
\end{align*}
Therefore, for any $T>0$,
\begin{align*}
{}&\sup_{(t,x)\in(0,T)\times X} \int_0^t \int_X K(t-s, x, y)|\nabla_Au(y)|^2\,d\vol_g(y) + \frac{1}{2}\sup_{(t,x)\in(0,T)\times X}|u(t, x)|^2
\\
&\leq \sup_{(t,x)\in(0,T)\times X} \int_0^t\int_X K(t-s, x, y)|L_Au(s, y)| |u(s, y)|\,d\vol_g(y) \,ds
\\
&\quad + \frac{1}{2}\sup_{(t,x)\in(0,T)} \int_X K(t, x, y)|u(0, y)|^2\,d\vol_g(y)
\\
&\leq \left(\sup_{(t,x)\in(0,T)\times X} \int_0^t\int_X K(t-s, x, y)|L_Au(s, y)|\,d\vol_g(y) \,ds\right)\left(\sup_{(s, y)\in(0,T)\times X}|u(s,y)|\right)
\\
&\quad + \frac{1}{2}\sup_{(t,x)\in(0,T)} \int_X K(t, x, y)|u(0, y)|^2\,d\vol_g(y).
\end{align*}
Consequently, using $ab \leq \eps a^2 + \eps^{-1}b^2$, for any $a, b \geq 0$ and $\eps > 0$ \cite[Equation (7.6)]{GilbargTrudinger} with $\eps = 1/4$ and rearrangement, we see that
\begin{align*}
{}&\sup_{(t,x)\in(0,T)\times X} \int_0^t \int_X K(t-s, x, y)|\nabla_Au(y)|^2\,d\vol_g(y) + \frac{1}{4}\sup_{(t,x)\in(0,T)\times X}|u(t, x)|^2
\\
&\leq \sup_{(t,x)\in(0,T)\times X} 4\left(\int_0^t\int_X K(t-s, x, y)|L_Au(s, y)|\,d\vol_g(y) \,ds\right)^2
\\
&\quad + \frac{1}{2}\sup_{(t,x)\in(0,T)} \int_X K(t, x, y)|u(0, y)|^2\,d\vol_g(y).
\end{align*}
Taking square roots and using $(a + b)/2 \leq (a^2 + b^2)^{1/2}$ and $(a^2 + b^2)^{1/2} \leq a + b$, for any $a, b \geq 0$, now yields the conclusion.
\end{proof}

We recall the following standard \apriori estimate \eqref{eq:Sell_You_42-30_heat_equation_alpha_is_one} from Theorem \ref{thm:Sell_You_42-12_heat_equation_alpha_is_one},
\begin{multline}
\label{eq:Standard_H2_parabolic_linear_apriori_estimate}
\|u\|_{C([0,T]; H_A^1(X; E))} + \|u\|_{L^2(0,T; H_A^2(X;E)} + \|\partial_tu\|_{L^2(0,T; L^2(X; E))}
\\
\leq 4\|u(0,\cdot)\|_{H_A^1(X;E)} + 6\|L_A u\|_{L^2(0,T; L^2(X; E))}.
\end{multline}
In view of Lemma \ref{lem:Feehan_5-3_heat_operator} and the
the preceding standard \apriori estimate, we obtain the following strengthened \apriori estimate, valid on $X$ of any dimension $d \geq 2$,
\begin{multline}
\label{eq:Critical_exponent_plus_standard_H2_parabolic_linear_apriori_estimate}
\|u\|_{C([0, T]; C(X; E))} + \|u\|_{C([0, T]; H_A^1(X; E))}
\\
+ \|\nabla _Au\|_{L^{2\lozenge}((0, T)\times X)} + \|u\|_{L^2(0, T; H_A^2(X; E))}
+ \|\partial_t u\|_{L^2(0, T; L^2(X; E))}
\\
\leq 8\|L_Au\|_{L^2(0, T; L^2(X; E))} + 8\|L_Au\|_{L^\lozenge((0, T)\times X; E)}
\\
+ 4\|u(0,\cdot)\|_{H_A^1(X;E)} + \sup_{(t, x) \in (0,T)\times X} 4\|u(0,\cdot)\|_{L^2(X, K_{t,x}; E)}.
\end{multline}
We now define the following parabolic analogues of our replacements for $L^2(X; E)$ and $H_A^2(X; E) \equiv W^{2, 2}_A(X; E)$ by the Banach spaces $L^{2+\sharp}(X; E)$ and $W^{2, 2+\sharp}_A(X; E)$. We let $\fW$ and $\fV$ denote the Banach spaces obtained by completing $C^\infty([0, T]\times X; E)$ with respect to the \emph{hybrid critical-exponent parabolic Sobolev space} norms,
\begin{subequations}
\label{eq:Parabolic_Sobolev_norms_hybrid_standard_and_critical_exponent}
\begin{align}
\label{eq:Parabolic_Sobolev_norm_hybrid_standard_L2_and_critical_exponent_range}
\|w\|_\fW &= \|w\|_{L^2(0, T; L^2(X; E))} +  \|w\|_{L^\lozenge((0, T)\times X; E)}
\\
\label{eq:Parabolic_Sobolev_norm_hybrid_standard_H2_and_critical_exponent_domain}
\|v\|_\fV &= \|u\|_{C([0, T]; C(X; E))} + \|u\|_{C([0, T]; H_A^1(X; E))} + \|\nabla _Au\|_{L^{2\lozenge}((0, T)\times X; E)}
\\
\notag
&\qquad  + \|L_Au\|_{L^2(0, T; L^2(X; E))} + \|L_Au\|_{L^\lozenge((0, T)\times X; E)}.
\end{align}
\end{subequations}
We recall that the heat kernel, $K$, belongs to $C^\infty((0,\infty)\times X \times X)$ and is everywhere positive, $K(t,x,y) > 0$ for all $(t,x,y) \in (0,\infty)\times X \times X$, by \cite[Theorem 5.2.1]{Davies_1989}, so the terms involving $K$ on the right-hand side in \eqref{eq:Parabolic_Sobolev_norms_hybrid_standard_and_critical_exponent} themselves define norms. Our \apriori estimate \eqref{eq:Critical_exponent_plus_standard_H2_parabolic_linear_apriori_estimate} now takes the more compact form,
\begin{equation}
\label{eq:Critical_exponent_plus_standard_H2_parabolic_linear_apriori_estimate_streamlined}
\|u\|_\fV \leq 8\|L_Au\|_\fW + 4\|u(0,\cdot)\|_{H_A^1(X;E)} + \sup_{(t, x) \in (0,T)\times X} 4\|u(0,\cdot)\|_{L^2(X, K_{t,x}; E)}.
\end{equation}
We recall from \cite[Theorem 5.2.6]{Davies_1989} that $e^{-\Delta t} 1 = 1$ on $X$ for all $t \geq 0$ (the `conservation of probability'), while $e^{-(\Delta + 2)t} = e^{-2t} e^{-\Delta t}$ on $C(X)$ or $L^p(X)$, where $1 \leq p < \infty$. In particular, $e^{-(\Delta + 2)t} 1 = e^{-2t}$ on $X$. If $K_0(t,x,y)$ is the kernel for $e^{-\Delta t}$, then $K(t,x,y) = e^{-2t} K_0(t,x,y)$ and
$$
\int_X K(t, x, y)\,d\vol_g(y) =  e^{-2t}\int_X K_0(t, x, y)\,d\vol_g(y) = e^{-2t}, \quad\forall\, t \geq 0 \hbox{ and } x \in X.
$$
These observations leads to useful inequalities for the norms in \eqref{eq:Critical-exponent_parabolic_Sobolev_spaces}, which are direct consequences of their definitions. We calculate
\begin{align*}
\|u\|_{L^\lozenge((0,T)\times X; E)}
&\leq
\sup_{(t,x)\in (0,T)\times X} \left(\int_0^t\int_X K(t-s, x, y)\,d\vol_g(y) \,ds\right)
\|u\|_{C([0,T]\times X; E)}
\\
&= \sup_{t\in(0,T)} \left(\int_0^t e^{s-t}\,ds\right) \|u\|_{C([0,T]\times X; E)},
\\
\|u\|_{L^{2\lozenge}((0, T)\times X; E)}
&\leq
\sup_{(t,x)\in(0,T)\times X} \left(\int_0^t\int_X K(t-s, x, y)\,d\vol_g(y) \right)^{1/2}
\|u\|_{C([0,T]\times X; E)}
\\
&=  \sup_{t\in(0,T)}\left(\int_0^t e^{s-t}\,ds \right)^{1/2} \|u\|_{C([0,T]\times X; E)},
\end{align*}
and therefore,
\begin{subequations}
\label{eq:Sobolev_embedding_Linfty_time_and_space_Llozenge_or_L2lozenge}
\begin{align}
\label{eq:Sobolev_embedding_Linfty_time_and_space_Llozenge}
\|u\|_{L^\lozenge((0,T)\times X; E)} &\leq T\|u\|_{C([0,T]\times X; E)},
\\
\label{eq:Sobolev_embedding_Linfty_time_and_space_into_L2lozenge}
\|u\|_{L^{2\lozenge}((0, T)\times X; E)} &\leq \sqrt{T}\|u\|_{C([0,T]\times X; E)}.
\end{align}
\end{subequations}
The norm on the initial data $u(0,\cdot) = u_0$ on $X$ in the \apriori estimate \eqref{eq:Critical_exponent_plus_standard_H2_parabolic_linear_apriori_estimate} can be bounded by its $C^0$ norm as follows:
\begin{align*}
\sup_{(t, x) \in (0,T)\times X}\|u_0\|_{L^2(X, K_{t,x}; E)}
&\equiv
\sup_{(t, x) \in (0,T)\times X} \left( \int_X K(t, x, y)|u_0(y)|^2\,d\vol_g(y) \right)^{1/2}
\\
&\quad \leq \sup_{(t, x) \in (0,T)\times X} \left( \int_X K(t, x, y)\,d\vol_g(y) \right)^{1/2}
\|u_0\|_{C(X; E)}
\\
&\quad = e^{-2T}\|u_0\|_{C(X; E)},
\end{align*}
that is,
\begin{equation}
\label{eq:Sobolev_embedding_Linfty_time_and_space_into_sup_time_and_space_L2space_dKtx}
\sup_{(t, x) \in (0,T)\times X}\|u_0\|_{L^2(X, K_{t,x}; E)} \leq e^{-2T}\|u_0\|_{C(X; E)}.
\end{equation}
Recall from \cite[Theorem 4.12]{AdamsFournier} that $W^{2\beta,p}(X) \hookrightarrow C(X)$ when $2\beta p > d$, where $d \geq 2$ and $1 < p < \infty$ and $\beta \in (0, 1)$. Therefore,
\begin{align*}
\|u\|_{L^\lozenge((0,T)\times X; E)}
&=
\sup_{(t,x)\in (0,T)\times X} \int_0^t e^{(\Delta + 2)(t-s)}|u|(s,x) \,ds
\\
&\leq \sup_{t \in (0,T)} \int_0^t \| e^{(\Delta + 2)(t-s)}|u|(s,\cdot) \|_{L^\infty(X)} \,ds
\\
&\leq \sup_{t \in (0,T)} C_p\int_0^t \| e^{(\Delta + 2)(t-s)}|u|(s,\cdot) \|_{W^{2\beta,p}(X)} \,ds
\\
&\leq \sup_{t \in (0,T)} C_p\int_0^t \| (\Delta + 2)^\beta e^{(\Delta + 2)(t-s)}|u|(s,\cdot) \|_{L^p(X)} \,ds
\\
&\leq \sup_{t \in (0,T)} C_{p,\beta} \int_0^t (t-s)^{-\beta} e^{2(t-s)} \|u(s,\cdot)\|_{L^p(X; E)} \,ds
\quad\hbox{(by Theorem \ref{thm:Sell_You_37-5})}
\\
&\leq \frac{C_{p,\beta} T^\beta}{1-\beta} \sup_{t \in (0,T)} \|u(t,\cdot)\|_{L^p(X; E)},
\end{align*}
with a similar calculation for the $L^{2\lozenge}$ norm. We conclude that, for all $p \in (1, \infty)$ and $\beta \in (0, 1)$ with $2\beta p > d$, we have
\begin{align}
\label{eq:Sobolev_embedding_Linfty_time_Lp_space_into_Llozenge}
\|u\|_{L^\lozenge((0,T)\times X; E)} &\leq \frac{C_{p,\beta} T^\beta}{1-\beta} \|u\|_{C([0, T]; L^p(X; E))},
\\
\label{eq:Sobolev_embedding_Linfty_time_L2p_space_L2lozenge}
\|u\|_{L^{2\lozenge}((0,T)\times X; E)} &\leq \left(\frac{C_{p,\beta} T^\beta}{1-\beta}\right)^{1/2} \|u\|_{C([0, T]; L^{2p}(X; E))}.
\end{align}
Note that, because $\beta \in (0, 1)$, we \emph{must} select $p > d/2$ for these estimates to be valid.

\section[Local well-posedness for the Yang-Mills heat equation]{Local well-posedness for the Yang-Mills heat equation with initial data of minimal regularity}
\label{sec:Struwe_3and4}
In Section \ref{sec:Local_well-posedness_yang_mills_heat_equation}, we obtained well-posedness for strong solutions in $W^{2\beta, p}(X; \Lambda^1\otimes\ad P)$ to the Yang-Mills heat equation \eqref{eq:Yang-Mills_heat_equation_as_perturbation_rough_Laplacian_plus_one_heat_equation}
over a closed Riemannian manifold, $X$, of arbitrary dimension $d\geq 2$ by appealing to the abstract theory of nonlinear evolution equations and semigroup theory via the fractional-order Banach spaces, $\calV^{2\beta}$, for $\beta \in (0,1)$ described, for example, by Sell and You in \cite[Sections 4.6 and 4.8]{Sell_You_2002}. Kozono, Maeda, and Naito \cite{Kozono_Maeda_Naito_1995} also use semigroup theory, but we obtain our results more easily by relying on the more developed treatment in \cite{Sell_You_2002}, although our more abstract treatment leads to the requirement that the initial data, $a_0$, belong to $W^{2\beta, p}(X; \Lambda^1\otimes\ad P)$ rather than $L^d(X; \Lambda^1\otimes\ad P)$ as in \cite{Kozono_Maeda_Naito_1995}. For example, when $\calV = W^{2, p}(X; \Lambda^1\otimes\ad P)$, then Theorems
\ref{thm:Existence_uniqueness_mild_solution_Yang-Mills_heat_equation_in_W_2beta_p_initial_data_in_W_2beta_p} and
\ref{thm:Existence_uniqueness_strong_solution_Yang-Mills_heat_equation_in_W_2beta_p_initial_data_in_W_2beta_p}
require us to choose $\beta \in [1/2, 1)$ when $p > d$ or $\beta \in (1/2, 1)$ when $p < d$.

An alternative and arguably more elementary approach, which does not require semigroup theory but is still flexible enough to handle the Yang-Mills heat equation nonlinearity \eqref{eq:Yang-Mills_heat_equation_nonlinearity_relative_rough_Laplacian_plus_one},
can also be based on the $L^p$ theory for linear parabolic equations and the contraction mapping principle to handle the nonlinearity. Existence, uniqueness, and \apriori $L^p$ estimates are provided in modern sources such as those of Maugeri, Palagachev, and Softova \cite{Maugeri_Palagachev_Softova_2000}, for parabolic scalar equations, or Dong and Kim \cite{Dong_Kim_2011arma, Dong_Kim_2011cvpde} for parabolic scalar systems, all in the case of domains $\Omega \subseteqq \RR^d$, while
Lady{\v{z}}enskaja, Solonnikov, and Ural$'$ceva \cite{LadyzenskajaSolonnikovUralceva} provides a classic, if older source. For example, \cite[Theorem 2.5.1]{Maugeri_Palagachev_Softova_2000} gives (we temporarily adopt the notation of \cite{Maugeri_Palagachev_Softova_2000} for ease of reference)
\begin{equation}
\label{eq:Maugeri_Palagachev_Softova_2-93}
\|u\|_{W_p^{2,1}(Q_T)} \leq C\left(\|\sP u\|_{L^p(Q_T)} + \|\sB u\|_{W_p^{1 - 1/p, 1/2 - 1/2p}(S_T)} + \|\sI u\|_{W_p^{2-2/p}(\Omega)}\right),
\end{equation}
where $p \in (1, \infty)$ and $Q_T = (0,T)\times\Omega$ and $\Omega \Subset \RR^d$ is bounded domain with $C^{1,1}$ boundary and lateral boundary $S_T = (0,T)\times\partial\Omega $, and $W_p^{2,1}(Q_T)$ denotes the parabolic Sobolev space \cite[Equation (2.90)]{Maugeri_Palagachev_Softova_2000}
$$
\|u\|_{W_p^{2,1}(Q_T)} := \|u\|_{L^p(\Omega_T)} + \|\partial_t u\|_{L^p(\Omega_T)} + \sum_{i=1}^d\|u_{x_i}\|_{L^p(\Omega_T)} + \sum_{i,j=1}^d\|u_{x_ix_j}\|_{L^p(\Omega_T)},
$$
and scalar, uniformly and strictly elliptic, symmetric coefficients $a^{ij} \in \VMO(Q_T)$ and coefficients $\ell_i \in \Lip(\bar S_T)$ determine the operators \cite[Equations (2.89), (2.91), and (2.92)]{Maugeri_Palagachev_Softova_2000}
\begin{align*}
\sP u &:= \frac{\partial u}{\partial t} - \sum_{i,j=1}^d  a^{ij}u_{x_ix_j},
\\
\sB u &:= \sum_{i=1}^d \ell_i u_{x_i},
\\
\sI u &:= u(0, \cdot).
\end{align*}
The oblique derivative boundary and initial data, $\sB u$ and $\sI u$, belong to the (fractional derivative) \emph{Besov} spaces, $W_p^{1 - 1/p, 1/2 - 1/2p}(S_T)$ and $W_p^{2-2/p}(\Omega)$, respectively \cite[Section 2.5]{Maugeri_Palagachev_Softova_2000}. When $p = 2$, then the \apriori estimate \eqref{eq:Maugeri_Palagachev_Softova_2-93} is similar to the standard \apriori estimate \cite[Equation (7.46)]{Evans2} in \cite[Theorem 7.1.5]{Evans2} for a strong solution to a scalar parabolic equation (in divergence form), given homogeneous Dirichlet data on $S_T$ and initial data $u(\cdot, 0)$ in $H^1_0(\Omega)$.

However, as we shall explain, when $X$ has dimension four, an \apriori estimate such as \eqref{eq:Maugeri_Palagachev_Softova_2-93} with $p=2$ or \cite[Equation (7.46)]{Evans2} just falls short of what is required to prove existence of a strong solution
\begin{multline*}
a \in L^2(0,\tau; H^2_{A_1}(X; \Lambda^1\otimes\ad P)) \cap L^\infty(0,\tau; H^1_{A_1}(X; \Lambda^1\otimes\ad P))
\\
\cap H^1(0,\tau; L^2(X; \Lambda^1\otimes\ad P)).
\end{multline*}
to the Yang-Mills heat equation \eqref{eq:Yang-Mills_heat_equation_as_perturbation_rough_Laplacian_plus_one_heat_equation}, for initial data $a_0 \in H^1_{A_1}(X; \Lambda^1\otimes\ad P)$ and (small enough) $\tau > 0$, contrary to the assertion by Struwe \cite[Section 4.3]{Struwe_1994}. This failure is linked to the criticality of the Yang-Mills heat equation nonlinearity \eqref{eq:Yang-Mills_heat_equation_nonlinearity_relative_rough_Laplacian_plus_one} in dimension four. If the initial data is assumed to be more regular than supposed in \cite[Section 4]{Struwe_1994}, for example $a_0 \in H^2_{A_1}(X; \Lambda^1\otimes\ad P)$, we shall find that the contraction mapping argument then succeeds, as we can then appeal to an analogue of the more powerful \apriori estimate \cite[Equation (7.47)]{Evans2} in \cite[Theorem 7.1.5]{Evans2} for a strong solution
\begin{multline*}
a \in L^\infty(0,\tau; H^2_{A_1}(X; \Lambda^1\otimes\ad P)) \cap H^1(0,\tau; H^1_{A_1}(X; \Lambda^1\otimes\ad P))
\\
\cap W^{1,\infty}(0,\tau; L^2(X; \Lambda^1\otimes\ad P)) \cap H^2(0,\tau; H^{-1}_{A_1}(X; \Lambda^1\otimes\ad P)).
\end{multline*}
We shall describe an intermediate approach, based on our treatment for the model linear parabolic equation in Section
\ref{sec:Critical-exponent_parabolic_Sobolev_spaces_linear_parabolic_operator_vector_bundle_manifold} using critical-exponent parabolic Sobolev norms, with a condition on the regularity of the initial data which lies in between the assumption $a_0 \in H^1_{A_1}(X; \Lambda^1\otimes\ad P)$ of Struwe and the stronger assumption $a_0 \in H^2_{A_1}(X; \Lambda^1\otimes\ad P)$.

See \cite[p. 159]{Sell_You_2002} and \cite{Temam_1982} for a discussion of the potentially singular behavior of mild and strong solutions to linear and semilinear evolution equations at $t = 0$, due to mismatches between the regularity of the initial data and the desired regularity for positive time.

We also provide a modified version of Struwe's argument, yielding a slightly more regular solution (though at the cost of more regular initial data) using our system of critical-exponent parabolic Sobolev norms developed in Section \ref{sec:Critical-exponent_parabolic_Sobolev_spaces_linear_parabolic_operator_vector_bundle_manifold} or temporal weighting as in Section \ref{subsec:Struwe_page_137_contraction_mapping_initial_data_in_Ld}.

\subsection{$L^2(0,T;L^2(X))$ estimates for quadratic and cubic terms arising in the Yang-Mills heat equation}
\label{subsec:Struwe_section_4-3_L2_estimates_quadratic_and_cubic_terms_Yang-Mills_heat-equation}
We develop some preliminary estimates that we shall need in Section \ref{subsec:Struwe_page_137_contraction_mapping_small_initial_data_in_H1} and in the sequel. We adopt the notation and conventions of \cite[Section 4.3]{Struwe_1994}. We use a fixed $C^\infty$ reference connection, $A_1$ on $P$, in defining Sobolev spaces of sections of the vector bundle $\Lambda^1\otimes \ad P$ over a closed, Riemannian manifold, $X$, of dimension $d\geq 2$. Following \cite[p. 129]{Struwe_1994}, we consider the Hilbert space,
\begin{equation}
\label{eq:Struwe_page_129_Hilbert_space}
\sV := L^2(0,T; H^2_{A_1}(X; \Lambda^1\otimes\ad P)) \cap H^1(0,T; L^2(X; \Lambda^1\otimes\ad P)),
\end{equation}
and recall \cite[p. 130]{Struwe_1994} that there is a continuous embedding,
\begin{equation}
\label{eq:Struwe_page_130_L^2_time_H2_space_cap_H1_time_L2_space_embeds_into_C0_time_H1_space}
\sV \hookrightarrow C([0,T]; H^1_{A_1}(X; \Lambda^1\otimes\ad P)),
\end{equation}
with a universal embedding constant given by \cite[Equation (9)]{Struwe_1994},
\begin{equation}
\label{eq:Struwe_9}
\|a\|_{C([0,T]; H^1_{A_1}(X))}^2 \leq \|a(0)\|_{H_{A_1}^1(X)}^2 + 2\|a\|_\sV^2.
\end{equation}
Thus, it is convenient to consider the Banach space,
\begin{multline}
\label{eq:Struwe_page_129_Banach_space}
\fV := L^2(0,T; H^2_{A_1}(X; \Lambda^1\otimes\ad P)) \cap H^1(0,T; L^2(X; \Lambda^1\otimes\ad P))
\\
\cap C([0,T]; H^1_{A_1}(X; \Lambda^1\otimes\ad P)).
\end{multline}
We first estimate the quadratic terms appearing in \cite[p. 137]{Struwe_1994}, such as Struwe's terms VI and VII. For clarity, we write Struwe's $L^{2,2}$ in full as $L^2(0,T; L^2(X))$, namely
$$
\|f\|_{L^2(0,T; L^2(X))}^2 = \int_0^T \|f(t)\|_{L^2(X)}^2 \,dt = \int_0^T \int_X |f(t,x)|^2 \,dx\,dt.
$$
For a quadratic term like Struwe's term VI, we have
\begin{align*}
\|\nabla_{A_1}a_1\times a_2\|_{L^2(0,T; L^2(X))}^2
&= \int_0^T \|\nabla_{A_1}a_1(t)\times a_2(t)\|_{L^2(X)}^2 \,dt
\\
&\leq c\int_0^T \|\nabla_{A_1} a_1(t)\|_{L^4(X)}^2 \|a_2(t)\|_{L^4(X)}^2 \,dt
\\
&\leq c\left(\sup_{t \in (0,T)} \|a_2(t)\|_{L^4(X)}^2 \right)
\int_0^T \|\nabla_{A_1} a_1(t)\|_{L^4(X)}^2 \,dt,
\\
&\leq c\left(\sup_{t \in (0,T)} \|a_2(t)\|_{L^4(X)}^2 \right)
\int_0^T \|a_1(t)\|_{W_{A_1}^{1,4}(X)}^2 \,dt,
\end{align*}
that is, taking square roots,
\begin{equation}
\label{eq:Struwe_Yang-Mills_heat_equation_estimate_term_VI_W14}
\|\nabla_{A_1}a_1\times a_2\|_{L^2(0,T; L^2(X))}
\leq
c\|a_2\|_{L^\infty(0,T; L^4(X))} \|a_1\|_{L^2(0,T; W_{A_1}^{1,4}(X))}, \quad\hbox{for }d\geq 2.
\end{equation}
For $d\leq 4$, we have a Sobolev embedding $H^1(X) \hookrightarrow L^4(X)$ \cite[Theorem 4.12]{AdamsFournier} and so the Kato Inequality \eqref{eq:FU_6-20_first-order_Kato_inequality} gives
\begin{equation}
\label{eq:Struwe_Yang-Mills_heat_equation_estimate_term_VI_H2}
\|\nabla_{A_1}a_1\times a_2\|_{L^2(0,T; L^2(X))}
\leq
c\|a_2\|_{L^\infty(0,T; H_{A_1}^1(X))} \|a_1\|_{L^2(0,T; H_{A_1}^2(X))},
\quad\hbox{for }2\leq d\leq 4,
\end{equation}
and thus,
$$
\|\nabla_{A_1}a_1\times a_2\|_{L^2(0,T; L^2(X))} \leq c\|a_1\|_\fV \|a_2\|_\fV,
\quad\hbox{for }2\leq d\leq 4,
$$
where $c$ is a universal positive constant, depending at most on the Riemannian metric.

For quadratic term like Struwe's term VII, we have
\begin{align*}
\|a_1\times a_2\|_{L^2(0,T; L^4(X))}^2
&=
\int_0^T \|a_1(t)\times a_2(t)\|_{L^4(X)}^2 \,dt
\\
&\leq c\int_0^T \|a_1(t)\times a_2(t)\|_{H_{A_1}^1(X)}^2 \,dt
\\
&\leq c\int_0^T \left(\|\nabla_{A_1} (a_1(t)\times a_2(t))\|_{L^2(X)}^2 + \|a_1(t)\times a_2(t)\|_{L^2(X)}^2\right) \,dt
\\
&\leq c\int_0^T \left(\|\nabla_{A_1} a_1(t)\times a_2(t))\|_{L^2(X)}^2 + \|a_1(t)\times \nabla_{A_1}a_2(t))\|_{L^2(X)}^2
\right.
\\
&\quad + \left. \|a_1(t)\times a_2(t)\|_{L^2(X)}^2\right) \,dt,
\end{align*}
and thus,
\begin{align*}
{}&\|a_1\times a_2\|_{L^2(0,T; L^4(X))}^2
\\
&\quad \leq c\int_0^T \left(\|\nabla_{A_1} a_1(t)\|_{L^4(X)}^2 \|a_2(t)\|_{L^4(X)}^2
+ \|a_1(t)\|_{L^4(X)}^2 \|\nabla_{A_1} a_2(t)\|_{L^4(X)}^2  \right.
\\
&\qquad + \left. \|a_1(t)\|_{L^4(X)}^2 \|a_2(t)\|_{L^4(X)}^2\right) \,dt
\\
&\quad \leq c\left(\sup_{t \in (0,T)} \|a_1(t)\|_{L^4(X)}^2 \right)
\int_0^T \left(\|\nabla_{A_1} a_2(t)\|_{L^4(X)}^2 + \|a_2(t)\|_{L^4(X)}^2 \right)\,dt
\\
&\qquad + c\left(\sup_{t \in (0,T)} \|a_2(t)\|_{L^4(X)}^2 \right)
\int_0^T \|\nabla_{A_1} a_1(t)\|_{L^4(X)}^2\,dt
\\
&\leq c\left(\sup_{t \in (0,T)} \|a_1(t)\|_{L^4(X)}^2 \right)
\int_0^T \|a_2(t)\|_{W_{A_1}^{1,4}(X)}^2\,dt
\\
&\qquad + c\left(\sup_{t \in (0,T)} \|a_2(t)\|_{L^4(X)}^2 \right)
\int_0^T \|a_1(t)\|_{W_{A_1}^{1,4}(X)}^2\,dt
\end{align*}
that is, taking square roots, for any $d\geq 2$,
\begin{multline}
\label{eq:Struwe_Yang-Mills_heat_equation_estimate_term_VII_W14}
\|a_1\times a_2\|_{L^2(0,T; L^4(X))}
\\
\leq
c\left(\|a_1\|_{L^\infty(0,T; L^4(X))} \|a_2\|_{L^2(0,T; W_{A_1}^{1,4}(X))}
+ \|a_2\|_{L^\infty(0,T; L^4(X))} \|a_1\|_{L^2(0,T; W_{A_1}^{1,4}(X))}\right),
\end{multline}
and hence, using $H^1(X) \hookrightarrow L^4(X)$ and the Kato Inequality, for $2\leq d\leq 4$,
\begin{multline}
\label{eq:Struwe_Yang-Mills_heat_equation_estimate_term_VII_H2}
\|a_1\times a_2\|_{L^2(0,T; L^4(X))}
\\
\leq
c\left(\|a_1\|_{L^\infty(0,T; H_{A_1}^1(X))} \|a_2\|_{L^2(0,T; H_{A_1}^2(X))}
+ \|a_2\|_{L^\infty(0,T; H_{A_1}^1(X))} \|a_1\|_{L^2(0,T; H_{A_1}^2(X))}\right),
\end{multline}
and thus,
$$
\|a_1\times a_2\|_{L^2(0,T; L^4(X))} \leq c\|a_1\|_\fV \|a_2\|_\fV,
\quad\hbox{for }2\leq d\leq 4,
$$
where $c$ is a universal positive constant, depending at most on the Riemannian metric.

Next we consider the estimate for a cubic term like Struwe's term VIII appearing in \cite[p. 137]{Struwe_1994}. We have
\begin{align*}
\|a_1\times a_2\times a_3\|_{L^2(0,T; L^2(X))}^2
&= \int_0^T \|a_1(t)\times a_2(t)\times a_3(t)\|_{L^2(X)}^2 \,dt
\\
&\leq c\int_0^T \|a_1(t)\times a_2(t)\|_{L^4(X)}^2 \|a_3(t)\|_{L^4(X)}^2 \,dt
\\
&\leq c\left(\sup_{t \in (0,T)} \|a_3(t)\|_{L^4(X)}^2 \right)
\int_0^T \|a_1(t)\times a_2(t)\|_{L^4(X)}^2 \,dt,
\end{align*}
that is, taking square roots,
$$
\|a_1\times a_2\times a_3\|_{L^2(0,T; L^2(X))}
\leq
c\|a_3\|_{L^\infty(0,T; L^4(X))} \|a_1\times a_2\|_{L^2(0,T; L^4(X))}.
$$
Therefore, by combining the preceding inequality with the bound \eqref{eq:Struwe_Yang-Mills_heat_equation_estimate_term_VII_W14}, we have
\begin{multline}
\label{eq:Struwe_Yang-Mills_heat_equation_estimate_term_VIII_W14}
\|a_1\times a_2\times a_3\|_{L^2(0,T; L^2(X))}
\\
\leq
c\|a_3\|_{L^\infty(0,T; L^4(X))}
\left(\|a_1\|_{L^\infty(0,T; L^4(X))} \|a_2\|_{L^2(0,T; W_{A_1}^{1,4}(X))} \right.
\\
+ \left. \|a_2\|_{L^\infty(0,T; L^4(X))} \|a_1\|_{L^2(0,T; W_{A_1}^{1,4}(X))}\right),
\quad\hbox{for }d\geq 2,
\end{multline}
and hence, using $H^1(X) \hookrightarrow L^4(X)$ and the Kato Inequality,
\begin{multline}
\label{eq:Struwe_Yang-Mills_heat_equation_estimate_term_VIII_H2}
\|a_1\times a_2\times a_3\|_{L^2(0,T; L^2(X))}\
\\
\leq
c\|a_3\|_{L^\infty(0,T; H_{A_1}^1(X))}
\left(\|a_1\|_{L^\infty(0,T; H_{A_1}^1(X))} \|a_2\|_{L^2(0,T; H_{A_1}^2(X))} \right.
\\
+ \left. \|a_2\|_{L^\infty(0,T; H_{A_1}^1(X))} \|a_1\|_{L^2(0,T; H_{A_1}^2(X))}\right),
\quad\hbox{for }2\leq d\leq 4,
\end{multline}
and thus,
$$
\|a_1\times a_2\times a_3\|_{L^2(0,T; L^2(X))}
\leq
c\|a_1\|_\fV \|a_2\|_\fV \|a_3\|_\fV, \quad\hbox{for }2\leq d\leq 4,
$$
where $c$ is a universal positive constant, depending at most on the Riemannian metric.

\subsection{Local existence of solutions to the Yang-Mills heat equation over a closed manifold with dimension less than or equal to four and small initial data in $H^1$}
\label{subsec:Struwe_page_137_contraction_mapping_small_initial_data_in_H1}
In this subsection, we describe a simple and elegant approach due to Struwe \cite[Section 4.3]{Struwe_1994} for establishing local existence of solutions for the Yang-Mills heat equation with initial data in $H^1_{A_1}(X; \Lambda^1\otimes\ad P)$ when the dimension of $X$ is less than or equal to four. While we follow Struwe's idea, we simplify his development using our observation that it is not necessary to employ a time-varying family of reference connections for regularity reasons; instead, a \emph{fixed} $C^\infty$ reference connection, $A_1$ on $P$, will do when the initial data, $a_0 = A_0-A_1 \in H^1_{A_1}(X; \Lambda^1\otimes\ad P)$, has small norm, $\|a_0\|_{H^1_{A_1}(X)}$. In Section \ref{subsec:Struwe_page_137_contraction_mapping_arbitrary_initial_data_in_H1}, we remove the small $\|a_0\|_{H^1_{A_1}(X)}$ constraint using Struwe's idea in \cite[Section 4.2]{Struwe_1994}, but emphasize that we do not appeal to it for the regularity reasons described in \cite[Section 4.1]{Struwe_1994}. Struwe's goal, as is ours, was to treat the case when $X$ has dimension $d=4$ but here we make the simple observation that of course his method applies when $2\leq d\leq 4$.

The goal is to solve the Yang-Mills heat equation \eqref{eq:Yang-Mills_heat_equation_as_perturbation_rough_Laplacian_plus_one_heat_equation} by a contraction mapping in the Banach space $\fV$ given by \eqref{eq:Struwe_page_129_Banach_space}, with $\tau$ in place of $T$ and for small enough $\tau > 0$.

Given any $\tau>0$ and
$$
a_0 \in H^1_{A_1}(X; \Lambda^1\otimes\ad P) \quad\hbox{and}\quad f \in L^2(0, \tau; L^2(X; \Lambda^1\otimes\ad P)),
$$
then Theorem \ref{thm:Sell_You_42-12_heat_equation_alpha_is_one} and Corollary \ref{cor:Sell_You_42-13_heat_equation} imply that there is a unique strong solution $a \in \fV$ to \eqref{eq:Linear_heat_equation_with_rough_Laplacian_plus_one_on_Omega_1_adP} with initial data $a(0) = a_0$, that is,
$$
\frac{\partial a}{\partial t} + (\nabla_{A_1}^*\nabla_{A_1} + 1)a = f \quad\hbox{a.e. on } (0, \tau), \quad a(0) = a_0.
$$
As suggested by Struwe \cite[p. 137]{Struwe_1994}, we shall use a contraction mapping approach to solving \eqref{eq:Yang-Mills_heat_equation_as_perturbation_rough_Laplacian_plus_one_heat_equation} in $\fV$. With that in mind, we claim that, for any $w \in \fV$,
\begin{equation}
\label{eq:Assumption_sF(w)_in_L2from_0_to_T_into_L^2(X)_given_w_infH}
f := \sF(w) \in L^2(0, \tau; L^2(X; \Lambda^1\otimes\ad P)).
\end{equation}
We now verify the claim \eqref{eq:Assumption_sF(w)_in_L2from_0_to_T_into_L^2(X)_given_w_infH}. Proceeding as in \cite[p. 137]{Struwe_1994} and appealing to our expression \eqref{eq:Yang-Mills_heat_equation_nonlinearity_relative_rough_Laplacian_plus_one} for $\sF(w)$ and applying the inequalities \eqref{eq:Struwe_Yang-Mills_heat_equation_estimate_term_VI_H2} and \eqref{eq:Struwe_Yang-Mills_heat_equation_estimate_term_VIII_H2} to bound the quadratic and cubic terms, we have
\begin{align*}
{}&\|\sF(w)\|_{L^2(0,\tau; L^2(X))}
\\
&\quad \leq \|d_{A_1}^*F_{A_1}\|_{L^2(0,\tau; L^2(X))} + \|\Ric_g\times w\|_{L^2(0,\tau; L^2(X))} + \|(F_{A_1} - 1)\times w\|_{L^2(0,\tau; L^2(X))}
\\
&\qquad + \|\nabla_{A_1}w\times w\|_{L^2(0,\tau; L^2(X))} + \|w\times w\times w\|_{L^2(0,\tau; L^2(X))}
\\
&\quad \leq \tau\|d_{A_1}^*F_{A_1}\|_{L^2(X)} + c\left(1 + \|\Ric_g\|_{C(X)} + \|F_{A_1}\|_{C(X)}\right)\|w\|_{L^2(0,\tau; L^2(X))}
\\
&\qquad + c\|w\|_{L^\infty(0,T; H_{A_1}^1(X))} \|w\|_{L^2(0,T; H_{A_1}^2(X))}
+ c\|w\|_{L^\infty(0,T; H_{A_1}^1(X))}^2 \|w\|_{L^2(0,T; H_{A_1}^2(X))}
\\
&\quad \leq \tau\|d_{A_1}^*F_{A_1}\|_{L^2(X)}
+ c\sqrt{\tau}\left(1 + \|\Ric_g\|_{C(X)} + \|F_{A_1}\|_{C(X)}\right)\|w\|_{L^\infty(0,\tau; L^2(X))}
\\
&\qquad + c\|w\|_{L^\infty(0,T; H_{A_1}^1(X))} \|w\|_{L^2(0,T; H_{A_1}^2(X))}
+ c\|w\|_{L^\infty(0,T; H_{A_1}^1(X))}^2 \|w\|_{L^2(0,T; H_{A_1}^2(X))},
\end{align*}
where the positive constant, $c$, depends at most on the Riemannian metric, $g$, on $X$. The preceding bound verifies the claim \eqref{eq:Assumption_sF(w)_in_L2from_0_to_T_into_L^2(X)_given_w_infH} and, by our definition \eqref{eq:Struwe_page_129_Banach_space} of the Banach space, $\fV$, yields the inequality,
\begin{equation}
\label{eq:Struwe_page_137_nonlinearity_L2_time_L2_space_bound}
\|\sF(w)\|_{L^2(0,\tau; L^2(X))} \leq C_0\tau + C_0\sqrt{\tau}\|w\|_\fV + c_0\|w\|_\fV^2 + c_0\|w\|_\fV^3, \quad\forall\, w \in \fV,
\end{equation}
where the positive constant, $c_0$, depends at most on $g$ and $C_0$ depends at most on $A_1$ and $g$. For each $w \in \fV$, we can thus let
\begin{equation}
\label{eq:Definition_Struwe_section_4-3_contraction_map}
a := \Phi(w) \in \fV
\end{equation}
be the unique solution to \eqref{eq:Linear_heat_equation_with_rough_Laplacian_plus_one_on_Omega_1_adP}, with initial data $a(0) = a_0$, provided by Theorem \ref{thm:Sell_You_42-12_heat_equation_alpha_is_one}, therefore defining a map $\Phi: \fV \to \fV$ thanks to \eqref{eq:Assumption_sF(w)_in_L2from_0_to_T_into_L^2(X)_given_w_infH}.

We could now appeal to Corollary \ref{cor:Sell_You_42-12_heat_equation_alpha_is_one_apriori_estimates} to bound $\|a\|_\fV$ in terms of $\|(\partial_t + \nabla_{A_1}^*\nabla_{A_1})a\|_{L^2(0,\tau; L^2(X))}$ but, as we are using the traditional gauge-theoretic definition of Sobolev spaces via covariant derivatives,
the constants in the \apriori estimates in Corollary \ref{cor:Sell_You_42-12_heat_equation_alpha_is_one_apriori_estimates} (whose Sobolev norms are defined via spectral theory) would acquire a non-explicit dependence on $A_1$, as we noted in Remark \ref{rmk:Sell_You_42-12_heat_equation_alpha_is_zero_or_one_or_real_apriori_estimate_constant_dependencies}. We shall instead digress to record an alternative due to Struwe \cite[Lemma 3.2]{Struwe_1994}, where it is possible to more easily trace the dependence on $A_1$ in the \apriori estimate constant.

\begin{lem}[\Apriori estimate for the connection Laplace operator]
\label{lem:Struwe_3-1}
\cite[Lemma 3.1]{Struwe_1994}
Let $G$ be a compact Lie group and $P$ be a principal $G$-bundle over a closed, smooth manifold, $X$, of dimension $d\geq 2$, with Riemannian metric, $g$, and $A_1$ a fixed reference connection of class $C^\infty$ on $P$. Then there is a positive constant, $C_1 = C_1(A_1,g)$, with the following significance. If $l \in \NN$ and $a \in H_{A_1}^2(X; \Lambda^l\otimes\ad P)$, then
\begin{equation}
\label{eq:Struwe_lemma_3-1}
\|a\|_{H_{A_1}^2(X)} \leq C_1\left(\|\nabla_{A_1}^*\nabla_{A_1}a\|_{L^2(X)} + \|a\|_{L^2(X)}\right).
\end{equation}
\end{lem}

\begin{rmk}[On the dependencies of the constant in Lemma \ref{lem:Struwe_3-1}]
\label{rmk:Lemma_Struwe_3-1_constant_dependencies}
We recall from our Lemma \ref{lem:L22InfinityEstu} that, when $d=4$, it is possible to write $C_1 = c(1+\|F_{A_1}\|_{L^2(X)})$, where $c$ depends at most on $g$, \emph{provided} the term $\|\nabla_{A_1}^*\nabla_{A_1}a\|_{L^2(X)}$ on the right-hand side of the inequality \eqref{eq:Struwe_lemma_3-1} is replaced by $\|\nabla_{A_1}^*\nabla_{A_1}a\|_{L^{2,\sharp}(X)}$. Improvements of this kind are easier when $d=2$ or $3$.
\end{rmk}

\begin{lem}[\Apriori estimate for the connection heat operator]
\label{lem:Struwe_3-2}
\cite[Lemma 3.2]{Struwe_1994}
Let $G$ be a compact Lie group and $P$ be a principal $G$-bundle over a closed, smooth manifold, $X$, of dimension $d\geq 2$, with Riemannian metric, $g$, and $A_1$ a fixed reference connection of class $C^\infty$ on $P$. Then there are positive constants, $C_2 = C_2(A_1,g) \geq 1$ and $\tau = \tau(A_1,g) \in (0, 1]$, with the following significance. If
$$
a \in L^2(0,\tau; H_{A_1}^2(X; \Lambda^1\otimes\ad P)) \cap H^1(0,\tau; L^2(X; \Lambda^1\otimes\ad P)),
$$
then
\begin{multline}
\label{eq:Struwe_lemma_3-2_small_time_interval}
\|a\|_{L^2(0,\tau; H_{A_1}^2(X))} + \|a\|_{H^1(0,\tau; L^2(X))}
\\
\leq C_2\left( \left\|\left(\partial_t + \nabla_{A_1}^*\nabla_{A_1}\right)a\right\|_{L^2(0,\tau; L^2(X))}
+ \|a(0)\|_{H_{A_1}^1(X))} \right).
\end{multline}
Moreover, for any $T > 0$, there is a positive constant, $C_2' = C_2'(A_1,g, T) \geq 1$ such that
\begin{multline}
\label{eq:Struwe_lemma_3-2_finite_time_interval}
\|a\|_{L^2(0, T; H_{A_1}^2(X))} + \|a\|_{H^1(0, T; L^2(X))}
\\
\leq C_2'\left( \left\|\left(\partial_t + \nabla_{A_1}^*\nabla_{A_1}\right)a\right\|_{L^2(0, T; L^2(X))}
+ \|a(0)\|_{H_{A_1}^1(X))} \right).
\end{multline}
\end{lem}

\begin{proof}
The inequality \eqref{eq:Struwe_lemma_3-2_small_time_interval} is provided by \cite[Lemma 3.2]{Struwe_1994}. To prove \eqref{eq:Struwe_lemma_3-2_finite_time_interval}, it suffices to partition $[0, T$ into small subintervals and then use the methods of proof of the \apriori estimate
\eqref{eq:Rade_11-3_and_11-4_finite_time_interval} in Lemma \ref{lem:Rade_inequalities_11-3_and_11-4} or the \apriori estimate
\eqref{eq:Rade_apriori_interior_estimate_lemma_7.3} in Lemma \ref{lem:Rade_7-3} to extend \eqref{eq:Struwe_lemma_3-2_small_time_interval} from a small time interval, $[0,\tau]$, to an arbitrary finite time interval, $[0,T]$.
\end{proof}

\begin{rmk}[On relating \apriori estimates involving the connection and Hodge Laplace operators]
\label{rmk:Relating_apriori_estimates_for_connection_and_Hodge_Laplace_operators}
Struwe's \cite[Lemmata 3.1 and 3.2]{Struwe_1994} use the Hodge Laplace operator \eqref{eq:Lawson_page_93_Hodge_Laplacian} rather than the simpler connection Laplace operator which suffices for our development. Of course, the two styles of \apriori estimates may be related by Bochner-Weitzenb\"ock formulae such as \eqref{eq:Lawson_corollary_II-2}.
\end{rmk}

\begin{rmk}[On the dependencies of the constant in Lemma \ref{lem:Struwe_3-2}]
\label{rmk:Lemma_Struwe_3-2_constant_dependencies}
When $d=4$, we recall from our Lemma \ref{lem:Feehan_5-3_heat_operator} that, provided one is willing to employ a stronger system of Sobolev norms than those appearing in the \apriori estimate \eqref{eq:Struwe_lemma_3-2_finite_time_interval}, it is possible to replace $C_1 = C_1(A_1,g)$ by a universal numerical constant, independent of $A_1$ or $g$.
\end{rmk}

By combining \eqref{eq:Struwe_9} and \eqref{eq:Struwe_lemma_3-2_small_time_interval}, we see that (after absorbing a factor of $\sqrt{2}$ in $C_2$),
\begin{multline}
\label{eq:Struwe_lemma_3-2_plus_Linfty_H1_small_time_interval}
\|a\|_{L^2(0,\tau; H_{A_1}^2(X))} + \|a\|_{L^\infty(0,\tau; H_{A_1}^1(X))} + \|a\|_{H^1(0,\tau; L^2(X))}
\\
\leq C_2\left( \left\|\left(\partial_t + \nabla_{A_1}^*\nabla_{A_1}\right)a\right\|_{L^2(0,\tau; L^2(X))}
+ \|a(0)\|_{H_{A_1}^1(X))} \right),
\end{multline}
and of course we can similarly extend \eqref{eq:Struwe_lemma_3-2_finite_time_interval}.

To check the boundedness property of the map $\Phi$ in \eqref{eq:Definition_Struwe_section_4-3_contraction_map}, we apply the \apriori estimate \eqref{eq:Struwe_lemma_3-2_plus_Linfty_H1_small_time_interval}. Letting $a := \Phi(w)$ be the unique solution to \eqref{eq:Linear_heat_equation_with_rough_Laplacian_plus_one_on_Omega_1_adP} with $f := \sF(w)$ for $w \in \fV$, and initial data $a(0) = a_0 \in H_{A_1}^1(X; \Lambda^1\otimes \ad P)$, we discover that
\begin{align*}
\|\Phi(w)\|_\fV &= \|a\|_\fV
\\
&\leq C_2\left(\left\|\left(\partial_t + \nabla_{A_1}^*\nabla_{A_1}\right)a\right\|_{L^2(0,\tau; L^2(X))}
+ \|a_0\|_{H_{A_1}^1(X)} \right)
\quad\hbox{(by \eqref{eq:Struwe_page_129_Banach_space} and \eqref{eq:Struwe_lemma_3-2_plus_Linfty_H1_small_time_interval})}
\\
&= C_2\left(\left\|\sF(w)\right\|_{L^2(0,\tau; L^2(X))} + \|a_0\|_{H_{A_1}^1(X)} \right)
\quad\hbox{(by \eqref{eq:Linear_heat_equation_with_rough_Laplacian_plus_one_on_Omega_1_adP})}
\\
&\leq  C_2\left(C_0\tau + C_0\sqrt{\tau}\|w\|_\fV + c_0\|w\|_\fV^2 + c_0\|w\|_\fV^3 + \|a_0\|_{H_{A_1}^1(X)} \right)
\quad\hbox{(by \eqref{eq:Struwe_page_137_nonlinearity_L2_time_L2_space_bound})}.
\end{align*}
Now suppose that $\|a_0\|_{H_{A_1}^1(X)} \leq \eps$, where $\eps \in (0, 1]$ is a constant to be determined, so the preceding inequality gives
$$
\|\Phi(w)\|_\fV \leq C_2\eps + C_0C_2\tau + C_2C_0\sqrt{\tau}\|w\|_\fV + c_0C_2\left(1 + \|w\|_\fV\right)\|w\|_\fV^2.
$$
Setting $R = 2C_2\eps$ and assuming that $\|w\|_\fV \leq R$, we have
$$
c_0C_2 \|w\|_\fV\left(1 + \|w\|_\fV\right)
\leq
2\eps c_0C_2^2(1+2C_2\eps)
\leq
2\eps c_0C_2^2(1+2C_2).
$$
Now choose $\eps = \eps(c_0,C_2) \in (0, 1]$ small enough that $2\eps c_0C_2^2(1+2C_2) \leq 1/6$, that is
\begin{equation}
\label{eq:Struwe_4-3_bounded_mapping_epsilon_constraint}
0 < \eps \leq \frac{1}{12c_0C_2^2(1+2C_2)}.
\end{equation}
Next choose $\tau = \tau(\eps,C_0,C_2) = \tau(c_0,C_0,C_2) \in (0, 1]$ small enough that $C_0C_2\tau \leq 2C_2\eps/6$, that is, $C_0\tau \leq \eps/3$, and $C_0C_2\sqrt{\tau} \leq 1/6$, that is, $C_0^2C_2^2\tau \leq 1/36$, and so
\begin{equation}
\label{eq:Struwe_4-3_bounded_mapping_tau_constraint}
0 < \tau \leq \min\left\{\frac{\eps}{3C_0}, \frac{1}{36C_0^2C_2^2}\right\}.
\end{equation}
Denote $\fB_R := \{v \in \fV: \|v\|_\fV \leq R\}$. By assembling the preceding inequalities, we see that
$$
\|\Phi(w)\|_\fV \leq C_2\eps + \frac{2C_2\eps}{6} + \frac{2C_2\eps}{6} + \frac{2C_2\eps}{6}
= 2C_2\eps = R, \quad \forall\, w \in \fB_R,
$$
and consequently, $\Phi$ is a well-defined map of the closed ball, $\fB_R$, to itself.

To check the contraction-mapping property of the map $\Phi$ in \eqref{eq:Definition_Struwe_section_4-3_contraction_map}, we apply the \apriori estimate \eqref{eq:Struwe_lemma_3-2_plus_Linfty_H1_small_time_interval}, keeping in mind our definition \eqref{eq:Struwe_page_129_Banach_space} of $\fV$. Letting $a_i := \Phi(w_i)$ be the unique solution to \eqref{eq:Linear_heat_equation_with_rough_Laplacian_plus_one_on_Omega_1_adP} with $f_i := \sF(w_i)$ for $w_i \in \fV$, with $i=1,2$, and initial data $a_1(0) = a_2(0) = a_0 \in H_{A_1}^1(X; \Lambda^1\otimes \ad P)$, we find that
\begin{align*}
\|\Phi(w_1) - \Phi(w_2)\|_\fV &= \|a_1 - a_2\|_\fV
\\
&\leq C_2\left\|\left(\partial_t + \nabla_{A_1}^*\nabla_{A_1}\right)(a_1-a_2)\right\|_{L^2(0,\tau; L^2(X))}
\quad\hbox{(by \eqref{eq:Struwe_page_129_Banach_space} and \eqref{eq:Struwe_lemma_3-2_plus_Linfty_H1_small_time_interval})}
\\
&= C_2\|\sF(w_1) - \sF(w_2)\|_{L^2(0,\tau; L^2(X))}
\quad\hbox{(by \eqref{eq:Linear_heat_equation_with_rough_Laplacian_plus_one_on_Omega_1_adP})}.
\end{align*}
To bound the term $\sF(w_1) - \sF(w_2)$, we observe that the schematic expression \eqref{eq:Yang-Mills_heat_equation_nonlinearity_relative_rough_Laplacian_plus_one} for the nonlinearity, $\sF$, yields
\begin{equation}
\label{eq:Yang-Mills_nonlinearity_difference}
\begin{aligned}
{}& \sF(w_1) - \sF(w_2)
\\
&\quad = (F_{A_1} - 1)\times (w_1 - w_2) + \Ric_g\times (w_1 - w_2)
\\
&\qquad + \nabla_{A_1}(w_1 - w_2)\times w_1 + \nabla_{A_1} w_2\times (w_1 - w_2)
\\
&\qquad + (w_1 - w_2)\times w_1\times w_1 + w_2\times (w_1 - w_2) \times w_1 + w_2\times w_2 \times (w_1 - w_2).
\end{aligned}
\end{equation}
Therefore, applying \eqref{eq:Struwe_Yang-Mills_heat_equation_estimate_term_VI_H2} and \eqref{eq:Struwe_Yang-Mills_heat_equation_estimate_term_VIII_H2} to bound the quadratic and cubic terms,
\begin{align*}
{}& \|\sF(w_1) - \sF(w_2)\|_{L^2(0,\tau; L^2(X))}
\\
&\quad \leq
c\left(1 + \|F_{A_1}\|_{C(X)} + \|\Ric_g\|_{C(X)} \right)\|w_1 - w_2\|_{L^2(0,\tau; L^2(X))}
\\
&\qquad + c\|w_1 - w_2\|_{L^2(0,T; H_{A_1}^2(X))} \|w_1\|_{L^\infty(0,T; H_{A_1}^1(X))}
\\
&\qquad + c\|w_2\|_{L^2(0,T; H_{A_1}^2(X))} \|w_1 - w_2\|_{L^\infty(0,T; H_{A_1}^1(X))}
\\
&\qquad + c\left(\|w_1\|_{L^\infty(0,T; H_{A_1}^1(X))}^2 + \|w_2\|_{L^\infty(0,T; H_{A_1}^1(X))}^2 \right)
\|w_1 - w_2\|_{L^2(0,T; H_{A_1}^2(X))},
\end{align*}
where $c$ is a positive constant which depends at most on the Riemannian metric, $g$. But
$$
\|w_1 - w_2\|_{L^2(0,\tau; L^2(X))}
\leq
\sqrt{\tau} \|w_1 - w_2\|_{L^\infty(0,\tau; L^2(X))}
\leq
\sqrt{\tau} \|w_1 - w_2\|_{L^\infty(0,\tau; H_{A_1}^1(X))}.
$$
Thus, by definition \eqref{eq:Struwe_page_129_Banach_space} of $\fV$ and by combining the preceding inequalities, we obtain,
\begin{multline*}
\|\sF(w_1) - \sF(w_2)\|_{L^2(0,\tau; L^2(X))}
\\
\leq
C_3\sqrt{\tau}\|w_1 - w_2\|_\fV
+ c_1\left( \|w_1\|_\fV + \|w_2\|_\fV + \|w_1\|_\fV^2 + \|w_2\|_\fV^2 \right) \|w_1 - w_2\|_\fV,
\end{multline*}
where $C_3$ is a positive constant which depends at most on $A_1$ and $g$, while $c_1$ depends at most on $g$. Therefore,
\begin{multline*}
\|\Phi(w_1) - \Phi(w_2)\|_{L^2(0,\tau; L^2(X))}
\\
\leq
C_2C_3\sqrt{\tau}\|w_1 - w_2\|_\fV
+ c_1C_2\left( \|w_1\|_\fV + \|w_2\|_\fV + \|w_1\|_\fV^2 + \|w_2\|_\fV^2 \right) \|w_1 - w_2\|_\fV.
\end{multline*}
We now further constrain $\eps$ and $\tau$ to ensure that $\Phi$ is a contraction mapping. For $w_1, w_2 \in \fB_R$, we have
$$
c_1C_2\left( \|w_1\|_\fV + \|w_2\|_\fV + \|w_1\|_\fV^2 + \|w_2\|_\fV^2 \right) \leq 2c_1C_2R(1+R).
$$
Recall that $R = 2C_2\eps$ and observe that $2c_1C_2R(1+R) = 4\eps c_1C_2^2(1+2\eps C_2) \leq 4\eps c_1C_2^2(1+2C_2)$. Choose $\eps = \eps(c_1,C_1,C_2) \in (0, 1]$ small enough that  $4\eps c_1C_2^2(1+2C_2) \leq 1/4$, that is,
\begin{equation}
\label{eq:Struwe_4-3_contraction_mapping_epsilon_constraint}
0 < \eps \leq \frac{1}{16c_1C_2^2(1+2C_2)},
\end{equation}
and $\eps = \eps(c_0,C_2)$ continues to obey \eqref{eq:Struwe_4-3_bounded_mapping_epsilon_constraint}.
Choose $\tau = \tau(C_2,C_3) \in (0, 1]$ small enough that $C_2C_3\sqrt{\tau} \leq 1/4$, namely,
\begin{equation}
\label{eq:Struwe_4-3_contraction_mapping_tau_constraint}
0 < \tau \leq \frac{1}{16C_2^2C_3^2},
\end{equation}
and $\tau = \tau(c_0,C_0,C_2) \in (0, 1]$ continues to obey \eqref{eq:Struwe_4-3_bounded_mapping_tau_constraint}. Consequently, for $w_1, w_2 \in \fB_R$, we see that,
$$
\|\Phi(w_1) - \Phi(w_2)\|_{L^2(0,\tau; L^2(X))} \leq \frac{1}{2}\|w_1 - w_2\|_\fV, \quad\forall\, w_1, w_2 \in \fB_R.
$$
The Banach Contraction Mapping Theorem now implies that the map $\Phi:\fB_R \to \fB_R$ has a unique fixed point, $a \in \fV$ obeying $\|a\|_\fV \leq R  = 2C_2\eps$, and thus a unique solution to the Yang-Mills heat equation
\eqref{eq:Yang-Mills_heat_equation_as_perturbation_rough_Laplacian_plus_one_heat_equation} which obeys the \apriori estimate \eqref{eq:Struwe_section_4-3_apriori_estimate_Yang-Mills_heat_equation} below with $C = 2C_2$. If $a_1, a_2$ are the unique solutions to \eqref{eq:Yang-Mills_heat_equation_as_perturbation_rough_Laplacian_plus_one_heat_equation} corresponding to initial data $a_1(0)$ and $a_2(0)$, then the continuity estimate \eqref{eq:Struwe_section_4-3_continuity_estimate_Yang-Mills_heat_equation} below with $C = 2C_2$ follows from the inequality \eqref{eq:Struwe_lemma_3-2_plus_Linfty_H1_small_time_interval}, the definition \eqref{eq:Struwe_page_129_Banach_space} of $\fV$, the Yang-Mills heat equation \eqref{eq:Yang-Mills_heat_equation_as_perturbation_rough_Laplacian_plus_one_heat_equation}, the contraction mapping property for $\Phi$, and the fact that $\Phi(a)=a$. In summary, we have proved

\begin{thm}[Local well-posedness for the Yang-Mills heat equation over a closed manifold with small initial data in $H^1$]
\label{thm:Struwe_section_4-3_local_wellposedness_Yang-Mills_heat_equation_small_initial_data_in_H1}
\cite[Section 4.3]{Struwe_1994}
Let $G$ be a compact Lie group and $P$ a principal $G$-bundle over a closed, connected smooth manifold, $X$, of dimension $2 \leq d \leq 4$ and Riemannian metric, $g$. If $A_1$ is a reference connection of class $C^\infty$ on $P$, then there are positive constants $C = C(A_1,g) \geq 1$ and $\eps = \eps(A_1,g) \in (0, 1]$ and $\tau = \tau(A_1,g) \in (0, 1]$, with the following significance. If
$$
a_0 := A_0 - A_1 \in H^1_{A_1}(X; \Lambda^1\otimes\ad P),
$$
obeys
\begin{equation}
\label{eq:Norm_a0_H1_leq_epsilon}
\|a_0\|_{H^1_{A_1}(X)} \leq \eps,
\end{equation}
then there is a unique strong solution $a(t)$ in $\fV$ to the Yang-Mills heat equation
\eqref{eq:Yang-Mills_heat_equation_as_perturbation_rough_Laplacian_plus_one_heat_equation} on the interval $(0, \tau)$ with initial data $a(0) = a_0$, where $\fV$ is the Banach space \eqref{eq:Struwe_page_129_Banach_space}. Moreover, the solution $a(t)$ obeys the \apriori estimate,
\begin{equation}
\label{eq:Struwe_section_4-3_apriori_estimate_Yang-Mills_heat_equation_small_initial_data_in_H1}
\|a\|_\fV \leq C\eps.
\end{equation}
If $a_1, a_2$ are the unique solutions corresponding to initial data $a_1(0)$ and $a_2(0)$, respectively, then
\begin{equation}
\label{eq:Struwe_section_4-3_continuity_estimate_Yang-Mills_heat_equation}
\|a_1 - a_2\|_\fV \leq C\|a_1(0) - a_2(0)\|_{H^1_{A_1}(X)}.
\end{equation}
\end{thm}

Local well-posedness results similar to that of Theorem \ref{thm:Struwe_section_4-3_local_wellposedness_Yang-Mills_heat_equation_small_initial_data_in_H1} can be established by other methods, but apparently only at the cost of assuming greater regularity for the initial data, $a_0$, or by relaxing our requirement of local well-posedness for the Yang-Mills heat equation to local existence and uniqueness, that is, by dropping the requirement that the solution depend continuously on the initial data. For example:
\begin{enumerate}
\item $a_0 \in H^\alpha_{A_1}(X; \Lambda^1\otimes\ad P)$ with $\alpha > 1$ in our application of Theorem \ref{thm:Sell_You_42-12}, so $\calV^{\alpha+1} = H^{\alpha+1}_{A_1}(X; \Lambda^1\otimes\ad P)$ and $H^{\alpha+1}(X)$ is a Banach algebra when $d = 4$; or

\item $a_0 \in H^1_{A_1}(X; \Lambda^1\otimes\ad P) \cap C(X; \Lambda^1\otimes\ad P)$ in our application of our results in Section \ref{sec:Critical-exponent_parabolic_Sobolev_spaces_linear_parabolic_operator_vector_bundle_manifold} to hybrid or pure critical-exponent parabolic Sobolev spaces; or

\item $a_0 \in L^d(X; \Lambda^1\otimes\ad P)$ and we seek a solution in a Banach space defined via temporal weighting, essentially following the idea of Kozono, Maeda, and Naito \cite{Kozono_Maeda_Naito_1995}.
\end{enumerate}
In the remainder of this section we shall develop these alternative approaches.

\begin{rmk}[On the constraint in Theorem \ref{thm:Struwe_section_4-3_local_wellposedness_Yang-Mills_heat_equation_small_initial_data_in_H1} on the norm of the initial data]
\label{rmk:Theorem_Struwe_section_4-3_local_wellposedness_Yang-Mills_heat_equation_initial_data_norm}
Because the Fr\'echet affine space of connections of class $C^\infty$ is dense in the Banach affine space of connections of class $H^1$, there is no significant loss of generality in the constraint that $A_0$ obey $\|A_0-A_1\|_{H_{A_1}^1(X)} \leq \eps$, with $\eps \in (0, 1]$ sufficiently small, where $a_0 = A_0 - A_1 \in H_{A_1}^1(X; \Lambda^1\otimes\ad P)$ serves as initial data for the Yang-Mills heat equation \eqref{eq:Yang-Mills_heat_equation_as_perturbation_rough_Laplacian_plus_one_heat_equation}. However, as we shall see in Sections \ref{subsec:Struwe_page_137_contraction_mapping_hybrid_critical-exponent_parabolic_Sobolev_space}, \ref{subsec:Struwe_page_137_contraction_mapping_pure_critical-exponent_parabolic_Sobolev_space}, and \ref{subsec:Struwe_page_137_contraction_mapping_initial_data_in_Ld}, the alternative approaches we listed in the preceding comments do not carry this restriction on the norm of the initial data. Moreover, the restriction on the norm of the initial data required for the proof of Theorem \ref{thm:Struwe_section_4-3_local_wellposedness_Yang-Mills_heat_equation_small_initial_data_in_H1} can be removed by using Struwe's idea \cite[Section 4.2]{Struwe_1994} of constructing the solution to the Yang-Mills heat equation \eqref{eq:Yang-Mills_heat_equation_as_perturbation_rough_Laplacian_plus_one_heat_equation} as a small perturbation of a solution to suitably chosen linear heat equation. We shall develop this approach in Section \ref{subsec:Struwe_page_137_contraction_mapping_arbitrary_initial_data_in_H1}.
\end{rmk}

\subsection{Local existence of solutions to the Yang-Mills heat equation over a closed manifold with dimension less than or equal to four and arbitrary initial data in $H^1$}
\label{subsec:Struwe_page_137_contraction_mapping_arbitrary_initial_data_in_H1}
In order to remove the small-norm constraint on the initial data, $a_0 \in H_{A_1}^1(X; \Lambda^1\otimes\ad P)$, in Theorem \ref{thm:Struwe_section_4-3_local_wellposedness_Yang-Mills_heat_equation_small_initial_data_in_H1} we may instead consider the solution, $a(t)$, to the Yang-Mills heat equation \eqref{eq:Yang-Mills_heat_equation_as_perturbation_rough_Laplacian_plus_one_heat_equation} to be a small perturbation, $a(t) = \alpha(t)+\xi(t)$ of the solution, $\alpha(t)$, to a linear heat equation with \emph{smooth} initial data,
\begin{equation}
\label{eq:Struwe_19_with_augmented_connection_Laplacian}
\frac{\partial \alpha}{\partial t} + (\nabla_{A_1}^*\nabla_{A_1}+1)\alpha = 0 \quad\hbox{on }(0,\infty)\times X,
\quad \alpha(0) = \alpha_0 \in \Omega^1(X;\ad P).
\end{equation}
Indeed, Theorem \ref{thm:Sell_You_42-12_heat_equation_alpha_is_one} and Corollary \ref{cor:Sell_You_42-14_heat_equation_smoothness_for_t_geq_zero} provide a unique classical solution, $\alpha \in C^\infty([0,\infty)\times X;\Lambda^1\otimes\ad P)$, to the linear Cauchy problem \eqref{eq:Struwe_19_with_augmented_connection_Laplacian}.

We shall now write
\begin{equation}
\label{eq:Decomposition_initial_data_a0_as_large_smooth_alpha0_plus_small_H1_xi0}
A_0 = A_1+a_0 = A_1 + \alpha_0 + \xi_0,
\end{equation}
where $\alpha_0 \in \Omega^1(X;\ad P) = C^\infty(X;\Lambda^1\otimes \ad P)$ may have norm, $\|\alpha_0\|_{H_{A_1}^1(X)}$ of arbitrary size, but $\xi_0 = A_0 - A_1 - \alpha_0 \in H_{A_1}^1(X;\Lambda^1\otimes\ad P)$ will have \emph{small} norm,
\begin{equation}
\label{eq:Yang-Mills_heat_equation_H1_small_initial_data_bound}
\|\xi_0\|_{H_{A_1}^1(X)} \leq \eps,
\end{equation}
with measure of smallness $\eps \in (0, 1]$ to be determined. We shall seek a solution to the Yang-Mills heat equation \eqref{eq:Yang-Mills_heat_equation_as_perturbation_rough_Laplacian_plus_one_heat_equation} of the form
$$
A(t) = A_1 + a(t) = A_1 + \alpha(t) + \xi(t), \quad t\in [0,\tau),
$$
where $\alpha(t) \in \Omega^1(X;\ad P)$ for $t \geq 0$ is the unique solution to the linear Cauchy problem \eqref{eq:Struwe_19_with_augmented_connection_Laplacian} and $\xi(t) \in H_{A_1}^1(X;\Lambda^1\otimes\ad P)$ for $t\in [0,\tau)$ is a perturbation with initial data,
\begin{equation}
\label{eq:Yang-Mills_heat_equation_H1_small_initial_data_condition}
\xi(0) = \xi_0 \in H_{A_1}^1(X;\Lambda^1\otimes\ad P),
\end{equation}
required to solve the Yang-Mills heat equation \eqref{eq:Yang-Mills_heat_equation_as_perturbation_rough_Laplacian_plus_one_heat_equation}.

By substituting $a(t) = A_1 + \alpha(t) + \xi(t)$ into the Yang-Mills heat equation \eqref{eq:Yang-Mills_heat_equation_as_perturbation_rough_Laplacian_plus_one_heat_equation} with nonlinearity \eqref{eq:Yang-Mills_heat_equation_nonlinearity_relative_rough_Laplacian_plus_one}, we see that $a(t)$ solves \eqref{eq:Yang-Mills_heat_equation_as_perturbation_rough_Laplacian_plus_one_heat_equation} if and only if $\xi(t)$ solves, a.e. on $(0, \tau)\times X$,
\begin{multline*}
\frac{\partial\alpha}{\partial t} + \left(\nabla_{A_1}\nabla_{A_1}^* + 1\right)\alpha
+ \frac{\partial\xi}{\partial t} + \left(\nabla_{A_1}\nabla_{A_1}^* + 1\right)\xi
\\
+ d_{A_1}^*F_{A_1} + (F_{A_1} - 1)\times \alpha + \Ric_g\times \alpha + (F_{A_1} - 1)\times \xi + \Ric_g\times \xi
\\
+ \nabla_{A_1}\alpha\times \alpha + \nabla_{A_1}\alpha\times \xi + \nabla_{A_1}\xi\times \alpha + \nabla_{A_1}\xi\times\xi
\\
+ \alpha\times \alpha \times \alpha + \alpha\times \alpha \times \xi + \alpha\times \xi \times \xi
+ \xi\times \xi \times \xi = 0.
\end{multline*}
Using the fact that $\alpha(t)$ solves the linear Cauchy problem \eqref{eq:Struwe_19_with_augmented_connection_Laplacian}, we see that $a(t)$ solves the Yang-Mills heat equation \eqref{eq:Yang-Mills_heat_equation_as_perturbation_rough_Laplacian_plus_one_heat_equation} with initial condition $a(0)=a_0$ if and only if $\xi(t)$ solves the residual quasi-linear parabolic equation,
\begin{multline}
\label{eq:Yang-Mills_heat_equation_as_perturbation_rough_Laplacian_plus_one_heat_equation_residual}
\frac{\partial\xi}{\partial t} + \left(\nabla_{A_1}\nabla_{A_1}^* + 1\right)\xi
+ d_{A_1}^*F_{A_1} + (F_{A_1} - 1)\times \alpha + \Ric_g\times \alpha
+ \nabla_{A_1}\alpha\times \alpha + \alpha\times \alpha \times \alpha
\\
+ \Ric_g\times\xi + (F_{A_1}-1) \times \xi + \nabla_{A_1}\alpha\times \xi + \alpha\times \nabla_{A_1}\xi + \alpha\times \alpha \times \xi
\\
+ \nabla_{A_1}\xi\times\xi + \alpha\times \xi \times \xi + \xi\times \xi \times \xi = 0
\quad\hbox{a.e. on } (0, \tau)\times X,
\end{multline}
with initial condition $\xi(0)=\xi_0$.

Because the initial data, $\alpha(0) = \alpha_0$, is $C^\infty$, the solution, $\alpha(t)$, to the linear Cauchy problem \eqref{eq:Struwe_19_with_augmented_connection_Laplacian} defining the source term and coefficients of the quasi-linear parabolic equation \eqref{eq:Yang-Mills_heat_equation_as_perturbation_rough_Laplacian_plus_one_heat_equation_residual} for $\xi(t)$ is bounded via the \apriori estimates \eqref{eq:Sell_You_42-30_heat_equation_alpha_is_real_compact_form} and \eqref{eq:Sell_You_42-30_heat_equation_L2-Halpha-1_dudt_compact_form_finite_T0} provided by Corollary \ref{cor:Sell_You_42-12_heat_equation_alpha_is_real_apriori_estimates},
\begin{equation}
\label{eq:Sell_You_42-30_heat_equation_nonnegative_integer_k_compact_form}
\|\alpha\|_{L^2(0, \tau; H_{A_1}^{k+2}(X))} + \|\alpha\|_{L^\infty(0, \tau; H_{A_1}^{k+1}(X))}
+ \|\alpha\|_{H^1(0, \tau; H_{A_1}^k(X))}
\leq
\bar C_2\|\alpha_0\|_{H_{A_1}^{k+1}(X))},
\end{equation}
where $k\in\NN$ and $\bar C_2 = \bar C_2(A_1,g,k)$ and we have used the fact that $\tau \in (0,1]$. Because $X$ has dimension $2\leq d\leq 4$, there are Sobolev embeddings, $H^k(X) \hookrightarrow C^{k-3}(X)$ for all integers $k \geq 3$, by \cite[Theorem 4.12]{AdamsFournier}. Consequently, there are embeddings,
\begin{equation}
\label{eq:Sobolev_embedding_HAk_into_Ck_minus_3_dimX_leq_4}
H_{A_1}^k(X; \Lambda^1\otimes\ad P) \hookrightarrow C_{A_1}^{k-3}(X; \Lambda^1\otimes\ad P), \quad\forall\, k \geq 3,
\end{equation}
with embedding constants depending at most on the connection, $A_1$, the integer, $k$, and the Riemannian metric, $g$, on $X$. Thus,
\begin{equation}
\label{eq:Sobolev_embedding_norm_HAk_alpha_into_Ck_minus_3_dimX_leq_4}
\|\alpha\|_{L^\infty(0, \tau; W_{A_1}^{\infty,k-3}(X))} \leq C\|\alpha\|_{L^\infty(0, \tau; H_{A_1}^k(X))},
\quad\forall\, k \geq 3,
\end{equation}
for a positive constant, $C = C(A_1,g,k)$. This concludes our development of the \apriori estimates we shall need for the coefficients, $\alpha(t)$, for $t\in [0, \tau]$.

The quasi-linear parabolic equation \eqref{eq:Yang-Mills_heat_equation_as_perturbation_rough_Laplacian_plus_one_heat_equation_residual} for $\xi(t)$ has the same structure as that of the Yang-Mills heat equation
\eqref{eq:Yang-Mills_heat_equation_as_perturbation_rough_Laplacian_plus_one_heat_equation} for $a(t)$ considered in Theorem \ref{thm:Struwe_section_4-3_local_wellposedness_Yang-Mills_heat_equation_small_initial_data_in_H1}, aside from the harmless addition of a quadratic term, $\alpha\times\xi\times\xi$, and the presence of time-varying coefficients, $\alpha(t)$. Hence, for the proof of Theorem \ref{thm:Struwe_section_4-3_local_wellposedness_Yang-Mills_heat_equation_small_initial_data_in_H1} to apply to \eqref{eq:Yang-Mills_heat_equation_as_perturbation_rough_Laplacian_plus_one_heat_equation_residual} in place of \eqref{eq:Yang-Mills_heat_equation_as_perturbation_rough_Laplacian_plus_one_heat_equation}, it suffices to check that the additional terms involving $\alpha(t)$ do not cause any new difficulties. By appealing to the Sobolev embedding $H^1(X) \hookrightarrow L^4(X)$ \cite[Theorem 4.12]{AdamsFournier} and the Kato Inequality \eqref{eq:FU_6-20_first-order_Kato_inequality} as needed, we find that the new inhomogeneous terms are bounded via:
\begin{align*}
\|\nabla_{A_1}\alpha\times \alpha\|_{L^2(0, \tau; L^2(X))}
&\leq
c\sqrt{\tau} \|\alpha\|_{L^\infty(0, \tau; H_{A_1}^1(X))} \|\alpha\|_{L^\infty(0, \tau; C(X))},
\\
\|\alpha\times \alpha \times \alpha\|_{L^2(0, \tau; L^2(X))}
&\leq
c\sqrt{\tau} \|\alpha\|_{L^\infty(0, \tau; C(X))}^3.
\end{align*}
The new linear terms are bounded via:
\begin{align*}
\|\nabla_{A_1}\alpha\times \xi\|_{L^2(0, \tau; L^2(X))}
&\leq
c\sqrt{\tau} \|\nabla_{A_1}\alpha\|_{L^\infty(0, \tau; L^4(X))} \|\xi\|_{L^\infty(0, \tau; L^4(X))}
\\
&\leq
c\sqrt{\tau} \|\alpha\|_{L^\infty(0, \tau; H_{A_1}^2(X))} \|\xi\|_{L^\infty(0, \tau; H_{A_1}^1(X))},
\\
\|\alpha\times\nabla_{A_1}\xi\|_{L^2(0, \tau; L^2(X))}
&\leq
c\sqrt{\tau} \|\alpha\|_{L^\infty(0, \tau; C(X))} \|\xi\|_{L^\infty(0, \tau; H_{A_1}^1(X))},
\\
\|\alpha\times \alpha \times \xi\|_{L^2(0, \tau; L^2(X))}
&\leq
c\sqrt{\tau} \|\alpha\|_{L^\infty(0, \tau; C(X))}^2 \|\xi\|_{L^\infty(0, \tau; L^2(X))}.
\end{align*}
Finally, the new quadratic term is bounded via:
\begin{align*}
\|\alpha\times \xi \times \xi\|_{L^2(0, \tau; L^2(X))}
&\leq
c\sqrt{\tau} \|\alpha\|_{L^\infty(0, \tau; C(X))} \|\xi\|_{L^\infty(0, \tau; L^4(X))}^2
\\
&\leq
c\sqrt{\tau} \|\alpha\|_{L^\infty(0, \tau; C(X))} \|\xi\|_{L^\infty(0, \tau; H_{A_1}^1(X))}^2.
\end{align*}
By definition \eqref{eq:Struwe_page_129_Banach_space} of the Banach space, $\fV$, and the preceding estimates for the coefficients $\alpha$ appearing in \eqref{eq:Yang-Mills_heat_equation_as_perturbation_rough_Laplacian_plus_one_heat_equation_residual}, we see that
\begin{equation}
\label{eq:Struwe_20_initial_data_alpha0_arbitrary_H1_norm}
\begin{aligned}
\|\nabla_{A_1}\alpha\times w\|_{L^2(0, \tau; L^2(X))}
&\leq C_5\sqrt{\tau}\|w\|_\fV,
\\
\|\alpha\times \nabla_{A_1}w\|_{L^2(0, \tau; L^2(X))}
&\leq C_5\sqrt{\tau}\|w\|_\fV,
\\
\|\alpha\times \alpha\times w\|_{L^2(0, \tau; L^2(X))}
&\leq C_5\sqrt{\tau}\|w\|_\fV,
\\
\|\alpha\times w\times w\|_{L^2(0, \tau; L^2(X))}
&\leq C_5\sqrt{\tau}\|w\|_\fV^2,
\\
\|\nabla_{A_1}\alpha\times \alpha\|_{L^2(0, \tau; L^2(X))}
&\leq C_5\sqrt{\tau},
\\
\|\alpha\times \alpha\times \alpha\|_{L^2(0, \tau; L^2(X))}
&\leq C_5\sqrt{\tau},
\end{aligned}
\end{equation}
for a positive constant $C_5 = C_5(A_1,g,\|\alpha_0\|_{H_{A_1}^4(X)}) \geq 1$.

Given $w \in \fV$, we shall now define $\xi := \Psi(w) \in \fV$ to be the unique strong solution to the \emph{linear} Cauchy problem provided by Theorem \ref{thm:Sell_You_42-12_heat_equation_alpha_is_one} and Corollary \ref{cor:Sell_You_42-13_heat_equation},
\begin{equation}
\label{eq:Yang-Mills_heat_equation_as_perturbation_rough_Laplacian_plus_one_heat_equation_residual_linear}
\frac{\partial\xi}{\partial t} + \nabla_{A_1}\nabla_{A_1}^*\xi = f_\alpha \quad\hbox{a.e. on } (0,\tau)\times X,
\quad \xi(0) = \xi_0,
\end{equation}
where
$$
f_\alpha(t) := \sF_\alpha(t,w(t)),  \quad\hbox{a.e. } t \in (0,\tau),
$$
and the nonlinearity is defined through \eqref{eq:Yang-Mills_heat_equation_as_perturbation_rough_Laplacian_plus_one_heat_equation_residual} by
\begin{equation}
\label{eq:Yang-Mills_heat_equation_as_perturbation_rough_Laplacian_plus_one_heat_equation_residual_nonlinearity}
\begin{aligned}
-\sF_\alpha(t,w(t)) &:= d_{A_1}^*F_{A_1} + \nabla_{A_1}\alpha(t)\times \alpha(t) + \alpha(t)\times \alpha(t) \times \alpha(t)
\\
&\quad + (F_{A_1} - 1)\times w(t) + \Ric_g\times w(t)
\\
&\quad + \nabla_{A_1}\alpha(t)\times w(t) + \alpha(t)\times \nabla_{A_1}w(t)
 + \alpha(t)\times \alpha(t) \times w(t)
\\
&\quad + \nabla_{A_1}w(t)\times w(t) + \alpha\times w(t) \times w(t) + w(t)\times w(t) \times w(t).
\end{aligned}
\end{equation}
The verification that
$$
f_\alpha = \sF_\alpha(\cdot, w) \in L^2(0, \tau; L^2(X; \Lambda^1\otimes\ad P))
$$
follows just as in the case of \eqref{eq:Assumption_sF(w)_in_L2from_0_to_T_into_L^2(X)_given_w_infH}.
We thus aim to solve the quasi-linear parabolic equation
\eqref{eq:Yang-Mills_heat_equation_as_perturbation_rough_Laplacian_plus_one_heat_equation_residual}
for $\xi$, with initial condition $\xi(0) = \xi_0$, as the fixed point of a map, $\Psi:\fB_R \to \fB_R$, for a suitably chosen radius, $R = R(A_1,g,\|\alpha_0\|_{H_{A_1}^4(X)}) \in (0,1]$. The new dependency of $R$ on $\|\alpha_0\|_{H_{A_1}^4(X)}$ arises from the estimates of the coefficients, $\alpha$, in \eqref{eq:Struwe_20_initial_data_alpha0_arbitrary_H1_norm}.

Given the estimates in \eqref{eq:Struwe_20_initial_data_alpha0_arbitrary_H1_norm} for the terms
in the residual nonlinearity \eqref{eq:Yang-Mills_heat_equation_as_perturbation_rough_Laplacian_plus_one_heat_equation_residual_nonlinearity}
involving the coefficients $\alpha(t)$, which were not present in the
nonlinearity \eqref{eq:Yang-Mills_heat_equation_nonlinearity_relative_rough_Laplacian_plus_one},
the proof of Theorem \ref{thm:Struwe_section_4-3_local_wellposedness_Yang-Mills_heat_equation_small_initial_data_in_H1} now applies to give a unique solution, $\xi \in \fV$, to \eqref{eq:Yang-Mills_heat_equation_as_perturbation_rough_Laplacian_plus_one_heat_equation_residual} for $\xi_0$ obeying \eqref{eq:Yang-Mills_heat_equation_H1_small_initial_data_bound} and possibly smaller constants,
$$
\eps = \eps(A_1,g) \in (0, 1] \quad\hbox{and}\quad \tau = \tau(A_0,A_1,\alpha_0,g) \in (0, 1],
$$
and possibly larger constant $C = C(A_1,g) \geq 1$, with
\begin{equation}
\label{eq:Struwe_section_4-3_apriori_estimate_Yang-Mills_heat_equation_small_initial_data_in_H1_xi}
\|\xi\|_\fV  \leq C\eps.
\end{equation}
Because $a = \alpha+\xi$ on $[0,\tau)\times X$, we find that, for $\bar C_2 = \bar C_2(A_1,g)$,
\begin{align*}
\|a\|_\fV &\leq \|\alpha\|_\fV + \|\xi\|_\fV
\\
&\leq \bar C_2\|\alpha_0\|_{H_{A_1}^1(X)} + C\eps
\quad\hbox{(by \eqref{eq:Struwe_section_4-3_apriori_estimate_Yang-Mills_heat_equation_small_initial_data_in_H1_xi},
\eqref{eq:Sell_You_42-30_heat_equation_nonnegative_integer_k_compact_form} with $k=0$, and definition \eqref{eq:Struwe_page_129_Banach_space} of $\fV$)}
\\
&\leq \bar C_2\|a_0\|_{H_{A_1}^1(X)} + \bar C_2\|\xi_0\|_{H_{A_1}^1(X)} + C\eps
\quad\hbox{(by \eqref{eq:Decomposition_initial_data_a0_as_large_smooth_alpha0_plus_small_H1_xi0})}
\\
&\leq \bar C_2\|a_0\|_{H_{A_1}^1(X)} + (\bar C_2 + C)\eps
\quad\hbox{(by \eqref{eq:Yang-Mills_heat_equation_H1_small_initial_data_bound})},
\end{align*}
which gives the desired \apriori estimate in Theorem \ref{thm:Struwe_section_4-3_local_wellposedness_Yang-Mills_heat_equation} below.

To verify the continuity of the solution, $a(t)$, with respect to the initial data, $a_0$, as asserted by Theorem \ref{thm:Struwe_section_4-3_local_wellposedness_Yang-Mills_heat_equation}, suppose we are given $A_{0,i} = A_1+a_{0,i}$ for $a_{0,i} \in H^1_{A_1}(X; \Lambda^1\otimes\ad P)$, with $i=1,2$. Then, following \eqref{eq:Decomposition_initial_data_a0_as_large_smooth_alpha0_plus_small_H1_xi0}, we write
$$
A_{0,i} = A_1+a_{0,i} = A_1 + \alpha_0 + \xi_{0,i}, \quad\hbox{for } i = 1,2,
$$
where $\alpha_0 \in \Omega^1(X;\ad P)$ is chosen so that the $\xi_{0,i} \in H^1_{A_1}(X; \Lambda^1\otimes\ad P)$ obey \eqref{eq:Yang-Mills_heat_equation_H1_small_initial_data_bound} for $i=1,2$. Let $\xi_i$, for $i=1,2$, be the unique solutions to
\eqref{eq:Yang-Mills_heat_equation_as_perturbation_rough_Laplacian_plus_one_heat_equation_residual}, with initial conditions $\xi_i(0) = \xi_{i,0}$, produced by the preceding analysis. Let $a_i = \alpha+\xi_i$, for $i=1,2$, be the corresponding unique solutions to \eqref{eq:Yang-Mills_heat_equation_as_perturbation_rough_Laplacian_plus_one_heat_equation}, with initial conditions $a_i(0) = a_{i,0}$. Then,
\begin{align*}
\|a_1 - a_2\|_\fV &= \|\xi_1 - \xi_2\|_\fV
\\
&\leq C\|\xi_1(0) - \xi_2(0)\|_{H^1_{A_1}(X)}
\quad\hbox{(by \eqref{eq:Struwe_section_4-3_continuity_estimate_Yang-Mills_heat_equation})}
\\
&= C\|a_1(0) - a_2(0)\|_{H^1_{A_1}(X)},
\end{align*}
where $C = C(A_1,g)$. In summary, we have proved

\begin{thm}[Local well-posedness for the Yang-Mills heat equation over a closed manifold with arbitrary initial data in $H^1$]
\label{thm:Struwe_section_4-3_local_wellposedness_Yang-Mills_heat_equation}
Let $G$ be a compact Lie group and $P$ a principal $G$-bundle over a closed, connected smooth manifold, $X$, of dimension $2 \leq d \leq 4$ and Riemannian metric, $g$. If $A_1$ is a reference connection of class $C^\infty$ on $P$ and $A_0$ is a connection of class $H^1$ on $P$, then there are positive constants, $C = C(A_1,g) \geq 1$ and $\eps = \eps(A_1,g) \in (0, 1]$, and if $\alpha_0 \in \Omega^1(X;\ad P)$ is such that
$$
\|A_0 - A_1 - \alpha_0\|_{H^1_{A_1}(X)} < \eps,
$$
there is a positive constant, $\tau = \tau(A_0,A_1,\alpha_0,g) \in (0, 1]$, with the following significance. If
$$
a_0 := A_0-A_1 \in H^1_{A_1}(X; \Lambda^1\otimes\ad P),
$$
then there is a unique strong solution, $a(t)$ in $\fV$, to the Yang-Mills heat equation
\eqref{eq:Yang-Mills_heat_equation_as_perturbation_rough_Laplacian_plus_one_heat_equation} on the interval $(0, \tau)$ with initial data $a(0) = a_0$, where $\fV$ is the Banach space \eqref{eq:Struwe_page_129_Banach_space}. Moreover, the solution $a(t)$ obeys the \apriori estimate,
\begin{equation}
\label{eq:Struwe_section_4-3_apriori_estimate_Yang-Mills_heat_equation}
\|a\|_\fV \leq C\left(\|a_0\|_{H^1_{A_1}(X)} + \eps\right).
\end{equation}
If $a_1, a_2$ are the unique solutions to
\eqref{eq:Yang-Mills_heat_equation_as_perturbation_rough_Laplacian_plus_one_heat_equation}
with initial data $a_1(0)$ and $a_2(0)$, respectively, and
$$
\|a_1(0) - a_2(0)\|_{H^1_{A_1}(X)} \leq \frac{\eps}{2},
$$
then $a_1 - a_2$ obeys the continuity estimate \eqref{eq:Struwe_section_4-3_continuity_estimate_Yang-Mills_heat_equation}.
\end{thm}

\subsection{Regularity of solutions to the Yang-Mills heat equation over a closed manifold with dimension less than or equal to four and initial data in $H^1$}
\label{subsec:Struwe_page_137_contraction_mapping_regularity_initial_data_in_H1}
We recall from Section \ref{subsec:Struwe_page_137_contraction_mapping_small_initial_data_in_H1} that --- given $a\in \fV$, where $\fV$ is the Banach space defined in \eqref{eq:Struwe_page_129_Banach_space} --- the nonlinearity $\sF(a)$ in
\eqref{eq:Yang-Mills_heat_equation_nonlinearity_relative_rough_Laplacian_plus_one}
belongs to $L^2(0, \tau; L^2(X; \Lambda^1\otimes\ad P))$ by \eqref{eq:Assumption_sF(w)_in_L2from_0_to_T_into_L^2(X)_given_w_infH}. Consequently, the question of regularity of the unique solution, $a(t)$, to the Yang-Mills heat equation
\eqref{eq:Linear_heat_equation_with_rough_Laplacian_plus_one_on_Omega_1_adP} provided by Theorem \ref{thm:Struwe_section_4-3_local_wellposedness_Yang-Mills_heat_equation} is more subtle, at least when $X$ has dimension $d=4$, then in the case of the linear heat equation since there is no improvement in the regularity of source function, in this case $\sF(a)$, on intervals $[t_0,\tau)$ with $0 < t_0 < \tau$.

To circumvent this problem, we first appeal to Corollary \ref{cor:Sell_You_42-13_heat_equation}, with $p=2$, $\alpha=1$, and $\beta \in [0,1/2)$, to conclude from \eqref{eq:Sell_You_42-36_heat_equation} that
$$
a \in C([0, \tau); H_{A_1}^1(X;\Lambda^1\otimes\ad P))
\cap C((0, \tau); H_{A_1}^{2(1+\beta)}(X;\Lambda^1\otimes\ad P)).
$$
(A stronger assertion, involving H\"older regularity with respect to $t\in [0,\tau)$ is also possible, as evident from \eqref{eq:Sell_You_42-36_heat_equation}, but we shall not need this property.) By restricting to $\beta \in (0, 1/2)$, we recover the fact that $H^{2(1+\beta)}(X)$ is a Banach algebra \cite[Theorem 4.39]{ AdamsFournier} when $2\leq d\leq 4$ and hence
$$
\sF(a) \in C((0, \tau); H_{A_1}^{2(1+\beta)}(X;\Lambda^1\otimes\ad P)).
$$
To obtain higher-order temporal regularity for $t>0$, we observe that the expression \eqref{eq:Yang-Mills_heat_equation_nonlinearity_relative_rough_Laplacian_plus_one} for $\sF(a)$ yields
\begin{multline*}
-\frac{\partial \sF(a)}{\partial t}
=
(F_{A_1} - 1)\times \dot a + \Ric_g\times \dot a + \nabla_{A_1}\dot a\times a
+ \nabla_{A_1}a\times\dot a
\\
+ \dot a\times a\times a + a\times \dot a\times a + a\times a\times \dot a.
\end{multline*}
Theorem \ref{thm:Struwe_section_4-3_local_wellposedness_Yang-Mills_heat_equation} already ensures that (by definition \eqref{eq:Struwe_page_129_Banach_space} of $\fV$)
$$
\dot a \in L^2(0, \tau; L^2(X; \Lambda^1\otimes\ad P)),
$$
and so, using the fact that $H^{2(1+\beta)}(X)$ is a Banach algebra and $L^2(X)$ is a $H^{2(1+\beta)}(X)$-module, we discover that
$$
\frac{\partial \sF(a)}{\partial t} \in L_{\loc}^2(0, \tau; L^2(X; \Lambda^1\otimes\ad P)),
$$
and thus
$$
\sF(a) \in H_{\loc}^1(0, \tau; L^2(X; \Lambda^1\otimes\ad P)).
$$
In particular, for any $t_0 \in (0,\tau)$, we have
$$
\sF(a) \in C([t_0, \tau); H_{A_1}^1(X;\Lambda^1\otimes\ad P))
\cap H^1(t_0, \tau; L^2(X; \Lambda^1\otimes\ad P)).
$$
(In fact, we also have $\sF(a) \in C([t_0, \tau); H_{A_1}^{2(1+\beta)}(X;\Lambda^1\otimes\ad P))$.) Therefore, applying Theorem \ref{thm:Sell_You_42-15_heat_equation} with interval $[t_0,\tau)$, initial data $a(t_0)\in H_{A_1}^1(X;\Lambda^1\otimes\ad P)$ (in fact, $a(t_0)\in H_{A_1}^{2(1+\beta)}(X;\Lambda^1\otimes\ad P)$), source function $f = \sF(a)$, and $\alpha=1$ yields
$$
a \in C^1((t_0, \tau); H_{A_1}^1(X;\Lambda^1\otimes\ad P))
\cap C((t_0, \tau); H_{A_1}^3(X; \Lambda^1\otimes\ad P)).
$$
It is now routine to iterate the preceding argument to obtain any desired spatial and temporal regularity for the solution, $a(t)$ for $t\in (0,\tau)$, and conclude with the

\begin{thm}[Higher-order spatial and temporal regularity for a solution to the Yang-Mills heat equation over a closed manifold with arbitrary initial data in $H^1$]
\label{thm:Struwe_section_4-3_regularity_Yang-Mills_heat_equation}
Assume the hypotheses of Theorem \ref{thm:Struwe_section_4-3_local_wellposedness_Yang-Mills_heat_equation}. If $a_0 := A_1-A_0 \in H_{A_1}^1(X;\Lambda^1\otimes\ad P)$ and $a\in\fV$ is a strong solution to the Yang-Mills heat equation \eqref{eq:Linear_heat_equation_with_rough_Laplacian_plus_one_on_Omega_1_adP}, where $\fV$ is as in \eqref{eq:Struwe_page_129_Banach_space}, then $a(t)$ is a classical solution and obeys
$$
a \in C([0,\tau); H_{A_1}^1(X;\Lambda^1\otimes\ad P)) \cap C^\infty((0,\tau)\times X; \Lambda^1\otimes\ad P)).
$$
\end{thm}

\subsection{Local well-posedness for the Yang-Mills heat equation in hybrid critical-exponent parabolic Sobolev spaces}
\label{subsec:Struwe_page_137_contraction_mapping_hybrid_critical-exponent_parabolic_Sobolev_space}
Our goal in this subsection is prove the following version of the local well-posedness for the Yang-Mills heat equation
\eqref{eq:Yang-Mills_heat_equation_as_perturbation_rough_Laplacian_plus_one_heat_equation} with initial data in $H^1_{A_1}(X; \Lambda^1\otimes\ad P) \cap L^{2\lozenge}(X; \Lambda^1\otimes\ad P)$. See Brezis and Nirenberg \cite{Brezis_Nirenberg_1983} for an early discussion in the literature of critical-exponent Sobolev problems and some of the difficulties therein.

\begin{thm}[Local well-posedness for the Yang-Mills heat equation in hybrid critical-exponent parabolic Sobolev spaces over a closed manifold]
\label{thm:Struwe_section_4-3_local_wellposedness_Yang-Mills_heat_equation_hybrid_critical-exponent_parabolic_Sobolev_space}
Let $G$ be a compact Lie group and $P$ a principal $G$-bundle over a closed, connected, orientable, smooth manifold, $X$, of dimension $d \geq 2$ and Riemannian metric, $g$. If $A_1$ is a reference connection of class $C^\infty$ on $P$, then there is a positive constant $C$ and, if $b$ is a positive constant, then there is also a constant, $\tau \in (0, 1]$, with the following significance. If
$$
a_0 \in H^1_{A_1}(X; \Lambda^1\otimes\ad P) \cap C(X; \Lambda^1\otimes\ad P),
$$
obeys
\begin{equation}
\label{eq:Norm_ao_H1_and_sup_tx_L2_a0_dKtx_leq_b}
\|a_0\|_{H^1_{A_1}(X)} + \sup_{(t, x) \in (0,\tau)\times X} \|a_0\|_{L^2(X, K_{t,x})} \leq b,
\end{equation}
then there is a unique strong solution $a(t)$ in $\fV$ to the Yang-Mills heat equation
\eqref{eq:Yang-Mills_heat_equation_as_perturbation_rough_Laplacian_plus_one_heat_equation} on the interval $(0, \tau)$ with initial data $a(0) = a_0$, where $\fV$ is the Banach space completion of $C^\infty([0,\tau]\times X; \Lambda^1\otimes \ad P))$ with respect to the norm \eqref{eq:Parabolic_Sobolev_norm_hybrid_standard_H2_and_critical_exponent_domain}. Moreover, the solution $a(t)$ obeys the \apriori estimate,
\begin{equation}
\label{eq:Struwe_section_4-3_apriori_estimate_Yang-Mills_heat_equation_hybrid_critical-exponent_parabolic_Sobolev_space}
\|a\|_\fV \leq 32\sqrt{\tau} \|d_{A_1}^*F_{A_1}\|_{C(X)}
+ 8\|a_0\|_{H^1_{A_1}(X)} + \sup_{(t, x) \in (0,\tau)\times X} 8 \|a_0\|_{L^2(X, K_{t,x})},
\end{equation}
and the minimal lifetime, $\tau$, is given by
\begin{equation}
\label{eq:Struwe_section_4-3_minimallifetime_Yang-Mills_heat_equation_hybrid_critical-exponent_parabolic_Sobolev_space}
\tau = Cb\left(\|d_{A_1}^*F_{A_1}\|_{C(X)} + \left(\|\Ric_g\|_{C(X)} + \|F_{A_1}\|_{C(X)} + 1\right)(1 + b^3)\right)^{-2}.
\end{equation}
If $a_1, a_2$ are the unique solutions corresponding to initial data $a_1(0)$ and $a_2(0)$, respectively, then
\begin{equation}
\label{eq:Struwe_section_4-3_continuity_estimate_Yang-Mills_heat_equation_hybrid_critical-exponent_parabolic_Sobolev_space}
\|a_1 - a_2\|_\fV \leq 8\|a_1(0) - a_2(0)\|_{H^1_{A_1}(X)} + \sup_{(t, x) \in (0,\tau)\times X} 8 \|a_1(0) - a_2(0)\|_{L^2(X, K_{t,x})}.
\end{equation}
\end{thm}

\begin{rmk}[Validity of Theorem \ref{thm:Struwe_section_4-3_local_wellposedness_Yang-Mills_heat_equation_hybrid_critical-exponent_parabolic_Sobolev_space} for manifolds of dimension $d \geq 2$]
While Theorem \ref{thm:Struwe_section_4-3_local_wellposedness_Yang-Mills_heat_equation_hybrid_critical-exponent_parabolic_Sobolev_space} is of most interest when $d = 4$, it is valid for a manifold $X$ of any dimension $d \geq 2$, since we do not need to appeal to the Sobolev Embedding Theorem \cite[Theorem 4.12]{AdamsFournier} explicitly or implicitly in its proof.
\end{rmk}

\begin{rmk}[Norms for the initial data]
By the inequality \eqref{eq:Sobolev_embedding_Linfty_time_and_space_into_sup_time_and_space_L2space_dKtx}, we could replace the norms $\sup_{(t, x) \in (0,\tau)\times X}\|\cdot\|_{L^2(X, K_{t,x})}$ in the statement of Theorem
\ref{thm:Struwe_section_4-3_local_wellposedness_Yang-Mills_heat_equation_hybrid_critical-exponent_parabolic_Sobolev_space} by the simpler $\|\cdot\|_{C(X;E)}$ norms when referring to bounds involving initial data.
\end{rmk}

\begin{proof}[Proof of Theorem \ref{thm:Struwe_section_4-3_local_wellposedness_Yang-Mills_heat_equation_hybrid_critical-exponent_parabolic_Sobolev_space}]
For completeness, we shall include some of the steps already verified in Section \ref{subsec:Struwe_page_137_contraction_mapping_small_initial_data_in_H1}, but now making use of the stronger system of norms to simplify the estimate of the Yang-Mills heat equation nonlinearity
\eqref{eq:Yang-Mills_heat_equation_nonlinearity_relative_rough_Laplacian_plus_one}.

\setcounter{step}{0}
\begin{step}[Estimates in $L^2(0,\tau; L^2(X; \Lambda^1\otimes\ad P))$ and $L^\lozenge((0,\tau)\times X; \Lambda^1\otimes\ad P)$ for the quadratic term in the Yang-Mills heat equation nonlinearity]
We observe that
\begin{align*}
\|\nabla_{A_1}w\times w\|_{L^2(0,\tau; L^2(X))} &\leq c\|\nabla_{A_1}w\|_{L^2(0,\tau; L^2(X))} \|w\|_{C([0,\tau]; C(X))}
\\
&\leq c\|w\|_{L^2(0,\tau; H^1_{A_1}(X))} \|w\|_{C([0,\tau]\times X)},
\end{align*}
and thus,
\begin{equation}
\label{eq:Quadratic_term_norm_L2time_and_L2space_dimension_4}
\|\nabla_{A_1}w\times w\|_{L^2(0,\tau; L^2(X))} \leq c\sqrt{\tau}\|w\|_{C([0,\tau]; H^1_{A_1}(X))} \|w\|_{C([0,\tau]\times X)}.
\end{equation}
Similarly, the H\"older inequality and the definitions \eqref{eq:Critical-exponent_parabolic_Sobolev_spaces} of the $L^\lozenge$ and $L^{2\lozenge}$ norms gives
$$
\|\nabla_{A_1}w\times w\|_{L^\lozenge((0,\tau)\times X)}
\leq c\|\nabla_{A_1}w\|_{L^{2\lozenge}((0,\tau)\times X)} \|w\|_{L^{2\lozenge}((0,\tau)\times X)},
$$
and so, by  \eqref{eq:Sobolev_embedding_Linfty_time_and_space_into_L2lozenge},
\begin{equation}
\label{eq:Quadratic_term_norm_Llozenge_time_and_space_dimension_4}
\|\nabla_{A_1}w\times w\|_{L^\lozenge((0,\tau)\times X)}
\\
\leq c\sqrt{\tau}\|\nabla_{A_1}w\|_{L^{2\lozenge}((0,\tau)\times X)} \|w\|_{C([0,\tau]\times X)}.
\end{equation}
Combining \eqref{eq:Quadratic_term_norm_L2time_and_L2space_dimension_4} and \eqref{eq:Quadratic_term_norm_Llozenge_time_and_space_dimension_4} yields
\begin{multline}
\label{eq:Quadratic_term_norm_L2time_and_L2_space_Llozenge_time_and_space_dimension_4}
\|\nabla_{A_1}w\times w\|_{L^2(0,\tau; L^2(X))} + \|\nabla_{A_1}w\times w\|_{L^\lozenge((0,\tau)\times X)}
\\
\leq c\sqrt{\tau}\left(\|w\|_{L^2(0,\tau; H^1_{A_1}(X))} + \|\nabla_{A_1}w\|_{L^{2\lozenge}((0,\tau)\times X)} \right)
\|w\|_{C([0,\tau]\times X)}.
\end{multline}
The constant $c$ in the inequality \eqref{eq:Quadratic_term_norm_L2time_and_L2_space_Llozenge_time_and_space_dimension_4} depends at most on the Riemannian metric $g$ on $X$ but is independent of the connection $A_1$ on $P$. This completes our estimate for the quadratic term in the Yang-Mills nonlinearity \eqref{eq:Yang-Mills_heat_equation_nonlinearity_relative_rough_Laplacian_plus_one}.
\end{step}

\begin{step}[Estimates in $L^2(0,\tau; L^2(X; \Lambda^1\otimes\ad P))$ and $L^\lozenge((0,\tau)\times X; \Lambda^1\otimes\ad P)$ for the cubic term in the Yang-Mills nonlinearity]
We can partially repair the argument in \cite[Section 4.3]{Struwe_1994} by instead availing of the system of critical-exponent parabolic Sobolev norms that we described in Section \ref{sec:Critical-exponent_parabolic_Sobolev_spaces_linear_parabolic_operator_vector_bundle_manifold} to estimate
$$
\|w\times w\times w\|_{L^2(0,\tau; L^2(X))} \leq c\sqrt{\tau}\|w\times w\times w\|_{C([0,\tau]; C(X))}
$$
and thus,
\begin{equation}
\label{eq:Cubic_term_norm_L2time_and_L2space_dimension_4}
\|w\times w\times w\|_{L^2(0,\tau; L^2(X))} \leq c\sqrt{\tau}\|w\|_{C([0,\tau]\times X)}^3.
\end{equation}
Similarly, by the definition \eqref{eq:Critical-exponent_parabolic_Sobolev_space_Llozenge_0_to_T_times_X} of the $L^\lozenge$ norm,
\begin{align*}
\|w\times w\times w\|_{L^\lozenge((0,\tau)\times X)}
&\leq
c\|w\|_{L^\lozenge((0,\tau)\times X)} \|w\times w\|_{C([0,\tau]\times X)}
\\
&\leq c\tau\|w\|_{C([0,\tau]\times X)} \|w\times w\|_{C([0,\tau]\times X)}
\quad \hbox{(by \eqref{eq:Sobolev_embedding_Linfty_time_and_space_Llozenge})},
\end{align*}
and thus,
\begin{equation}
\label{eq:Cubic_term_norm_Llozenge_time_and_space_dimension_4}
\|w\times w\times w\|_{L^\lozenge((0,\tau)\times X)} \leq c\tau\|w\|_{C([0,\tau]\times X)}^3.
\end{equation}
Combining \eqref{eq:Cubic_term_norm_L2time_and_L2space_dimension_4} and \eqref{eq:Cubic_term_norm_Llozenge_time_and_space_dimension_4} yields
\begin{equation}
\label{eq:Cubic_term_norm_L2time_and_L2_space_Llozenge_time_and_space_dimension_4}
\|w\times w\times w\|_{L^2(0,\tau; L^2(X))} + \|w\times w\times w\|_{L^\lozenge((0,\tau)\times X)}
\leq
c(\sqrt{\tau} + \tau)\|w\|_{C([0,\tau]\times X)}^3.
\end{equation}
The constant $c$ in the inequality \eqref{eq:Cubic_term_norm_L2time_and_L2_space_Llozenge_time_and_space_dimension_4} depends at most on the Riemannian metric $g$ on $X$ but is independent of the connection $A_1$ on $P$. This completes our estimate for the cubic term in the Yang-Mills nonlinearity \eqref{eq:Yang-Mills_heat_equation_nonlinearity_relative_rough_Laplacian_plus_one}.
\end{step}

Consequently, we can now hope to apply a contraction mapping argument to solve the Yang-Mills heat equation
\eqref{eq:Yang-Mills_heat_equation_as_perturbation_rough_Laplacian_plus_one_heat_equation} in the Banach space $\fV$ defined by the norm \eqref{eq:Parabolic_Sobolev_norm_hybrid_standard_H2_and_critical_exponent_domain}. We note that
\begin{subequations}
\label{eq:Critical_exponent_plus_standard_H2_parabolic_Sobolev_spaces}
\begin{align}
\label{eq:Critical_exponent_plus_standard_H2_parabolic_Sobolev_space_range}
\fW &= L^2(0,\tau; L^2(X; \Lambda^1\otimes\ad P)) \cap L^\lozenge((0,\tau)\times X; \Lambda^1\otimes\ad P),
\\
\label{eq:Critical_exponent_plus_standard_H2_parabolic_Sobolev_space_domain}
\fV &\subset C([0, \tau]\times X; \Lambda^1\otimes\ad P) \cap H^1(0, \tau; L^2(X; \Lambda^1\otimes\ad P))
\\
\notag
&\qquad \cap C([0, \tau]; H_A^1(X; \Lambda^1\otimes\ad P)) \cap L^2(0, \tau; H_A^2(X; \Lambda^1\otimes\ad P)),
\end{align}
\end{subequations}
with the norms defining $\fW$ and $\fV$ given by \eqref{eq:Parabolic_Sobolev_norms_hybrid_standard_and_critical_exponent}.

\begin{step}[Boundedness of the Yang-Mills nonlinearity in $\fW$]
Given $w \in \fV$, we first estimate the $L^2(0,\tau; L^2(X; \Lambda^1\otimes\ad P))$-norm of the nonlinearity $\sF(w)$ defined by
\eqref{eq:Yang-Mills_heat_equation_nonlinearity_relative_rough_Laplacian_plus_one} to give
\begin{align*}
{}&\|\sF(w)\|_{L^2(0,\tau; L^2(X))}
\\
&\quad \leq \|d_{A_1}^*F_{A_1}\|_{L^2(0,\tau; L^2(X))} + \|\Ric_g\times w\|_{L^2(0,\tau; L^2(X))} + \|(F_{A_1}-1)\times w\|_{L^2(0,\tau; L^2(X))}
\\
&\qquad + \|\nabla_{A_1}w\times w\|_{L^2(0,\tau; L^2(X))} + \|w\times w\times w\|_{L^2(0,\tau; L^2(X))}
\\
&\quad \leq \sqrt{\tau}\|d_{A_1}^*F_{A_1}\|_{L^2(X)}
+ c\sqrt{\tau}\left(\|\Ric_g\|_{L^2(X)} + \|F_{A_1}\|_{L^2(X)} + 1\right)\|w\|_{C([0,\tau]\times X)}
\\
&\qquad + c\sqrt{\tau}\|w\|_{C([0,\tau]; H^1_{A_1}(X))} \|w\|_{C([0,\tau]\times X)}
\\
&\qquad + c\sqrt{\tau}\|w\|_{C([0,\tau]\times X)}^3
\quad\hbox{(by \eqref{eq:Quadratic_term_norm_L2time_and_L2space_dimension_4} and \eqref{eq:Cubic_term_norm_L2time_and_L2space_dimension_4})},
\end{align*}
where the constant $c$ depends on the Riemannian metric $g$ on $X$ but is independent of the connection $A_1$ on $P$. By applying our definition \eqref{eq:Parabolic_Sobolev_norm_hybrid_standard_H2_and_critical_exponent_domain} of the Banach space norm on $\fV$ in the preceding inequality for $\|\sF(w)\|_{L^2(0,\tau; L^2(X))}$, we obtain
\begin{multline}
\label{eq:Struwe_nonlinear_map_bound_L2time_and_L2space_dimension_4}
\|\sF(w)\|_{L^2(0,\tau; L^2(X))} \leq c\sqrt{\tau}\left(\|\Ric_g\|_{L^2(X)} + \|F_{A_1}\|_{L^2(X)} + 1\right)\|w\|_\fV
\\
+ \sqrt{\tau}\|d_{A_1}^*F_{A_1}\|_{L^2(X)} + c\sqrt{\tau}\left(\|w\|_\fV^2 + \|w\|_\fV^3\right),
\end{multline}
where the constant $c$ depends on the Riemannian metric $g$ on $X$ but is independent of the connection $A_1$ on $P$. Similarly,
\begin{align*}
{}&\|\sF(w)\|_{L^\lozenge((0,\tau)\times X)}
\\
&\quad \leq \|d_{A_1}^*F_{A_1}\|_{L^\lozenge((0,\tau)\times X)}
+ \|\Ric_g\times w\|_{L^\lozenge((0,\tau)\times X)} + \|(F_{A_1}-1)\times w\|_{L^\lozenge((0,\tau)\times X)}
\\
&\qquad + \|\nabla_{A_1}w\times w\|_{L^\lozenge((0,\tau)\times X)} + \|w\times w\times w\|_{L^\lozenge((0,\tau)\times X)}
\\
&\quad \leq \tau\|d_{A_1}^*F_{A_1}\|_{C(X)}
+ c\tau\left(\|\Ric_g\|_{C(X)} + \|F_{A_1}\|_{C(X)} + 1\right)\|w\|_{C([0,\tau]\times X)}
\\
&\qquad + c\sqrt{\tau}\|\nabla_{A_1}w\|_{L^{2\lozenge}((0,\tau)\times X)} \|w\|_{C([0,\tau]\times X)}
\\
&\qquad + c\tau\|w\|_{C([0,\tau]\times X)}^3
\quad\hbox{(by \eqref{eq:Sobolev_embedding_Linfty_time_and_space_Llozenge},
\eqref{eq:Quadratic_term_norm_Llozenge_time_and_space_dimension_4}, and \eqref{eq:Cubic_term_norm_Llozenge_time_and_space_dimension_4})},
\end{align*}
where the constant $c$ depends on the Riemannian metric $g$ on $X$ but is independent of the connection $A_1$ on $P$. By applying our definition \eqref{eq:Parabolic_Sobolev_norm_hybrid_standard_H2_and_critical_exponent_domain} of the Banach space norm on $\fV$ in the preceding inequality for
$\|\sF(w)\|_{L^\lozenge((0,\tau)\times X)}$, we obtain
\begin{multline}
\label{eq:Struwe_nonlinear_map_bound_Llozenge_time_and_space_dimension_4}
\|\sF(w)\|_{L^\lozenge((0,\tau)\times X)}
\leq \tau\|d_{A_1}^*F_{A_1}\|_{C(X)} + c\tau\left(\|\Ric_g\|_{C(X)} + \|F_{A_1}\|_{C(X)} + 1\right)\|w\|_\fV
\\
+ c\left(\sqrt{\tau}\|w\|_\fV^2 + \tau\|w\|_\fV^3\right),
\end{multline}
where the constant $c$ depends on the Riemannian metric $g$ on $X$ but is independent of the connection $A_1$ on $P$. Combining \eqref{eq:Struwe_nonlinear_map_bound_L2time_and_L2space_dimension_4} and \eqref{eq:Struwe_nonlinear_map_bound_Llozenge_time_and_space_dimension_4} and recalling our definition of the norm \eqref{eq:Parabolic_Sobolev_norm_hybrid_standard_L2_and_critical_exponent_range} on $\fW$ yields
\begin{multline}
\label{eq:Struwe_nonlinear_map_bound_L2time_and_L2_space_and_Llozenge_time_and_space_dimension_4}
\|\sF(w)\|_\fW
\leq
(\sqrt{\tau} + \tau)\|d_{A_1}^*F_{A_1}\|_{C(X)}
\\
+ C_0(\sqrt{\tau} + \tau)\left(\|\Ric_g\|_{C(X)} + \|F_{A_1}\|_{C(X)} + 1\right)\|w\|_\fV
\\
+ C_0\sqrt{\tau}\left(\|w\|_\fV^2 + (1 + \sqrt{\tau})\|w\|_\fV^3\right),
\end{multline}
where the constant $C_0$ depends on the Riemannian metric $g$ on $X$ but is independent of the connection $A_1$ on $P$. This completes our verification of boundedness for the Yang-Mills nonlinearity \eqref{eq:Yang-Mills_heat_equation_nonlinearity_relative_rough_Laplacian_plus_one}.
\end{step}

\begin{step}[Lipschitz property the Yang-Mills nonlinearity in $\fW$]
Second, given $w_1, w_2 \in \fV$, we can estimate the expression $\sF(w_1) - \sF(w_2)$ defined by the nonlinearity \eqref{eq:Yang-Mills_heat_equation_nonlinearity_relative_rough_Laplacian_plus_one} to give
\begin{align*}
{}&\|\sF(w_1) - \sF(w_2)\|_{L^2(0,\tau; L^2(X))}
\\
&\quad \leq \|\Ric_g\times (w_1-w_2)\|_{L^2(0,\tau; L^2(X))} + \|(F_{A_1}-1)\times (w_1-w_2)\|_{L^2(0,\tau; L^2(X))}
\\
&\qquad + \|\nabla_{A_1}(w_1-w_2)\times w_1\|_{L^2(0,\tau; L^2(X))}
+ \|\nabla_{A_1}w_2\times (w_1-w_2)\|_{L^2(0,\tau; L^2(X))}
\\
&\qquad + \|(w_1-w_2)\times w_1\times w_1\|_{L^2(0,\tau; L^2(X))}
+ \|w_2\times (w_1-w_2)\times w_1\|_{L^2(0,\tau; L^2(X))}
\\
&\qquad + \|w_2\times w_2\times (w_1-w_2)\|_{L^2(0,\tau; L^2(X))}
\\
&\quad \leq c\left(\|\Ric_g\|_{L^2(0,\tau; L^2(X))} + \|F_{A_1}\|_{L^2(0,\tau; L^2(X))} + 1\right)\|w_1-w_2\|_{C([0,\tau] \times X)}
\\
&\qquad + c\|\nabla_{A_1}(w_1-w_2)\|_{L^2(0,\tau; L^2(X))} \|w_1\|_{C([0,\tau] \times X)}
\\
&\qquad + c\|\nabla_{A_1}w_2\|_{L^2(0,\tau; L^2(X))} \|w_1-w_2\|_{C([0,\tau] \times X)}
\\
&\qquad + c\|w_1-w_2\|_{L^2(0,\tau; L^2(X))}\|w_1\|_{C([0,\tau] \times X)}^2
\\
&\qquad + c\|w_1-w_2\|_{L^2(0,\tau; L^2(X))}\|w_1\|_{C([0,\tau] \times X)}
\|w_2\|_{C([0,\tau] \times X)}
\\
&\qquad + c\|w_1-w_2\|_{L^2(0,\tau; L^2(X))}\|w_2\|_{C([0,\tau] \times X)}^2,
\end{align*}
by the same reasoning as that which led to the estimates \eqref{eq:Quadratic_term_norm_L2time_and_L2space_dimension_4} and \eqref{eq:Cubic_term_norm_L2time_and_L2space_dimension_4} in $L^2(0,\tau; L^2(X))$ for the quadratic and cubic terms, respectively. The constant $c$ depends on the Riemannian metric $g$ on $X$ but is independent of the connection $A_1$ on $P$. Therefore,
\begin{align*}
{}&\|\sF(w_1) - \sF(w_2)\|_{L^2(0,\tau; L^2(X))}
\\
&\quad \leq c\sqrt{\tau}\left(\|\Ric_g\|_{L^2(X)} + \|F_{A_1}\|_{L^2(X)} + 1\right)\|w_1-w_2\|_{C([0,\tau] \times X)}
\\
&\qquad + c\sqrt{\tau}\|w_1-w_2\|_{C([0,\tau]; H^1_{A_1}(X))} \|w_1\|_{C([0,\tau] \times X)}
\\
&\qquad + c\sqrt{\tau}\|w_2\|_{C([0,\tau]; H^1_{A_1}(X))} \|w_1-w_2\|_{C([0,\tau] \times X)}
\\
&\qquad + c\sqrt{\tau}\|w_1-w_2\|_{C([0,\tau]; L^2(X))} \|w_1\|_{C([0,\tau] \times X)}^2
\\
&\qquad + c\sqrt{\tau}\|w_1-w_2\|_{C([0,\tau]; L^2(X))} \|w_1\|_{C([0,\tau] \times X)}
\|w_2\|_{C([0,\tau] \times X)}
\\
&\qquad + c\sqrt{\tau}\|w_1-w_2\|_{C([0,\tau]; L^2(X))} \|w_2\|_{C([0,\tau] \times X)}^2.
\end{align*}
Applying our definition \eqref{eq:Parabolic_Sobolev_norm_hybrid_standard_H2_and_critical_exponent_domain} of the Banach space norm on $\fV$ in the preceding inequality, we see that
\begin{multline}
\label{eq:Struwe_nonlinear_contraction_map_bound_L2time_and_L2_space_dimension_4}
\|\sF(w_1) - \sF(w_2)\|_{L^2(0,\tau; L^2(X))}
\\
\leq c\sqrt{\tau}\left(\|\Ric_g\|_{L^2(X)} + \|F_{A_1}\|_{L^2(X)} + 1\right)\|w_1-w_2\|_\fV
\\
+ c\sqrt{\tau}\left(\|w_1\|_\fV + \|w_2\|_\fV + \|w_1\|_\fV^2 + \|w_2\|_\fV^2\right)\|w_1-w_2\|_\fV,
\end{multline}
where the constant $c$ depends on the Riemannian metric $g$ on $X$ but is independent of the connection $A_1$ on $P$.

Similarly, by the reasoning which led to the estimates \eqref{eq:Quadratic_term_norm_Llozenge_time_and_space_dimension_4} and \eqref{eq:Cubic_term_norm_Llozenge_time_and_space_dimension_4} in $L^\lozenge((0,\tau)\times X; \Lambda^1\otimes\ad P)$ for the quadratic and cubic terms, respectively, the derivations of the inequalities \eqref{eq:Struwe_nonlinear_map_bound_Llozenge_time_and_space_dimension_4} and \eqref{eq:Struwe_nonlinear_contraction_map_bound_L2time_and_L2_space_dimension_4} lead, \emph{mutatis mutandis}, to
\begin{multline}
\label{eq:Struwe_nonlinear_contraction_map_bound_Llozenge_time_and_space_dimension_4}
\|\sF(w_1) - \sF(w_2)\|_{L^\lozenge((0,\tau)\times X)} \leq c\tau\left(\|\Ric_g\|_{C(X)} + \|F_{A_1}\|_{C(X)} + 1\right)\|w_1-w_2\|_\fV
\\
+ c\sqrt{\tau}\left(\|w_1\|_\fV + \|w_2\|_\fV + \sqrt{\tau}\|w_1\|_\fV^2 + \sqrt{\tau}\|w_2\|_\fV^2\right)\|w_1-w_2\|_\fV,
\end{multline}
where the constant $c$ depends on the Riemannian metric $g$ on $X$ but is independent of the connection $A_1$ on $P$. Combining \eqref{eq:Struwe_nonlinear_contraction_map_bound_L2time_and_L2_space_dimension_4} and \eqref{eq:Struwe_nonlinear_contraction_map_bound_Llozenge_time_and_space_dimension_4} yields
\begin{equation}
\begin{aligned}
\label{eq:Struwe_nonlinear_contraction_map_bound_L2time_and_L2_space_Llozenge_time_and_space_dimension_4}
{}&\|\sF(w_1) - \sF(w_2)\|_\fW
\\
&\quad \leq C_1(\sqrt{\tau} + \tau)\left(\|\Ric_g\|_{C(X)} + \|F_{A_1}\|_{C(X)} + 1\right)\|w_1-w_2\|_\fV
\\
&\qquad + C_1\sqrt{\tau}\left(\|w_1\|_\fV + \|w_2\|_\fV + (1 + \sqrt{\tau})\|w_1\|_\fV^2
+ (1 + \sqrt{\tau})\|w_2\|_\fV^2\right)\|w_1-w_2\|_\fV,
\end{aligned}
\end{equation}
where the constant $C_1$ depends on the Riemannian metric $g$ on $X$ but is independent of the connection $A_1$ on $P$. This completes our verification of the Lipschitz property for the Yang-Mills nonlinearity \eqref{eq:Yang-Mills_heat_equation_nonlinearity_relative_rough_Laplacian_plus_one}.
\end{step}

\begin{step}[Contraction mapping]
We now proceed by analogy with the proofs of \cite[Theorem 46.1 and Lemma 47.1]{Sell_You_2002}. Given
$$
w \in \{v \in \fV: \|v\|_\fV \leq R \},
$$
with $R$ a positive constant to be determined, we wish to define $a := \Phi(w) \in \fV$ by solving the linear heat equation, that is,
$$
L_{A_1}a \equiv \frac{\partial a}{\partial t} + (\nabla_{A_1}^*\nabla_{A_1} + 1)a =  f \quad\hbox{on } (0,\tau), \quad a(0) = a_0,
$$
with source function $f = \sF(w)$. By the inequality \eqref{eq:Struwe_nonlinear_map_bound_L2time_and_L2_space_and_Llozenge_time_and_space_dimension_4}, we have $\sF(w) \in \fW$, for any choice of $\tau \in (0,\infty)$. Because $C^\infty([0,\tau]\times X; \Lambda^1\otimes \ad P)$ is dense in $\fW$, by the definition of this Banach space by completion with respect to the norm \eqref{eq:Parabolic_Sobolev_norm_hybrid_standard_L2_and_critical_exponent_range} and because of the \apriori estimate \eqref{eq:Critical_exponent_plus_standard_H2_parabolic_linear_apriori_estimate_streamlined}, the preceding linear heat equation has a unique strong solution $a \in \fV$, given $f \equiv \sF(w)$ (for example, by Corollary \ref{cor:Sell_You_42-13} with $\alpha = 1$ and $\calV = H^1_{A_1}(X; \Lambda^1\ad P)$), which also obeys,
\begin{align*}
\|\Phi(w)\|_\fV = \|a\|_\fV
&\leq 8\|L_Aa\|_\fW +  4\|a_0\|_{H^1_{A_1}(X)} + \sup_{(t, x) \in (0,\tau)\times X} 4\|a_0\|_{L^2(X, K_{t,x})}
\\
&= 8\|\sF(w)\|_\fW + 4\|a_0\|_{H^1_{A_1}(X)} + \sup_{(t, x) \in (0,\tau)\times X} 4\|a_0\|_{L^2(X, K_{t,x})},
\end{align*}
Thus, by \eqref{eq:Struwe_nonlinear_map_bound_L2time_and_L2_space_and_Llozenge_time_and_space_dimension_4},
\begin{equation}
\label{eq:Struwe_contraction_map_solution_bound_sharp}
\begin{aligned}
\|\Phi(w)\|_\fV &\leq 8(\sqrt{\tau} + \tau)\|d_{A_1}^*F_{A_1}\|_{C(X)}
\\
&\quad + 8C_0(\sqrt{\tau} + \tau)\left(\|\Ric_g\|_{C(X)} + \|F_{A_1}\|_{C(X)} + 1\right)\|w\|_\fV
\\
&\quad + 8C_0\sqrt{\tau}\left(\|w\|_\fV^2 + (1 + \sqrt{\tau})\|w\|_\fV^3\right)
\\
&\quad + 4\|a_0\|_{H^1_{A_1}(X)} + \sup_{(t, x) \in (0,\tau)\times X} 4\|a_0\|_{L^2(X, K_{t,x})}.
\end{aligned}
\end{equation}
Recalling the bound on the initial data in terms of $b$ given by \eqref{eq:Norm_ao_H1_and_sup_tx_L2_a0_dKtx_leq_b}, if we set
\begin{equation}
\label{eq:Struwe_section_4-3_Yang-Mills_heat_equation_hybrid_critical-exponent_parabolic_Sobolev_space_radius}
R := 8b,
\end{equation}
then we obtain $\|\Phi(w)\|_\fV \le R$ provided $\|w\|_\fV \le R$ and, supposing without loss of generality that $0 < \tau \leq 1$ and so $\tau \leq \sqrt{\tau}$, the constant $\tau \in (0, 1]$ is chosen small enough that
\begin{multline}
\label{eq:Struwe_section_4-3_minimallifetime_Yang-Mills_heat_equation_hybrid_critical-exponent_parabolic_Sobolev_space_prelim1}
2\sqrt{\tau}\|d_{A_1}^*F_{A_1}\|_{C(X)}
\\
+ 2C_0\sqrt{\tau} \left(\|\Ric_g\|_{C(X)} + \|F_{A_1}\|_{C(X)} + 1\right)R
+ C_0\sqrt{\tau}\left(R^2 + 2R^3\right)
\leq b.
\end{multline}
Thus, for such $\tau > 0$, we have $\|\Phi(w)\|_\fV \leq R$ for all $w \in \fV$ obeying $\|w\|_\fV \leq R$.

Similarly, the \apriori estimates \eqref{eq:Critical_exponent_plus_standard_H2_parabolic_linear_apriori_estimate_streamlined} and \eqref{eq:Struwe_nonlinear_contraction_map_bound_L2time_and_L2_space_Llozenge_time_and_space_dimension_4} yield a bound for $a_1-a_2 = \Phi(w_1)-\Phi(w_2)$,
\begin{multline}
\label{eq:Struwe_section_4-3_continuity_estimate_Yang-Mills_heat_equation_hybrid_critical-exponent_parabolic_Sobolev_space_prelim}
\|a_1 - a_2\|_\fV \leq 8C_1(\sqrt{\tau} + \tau)\left(\|\Ric_g\|_{C(X)} + \|F_{A_1}\|_{C(X)} + 1\right)\|w_1-w_2\|_\fV
\\
+ 8C_1\sqrt{\tau}\left(\|w_1\|_\fV + \|w_2\|_\fV + (1 + \sqrt{\tau})\|w_1\|_\fV^2
+ (1 + \sqrt{\tau})\|w_2\|_\fV^2\right)\|w_1-w_2\|_\fV.
\end{multline}
Therefore, $\Phi$ is a contraction map on $\{w \in \fV:\|w\|_\fV \leq R\}$ with contraction coefficient less than or equal to $1/2$ provided $\tau \in (0, 1]$ is also chosen small enough that
\begin{equation}
\label{eq:Struwe_section_4-3_minimallifetime_Yang-Mills_heat_equation_hybrid_critical-exponent_parabolic_Sobolev_space_prelim2}
16C_1\sqrt{\tau}\left(\|\Ric_g\|_{C(X)} + \|F_{A_1}\|_{C(X)} + 1\right)
+
16C_1\sqrt{\tau}\left( R + 2R^2 \right) \leq \frac{1}{2}.
\end{equation}
An application of the contraction mapping principle now completes the proof of existence and uniqueness in Theorem \ref{thm:Struwe_section_4-3_local_wellposedness_Yang-Mills_heat_equation_hybrid_critical-exponent_parabolic_Sobolev_space}, completing this step.
\end{step}

The proof of the last step also provides the \apriori estimate for the solution and lower bound for $\tau$. Indeed, the inequality \eqref{eq:Struwe_contraction_map_solution_bound_sharp} and the fact that $a = \Phi(a)$ yields
\begin{align*}
\|a\|_\fV &\leq 8(\sqrt{\tau} + \tau) \|d_{A_1}^*F_{A_1}\|_{C(X)}
\\
&\quad + 8C_0(\sqrt{\tau} + \tau) \left(\|\Ric_g\|_{C(X)} + \|F_{A_1}\|_{C(X)} + 1\right) \|a\|_\fV
\\
&\quad + 8C_0\sqrt{\tau}\left(\|a\|_\fV^2 + (1 + \sqrt{\tau})\|a\|_\fV^3\right)
\\
&\quad + 4\|a_0\|_{H^1_{A_1}(X)} + \sup_{(t, x) \in (0,\tau)\times X} 4 \|a_0\|_{L^2(X, K_{t,x})}.
\end{align*}
We may further choose $\tau \in (0,1]$ small enough that it also obeys
\begin{equation}
\label{eq:Struwe_section_4-3_minimallifetime_Yang-Mills_heat_equation_hybrid_critical-exponent_parabolic_Sobolev_space_prelim3}
16C_0\sqrt{\tau}\left(\|\Ric_g\|_{C(X)} + \|F_{A_1}\|_{C(X)} + 1\right)
+
8C_0\sqrt{\tau}\left( R + 2R^2 \right) \leq \frac{1}{2}.
\end{equation}
Therefore,
\begin{align*}
\|a\|_\fV \leq 8(\sqrt{\tau} + \tau) \|d_{A_1}^*F_{A_1}\|_{C(X)}
+ 4\|a_0\|_{H^1_{A_1}(X)} + \sup_{(t, x) \in (0,\tau)\times X} 4 \|a_0\|_{L^2(X, K_{t,x})}
+ \frac{1}{2}\|a\|_\fV,
\end{align*}
and this yields the desired \apriori estimate
\eqref{eq:Struwe_section_4-3_apriori_estimate_Yang-Mills_heat_equation_hybrid_critical-exponent_parabolic_Sobolev_space}.

If $a_1, a_2$ are the unique solutions corresponding to initial data $a_1(0)$ and $a_2(0$, respectively, then the derivation of the contraction mapping inequality \eqref{eq:Struwe_section_4-3_continuity_estimate_Yang-Mills_heat_equation_hybrid_critical-exponent_parabolic_Sobolev_space_prelim} and the contraction mapping coefficient condition \eqref{eq:Struwe_section_4-3_minimallifetime_Yang-Mills_heat_equation_hybrid_critical-exponent_parabolic_Sobolev_space_prelim2} and the fact that $a_i = \Phi(a_i)$ for $i = 1,2$ yield
\begin{multline*}
\|a_1 - a_2\|_\fV \leq \frac{1}{2}\|a_1-a_2\|_\fV
\\
+ 4\|a_1(0) - a_2(0)\|_{H^1_{A_1}(X)} + \sup_{(t, x) \in (0,\tau)\times X} 4 \|a_1(0) - a_2(0)\|_{L^2(X, K_{t,x})},
\end{multline*}
and this yields the desired continuity estimate
\eqref{eq:Struwe_section_4-3_continuity_estimate_Yang-Mills_heat_equation_hybrid_critical-exponent_parabolic_Sobolev_space}.

To obtain a positive lower bound for the lifetime of the solution, it suffices to notice that we can choose $\tau \in (0, 1]$ to be the largest constant obeying
\eqref{eq:Struwe_section_4-3_minimallifetime_Yang-Mills_heat_equation_hybrid_critical-exponent_parabolic_Sobolev_space_prelim1},
\eqref{eq:Struwe_section_4-3_minimallifetime_Yang-Mills_heat_equation_hybrid_critical-exponent_parabolic_Sobolev_space_prelim2}, and \eqref{eq:Struwe_section_4-3_minimallifetime_Yang-Mills_heat_equation_hybrid_critical-exponent_parabolic_Sobolev_space_prelim3}.
Combining the preceding observation with our definition
\eqref{eq:Struwe_section_4-3_Yang-Mills_heat_equation_hybrid_critical-exponent_parabolic_Sobolev_space_radius} of $R$ yields the desired minimal value of $\tau$ in
\eqref{eq:Struwe_section_4-3_minimallifetime_Yang-Mills_heat_equation_hybrid_critical-exponent_parabolic_Sobolev_space}.
\end{proof}

We can now apply a bootstrapping argument to obtain higher-order regularity for the solution $a(t)$ obtained in Theorem \ref{thm:Struwe_section_4-3_local_wellposedness_Yang-Mills_heat_equation_hybrid_critical-exponent_parabolic_Sobolev_space} when $t > 0$. This can be done for a manifold of any dimension $d \geq 2$, in various ways, but we shall defer the proof of higher-order regularity until the end of Section \ref{subsec:Struwe_page_137_contraction_mapping_initial_data_in_Ld}, where the desired smoothness follows more easily from the results of that section. The new regularity feature in Corollary
\ref{cor:Struwe_section_4-3_higher_order_regularity_Yang-Mills_heat_equation_hybrid_critical-exponent_parabolic_Sobolev_space}
is the fact that $C^\infty((0, T) \times X; \Lambda^1\otimes\ad P))$; the remaining regularity property of the solution is implicit in the definition of the Banach space $\fV$ and the \apriori estimate \eqref{eq:Critical_exponent_plus_standard_H2_parabolic_linear_apriori_estimate}.

\begin{cor}[Higher-order regularity for a solution to the Yang-Mills heat equation in a hybrid critical-exponent parabolic Sobolev space]
\label{cor:Struwe_section_4-3_higher_order_regularity_Yang-Mills_heat_equation_hybrid_critical-exponent_parabolic_Sobolev_space}
Let $G$ be a compact Lie group and $P$ a principal $G$-bundle over a closed, connected, orientable, smooth manifold, $X$, of dimension $d \geq 2$ and Riemannian metric, $g$. Let $A_1$ be a reference connection of class $C^\infty$ on $P$. If $a_0 \in H^1_{A_1}(X; \Lambda^1\otimes\ad P) \cap C(X; \Lambda^1\otimes\ad P)$ and $a(t)$ is a strong solution in $\fV$ to the Yang-Mills heat equation \eqref{eq:Yang-Mills_heat_equation_as_perturbation_rough_Laplacian_plus_one_heat_equation} on an interval $(0, T)$ with initial data $a(0) = a_0$, for some $T > 0$, where $\fV$ is the Banach space completion of $C^\infty([0,T]\times X; \Lambda^1\otimes \ad P))$ with respect to the norm \eqref{eq:Parabolic_Sobolev_norm_hybrid_standard_H2_and_critical_exponent_domain}, then
\begin{multline*}
a \in C([0, T); C(X; \Lambda^1\otimes\ad P)) \cap C([0, T); H^1_{A_1}(X; \Lambda^1\otimes\ad P))
\\
\cap L^2(0, T; H^2_{A_1}(X; \Lambda^1\otimes\ad P)) \cap H^1(0, T; L^2(X; \Lambda^1\otimes\ad P)) \cap
\\
\cap C^\infty((0, T) \times X; \Lambda^1\otimes\ad P)),
\end{multline*}
and $a(t)$ is a classical solution on $(0, T)$. Furthermore, if $a_0 \in C^\infty(X; \Lambda^1\otimes\ad P)$, then
$$
a \in C^\infty([0, T) \times X; \Lambda^1\otimes\ad P)).
$$
\end{cor}

\subsection{Local well-posedness for the Yang-Mills heat equation in pure critical-exponent parabolic Sobolev spaces}
\label{subsec:Struwe_page_137_contraction_mapping_pure_critical-exponent_parabolic_Sobolev_space}
Although we introduced our hybrid critical-exponent parabolic Sobolev spaces as a way of repairing the argument of Struwe in \cite[Section 4.3]{Struwe_1994}, it is easier to use pure critical-exponent parabolic Sobolev spaces to obtain what appears to be an optimal result in terms of keeping the regularity requirement on the initial data, $a_0$, to the minimum needed to guarantee local well-posedness for the Yang-Mills heat equation \eqref{eq:Yang-Mills_heat_equation_as_perturbation_rough_Laplacian_plus_one_heat_equation}. Indeed, this is the approach taken by Taubes when solving elliptic quasilinear partial differential systems arising in the study of anti-self-dual conformal metrics or anti-self-dual connections \cite{TauGluing}.

Given a solution in a pure critical-exponent parabolic Sobolev space, we will be able to appeal to results on embeddings for pure critical-exponent parabolic Sobolev spaces into standard Sobolev spaces which become available when $t \geq t_0 > 0$ and obtain regularity for our solution in standard Sobolev spaces when $t > 0$. For example, the definitions of the $L^\lozenge$ and $L^{2\lozenge}$ norms in \eqref{eq:Critical-exponent_parabolic_Sobolev_spaces} for $u \in C^\infty((0, T)\times X; E)$ and the norm on $u(0,\cdot) \in C^\infty(X; E)$ in \eqref{eq:Feehan_5-3_heat_operator} and the fact that by Benjamini, Chavel, and Feldman \cite[Theorem 2]{Benjamini_Chavel_Feldman_1996}\footnote{Their result is stated for non-compact Riemannian manifolds of dimension $d \geq 2$ with bounded geometry. Related results for compact Riemannian manifolds may be extracted from Berger, Gauduchon, Mazet \cite{Berger_Gauduchon_Mazet_1971}, Chavel \cite{Chavel}, Davies \cite{Davies_1989}, and Grigor'yan \cite{Grigoryan_2009}, although not in the simple Gaussian form we need as far as the author can determine. For $t \in [0, T]$ with $T < \infty$, the pointwise Gaussian lower bound \eqref{eq:Heat_kernel_Gaussian_pointwise_lower_bound} follows from \cite[Equation (6.45)]{Chavel}, with constants $c_0 = c_0(g,T)$ and $\bar c_1 = 1/4$.},
\begin{equation}
\label{eq:Heat_kernel_Gaussian_pointwise_lower_bound}
c_0^{-1}t^{-d/2} e^{-\bar c_1\dist_g^2(x,y)/t} \leq K(t, x, y), \quad \forall\, (t, x, y) \in [t_0, \infty) \times X \times X,
\end{equation}
and, in particular,
$$
c_0^{-1}T^{-d/2} e^{-c_1/t_0} \leq K(t, x, y), \quad \forall\, (t, x, y) \in [t_0, T] \times X \times X,
$$
for positive constants $c_0, c_1$, and $\bar c_1 = c_1/\diam_g^2(X)$ depending at most on the Riemannian metric, $g$, on $X$, yield the inequalities,
\begin{subequations}
\label{eq:Sobolev_embedding_critical-exponent_parabolic_into_standard}
\begin{align}
\label{eq:Sobolev_embedding_Llozenge_into_L1_time_and_space}
\|u\|_{L^1((t_0, T)\times X; E)} &\leq c_0T^{d/2} e^{c_1/t_0} \|u\|_{L^\lozenge((0, T)\times X; E)},
\\
\label{eq:Sobolev_embedding_L2lozenge_into_L2_time_and_space}
\|u\|_{L^2((t_0, T)\times X; E)} &\leq c_0T^{d/4} e^{c_1/2t_0} \|u\|_{L^{2\lozenge}((0, T)\times X; E)},
\\
\label{eq:Sobolev_embedding_sup_time_and_space_L2space_dKtx_into_L2space}
\|u_0\|_{L^2(X; E))} &\leq c_0T^{d/4} e^{c_1/2t_0} \sup_{(t, x) \in (t_0,T)\times X}\|u_0\|_{L^2(X, K_{t,x}; E)},
\end{align}
\end{subequations}
analogous to the elliptic embedding results in Lemma \ref{lem:Feehan_4-2}. By combining the Sobolev embeddings \eqref{eq:Sobolev_embedding_critical-exponent_parabolic_into_standard} with the \apriori estimate \eqref{eq:Feehan_5-3_heat_operator} we see that
\begin{multline}
\label{eq:Feehan_5-3_heat_operator_time_geq_t0}
\|\nabla_Au\|_{L^2((t_0, T)\times X; E)} + \|u\|_{C([t_0, T]\times X; E)}
\\
\leq
c_0T^{d/2} e^{c_1/t_0} \|L_Au\|_{L^\lozenge((0, T)\times X; E)}
+ c_0T^{d/4} e^{c_1/2t_0} \sup_{(t, x) \in (0,T)\times X}\|u(0,\cdot)\|_{L^2(X, K_{t,x}; E))}.
\end{multline}
We now replace our previous choices for $\fW$ and $\fV$ defined by the hybrid critical-exponent parabolic Sobolev space norms in \eqref{eq:Parabolic_Sobolev_norms_hybrid_standard_and_critical_exponent} with Banach spaces defined by completing $C^\infty([0, T]\times X; E)$ with respect to the \emph{pure critical-exponent parabolic Sobolev space} norms,
\begin{subequations}
\label{eq:Parabolic_Sobolev_norms_pure_critical_exponent}
\begin{align}
\label{eq:Parabolic_Sobolev_norm_pure_critical_exponent_range}
\|w\|_\fW &= \|w\|_{L^\lozenge((0, T)\times X; E)}
\\
\label{eq:Parabolic_Sobolev_norm_pure_critical_exponent_domain}
\|v\|_\fV &= \|u\|_{C([0, T]\times X)} + \|\nabla _Au\|_{L^{2\lozenge}((0, T)\times X; E)}
+ \|L_Au\|_{L^\lozenge((0, T)\times X; E)},
\end{align}
\end{subequations}
where we recall from \eqref{eq:Augmented_heat_operator_on_sections_vectorbundle_over_manifold} that $L_A \equiv \partial_t + \nabla_A^*\nabla_A + 1$. We can use the \apriori estimate \eqref{eq:Feehan_5-3_heat_operator_time_geq_t0} to convert a solution in the pure critical-exponent parabolic Sobolev space, $\fV$ in \eqref{eq:Parabolic_Sobolev_norm_pure_critical_exponent_domain}, to the Yang-Mills heat equation \eqref{eq:Yang-Mills_heat_equation_as_perturbation_rough_Laplacian_plus_one_heat_equation} on $(0,T)\times X$ into a classical solution on $(t_0,T)\times X$ for any $t_0 > 0$.

Minor changes in our proof of Theorem \ref{thm:Struwe_section_4-3_local_wellposedness_Yang-Mills_heat_equation_hybrid_critical-exponent_parabolic_Sobolev_space}, but now using the choices in \eqref{eq:Parabolic_Sobolev_norms_pure_critical_exponent} for $\fW$ and $\fV$ instead of those in \eqref{eq:Parabolic_Sobolev_norms_hybrid_standard_and_critical_exponent}, yield the proof of

\begin{thm}[Local well-posedness for the Yang-Mills heat equation in pure critical-exponent parabolic Sobolev spaces over a closed manifold]
\label{thm:Struwe_section_4-3_local_wellposedness_Yang-Mills_heat_equation_pure_critical-exponent_parabolic_Sobolev_space}
Let $G$ be a compact Lie group and $P$ a principal $G$-bundle over a closed, connected, orientable, smooth manifold, $X$, of dimension $d \geq 2$ and Riemannian metric, $g$. Then there is a positive constant $C = C(d, g)$ and, if $b$ is a positive constant, there is a positive constant,
$$
\tau = \tau(b, A_1, d, g) \in (0, 1],
$$
with the following significance. If $A_1$ is a reference connection of class $C^\infty$ on $P$ and
$$
a_0 \in C(X; \Lambda^1\otimes\ad P)
$$
obeys
\begin{equation}
\label{eq:Norm_ao_sup_tx_L2_a0_dKtx_leq_b}
\sup_{(t,x) \in (0,\tau)\times X} \|a_0\|_{L^2(X, K_{t,x})} \leq b,
\end{equation}
then there is a unique strong solution $a(t)$ in $\fV$ to the Yang-Mills heat equation
\eqref{eq:Yang-Mills_heat_equation_as_perturbation_rough_Laplacian_plus_one_heat_equation} on the interval $(0, \tau)$ with initial data $a(0) = a_0$, where $\fV$ is the Banach space completion of $C^\infty([0,\tau]\times X; \Lambda^1\otimes \ad P))$ with respect to the norm \eqref{eq:Parabolic_Sobolev_norm_pure_critical_exponent_domain}. Moreover, the solution $a(t)$ obeys the \apriori estimate,
\begin{equation}
\label{eq:Apriori_estimate_Yang-Mills_heat_equation_pure_critical-exponent_parabolic_Sobolev_space}
\|a\|_\fV \leq 16\sqrt{\tau} \|d_{A_1}^*F_{A_1}\|_{C(X)} + \sup_{(t, x) \in (0,\tau)\times X} 8\|a_0\|_{L^2(X, K_{t,x})},
\end{equation}
and the minimal lifetime, $\tau$, is given by
\begin{equation}
\label{eq:Minimallifetime_Yang-Mills_heat_equation_pure_critical-exponent_parabolic_Sobolev_space}
\tau = Cb\left(\|d_{A_1}^*F_{A_1}\|_{C(X)} + \left(\|\Ric_g\|_{C(X)} + \|F_{A_1}\|_{C(X)} + 1\right)(1 + b^3)\right)^{-2}.
\end{equation}
If $a_1, a_2$ are the unique solutions corresponding to initial data $a_1(0)$ and $a_2(0)$, respectively, then
\begin{equation}
\label{eq:Struwe_section_4-3_continuity_estimate_Yang-Mills_heat_equation_pure_critical-exponent_parabolic_Sobolev_space}
\|a_1 - a_2\|_\fV \leq \sup_{(t, x) \in (0,T)\times X} 8 \|a_1(0) - a_2(0)\|_{L^2(X, K_{t,x})}.
\end{equation}
\end{thm}

While Theorem \ref{thm:Struwe_section_4-3_local_wellposedness_Yang-Mills_heat_equation_pure_critical-exponent_parabolic_Sobolev_space} has greatest value when $d = 4$, it is worth noting that, like Theorem \ref{thm:Struwe_section_4-3_local_wellposedness_Yang-Mills_heat_equation_hybrid_critical-exponent_parabolic_Sobolev_space},  it is valid for a manifold $X$ of any dimension $d \geq 2$, since we do not need to appeal to the Sobolev Embedding Theorem \cite[Theorem 4.12]{AdamsFournier} explicitly or implicitly in its proof.

We shall again defer the proof of the higher-order regularity property of the solution, $a(t)$, produced by
Theorem \ref{thm:Struwe_section_4-3_local_wellposedness_Yang-Mills_heat_equation_pure_critical-exponent_parabolic_Sobolev_space}
until the end of Section \ref{subsec:Struwe_page_137_contraction_mapping_initial_data_in_Ld}, where the desired smoothness for $t > 0$ follows more easily from the results of that section.

\begin{cor}[Higher-order regularity for a solution to the Yang-Mills heat equation in a pure critical-exponent parabolic Sobolev space]
\label{cor:Struwe_section_4-3_higher_order_regularity_Yang-Mills_heat_equation_pure_critical-exponent_parabolic_Sobolev_space}
Let $G$ be a compact Lie group and $P$ a principal $G$-bundle over a closed, connected, orientable, smooth manifold, $X$, of dimension $d \geq 2$ and Riemannian metric, $g$. Let $A_1$ be a reference connection of class $C^\infty$ on $P$. If $a_0 \in C(X; \Lambda^1\otimes\ad P)$ and $a(t)$ is a strong solution in $\fV$ to the Yang-Mills heat equation
\eqref{eq:Yang-Mills_heat_equation_as_perturbation_rough_Laplacian_plus_one_heat_equation} on an interval $(0, T)$ with initial data $a(0) = a_0$, for some $T > 0$, where $\fV$ is the Banach space completion of $C^\infty([0,T]\times X; \Lambda^1\otimes \ad P))$ with respect to the norm \eqref{eq:Parabolic_Sobolev_norm_pure_critical_exponent_domain}, then
$$
a \in C([0, T); C(X; \Lambda^1\otimes\ad P)) \cap C^\infty((0, T) \times X; \Lambda^1\otimes\ad P)),
$$
and $a(t)$ is a classical solution on $(0, T)$. Furthermore, if $a_0 \in C^\infty(X; \Lambda^1\otimes\ad P)$, then
$$
a \in C^\infty([0, T) \times X; \Lambda^1\otimes\ad P)).
$$
\end{cor}

\subsection{Local existence and uniqueness of solutions to the Yang-Mills heat equation with initial data in $L^d$}
\label{subsec:Struwe_page_137_contraction_mapping_initial_data_in_Ld}
Our abstract Theorem \ref{thm:Kozono_Maeda_Naito_lemma_3-4_plus_uniqueness} gives another approach to local existence and uniqueness of a solution, $a(t)$ for $t\in [0, \tau)$, to the Yang-Mills heat equation over a closed manifold of dimension $d \geq 2$ and initial data $a_0$ in $L^d(X; \Lambda^1\otimes\ad P)$. However, as we explained in Remark \ref{rmk:Kozono_Maeda_Naito_lemma_3-4_plus_uniqueness_initial_data}, that result does not yield local well-posedness since the solution need not depend continuously on initial data, $a_0$, with only that regularity.

\begin{thm}[Local existence and uniqueness of mild solutions to the Yang-Mills heat equation with initial data in $L^d$]
\label{thm:Struwe_section_4-3_local_existence_uniqueness_Yang-Mills_heat_equation_initial_data_in_Ld}
Let $G$ be a compact Lie group and $P$ a principal $G$-bundle over a closed, connected, orientable, smooth manifold, $X$, of dimension $d \geq 2$ and Riemannian metric, $g$. If $A_1$ is a reference connection of class $C^\infty$ on $P$ and $a_0 \in L^d(X; \Lambda^1\otimes\ad P)$ and $\delta_0 \in [\frac{1}{2}, \frac{3}{4})$, then there are positive constants,
$$
\tau = \tau\left(a_0, A_1, d, g, \delta_0\right)
\quad\hbox{and}\quad
C_0 = C_0(A_1, d, g, \delta_0),
$$
with the following significance. There is a unique mild solution,
\begin{equation}
\label{eq:Local_existence_uniqueness_Yang-Mills_heat_equation_initial_data_in_Ld}
a \in C([0, \tau]; L^d(X; \Lambda^1\otimes\ad P) \cap C((0, \tau]; W_{A_1}^{2\delta_0, d}(X; \Lambda^1\otimes\ad P)),
\end{equation}
to the Yang-Mills heat equation \eqref{eq:Yang-Mills_heat_equation_as_perturbation_rough_Laplacian_plus_one_heat_equation} on the interval $[0, \tau)$ with initial data $a(0) = a_0$.
Moreover, for every $\delta \in [0, \delta_0]$, the solution $u$ obeys the \apriori estimate,
\begin{equation}
\label{eq:Apriori_estimate_Yang-Mills_heat_equation_initial_data_in_Ld}
\sup_{t\in (0, \tau)}t^\delta\|a(t)\|_{W_{A_1}^{2\delta, d}(X)} \leq C_0.
\end{equation}
\end{thm}

\begin{rmk}[On the choice of $\cW$ in Theorem \ref{thm:Kozono_Maeda_Naito_lemma_3-4_plus_uniqueness} in the proof of Theorem
\ref{thm:Struwe_section_4-3_local_existence_uniqueness_Yang-Mills_heat_equation_initial_data_in_Ld}]
\label{rmk:Struwe_section_4-3_local_existence_uniqueness_Yang-Mills_heat_equation_initial_data_in_Ld}
Because of the Sobolev embedding $W^{1, d/2}(X) \hookrightarrow L^d(X)$ \cite[Theorem 4.12]{AdamsFournier}, one might be tempted to choose $\cW = W_{A_1}^{1, d/2}(X; \Lambda^1\otimes\ad P)$ in our application proof of Theorem \ref{thm:Kozono_Maeda_Naito_lemma_3-4_plus_uniqueness} to the proof of Theorem
\ref{thm:Struwe_section_4-3_local_existence_uniqueness_Yang-Mills_heat_equation_initial_data_in_Ld}. Unfortunately, as is clear from the Sobolev multiplication result \cite[Lemma 3.2]{Kozono_Maeda_Naito_1995}, one cannot expect the inequalities \eqref{eq:Kozono_Maeda_Naito_3-7} and \eqref{eq:Kozono_Maeda_Naito_3-8} for the nonlinearity $\cF$ assumed in the hypotheses of Theorem \ref{thm:Kozono_Maeda_Naito_lemma_3-4_plus_uniqueness} to hold when choosing $\cW = W_{A_1}^{1, d/2}(X; \Lambda^1\otimes\ad P)$ and $\cF$ to be the Yang-Mills heat equation nonlinearity \eqref{eq:Yang-Mills_heat_equation_nonlinearity_relative_rough_Laplacian_plus_one}.
\end{rmk}

Theorem \ref{thm:Struwe_section_4-3_local_existence_uniqueness_Yang-Mills_heat_equation_initial_data_in_Ld} will follow almost immediately from Theorem \ref{thm:Kozono_Maeda_Naito_lemma_3-4_plus_uniqueness} and the following Sobolev multiplication lemma, due to Kozono, Maeda, and Naito, but whose proof we include for completeness and the fact that (because of the critical dimension) it does not follow from standard Sobolev multiplication results, such as those in \cite[pp. 95--96]{FU}. Let $\cA = \sA_d$ be the realization of the augmented connection Laplace operator, $\sA = \nabla_{A_1}^*\nabla_{A_1} + 1$, on $\sD(\sA_d) \subset L^d(X; \Lambda^1\otimes\ad P)$, and recall that $\cA$ satisfies Hypothesis \ref{hyp:Sell_You_4_standing_hypothesis_A} by the discussion in Section \ref{subsec:Standing_hypothesis_A_and_the_augmented_connection_Laplacian}.

\begin{lem}[Estimates for the quadratic and cubic terms in the Yang-Mills heat equation nonlinearity]
\label{lem:Kozono_Maeda_Naito_3-3}
\cite[Lemma 3.3]{Kozono_Maeda_Naito_1995}
Let $G$ be a compact Lie group and $P$ a principal $G$-bundle over a closed, connected, orientable, Riemannian, smooth manifold, $X$, of dimension $d \geq 2$. If $A_1$ is a connection of class $C^\infty$ on $P$, then there is a positive constant, $C$, with the following significance. If $a_, a_1, a_2 \in W^{1, d}(X; \Lambda^1\otimes\ad P)$, then
\begin{subequations}
\label{eq:Kozono_Maeda_Naito_3-7_Yang-Mills}
\begin{align}
\label{eq:Kozono_Maeda_Naito_3-7_Yang-Mills_quadratic}
\|\cA^{-\frac{1}{4}}(\nabla_{A_1}a \times a)\|_{L^d(X)} &\leq C\|\cA^{\frac{1}{2}} a\|_{L^d(X)} \|\cA^{\frac{1}{4}}a\|_{L^d(X)},
\\
\label{eq:Kozono_Maeda_Naito_3-7_Yang-Mills_cubic}
\|\cA^{-\frac{1}{4}}(a\times a \times a)\|_{L^d(X)} &\leq C\|\cA^{\frac{1}{4}}a\|_{L^d(X)}^3,
\end{align}
\end{subequations}
and
\begin{subequations}
\label{eq:Kozono_Maeda_Naito_3-8_Yang-Mills}
\begin{multline}
\label{eq:Kozono_Maeda_Naito_3-8_Yang-Mills_quadratic}
\|\cA^{-\frac{1}{4}}(\nabla_{A_1}a_1 \times a_1 - \nabla_{A_1}a_2 \times a_2)\|_{L^d(X)}
\\
\leq C\|\cA^{\frac{1}{4}} a_1\|_{L^d(X)}\|\cA^{\frac{1}{2}}(a_1-a_2)\|_{L^d(X)}
+ C\|\cA^{\frac{1}{2}}a_2\|_{L^d(X)}\|\cA^{\frac{1}{4}}(a_1-a_2)\|_{L^d(X)},
\end{multline}
\begin{multline}
\label{eq:Kozono_Maeda_Naito_3-8_Yang-Mills_cubic}
\|\cA^{-\frac{1}{4}}(a_1\times a_1 \times a_1 - a_2\times a_2 \times a_2)\|_{L^d(X)}
\\
\leq C\left(\|\cA^{\frac{1}{4}} a_1\|_{L^d(X)}^2 + \|\cA^{\frac{1}{4}}a_2\|_{L^d(X)}^2 \right) \|\cA^{\frac{1}{4}}(a_1-a_2)\|_{L^d(X)}.
\end{multline}
\end{subequations}
\end{lem}

\begin{proof}
By the reverse H\"older inequality (duality) \cite[Lemma 2.7]{AdamsFournier} and letting $d' \in (1, 2]$ denote the dual exponent to $d \geq 2$, so $1 = 1/d + 1/d'$, we have (writing $L^{d'}(X)$ for $L^{d'}(X; \Lambda^1\otimes\ad P)$ in subscripts for brevity)
\begin{align*}
\|\cA^{-\frac{1}{4}}(\nabla_{A_1}a \times a)\|_{L^d(X)}
&= \sup_{\begin{subarray}{c} b \in L^{d'}(X), \\ \|b\|_{L^{d'}(X)} \leq 1 \end{subarray}} (\cA^{-\frac{1}{4}}(\nabla_{A_1}a \times a), b)_{L^2(X)}
\\
&= \sup_{\|b\|_{L^{d'}(X)} \leq 1} (\nabla_{A_1}a \times a, \cA^{-\frac{1}{4}}b)_{L^2(X)}
\\
&\leq \sup_{\|b\|_{L^{d'}(X)} \leq 1} \|\nabla_{A_1}a \times a\|_{L^{q'}(X)} \|\cA^{-\frac{1}{4}}b\|_{L^q(X)},
\end{align*}
for $q \in (1, \infty)$ and dual exponent $q'$, defined via $1 = 1/q + 1/q'$, still to be chosen.

The Sobolev Embedding Theorem yields $W^{s, p}(X) \hookrightarrow L^q(X)$ for $s \geq 0$ and $p \geq 1$ obeying $sp < d$ and $p \leq q \leq dp/(d-sp)$ by \cite[Theorem 4.12]{AdamsFournier}. Setting $p = d' = d/(d-1)$ and $s = 1/2$ allows us to choose
$$
q = \frac{dp}{d-sp} = \frac{dd'}{d - d'/2} = \frac{d^2/(d-1)}{d - d/2(d-1)} = \frac{d^2}{d^2 - d - d/2} = \frac{d}{d - 3/2}.
$$
Therefore, by the embedding $W^{\frac{1}{2}, d'}(X) \hookrightarrow L^{d/(d-3/2)}(X)$ and our choice $q = d/(d-3/2)$,
$$
\|\cA^{-\frac{1}{4}}b\|_{L^q(X)} = \|\cA^{-\frac{1}{4}}b\|_{L^{d/(d-3/2)}(X)} \leq C\|\cA^{-\frac{1}{4}}b\|_{W_{A_1}^{1/2, d'}(X)} \leq C\|b\|_{L^{d'}(X)} \leq C,
$$
using the fact that $\cA^{-\frac{1}{4}}: L^{d'}(X; \Lambda^1\otimes\ad P) \to W_{A_1}^{1/2, d'}(X; \Lambda^1\otimes\ad P)$ is a bounded operator. The dual exponent, $q'$, is given by
$$
q' = \frac{q}{q - 1} = \frac{d/(d - 3/2)}{d/(d - 3/2) - 1} = \frac{2d}{3}.
$$
Thus, substituting $q = d/(d-3/2)$ and $q' = 2d/3$,
\begin{align*}
\|\cA^{-\frac{1}{4}}(\nabla_{A_1}a \times a)\|_{L^d(X)}
&\leq
\sup_{\|b\|_{L^{d'}(X)} \leq 1} \|\nabla_{A_1}a \times a\|_{L^{2d/3}(X)} \|\cA^{-\frac{1}{4}}b\|_{L^{d/(d-3/2)}(X)}
\\
&\leq C\|\nabla_{A_1}a \times a\|_{L^{2d/3}(X)}
\\
&\leq C\|\nabla_{A_1}a\|_{L^d(X)} \|a\|_{L^{2d}(X)},
\\
&\leq C\|\nabla_{A_1} a\|_{L^d(X)} \|a\|_{W_{A_1}^{1/2, d}(X)},
\\
&\leq C\|\cA^{\frac{1}{2}} a\|_{L^d(X)} \|\cA^{\frac{1}{4}}a\|_{L^d(X)},
\end{align*}
where we used $3/2d = 1/d + 1/2d$ and the associated H\"older inequality, together with the Sobolev embedding,  $W^{\frac{1}{2}, d}(X) \hookrightarrow L^{2d}(X)$, from \cite[Theorem 4.12]{AdamsFournier}, and the fact that
$$
\|a\|_{W_{A_1}^{1/2, d}(X)} = \|\cA^{\frac{1}{4}}a\|_{L^d(X)}.
$$
This establishes the inequality \eqref{eq:Kozono_Maeda_Naito_3-7_Yang-Mills_quadratic}. Similarly, using $3/2d =  1/2d + 1/2d + 1/2d$ and the associated H\"older inequality, we find that
\begin{align*}
\|\cA^{-1/4}(a\times a \times a)\|_{L^d(X)}
&\leq
C\|a\times a \times a\|_{L^{2d/3}(X)}
\\
&\leq C\|a\|_{L^{2d}(X)}^3
\\
&\leq C\|a\|_{W_{A_1}^{1/2, d}(X)}^3
\\
&= C\|\cA^{\frac{1}{4}}a\|_{L^d(X)}^3,
\end{align*}
which gives the inequality \eqref{eq:Kozono_Maeda_Naito_3-7_Yang-Mills_cubic}.

For the Lipschitz inequalities, we write
\begin{align*}
\nabla_{A_1}a_1 \times a_1 - \nabla_{A_1}a_2 \times a_2 &= \nabla_{A_1}(a_1 - a_2) \times a_1 + \nabla_{A_1}a_2 \times (a_1 - a_2)
\\
a_1\times a_1 \times a_1 - a_2\times a_2 \times a_2 &= (a_1 - a_2)\times a_1 \times a_1 + a_2\times (a_1 - a_2) \times a_1
\\
&\quad + a_2\times a_2 \times (a_1 - a_2),
\end{align*}
and now estimate each product term, just as before, to obtain \eqref{eq:Kozono_Maeda_Naito_3-8_Yang-Mills}.
\end{proof}

We are now ready to complete the

\begin{proof}[Proof of Theorem \ref{thm:Struwe_section_4-3_local_existence_uniqueness_Yang-Mills_heat_equation_initial_data_in_Ld}]
Lemma \ref{lem:Kozono_Maeda_Naito_3-3} ensures that the Yang-Mills heat equation nonlinearity $\sF$ in \eqref{eq:Yang-Mills_heat_equation_nonlinearity_relative_rough_Laplacian_plus_one} satisfies the hypotheses \eqref{eq:Kozono_Maeda_Naito_3-7}, \eqref{eq:Kozono_Maeda_Naito_3-2_nonlinearity}, and \eqref{eq:Kozono_Maeda_Naito_3-8} of Theorem \ref{thm:Kozono_Maeda_Naito_lemma_3-4_plus_uniqueness}, with $f_0 = d_{A_1}^*F_{A_1} \in C^\infty(X; \Lambda^1\otimes\ad P))$ and taking $\cW$ to be $L^d(X; \Lambda^1\otimes\ad P)$. Hence, Theorem \ref{thm:Kozono_Maeda_Naito_lemma_3-4_plus_uniqueness} now provides a unique mild solution, $a(t)$, to the Yang-Mills heat equation
\eqref{eq:Yang-Mills_heat_equation_as_perturbation_rough_Laplacian_plus_one_heat_equation} on an interval $[0, \tau)$ with initial data $a(0) = a_0$ and the regularity
\eqref{eq:Local_existence_uniqueness_Yang-Mills_heat_equation_initial_data_in_Ld}.
\end{proof}

\begin{cor}[Strong solutions and higher-order regularity for mild solutions to the Yang-Mills heat equation with initial data in $L^d$]
\label{cor:Struwe_section_4-3_higher_order_regularity_Yang-Mills_heat_equation_initial_data_in_Ld}
Let $d \geq 2$. Then there is a constant $\beta = \beta(d) \in [\frac{1}{2}, \frac{3}{4})$ with the following significance. Let $G$ be a compact Lie group and $P$ a principal $G$-bundle over a closed, connected, orientable, smooth manifold, $X$, of dimension $d$ and Riemannian metric, $g$. Let $A_1$ be a reference connection of class $C^\infty$ on $P$. If $a_0 \in L^d(X; \Lambda^1\otimes\ad P)$ and
$$
a \in C([0, T); L^d(X; \Lambda^1\otimes\ad P) \cap C((0, T); W_{A_1}^{2\beta, d}(X; \Lambda^1\otimes\ad P)),
$$
is a mild solution to the Yang-Mills heat equation
\eqref{eq:Yang-Mills_heat_equation_as_perturbation_rough_Laplacian_plus_one_heat_equation} on an interval $[0, T)$ with initial data $a(0) = a_0$, for some $T > 0$, then
$$
a \in C([0, T); L^d(X; \Lambda^1\otimes\ad P)) \cap C^\infty((0, T) \times X; \Lambda^1\otimes\ad P)),
$$
and $a(t)$ is a classical solution on $(0, T)$. Furthermore, if $a_0 \in C^\infty(X; \Lambda^1\otimes\ad P)$, then
$$
a \in C^\infty([0, T) \times X; \Lambda^1\otimes\ad P)).
$$
\end{cor}

\begin{proof}
Set $\beta = 5/8$ and choose $p$ so that the hypotheses of Theorems
\ref{thm:Existence_uniqueness_mild_solution_Yang-Mills_heat_equation_in_W_2beta_p_initial_data_in_W_2beta_p} and
\ref{thm:Existence_uniqueness_strong_solution_Yang-Mills_heat_equation_in_W_2beta_p_initial_data_in_W_2beta_p} are obeyed:
\begin{inparaenum}[\itshape i\upshape)]
\item $d \geq 3$ and $p > d/2$ and $d/2p < 5/8$, that is, $p > 8d/10$, for example, $p = 9d/10$; or
\item $d \geq 2$ and $d/3 < p < d$ and $d/4p + 1/4 \leq  5/8 < d/2p$, that is, $p$ obeys $d/4p \leq 3/8$ and $p < 8d/10$, namely $p$ obeys $8d/12 \leq p < 8d/10$, for example, $p = 2d/3$.
\end{inparaenum}

Theorem \ref{thm:Existence_uniqueness_strong_solution_Yang-Mills_heat_equation_in_W_2beta_p_initial_data_in_W_2beta_p} now ensures that, on an interval $[t_0, T)$ for any $t_0 \in (0, T)$, the function $a(t)$ is a strong solution to the Yang-Mills heat equation
\eqref{eq:Yang-Mills_heat_equation_as_perturbation_rough_Laplacian_plus_one_heat_equation} with the regularity
\ref{eq:Sell_You_47-7_Yang-Mills_strong_solution_in_W_2beta_p_initial_data_in_W_2beta_p} on $[t_0, T)$. Furthermore, Theorem \ref{thm:Smoothness_strong_solution_Yang-Mills_heat_equation_initial_data_in_W_2beta_p} implies that $a \in C^\infty((t_0, T) \times X; \Lambda^1\otimes\ad P))$ and since $t_0 \in (0, T)$ was arbitrary, we have $a \in C^\infty((0, T) \times X; \Lambda^1\otimes\ad P))$.

If $a_0 \in C^\infty(X; \Lambda^1\otimes\ad P)$, then $C^\infty$ regularity on $[0, T)\times X$ follows from Theorem \ref{thm:Smoothness_strong_solution_Yang-Mills_heat_equation_smooth_initial_data}.
\end{proof}

We can now conclude the

\begin{proof}[Proof of Corollaries
\ref{cor:Struwe_section_4-3_higher_order_regularity_Yang-Mills_heat_equation_hybrid_critical-exponent_parabolic_Sobolev_space} and
\ref{cor:Struwe_section_4-3_higher_order_regularity_Yang-Mills_heat_equation_pure_critical-exponent_parabolic_Sobolev_space}]
Since $a_0 \in C(X; \Lambda^1\otimes\ad P)$ by hypothesis, then $a_0 \in L^d(X; \Lambda^1\otimes\ad P)$ and the conclusions follow from Corollary \ref{cor:Struwe_section_4-3_higher_order_regularity_Yang-Mills_heat_equation_initial_data_in_Ld}.
\end{proof}

\section{Local existence and uniqueness for Yang-Mills gradient flow}
\label{sec:Local_existence_yang_mills_gradient_coulomb_gauge_flow}
We begin in Section \ref{subsec:Struwe_4-4_Donaldson-DeTurck trick} by reviewing Donaldson's version \cite{DonASD} of the DeTurck trick \cite{DeTurck_1983} (for Ricci flow) used to relate solutions of the Yang-Mills gradient flow and Yang-Mills heat equations. In Section \ref{subsec:Struwe_5}, we recall Struwe's proof of existence of a weak solution to the Cauchy problem for a gauge-equivalent version of Yang-Mills gradient flow that obeys a Coulomb gauge condition. We digress in Section \ref{subsubsec:Struwe_5-1} to discuss irreducible connections, as required for Struwe's approach, while Section \ref{subsubsec:Struwe_5-2} reviews the main existence result (Proposition \ref{prop:Struwe_5-2}). In Section \ref{subsec:Struwe_6}, we describe Struwe's  uniqueness result (Proposition \ref{prop:Struwe_6-1}) for weak solutions the Cauchy problem for a gauge-equivalent version of Yang-Mills gradient flow coupled with a Coulomb gauge equation. We conclude in Section \ref{subsec:Kozono_Maeda_Naito_6} by presenting an alternative approach, due to Kozono, Maeda, and Naito \cite{Kozono_Maeda_Naito_1995}, for uniqueness modulo gauge transformations of a solution to the Cauchy problem for Yang-Mills gradient flow.

\subsection{Gauge transformations and the Donaldson-DeTurck trick for Yang-Mills gradient flow}
\label{subsec:Struwe_4-4_Donaldson-DeTurck trick}
The adaptation of the DeTurck trick \cite{DeTurck_1983} (for the Ricci flow of a Riemannian metric on a three-dimensional manifold) to establish the gauge-equivalence of a solution to the Yang-Mills gradient flow and a solution to the \emph{modified Yang-Mills gradient flow} (in the terminology of \cite[p. 96]{Kozono_Maeda_Naito_1995}) or \emph{Yang-Mills heat equation\footnote{Our terminology differs from that of R\r{a}de \cite[p. 123]{Rade_1992}, who applies the term to the Yang-Mills gradient flow equation itself, although that equation is not parabolic.}} (as we call it here) is due to Donaldson \cite{DonASD}. We shall review Struwe's treatment of the DeTurck trick here in order to confirm the minimum regularity required for the initial data, $A_0$, and path $A(t)$, for $t\in [0,T)$, and motivate our development of the regularity theory for the Yang-Mills gradient flow. As Struwe notes in \cite[Section 4.1]{Struwe_1994}, the DeTurck trick imposes a regularity requirement on a path $A(t)$ and initial data which is more subtle than it might appear at first glance; introductory descriptions (for example, \cite[Section 6.3.1]{DK}) assume for simplicity that $A_0$ is smooth and $A(t)$ are smooth with respect to the spatial variables and that $A(t)$ is smooth with respect to $t \in [0,T)$, but that obscures issues which will be important in our application.

Let $P$ be a principal $G$-bundle a closed, Riemannian, smooth manifold, $X$, of dimension $d\geq 2$. Recall that the group of (continuous or smooth) gauge transformations (automorphisms) of $P$ may be interpreted as the group of sections of the bundle $\Ad P$ over $X$ with fiberwise multiplication, that is, $\Aut P \cong C^\infty(X;\Ad P)$ by \cite[Section 3.1]{FrM}. To define \emph{Sobolev gauge transformations} of $P$, we follow Uhlenbeck \cite[Section 1]{UhlRem}; see also Freed and Uhlenbeck \cite[Appendix A]{FU} and Struwe \cite[Section 1.3]{Struwe_1994}. One considers a representation $\rho:G \subset \SO(n)$, where the metric on $G$ is induced by that on $\SO(n)$, and views $\SO(n) \subset \RR^{n\times n}$ and $\Ad P \subset \End V$, where $V = P\times_\rho\RR^n\to X$ is a real, Riemannian vector bundle of rank $n$. For $s \geq 1$ and $p\geq 2$ obeying $sp > d$, there is a Sobolev embedding $W^{s,p}(X) \hookrightarrow C(X)$ \cite[Theorem 4.12]{AdamsFournier} and $W^{s,p}(X)$ is a Banach algebra
\cite[Theorem 4.39]{AdamsFournier}. For $s \geq 1$ and $p\geq 2$ obeying $sp > d$, one defines the set of gauge transformations of $P$ of class $W^{s,p}$ by the inclusion
$$
W^{s,p}(X; \Ad P) \subset W^{s,p}(X; \End V),
$$
and, when $p = 2$ and $s>d/2$, abbreviate $H^s(X; \Ad P) = W^{s,2}(X; \Ad P)$, so that
$$
H^s(X; \Ad P) \subset H^s(X; \End V).
$$
According to \cite[Proposition A.2]{FU}, one has that $W^{s,p}(X; \Ad P)$ (respectively, $H^s(X; \Ad P)$) is a Banach (respectively, Hilbert) Lie group with Lie algebra $W^{s,p}(X; \ad P)$ (respectively, $H^s(X; \ad P)$). (The statement and proof of \cite[Proposition A.2]{FU} assume that $d=4$ and $s$ is integer and $G=\SU(2)$, but extensions to more general situations present no new difficulties.) Bethuel \cite{Bethuel_1991} provides a careful discussion of Sobolev spaces of maps between manifolds relevant to the definition of $W^{s,p}(X; \Ad P)$.

To apply the DeTurck trick \cite{DeTurck_1983} to the Yang-Mills gradient flow equation \eqref{eq:Yang-Mills_gradient_flow_equation}, following Donaldson \cite[p. 7]{DonASD} or Donaldson and Kronheimer \cite[Section 6.3.1]{DK}, one must solve a certain ordinary differential equation for a one-parameter family of gauge transformations, $u(t) \in \Aut P$ for $t\in[0,\infty)$. The required result is proved by Nagasawa \cite{Nagasawa_1989}.

\begin{lem}[Existence and uniqueness of a family of gauge transformations]
\label{lem:Kozono_Maeda_Naito_6-2}
\cite[Lemma 6.2]{Kozono_Maeda_Naito_1995}, \cite[Theorem 3.2.1 and p. 514]{Nagasawa_1989}
Let $G$ be a compact Lie group and $P$ be a principal $G$-bundle over a closed, Riemannian, smooth manifold, $X$, of dimension $d\geq 2$.
\begin{enumerate}
\item
\label{item:Kozono_Maeda_Naito_6-2_Linfty}
If $0 < T \leq \infty$ and
$$
\zeta \in L_{\loc}^\infty([0,T); L^\infty(X;\ad P)),
$$
then there exists a unique family of gauge transformations,
$$
u \in C([0,T); L^\infty(X;\Ad P)) \cap C^1((0,T); L^\infty(X;\Ad P)),
$$
such that
\begin{equation}
\label{eq:Kozono_Maeda_Naito_lemma_6-2_ode}
u^{-1}\frac{\partial u}{\partial t}(t) = \zeta(t)
\quad\forall\, t \in (0, T) \hbox{ and a.e. on } X, \quad u(0) = \id_P.
\end{equation}

\item
\label{item:Kozono_Maeda_Naito_6-2_WA1_sp}
If $A_1$ is a $C^\infty$ reference connection on $P$ and $s\in \RR$ and $p\in[1,\infty]$ obey $sp > d$ and
$$
\zeta \in L_{\loc}^\infty([0,T); W_{A_1}^{s,p}(X;\Ad P)),
$$
then
$$
u \in C([0,T); W_{A_1}^{s,p}(X;\Ad P)) \cap C^1((0,T); W_{A_1}^{s,p}(X;\Ad P)).
$$

\item
\label{item:Kozono_Maeda_Naito_6-2_smooth}
If
$$
\zeta \in C([0,T)\times X; \ad P) \cap C^\infty((0,T)\times X; \ad P)
\quad (\hbox{respectively}, \quad C^\infty([0,T)\times X; \ad P)),
$$
then
$$
u \in C([0,T)\times X; \Ad P) \cap C^\infty((0,T)\times X; \Ad P)
\quad (\hbox{respectively}, \quad C^\infty([0,T)\times X; \Ad P)).
$$
\end{enumerate}
\end{lem}


\begin{rmk}[On the regularity of the gauge transformations in Lemma \ref{lem:Kozono_Maeda_Naito_6-2}]
\label{rmk:Kozono_Maeda_Naito_lemma_6-2_regularity}
Assertion \ref{item:Kozono_Maeda_Naito_6-2_Linfty} of Lemma \ref{lem:Kozono_Maeda_Naito_6-2} is proved by Nagasawa as \cite[Theorem 3.2.1]{Nagasawa_1989}. Assertion \ref{item:Kozono_Maeda_Naito_6-2_smooth} (essentially in this form) is stated in \cite[p. 514]{Nagasawa_1989} and this, together with Assertion \ref{item:Kozono_Maeda_Naito_6-2_WA1_sp}, follow by standard results for regularity and continuous dependence on parameters of solutions to ordinary differential equations.
\end{rmk}


Before proceeding to review the more difficult case where $A_0$ is a connection of class $H^1$ on $P$, due to Struwe in \cite[Section 4.4]{Struwe_1994}, let us first combine Lemma \ref{lem:Kozono_Maeda_Naito_6-2} with the main result of the analysis due to Donaldson and Kronheimer in \cite[Section 6.3.1]{DK} in the simpler case when $A_0$ is of class $W^{s+1,p}$ with $sp>d$.

\begin{lem}[Gauge equivalence of a solution to the Yang-Mills heat and gradient flows for an initial connection of class $W^{s+1,p}$ with $sp>d$]
\cite[p. 7]{DonASD}, \cite[Section 6.3.1]{DK}
(See also \cite[Section 2]{Kono_Nagasawa_1988}, \cite[Equations (1.3) and (3.1)]{Kozono_Maeda_Naito_1995}, \cite[Section 4.1]{Jost_1991}, and \cite[Section 4.2]{Struwe_1994})
\label{lem:Donaldson_DeTurck_trick}
Let $G$ be a compact Lie group and $P$ be a principal $G$-bundle over a closed, Riemannian, smooth manifold, $X$, of dimension $d\geq 2$. Let $A_0$ be a connection of class $W^{s+1,p}$ on $P$ such that $s\geq 1$ and $p\geq 2$ obey $sp>d$ and $T \in (0,\infty]$. Assume that $A=A_1+a$ is a strong solution to the Cauchy problem for the \emph{Yang-Mills heat equation},
\begin{subequations}
\label{eq:Struwe_17_initial_value_problem}
\begin{align}
\label{eq:Struwe_17}
\frac{\partial A}{\partial t} + d_A^*F_A + d_Ad_A^*a &= 0
\quad \hbox{on } (0, T)\times X,
\\
\label{eq:Struwe_17_initial_data}
A(0) &= A_0,
\end{align}
\end{subequations}
such that
\begin{equation}
\label{eq:Struwe_17_regularity}
a \in C([0,T); W^{s+1,p}(X; \Lambda^1\otimes\ad P))
\cap C^1((0,T); W^{s+1,p}(X; \Lambda^1\otimes\ad P)).
\end{equation}
Then there is a unique solution, $u \in C([0,T); W^{s,p}(X; \Lambda^1\otimes\Ad P))
\cap C^1((0,T); W^{s,p}(X; \Lambda^1\otimes\Ad P))$, to the ordinary differential equation,
\begin{equation}
\label{eq:Struwe_18}
u^{-1}\frac{\partial u}{\partial t} = -d_A^*a
\quad\hbox{on } (0,T)\times X, \quad u(0) = \id_P,
\end{equation}
with the following significance. If\footnote{We apply the DeTurck trick in order to convert a solution of the Yang-Mills heat equation to a solution to the Yang-Mills gradient flow, which is the opposite direction to the application in \cite[Section 4]{Struwe_1994}.}
\begin{equation}
\label{eq:tildeA_gauge_transformation_A}
\tilde A := (u^{-1})^*A,
\end{equation}
then $\tilde A$ is a strong solution to the Cauchy problem for the \emph{Yang-Mills gradient flow equation},
\begin{subequations}
\label{eq:Struwe_3_and_4}
\begin{align}
\label{eq:Struwe_3}
\frac{\partial\tilde A}{\partial t} + d_{\tilde A}^*F_{\tilde A} &= 0
\quad \hbox{on } (0, T)\times X,
\\
\label{eq:Struwe_4}
\tilde A(0) &= A_0,
\end{align}
\end{subequations}
such that
\begin{equation}
\label{eq:Struwe_3_regularity}
\tilde a := \tilde A - A_1 \in C([0,T); W^{s-1,p}(X; \Lambda^1\otimes\ad P))
\cap C^1((0,T); W^{s-1,p}(X; \Lambda^1\otimes\ad P)).
\end{equation}
Conversely, if $\tilde A$ is a solution to the Yang-Mills gradient flow equation \eqref{eq:Struwe_3_and_4} with regularity \eqref{eq:Struwe_17_regularity}, then $A=u^*\tilde A$ is a solution to the Yang-Mills heat equation \eqref{eq:Struwe_17_initial_value_problem} with regularity \eqref{eq:Struwe_3_regularity}.

A solution, $A-A_1$, belongs to $C([0,T)\times X; \Lambda^1\otimes\ad P)
\cap C^\infty((0,T)\times X; \Lambda^1\otimes\ad P)$ if and only if the same holds for $\tilde A-A_1$. If $A_0$ is a connection of class $C^\infty$ on $P$, then $A-A_1$ belongs to $C^\infty([0,T)\times X; \Lambda^1\otimes\ad P)$ if and only if the same holds for $\tilde A-A_1$.
\end{lem}

\begin{rmk}[Variants of Lemma \ref{lem:Donaldson_DeTurck_trick}]
\label{lem:Donaldson_DeTurck_trick_variants}
It should be possible to state and prove other versions of Lemma \ref{lem:Donaldson_DeTurck_trick} (the `Donaldson-DeTurck trick') involving different solution concepts (weak, strong, or classical) or different combinations of spatial or temporal regularities. However, in our applications, we shall generally only need Lemma \ref{lem:Donaldson_DeTurck_trick} in the simplest setting where $A-A_1$ (or $\tilde A-A_1$) belongs to $C([0,T)\times X; \Lambda^1\otimes\ad P)
\cap C^\infty((0,T)\times X; \Lambda^1\otimes\ad P)$ or $C^\infty([0,T)\times X; \Lambda^1\otimes\ad P)$.
\end{rmk}

\begin{proof}
Let us first note that the source term, $d_A^*a$, in the ordinary differential equation \eqref{eq:Struwe_18} belongs to $C([0,T); W^{s,p}(X; \Lambda^1\otimes\Ad P))$ by hypothesis on $A$ and so a unique solution, $u$, to \eqref{eq:Struwe_18} exists by Lemma \ref{lem:Kozono_Maeda_Naito_6-2}. Since $u(0)=\id_P$, then $A$ obeys the initial condition \eqref{eq:Struwe_17_initial_data} if and only if $\tilde A$ obeys the initial condition \eqref{eq:Struwe_4}. To see that $\tilde A = (u^{-1})^*A$ is a solution to \eqref{eq:Struwe_3} if and only if $A = u^*\tilde A$ is a solution to \eqref{eq:Struwe_17}, we recall that $u$ acts on $\tilde A = A_1 + \tilde a$ by
$$
u^*\tilde A = A_1 + u^{-1}\tilde au + u^{-1}d_{A_1}u.
$$
Compare \cite[proof of Proposition A.3]{FU}. Note also that
$$
d_{A_1}\left(u^{-1} \frac{\partial u}{\partial t}\right)
=
-u^{-1}(d_{A_1}u)u^{-1}\frac{\partial u}{\partial t}
+ u^{-1} d_{A_1}\frac{\partial u}{\partial t}.
$$
Thus,
\begin{align*}
\frac{\partial u^*\tilde A}{\partial t}
&= -u^{-1}\frac{\partial u}{\partial t}u^{-1}\tilde au
+ u^{-1}\frac{\partial\tilde a}{\partial t}u
+ u^{-1}\tilde a\frac{\partial u}{\partial t}
- u^{-1}\frac{\partial u}{\partial t}u^{-1}d_{A_1}u
+ u^{-1}d_{A_1}\frac{\partial u}{\partial t}
\\
&= u^{-1}\frac{\partial\tilde a}{\partial t}u
+ u^{-1}d_{A_1}\frac{\partial u}{\partial t}
+ \left[u^{-1}\tilde au, u^{-1}\frac{\partial u}{\partial t} \right]
- u^{-1}\frac{\partial u}{\partial t}u^{-1}d_{A_1}u
\\
&= u^{-1}\frac{\partial\tilde a}{\partial t}u
+ u^{-1}d_{A_1}\frac{\partial u}{\partial t}
+ \left[u^{-1}\tilde au, u^{-1}\frac{\partial u}{\partial t} \right]
+ \left[u^{-1}d_{A_1}u, u^{-1}\frac{\partial u}{\partial t}\right]
- u^{-1}(d_{A_1}u)u^{-1}\frac{\partial u}{\partial t}
\\
&= u^{-1}\frac{\partial\tilde a}{\partial t}u
+ d_{A_1}\left(u^{-1} \frac{\partial u}{\partial t}\right)
+ \left[u^{-1}\tilde au, u^{-1}\frac{\partial u}{\partial t} \right]
+ \left[u^{-1}d_{A_1}u, u^{-1}\frac{\partial u}{\partial t}\right]
\\
&= u^{-1}\frac{\partial\tilde a}{\partial t}u
+ d_{A_1 + u^{-1}\tilde au + u^{-1}d_{A_1}u}\left(u^{-1} \frac{\partial u}{\partial t}\right)
\\
&= u^{-1}\frac{\partial\tilde a}{\partial t}u
+ d_{u^*\tilde A}\left(u^{-1} \frac{\partial u}{\partial t}\right).
\end{align*}
Therefore, $A = u^*\tilde A$ obeys (compare \cite[Equation (16)]{Struwe_1994})
\begin{equation}
\label{eq:Struwe_16}
\frac{\partial A}{\partial t} = u^{-1}\frac{\partial\tilde A}{\partial t}u
+ d_A\left(u^{-1} \frac{\partial u}{\partial t}\right),
\end{equation}
and hence,
\begin{align*}
\frac{\partial A}{\partial t} &= u^{-1}\frac{\partial\tilde A}{\partial t}u - d_Ad_A^*a
\quad\hbox{(by \eqref{eq:Struwe_18})}
\\
&= u^{-1}\left(-d_{\tilde A}^*F_{\tilde A}\right)u - d_Ad_A^*a \quad\hbox{(by \eqref{eq:Struwe_3})}
\\
&= -d_A^*F_A - d_Ad_A^*a \quad\hbox{(by \eqref{eq:tildeA_gauge_transformation_A})},
\end{align*}
which is \eqref{eq:Struwe_17}. The preceding steps are clearly reversible, so $A$ is a solution to \eqref{eq:Struwe_17} if and only if $\tilde A$ is a solution to \eqref{eq:Struwe_3}. The regularity of $\tilde A$ is implied by that of $A$ and the gauge transformation, $u$.
\end{proof}

The significance of Lemma \ref{lem:Donaldson_DeTurck_trick} is that one can now expand \eqref{eq:Struwe_17} around the linear heat equation defined by a $C^\infty$ connection, $A_1$, to give the following equation for $a = A-A_1$, expressed in schematic form,
\begin{equation}
\label{eq:Struwe_21}
\frac{\partial a}{\partial t} + \Delta_{A_1}a + d_{A_1}^*F_{A_1} + F_{A_1}\times a
+ \nabla_{A_1}a\times a + a\times a\times a = 0 \quad\hbox{a.e. on } (0, T)\times X.
\end{equation}
Observe that \eqref{eq:Struwe_21} is nonlinear parabolic equation for $a(t)$ with coefficients and source function which are $C^\infty$ in space and independent of time $t\geq 0$ when $A_1$ is a fixed connection of class $C^\infty$ on $P$.

As one might expect, that existence problem becomes more delicate when $A_0$ is of class $W^{s+1,p}$ with $s\geq 1$ and $p\geq 2$ such that $sp \leq d$ and it is no longer possible to assure the existence of continuous gauge transformations due to the lack of a Sobolev embedding $W^{s,p}(X) \hookrightarrow C(X)$ for this range of $s$ and $p$. This is the issue addressed by Struwe in \cite[Section 4.4]{Struwe_1994}, where it is assumed only that $A_0$ is a connection of class $H^1$ on $P$.

\begin{lem}[Gauge equivalence of a solution to Yang-Mills heat flow to that of Yang-Mills gradient flow for an initial connection of class $H^1$ over a manifold of dimension at most four]
\label{lem_Struwe_sections_4-2_and_4-4}
\cite[Sections 4.2 and 4.4]{Struwe_1994}
Let $G$ be a compact Lie group and $P$ be a principal $G$-bundle over a closed, Riemannian, smooth manifold, $X$, of dimension $2 \leq d\leq 4$. Suppose that $T \in (0,\infty]$ and $A_0$ is a connection on $P$ of class $H^1$ and $A(t)$, for $t \in [0,T)$, is a classical solution to the Yang-Mills heat equation \eqref{eq:Struwe_17_initial_value_problem} in the sense that, for any fixed $C^\infty$ reference connection, $A_1$, on $P$,
$$
A-A_1
\in
C([0,T);H_{A_1}^1(X;\Lambda^1\otimes\ad P))
\cap C^\infty((0,T)\times X; \Lambda^1\otimes\ad P).
$$
Let $\{t_k\}_{k\in\NN} \subset (0,T]$ be a sequence such that $t_k\searrow 0$ as $k \to \infty$ and let $\{u_k\}_{k\in\NN} \subset C^\infty([0,T)\times X; \Ad P)$ be a sequence of gauge transformations defined by solving \eqref{eq:Struwe_18} on $[t_k,T)$ with initial condition $u(t_k) = \id_P$ and setting $u_k(t) := \id_P$ for $t\in [0,t_k)$, for each $k\in\NN$. Define
$$
A_k := u_k^*A \quad\hbox{on } [0, T) \times P, \quad\forall\, k \in \NN.
$$
Then $\{A_k\}_{k\in\NN}$ converges to a \emph{Struwe weak} solution, $\tilde A$, to the Cauchy problem for the Yang-Mills gradient flow equation \eqref{eq:Struwe_3_and_4} in the sense of Definition \ref{defn:Struwe_weak_solution} with, as $k\to\infty$,
\begin{align*}
A_k &\to \tilde A \in C([0,T); L^2(X; \Lambda^1\otimes\ad P)),
\\
F_{A_k} &\to F_{\tilde A} \in C([0,T); L^2(X; \Lambda^2\otimes\ad P)).
\end{align*}
Furthermore, the \emph{energy}, $\sE(\tilde A(t)) = (1/2)\|F_{\tilde A}(t)\|_{L^2(X)}^2$, obeys
\begin{equation}
\label{eq:Struwe_12}
\sE(\tilde A(0))
=
\int_0^t \|d_{\tilde A}^*F_{\tilde A}(s)\|_{L^2(X)}^2\,ds
+
\sE(\tilde A(t)), \quad\forall\, t\in [0,T).
\end{equation}
\end{lem}

The following useful property of a solution to a Yang-Mills gradient flow equation was noted by Donaldson in \cite[p. 7]{DonASD} and explained further by Jost in \cite[Section 4.1]{Jost_1991}.

\begin{lem}[Energy-decreasing property of a solution to modified Yang-Mills gradient flow]
\cite[p. 7]{DonASD}, \cite[Section 4.1]{Jost_1991}
\label{lem:Jost_equations_4-1-2_and_3}
Let $G$ be a compact Lie group and $P$ be a principal $G$-bundle over a closed, Riemannian, smooth manifold, $X$, of dimension $d\geq 2$. Suppose that $T \in (0,\infty]$ and
$$
\zeta \in L_{\loc}^2(0,T;H_{A_1}^1(X;\Lambda^1\otimes\ad P)).
$$
If $A(t)$, for $t \in [0,T)$, is a strong solution to the \emph{modified Yang-Mills gradient flow equation},
\begin{equation}
\label{eq:Jost_4-1-3}
\frac{\partial A}{\partial t} = -d_A^*F_A + d_A\zeta,
\end{equation}
in the sense that, for any fixed $C^\infty$ reference connection, $A_1$, on $P$,
$$
A-A_1
\in
L_{\loc}^2(0,T;H_{A_1}^2(X;\Lambda^1\otimes\ad P))
\cap H_{\loc}^1(0,T;L^2(X;\Lambda^1\otimes\ad P)),
$$
then $A(t)$ obeys
\begin{equation}
\label{eq:Jost_4-1-2}
\frac{\partial}{\partial t} \|F_A\|_{L^2(X)}^2 = -2\|d_A^*F_A\|_{L^2(X)}^2
\quad\hbox{a.e. on }(0,T)\times X,
\end{equation}
and so the energy, $\sE(A(t)) = (1/2)\|F_A(t)\|_{L^2(X)}^2$, is non-increasing on $[0,T)$.
\end{lem}

In particular, Lemma \ref{lem:Jost_equations_4-1-2_and_3} shows that passage from pure Yang-Mills gradient flow \eqref{eq:Struwe_3} to the Yang-Mills heat equation \eqref{eq:Struwe_17} via the Donaldson-DeTurck trick in Lemma \ref{lem:Donaldson_DeTurck_trick} preserves the non-increasing energy property of Yang-Mills gradient flow.

\subsection{Gauge normalization for weak solutions to the Cauchy problem for Yang-Mills gradient flow}
\label{subsec:Struwe_5}
Next we review the question of uniqueness of solutions to Yang-Mills gradient flow. There are two approaches to this problem. In the first, due to Struwe \cite[Sections 5 and 6]{Struwe_1994}, one appends a Coulomb gauge condition to the Yang-Mills gradient flow and considers the uniqueness of the gauge transformation used to pass between Yang-Mills gradient and heat flows via the Donaldson-DeTurck trick. In the second, due to Kozono, Maeda, and Naito \cite[Section 6]{Kozono_Maeda_Naito_1995}, one shows that any two solutions of Yang-Mills gradient flow must be gauge equivalent.

We first review the main results of \cite[Section 5]{Struwe_1994}. For the proof of uniqueness of solutions to Yang-Mills gradient flow, we follow Struwe and consider the gauge-equivalent version \eqref{eq:Struwe_16} of the Yang-Mills gradient flow \eqref{eq:Struwe_3}. For that purpose, one needs to specify a gauge condition. Before getting into details, observe that uniqueness of $s(t) \in \Omega^0(X;\ad P)$, for $t\in [0,T)$, and hence of the time-varying family gauge transformations, $S(t) \in \Aut P$, determined by \eqref{eq:Struwe_18} can only hold if the operators,
$$
d_{A(t)}: \Omega^0(X;\ad P) \to \Omega^1(X;\ad P), \quad\hbox{for } t\in [0,T),
$$
are invertible. For smooth $A(t)$, this condition is equivalent to a topological condition on the connection.

See Propositions \ref{prop:Struwe_5-2} and \ref{prop:Struwe_6-1} for Struwe's construction of a unique solution, $\bar A(t) = A_1+\alpha(t)+\bar(t)$, to a gauge-equivalent version of Yang-Mills gradient flow augmented with Coulomb gauge condition, $d_{\bar A}^*\bar a = 0$ a.e. on $(0,T)\times X$, as an alternative to his construction of a strong solution to the Yang-Mills heat equation.

\subsubsection{Irreducible connections}
\label{subsubsec:Struwe_5-1}
We briefly review the definition of an irreducible connection from \cite[Section 4.2.2]{DK}, essentially following the abbreviated discussion in \cite[Section 5.1]{Struwe_1994}. Given a Lie group $G$ and a smooth connection $A$ on a principal $G$-bundle $P$, we let
\begin{equation}
\label{eq:Donaldson_Kronheimer_4-2-7}
\Gamma_A := \{u \in \Aut P: u^*A = A\}
\end{equation}
denote the isotropy subgroup of $A$ and recall that $\Gamma_A \subset G$ is a closed Lie subgroup with Lie algebra,
$$
\gamma_A := \{\zeta \in \Omega^0(X;\ad P): d_A\zeta = 0\}.
$$
Recall that $\Gamma_A$ always contains the center, $C(G)$, of $G$ \cite[p. 132]{DK}. By the customary abuse of terminology \cite[pp. 132--133]{DK}, we shall a connection $A$ \emph{irreducible} if
$$
\Gamma_A = C(G),
$$
where $C(G)$ is the center of $G$. We recall the holonomy definition of reducible connection in Remark below and the fact that, when $G = \SU(2)$ or $\SO(3)$, the two definitions of `irreducible' agree \cite[p. 133]{DK}.

\begin{rmk}
\label{rmk:Donaldson_Kronheimer_page_131}
A smooth connection $A$ on a principle $G$-bundle, $P\to X$, is \emph{reducible} if for each point $x\in X$, the holonomy maps, $T_\ell$, of all loops $\ell$ based at $x$ lie in some proper subgroup of the automorphism group, $\Aut P_x \cong G$ \cite[Section 4.2.2]{DK}. If the base space, $X$ is connected one can restrict attention to a single fibre, $P_{x_0}$, and obtain a holonomy group, $H_A \subset G$, or more precisely a conjugacy class of subgroups. It can be shown that $H_A$ is a closed Lie subgroup of $G$ and, when $X$ is connected, that $\Gamma_A$ is isomorphic to the centralizer of $H_A$ in $G$ \cite[Lemma 4.2.8]{DK}.
\end{rmk}

For a connection $A$ of class $H^1$ on $P$, we shall require \emph{irreducibility} in the sense that \cite[Equation (26)]{Struwe_1994},
\begin{equation}
\label{eq:Struwe_26}
\|\zeta\|_{H_A^1(X)} \leq C\|d_A\zeta\|_{L^2(X)},
\quad\forall\, \zeta \in H_A^1(X;\ad P),
\end{equation}
for some positive constant, $C = C(A,g)$. The following lemma asserts that this constant, $C$, can be chosen locally  uniformly.

\begin{lem}[Openness of the irreducibility condition]
\label{lem:Struwe_5-1}
\cite[Lemma  5.1]{Struwe_1994}
Suppose a connection $A_0$ of class $H^1$ on $P$ satisfies \eqref{eq:Struwe_26} with constant $C_0 = C(A_0,g)$.  Then there exists an open neighborhood $\sU$ of the origin in $H_{A_0}^1(X;\Lambda^1\otimes\ad P)$ and a positive constant, $C$, such that any $A \in A_0 + \sU$  is irreducible in the sense that
$$
\|\zeta\|_{H_A^1(X)} \leq C\|d_A\zeta\|_{L^2(X)},
\quad\forall\, \zeta \in H_A^1(X;\ad P).
$$
\end{lem}


\subsubsection{Gauge-fixing for a weak solution to the Cauchy problem for Yang-Mills gradient flow}
\label{subsubsec:Struwe_5-2}
We next recall Struwe's global analogue of Uhlenbeck's theorem \cite[Theorem 2.1 and Corollary 2.2]{UhlLp} on the existence of local Coulomb gauges, depending smoothly on the connection.

\begin{prop}[Gauge-fixing for a weak solution to the Cauchy problem for Yang-Mills gradient flow]
\label{prop:Struwe_5-2}
(Compare \cite[Proposition 5.2 and Equations (29), (30), (31), and (32)]{Struwe_1994}.)
Let $G$ be a compact Lie group and $P$ a principal $G$-bundle over a closed, connected smooth manifold, $X$, of dimension $2 \leq d \leq 4$ and Riemannian metric, $g$. If $A_1$ is a reference connection of class $C^\infty$ on $P$ and $A_0$ is a connection of class $H^1$ on $P$ that is irreducible in the sense of \eqref{eq:Struwe_26}, then there is a constant, $\eps_0 = \eps_0(A_1,g) \in (0, 1]$ and, if $\alpha_0 \in \Omega^1(X;\ad P)$ is such that
$$
a_0 := A_0-A_1-\alpha_0 \in H_{A_1}^1(X; \Lambda^1\otimes\ad P)
$$
obeys
$$
\|a_0\|_{H^1_{A_1}(X)} < \eps_0,
$$
there are a constant, $\tau_0 = \tau_0(A_0,A_1,\alpha_0,g) \in (0, 1]$, and a sequence of gauge transformations,
$$
\{u_k\}_{k\in\NN} \subset C^\infty([0,\tau_0)\times X;\Ad P),
$$
with the following significance.

Let $\alpha \in C^\infty([0,\infty)\times X; \Lambda^1\otimes\ad P)$ be a solution to the Cauchy problem \eqref{eq:Struwe_19_with_augmented_connection_Laplacian} for the linear heat equation with initial data $\alpha(0) = \alpha_0$. Let $A = A_1+a$ be a weak solution (in the sense of Definition \ref{defn:Struwe_weak_solution}) to the Cauchy problem \eqref{eq:Struwe_3_and_4} on $[0, \tau_0)\times P$ for Yang-Mills gradient flow with initial data $A(0)=A_0$. If
$$
A_k := u_k^*A \quad\hbox{and}\quad a_k := A_k - A_1 - \alpha, \quad\forall\, k \in \NN,
$$
then
\begin{align*}
a_k &\to \bar a \quad\hbox{in }
L^\infty(0,\tau_0; L^2(X; \Lambda^1\otimes\ad P))
\cap H^1(0,\tau_0; L^2(X; \Lambda^1\otimes\ad P)), \quad\hbox{as } k \to \NN,
\\
u_k^{-1}\frac{\partial u_k}{\partial t}
&\to \zeta \quad\hbox{in } L^2(0,\tau_0; H^1_{A_1}(X; \Lambda^1\otimes\ad P)),
\quad\hbox{as } k \to \NN,
\end{align*}
and $\bar A := A_1 + \alpha + \bar a$ and $\zeta$ obey
\begin{align*}
{}& \bar a \in L^\infty(0,\tau_0; H^1_{A_1}(X; \Lambda^1\otimes\ad P))
\cap H^1(0,\tau_0; L^2(X; \Lambda^1\otimes\ad P)),
\\
{}& \bar a(t) \to  a_0 \quad\hbox{in } H_{A_1}^1(X; \Lambda^1\otimes\ad P)
\quad\hbox{as } t \searrow 0,
\\
{}& d_{\bar A}^*\bar a = 0 \quad\hbox{a.e. on } (0,\tau_0)\times X,
\end{align*}
and $\bar A$ is a weak solution (in the variational sense of Equation \eqref{eq:Struwe_Definition_2-1_variational_equation}) to the Cauchy problem,
$$
\frac{\partial \bar A}{\partial t} = -d_{\bar A}^*F_{\bar A} + d_{\bar A}\zeta,
\quad \bar A(0) = A_0.
$$
\end{prop}

\begin{rmk}[On the hypotheses and statement of Proposition \ref{prop:Struwe_5-2}]
\label{rmk:Struwe_proposition_5-2_hypotheses_and_statement}
Our version of Proposition \ref{prop:Struwe_5-2} differs in several respects from \cite[Proposition 5.2]{Struwe_1994}.
\begin{enumerate}
\item \emph{Dimension of the base manifold.} In keeping with our desire to provide a framework that includes R\r{a}de's results when $X$ has dimension two or three \cite{Rade_1992}, we allow $X$ to have dimension $d$ with $2\leq d \leq 4$; naturally, the proof of Proposition \ref{prop:Struwe_5-2} simplifies when $d<4$.

\item \emph{Gauge transformations.} Struwe allows the gauge transformation, $u$, to belong to the $H_{A_1}^1([0,\tau_0)\times X;\Ad P)$-closure of $C^\infty([0,\tau_0)\times X;\Ad P)$, the family of smooth gauge transformations of $P$ varying smoothly with $t\in [0,\tau_0)$. Because gauge transformations of class $W^{s,p}$ with $s\geq 1$ and $p\geq 2$ obeying $sp\leq d$ need not be continuous (see the discussion in \cite[Appendix A]{FU}), we prefer to state Proposition \ref{prop:Struwe_5-2} in terms of sequences of smooth families of gauge transformations, $\{u_k\}_{k\in\NN} \subset C^\infty([0,\tau_0)\times X;\Ad P)$, as in Struwe's proof of his \cite[Proposition 5.2]{Struwe_1994}.

\item \emph{Reference connection and initial data.} In \cite[Section 5.2]{Struwe_1994}, Struwe continues his choice of a time-varying family of background connections, with (after translation to our notation) $A_{\textrm{bg}}(t) = A_1+a_{\textrm{bg}}(t)$ for $t\in[0,\infty)$, where $A_1$ is a fixed $C^\infty$ reference connection and
$$
a_{\textrm{bg}} \in C([0,\infty); H_{A_1}^1(X; \Lambda^1\otimes\ad P))
\cap C^\infty((0,\infty)\times X; \Lambda^1\otimes\ad P)
$$
    is determined as the classical solution to a Cauchy problem for the linear heat equation \cite[Equation (19)]{Struwe_1994},
$$
\frac{\partial a_{\textrm{bg}}}{\partial t} + \Delta_{A_1}a_{\textrm{bg}} = 0
\quad\hbox{on } (0,\infty)\times X,
\quad a_{\textrm{bg}}(0) = A_0-A_1,
$$
and thus $A_{\textrm{bg}}(0) = A_0$, the initial data of class $H^1$, while $\|A_0-A_1\|_{H_{A_1}^1(X)} < \eps$ and $\eps \in (0,1]$ is arbitrarily small. Our choice of time-varying family, $\alpha(t)$, with initial data $\alpha(0)=\alpha_0$ of arbitrary norm $\|\alpha_0\|_{H_{A_1}^1(X)}$ removes the need for $A_1$ to be arbitrarily $H_{A_1}^1(X)$-close to $A_0$, albeit at the cost of introducing the additional dependency on $\alpha_0$. This also serves to explain the fact that our family, $\bar a(t)$, has non-zero initial data, $\bar a(0) = A_0-A_1-\alpha_0$, unlike in \cite[Proposition 5.2]{Struwe_1994}.

\item \emph{Estimates for gauge transformations.} \Apriori estimates for the sequence of gauge transformations, $\{u_k\}_{k\in\NN}$, and the limit $\zeta \in H_{A_1}^1(X; \ad P)$ of $u_k^{-1}\partial_t u_k$ can be extracted from Struwe's
    \cite[Lemma 5.3]{Struwe_1994}.
\end{enumerate}
The proof of Proposition \ref{prop:Struwe_5-2} follows Struwe's proof of \cite[Proposition 5.2]{Struwe_1994} \mutatis just as the methods of Sections \ref{subsec:Struwe_page_137_contraction_mapping_small_initial_data_in_H1} and \ref{subsec:Struwe_page_137_contraction_mapping_arbitrary_initial_data_in_H1} adapt those of Struwe in \cite[Section 4.4]{Struwe_1994}.
\end{rmk}


\subsection{Uniqueness of a solution to the Cauchy problem for Yang-Mills gradient flow}
\label{subsec:Struwe_6}
We review the main result of \cite[Section 6]{Struwe_1994}. We continue the setup of Section \ref{subsubsec:Struwe_5-2}, namely, $G$ is a compact Lie group, $P$ is a principal $G$-bundle over a closed, connected smooth manifold, $X$, of dimension $2 \leq d \leq 4$ and Riemannian metric, $g$, and $A_1$ is a reference connection of class $C^\infty$ on $P$, and $A_0$ is a connection of class $H^1$ on $P$ that is irreducible in the sense of \eqref{eq:Struwe_26}.

Suppose that $T > 0$ and $A = A_1 + a$, with
$$
a \in C([0,T); L^2(X;\Lambda^1\otimes\ad P))
\cap H^1(0,T; L^2(X;\Lambda^1\otimes\ad P)),
$$
is a weak solution (in the sense of Definition \ref{defn:Struwe_weak_solution}) to the Cauchy problem for Yang-Mills gradient flow on $P$ with initial data $A(0) = A_0$. Let $\bar A(t) = A_1 + \alpha(t) + \bar a(t)$, for $t \in [0,\tau_0)$, be the corresponding family of normalized connections produced by Proposition \ref{prop:Struwe_5-2}. The family $\bar A$ is a weak solution to the initial-value problem,
\begin{subequations}
\label{eq:Struwe_29_and_30_and_31}
\begin{align}
\label{eq:Struwe_29}
\frac{\partial\bar A}{\partial t} &= - d_{\bar A}^*F_{\bar A} + d_{\bar A}\zeta
\quad\hbox{a.e. on } (0,\tau_0)\times X,
\\
\label{eq:Struwe_30}
d_{\bar A}^*\bar a &= 0 \quad\hbox{a.e. on } (0,\tau_0)\times X,
\\
\label{eq:Struwe_31}
\bar a(0) &= A_0-A_1-\alpha_0,
\end{align}
\end{subequations}
where \eqref{eq:Struwe_29} is obeyed in variational sense of Equation \eqref{eq:Struwe_Definition_2-1_variational_equation} and
\begin{equation}
\label{eq:Struwe_32}
\begin{aligned}
\bar a &\in L^\infty(0,\tau_0; H_{A_1}^1(X; \Lambda^1\otimes\ad P))
\cap H^1(0,\tau_0; L^2(X; \Lambda^1\otimes\ad P)),
\\
F_{\bar A} &\in C([0,\tau_0]; L^2(X; \Lambda^2\otimes\ad P)),
\\
\zeta &\in L^2(0,\tau_0; H_{A_1}^1(X; \ad P)),
\end{aligned}
\end{equation}
and in \eqref{eq:Struwe_31} the solution, $\bar a$, attains its initial value in the sense of $H_0^1(0,\tau_0; L^2(X; \Lambda^1\otimes\ad P))$. The following result shows that the solution, $\bar A$, to \eqref{eq:Struwe_29_and_30_and_31} is unique when $A_0$ is irreducible in the sense of \eqref{eq:Struwe_26}.

\begin{prop}
\label{prop:Struwe_6-1}
\cite[Proposition 6.1]{Struwe_1994}
Let $G$ be a compact Lie group and $P$ a principal $G$-bundle over a closed, connected smooth manifold, $X$, of dimension $2 \leq d \leq 4$ and Riemannian metric, $g$. If $A_1$ is a reference connection of class $C^\infty$ on $P$ and $A_0$ is a connection of class $H^1$ on $P$ that is irreducible in the sense of \eqref{eq:Struwe_26}, then there is a constant, $\eps_0 = \eps_0(A_1,g) \in (0, 1]$ and, if $\alpha_0 \in \Omega^1(X;\ad P)$ is such that
$$
a_0 := A_0-A_1-\alpha_0 \in H_{A_1}^1(X; \Lambda^1\otimes\ad P)
$$
obeys
$$
\|a_0\|_{H^1_{A_1}(X)} < \eps_0,
$$
there is a constant, $\tau_0 = \tau_0(A_0,A_1,\alpha_0,g) \in (0, 1]$, with the following significance. There exists a \emph{unique} solution, $(\bar a, \zeta)$, to \eqref{eq:Struwe_29_and_30_and_31} on $[0, \tau_0)$ satisfying \eqref{eq:Struwe_32}. In addition,
\begin{gather*}
\bar a \in L^2(0,\tau_0; H_{A_1}^2(X;\Lambda^1\otimes\ad P))
\cap C^\infty((0,\tau_0)\times X; \Lambda^1\otimes\ad P),
\\
\zeta \in C^\infty((0,\tau_0)\times X; \ad P).
\end{gather*}
If $A_0$ is of class $C^\infty$, then
$$
a\in C^\infty([0,\tau_0)\times X; \Lambda^1\otimes\ad P)
\quad\hbox{and}\quad
\zeta \in C^\infty([0,\tau_0)\times X; \ad P).
$$
\end{prop}

Again, there are small differences between our version of Proposition \ref{prop:Struwe_6-1} and Struwe's \cite[Proposition 6.1]{Struwe_1994}; those differences follow the pattern described in Remark \ref{rmk:Struwe_proposition_5-2_hypotheses_and_statement}, comparing our Proposition \ref{prop:Struwe_5-2} and Struwe's \cite[Proposition 5.2]{Struwe_1994}.


\subsection{Uniqueness modulo gauge transformations of a solution to the Cauchy problem for Yang-Mills gradient flow}
\label{subsec:Kozono_Maeda_Naito_6}
Before we proceed to discuss the approach of Kozono, Maeda, and Naito \cite[Section 6]{Kozono_Maeda_Naito_1995} to uniqueness of solutions to the Cauchy problem for Yang-Mills gradient flow, we digress to review why the question of irreducibility of the initial data, $A_0$, does not arise in R\r{a}de's approach to that problem when $X$ has dimension two or three. Rather than apply the Donaldson-DeTurck trick to convert the non-parabolic Yang-Mills gradient flow equation on $(0, T)\times P$ to a parabolic Yang-Mills heat equation, he instead notes that if $A(t)$ is a solution to the Yang-Mills gradient flow equation, then $(A(t), \Omega(t)) = (A(t), F_A(t))$ is a solution to the system \cite[Equation (4.4)]{Rade_1992},
\begin{equation}
\label{eq:Rade_1992}
\begin{aligned}
\frac{\partial A}{\partial t} + d_A^*\Omega &= 0, \quad A(0) = A_0,
\\
\frac{\partial \Omega}{\partial t} + \Delta_A\Omega &= 0, \quad \Omega(0) = F_{A_0},
\end{aligned}
\end{equation}
where $\Delta_A := d_A^*d_A + d_Ad_A^*$ is the Hodge Laplace operator \eqref{eq:Lawson_page_93_Hodge_Laplacian}. The system \eqref{eq:Rade_1992} is strong enough to give uniqueness without appeal to the Donaldson-DeTurck trick and hence no need to appeal to irreducibility of $A_0$.

The following uniqueness result, due to Kozono, Maeda, and Naito, complements Struwe's \cite[Theorem 6.1]{Struwe_1994} and does not require $A_0$ to be irreducible. However, it imposes a stronger regularity requirement on the solutions and thus, implicitly, a stronger regularity requirement on the initial data than is explicitly stated in Theorem \ref{thm:Kozono_Maeda_Naito_6-1}.

\begin{thm}[Uniqueness up to gauge transformation for weak solutions to the Cauchy problem for Yang-Mills gradient flow]
\label{thm:Kozono_Maeda_Naito_6-1}
\cite[Theorem 6.1]{Kozono_Maeda_Naito_1995}
Let $G$ be a compact Lie group and $P$ a principal $G$-bundle over a closed, connected smooth manifold, $X$, of dimension $2 \leq d \leq 4$ and Riemannian metric, $g$, and $A_1$ a reference connection of class $C^\infty$ on $P$ and $A_0$ a connection of class $H^1$ on $P$, and $T \in (0, \infty]$. Let $A^i(t) = A_1 + a_i(t)$ be two weak solutions (in the sense of Definition \ref{defn:Weak_solution_Yang-Mills_gradient_flow}) to the Cauchy problem for Yang-Mills gradient flow with initial data, $A^i(0) = A_0$ for $i=1,2$. If in addition,
$$
a_i \in L^q(0, T; L^r(X;\Lambda^1\otimes\ad P)), \quad i = 1, 2,
$$
for $q \geq 2$ and $r > 4$ with $2/q+4/r \leq 1$ and
$$
d_{A_1}^*a_i \in L^\infty(0,T; W_{A_1}^{1,\infty}(X;\ad P)), \quad i = 1, 2,
$$
then there exist gauge transformations,
$$
u_i \in C_b([0,T); W_{A_1}^{1,\infty}(X;\Ad P)) \cap C^1((0,T); W_{A_1}^{1,\infty}(X;\Ad P)),
\quad i = 1,2,
$$
solving
$$
u_i^{-1}\frac{\partial u_i}{\partial t} = d_{A_1}^*a_i \quad\hbox{on } (0,T)\times X,
\quad u_i(0) = \id_P \quad\hbox{for } i=1,2,
$$
such that
$$
u_1^*A^1 = u_2^*A^2 \quad\hbox{a.e. on } (0,T)\times P.
$$
\end{thm}

\begin{rmk}[On the hypotheses and statement of Theorem \ref{thm:Kozono_Maeda_Naito_6-1}]
\label{rmk:Kozono_Maeda_Naito_theorem_6-1_hypotheses_and_statement}
Our version of Theorem \ref{thm:Kozono_Maeda_Naito_6-1} differs slightly from that of \cite[Theorem 6.1]{Kozono_Maeda_Naito_1995}. We clarify what appear to be ambiguities in the statements of regularity for the terms $d_{A_1}^*a_i$ and gauge transformations $u_i$ for $i=1,2$. The regularity of the gauge transformations in Theorem \ref{thm:Kozono_Maeda_Naito_6-1} follows from Lemma \ref{lem:Kozono_Maeda_Naito_6-2}, which is due to Nagasawa \cite[Theorem 3.2.1 and p. 514]{Nagasawa_1989}.

The \cite[Theorem 6.1]{Kozono_Maeda_Naito_1995} is stated and proved for the case $d=4$ but, as usual, the proof simplifies when $d < 4$. The statement of the Cauchy problems for the gauge transformations, $u_i$, is given in the \cite[proof of Theorem 6.1]{Kozono_Maeda_Naito_1995}.
\end{rmk}

\chapter[{\L}ojasiewicz-Simon inequality and convergence for gradient systems]{The {\L}ojasiewicz-Simon gradient inequality, stability, and convergence for gradient systems}
\label{chapter:Lojasiewicz-Simon_gradient_inequality_and_stability_and_convergence}

\section[Finite-dimensional {\L}ojasiewicz gradient inequality]{{\L}ojasiewicz gradient inequality and finite-dimensional dynamical systems}
\label{sec:Lojasiewicz_gradient_inequality_and_finite-dimensional_dynamical_systems}
While the application of the infinite-dimensional {\L}ojasiewicz-Simon gradient inequality to prove global existence and convergence of solutions to gradient systems can appear rather technical at first glance (for example, see Simon \cite{Simon_1983}), one can gain a useful understanding of the fundamental ideas by restricting attention to the far simpler setting of finite-dimensional gradient systems. Thus, by way of introduction to this chapter, we shall review an elementary result due to {\L}ojasiewicz (namely, Theorem \ref{thm:Lojasiewicz_1984} below), keeping in mind that this appeared around the same time as Simon's development of his infinite-dimensional gradient inequality and its application to certain infinite-dimensional gradient systems in geometric analysis in \cite{Simon_1983}.

\begin{thm}[Global existence and convergence of a solution to a gradient system in $\RR^n$]
\label{thm:Lojasiewicz_1984}
\cite[Theorem 1]{Lojasiewicz_1984}
\cite[p. 1592]{Lojasiewicz_1993}
Let $\sE$ be an analytic, non-negative function on a neighborhood of the origin in $\RR^n$ such that $\sE(0) = 0$. Then there exists a neighborhood, $U = \{x\in\RR^n: |x| < \sigma\}$, of the origin such that each trajectory, $u_{x_0}(t)$, with $u_{x_0}(0) = x_0 \in U$, of the system,
\begin{equation}
\label{eq:Lojasiewicz_1984_gradient_system}
\dot u(t) = -\sE'(u(t)),
\end{equation}
is defined on $[0, \infty)$, has finite length, and converges uniformly to a point $u_{x_0}(\infty) \in Z := \{z \in U: \sE'(z) = 0\}$ as $t \to \infty$. For a constant $\theta \in (0,1)$ depending only on $\sE$, one has
\begin{align*}
|u_{x_0}(t) - x_0| &\leq \int_0^t |\dot u_{x_0}(s)|\,ds \leq \frac{\sE(x_0)^{1-\theta}}{1-\theta},
\\
|u_{x_0}(\infty) - u_{x_0}(t)| &\leq \int_t^\infty |\dot u_{x_0}(s)|\,ds \leq \frac{(1+t)^{\theta-1}}{1-\theta},
\quad\text{for } 0 \leq t < \infty.
\end{align*}
\end{thm}

Note that because $\sE(x) \geq 0$ for all $x$ in its domain and $\sE(0) = 0$, then $\sE$ achieves its absolute minimum at $x = 0$ and thus $\sE'(0) = 0$, so $0 \in Z$. Assuming global existence of a solution, $u(t)$ for $t \in [0,\infty)$, to \eqref{eq:Lojasiewicz_1984_gradient_system}, the convergence assertion in Theorem \ref{thm:Lojasiewicz_1984} (though not the convergence rate) is also proved by Chill and Jendoubi as \cite[Theorem 2.2]{Chill_Jendoubi_2003}, in the abstract, infinite-dimensional setting. The convergence rate estimate asserted by Theorem \ref{thm:Lojasiewicz_1984} can be improved and, indeed, we shall pursue such improvements when we consider its infinite-dimensional analogues.

To prove Theorem \ref{thm:Lojasiewicz_1984}, {\L}ojasiewicz applied the following version of his gradient inequality \cite{Lojasiewicz_1965}:

\begin{thm}[Finite-dimensional {\L}ojasiewicz and Simon gradient inequalities]
\label{thm:Huang_2-3-1}
\cite[Theorem 2.3.1]{Huang_2006}
\footnote{There is a typographical error in the statement of \cite[Theorem 2.3.1 (i)]{Huang_2006}, as Huang omits the hypothesis that $\sE'(z) = 0$; also our statement differs slightly from that of \cite[Theorem 2.3.1 (i)]{Huang_2006}, but is based on original sources.}
Let $U \subset \RR^n$ be an open subset, $z \in U$, and let $\sE: U \to \RR$ be a real-valued function.
\begin{enumerate}
\item
\label{item:Huang_theorem_2-3-1_i}
If $\sE$ is real analytic on a neighborhood of $z$ and $\sE'(z) = 0$, then there exist constants $\theta \in (0,1)$ and $\sigma > 0$ such that
\begin{equation}
\label{eq:Lojasiewicz_1984_star}
|\sE'(x)| \geq |\sE(x) - \sE(z)|^\theta,
\quad\forall\, x \in \RR^n, \ |x - z| < \sigma.
\end{equation}

\item
\label{item:Huang_theorem_2-3-1_ii}
Assume that $\sE$ is a $C^2$ function and $\sE'(z) = 0$. If the connected component, $C$, of the critical point set, $\{x \in U : \sE'(x) = 0\}$, that contains $z$ has the same dimension as the kernel of the Hessian matrix $\Hess_\sE(z)$ of $\sE$ at $z$ locally near $z$, and $z$ lies in the interior of the component, $C$, then there are positive constants, $c$ and $\sigma$, such that
\begin{equation}
\label{eq:Simon_1996_lemma_1_page_80}
|\sE'(x)| \geq c|\sE(x) - \sE(z)|^{1/2},
\quad\forall\, x \in \RR^n, \ |x - z| < \sigma.
\end{equation}
\end{enumerate}
\end{thm}

Theorem \ref{thm:Huang_2-3-1} \eqref{item:Huang_theorem_2-3-1_i} is well known and was stated by {\L}ojasiewicz in \cite{Lojasiewicz_1963} and proved by him as \cite[Proposition 1, p. 92]{Lojasiewicz_1965} and Bierstone and Milman as \cite[Proposition 6.8]{BierstoneMilman}; see also the statements by Chill and Jendoubi \cite[Proposition 5.1 (i)]{Chill_Jendoubi_2003} and by {\L}ojasiewicz \cite[p. 1592]{Lojasiewicz_1993}.

Theorem \ref{thm:Huang_2-3-1} \eqref{item:Huang_theorem_2-3-1_ii} was proved by Simon as \cite[Lemma 1, p. 80]{Simon_1996} and Haraux and Jendoubi as \cite[Theorem 2.1]{Haraux_Jendoubi_2007}; see also the statement by Chill and Jendoubi \cite[Proposition 5.1 (ii)]{Chill_Jendoubi_2003}.

{\L}ojasiewicz used methods of \emph{semi-analytic sets} \cite{Lojasiewicz_1965} to prove Theorem \ref{thm:Huang_2-3-1} \eqref{item:Huang_theorem_2-3-1_i}. For the inequality \eqref{eq:Lojasiewicz_1984_star}, unlike \eqref{eq:Simon_1996_lemma_1_page_80}, the constant, $c$, is equal to one while $\theta \in (0,1)$. In general, so long as $c$ is positive, its actual value is irrelevant to applications; the value of $\theta$ in the infinite-dimensional setting \cite[Theorem 2.4.2 (i)]{Huang_2006}, at least, is restricted to the range $[1/2, 1)$ and $\theta=1/2$ is optimal \cite[Theorem 2.7.1]{Huang_2006}.

We recall the following well-known local existence and uniqueness results from the classical theory of ordinary differential equations \cite{Hartman_2002}.

\begin{thm}[Peano existence]
\label{thm:Hartman_II-2-1}
\cite[Theorem II.2.1]{Hartman_2002}
Let $t_0 \in \RR$, $a > 0$, $b > 0$, $x_0 \in \RR^n$ (for $n \geq 1$) and
\[
\bar R := \{(t,x) \in \RR\times\RR^n: t_0 \leq t \leq t_0+a, \ |x - x_0| \leq b\}.
\]
Let $f:\bar R\to \RR^n$ be continuous and $\alpha := \min\{a, b/M\}$, where $M > 0$ is a constant such that $|f| \leq M$ on $\bar R$. Then there exists a solution, $u(t)$ for $t \in [t_0, t_0 + \alpha]$, to the initial value problem,
\begin{equation}
\label{eq:Hartman_II-1-1}
\dot u(t) = f(t, u(t)), \quad u(t_0) = x_0.
\end{equation}
\end{thm}

\begin{thm}[Picard-Lindel\"of existence and uniqueness]
\label{thm:Hartman_II-1-1}
\cite[Theorem II.1.1]{Hartman_2002}
Assume the hypotheses of Theorem \ref{thm:Hartman_II-2-1} and, in addition, that $f$ is uniformly Lipschitz continuous with respect to $x$, that is,
\[
|f(t,x_1) - f(t,x_2)| \leq K|x_1-x_2|, \quad\forall\, (t,x_1), (t,x_2) \in \bar R.
\]
Then there exists a unique solution, $u(t)$ for $t \in [t_0, t_0 + \alpha]$, to the initial value problem \eqref{eq:Hartman_II-1-1}.
\end{thm}

We now turn to the proof of Theorem \ref{thm:Lojasiewicz_1984} following \cite{Lojasiewicz_1984}, but provide a few additional details.

\begin{proof}[Proof of Theorem \ref{thm:Lojasiewicz_1984}]
Let $G \subset \RR^d$ (an open subset) be the domain of $\sE$ and let $[0,\tau_{x_0})$ be the maximal interval of existence for a solution $u_{x_0}(t)$ to \eqref{eq:Lojasiewicz_1984_gradient_system}, for each $x_0 \in G$. The existence of solution, $u_{x_0}(t)$ with $t\in [0,\tau_{x_0})$, for some $\tau_{x_0} > 0$ is ensured by Theorem \ref{thm:Hartman_II-2-1} and, in addition, uniqueness is assured by Theorem \ref{thm:Hartman_II-1-1}, since $\sE$ is analytic.

Thanks to uniqueness, one has the alternative,
\[
u_{x_0} \equiv x_0 \quad\text{or}\quad u_{x_0}(t) \in G \less Z, \quad \forall\, t \in [0, \tau_{x_0}),
\]
according to whether $x_0 \in Z = \{z \in G: \sE'(z) = 0\}$ or $x_0 \in G\less Z$. From the {\L}ojasiewicz gradient inequality \eqref{eq:Lojasiewicz_1984_star}, the continuity of $\sE$, and the facts that $\sE \geq 0$ on $G$ and $\sE(0) = 0$, there is an open neighborhood of the origin, $Q \subset G$, such that
\begin{equation}
\label{eq:Lojasiewicz_1984_1}
0 \leq \sE(x) < 1 \quad\text{and}\quad |\sE'(x)| \geq |\sE(x)|^\theta, \quad\forall\, x \in Q.
\end{equation}
Define
\begin{equation}
\label{eq:Huang_proposition_2-3-2_finite_dimensions_proof_line_one}
H_x(t) := \sE(u_x(t)), \quad\forall \, (t,x) \in [0, \tau_x)\times G.
\end{equation}
One therefore has, by \eqref{eq:Lojasiewicz_1984_gradient_system},
\begin{equation}
\label{eq:Huang_2-3-16a_finite_dimensions}
\dot H_{x_0}(t)
=
\langle \sE'(u_{x_0}(t)), \dot u_{x_0}(t) \rangle
=
-|\dot u_{x_0}(t)||\sE'(u_{x_0}(t))|
=
-|\sE'(u_{x_0}(t))|^2, \quad\forall \, x_0 \in G,
\end{equation}
and
\begin{multline}
\label{eq:Lojasiewicz_1984_2}
H_{x_0}(t) = 0 \quad\text{or}\quad 0 < H_{x_0}(t) < 1, \quad\forall\, t\in [0,\tau_{x_0}),
\\
\text{according to whether } x_0 \in Z \text{ or } Q\less Z.
\end{multline}
To prove Theorem \ref{thm:Lojasiewicz_1984}, it suffices to show that there exists an open neighborhood, $U \subset Q$, of the origin such that if $x_0 \in U \less Z$, then one has
\begin{align}
\label{eq:Lojasiewicz_1984_3}
\tau_{x_0} &= \infty \quad\text{and}\quad u_{x_0}([0,\infty)) \subset Q,
\\
\label{eq:Lojasiewicz_1984_4}
|u_{x_0}(t_1) - u_{x_0}(t_0)| &\leq \int_{t_0}^{t_1} |\dot u_{x_0}(s)|\,ds \leq \frac{H_{x_0}(t_0)^{1-\theta}}{1-\theta},
\quad\text{for } 0 \leq t_0 \leq t_1 < \infty,
\\
\label{eq:Lojasiewicz_1984_5}
H_{x_0}(t) &\leq (1+t)^{-1}, \quad\text{for } 0 \leq t < \infty.
\end{align}
Properties \eqref{eq:Lojasiewicz_1984_3}, \eqref{eq:Lojasiewicz_1984_4}, and \eqref{eq:Lojasiewicz_1984_5}, hold trivially for $x_0 \in Z$.

Let us first prove \eqref{eq:Lojasiewicz_1984_4} under the hypothesis that $u_{x_0}([t_0,t_1]) \subset Q$ (where $x_0 \in Q\less Z$ and $0\leq t_0 \leq t_1 < \tau_{x_0}$). One has, for $t \in (t_0, t_1)$,
\begin{align*}
|\dot u_{x_0}(t)|
&= -\frac{\dot H_{x_0}(t)}{|\sE'(u_{x_0}(t))|}
\quad\text{(by \eqref{eq:Huang_2-3-16a_finite_dimensions})}
\\
&\leq
-\frac{\dot H_{x_0}(t)}{\sE(u_{x_0}(t))^\theta}
\quad\text{by \eqref{eq:Lojasiewicz_1984_1}}
\\
&= -\frac{\dot H_{x_0}(t)}{H_{x_0}(t)^\theta}
\quad\text{(by \eqref{eq:Huang_proposition_2-3-2_finite_dimensions_proof_line_one})}
\\
&=
-\frac{\left(H_{x_0}(t)^{1-\theta}\right)'}{1-\theta}
\quad\text{(as $H_{x_0}(t) > 0$ by \eqref{eq:Lojasiewicz_1984_2}).}
\end{align*}
In particular, $(H_{x_0}(t)^{1-\theta})' \leq 0$ on $(t_0,t_1)$ and
$H_{x_0}(t)^{1-\theta}$ is decreasing function of $t \in [t_0, t_1]$ and we obtain the inequality \eqref{eq:Lojasiewicz_1984_4}, since
\begin{align*}
\int_{t_0}^{t_1} |\dot u_{x_0}(t)|\,dt
&\leq
-\frac{1}{1-\theta} \int_{t_0}^{t_1} \left(H_{x_0}(t)^{1-\theta}\right)'\,dt
\\
&=
-\frac{1}{1-\theta} \left( H_{x_0}(t_1)^{1-\theta} - H_{x_0}(t_0)^{1-\theta}\right)
\\
&= \frac{1}{1-\theta} \left( H_{x_0}(t_0)^{1-\theta} - H_{x_0}(t_1)^{1-\theta}\right)
\\
&\leq \frac{1}{1-\theta} H_{x_0}(t_0)^{1-\theta}.
\end{align*}
Choose $\eps>0$ small enough that $B_\eps(0) \subset Q$. We have
\begin{equation}
\label{eq:Lojasiewicz_1984_0_leq_energy_leq_eps_over_2}
0 \leq \frac{\sE(x)^{1-\theta}}{1-\theta} < \frac{\eps}{2},
\end{equation}
for all $x \in U := B_\sigma(0)$ and small enough $\sigma$ with $0 < \sigma < \eps/2$, since $\sE(0) = 0$ and $\sE$ is continuous and non-negative on its domain, $G$.

Let us prove that \eqref{eq:Lojasiewicz_1984_3} holds for this open set, $U$, and which then implies \eqref{eq:Lojasiewicz_1984_4} in the general case --- without the temporary restriction that $u_{x_0}([t_0,t_1]) \subset Q$. Property \eqref{eq:Lojasiewicz_1984_3} holds because, if $x_0 \in U \less Z$, we find that
\begin{equation}
\label{eq:Lojasiewicz_1984_path_in_Q_implies_path_in_B_eps}
u_{x_0}([0,t]) \subset Q \quad\text{for } 0 < t < \tau_{x_0}
\implies
u_{x_0}([0,t]) \subset B_\eps(0),
\end{equation}
since, if $0 \leq s \leq t$, we have
\begin{align*}
|u_{x_0}(s)|
&\leq
|u_{x_0}(s) - u_{x_0}(0)| + |u_{x_0}(0)|
\\
&=
|u_{x_0}(s) - u_{x_0}(0)| + |x_0| \quad\text{(as $u_{x_0}(0)=x_0$)}
\\
&\leq
\frac{H_{x_0}(0)^{1-\theta}}{1-\theta} + |x_0|
\quad\text{(by \eqref{eq:Lojasiewicz_1984_4} when $u_{x_0}([0,s]) \subset Q$)}
\\
&=
\frac{\sE(x_0)^{1-\theta}}{1-\theta} + |x_0|
\quad\text{(by \eqref{eq:Huang_proposition_2-3-2_finite_dimensions_proof_line_one})}
\\
&< \frac{\eps}{2} + \frac{\eps}{2}
= \eps \quad\text{(by \eqref{eq:Lojasiewicz_1984_0_leq_energy_leq_eps_over_2}
and $x_0\in B_\sigma(0)$ with $\sigma<\eps/2$)},
\end{align*}
and thus
\[
u_{x_0}(s) \in Q \quad\forall\, s \in [0, t] \implies |u_{x_0}(s)| < \eps \quad\forall\, s \in [0, t].
\]
This verifies the claim \eqref{eq:Lojasiewicz_1984_path_in_Q_implies_path_in_B_eps} and this in turn leads to \eqref{eq:Lojasiewicz_1984_3}.

To better understand why \eqref{eq:Lojasiewicz_1984_3}, namely the assertion that $\tau_{x_0}=\infty$, indeed follows from \eqref{eq:Lojasiewicz_1984_path_in_Q_implies_path_in_B_eps}, suppose $\hat\tau_{x_0} \in (0, \tau_{x_0})$ is the smallest (finite) time such that
\begin{gather*}
u_{x_0}(s) \in Q \quad\forall\, s \in [0, \hat\tau_{x_0}] \quad\text{and}\quad
|u_{x_0}(t)| < \eps \quad\forall\, t \in [0, \hat\tau_{x_0}),
\\
\text{but}\quad |u_{x_0}(\hat\tau_{x_0})| \geq \eps.
\end{gather*}
But property \eqref{eq:Lojasiewicz_1984_path_in_Q_implies_path_in_B_eps} implies that $|u_{x_0}(\hat\tau_{x_0})| < \eps$, a contradiction, so we must have $\hat\tau_{x_0} = \infty$ and thus $\tau_{x_0}=\infty$.

Regarding \eqref{eq:Lojasiewicz_1984_5}, let $x_0 \in U \less Z$. One has
\[
\dot H_{x_0} \leq -H_{x_0}^{2\theta} \leq -H_{x_0}^2
\quad\text{(as $0<H_{x_0}<1$ by \eqref{eq:Lojasiewicz_1984_2}),}
\]
which gives $\dot H_{x_0}/H_{x_0}^2 \leq -1$ and
\[
(1/H_{x_0} - t)' = -\dot H_{x_0}/H_{x_0}^2 - 1 \geq 0
\]
and
\[
1/H_{x_0} - t \geq 1,
\]
which is equivalent to \eqref{eq:Lojasiewicz_1984_5}.

We can now reinterpret \eqref{eq:Lojasiewicz_1984_4} and \eqref{eq:Lojasiewicz_1984_5} in terms of the length of the trajectory, $\{u(t): 0 \leq t < \infty\}$, distance between the point $u_{x_0}(t)$ and the initial data $x_0$, and energy decay along the trajectory. Setting $t_0 = 0$ and $t_1=t$ in \eqref{eq:Lojasiewicz_1984_4} and recalling that $u_{x_0}(0) = x_0$ and $H_{x_0}(t) = \sE(u_{x_0}(t))$ by \eqref{eq:Huang_proposition_2-3-2_finite_dimensions_proof_line_one} yields
\[
|u_{x_0}(t) - x_0| \leq \int_0^t |\dot u_{x_0}(s)|\,ds \leq \frac{\sE(x_0)^{1-\theta}}{1-\theta}, \quad\text{for } 0 \leq t < \infty.
\]
Similarly, setting $t_0 = t$ and $t_1=\infty$ in \eqref{eq:Lojasiewicz_1984_4} yields
\[
|u_{x_0}(\infty) - u_{x_0}(t)| \leq \int_t^\infty |\dot u_{x_0}(s)|\,ds \leq \frac{\sE(u_{x_0}(t))^{1-\theta}}{1-\theta}, \quad\text{for } 0 \leq t < \infty,
\]
and hence the convergence rate estimate, since \eqref{eq:Lojasiewicz_1984_5} gives
\[
\sE(u_{x_0}(t)) \leq (1+t)^{-1}, \quad\text{for } 0 \leq t < \infty.
\]
This completes the proof of Theorem \ref{thm:Lojasiewicz_1984}.
\end{proof}

\section{Abstract gradient inequalities}
\label{sec:Huang_2}
In this section, we review abstract gradient inequalities for functions on Banach spaces and, in particular, the infinite-dimensional {\L}ojasiewicz-Simon gradient inequality. The monograph of Huang \cite{Huang_2006} provides a comprehensive introduction to gradient inequalities and we refer to it for further details and background --- see, in particular, \cite[Chapter 2]{Huang_2006}.

\subsection{Basic properties of gradient maps}
\label{subsec:Huang_2-1}
We refer to \cite[Section 2.1]{Huang_2006}.

\subsubsection{Differentiable and analytic maps}
\label{subsubsec:Huang_2-1A}
We refer to \cite[Section 2.1A]{Huang_2006}; see also \cite[Section 2.3]{Berger_1977}. Let $\sX, \sY$ be two Banach spaces and $\sL(\sX,\sY)$ denote the Banach space of bounded, linear operators from $\sX$ to $\sY$. Let $\sU\subset\sX$ be an open subset and $\sF:\sU \to \sY$ be a map. Recall that $\sF$ is \emph{G\^ateaux differentiable} at a point $u \in \sU$ with a G\^ateaux derivative, $\sF'(u) \in \sL(\sX,\sY)$, if
$$
\lim_{t\to 0} \frac{1}{t} \|\sF(u + tv) - \sF(u) - \sF'(u)tv\|_\sY = 0, \quad \forall\, v \in \sX.
$$
Furthermore, if $\sF$ is G\^ateaux differentiable at $u \in \sU$ and
$$
\lim_{v\to 0} \frac{1}{\|v\|_\sX}\|\sF(u + v) - \sF(u) - \sF'(u)v\|_\sY = 0,
$$
then $\sF$ is said to be \emph{Fr\'echet differentiable} at $u \in \sU$. If $\sF$ is G\^ateaux differentiable near $u$ and the G\^ateaux derivative is continuous at $u$, then $\sF$ is \emph{Fr\'echet differentiable} at $u$ \cite[Proposition 2.7.5]{Deimling_1985}.

Recall that $\sF$ is (real) \emph{analytic} at $u \in \sU$ if there exists a constant $r > 0$ and a sequence of continuous symmetric $n$-linear forms, $L_n:\sX \times  \cdots \times \sX \to \sY$, such that $\sum_{n\geq 1} \|L_n\| r^n < \infty$ and there is a positive constant $\eps = \eps(u)$ such that
$$
\sF(u + v) = \sF(v) + \sum_{n\geq 1} L_n(v^n), \quad \|v\|_\sX < \eps,
$$
where $v^n \equiv (v,\ldots,v) \in \sX \times \cdots \times \sX$ ($n$-fold product). If $\sF$ is differentiable (respectively, analytic) at every point $u \in \sU$, then $\sF$ is differentiable (respectively, analytic) on $\sU$.

\subsubsection{Smooth and analytic inverse and implicit function theorems for maps on Banach spaces}
\label{subsubsec:Smooth_and_analytic_inverse_and_implicit function theorems}
Statements and proofs of the Inverse Function Theorem for $C^k$ maps of Banach spaces are provided by Abraham, Marsden, and Ratiu \cite[Theorem 2.5.2]{AMR}, Deimling \cite[Theorem 4.15.2]{Deimling_1985}, Zeidler \cite[Theorem 4.F]{Zeidler_nfaa_v1}; statements and proofs of the Inverse Function Theorem for \emph{analytic} maps of Banach spaces are provided by Berger \cite[Corollary 3.3.2]{Berger_1977} (complex), Deimling \cite[Theorem 4.15.3]{Deimling_1985} (real or complex), and Zeidler \cite[Corollary 4.37]{Zeidler_nfaa_v1} (real or complex). The corresponding $C^k$ or Analytic Implicit Function Theorems are proved in the standard way as corollaries, for example \cite[Theorem 2.5.7]{AMR} and \cite[Theorem 4.H]{Zeidler_nfaa_v1}.

\subsubsection{Gradient maps}
\label{subsubsec:Huang_2-1B}
We refer to \cite[Section 2.1B]{Huang_2006}; see also \cite[Section 2.5]{Berger_1977}. Let $\sX$ be a Banach space with norm $\|\cdot\|_\sX$ and let $\sY$ be a linear subspace of the dual space $\sX'$ which becomes a Banach space under its own norm $\|\cdot\|_\sY$, so the embedding $\sY \hookrightarrow \sX'$ is continuous.

\begin{defn}[Gradient map]
\cite[Definition 2.1.1]{Huang_2006}
\label{defn:Huang_2-1-1}
Let $\sU\subset \sX$ be an open subset of a Banach space, $\sX$, and let $\sY$ be a Banach space with continuous embedding, $\sY \hookrightarrow \sX'$. A continuous map, $\sM:\sU\to \sY$, is called a \emph{gradient map} if there exists a $C^1$ function, $\sE:\sU\to\RR$, such that $\sM(u) = \sE'(u)$ for all $u \in \sU$ in the sense that,
$$
\sE'(u)h = \langle \sM(u), h \rangle_{\sX'\times\sX}, \quad \forall\, u \in \sU, \quad h \in \sX,
$$
where $\langle \cdot , \cdot \rangle_{\sX'\times\sX}$ is the canonical bilinear form on $\sX'\times\sX$. The real-valued function, $\sE$, is called a \emph{potential} for the map $\sM$.
\end{defn}

We shall ultimately choose $\sY$ to be a Hilbert space and denote $\sY \equiv \sH$. We recall the following basic facts concerning gradient maps.

\begin{prop}[Properties of gradient maps]
\cite[Proposition 2.1.2]{Huang_2006}
\label{prop:Huang_2-1-2}
Let $\sU\subset \sX$ be an open subset of a Banach space, $\sX$, and let $\sM:\sU \to \sY \subset \sX'$ be a continuous map. Then the following hold.
\begin{enumerate}
\item If $\sM$ is of class $C^1$, then $\sM$ is a gradient map if and only if all of its Fr\'echet derivatives, $\sM'(u)$ for $u \in \sU$, are symmetric in the sense that
$$
\langle \sM'(u)v, w \rangle_{\sX'\times\sX} = \langle \sM'(u)w, v \rangle_{\sX'\times\sX}, \quad\forall\, u, v, w \in \sU.
$$
\item A bounded linear operator $\sA: \sX \to \sX'$ is a gradient operator if and only if $\sA$ is symmetric, in which case a potential for $\sA$ is given by $\sE(u) = \frac{1}{2}\langle \sA u, u \rangle_{\sX'\times\sX}$, for $u \in \sX$.
\item If $\sM$ is an analytic gradient map, then any potential $\sE:\sU\to\RR$ such that $\sM = \sE'$ is analytic as well.
\end{enumerate}
\end{prop}

\subsection{Gradient inequalities}
\label{subsubsec:Huang_2-2}
The goal of this section is to describe analogues of Theorem \ref{thm:Huang_2-3-1} \eqref{item:Huang_theorem_2-3-1_i} for analytic functionals on open subsets of Banach spaces. We closely follow the development in \cite[Section 2.2]{Huang_2006}.

\begin{defn}[{\L}ojasiewicz-Simon gradient inequality]
\cite[Definition 2.2.2 (ii)]{Huang_2006}
\label{defn:Huang_2-2-2}
Let $\sM:\sU\to \sY$ be a gradient map in the sense of Definition \ref{defn:Huang_2-1-1}. The map $\sM$ satisfies a \emph{{\L}ojasiewicz-Simon gradient inequality} near a given point $\varphi \in \sU$ if it has the form,
\begin{equation}
\label{eq:Definition_Huang_2-2-2}
\|\sE'(u)\|_\sY \geq c|\sE(u) - \sE(\varphi)|^\theta, \quad \forall\, u \in \sU, \quad \|u-\varphi\|_\sX < \sigma,
\end{equation}
for some positive constants $c$, $\sigma$, and $\theta \in [1/2,1)$.
\end{defn}

\begin{rmk}[Generalized gradient inequalities]
More generally, given a subset $\sV \subset \sU$, one says (see \cite[Definitions 2.2.1 and 2.2.2 (i)]{Huang_2006}) that $\sM$ satisfies a \emph{gradient inequality} in $\sV$ if there exists a Lebesgue-measurable function $\phi:\RR \to [0,\infty)$ such that $1/\phi \in L^1_{\loc}(\RR)$ and
$$
\|\sM(v)\|_\sY \geq \phi(\sE(v)), \quad \forall\, v \in \sV.
$$
Accordingly, one says that $\sM$ satisfies a gradient inequality near a given point in $\sU$ if it satisfies a gradient inequality in an open neighborhood of that point. One obtains a {\L}ojasiewicz-Simon gradient inequality near a point $\varphi \in \sU$ when $\phi(x) = c|x - \sE(\varphi)|^\theta$, for $x \in \RR$, and $\sV$ is the set of $v \in \sX$ such that $\|v - \varphi\|_\sX < \sigma$, for some positive constants $c$, $\sigma$, and $\theta \in [1/2,1)$. (One may write $\sE'(v)$ or $\sM(v)$ in the gradient inequality.) As explained in \cite[Chapters 3, 4, and 5]{Huang_2006}, many of the convergence and stability results one can deduce from the {\L}ojasiewicz-Simon gradient inequality also follow from more general gradient inequalities. We shall restrict our attention to the {\L}ojasiewicz-Simon gradient inequality and its application in this article since more general gradient inequalities are not necessarily easier to prove.
\end{rmk}

Recall that the set of \emph{critical points} of $\sE$ is defined by
\begin{equation}
\label{eq:Huang_page_25_defn_critical_point_set}
\sS := \{ \varphi \in \sU: \sE'(\varphi) = 0\} = \{ \varphi \in \sU: \sM(\varphi) = 0\},
\end{equation}
while points in $\sU \less \sS$ are called \emph{regular}. We note the useful

\begin{lem}[Discreteness of critical values]
\label{lem:Huang_page_25}
\cite[p. 25]{Huang_2006}
If a gradient map satisfies a {\L}ojasiewicz-Simon gradient inequality near a critical point $\varphi$, then there is only one critical value, $\sE(\varphi)$, in some neighborhood of $\varphi$, that is, we must have $\sE(\varphi_1) = \sE(\varphi)$ for any critical point $\varphi_1$ that is sufficiently near $\varphi$.
\end{lem}

\begin{proof}
If $\|\varphi_1 - \varphi\|_\sX < \sigma$, then \eqref{eq:Definition_Huang_2-2-2} implies that $c|\sE(\varphi_1) - \sE(\varphi)|^\theta \leq \|\sM(\varphi_1)\|_\sY$ and so $\sE(\varphi_1) = \sE(\varphi)$, since $\sM(\varphi_1) = 0$ because $\varphi_1$ is a critical point.
\end{proof}

\subsection{Abstract {\L}ojasiewicz-Simon gradient inequalities}
\label{subsec:Huang_2-4}
The {\L}ojasiewicz gradient inequality \cite[Th\'eor\`eme 4]{Lojasiewicz_1963}, \cite[Proposition 1, p. 92]{Lojasiewicz_1965} in finite dimensions was extended by Simon \cite[Theorem 3]{Simon_1983} to an infinite-dimensional setting, where the Banach and Hilbert spaces in Theorem \ref{thm:Huang_2-4-2}, respectively, are $C^{2,\alpha}$ H\"older and $L^2$ spaces of sections of a Riemannian vector bundle over a closed, finite-dimensional Riemannian manifold. However, that setting is too restrictive for the analysis of the Yang-Mills energy functional. Rather than directly adapt Simon's proof of \cite[Theorem 3]{Simon_1983}, as carried out, for example, in \cite{MMR, Rade_1992, RadeThesis, Yang_2003aim, Wilkin_2008}, we shall instead derive our versions of the {\L}ojasiewicz-Simon gradient inequality for the Yang-Mills energy functional from generalizations of \cite[Theorem 3]{Simon_1983} to the setting of abstract Banach and Hilbert spaces that we review in this subsection. Bierstone and Milman \cite{BierstoneMilman} provide a more recent exposition of the difficult proof of the finite-dimensional {\L}ojasiewicz gradient inequality. An application to gradient flow in finite dimensions was given by {\L}ojasiewicz in \cite{Lojasiewicz_1984}.

We shall first describe the hypotheses for Theorem \ref{thm:Huang_2-4-2}.

\begin{hyp}[Hypotheses for the abstract {\L}ojasiewicz-Simon gradient inequality with Hilbert space gradient norm]
\label{hyp:Huang_2-4_H1_H2_H3}
Assume the following conditions.
\begin{enumerate}
\item Let $\sH$ be a Hilbert space, $\sA:\sD(\sA)\subset \sH \to \sH$ a linear, positive definite, self-adjoint operator, and $\sH_\sA := (\sD(\sA), (\cdot,\cdot)_\sA)$ be the Hilbert space with inner product,
$$
(u,v)_\sA := (\sA u, \sA v)_\sH, \quad\forall\, u, v \in \sD(\sA),
$$
where $(\cdot, \cdot)_\sH$ is the inner product on $\sH$.

\item Let $\sX \subset \tilde\sX$ be Banach spaces such the following embeddings are continuous,
$$
\sX \hookrightarrow \sH_\sA, \quad \tilde \sX \hookrightarrow \sH.
$$
\item Let $\sE:\sX\to\RR$ be a function with an analytic gradient map, $\sM \equiv \sE':\sU\subset \sX\to \tilde\sX$, where $\sU\subset \sX$ is an open subset, and having the following properties:
\begin{enumerate}
\item $\sM$ is a Fredholm map of index zero, that is, for each $u \in \sX$, the bounded linear operator,
$$
\sM'(u): \sX \to \tilde \sX,
$$
is a Fredholm operator of index zero.

\item For each $u \in \sU$, the bounded, linear operator,
$$
\sM'(u): \sX \to \tilde \sX,
$$
has an extension
$$
\sM_1(u): \sH_A \to \sH
$$
which is symmetric\footnote{See \cite[Section 7.4]{Brezis}, \cite[Section 5.3.3]{Kato}, or \cite[Section 7.3]{Yosida}.} and also a Fredholm operator of index zero and such that the map
$$
\sU \ni u \mapsto \sM_1(u) \in \sL(\sH_\sA, \sH) \quad\hbox{is continuous},
$$
or, equivalently, the map $\sU \ni u \mapsto \sM_1(u)\sA^{-1} \in \sL(\sH)$ is continuous.
\end{enumerate}
\end{enumerate}
\end{hyp}

Recall that to say $P \in \sL(\sX, \tilde\sX)$ is a Fredholm operator of index zero means that $P$ has closed range, finite-dimensional kernel and cokernel, and $\dim\Ker P = \dim\Coker P$.

Under the preceding conditions, we consider a critical point $\varphi \in \sU$ of $\sE$. Without loss of generality we may assume $\varphi=0$ and thus $\sM(0)=0$ when convenient; indeed, $\sE$ has critical point $\varphi$ if and only if $\sE_\varphi := \sE(\cdot - \varphi)$ has critical point zero.

\begin{thm}[Abstract {\L}ojasiewicz-Simon gradient inequality with Hilbert space gradient norm]
\cite[Theorem 2.4.2 (i)]{Huang_2006}
\label{thm:Huang_2-4-2}
Assume Hypothesis \ref{hyp:Huang_2-4_H1_H2_H3} on $\sE$, $\sM$, $\sM_1$, $\sH$, $\sU$, $\sX$, and $\tilde \sX$. If $\varphi \in \sU$ is a critical point of $\sE$, that is, $\sM(\varphi) = 0$, then there are constants, $c \in (0,,\infty)$ and $\sigma \in (0,1]$ and $\theta \in [1/2,1)$,
such that
\begin{equation}
\label{eq:Simon_2-2}
\|\sE'(u)\|_\sH \geq c|\sE(u) - \sE(\varphi)|^\theta, \quad \forall\, u \in \sU \hbox{ such that } \|u-\varphi\|_\sX < \sigma.
\end{equation}
\end{thm}

\begin{rmk}[Related results]
\label{rmk:Theorem_Huang_2-4-2_related_results}
Related results are due to Chill \cite[Theorem 3.10]{Chill_2003}, Chill, Haraux, and Jendoubi \cite{Chill_Haraux_Jendoubi_2009}, and Jendoubi \cite[Proposition 1.3]{Jendoubi_1998jfa}. Haraux provides a recent review in \cite{Haraux_2012}. Huang notes \cite[p. 41]{Huang_2006} that the original result of Simon \cite[Theorem 3]{Simon_1983} was improved by Feireisl and Simondon \cite[Proposition 6.1]{Feireisl_Simondon_2000}, Rybka and Hoffmann \cite[Theorem 3.2]{Rybka_Hoffmann_1998}, \cite[Theorem 3.2]{Rybka_Hoffmann_1999}, and Takac \cite[Proposition 8.1]{Takac_2000} by replacing the $L^2(X)$ norm used by Simon with dual Sobolev norms, such as $H^{-1}(X) \equiv W^{-1,2}(X)$, and replacing the $C^{2,\alpha}$ H\"older norm used by Simon to define the neighborhood of the critical point with a Sobolev $H^1(X) \equiv W^{1,2}(X)$ norm. However, as far as we can tell, R\r{a}de \cite{Rade_1992, RadeThesis} was the first to make such an improvement of Simon's result, though R\r{a}de's work does not appear to be well-known among the broader applied mathematics and mathematical physics communities. Jendoubi \cite{Jendoubi_1998jfa} gives a simplified approach to Simon's method and proves gradient inequalities for gradient maps of the semilinear form, $\sM(u) = Au + f(x,u)$, for $x \in \Omega \subset \RR^d$.
\end{rmk}

The following result will be more useful, however, in our application to the Yang-Mills energy functional.

\begin{thm}[Abstract {\L}ojasiewicz-Simon gradient inequality with dual Banach space gradient norm]
\cite[Theorem 1]{Feehan_Maridakis_Lojasiewicz-Simon_harmonic_maps}, \cite[Theorem 2.4.5]{Huang_2006}
\label{thm:Huang_2-4-5}
Let $\sX$ be a Banach space that is continuously embedded in a Hilbert space $\sH$. Let $\sU \subset \sX$ be an open subset, $\sE:\sU\to\RR$ be an analytic function, and $\varphi\in\sU$ be a critical point of $\sE$, that is, $\sE'(\varphi) = 0$. Assume that $\sE''(\varphi):\sX\to \sX'$ is a Fredholm operator with index zero. Then there are constants, $c \in (0,,\infty)$ and $\sigma \in (0,1]$ and $\theta \in [1/2,1)$, such that
\begin{equation}
\label{eq:Simon_2-2_dualspacenorm}
\|\sE'(u)\|_{\sX'} \geq c|\sE(u) - \sE(\varphi)|^\theta, \quad \forall\, u \in \sU \hbox{ such that } \|u-\varphi\|_\sX < \sigma.
\end{equation}
\end{thm}

\begin{rmk}[History and related results]
\label{rmk:Theorem_Huang_2-4-5_related_results}
Theorem \ref{thm:Huang_2-4-5} was stated by Huang as \cite[Theorem 2.4.5]{Huang_2006}, but not proved and its hypotheses do not match those of the cited \cite[Proposition 3.3]{Huang_Takac_2001}, where it is assumed that $\sX$ is a Hilbert space --- see Feehan and Maridakis \cite{Feehan_Maridakis_Lojasiewicz-Simon_harmonic_maps} for a detailed discussion and more general statements. See also Haraux, and Jendoubi \cite[Theorem 2.1]{Haraux_Jendoubi_2007} and, in particular, \cite[Theorem 4.1]{Haraux_Jendoubi_2011} for a result due to Haraux and Jendoubi which is similar to Theorem \ref{thm:Huang_2-4-5} and which they argue is optimal based on examples that they discuss in \cite[Section 3]{Haraux_Jendoubi_2011}. Additional related results include those of Haraux, Jendoubi, and Kavian \cite[Proposition 1.1 and Theorem 2.1]{Haraux_Jendoubi_Kavian_2003}, Huang and Takac, \cite[Proposition 3.3]{Huang_Takac_2001} and Takac \cite[Proposition 8.1]{Takac_2000}.
\end{rmk}

\begin{rmk}[Applications of abstract {\L}ojasiewicz-Simon gradient inequalities]
\label{rmk:Theorems_Huang_2-4-2_and_2-4-5_applications}
Applications of abstract {\L}ojasiewicz-Simon gradient inequalities to asymptotic convergence questions are described by Chill and Fiorenza, and Jendoubi \cite{Chill_Fiorenza_2006}, Chill and Jendoubi \cite{Chill_Jendoubi_2003, Chill_Jendoubi_2007}, Feireisl and Simondon \cite{Feireisl_Simondon_2000}, Huang and Takac \cite{Huang_Takac_2001}, Jendoubi \cite{Jendoubi_1998jde}, Rybka and Hoffmann \cite{Rybka_Hoffmann_1998, Rybka_Hoffmann_1999}, and many other authors.
\end{rmk}

One sufficient condition for the gradient inequality to hold with the optimal exponent, $\theta=1/2$, is given by

\begin{prop}[Optimal exponent in the abstract {\L}ojasiewicz-Simon gradient inequality]
\cite[Proposition 2.7.1]{Huang_2006}
Let $\sH$ be a real Hilbert space, $\sE:\sU \subset \sH \to \sH$ be a $C^2$ functional, and $\varphi \in \sU$ be a critical point of $\sE$. If the self-adjoint operator $\sE''(\varphi)$ is invertible, then there are positive constants $c$ and $\sigma$ such that \eqref{eq:Simon_2-2} holds with $\theta = 1/2$.
\end{prop}

We have the following slight generalization of Theorem \ref{thm:Huang_2-4-2}.

\begin{cor}[Abstract {\L}ojasiewicz-Simon gradient inequality with Hilbert space gradient norm]
\label{cor:Huang_2-4-2_generalized}
Assume the hypotheses of Theorem \ref{thm:Huang_2-4-2}, except that now the operators $\sM'(u):\sX\to\tilde\sX$ and $\sM_1(u): \sH_A \to \sH$ are only required to be Fredholm of index zero at $u=\varphi$ rather than all $u \in \sU$. Then the conclusions of Theorem \ref{thm:Huang_2-4-2} continue to hold.
\end{cor}

\begin{proof}
We shall describe the minor modifications required to Huang's \cite[Proof of Proposition 2.4.1, pp. 35--41]{Huang_2006} (Lyapunov-Schmidt reduction) and \cite[Proof of Theorem 2.4.2 (i), pp. 41-41]{Huang_2006} ({\L}ojasiewicz-Simon gradient inequality). As in \cite{Huang_2006}, we may assume without loss of generality that $\varphi = 0$ by replacing $\sE:\sU\to\RR$ by $\sF:\sU_0\to\RR$, where $\sF := \sE(\cdot+\varphi)$ and $\sU_0 := \sU - \varphi \subset \sX$.

For the definition of the operators \cite[Equation (2.4.2a) and (2.4.3a)]{Huang_2006},
\[
G := \sM_1(0) \in \sL(\sH_\sA,\sH) \quad\text{and}\quad L := \sM(0) \in \sL(\sX,\tilde\sX),
\]
the proof of \cite[Proposition 2.4.1]{Huang_2006} requires only that $\sM_1(u) \in \sL(\sH_\sA,\sH)$ and $\sM(u) \in \sL(\sX,\tilde\sX)$ be Fredholm operators with index zero when $u = 0$. Following \cite{Huang_2006}, we denote $N(L) := \Ker(L:\sX\to\tilde\sX)$ and $R(L) := \{Lu:\forall\, u \in \sX\}$ together with $N(G) := \Ker(G:\sH_\sA\to\sH)$ and $R(G) := \{Gh:\forall\, h \in \sH\}$. As in \cite[Equation (2.4.4a)]{Huang_2006}, we have
\[
\sX = N(L) \oplus (\sX\cap R(G)) \quad\text{and}\quad \tilde\sX = N(L) \oplus R(L),
\]
direct sums of Banach subspaces of $\sX \subset \sH$ and $\tilde\sX \subset \sH$, respectively, that are $\sH$-orthogonal. Also as in \cite[Equation (2.4.4c)]{Huang_2006}, we define the finite-rank projection operator, $\Pi: \sX \to \tilde\sX$, by
\begin{equation}
\label{eq:Pi_X}
\Pi u
:=
\begin{cases}
u, &\quad\forall\, u \in N(L),
\\
0, &\quad\forall\, u \in \sX\cap R(G),
\end{cases}
\end{equation}
The map defined in \cite[Equation (2.4.5a)]{Huang_2006},
\[
\sN \equiv \Pi + \sM:\sU \to \tilde\sX,
\]
is real analytic by Hypothesis \ref{hyp:Huang_2-4_H1_H2_H3} and $\sN'(0) = \Pi + \sM'(0):\sX\to\tilde\sX$ is an isomorphism of Banach spaces. As in \cite[p. 36]{Huang_2006}, the Inverse Function Theorem for real analytic maps of Banach spaces (see Section \ref{subsubsec:Smooth_and_analytic_inverse_and_implicit function theorems}) yields convex neighborhoods of the origin, $B_1 \subset \sX$ and $B_2 \subset \tilde\sX$, and a real analytic map, $\Psi:B_2\to B_1$, such that \cite[Equation (2.4.5b)]{Huang_2006} holds, namely
\begin{align*}
\Psi(\sN(u)) &= u, \quad\forall\, u \in B_1,
\\
\sN(\Psi(v)) &= v, \quad\forall\, v \in B_2.
\end{align*}
The remainder of the proof of \cite[Theorem 2.4.2 (i)]{Huang_2006} remain unchanged. This completes the proof of Corollary \ref{cor:Huang_2-4-2_generalized}.
\end{proof}

\section{{\L}ojasiewicz-Simon gradient inequality for the Yang-Mills energy functional}
\label{sec:Huang_2_yangmills}
In this section, we shall apply the {\L}ojasiewicz-Simon gradient inequalities provided by Theorems \ref{thm:Huang_2-4-2} and \ref{thm:Huang_2-4-5} for an analytic potential function on an open subset of an abstract Banach space, $\sE:\sU\subset\sX \to \RR$, to the Yang-Mills energy functional,
$$
\sE:\sB^*(P,g) \to \RR, \quad [A] \mapsto \frac{1}{2}\int_X |F_A|_g^2 \, d\vol_g,
$$
where $G$ is a compact Lie group and $P$ is a principal $G$-bundle over a closed, smooth manifold, $X$, of dimension $d \geq 2$ with Riemannian metric $g$, and $\sB^*(P,g)$ is the Banach manifold of $W^{k,q}$ connections, $A$, modulo the action of the Banach Lie group $\Aut P$ of $W^{k+1,q}$ automorphisms, $u$, of the bundle, $P$, where $(k+1)q>d$, so $W^{k+1,q}(X) \hookrightarrow C(X)$. To obtain our first version of the {\L}ojasiewicz-Simon gradient inequality for the Yang-Mills energy functional, Theorem \ref{thm:Rade_proposition_7-2_L2}, we apply Theorem \ref{thm:Huang_2-4-2} with $\sH = L^2(X;\Lambda^1\otimes\ad P)$ and $\sX = W^{2,p}_{A_1}(X;\Lambda^1\otimes\ad P) \cap \Ker d_{A_1}^*$, a slice for the action of the Banach Lie group $\Aut P$ of $W^{3,p}$ bundle automorphisms, where $d \geq 2$ and $p \in [2,\infty)$ obeys $p > d/3$ and $A_1$ is a $C^\infty$ reference connection on $P$.

To derive our second version of the {\L}ojasiewicz-Simon gradient inequality for the Yang-Mills energy functional, Theorem \ref{thm:Rade_proposition_7-2}, we apply Theorem \ref{thm:Huang_2-4-5}, which allows us to use a much weaker system of norms on spaces of connections when $2\leq d \leq 4$, as we can then choose $\sX = H^1_{A_1}(X;\Lambda^1\otimes\ad P) \cap \Ker d_{A_1}^*$, so $\sX' = H^{-1}_{A_1}(X;\Lambda^1\otimes\ad P) \cap \Ker d_{A_1}^*$.

Consequently, by appealing to Theorem \ref{thm:Rade_proposition_7-2} one finds that the {\L}ojasiewicz-Simon gradient inequality for the Yang-Mills energy functional, when $d=4$, holds for connections $A$ in a neighborhood of a Yang-Mills connection of the form,
$$
\|A - A_\infty\|_{H^1_{A_1}(X)} < \sigma,
$$
rather than,
$$
\|A - A_\infty\|_{H^2_{A_1}(X)} < \sigma,
$$
resulting from an application of Theorem \ref{thm:Rade_proposition_7-2_L2}. If $A$ already obeys $d_{A_\infty}^*(A-A_\infty) = 0$ and $A_\infty$ is of class $C^\infty$, one can choose $A$ to be of class $H^1$ in Theorem \ref{thm:Rade_proposition_7-2} rather than of class $H^2$ as would be required by Theorem \ref{thm:Rade_proposition_7-2_L2}.

The {\L}ojasiewicz-Simon gradient inequality was first introduced in Yang-Mills gauge theory by Morgan, Mrowka, and Ruberman in their adaptation of the proof of \cite[Theorem 3]{Simon_1983} (see \cite[Proposition 4.2.1 and Corollary 4.2.3]{MMR}) to the \emph{Chern-Simons functional} \cite[Section 2.1]{MMR},
\begin{equation}
\label{eq:Chern_simons_functional}
\CS(A) := \int_X \tr\left(2a\wedge F_\Gamma + a\wedge d_\Gamma a + \frac{2}{3}a\wedge a\wedge a\right),
\end{equation}
for connections on a principal $G$-bundle, $P$, over a closed, Riemannian, smooth, three-manifold, $X$, where $A = \Gamma + a$ and $\Gamma$ is a fixed $C^\infty$ reference connection on $P$, so $a \in \Omega^1(X;\ad P)$. The critical points of the Chern-Simons functional are \emph{flat connections} on $P$. It is important to note that in the statement of the original of the {\L}ojasiewicz-Simon gradient inequality due to Leon Simon (see \cite[Theorem 3]{Simon_1983}), he takes $\sH = L^2(X;V)$ and $\sX = C^{2,\alpha}(X;V)$ (for $\alpha \in (0,1)$), where $V$ is a Riemannian vector bundle over a closed, Riemannian, smooth manifold, $X$, of dimension $d \geq 2$. Hence, there is slightly more work involved in the translation of the proof of \cite[Theorem 3]{Simon_1983} than perhaps might be evident in the short argument in \cite[pp. 62-63]{MMR} establishing \cite[Proposition 4.2.1]{MMR}, where $\sH = L^2(X;\Lambda^1\otimes\ad P)$ and $\sX = W^{2,k}_\Gamma(X;\Lambda^1\otimes\ad P)$ (for suitable $k \geq 1$) because, for example, the analyticity of the functional \eqref{eq:Chern_simons_functional} depends on the choice of $\sX$.

In his doctoral dissertation \cite{Rade_1992, RadeThesis}, R\r{a}de also adapted the proof of the {\L}ojasiewicz-Simon gradient inequality in \cite[Theorem 3]{Simon_1983} to derive gradient inequalities \cite[Proposition 7.2 and Equation (9.1)]{Rade_1992} for the Yang-Mills energy functional, yielding Theorem \ref{thm:Rade_proposition_7-2} when $X$ has dimension $d=2, 3$, with $\sH = L^2(X;\Lambda^1\otimes\ad P)$ and $\sX = H^1_{A_1}(X;\Lambda^1\otimes\ad P)$.

Wilkin adapted the proof of the {\L}ojasiewicz-Simon gradient inequality in \cite[Theorem 3]{Simon_1983} to derive gradient inequalities for the Yang-Mills-Higgs functional over a Riemann surface \cite[Proposition 3.5 and Theorem 3.19]{Wilkin_2008}, with choices of $\sH$ and $\sX$ for Sobolev spaces of pairs analogous to those of R\r{a}de in \cite[Proposition 7.2 and Equation (9.1)]{Rade_1992} for Sobolev spaces of connections. B. Yang established a version of the {\L}ojasiewicz-Simon gradient inequality for the Yang-Mills energy functional on a principal $G$-bundle over a four-dimensional manifold \cite[Lemma 12]{Yang_2003aim}. However, Yang's result is a direct analogue of \cite[Theorem 3]{Simon_1983} in that $\sH = L^2(X;\Lambda^1\otimes\ad P)$ and $\sX = C^{2,\alpha}(X;\Lambda^1\otimes\ad P)$ and so cannot be used in our application to establish global existence of smooth solutions to the gradient flow for the Yang-Mills energy functional.

\subsection{A $W^{2,p}$ {\L}ojasiewicz-Simon gradient inequality for the Yang-Mills energy functional}
\label{subsec:W2p_Lojasiewicz-Simon_gradient_inequality_Yang-Mills}
Suppose $G$ is a compact Lie group and that $X$ is a closed, Riemannian, smooth manifold of dimension $d \geq 2$. According to \cite[Corollary 1.4]{UhlLp}, when $p\geq d/2$, we may choose local trivializations of a principal $G$-bundle $P$ such that the local representatives of a $W^{2,p}$ Yang-Mills connection, $A_\infty$, are $C^\infty$, that is, $A_\infty$ is a connection of class $C^\infty$ on $P$. Thus, when $2p \geq d$ in Theorems \ref{thm:Rade_proposition_7-2_L2} or \ref{thm:Rade_proposition_7-2}, we may choose $A_1=A_\infty$ as the $C^\infty$ reference connection.

We use the $C^\infty$ reference connection, $A_1$, on $P$ to define Sobolev spaces, $W^{s,p}_{A_1}(X;\Lambda^1\otimes\ad P)$, and their norms, $\|\cdot\|_{W^{s,p}_{A_1}(X)}$. When no confusion can arise we may suppress $A_1$ from our notation, but we caution the reader that Sobolev embedding and multiplication constants will generally depend on $A_1$, though this dependence may not always be explicitly noted.

As a consequence of Theorem \ref{thm:Huang_2-4-2}, we obtain our first version of \cite[Proposition 7.2]{Rade_1992}, without the restriction on the dimension, $d$, of $X$ that $d=2,3$ but with the $H^1$-norm replaced by the $W^{2,p}$-norm.

\begin{thm}[A $W^{2,p}$ {\L}ojasiewicz-Simon gradient inequality for the Yang-Mills energy functional and $L^2$ gradient norm]
\label{thm:Rade_proposition_7-2_L2}
Let $X$ be a closed, Riemannian, smooth manifold of dimension $d$, let $G$ be a compact Lie group, and assume that
\begin{enumerate}
\item $d\geq 2$ and $p \in [2,\infty)$ obeys $p>d/2$, or
\item $2 \leq d \leq 5$ and $p = 2$.
\end{enumerate}
If $A_1$ is a connection of class $C^\infty$ and $A_\infty$ is a Yang-Mills connection of class $W^{2,p}$ on a principal $G$-bundle $P$ over $X$, then there are constants, $c\in (0,\infty)$ and $\sigma \in (0,1]$ and $\theta \in [1/2,1)$ depending on $A_1$, $A_\infty$, $g$, $G$, and $p$, with the following significance.  If $A$ is a connection of class $W^{2,p}$ on $P$ and
\begin{equation}
\label{eq:Rade_7-1_L2neighborhood}
\|A - A_\infty\|_{W^{2,p}_{A_1}(X)} < \sigma,
\end{equation}
then
\begin{equation}
\label{eq:Rade_7-1_L2}
\|d_A^*F_A\|_{L^2(X)} \geq c|\sE(A) - \sE(A_\infty)|^\theta,
\end{equation}
where
\begin{equation}
\label{eq:Potential_yang_mills}
\sE(A) := \frac{1}{2}\|F_A\|_{L^2(X)}^2.
\end{equation}
\end{thm}

\begin{proof}
The inequality \eqref{eq:Rade_7-1_L2} is gauge-invariant, that is, if $u \in \Aut P$ is a gauge transformation of class $W^{3,p}$ (and thus continuous under our hypotheses on $d$ and $p$), then \eqref{eq:Rade_7-1_L2} is equivalent to
$$
\|d_{u(A)}^*F_{u(A)}\|_{L^2(X)} \geq c|\sE(u(A)) - \sE(u(A_\infty))|^\theta.
$$
When $d=4$ and $k\geq 2$, the Slice Theorem\footnote{More refined versions of the Slice Theorem in dimension four are provided by \cite[Theorem 1.1]{FeehanSlice}, where the neighborhood of the reference connection, $A_1$, is shown to depend at most on $\|F_{A_1}\|_{L^2(X)}$ and the least positive eigenvalue of the Laplace operator $d_{A_1}^*d_{A_1}$ on $C^\infty(X;\ad P)$, and in arbitrary dimensions by \cite[Theorem 10]{Feehan_Maridakis_Lojasiewicz-Simon_coupled_Yang-Mills}.}
for the quotient space of $W^{k,2}$ connections modulo $W^{k+1,2}$ gauge transformations \cite[Proposition 4.2.9]{DK}, \cite[Theorem 3.2]{FU}, \cite[Theorem 2.10.4]{Lawson} yields constants $\eps_0 = \eps_0(A_1,A_\infty,g,G,k) \in (0,1]$ and $C_0 = C_0(A_1,A_\infty,g,G,k) \in [1,\infty)$ and a $W^{k+1,2}$ gauge transformation $u \in \Aut P$ such that
\begin{gather*}
d_{A_\infty}^*(u(A) - A_1) = 0 \quad\hbox{if } \|A - A_\infty\|_{W^{k,2}_{A_1}(X)} < \eps_0,
\\
\|u(A) - A_\infty\|_{W^{k,2}_{A_1}(X)} \leq C_0\|A -A_\infty\|_{W^{k,2}_{A_1}(X)}.
\end{gather*}
The proof of the preceding gauge-fixing result is a consequence of the Implicit Function Theorem for maps of Banach manifolds and depends on the dimension, $d$, of the base manifold, $X$, only insofar as that the gauge transformations $u \in \Aut P$ should be continuous; in particular, it continues to hold under the hypotheses of Theorem \ref{thm:Rade_proposition_7-2_L2} on $d$ and $p$ for $W^{2,p}$ connections and $W^{3,p}$ gauge transformations to give, for constants $\eps_0 = \eps_0(A_1,A_\infty,g,G,p) \in (0,1]$ and $C_0 = C_0(A_1,A_\infty,g,G,p) \in [1,\infty)$, a $W^{3,p}$ gauge transformation $u \in \Aut P$ such that
\begin{gather}
\label{eq:Coulomb_gauge_L2}
d_{A_1}^*(u(A) - A_\infty) = 0 \quad\hbox{if } \|A - A_\infty\|_{W^{2,p}_{A_1}(X)} < \eps_0,
\\
\label{eq:Coulomb_gauge_transformed_connection_near_original_L2}
\|u(A) - A_\infty\|_{W^{2,p}_{A_1}(X)} \leq C_0\|A - A_\infty\|_{W^{2,p}_{A_1}(X)}.
\end{gather}
We divide the remainder of the proof of Theorem \ref{thm:Rade_proposition_7-2_L2} into the two cases suggested by the hypotheses.

\setcounter{case}{0}
\begin{case}[$d\geq 2$ and $p>d/2$]
\label{case:p_greaterthan_dover2}
For this case, the embedding $W^{2,p}(X) \hookrightarrow C(X)$ is continuous and $W^{2,p}(X)$ is a Banach algebra \cite[Theorem 4.39]{AdamsFournier}, \cite[p. 96]{FU}; in particular, the following embedding is continuous,
\begin{equation}
\label{eq:Sobolev_embedding_W2p_into_C}
W^{2,p}_{A_1}(X;\Lambda^1\otimes\ad P) \hookrightarrow C(X;\Lambda^1\otimes\ad P), \quad\hbox{for } d \geq 2 \hbox{ and } p > d/2.
\end{equation}
We assume first that $A_\infty$ is $C^\infty$ and remove that restriction in the forthcoming Lemma \ref{lem:Proof_Rade_7-1_L2_when_A_Coulomb_p>dover2} at the end of our discussion of Case \ref{case:p_greaterthan_dover2}. We shall apply Corollary \ref{cor:Huang_2-4-2_generalized} with the linear, positive, self-adjoint operator,
$$
\sA := \Delta_{A_\infty} + 1 = d_{A_\infty}^*d_{A_\infty} + d_{A_\infty}d_{A_\infty}^* + 1: \sD(\sA) \subset \sH \to \sH,
$$
where
\begin{gather*}
\sD(\sA) \equiv \sH_\sA \cong \Ker d_{A_\infty}^*\cap H^2_{A_1}(X;\Lambda^1\otimes\ad P),
\\
\sH := \Ker d_{A_\infty}^*\cap L^2(X;\Lambda^1\otimes\ad P),
\end{gather*}
with Banach spaces and continuous embeddings (when $p\geq 2$),
\begin{gather*}
\sX :=  \Ker d_{A_\infty}^*\cap W^{2,p}_{A_1}(X;\Lambda^1\otimes\ad P) \hookrightarrow  \sH_\sA,
\\
\tilde\sX :=  \Ker d_{A_\infty}^*\cap L^p(X;\Lambda^1\otimes\ad P) \hookrightarrow  \sH.
\end{gather*}
We choose $\sU = \sX$ and define $\sE:A_\infty+\sX \to \RR$ according to \eqref{eq:Potential_yang_mills}. For all $a, b \in W^{2,p}_{A_1}(X;\Lambda^1\otimes\ad P)$, the differential, Hessian, pre-gradient, and pre-Hessian operators are given by
\begin{align}
\label{eq:Pre_gradient}
\sE'(A)a &= (F_A, d_Aa)_{L^2(X)} = (d_A^*F_A, a)_{L^2(X)} =: (\tilde\sM(A), a)_{L^2(X)},
\\
\label{eq:Pre_Hessian}
\sE''(A)(a,b) &= (d_Aa, d_Ab)_{L^2(X)} +
(F_A, a \wedge b + b \wedge a)_{L^2(X)} =: (\tilde\sM'(A)a, b)_{L^2(X)}.
\end{align}
For $k\in\ZZ$ and $p \in (1,\infty)$, denote the $L^2$-orthogonal projection onto the Coulomb-gauge slice through $A_\infty$ by
\begin{equation}
\label{eq:L2-orthogonal_projection_onto_slice_Wkp}
\Pi_{A_\infty}:W_{A_1}^{k,p}(X;\Lambda^1\otimes \ad P) \to \Ker d_{A_\infty}^* \cap W_{A_1}^{k,p}(X;\Lambda^1\otimes \ad P).
\end{equation}
Let $G_{A_\infty}$ denote the Green's operator for $d_{A_\infty}^*d_{A_\infty}:C^\infty(X;\ad P)\to C^\infty(X;\ad P)$. Because
\begin{align*}
W_{A_1}^{k,p}(X;\Lambda^1\otimes \ad P)
&= \Ker d_{A_\infty}^* \cap W_{A_1}^{k,p}(X;\Lambda^1\otimes \ad P)
\\
&\quad \oplus \Ran\left(d_{A_\infty}:W_{A_1}^{k+1,p}(X;\ad P)
\to W_{A_1}^{k,p}(X;\Lambda^1\otimes \ad P) \right)
\end{align*}
is an $L^2$-orthogonal direct sum, we see that
\begin{equation}
\label{eq:L2-orthogonal_projection_onto_slice_intermsof_Greens_operator}
\Pi_{A_\infty} = \id - d_{A_\infty}G_{A_\infty}d_{A_\infty}^*.
\end{equation}

The corresponding gradient and Hessian maps,
\begin{gather}
\label{eq:Gradient_map_yang_mills_A1_plus_W2p_cap_kerdA1star_to_Lp}
\begin{split}
\sM:&A_\infty + \Ker d_{A_\infty}^*\cap W^{2,p}_{A_1}(X;\Lambda^1\otimes\ad P) \ni A
\\
&\quad \mapsto \sM(A) \in \Ker d_{A_\infty}^*\cap L^p(X;\Lambda^1\otimes\ad P),
\end{split}
\\
\label{eq:Hessian_map_yang_mills_A1_plus_W2p_cap_kerdA1star_to_operators_from_W2p_cap_kerdA1star_to_Lp}
\begin{split}
\sM':&A_\infty + \Ker d_{A_\infty}^*\cap W^{2,p}_{A_1}(X;\Lambda^1\otimes\ad P) \ni A
\\
&\quad \mapsto \sM'(A) \in \sL\left(\Ker d_{A_\infty}^*\cap W^{2,p}_{A_1}(X;\Lambda^1\otimes\ad P), \Ker d_{A_\infty}^*\cap L^p(X;\Lambda^1\otimes\ad P)\right),
\end{split}
\end{gather}
are given formally by
\begin{align}
\label{eq:Gradient_map_and_potential_yang_mills}
\sE'(A)a &= (\sM(A),a)_{L^2(X)}, \quad\forall\, a \in \Ker d_{A_\infty}^*\cap W_{A_1}^{2,p}(X;\Lambda^1\otimes\ad P),
\\
\label{eq:Hessian_map_and_potential_yang_mills}
\sE''(A)(a,b) &= (\sM'(A)a,b)_{L^2(X)}, \quad\forall\, a,b \in \Ker d_{A_\infty}^*\cap W_{A_1}^{2,p}(X;\Lambda^1\otimes\ad P),
\end{align}
we obtain, for all $A \in A_\infty + \Ker d_{A_\infty}^*\cap W^{2,p}_{A_1}(X;\Lambda^1\otimes\ad P)$,
\begin{align}
\label{eq:Gradient_map_yang_mills_short_expression}
\sM(A) &= \Pi_{A_\infty}\tilde\sM(A) = \Pi_{A_\infty}d_A^*F_A,
\\
\label{eq:Hessian_map_yang_mills_short_expression}
\begin{split}
\sM'(A)b &= \Pi_{A_\infty}\tilde\sM'(A) = \Pi_{A_\infty}(d_A^*d_Ab + F_A\times b),
\\
&\qquad \forall\, b \in \Ker d_{A_\infty}^*\cap W^{2,p}_{A_1}(X;\Lambda^1\otimes\ad P).
\end{split}
\end{align}
We recall from R\r{a}de \cite[p. 148]{Rade_1992} that
\begin{align*}
d_A^*\tilde\sM(A) &= 0,
\\
d_A^*\tilde\sM'(A)a &= 0, \quad\forall\, a \in W^{2,p}_{A_1}(X;\Lambda^1\otimes\ad P),
\\
\tilde\sM'(A)d_A\xi &= 0, \quad\forall\, \xi \in W^{3,p}_{A_1}(X;\ad P),
\end{align*}
as a consequence of the fact that $\sE(u(A)) = \sE(A)$, for all $W^{2,p}$ connections $A$ and $W^{3,p}$ gauge transformations $u \in \Aut P$. In particular, we have
\[
\sM'(A_\infty) = \Pi_{A_\infty}\tilde\sM'(A_\infty) = \tilde\sM'(A_\infty).
\]
We now use the formal expressions \eqref{eq:Gradient_map_yang_mills_short_expression} and \eqref{eq:Hessian_map_yang_mills_short_expression} to verify that the maps \eqref{eq:Gradient_map_yang_mills_A1_plus_W2p_cap_kerdA1star_to_Lp} and \eqref{eq:Hessian_map_yang_mills_A1_plus_W2p_cap_kerdA1star_to_operators_from_W2p_cap_kerdA1star_to_Lp} have the desired properties for suitable $p$. By defining
\begin{equation}
\label{eq:Laplacian_omega1adP}
\Delta_A := d_A^*d_A + d_Ad_A^*,
\end{equation}
and writing $A = A_1 + a$ to give formal expressions,
\begin{align*}
F_A &= F_{A_1} + d_{A_1}a + a\times a, \quad \forall\, a \in W^{2,p}_{A_1}(X;\Lambda^1\otimes\ad P),
\\
d_A^*d_Ab &= d_{A_1}^*d_{A_1}b + a\times \nabla_{A_1}b + \nabla_{A_1}a\times b + a\times a \times b,
\quad \forall\, b \in W^{2,p}_{A_1}(X;\Lambda^1\otimes\ad P).
\end{align*}
If $b \in \Ker d_{A_\infty}^*$, then $\Delta_{A_\infty}b = d_A^*d_Ab$ and this motivates the formal expression,
\begin{equation}
\label{eq:Hessian_map_yang_mills_long_expression}
\begin{gathered}
\tilde\sM'(A)b = \Delta_{A_\infty}b + (\Delta_{A_1}-\Delta_{A_\infty})b + a\times \nabla_{A_1}b + F_{A_1}\times b + \nabla_{A_1}a\times b + a\times a \times b,
\\
\quad\forall\, b \in  W^{2,p}_{A_1}(X;\Lambda^1\otimes\ad P).
\end{gathered}
\end{equation}
Clearly, $\Delta_{A_\infty}$ defines a bounded linear operator,
\begin{equation}
\label{eq:Laplacian_yang_mills_W2p_to_Lp}
\Delta_{A_\infty}: W_{A_1}^{2,p}(X;\Lambda^1\otimes\ad P)
\to L^p(X;\Lambda^1\otimes\ad P),
\end{equation}
while $\Delta_{A_1}-\Delta_{A_\infty}$ is a first-order partial differential operator with $C^\infty$ coefficients and thus defines a bounded operator,
\[
\Delta_{A_1}-\Delta_{A_\infty}: W_{A_1}^{2,p}(X;\Lambda^1\otimes\ad P)
\to  W_{A_1}^{1,p}(X;\Lambda^1\otimes\ad P).
\]
We can estimate the remaining terms in the expression \eqref{eq:Hessian_map_yang_mills_long_expression} by
\begin{equation}
\label{eq:Hessian_map_minus_laplacian_yang_mills_Lp_estimates}
\begin{aligned}
\|a\times \nabla_{A_1}b\|_{L^p(X)} &\leq \|a\|_{C(X)}\|\nabla_{A_1}b\|_{L^p(X)}
\\
&\leq C\|a\|_{W^{2,p}_{A_1}(X)}\|b\|_{W^{1,p}_{A_1}(X)} \quad\hbox{(by \eqref{eq:Sobolev_embedding_W2p_into_C})},
\\
\|F_{A_1}\times b\|_{L^p(X)} &\leq \|F_{A_1}\|_{L^p(X)}\|b\|_{C(X)}
\\
&\leq C\|F_{A_1}\|_{L^p(X)}\|b\|_{W^{2,p}_{A_1}(X)} \quad\hbox{(by \eqref{eq:Sobolev_embedding_W2p_into_C})},
\\
\|\nabla_{A_1}a\times b\|_{L^p(X)} &\leq \|\nabla_{A_1}a\|_{L^p(X)}\|b\|_{C(X)}
\\
&\leq C\|a\|_{W^{1,p}_{A_1}(X)}\|b\|_{W^{2,p}_{A_1}(X)} \quad\hbox{(by \eqref{eq:Sobolev_embedding_W2p_into_C})},
\\
\|a\times a\times b\|_{L^p(X)} &\leq \|a\|_{C(X)}^2\|b\|_{L^p(X)}
\\
&\leq C\|a\|_{W^{2,p}_{A_1}(X)}^2\|b\|_{W^{2,p}_{A_1}(X)} \quad\hbox{(by \eqref{eq:Sobolev_embedding_W2p_into_C})}.
\end{aligned}
\end{equation}
Hence, the formal expression \eqref{eq:Hessian_map_yang_mills_short_expression} for $\tilde\sM'(A)$ defines a bounded linear operator,
\begin{equation}
\label{eq:Pre_Hessian_map_yang_mills_fixed_A_operator_W2p_cap_ker_dA1star_to_Lp}
\tilde\sM'(A): W_{A_1}^{2,p}(X;\Lambda^1\otimes\ad P) \to L^p(X;\Lambda^1\otimes\ad P).
\end{equation}
We now make the

\begin{claim}[Bounded operator and extension]
\label{claim:Hessian_map_yang_mills_H2_extension}
Assume the hypotheses of Theorem \ref{thm:Rade_proposition_7-2_L2} but require that $A_\infty$ be $C^\infty$ and $p > d/2$. Then the Hessian map operator,
\begin{equation}
\label{eq:Hessian_map_yang_mills_fixed_A_operator_W2p_cap_ker_dA1star_to_Lp}
\sM'(A):\Ker d_{A_\infty}^* \cap W_{A_1}^{2,p}(X;\Lambda^1\otimes\ad P)
\to \Ker d_{A_\infty}^* \cap L^p(X;\Lambda^1\otimes\ad P),
\end{equation}
is bounded and has an extension as a bounded, linear operator,
\begin{equation}
\label{eq:Hessian_map_yang_mills_fixed_A_operator_H2_cap_ker_dA1star_to_L2}
\sM_1(A): \Ker_{A_\infty}^* \cap H^2_{A_1}(X;\Lambda^1\otimes\ad P) \to \Ker_{A_\infty}^* \cap L^2(X;\Lambda^1\otimes\ad P).
\end{equation}
\end{claim}

For this purpose, we shall need the

\begin{lem}[Continuous Sobolev multiplication maps]
\label{lem:Sobolev_multiplication_W1p_times_H1_into_L2}
Let $(X,g)$ be a closed, Riemannian, smooth manifold of dimension $d \geq 2$. Then the following Sobolev multiplication maps are continuous:
\begin{align}
\label{eq:Sobolev_multiplication_W1p_times_H1_into_L2}
W^{1,p}(X) \times H^1(X) &\to L^2(X), \quad\hbox{if } d \geq 2 \hbox{ and } p \geq d/2,
\\
\label{eq:Sobolev_multiplication_H1_times_H1_into_L2}
H^1(X) \times H^1(X) &\to L^2(X), \quad\hbox{if } 2 \leq d \leq 4.
\end{align}
\end{lem}

\begin{proof}[Proof of Claim \ref{claim:Hessian_map_yang_mills_H2_extension}]
Boundedness of the operator $\sM'(A) = \Pi_{A_\infty}\tilde\sM'(A)$ in \eqref{eq:Hessian_map_yang_mills_fixed_A_operator_W2p_cap_ker_dA1star_to_Lp} follows from boundedness of the operator $\tilde\sM'(A)$ in \eqref{eq:Pre_Hessian_map_yang_mills_fixed_A_operator_W2p_cap_ker_dA1star_to_Lp} and boundedness of the $L^2$-orthogonal projection $\Pi_{A_\infty}$ in \eqref{eq:L2-orthogonal_projection_onto_slice_Wkp}.

The Laplace operator, $\Delta_{A_\infty}:H^2_{A_1}(X;\Lambda^1\otimes\ad P)\to L^2(X;\Lambda^1\otimes\ad P)$, is clearly bounded, while $\Delta_{A_1}-\Delta_{A_\infty}$ defines a bounded operator,
\[
\Delta_{A_1}-\Delta_{A_\infty}: H_{A_1}^2(X;\Lambda^1\otimes\ad P)
\to H_{A_1}^1(X;\Lambda^1\otimes\ad P).
\]
Furthermore, the following operator is bounded,
\begin{equation}
\label{eq:Hessian_map_minus_laplacian_from_H1_cap_ker_dA1star_to_L2}
\tilde\sM'(A)-\Delta_1:H^1_{A_1}(X;\Lambda^1\otimes\ad P) \to \Ker d_{A_\infty}^*\cap L^2(X;\Lambda^1\otimes\ad P),
\end{equation}
because, recalling the expression \eqref{eq:Hessian_map_yang_mills_long_expression} for $\tilde\sM'(A)$, we have
\begin{equation}
\label{eq:Hessian_map_minus_laplacian_yang_mills_L2_estimates}
\begin{aligned}
\|a\times \nabla_{A_1}b\|_{L^2(X)} &\leq \|a\|_{C(X)}\|\nabla_{A_1}b\|_{L^2(X)}
\\
&\leq C\|a\|_{W^{2,p}_{A_1}(X)}\|b\|_{H^1_{A_1}(X)} \quad\hbox{(by \eqref{eq:Sobolev_embedding_W2p_into_C})},
\\
\|F_{A_1}\times b\|_{L^2(X)} &\leq \|F_{A_1}\|_{C(X)}\|b\|_{L^2(X)},
\\
\|\nabla_{A_1}a\times b\|_{L^2(X)} &\leq C\|\nabla_{A_1}a\|_{W^{1,p}(X)}\|b\|_{H^1_{A_1}(X)}
\quad\hbox{(by \eqref{eq:Sobolev_multiplication_W1p_times_H1_into_L2})}
\\
&\leq C\|a\|_{W^{2,p}_{A_1}(X)}\|b\|_{H^1_{A_1}(X)},
\\
\|a\times a\times b\|_{L^2(X)} &\leq \|a\|_{C(X)}^2\|b\|_{L^2(X)}
\\
&\leq C\|a\|_{W^{2,p}_{A_1}(X)}^2\|b\|_{L^2(X)} \quad\hbox{(by \eqref{eq:Sobolev_embedding_W2p_into_C})}.
\end{aligned}
\end{equation}
Boundedness of the operator $\sM_1(A)$ in \eqref{eq:Hessian_map_yang_mills_fixed_A_operator_H2_cap_ker_dA1star_to_L2} follows from boundedness of the operator $\tilde\sM'(A)$ in \eqref{eq:Pre_Hessian_map_yang_mills_fixed_A_operator_W2p_cap_ker_dA1star_to_Lp} (with $p=2$) and boundedness of the $L^2$-orthogonal projection $\Pi_{A_\infty}$ in \eqref{eq:L2-orthogonal_projection_onto_slice_Wkp} (with $p=2$).
This completes the proof of Claim \ref{claim:Hessian_map_yang_mills_H2_extension}.
\end{proof}

We now turn to the

\begin{proof}[Proof of Lemma \ref{lem:Sobolev_multiplication_W1p_times_H1_into_L2}]
We recall by \cite[Equation (6.34)]{FU} that \eqref{eq:Sobolev_multiplication_W1p_times_H1_into_L2} is continuous when
$$
(1-d/p) + (s - d/2) \geq -d/2 \quad \hbox{and}\quad s \geq 0 \quad \hbox{and}\quad \max\{2s, p\} < d,
$$
that is, when $p \geq d/(1+s)$ and $s \geq 0$ and $\max\{2s, p\} < d$. For $d>2$, we can take $s=1$ and this yields \eqref{eq:Sobolev_multiplication_W1p_times_H1_into_L2} when $p\geq d/2$.

H\"older's inequality yields a continuous map $L^4(X) \times L^4(X) \to L^2(X)$ and so, for $2\leq d \leq 4$, there are continuous embeddings
\begin{inparaenum}[\itshape a\upshape)]
\item $H^1(X) \hookrightarrow L^q(X)$ for $d=2$ and $2 \leq q < \infty$ by \cite[Theorem 4.12, Part I (B)]{AdamsFournier},
\item $H^1(X) \hookrightarrow L^6(X)$ for $d=3$ by \cite[Theorem 4.12, Part I (C)]{AdamsFournier}, and
\item $H^1(X) \hookrightarrow L^4(X)$ for $d=4$ by \cite[Theorem 4.12, Part I (C)]{AdamsFournier}.
\end{inparaenum}
In particular, there is a continuous embedding $H^1(X) \hookrightarrow L^4(X)$, so \eqref{eq:Sobolev_multiplication_H1_times_H1_into_L2} holds when $2\leq d \leq 4$. This completes the proof of Lemma \ref{lem:Sobolev_multiplication_W1p_times_H1_into_L2}.
\end{proof}

\begin{prop}[Continuity of the extended Hessian operator map for the Yang-Mills energy functional on $W^{2,p}$ when $p>d/2$]
\label{prop:Continuity_extended_Hessian_operator_map_Yang-Mills_energy_W2p_and_p>dover2}
Assume the hypotheses of Theorem \ref{thm:Rade_proposition_7-2_L2} but require that $A_\infty$ is $C^\infty$ and $p>d/2$. Then the following map is continuous,
\begin{multline}
\label{eq:Hessian_map_yang_mills_from_A1_plus_W2p_to_extended_operators_from_H2_to_L2}
A_\infty + \Ker d_{A_\infty}^*\cap W^{2,p}_{A_1}(X;\Lambda^1\otimes\ad P) \ni A
\\
\mapsto \sM_1(A) \in \sL\left(\Ker d_{A_\infty}^*\cap H^2_{A_1}(X;\Lambda^1\otimes\ad P), \Ker d_{A_\infty}^*\cap L^2(X;\Lambda^1\otimes\ad P)\right).
\end{multline}
\end{prop}

\begin{proof}
The fact that each linear operator,
\[
\tilde\sM'(A):H^2_{A_1}(X;\Lambda^1\otimes\ad P) \to L^2(X;\Lambda^1\otimes\ad P),
\]
is a bounded follows from Claim \ref{claim:Hessian_map_yang_mills_H2_extension}.

By \eqref{eq:Hessian_map_yang_mills_long_expression} and writing $A = A_1 + a_1$, we see that for all $a, a_1 \in W^{2,p}_{A_1}(X;\Lambda^1\otimes\ad P)$ and $b \in H^2(X;\Lambda^1\otimes\ad P)$ we have
$$
\tilde\sM'(A+a)b = \Delta_{A_1} b + (a_1+a)\times \nabla_{A_1} b + F_{A_1}\times b + \nabla_{A_1}(a_1+a)\times b
+ (a_1+a)\times (a_1+a) \times b,
$$
which gives
\begin{equation}
\label{eq:Hessian_map_yang_mills_at_Aplusa_minus_at_A_on_b}
\tilde\sM'(A+a)b - \tilde\sM'(A)b = a\times \nabla_{A_1} b + \nabla_{A_1} a\times b + (a_1 + a) \times a \times b.
\end{equation}
Therefore,
\begin{align*}
{}&\|(\tilde\sM'(A+a) - \tilde\sM'(A))b\|_{L^2(X)}
\\
&\quad \leq \|a\times \nabla_{A_1} b\|_{L^2(X)} + \|\nabla_{A_1} a\times b\|_{L^2(X)}
+ \|(a_1+a)\times a \times b\|_{L^2(X)}
\\
&\quad \leq C\|a\|_{W^{1,p}_{A_1}(X)}\|\nabla_{A_1} b\|_{H^1_{A_1}(X)}
+ C\|\nabla_{A_1} a\|_{W^{1,p}_{A_1}(X)}\|b\|_{H^1_{A_1}(X)}
\\
&\qquad + \|a_1+a\|_{C(X)}\|a\|_{C(X)}\|b\|_{L^2(X)} \quad\hbox{(by \eqref{eq:Sobolev_multiplication_W1p_times_H1_into_L2})}
\\
&\quad \leq C\left(1 + \|a_1\|_{W^{2,p}_{A_1}(X)} + \|a\|_{W^{2,p}_{A_1}(X)}\right)\|a\|_{W^{2,p}_{A_1}(X)}\|b\|_{H^2_{A_1}(X)}
\quad\hbox{(by \eqref{eq:Sobolev_embedding_W2p_into_C})}.
\end{align*}
Hence, the pre-gradient map,
\[
\tilde\sM: A_\infty + W^{2,p}_{A_1}(X;\Lambda^1\otimes\ad P) \ni A
\mapsto \tilde\sM(A) \in L^p(X;\Lambda^1\otimes\ad P),
\]
is Lipschitz continuous. Boundedness of the $L^2$-orthogonal projection $\Pi_{A_\infty}$ in \eqref{eq:L2-orthogonal_projection_onto_slice_Wkp} implies that the map,
\[
A \mapsto \sM_1(A) = \Pi_{A_\infty}\tilde\sM'(A),
\]
in \eqref{eq:Hessian_map_yang_mills_from_A1_plus_W2p_to_extended_operators_from_H2_to_L2} is also Lipschitz continuous. This completes the proof of Proposition \ref{prop:Continuity_extended_Hessian_operator_map_Yang-Mills_energy_W2p_and_p>dover2}.
\end{proof}

Next, we verify that the gradient map is real analytic.

\begin{prop}[Real analyticity of the gradient map for the Yang-Mills energy functional on $W^{2,p}$ when $p>d/2$]
\label{prop:Real_analyticity_gradient_map_Yang-Mills_energy_W2p_and_p>dover2}
Assume the hypotheses of Theorem \ref{thm:Rade_proposition_7-2_L2} but require that $A_\infty$ is $C^\infty$ and $p>d/2$. Then the gradient map $\sM$ in \eqref{eq:Gradient_map_yang_mills_A1_plus_W2p_cap_kerdA1star_to_Lp} is real analytic.
\end{prop}

\begin{proof}
By the expression \eqref{eq:Pre_gradient} for $\tilde\sM(A)$, we see that
\begin{align*}
\tilde\sM(A_1+a) &= d_{A_1+a}^*F_{A_1+a}
\\
&= d_{A_1}^*d_{A_1}a + F_{A_1}\times a + \nabla_{A_1} a\times a + a\times a \times a + d_{A_1}^*F_{A_1},
\end{align*}
for all $a \in W^{2,p}_{A_1}(X;\Lambda^1\otimes\ad P)$, and so, restricting to $a \in W^{2,p}_{A_1}(X;\Lambda^1\otimes\ad P)\cap\Ker d_{A_1}^*$,
\begin{equation}
\label{eq:MA1_plus_a1_plus_a_minus_MA1_plus_a1}
\begin{aligned}
\tilde\sM(A_1+a_1+a) - \tilde\sM(A_1+a_1) &= d_{A_1}^*d_{A_1}a + F_{A_1}\times a
\\
&\quad + a\times \nabla_{A_1} a + \nabla_{A_1} a_1\times a + a_1\times \nabla_{A_1} a
\\
&\quad + a\times a \times a + a_1\times a \times a + a_1\times a_1 \times a.
\end{aligned}
\end{equation}
Note that it is immaterial whether we include the term $d_{A_1}d_{A_1}^*a$ or not when we restrict to $\Ker d_{A_1}^*$. Using the embedding \eqref{eq:Sobolev_embedding_W2p_into_C}, the fact that $W^{2,p}(X)$ is a Banach algebra when $p > d/2$, and writing $A = A_1+a_1$, we see that, for all $a_1, a \in W^{2,p}_{A_1}(X;\Lambda^1\otimes\ad P)$,
\begin{align*}
{}&\|\tilde\sM(A+a) - \tilde\sM(A)\|_{L^p(X)}
\\
&\quad = \|\tilde\sM(A_1+a_1+a) - \tilde\sM(A_1+a_1)\|_{L^p(X)}
\\
&\quad \leq \|d_{A_1}^*d_{A_1}a\|_{L^p(X)} + \|F_{A_1}\times a\|_{L^p(X)}
\\
&\qquad + \|\nabla_{A_1} a\times a\|_{L^p(X)} + \|\nabla_{A_1} a_1\times a\|_{L^p(X)} + \|a_1\times \nabla_{A_1} a\|_{L^p(X)}
\\
&\qquad + \|a\times a \times a\|_{L^p(X)} + \|a_1\times a \times a\|_{L^p(X)} + \|a_1\times a_1 \times a\|_{L^p(X)}
\\
&\quad \leq \|d_{A_1}^*d_{A_1}a\|_{L^p(X)} + \|F_{A_1}\|_{L^p(X)}\|a\|_{C(X)}
\\
&\qquad + \|\nabla_{A_1} a\|_{L^p(X)}\|a\|_{C(X)} + \|\nabla_{A_1} a_1\|_{L^p(X)}\|a\|_{C(X)}
+ \|a_1\|_{C(X)}\|\nabla_{A_1} a\|_{L^p(X)}
\\
&\qquad + \|a\|_{L^p(X)}\|a\|_{C(X)}^2 + \|a_1\|_{L^p(X)}\|a\|_{C(X)}^2 + \|a_1\|_{C(X)}^2\|a\|_{L^p(X)}
\\
&\quad \leq C\left(1 + \|F_{A_1}\|_{L^p(X)} + \|a_1\|_{W^{2,p}_{A_1}(X)} + \|a\|_{W^{2,p}_{A_1}(X)} + \|a_1\|_{W^{2,p}_{A_1}(X)}^2 + \|a\|_{W^{2,p}_{A_1}(X)}^2 \right)\|a\|_{W^{2,p}_{A_1}(X)}.
\end{align*}
Hence, the pre-gradient map,
\[
\tilde\sM: A_\infty + W^{2,p}_{A_1}(X;\Lambda^1\otimes\ad P) \ni A
\mapsto \tilde\sM(A) \in L^p(X;\Lambda^1\otimes\ad P),
\]
is Lipschitz continuous
and thus real analytic since the function $\tilde\sM(A_1+a)$ is a polynomial of degree three in $a \in W^{2,p}_{A_1}(X;\Lambda^1\otimes\ad P)$.

Boundedness of the $L^2$-orthogonal projection $\Pi_{A_\infty}$ in \eqref{eq:L2-orthogonal_projection_onto_slice_Wkp} and the expression \eqref{eq:Gradient_map_yang_mills_short_expression} for $\sM(A) = \Pi_{A_\infty}\tilde\sM(A)$ thus imply that the gradient map $\sM$ in \eqref{eq:Gradient_map_yang_mills_A1_plus_W2p_cap_kerdA1star_to_Lp} is real analytic. This completes the proof of Proposition \ref{prop:Real_analyticity_gradient_map_Yang-Mills_energy_W2p_and_p>dover2}.
\end{proof}

Finally, we establish the required Fredholm and index zero properties for the Hessian operator, $\sM'(A_\infty)$, and its extension, $\sM_1(A_\infty)$.

\begin{lem}[Fredholm and index zero properties of the Hessian operator for the Yang-Mills $L^2$-energy functional on a Coulomb-gauge slice]
\label{lem:fred_Lp}
Let $(X,g)$ be a closed, Riemannian, smooth manifold of dimension $d \geq 2$, and $G$ a compact Lie group, $P$ a smooth principal $G$-bundle over $X$, and $A_1$ a $C^\infty$ reference connection on $P$. If $A_\infty$ is a $C^\infty$ connection on $P$ and $1 < p < \infty$, then the following operator is Fredholm with index zero,
\[
\sM'(A_\infty): \Ker d_{A_\infty}^* \cap W_{A_1}^{2,p}(X; \Lambda^1 \otimes \ad P)
\to \Ker d_{A_\infty}^* \cap L^p(X; \Lambda^1 \otimes \ad P).
\]
\end{lem}

\begin{proof}
Lemma \ref{lem:Gilkey_1-4-5_Sobolev} implies that the operator,
\[
\sM'(A_\infty): W_{A_1}^{2,p}(X; \Lambda^1 \otimes \ad P) \to L^p(X; \Lambda^1 \otimes \ad P),
\]
is Fredholm with index zero. The argument of R\r{a}de \cite[p. 148]{Rade_1992} for the case $p=2$ (and $d=2,3$) now adapts to prove that the restriction of the preceding operator to a slice domain (with slice range) is also Fredholm with index zero. (See the proof of \cite[Proposition 3.6]{Feehan_Maridakis_Lojasiewicz-Simon_coupled_Yang-Mills} for details.)
\end{proof}

Just as with $\sM'(A)$ in \eqref{eq:Hessian_map_yang_mills_fixed_A_operator_W2p_cap_ker_dA1star_to_Lp}, the operator $\sM_1(A)$ in \eqref{eq:Hessian_map_yang_mills_fixed_A_operator_H2_cap_ker_dA1star_to_L2} is symmetric (in fact, $L^2$-self-adjoint).
The Fredholm and index zero properties for the extension, $\sM_1(A_\infty)$, of $\sM'(A_\infty)$,
\[
\sM_1(A_\infty): \Ker d_{A_\infty}^* \cap H_{A_1}^2(X; \Lambda^1 \otimes \ad P)
\to \Ker d_{A_\infty}^* \cap L^2(X; \Lambda^1 \otimes \ad P).
\]
are an immediate consequence of Lemma \ref{lem:fred_Lp} with $p=2$.

\begin{lem}[Completion of the proof of the gradient inequality \eqref{eq:Rade_7-1_L2} when $A_\infty$ is $C^\infty$ and $A$ is in Coulomb gauge relative to $A_\infty$]
\label{lem:Proof_Rade_7-1_L2_when_Ainfty_smooth_A_Coulomb_p>dover2}
Assume the hypotheses of Theorem \ref{thm:Rade_proposition_7-2_L2} but require that $A_\infty$ is $C^\infty$ and $p>d/2$ and $d_{A_\infty}^*(A-A_\infty) = 0$. Then the gradient inequality \eqref{eq:Rade_7-1_L2} holds.
\end{lem}

\begin{proof}
Proposition \ref{prop:Continuity_extended_Hessian_operator_map_Yang-Mills_energy_W2p_and_p>dover2} guarantees that the map, $A_\infty+\sX \ni A \mapsto \sM_1(A) \in \sL(\sH_\sA,\sH)$, is continuous. Proposition \ref{prop:Real_analyticity_gradient_map_Yang-Mills_energy_W2p_and_p>dover2}
ensures that the gradient map, $\sM:A_\infty+\sX \to \tilde\sX$, is real analytic. Lemma \ref{lem:fred_Lp} implies that $\sM'(A_\infty) \in \sL(\sX,\tilde\sX)$ and $\sM_1(A_\infty) \in \sL(\sH_\sA,\sH)$ are Fredholm with index zero. Hence, Corollary \ref{cor:Huang_2-4-2_generalized} yields the gradient inequality \eqref{eq:Rade_7-1_L2}.
\end{proof}

\begin{lem}[Completion of the proof of the gradient inequality \eqref{eq:Rade_7-1_L2} when $A$ is in Coulomb gauge relative to $A_\infty$]
\label{lem:Proof_Rade_7-1_L2_when_A_Coulomb_p>dover2}
Assume the hypotheses of Theorem \ref{thm:Rade_proposition_7-2_L2} but require that $p>d/2$ and $d_{A_\infty}^*(A-A_\infty) = 0$. Then the gradient inequality \eqref{eq:Rade_7-1_L2} holds.
\end{lem}

\begin{proof}
By hypothesis, $A_\infty$ is a $W^{2,p}$ connection that is a critical point for the functional $\sE$ in \eqref{eq:Yang-Mills_energy_functional}. When $p\geq d/2$, the regularity result \cite[Corollary 1.4]{UhlLp} implies that there exists a $W^{3,p}$ gauge transformation, $u_\infty \in \Aut P$, such that $u_\infty(A_\infty)$ is a $C^\infty$ connection. Moreover, $u_\infty(A)$ is in Coulomb gauge relative to $u_\infty(A_\infty)$ and (for example, see \cite[Lemma 3.10]{Feehan_Maridakis_Lojasiewicz-Simon_coupled_Yang-Mills}), there is a constant, $C_1 = C_1(A_1,A_\infty,g,G,p) \in [1,\infty)$ such that
\[
\|u_\infty(A) - u_\infty(A_\infty)\|_{W_{A_1}^{2,p}(X)}
\leq
C_1\|A - A_\infty\|_{W_{A_1}^{2,p}(X)}.
\]
For $\sigma_2 = \sigma_1/C_1 \in (0,1]$ and constant $\sigma_1 \in (0,1]$ given by Lemma \ref{lem:Proof_Rade_7-1_L2_when_Ainfty_smooth_A_Coulomb_p>dover2}, the condition
\[
\|A - A_\infty\|_{W_{A_1}^{2,p}(X)} < \sigma_2
\]
ensures that
\[
\|u_\infty(A) - u_\infty(A_\infty)\|_{W_{A_1}^{2,p}(X)} < \sigma_1.
\]
Lemma \ref{lem:Proof_Rade_7-1_L2_when_Ainfty_smooth_A_Coulomb_p>dover2} now yields, with $A$ and $A_\infty$ replaced by $u_\infty(A)$ and $u_\infty(A_\infty)$, respectively,
\[
\|d_{u_\infty(A)}^*F_{u_\infty(A)}\|_{L^2(X)}
\geq c|\sE(u_\infty(A)) - \sE(u_\infty(A_\infty))|^\theta,
\]
By gauge invariance, the preceding inequality is equivalent to \eqref{eq:Rade_7-1_L2}.
\end{proof}

\begin{lem}[Completion of the proof of the gradient inequality \eqref{eq:Rade_7-1_L2} when $A_\infty$ is $W^{2,p}$ and $A$ is not in Coulomb gauge relative to $A_\infty$]
\label{lem:Proof_Rade_7-1_L2_when_p>dover2}
Assume the hypotheses of Theorem \ref{thm:Rade_proposition_7-2_L2} but require that $p > d/2$. Then the gradient inequality \eqref{eq:Rade_7-1_L2} holds.
\end{lem}

\begin{proof}
By the remarks at the beginning of the proof of Theorem \ref{thm:Rade_proposition_7-2_L2}, there are constants $\eps_0 = \eps_0(A_1,A_\infty,g,G,p) \in (0,1]$ and $C_0 = C_0(A_1,A_\infty,g,G,p) \in [1,\infty)$ and a $W^{3,p}$ gauge transformation $u \in \Aut P$ such that
\begin{gather*}
d_{A_\infty}^*(u(A) - A_\infty) = 0 \quad\hbox{if } \|A - A_\infty\|_{W^{2,p}_{A_1}(X)} < \eps_0,
\\
\|u(A) - A_\infty\|_{W^{2,p}_{A_1}(X)} \leq C_0\|A - A_\infty\|_{W^{2,p}_{A_1}(X)}.
\end{gather*}
By choosing $\sigma$ in the hypotheses of Theorem \ref{thm:Rade_proposition_7-2_L2} small enough that $\sigma = \sigma_2/C_0\in (0,1]$ obeys $\sigma \leq \eps_0$ (for constant $\sigma_2$ given by Lemma \ref{lem:Proof_Rade_7-1_L2_when_A_Coulomb_p>dover2}), then
\[
\|u(A) - A_\infty\|_{W^{2,p}_{A_1}(X)} < \sigma_2.
\]
Lemma \ref{lem:Proof_Rade_7-1_L2_when_A_Coulomb_p>dover2}, with $A$ replaced by $u(A)$, thus yields
\[
\|d_{u(A)}^*F_{u(A)}\|_{L^2(X)}
\geq c|\sE(u(A)) - \sE(A_\infty)|^\theta.
\]
By gauge invariance, the preceding inequality is equivalent to \eqref{eq:Rade_7-1_L2}.
\end{proof}

This completes the proof of Theorem \ref{thm:Rade_proposition_7-2_L2} for Case \ref{case:p_greaterthan_dover2}.
\end{case}

\begin{case}[$2 \leq d \leq 5$ and $p = 2$]
\label{case:2_leq_d_leq_5_and_p_is_2}
The only difference from Case \ref{case:p_greaterthan_dover2} is that we choose
\begin{gather*}
\sX :=  \Ker d_{A_\infty}^*\cap H^2_{A_1}(X;\Lambda^1\otimes\ad P) = \sH_\sA,
\\
\tilde\sX :=  \Ker d_{A_\infty}^*\cap L^2(X;\Lambda^1\otimes\ad P) = \sH.
\end{gather*}
We shall again apply Corollary \ref{cor:Huang_2-4-2_generalized}, but before proceeding to verify its hypotheses, we need to establish continuity of certain Sobolev multiplication maps. The failure of the embedding $H^2(X) \hookrightarrow C(X)$ when $d=4$ indicates that care is needed in the case of a critical exponent when establishing continuity of related Sobolev multiplication maps. However, we have the

\begin{lem}[Continuous Sobolev multiplication map]
\label{lem:Sobolev_multiplication_H2_times_H2_into_H1}
Let $X$ be a closed, Riemannian, smooth manifold of dimension $d \geq 2$. Then the following Sobolev multiplication map is continuous:
\begin{equation}
\label{eq:Sobolev_multiplication_H2_times_H2_into_H1}
H^2(X) \times H^2(X) \to H^1(X), \quad\hbox{if } 2 \leq d \leq 6.
\end{equation}
\end{lem}

\begin{proof}
From \cite[Equation (6.34)]{FU}, one knows that \eqref{eq:Sobolev_multiplication_H2_times_H2_into_H1} holds when $d$ obeys $(2-d/2) + (2-d/2) \geq 1-d/2$ and $d>4$, that is, when $d=5,6$. To see that \eqref{eq:Sobolev_multiplication_H2_times_H2_into_H1} also holds when $d=4$, consider $f_1, f_2 \in H^2(X)$ and observe that
\begin{align*}
\|f_1f_2\|_{L^2(X)} &\leq \|f_1\|_{L^4(X)}\|f_2\|_{L^4(X)}
\\
&\leq \|f_1\|_{H^1(X)}\|f_2\|_{H^1(X)} \quad\hbox{(when $d=4$)}
\\
&\leq \|f_1\|_{H^1(X)}\|f_2\|_{H^2(X)},
\end{align*}
and
\begin{align*}
\|\nabla(f_1f_2)\|_{L^2(X)} &= \|(\nabla f_1)f_2 + f_1\nabla f_2\|_{L^2(X)}
\\
&\leq \|\nabla f_1\|_{L^4(X)}\|f_2\|_{L^4(X)} + \|f_1\|_{L^4(X)}\|\nabla f_2\|_{L^4(X)}
\\
&\leq C\|\nabla f_1\|_{H^1(X)}\|f_2\|_{H^1(X)} + C\|f_1\|_{H^1(X)}\|\nabla f_2\|_{H^1(X)} \quad\hbox{(when $d=4$)}
\\
&\leq C\|f_1\|_{H^2(X)}\|f_2\|_{H^2(X)},
\end{align*}
and this yields the desired continuity of the Sobolev multiplication map \eqref{eq:Sobolev_multiplication_H2_times_H2_into_H1} for the case $d=4$. If $d=2,3$, then $H^2(X)$ is a Banach algebra and so \eqref{eq:Sobolev_multiplication_H2_times_H2_into_H1} holds in this case as well.
\end{proof}

\begin{lem}[Continuous Sobolev multiplication maps]
\label{lem:Sobolev_multiplication_H2_or_H1_times_Hs_into_L2}
Let $X$ be a closed, Riemannian, smooth manifold of dimension $d \geq 2$. Then the following Sobolev multiplication maps are continuous:
\begin{align}
\label{eq:Sobolev_multiplication_H2_times_Hs_into_L2}
H^2(X) \times H^s(X) &\to L^2(X), \quad\hbox{if } d \geq 2 \hbox{ and } \max\{0, d/2 - 2\} \leq s < d/2,
\\
\label{eq:Sobolev_multiplication_H1_times_Hs_into_L2}
H^1(X) \times H^s(X) &\to L^2(X), \quad\hbox{if } d \geq 2 \hbox{ and } \max\{0, d/2 - 1\} \leq s < d/2,
\\
\label{eq:Sobolev_multiplication_H1_times_H2_into_L2}
H^1(X) \times H^2(X) &\to L^2(X), \quad\hbox{if } 2 \leq d \leq 6.
\end{align}
\end{lem}

\begin{proof}
We recall by \cite[Equation (6.34)]{FU} that \eqref{eq:Sobolev_multiplication_H2_times_Hs_into_L2} is continuous when
$$
(2-d/2) + (s - d/2) \geq -d/2 \quad \hbox{and}\quad s \geq 0 \quad \hbox{and}\quad 2s < d,
$$
that is, when $s \geq d/2-2$ and $0 \leq s < d/2$, and this yields \eqref{eq:Sobolev_multiplication_H2_times_Hs_into_L2}.

By \cite[Equation (6.34)]{FU},
the Sobolev multiplication map \eqref{eq:Sobolev_multiplication_H1_times_Hs_into_L2} is continuous when $(1-d/2) + (s-d/2) \geq -d/2$ and $s \geq 0$ and $\max\{2, 2s\} < d$, that is, when $\max\{0, d/2 - 1\} \leq s < d/2$.

The continuity of the Sobolev multiplication map \eqref{eq:Sobolev_multiplication_H1_times_H2_into_L2} is implied by \eqref{eq:Sobolev_multiplication_H1_times_Hs_into_L2} when $s=2$ and $d=5,6$ and by \eqref{eq:Sobolev_multiplication_H1_times_H1_into_L2} and the embedding $H^2(X) \hookrightarrow H^1(X)$ when $2 \leq d \leq 4$.
\end{proof}

The embedding $H^2(X) \hookrightarrow H^s(X)$ is compact for any $d\geq 2$ and $s < 2$ by
\cite[Theorem 6.3]{AdamsFournier}.
Thus, choosing $s=1/2$ in \eqref{eq:Sobolev_multiplication_H2_times_Hs_into_L2}, we obtain a compact embedding $H^2(X) \hookrightarrow H^{3/2}(X)$ and a continuous Sobolev multiplication map,
\begin{equation}
\label{eq:Sobolev_multiplication_H2_times_H1over2_into_L2_when_2_leq_d_leq_5}
H^2(X) \times H^{1/2}(X) \to L^2(X), \quad\hbox{for } 2 \leq d \leq 5.
\end{equation}
Choosing $s=3/2$ in \eqref{eq:Sobolev_multiplication_H1_times_Hs_into_L2} gives a continuous Sobolev multiplication map when $d=4,5$, for $s = 1/2$ when $d=3$, and for $s = 0$ when $d=2$. Therefore, we obtain a continuous Sobolev multiplication map,
\begin{equation}
\label{eq:Sobolev_multiplication_H1_times_H3over2_into_L2_when_2_leq_d_leq_5}
H^1(X) \times H^{3/2}(X) \to L^2(X), \quad\hbox{for } 2 \leq d \leq 5.
\end{equation}
In the case $d=3$, the map \eqref{eq:Sobolev_multiplication_H1_times_H3over2_into_L2_when_2_leq_d_leq_5} is obtained
by composing \eqref{eq:Sobolev_multiplication_H1_times_Hs_into_L2} for $s = 1/2$ with the embedding $H^{3/2}(X)\hookrightarrow H^{1/2}(X)$; in the case $d=2$, the map \eqref{eq:Sobolev_multiplication_H1_times_H3over2_into_L2_when_2_leq_d_leq_5} is obtained by composing \eqref{eq:Sobolev_multiplication_H1_times_Hs_into_L2} for $s = 0$ with the embedding $H^{3/2}(X)\hookrightarrow L^2(X)$.

\begin{claim}[Bounded operator]
\label{claim:Hessian_map_yang_mills_H2_extension_2_leq_d_leq_6}
Assume the hypotheses of Theorem \ref{thm:Rade_proposition_7-2_L2} but require that $A_\infty$ be $C^\infty$ and $2 \leq d \leq 5$ with $p = 2$. Then the following linear operator is bounded,
\begin{equation}
\label{eq:Hessian_map_from_H2_to_L2}
\tilde\sM'(A):H^2_{A_1}(X;\Lambda^1\otimes\ad P) \to L^2(X;\Lambda^1\otimes\ad P).
\end{equation}
\end{claim}

\begin{proof}
From the expression \eqref{eq:Hessian_map_yang_mills_long_expression} for $\tilde\sM'(A)$, we observe that
\begin{align*}
\|a\times \nabla_{A_1}b\|_{L^2(X)} &\leq C\|a\|_{H^2_{A_1}(X)}\|\nabla_{A_1}b\|_{H^{1/2}(X)}
\quad\hbox{(by \eqref{eq:Sobolev_multiplication_H2_times_H1over2_into_L2_when_2_leq_d_leq_5})}
\\
&\leq C\|a\|_{H^2_{A_1}(X)}\|b\|_{H^{3/2}_{A_1}(X)},
\\
\|F_{A_1}\times b\|_{L^2(X)} &\leq \|F_{A_1}\|_{C(X)}\|b\|_{L^2(X)},
\\
\|\nabla_{A_1}a\times b\|_{L^2(X)} &\leq C\|\nabla_{A_1}a\|_{H^1_{A_1}(X)}\|b\|_{H^{3/2}_{A_1}(X)}
\quad\hbox{(by \eqref{eq:Sobolev_multiplication_H1_times_H3over2_into_L2_when_2_leq_d_leq_5})}
\\
&\leq C\|a\|_{H^2_{A_1}(X)}\|b\|_{H^{3/2}_{A_1}(X)},
\\
\|a\times a\times b\|_{L^2(X)} &\leq C\|a\times a\|_{H^1_{A_1}(X)}\|b\|_{H^{3/2}_{A_1}(X)}
\quad\hbox{(by \eqref{eq:Sobolev_multiplication_H1_times_H3over2_into_L2_when_2_leq_d_leq_5})}
\\
&\leq C\|a\|_{H^2_{A_1}(X)}^2\|b\|_{H^{3/2}_{A_1}(X)}
\quad\hbox{(by \eqref{eq:Sobolev_multiplication_H2_times_H2_into_H1})}.
\end{align*}
Consequently, the following operator is bounded,
$$
\sM'(A) - \Delta_{A_1}:H^2_{A_1}(X;\Lambda^1\otimes\ad P) \to L^2(X;\Lambda^1\otimes\ad P).
$$
Clearly, the operator $\Delta_{A_1}:H^2_{A_1}(X;\Lambda^1\otimes\ad P) \to L^2(X;\Lambda^1\otimes\ad P)$ is bounded and combining these observations yields the conclusion. This completes the proof of Claim \ref{claim:Hessian_map_yang_mills_H2_extension_2_leq_d_leq_6}.
\end{proof}

\begin{prop}[Continuity of the extended Hessian operator map for the Yang-Mills energy functional on $H^2$ when $2\leq d \leq 6$]
\label{prop:Continuity_extended_Hessian_operator_map_Yang-Mills_energy_H2_and_2_leq_d_leq_6}
Assume the hypotheses of Theorem \ref{thm:Rade_proposition_7-2_L2} but require that $A_\infty$ be $C^\infty$ and $2\leq d \leq 6$ with $p=2$. Then the following map is continuous,
\begin{multline}
\label{eq:Hessian_map_yang_mills_from_A1_plus_H2_cap_kerdA1star_to_operators_from_H2_cap_kerdA1star_to_L2}
A_\infty + \Ker d_{A_\infty}^*\cap H^2_{A_1}(X;\Lambda^1\otimes\ad P) \ni A \mapsto \sM'(A) \equiv \sM_1(A)
\\
\in \sL\left(\Ker d_{A_\infty}^*\cap H^2_{A_1}(X;\Lambda^1\otimes\ad P), \Ker d_{A_\infty}^*\cap L^2(X;\Lambda^1\otimes\ad P)\right).
\end{multline}
\end{prop}

\begin{proof}
The fact that the linear operator,
\[
\tilde\sM'(A):H^2_{A_1}(X;\Lambda^1\otimes\ad P) \to L^2(X;\Lambda^1\otimes\ad P),
\]
is a bounded follows from Claim \ref{claim:Hessian_map_yang_mills_H2_extension_2_leq_d_leq_6}. The expression \eqref{eq:Hessian_map_yang_mills_at_Aplusa_minus_at_A_on_b}, for $a, a_1, b \in H^2_{A_1}(X;\Lambda^1\otimes\ad P)$ and writing $A=A_1+a_1$, gives
\begin{align*}
{}&\|(\tilde\sM'(A+a) - \tilde\sM'(A))b\|_{L^2(X)}
\\
&\quad \leq \|a\times \nabla_{A_1} b\|_{L^2(X)} + \|\nabla_{A_1} a\times b\|_{L^2(X)}
+ \|(a_1+a)\times a \times b\|_{L^2(X)}
\\
&\quad \leq C\|a\|_{H^2_{A_1}(X)}\|\nabla_{A_1} b\|_{H^1_{A_1}(X)}
+ C\|\nabla_{A_1} a\|_{H^1_{A_1}(X)}\|b\|_{H^2_{A_1}(X)}
\\
&\qquad + \|(a_1+a)\times a\|_{H^1_{A_1}(X)}\|b\|_{H^2_{A_1}(X)}
\quad \hbox{(by \eqref{eq:Sobolev_multiplication_H1_times_H2_into_L2} for $2 \leq d \leq 6$)}
\\
&\quad \leq C\left(1 + \|a_1\|_{H^2_{A_1}(X)} + \|a\|_{H^2_{A_1}(X)}\right)\|a\|_{H^2_{A_1}(X)}\|b\|_{H^2_{A_1}(X)}
\quad\hbox{(by \eqref{eq:Sobolev_multiplication_H2_times_H2_into_H1})}.
\end{align*}
Hence, the pre-gradient map,
\[
\tilde\sM: A_\infty + H^2_{A_1}(X;\Lambda^1\otimes\ad P) \ni A
\mapsto \tilde\sM(A) \in L^2(X;\Lambda^1\otimes\ad P),
\]
is Lipschitz continuous. Boundedness of the $L^2$-orthogonal projection $\Pi_{A_\infty}$ in \eqref{eq:L2-orthogonal_projection_onto_slice_Wkp} implies that the map,
\[
A \mapsto \sM'(A) \equiv \sM_1(A) = \Pi_{A_\infty}\tilde\sM'(A),
\]
in \eqref{eq:Hessian_map_yang_mills_from_A1_plus_H2_cap_kerdA1star_to_operators_from_H2_cap_kerdA1star_to_L2} is also Lipschitz continuous.
\end{proof}

It remains to verify real analyticity of the gradient map.

\begin{prop}[Real analyticity of the gradient map for the Yang-Mills energy functional on $H^2$ when $2\leq d \leq 6$]
\label{prop:Real_analyticity_gradient_map_Yang-Mills_energy_H2_and_2_leq_d_leq_6}
Assume the hypotheses of Theorem \ref{thm:Rade_proposition_7-2_L2} but require that $A_\infty$ be $C^\infty$ and $2 \leq d \leq 6$ with $p=2$. Then the gradient map $\sM$ in \eqref{eq:Gradient_map_yang_mills_A1_plus_W2p_cap_kerdA1star_to_Lp} is real analytic.
\end{prop}

\begin{proof}
By the expression \eqref{eq:Pre_gradient} for $\tilde\sM(A)$ and writing $A = A_1+a_1$ we see that, for all $a_1, a \in H^2_{A_1}(X;\Lambda^1\otimes\ad P)$,
\begin{align*}
{}&\|\tilde\sM(A+a) - \tilde\sM(A)\|_{L^2(X)}
\\
&\quad = \|\tilde\sM(A_1+a_1+a) - \tilde\sM(A_1+a_1)\|_{L^2(X)}
\\
&\quad \leq \|d_{A_1}^*d_{A_1}a\|_{L^2(X)} + \|F_{A_1}\times a\|_{L^2(X)}
\\
&\qquad + \|\nabla_{A_1} a\times a\|_{L^2(X)} + \|\nabla_{A_1} a_1\times a\|_{L^2(X)} + \|a_1\times \nabla_{A_1} a\|_{L^2(X)}
\\
&\qquad + \|a\times a \times a\|_{L^2(X)} + \|a_1\times a \times a\|_{L^2(X)} + \|a_1\times a_1 \times a\|_{L^2(X)}
\\
&\quad \leq C\left(\|a\|_{H^2_{A_1}(X)} + \|F_{A_1}\|_{H^1_{A_1}(X)}\|a\|_{H^2_{A_1}(X)} \right.
\\
&\qquad + \|\nabla_{A_1} a\|_{H^1_{A_1}(X)}\|a\|_{H^2_{A_1}(X)} + \|\nabla_{A_1} a_1\|_{H^1_{A_1}(X)}\|a\|_{H^2_{A_1}(X)}
+ \|a_1\|_{H^2_{A_1}(X)}\|\nabla_{A_1} a\|_{H^1_{A_1}(X)}
\\
&\qquad + \left. \|a\times a\|_{H^1_{A_1}(X)}\|a\|_{H^2_{A_1}(X)}
+ \|a_1\|_{H^2_{A_1}(X)}\|a\times a\|_{H^1_{A_1}(X)} + \|a_1\times a_1\|_{H^1_{A_1}(X)}\|a\|_{H^2_{A_1}(X)} \right),
\end{align*}
where we applied \eqref{eq:Sobolev_multiplication_H1_times_H2_into_L2} to obtain the last inequality. Thus,
\begin{align*}
{}&\|\tilde\sM(A+a) - \tilde\sM(A)\|_{L^2(X)}
\\
&\quad \leq C\left( \|a\|_{H^2_{A_1}(X)} + \|F_{A_1}\|_{H^1_{A_1}(X)}\|a\|_{H^2_{A_1}(X)} \right.
\\
&\qquad + \|a\|_{H^2_{A_1}(X)}^2 + \|a_1\|_{H^2_{A_1}(X)}\|a\|_{H^2_{A_1}(X)}
+ \|a_1\|_{H^2_{A_1}(X)}\|a\|_{H^2_{A_1}(X)}
\\
&\qquad + \left. \|a\|_{H^2_{A_1}(X)}^3 + \|a_1\|_{H^2_{A_1}(X)}\|a\|_{H^1_{A_1}(X)}^2
+ \|a_1\|_{H^2_{A_1}(X)}^2\|a\|_{H^2_{A_1}(X)} \right)
\quad\hbox{(by \eqref{eq:Sobolev_multiplication_H2_times_H2_into_H1})}
\\
&\quad \leq C\left(1 + \|F_{A_1}\|_{H^1_{A_1}(X)} + \|a_1\|_{H^2_{A_1}(X)} + \|a\|_{H^2_{A_1}(X)} + \|a_1\|_{H^2_{A_1}(X)}^2 + \|a\|_{H^2_{A_1}(X)}^2 \right)\|a\|_{H^2_{A_1}(X)}.
\end{align*}
Hence, the pre-gradient map,
\[
\tilde\sM: A_\infty + H^2_{A_1}(X;\Lambda^1\otimes\ad P) \ni A
\mapsto \tilde\sM(A) \in L^2(X;\Lambda^1\otimes\ad P),
\]
is Lipschitz continuous
and thus real analytic since the function $\tilde\sM(A_1+a)$ is a polynomial of degree three in $a \in H^2_{A_1}(X;\Lambda^1\otimes\ad P)$.

Boundedness of the $L^2$-orthogonal projection $\Pi_{A_\infty}$ in \eqref{eq:L2-orthogonal_projection_onto_slice_Wkp} and the expression \eqref{eq:Gradient_map_yang_mills_short_expression} for $\sM(A) = \Pi_{A_\infty}\tilde\sM(A)$ thus imply that the gradient map $\sM$ in \eqref{eq:Gradient_map_yang_mills_A1_plus_W2p_cap_kerdA1star_to_Lp} is real analytic when $2\leq d\leq 6$ and $p=2$. This completes the proof of Proposition \ref{prop:Real_analyticity_gradient_map_Yang-Mills_energy_H2_and_2_leq_d_leq_6}.
\end{proof}

\begin{lem}[Completion of the proof of the gradient inequality \eqref{eq:Rade_7-1_L2} when $A_\infty$ is $C^\infty$ and $A$ is in Coulomb gauge relative to $A_\infty$]
\label{lem:Proof_Rade_7-1_L2_when_Ainfty_smooth_A_Coulomb_2_leq_d_leq_6_and_p=2}
Assume the hypotheses of Theorem \ref{thm:Rade_proposition_7-2_L2} but require that $A_\infty$ be $C^\infty$ and $2\leq d\leq 5$ with $p=2$ and $d_{A_\infty}^*(A-A_\infty) = 0$. Then the gradient inequality \eqref{eq:Rade_7-1_L2} holds.
\end{lem}

\begin{proof}
Just as in Case \ref{case:p_greaterthan_dover2}, Proposition \ref{prop:Continuity_extended_Hessian_operator_map_Yang-Mills_energy_W2p_and_p>dover2} guarantees that the map, $A_\infty+\sX \ni A \mapsto \sM_1(A) \in \sL(\sH_\sA,\sH)$, is continuous. Proposition \ref{prop:Real_analyticity_gradient_map_Yang-Mills_energy_H2_and_2_leq_d_leq_6}
ensures that the gradient map, $\sM:A_\infty+\sX \to \tilde\sX$, is real analytic. Lemma \ref{lem:fred_Lp} implies that $\sM'(A_\infty) \in \sL(\sX,\tilde\sX) = \sL(\sH_\sA,\sH)$ is Fredholm with index zero. Hence, Corollary \ref{cor:Huang_2-4-2_generalized} yields the gradient inequality \eqref{eq:Rade_7-1_L2}.
\end{proof}

The additional technical conditions in the hypothesis of Lemma \ref{lem:Proof_Rade_7-1_L2_when_Ainfty_smooth_A_Coulomb_2_leq_d_leq_6_and_p=2} (that $A_\infty$ is $C^\infty$ and $d_{A_\infty}^*(A-A_\infty) = 0$) are removed exactly as in Case \ref{case:p_greaterthan_dover2}, so this completes the proof of Theorem \ref{thm:Rade_proposition_7-2_L2} for Case \ref{case:2_leq_d_leq_5_and_p_is_2}.
\end{case}

Combining Cases \ref{case:p_greaterthan_dover2} and \ref{case:2_leq_d_leq_5_and_p_is_2} concludes the proof of Theorem \ref{thm:Rade_proposition_7-2_L2}.
\end{proof}

\begin{rmk}[Generalization to the case $d \geq 2$ and $p \in [2,\infty)$ obeying $p \geq d/3$]
\label{rmk:Proposition_Rade_7-2_L2_for_d_geq_6_and_2_lessthan_p_leq_d}
It should be possible to adapt the arguments of Feehan and Maridakis \cite{Feehan_Maridakis_Lojasiewicz-Simon_coupled_Yang-Mills} to prove a refinement of Theorem \ref{thm:Rade_proposition_7-2_L2} that holds for all $d \geq 2$ and $p \in [2,\infty)$ obeying $p \geq d/3$, with $A$ and $A_\infty$ of class $W^{2,q}$ and $q > d/3$.
\end{rmk}

\subsection{A $W^{1,2}$ {\L}ojasiewicz-Simon gradient inequality for the Yang-Mills energy functional}
\label{subsec:W1p_Lojasiewicz-Simon_gradient_inequality_Yang-Mills}
We now use Theorem \ref{thm:Huang_2-4-5} to obtain our second version of \cite[Proposition 7.2]{Rade_1992}, with the restriction on the dimension, $d$, of $X$ that $d=2,3$ relaxed to $2\leq d\leq 4$.

We digress in order to identify a dual space arising in the proof of Theorem \ref{thm:Rade_proposition_7-2} below. If $W$ is a subspace of a Banach space, $V$, with dual space, $V'$, one denotes $W^\perp : = \{\ell \in V': \ell(W) = \{0\}\} \subset W'$. If $V$ is reflexive, so $V'' = V$, then $(W^\perp)^\perp = W$.
The map
$$
V'/W^\perp \cong W'
$$
induced by $\ell \mapsto \ell\restriction W$ is an isometric isomorphism.

Given a bounded, linear operator $Q:V_1 \to V_2$ between Banach spaces $V_1$ and $V_2$, then $Q^*:V_2'\to V_1'$ is defined by $Q^*\ell_2(v_1) = \ell_2(Qv_1)$ for all $v_1 \in V_1$ and $\ell_2 \in V_2'$. In particular, $\Ker Q^* = (\Imag Q)^\perp$ since
$$
\ell_2 \in \Ker Q^* \iff Q^*\ell_2(v_1) = \ell_2(Qv_1) = 0 \quad\forall\, v_1 \in V_1.
$$
Applying the preceding observation to $Q^*:V_2' \to V_1'$ yields $\Ker Q^* = \Ker Q^{**} = (\Imag Q^*)^\perp$. If $W := \Ker Q \subset V_1$ and $V_1$ is reflexive, then $W^\perp = (\Ker Q)^\perp = ((\Imag Q^*)^\perp)^\perp = \Imag Q^*$ and
$$
(\Ker Q)' \cong V_1'/(\Ker Q)^\perp = V_1'/\Imag Q^*.
$$
As a consequence of Theorem \ref{thm:Huang_2-4-5}, we obtain the following version of \cite[Proposition 7.2]{Rade_1992}, including R\r{a}de's stronger assertion \cite[Equation (9.1)]{Rade_1992} analogous to \eqref{eq:Simon_2-2_dualspacenorm}, but with the restriction on the dimension, $d$, of $X$ relaxed from $d = 2,3$ relaxed to $2 \leq d \leq 4$ when $p=2$.

\begin{thm}[An $H^1$ {\L}ojasiewicz-Simon gradient inequality for the Yang-Mills energy functional with $H^{-1}$ gradient norm]
\label{thm:Rade_proposition_7-2}
Let $(X,g)$ be a closed, Riemannian, smooth manifold of dimension $d$, and $G$ a compact Lie group, $A_1$ a connection of class $C^\infty$, and $A_\infty$ a Yang-Mills connection of class $W^{1,q}$, with $q>d/2$, on a principal $G$-bundle $P$ over $X$. If $2 \leq d\leq 4$ and $p = 2$, then there are constants, $c\in (0,\infty)$ and $\sigma \in (0,1]$ and $\theta \in [1/2,1)$ depending on $A_1$, $A_\infty$, $g$, $G$, and $p$, with the following significance.  If $A$ is a connection of class $W^{1,q}$ on $P$ and
\begin{equation}
\label{eq:Rade_7-1_neighborhood}
\|A - A_\infty\|_{H^1_{A_1}(X)} < \sigma,
\end{equation}
then
\begin{equation}
\label{eq:Rade_7-1}
\|d_A^*F_A\|_{H^{-1}_{A_1}(X)} \geq c|\sE(A) - \sE(A_\infty)|^\theta.
\end{equation}
where $\sE(A)$ is given by \eqref{eq:Potential_yang_mills}.
\end{thm}

\begin{proof}
The inequality \eqref{eq:Rade_7-1} is gauge-invariant, that is, if $u \in \Aut P$ is a gauge transformation of class $W^{2,q}$ (and thus continuous under our hypotheses on $d$ and $q$), then \eqref{eq:Rade_7-1} is equivalent to
$$
\|d_{u(A)}^*F_{u(A)}\|_{W_{u(A_1)}^{-1,p'}(X)} \geq c|\sE(u(A)) - \sE(u(A_\infty))|^\theta.
$$
By analogy with the proof of Theorem \ref{thm:Rade_proposition_7-2_L2}, we shall first apply a $W^{2,q}$ gauge transformation $u \in \Aut P$
to yield a global Coulomb gauge condition, $d_{A_\infty}^*(u(A)-A_\infty)=0$, when $A$ is $W_{A_1}^{1,p}$-close to $A_\infty$. Because we allow $p=2$ when $d=4$ (the critical exponent), the question of existence of the gauge transformation is more delicate. In particular, it does \emph{not} follow from standard references as in the proof of Theorem \ref{thm:Rade_proposition_7-2_L2}. However, we can apply our \cite[Theorem 10]{Feehan_Maridakis_Lojasiewicz-Simon_coupled_Yang-Mills} to give, for constants $\eps_0 = \eps_0(A_1,A_\infty,g,p,q) \in (0,1]$ and $C_0 = C_0(A_1,A_\infty,g,p,q) \in [1,\infty)$, a $W^{2,q}$ gauge transformation $u \in \Aut P$ such that
\begin{gather}
\label{eq:Coulomb_gauge}
d_{A_\infty}^*(u(A) - A_\infty) = 0 \quad\hbox{if } \|A - A_\infty\|_{W_{A_1}^{1,p}(X)} < \eps_0,
\\
\label{eq:Coulomb_gauge_transformed_connection_near_original}
\|u(A) - A_\infty\|_{W_{A_1}^{1,p}(X)} \leq C_0\|A - A_\infty\|_{W_{A_1}^{1,p}(X)},
\end{gather}
where $u$ is $W_{A_1}^{2,q}$-close to the identity automorphism of $P$ and $p$ obeys $d/2\leq p \leq q$.

We shall apply Theorem \ref{thm:Huang_2-4-5} with
$$
\sX := \Ker d_{A_\infty}^*\cap H^1_{A_1}(X;\Lambda^1\otimes\ad P) \quad\hbox{and}\quad
\sH := \Ker d_{A_\infty}^*\cap L^2(X;\Lambda^1\otimes\ad P).
$$
From the digression preceding the statement of Theorem \ref{thm:Rade_proposition_7-2}, we see that
$$
\sX' = H^{-1}_{A_1}(X;\Lambda^1\otimes\ad P)/\Imag d_{A_1},
$$
by choosing $V_1 := H^1_{A_1}(X;\Lambda^1\otimes\ad P)$ and $V_2 = L^2(X;\Lambda^1\otimes\ad P)$ and
$$
Q := d_{A_\infty}^*: H^1_{A_1}(X;\Lambda^1\otimes\ad P) \to L^2(X;\Lambda^1\otimes\ad P),
$$
so that
$$
Q^* = d_{A_\infty}^{**} = d_{A_1}: L^2(X;\Lambda^1\otimes\ad P) \to H^{-1}_{A_1}(X;\Lambda^1\otimes\ad P).
$$
But $H^{-1}_{A_1}(X;\Lambda^1\otimes\ad P)/\Imag d_{A_\infty} = \Ker d_{A_\infty}^*\cap H^{-1}_{A_1}(X;\Lambda^1\otimes\ad P)$, where $d_{A_\infty}^*:H^{-1}_{A_1}(X;\Lambda^1\otimes\ad P) \to H^{-2}_{A_1}(X;\Lambda^1\otimes\ad P)$ is a bounded operator. Therefore, we also have
$$
\sX' = \Ker d_{A_\infty}^*\cap H^{-1}_{A_1}(X;\Lambda^1\otimes\ad P).
$$
The formal expressions for the differential map and Hessian operator,
\begin{multline}
\label{eq:Gradient_map_yang_mills_A1_plus_H1_cap_kerdA1star_to_Hminus1}
\sE':A_\infty + \Ker d_{A_\infty}^*\cap H^1_{A_1}(X;\Lambda^1\otimes\ad P) \ni A
\\
\mapsto \sE'(A) \in \Ker d_{A_\infty}^*\cap H^{-1}_{A_1}(X;\Lambda^1\otimes\ad P),
\end{multline}
\begin{multline}
\label{eq:Hessian_map_yang_mills_A1_plus_H1_cap_kerdA1star_to_operators_from_H1_cap_kerdA1star_to_Hminus1}
\sE''(A_\infty): \Ker d_{A_\infty}^*\cap H^1_{A_1}(X;\Lambda^1\otimes\ad P)  \ni a
\\
\mapsto \sE''(A_\infty)a \in \Ker d_{A_\infty}^*\cap H^{-1}_{A_1}(X;\Lambda^1\otimes\ad P),
\end{multline}
are given by \eqref{eq:Pre_gradient} and \eqref{eq:Pre_Hessian}. We now establish the required Fredholm and index zero properties for the Hessian, $\sE''(A_\infty)$.

\begin{lem}[Fredholm and index zero properties of the Hessian on $H^1$ for the Yang-Mills $L^2$-energy functional on a Coulomb-gauge slice]
\label{lem:fred_L2}
Let $(X,g)$ be a closed, Riemannian, smooth manifold of dimension $d \geq 2$, and $G$ a compact Lie group, $P$ a smooth principal $G$-bundle over $X$, and $A_1$ a $C^\infty$ reference connection on $P$. If $A_\infty$ is a $C^\infty$ connection on $P$, then the following operator is Fredholm with index zero,
\[
\sE''(A_\infty): \Ker d_{A_\infty}^* \cap H_{A_1}^1(X; \Lambda^1 \otimes \ad P)
\to \Ker d_{A_\infty}^* \cap H_{A_1}^{-1}(X; \Lambda^1 \otimes \ad P).
\]
\end{lem}

\begin{proof}
By H\"ormander \cite[Theorem 19.2.1]{Hormander_v3}, the operator,
\[
\sE''(A_\infty): H_{A_1}^1(X; \Lambda^1 \otimes \ad P) \to H_{A_1}^{-1}(X; \Lambda^1 \otimes \ad P),
\]
is Fredholm with index zero. The argument of R\r{a}de \cite[p. 148]{Rade_1992} for the case $p=2$ (and $d=2,3$) again adapts to prove that the restriction of the preceding operator to a slice domain (with slice range) is also Fredholm with index zero. (See the proof of \cite[Proposition 3.6]{Feehan_Maridakis_Lojasiewicz-Simon_coupled_Yang-Mills} for details.)
\end{proof}

It remains to verify that the potential function $\sE$ is analytic.

\begin{prop}[Real analyticity of the Yang-Mills energy functional on $H^1$ when $2\leq d \leq 4$]
\label{prop:Real_analyticity_gradient_map_Yang-Mills_energy_H1_and_2_leq_d_leq_4}
Assume the hypotheses of Theorem \ref{thm:Rade_proposition_7-2} but require that $A_\infty$ be $C^\infty$. Then the following potential function is real analytic,
\begin{equation}
\label{eq:Potential_function_yang_mills_from_A1_plus_H1_cap_kerdA1star_to_reals}
\sE: A_\infty + \Ker d_{A_\infty}^*\cap H^1_{A_1}(X;\Lambda^1\otimes\ad P) \ni A \mapsto \sE(A) \in \RR.
\end{equation}
\end{prop}

\begin{proof}
Writing $A = A_1+a_1$ and using $F_{A+a} = F_A + d_Aa + a\times a$ and
\begin{align*}
2\sE(A+a) &= \|F_{A+a}\|_{L^2(X)}^2
\\
&= \|F_A\|_{L^2(X)}^2 + \|d_Aa\|_{L^2(X)}^2 + \|a\times a\|_{L^2(X)}^2
\\
&\quad + 2(F_A, d_Aa)_{L^2(X)} + 2(F_A, a \times a)_{L^2(X)} + 2(d_Aa, a \times a)_{L^2(X)},
\end{align*}
we see that, for all $a_1, a \in H^1_{A_1}(X;\Lambda^1\otimes\ad P)$,
\begin{align*}
{}&2\sE(A+a) - 2\sE(A)
\\
&\quad = 2\sE(A_1+a_1+a) - 2\sE(A_1+a_1)
\\
&\quad = \|d_{A_1}(a_1+a)\|_{L^2(X)}^2 + \|(a_1+a)\times (a_1+a)\|_{L^2(X)}^2 + 2(F_{A_1}, d_{A_1}(a_1+a))_{L^2(X)}
\\
&\qquad + 2(F_{A_1}, (a_1+a) \times (a_1+a))_{L^2(X)} + 2(d_{A_1}(a_1+a), (a_1+a) \times (a_1+a))_{L^2(X)}
\\
&\qquad - \|d_{A_1}a_1\|_{L^2(X)}^2 - \|a_1\times a_1\|_{L^2(X)}^2 - 2(F_{A_1}, d_{A_1}a_1)_{L^2(X)}
\\
&\qquad - 2(F_{A_1}, a_1 \times a_1)_{L^2(X)} - 2(d_{A_1}a_1, a_1 \times a_1)_{L^2(X)}
\\
&\quad = \|d_{A_1}a\|_{L^2(X)}^2 + 2(d_{A_1}a_1,d_{A_1}a)_{L^2(X)} + \|a\times a\|_{L^2(X)}^2
+ 2\|a_1\times a\|_{L^2(X)}^2
\\
&\qquad + 2(F_{A_1}, d_{A_1}a)_{L^2(X)} + 2(F_{A_1}, a \times a)_{L^2(X)}
+ 4(F_{A_1}, a_1 \times a)_{L^2(X)}
\\
&\qquad + 2(d_{A_1}a, (a_1+a) \times (a_1+a))_{L^2(X)}
\\
&\qquad + 2(d_{A_1}a_1, a \times a)_{L^2(X)} + 2(d_{A_1}a_1, a \times a_1)_{L^2(X)} + 2(d_{A_1}a_1, a_1 \times a)_{L^2(X)},
\end{align*}
and therefore, employing \eqref{eq:Sobolev_multiplication_H1_times_H1_into_L2} for $2\leq d\leq 4$,
\begin{align*}
{}&|\sE(A+a) - \sE(A)|
\\
&\quad \leq C\left(\|a\|_{H^1_{A_1}(X)}^2 + \|a_1\|_{H^1_{A_1}(X)}\|a\|_{H^1_{A_1}(X)} + \|a\|_{H^1_{A_1}(X)}^4
+ \|a_1\|_{H^1_{A_1}(X)}^2\|a\|_{H^1_{A_1}(X)}^2 \right.
\\
&\qquad + \|F_{A_1}\|_{L^2(X)}\|a\|_{H^1_{A_1}(X)} + \|F_{A_1}\|_{L^2(X)}\|a\|_{H^1_{A_1}(X)}^2
+ \|F_{A_1}\|_{L^2(X)}\|a_1\|_{H^1_{A_1}(X)}\|a\|_{H^1_{A_1}(X)}
\\
&\qquad + \left. \|a\|_{H^1_{A_1}(X)}^3 + \|a_1\|_{H^1_{A_1}(X)}\|a\|_{H^1_{A_1}(X)}^2
+ \|a_1\|_{H^1_{A_1}(X)}^2\|a\|_{H^1_{A_1}(X)}\right).
\end{align*}
Consequently, the function,
\[
\sE: A_\infty + H^1_{A_1}(X;\Lambda^1\otimes\ad P) \ni A \mapsto \sE(A) \in H^{-1}_{A_1}(X;\Lambda^1\otimes\ad P),
\]
is Lipschitz continuous
and thus real analytic since the function $\sE(A_1+a)$ is a polynomial of degree four in $a \in H^1_{A_1}(X;\Lambda^1\otimes\ad P)$. Hence, the potential function \eqref{eq:Potential_function_yang_mills_from_A1_plus_H1_cap_kerdA1star_to_reals} is real analytic, as it is obtained by restriction from $H^1_{A_1}(X;\Lambda^1\otimes\ad P)$ to the closed subspace $\Ker d_{A_\infty}^* \cap H^1_{A_1}(X;\Lambda^1\otimes\ad P)$.
\end{proof}

\begin{lem}[Completion of the proof of the gradient inequality \eqref{eq:Rade_7-1} when $A_\infty$ is $C^\infty$ and $A$ is in Coulomb gauge relative to $A_\infty$]
\label{lem:Proof_Rade_7-1_when_Ainfty_smooth_A_Coulomb}
Assume the hypotheses of Theorem \ref{thm:Rade_proposition_7-2_L2} but require that $A_\infty$ be $C^\infty$ and $d_{A_\infty}^*(A-A_\infty) = 0$. Then the gradient inequality \eqref{eq:Rade_7-1} holds.
\end{lem}

\begin{proof}
Proposition \ref{prop:Real_analyticity_gradient_map_Yang-Mills_energy_H1_and_2_leq_d_leq_4} ensures that the potential function, $\sE:A_\infty+\sX \to \RR$, is real analytic. Lemma \ref{lem:fred_L2} implies that the Hessian, $\sE''(A_\infty) \in \sL(\sX,\sX^*)$, is Fredholm with index zero. Hence, Theorem \ref{thm:Huang_2-4-5} yields the gradient inequality \eqref{eq:Rade_7-1}.
\end{proof}

The additional technical conditions in the hypothesis of Lemma \ref{lem:Proof_Rade_7-1_when_Ainfty_smooth_A_Coulomb} (that $A_\infty$ is $C^\infty$ and $d_{A_\infty}^*(A-A_\infty) = 0$) are removed exactly as in Case \ref{case:p_greaterthan_dover2} of the proof of Theorem \ref{thm:Rade_proposition_7-2_L2}. This completes the proof of Theorem \ref{thm:Rade_proposition_7-2}.
\end{proof}

Feehan and Maridakis have established a generalization of Theorem \ref{thm:Rade_proposition_7-2} for arbitrary $d \geq 2$ and coupled Yang-Mills energy functionals. We record a special case of their result in the following

\begin{thm}[A $W^{1,p}$ {\L}ojasiewicz-Simon gradient inequality for the Yang-Mills energy functional with $W^{-1,p}$ gradient norm]
\label{thm:Rade_proposition_7-2_W-1p}
(See \cite[Theorem 4]{Feehan_Maridakis_Lojasiewicz-Simon_coupled_Yang-Mills}.)
Let $(X,g)$ be a closed, Riemannian, smooth manifold of dimension $d\geq 2$, and $G$ be a compact Lie group, and $P$ be a smooth principal $G$-bundle over $X$. Let $A_1$ be a $C^\infty$ reference connection on $P$, and $A_\infty$ a Yang-Mills connection on $P$ for $g$ of class $W^{1,q}$, with $q \in [2,\infty)$ obeying $q > d/2$. If $p \in [2,\infty)$ obeys $d/2 \leq p \leq q$, then there are constants $c \in (0, \infty)$, and $\sigma \in (0,1]$, and $\theta \in [1/2,1)$, depending on $A_1$, $A_\infty$, $g$, $G$, $p$, and $q$ with the following significance. If $A$ is a $W^{1,q}$ connection on $P$ obeying the \emph{{\L}ojasiewicz-Simon neighborhood} condition,
\begin{equation}
\label{eq:Lojasiewicz-Simon_gradient_inequality_Yang-Mills_neighborhood}
\|A - A_\infty\|_{W^{1,p}_{A_1}(X)} < \sigma,
\end{equation}
then the Yang-Mills energy functional \eqref{eq:Yang-Mills_energy_functional} obeys the \emph{{\L}ojasiewicz-Simon gradient inequality}
\begin{equation}
\label{eq:Lojasiewicz-Simon_gradient_inequality_Yang-Mills_energy_functional}
\|d_A^*F_A\|_{W^{-1,p}_{A_1}(X)}
\geq
c|\sE(A) - \sE(A_\infty)|^\theta.
\end{equation}
\end{thm}

\begin{rmk}[Range of $d$ in Theorem \ref{thm:Rade_proposition_7-2_W-1p} allowing replacement of $W^{-1,p}$ by $L^2$]
While Theorem \ref{thm:Rade_proposition_7-2_W-1p} is valid for arbitrary $d \geq 2$, our application to Yang-Mills gradient flow requires an inequality of the form
\begin{equation}
\label{eq:Lojasiewicz-Simon_gradient_inequality_Yang-Mills_energy_functional}
\|d_A^*F_A\|_{L^2(X)}
\geq
c|\sE(A) - \sE(A_\infty)|^\theta,
\end{equation}
and thus a continuous embedding, $L^2(X) \subset W^{-1,p}(X)$, or equivalently, $W^{1,p'}(X) \subset L^2(X)$, where $p' = p/(p-1) \in (1,2]$ is the H\"older exponent dual to $p \in [2,\infty)$. The largest possible value of $p'$ is attained with the smallest allowable value of $p$, namely, $p=2$ for $d=2,3$ and $p = d/2$ for $d \geq 4$; in the latter case, we obtain $p' = (d/2)/((d/2) - 1) = d/(d-2) < d$. According to \cite[Theorem 4.12]{AdamsFournier}, when $p'<d$, the embedding, $W^{1,p'}(X) \subset L^2(X)$, is continuous if $(p')^* = dp'/(d-p') = dp/(d(p-1)-p) \geq 2$ and substituting $p=d/2$, this yields the constraint
\[
d^2/2 \geq 2(d^2/2 - 3d/2),
\]
or $d \geq 2d - 6$, that is, $d \leq 6$. In other words, Theorem \ref{thm:Rade_proposition_7-2_W-1p} implies that the inequality \eqref{eq:Lojasiewicz-Simon_gradient_inequality_Yang-Mills_energy_functional} holds for $2 \leq d \leq 6$.
\end{rmk}

\subsection{R\r{a}de's proof of the {\L}ojasiewicz-Simon gradient inequality for the Yang-Mills energy functional in dimensions two and three}
R\r{a}de gives two proofs of his \cite[Proposition 7.2]{Rade_1992} when $\dim X = 2$ or $3$; more precisely, he proves \eqref{eq:Rade_7-1} when $p=2$ and $\dim X = 2$ or $3$ in \cite[Section 9]{Rade_1992} while he proves \eqref{eq:Rade_7-1_L2} when $G=U(n)$ and $\dim X = 2$ in \cite[Section 10]{Rade_1992}, with $\theta = 1/2$ when $A_\infty$ is irreducible and $\theta = 3/4$ when $A_\infty$ is reducible.

R\r{a}de's first proof of \cite[Proposition 7.2]{Rade_1992} in \cite[Section 9]{Rade_1992} is modeled closely on that of Simon's proof of \cite[Theorem 3]{Simon_1983}, as R\r{a}de notes in \cite[p. 146]{Rade_1992}. His proof primarily depends on the dimension of $X$ through an appeal to the Sobolev multiplication map \cite[p. 147]{Rade_1992},
\begin{equation}
\label{eq:H1_times_Hminus1_to_Hminus2_when_d_is_2_or_3}
H^1(X) \times H^{-1}(X) \to H^{-2}(X).
\end{equation}
To understand for which dimensions \eqref{eq:H1_times_Hminus1_to_Hminus2_when_d_is_2_or_3} is continuous, we recall from \cite[p. 96]{FU}, with $k_1=k=1$, $k_2=2$, and $p_1=q=q_2=2$, that \eqref{eq:H1_times_Hminus1_to_Hminus2_when_d_is_2_or_3} is continuous if and only if the same is true for \cite[Equation (6.34)]{FU} with $p_2=p=2$ and so $1 = 1/p+1/q$ and $1 = 1/p_2+1/q_2$, that is,
\begin{equation}
\label{eq:H1_times_H2_to_H1_when_d_is_2_or_3}
H^1(X) \times H^2(X) \to H^1(X).
\end{equation}
The difficulty here is that, when $d=2$ or $3$, one has $d>4$ and hence $H^2(X) \hookrightarrow C(X)$ by \cite[Theorem 4.12, Part I (A)]{AdamsFournier}, but when $d=4$ one only has $H^2(X) \hookrightarrow L^q(X)$ for all $q<\infty$ by \cite[Theorem 4.12, Part I (B)]{AdamsFournier}. Restricting temporarily to $d=2$ or $3$, we see that
\begin{align*}
\|f_1f_2\|_{L^2(X)} &\leq \|f_1\|_{L^4(X)}\|f_2\|_{L^4(X)}
\\
&\leq \|f_1\|_{H^1(X)}\|f_2\|_{H^1(X)} \quad\hbox{(when $d\leq 4$)}
\\
&\leq \|f_1\|_{H^1(X)}\|f_2\|_{H^2(X)},
\end{align*}
and
\begin{align*}
\|\nabla(f_1f_2)\|_{L^2(X)} &= \|(\nabla f_1)f_2 + f_1\nabla f_2\|_{L^2(X)}
\\
&\leq \|\nabla f_1\|_{L^2(X)}\|f_2\|_{L^\infty(X)} + \|f_1\|_{L^4(X)}\|\nabla f_2\|_{L^4(X)}
\\
&\leq C\|f_1\|_{H^1(X)}\|f_2\|_{H^2(X)} + C\|f_1\|_{H^1(X)}\|\nabla f_2\|_{H^1(X)} \quad\hbox{(when $d=2,3$)}
\\
&\leq C\|f_1\|_{H^1(X)}\|f_2\|_{H^2(X)},
\end{align*}
and so we obtain \eqref{eq:H1_times_H2_to_H1_when_d_is_2_or_3}. To derive \eqref{eq:H1_times_Hminus1_to_Hminus2_when_d_is_2_or_3}, we observe that
\begin{align*}
\|f_1f_2\|_{H^{-2}(X)}
&= \sup_{\begin{subarray}{c}f \in H^2(X), \\ \|f\|_{H^2(X)} \leq 1\end{subarray}} (f_1f_2, f)_{L^2(X)}
\\
&\leq \sup_{\begin{subarray}{c}f \in H^2(X), \\ \|f\|_{H^2(X)} \leq 1\end{subarray}}
\|f_1f\|_{H^1(X)}\|f_2\|_{H^{-1}(X)} \quad\hbox{(by duality)}
\\
&\leq \sup_{\begin{subarray}{c}f \in H^2(X), \\ \|f\|_{H^2(X)} \leq 1\end{subarray}}
\|f_1\|_{H^1(X)}\|f\|_{H^2(X)}\|f_2\|_{H^{-1}(X)} \quad\hbox{(by \eqref{eq:H1_times_H2_to_H1_when_d_is_2_or_3})}
\\
&\leq \|f_1\|_{H^1(X)}\|f_2\|_{H^{-1}(X)},
\end{align*}
as desired. However, when $d=4$, the preceding argument breaks down due to the failure of the embedding $H^2(X) \hookrightarrow C(X)$ when $d=4$. Indeed, by \eqref{eq:Sobolev_multiplication_H2_times_H2_into_H1}, one only has the embedding
$$
H^2(X) \times H^2(X) \to H^1(X).
$$
Of course, one could envisage replacing the Sobolev multiplication map \eqref{eq:Sobolev_multiplication_H2_times_H2_into_H1}, and thus modifying R\r{a}de's proof of \eqref{eq:Rade_7-1} in \cite[Section 9]{Rade_1992}, by exploiting the critical-exponent Sobolev embeddings and norms introduced by Taubes in \cite{TauFrame, TauStable, TauConf}. (See also \cite{FeehanSlice} for an application due to the author of Taubes' norms to the development of refined slice theorems for the quotient space of connections.) However, we shall see that Theorem \ref{thm:Rade_proposition_7-2} may instead be quickly deduced from the abstract Theorem \ref{thm:Huang_2-4-5}.

While R\r{a}de's proof \cite[Section 10]{Rade_1992} of \eqref{eq:Rade_7-1_L2} also makes use of a dimension-dependent Sobolev embedding \cite[p. 151]{Rade_1992}, namely $H^1(X) \hookrightarrow L^4(X)$, that holds when $\dim X = 2,3,4$ although not $\dim X \geq 5$ (see \cite[Theorem 4.12, Part I (C)]{AdamsFournier}. However, the argument in \cite[p. 153]{Rade_1992} makes explicit use of the assumption that $\dim X = 2$ (and $G = U(n)$).

\section{Convergence and stability for abstract gradient systems}
\label{sec:Huang_3_and_5_gradient-system}
The broad outline of our proofs of the abstract convergence and stability theorems follows the scheme described by Huang in \cite[Sections 3.3 and 5.1]{Huang_2006}, but we simplify his treatment by
\begin{inparaenum}[\itshape a\upshape)]
\item occasionally assuming greater regularity for the flow, $u(t)$, since this is normally a consequence of regularity for solutions to nonlinear, parabolic, quasilinear equations,
\item restricting our attention initially to pure gradient flow, in particular $F\equiv 0$ and $G\equiv 0$ in \cite[Equation (3.10a)]{Huang_2006}, and only later considering gradient-like (or pseudogradient) flow as needed,
\item replacing Huang's rather strong growth condition \cite[Equation (3.10a)]{Huang_2006} with a simpler `energy' counterpart such as \cite[Equation (3.5)]{Simon_1983}, where $\sX$ is a real, reflexive Banach space that is continuously embedded and dense in a Hilbert space, $\sH$, and (following \cite[p. 76]{Huang_2006}) we choose $\sY=\sH$, so $\sY' = \sH'$, and we identify $\sH' \cong \sH$ in \cite[Section 3.3]{Huang_2006}, and
\item Assuming \apriori interior estimates for $\|\dot u\|_{L^1((S+\delta,T); \sX)}$ in terms of $\|\dot u\|_{L^1((S,T); \sH)}$ with $\delta > 0$ and $S \geq 0$ and $T \leq \infty$ obeying $S + \delta \leq T$.
\end{inparaenum}
The alternative framework we develop is still quite general but closer in spirit to our applications to Yang-Mills gradient and gradient-like flows.

\subsection{Abstract gradient systems}
\label{sec:Huang_3-1}
Let $\sX$ be a real, reflexive Banach space that is continuously embedded and dense in a Hilbert space, $\sH$, and identify $\sH \cong \sH'$ by the Riesz map, so that $\sX\hookrightarrow \sH \hookrightarrow \sX'$ with
\begin{equation}
\label{eq:Riesz_pairing_relations}
\langle \rho, v \rangle_{\sX'\times \sX} = \langle \rho, v \rangle_{\sH\times \sH}, \quad \forall \, \rho \in \sH \subset \sX' \hbox{ and } v \in \sX.
\end{equation}
Let $\sU$ denote an open subset of $\sX$. We begin with a notion of trajectory which specializes that of \cite[Definition 3.1.1]{Huang_2006} by taking $\sY = \sH \cong \sH' = \sY'$.

\begin{defn}[Trajectory]
\label{defn:Huang_3-1-1}
By a \emph{trajectory} in $\sU$ we mean a map $u:[0,T) \to \sU$, for $0<T\leq\infty$, obeying the following conditions:
\begin{enumerate}
\item (Strong continuity) The map $u:[0,T) \to \sU$ is continuous with respect to the norm of $\sX$;
\item (Weak-star differentiability) There is a (weakly measurable\footnote{See \cite[Appendix C.1]{Sell_You_2002}}) map, $\dot u:[0,T)\to \sH$, called the \emph{weak-star derivative of $u$}, such that for each $\rho \in \sH$, the function $[0,T) \ni t \to \langle \rho, \dot u(t) \rangle_{\sX'\times \sX}$ is continuous and such that
\begin{equation}
\label{eq:Huang_3-1d}
\langle \rho, u(a) - u(b) \rangle_{\sX'\times \sX} = \int_a^b \langle \rho, \dot u(t) \rangle_{\sX'\times \sX} \,dt,
\quad\forall\, a, b \in [0,T) \hbox{ and } \rho \in \sH.
\end{equation}
\end{enumerate}
In the sequel, we use
$$
O(u) := \{u(t): t \in [0,T)\},
$$
to denote the \emph{orbit} of the trajectory $u:[0,T)\to \sU$.
\end{defn}

The weak-star continuity of the derivative, $\dot u:[0,T)\to \sH$, implies that for each compact subset, $K \subset [0,T)$,
$$
\sup_{t\in K} |\langle \rho, \dot u(t) \rangle_{\sX'\times \sX}| < \infty, \quad\forall\, \rho \in \sH.
$$
An application of the classical Resonance Theorem \cite[Corollary 2.1.1]{Yosida} ensures that the norms $\|\dot u(t)\|_\sH$ are uniformly bounded on $K$. It follows that the function $[0,T) \ni t \to \|\dot u(t)\|_\sH$ is lower semicontinuous. Note that since the map $[0, T] \ni t\mapsto \|\dot u(t)\|_\sH \in \RR$ is lower semicontinuous, then it is also Lebesgue measurable \cite[Definition 2.7.20 and Corollary 3.11.5]{HewittStromberg}.

Since $\sX$ is reflexive, the canonical map $\varphi:\sX\cong \sX''$ given by $x\mapsto \varphi(x)$, where $\varphi(x)\rho := \rho(x)$ for all $\rho \in \sX'$, is an isometric isomorphism and thus
$$
\rho(u(a) - u(b))
\equiv
\langle \rho, u(a) - u(b) \rangle_{\sX'\times \sX}
\equiv
\langle \rho, \varphi(u(a) - u(b)) \rangle_{\sX'\times \sX''}, \quad \forall\, \rho \in \sX'.
$$
Moreover, $\sX'$ is necessarily reflexive as well and so by James' Theorem \cite[p. 120]{Megginson_1998}, there exists a $\rho_0 \in \sX'$ such that $\|\rho_0\|_{\sX'} \leq 1$ and
$$
\langle \rho_0, \varphi(u(a) - u(b)) \rangle_{\sX'\times \sX''} = \|\varphi(u(a) - u(b))\|_{\sX''},
$$
Therefore, since $\varphi:\sX\cong \sX''$ is an isometry,
$$
\langle \rho_0, u(a) - u(b) \rangle_{\sX'\times \sX} = \|u(a) - u(b)\|_\sX,
$$
and hence \eqref{eq:Huang_3-1d} implies that, if we further assume $\dot u:[0, T) \to \sX$ is Bochner integrable \cite[Appendix C.1]{Sell_You_2002},
\begin{equation}
\label{eq:Huang_3-1d_inequality}
\|u(a) - u(b)\|_\sX \leq \int_a^b \|\dot u(t)\|_\sX \,dt,
\quad\forall\, a, b \in [0,T),
\end{equation}
noting that $\langle \rho, \dot u(t) \rangle_{\sX'\times \sX} \leq \|\rho\|_{\sX'}\|\dot u(t)\|_\sX$ for all $\rho \in \sX'$.

\begin{rmk}[Absolutely continuous paths in Banach spaces]
\label{rmk:Absolutely_continuous_paths_Banach_spaces}
Following \cite[p. 104]{Showalter}, we let $W^{1,p}(0,T;\sX)$ be the set of functions $u:[0,T]\to \sX$ such that for some $v \in L^p(0,T;\sX)$,
$$
u(t) = u(0) + \int_0^t v(s)\,ds, \quad\forall\, t \in [0, T],
$$
where $L^p(S;\sX)$, for $(S,\Sigma,\mu)$ a measure space and $1\leq p \leq \infty$, is the Banach space of (equivalence classes of) measurable functions $v:S\to \sX$ such that $\|v(\cdot)\|_\sX \in L^p(S;\RR)$, with norm
$$
\|v\|_{L^p(S;\sX)}
:=
\begin{cases}
\displaystyle
\left(\int_S\|v(s)\|_\sX\,d\mu\right)^{1/p}, & 1 \leq p < \infty,
\\
\displaystyle
\esssup_{s\in S}\|v(s)\|_\sX, & p = \infty.
\end{cases}
$$
If $u \in W^{1,1}(0,T;\sX)$, for $0<T\leq\infty$, then $u$ is strongly differentiable a.e. on $(0,T)$ with $\dot u(t) = v(t)$ for a.e. $t\in (0,T)$ \cite[Theorem III.1.6]{Showalter} and
$$
\|u(b) - u(a)\|_\sX \leq \int_a^b\|\dot u(t)\|_\sX\,dt, \quad\forall\, a, b \in [0, T],
$$
so $u:[0,T]\to \sX$ is absolutely continuous and, in particular, continuous. Conversely, if $\sX$ is reflexive and $u:[0,T]\to \sX$ is absolutely continuous, then $u$ is strongly differentiable a.e. on $(0,T)$, $\dot u \in L^1(0,T;\sX)$ and $u(t) = u(0) + \int_0^t \dot u(s)\,ds$ for $0\leq t \leq T$ \cite[Theorem III.1.7]{Showalter}.
\end{rmk}

Recall the

\begin{prop}[Differentiability of an energy functional along a trajectory]
\label{prop:Huang_3-1-2}
\cite[Proposition 3.1.2]{Huang_2006}
Let $\sE':\sU\subset \sX \to \sH$ be a gradient map associated with a $C^1$ functional, $\sE:\sU \to \RR$, where $\sU$ is an open subset of a real, reflexive Banach space, $\sX$, that is continuously embedded and dense in a Hilbert space, $\sH$. If $u:[0,T)\to \sU$ is a trajectory in $\sU$ in the sense of Definition \ref{defn:Huang_3-1-1}, then the function $[0,T) \ni t \mapsto \sE(u(t)) \in \RR$ is continuously differentiable and
$$
\frac{d}{dt}\sE(u(t)) = \langle \sE'(u(t)), \dot u(t) \rangle_{\sX'\times \sX}, \quad\forall\, t\in [0, T).
$$
\end{prop}

We have the following notions of solution to a gradient system.

\begin{defn}[Strong solution]
\label{defn:Strong_solution_to_gradient_system}
Let $\sU$ be an open subset of a real Banach space, $\sX$, that is continuously embedded and dense in a Hilbert space, $\sH$. Let $\sE:\sU\subset \sX\to\RR$ be a $C^1$ functional with gradient map $\sE':\sU\subset \sX \to \sH$ in the sense of Definition \ref{defn:Huang_2-1-1}. We call a trajectory, $u:[0,T)\to \sU$, a \emph{strong solution} of the gradient system for $\sE$ if
\begin{equation}
\label{eq:Huang_3-3a}
\dot u(t) = -\sE'(u(t)), \quad\hbox{a.e. } t \in (0,T), \quad u(0) = u_0 \in \sU,
\end{equation}
as an equation in $\sH$.
\end{defn}

Compare \eqref{eq:Huang_3-3a} with \cite[Equations (3.4A*) or (3.28a)]{Huang_2006}. We have the following specialization of \cite[Definition 3.1.3]{Huang_2006} (compare \cite[Proposition III.2.1 and p. 108]{Showalter}).

\begin{defn}[Mild and weak solutions]
\label{defn:Huang_3-1-3}
Let $\sU$ be an open subset of a real Banach space, $\sX$, which is continuously embedded and dense in a Hilbert space, $\sH$.
\begin{enumerate}
\item A continuous map $u:[0,T)\to \sU \subset \sX$ is called a \emph{mild solution} of \eqref{eq:Huang_3-3a} if $u$ is a solution to the integral equation in $\sH$,
\begin{equation}
\label{eq:Huang_3-3b}
u(t)= u_0 - \int_0^t \sE'(u(s))\, ds, \quad\forall\, t\in [0, T).
\end{equation}
\item A map $u:[0,T)\to \sU \subset \sX$ is called a \emph{weak solution} of \eqref{eq:Huang_3-3a} if $u$ is a trajectory in $\sU$ in the sense of Definition \ref{defn:Huang_3-1-1} and obeys
\begin{equation}
\label{eq:Huang_3-4}
\langle \dot u(t), v \rangle_{\sH\times \sH} = \langle -\sE'(u(t)), v \rangle_{\sX'\times \sX},
\quad\forall\, t\in [0, T), \quad\forall\, v \in \sX.
\end{equation}
\end{enumerate}
\end{defn}

Note that a strong solution to \eqref{eq:Huang_3-3a} is also a mild solution and a mild solution is also a weak solution \cite[p. 61]{Huang_2006}. If $u$ is a weak solution, then
\begin{equation}
\label{eq:Huang_3-3c}
\langle u(t) - u_0, v \rangle_{\sH\times \sH} = -\int_0^t \langle \sE'(u(t)), v \rangle_{\sX'\times \sX},
\quad\forall\, t\in [0, T) \hbox{ and } v \in \sH,
\end{equation}
where the canonical pairing $\langle \cdot , \cdot \rangle_{\sX\times \sX'}$ may be replaced in the preceding discussion by the inner product $\langle \cdot , \cdot \rangle_{\sH\times \sH}$ using \eqref{eq:Riesz_pairing_relations}.

Consider a strong solution, $u:[0,T)\to \sU$, of \eqref{eq:Huang_3-3a}. It follows from Proposition \ref{prop:Huang_3-1-2} that the function $[0,T) \ni t \mapsto \sE(u(t)) \in \RR$ is differentiable and
\begin{equation}
\label{eq:Huang_3-5}
\frac{d}{dt}\sE(u(t))
=
\langle \sE'(u(t)), \dot u(t) \rangle_{\sH\times \sH}
=
-\|\sE'(u(t))\|_\sH^2, \quad\forall\, t\in [0, T],
\end{equation}
where we use \eqref{eq:Huang_3-3a} to obtain the second equality. The function $\sE(u(t))$ is called the \emph{energy} of the solution, $u$.

\subsection{Huang's three technical lemmata}
\label{subsec:Huang_3-2}
We review three useful technical lemmata described by Huang in \cite[Section 3.2]{Huang_2006}.

\begin{lem}
\label{lem:Huang_3-2-1}
\cite[Lemma 3.2.1]{Huang_2006}
For $-\infty < t_0 < t_1 < \infty$, let $H:[t_0, t_1] \to \RR$ be an absolutely continuous, monotone function.\footnote{Hence, $H' \in L^1(t_0, t_1; \RR)$.} Let $\Phi: [H(t_0) \wedge H(t_1), H(t_0) \vee H(t_1)] \to \RR$ be an absolutely continuous function. Then the composition $G := \Phi \circ H$ defined by $G(t) = \Phi(H(t))$ for $t\in [t_0, t_1]$  is also absolutely continuous and obeys
$$
G'(t) = \Phi'(H(t))H'(t), \quad \hbox{a.e. } t \in [t_0, t_1].
$$
\end{lem}

Lemma \ref{lem:Huang_3-2-2} below will not directly be used in the sequel, but is included for completeness.

\begin{lem}
\label{lem:Huang_3-2-2}
\cite[Lemma 3.2.2]{Huang_2006}
Let $h:[0, \infty) \to \RR$ be an absolutely continuous, bounded, monotone function. Assume that the following inequality holds,
\begin{equation}
\label{eq:Huang_3-8}
\phi(h(t)) \leq \sqrt{|h'(t)|} + g(t), \quad \forall\, t \in \Lambda,
\end{equation}
where $\Lambda \subset [0, \infty)$ is a measurable subset and $g \in L^1(0, \infty; \RR)$ is non-negative and $\phi:\RR \to [0, \infty)$ is a (non-negative) measurable function with $1/\phi \in L^1_{\loc}(\RR)$. Then $\sqrt{h'} \in L^1(\Lambda; \RR)$ and
$$
\int_\Lambda \sqrt{h'(t)} \, dt \leq \|g\|_{L^1(0, \infty)} + \left| \int_{h(0)}^{h(\infty)} \frac{ds}{\phi(s)}\right|.
$$
where $h(\infty) := \lim_{t \to \infty} h(t)$. As a consequence, if $f \in L^2(0, \infty; \RR)$ is non-negative and obeys
\begin{equation}
\label{eq:Huang_3-8prime}
\tag{\ref*{eq:Huang_3-8}$'$}
\phi\left(\int_t^\infty f^2(s)\,ds\right) \leq f(t) + g(t), \quad \forall\, t \in \Lambda,
\end{equation}
then $f \in L^1(\Lambda; \RR)$ and
$$
\int_\Lambda f(t)\,dt \leq \|g\|_{L^1(0, \infty)} + \int_0^{\|f\|_{L^2(0,\infty)}^2} \frac{ds}{\phi(s)}.
$$
\end{lem}

Recall that a subset $O$ of a Banach space $\sX$ is \emph{precompact} if every sequence in $O$ has a subsequence which converges to a limit in $\sX$ \cite[Section D.5]{Evans2}. If $O\subset \sX$ is bounded and $\sX$ is reflexive, then every sequence in $O$ contains a weakly convergent subsequence in $\sX$ \cite[Theorem D.3]{Evans2}. If $O\subset \sX$ is bounded and $\sX\hookrightarrow \sY$ is a compact embedding of $\sX$ in a Banach space $\sY$, then every sequence in $O$ contains a convergent subsequence in $\sY$ \cite[Section D.5]{Evans2}. One says that $\varphi \in \sX$ is a \emph{cluster point} of the orbit $O(u) = \{u(t): t \geq 0\}$ of a trajectory $u:[0,\infty)\to \sX$ \cite[p. 65]{Huang_2006} if there exists an unbounded sequence $\{t_n\}_{n\geq 1} \subset [0, \infty)$ such that
\begin{equation}
\label{eq:Orbit_cluster_point}
\|\varphi(t_n) - \varphi\|_\sX \to 0 \quad\hbox{as } n \to \infty.
\end{equation}
We have the following simpler version of \cite[Lemma 3.2.3]{Huang_2006}, obtained by choosing $\sY = \sH$ and $\sY' \cong \sH$.

\begin{lem}[Convergence of trajectories]
\label{lem:Huang_3-2-3}
Let $\sX$ be a real Banach space, $\sX$, that is continuously embedded and dense in a Hilbert space, $\sH$, and let $u$ be a trajectory in the sense of Definition \ref{defn:Huang_3-1-1} whose orbit, $O(u)$, is precompact in $\sX$. Moreover, assume that there exists a cluster point $\varphi \in \sX$ of $O(u)$ and a constant $\sigma > 0$ such that
\begin{subequations}
\label{eq:Huang_3-9}
\begin{gather}
\label{eq:Huang_3-9a}
\int_W \|\dot u\|_\sH\,dt < \infty,
\\
\label{eq:Huang_3-9b}
\hbox{where } W := \{t\geq 0 : \|u(t) - \varphi\|_\sX < \sigma \},
\end{gather}
\end{subequations}
with $\dot u: [0, \infty) \to \sH$ locally Bochner integrable. Then $u(t)$ converges to $\varphi$ as $t\to\infty$ in the sense that
\begin{equation}
\label{eq:Huang_3-9c}
\tag{\ref*{eq:Huang_3-9}c}
\lim_{t\to\infty}\|u(t)-\varphi\|_\sX = 0 \quad\hbox{and}\quad \int_0^\infty \|\dot u\|_\sH\,dt < \infty.
\end{equation}
\end{lem}

\begin{proof}
By continuity of the trajectory $u:[0,\infty)\to \sX$, the set $W\subset [0,\infty)$ in \eqref{eq:Huang_3-9b} is open and also non-empty since $\varphi$ is a cluster point of $O(u)$. We choose an unbounded sequence $\{t_n\}_{n\geq 1} \subset W$ such that
\begin{equation}
\label{eq:Orbit_cluster_point_lemma_Huang_3-2-3_proof}
\|u(t_n) - \varphi\|_\sX \to 0 \quad\hbox{as } n \to \infty.
\end{equation}
For each $t_n$, we let
$$
\tilde t_n := \sup\{t \geq t_n: [t_n, t] \subset W\}.
$$
We claim that there exists some integer $N\in\NN$ such that $\tilde t_N = \infty$. We give an argument by contradiction. Suppose that for each $n \in N$ we have $\tilde t_n < \infty$. Then $[t_n, \tilde t_n) \subset W$, but the point $\tilde t_n$ is located in the boundary of $W$, that is,
\begin{equation}
\label{eq:tilde_tn_boundary_W}
\|u(\tilde t_n) - \varphi\|_\sX = \sigma, \quad \forall\, n \in \NN.
\end{equation}
By hypothesis, $\dot u: [0, \infty) \to \sH$ is locally Bochner integrable and the Newton-Leibnitz formula (Theorem \ref{thm:Sell_You_C-9}, with $\sX$ replaced by $\sH$) gives
$$
\|u(\tilde t_n) - u(t_n)\|_\sH \leq \int_{t_n}^{\tilde t_n} \|\dot u(t)\|_\sH\,dt,
$$
and thus
$$
\|u(\tilde t_n) - \varphi\|_\sH \leq \|u(t_n) - \varphi\|_\sH + \int_{t_n}^{\tilde t_n}  \|\dot u(t)\|_\sH\,dt, \quad \forall\, n \in \NN.
$$
By \eqref{eq:Huang_3-9a} we have
$$
\sum_{n=1}^\infty \int_{t_n}^{\tilde t_n}  \|\dot u(t)\|_\sH\,dt \leq \int_W  \|\dot u(t)\|_\sH\,dt < \infty,
$$
and so
$$
\int_{t_n}^{\tilde t_n}  \|\dot u(t)\|_\sH\,dt \to 0  \quad\hbox{as } n \to \infty,
$$
and \eqref{eq:Orbit_cluster_point_lemma_Huang_3-2-3_proof} therefore yields,
\begin{equation}
\label{eq:limit_norm_H_u(tilde_tn)_minus_varphi_is_zero}
\lim_{t\to\infty} \|u(\tilde t_n) - \varphi\|_\sH = 0.
\end{equation}
This implies that if $\tilde\varphi \in \sX$ is any cluster point of the precompact sequence $\{u(\tilde t_n)\}_{n \in \NN} \subset O(u)$, then we must have $\|\tilde\varphi - \varphi\|_\sH = 0$ and thus $\tilde\varphi = \varphi$. Hence, $\varphi \in \sX$ is the unique cluster point of the precompact sequence $\{u(\tilde t_n)\}_{n \in \NN}$ in $\sX$ and thus
$$
\|u(\tilde t_n) - \varphi\|_\sX  \to 0  \quad\hbox{as } n \to \infty,
$$
contradicting \eqref{eq:tilde_tn_boundary_W}.

Hence, there must be some $N \in \NN$ such that $\tilde t_N = \infty$, that is, $[t_N, \infty) \subset W$ and so \eqref{eq:Huang_3-9a} implies that
$$
\int_{t_N}^\infty  \|\dot u(t)\|_\sH \,dt < \infty.
$$
This yields the integral convergence in \eqref{eq:Huang_3-9c}, since
$$
\int_0^\infty  \|\dot u(t)\|_\sH\,dt = \int_0^{t_N}  \|\dot u(t)\|_\sH\,dt + \int_{t_N}^\infty  \|\dot u(t)\|_\sH\,dt < \infty,
$$
where the finiteness of the integral $\int_0^{t_N}  \|\dot u(t)\|_\sH\,dt$ is assured by our hypothesis that $\dot u: [0, \infty) \to \sH$ is locally Bochner integrable and the fact that $t_N<\infty$. As a consequence, we also have $\|u(t) - \varphi\|_\sH \to 0$ as $t \to \infty$.

Moreover, if $\hat\varphi \in \sX$ is any cluster point of the precompact orbit $O(u) \subset \sX$, then we must have $\|u(t) - \hat\varphi\|_\sH \to 0$ as $t \to \infty$ and consequently $\|\hat\varphi - \varphi\|_\sH = 0$ and so $\hat\varphi = \varphi$. Hence, $\varphi \in \sX$ is the unique cluster point of the precompact orbit $O(u) \subset \sX$ and thus
$$
\|u(t) - \varphi\|_\sX  \to 0  \quad\hbox{as } n \to \infty,
$$
yielding the desired convergence in the norm of $\sX$ in \eqref{eq:Huang_3-9c}.
\end{proof}

Lemma \ref{lem:Huang_3-2-3_X} below will not directly be used in the sequel, but is included for completeness.

\begin{lem}[Convergence of trajectories]
\label{lem:Huang_3-2-3_X}
Let $\sX$ be a real Banach space, $\sX$, that is continuously embedded and dense in a Hilbert space, $\sH$, and let $u$ be a trajectory in the sense of Definition \ref{defn:Huang_3-1-1} whose orbit, $O(u)$, is precompact in $\sX$. Moreover, assume that there exists a cluster point $\varphi \in \sX$ of $O(u)$ and a constant $\sigma > 0$ such that
\begin{subequations}
\label{eq:Huang_3-9_X}
\begin{gather}
\label{eq:Huang_3-9a_X}
\int_W \|\dot u\|_\sX\,dt < \infty,
\\
\label{eq:Huang_3-9b_X}
\hbox{where } W := \{t\geq 0 : \|u(t) - \varphi\|_\sX < \sigma \},
\end{gather}
\end{subequations}
with $\dot u: (0, \infty) \to \sX$ locally Bochner integrable. Then $u(t)$ converges to $\varphi$ as $t\to\infty$ in the sense that
\begin{equation}
\label{eq:Huang_3-9c_X}
\tag{\ref*{eq:Huang_3-9_X}c}
\lim_{t\to\infty}\|u(t)-\varphi\|_\sX = 0 \quad\hbox{and}\quad \int_1^\infty \|\dot u\|_\sX\,dt < \infty.
\end{equation}
\end{lem}

\begin{proof}
The proof is a trivial modification of that of Lemma \ref{lem:Huang_3-2-3}. Indeed, in the statement of \cite[Lemma 3.2.3]{Huang_2006}, the two results are essentially combined.
\end{proof}

\subsection{Convergence in abstract gradient systems}
\label{subsec:Huang_3-3}
We simplify \cite[Section 3.3]{Huang_2006}, restricting our attention to gradient systems defined by energy functionals obeying the {\L}ojasiewicz-Simon gradient inequality. We assume that $\sU$ is an open subset of a Banach space $\sX$ that is continuously, dense embedding, $\sX \hookrightarrow \sH$, in a Hilbert space, $\sH$, and thus $\sH \cong \sH' \hookrightarrow \sX'$, where we identify $\sH$ with $\sH'$ by the Riesz map \eqref{eq:Riesz_pairing_relations}. We also assume that $\sE:\sU\to\RR$ is a $C^1$ functional and that the gradient map, $\sE':\sU\to \sH$, is continuous. We wish to study the convergence of gradient trajectories $u:[0,\infty)\to\sU$, with respect to the gradient map, $\sE':\sU\to \sH$, which locate themselves in a neighborhood of the set of stationary points \eqref{eq:Huang_page_25_defn_critical_point_set}. We shall need the following technical hypothesis to allow us to relate convergence of $u(t)$ in the Hilbert space $\sH$ to the stronger convergence of $u(t)$ in the Banach space $\sX$ as $t \to \infty$.

\begin{hyp}[Regularity and \apriori interior estimate for a trajectory]
\label{hyp:Abstract_apriori_interior_estimate_trajectory}
Let $C_1$ and $\rho$ be positive constants and let $T \in (0, \infty]$. Given a trajectory $u:[0,T)\to \sU$, in the sense of Definition \ref{defn:Huang_3-1-1}, we say that $\dot u:[0,T) \to \sH$ obeys an \apriori \emph{interior estimate on $(0, T]$} if, for every $S \geq 0$ and $\delta > 0$ obeying $S+\delta \leq T$, the map $\dot u:[S+\delta,T) \to \sX$ is Bochner integrable and there holds
\begin{equation}
\label{eq:Abstract_apriori_interior_estimate_trajectory}
\int_{S+\delta}^T \|\dot u(t)\|_\sX\,dt \leq C_1(1+\delta^{-\rho})\int_S^T \|\dot u(t)\|_\sH\,dt.
\end{equation}
\end{hyp}

\begin{rmk}[Regularity and \apriori interior estimates for solutions to Yang-Mills gradient flow]
\label{rmk:Abstract_apriori_interior_estimate_trajectory}
The Hypothesis \ref{hyp:Abstract_apriori_interior_estimate_trajectory} is verified for nonlinear evolution equations on Banach spaces by Lemma \ref{lem:Rade_7-3_abstract_interior_L1_in_time_V2beta_space_time_derivative_interior} and for solutions to Yang-Mills gradient flow with $\sH = L^2(X; \Lambda^1\otimes\ad P)$ in Lemma \ref{lem:Rade_7-3}, where $\sX = H_{A_1}^1(X; \Lambda^1\otimes\ad P)$ and $\rho = 1/2$ and $X$ has dimension $d$ in the range $2\leq d\leq 4$, in Lemma \ref{lem:Rade_7-3_L1_in_time_H2beta_in_space_apriori_estimate_by_L1_in_time_L2_in_space}, where $\sX = H_{A_1}^{2\beta}(X; \Lambda^1\otimes\ad P)$ with $\beta \in [1/4+d/8, 1)$ when $2\leq d\leq 5$ and $\rho = 1$, in Lemma \ref{lem:Rade_7-3_arbitrary_dimension}, where $\sX = W_{A_1}^{2\beta,r}(X; \Lambda^1\otimes\ad P)$ for suitable $\beta \in (1/2, 1)$ and $r \in (1,\infty)$ when $d \geq 2$ and $\rho = 1$, and in Corollary \ref{cor:Rade_7-3_arbitrary_dimension_L1_time_W1p_space}, where $\sX = W_{A_1}^{1,p}(X; \Lambda^1\otimes\ad P)$ for $p \in (2\vee d/2, \infty)$ when $d \geq 2$ and $\rho = 1$.
\end{rmk}

We then have the following simpler version of \cite[Proposition 3.3.2]{Huang_2006}, which we include for completeness though it is not explicitly used in the sequel.

\begin{prop}[Convergence under the validity of the {\L}ojasiewicz-Simon gradient inequality]
\label{prop:Huang_3-3-2}
Let $\sU$ be an open subset of a real Banach space, $\sX$, that is continuously embedded and dense in a Hilbert space, $\sH$. Let $\sE:\sU\subset \sX \to \RR$ be a $C^1$ functional on an open subset $\sU$ of a Banach space $\sX$ with gradient map $\sE':\sU\subset \sX \to \sH$ and let $u:[0,\infty)\to \sU$ be a strong solution to \eqref{eq:Huang_3-3a} in the sense of Definition \ref{defn:Strong_solution_to_gradient_system} such that\footnote{We correct a small typographical error in the statement of \cite[Proposition 3.3.2]{Huang_2006}, where a supremum is indicated rather than the required infimum in \cite[Equation (3.14)]{Huang_2006}}
\begin{equation}
\label{eq:Huang_3-14}
\inf\{|\sE(u(t))| : t \geq 0\} > - \infty.
\end{equation}
If the gradient map $\sE'$ satisfies a {\L}ojasiewicz-Simon gradient inequality \eqref{eq:Simon_2-2} in the orbit $O(u)$, that is,
\begin{equation}
\label{eq:Huang_3-15a}
\|\sE'(u(t))\|_\sH \geq c|\sE(u(t)) - \sE(\varphi)|^\theta, \quad \forall\, t \geq 0,
\end{equation}
for positive constants $c$ and $\theta \in [1/2, 1)$, then
\begin{equation}
\label{eq:Huang_3-15b_H}
\int_0^\infty \|\dot u(t)\|_\sH\,dt
\leq
\int_{\sE_\infty}^{\sE(u(0))} \frac{1}{c|s - \sE(\varphi)|^\theta} \,ds < \infty,
\end{equation}
where $\sE_\infty := \lim_{t\to\infty}\sE(u(t))$. If $u$ obeys Hypothesis \ref{hyp:Abstract_apriori_interior_estimate_trajectory}, then, for all $\delta > 0$,
\begin{equation}
\label{eq:Huang_3-15b_X}
\int_\delta^\infty \|\dot u\|_\sX\,dt
\leq
C_1(1+\delta^{-\rho}) \int_{\sE_\infty}^{\sE(u(0))} \frac{1}{c|s - \sE(\varphi)|^\theta} \,ds < \infty.
\end{equation}
If in addition $\dot u: [0,\infty)\to \sX$ is locally Bochner integrable on $[0,\infty)$, then
$$
\int_0^\infty \|\dot u\|_\sX\,dt < \infty.
$$
\end{prop}

\begin{proof}
The function $[0, \infty) \ni t \mapsto \sE(u(t)) \in \RR$ is $C^1$ by Proposition \ref{prop:Huang_3-1-2} and obeys
\begin{align*}
-\frac{d}{dt}\sE(u(t)) &= \langle -\sE'(u(t)), \dot u(t) \rangle_{\sX'\times \sX} \quad\hbox{(by Proposition \ref{prop:Huang_3-1-2})}
\\
&= (-\sE'(u(t)), \dot u(t))_\sH  \quad\hbox{(by \eqref{eq:Riesz_pairing_relations})}
\\
&= \|\sE'(u(t))\|_\sH\|\dot u(t)\|_\sH, \quad\forall\, t \in [0, \infty) \quad\hbox{(by \eqref{eq:Huang_3-3a})}.
\end{align*}
Hence, $\sE(u(t))$ is a nonincreasing and uniformly bounded function of $t \in [0,\infty)$ by \eqref{eq:Huang_3-14}, so $a \equiv \sE_\infty = \lim_{t\to\infty}\sE(u(t))$ exists, as asserted by the proposition. Set $H(t) := \sE(u(t))$, for all $t \in [0, \infty)$, and observe that $H(t)$ is monotone and absolutely continuous on $[0, \infty)$ and obeys, by the preceding equality,
\begin{equation}
\label{eq:Huang_3-16a}
-\frac{d}{dt}H(t) = \|\sE'(u(t))\|_\sH\|\dot u(t)\|_\sH, \quad\forall\, t \in [0, \infty).
\end{equation}
Let $\phi:\RR\to\RR$ be the function defined by $\phi(s) = c|s - \sE(\varphi)|^\theta$, for all $s \in \RR$, and let $\Phi:\RR \to \RR$ be the absolutely continuous function given by
$$
\Phi(x) := \int_a^x \frac{1}{\phi(s)} \,ds = \int_a^x \frac{1}{c|s - \sE(\varphi)|^\theta} \,ds, \quad \forall\, x \in \RR,
$$
where the limit $a := \lim_{t\to\infty}H(t) = \sE_\infty > -\infty$ exists. The function $\Phi$ is differentiable a.e. on $\RR$ with $\Phi'(x) = 1/\phi(x)$ for a.e. $x \in \RR$. According to Lemma \ref{lem:Huang_3-2-1}, the composition $\Phi\circ H$ is absolutely continuous on $[0, \infty)$ and there holds
\begin{equation}
\label{eq:Huang_3-16b}
\frac{d}{dt}\Phi(H(t)) = \frac{H'(t)}{\phi(H(t))}, \quad \forall\, t \in \Lambda,
\end{equation}
where $\Lambda \subset [0, \infty)$ is such that the complement, $[0, \infty) \less \Lambda$, has zero Lebesgue measure.

For any $t \in \Lambda$, we have two possibilities: either
\begin{inparaenum}[\itshape i\upshape)]
\item $\|\sE'(u(t))\|_\sH = 0$, or
\item $\|\sE'(u(t))\|_\sH > 0$.
\end{inparaenum}
For Case (ii), we observe that the {\L}ojasiewicz-Simon gradient inequality \eqref{eq:Huang_3-15a} takes the shape,
\begin{equation}
\label{eq:Huang_3-16c}
\phi(H(t)) = c|\sE(u(t)) - \sE(\varphi)|^\theta \leq \|\sE'(u(t))\|_\sH ,
\end{equation}
and so
\begin{align*}
-\frac{d}{dt}\Phi(H(t)) &= -\frac{H'(t)}{\phi(H(t))} \quad\text{(by \eqref{eq:Huang_3-16b})}
\\
&= \frac{\|\sE'(u(t))\|_\sH\|\dot u(t)\|_\sH}{\phi(H(t))} \quad\text{(by  \eqref{eq:Huang_3-16a})}
\\
&\geq \frac{\|\sE'(u(t))\|_\sH\|\dot u(t)\|_\sH}{\|\sE'(u(t))\|_\sH}  \quad\text{(by \eqref{eq:Huang_3-16c})}
\\
&= \|\dot u(t)\|_\sH,
\end{align*}
that is,
\begin{equation}
\label{eq:Huang_3-16d}
-\frac{d}{dt}\Phi(H(t))  \geq \|\dot u(t)\|_\sH.
\end{equation}
Therefore, by the non-negativity of the function $-d\Phi(H(t))/dt$, combined with the fact that $[0, \infty) \less \Lambda$ has Lebesgue measure zero, we obtain the estimate,
\begin{equation}
\label{eq:Huang_3-18a_Hprelim}
-\frac{d}{dt}\Phi(H(t))  \geq \|\dot u(t)\|_\sH, \quad \hbox{a.e. } t \in [0, \infty),
\end{equation}
for both Cases (i) and (ii). Integration and the fact that $\lim_{t\to\infty}H(t) = a$ yields
$$
\int_0^\infty \|\dot u(t)\|_\sH\,dt \leq \Phi(H(0)) - \lim_{t\to\infty}\Phi(H(t)) = \Phi(H(0)) - \Phi(a) = \Phi(H(0)).
$$
By the definitions of $\Phi(x)$ and $H(t)$, this is \eqref{eq:Huang_3-15b_H}, since
$$
\Phi(H(0))
=
\int_a^{H(0)} \frac{1}{\phi(s)} \,ds
=
\int_{\sE_\infty}^{\sE(u(0))} \frac{1}{c|s - \sE(\varphi)|^\theta} \,ds.
$$
Furthermore, for any $\delta >0$ we have
$$
\int_\delta^\infty \|\dot u(t)\|_\sX\,dt \leq C_1(1+\delta^{-\rho})\int_0^\infty \|\dot u(t)\|_\sH\,dt
\quad\hbox{(by \eqref{eq:Abstract_apriori_interior_estimate_trajectory}),}
$$
and combining the preceding inequality with \eqref{eq:Huang_3-15b_H} yields \eqref{eq:Huang_3-15b_X}. Finiteness of the integrals $\int_0^\delta \|\dot u(t)\|_\sX\,dt$ and $\int_0^\infty \|\dot u(t)\|_\sX\,dt$ follow from the local Bochner integrability condition.
\end{proof}

\begin{rmk}[Evaluation of the integral in \eqref{eq:Huang_3-15b_H} and \eqref{eq:Huang_3-15b_X}]
We note that $\sE_\infty - \sE(\varphi) \leq \sE(u(0)) - \sE(\varphi)$, since $\sE(u(t))$ is a non-increasing function of $t \in [0, \infty)$, and separately consider each of the cases:
\begin{inparaenum}[\itshape i\upshape)]
\item $\sE_\infty - \sE(\varphi) \geq 0$;
\item $\sE(u(0)) - \sE(\varphi) \geq 0 \geq \sE_\infty - \sE(\varphi)$; and
\item $0 \geq \sE(u(0)) - \sE(\varphi)$.
\end{inparaenum}
The integral on the right-hand side of the inequality \eqref{eq:Huang_3-15b_H} is
$$
\int_{\sE_\infty}^{\sE(u(0))} \frac{1}{c|s - \sE(\varphi)|^\theta} \,ds
=
\int_{\sE_\infty - \sE(\varphi)}^{\sE(u(0)) - \sE(\varphi)} \frac{1}{c|r|^\theta} \,dr.
$$
Integration then yields
\begin{multline*}
\int_{\sE_\infty}^{\sE(u(0))} \frac{1}{c|s - \sE(\varphi)|^\theta} \,ds
\\
=
\frac{1}{c(1-\theta)}
\times
\begin{cases}
|\sE(u(0)) - \sE(\varphi)|^{1-\theta} - |\sE_\infty - \sE(\varphi)|^{1-\theta},
&\hbox{if } \sE_\infty - \sE(\varphi) \geq 0,
\\
|\sE(u(0)) - \sE(\varphi)|^{1-\theta} + |\sE_\infty - \sE(\varphi)|^{1-\theta},
&\hbox{if } \sE(u(0)) - \sE(\varphi) \geq 0 \geq \sE_\infty - \sE(\varphi),
\\
|\sE_\infty - \sE(\varphi)|^{1-\theta} -|\sE(u(0)) - \sE(\varphi)|^{1-\theta},
&\hbox{if } 0 \geq \sE(u(0)) - \sE(\varphi).
\end{cases}
\end{multline*}
Naturally, the calculation simplifies when it is known that $\sE(u(t)) \geq \sE(\varphi)$ for all $t \in [0, \infty)$.
\end{rmk}

\subsection{An abstract generalization of Simon's theorem on convergence of a subsequence implying convergence}
\label{subsec:Simon_convergence_subsequence_implies_convergence_in_time}
We have the following analogue of \cite[Theorems 3.3.3 and 3.3.6]{Huang_2006}, \cite[Corollary 2]{Simon_1983}.

\begin{thm}[Convergence of a subsequence implies convergence for a strong solution to a gradient system under the validity of the {\L}ojasiewicz-Simon gradient inequality]
\label{thm:Simon_corollary_2}
Let $\sU$ be an open subset of a real Banach space, $\sX$, that is continuously embedded and dense in a Hilbert space, $\sH$. Let $\sE:\sU\subset \sX\to\RR$ be a $C^1$ functional with gradient map $\sE':\sU\subset \sX \to \sH$. Assume that $\varphi \in \sU$ is a critical point of $\sE$, that is $\sE'(\varphi)=0$, and that the gradient map $\sE':\sU_\sigma\subset \sX \to \sX'$ satisfies a {\L}ojasiewicz-Simon gradient inequality \eqref{eq:Simon_2-2}, for positive constants $c$, $\sigma$, and $\theta \in [1/2,1)$. If $u:[0, \infty)\to \sU$ is a strong solution to \eqref{eq:Huang_3-3a} in the sense of Definition \ref{defn:Strong_solution_to_gradient_system}, the orbit $O(u) = \{u(t): t\geq 0\} \subset \sX$ is precompact\footnote{Recall that \emph{precompact} (or \emph{relatively compact}) subspace $Y$ of a topological space $X$ is a subset whose closure is compact. If the topology on $X$ is metrizable, then a subspace $Z \subset X$ is compact if and only if $Z$ is \emph{sequentially compact} \cite[Theorem 28.2]{Munkres2}, that is, every infinite sequence in $Z$ has a convergent subsequence in $Z$ \cite[Definition, p. 179]{Munkres2}.}, and $\varphi$ is a cluster point of $O(u)$, then the trajectory $u(t)$ converges to $\varphi$ as $t\to\infty$ in the sense that
\[
\lim_{t\to\infty}\|u(t)-\varphi\|_\sX = 0
\quad\hbox{and}\quad
\int_0^\infty \|\dot u\|_\sH\,dt < \infty.
\]
Furthermore, if $u$ satisfies Hypothesis \ref{hyp:Abstract_apriori_interior_estimate_trajectory} on $(0, \infty)$, then
\[
\int_1^\infty \|\dot u\|_\sX\,dt < \infty.
\]
\end{thm}

\begin{proof}
We shall adapt the proofs of \cite[Theorems 3.3.3 and 3.3.6]{Huang_2006}. Let $H:[0,\infty) \to \RR$ and $\Phi:\RR \to \RR$ be as in the proof of Proposition \ref{prop:Huang_3-3-2}. Note that the validity of the differential inequality \eqref{eq:Huang_3-16d} follows by the application of the {\L}ojasiewicz-Simon gradient inequality \eqref{eq:Huang_3-16c} at $u(t)$. Here, we assume that the {\L}ojasiewicz-Simon gradient inequality holds at $u(t)$ for all $t$ in the set
\[
W := \left\{t\geq 0: \|u(t) - \varphi\|_\sX < \sigma\right\}.
\]
Therefore, we have
\[
-\frac{d}{dt}\Phi(H(t)) \geq \|\dot u(t)\|_\sX, \quad\text{a.e. } t \in W.
\]
Integrating the preceding differential inequality yields
\[
\int_W \|\dot u(t)\|_\sX\,dt \leq \Phi(H(0)) - \lim_{t\to\infty}\Phi(H(t)) = \Phi(H(0)) < \infty,
\]
where we use the facts that $H(t) = \sE(u(t)) \leq \sE(u(0)) = H(0)$ (since the energy $\sE(u(t))$ is non-increasing when $u(t)$ is a solution to the gradient flow for $\sE$) and $\lim_{t\to\infty}\Phi(H(t)) = \Phi(a) = 0$, by definition of $\Phi(x)$ and $a = \lim_{t\to\infty} H(t)$. We can now apply Lemma \ref{lem:Huang_3-2-3} to show that the trajectory $u(t)$ converges to $\varphi$ as $t\to\infty$ in the sense that
\[
\lim_{t\to\infty}\|u(t)-\varphi\|_\sX = 0
\quad\hbox{and}\quad
\int_0^\infty \|\dot u\|_\sH\,dt < \infty.
\]
Finally, under Hypothesis \ref{hyp:Abstract_apriori_interior_estimate_trajectory} on $(0, \infty)$, we obtain
\[
\int_1^\infty \|\dot u\|_\sX\,dt < \infty,
\]
as desired.
\end{proof}

\subsection{Convergence for a solution to an abstract gradient system under the validity of the {\L}ojasiewicz-Simon gradient inequality}
\label{subsec:Convergence_gradient_system_orbit_within_Lojasiewicz-Simon_neighborhood}
We have the following simpler version of \cite[Lemma 3.3.4]{Huang_2006}; the result may be viewed as an abstract version of \cite[Lemma 1]{Simon_1983}.

\begin{lem}[Growth estimate for a strong solution to a gradient system under the validity of the {\L}ojasiewicz-Simon gradient inequality]
\label{lem:Huang_3-3-4}
Let $\sU$ be an open subset of a real Banach space, $\sX$, that is continuously embedded and dense in a Hilbert space, $\sH$. Let $\sE:\sU\subset \sX\to\RR$ be a $C^1$ functional with gradient map $\sE':\sU\subset \sX \to \sH$. Assume that the gradient map $\sE'$ satisfies a {\L}ojasiewicz-Simon gradient inequality \eqref{eq:Simon_2-2} with positive constants $c$, $\sigma$, and $\theta \in [1/2, 1)$, and critical point $\varphi \in \sU$. If $T>0$ and $u \in C([0, T]; \sU)$ is a strong solution to the Cauchy problem \eqref{eq:Huang_3-3a} in the sense of Definition \ref{defn:Strong_solution_to_gradient_system} and the gradient map $\sE'$ satisfies the {\L}ojasiewicz-Simon gradient inequality in the orbit $O(u)$ as in \eqref{eq:Huang_3-15a}, then
\begin{subequations}
\label{eq:Huang_3-18_H}
\begin{align}
\label{eq:Huang_3-18a_H}
\|\dot u(t)\|_\sH &\leq -\frac{d}{dt} \int_{\sE(u(T))}^{\sE(u(t))} \frac{1}{c|s - \sE(\varphi)|^\theta} \,ds,
 \quad\hbox{a.e. } t \in (0, T),
\\
\label{eq:Huang_3-18b_H}
\int_0^T\|\dot u(t)\|_\sH\,dt &\leq \int_{\sE(u(T))}^{\sE(u(0))} \frac{1}{c|s - \sE(\varphi)|^\theta} \,ds.
\end{align}
\end{subequations}
If in addition $u$ obeys Hypothesis \ref{hyp:Abstract_apriori_interior_estimate_trajectory}, then, for all $\delta \in (0, T]$,
\begin{equation}
\label{eq:Huang_3-18b_X_interior}
\int_\delta^T\|\dot u(t)\|_\sX
\leq
C_1(1+\delta^{-\rho}) \int_{\sE(u(T))}^{\sE(u(0))} \frac{1}{c|s - \sE(\varphi)|^\theta} \,ds.
\end{equation}
\end{lem}

\begin{proof}
The proof is the same as that of Proposition \ref{prop:Huang_3-3-2}, the only difference being the assumption that $0 < T < \infty$ and replacement of the role played by $a = \sE_\infty$ by that of $\sE(u(T))$. The inequality \eqref{eq:Huang_3-18a_H} follows from the inequality \eqref{eq:Huang_3-18a_Hprelim} and the definitions of $\Phi(x)$ and $H(t)$ in the proof of Proposition \ref{prop:Huang_3-3-2}.
\end{proof}

An immediate consequence of Lemma \ref{lem:Huang_3-3-4} is a simpler version of \cite[Theorem 3.3.5]{Huang_2006} which we again include for completeness, though it is not used in the sequel. Theorem \ref{thm:Huang_3-3-5} gives a simple convergence result for a solution $u$ under the strong hypothesis that the orbit $O(u)$ belongs to a {\L}ojasiewicz-Simon neighborhood.

\begin{thm}[Convergence for a strong solution to a gradient system under the validity of the {\L}ojasiewicz-Simon gradient inequality]
\label{thm:Huang_3-3-5}
Let $\sU$ be an open subset of a real Banach space, $\sX$, that is continuously embedded and dense in a Hilbert space, $\sH$. Let $\sE:\sU\subset \sX\to\RR$ be a $C^1$ functional with gradient map $\sE':\sU\subset \sX \to \sH$ and $u:[0, \infty)\to \sU$ be a strong solution to the gradient system \eqref{eq:Huang_3-3a} in the sense of Definition \ref{defn:Strong_solution_to_gradient_system}. If the gradient map $\sE'$ satisfies a {\L}ojasiewicz-Simon gradient inequality \eqref{eq:Simon_2-2} in the orbit $O(u)$, then $u(t)$ converges in $\sH$ as $t\to\infty$ in the sense that
$$
\int_0^\infty \|\dot u(t)\|_\sH\,dt < \infty.
$$
If in addition $u$ obeys Hypothesis \ref{hyp:Abstract_apriori_interior_estimate_trajectory}, then $u(t)$ converges in $\sX$ as $t\to\infty$ in the sense that
$$
\int_1^\infty \|\dot u(t)\|_\sX\,dt < \infty.
$$
\end{thm}

\subsection{Simon Alternative and convergence for an abstract gradient system}
\label{subsec:Simon_alternative}
We next have the following abstract analogue of R\r{a}de's \cite[Proposition 7.4]{Rade_1992}, in turn a variant the \emph{Simon Alternative}, namely \cite[Theorem 2]{Simon_1983}.

\begin{thm}[Simon Alternative for convergence for a strong solution to a gradient system]
\label{thm:Huang_3-3-6}
Let $\sU$ be an open subset of a real Banach space, $\sX$, that is continuously embedded and dense in a Hilbert space, $\sH$. Let $\sE:\sU\subset \sX\to\RR$ be a $C^1$ functional with gradient map $\sE':\sU\subset \sX \to \sH$. Assume that
\begin{enumerate}
\item $\varphi \in \sU$ is a critical point of $\sE$, that is $\sE'(\varphi)=0$, and that the gradient map $\sE':\sU_\sigma\subset \sX \to \sX'$ satisfies a {\L}ojasiewicz-Simon gradient inequality \eqref{eq:Simon_2-2}, for positive constants $c$, $\sigma$, and $\theta \in [1/2,1)$;

\item Given positive constants $b$, $\eta$, and $\tau$, there is a constant $\delta = \delta(\eta, \tau, b) \in (0, \tau]$ such that if $v$ is a solution to the gradient system \eqref{eq:Huang_3-3a} on $[t_0, t_0 + \tau)$ with $\|v(t_0)\|_\sX \leq b$, then
\begin{equation}
\label{eq:Gradient_solution_near_initial_data_at_t0_for_short_enough_time}
\sup_{t\in [t_0, t_0+\delta]}\|v(t) - v(t_0)\|_\sX < \eta.
\end{equation}
\end{enumerate}
Then there is a constant $\eps = \eps(c, C_1, \delta, \theta, \rho, \sigma, \tau, \varphi) \in (0, \sigma/4)$ with the following significance.  If $u:[0, \infty)\to \sU$ is a strong solution to \eqref{eq:Huang_3-3a} in the sense of Definition \ref{defn:Strong_solution_to_gradient_system} that satisfies Hypothesis \ref{hyp:Abstract_apriori_interior_estimate_trajectory} on $(0, \infty)$ and there is a $T \geq 0$ such that
\begin{equation}
\label{eq:Rade_7-2_banach}
\|u(T) - \varphi\|_\sX < \eps,
\end{equation}
then either
\begin{enumerate}
\item
\label{item:Theorem_3-3-6_energy_u_at_time_t_below_energy_critical_point}
$\sE(u(t)) < \sE(\varphi)$ for some $t>T$, or
\item
\label{item:Theorem_3-3-6_u_converges_to_limit_u_at_infty}
the trajectory $u(t)$ converges in $\sX$ to a limit $u_\infty \in \sX$ as $t\to\infty$ in the sense that
$$
\lim_{t\to\infty}\|u(t)-u_\infty\|_\sX =0
\quad\hbox{and}\quad
\int_1^\infty \|\dot u\|_\sX\,dt < \infty.
$$
If $\varphi$ is a cluster point of the orbit $O(u) = \{u(t): t\geq 0\}$, then $u_\infty = \varphi$.
\end{enumerate}
\end{thm}

\begin{rmk}[On the role of the constant $\tau$]
\label{rmk:Huang_3-3-6}
In our applications, the constant, $\tau$, in the hypotheses of Theorem \ref{thm:Huang_3-3-6} will depend on $u(t_0)$ through at most an upper bound, $b > 0$, for $\|u(t_0)\|_\sX$. In the present situation, we have $\|u(T)\|_\sX \leq \|u(T) - \varphi\|_\sX + \|\varphi\|_\sX < \eps + \|\varphi\|_\sX < \sigma/4 + \|\varphi\|_\sX$,
where we use the fact that the positive constant, $\eps$, in Theorem \ref{thm:Huang_3-3-6} belongs to $(0, \sigma/4)$.
\end{rmk}

\begin{rmk}[Application to Yang-Mills gradient flow over a closed manifold of dimension two or three]
\label{rmk:Rade_propositions_7-1_and_7-4}
In the context of Yang-Mills gradient flow $A(t)$ on a principal $G$-bundle $P$ over a closed manifold $X$ of dimension two or three, Theorem \ref{thm:Simon_corollary_2} is useful because one can appeal to Uhlenbeck convergence without bubbling and the theory of parabolic partial differential equations, as in R\r{a}de's \cite[Proposition 7.1]{Rade_1992} and its proof, to obtain precompactness modulo gauge transformations and thus an analogue of \eqref{eq:Rade_7-2_banach} in R\r{a}de's \cite[Proposition 7.4]{Rade_1992}. Indeed, given any  unbounded sequence, $\{t_m\}_{m=0}^\infty \subset [0,\infty)$, there is subsequence relabelled as $\{t_m\}_{m=0}^\infty$, a sequence of gauge transformations, $\{\Phi_m\}_{m=0}^\infty$ of $P$, of Sobolev class $H^2$, and a $C^\infty$ Yang-Mills connection $A_\infty$ on $P$, such that $t_m \to \infty$ and $\Phi_m^*A(t_m) \to A_\infty$ strongly in the $H^1$ topology as $m\to\infty$. Hence, for large enough $m$ and $T:=t_m$, we have
\[
\|A(T) - \tilde A_\infty\|_{H_{A_1}^1(X)} < \eps,
\]
where $\tilde A_\infty := (\Phi_m^{-1})^*A_\infty$ and $A_1$ is a $C^\infty$ reference connection on $P$.
\end{rmk}

\begin{proof}[Proof of Theorem \ref{thm:Huang_3-3-6}]
We shall assume that $\sE(u(t)) \geq \sE(\varphi)$ for all $t \in [T, \infty)$ and thus aim to establish the Alternative \eqref{item:Theorem_3-3-6_u_converges_to_limit_u_at_infty} in the statement of the proposition. Recall that $\sE:\sU\subset \sX\to \RR$ is $C^1$ map and that $\varphi$ is a critical point of $\sE$. Consequently, by definition of the Fr\'echet derivative,
\begin{align*}
\sE(v) &= \sE(\varphi) + \langle\sE'(\varphi), v-\varphi\rangle_{\sX'\times\sX} + o(\|v - \varphi\|_\sX)
\\
&= \sE(\varphi) + o(\|v - \varphi\|_\sX), \quad \forall\, v \in \sU,
\end{align*}
and so, for a small enough positive constant $\eps_\varphi$, we have\footnote{If we assume that $\sE$ is $C^2$, then we can apply the Taylor formula with explicit form of the remainder \cite[Proposition 2.1.33]{Berger_1977} to write $o(\|v - \varphi\|_\sX) = M\|v - \varphi\|_\sX^2$, where $M := \sup_{s\in [0,1]}\|\sE''((1-s)\varphi + sv)\|$.}
\begin{equation}
\label{eq:Rade_7-4_banach}
|\sE(v) - \sE(\varphi)| \leq \|v - \varphi\|_\sX, \quad\forall\, v \in \sX \hbox{ such that } \|v-\varphi\|_\sX < \eps_\varphi.
\end{equation}
We make the

\begin{claim}[Confinement of solution to a {\L}ojasiewicz-Simon gradient inequality neighborhood]
\label{claim:norm_X_of_u_minus_varphi_lessthan_sigma_all_t_geq_T}
There is a constant
$\eps = \eps(c,C_1,\delta, \theta, \rho, \sigma, \tau, \varphi) \in (0, \sigma/4)$ in \eqref{eq:Rade_7-2_banach} such that
\begin{equation}
\label{eq:norm_X_of_u_minus_varphi_lessthan_sigma_all_t_geq_T}
\|u(t) - \varphi\|_\sX < \frac{\sigma}{2}, \quad\forall\, t \in [T, \infty).
\end{equation}
\end{claim}

\begin{proof}[Proof of Claim \ref{claim:norm_X_of_u_minus_varphi_lessthan_sigma_all_t_geq_T}]
We suppose the opposite and obtain a contradiction. Let $\widehat T > T$ be the smallest number such that
\begin{equation}
\label{eq:Rade_7-6_banach}
\|u(t) - \varphi\|_\sX < \frac{\sigma}{2} \quad\forall\, t \in [T, \widehat T),
\quad\hbox{but}\quad
\|u(\widehat T) - \varphi\|_\sX \geq \frac{\sigma}{2}.
\end{equation}
By \eqref{eq:Gradient_solution_near_initial_data_at_t0_for_short_enough_time} with $t_0 = T$ and $\eta = \sigma/4$, we may choose $\delta = \delta(\sigma, \tau, \|\varphi\|_\sX) \in (0, \tau]$ so that
\begin{equation}
\label{eq:Rade_7-5_banach}
\sup_{t\in [T, T+\delta]}\|u(t) - u(T)\|_\sX < \frac{\sigma}{4},
\end{equation}
and, in particular,
$$
\|u(T+\delta) - u(T)\|_\sX < \frac{\sigma}{4}.
$$
By hypothesis \eqref{eq:Rade_7-2_banach} and the fact that we seek $\eps \in (0, \sigma/4)$, we also have
$$
\|u(T) - \varphi\|_\sX < \frac{\sigma}{4},
$$
and combining this inequality with \eqref{eq:Rade_7-5_banach} gives
$$
\sup_{t\in [T, T+\delta]}\|u(t) - \varphi\|_\sX < \frac{\sigma}{2}.
$$
Therefore, by definition of $\widehat T$ in \eqref{eq:Rade_7-6_banach}, we must have
$$
T > T + \delta.
$$
Inequality \eqref{eq:Huang_3-18a_H}, with $[0,T]$ replaced by $[T,\widehat T]$, gives
$$
\|\dot u(t)\|_\sH \leq -\frac{d}{dt}\frac{1}{c(1-\theta)}(\sE(u(t)) - \sE(\varphi))^{1-\theta},
\quad\hbox{a.e. } t \in [T,\widehat T].
$$
Integrating this inequality yields \eqref{eq:Huang_3-18b_H}, with $[0,T]$ replaced by $[T,\widehat T]$,
\begin{equation}
\label{eq:Huang_3-18b_H_with_0toT_replaced_by_TtotildeT}
\int_T^{\widehat T}\|\dot u(t)\|_\sH\,dt
\leq
\frac{1}{c(1-\theta)}\left((\sE(u(T)) - \sE(\varphi))^{1-\theta} - (\sE(u(\widehat T)) - \sE(\varphi))^{1-\theta}\right).
\end{equation}
Therefore, discarding the negative term on the right-hand side of the preceding inequality,
\begin{align*}
\int_T^{\widehat T}\|\dot u(t)\|_\sH\,dt &\leq \frac{1}{c(1-\theta)}\left(\sE(u(T)) - \sE(\varphi)\right)^{1-\theta}
\\
&\leq \frac{1}{c(1-\theta)}\|u(T) - \varphi\|_\sX^{1-\theta} \quad\hbox{(applying \eqref{eq:Rade_7-4_banach} with $v=u(T)$),}
\end{align*}
and thus, by \eqref{eq:Rade_7-2_banach} and $\eps\leq\eps_\varphi$,
\begin{equation}
\label{eq:Bound_on_integral_from_T_to_tildeT_of_Hnorm_dotu}
\int_T^{\widehat T}\|\dot u(t)\|_\sH\,dt \leq \frac{\eps^{1-\theta}}{c(1-\theta)}.
\end{equation}
On the other hand, noting that $T+\delta < \widehat T$,
\begin{align*}
\int_{T+\delta}^{\widehat T}\|\dot u(t)\|_\sX\,dt &\geq \|u(\widehat T) - u(T+\delta)\|_\sX
\\
&\geq \|u(\widehat T) - \varphi\|_\sX - \|u(T+\delta) - u(T)\|_\sX - \|u(T) - \varphi\|_\sX
\\
&> \frac{\sigma}{2} - \frac{\sigma}{4} - \eps  \quad\hbox{(by \eqref{eq:Rade_7-2_banach}, \eqref{eq:Rade_7-6_banach}, and \eqref{eq:Rade_7-5_banach})}.
\end{align*}
The preceding inequality implies that
\begin{align*}
\frac{\sigma}{4} - \eps &< \int_{T+\delta}^{\widehat T}\|\dot u(t)\|_\sX\,dt
\\
&\leq C_1(1+\delta^{-\rho})\int_T^{\widehat T}\|\dot u(t)\|_\sH\,dt  \quad\hbox{(by \eqref{eq:Abstract_apriori_interior_estimate_trajectory})}
\\
&\leq \frac{C_1(1+\delta^{-\rho})}{c(1-\theta)}\eps^{1-\theta}  \quad\hbox{(by \eqref{eq:Bound_on_integral_from_T_to_tildeT_of_Hnorm_dotu})}.
\end{align*}
Choosing $\eps = \eps(c, C_1, \delta, \theta, \rho, \sigma, \tau, \varphi) \in (0, \sigma/4)$ small enough that $\eps + (C_1(1+\delta^{-\rho})/c(1-\theta))\eps^{1-\theta} < \sigma/4$ leads to a contradiction and this completes the proof of Claim \ref{claim:norm_X_of_u_minus_varphi_lessthan_sigma_all_t_geq_T}.
\end{proof}

Therefore, because we may now set $\widehat T = \infty$ in \eqref{eq:Huang_3-18b_H_with_0toT_replaced_by_TtotildeT} and discard the negative term, we obtain the inequality,
\begin{equation}
\label{eq:Rade_7-7H_banach}
\int_T^\infty\|\dot u(t)\|_\sH\,dt
\leq
\frac{1}{c(1-\theta)}\left(\sE(u(T)) - \sE(\varphi)\right)^{1-\theta}.
\end{equation}
From the \apriori estimate \eqref{eq:Abstract_apriori_interior_estimate_trajectory}, we have
\begin{equation}
\label{eq:Rade_7-7X_banach}
\int_{T+1}^\infty\|\dot u(t)\|_\sX\,dt \leq 2C_1\int_T^\infty\|\dot u(t)\|_\sH\,dt.
\end{equation}
Combining \eqref{eq:Rade_7-7H_banach} and \eqref{eq:Rade_7-7X_banach} yields
\begin{equation}
\label{eq:Rade_7-7_banach}
\int_{T+1}^\infty\|\dot u(t)\|_\sX\,dt \leq \frac{2C_1}{c(1-\theta)}\left(\sE(u(T)) - \sE(\varphi)\right)^{1-\theta}.
\end{equation}
In particular,
$$
\int_{T+1}^\infty\|\dot u(t)\|_\sX\,dt < \infty.
$$
Therefore, $u(t)$ converges in $\sX$ to a limit $u_\infty \in \sX$ as $t \to \infty$. Indeed, for $\eps>0$ and any sequence of times $\{t_n\}_{n=1}^\infty \subset [0,\infty)$ with $t_n\to\infty$ as $n\to\infty$, we have
\begin{align*}
\|u(t_n) - u(t_m)\|_\sX &\leq \int_{t_m\wedge t_n}^{t_m\vee t_n} \|\dot u(t)\|_\sX\,dt
\\
&\leq \int_{t_m\wedge t_n}^\infty \|\dot u(t)\|_\sX\,dt < \eps, \quad\forall\, n, m \geq N,
\end{align*}
for some $N=N(\eps)$. Hence, the sequence $\{u(t_n)\}_{n=1}^\infty$ is Cauchy in $\sX$ and converges to a limit $u_\infty \in \sX$ as $n\to\infty$. Consequently,
$$
\|u(t) - u_\infty\|_\sX \leq \int_t^\infty \|\dot u(t)\|_\sX\,dt, \quad\forall\, t \geq 0,
$$
and so $u(t) \to u_\infty$ in $\sX$ as $t\to\infty$. If $\varphi$ is a cluster point of $O(u)$, then there is an unbounded sequence $\{s_n\}_{n\geq 0} \subset [0,\infty)$ such that $u(s_n) \to \varphi$ in $\sX$ as $n\to\infty$ and therefore $u_\infty = \varphi$.
\end{proof}

\subsection{Convergence rate}
\label{subsec:Huang_3-4C}
We essentially follow \cite[Section 3.4C]{Huang_2006}, with one enhancement to provide a convergence rate in $\sX$ as well as in $\sH$. We consider the gradient system,
\begin{equation}
\label{eq:Huang_3-42}
\dot u(t) = -\sE'(u(t)),
\end{equation}
for $u :[0, \infty) \to \sX$, where $\sX$ is a Banach space and $\sH$ a Hilbert space with continuous embeddings, $\sX \hookrightarrow \sH \hookrightarrow \sX'$. We assume that $\sU \subset \sX$ is an open subset and $a, \gamma$ are constants such that $\gamma > a$ and
\begin{equation}
\label{eq:Huang_3-43a}
a \leq \sE(v)  \leq \gamma, \quad\forall\, v \in \sU.
\end{equation}
Moreover, we assume that $\sE$ satisfies a {\L}ojasiewicz-Simon gradient inequality \eqref{eq:Simon_2-2} in $\sU$ with positive constants $c$ and $\theta \in [1/2,1)$. For example, this holds when $\sE:\sU\to\RR$ is analytic (by Theorem \ref{thm:Huang_2-4-2} or \ref{thm:Huang_2-4-5}) and $\varphi \in \sU = \sU_\sigma$ is a critical point with critical value $a$, that is, $\sE'(\varphi) = 0$ and $\sE(\varphi) = a$. It follows that (compare \cite[Example 3.4.9]{Huang_2006}, where $a=0$)
$$
\int_x^\gamma \frac{ds}{c^2(s-a)^{2\theta}}
=
\begin{cases}
\displaystyle\frac{1}{c^2(2\theta-1)}\left((x-a)^{1-2\theta} - (\gamma-a)^{1-2\theta}\right), & 1/2 < \theta \leq 1,
\\
\displaystyle\frac{1}{c^2}(\ln (\gamma-a) - \ln (x-a)), & \theta = 1/2,
\end{cases}
$$
defines a strictly \emph{decreasing} function of $x \in (a, \gamma]$. When $1/2 < \theta \leq 1$, solving
$$
t =\frac{1}{c^2(2\theta-1)}\left((x-a)^{1-2\theta} - (\gamma-a)^{1-2\theta}\right)
$$
for $x = g(t)$ yields
\begin{equation}
\label{eq:Huang_3-44a_theta_less_than_half}
g(t) = a + \left(c^2(2\theta-1)t + (\gamma-a)^{1-2\theta}\right)^{-1/(2\theta-1)}, \quad t \in [0, \infty).
\end{equation}
When $\theta = 1/2$, solving
$$
t =\frac{1}{c^2}(\ln (\gamma-a) - \ln (x-a))
$$
for $x = g(t)$ yields $g(t) = a + \exp(\ln(\gamma-a) - c^2t)$, that is,
\begin{equation}
\label{eq:Huang_3-44a_theta_is_half}
g(t) = a + (\gamma-a)\exp(-c^2t), \quad t \in [0, \infty).
\end{equation}
Observe that, for either case, $g(0) = \gamma$ and $g(t) \to a$ as $t \to \infty$. For $1/2 \leq \theta < 1$, let
\begin{equation}
\label{eq:Huang_3-44b}
\Phi(x) := \int_a^x \frac{ds}{c (s-a)^\theta} = \frac{(x-a)^{1-\theta}}{c(1-\theta)},
\quad x \in [a, \infty).
\end{equation}
This defines a strictly increasing function of $x\geq a$.

We have the following enhancement of \cite[Theorem 3.4.8]{Huang_2006}.

\begin{thm}[Convergence rate under the validity of a {\L}ojasiewicz-Simon gradient inequality]
\label{thm:Huang_3-4-8}
Let $\sU$ be an open subset of a real Banach space, $\sX$, that is continuously embedded and dense in a Hilbert space, $\sH$. Let $\sE:\sU\subset \sX\to\RR$ be a $C^1$ functional with gradient map $\sE':\sU\subset \sX \to \sH$ and suppose that $\sE'$ obeys a {\L}ojasiewicz-Simon gradient inequality \eqref{eq:Simon_2-2} with positive constants $c$, $\sigma$, and $\theta \in [1/2, 1)$. Let $u:[0,\infty) \to \sX$ be a strong solution to the gradient system \eqref{eq:Huang_3-42}, in the sense of Definition \ref{defn:Strong_solution_to_gradient_system}, and assume that $O(u) \subset \sU_\sigma \subset \sU$. Then there exists a $u_\infty \in \sH$ such that
\begin{equation}
\label{eq:Huang_3-45_H}
\|u(t) - u_\infty\|_\sH \leq \Phi(g(t)), \quad t\geq 0,
\end{equation}
where
\begin{equation}
\label{eq:Huang_3-45_growth_rate}
\Phi(g(t))
=
\begin{cases}
\displaystyle
\frac{1}{c(1-\theta)}\left(c^2(2\theta-1)t + (\gamma-a)^{1-2\theta}\right)^{-(1-\theta)/(2\theta-1)},
& 1/2 < \theta < 1,
\\
\displaystyle
\frac{2}{c}\sqrt{\gamma-a}\exp(-c^2t/2),
&\theta = 1/2.
\end{cases}
\end{equation}
If in addition $u$ obeys Hypothesis \ref{hyp:Abstract_apriori_interior_estimate_trajectory}, then $u_\infty \in \sX$ and
\begin{equation}
\label{eq:Huang_3-45_X}
\|u(t+1) - u_\infty\|_\sX \leq 2C_1\Phi(g(t)), \quad t\geq 0,
\end{equation}
where $\Phi(g(t))$ is as in \eqref{eq:Huang_3-45_growth_rate} and $C_1$ is the positive constant in \eqref{eq:Abstract_apriori_interior_estimate_trajectory}.
\end{thm}

\begin{proof}
The existence of $u_\infty \in \sH$ and validity of the estimate \eqref{eq:Huang_3-45_H} is given by \cite[Theorem 3.4.8]{Huang_2006}. The \apriori estimate \eqref{eq:Abstract_apriori_interior_estimate_trajectory} asserts that, for $\delta=1$, $S=t$, and $T=\infty$,
$$
\int_{t+1}^\infty \|\dot u(s)\|_\sX \,ds \leq 2C_1\int_t^\infty \|\dot u(s)\|_\sH \,ds, \quad t \geq 0,
$$
while \cite[Inequality (3.46c)]{Huang_2006} and the penultimate line of the proof of \cite[Theorem 3.4.8]{Huang_2006} gives
$$
\int_t^\infty \|\dot u(s)\|_\sH \,ds \leq \Phi(g(t)), \quad t \geq 0.
$$
Combining the preceding two inequalities yields
$$
\int_{t+1}^\infty \|\dot u(s)\|_\sX \,ds \leq 2C_1\Phi(g(t)), \quad t \geq 0.
$$
As in the proof of \cite[Theorem 3.4.8]{Huang_2006}, finiteness of the integral $\int_0^\infty \|\dot u(s)\|_\sH \,ds$ implies that $\|u(t) - u_\infty\|_\sH \to 0$ as $t \to \infty$, for some $u_\infty \in \sH$. Similarly, finiteness of the integral $\int_0^\infty \|\dot u(s)\|_\sX \,ds$ implies that $\|u(t) - v_\infty\|_\sX \to 0$ as $t \to \infty$, for some $v_\infty \in \sX$. Since $\sX\hookrightarrow \sH$, we must have $u_\infty = v_\infty$ and so $u_\infty \in \sX$. Moreover,
$$
\|u(t+1) - u_\infty\|_\sX \leq \int_{t+1}^\infty \|\dot u(s)\|_\sX \,ds, \quad t \geq 0,
$$
and combining the preceding inequalities gives \eqref{eq:Huang_3-45_X}.
\end{proof}

\subsection{Existence and convergence of a global solution to an abstract gradient system started near a local minimum}
\label{subsec:Huang_5-1_existence_and_convergence_global_solution_near_local_minimum}
In this subsection, we present a simplification of the first part of \cite[Section 5.1]{Huang_2006}, concerning existence and convergence of global solutions to an abstract gradient system. We restrict our attention to the case of a $C^1$ potential function, $\sE:\sU\subset \sX\to\RR$, which obeys a \L ojasiewicz-Simon gradient inequality rather than the more general gradient inequalities allowed by Huang in \cite{Huang_2006}. We also present a modification of Huang's treatment that will be important in our application to Yang-Mills gradient flow and its perturbations. Rather than consider the gradient system \cite[Equation (5.1a)]{Huang_2006},
$$
\dot u(t) = -\sN(t, u(t)), \quad t > 0, \quad u(0) = u_0,
$$
for a solution $u:[0, T) \to \sX$, where $\sX$ is a real Banach space and $\sN:[0, \infty) \times \sX \to \sX$ is a continuous map that obeys a local Lipschitz property \cite[Equation (5.1b)]{Huang_2006} and local growth conditions specified by \cite[Equation (5.3) and (5.5)]{Huang_2006}, we shall instead consider the pure gradient flow equation,
$$
\dot u(t) = -\sE'(u(t)), \quad t > 0, \quad u(0) = u_0,
$$
and assume that any solution obeys a suitable \apriori estimate. This modification is important since the growth condition \cite[Equation (5.3)]{Huang_2006} appears to us to be rather difficult to verify in practice whereas, in the case of the Yang-Mills gradient flow, we shall obtain suitable \apriori estimates in the sequel in Lemmata \ref{lem:Rade_7-3} and \ref{lem:Rade_7-3_L1_in_time_H2beta_in_space_apriori_estimate_by_L1_in_time_L2_in_space}. In the sequel, we shall also consider a situation where $u$ is a solution to a small perturbation of the abstract pure gradient flow equation and then apply those abstract results to the case of a small perturbation of the Yang-Mills gradient flow equation.

Our interest in \cite[Section 5.1]{Huang_2006} is due to the fact that \cite[Theorems 5.1.1 and 5.1.2]{Huang_2006} provide global existence, Lyapunov stability, and asymptotic convergence to ground states under certain hypotheses. These provide the abstract results closest to those obtained by R\r{a}de \cite{Rade_1992} for the Yang-Mills gradient flow over a closed, Riemannian, smooth manifold of dimension two or three.

Recall that $\sS$ in \eqref{eq:Huang_page_25_defn_critical_point_set} is the set of critical points of $\sE$. We call a critical point $\varphi \in \sU$ of $\sE$ a \emph{ground state} if $\sE$ attains its minimum on $\sU$ at this point, that is,
\begin{equation}
\label{eq:Huang_5-2}
\sE(\varphi) = \inf_{u\in \sU}\sE(u).
\end{equation}
For $r>0$, we define an open neighborhood of a critical point, $\varphi$, by
\begin{equation}
\label{eq:Huang_page_164_Usigma}
\sU_r := \{u\in \sX: \|u-\varphi\|_\sX < r\}.
\end{equation}
We have the following analogue of \cite[Theorem 5.1.1]{Huang_2006}; our proof is similar to that of Theorem \ref{thm:Huang_3-3-6}.

\begin{thm}[Existence and convergence of a global solution to a gradient system near a local minimum]
\label{thm:Huang_5-1-1}
Let $\sU$ be an open subset of a real Banach space, $\sX$, that is continuously embedded and dense in a Hilbert space, $\sH$. Let $\sE:\sU\subset \sX\to\RR$ be a $C^1$ functional with gradient map $\sE':\sU\subset \sX \to \sH$. Let $\varphi \in \sU$ be a ground state of $\sE$ on $\sU$ and suppose that $\sE'$ obeys a {\L}ojasiewicz-Simon gradient inequality \eqref{eq:Simon_2-2} with positive constants $c$, $\sigma$, and $\theta \in [1/2, 1)$.
Assume that
\begin{enumerate}
\item For each $u_0 \in \sU$, there exists a unique strong solution to the Cauchy problem \eqref{eq:Huang_3-3a}, in the sense of Definition \ref{defn:Strong_solution_to_gradient_system}, on a time interval $[0, \tau)$ for some positive constant, $\tau$;

\item Hypothesis \ref{hyp:Abstract_apriori_interior_estimate_trajectory} holds with positive constants $C_1$ and $\rho$ for strong solutions to the gradient system \eqref{eq:Huang_3-3a}; and

\item Given positive constants $b$ and $\eta$, there is a constant $\delta = \delta(\eta, \tau, b) \in (0, \tau]$ such that if $v$ is a solution to the gradient system \eqref{eq:Huang_3-3a} on $[0, \tau)$ with $\|v_0\|_\sX \leq b$, then
\begin{equation}
\label{eq:Gradient_solution_near_initial_data_at_time_zero_for_short_enough_time}
\sup_{t\in [0, \delta]}\|v(t) - v_0\|_\sX < \eta.
\end{equation}
\end{enumerate}
Then there is a constant $\eps = \eps(c,C_1,\delta, \theta, \rho, \sigma, \tau, \varphi) \in (0, \sigma/4)$ with the following significance. For each $u_0 \in \sU_\eps$, the Cauchy problem \eqref{eq:Huang_3-3a} admits a global strong solution, $u:[0,\infty) \to \sU_{\sigma/2}$, that converges to a limit $u_\infty \in \sX$ as $t\to\infty$ with respect to the $\sX$ norm in the sense that
$$
\lim_{t \to \infty} \|u(t) - u_\infty\|_\sX = 0 \quad\hbox{and}\quad \int_1^\infty\|\dot u(t)\|_\sX\,dt < \infty.
$$
\end{thm}

\begin{rmk}[On the role of the minimal lifetime, $\tau$, of a solution with initial data, $u_0$]
\label{rmk:Huang_5-1-1}
In our applications, the minimal lifetime, $\tau$, of a solution $u$ to the gradient system \eqref{eq:Huang_3-3a} will depend on $u_0$ through at most an upper bound, $b > 0$, for $\|u_0\|_\sX$. In the present situation, we have $\|u_0\|_\sX \leq \|u_0 - \varphi\|_\sX + \|\varphi\|_\sX < \eps + \|\varphi\|_\sX < \sigma/4 + \|\varphi\|_\sX$,
where we use the fact that the positive constant, $\eps$, in Theorem \ref{thm:Huang_5-1-1} belongs to $(0, \sigma/4)$.
\end{rmk}

It will be useful at this point to recall a special case of the general properties of vector-valued Sobolev spaces of functions and vector-valued absolutely continuous functions on bounded intervals $(a, b) \subset \RR$ \cite[Appendix C]{Sell_You_2002}, with proofs provided in \cite{Barbu_Precupanu_2012}, for example. Closely related results are provided in \cite[Section 3.1]{Showalter} --- see Remark \ref{rmk:Absolutely_continuous_paths_Banach_spaces}.

\begin{thm}[Vector-valued Sobolev spaces of functions and vector-valued absolutely continuous functions on bounded intervals]
\label{thm:Sell_You_C-9}
Let $\sX$ be a Banach space and $(a, b) \subset \RR$ a bounded interval. Then $u \in W^{1,1}(a,b;\sX)$ if and only if there is an absolutely continuous function $v \in C([a,b]; \sX)$ such that $\dot v \in L^1(a,b; \sX)$ and $u = v$ a.e. on $[a, b]$. Moreover, if $u \in W^{1,1}(a,b;\sX)$, then\footnote{After identifying $u$ and $v$ on $[a,b]$.} the \emph{Newton-Leibnitz formula} holds, that is,
\begin{equation}
\label{eq:Sell_You_93-1}
u(t) = u(a) + \int_a^t \dot u(s) \,ds, \quad\forall\, t \in [a, b].
\end{equation}
\end{thm}

\begin{proof}
The first assertion is provided by \cite[Theorem C.9]{Sell_You_2002} with $m = 1$ and $p = 1$. The second assertion is provided by \cite[Lemma C.5]{Sell_You_2002}.
\end{proof}

We will need the following version of a well-known idea used in the proof, for example, of \cite[Theorem 4.1]{Neuberger_2010}. Compare also Theorem \ref{thm:Sell_You_lemma_47-4}.

\begin{lem}[Global existence of a strong solution to a gradient flow equation]
\label{lem:Neuberger_theorem_4-1}
Let $\sU$ be an open subset of a real Banach space, $\sX$, that is continuously embedded and dense in a Hilbert space, $\sH$. Let $\sE:\sU\subset \sX\to\RR$ be a $C^1$ functional with gradient map $\sE':\sU\subset \sX \to \sH$. Assume that for each $v_0 \in \sU$, there is a maximal time, $T(v_0)$, such that there exists a strong solution, $v$, to \eqref{eq:Huang_3-3a}, in the sense of Definition \ref{defn:Strong_solution_to_gradient_system}, on the interval, $[0, T(v_0))$. Let $u_0 \in \sU$ and let $u$ be a strong solution to \eqref{eq:Huang_3-3a} on $[0, T(u_0))$. If there is an open subset $\sV \subset \sU$ such that $\bar\sV \subset \sU$ and
$$
\{u(t): t \in [0, T(u_0))\} \subset \sV,
$$
and $u$ is uniformly continuous on $[0, T(u_0))$, then $T(u_0) = \infty$.
\end{lem}

\begin{proof}
Suppose $T := T(u_0) < \infty$. By hypothesis, $u$ is uniformly continuous on $[0, T)$ and thus has a unique extension to a continuous function, also denoted $u$, belonging to $C([0, T]; \sX)$ with $u(T) := \lim_{t\to T^-} u(t)\in \bar\sV \subset \sU$. Also by our hypothesis, there is a strong solution, $v$, on a time interval $[T, T_1)$, for some $T_1 := T_1(u(T)) > T$, to the Cauchy problem,
$$
\dot v(t) = -\sE(v(t)), \quad\hbox{a.e. } t \in [T, T_1), \quad v(T) = u(T).
$$
Now define $w:[0, T_1) \to \sX$ by setting
$$
w(t) = \begin{cases} u(t) & \hbox{if } 0 \leq t < T, \\ v(t) & \hbox{if } T \leq t < T_1,\end{cases}
$$
and observe that $w$ is a strong solution to the Cauchy problem,
$$
\dot w(t) = -\sE(w(t)), \quad\hbox{a.e. } t \in [0, T_1), \quad w(0) = u_0.
$$
But this contradicts the maximality of $T$, since $T_1 > T$, and so we must have $T = \infty$.
\end{proof}

\begin{cor}[Global existence of a strong solution to a gradient flow equation]
\label{cor:Neuberger_theorem_4-1}
Assume the hypotheses of Lemma \ref{lem:Neuberger_theorem_4-1}, except for the uniform continuity of $u$ on $[0, T(v_0))$. If there is a positive constant $\delta(u_0)$ such that $u \in W^{1,1}(\delta(u_0), T(u_0))$, then $T(u_0) = \infty$.
\end{cor}

\begin{proof}
Suppose $T := T(u_0) < \infty$ and denote $\delta = \delta(u_0)$. Because $u \in W^{1,1}(\delta, T)$, Theorem \ref{thm:Sell_You_C-9} implies that $u \in C([\delta, T]; \sX)$ and the result now follows from Lemma \ref{lem:Neuberger_theorem_4-1}.
\end{proof}

We can now turn to the

\begin{proof}[Proof of Theorem \ref{thm:Huang_5-1-1}]
By hypothesis, for $\eps > 0$ still to be determined and any $u_0 \in \sU_\eps$, the Cauchy problem \eqref{eq:Huang_3-3a} admits a strong solution, $u:[0,\widetilde T)\to \sX$, where we let $\widetilde T$ denote its maximal lifetime. Clearly, $\widetilde T \geq \tau$, where $\tau > 0$ is the minimal lifetime of $u$ provided by our hypotheses.

\begin{claim}[Confinement of solution to a {\L}ojasiewicz-Simon gradient inequality neighborhood]
\label{claim:Equation_Huang_5-7}
There is a constant
$\eps = \eps(c,C_1,\delta, \theta, \rho, \sigma, \tau, \varphi) \in (0, \sigma/4)$ such that, for any $u_0 \in \sU_\eps$ and solution, $u:[0,\widetilde T)\to \sX$, to \eqref{eq:Huang_3-3a}, there holds
\begin{equation}
\label{eq:Huang_5-7}
\|u(t) - \varphi\|_\sX < \frac{\sigma}{2}, \quad\forall\, t \in [0,\widetilde T).
\end{equation}
\end{claim}

\begin{proof}[Proof of Claim \ref{claim:Equation_Huang_5-7}]
We provide an argument by contradiction. Suppose there exists $T \in (0,\widetilde T)$ such that
\begin{equation}
\label{eq:Huang_5-8}
\|u(t) - \varphi\|_\sX < \frac{\sigma}{2}, \quad\forall\, t \in [0,T),
\quad\hbox{but}\quad \|u(T) - \varphi\|_\sX \geq \frac{\sigma}{2}.
\end{equation}
By the inequality \eqref{eq:Gradient_solution_near_initial_data_at_time_zero_for_short_enough_time} with $\eta = \sigma/8$, we may choose $\delta = \delta(\sigma, \tau, \|\varphi\|_\sX) \in (0,\widetilde T)$ small enough that
\begin{equation}
\label{eq:Huang_theorem_5-1-1_proof_norm_X_udelta_minus_u0_inequality}
\sup_{t \in [0, \delta]}\|u(t)-u_0\|_\sX < \frac{\sigma}{8},
\end{equation}
and, in particular,
$$
\|u(\delta)-u_0\|_\sX < \frac{\sigma}{8}.
$$
Note also that, as we can see by comparing \eqref{eq:Huang_5-8} and \eqref{eq:Huang_theorem_5-1-1_proof_norm_X_udelta_minus_u0_inequality}, we must have
$$
T > \delta.
$$
We can thus apply the growth estimate \eqref{eq:Huang_3-18b_X_interior} in Lemma \ref{lem:Huang_3-3-4} for the interval $[0, T)$ to provide
\begin{align*}
\int_\delta^T \|\dot u(t)\|_\sX\,dt
&\leq
C_1(1+\delta^{-\rho}) \int_{\sE(u(T))}^{\sE(u_0)} \frac{1}{c|s - \sE(\varphi)|^\theta} \,ds
\\
&= \frac{C_1(1+\delta^{-\rho})}{c(1-\theta)}
\left((\sE(u_0) - \sE(\varphi))^{1-\theta} - (\sE(u(T)) - \sE(\varphi))^{1-\theta}\right),
\end{align*}
where to obtain the equality we use the fact that $\sE(v) \geq \sE(\varphi)$ for all $v \in \sU$ by our hypothesis that $\varphi$ is a ground state for $\sE$ on $\sU$ and, in particular, $\sE(u_0) \geq \sE(u(T)) \geq \sE(\varphi)$. By discarding the negative term in the preceding inequality and recalling that \eqref{eq:Rade_7-4_banach} yields, for $\|u_0-\varphi\|_\sX < \eps_\varphi$ and using the fact that $\sE:\sU \to \RR$ is $C^1$,
$$
|\sE(u_0) - \sE(\varphi)| \leq \|u_0 - \varphi\|_\sX,
$$
we obtain
$$
\int_\delta^T \|\dot u(t)\|_\sX\,dt
\leq
\frac{C_1(1+\delta^{-\rho})}{c(1-\theta)} \|u_0 - \varphi\|_\sX^{1-\theta}
<
\frac{C_1(1+\delta^{-\rho})}{c(1-\theta)} \eps^{1-\theta},
$$
where, to obtain the last inequality, we apply our hypothesis that $u_0 \in \sU_\eps$ and the definition \eqref{eq:Huang_page_164_Usigma} of $\sU_\eps$. Hence, for $\eps = \eps(c,C_1,\delta, \theta, \rho, \sigma, \tau, \varphi) \in (0, \sigma/4)$ small enough that $\eps \leq \eps_\varphi$ and
$$
\frac{C_1(1+\delta^{-\rho})}{c(1-\theta)} \eps^{1-\theta} \leq \frac{\sigma}{8},
$$
we find that
$$
\int_\delta^T \|\dot u(t)\|_\sX\,dt < \frac{\sigma}{8},
$$
and thus,
\begin{equation}
\label{eq:Huang_theorem_5-1-1_proof_norm_X_uT_minus_udelta_inequality}
\|u(T)-u(\delta)\|_\sX < \frac{\sigma}{8}.
\end{equation}
Combining the inequalities \eqref{eq:Huang_theorem_5-1-1_proof_norm_X_udelta_minus_u0_inequality} and \eqref{eq:Huang_theorem_5-1-1_proof_norm_X_uT_minus_udelta_inequality} gives
$$
\|u(T)-u_0\|_\sX \leq \|u(\delta)-u_0\|_\sX + \|u(T)-u(\delta)\|_\sX < \frac{\sigma}{8} + \frac{\sigma}{8} = \frac{\sigma}{4},
$$
that is,
\begin{equation}
\label{eq:Huang_theorem_5-1-1_proof_u0_uT_X_inequality}
\|u_0 - u(T)\|_\sX < \frac{\sigma}{4},
\end{equation}
and thus\footnote{Here, we correct a small typographical error in Huang's proof of \cite[Theorem 5.1.1]{Huang_2006}, where \cite[Inequality (5.6)]{Huang_2006} only yields $\eps>0$ in the last displayed equation in \cite[p. 164]{Huang_2006} and not $\eps>\sigma/2$ as required to yield the desired contradiction.},
\begin{align*}
\eps &> \|u_0 - \varphi\|_\sX
\quad\hbox{(by definition \eqref{eq:Huang_page_164_Usigma} of $\sU_\eps$ and the fact that $u_0 \in \sU_\eps$)}
\\
&\geq \|u(T) - \varphi\|_\sX - \|u_0 - u(T)\|_\sX
\\
&> \frac{\sigma}{2} - \frac{\sigma}{4} = \frac{\sigma}{4} \quad \hbox{(by \eqref{eq:Huang_5-8} and \eqref{eq:Huang_theorem_5-1-1_proof_u0_uT_X_inequality})},
\end{align*}
contradicting the choice of $\eps \in (0, \sigma/4)$ in the hypotheses of Theorem \ref{thm:Huang_5-1-1}. This completes the proof of Claim \ref{claim:Equation_Huang_5-7}.
\end{proof}

By virtue of Claim \ref{claim:Equation_Huang_5-7}, we can apply the growth estimate \eqref{eq:Huang_3-18b_X_interior} in Lemma \ref{lem:Huang_3-3-4} for the maximal interval $[0, \widetilde T)$ to provide
$$
\int_\delta^{\widetilde T} \|\dot u(t)\|_\sX\,dt
\leq
\frac{C_1(1+\delta^{-\rho})}{c(1-\theta)} (\sE(u_0) - \sE(\varphi))^{1-\theta} < \infty,
$$
and so $u \in W^{1,1}(\delta, \widetilde T; \sX)$. Corollary \ref{cor:Neuberger_theorem_4-1} now implies that $\widetilde T = \infty$ and so the solution $u$ exists globally on $[0,\infty)$ and \eqref{eq:Huang_5-7} holds with $\widetilde T = \infty$, that is,
$$
\|u(t) - \varphi\|_\sX < \frac{\sigma}{2}, \quad\forall\, t \in [0,\infty).
$$
The orbit $O(u) = \{u(t): t\geq 0\}$ is contained in the open set $\sU_\sigma$ (actually, $\sU_{\sigma/2}$ by \eqref{eq:Huang_5-7} because $\widetilde T = \infty$) and by hypothesis the {\L}ojasiewicz-Simon gradient inequality \eqref{eq:Simon_2-2} holds on $\sU_\sigma$. Therefore, Theorem \ref{thm:Huang_3-3-6} ensures the convergence of the integral $\int_1^\infty\|\dot u(t)\|_\sX\,dt$ and the convergence of $u(t)$ in the norm of $\sX$ to a limit $u_\infty \in \sX$.
\end{proof}

\subsection{Stability of ground states in abstract gradient systems}
\label{subsec:Huang_5-1_stability_ground_state}
In this subsection, we provide a simplification of the second part of \cite[Section 5.1]{Huang_2006}, concerning the question of stability of ground states of an abstract gradient system.

We digress to discuss concepts of stability \cite{Sell_You_2002}. Let $W$ denote a complete metric space. The distance between two points $u$ and $v$ in $W$ will be denoted by $d(u, v) =
d_W(u, v)$, where $d = d_W$ is a metric on $W$; recall that if $W$ is a Banach space, then the standard metric on $W$ is given by $d(u, v) = \|u - v\|_W$, where $\|\cdot\|_W$ is the norm on $W$  \cite[p. 12]{Sell_You_2002}.

\begin{defn}[Semiflow]
\label{defn:Sell_You_page_12_semiflow}
\cite[Section 2.1]{Sell_You_2002}
Let $M$ be a subset of a complete metric space $W$. A mapping $[0, \infty) \times M \ni (t, u) \mapsto \bsigma(t, u) \in M$ is called a \emph{semiflow} on $M$, provided the following hold:

\begin{enumerate}
\item $\bsigma(0, w) = w$, for all $w \in M$;

\item The semigroup property holds, that is,
\begin{equation}
\label{defn:Sell_You_21-1}
\bsigma(t, \bsigma(s, w)) = \bsigma(s + t, w), \quad\forall\, w \in M
\hbox{ and }s, t \in \RR^+.
\end{equation}

\item The mapping $\bsigma : (0,\infty) \times M \to M$ is continuous.
\end{enumerate}
If in addition, the mapping $\bsigma: [0,\infty) \times M \to M$ is continuous, we will
say that the \emph{semiflow $\bsigma$ is continuous at $t = 0$}.
\end{defn}

Recall that a set $K \subset M$ is said to be \emph{positively invariant} if $S(t)K \subset K$, for all $t \geq 0$, and $K$ is said to be an \emph{invariant set} if $S(t)K = K$, for all $t \geq 0$ \cite[Section 2.1.1]{Sell_You_2002}.

\begin{defn}[Lyapunov and uniform asymptotic stability]
\label{defn:Sell_You_page_32_Lyapunov_and_uniform_asymptotic_stability}
\cite[Section 2.3.3]{Sell_You_2002}
Let $\bsigma$ be a semiflow on $W$ and let $A \subset W$. The set $A$ is said to be \emph{Lyapunov stable} provided that $A$ is positively invariant and for every $T > 0$ and every open neighborhood $V$ of $A$, there is an open neighborhood $U$ of $A$ such that
\begin{equation}
\label{defn:Sell_You_23-5}
S(t)U \subset V, \quad\forall\, t \geq T.
\end{equation}
The set $A$ is said to be \emph{uniformly asymptotically stable} if it is Lyapunov
stable and there is a neighborhood $U_0$ of $A$ such that $A$ attracts $U_0$, that is,
\begin{equation}
\label{defn:Sell_You_23-6}
d(S(t)U_0, A) \to 0, \quad\hbox{as } t \to \infty.
\end{equation}
\end{defn}

Suppose we are in the setting of Hypothesis \ref{hyp:Sell_You_4_standing_hypothesis_A} (`Standing Hypothesis A') with $\cF \in C_{\Lip}(\calV^{2\beta}, \cW)$, where $0 \leq \beta < 1$ and $\cW$ and $\calV$ are Banach spaces with $\calV = \sD(\cA)$ and that $u$ is a solution to \eqref{eq:Sell_You_47-1}, that is,
\begin{equation}
\label{eq:Sell_You_47-12}
\frac{\partial u}{\partial t} + \cA u = \cF(u(t)), \quad u(0) = u_0 \in \calV^{2\beta}.
\end{equation}
In this setting, we define
$$
M \equiv M(\cF) := \{u_0 \in \calV^{2\beta}: T(u_0, \cF) = \infty\},
$$
and for $u_0 \in M$, we set
\begin{equation}
\label{defn:Sell_You_47-13}
S(t)u_0 \equiv \bsigma(t, u_0) := \phi(t, u_0, \cF).
\end{equation}
According to \cite[Sections 4.6.4 and 4.7.5]{Sell_You_2002}, $\bsigma$ is a semiflow on $M$.

An \emph{equilibrium point} (or \emph{stationary solution}) for the evolutionary equation \eqref{eq:Sell_You_47-12} is a point $u_0 \in \sD(\cA)$, with the property that $\cA u_0 = \cF(u_0)$, that is, $\partial u/\partial t = 0$ and on $(0, \infty)$ when $u$ solves \eqref{eq:Sell_You_47-12}. See \cite[Section 7.1]{Sell_You_2002} for further discussion of the local dynamics of an evolutionary equation near an equilibrium point.

We now end the digression and return to the setting of solutions to gradient systems, keeping in mind that, in our application to the Yang-Mills gradient flow, our gradient systems can be given the structure of an evolutionary equation \eqref{eq:Sell_You_47-12}. We have the following analogue of \cite[Theorem 5.1.2]{Huang_2006}.

\begin{thm}[Convergence to a critical point and stability of a ground state]
\label{thm:Huang_5-1-2}
Assume the hypotheses of Theorem \ref{thm:Huang_5-1-1}. Then, for each $u_0 \in \sU_\eps$, the autonomous Cauchy problem \eqref{eq:Huang_3-3a} admits a global strong solution $u:[0,\infty) \to \sU_{\sigma/2}$, in the sense of Definition \ref{defn:Strong_solution_to_gradient_system}, that converges in $\sX$ as $t\to\infty$ to some critical point $u_\infty \in \sU_\sigma$. The critical point, $u_\infty$, satisfies $\sE(u_\infty) = \sE(\varphi)$. As an equilibrium of \eqref{eq:Huang_3-3a}, the point $\varphi$ is Lyapunov stable (see Definition \ref{defn:Sell_You_page_32_Lyapunov_and_uniform_asymptotic_stability}). If $\varphi$ is isolated or a cluster point of the orbit $O(u)$, then $\varphi$ is uniformly asymptotically stable (see Definition \ref{defn:Sell_You_page_32_Lyapunov_and_uniform_asymptotic_stability}).
\end{thm}

\begin{proof}
The existence of a global strong solution, $u:[0,\infty) \to \sU_{\sigma/2}$, for each $u_0 \in \sU_\sigma$ follows from Theorem \ref{thm:Huang_5-1-1}.

The fact that a solution $u$ converges to a limit $u_\infty \in \sX$ in the sense that $\|u(t) - u_\infty\|_\sX \to 0$ as $t \to \infty$ and $\int_1^\infty \|\dot u(t)\|_\sX \,dt < \infty$ is also implied Theorem \ref{thm:Huang_5-1-1}.

The fact that $u_\infty$ is a critical point of $\sE$ is due to Lemma \ref{lem:Limit_is_critical_point}.\footnote{Compare \cite[Theorem 4.2]{Neuberger_2010}.}

The {\L}ojasiewicz-Simon gradient inequality (Theorem \ref{thm:Huang_2-4-2}) and the facts that $\sE'(\varphi) = 0$ and $u_\infty \in \sU_\sigma$ (in fact, $\sU_{\sigma/2}$) imply that $\sE(u_\infty) = \sE(\varphi)$.

To prove the Lyapunov stability of $\varphi$, it suffices to apply Theorem \ref{thm:Huang_5-1-1} for each fixed $\varsigma \in (0,\sigma)$. Thus, for each $\varsigma \in (0, \sigma)$ Theorem \ref{thm:Huang_5-1-1} implies that there is an $\eps = \eps_\varsigma = \eps(c,C_1,T,\varphi,\varsigma,\theta) \in (0, \varsigma/4)$ such that for each $u_0 \in \sU_{\eps_\varsigma}$, the gradient system \eqref{eq:Huang_3-3a} has a global solution $u$ such that
$$
\|u(t) - \varphi\|_\sX < \frac{\varsigma}{2}, \quad \forall\, t \geq 0.
$$
This proves the desired Lyapunov stability of $\varphi$ in the sense of Definition \ref{defn:Sell_You_page_32_Lyapunov_and_uniform_asymptotic_stability}.

If $\varphi$ is isolated, then $u_\infty = \varphi$; alternatively, if $\varphi$ is a cluster point of $O(u)$, then $u_\infty = \varphi$ just as in the proof of the final assertion of Theorem \ref{thm:Huang_3-3-6}. When $u_\infty = \varphi$, uniform asymptotic stability of $\varphi$ follows from the convergence rate estimate \eqref{eq:Huang_3-45_X} for $u(t)$ to $u_\infty$, as $t \to \infty$, given in Theorem \ref{thm:Huang_3-4-8}. Indeed, the inequality \eqref{eq:Huang_3-45_X} and the expression \eqref{eq:Huang_3-45_growth_rate} for $\Phi(g(t))$ yield
$$
\|u(t) - u_\infty\|_\sX \to 0 \quad\hbox{as } t \to \infty,
$$
and so $\varphi$ is uniformly asymptotic stable by Definition \ref{defn:Sell_You_page_32_Lyapunov_and_uniform_asymptotic_stability}.
\end{proof}

\begin{rmk}[On hypotheses for uniform asymptotic stability of $\varphi$ in Theorem \ref{thm:Huang_5-1-2}]
\label{rmk:Huang_5-1-2}
It is worth mentioning that, in the setting of the Yang-Mills gradient flow, critical points are rarely known to be isolated. Huang obtains uniform asymptotic stability of $\varphi$ in his \cite[Theorem 5.1.2]{Huang_2006} under the conditions that, for all sufficiently small $r > 0$,
\begin{subequations}
\label{eq:Huang_5-15}
\begin{align}
\label{eq:Huang_5-15a}
\inf\{\|\sE'(u)\|_{\sX'}: u \in \sU_\sigma \hbox{ and } \sE(u) = r\} > 0,
\\
\label{eq:Huang_5-15b}
\frac{1}{\alpha}\|\sE'(u)\|_\sX \geq \|\sE'(u)\|_{\sX'} \geq \alpha \|\sE'(u)\|_\sX, \quad \forall\, u \in \sU.
\end{align}
\end{subequations}
However, it is unclear that these conditions can be verified when $\sE$ is the Yang-Mills energy functional even when the Banach space, $\sX$, in \eqref{eq:Huang_5-15} is replaced by a Hilbert space, $\sH$, so $\sH' \cong \sH$.
\end{rmk}

\chapter[Yang-Mills gradient flow near a local minimum]{Global existence and convergence for Yang-Mills gradient flow near a local minimum}
\label{chapter:Huang_5-1_application_pure_Yang-Mills_gradient_flow}
In this chapter, we derive several consequences of our results thus far when considering the asymptotic behavior of solutions to Yang-Mills gradient flow near a local minimum of the Yang-Mills energy functional. These consequences will follow in a straightforward manner from the general analytical theory that we have established for gradient flow and Yang-Mills gradient flow, in particular. In a later chapter, we shall consider the extensions of these results that can be achieved with the aid of an analysis of solutions to a pseudogradient flow, that is, gradient flow plus a suitably small, time-varying perturbation. In keeping with the pattern of our previous analysis, we first describe results that are valid for Yang-Mills gradient flow over a closed, Riemannian, smooth manifold, $X$, of any dimension $d \geq 2$ and then specialize, as needed, to the case of $d = 4$.

The simplest case to consider is when the Yang-Mills gradient flow, $A(t)$, is started at a connection, $A_0$, on $P$ that is sufficiently close, in the sense of the usual norm on $H_{A_1}^1(X; \Lambda^1\otimes\ad P)$, to a local minimum, say $A_1$, of the Yang-Mills energy functional. In this situation, our conclusions for global existence, convergence, convergence rate, and stability of solutions to Yang-Mills gradient flow will follow quickly from the results of Sections \ref{subsec:Huang_3-4C}, \ref{subsec:Huang_5-1_existence_and_convergence_global_solution_near_local_minimum}, and \ref{subsec:Huang_5-1_stability_ground_state}. More complex situations specific to dimension four, involving initial data, $A_0$, with small $\|F_{A_0}^+\|_{L^p(X)}$ or even initial data, $A_0$, with arbitrary energy, $\sE(A_0) = \frac{1}{2}\|F_{A_0}\|_{L^2(X)}^2$, are treated in later chapters.

\section[Energy inequalities and a priori estimates]{Energy inequalities and \apriori estimates for a variational solution to a linear heat equation}
\label{sec:Simon_4_linear_theory}
In order to extend the proof of \cite[Lemma 7.3]{Rade_1992} to the case where $\dim X = 4$ (rather than $2$ or $3$ as in \cite{Rade_1992}), we shall require certain \apriori estimates for variational solutions to a linear heat equation. The estimates we require are essentially identified by Simon in \cite[pp. 543--544]{Simon_1983}. We also include some useful extensions due to R\r{a}de \cite{Rade_1992}. Estimates for a nonlinear, second-order, parabolic system are given by Simon in \cite[p. 545]{Simon_1983}. The estimates in Lemmata \ref{lem:Rade_inequalities_11-3_and_11-4} and \ref{lem:Evans_secondedition_equation_7-1-20_L1_and_L2} complement those obtained by semigroup methods in Corollary \ref{cor:Sell_You_42-12_heat_equation_alpha_is_zero_apriori_estimates}.

We begin with the following basic \emph{energy estimate} (compare, for example, \cite[p. 543]{Simon_1983} but note that our sign convention for the Laplace operator is opposite to Simon's choice).

\begin{lem}[\Apriori estimate for a solution to a variational solution to a linear heat equation]
\label{lem:Rade_inequalities_11-3_and_11-4}
Let $V$ be a finite-rank, Riemannian, smooth vector bundle with
Riemannian smooth covariant derivative, $\nabla$, over a closed, Riemannian, smooth manifold, $X$, of dimension $d \geq 2$. For $-\infty < t_1 < t_2 \leq \infty$, let $f_1 \in L^2(t_1,t_2; H^{-1}(X;V))$, and $f_2 \in L^1(t_1,t_2; L^2(X;V))$. If
$$
u \in L^2(t_1,t_2; H^1(X;V)) \quad\hbox{with }\dot u \in L^2(t_1,t_2; H^{-1}(X;V)),
$$
obeys
\begin{equation}
\label{eq:Simon_linear_second-order_parabolic_variational_equation}
\begin{aligned}
{}& (\dot u(t), v(t))_{L^2(X;V)} + (\cov u(t), \cov v(t))_{L^2(X;V)}
\\
&\quad = (f_1(t), v(t))_{L^2(X;V)} + (f_2(t), v(t))_{L^2(X;V)}, \quad\hbox{a.e. } t \in (t_1,t_2),
\end{aligned}
\end{equation}
for all $v \in L^2(t_1,t_2; H^1(X;V))$, then $u$ satisfies
\begin{multline}
\label{eq:Energy_inequality_compact_form}
\frac{1}{2\sqrt{2}}\|u\|_{L^\infty(t_1,t_2; L^2(X;V))} + \frac{1}{2}\|u\|_{L^2(t_1,t_2; H^1(X;V))}
\\
\leq \|u(t_1)\|_{L^2(X;V)} + \sqrt{2}\|u\|_{L^2(t_1,t_2; L^2(X;V))}
\\
+ \|f_1\|_{L^2(t_1,t_2; H^{-1}(X;V))} + \sqrt{2}\|f_2\|_{L^1(t_1,t_2;  L^2(X;V))}.
\end{multline}
If $|t_2-t_1| \leq 1/8$, then
\begin{equation}
\label{eq:Rade_11-3_and_11-4}
\begin{aligned}
{}&\frac{1}{4}\|u\|_{L^\infty(t_1,t_2; L^2(X;V))} + \frac{1}{2}\|u\|_{L^2(t_1,t_2; H^1(X;V))}
\\
&\quad\leq \|u(t_1)\|_{L^2(X;V)} + \|f_1\|_{L^2(t_1,t_2; H^{-1}(X;V))} + \sqrt{2}\|f_2\|_{L^1(t_1,t_2;  L^2(X;V))}.
\end{aligned}
\end{equation}
If $|t_2-t_1| < \infty$, then there is a positive constant, $C = C(t_2-t_1)$, such that
\begin{equation}
\label{eq:Rade_11-3_and_11-4_finite_time_interval}
\begin{aligned}
{}&\|u\|_{L^\infty(t_1,t_2; L^2(X;V))} + \|u\|_{L^2(t_1,t_2; H^1(X;V))}
\\
&\quad\leq C\left(\|u(t_1)\|_{L^2(X;V)} + \|f_1\|_{L^2(t_1,t_2; H^{-1}(X;V))}
+ \|f_2\|_{L^1(t_1,t_2;  L^2(X;V))} \right).
\end{aligned}
\end{equation}
\end{lem}

\begin{proof}
By analogy with \cite[Theorem 5.9.3]{Evans2} in the scalar case, we have $u \in C([t_1,t_2]; L^2(X;V))$. By choosing $v = u$ and substituting (by analogy with \cite[Theorem 5.9.3]{Evans2} in the scalar case),
$$
\frac{d}{dt}\|u(t)\|_{L^2(X;V)}^2 = 2(u(t), \dot u(t))_{L^2(X;V)}, \quad\hbox{a.e. } t \in (t_1,t_2),
$$
into \eqref{eq:Simon_linear_second-order_parabolic_variational_equation} to give
$$
\frac{d}{dt}\|u(t)\|_{L^2(X;V)}^2 + 2\|\cov u(t)\|_{L^2(X;V)}^2 = 2(f_1(t), u(t))_{L^2(X;V)} + 2(f_2(t), u(t))_{L^2(X;V)},
$$
and integrating over $[t_1,\tau]$, where $\tau \in [t_1, t_2]$, one finds that (compare \cite[p. 543]{Simon_1983})
\begin{align*}
{}&\|u(\tau)\|_{L^2(X;V)}^2 - \|u(t_1)\|_{L^2(X;V)}^2 + 2\int_{t_1}^\tau \|\cov u(t)\|_{L^2(X;V)}^2 \,dt
\\
&\quad = 2\int_{t_1}^\tau (f_1(t), u(t))_{L^2(X;V)} \,dt + 2\int_{t_1}^\tau (f_2(t), u(t))_{L^2(X;V)} \,dt
\\
&\quad \leq 2\int_{t_1}^{t_2} \|f_1(t)\|_{H^{-1}(X;V)} \|u(t)\|_{H^1(X;V)} \,dt
+ 2\sup_{t\in[t_1, t_2]}\|u(t)\|_{L^2(X;V)} \int_{t_1}^{t_2} \|f_2(t)\|_{L^2(X;V)} \,dt
\\
&\quad \leq \int_{t_1}^{t_2} \|f_1(t)\|_{H^{-1}(X;V)}^2 \,dt
+ \int_{t_1}^{t_2} \|u(t)\|_{H^1(X;V)}^2 \,dt
+ \frac{1}{2}\sup_{t\in[t_1, t_2]}\|u(t)\|_{L^2(X;V)}^2
\\
&\qquad + 2\left(\int_{t_1}^{t_2} \|f_2(t)\|_{L^2(X;V)} \,dt\right)^2, \quad\forall\, \tau \in [t_1, t_2].
\end{align*}
Thus, using $\|u(t)\|_{H^1(X;V)}^2 = \|u(t)\|_{L^2(X;V)}^2 + \|\cov u(t)\|_{L^2(X;V)}^2$, taking the supremum over $\tau \in [t_1,t_2]$ on the left-hand side, and rearranging, we see that
\begin{align*}
{}&\frac{1}{2}\sup_{t\in[t_1, t_2]}\|u(t)\|_{L^2(X;V)}^2 + \int_{t_1}^{t_2} \|\cov u(t)\|_{L^2(X;V)}^2 \,dt
\\
&\quad \leq \|u(t_1)\|_{L^2(X;V)}^2 + \int_{t_1}^{t_2} \|f_1(t)\|_{H^{-1}(X;V)}^2 \,dt + \int_{t_1}^{t_2} \|u(t)\|_{L^2(X;V)}^2 \,dt
\\
&\qquad + 2\left(\int_{t_1}^{t_2} \|f_2(t)\|_{L^2(X;V)} \,dt\right)^2.
\end{align*}
Therefore, adding $\int_{t_1}^{t_2} \|u(t)\|_{L^2(X;V)}^2 \,dt$ to both sides yields the inequality,
\begin{equation}
\label{eq:Energy_inequality}
\begin{aligned}
{}&\frac{1}{2}\sup_{t\in[t_1, t_2]}\|u(t)\|_{L^2(X;V)}^2 + \int_{t_1}^{t_2}\|u(t)\|_{H^1(X;V)}^2 \,dt
\\
&\quad \leq  \|u(t_1)\|_{L^2(X;V)}^2 + 2\int_{t_1}^{t_2} \|u(t)\|_{L^2(X;V)}^2 \,dt
\\
&\qquad + \int_{t_1}^{t_2} \|f_1(t)\|_{H^{-1}(X;V)}^2 \,dt + 2\left(\int_{t_1}^{t_2} \|f_2(t)\|_{L^2(X;V)} \,dt\right)^2.
\end{aligned}
\end{equation}
By taking square roots and using $\frac{1}{2}(a + b) \leq (a^2 + b^2)^{1/2} \leq a + b$ (for $a, b \geq 0$), we obtain the desired inequality \eqref{eq:Energy_inequality_compact_form} from \eqref{eq:Energy_inequality}.

Provided $|t_2-t_1| \leq 1/8$, we have
$$
\int_{t_1}^{t_2} \|u(t)\|_{L^2(X;V)}^2 \,dt \leq \frac{1}{8}\sup_{t\in[t_1, t_2]}\|u(t)\|_{L^2(X;V)}^2,
$$
and we obtain the following extension of \cite[Equation (4.1)]{Simon_1983},
\begin{equation}
\label{eq:Energy_inequality_improved}
\begin{aligned}
{}&\frac{1}{4}\sup_{t\in[t_1, t_2]}\|u(t)\|_{L^2(X;V)}^2  + \int_{t_1}^{t_2}\|u(t)\|_{H^1(X;V)}^2 \,dt
\\
&\quad \leq \|u(t_1)\|_{L^2(X;V)}^2 + \int_{t_1}^{t_2} \|f_1(t)\|_{H^{-1}(X;V)}^2 \,dt
+ 2\left(\int_{t_1}^{t_2} \|f_2(t)\|_{L^2(X;V)} \,dt\right)^2.
\end{aligned}
\end{equation}
By writing the \apriori estimate \eqref{eq:Energy_inequality_improved} more succinctly,
\begin{align*}
{}&\frac{1}{4}\|u\|_{L^\infty(t_1,t_2; L^2(X;V))}^2 + \|u\|_{L^2(t_1,t_2; H^1(X;V))}^2
\\
&\quad \leq \|u(t_1)\|_{L^2(X;V)}^2 + \|f_1\|_{L^2(t_1,t_2; H^{-1}(X;V))}^2 + 2\|f_2\|_{L^1(t_1,t_2;  L^2(X;V))}^2,
\end{align*}
and taking square roots and using $\frac{1}{2}(a + b) \leq (a^2 + b^2)^{1/2} \leq a + b$ (for $a, b \geq 0$), we obtain the desired inequality \eqref{eq:Rade_11-3_and_11-4}.

More generally, provided only $|t_2 - t_1| < \infty$ and partitioning the interval $[t_1, t_2]$ into $k$ subintervals $[s_j, s_{j+1}]$ of length at most $1/8$, for $1 \leq j \leq k$, where $t_1 = s_1 < s_2 < \cdots < s_k < s_{k+1} = t_2$, and using
\begin{equation}
\label{eq:Geometric_sum_inequality}
\frac{1}{k}(a_1 + \cdots + a_k) \leq (a_1^2 + \cdots + a_k^2)^{1/2} \leq a_1 + \cdots + a_k
\quad \hbox{(with $a_j \geq 0$ for $1 \leq j \leq k$)},
\end{equation}
we see that, via $[t_1, t_2] = \cup_{j=1}^k[s_j, s_{j+1}] = [s_1, s_{k+1}]$,
\begin{align*}
{}&\frac{1}{4}\|u\|_{L^\infty(t_1,t_2; L^2(X;V))} + \frac{1}{2}\|u\|_{L^2(t_1,t_2; H^1(X;V))}
\\
&\quad \leq \sum_{j=1}^k \left( \frac{1}{4}\|u\|_{L^\infty(s_j,s_{j+1}; L^2(X;V))}
+ \frac{1}{2}\|u\|_{L^2(s_j,s_{j+1}; H^1(X;V))} \right)
\\
&\quad \leq \sum_{j=1}^k \left(\|u(s_j)\|_{L^2(X;V)} + \|f_1\|_{L^2(s_j,s_{j+1}; H^{-1}(X;V))}
+ \sqrt{2}\|f_2\|_{L^1(s_j,s_{j+1};  L^2(X;V))}\right)
\\
&\quad \leq k\|f_1\|_{L^2(s_1,s_{k+1}; H^{-1}(X;V))} + \sqrt{2}\|f_2\|_{L^1(s_1,s_{k+1};  L^2(X;V))}
+ \sum_{j=1}^k \|u(s_j)\|_{L^2(X;V)},
\end{align*}
and so
\begin{multline}
\label{eq:Rade_11-3_and_11-4_finite_time_interval_raw}
\frac{1}{4}\|u\|_{L^\infty(t_1,t_2; L^2(X;V))} + \frac{1}{2}\|u\|_{L^2(t_1,t_2; H^1(X;V))}
\\
\leq k\|f_1\|_{L^2(t_1,t_2; H^{-1}(X;V))} + \sqrt{2}\|f_2\|_{L^1(s_1,s_{k+1};  L^2(X;V))}
+ \sum_{j=1}^k \|u(s_j)\|_{L^2(X;V)}.
\end{multline}
On the other hand, for $2 \leq j \leq k$, the inequality \eqref{eq:Rade_11-3_and_11-4} gives
\begin{align*}
\|u(s_j)\|_{L^2(X;V)} &\leq \|u\|_{L^\infty(s_{j-1},s_j; L^2(X;V))}
\\
&\leq 4\left(\|u(s_{j-1})\|_{L^2(X;V)} + \|f_1\|_{L^2(s_{j-1},s_j; H^{-1}(X;V))}
+ \sqrt{2}\|f_2\|_{L^1(s_{j-1},s_j;  L^2(X;V))}\right),
\end{align*}
and thus, by induction on $j \geq 2$ and using
$$
\sum_{i=2}^j\|f_1\|_{L^2(s_{i-1},s_i; H^{-1}(X;V))} \leq (j-1)\|f_1\|_{L^2(s_1,s_j; H^{-1}(X;V))},
$$
we obtain
\begin{multline*}
\|u(s_j)\|_{L^2(X;V)} \leq 4^{j-1}\left(\|u(s_1)\|_{L^2(X;V)} + (j-1)\|f_1\|_{L^2(s_1,s_j; H^{-1}(X;V))} \right.
\\
+ \left. \sqrt{2}\|f_2\|_{L^1(s_1,s_j;  L^2(X;V))}\right), \quad 2 \leq j \leq k.
\end{multline*}
Thus,
\begin{multline*}
\sum_{j=2}^k \|u(s_j)\|_{L^2(X;V)}
\leq
4^{k-1}(k-1)\left(\|u(s_1)\|_{L^2(X;V)} + (k-1)\|f_1\|_{L^2(s_1,s_k; H^{-1}(X;V))} \right.
\\
+ \left. \sqrt{2}\|f_2\|_{L^1(s_1,s_k;  L^2(X;V))}\right).
\end{multline*}
Combining the preceding bound with \eqref{eq:Rade_11-3_and_11-4_finite_time_interval_raw} gives
\begin{align*}
{}&\frac{1}{4}\|u\|_{L^\infty(t_1,t_2; L^2(X;V))} + \frac{1}{2}\|u\|_{L^2(t_1,t_2; H^1(X;V))}
\\
&\quad \leq \left(k + 4^{k-1}(k-1)^2\right)\|f_1\|_{L^2(s_1,s_{k+1}; H^{-1}(X;V))}
\\
&\qquad + \sqrt{2}\left(1 + 4^{k-1}(k-1)\right)\|f_2\|_{L^1(s_1,s_{k+1};  L^2(X;V))}
+ \left(1 + 4^{k-1}(k-1)\right)\|u(s_1)\|_{L^2(X;V)}.
\end{align*}
But $[t_1, t_2] = \cup_{j=1}^k[s_j, s_{j+1}]$ and we may assume without loss of generality that $|s_{j+1} - s_j| = 1/8$ for $1\leq j\leq k$, so $|t_2 - t_1| = k/8$ and the conclusion \eqref{eq:Rade_11-3_and_11-4_finite_time_interval} follows for some constant $C = C(t_2-t_1)$.
\end{proof}

\begin{rmk}[Alternative proofs of Lemma \ref{lem:Rade_inequalities_11-3_and_11-4} and related results]
\label{rmk:Rade_inequalities_11-3_and_11-4_alternative_proofs}
When $f_2=0$, the \apriori estimate \eqref{eq:Rade_11-3_and_11-4_finite_time_interval} can be compared with that provided by \cite[Theorem 7.1.2 and Equation (7.1.20)]{Evans2} for the scalar parabolic equation over a bounded cylinder, $(t_1,t_2)\times U$ with $U\subset\RR^d$, and proved using the differential form of Gronwall's Inequality \cite[Appendix B.2 (j)]{Evans}.

When $f_2 = 0$, the inequalities \eqref{eq:Rade_11-3_and_11-4} and \eqref{eq:Rade_11-3_and_11-4_finite_time_interval} also follow (with different constants) from the \apriori estimates \eqref{eq:Sell_You_42-30_heat_equation_Linfty-L2_u_compact_form} and \eqref{eq:Sell_You_42-30_heat_equation_L2-H1_u_compact_form} in Corollary \ref{cor:Sell_You_42-12_heat_equation_alpha_is_zero_apriori_estimates}, though the latter inequalities are stronger since they hold for an unbounded time interval and universal numerical constants independent of the interval length.
\end{rmk}

\begin{rmk}[Application of Lemma \ref{lem:Rade_inequalities_11-3_and_11-4}]
\label{rmk:Rade_inequalities_11-3_and_11-4_applications}
We shall use the \apriori estimate \eqref{eq:Rade_11-3_and_11-4} in place of R\r{a}de's \cite[Inequalities (11.3) and (11.4)]{Rade_1992} in our proof of our version, Lemma \ref{lem:Rade_7-3} (for $\dim X=4$) in the sequel, of his \cite[Lemma 7.3]{Rade_1992} (for $\dim X=2$ or $3$).
\end{rmk}

\begin{lem}[\Apriori estimate for the time derivative of a variational solution to a linear heat equation]
\label{lem:Evans_secondedition_equation_7-1-20_L1_and_L2}
Assume the hypotheses of Lemma \ref{lem:Rade_inequalities_11-3_and_11-4}, except for the conditions on $t_1, t_2$ and $f_1,f_2$. Then the following hold.
\begin{enumerate}
\item If $-\infty \leq t_1 < t_2 \leq \infty$ and $f_1 \in L^1(t_1,t_2; H^{-1}(X;V))$ and $f_2 \in L^1(t_1,t_2; L^2(X;V))$, then
\begin{multline}
\label{eq:Evans_secondedition_equation_7-1-20_L1}
\|\dot u\|_{L^1(t_1,t_2; H^{-1}(X;V))} \leq \|u\|_{L^1(t_1,t_2; H^1(X;V))} + \|f_1\|_{L^1(t_1,t_2; H^{-1}(X;V))}
\\
+ \|f_2\|_{L^1(t_1,t_2; L^2(X;V))}.
\end{multline}
\item If $-\infty \leq t_1 < t_2 \leq \infty$ and $f_1 \in L^2(t_1,t_2; H^{-1}(X;V))$ and $f_2 \in L^2(t_1,t_2;L^2(X;V))$, then
\begin{multline}
\label{eq:Evans_secondedition_equation_7-1-20_L2}
\frac{1}{\sqrt{3}}\|\dot u\|_{L^2(t_1,t_2;H^{-1}(X;V))} \leq \|u\|_{L^2(t_1,t_2;H^1(X;V))} + \|f_1\|_{L^2(t_1,t_2;H^{-1}(X;V))}
\\
+ \|f_2\|_{L^2(t_1,t_2;L^2(X;V))}.
\end{multline}
\item If $|t_2-t_1|<\infty$ and $f_1 \in L^2(t_1,t_2; H^{-1}(X;V))$ and $f_2 \in L^1(t_1,t_2; L^2(X;V))$, then there is a positive constant, $C=C(t_2-t_1)$, such that
\begin{multline}
\label{eq:Evans_secondedition_equation_7-1-20_L1_finite_time_interval}
\|\dot u\|_{L^1(t_1,t_2; H^{-1}(X;V))}
\\
\leq C\left(\|u(t_1)\|_{L^2(X;V)} + \|f_1\|_{L^2(t_1,t_2; H^{-1}(X;V))}
+ \|f_2\|_{L^1(t_1,t_2;  L^2(X;V))}\right).
\end{multline}
\end{enumerate}
\end{lem}

\begin{proof}
Returning to \eqref{eq:Simon_linear_second-order_parabolic_variational_equation} and choosing $v$ to be constant on $[t_1,t_2]$ with $\|v\|_{H^1(X;V)} \leq 1$, we obtain, for a.e. $t \in (t_1,t_2)$,
\begin{align*}
\left|(\dot u(t), v)_{L^2(X;V)}\right|
& \leq \|\cov u(t)\|_{L^2(X;V)}\|\cov v\|_{L^2(X;V)} + \|f_1(t)\|_{H^{-1}(X;V))}\|v\|_{H^1(X;V))}
\\
&\quad + \|f_2(t)\|_{L^2(X;V)}\|v\|_{L^2(X;V)},
\end{align*}
and thus
\begin{equation}
\label{eq:Evans_secondedition_equation_7-1-20_L1_raw}
\|\dot u(t)\|_{H^{-1}(X;V)}
\leq \|u(t)\|_{H^1(X;V)} + \|f_1(t)\|_{H^{-1}(X;V))} + \|f_2(t)\|_{L^2(X;V)}.
\end{equation}
Integration over $[t_1, t_2]$ yields \eqref{eq:Evans_secondedition_equation_7-1-20_L1}. If $|t_2-t_1| < \infty$, then
\begin{align*}
\|u\|_{L^1(t_1,t_2; H^1(X;V))} &\leq |t_2-t_1|^{1/2} \|u\|_{L^2(t_1,t_2;H^1(X;V))},
\\
\|f_1\|_{L^1(t_1,t_2; H^{-1}(X;V))} &\leq |t_2-t_1|^{1/2} \|f_1\|_{L^2(t_1,t_2; H^{-1}(X;V))}.
\end{align*}
Combining the preceding inequalities with \eqref{eq:Rade_11-3_and_11-4_finite_time_interval} yields \eqref{eq:Evans_secondedition_equation_7-1-20_L1_finite_time_interval}.

When $f_1 \in L^2(t_1,t_2; H^{-1}(X;V))$ and $f_2 \in L^2(t_1,t_2;L^2(X;V))$, then \eqref{eq:Evans_secondedition_equation_7-1-20_L1_raw} also leads to \eqref{eq:Evans_secondedition_equation_7-1-20_L2} using \eqref{eq:Geometric_sum_inequality}.
\end{proof}

\begin{rmk}[Alternative proofs of Lemma \ref{lem:Evans_secondedition_equation_7-1-20_L1_and_L2} and related results]
\label{rmk:Evans_secondedition_equation_7-1-20_L1_and_L2_alternative_proofs}
When $f_2=0$, the \apriori estimate \eqref{eq:Evans_secondedition_equation_7-1-20_L1_finite_time_interval} corresponds to \cite[Equation (7.1.20)]{Evans2} in the case of a scalar parabolic equation over a bounded cylinder, $(t_1,t_2)\times U$ with $U\subset\RR^d$, and proved using the differential form of Gronwall's Inequality \cite[Appendix B.2 (j)]{Evans}. When $f_2=0$, the inequality \eqref{eq:Evans_secondedition_equation_7-1-20_L1_finite_time_interval} also follows (with different constants) from the \apriori estimate \eqref{eq:Sell_You_42-30_heat_equation_L2-Hminus1_dudt_compact_form_finite_T0} in Corollary \ref{cor:Sell_You_42-12_heat_equation_alpha_is_zero_apriori_estimates}.
\end{rmk}

\begin{rmk}[Differential form of Gronwall's inequality]
\label{rmk:Gronwall_inequality_differential_form}
It is convenient at this point to recall the differential form of \emph{Gronwall's Inequality} \cite[Appendix B.2 (j)]{Evans2}. Let $\eta(t)$ be a nonnegative, absolutely continuous function on $[0,T]$, which satisfies for a.e. $t \in (0,T)$ the differential inequality
\begin{equation}
\label{eq:Evans_B-15}
\dot \eta(t) \leq \phi(t)\eta(t) + \psi(t),
\end{equation}
where $\phi(t)$ and $\psi(t)$ are nonnegative, summable functions on $[0,T]$. Then
\begin{equation}
\label{eq:Evans_B-16}
\eta(t) \leq \left(\eta(0) +  \int_0^t \psi(s)\,ds\right)\exp\left(\int_0^t \phi(s)\,ds\right),
\end{equation}
for all $0 < t < T$.
\end{rmk}

\section[A priori estimates for lengths of Yang-Mills gradient-like flow lines]{\Apriori estimates for lengths of Yang-Mills gradient-like flow lines}
\label{sec:Rade_Lemma_7-3_generalization}
In this section, we develop two extensions to R\r{a}de's key \apriori $L^1$-in-time interior estimate for a solution to Yang-Mills gradient flow \cite[Lemma 7.3]{Rade_1992}, by relaxing the conditions on both the dimension, $d$, of the base manifold, $X$, and the spatial regularity of the solution. We develop a similar generalization of R\r{a}de's \cite[Lemma 8.1]{Rade_1992}.

\subsection{A generalization of R\r{a}de's Lemma 7.3}
\label{subsec:Rade_Lemma_7-3_generalization_dimension_2_leq_d_leq_4}
In this subsection, we shall prove a generalization of \cite[Lemma 7.3]{Rade_1992}, when $X$ has dimension $d$ obeying $2 \leq d \leq 4$ rather than $d = 2$ or $3$, and extending from consideration of Yang-Mills gradient flow to Yang-Mills gradient-like flow by allowing a possibly non-zero perturbation term, $R(t)$, on the right-hand side.

In preparation for the proof of Lemma \ref{lem:Rade_7-3}, it is useful to recall the following fact noted in \cite[p. 235]{DK}, reminiscent of the better-known identities, $d_Ad_A = F_A$ \cite[Equation (2.1.3)]{DK} or $d_AF_A = 0$ \cite[Equation (2.1.21)]{DK}, when $A$ is a connection on a principal $G$-bundle with curvature $F_A$.

\begin{lem}
\label{lem:Donaldson_Kronheimer_page_235}
\cite[p. 235]{DK}
\cite[p. 577]{ParkerGauge}
Let $G$ be a Lie group and $P$ a principal $G$-bundle over a Riemannian, smooth manifold, $X$, of dimension $d$ with $d \geq 2$. If $A$ is a $C^\infty$ connection on $P$ with curvature $F_A \in \Omega^2(X; \ad P)$, then
\begin{equation}
\label{eq:Donaldson_Kronheimer_page_235}
d_A^*d_A^*F_A = 0.
\end{equation}
\end{lem}

\begin{proof}
We include the proof provided in \cite[p. 235]{DK} for the sake of completeness. (See also \cite[p. 577]{ParkerGauge}.) We have
$$
d_A^*d_A^*F_A = (d_Ad_A)^*F_A = F_A^*F_A = \{F_A, F_A\},
$$
where $\{,\}$ denotes the tensor product of the symmetric inner product on two-forms and the skew symmetric bracket on the Lie algebra of $G$. Thus, $\{\cdot, \cdot\}$ is skew symmetric and therefore $d_A^*d_A^* F_A = 0$.
\end{proof}

\begin{lem}[An \apriori $L^1$-in-time interior estimate for a solution to Yang-Mills gradient-like flow]
\label{lem:Rade_7-3}
Let $G$ be a compact Lie group and $P$ a principal $G$-bundle over a closed, connected, orientable,
smooth manifold, $X$, of dimension $d$ with $2 \leq d\leq 4$ and Riemannian metric, $g$. Let $A_1$ be a reference connection of class\footnote{This regularity assumption for $A_1$ is stronger than necessary but simplifies the discussion.}
$C^\infty$ on $P$. Then there are positive constants, $C = C(A_1, d, g)$ and $\eps_1 = \eps_1(A_1, d, g) \in (0, 1]$, such that if $A(t)$ is a strong solution to the Yang-Mills gradient-like flow on $P$ over an interval $(S, T)$,
\begin{equation}
\label{eq:Yang-Mills_gradientlike_flow_equation}
\dot A = -d_A^*F_A + R \quad\hbox{a.e. on }(S,T)\times X,
\end{equation}
for a source term, $R$, satisfying
\begin{multline*}
R \in L^1_{\loc}(S, T; L^\infty(X; \Lambda^1\otimes\ad P))
\cap L^1_{\loc}(S, T; H^2_{A_1}(X; \Lambda^1\otimes\ad P))
\\
\cap W^{1,1}_{\loc}(S, T; L^2(X; \Lambda^1\otimes\ad P)),
\end{multline*}
and regularity,
\begin{multline*}
A - A_1 \in L^1_{\loc}(S, T; H^4_{A_1}(X; \Lambda^1\otimes\ad P))
\cap W^{1,2}_{\loc}(S, T; H^1_{A_1}(X; \Lambda^1\otimes\ad P))
\\
\cap W^{1,1}_{\loc}(S, T; H^2_{A_1}(X; \Lambda^1\otimes\ad P))
\cap W^{2,1}_{\loc}(S, T; L^2(X; \Lambda^1\otimes\ad P)),
\end{multline*}
where $S \in \RR$ and $\delta \in (0, 1/16]$ and $T$ obey $S + 2\delta \leq T \leq \infty$, and
\begin{equation}
\label{eq:Linfinity_in_time_H1_in_space_small_norm_At_minus_A_1_condition_lemma_7.3}
\|A(t) - A_1\|_{H^1_{A_1}(X)} \leq \eps_1 \quad\forall\, t \in (S, T),
\end{equation}
then
\begin{multline}
\label{eq:Rade_apriori_interior_estimate_lemma_7.3}
\int_{S+\delta}^T \|\dot A(t)\|_{H_{A_1}^1(X)}\,dt \leq C(1 + \delta^{-1/2})\int_S^T \|\dot A(t)\|_{L^2(X)}\,dt
\\
+ C\sqrt{\delta}\int_{S+\delta/2}^T \left(\|R(t)\|_{H^2_{A_1}(X)} + \|R(t)\|_{L^\infty(X)} + \|\dot R\|_{L^2(X)}\right)\,dt.
\end{multline}
\end{lem}

\begin{rmk}[On the choice of connection, $A_1$, in Lemma \ref{lem:Rade_7-3}]
\label{rmk:Rade_7-3}
In R\r{a}de's statement of his \cite[Lemma 7.3]{Rade_1992}, he assumes that the connection $A_1$ in the expression $A(t)-A_1$ in \eqref{eq:Linfinity_in_time_H1_in_space_small_norm_At_minus_A_1_condition_lemma_7.3} (denoted $A_\infty$ in
\cite[Lemma 7.3]{Rade_1992}) is Yang-Mills but an examination of his proof reveals that this hypothesis is never used. We may choose any $C^\infty$ connection, $A_\infty$, to define a difference, $a(t) = A(t)-A_\infty$, in the proof of Lemma \ref{lem:Rade_7-3} and a different $C^\infty$ connection, $A_1$, to define a norm for the Sobolev space, $H^1_{A_1}(X; \Lambda^1\otimes\ad P)$ (that choice is suppressed in the statement of \cite[Lemma 7.3]{Rade_1992}). R\r{a}de's assumption that his connection, $A_\infty$, is Yang-Mills is made only for the sake of notational consistency with the statement of his version \cite[Proposition 7.2]{Rade_1992} of the {\L}ojasiewicz-Simon inequality, which does require a Yang-Mills connection. In our application of Lemma \ref{lem:Rade_7-3}, just as in \cite{Rade_1992}, we will need to assume that the hypothesis \eqref{eq:Linfinity_in_time_H1_in_space_small_norm_At_minus_A_1_condition_lemma_7.3} holds when $A_1$ in the expression $A(t)-A_1$ in \eqref{eq:Linfinity_in_time_H1_in_space_small_norm_At_minus_A_1_condition_lemma_7.3} is a Yang-Mills connection, although the Sobolev norms defined by $A_1$ appearing in Lemma \ref{lem:Rade_7-3} could be replaced by equivalent norms defined by a different $C^\infty$ reference connection, $A_{\myref}$.
\end{rmk}

\begin{proof}
We may assume without loss of generality that $S = 0$. The solution, $A(t)$ to the Yang-Mills gradient-like flow \eqref{eq:Yang-Mills_gradientlike_flow_equation} satisfies the estimate \eqref{eq:Linfinity_in_time_H1_in_space_small_norm_At_minus_A_1_condition_lemma_7.3} on for all $t\in [0, T]$ and in particular for all $t \in [0, 2t_0]$, for any $t_0$ obeying $0 < 2t_0 \leq T$ and $t_0 \geq \delta$, that is,
\begin{equation}
\label{eq:Linfinity_in_time_H1_in_space_small_norm_At_minus_A_1_0_to_2tzero}
\|A(t) - A_1\|_{H^1_{A_1}(X)} \leq \eps_1 \quad\forall\, t \in [0, 2t_0].
\end{equation}
For the remainder of the proof, we assume without loss of generality that $t_0 = \delta$.

For $A(t)$ obeying \eqref{eq:Yang-Mills_gradientlike_flow_equation}, we have
$$
\frac{\partial F_A}{\partial t} = d_A \dot A = -d_Ad_A^*F_A + d_A R \quad\hbox{(by \eqref{eq:Yang-Mills_gradientlike_flow_equation})},
$$
and so, by the Bianchi identity \eqref{eq:Freed_Uhlenbeck_2-16_Bianchi_identity},
one obtains a gradient-like analogue of the familiar \cite[Equation (4.2)]{Rade_1992},
\begin{equation}
\label{eq:Rade_4-2_gradientlike}
\frac{\partial F_A}{\partial t} = -\Delta_AF_A + d_A R
\end{equation}
where $\Delta_A = d_Ad_A^* + d_A^*d_A$ is the Hodge Laplace operator \eqref{eq:Lawson_page_93_Hodge_Laplacian} on $\Omega^2(X;\ad P)$. We seek a similar parabolic equation for $\dot A$, that is, an equality between $\partial\dot A/\partial t + \Delta_A \dot A$ and a source term involving $\dot R$, spatial derivatives of $R$, and lower-order spatial derivatives of $\dot A$.

Recalling that (due to \cite[Equation (6.2)]{Warner}),
\begin{equation}
\label{eq:Warner_6-2}
d_A^* = (-1)^{-d(p+1)+1}*d_A* \quad\hbox{on } \Omega^p(X;\ad P),
\end{equation}
we see that $b(t) := \dot A(t)$ obeys the partial differential equation,
\begin{align*}
\dot b &= \frac{\partial \dot A}{\partial t} = \frac{\partial }{\partial t}(-d_A^*F_A + R) \quad\hbox{(by \eqref{eq:Yang-Mills_gradientlike_flow_equation})}
\\
&= - d_A^*\frac{\partial F_A}{\partial t} - (-1)^{-3d+1}*\left[\frac{\partial A}{\partial t}, *F_A\right] + \dot R
\\
&= - d_A^*d_A\dot A + (-1)^{-3d} *[\dot A, *F_A] + \dot R,
\end{align*}
and thus, as $b = \dot A$,
\begin{equation}
\label{eq:Parabolic_equation_for_time_derivative_of_solution_to_Yang-Mills_gradientlike_flow_prelim}
\dot b = - d_A^*d_A b + (-1)^{-3d} *[b, *F_A] + \dot R.
\end{equation}
Recalling that $\Delta_A = d_A^*d_A + d_Ad_A^*$ denotes the covariant Hodge Laplace operator \eqref{eq:Lawson_page_93_Hodge_Laplacian} on $\Omega^1(X,\ad P)$, we have
\begin{align*}
\Delta_A b &= (d_A^*d_A + d_Ad_A^*)b
\\
&= d_A^*d_A b + d_Ad_A^*(-d_A^*F_A + R) \quad\hbox{(by \eqref{eq:Yang-Mills_gradientlike_flow_equation})}
\\
&= d_A^*d_A b - d_Ad_A^*d_A^*F_A + d_Ad_A^*R
\\
&= d_A^*d_A b + d_Ad_A^*R \quad\hbox{(by \eqref{eq:Donaldson_Kronheimer_page_235})},
\end{align*}
and so
$$
d_A^*d_A b = \Delta_A b - d_Ad_A^*R.
$$
By combining the preceding identity with \eqref{eq:Parabolic_equation_for_time_derivative_of_solution_to_Yang-Mills_gradientlike_flow_prelim}, we see that $b$ obeys the parabolic equation,
$$
\dot b = - \Delta_A b + (-1)^{-3d}*[b, *F_A] + d_Ad_A^*R + \dot R  \quad\hbox{on } (0,T)\times X,
$$
It follows from the Bochner-Weitzenb\"ock formula \eqref{eq:Lawson_corollary_II-2} (compare \cite[Equation (2.3)]{Rade_1992}) and by absorbing negative signs into the implied definitions of the bilinear operations, $\times$, that
\begin{equation}
\label{eq:Parabolic_equation_for_time_derivative_of_solution_to_Yang-Mills_gradientlike_flow_simple}
\dot b + \nabla_A^*\nabla_A b = F_A\times b + \Ric_g\times b + d_Ad_A^*R + \dot R  \quad\hbox{on }(0,T)\times X.
\end{equation}
Writing $A(t) = A_1 + a(t)$ on $(0,T)\times X$ and $F_A = F_{A_1} + d_{A_1}a + [a, a]$, we can formally expand this equation on $(0,T)\times X$ as\footnote{Compare R\r{a}de's equation for $\partial G/dt + \nabla_{A_\infty}^*\nabla_{A_\infty} G$ in \cite[page 156]{Rade_1992}, where $G = \dot A$.}
\begin{multline}
\label{eq:Parabolic_equation_for_time_derivative_of_solution_to_Yang-Mills_gradientlike_flow}
\dot b + \nabla_{A_1}^*\nabla_{A_1} b = a\times\nabla_{A_1}b + \nabla_{A_1} a\times b + a\times a \times b + F_{A_1}\times b + \Ric_g\times b
\\
+ \nabla_{A_1}^2 R + a\times \nabla_{A_1} R + \nabla_{A_1} a\times R +  a\times a\times R + \dot R.
\end{multline}
Let $\eta \in C^\infty(\RR)$ be a cut-off function with $\eta = 0$ on $(-\infty, t_0/2]$ and $\eta = 1$ on $[t_0, \infty)$ and $0 \leq \eta \leq 1$ on $\RR$. Then, on $(0,T)\times X$,
\begin{multline}
\label{eq:Parabolic_equation_for_cut-off_time_derivative_of_solution_to_Yang-Mills_gradientlike_flow}
\frac{\partial (\eta b)}{\partial t} + \nabla_{A_1}^*\nabla_{A_1} (\eta b) = a\times\nabla_{A_1}(\eta b)
+ \nabla_{A_1} a\times \eta b + a\times a \times \eta b
+ F_{A_1}\times \eta b + \Ric_g\times \eta b
\\
+ \dot\eta b
+ \nabla_{A_1}^2 \eta R + a\times \nabla_{A_1} \eta R + \nabla_{A_1} a\times \eta R +  a\times a\times \eta R + \eta \dot R.
\end{multline}
It follows from \eqref{eq:Rade_11-3_and_11-4}, for $t_0 \leq 1/16$ (this is permissible because $t_0 = \delta$ by assumption and $\delta \in (0, 1/16]$ by hypothesis) that
\begin{equation}
\label{eq:Rade_final_inequality_page_156}
\begin{aligned}
\|\eta b\|_{L^2(0, 2t_0; H_{A_1}^1(X))}
&\leq 2\left(\|a\times\nabla_{A_1}(\eta b)\|_{L^2(0, 2t_0; H_{A_1}^{-1}(X))} + \|\nabla_{A_1} a\times \eta b\|_{L^2(0, 2t_0; H_{A_1}^{-1}(X))} \right.
\\
&\quad + \left. \|a\times a \times \eta b\|_{L^2(0, 2t_0; H_{A_1}^{-1}(X))} \right)
+ 2\sqrt{2}\left(\|F_{A_1}\times \eta b\|_{L^1(0, 2t_0; L^2(X))} \right.
\\
&\quad + \|\Ric_g\times \eta b\|_{L^1(0, 2t_0; L^2(X))} + \|\dot\eta b\|_{L^1(0, 2t_0; L^2(X))}
\\
&\quad + \|\nabla_{A_1}^2 \eta R\|_{L^1(0, 2t_0; L^2(X))} + \|a\times \nabla_{A_1} \eta R\|_{L^1(0, 2t_0; L^2(X))}
\\
&\quad + \|\nabla_{A_1} a\times \eta R\|_{L^1(0, 2t_0; L^2(X))} + \|a\times a\times \eta R\|_{L^1(0, 2t_0; L^2(X))}
\\
&\quad + \left. \|\eta\dot R\|_{L^1(0, 2t_0; L^2(X))} \right).
\end{aligned}
\end{equation}
Recall from Lemma \ref{lem:Freed_Uhlenbeck_equation_6-34_minus_k_k2}, with $k = 0$, $k_1=k_2=1$, $p=p_1=p_2=2$, and $2\leq d\leq 4$, that the following Sobolev multiplication map is continuous, in the case of pointwise multiplication of real or complex-valued functions,
\begin{equation}
\label{eq:Sobolev_multiplication_H1_times_L2_into_Hminus1}
H^1(X) \times L^2(X) \to H^{-1}(X),
\end{equation}
and similarly in the case of pointwise tensor product of sections of vector bundles. For the first term on the right-hand side of \eqref{eq:Rade_final_inequality_page_156}, the Sobolev multiplication result \eqref{eq:Sobolev_multiplication_H1_times_L2_into_Hminus1} (see also the \cite[Appendix]{Rade_1992}) with Sobolev multiplication constant, $C>0$, for $2 \leq d \leq 4$, implies that
\begin{align*}
\|a\times\nabla_{A_1}(\eta b)\|_{L^2(0, 2t_0; H_{A_1}^{-1}(X))}
&\leq C\|a\|_{L^\infty(0, 2t_0; H_{A_1}^1(X))} \|\nabla_{A_1}(\eta b)\|_{L^2(0, 2t_0; L^2(X))}
\\
&\leq C\eps_1 \|\eta b\|_{L^2(0, 2t_0; H_{A_1}^1(X))} \quad\hbox{(by \eqref{eq:Linfinity_in_time_H1_in_space_small_norm_At_minus_A_1_0_to_2tzero})}.
\end{align*}
For the second term on the right-hand side of \eqref{eq:Rade_final_inequality_page_156}, we have
\begin{align*}
\|\nabla_{A_1} a\times \eta b\|_{L^2(0, 2t_0; H_{A_1}^{-1}(X))}
&\leq C\|\nabla_{A_1} a\|_{L^\infty(0, 2t_0; L^2(X))} \|\eta b\|_{L^2(0, 2t_0; H_{A_1}^1(X))}
\\
&\leq C\eps_1 \|\eta b\|_{L^2(0, 2t_0; H_{A_1}^1(X))} \quad\hbox{(by \eqref{eq:Linfinity_in_time_H1_in_space_small_norm_At_minus_A_1_0_to_2tzero})}.
\end{align*}
For the third term on the right-hand side of \eqref{eq:Rade_final_inequality_page_156}, applying in addition the Kato Inequality \eqref{eq:FU_6-20_first-order_Kato_inequality}
and the Sobolev embedding \cite[Theorem 4.12]{AdamsFournier}, $H^1(X) \hookrightarrow L^4(X)$ for $2 \leq d \leq 4$, we obtain
\begin{align*}
\|a\times a \times \eta b\|_{L^2(0, 2t_0; H_{A_1}^{-1}(X))}
&\leq C\|a\times a\|_{L^\infty(0, 2t_0; L^2(X))} \|\eta b\|_{L^2(0, 2t_0; H_{A_1}^1(X))}
\\
&\leq C\|a\|_{L^\infty(0, 2t_0; L^4(X))}^2 \|\eta b\|_{L^2(0, 2t_0; H_{A_1}^1(X))}
\\
&\leq C\|a\|_{L^\infty(0, 2t_0; H_{A_1}^1(X))}^2 \|\eta b\|_{L^2(0, 2t_0; H_{A_1}^1(X))}
\\
&\leq C\eps_1^2 \|\eta b\|_{L^2(0, 2t_0; H_{A_1}^1(X))} \quad\hbox{(by \eqref{eq:Linfinity_in_time_H1_in_space_small_norm_At_minus_A_1_0_to_2tzero})}.
\end{align*}
For the fourth term on the right-hand side of \eqref{eq:Rade_final_inequality_page_156}, we see that
\begin{align*}
\|F_{A_1}\times \eta b\|_{L^1(0, 2t_0; L^2(X))}
&\leq C\|F_{A_1}\|_{L^\infty(X)} \|\eta b\|_{L^1(0, 2t_0; L^2(X))}
\\
&\leq C\|\eta b\|_{L^1(0, 2t_0; L^2(X))},
\end{align*}
where $C$ depends on $F_{A_1}$ and $g$, and similarly for the fifth term on the right-hand side of \eqref{eq:Rade_final_inequality_page_156}, we find that
\begin{align*}
\|\Ric_g\times \eta b\|_{L^1(0, 2t_0; L^2(X))}
&\leq C\|\Ric_g\|_{L^\infty(X)} \|\eta b\|_{L^1(0, 2t_0; L^2(X))}
\\
&\leq C\|\eta b\|_{L^1(0, 2t_0; L^2(X))},
\end{align*}
where $C$ depends on $g$.

For the seventh through tenth terms on the right-hand side of \eqref{eq:Rade_final_inequality_page_156}, we can again exploit the \apriori bound \eqref{eq:Linfinity_in_time_H1_in_space_small_norm_At_minus_A_1_0_to_2tzero} and our assumption that $0 < \eps_1 \leq 1$ to give
\begin{align*}
{}&\|\nabla_{A_1}^2 \eta R\|_{L^1(0, 2t_0; L^2(X))}
+ \|a\times \nabla_{A_1} \eta R\|_{L^1(0, 2t_0; L^2(X))} + \|\nabla_{A_1} a\times \eta R\|_{L^1(0, 2t_0; L^2(X))}
\\
&\qquad + \|a\times a\times \eta R\|_{L^1(0, 2t_0; L^2(X))}
\\
&\quad \leq C\left( \|\nabla_{A_1}^2 \eta R\|_{L^1(0, 2t_0; L^2(X))}
 + \|a\|_{L^\infty(0, 2t_0; L^4(X))}\|\nabla_{A_1} \eta R\|_{L^1(0, 2t_0; L^4(X))} \right.
\\
&\qquad + \left. \|\nabla_{A_1} a\|_{L^\infty(0, 2t_0; L^2(X))}\|\eta R\|_{L^1(0, 2t_0; L^\infty(X))}
+ \|a\|_{L^\infty(0, 2t_0; L^4(X))}^2 \|\eta R\|_{L^1(0, 2t_0; L^\infty(X))} \right)
\\
&\quad \leq C\left(\|\eta R\|_{L^1(0, 2t_0; H^2_{A_1}(X))} + \|\eta R\|_{L^1(0, 2t_0; L^\infty(X))} \right).
\end{align*}
By combining the preceding inequalities with \eqref{eq:Rade_final_inequality_page_156}, we obtain
\begin{multline*}
\|\eta b\|_{L^2(0, 2t_0; H_{A_1}^1(X))}
\leq
C\eps_1\|\eta b\|_{L^2(0, 2t_0; H_{A_1}^1(X))} + 2\sqrt{2}\|\dot\eta b\|_{L^1(0, 2t_0; L^2(X))}
+ C\|\eta b\|_{L^1(0, 2t_0; L^2(X))}
\\
+ C\left(\|\eta R\|_{L^1(0, 2t_0; H^2_{A_1}(X))} + \|\eta R\|_{L^1(0, 2t_0; L^\infty(X))} + \|\eta\dot R\|_{L^1(0, 2t_0; L^2(X))}\right).
\end{multline*}
For $0<\eps_1 \leq 1/(2C)$, rearrangement in the preceding inequality yields
\begin{multline*}
\|\eta b\|_{L^2(0, 2t_0; H_{A_1}^1(X))}
\leq
4\sqrt{2}\|\dot\eta b\|_{L^1(0, 2t_0; L^2(X))} + 2C\|\eta b\|_{L^1(0, 2t_0; L^2(X))}
\\
+ 2C\left(\|\eta R\|_{L^1(0, 2t_0; H^2_{A_1}(X))} + \|\eta R\|_{L^1(0, 2t_0; L^\infty(X))} + \|\eta\dot R\|_{L^1(0, 2t_0; L^2(X))}\right).
\end{multline*}
We conclude, for $|\dot\eta| \leq 4t_0^{-1}$ on $\RR$, that
\begin{align*}
\|b\|_{L^2(t_0, 2t_0; H_{A_1}^1(X))} &\leq \|\eta b\|_{L^2(0, 2t_0; H_{A_1}^1(X))}
\\
&\leq 4\sqrt{2}\|\dot\eta b\|_{L^1(0, 2t_0; L^2(X))} + 2C\|\eta b\|_{L^1(0, 2t_0; L^2(X))}
\\
&\quad + 2C\left(\|\eta R\|_{L^1(0, 2t_0; H^2_{A_1}(X))} + \|\eta R\|_{L^1(0, 2t_0; L^\infty(X))} + \|\eta\dot R\|_{L^1(0, 2t_0; L^2(X))} \right)
\\
&\leq \left(16\sqrt{2}\,t_0^{-1} + 2C\right)\|b\|_{L^1(0, 2t_0; L^2(X))}
\\
&\quad + 2C\left(\|R\|_{L^1(t_0/2, 2t_0; H^2_{A_1}(X))} + \|R\|_{L^1(t_0/2, 2t_0; L^\infty(X))} + \|\dot R\|_{L^1(t_0/2, 2t_0; L^2(X))} \right).
\end{align*}
Consequently, assuming without loss of generality that $0 < t_0 \leq 1$, we have
\begin{align*}
{}&\|b\|_{L^1(t_0, 2t_0; H_{A_1}^1(X))}
\\
&\leq t_0^{1/2}\|b\|_{L^2(t_0, 2t_0; H_{A_1}^1(X))}
\\
&\leq \left(16\sqrt{2}\,t_0^{-1/2} + 2Ct_0^{1/2}\right)\|b\|_{L^1(0, 2t_0; L^2(X))}
\\
&\quad + 2t_0^{1/2}C\left(\|R\|_{L^1(t_0/2, 2t_0; H^2_{A_1}(X))} + \|R\|_{L^1(t_0/2, 2t_0; L^\infty(X))} + \|\dot R\|_{L^1(t_0/2, 2t_0; L^2(X))} \right)
\\
&\leq \left(16\sqrt{2}\,t_0^{-1/2} + 2C\right)\|b\|_{L^1(0, 2t_0; L^2(X))}
\\
&\quad + 2t_0^{1/2}C\left(\|R\|_{L^1(t_0/2, 2t_0; H^2_{A_1}(X))} + \|R\|_{L^1(t_0/2, 2t_0; L^\infty(X))} + \|\dot R\|_{L^1(t_0/2, 2t_0; L^2(X))} \right).
\end{align*}
To finish the proof of Lemma \ref{lem:Rade_7-3}, we divide the interval, $[0, T]$, into (possibly infinitely many) subintervals of length $t_0$, and we apply the last inequality on each pair of consecutive subintervals. In writing the interior \apriori estimate \eqref{eq:Rade_apriori_interior_estimate_lemma_7.3}, we recall our assumption that $t_0 = \delta$.
\end{proof}

\subsection{A generalization of R\r{a}de's Lemma 8.1}
\label{subsec:Rade_Lemma_8-1_generalization_dimension_2_leq_d_leq_4}
Next we observe that a slight modification of the argument in \cite{Rade_1992}, for $\dim X=2$ or $3$ (using a different Sobolev embedding) yields the following version of \cite[Lemma 8.1]{Rade_1992}, with a weaker conclusion valid when $\dim X=4$ (though not $\dim X \geq 5$). R\r{a}de noted in \cite{Rade_1992} that estimates similar to that in his \cite[Lemma 8.1]{Rade_1992} had been obtained by Sadun (for example, see \cite[Section 6.2]{Sadun_1987thesis}).
We include the details of the proof for completeness since Lemma \ref{lem:Rade_8-1} enhances \cite[Lemma 8.1]{Rade_1992} in several ways. A similar analysis was outlined in \cite[p. 221]{DK}, but not pursued, presumably due to the much weaker conclusion one can obtain in general when $\dim X = 4$ has dimension four relative to the case of $\dim X=2$ or $3$.

To motivate our Lemma \ref{lem:Rade_8-1}, suppose first that $v \in C^\infty([0,\infty)\times X; \Lambda^2\otimes\ad P)$ is a solution to a homogeneous \emph{linear} second-order parabolic equation,
\begin{equation}
\label{eq:Evans_theorem_7-1-5_linear_second-order_parabolic_system}
\frac{\partial v}{\partial t} + \Delta_{A_1} v = 0 \quad\hbox{on }(0,\infty)\times X,
\end{equation}
analogous to the a \emph{nonlinear} second-order parabolic equation \eqref{eq:Rade_4-3} for $F_A(t)$, arising in the proof of Lemma \ref{lem:Rade_8-1}, and where $A_1$ is a fixed, $C^\infty$ reference connection on $P$. The standard \apriori estimates in \cite[Theorems 7.1.2 and 7.1.5]{Evans2} for a solution $v$ to a linear, scalar second-order parabolic equation on a bounded open subset $U \subset \RR^d$ with $C^\infty$ boundary, $\partial U$, and a homogeneous Dirichlet boundary condition along $\partial U$, suggests that $v$ should obey the \apriori estimates,
\begin{align}
\label{eq:Evans_7-1-20}
\|v\|_{L^\infty(0,T;L^2(X))} + \|v\|_{L^2(0,T;H^1_{A_1}(X))} + \|\dot v\|_{L^2(0,T;H^{-1}(X))}
&\leq C\|v\|_{L^2(X)},
\\
\label{eq:Evans_7-1-46}
\|v\|_{L^\infty(0,T;H^1_{A_1}(X))} + \|v\|_{L^2(0,T;H^2_{A_1}(X))} + \|\dot v\|_{L^2(0,T;L^2(X))}
&\leq C\|v\|_{H^1_{A_1}(X)},
\\
\label{eq:Evans_7-1-47}
\|v\|_{L^\infty(0,T;H^2_{A_1}(X))} + \|\dot v\|_{L^2(0,T;H^1_{A_1}(X))} + \|\dot v\|_{L^\infty(0,T;L^2(X))}
&\leq C\|v\|_{H^2_{A_1}(X)},
\end{align}
for a positive constant $C$ depending on the Riemannian metric $g$, the time interval length, $0<T<\infty$, and, possibly, the connection, $A_1$, respectively for $v$ with the following minimal regularities,
\begin{align*}
v &\in L^2(0,T; H_{A_1}^1(X; \Lambda^2\otimes\ad P)) \cap C([0,T]; L^2(X; \Lambda^2\otimes\ad P)),
\quad \dot v \in L^2(0,T; H_{A_1}^{-1}(X; \Lambda^2\otimes\ad P)),
\\
v &\in L^2(0,T; H_{A_1}^2(X; \Lambda^2\otimes\ad P)) \cap C([0,T; H_{A_1}^1(X; \Lambda^2\otimes\ad P))
\quad \dot v \in L^2(0,T; L^2(X; \Lambda^2\otimes\ad P)),
\\
v &\in C([0,T]; H_{A_1}^2(X; \Lambda^2\otimes\ad P)),
\quad \dot v \in L^2(0,T; H_{A_1}^1(X; \Lambda^2\otimes\ad P)) \cap C([0,T; L^2(X; \Lambda^2\otimes\ad P)).
\end{align*}
Provided $T$ is sufficiently small, as determined by $\|F_A(0)\|_{L^2(X)}$, we shall see in Lemma \ref{lem:Rade_8-1} that this is indeed the type of \apriori estimate obeyed by $F_A(t)$, but with positive constant $C$ depending only on the Riemannian metric $g$ and \emph{independent} of $A(t)$ or any fixed reference connection, $A_1$, on $P$. We shall use the well-known pointwise \emph{Kato Inequality} \cite[Equation (6.20)]{FU},
\begin{equation}
\label{eq:FU_6-20_first-order_Kato_inequality}
|d|v|| \leq |\nabla v| \quad\hbox{a.e. on } X,
\end{equation}
for a section $v \in W^{1, p}(X; E)$ of a Hermitian or Riemannian vector bundle, $E$, over a Riemannian manifold, $X$, of dimension $d \geq 2$ and covariant derivative, $\nabla$, defined by a compatible connection.

\begin{lem}[\Apriori estimate for the curvature of a solution to Yang-Mills gradient flow]
\label{lem:Rade_8-1}
Let $G$ be a compact Lie group and $P$ a principal $G$-bundle over a closed, connected, orientable, smooth manifold, $X$, of dimension $d$ with $2 \leq d\leq 4$ and Riemannian metric, $g$. Then there is a positive constant, $c = c(d, g)$, such that if $A(t)$ is a solution to the Yang-Mills gradient flow on $(0, T)$, in the sense that
$$
A - A_1 \in L^2_{\loc}(0, T; H^3_{A_1}(X; \Lambda^1\otimes\ad P)),
$$
for a $C^\infty$ reference connection, $A_1$, on $P$, then, for $d = 4$,
\begin{multline}
\label{eq:lemma_Rade_8-1_apriori_estimate_nabla_A_FA_d_is_4}
\int_0^T \int_X |\nabla_{A(t)}F_A(t)|^2\,d\vol_g\,dt
\\
\leq
cT\left(\int_X |F_A(0)|^2\,d\vol_g\right)^{3/2}
+ c\left(\int_X |F_A(0)|^2\,d\vol_g\right)^{1/2}\int_0^T\int_X|\nabla_{A(t)} F_A(t)|^2\,d\vol_g\,dt,
\end{multline}
and\footnote{There is a typographical error in R\r{a}de's application of the Sobolev embedding $W^{1,2}(X) \hookrightarrow L^6(X)$ when $d = 2$ or $3$, which we correct here, and that explains the discrepancy between the bound \eqref{eq:lemma_Rade_8-1_apriori_estimate_nabla_A_FA_d_is_2_or_3} and that in \cite[Lemma 8.1]{Rade_1992}},
for $d = 2$ or $3$,
\begin{equation}
\label{eq:lemma_Rade_8-1_apriori_estimate_nabla_A_FA_d_is_2_or_3}
\int_0^T \int_X |\nabla_{A(t)}F_A(t)|^2\,d\vol_g\,dt
\leq
cT\left(\left(\int_X |F_A(0)|^2\,d\vol_g\right)^3 + \int_X |F_A(0)|^2\,d\vol_g\right).
\end{equation}
\end{lem}

\begin{rmk}[Alternative proofs of Lemma \ref{lem:Rade_8-1}]
\label{rmk:Alternative_proof_Rade_Lemma_8-1}
The \apriori estimates \eqref{eq:lemma_Rade_8-1_apriori_estimate_nabla_A_FA_d_is_4} and \eqref{eq:lemma_Rade_8-1_apriori_estimate_nabla_A_FA_d_is_2_or_3} could also be deduced as a consequence of the \apriori estimate \eqref{eq:Rade_11-3_and_11-4} (and its proof) or the \apriori estimates \eqref{eq:Sell_You_42-30_Linfty_time_Valpha_space_u} and \eqref{eq:Sell_You_42-30_L2_time_Valpha+1_space_u}  (and their proofs) in Theorem \ref{thm:Sell_You_42-12} by taking $\alpha = 0$, together with the Bochner-Weitzenb\"ock formula \eqref{eq:Lawson_corollary_II-3}.
\end{rmk}

\begin{proof}[Proof of Lemma \ref{lem:Rade_8-1}]
Recall that $\|F_A(t)\|_{L^2(X)}^2$ is a non-increasing function of $t\geq 0$ when $A(t)$ is a solution to Yang-Mills gradient flow. (This a general property of gradient flow --- see the proof of Proposition \ref{prop:Huang_3-3-2}.)  The curvature, $F_A$, satisfies \cite[Equation (4.3)]{Rade_1992} (see also \cite[p. 221]{DK} or \cite[Equation (11)]{Struwe_1994}), namely
\begin{equation}
\label{eq:Rade_4-3}
\frac{\partial F_A}{\partial t} + \Delta_A F_A = 0, \quad\hbox{a.e. on } (0,T)\times X,
\end{equation}
where we denote $\Delta_A = d_Ad_A^* + d_A^*d_A$ on $\Omega^2(X,\ad P)$ here. Using the Bochner-Weitzenb\"ock formula \eqref{eq:Lawson_corollary_II-3} (compare \cite[Equation (2.3)]{Rade_1992}), the preceding equation can be written
\begin{equation}
\label{eq:Rade_4-3_Bochner-Weitzenbock}
\frac{\partial F_A}{\partial t} + \nabla_A^*\nabla_A F_A = - F_A\times F_A - \Riem_g \times F_A,
\quad\hbox{a.e. on } (0,T)\times X.
\end{equation}
We take the $L^2$ inner product of the preceding equation with $F_A$ to give
\begin{multline*}
\int_0^T\int_X \left\langle \frac{\partial F_A}{\partial t}, F_A \right\rangle \,d\vol_g\,dt
+ \int_0^T\int_X \langle \nabla_A^*\nabla_A F_A, F_A \rangle \,d\vol_g\,dt
\\
= - \int_0^T\int_X \langle F_A\times F_A, F_A \rangle \,d\vol_g\,dt
- \int_0^T\int_X \langle \Riem_g \times F_A, F_A \rangle \,d\vol_g\,dt.
\end{multline*}
We now integrate by parts over $X$ and use the pointwise identity
$$
\left\langle \frac{\partial F_A}{\partial t}, F_A \right\rangle
=
\frac{1}{2}\frac{\partial}{\partial t} |F_A|^2, \quad\hbox{a.e. on } (0,T)\times X,
$$
to get
\begin{multline*}
\frac{1}{2}\int_0^T\int_X \frac{\partial}{\partial t} |F_A|^2 \,d\vol_g\,dt
+ \int_0^T\int_X |\nabla_A F_A|^2 \,d\vol_g\,dt
\\
= - \int_0^T\int_X \langle F_A\times F_A, F_A \rangle \,d\vol_g\,dt
- \int_0^T\int_X \langle \Riem_g \times F_A, F_A \rangle \,d\vol_g\,dt,
\end{multline*}
and therefore, with $c_0$ denoting a positive constant depending at most on the Riemannian metric, $g$, on $X$,
\begin{multline*}
\frac{1}{2}\int_X |F_A(T)|^2 \,d\vol_g - \frac{1}{2}\int_X |F_A(0)|^2 \,d\vol_g
+ \int_0^T\int_X |\nabla_A F_A|^2 \,d\vol_g\,dt
\\
\leq c_0 \int_0^T\int_X \left(|F_A|^3 + |\Riem_g| |F_A|^2\right)\,d\vol_g\,dt.
\end{multline*}
Thus,
\begin{align*}
\int_0^T\int_X |\nabla_A F_A|^2 \,d\vol_g\,dt
&\leq \frac{1}{2}\int_X |F_A(0)|^2 \,d\vol_g - \frac{1}{2}\int_X |F_A(T)|^2 \,d\vol_g
\\
&\quad + c_0 \int_0^T\int_X \left(|F_A|^3 + |\Riem_g| |F_A|^2\right)\,d\vol_g\,dt
\\
&\leq \frac{1}{2}\int_X |F_A(0)|^2 \,d\vol_g + c_0 \int_0^T\int_X \left(|F_A|^3 + |\Riem_g| |F_A|^2\right)\,d\vol_g\,dt.
\end{align*}
Again since $\int_X |F_A(t)|^2\,d\vol_g$ is a non-increasing function of $t \geq 0$, we have
$$
\int_0^T\int_X |\Riem_g||F_A|^2\,d\vol_g\,dt \leq c_1T\int_X |F_A(0)|^2\,d\vol_g,
$$
where $c_1$ is a positive constant depending at most on the Riemannian metric, $g$, on $X$.

We estimate the cubic term by applying the Cauchy-Schwarz inequality to $|F_A|^3 = |F_A||F_A|^2$ and the Kato Inequality \eqref{eq:FU_6-20_first-order_Kato_inequality} to give $|\nabla|F_A|| \leq |\nabla_A F_A|$ (a.e. on $X$). Suppose first that $d=4$, so the Sobolev embedding $W^{1,2}(X;\RR) \hookrightarrow L^4(X;\RR)$ \cite[Theorem 4.12]{AdamsFournier} with Sobolev embedding constant $c_2$, depending at most on the Riemannian metric, $g$, on $X$, yields
$$
\left(\int_X |f|^4 \,d\vol_g\right)^{1/4}
\leq
c_2\left(\int_X (|\nabla f|^2 + |f|^2)\,d\vol_g\right)^{1/2}, \quad \forall\, f \in W^{1,2}(X;\RR).
$$
Then we get
\begin{align*}
{}&\int_0^T\int_X |F_A|^3\,d\vol_g\,dt
\\
&\quad \leq \int_0^T \left(\int_X |F_A|^2\,d\vol_g\right)^{1/2}\left(\int_X |F_A|^4\,d\vol_g\right)^{1/2}\,dt
\\
&\quad \leq c_2^2\int_0^T\left(\int_X |F_A|^2\,d\vol_g\right)^{1/2}
\left(\int_X \left(|F_A|^2 + |\nabla|F_A||^2\right)\,d\vol_g\right)\,dt
\\
&\quad \leq c_2^2\left(\int_X |F_A(0)|^2\,d\vol_g\right)^{1/2}
\left(T\int_X |F_A(0)|^2\,d\vol_g + \int_0^T\int_X|\nabla_A F_A|^2\,d\vol_g\,dt\right)
\\
&\quad = c_2^2T\left(\int_X |F_A(0)|^2\,d\vol_g\right)^{3/2}
+ c_2^2\left(\int_X |F_A(0)|^2\,d\vol_g\right)^{1/2}\int_0^T\int_X|\nabla_A F_A|^2\,d\vol_g\,dt.
\end{align*}
For $d=4$, the conclusion now follows by combining the preceding inequalities.

For the case $d = 2$ or $3$, we may use the Sobolev embedding $W^{1,2}(X;\RR) \hookrightarrow L^6(X;\RR)$ \cite[Theorem 4.12]{AdamsFournier} with Sobolev embedding constant $c_3$, depending at most on the Riemannian metric, $g$, on $X$, to give
$$
\left(\int_X |f|^6 \,d\vol_g\right)^{1/6}
\leq
c_3\left(\int_X (|\nabla f|^2 + |f|^2)\,d\vol_g\right)^{1/2}, \quad \forall\, f \in W^{1,2}(X;\RR).
$$
Then, applying the H\"older inequality with exponents $p = 4/3$ and $q=4$ to $|F_A|^3 = |F_A|^{3/2} |F_A|^{3/2}$,
$$
\int_X |F_A|^{3/2} |F_A|^{3/2} \,d\vol_g
\leq
\left(\int_X |F_A|^2 \,d\vol_g \right)^{3/4}
\left(\int_X |F_A|^6  \,d\vol_g\right)^{1/4},
$$
we get
$$
\int_0^T\int_X |F_A|^3\,d\vol_g\,dt
\leq
\int_0^T \left(\int_X |F_A|^2 \,d\vol_g \right)^{3/4} \left(\int_X |F_A|^6  \,d\vol_g\right)^{1/4}\,dt.
$$
We now apply Young's Inequality \cite[Equation (7.5)]{GilbargTrudinger}, $ab \leq a^p/p + b^q/q$ with $a, b > 0$ and $1/p + 1/q = 1$, thus $q = 4/3$ and $p = 4$, to give
\begin{multline*}
\left(\int_X |F_A|^2 \,d\vol_g \right)^{3/4} \left(\int_X |F_A|^6  \,d\vol_g\right)^{1/4}
\\
\leq
\frac{1}{4\eps^4}\left(\int_X |F_A|^2 \,d\vol_g \right)^3
+ \frac{3\eps^{4/3}}{4}\left(\int_X |F_A|^6  \,d\vol_g\right)^{1/3}.
\end{multline*}
Hence,
\begin{align*}
{}&\int_0^T\int_X |F_A|^3\,d\vol_g\,dt
\\
&\quad\leq \frac{1}{4\eps^4}\int_0^T \left(\int_X |F_A|^2 \,d\vol_g \right)^3 \,dt
+ \frac{3\eps^{4/3}}{4}\int_0^T \left(\int_X |F_A|^6  \,d\vol_g\right)^{1/3} \,dt
\\
&\quad\leq \frac{1}{4\eps^4}\int_0^T \left(\int_X |F_A|^2 \,d\vol_g \right)^3 \,dt
+ \frac{3\eps^{4/3}}{4}\int_0^T \left(|F_A|^2 + |\nabla|F_A||^2\right)\,d\vol_g \,dt
\\
&\quad\leq \frac{1}{4\eps^4}T \left(\int_X |F_A(0|^2 \,d\vol_g \right)^3
+ \frac{3\eps^{4/3}}{4}\left(T\int_X |F_A(0|^2 \,d\vol_g + \int_0^T\int_X |\nabla_AF_A|^2\,d\vol_g \,dt\right).
\end{align*}
For $d=2$ or $3$, the conclusion now follows by combining the preceding inequalities and rearrangement.
\end{proof}

Lemma \ref{lem:Rade_8-1} provides a good illustration of the impact of the dimension, $d$, of the base manifold, $X$, on the available strength of integral-norm bounds on the curvature, $F_A(t)$, for a solution, $A(t)$, to Yang-Mills gradient flow. For example, when $d = 2$ or $3$, the \apriori estimate \eqref{eq:lemma_Rade_8-1_apriori_estimate_L4_FA_d_is_2_or_3} below can be used to show that bubbling cannot occur, irrespective of the energy of the initial connection, $A(0)$, whereas when $d = 4$, the \apriori estimate \eqref{eq:lemma_Rade_8-1_apriori_estimate_L4_FA_d_is_4} below can only be used to show that bubbling cannot occur when the energy of the initial connection, $A(0)$, is sufficiently small.

\begin{cor}[\Apriori estimate for the curvature of a solution to Yang-Mills gradient flow]
\label{cor:Rade_8-1}
Assume the hypotheses of Lemma \ref{lem:Rade_8-1}. If $d = 2$ or $3$, then there is a positive constant $c$, depending at most on the Riemannian metric, $g$, on $X$, such that
\begin{equation}
\label{eq:lemma_Rade_8-1_apriori_estimate_L4_FA_d_is_2_or_3}
\|F_A\|_{L^2(0,T; L^6(X))} \leq c\sqrt{T}\left(\|F_A(0)\|_{L^2(X)}^3 + \|F_A(0)\|_{L^2(X)}\right).
\end{equation}
If $d=4$, then there are positive constants $c$ and $\eps$, depending at most on the Riemannian metric, $g$, on $X$, with the following significance: If $\|F_A(0)\|_{L^2(X)} < \eps$, then
\begin{equation}
\label{eq:lemma_Rade_8-1_apriori_estimate_L4_FA_d_is_4}
\|F_A\|_{L^2(0,T; L^4(X))} \leq c\sqrt{T}\|F_A(0)\|_{L^2(X)}^{3/2}.
\end{equation}
\end{cor}

\begin{proof}
The estimate \eqref{eq:lemma_Rade_8-1_apriori_estimate_L4_FA_d_is_2_or_3} follows from \eqref{eq:lemma_Rade_8-1_apriori_estimate_nabla_A_FA_d_is_2_or_3} in Lemma \ref{lem:Rade_8-1}, the Kato Inequality, and the Sobolev embedding $W^{1,2}(X;\RR) \hookrightarrow L^6(X)$. Similarly, the estimate \eqref{eq:lemma_Rade_8-1_apriori_estimate_L4_FA_d_is_4} follows from \eqref{eq:lemma_Rade_8-1_apriori_estimate_nabla_A_FA_d_is_4} in Lemma \ref{lem:Rade_8-1}, the Kato Inequality, and the Sobolev embedding $W^{1,2}(X;\RR) \hookrightarrow L^4(X)$.
\end{proof}

\subsection{An alternative to Lemma \ref{lem:Rade_7-3} via Lemma \ref{lem:Rade_7-3_abstract_L1_in_time_V2beta_space}}
\label{subsec:Rade_Lemma_7-3_generalization_H2beta_dimension_2_leq_d_leq_5}
By employing Sobolev spaces with fractional derivative index, we can apply Lemma \ref{lem:Rade_7-3_abstract_L1_in_time_V2beta_space} to give a result which complements Lemma \ref{lem:Rade_7-3}; see Remark \ref{rmk:Rade_7-3} for a discussion of the two roles played by the connection, $A_1$, in the statement of Lemma \ref{lem:Rade_7-3_L1_in_time_H2beta_in_space_apriori_estimate_by_L1_in_time_L2_in_space}.

\begin{lem}[An \apriori $L^1$-in-time interior estimate for a solution to Yang-Mills gradient-like flow]
\label{lem:Rade_7-3_L1_in_time_H2beta_in_space_apriori_estimate_by_L1_in_time_L2_in_space}
Let $G$ be a compact Lie group and $P$ a principal $G$-bundle over a closed, connected, orientable, smooth manifold, $X$, of dimension $d$ with $2 \leq d\leq 5$ and Riemannian metric, $g$. Let $A_1$ be a reference connection of class $C^\infty$ on $P$ and $\beta \in [1/4 + d/8, 1)$. Then there are positive constants, $C = C(A_1, d, g, \beta)$ and $\eps_1 = \eps_1(A_1, d, g, \beta)$, such that if $A(t)$ is a strong solution to a Yang-Mills gradient-like flow \eqref{eq:Yang-Mills_gradientlike_flow_equation} on $P$ over an interval $(S, T)$, with regularity,
$$
A - A_1 \in C((S, T); H^1_{A_1}(X; \Lambda^1\otimes\ad P))
\cap W^{1,1}_{\loc}(S, T; H^{2\beta}_{A_1}(X; \Lambda^1\otimes\ad P)),
$$
where $S \in \RR$ and $\delta > 0$ and $T$ obey $S + 2\delta \leq T \leq \infty$, and
\begin{equation}
\label{eq:Linfinity_in_time_H2beta_in_space_small_norm_At_minus_A_1_condition_lemma_7.3}
\|A(t) - A_1\|_{H^{2\beta}_{A_1}(X)} \leq \eps_1 \quad\forall\, t \in (S, T),
\end{equation}
then
\begin{multline}
\label{eq:Rade_7-3_L1_in_time_H2beta_in_space_apriori_estimate_by_L1_in_time_L2_in_space}
\int_{S+\delta}^T \|\dot A\|_{H_{A_1}^{2\beta}(X)}\,dt \leq C(1 + \delta^{-1})\int_S^T \|\dot A\|_{L^2(X)}\,dt
\\
+ C\int_{S+\delta/2}^T\left(\|R(t)\|_{H^2_{A_1}(X)} + \|R(t)\|_{L^\infty(X)} + \|\dot R(t)\|_{L^2(X)} \right)\,dt.
\end{multline}
\end{lem}

\begin{proof}
Write $a(t) = A(t) - A_1$ and $b(t) = \dot A(t)$ as in the proof of Lemma \ref{lem:Rade_7-3}. Adding $b$ to both sides of the evolution equation \eqref{eq:Parabolic_equation_for_time_derivative_of_solution_to_Yang-Mills_gradientlike_flow} yields the following identity on $(S,T)\times X$,
\begin{multline}
\label{eq:Augmented_parabolic_equation_for_time_derivative_of_solution_to_Yang-Mills_gradientlike_flow}
\dot b + \left(\nabla_{A_1}^*\nabla_{A_1} + 1\right)b
\\
=
a\times\nabla_{A_1}b + \nabla_{A_1} a\times b + a\times a \times b + F_{A_1}\times b + \Ric_g\times b + b
\\
+ \nabla_{A_1}^2 R + a\times \nabla_{A_1} R + \nabla_{A_1} a\times R +  a\times a\times R + \dot R.
\end{multline}
Recall that the realization $\cA = \sA_2$ of $\sA = \nabla_{A_1}^*\nabla_{A_1} + 1$ on $\cW = L^2(X; \Lambda^1\otimes\ad P)$ is a positive, sectorial operator, that is, $\cA$ satisfies Hypothesis \ref{hyp:Sell_You_4_standing_hypothesis_A} by the discussion in Section \ref{subsec:Standing_hypothesis_A_and_the_augmented_connection_Laplacian} and that $\calV^2 = H_{A_1}^2(X; \Lambda^1\otimes\ad P)$.

We aim first to apply Lemma \ref{lem:Rade_7-3_abstract_L1_in_time_V2beta_space} to the nonlinear evolution equation \eqref{eq:Augmented_parabolic_equation_for_time_derivative_of_solution_to_Yang-Mills_gradientlike_flow} for $b$, so we check its hypotheses for the nonlinearity, that is
\begin{multline}
\cF(a; b) := a\times\nabla_{A_1}b + \nabla_{A_1} a\times b + a\times a \times b
+ F_{A_1}\times b + \Ric_g\times b + b
\\
+ \nabla_{A_1}^2 R + a\times \nabla_{A_1} R + \nabla_{A_1} a\times R +  a\times a\times R + \dot R.
\end{multline}
so that, in the notation of Lemma \ref{lem:Rade_7-3_abstract_L1_in_time_V2beta_space},
\begin{subequations}
\label{eq:Rade_7-3_Yang-Mills_gradientlike_abstract_nonlinearity_formulae}
\begin{align}
\label{eq:Rade_7-3_Yang-Mills_gradientlike_abstract_nonlinearity_f0}
f_0 &= \nabla_{A_1}^2 R + a\times \nabla_{A_1} R + \nabla_{A_1} a\times R +  a\times a\times R + \dot R,
\\
\cF_1(b) &= F_{A_1}\times b + \Ric_g\times b + b,
\\
\cF_2(a; b) &= a\times\nabla_{A_1}b + \nabla_{A_1} a\times b + a\times a \times b.
\end{align}
\end{subequations}
We observe that Lemma \ref{lem:Freed_Uhlenbeck_equation_6-34_nonnegative_real} yields a continuous Sobolev multiplication map,
$$
H^{2\beta-1}(X) \times H^{2\beta}(X) \to L^2(X),
$$
provided $2\beta-1 \geq 0$, that is, $\beta \geq 1/2$, and $2\beta < d$, and
$$
(2\beta-1 - d/2) + (2\beta - d/2) \geq 0 - d/2,
$$
and thus, $4\beta \geq 1 + d/2$, that is, $\beta \geq 1/4 + d/8$; in other words, noting that $d \geq 2$, we must restrict to $\beta \in [1/4 + d/8, 1)$. Similarly, Lemma \ref{lem:Freed_Uhlenbeck_equation_6-34_nonnegative_real_cubic} yields a continuous Sobolev multiplication map
$$
H^{2\beta}(X) \times H^{2\beta}(X) \times H^{2\beta}(X) \to L^2(X),
$$
provided $\beta \geq 0$ and $2\beta < d$, and
$$
(2\beta - d/2) + (2\beta - d/2) + (2\beta - d/2) \geq 0 - d/2,
$$
and thus, $6\beta \geq d$, that is, $\beta \geq d/6$; in other words, we must restrict to $\beta \in [d/6, 1)$. Consequently, for both quadratic and cubic Sobolev multiplication maps to be continuous, $\beta$ and $d$ must obey
$$
1/4 + d/8 \leq \beta < 1, \quad\hbox{with } 2 \leq d \leq 5,
$$
since $1/4 + d/8 \geq d/6$ for $2 \leq d \leq 6$ and we must already restrict $\beta \in [1/2, 1)$.

With those conditions on $\beta$ and $d$ understood, we have
\begin{align*}
\|\nabla_{A_1}a\times b\|_{L^2(X)} &\leq C\|\nabla_{A_1}a\|_{H_{A_1}^{2\beta-1}(X)} \|b\|_{H_{A_1}^{2\beta}(X)}
\\
&\leq C\|a\|_{H_{A_1}^{2\beta}(X)} \|b\|_{H_{A_1}^{2\beta}(X)},
\end{align*}
and
\begin{align*}
\|a\times \nabla_{A_1}b\|_{L^2(X)} &\leq C\|a\|_{H_{A_1}^{2\beta}(X)} \|\nabla_{A_1}b\|_{H_{A_1}^{2\beta-1}(X)}
\\
&\leq C\|a\|_{H_{A_1}^{2\beta}(X)} \|b\|_{H_{A_1}^{2\beta}(X)},
\end{align*}
and
$$
\|a\times a\times b\|_{L^2(X)} \leq C\|a\|_{H_{A_1}^{2\beta}(X)}^2 \|b\|_{H_{A_1}^{2\beta}(X)}.
$$
By definition of $\cF_1(b)$ we have,
$$
\|\cF_1(b)\|_{L^2(X)} \leq K\|b\|_{L^2(X)}, \quad\hbox{for } b \in L^2(X; \Lambda^1\otimes\ad P),
$$
where
$$
K := c\left(1 + \|F_{A_1}\|_{C(X)} + \|\Ric_g\|_{C(X)}\right),
$$
and similarly, by definition of $\cF_2(a; b)$, we see that
$$
\|\cF_2(a; b)\|_{L^2(X)} \leq c\|a\|_{H_{A_1}^{2\beta}(X)}\left(1 + \|a\|_{H_{A_1}^{2\beta}(X)}\right)\|b\|_{H_{A_1}^{2\beta}(X)},
\quad\hbox{for } b \in H_{A_1}^{2\beta}(X; \Lambda^1\otimes\ad P).
$$
Consequently, if we choose $\eps_1 > 0$ small enough to ensure that $c\eps_1(1 + \eps_1) \leq \eps$, where $\eps > 0$ is the constant in the hypotheses of Lemma \ref{lem:Rade_7-3_abstract_L1_in_time_V2beta_space}, the \apriori estimate,
\begin{equation}
\label{eq:Rade_7-3_L1_in_time_H2beta_in_space_apriori_estimate_by_L1_in_time_L2_in_space_global}
\int_S^T \|b(t)\|_{H_{A_1}^{2\beta}(X)} \,dt
\leq C \left(\|b(S)\|_{L^2(X)} + \int_S^T \|f_0(t)\|_{L^2(X)}\,dt + \int_S^T \|b(t)\|_{L^2(X)} \,dt\right).
\end{equation}
follows from the abstract version \eqref{eq:Rade_7-3_abstract_L1_in_time_V2beta_space_apriori_estimate} with $\cW = L^2(X; \Lambda^1\otimes\ad P)$ and
$\calV^2 = H^2_{A_1}(X; \Lambda^1\otimes\ad P)$. (It is clear from the proof of Lemma \ref{lem:Rade_7-3_abstract_L1_in_time_V2beta_space} that the regularity requirement on the mild solution $u$ in its hypotheses can be relaxed.)

In order to remove the dependency on the initial data, $b(S)$, and convert \eqref{eq:Rade_7-3_L1_in_time_H2beta_in_space_apriori_estimate_by_L1_in_time_L2_in_space_global} to an \apriori interior estimate, we proceed similarly to the proof of Lemma \ref{lem:Rade_7-3}, and introduce a cut-off function, $\eta \in C^\infty(\RR)$ with $\eta = 0$ on $(-\infty, \delta/2]$ and $\eta = 1$ on $[\delta, \infty)$ and $0 \leq \eta \leq 1$ on $\RR$ and $|\dot\eta| \leq 4/\delta$ on $\RR$. We now consider the augmented version of the nonlinear evolution equation
\eqref{eq:Parabolic_equation_for_cut-off_time_derivative_of_solution_to_Yang-Mills_gradientlike_flow} on $(S,T)\times X$, namely,
\begin{multline}
\label{eq:Augmented_parabolic_equation_for_cut-off_time_derivative_of_solution_to_Yang-Mills_gradientlike_flow}
\frac{\partial (\eta b)}{\partial t} + \left(\nabla_{A_1}^*\nabla_{A_1} + 1\right)(\eta b)
\\
=
a\times\nabla_{A_1}(\eta b) + \nabla_{A_1} a\times \eta b + a\times a \times \eta b
\\
+ F_{A_1}\times \eta b + \Ric_g\times \eta b + \eta b + \dot\eta b
\\
+ \nabla_{A_1}^2 \eta R
+ a\times \nabla_{A_1} \eta R + \nabla_{A_1} a\times \eta R +  a\times a\times \eta R + \eta \dot R.
\end{multline}
Proceeding exactly as in the case of the nonlinear evolution equation \eqref{eq:Augmented_parabolic_equation_for_time_derivative_of_solution_to_Yang-Mills_gradientlike_flow} for $b$, we have
\begin{align*}
\|\cF_1(\eta b) + \dot\eta b\|_{L^2(X)} &\leq \left(K + \delta^{-1}\right)\|b\|_{L^2(X)},
\\
\|\cF_2(a; \eta b)\|_{L^2(X)} &\leq \eps \|\eta b\|_{H_{A_1}^{2\beta}(X)}.
\end{align*}
Therefore,
\begin{align*}
\int_{S+\delta}^T \|b(t)\|_{H_{A_1}^{2\beta}(X)}\,dt &\leq \int_S^T \|\eta b(t)\|_{H_{A_1}^{2\beta}(X)}\,dt
\\
&\leq C(1 + \delta^{-1})\int_S^T \|b(t)\|_{L^2(X)} \,dt + C\int_S^T \|\eta f_0(t)\|_{L^2(X)}\,dt
\quad\hbox{(by Lemma \ref{lem:Rade_7-3_abstract_L1_in_time_V2beta_space})}.
\end{align*}
It remains to bound $\|f_0(t)\|_{L^2(X)}$ in terms of norms of $R$. From the expression \eqref{eq:Rade_7-3_Yang-Mills_gradientlike_abstract_nonlinearity_f0} for $f_0$, our hypothesis \eqref{eq:Linfinity_in_time_H2beta_in_space_small_norm_At_minus_A_1_condition_lemma_7.3} that $\|a(t)\|_{H^{2\beta}_{A_1}(X)} \leq \eps_1$ for all $t \in (S, T)$, with $\eps_1 \in (0, 1]$ and $\beta \in [1/4 + d/8, 1)$ and so $2\beta \geq 1$, and the Kato Inequality \eqref{eq:FU_6-20_first-order_Kato_inequality}, we have
\begin{align*}
\|f_0(t)\|_{L^2(X)}
&\leq
C\left(\|\nabla_{A_1}^2 R(t)\|_{L^2(X)}
+ \|a(t)\|_{L^4(X)} \|\nabla_{A_1} R(t)\|_{L^4(X)} \right.
\\
&\quad + \left. \|\nabla_{A_1} a(t)\|_{L^2(X)} \|R(t)\|_{L^\infty(X)}
+ \|a(t)\|_{L^4(X)}^2 \|R(t)\|_{L^\infty(X)} + \|\dot R(t)\|_{L^2(X)} \right)
\\
&\leq C\left(\|R(t)\|_{H^2_{A_1}(X)} + \|R(t)\|_{L^\infty(X)} + \|\dot R(t)\|_{L^2(X)} \right),
\end{align*}
where $C$ is a positive constant depending at most on $A_1$ and $g$. Combining the preceding inequalities yields \eqref{eq:Rade_7-3_L1_in_time_H2beta_in_space_apriori_estimate_by_L1_in_time_L2_in_space}, recalling that $b = \dot A$ and $\eta(t) = 0$ for all $t\in [S, S+\delta/2]$.
\end{proof}


\subsection{Higher-order \apriori estimates for lengths of Yang-Mills gradient flow lines over manifolds of arbitrary dimension}
\label{subsec:Rade_Lemma_7-3_higher-order_arbitrary_dimension}
In Section \ref{subsec:Rade_Lemma_7-3_generalization_dimension_2_leq_d_leq_4}, we established Lemma \ref{lem:Rade_7-3}, a generalization of R\r{a}de's key \apriori $L^1$-in-time interior estimate \cite[Lemma 7.3]{Rade_1992} for a solution to Yang-Mills gradient flow \cite[Lemma 7.3]{Rade_1992} from the case of manifolds of dimension $d=2$ or $3$ to allow $d=4$ and Yang-Mills gradient-like rather than pure gradient flow. Section \ref{subsec:Rade_Lemma_7-3_generalization_H2beta_dimension_2_leq_d_leq_5} contained a further generalization, Lemma \ref{lem:Rade_7-3_L1_in_time_H2beta_in_space_apriori_estimate_by_L1_in_time_L2_in_space}, now allowing $d=5$ and replacing the Sobolev space $H_{A_1}^1(X;\Lambda^1\otimes\ad P)$ by $H_{A_1}^{2\beta}(X;\Lambda^1\otimes\ad P)$, with $\beta \in [1/4+d/8,1)$ (the implicit constraint on $d$ here allows $2 \leq d \leq 5$ but not $d \geq 6$). However, to remove the constraint $d \leq 5$, we need to replace $H_{A_1}^1(X;\Lambda^1\otimes\ad P)$ by  $W_{A_1}^{1,p}(X;\Lambda^1\otimes\ad P)$ for suitably large $p$ depending on $d \geq 2$.

The proof of Lemma \ref{lem:Rade_7-3_L1_in_time_H2beta_in_space_apriori_estimate_by_L1_in_time_L2_in_space} relied the abstract Lemma \ref{lem:Rade_7-3_abstract_L1_in_time_V2beta_space} and we now apply Lemma \ref{lem:Rade_7-3_abstract_L1_in_time_V2beta_space}, together with finite Moser iteration and the Sobolev Embedding Theorem, to obtain \apriori estimates for lengths of Yang-Mills gradient flow lines over manifolds of arbitrary dimension $d \geq 2$. When $d=4$, Lemma \ref{lem:Rade_7-3} provides a stronger conclusion than those in Lemma \ref{lem:Rade_7-3_arbitrary_dimension} or Corollary \ref{cor:Rade_7-3_arbitrary_dimension_L1_time_W1p_space} because in the latter statements we cannot choose $p=2$ (as we do in Lemma \ref{lem:Rade_7-3}) but must restrict to $p > 2$. See Remark \ref{rmk:Rade_7-3} for a discussion of the two roles played by the connection, $A_1$, in the statement of Lemma \ref{lem:Rade_7-3_L1_in_time_H2beta_in_space_apriori_estimate_by_L1_in_time_L2_in_space}.

\begin{lem}[\Apriori $L^1$-in-time-$W^{2\beta,r}$-in-space interior estimate for a solution to Yang-Mills gradient flow over base manifolds of arbitrary dimension]
\label{lem:Rade_7-3_arbitrary_dimension}
Let $G$ be a compact Lie group and $P$ a principal $G$-bundle over a closed, connected, orientable,
smooth manifold, $X$, of dimension $d \geq 2$ and Riemannian metric, $g$. Let $A_1$ be a reference connection of class $C^\infty$ on $P$, and $p \in (d/2,\infty)$, and $r \in (1,\infty)$, and $\beta \in (1/2, 1)$ obey one of the following conditions:
\begin{subequations}
\label{eq:Rade_7-3_arbitrary_dimension_r_beta_conditions}
\begin{gather}
\label{eq:Rade_7-3_arbitrary_dimension_r_beta_supercritical}
r \leq p \text{ and } (2\beta-1)r\geq d, \quad\text{or}
\\
\label{eq:Rade_7-3_arbitrary_dimension_r_beta_subcritical}
r < p \text{ and } \beta \geq d/4p + 1/2 \text{ and } (2\beta-1)r<d.
\end{gather}
\end{subequations}
Then there are positive constants, $C = C(A_1, d, g, p, r, \beta) \in [1,\infty)$ and $\eps_1 = \eps_1(A_1, d, g, p, r, \beta) \in (0, 1]$, such that if $A(t)$ is a strong
solution to the Yang-Mills gradient flow \eqref{eq:Yang-Mills_gradient_flow} on $P$ over an interval $(S, T)$ with regularity,
\[
A - A_1 \in L^\infty(S, T; W_{A_1}^{1,p}(X; \Lambda^1\otimes\ad P)) \cap W^{1,1}_{\loc}(S, T; W_{A_1}^{2\beta,r}(X; \Lambda^1\otimes\ad P)),
\]
where $S \in \RR$ and $\delta > 0$ and $T$ obey $S + 2\delta \leq T \leq \infty$, and
\begin{equation}
\label{eq:Linfinity_in_time_W1p_in_space_small_norm_At_minus_A_1_condition_lemma_7-3}
\|A(t) - A_1\|_{W_{A_1}^{1,p}(X)} \leq \eps_1, \quad\text{a.e. } t \in (S, T),
\end{equation}
then
\begin{equation}
\label{eq:Rade_apriori_interior_estimate_lemma_7-3_arbitrary_dimension_W2beta_r_Lr}
\int_{S+\delta}^T \|\dot A(t)\|_{W_{A_1}^{2\beta,r}(X)}\,dt
\leq
C\left(1 + \delta^{-1}\right)\int_S^T \|\dot A(t)\|_{L^r(X)}\,dt.
\end{equation}
If in addition,
\begin{subequations}
\label{eq:Rade_7-3_arbitrary_dimension_p_greaterthan_2_and_2beta p_geq_d}
\begin{align}
\label{eq:Rade_7-3_arbitrary_dimension_p_greaterthan_2}
p &> 2,
\\
\label{eq:Rade_7-3_arbitrary_dimension_2beta p_geq_d}
2\beta p &\geq d,
\end{align}
\end{subequations}
then there is an integer $n = n(d,r,\beta) \geq 1$ such that, for a possibly larger $C \in [1,\infty)$,
\begin{equation}
\label{eq:Rade_apriori_interior_estimate_lemma_7-3_arbitrary_dimension_W2betar_L2}
\int_{S+\delta}^T \|\dot A(t)\|_{W_{A_1}^{2\beta,r}(X)}\,dt
\leq
C\left(1 + \delta^{-n}\right)\int_S^T \|\dot A(t)\|_{L^2(X)}\,dt,
\end{equation}
and, in particular, for a possibly larger $C \in [1,\infty)$,
\begin{equation}
\label{eq:Rade_apriori_interior_estimate_lemma_7-3_arbitrary_dimension_W1p_L2}
\int_{S+\delta}^T \|\dot A(t)\|_{W_{A_1}^{1,p}(X)}\,dt
\leq
C\left(1 + \delta^{-n}\right)\int_S^T \|\dot A(t)\|_{L^2(X)}\,dt.
\end{equation}
\end{lem}

\begin{proof}
We shall adapt and extend the proof\footnote{While we have assumed that $R\equiv 0$ for simplicity, the argument could easily be extended to allow for non-zero $R$ and gradient-like flow.}
of Lemma \ref{lem:Rade_7-3_L1_in_time_H2beta_in_space_apriori_estimate_by_L1_in_time_L2_in_space}, by first replacing the choices $\calV^2 = H_{A_1}^2(X;\Lambda^1\otimes\ad P)$ and $\cW = L^2(X;\Lambda^1\otimes\ad P)$ (which had limited the admissible dimension of $X$ to $2 \leq d \leq 5$), with $\calV^2 = W_{A_1}^{2,r}(X;\Lambda^1\otimes\ad P)$ and $\cW = L^r(X;\Lambda^1\otimes\ad P)$ to give \eqref{eq:Rade_apriori_interior_estimate_lemma_7-3_arbitrary_dimension_W2beta_r_Lr}, then applying finite Moser iteration to obtain \eqref{eq:Rade_apriori_interior_estimate_lemma_7-3_arbitrary_dimension_W2betar_L2}, and finally appealing to the Sobolev Embedding to deduce \eqref{eq:Rade_apriori_interior_estimate_lemma_7-3_arbitrary_dimension_W1p_L2}. Observe that $\calV^{2\beta} = W_{A_1}^{2\beta,r}(X;\Lambda^1\otimes\ad P)$.

We first prove \eqref{eq:Rade_apriori_interior_estimate_lemma_7-3_arbitrary_dimension_W2beta_r_Lr}. Write $a(t) = A(t) - A_1$ and $b(t) = \dot A(t)$ as in the proof of Lemma \ref{lem:Rade_7-3_L1_in_time_H2beta_in_space_apriori_estimate_by_L1_in_time_L2_in_space} and observe that, to extend that proof, we shall need to verify the bounds on nonlinearities asserted in the following

\begin{claim}[$L^r$ bounds on nonlinearities]
\label{claim:Rade_7-3_arbitrary_dimension_Lr_bounds_nonlinearities}
Given the hypotheses of Lemma \ref{lem:Rade_7-3_arbitrary_dimension}, excluding \eqref{eq:Rade_7-3_arbitrary_dimension_2beta p_geq_d}, there is a constant $C = C(A_1,g,G,r,\beta) \in [1,\infty)$ such that
\begin{subequations}
\label{eq:Rade_7-3_arbitrary_dimension_Lr_bounds_nonlinearities}
\begin{align}
\label{eq:Rade_7-3_arbitrary_dimension_Lr_bound_a_times_nabla_b}
\|a\times\nabla_{A_1}b\|_{L^r(X)} &\leq C\eps_1\|b\|_{W_{A_1}^{2\beta,r}(X)}
\\
\label{eq:Rade_7-3_arbitrary_dimension_Lr_bound_nabla_a_times_b}
\|\nabla_{A_1}a\times b\|_{L^r(X)} &\leq C\eps_1\|b\|_{W_{A_1}^{2\beta,r}(X)},
\\
\label{eq:Rade_7-3_arbitrary_dimension_Lr_bound_a_times_a_times_a}
\|a\times a\times b\|_{L^r(X)} &\leq C\eps_1\|b\|_{W_{A_1}^{2\beta,r}(X)}.
\end{align}
\end{subequations}
\end{claim}

\begin{proof}[Proof of Claim \ref{claim:Rade_7-3_arbitrary_dimension_Lr_bounds_nonlinearities}]
For $1 < r \leq p < \infty$, we define $q \in (r,2p]$ by $1/r = 1/2p + 1/q$, so that by \cite[Theorem 4.12]{AdamsFournier} we have
\begin{equation}
\label{eq:Sobolev_embedding_W2beta-1_r_into_Lu}
W^{2\beta-1,r}(X)
\subset
\begin{cases}
L^{r^*}(X), &\text{if } 0 \leq (2\beta-1)r < d \text{ and } r^* = dr/(d - (2\beta-1)r),
\\
L^u(X), &\text{if } (2\beta-1)r = d \text{ and } 1 \leq u < \infty,
\\
L^\infty(X), &\text{if } (2\beta-1)r > d.
\end{cases}
\end{equation}
It will be convenient to separately consider the two cases in \eqref{eq:Rade_7-3_arbitrary_dimension_r_beta_conditions}.

\setcounter{case}{0}
\begin{case}[$r$ and $\beta$ obey \eqref{eq:Rade_7-3_arbitrary_dimension_r_beta_supercritical}]
We apply the H\"older inequalities defined by $1/r = 1/2p + 1/q$ and $1/r = 1/p + 1/s$ determining $s\in (r,\infty]$, so that
\begin{align*}
\|a\times\nabla_{A_1}b\|_{L^r(X)} &\leq c\|a\|_{L^{2p}(X)} \|\nabla_{A_1}b\|_{L^q(X)} \leq C\|a\|_{W_{A_1}^{1,p}(X)} \|\nabla_{A_1}b\|_{W_{A_1}^{2\beta-1,r}(X)},
\\
\|\nabla_{A_1}a\times b\|_{L^r(X)} &\leq c\|\nabla_{A_1}a\|_{L^p(X)} \|b\|_{L^s(X)} \leq C\|\nabla_{A_1}a\|_{L^p(X)} \|b\|_{W_{A_1}^{2\beta,r}(X)},
\\
\|a\times a\times b\|_{L^r(X)} &\leq c\|a\|_{L^{2p}(X)}^2 \|b\|_{L^s(X)} \leq C\|a\|_{W_{A_1}^{1,p}(X)}^2 \|b\|_{W_{A_1}^{2\beta,r}(X)},
\end{align*}
where we appeal to \cite[Theorem 4.12]{AdamsFournier} for continuity of the Sobolev embeddings, $W^{1,p}(X) \subset L^{2p}(X)$ (for any $p \geq d/2$) and $W^{2\beta,r}(X) \subset L^\infty(X)$ (using the fact that because $(2\beta-1)r \geq d$, we clearly have $2\beta r > d$) and appeal to \eqref{eq:Sobolev_embedding_W2beta-1_r_into_Lu} for continuity of the Sobolev embedding, $W^{2\beta-1,r}(X) \subset L^q(X)$ (noting that $(2\beta-1)r \geq d$ by hypothesis \eqref{eq:Rade_7-3_arbitrary_dimension_r_beta_supercritical}).
\end{case}

\begin{case}[$r$ and $\beta$ obey \eqref{eq:Rade_7-3_arbitrary_dimension_r_beta_subcritical}]
Considering the first nonlinear term, $a\times\nabla_{A_1}b$, we observe that $q \leq r^* = dr/(d - (2\beta-1)r)$ if and only if
\[
\frac{1}{q} = \frac{1}{r} - \frac{1}{2p} \geq \frac{1}{r^*} = \frac{1}{r} - \frac{(2\beta-1)}{d},
\]
namely, $(2\beta-1)/d \geq 1/2p$ or $2\beta \geq d/2p + 1$ or $\beta \geq d/4p + 1/2$ (as assured by hypothesis \eqref{eq:Rade_7-3_arbitrary_dimension_r_beta_subcritical}). When $p > d/2$, as we assume in the hypotheses of Lemma \ref{lem:Rade_7-3_arbitrary_dimension}, then $d/4p < 1/2$ and we may choose $\beta \in [d/4p + 1/2, 1)$. Therefore, we have
\begin{align*}
\|a\times\nabla_{A_1}b\|_{L^r(X)} &\leq c\|a\|_{L^{2p}(X)}\|\nabla_{A_1}b\|_{L^q(X)}
\\
&\leq  C\|a\|_{W_{A_1}^{1,p}(X)}\|\nabla_{A_1}b\|_{W_{A_1}^{2\beta-1,r}(X)}
\quad\text{(by \eqref{eq:Sobolev_embedding_W2beta-1_r_into_Lu})}
\\
&\leq  C\|a\|_{W_{A_1}^{1,p}(X)}\|b\|_{W_{A_1}^{2\beta,r}(X)}.
\end{align*}
For the second nonlinear term, $\nabla_{A_1}a\times b$, we use $1/r = 1/p + 1/s$ and $W^{2\beta,r}(X) \subset L^s(X)$ when $s \leq r^* = dr/(d - 2\beta r)$ by \cite[Theorem 4.12]{AdamsFournier}, that is,
\[
\frac{1}{s} = \frac{1}{r} - \frac{1}{p} \geq \frac{1}{r^*} = \frac{1}{r} - \frac{2\beta}{d},
\]
namely, $2\beta/d \geq 1/p$ or $\beta \geq d/2p$. When $p > d/2$, as we suppose, then $d/2p < 1$ and we may choose $\beta \in [d/2p, 1)$. Therefore, we have
\begin{align*}
\|\nabla_{A_1}a\times b\|_{L^r(X)} &\leq c\|\nabla_{A_1}a\|_{L^p(X)}\|b\|_{L^s(X)}
\\
&\leq  C\|a\|_{W_{A_1}^{1,p}(X)}\|b\|_{W_{A_1}^{2\beta,r}(X)},
\end{align*}
as desired. Note that $d/2p = d/4p + d/4p < d/4p + 1/2$ and so a choice of $\beta \in [d/4p + 1/2, 1)$ is valid for both the first and second nonlinear terms.

For the third nonlinear term, $a\times a\times b$, using $1/r = 1/p + 1/s = 1/2p + 1/2p + 1/s$, for $s$ and $\beta$ chosen as in our analysis of the second term, we have
\begin{align*}
\|a\times a\times b\|_{L^r(X)} &\leq c\|a\|_{L^{2p}(X)}^2 \|b\|_{L^s(X)}
\\
&\leq  C\|a\|_{W_{A_1}^{1,p}(X)}^2 \|b\|_{W_{A_1}^{2\beta,r}(X)},
\end{align*}
as desired.
\end{case}
This completes the proof of Claim \ref{claim:Rade_7-3_arbitrary_dimension_Lr_bounds_nonlinearities}.
\end{proof}

The inequality \eqref{eq:Rade_apriori_interior_estimate_lemma_7-3_arbitrary_dimension_W2beta_r_Lr} follows from Lemma \ref{lem:Rade_7-3_abstract_L1_in_time_V2beta_space}, by adapting \mutatis the proof of Lemma \ref{lem:Rade_7-3_L1_in_time_H2beta_in_space_apriori_estimate_by_L1_in_time_L2_in_space}, for $r$ obeying either of the two cases in \eqref{eq:Rade_7-3_arbitrary_dimension_r_beta_conditions}.

We now apply a finite Moser iteration argument to deduce the inequality \eqref{eq:Rade_apriori_interior_estimate_lemma_7-3_arbitrary_dimension_W2betar_L2} from \eqref{eq:Rade_apriori_interior_estimate_lemma_7-3_arbitrary_dimension_W2beta_r_Lr}. It is convenient to consider three separate cases,
\begin{inparaenum}[\itshape i\upshape)]
\item $1 < r \leq 2$,
\item $1 < r < \infty$ and $4\beta \geq d$,
\item $1 < r < \infty$ and $4\beta < d$,
\end{inparaenum}
where $r,\beta$ also satisfy one of the conditions in \eqref{eq:Rade_7-3_arbitrary_dimension_r_beta_conditions}.

\setcounter{case}{0}
\begin{case}[$1 < r \leq 2$]
\label{case:1_lessthan_r_lessthanequal_2}
In this situation, the inequality \eqref{eq:Rade_apriori_interior_estimate_lemma_7-3_arbitrary_dimension_W2betar_L2} (with $n=2$) follows immediately from \eqref{eq:Rade_apriori_interior_estimate_lemma_7-3_arbitrary_dimension_W2beta_r_Lr}.
\end{case}

\begin{case}[$2 < r < \infty$ and $4\beta \geq d$]
\label{case:2_lessthan_r_lessthan_infty_and_4beta_geq_d}
Because $\beta < 1$, this case can only occur when $d = 2$ or $3$. Continuity of the Sobolev embedding $W^{2\beta,2}(X) \subset L^r(X)$ from \cite[Theorem 4.12]{AdamsFournier} ensures that
\begin{align*}
\int_{S+\delta}^T \|\dot A(t)\|_{W_{A_1}^{2\beta,r}(X)}\,dt
&\leq
C\left(1 + (\delta/2)^{-1}\right)\int_{S+\delta/2}^T \|\dot A(t)\|_{L^r(X)}\,dt \quad\text{(by \eqref{eq:Rade_apriori_interior_estimate_lemma_7-3_arbitrary_dimension_W2beta_r_Lr})}
\\
&\leq
C\left(1 + (\delta/2)^{-1}\right)\int_{S+\delta/2}^T \|\dot A(t)\|_{W_{A_1}^{2\beta,2}(X)}\,dt \quad\text{(by \cite[Theorem 4.12]{AdamsFournier})}
\\
&\leq
C\left(1 + (\delta/2)^{-1}\right)^2\int_S^T \|\dot A(t)\|_{L^2(X)}\,dt  \quad\text{(by \eqref{eq:Rade_apriori_interior_estimate_lemma_7-3_arbitrary_dimension_W2beta_r_Lr})},
\end{align*}
where in the last inequality we apply \eqref{eq:Rade_apriori_interior_estimate_lemma_7-3_arbitrary_dimension_W2beta_r_Lr} with $r$ replaced by $r_0=2$ and observe that for a given triple $d,p,\beta$ and $p>2$ by \eqref{eq:Rade_7-3_arbitrary_dimension_p_greaterthan_2}, any choice of $r \in (1,p)$ will obey one of the two conditions in \eqref{eq:Rade_7-3_arbitrary_dimension_r_beta_conditions}. Hence, we obtain \eqref{eq:Rade_apriori_interior_estimate_lemma_7-3_arbitrary_dimension_W2betar_L2} with $n=2$ for this case.
\end{case}

\begin{case}[$2 < r < \infty$ and $4\beta < d$]
\label{case:2_lessthan_r_lessthan_infty_and_4beta_lessthan_d}
Recall again that for a given triple $d,p,\beta$ and $p>2$ by \eqref{eq:Rade_7-3_arbitrary_dimension_p_greaterthan_2}, any choice of $r \in (1,p)$ will obey one of the two conditions in \eqref{eq:Rade_7-3_arbitrary_dimension_r_beta_conditions}. Set $r_0=2$ and observe that \eqref{eq:Rade_apriori_interior_estimate_lemma_7-3_arbitrary_dimension_W2betar_L2} holds with $r$ replaced by $r_0$ by Case \ref{case:1_lessthan_r_lessthanequal_2}, so
\[
\int_{S+\delta_0}^T \|\dot A(t)\|_{W_{A_1}^{2\beta,r_0}(X)}\,dt \leq C_0\left(1 + \delta_0^{-1}\right)\int_S^T \|\dot A(t)\|_{L^2(X)}\,dt,
\]
for $\delta_0 = \delta/2$ and $C_0 \in [1,\infty)$. We observe that $2\beta r_0 < d$ by the assumptions for this case and now define a (finite) sequence $\{r_i\} \subset [2, d/2\beta)$ by induction using
\[
r_{i+1} := r_i^* = \frac{dr_i}{d - 2\beta r_i} = \frac{r_i(d - 2\beta r_i) + 2\beta r_i^2}{d - 2\beta r_i} = r_i + \frac{2\beta r_i^2}{d - 2\beta r_i}, \quad i=0,1,2,\ldots
\]
so that, for $\delta_i := \delta/2^{i+1}$ and each $i\geq 0$,
\begin{align*}
\int_{S+\delta_i+\delta_{i+1}}^T \|\dot A(t)\|_{W_{A_1}^{2\beta,r_{i+1}}(X)}\,dt
&\leq
C_{i+1}\left(1 + \delta_{i+1}^{-1}\right)\int_{S+\delta_i}^T \|\dot A(t)\|_{L^{r_{i+1}}(X)}\,dt
\quad\text{(by \eqref{eq:Rade_apriori_interior_estimate_lemma_7-3_arbitrary_dimension_W2beta_r_Lr})}
\\
&\leq \kappa_{i+1} C_{i+1}\left(1 + \delta_{i+1}^{-1}\right)\int_{S+\delta_i}^T \|\dot A(t)\|_{W_{A_1}^{2\beta,r_i}(X)}\,dt,
\end{align*}
where $\kappa_{i+1} \in [1,\infty)$ is the norm of the continuous Sobolev embedding $W^{2\beta,r_i}(X) \subset L^{r_{i+1}}(X)$ provided by \cite[Theorem 4.12]{AdamsFournier}, and thus
\begin{equation}
\label{eq:Interior_L1_in_time_W2beta_in_space_estimate_step_i}
\int_{S+\Delta_i}^T \|\dot A(t)\|_{W_{A_1}^{2\beta,r_i}(X)}\,dt
\leq
K_i\left(1 + \delta^{-1}\right)^{i+1}\int_{S}^T \|\dot A(t)\|_{L^2(X)}\,dt,
\end{equation}
where $\Delta_i := (\delta/2)\sum_{j=0}^i 2^{-j} \in (0,\delta)$ and $K_i \in [1,\infty)$.

But $h(r) := 2\beta r^2/(d-2\beta r)$ is a strictly increasing function of $r \in (1,d/2\beta)$ and $h(r) \nearrow \infty$ as $r \nearrow d/2\beta$ and thus, after finitely many steps, we obtain an index $m$ and a large enough exponent $r_m \in (1,d/2\beta)$ such that $r_{m+1} \in (r_m,\infty)$ obeys $2\beta r_{m+1} \geq d$. Choose $\hat r_{m+1} = r_{m+1}\wedge p$, noting that $2\beta p \geq d$ by hypothesis \eqref{eq:Rade_7-3_arbitrary_dimension_2beta p_geq_d}. By \cite[Theorem 4.12]{AdamsFournier}, we have $W^{2\beta,\hat r_{m+1}}(X) \subset L^q(X)$, for $2\beta \hat r_{m+1} = d$ and $1 \leq q < \infty$, and $W^{2\beta,\hat r_{m+1}}(X) \subset L^\infty(X)$, for $2\beta \hat r_{m+1} > d$. Therefore,
\begin{align*}
\int_{S+\delta}^T \|\dot A(t)\|_{W_{A_1}^{2\beta,r}(X)}\,dt
&\leq
C\left(1 + (\delta/2)^{-1}\right)\int_{S+\delta/2}^T \|\dot A(t)\|_{L^r(X)}\,dt \quad\text{(by \eqref{eq:Rade_apriori_interior_estimate_lemma_7-3_arbitrary_dimension_W2beta_r_Lr})}
\\
&\leq
C\left(1 + (\delta/2)^{-1}\right)\int_{S+\delta/2}^T \|\dot A(t)\|_{W_{A_1}^{2\beta,\hat r_{m+1}}(X)}\,dt
\\
&\qquad\text{(as $r<\infty$ by hypothesis and applying \cite[Theorem 4.12]{AdamsFournier})}
\\
&\leq
C\left(1 + (\delta/2)^{-1}\right)\int_{S+\Delta_{m+1}}^T \|\dot A(t)\|_{W_{A_1}^{2\beta,\hat r_{m+1}}(X)}\,dt
\quad\text{(by $\Delta_{m+1} < \delta/2$)}
\\
&\leq
C\left(1 + \delta^{-1}\right)^{m+2}\int_S^T \|\dot A(t)\|_{L^2(X)}\,dt
\\
&\qquad\text{(because $\hat r_{m+1}$ obeys \eqref{eq:Rade_7-3_arbitrary_dimension_r_beta_supercritical} and applying \eqref{eq:Interior_L1_in_time_W2beta_in_space_estimate_step_i})}
\end{align*}
which yields \eqref{eq:Rade_apriori_interior_estimate_lemma_7-3_arbitrary_dimension_W2betar_L2} with $n=m+2$.
\end{case}

We may choose $r = p$ in \eqref{eq:Rade_apriori_interior_estimate_lemma_7-3_arbitrary_dimension_W2betar_L2} and note that because $\beta \geq 1/2$, we necessarily have a continuous Sobolev embedding, $W^{2\beta,r}(X) \subset W^{1,p}(X)$, and this yields \eqref{eq:Rade_apriori_interior_estimate_lemma_7-3_arbitrary_dimension_W1p_L2}. This completes the proof of Lemma \ref{lem:Rade_7-3_arbitrary_dimension}.
\end{proof}

It will be convenient to isolate the following special case of Lemma \ref{lem:Rade_7-3_arbitrary_dimension} and one that we shall most often use.

\begin{cor}[\Apriori $L^1$-in-time-$W^{1,p}$-in-space interior estimate for a solution to Yang-Mills gradient flow over base manifolds of arbitrary dimension]
\label{cor:Rade_7-3_arbitrary_dimension_L1_time_W1p_space}
Let $G$ be a compact Lie group and $P$ a principal $G$-bundle over a closed, connected, orientable,
smooth manifold, $X$, of dimension $d \geq 2$ and Riemannian metric, $g$. Let $A_1$ be a reference connection of class $C^\infty$ on $P$, and $p \in (d/2,\infty)$ obey $p > 2$. Then there are positive constants, $C = C(A_1, d, g, p) \in [1,\infty)$ and $\eps_1 = \eps_1(A_1, d, g, p) \in (0, 1]$, such that if $A(t)$ is a strong
solution to the Yang-Mills gradient flow \eqref{eq:Yang-Mills_gradient_flow} on $P$ over an interval $(S, T)$ with regularity,
\[
A - A_1 \in L^\infty(S, T; W_{A_1}^{1,p}(X; \Lambda^1\otimes\ad P)) \cap W^{1,1}_{\loc}(S, T; W_{A_1}^{2,p}(X; \Lambda^1\otimes\ad P)),
\]
where $S \in \RR$ and $\delta > 0$ and $T$ obey $S + 2\delta \leq T \leq \infty$, and
\begin{equation}
\label{eq:Linfinity_in_time_W1p_in_space_small_norm_At_minus_A_1_condition_lemma_7-3_corollary}
\|A(t) - A_1\|_{W_{A_1}^{1,p}(X)} \leq \eps_1, \quad\text{a.e. } t \in (S, T),
\end{equation}
then there is an integer $n = n(d,p) \geq 1$ such that
\begin{equation}
\label{eq:Rade_apriori_interior_estimate_lemma_7-3_arbitrary_dimension_W1p_L2_corollary}
\int_{S+\delta}^T \|\dot A(t)\|_{W_{A_1}^{1,p}(X)}\,dt
\leq
C\left(1 + \delta^{-n}\right)\int_S^T \|\dot A(t)\|_{L^2(X)}\,dt.
\end{equation}
\end{cor}

We next establish two simple extensions of Corollary \ref{cor:Rade_7-3_arbitrary_dimension_L1_time_W1p_space} from the case of an $L^1$-in-time-$W^{1,p}$-in-space \apriori estimate to those of $L^1$-in-time-$W^{k,q}$-in-space and $L^1$-in-time-$C^{l,\alpha}$-in-space \apriori interior estimates for arbitrary $q \in (1,\infty)$ and $\alpha \in (0,1)$ and integers $k\geq 1$ and $l \geq 0$.

\begin{cor}[\Apriori $L^1$-in-time-$W^{k,q}$-in-space interior estimate for a solution to Yang-Mills gradient flow over base manifolds of arbitrary dimension]
\label{cor:Rade_7-3_arbitrary_dimension_L1_time_Wkp_space}
Let $G$ be a compact Lie group and $P$ a principal $G$-bundle over a closed, connected, orientable,
smooth manifold, $X$, of dimension $d \geq 2$ and Riemannian metric, $g$. Let $A_1$ be a reference connection of class $C^\infty$ on $P$, and $k \geq 1$ an integer, and $q \in (1,\infty)$ such that $kq > d$. Then there are positive constants, $C = C(A_1, d, g, k, q) \in [1,\infty)$ and $\eps_1 = \eps_1(A_1, d, g, k, q) \in (0, 1]$, such that if $A(t)$ is a strong
solution to the Yang-Mills gradient flow \eqref{eq:Yang-Mills_gradient_flow} on $P$ over an interval $(S, T)$ with regularity,
\[
A - A_1 \in L^\infty(S, T; W_{A_1}^{k,q}(X; \Lambda^1\otimes\ad P)) \cap W^{1,1}_{\loc}(S, T; W_{A_1}^{k+2,q}(X; \Lambda^1\otimes\ad P)),
\]
where $S \in \RR$ and $\delta > 0$ and $T$ obey $S + 2\delta \leq T \leq \infty$, and
\begin{equation}
\label{eq:Linfinity_in_time_Wkq_in_space_small_norm_At_minus_A_1_condition_lemma_7-3_corollary}
\|A(t) - A_1\|_{W_{A_1}^{k,q}(X)} \leq \eps_1, \quad\text{a.e. } t \in (S, T),
\end{equation}
then there is an integer $n = n(d,k,q) \geq 1$ such that
\begin{equation}
\label{eq:Rade_apriori_interior_estimate_lemma_7-3_arbitrary_dimension_Wkq_L2_corollary}
\int_{S+\delta}^T \|\dot A(t)\|_{W_{A_1}^{k,q}(X)}\,dt
\leq
C\left(1 + \delta^{-n}\right)\int_S^T \|\dot A(t)\|_{L^2(X)}\,dt.
\end{equation}
\end{cor}

\begin{proof}
We proceed exactly as in the proof of Lemma \ref{lem:Rade_7-3_L1_in_time_H2beta_in_space_apriori_estimate_by_L1_in_time_L2_in_space}, but we now apply Lemma \ref{lem:Rade_7-3_abstract_L1_in_time_V2beta_space} with $\cW = W_{A_1}^{k-1,q}(X; \Lambda^1\otimes\ad P)$ and $\calV^2 = W_{A_1}^{k+1,q}(X; \Lambda^1\otimes\ad P)$ and exploit the fact that because $W^{k,q}(X)$ is a Banach algebra by our hypothesis on $(k,q)$, all calculations involving nonlinearities in the proof of Lemma \ref{lem:Rade_7-3_L1_in_time_H2beta_in_space_apriori_estimate_by_L1_in_time_L2_in_space} simplify considerably and yield, for $\beta \in [0,1)$,
\[
\int_{S+\delta}^T \|\dot A(t)\|_{W_{A_1}^{k-1+2\beta,q}(X)}\,dt
\leq
C\left(1 + \delta^{-n}\right)\int_S^T \|\dot A(t)\|_{W_{A_1}^{k-1,q}(X)}\,dt.
\]
Thus, in particular we obtain (taking $\beta = 1/2$),
\[
\int_{S+\delta}^T \|\dot A(t)\|_{W_{A_1}^{k,q}(X)}\,dt
\leq
C\left(1 + \delta^{-n}\right)\int_S^T \|\dot A(t)\|_{W_{A_1}^{k-1,q}(X)}\,dt.
\]
We now take $\cW = W_{A_1}^{k-j,q}(X; \Lambda^1\otimes\ad P)$ and $\calV^2 = W_{A_1}^{k+2-j,q}(X; \Lambda^1\otimes\ad P)$ for each $j$ obeying $1 \leq j \leq k$ to give (for possibly larger $C$ and $n$)
\[
\int_{S+\delta}^T \|\dot A(t)\|_{W_{A_1}^{k+1-j,q}(X)}\,dt
\leq
C\left(1 + \delta^{-n}\right)\int_S^T \|\dot A(t)\|_{W_{A_1}^{k-j,q}(X)}\,dt, \quad\text{for } 1 \leq j \leq k.
\]
Combining the estimates for $1 \leq j \leq k$ by replacing  $(S+\delta, T)$ by $(S+j\delta/k, T)$ on the left-hand side and replacing $(S, T)$ by $(S+(j-1)\delta/k, T)$ on the right-hand side yields (for possibly larger $C$ and $n$)
\[
\int_{S+\delta}^T \|\dot A(t)\|_{W_{A_1}^{k,q}(X)}\,dt
\leq
C\left(1 + \delta^{-n}\right)\int_S^T \|\dot A(t)\|_{L^q(X)}\,dt.
\]
Combining the preceding inequality with Corollary \ref{cor:Rade_7-3_arbitrary_dimension_L1_time_W1p_space} and $p \in (2\vee d/2, \infty)$ chosen large enough that $W^{1,p}(X) \subset L^q(X)$ yields the desired conclusion.
\end{proof}

By choosing a large enough integer $k = k(d,l,p,\alpha) \geq 1$ to obtain a continuous Sobolev embedding $W^{k,p}(X) \subset C^{l,\alpha}$ assured by \cite[Theorem 4.12]{AdamsFournier}, Corollary \ref{cor:Rade_7-3_arbitrary_dimension_L1_time_Wkp_space} immediately yields

\begin{cor}[\Apriori $L^1$-in-time-$C^{l,\alpha}$-in-space interior estimate for a solution to Yang-Mills gradient flow over base manifolds of arbitrary dimension]
\label{cor:Rade_7-3_arbitrary_dimension_L1_time_Clalpha_space}
Let $G$ be a compact Lie group and $P$ a principal $G$-bundle over a closed, connected, orientable,
smooth manifold, $X$, of dimension $d \geq 2$ and Riemannian metric, $g$. Let $l \geq 0$ be an integer and $\alpha \in (0,1)$. Then there are positive constants, $C = C(d, g, l, \alpha) \in [1,\infty)$ and $\eps_1 = \eps_1(d, g, l, \alpha) \in (0, 1]$, such that if $A(t)$ is a $C^\infty$ classical
solution to the Yang-Mills gradient flow \eqref{eq:Yang-Mills_gradient_flow} on $P$ over an interval $(S, T)$, where $S \in \RR$ and $\delta > 0$ and $T$ obey $S + 2\delta \leq T \leq \infty$, and
\begin{equation}
\label{eq:Linfinity_in_time_Clalpha_in_space_small_norm_At_minus_A_1_condition_lemma_7-3_corollary}
\|A(t) - A_1\|_{C^{l,\alpha}(X)} \leq \eps_1, \quad\forall\, t \in (S, T),
\end{equation}
then there is an integer $n = n(d,l,\alpha) \geq 1$ such that
\begin{equation}
\label{eq:Rade_apriori_interior_estimate_lemma_7-3_arbitrary_dimension_Clalpha_L2_corollary}
\int_{S+\delta}^T \|\dot A(t)\|_{C^{l,\alpha}(X)}\,dt
\leq
C\left(1 + \delta^{-n}\right)\int_S^T \|\dot A(t)\|_{L^2(X)}\,dt.
\end{equation}
\end{cor}

\section[Application to the Yang-Mills gradient system]{Application to the Yang-Mills gradient system and proofs of main results for Yang-Mills gradient flow near a minimum}
\label{sec:Application_abstract_gradient_system_results_Yang-Mills_energy_functional}
In this section, we shall apply to the Yang-Mills gradient system our results for an abstract gradient system defined by a potential function obeying a {\L}ojasiewicz-Simon gradient inequality and bearing on convergence and convergence rates together with global existence and stability of solutions started near a local minimum. Collectively, the results of this section prove Theorem \ref{mainthm:Yang-Mills_gradient_flow_global_existence_and_convergence_started_near_minimum}.

As we shall see in this section, Theorem \ref{mainthm:Yang-Mills_gradient_flow_global_existence_and_convergence_started_near_minimum} is a consequence of short-time existence results for Yang-Mills gradient flow \eqref{eq:Yang-Mills_gradient_flow}, which are valid for a base manifold, $X$, of any dimension $d\geq 2$, and the results we have developed for an abstract gradient system defined by a potential function obeying a {\L}ojasiewicz-Simon gradient inequality.

\subsection{Short-time well-posedness, \apriori estimate, regularity, and minimal lifetime of a solution to the Yang-Mills heat equation}
\label{subsec:Short-time_well-posedness_regularity_minimal_lifetime_solution_Yang-Mills_heat_equation}
Short-time existence, uniqueness, an \apriori estimate, and regularity of a solution, $A(t)$ for $t\in [0,\tau)$ with minimal lifetime $\tau \in (0,\infty]$, to the Yang-Mills heat equation \eqref{eq:Yang-Mills_heat_equation}, and continuity with respect to the initial data are provided by many different methods in Sections \ref{sec:Local_well-posedness_yang_mills_heat_equation} and \ref{sec:Struwe_3and4}, depending in part on the dimension, $d$, of the manifold, $X$:
\begin{itemize}
\item Theorems \ref{thm:Existence_uniqueness_strong_solution_Yang-Mills_heat_equation_in_W_2beta_p_initial_data_in_W_2beta_p}, \ref{thm:Smoothness_strong_solution_Yang-Mills_heat_equation_initial_data_in_W_2beta_p} and \ref{thm:Smoothness_strong_solution_Yang-Mills_heat_equation_smooth_initial_data}, for any $d \geq 2$;

\item Theorems \ref{thm:Struwe_section_4-3_local_wellposedness_Yang-Mills_heat_equation} and \ref{thm:Struwe_section_4-3_regularity_Yang-Mills_heat_equation}, for $d \leq d \leq 4$;

\item Theorem \ref{thm:Struwe_section_4-3_local_wellposedness_Yang-Mills_heat_equation_hybrid_critical-exponent_parabolic_Sobolev_space} and Corollary \ref{cor:Struwe_section_4-3_higher_order_regularity_Yang-Mills_heat_equation_hybrid_critical-exponent_parabolic_Sobolev_space}, for any $d \geq 2$;

\item Theorem \ref{thm:Struwe_section_4-3_local_wellposedness_Yang-Mills_heat_equation_pure_critical-exponent_parabolic_Sobolev_space} and Corollary \ref{cor:Struwe_section_4-3_higher_order_regularity_Yang-Mills_heat_equation_pure_critical-exponent_parabolic_Sobolev_space}, for any $d \geq 2$;

\item Theorem \ref{thm:Struwe_section_4-3_local_existence_uniqueness_Yang-Mills_heat_equation_initial_data_in_Ld} and Corollary \ref{cor:Struwe_section_4-3_higher_order_regularity_Yang-Mills_heat_equation_initial_data_in_Ld}, for any $d \geq 2$.
\end{itemize}
The preceding results each allow initial connections, $A_0$, of different regularity while others (namely, Theorems \ref{thm:Struwe_section_4-3_local_wellposedness_Yang-Mills_heat_equation} and \ref{thm:Struwe_section_4-3_regularity_Yang-Mills_heat_equation}) restrict the dimension $d$ of $X$ to the range $2\leq d\leq 4$ and so they are not all equivalent. However, when the initial data, $A_0$, is a $C^\infty$ connection on $P$, as assumed in Section \ref{subsec:Yang-Mills_gradient_flow_near_local_minimum}, each of these methods yield a classical solution, $A(t) = A_0 + a(t)$, with the property that
$$
a \in C^\infty([0,\infty)\times X; \Lambda^1\otimes\ad P),
$$
together with solutions with different initial regularities,
$$
a \in C([0,\infty); W_{A_1}^{s,p}(X; \Lambda^1\otimes\ad P) \cap C^\infty((0,\infty)\times X; \Lambda^1\otimes\ad P),
$$
when $A_0$ is a connection of class $W^{s,p}$, for suitable $s\geq 0$ and $p>1$.

\subsection{Short-time existence, regularity, and minimal lifetime of a solution to Yang-Mills gradient flow}
\label{subsec:Short-time_existence_regularity_minimal_lifetime_solution_Yang-Mills gradient flow}
Short-time existence, an \apriori estimate, minimal lifetime, and regularity of a solution, $A(t)$ for $t\in [0,\tau)$ and minimal lifetime $\tau \in (0,\infty]$, to the Yang-Mills gradient flow \eqref{eq:Yang-Mills_gradient_flow} follows from Lemma \ref{lem:Donaldson_DeTurck_trick} (the Donaldson DeTurck trick), in conjunction with any of the results in Section \ref{subsec:Short-time_well-posedness_regularity_minimal_lifetime_solution_Yang-Mills_heat_equation} providing short-time existence and regularity of a solution to the Yang-Mills heat equation for initial data, $A_0$, of class $W^{s+1,p}$ with $sp>d$, for a manifold, $X$, of any dimension $d \geq 2$.

A weaker result, again for $X$ with dimension $d$ in the range $2 \leq d \leq 4$ but when $A_0$ is only of class $H^1$, is provided by Lemma \ref{lem_Struwe_sections_4-2_and_4-4}.

\subsection{{\L}ojasiewicz-Simon gradient inequality for the Yang-Mills energy functional}
Theorems \ref{thm:Rade_proposition_7-2_L2} and \ref{thm:Rade_proposition_7-2} provide {\L}ojasiewicz-Simon gradient inequalities for the Yang-Mills energy functional for $d \geq 2$, with Theorem \ref{thm:Rade_proposition_7-2} applied to the proof of Theorem \ref{mainthm:Yang-Mills_gradient_flow_global_existence_and_convergence_started_near_minimum} for $2 \leq d \leq 4$ and Theorem \ref{thm:Rade_proposition_7-2_L2} for $d \geq 5$.

\subsection{\Apriori estimates for lengths of gradient flow lines}
The crucial Hypothesis \ref{hyp:Abstract_apriori_interior_estimate_trajectory} for abstract gradient flow, required for Theorems \ref{thm:Huang_3-3-6}, \ref{thm:Huang_3-4-8}, \ref{thm:Huang_5-1-1}, and \ref{thm:Huang_5-1-2}, is verified for Yang-Mills gradient flow \eqref{eq:Yang-Mills_gradient_flow} by Lemma \ref{lem:Rade_7-3} for $X$ with dimension $d$ in the range $2\leq d \leq 4$ by Lemma \ref{lem:Rade_7-3_L1_in_time_H2beta_in_space_apriori_estimate_by_L1_in_time_L2_in_space} and by Corollary \ref{cor:Rade_7-3_arbitrary_dimension_L1_time_Wkp_space} for $d \geq 5$.

\subsection{Convergence alternative for a global solution to Yang-Mills gradient flow near a critical point}
Before proceeding to apply Theorem \ref{thm:Huang_3-3-6}, we provide the following characterization of the limit point in the that result. R\r{a}de observes that it follows as in the proof of his \cite[Proposition 7.1]{Rade_1992} that the limit point in his \cite[Proposition 7.4]{Rade_1992} is Yang-Mills. We give a simpler and more general proof here, which neither requires the estimate \cite[Proposition 8.1]{Rade_1992} nor an appeal to Uhlenbeck's compactness theorem \cite{UhlLp}.

\begin{lem}[Criterion for an orbit cluster point to be a critical point]
\label{lem:Limit_is_critical_point}
Let $\sH$ be a real Hilbert space and let $\sX\hookrightarrow \sH$ be a subspace that becomes a Banach space under its own norm (as in the introduction to \cite[Section 3.4]{Huang_2006} and setup for \cite[Equation (3.27a)]{Huang_2006}), so $\sX \hookrightarrow \sH \hookrightarrow \sX'$, where we identify $\sH'\cong \sH$. Suppose that $u:[0,\infty)\to \sX$ is a weak solution to the gradient flow for a $C^1$ potential function, $\sE:\sX\to\RR$, namely
$$
\langle\dot u(t), v\rangle = -\langle \sE'(u(t)), v \rangle, \quad \forall\, v \in \sX, \quad t\geq 0,
$$
where $\langle \cdot , \cdot\rangle:\sX'\times \sX \to \RR$ is the canonical bilinear pairing. If $u(t)$ converges to $\varphi\in \sX$ as $t\to\infty$ in the sense that \cite[Equation (3.9c)]{Huang_2006}
$$
\lim_{t\to\infty}\|u(t) - \varphi\|_\sX = 0
\quad\hbox{and}\quad
\int_0^\infty \|\dot u\|_\sX\,dt < \infty,
$$
then $\varphi$ is a critical point of $\sE$, namely, $\sE'(\varphi) = 0$.
\end{lem}

\begin{proof}
Since $\int_0^\infty\|\dot u(t)\|_\sX\,dt < \infty$, it follows that $\liminf_{t\to\infty} \|\dot u(t)\|_\sX = 0$ and hence there is an unbounded sequence $\{t_n\}_{n=1}^\infty\subset[0,\infty)$ such that $\|\dot u(t_n)\|_\sX \to 0$ as $n\to\infty$, that is, $\dot u(t_n) \to 0$ in $\sX$ as $n\to\infty$, while $u(t_n) \to \varphi$ in $\sX$ as $n\to\infty$ by Lemma \ref{lem:Huang_3-2-3}. But
$$
\langle\dot u(t_n), v\rangle = -\langle \sE'(u(t_n)), v \rangle, \quad \forall\, v \in \sX,
$$
and taking the limit as $n \to \infty$ yields
$$
\lim_{n\to\infty}\langle \sE'(u(t_n)), v \rangle = 0, \quad \forall\, v \in \sX.
$$
But $\sE':\sX\to \sX'$ is continuous and so
$$
\langle \sE'(\varphi), v \rangle = \lim_{n\to\infty}\langle \sE'(u(t_n)), v \rangle = 0, \quad \forall\, v \in \sX.
$$
Thus, $\sE'(\varphi) = 0$, as claimed.
\end{proof}

By virtue of the preceding results, we can now apply Theorem \ref{thm:Huang_3-3-6}, our convergence result for an abstract gradient system obeying a {\L}ojasiewicz-Simon inequality, to Yang-Mills gradient flow \eqref{eq:Yang-Mills_gradient_flow} to give the following analogue of the convergence alternative R\r{a}de's \cite[Proposition 7.4]{Rade_1992} and Simon's \cite[Theorem 2]{Simon_1983}:

\begin{thm}[Convergence alternative for a global solution to Yang-Mills gradient flow near a critical point]
\label{thm:Huang_3-3-6_Yang-Mills}
Let $G$ be a compact Lie group, $P$ a principal $G$-bundle over a closed, Riemannian, smooth manifold, $X$, of dimension $d \geq 2$, and $A_1$ a $C^\infty$ reference connection on $P$, and $A_\ym$ a $C^\infty$ Yang-Mills connection on $P$. If $\sigma>0$ is the {\L}ojasiewicz-Simon constant, then there is a constant $\eps \in (0,\sigma/4)$ with the following significance. If $A(t)$ is a strong solution to the Yang-Mills gradient flow \eqref{eq:Yang-Mills_gradient_flow} over $[0,\infty)\times X$
and there is a $T \geq 0$ such that
\begin{subequations}
\label{eq:Rade_7-2}
\begin{align}
\label{eq:Rade_7-2_2_leq_d_leq_4}
\|A(T) - A_\ym\|_{H_{A_1}^1(X)} &< \eps, \quad\hbox{for } 2 \leq d \leq 4, \text{ or}
\\
\label{eq:Rade_7-2_d_geq_5}
\|A(T) - A_\ym\|_{W_{A_1}^{2,p}(X)} &< \eps, \quad\hbox{for } d\geq 5,
\end{align}
\end{subequations}
where $p \in (d/2,\infty)$ and $\eps \in (0,1]$ depends in addition on $p$ in this case, then either
\begin{enumerate}
\item $\sE(A(t)) < \sE(A_\ym)$ for some $t>T$, or

\item the trajectory $A(t)$ converges in $H_{A_1}^1(X;\Lambda^1\otimes\ad P)$, for $2\leq d\leq 4$, to a limit, $A_\infty$ on $P$, as $t \to \infty$ in the sense that
    $$
    \lim_{t\to\infty}\|A(t)-A_\infty\|_{H_{A_1}^1(X)} =0
    \quad\hbox{and}\quad
    \int_1^\infty \|\dot A\|_{H_{A_1}^1(X)}\,dt < \infty;
    $$
    if $d \geq 5$, then the analogous convergence results hold with the norm on $H_{A_1}^1(X;\Lambda^1\otimes\ad P)$ replaced by that on $W_{A_1}^{2,p}(X;\Lambda^1\otimes\ad P)$. The connection, $A_\infty$, is Yang-Mills and has energy $\sE(A_\infty) = \sE(A_\ym)$. If $A_\ym$ is a cluster point of the orbit, $O(A) = \{A(t):t\geq 0\}$, then $A_\infty = A_\ym$.
\end{enumerate}
\end{thm}

\begin{cor}[Convergence for a global solution to Yang-Mills gradient flow near a local minimum]
\label{cor:Huang_3-3-6_Yang-Mills}
Assume the hypotheses of Theorem \ref{thm:Huang_3-3-6_Yang-Mills}. If $A_\infty$ in Theorem \ref{thm:Huang_3-3-6_Yang-Mills} is a local minimum for the Yang-Mills energy functional, then the second alternative in Theorem \ref{thm:Huang_3-3-6_Yang-Mills} necessarily holds.
\end{cor}

\subsection{Convergence rate for solution to the Yang-Mills gradient system}
Again by virtue of the preceding results, we can apply Theorem \ref{thm:Huang_3-4-8}, our enhancement of Huang's   \cite[Theorem 3.4.8]{Huang_2006}, to give the following analogue of the convergence rate results in R\r{a}de's \cite[Proposition 7.4]{Rade_1992} and Simon's   \cite[Theorem 2]{Simon_1983}:

\begin{thm}[Convergence rate for a global solution to Yang-Mills gradient flow near a critical point]
\label{thm:Huang_3-4-8_Yang-Mills}
Assume the hypotheses of Theorem \ref{thm:Huang_3-3-6_Yang-Mills} and, in addition, that the second alternative holds. If $c$, $\sigma$, and $\theta \in [1/2,1)$ denote the {\L}ojasiewicz-Simon constants for the Yang-Mills energy functional in Theorem \ref{thm:Huang_3-3-6_Yang-Mills} then, for $2\leq d\leq 4$ and all $t \geq T+1$,
\begin{multline}
\label{eq:Rade_Proposition_7-4_convergence_rate}
\|A(t) - A_\infty\|_{H_{A_1}^1(X)}
\\
\leq
\begin{cases}
\displaystyle
\frac{1}{c(1-\theta)}\left(c^2(2\theta-1)(t-T-1) + (\sE_T-\sE_\infty)^{1-2\theta}\right)^{-(1-\theta)/(2\theta-1)},
& 1/2 < \theta < 1,
\\
\displaystyle
\frac{2}{c}\sqrt{\sE_T-\sE_\infty}\exp(-c^2(t-T-1)/2),
&\theta = 1/2,
\end{cases}
\end{multline}
where $\sE_T := \sE(A(T))$ and $\sE_\infty := \sE(A_\infty)$; if $d \geq 5$, then \eqref{eq:Rade_Proposition_7-4_convergence_rate} holds with the norm on $H_{A_1}^1(X;\Lambda^1\otimes\ad P)$ replaced by that on $W_{A_1}^{2,p}(X;\Lambda^1\otimes\ad P)$ with $p \in (d/2,\infty)$.
\end{thm}

\begin{rmk}
Our convergence rate estimate \eqref{eq:Rade_Proposition_7-4_convergence_rate} is similar, though not identical, to R\r{a}de's estimate in his \cite[Proposition 7.4]{Rade_1992}:
$$
\|A(t) - A_\infty\|_{H_{A_1}^1(X)}
\leq
\begin{cases}
C(t-T)^{-(1-\theta)/(2\theta-1)}, &\hbox{if } 1/2 < \theta < 1,
\\
Ce^{-k(t-T)/2}, &\hbox{if } \theta = 1/2,
\end{cases}
$$
for positive constants, $C$ and $k$.
\end{rmk}

\subsection{Global existence of a solution to Yang-Mills gradient flow near a local minimum}
We now proceed to apply Theorem \ref{thm:Huang_5-1-1}, our version of Huang's \cite[Theorem 5.1.1]{Huang_2006}, giving the following analogue of Simon's \cite[Corollary 1]{Simon_1983}.

\begin{thm}[Existence and convergence of a global solution to Yang-Mills gradient flow near a local minimum]
\label{thm:Huang_5-1-1_Yang-Mills}
Let $G$ be a compact Lie group, $P$ a principal $G$-bundle over a closed, Riemannian, smooth manifold, $X$, of dimension $d \geq 2$, and $A_1$ a $C^\infty$ reference connection on $P$, and $A_{\min}$ a $C^\infty$ connection on $P$ that is a local minimum for the Yang-Mills energy functional. If $\sigma>0$ is the {\L}ojasiewicz-Simon constant, then there is a constant $\eps \in (0,\sigma/4)$ with the following significance. If $A_0$ is a connection of class $W^{s+1,q}$ on $P$ such that $s\geq 1$ and $q\geq 2$ obey $sq>d$ and satisfies
\begin{subequations}
\label{eq:Initial_connection_Yang-Mills_gradient_flow_in_LS_quarter_neighborhood}
\begin{align}
\label{eq:Initial_connection_Yang-Mills_gradient_flow_in_LS_quarter__neighborhood_2_leq_d_leq_4}
\|A_0 - A_{\min}\|_{H_{A_1}^1(X)} &< \eps, \quad\hbox{for } 2 \leq d \leq 4, \text{ or}
\\
\label{eq:Initial_connection_Yang-Mills_gradient_flow_in_LS_quarter__neighborhood_d_geq_5}
\|A_0 - A_{\min}\|_{W_{A_1}^{2,p}(X)} &< \eps, \quad\hbox{for } d\geq 5,
\end{align}
\end{subequations}
where $p \in (d/2,\infty)$ and $\eps \in (0,1]$ depends in addition on $p$ in this case, then there exists a global strong solution, $A(t)$, to the Yang-Mills gradient flow \eqref{eq:Yang-Mills_gradient_flow} over $[0,\infty)\times X$
such that, for all $t \in [0,\infty)$,
\begin{subequations}
\label{eq:Yang-Mills_gradient_flow_stays_in_LS_half_neighborhood}
\begin{align}
\label{eq:Yang-Mills_gradient_flow_stays_in_LS_half_neighborhood_2_leq_d_leq_4}
\|A(t) - A_{\min}\|_{H_{A_1}^1(X)} &< \sigma/2, \quad\hbox{for } 2 \leq d \leq 4, \text{ or}
\\
\label{eq:Yang-Mills_gradient_flow_stays_in_LS_half_neighborhood_d_geq_5}
\|A(t) - A_{\min}\|_{W_{A_1}^{2,p}(X)} &< \sigma/2, \quad\hbox{for } d\geq 5,
\end{align}
\end{subequations}
and $A(t)$ converges to a Yang-Mills connection, $A_\infty$ on $P$, as $t\to\infty$ in the sense that, for $2 \leq d \leq 4$,
$$
\lim_{t \to \infty} \|A(t) - A_\infty\|_{H_{A_1}^1(X)} = 0
\quad\hbox{and}\quad
\int_1^\infty\|\dot A(t)\|_{H_{A_1}^1(X)}\,dt < \infty,
$$
and for $d \geq 5$, then the preceding convergence results hold with the norm on $H_{A_1}^1(X;\Lambda^1\otimes\ad P)$ replaced by that on $W_{A_1}^{2,p}(X;\Lambda^1\otimes\ad P)$.
\end{thm}

\subsection{Stability of a solution to Yang-Mills gradient flow}
Furthermore, we can apply Theorem \ref{thm:Huang_5-1-2}, our version of Huang's \cite[Theorem 5.1.2]{Huang_2006}, on convergence to a critical point and stability of a local minimum, to the Yang-Mills energy functional and obtain the

\begin{thm}[Convergence to a Yang-Mills connection and stability of a local minimum]
\label{thm:Huang_5-1-2_Yang-Mills}
Assume the hypotheses of Theorem \ref{thm:Huang_5-1-1_Yang-Mills}. Then, for each connection $A_0$ on $P$ obeying \eqref{eq:Initial_connection_Yang-Mills_gradient_flow_in_LS_quarter_neighborhood}, there exists a global strong solution, $A(t)$, to the Yang-Mills gradient flow \eqref{eq:Yang-Mills_gradient_flow} over $[0,\infty)\times X$ that obeys \eqref{eq:Yang-Mills_gradient_flow_stays_in_LS_half_neighborhood} and converges as $t\to\infty$ to a Yang-Mills connection, $A_\infty$ on $P$, satisfying
\begin{subequations}
\label{eq:Yang-Mills_flow_limit_point_in_LS_neighborhood}
\begin{align}
\label{eq:Yang-Mills_flow_limit_point_in_LS_neighborhood_2_leq_d_leq_4}
\|A_\infty - A_{\min}\|_{H_{A_1}^1(X)} &< \sigma, \quad\hbox{for } 2 \leq d \leq 4, \text{ or}
\\
\label{eq:Yang-Mills_flow_limit_point_in_LS_neighborhood_d_geq_5}
\|A_\infty - A_{\min}\|_{W_{A_1}^{2,p}(X)} &< \sigma, \quad\hbox{for } d\geq 5,
\end{align}
\end{subequations}
where $p \in (d/2,\infty)$. The Yang-Mills connection, $A_\infty$, satisfies $\sE(A_\infty) = \sE(A_{\min})$. As an equilibrium of the Yang-Mills gradient flow \eqref{eq:Yang-Mills_gradient_flow}, the point $A_\infty$ is Lyapunov stable (see Definition \ref{defn:Sell_You_page_32_Lyapunov_and_uniform_asymptotic_stability}). If $A_\infty$ is isolated or a cluster point of the orbit $O(A) = \{A(t):t\geq 0\}$, then $A_\infty$ is uniformly asymptotically stable (see Definition \ref{defn:Sell_You_page_32_Lyapunov_and_uniform_asymptotic_stability}).
\end{thm}

\begin{rmk}[Dynamical systems that are Lyapunov but not asymptotically stable]
\label{rmk:Lyapunov_but_not_asymptotically_stable}
For some examples of dynamical systems that are stable in the sense of Lyapunov but not asymptotically stable, we refer the reader to \cite{Grillakis_Shatah_Strauss_1987, Grillakis_Shatah_Strauss_1990, Weinstein_1987}.
\end{rmk}



\subsection{Uniqueness of solutions to Yang-Mills gradient flow}
Uniqueness of the solution modulo a $C^\infty$ path of $C^\infty$ gauge transformations is provided by Theorem \ref{thm:Kozono_Maeda_Naito_6-1} when $d=4$.

\section[Proofs of R\r{a}de's results on Yang-Mills gradient flow]{Proofs of R\r{a}de's results on global existence and convergence of Yang-Mills gradient flow over base manifolds of dimensions two or three}
\label{sec:Yang-Mills_gradient_flow_global_existence_and_convergence_rade_proofs}
In this section, we provide a proof of Item \eqref{item:Rade_theorem_1_and_2_global_existence} in Theorem \ref{mainthm:Rade_1_and_2} that is far simpler than R\r{a}de's argument in \cite[Sections 4, 5, and 6]{Rade_1992}. R\r{a}de's idea is to observe that if $A(t)$ is \emph{any} $C^\infty$ path of connections, so $A-A_0 \in C^\infty((0,T)\times X;\Lambda^1\otimes\ad P)$, then $F_A$ obeys \cite[Equation (4.1)]{Rade_1992}, that is,
\[
\frac{dF_A}{dt} = d_A\left(\frac{dA}{dt}\right), \quad t > 0,
\]
and thus, if $A(t)$ is a solution to \eqref{eq:Yang-Mills_gradient_flow}, then $F_A$ obeys \cite[Equation (4.2)]{Rade_1992}, that is,
\[
\frac{dF_A}{dt} + d_Ad_A^*F_A = 0, \quad t > 0.
\]
As usual, the Bianchi identity \eqref{eq:Freed_Uhlenbeck_2-16_Bianchi_identity} implies that $F_A \in C^\infty((0,T)\times X;\Lambda^2\otimes\ad P)$ solves the parabolic equation \cite[Equation (4.3)]{Rade_1992}, namely
\[
\frac{dF_A}{dt} + \Delta_AF_A = 0, \quad t > 0,
\]
where $\Delta_A = d_Ad_A^* + d_A^*d_A$ is the Hodge Laplace operator \eqref{eq:Lawson_page_93_Hodge_Laplacian}. R\r{a}de's strategy is to view $(A,F_A)$ as a solution $(A,\Omega) \in C^\infty((0,T)\times X;(\Lambda^1\oplus\Omega^2)\otimes\ad P)$ to the \emph{parabolic} system \cite[Equation (4.4)]{Rade_1992},
\begin{align*}
\frac{dA}{dt} + d_A^*\Omega &= 0,
\\
\displaystyle
\frac{d\Omega}{dt} + \Delta_A\Omega &= 0, \quad t > 0,
\end{align*}
with initial data, $(A(0),\Omega(0)) = (A_0,F_{A_0})$, and prove global existence for this system. However, the argument is lengthy and technical.

\begin{proof}[Proof of Theorem \ref{mainthm:Rade_1_and_2}]
Consider Item \eqref{item:Rade_theorem_1_and_2_global_existence}. We have proved regularity and short-time existence of solutions to the Yang-Mills heat equation \eqref{eq:Yang-Mills_heat_equation} by several methods, for low dimension $d=2,3,4$ and for arbitrary dimension $d \geq 2$, as part of the proof of Theorem \ref{mainthm:Yang-Mills_gradient_flow_global_existence_and_convergence_started_near_minimum}; those results are summarized in Section \ref{subsec:Short-time_well-posedness_regularity_minimal_lifetime_solution_Yang-Mills_heat_equation}. We have also proved regularity and short-time existence of solutions to the Yang-Mills gradient system \eqref{eq:Yang-Mills_gradient_flow} by several methods, for arbitrary dimension $d \geq 2$, again as part of the proof of Theorem \ref{mainthm:Yang-Mills_gradient_flow_global_existence_and_convergence_started_near_minimum}; those methods are summarized in Section \ref{subsec:Short-time_existence_regularity_minimal_lifetime_solution_Yang-Mills gradient flow}.

Hence, there exists $T_1 \in (0,\infty]$ such that $A(t) = A_0 + a(t)$ for $t\in [0,T_1)$, with
$$
a \in C^\infty([0,\infty)\times X; \Lambda^1\otimes\ad P),
$$
is a solution to the Yang-Mills gradient flow \eqref{eq:Yang-Mills_gradient_flow} with initial data, $A(0) = A_0$. If $T_1 = \infty$, we are done, so we suppose that $T_1 < \infty$. Because $A$ is a solution to the Yang-Mills gradient system \eqref{eq:Yang-Mills_gradient_flow}, then Lemma \ref{lem:Kozono_Maeda_Naito_4-1} yields
\[
\sup_{t\in[0,T_1)}\|F_{A(t)}\|_{L^2(X)} = \|F_{A(0)}\|_{L^2(X)}.
\]
We now apply the argument used in the proof of Theorem \ref{thm:Kozono_Maeda_Naito_5-1}, where $X$ has dimension four, to the present situation where $X$ has dimension two or three. Since $(X,g)$ is a closed Riemannian manifold, for every constant $r \in (0,1\wedge\Inj(X,g))$ there exists a finite integer $N = N(g,r)$ such that $X$ is covered by the union of at most $N$ geodesic balls, $B_r(x_i) \subset X$.

For $d=2$ or $3$ and any $p \in [\frac{d}{2}, 2)$, with $q \in [2d/(4-d), \infty)$ defined by $1/p = 1/2 + 1/q$, we have
\begin{align*}
\|F_{A(t)}\|_{L^p(B_r(x_0))}
&\leq
\|F_{A(t)}\|_{L^2(B_r(x_0))} \|1\|_{L^q(B_r(x_0))}
\\
&= cr^{d/q}\|F_{A(t)}\|_{L^2(B_r(x_0))}
\\
&\leq cr^{d/q}\|F_{A(t)}\|_{L^2(X)},
\end{align*}
for $c = c(d,g) \in [1,\infty)$ and all $t \in [0,T_1)$. Given $\eps \in (0,1]$, then for small enough $r = r(d,g,\eps) \in (0, 1\wedge\Inj(X,g))$, we obtain
\begin{align*}
\|F_{A(t)}\|_{L^{d/2}(B_r(x_0))}
&\leq
cr^{d/q}\|F_{A(t)}\|_{L^2(X)}
\\
&\leq cr^{d/q}\|F_{A(0)}\|_{L^2(X)}
\\
&\leq \eps, \quad\forall\,t \in [0,T_1).
\end{align*}
But $X$ is covered by the union of $N = N(g,r)$ geodesic balls, $B_r(x_i) \subset X$, with no exceptional points. Hence, the proof of Theorem \ref{thm:Kozono_Maeda_Naito_5-1} implies that the solution, $A(t)$, to \eqref{eq:Yang-Mills_gradient_flow} extends to $t \in [0,T_2)$ for some $T_2 > T_1$, contradicting the maximality of $T_1$. Hence, $A(t)$ exists for all $t \in [0,\infty)$. This completes the proof of Item \eqref{item:Rade_theorem_1_and_2_global_existence}.

R\r{a}de already gives a concise proof of Item \eqref{item:Rade_theorem_1_and_2_convergence} in \cite[Section 7]{Rade_1992} using the {\L}ojasiewicz-Simon gradient inequality and Uhlenbeck compactness \cite{UhlLp}. In particular, by combining Uhlenbeck compactness in the form of his \cite[Proposition 7.1]{Rade_1992} and Theorem \ref{mainthm:Yang-Mills_gradient_flow_global_with_critical_point_in_orbit_closure}, we obtain Item \eqref{item:Rade_theorem_1_and_2_convergence}. The remaining assertions of Theorem \ref{mainthm:Rade_1_and_2} are provided by Theorem \ref{mainthm:Yang-Mills_gradient_flow_global_existence_and_convergence_started_near_minimum}. This completes the proof of Theorem \ref{mainthm:Rade_1_and_2}.
\end{proof}

\chapter[Yang-Mills gradient flow with arbitrary initial energy]{Global existence and convergence of smooth solutions to Yang-Mills gradient flow with arbitrary initial energy}
\label{chapter:Gradient_flow_arbitrary_initial_energy}

\section{Convergence and stability for abstract gradient-like systems}
\label{sec:Huang_3_and_5_gradientlike_system}
We return to the abstract setting considered in Section \ref{sec:Huang_3_and_5_gradient-system}, but now develop growth estimates for `pseudogradient' or `gradient-like' systems, rather than the pure gradient systems discussed in Section \ref{sec:Huang_3_and_5_gradient-system}. The generality we allow lies somewhere between that of Section \ref{sec:Huang_3_and_5_gradient-system} and that considered by Huang in \cite[Section 3.3]{Huang_2006} as results in the latter source, in our view, require hypotheses which appear to us more difficult to apply in practice than those required in this section.

\subsection{Gradient-like flow}
\label{sec:Huang_3_gradientlike_system}
We introduce a modification of the pure gradient flow for $\sE$ previously considered and modeled on that of \cite[Equations (0.1) or (3.1)]{Simon_1983}.

\begin{defn}[Strong solution to a gradient-like system]
\label{defn:Strong_solution_to_gradientlike_system}
Let $\sE':\sU\subset \sX \to \sH$ be a gradient map associated with a $C^1$ functional, $\sE:\sU \to \RR$, where $\sU$ is an open subset of a real, reflexive Banach space, $\sX$, that is continuously embedded and dense in a Hilbert space, $\sH$, and let $R:[0, T) \to \sH$ be a continuous map for some $0 < T \leq \infty$. We call a trajectory, $u:[0,T)\to \sU$ (in the sense of Definition \ref{defn:Huang_3-1-1}), a \emph{strong solution} of a \emph{pseudogradient} or \emph{gradient-like system} for $\sE$ if
\begin{equation}
\label{eq:Simon_0-1}
\dot u(t) = -\sE'(u(t)) + R(t), \quad\hbox{a.e. } t \in (0,T), \quad u(0) = u_0 \in \sU,
\end{equation}
as an equation in $\sH$.
\end{defn}

Our Definition \ref{defn:Huang_3-1-3} of mild and weak solutions to the Cauchy problem \eqref{eq:Simon_0-1} for a pure gradient system easily extend to those of mild and weak solutions to the gradient-like system \eqref{eq:Simon_0-1}, but, as with pure gradient systems, we shall make little use of those solution concepts.

Naturally, one requires conditions on the perturbation, $R(t)$, in \eqref{eq:Simon_0-1} in order to develop properties of the flow. However, we shall \emph{not} require that the autonomous term, $R(t)$, obey the exponential decay condition in \cite[Equation (0.1)]{Simon_1983} or the inequality constraint in \cite[Equation (3.1)]{Simon_1983}, but rather more relaxed conditions consistent with our Hypothesis \ref{hyp:Huang_3-10} below and which are simpler versions of those given by Huang in \cite[Equations (3.10a) and (3.10b) or (3.10b$'$)]{Huang_2006}, \cite[Equations (3.27a) and (3.27b) or (3.27b$'$)]{Huang_2006}, or \cite[(5.3) and (5.5)]{Huang_2006}.

\begin{hyp}[Conditions on the perturbation in a gradient-like system]
\label{hyp:Huang_3-10}
The perturbation, $R\in C([0, T); \sH)$, in Definition \ref{defn:Strong_solution_to_gradientlike_system} is such that
\begin{subequations}
\label{eq:Huang_3-10a}
\begin{align}
\label{eq:Huang_3-10a_gradE_innerproduct_dotu}
(-\sE'(u(t)), \dot u(t))_\sH &\geq \|\sE'(u(t))\|_\sH \|\dot u(t)\|_\sH + F'(t),
\\
\label{eq:Huang_3-10a_gradE_norm_H}
\|\sE'(u(t))\|_\sH &\geq \|\dot u(t)\|_\sH + G'(t),  \quad\hbox{for a.e. } t \geq 0,
\end{align}
\end{subequations}
where the functions $F, G: [0, \infty) \to [0, \infty)$ are absolutely continuous, non-increasing (non-negative) functions satisfying
\begin{subequations}
\label{eq:Huang_3-10b}
\begin{align}
\label{eq:Huang_3-10b_limit_F+G_zero}
\lim_{t\to \infty} F(t) + G(t) &= 0,
\\
\label{eq:Huang_3-10b_integral_zero_to_infinity_of_phiF_finite}
\int_0^\infty |F(t)|^\theta \,dt &< \infty,
\end{align}
\end{subequations}
where $\phi(x) := c|x - a|^\theta$, $x \in \RR$, for constants $c > 0$ and $a \in \RR$ and $\theta \in [1/2, 1)$.
\end{hyp}

See Section \ref{sec:Huang_3_and_5_refined_hypotheses} for a discussion of generalizations and variants of the conditions on $F$, $G$, and (implicitly, $R$) in Hypothesis \ref{hyp:Huang_3-10}.

To give some examples of constraints on $R(t)$ that are sufficient to ensure that conditions \eqref{eq:Huang_3-10b_limit_F+G_zero} and \eqref{eq:Huang_3-10b_integral_zero_to_infinity_of_phiF_finite} or \eqref{eq:Huang_3-10b_integral_one_over_rho_plus_rho_times_Fprime} are satisfied, we proceed as in the proof of \cite[Lemma 1]{Simon_1983} and expand
$$
\|R(t)\|_\sH^2 = \|\dot u(t) + \sE'(u(t))\|_\sH^2 = \|\dot u(t)\|_\sH^2 + 2(\sE'(u(t)), \dot u(t))_\sH + \|\sE'(u(t))\|_\sH^2,
$$
to give
\begin{align*}
(-\sE'(u(t)), \dot u(t))_\sH &= \frac{1}{2}\left(\|\dot u(t)\|_\sH^2 + \|\sE'(u(t))\|_\sH^2\right) - \frac{1}{2}\|R(t)\|_\sH^2,
\\
&\geq \|\dot u(t)\|_\sH\|\sE'(u(t))\|_\sH  - \frac{1}{2}\|R(t)\|_\sH^2,
\end{align*}
while
$$
\|\sE'(u(t))\|_\sH \geq \|\dot u(t)\|_\sH - \|R(t)\|_\sH.
$$
Thus, we identify
\begin{subequations}
\label{eq:Huang_Time_derivatives_of_F_and_G_in_terms_of_Simon_R}
\begin{align}
\label{eq:Huang_Time_derivative_of_F_in_terms_of_Simon_R}
\dot F(t) &= -\frac{1}{2}\|R(t)\|_\sH^2,
\\
\label{eq:Huang_Time_derivative_of_G_in_terms_of_Simon_R}
\dot G(t) &= -\|R(t)\|_\sH, \quad\hbox{for a.e. } t > 0.
\end{align}
\end{subequations}
and hence
$$
F(t) = F(0) - \frac{1}{2}\int_0^t\|R(s)\|_\sH^2\,ds
\quad\hbox{and}\quad
G(t) = G(0) - \int_0^t\|R(s)\|_\sH\,ds, \quad\forall\, t \geq 0.
$$
The condition \eqref{eq:Huang_3-10b} implies that
$$
F(0) = \frac{1}{2}\int_0^\infty\|R(s)\|_\sH^2\,ds
\quad\hbox{and}\quad
G(0) = \int_0^\infty\|R(s)\|_\sH\,ds,
$$
and this also ensures that $F$ and $G$ are non-negative on $[0, \infty)$. We can thus write,
\begin{subequations}
\label{eq:Huang_F_and_G_in_terms_of_Simon_R}
\begin{align}
\label{eq:Huang_F_in_terms_of_Simon_R}
F(t) &= \frac{1}{2}\int_t^\infty\|R(s)\|_\sH^2\,ds,
\\
\label{eq:Huang_G_in_terms_of_Simon_R}
G(t) &= \int_t^\infty\|R(s)\|_\sH\,ds, \quad\forall\, t \geq 0.
\end{align}
\end{subequations}
We illustrate that the conditions on $F$ and $G$ can be obeyed with an example of a decay condition for $\|R(t)\|_\sH$ as $t \infty$, reminiscent of an explicit assumption by Simon in his \cite[Theorem 2]{Simon_1983}.

\begin{exmp}[Exponential decay for the perturbation term in a gradient-like system]
If $\|R(t)\|_\sH \leq Ct^{-\beta}$, for some $C > 0$ and $\beta > 1$ as $t\to\infty$, then
$$
G(t) = \int_t^\infty\|R(s)\|_\sH\,ds \leq C\int_t^\infty s^{-\beta}\,ds = \frac{Ct^{1-\beta}}{\beta-1} < \infty,
$$
for all $t \geq 0$, while
$$
F(t) = \frac{1}{2}\int_t^\infty\|R(s)\|_\sH^2\,ds \leq \frac{C^2}{2}\int_t^\infty s^{-2\beta}\,ds
= \frac{C^2t^{1-2\beta}}{2(2\beta-1)} < \infty.
$$
Hence,
$$
\int_1^\infty \sqrt{F(t)}\,dt \leq \frac{C}{\sqrt{4\beta-2}}\int_1^\infty t^{1/2 - \beta}\,dt < \infty,
$$
provided $\beta > 3/2$; finiteness of the integral $\int_0^1 \sqrt{F(t)}\,dt$ is assured by the fact that $F$ is (absolutely) continuous on $[0,\infty)$. Therefore, it will suffice to show that the autonomous term, $R(t)$, in \eqref{eq:Simon_0-1} obeys
\begin{equation}
\label{eq:Simon_3-1_bridge_Huang_3-10b}
\|R(t)\|_\sH \leq Ct^{-\beta}, \quad t \geq T,
\end{equation}
for positive constants $C$, $T$, and $\beta > 3/2$.
\end{exmp}

We shall also consider situations where it is only known \apriori that $0<T<\infty$ in \eqref{eq:Simon_0-1} and so the solution, $u(t)$, and perturbation, $R(t)$, are only defined \apriori for $t \in [0, T)$. In that case, it is natural to choose
\begin{subequations}
\label{eq:Huang_F_and_G_in_terms_of_Simon_R_finite_interval}
\begin{align}
\label{eq:Huang_F_in_terms_of_Simon_R_finite_interval}
F(t) &= \frac{1}{2}\int_t^T \|R(s)\|_\sH^2\,ds,
\\
\label{eq:Huang_G_in_terms_of_Simon_R_finite_interval}
G(t) &= \int_t^T \|R(s)\|_\sH\,ds, \quad\forall\, t \in [0, T).
\end{align}
\end{subequations}
The functions $F$ and $G$ defined in \eqref{eq:Huang_F_and_G_in_terms_of_Simon_R_finite_interval} also obey the conditions in Hypothesis \ref{hyp:Huang_3-10}, but with $T = \infty$ replaced by $T < \infty$.

\subsection{Convergence in abstract gradient-like systems}
\label{subsec:Huang_3-3_gradientlike_system}
We have the following version of \cite[Lemma 3.3.4]{Huang_2006}, an extension of our Lemma \ref{lem:Huang_3-3-4} to the case of a gradient-like system.


\begin{lem}[Growth estimate in $\sH$ for a strong solution to a gradient-like system under the validity of the {\L}ojasiewicz-Simon gradient inequality]
\label{lem:Huang_3-3-4_gradient-like}
Let $\sU$ be an open subset of a real Banach space, $\sX$, that is continuously embedded and dense in a Hilbert space, $\sH$. Let $\sE:\sU\subset \sX\to\RR$ be a $C^1$ functional with gradient map $\sE':\sU\subset \sX \to \sH$, and $T>0$, and $R \in C([0, T]; \sH)$. Assume that the gradient map $\sE'$ satisfies a {\L}ojasiewicz-Simon gradient inequality \eqref{eq:Simon_2-2} with positive constants $c$, $\sigma$, and $\theta \in [1/2, 1)$, and critical point $\varphi \in \sU$, and that the perturbation, $R$, obeys Hypothesis \ref{hyp:Huang_3-10}. If $u \in C([0, T]; \sU)$ is a strong solution to the Cauchy problem \eqref{eq:Simon_0-1} in the sense of Definition \ref{defn:Strong_solution_to_gradientlike_system} and the gradient map $\sE'$ satisfies the {\L}ojasiewicz-Simon gradient inequality in the orbit $O(u)$ as in \eqref{eq:Huang_3-15a}, then,
\begin{subequations}
\label{eq:Huang_3-18_H_gradient-like}
\begin{gather}
\label{eq:Huang_3-18a_H_gradient-like}
\|\dot u(t)\|_\sH
\leq
-\frac{d}{dt} \Phi(H(t)) + |F(t)|^\theta - G'(t), \quad\hbox{a.e. } t \in [0,T],
\\
\label{eq:Huang_3-18b_H_gradient-like}
\int_0^T\|\dot u(t)\|_\sH\,dt \leq \int_{H(T)}^{H(0)} \frac{1}{\phi(s)} \,ds + G(0) + \int_0^T |F(t)|^\theta \,dt,
\end{gather}
\end{subequations}
where
\begin{subequations}
\begin{align}
\label{eq:Huang_3-18c}
H(t) &:= \sE(u(t)) + F(t), \quad\forall\, t \in [0, T],
\\
\label{eq:Huang_3-18d}
\Phi(x) &:= \int_0^x \frac{ds}{\phi(s)}, \quad\forall\, x \in \RR,
\end{align}
\end{subequations}
and, following \eqref{eq:Simon_2-2},
\begin{equation}
\label{eq:Simon_2-2_gradient_inequality_function_phi}
\phi(x) := c|x - \sE(\varphi)|^\theta, \quad \forall\, x \in \RR.
\end{equation}
\end{lem}

\begin{proof}
We adapt the proof of \cite[Lemma 3.3.4]{Huang_2006} and our Lemma \ref{lem:Huang_3-3-4}. The function $[0, \infty) \ni t \mapsto \sE(u(t)) \in \RR$ is $C^1$ by Proposition \ref{prop:Huang_3-1-2} and therefore $H(t)$ obeys
\begin{align*}
-H'(t) &= \langle -\sE'(u(t)), \dot u(t) \rangle_{\sX'\times \sX} - F'(t) \quad\hbox{(by Proposition \ref{prop:Huang_3-1-2} and \eqref{eq:Huang_3-18c})}
\\
&= (-\sE'(u(t)), \dot u(t))_\sH - F'(t)  \quad\hbox{(by \eqref{eq:Riesz_pairing_relations}),}
\end{align*}
and thus it follows from \eqref{eq:Huang_3-10a_gradE_innerproduct_dotu} that
\begin{equation}
\label{eq:Huang_3-19}
-H'(t) \geq \|\dot u(t)\|_\sH\|\sE'(u(t))\|_\sH, \quad\hbox{a.e. } t \in [0,T].
\end{equation}
Hence, the absolutely continuous function, $H$, is non-increasing. Therefore, an application of Lemma \ref{lem:Huang_3-2-1} shows that the composition $\Phi\circ H$ of $H$ with the absolutely continuous function $\Phi$ defined by \eqref{eq:Huang_3-18d} is again absolutely continuous and
$$
\frac{d}{dt}\Phi(H(t)) = \Phi'(H(t)) H'(t) = \frac{H'(t)}{\phi(H(t))}, \quad \hbox{a.e. } t \in [0,T].
$$
(The arguments in the remainder of the proof show that the preceding ratio is well-defined a.e. on $[0, T]$.) To estimate $\phi(H(t))$, we can use the elementary inequality (see Remark \ref{rmk:Proof_power_inequality_nonnegative_x_and_y}), for $p \in [0, 1]$,
\begin{equation}
\label{eq:Power_inequality_nonnegative_x_and_y}
(x + y)^p \leq x^p + y^p, \quad\forall\, x, y \geq 0,
\end{equation}
noting that, if $x < 0$ and $y \geq 0$, then we also have
$$
|x + y|^p \leq |-x + y|^p \leq (-x)^p + y^p.
$$
Therefore,
\begin{equation}
\label{eq:Power_inequality_real_x_and_y}
|x + y|^p \leq |x|^p + |y|^p, \quad\forall\, x, y \in \RR.
\end{equation}
It follows that, by our definition $\phi(x) = c|x - \sE(\varphi)|^\theta$, for $x \in \RR$, and noting that $\theta \in [1/2, 1)$ by hypothesis,
\begin{align*}
\phi(H(t)) &= \phi(\sE(u(t)) + F(t))
\\
&= c|\sE(u(t)) - \sE(\varphi) + F(t)|^\theta
\\
&\leq c|\sE(u(t)) - \sE(\varphi)|^\theta + |F(t)|^\theta
\\
&= \phi(\sE(u(t))) + |F(t)|^\theta, \quad \forall\, t \in [0, T].
\end{align*}
(The preceding inequality replaces the alternative bound for $\phi(H(t))$ given in \cite[p. 72]{Huang_2006}, which in turn relies on \cite[Definition 2.2.1]{Huang_2006}.) Combining the preceding inequality with \eqref{eq:Huang_3-19} yields
$$
-\frac{d}{dt}\Phi(H(t)) \geq \frac{\|\dot u(t)\|_\sH\|\sE'(u(t))\|_\sH}{\phi(\sE(u(t))) + |F(t)|^\theta},
$$
and thus, by the {\L}ojasiewicz-Simon gradient inequality \eqref{eq:Simon_2-2},
\begin{equation}
\label{eq:Huang_3-20}
-\frac{d}{dt}\Phi(H(t)) \geq \frac{\|\dot u(t)\|_\sH\|\sE'(u(t))\|_\sH}{\|\sE'(u(t))\|_\sH + |F(t)|^\theta},
 \quad\hbox{a.e. } t \in [0,T].
\end{equation}
Hence,
$$
-\frac{d}{dt}\Phi(H(t)) \geq \|\dot u(t)\|_\sH - \frac{|F(t)|^\theta \|\dot u(t)\|_\sH}{\|\sE'(u(t))\|_\sH + |F(t)|^\theta},
$$
which can be further estimated as, for a.e. $t \in [0,T]$,
\begin{equation}
\label{eq:Huang_3-21a}
\|\dot u(t)\|_\sH \leq -\frac{d}{dt}\Phi(H(t)) + |F(t)|^\theta \frac{\|\dot u(t)\|_\sH}{\|\sE'(u(t))\|_\sH + |F(t)|^\theta}.
\end{equation}
To estimate the second term in the right-hand side of \eqref{eq:Huang_3-21a}, we use the inequality \eqref{eq:Huang_3-10a_gradE_norm_H},
$$
\frac{\|\dot u(t)\|_\sH}{\|\sE'(u(t))\|_\sH + |F(t)|^\theta}
\leq
\frac{\|\sE'(u(t))\|_\sH - G'(t)}{\|\sE'(u(t))\|_\sH + |F(t)|^\theta},
$$
and thus, noting that $-G'(t) \geq 0$ for a.e. $t \in [0,T]$ (because $G$ is non-increasing on $[0,T]$), we find that
\begin{equation}
\label{eq:Huang_3-21b}
\frac{\|\dot u(t)\|_\sH}{\|\sE'(u(t))\|_\sH + |F(t)|^\theta} \leq 1 - \frac{G'(t)}{|F(t)|^\theta},
\quad\hbox{a.e. } t \in [0,T].
\end{equation}
Combining \eqref{eq:Huang_3-21a} and \eqref{eq:Huang_3-21b} gives
$$
\|\dot u(t)\|_\sH \leq -\frac{d}{dt}\Phi(H(t)) + |F(t)|^\theta - G'(t), \quad\hbox{a.e. } t \in [0,T],
$$
which is the desired estimate \eqref{eq:Huang_3-18a_H_gradient-like}.

Integrating \eqref{eq:Huang_3-18a_H_gradient-like} over $[0,T]$ gives
$$
\int_0^T\|\dot u(t)\|_\sH \,dt
\leq
\Phi(H(0)) - \Phi(H(T)) + \int_0^T |F(t)|^\theta \,dt + G(0) - G(T).
$$
Noting that $G(T) \geq 0$ by Hypothesis \ref{hyp:Huang_3-10} and recalling the definition \eqref{eq:Huang_3-18d} of $\Phi$ yields the estimate \eqref{eq:Huang_3-18b_H_gradient-like}.
\end{proof}

\begin{rmk}[Proof of Inequality \eqref{eq:Power_inequality_nonnegative_x_and_y}]
\label{rmk:Proof_power_inequality_nonnegative_x_and_y}
First observe that it suffices to assume $x+y=1$ because, if the Inequality \eqref{eq:Power_inequality_nonnegative_x_and_y} is true under that constraint then, in general, we have
$$
(x + y)^p
=
\left(\frac{x}{x+y} + \frac{y}{x+y}\right)^p (x + y)^p
\leq
\left(\frac{x}{x+y}\right)^p (x + y)^p + \left(\frac{y}{x+y}\right)^p (x + y)^p
=
x^p + y^p,
$$
assuming that one of $x$ or $y$ is positive. Thus, assuming $x+y=1$, we have $0 \leq x, y \leq 1$ and
$$
(x + y)^p = 1 = x + y \leq x^p + y^p,
$$
noting that $0 \leq p \leq 1$, and this completes the proof.
\end{rmk}

In the case of a strong solution, $u:[0,T)\to \sX$, to a pure gradient system \eqref{eq:Huang_3-3a}, one makes frequent use of the fact that the energy, $\sE(u(t))$, is a non-increasing function of $t \in [0, T)$. For a gradient-like system, we have the following corollary of the proof of Lemma \ref{lem:Huang_3-3-4_gradient-like}.

\begin{cor}[Approximate monotonicity of the energy of a strong solution to a gradient-like system]
\label{cor:Huang_3-3-4_gradient-like_approximate_energy_monotonicity}
Assume the hypotheses and notation of Lemma \ref{lem:Huang_3-3-4_gradient-like}. Then the function, $H(t) = \sE(u(t)) + F(t)$ for $t \in [0, T)$, in \eqref{eq:Huang_3-18c} is absolutely continuous and non-increasing on $[0, T)$.
\end{cor}

We shall need a technical hypothesis to allow us to relate convergence of $u(t)$ in the Hilbert space $\sH$ to the stronger convergence of $u(t)$ in the Banach space $\sX$ as $t \to \infty$, extending our previous Hypothesis \ref{hyp:Abstract_apriori_interior_estimate_trajectory} appropriate for a solution to a pure gradient system \eqref{eq:Huang_3-3a} to one appropriate for a solution to a gradient-like system such as \eqref{eq:Simon_0-1}. The shape of that \apriori estimate will depend on how \eqref{eq:Simon_0-1} is expressed in the form of a parabolic equation --- for example, Yang-Mills gradient-like flow together with a Coulomb gauge condition and how the Coulomb gauge condition is introduced.

In practice, gradient-like systems such as \eqref{eq:Simon_0-1} arise because of the application of simple operations to a pure gradient system \eqref{eq:Huang_3-3a} --- for example, application of a cut off function or changing a Riemannian metric on a base manifold --- and thus we can make use of the \apriori estimate obeyed by solutions to the original pure gradient system \eqref{eq:Huang_3-3a}, after allowing for a small error. With that background in mind, we introduce the

\begin{hyp}[Regularity and \apriori interior estimate for a trajectory]
\label{hyp:Abstract_apriori_interior_estimate_trajectory_with_perturbation}
Let $C_1$ and $\mu$ be positive constants and let $T \in (0, \infty]$. Given a trajectory $u:[0,T)\to \sU$, in the sense of Definition \ref{defn:Huang_3-1-1}, and a positive constant $\eps_0$, we say that $\dot u:[0,T) \to \sH$ obeys an \apriori \emph{interior estimate on $(0, T]$} if, for every $S \geq 0$ and $\delta > 0$ obeying $S+\delta \leq T$, the map $\dot u:[S+\delta,T) \to \sX$ is Bochner integrable and there holds
\begin{equation}
\label{eq:Abstract_apriori_interior_estimate_trajectory_with_small_constant_error}
\int_{S+\delta}^T \|\dot u(t)\|_\sX\,dt \leq C_1(1+\delta^{-\mu})\int_S^T \|\dot u(t)\|_\sH\,dt + \eps_0.
\end{equation}
\end{hyp}

Often it will suffice that the constant, $\eps_0$, in \eqref{eq:Abstract_apriori_interior_estimate_trajectory_with_small_constant_error} is merely finite but occasionally we shall also require that $\eps_0$ be small. An immediate consequence of Lemma \ref{lem:Huang_3-3-4_gradient-like} and Hypothesis \ref{hyp:Abstract_apriori_interior_estimate_trajectory_with_perturbation} is the following refined growth estimate.

\begin{cor}[Growth estimate in $\sX$ for a strong solution to a gradient-like system under the validity of the {\L}ojasiewicz-Simon gradient inequality]
\label{cor:Huang_3-3-4_gradient-like}
Assume the hypotheses of Lemma \ref{lem:Huang_3-3-4_gradient-like} and, in addition, that $u$ satisfies Hypothesis \ref{hyp:Abstract_apriori_interior_estimate_trajectory_with_perturbation} with $R \in W^{1,1}(0,T;\sH)$ in \eqref{eq:Simon_0-1}. Then, for every $\delta \in (0, T]$, the solution $u$ obeys
\begin{equation}
\label{eq:Huang_3-18b_X_gradient-like}
\int_\delta^T \|\dot u(t)\|_\sX\,dt \leq C_1(1+\delta^{-\mu})\left(\int_{H(0)}^{H(T)} \frac{1}{\phi(s)} \,ds + G(0)
+ \int_0^T |F(t)|^\theta \,dt\right) + \eps_0.
\end{equation}
\end{cor}

In addition, we obtain convergence for a strong solution to a gradient-like system, analogous to \cite[Theorem 3.3.5]{Huang_2006}, and extending the previous convergence result for a strong solution to a gradient system in Theorem \ref{thm:Huang_3-3-5}.

\begin{thm}[Convergence for a strong solution to a gradient-like system under the validity of the {\L}ojasiewicz-Simon gradient inequality]
\label{thm:Huang_3-3-5_gradient-like}
Let $\sU$ be an open subset of a real Banach space, $\sX$, that is continuously embedded and dense in a Hilbert space, $\sH$. Let $\sE:\sU\subset \sX\to\RR$ be a $C^1$ functional with gradient map $\sE':\sU\subset \sX \to \sH$ and $u:[0, \infty)\to \sU$ be a strong solution to the gradient-like system \eqref{eq:Simon_0-1} in the sense of Definition \ref{defn:Strong_solution_to_gradientlike_system}. If the gradient map $\sE'$ satisfies a {\L}ojasiewicz-Simon gradient inequality \eqref{eq:Simon_2-2} in the orbit $O(u)$, then $u(t)$ converges in $\sH$ as $t\to\infty$ in the sense that
$$
\int_0^\infty \|\dot u(t)\|_\sH\,dt < \infty.
$$
If in addition $u$ obeys Hypothesis \ref{hyp:Abstract_apriori_interior_estimate_trajectory_with_perturbation}, then $u(t)$ converges in $\sX$ as $t\to\infty$ in the sense that
$$
\int_1^\infty \|\dot u(t)\|_\sX\,dt < \infty.
$$
\end{thm}

\subsection{Simon Alternative and convergence for an abstract gradient-like system}
\label{subsec:Simon_alternative_gradientlike}
We next have the following abstract, gradient-like analogue of R\r{a}de's \cite[Proposition 7.4]{Rade_1992}, in turn a variant the \emph{Simon Alternative}, namely \cite[Theorem 2]{Simon_1983}. Theorem \ref{thm:Huang_3-3-6_and_7_gradientlike} below extends our previous Theorem \ref{thm:Huang_3-3-6} from the context of gradient to gradient-like systems.

\begin{thm}[Simon Alternative for convergence for a strong solution to a gradient system]
\label{thm:Huang_3-3-6_and_7_gradientlike}
Let $\sU$ be an open subset of a real Banach space, $\sX$, that is continuously embedded and dense in a Hilbert space, $\sH$. Let $\sE:\sU\subset \sX\to\RR$ be a $C^1$ functional with gradient map $\sE':\sU\subset \sX \to \sH$. Assume that
\begin{enumerate}
\item $\varphi \in \sU$ is a critical point of $\sE$, that is $\sE'(\varphi)=0$, and that the gradient map $\sE':\sU_\sigma\subset \sX \to \sX'$ satisfies a {\L}ojasiewicz-Simon gradient inequality \eqref{eq:Simon_2-2}, for positive constants $c$, $\sigma$, and $\theta \in [1/2,1)$;

\item Given positive constants $b$, $\eta$, and $\tau$, there is a constant $\delta = \delta(\eta, \tau, b) \in (0, \tau]$ such that if $v$ is a solution to the gradient-like system \eqref{eq:Simon_0-1} on $[t_0, t_0 + \tau)$ with $\|v(t_0)\|_\sX \leq b$, then
\begin{equation}
\label{eq:Gradientlike_solution_near_initial_data_at_t0_for_short_enough_time}
\sup_{t\in [t_0, t_0+\delta]}\|v(t) - v(t_0)\|_\sX < \eta.
\end{equation}
\end{enumerate}
Then there is a constant $\eps = \eps(c, C_1, \delta, \theta, \mu, \sigma, \tau, \varphi) \in (0, \sigma/4)$ with the following significance.  If $u:[0, \infty)\to \sU$ is a strong solution to \eqref{eq:Simon_0-1} in the sense of Definition \ref{defn:Strong_solution_to_gradientlike_system} that satisfies Hypothesis \ref{hyp:Abstract_apriori_interior_estimate_trajectory_with_perturbation} on $(0, \infty)$ and there is a $T \geq 0$ such that
\begin{equation}
\label{eq:Rade_7-2_banach_gradientlike}
\|u(T) - \varphi\|_\sX < \eps,
\end{equation}
and, denoting $\phi(x) = c|x - \sE(\varphi)|^\theta$ as in \eqref{eq:Simon_2-2_gradient_inequality_function_phi}, for all $x \in \RR$,
\begin{subequations}
\label{eq:Huang_5-6_Tfinite}
\begin{align}
\label{eq:Huang_5-6_Tfinite_F}
F(T) &< \eps,
\\
\label{eq:Huang_5-6_Tfinite_G_plus_integral_phiF}
G(T) + \int_T^\infty |F(t)|^\theta \,dt &< \eps,
\\
\label{eq:Huang_5-6_Tfinite_R}
\eps_0 &< \eps,
\end{align}
\end{subequations}
where $F$ and $G$ are as in Hypothesis \ref{hyp:Huang_3-10} and $\eps_0$ is the positive constant in Hypothesis \ref{hyp:Abstract_apriori_interior_estimate_trajectory_with_perturbation}, then either
\begin{enumerate}
\item
\label{item:Theorem_3-3-6_and_7_gradientlike_energy_u_at_time_t_below_energy_critical_point}
$\sE(u(t)) \leq \sE(\varphi) - \eps$ for some $t>T$, or
\item
\label{item:Theorem_3-3-6_and_7_gradientlike_u_converges_to_limit_u_at_infty}
the trajectory $u(t)$ converges in $\sX$ to a limit $u_\infty \in \sX$ in the sense that
$$
\lim_{t\to\infty}\|u(t)-u_\infty\|_\sX =0
\quad\hbox{and}\quad
\int_1^\infty \|\dot u\|_\sX\,dt < \infty.
$$
If $\varphi$ is a cluster point of the orbit $O(u) = \{u(t): t\geq 0\}$, then $u_\infty = \varphi$.
\end{enumerate}
\end{thm}

\begin{rmk}[Comparison between Theorem \ref{thm:Huang_3-3-6_and_7_gradientlike} and Huang's Theorems 3.3.6 and 3.3.7 in \cite{Huang_2006}]
\label{rmk:Huang_incomplete_proof_of_Theorem_3-3-6}
In his \cite[Theorem 3.3.6 and 3.3.7]{Huang_2006}, Huang also asserts strong convergence in $\sX$ of $u(t)$ as $t\to\infty$ but (in his more general setting) he hypothesizes that $\varphi$ is a cluster point of $O(u)$ in $\sX$ and that either
\begin{inparaenum}[\itshape i\upshape)]
\item the orbit $O(u)$ is precompact in $\sX$; or
\item $\sY = \sX'$, and thus $\sY' = \sX''$, and the operator $\gamma$ appearing in his \cite[Equation (3.10a)]{Huang_2006} is a homeomorphism on $\sY'$.
\end{inparaenum}
An essential difference between Theorem \ref{thm:Huang_3-3-6_and_7_gradientlike} and \cite[Theorems 3.3.6 or 3.3.7]{Huang_2006} is that we do \emph{not} assume that $\varphi$ is a cluster point of the orbit $O(u)$ or that the orbit $O(u)$ is precompact in $\sX$. Indeed, our Theorem \ref{thm:Huang_3-3-6_and_7_gradientlike} is closer in spirit to Simon's \cite[Theorem 2]{Simon_1983}.
\end{rmk}

\begin{proof}[Proof of Theorem \ref{thm:Huang_3-3-6_and_7_gradientlike}]
We shall apply an argument similar to the one we used to prove Theorem \ref{thm:Huang_3-3-6}. We assume that Alternative \eqref{item:Theorem_3-3-6_and_7_gradientlike_energy_u_at_time_t_below_energy_critical_point} does not occur and thus,
\begin{equation}
\label{eq:Simon_3-2}
\sE(u(t)) > \sE(\varphi) - \eps, \quad\forall\, t \geq T.
\end{equation}
Note also that $\lim_{t\to\infty}\sE(u(t)) \geq \sE(\varphi) - \eps$, where the limit exists because the function $H(t) = \sE(u(t)) + F(t)$ is non-increasing by the proof of Lemma \ref{lem:Huang_3-3-4_gradient-like}. We seek to establish Alternative \ref{item:Theorem_3-3-6_and_7_gradientlike_u_converges_to_limit_u_at_infty} and begin with the following analogue of Claim \ref{claim:norm_X_of_u_minus_varphi_lessthan_sigma_all_t_geq_T}.

\begin{claim}[Confinement of solution to a {\L}ojasiewicz-Simon gradient inequality neighborhood]
\label{claim:norm_X_of_u_minus_varphi_lessthan_sigma_all_t_geq_T_gradientlike}
There is a constant
$\eps = \eps(c,C_1,\delta, \theta, \mu, \sigma, \tau, \varphi) \in (0, \sigma/4)$ in \eqref{eq:Rade_7-2_banach_gradientlike} such that
\begin{equation}
\label{eq:norm_X_of_u_minus_varphi_lessthan_sigma_all_t_geq_T_gradientlike}
\|u(t) - \varphi\|_\sX < \frac{\sigma}{2}, \quad\forall\, t \in [T, \infty).
\end{equation}
\end{claim}

\begin{proof}[Proof of Claim \ref{claim:norm_X_of_u_minus_varphi_lessthan_sigma_all_t_geq_T_gradientlike}]
We adapt the proof of Claim \ref{claim:norm_X_of_u_minus_varphi_lessthan_sigma_all_t_geq_T}. We suppose the conclusion is false and obtain a contradiction. Let $\widehat T > T$ be the smallest number such that
\begin{equation}
\label{eq:Rade_7-6_banach_gradientlike}
\|u(t) - \varphi\|_\sX < \frac{\sigma}{2} \quad\forall\, t \in [T, \widehat T),
\quad\hbox{but}\quad
\|u(\widehat T) - \varphi\|_\sX \geq \frac{\sigma}{2}.
\end{equation}
By \eqref{eq:Gradientlike_solution_near_initial_data_at_t0_for_short_enough_time} with $t_0 = T$ and $\eta = \sigma/4$, we may choose $\delta = \delta(\sigma, \tau, \|\varphi\|_\sX) \in (0, \tau]$ so that
\begin{equation}
\label{eq:Rade_7-5_banach_gradientlike}
\sup_{t\in [T, T+\delta]}\|u(t) - u(T)\|_\sX < \frac{\sigma}{4},
\end{equation}
and, in particular,
$$
\|u(T+\delta) - u(T)\|_\sX < \frac{\sigma}{4}.
$$
By hypothesis \eqref{eq:Rade_7-2_banach_gradientlike} and the fact that we seek $\eps \in (0, \sigma/4)$, we also have
$$
\|u(T) - \varphi\|_\sX < \frac{\sigma}{4},
$$
and combining this inequality with \eqref{eq:Rade_7-5_banach_gradientlike} gives
$$
\sup_{t\in [T, T+\delta]}\|u(t) - \varphi\|_\sX < \frac{\sigma}{2}.
$$
Therefore, by definition of $\widehat T$ in \eqref{eq:Rade_7-6_banach_gradientlike}, we must have
$$
T > T + \delta.
$$
Inequality \eqref{eq:Huang_3-18a_H_gradient-like}, with $[0,T]$ replaced by $[T,\widehat T]$, gives
$$
\|\dot u(t)\|_\sH \leq -\frac{d}{dt} \Phi(H(t)) + |F(t)|^\theta - G'(t), \quad\hbox{a.e. } t \in [T,\widehat T],
$$
where $H(t) = \sE(u(t)) + F(t)$ from \eqref{eq:Huang_3-18c} and $\Phi$ is defined by \eqref{eq:Huang_3-18d}. Integrating this inequality over $[T, \widehat T]$, noting that $-\int_{[T, \widehat T]} G'(t) \,dt = G(T) - G(\widehat T)$ and discarding the term $-G(\widehat T) \leq 0$ (the function $G(t)$ is non-negative and non-increasing by Hypothesis \ref{hyp:Huang_3-10}), yields
\begin{equation}
\label{eq:Huang_3-18b_H_with_0toT_replaced_by_TtotildeT_gradientlike}
\int_T^{\widehat T}\|\dot u(t)\|_\sH\,dt
\leq
\int_{H(\widehat T)}^{H(T)} \frac{1}{\phi(s)} \,ds + G(T) + \int_T^{\widehat T} |F(t)|^\theta \,dt.
\end{equation}
Here, we note that
$$
H(T) = \sE(u(T)) + F(T) \quad\hbox{and}\quad \phi(s) = c|s - \sE(\varphi)|^\theta,
$$
The terms $G(T)$ and $\int_T^{\widehat T} |F(t)|^\theta \,dt$ may be assumed as small as desired for large enough $T \geq 0$ by \eqref{eq:Huang_3-10b_limit_F+G_zero}. To examine the first integral term on the right-hand side of \eqref{eq:Huang_3-18b_H_with_0toT_replaced_by_TtotildeT_gradientlike}, recall from the proof of Lemma \ref{lem:Huang_3-3-4_gradient-like} that $H(t) = \sE(u(t)) + F(t)$ is continuous and non-increasing and so there are times $t_*$ and $t_{**}$ such that $T \leq t_* \leq t_{**} \leq \widehat T$ with $H(t) \geq \sE(\varphi)$ on $[T, t_*]$, and $H(t) = \sE(\varphi)$ on $[t_*, t_{**}]$, and $H(t) \leq \sE(\varphi)$ on $[t_{**}, \widehat T]$. (By the assumption of the second alternative in our hypotheses and the fact that $\lim_{t\to\infty} F(t) = 0$ by \eqref{eq:Huang_3-10b_limit_F+G_zero}, we have $\lim_{t\to\infty} H(t) \geq \sE(\varphi)-\eps$.) Therefore, keeping in mind the possibility that $t_{**} = T$ or $t_* = \widehat T$ (compare the proof of \cite[Lemma 1]{Simon_1983}),
\begin{align*}
\int_{H(\widehat T)}^{H(T)} \frac{1}{\phi(s)} \,ds
&=
1_{\{t_{**} < \widehat T\}}\int_{H(\widehat T)}^{H(t^{**})} \frac{1}{c(\sE(\varphi) - s)^\theta} \,ds
+ 1_{\{t_* > T\}}\int_{H(t^*)}^{H(T)} \frac{1}{c(s - \sE(\varphi))^\theta} \,ds
\\
&= \frac{1}{c(1-\theta)}\left( (\sE(\varphi) - H(\widehat T))^{1-\theta} - (\sE(\varphi) - H(t^{**}))^{1-\theta}\right) 1_{\{t_{**} < \widehat T\}}
\\
&\quad + \frac{1}{c(1-\theta)}\left( (H(T) - \sE(\varphi))^{1-\theta} - (H(t^*) - \sE(\varphi))^{1-\theta} \right) 1_{\{t_* > T\}}
\\
&= \frac{1}{c(1-\theta)}\left( 1_{\{t_* > T\}}(H(T) - \sE(\varphi))^{1-\theta} + 1_{\{t_{**} < \widehat T\}}(\sE(\varphi) - H(\widehat T))^{1-\theta} \right),
\end{align*}
noting that $H(t^*) = \sE(\varphi) = H(t^{**})$. Thus, using $H(T) = \sE(u(T)) + F(T)$, the fact that $F \geq 0$ on $[0, \infty)$ by Hypothesis \ref{hyp:Huang_3-10}, the condition \eqref{eq:Simon_3-2} on $\sE(u(t))$, and
$$
H(\widehat T) = \sE(u(\widehat T)) + F(\widehat T) \geq \sE(u(\widehat T)) > \sE(\varphi) - \eps,
$$
we find that
\begin{equation}
\label{eq:Integral_HwidehatT_to_HT_of_1_over_phi}
\int_{H(\widehat T)}^{H(T)} \frac{1}{\phi(s)} \,ds \leq \frac{1}{c(1-\theta)}|\sE(u(T)) + F(T) - \sE(\varphi)|^{1-\theta} + \frac{\eps^{1-\theta}}{c(1-\theta)}.
\end{equation}
The term $F(T)$ may be assumed as small as desired for large enough $T \geq 0$ by \eqref{eq:Huang_3-10b_limit_F+G_zero}, while \eqref{eq:Rade_7-4_banach} with $v=u(T)$ and $\|u(T) - \varphi\|_\sX < \eps_\varphi$ (assured by \eqref{eq:Rade_7-2_banach_gradientlike} with $\eps\leq\eps_\varphi$) gives
$$
|\sE(\varphi) - \sE(u(T))| \leq \|u(T) - \varphi\|_\sX.
$$
Therefore, applying the inequality \eqref{eq:Power_inequality_nonnegative_x_and_y} to give
$$
|\sE(u(T)) + F(T) - \sE(\varphi)|^{1-\theta} \leq |\sE(u(T)) - \sE(\varphi)|^{1-\theta} + (F(T))^{1-\theta},
$$
and combining the preceding observations with \eqref{eq:Huang_3-18b_H_with_0toT_replaced_by_TtotildeT_gradientlike} and \eqref{eq:Integral_HwidehatT_to_HT_of_1_over_phi} yields
\begin{multline}
\label{eq:Bound_on_integral_T_to_widehatT_of_dotu_Hnorm}
\int_T^{\widehat T}\|\dot u(t)\|_\sH\,dt
\leq
\frac{1}{c(1-\theta)}\left( \|u(T) - \varphi\|_\sX^{1-\theta} + (F(T))^{1-\theta} \right)
+ \frac{\eps^{1-\theta}}{c(1-\theta)}
\\
+ G(T) + \int_T^\infty |F(t)|^\theta \,dt.
\end{multline}
Thus, by \eqref{eq:Rade_7-2_banach_gradientlike} and choosing $\eps\leq\eps_\varphi$ (given in \eqref{eq:Rade_7-4_banach}), and supposing that $T$ is large enough that \eqref{eq:Huang_5-6_Tfinite_F} and \eqref{eq:Huang_5-6_Tfinite_G_plus_integral_phiF} hold, we obtain
\begin{equation}
\label{eq:Bound_on_integral_from_T_to_tildeT_of_Hnorm_dotu_gradientlike}
\int_T^{\widehat T}\|\dot u(t)\|_\sH\,dt \leq \frac{3\eps^{1-\theta}}{c(1-\theta)} + \eps.
\end{equation}
On the other hand, noting that $T+\delta < \widehat T$,
\begin{align*}
\int_{T+\delta}^{\widehat T}\|\dot u(t)\|_\sX\,dt &\geq \|u(\widehat T) - u(T+\delta)\|_\sX
\\
&\geq \|u(\widehat T) - \varphi\|_\sX - \|u(T+\delta) - u(T)\|_\sX - \|u(T) - \varphi\|_\sX
\\
&> \frac{\sigma}{2} - \frac{\sigma}{4} - \eps  \quad\hbox{(by \eqref{eq:Rade_7-2_banach_gradientlike}, \eqref{eq:Rade_7-6_banach_gradientlike},
and \eqref{eq:Rade_7-5_banach_gradientlike})}.
\end{align*}
The preceding inequality implies that
\begin{align*}
\frac{\sigma}{4} - \eps &< \int_{T+\delta}^{\widehat T}\|\dot u(t)\|_\sX\,dt
\\
&\leq C_1(1+\delta^{-\mu})\int_T^{\widehat T}\|\dot u(t)\|_\sH\,dt + \eps_0
\quad\hbox{(by \eqref{eq:Simon_0-1} and \eqref{eq:Abstract_apriori_interior_estimate_trajectory_with_small_constant_error})}
\\
&\leq C_1(1+\delta^{-\mu})\left(\frac{3\eps^{1-\theta}}{c(1-\theta)} + \eps\right) + \eps_0
\quad\hbox{(by \eqref{eq:Bound_on_integral_from_T_to_tildeT_of_Hnorm_dotu_gradientlike})},
\end{align*}
and because \eqref{eq:Huang_5-6_Tfinite_R} holds, that is, $\eps_0 < \eps$, we obtain
\begin{equation}
\label{eq:Rade_gradientlike_epsiloncombination_greater_than_sigma_over_4}
\frac{\sigma}{4} - \eps
<
C_1(1+\delta^{-\mu})\left(\frac{3\eps^{1-\theta}}{c(1-\theta)} + \eps\right) + \eps.
\end{equation}
But choosing $\eps = \eps(c, C_1, \delta, \theta, \mu, \sigma, \tau, \varphi) \in (0, \sigma/4)$ small enough (and also requiring $T \geq 0$ large enough for
\eqref{eq:Huang_5-6_Tfinite} to hold) that
\begin{equation}
\label{eq:norm_X_of_u_minus_varphi_lessthan_sigma_all_t_geq_T_gradientlike_contradiction}
C_1(1+\delta^{-\mu})\left(\frac{3\eps^{1-\theta}}{c(1-\theta)} + \eps\right) + 2\eps
< \frac{\sigma}{4},
\end{equation}
contradicts \eqref{eq:Rade_gradientlike_epsiloncombination_greater_than_sigma_over_4}, and this completes the proof of Claim \ref{claim:norm_X_of_u_minus_varphi_lessthan_sigma_all_t_geq_T_gradientlike}.
\end{proof}

Therefore, because we may now set $\widehat T = \infty$ in \eqref{eq:Huang_3-18b_H_with_0toT_replaced_by_TtotildeT_gradientlike} and \eqref{eq:Integral_HwidehatT_to_HT_of_1_over_phi}, we obtain the inequality,
\begin{multline}
\label{eq:Rade_7-7H_banach_gradientlike}
\int_T^\infty\|\dot u(t)\|_\sH\,dt
\leq
\frac{1}{c(1-\theta)}|\sE(u(T)) + F(T) - \sE(\varphi)|^{1-\theta}
+ \frac{\eps^{1-\theta}}{c(1-\theta)}
\\
+ G(T) + \int_T^\infty |F(t)|^\theta \,dt.
\end{multline}
From \eqref{eq:Simon_0-1} and the \apriori estimate \eqref{eq:Abstract_apriori_interior_estimate_trajectory_with_small_constant_error}, we have
\begin{equation}
\label{eq:Rade_7-7X_banach_gradientlike}
\int_{T+1}^\infty\|\dot u(t)\|_\sX\,dt \leq 2C_1\int_T^\infty\|\dot u(t)\|_\sH\,dt + \eps_0.
\end{equation}
Combining \eqref{eq:Rade_7-7H_banach_gradientlike} and \eqref{eq:Rade_7-7X_banach_gradientlike} yields
\begin{multline}
\label{eq:Rade_7-7_banach_gradientlike}
\int_{T+1}^\infty\|\dot u(t)\|_\sX\,dt \leq
\frac{2C_1}{c(1-\theta)}|\sE(u(T)) + F(T) - \sE(\varphi)|^{1-\theta}
+ \frac{2C_1 \eps^{1-\theta}}{c(1-\theta)}
\\
+ 2C_1G(T) + 2C_1\int_T^\infty |F(t)|^\theta \,dt + \eps_0.
\end{multline}
In particular,
$$
\int_{T+1}^\infty\|\dot u(t)\|_\sX\,dt < \infty.
$$
The remainder of the proof is the same as that of Theorem \ref{thm:Huang_3-3-6}.
\end{proof}

\subsection{Existence and convergence of global solutions to abstract gradient-like systems started near a local minimum}
\label{subsec:Huang_5-1_existence_and_convergence_global_solution_gradientlike}
In this subsection, we extend our previous results in Section \ref{subsec:Huang_5-1_existence_and_convergence_global_solution_near_local_minimum} on existence and convergence of global solutions from abstract gradient to gradient-like systems. We have the following analogue of \cite[Theorem 5.1.1]{Huang_2006}. Theorem \ref{thm:Huang_5-1-1_gradientlike} below extends our previous Theorem \ref{thm:Huang_5-1-1} from the context of gradient to gradient-like systems.

\begin{thm}[Existence and convergence of a global solution to a gradient-like system near a local minimum]
\label{thm:Huang_5-1-1_gradientlike}
Let $\sU$ be an open subset of a real Banach space, $\sX$, that is continuously embedded and dense in a Hilbert space, $\sH$. Let $\sE:\sU\subset \sX\to\RR$ be a $C^1$ functional with gradient map $\sE':\sU\subset \sX \to \sH$. Let $\varphi \in \sU$ be a ground state of $\sE$ on $\sU$ and suppose that $\sE'$ obeys a {\L}ojasiewicz-Simon gradient inequality \eqref{eq:Simon_2-2} with positive constants $c$, $\sigma$, and $\theta \in [1/2, 1)$.
Assume that
\begin{enumerate}
\item For each $u_0 \in \sU$, there exists a unique strong solution to the Cauchy problem \eqref{eq:Simon_0-1}, in the sense of Definition \ref{defn:Strong_solution_to_gradientlike_system}, on a time interval $[0, \tau)$ for some positive constant, $\tau$;

\item Hypothesis \ref{hyp:Abstract_apriori_interior_estimate_trajectory_with_perturbation} holds with positive constants $C_1$ and $\mu$ for strong solutions to the gradient-like system \eqref{eq:Simon_0-1}; and

\item Given positive constants $b$ and $\eta$, there is a constant $\delta = \delta(\eta, \tau, b) \in (0, \tau]$ such that if $v$ is a solution to the gradient-like system \eqref{eq:Simon_0-1} on $[0, \tau)$ with $\|v_0\|_\sX \leq b$, then
\begin{equation}
\label{eq:Gradientlike_solution_near_initial_data_at_time_zero_for_short_enough_time}
\sup_{t\in [0, \delta]}\|v(t) - v_0\|_\sX < \eta.
\end{equation}
\end{enumerate}
Then there is a constant $\eps = \eps(c,C_1,\delta, \theta, \mu, \sigma, \tau, \varphi) \in (0, \sigma/4)$ with the following significance. If $F$ and $G$ in Hypothesis \ref{hyp:Huang_3-10} and the positive constant $\eps_0$ in Hypothesis \ref{hyp:Abstract_apriori_interior_estimate_trajectory_with_perturbation} obey
\begin{subequations}
\label{eq:Huang_5-6}
\begin{align}
\label{eq:Huang_5-6_F0_lessthan_epsilon}
F(0) &< \eps,
\\
\label{eq:Huang_5-6_G0_plus_integral_0_to_infinity_phiF_lessthan_epsilon}
G(0) + \int_0^\infty |F(t)|^\theta \,dt &< \eps,
\\
\label{eq:Huang_5-6_integral_0_to_infinity_dotR_norm_H_lessthan_epsilon}
\eps_0 &< \eps,
\end{align}
\end{subequations}
where $\phi(x) = c|x - \sE(\varphi)|^\theta$ as in \eqref{eq:Simon_2-2_gradient_inequality_function_phi}, for all $x \in \RR$, then, for each $u_0 \in \sU_\eps$, the Cauchy problem \eqref{eq:Simon_0-1} admits a global strong solution, $u:[0,\infty) \to \sU_{\sigma/2}$, that converges to a limit $u_\infty \in \sX$ as $t\to\infty$ with respect to the $\sX$ norm in the sense that
$$
\lim_{t \to \infty} \|u(t) - u_\infty\|_\sX = 0 \quad\hbox{and}\quad \int_1^\infty\|\dot u(t)\|_\sX\,dt < \infty.
$$
\end{thm}

\begin{proof}
We adapt the proof of Theorem \ref{thm:Huang_5-1-1}. By hypothesis, for $\eps > 0$ still to be determined and any $u_0 \in \sU_\eps$, the Cauchy problem \eqref{eq:Simon_0-1} admits a strong solution, $u:[0,\widetilde T)\to \sX$, where we let $\widetilde T$ denote its maximal lifetime. Clearly, $\widetilde T \geq \tau$, where $\tau > 0$ is the minimal lifetime of $u$ provided by our hypotheses. We make the

\begin{claim}[Confinement of solution to a {\L}ojasiewicz-Simon gradient inequality neighborhood]
\label{claim:Equation_Huang_5-7_gradientlike}
There is a constant
$\eps = \eps(c,C_1,\delta, \theta, \mu, \sigma, \tau, \varphi) \in (0, \sigma/4)$ such that, for any $u_0 \in \sU_\eps$ and solution, $u:[0,\widetilde T)\to \sX$, to \eqref{eq:Simon_0-1}, there holds
\begin{equation}
\label{eq:Huang_5-7_gradientlike}
\|u(t) - \varphi\|_\sX < \frac{\sigma}{2}, \quad\forall\, t \in [0,\widetilde T).
\end{equation}
\end{claim}

\begin{proof}[Proof of Claim \ref{claim:Equation_Huang_5-7_gradientlike}]
We provide an argument by contradiction. Suppose there exists $T \in (0,\widetilde T)$ such that
\begin{equation}
\label{eq:Huang_5-8_gradientlike}
\|u(t) - \varphi\|_\sX < \frac{\sigma}{2}, \quad\forall\, t \in [0,T),
\quad\hbox{but}\quad \|u(T) - \varphi\|_\sX \geq \frac{\sigma}{2}.
\end{equation}
By the inequality \eqref{eq:Gradientlike_solution_near_initial_data_at_time_zero_for_short_enough_time} with $\eta = \sigma/8$, we may choose $\delta = \delta(\sigma, \tau, \|\varphi\|_\sX) \in (0,\widetilde T)$ small enough that
\begin{equation}
\label{eq:Huang_theorem_5-1-1_proof_norm_X_udelta_minus_u0_inequality_gradientlike}
\sup_{t \in [0, \delta]}\|u(t)-u_0\|_\sX < \frac{\sigma}{8},
\end{equation}
and, in particular,
$$
\|u(\delta)-u_0\|_\sX < \frac{\sigma}{8}.
$$
Note also that, as we can see by comparing \eqref{eq:Huang_5-8_gradientlike} and \eqref{eq:Huang_theorem_5-1-1_proof_norm_X_udelta_minus_u0_inequality_gradientlike}, we must have
$$
T > \delta.
$$
We can thus apply the growth estimate \eqref{eq:Huang_3-18b_H_gradient-like} in Lemma \ref{lem:Huang_3-3-4_gradient-like} on the interval $[0, T)$ to provide
$$
\int_0^T\|\dot u(t)\|_\sH\,dt \leq \int_{H(T)}^{H(0)} \frac{1}{\phi(s)} \,ds + G(0) + \int_0^T |F(t)|^\theta \,dt
$$
where $H(t) = \sE(u(t)) + F(t)$ by \eqref{eq:Huang_3-18c} and $\phi(s) = c|s - \sE(\varphi)|^\theta$ by \eqref{eq:Simon_2-2_gradient_inequality_function_phi}, with $F$ and $G$ as in Hypothesis \ref{hyp:Huang_3-10}. Combining the preceding inequality with the \apriori interior estimate \eqref{eq:Abstract_apriori_interior_estimate_trajectory_with_small_constant_error} in Hypothesis \ref{hyp:Abstract_apriori_interior_estimate_trajectory_with_perturbation} for a solution $u$ to \eqref{eq:Simon_0-1} yields
\begin{align*}
\int_\delta^T \|\dot u(t)\|_\sX\,dt
&\leq C_1(1+\delta^{-\mu})\left(\int_{\sE(u(T)) + F(T)}^{\sE(u_0) + F(0)} \frac{1}{c(s - \sE(\varphi))^\theta} \,ds + G(0) + \int_0^T |F(t)|^\theta \,dt\right)
+ \eps_0
\\
&= \frac{C_1(1+\delta^{-\mu})}{c(1-\theta)}
\left((\sE(u_0) - \sE(\varphi) + F(0))^{1-\theta} - (\sE(u(T)) - \sE(\varphi) + F(T))^{1-\theta}\right.
\\
&\qquad + \left. G(0) + \int_0^T |F(t)|^\theta \,dt\right) + \eps_0,
\end{align*}
where to obtain the equality we use the fact that $\sE(v) \geq \sE(\varphi)$ for all $v \in \sU$ by our hypothesis that $\varphi$ is a ground state for $\sE$ on $\sU$ and, in particular, $\sE(u_0) \geq \sE(u(T)) \geq \sE(\varphi)$. By discarding the negative term in the preceding inequality and recalling that \eqref{eq:Rade_7-4_banach} yields, for $\|u_0-\varphi\|_\sX < \eps_\varphi$ and using the fact that $\sE:\sU \to \RR$ is $C^1$,
$$
|\sE(u_0) - \sE(\varphi)| \leq \|u_0 - \varphi\|_\sX,
$$
we obtain
\begin{align*}
\int_\delta^T \|\dot u(t)\|_\sX\,dt
&\leq
\frac{C_1(1+\delta^{-\mu})}{c(1-\theta)} \left( \left(\|u_0 - \varphi\|_\sX  + F(0) \right)^{1-\theta} + G(0)
+ \int_0^T |F(t)|^\theta \,dt \right)
+ \eps_0
\\
&< \frac{C_1(1+\delta^{-\mu})}{c(1-\theta)} (2\eps^{1-\theta} + \eps) + \eps,
\end{align*}
where, to obtain the last inequality, we apply our hypothesis that $u_0 \in \sU_\eps$, the definition \eqref{eq:Huang_page_164_Usigma} of $\sU_\eps$, the algebraic inequality \eqref{eq:Power_inequality_nonnegative_x_and_y}, and our hypotheses \eqref{eq:Huang_5-6} on $F(0)$, and $G(0)$ and $|F(t)|^\theta$, and $\eps_0$. Hence, for $\eps = \eps(c,C_1,\delta, \theta, \mu, \sigma, \tau, \varphi) \in (0, \sigma/4)$ small enough that $\eps \leq \eps_\varphi$ and
$$
\frac{C_1(1+\delta^{-\mu})}{c(1-\theta)} (2\eps^{1-\theta} + \eps) + \eps \leq \frac{\sigma}{8},
$$
we find that
$$
\int_\delta^T \|\dot u(t)\|_\sX\,dt < \frac{\sigma}{8},
$$
and thus,
\begin{equation}
\label{eq:Huang_theorem_5-1-1_proof_norm_X_uT_minus_udelta_inequality_gradientlike}
\|u(T)-u(\delta)\|_\sX < \frac{\sigma}{8}.
\end{equation}
Combining the inequalities \eqref{eq:Huang_theorem_5-1-1_proof_norm_X_udelta_minus_u0_inequality_gradientlike} and \eqref{eq:Huang_theorem_5-1-1_proof_norm_X_uT_minus_udelta_inequality_gradientlike} gives
$$
\|u(T)-u_0\|_\sX \leq \|u(\delta)-u_0\|_\sX + \|u(T)-u(\delta)\|_\sX < \frac{\sigma}{8} + \frac{\sigma}{8} = \frac{\sigma}{4},
$$
that is,
\begin{equation}
\label{eq:Huang_theorem_5-1-1_proof_u0_uT_X_inequality_gradientlike}
\|u_0 - u(T)\|_\sX < \frac{\sigma}{4},
\end{equation}
and thus\footnote{Here again we correct a small typographical error in Huang's proof of \cite[Theorem 5.1.1]{Huang_2006}, where \cite[Inequality (5.6)]{Huang_2006} only yields $\eps>0$ in the last displayed equation in \cite[p. 164]{Huang_2006} and not $\eps>\sigma/2$ as required to yield the desired contradiction.},
\begin{align*}
\eps &> \|u_0 - \varphi\|_\sX
\quad\hbox{(by definition \eqref{eq:Huang_page_164_Usigma} of $\sU_\eps$ and the fact that $u_0 \in \sU_\eps$)}
\\
&\geq \|u(T) - \varphi\|_\sX - \|u_0 - u(T)\|_\sX
\\
&> \frac{\sigma}{2} - \frac{\sigma}{4} = \frac{\sigma}{4} \quad \hbox{(by \eqref{eq:Huang_5-8_gradientlike} and \eqref{eq:Huang_theorem_5-1-1_proof_u0_uT_X_inequality_gradientlike})},
\end{align*}
contradicting the choice of $\eps \in (0, \sigma/4)$ in the hypotheses of Theorem \ref{thm:Huang_5-1-1_gradientlike}. This completes the proof of Claim \ref{claim:Equation_Huang_5-7_gradientlike}.
\end{proof}

By virtue of Claim \ref{claim:Equation_Huang_5-7_gradientlike}, we can apply the growth estimate \eqref{eq:Huang_3-18b_H_gradient-like} in Lemma \ref{lem:Huang_3-3-4_gradient-like} on the maximal interval $[0, \widetilde T)$ to provide (just as in the proof of Claim \ref{claim:Equation_Huang_5-7_gradientlike})
\begin{align*}
\int_\delta^{\widetilde T} \|\dot u(t)\|_\sX\,dt
&\leq
\frac{C_1(1+\delta^{-\mu})}{c(1-\theta)}
\left((\sE(u_0) - \sE(\varphi) + F(0))^{1-\theta} + G(0) + \int_0^{\widetilde T} |F(t)|^\theta \,dt\right)
+ \eps_0
\\
&< \infty,
\end{align*}
and so $u \in W^{1,1}(\delta, \widetilde T; \sX)$. Corollary \ref{cor:Neuberger_theorem_4-1} now implies that $\widetilde T = \infty$ and so the solution $u$ exists globally on $[0,\infty)$ and \eqref{eq:Huang_5-7_gradientlike} holds with $\widetilde T = \infty$, that is,
$$
\|u(t) - \varphi\|_\sX < \frac{\sigma}{2}, \quad\forall\, t \in [0,\infty).
$$
The orbit $O(u) = \{u(t): t\geq 0\}$ is contained in the open set $\sU_\sigma$ (actually, $\sU_{\sigma/2}$ by \eqref{eq:Huang_5-7_gradientlike} because $\widetilde T = \infty$) and by hypothesis the {\L}ojasiewicz-Simon gradient inequality \eqref{eq:Simon_2-2} holds on $\sU_\sigma$. Therefore, noting that the hypotheses \eqref{eq:Huang_5-6} required by Theorem \ref{thm:Huang_5-1-1_gradientlike} imply those of \eqref{eq:Huang_5-6_Tfinite} in Theorem \ref{thm:Huang_3-3-6_and_7_gradientlike}, we see that Theorem \ref{thm:Huang_3-3-6_and_7_gradientlike} provides the convergence of the integral $\int_1^\infty\|\dot u(t)\|_\sX\,dt$ and the convergence of $u(t)$ in the norm of $\sX$ to a limit $u_\infty \in \sX$.
\end{proof}

For the sake of clarity in the statement and proof of Theorem \ref{thm:Huang_5-1-1_gradientlike}, and consistency with \cite[Theorem 5.1.1]{Huang_2006}, we assumed that the critical point, $\varphi \in \sX$, was a minimum of $\sE$ on the open neighborhood $\sU \subset \sX$ of $\varphi$. However, as we saw in the statement and proof of Theorem \ref{thm:Huang_3-3-6_and_7_gradientlike}, it one can draw useful conclusions when $\varphi$ is only an `approximate local minimum' for a solution $u(t)$ on a maximal interval, $[0, \widetilde T)$, essentially as allowed by Simon in his \cite[Theorem 2]{Simon_1983}. The following result, a corollary of the proofs of Theorems \ref{thm:Huang_3-3-6_and_7_gradientlike} and \ref{thm:Huang_5-1-1_gradientlike}, provides such an extension.

\begin{cor}[Existence and convergence of a global solution to a gradient-like system near a critical point that is an approximate local minimum]
\label{cor:Huang_5-1-1_gradientlike_approximate_minimum}
Assume the hypotheses of Theorem \ref{thm:Huang_5-1-1_gradientlike}, but relax the requirement that $\varphi \in \sX$ is a ground state for $\sE$ on the open neighborhood $\sU \subset \sX$ of $\varphi$ to the requirement that $\varphi \in \sX$ is a critical point for $\sE$ on $\sU$. For each $u_0 \in \sU_\eps$, the following alternative holds. If $\widetilde T > 0$ is the maximal lifetime of a strong solution, $u$, to the Cauchy problem \eqref{eq:Simon_0-1} with $u(0) = u_0$, then either
\begin{enumerate}
\item
\label{item:Huang_5-1-1_gradientlike_approximate_minimum_u_at_time_t_below_energy_critical_point}
$\sE(u(t)) \leq \sE(\varphi) - \eps$ for some $t \in [0, \widetilde T)$, or
\item
\label{item:Huang_5-1-1_gradientlike_approximate_minimum_u_converges_to_limit_u_at_infty}
$\widetilde T = \infty$ and the trajectory $u(t)$ converges in $\sX$ to a limit $u_\infty \in \sX$ in the sense that
$$
\lim_{t\to\infty}\|u(t)-u_\infty\|_\sX =0
\quad\hbox{and}\quad
\int_1^\infty \|\dot u\|_\sX\,dt < \infty.
$$
If $\varphi$ is a cluster point of the orbit $O(u) = \{u(t): t\geq 0\}$, then $u_\infty = \varphi$.
\end{enumerate}
\end{cor}

\begin{proof}
We just indicate the changes to the proof of Theorem \ref{thm:Huang_5-1-1_gradientlike}. We assume that Alternative \eqref{item:Huang_5-1-1_gradientlike_approximate_minimum_u_at_time_t_below_energy_critical_point} does not occur and thus,
\begin{equation}
\label{eq:Simon_3-2_possibly_finite_maximal_lifetime}
\sE(u(t)) > \sE(\varphi) - \eps, \quad\forall\, t \in [0, \widetilde T).
\end{equation}
The proof of Claim \ref{claim:Equation_Huang_5-7_gradientlike} for $u$ on $[0,\widetilde T)$, under the hypotheses of Corollary \ref{cor:Huang_5-1-1_gradientlike_approximate_minimum} rather than Theorem \ref{thm:Huang_5-1-1_gradientlike}, is identical to that of Claim \ref{claim:norm_X_of_u_minus_varphi_lessthan_sigma_all_t_geq_T_gradientlike} on $[T, \infty)$ in the proof of Theorem \ref{thm:Huang_3-3-6_and_7_gradientlike}, except that the roles of $T$ and $t = \infty$ are replaced by those of $t = 0$ and $t = \widetilde T$, respectively. The remainder of the proof of Corollary \ref{cor:Huang_5-1-1_gradientlike_approximate_minimum} is now identical to that of Theorem \ref{thm:Huang_5-1-1_gradientlike}.
\end{proof}

\subsection{Uniform continuity of solutions to abstract gradient-like systems}
\label{subsec:Huang_5-1_uniform continuity}
Corollary \ref{cor:Huang_5-1-1_gradientlike_approximate_minimum} is close to the type of result we seek in connection with our application to the Yang-Mills gradient flow over closed four-dimensional manifolds, but its hypotheses on the perturbation term, $R$, and $F$ and $G$ (which are typically defined in terms of $R$ as discussed in Section \ref{sec:Huang_3_gradientlike_system}), are slightly stronger than we would ideally like since they are assumed --- through the hypotheses of Theorem \ref{thm:Huang_5-1-1_gradientlike} --- to obey certain conditions on $[0, \infty)$.

In our application, the perturbation term, $R$, and hence also $F$ and $G$, are only known to be defined \apriori on $[0, \widetilde T)$, where $\widetilde T$ is the maximal, but possibly finite lifetime of a strong solution, $u$, to the gradient-like system \eqref{eq:Simon_0-1}, in which case the expressions \eqref{eq:Huang_F_and_G_in_terms_of_Simon_R} for $F$ and $G$ in terms of $R$ on $[0, \infty)$ would be replaced by the expressions \eqref{eq:Huang_F_and_G_in_terms_of_Simon_R_finite_interval} (with $T$ replaced by $\widetilde T$) in terms of $R$ on $[0, \widetilde T)$.

On the other hand, as we know from Sections \ref{subsec:Huang_5-1_existence_and_convergence_global_solution_near_local_minimum} and \ref{subsec:Huang_5-1_existence_and_convergence_global_solution_gradientlike}, a crucial ingredient in our proofs of global existence for solutions to gradient or gradient-like systems, is to show that the solution $u$ is continuous on $[0, \widetilde T]$, not just $[0, \widetilde T)$, as continuity of $u(t)$ at $t = \widetilde T$ allows to continue the flow for $t > \widetilde T$ via short-time existence of solutions to the gradient or gradient-like system. In other words, we can conclude that $u(t)$ is defined for $t\in [0, \infty)$ and $u \in C([0, \infty); \sX)$.
Therefore, keeping in mind the preceding remarks, we have the following refinement of Corollary \ref{cor:Huang_5-1-1_gradientlike_approximate_minimum}.

\begin{thm}[Uniform continuity of solutions to gradient-like systems near a critical point that is an approximate local minimum]
\label{thm:Huang_5-1-1_gradientlike_uniform_continuity}
Let $\sU$ be an open subset of a real Banach space, $\sX$, that is continuously embedded and dense in a Hilbert space, $\sH$. Let $\sE:\sU\subset \sX\to\RR$ be a $C^1$ functional with gradient map $\sE':\sU\subset \sX \to \sH$, let $\varphi \in \sU$ be a critical point of $\sE$, and suppose that $\sE'$ obeys a {\L}ojasiewicz-Simon gradient inequality \eqref{eq:Simon_2-2} with positive constants $c$, $\sigma$, and $\theta \in [1/2, 1)$.
Assume that
\begin{enumerate}
\item For each $u_0 \in \sU$, there exists a unique strong solution to the Cauchy problem \eqref{eq:Simon_0-1}, in the sense of Definition \ref{defn:Strong_solution_to_gradientlike_system}, on a time interval $[0, \tau)$ for some positive constant, $\tau$;

\item Hypothesis \ref{hyp:Abstract_apriori_interior_estimate_trajectory_with_perturbation} holds with positive constants $C_1$ and $\mu$ for solutions to the gradient-like system \eqref{eq:Simon_0-1}; and

\item Given positive constants $b$ and $\eta$, there is a constant $\delta = \delta(\eta, \tau, b) \in (0, \tau]$ such that if $v$ is a solution to the gradient-like system \eqref{eq:Simon_0-1} on $[0, \tau)$ with $\|v_0\|_\sX \leq b$, then
\begin{equation}
\label{eq:Gradientlike_solution_near_initial_data_at_time_zero_for_short_enough_time_uniform continuity}
\sup_{t\in [0, \delta]}\|v(t) - v_0\|_\sX < \eta.
\end{equation}
\end{enumerate}
Then there is a constant $\eps = \eps(c,C_1,\delta, \theta, \mu, \sigma, \tau, \varphi) \in (0, \sigma/4)$ with the following significance. If $u_0 \in \sU_\eps$ and $u(t)$, for $t \in [0,\widetilde T)$, is a strong solution to the Cauchy problem \eqref{eq:Simon_0-1} with $u(0) = u_0$ and
\begin{equation}
\label{eq:Simon_3-2_possibly_finite_maximal_lifetime_uniform continuity}
\sE(u(t)) > \sE(\varphi) - \eps, \quad\forall\, t \in [0, \widetilde T),
\end{equation}
and $F$ and $G$ in Hypothesis \ref{hyp:Huang_3-10} and $\eps_0$ in Hypothesis \ref{hyp:Abstract_apriori_interior_estimate_trajectory_with_perturbation} are defined on $[0, \widetilde T)$ and obey
\begin{subequations}
\label{eq:Huang_5-6_uniform continuity}
\begin{align}
\label{eq:Huang_5-6_F0_lessthan_epsilon_uniform continuity}
F(0) &< \eps,
\\
\label{eq:Huang_5-6_G0_plus_integral_0_to_infinity_phiF_lessthan_epsilon_uniform continuity}
G(0) + \int_0^{\widetilde T} |F(t)|^\theta \,dt &< \eps,
\\
\label{eq:Huang_5-6_integral_0_to_infinity_dotR_norm_H_lessthan_epsilon_uniform continuity}
\eps_0 &< \eps,
\end{align}
\end{subequations}
where $\phi(x) = c|x - \sE(\varphi)|^\theta$ as in \eqref{eq:Simon_2-2_gradient_inequality_function_phi}, for all $x \in \RR$, then
$$
\sup_{t \in [0, \widetilde T)}\|u(t) - \varphi\|_\sX < \frac{\sigma}{2} \quad\hbox{and}\quad u \in C([0, \widetilde T]; \sX).
$$
\end{thm}

\begin{proof}
It is enough to observe that the proof of Corollary \ref{cor:Huang_5-1-1_gradientlike_approximate_minimum} only requires the conditions \eqref{eq:Huang_5-6_G0_plus_integral_0_to_infinity_phiF_lessthan_epsilon_uniform continuity} and \eqref{eq:Huang_5-6_integral_0_to_infinity_dotR_norm_H_lessthan_epsilon_uniform continuity}, for $F$, $G$, and $\eps_0$ on $[0, \widetilde T)$, respectively, rather than the conditions \eqref{eq:Huang_5-6_G0_plus_integral_0_to_infinity_phiF_lessthan_epsilon} and \eqref{eq:Huang_5-6_integral_0_to_infinity_dotR_norm_H_lessthan_epsilon}, for $F$, $G$, and $\eps_0$ on $[0, \infty)$. More specifically, see the proof of Theorem \ref{thm:Huang_5-1-1_gradientlike}, especially Claim \ref{claim:Equation_Huang_5-7_gradientlike}, in the case that $\varphi$ is a ground state and, more generally when $\varphi$ is a critical point that is only an approximate minimum in the sense of \eqref{eq:Simon_3-2_possibly_finite_maximal_lifetime_uniform continuity}, see the proof of Corollary \ref{cor:Huang_5-1-1_gradientlike_approximate_minimum}, especially Claim \ref{claim:norm_X_of_u_minus_varphi_lessthan_sigma_all_t_geq_T_gradientlike}, but for the interval $[0,\widetilde T)$ rather than $[T, \infty)$. The proof of Corollary \ref{cor:Huang_5-1-1_gradientlike_approximate_minimum} then shows that $u$ satisfies $\sup_{t \in [0, \widetilde T)}\|u(t) - \varphi\|_\sX < \sigma/2$ and $u \in W^{1,1}(0, \widetilde T); \sX)$, and hence $u \in C([0, \widetilde T]; \sX)$ by Theorem \ref{thm:Sell_You_C-9}.
\end{proof}

\subsection{Further refinements}
\label{sec:Huang_3_and_5_refined_hypotheses}
The results of Sections \ref{sec:Huang_3_gradientlike_system}, \ref{subsec:Huang_3-3_gradientlike_system}, \ref{subsec:Huang_5-1_existence_and_convergence_global_solution_gradientlike}, and \ref{subsec:Huang_5-1_uniform continuity} may be extended with very minor modifications to their proofs and statements by replacing the hypotheses with one or more of the following variants.

In Hypothesis \ref{hyp:Huang_3-10}, the conditions \eqref{eq:Huang_3-10a} may be replaced by the variants
\begin{subequations}
\label{eq:Huang_3-10a_prime}
\begin{align}
\label{eq:Huang_3-10a_gradE_innerproduct_dotu_prime}
\tag{\ref*{eq:Huang_3-10a_gradE_innerproduct_dotu}$'$}
(-\sE'(u(t)), \dot u(t))_\sH &\geq \alpha_0\|\sE'(u(t))\|_\sH \|\dot u(t)\|_\sH + F'(t),
\\
\label{eq:Huang_3-10a_gradE_norm_H_prime}
\tag{\ref*{eq:Huang_3-10a_gradE_norm_H}$'$}
\|\sE'(u(t))\|_\sH &\geq \alpha_1\|\dot u(t)\|_\sH + G'(t),  \quad\hbox{for a.e. } t \geq 0,
\end{align}
\end{subequations}
for positive constants $\alpha_0$ and $\alpha_1$; compare \cite[Equation (3.10a)]{Huang_2006}.

Furthermore, in Hypothesis \ref{hyp:Huang_3-10}, it is possible to replace the condition \eqref{eq:Huang_3-10b_integral_zero_to_infinity_of_phiF_finite} by the weaker variant,
\begin{equation}
\label{eq:Huang_3-10b_integral_one_over_rho_plus_rho_times_Fprime}
\tag{\ref*{eq:Huang_3-10b_integral_zero_to_infinity_of_phiF_finite}$'$}
\int_0^\infty \left(\frac{1}{\rho(t)} + \rho(t)|F'(t)|\right)\,dt < \infty,
\end{equation}
for \emph{some} increasing, positive function, $\rho:[0, \infty) \to (0, \infty)$. See, for example, \cite[Theorem 3.3.7]{Huang_2006}.

The convergence and stability results in \cite[Sections 3.3 and 5.1]{Huang_2006} follow not only from the {\L}ojasiewicz-Simon gradient inequality \eqref{eq:Simon_2-2}, but more generally from a far weaker gradient inequality of the form,
\begin{equation}
\label{eq:Huang_3-13a}
\|\sE'(v)\|_\sH \geq \phi(\sE(v)), \quad\forall\, v \in \sU,
\end{equation}
where $\phi$ is a function $\phi$ belonging to the class $\sG$ (see \cite[Definition 2.2.1]{Huang_2006}), that is, $\phi:\RR \to [0, \infty)$ is a Lebesgue-measurable function such that $1/\phi \in L^1_{\loc}(\RR)$ --- see \cite[Propositions 3.3.1 and 3.3.2 and Theorems 3.3.3 and 5.1.2]{Huang_2006}. The key results in \cite[Sections 3.3 and 5.1]{Huang_2006} require in addition that $\phi$ belong to the class $\sG_k$ for some constant $k \geq 1$ (see \cite[Definition 2.2.1 (i)]{Huang_2006}), that is $\phi \in \sG$ and satisfies the inequality,
\begin{equation}
\label{eq:Huang_definition_2-1-i}
\phi(x + y) \leq \phi(x) + k|y|^{1/2}, \quad\forall\, x, y \in \RR.
\end{equation}
See \cite[Lemma 3.3.4 and Theorems 3.3.5, 3.3.6, and 5.1.1]{Huang_2006} for results of this kind. We could also have derived the results in Sections \ref{sec:Huang_3_gradientlike_system}, \ref{subsec:Huang_3-3_gradientlike_system}, \ref{subsec:Huang_5-1_existence_and_convergence_global_solution_gradientlike}, and \ref{subsec:Huang_5-1_uniform continuity} in a setting of this generality; the required modifications to the statements and proofs of these results are again very minor.

\section{Local a priori estimates for Yang-Mills heat and gradient flows}
\label{sec:Local_apriori_estimates_Yang-Mills_heat_and_gradient_flows}
In order to apply the results of Section \ref{sec:Huang_3_and_5_gradientlike_system} to a Yang-Mills gradient-like flow, we shall require \apriori estimates for solutions to the Yang-Mills heat and gradient flow equations over an annulus in the base four-dimensional manifold, $X$, when the energy (or the $L^2$ norm of the curvature) of the flow, restricted to that annulus, is sufficiently small. These local-in-space estimates will be used to show that the perturbation, $R(t)$, of pure Yang-Mills gradient flow over $S^4$, that arises when cutting off a pure Yang-Mills gradient flow over an annulus in $X$ and then transferring the resulting flow (over the annulus and the ball it encloses) to $S^4$ (with its standard round metric of radius one) obeys conditions similar to those in Hypothesis \ref{hyp:Huang_3-10}. Fortunately, \apriori estimates of this kind were developed independently by Kozono, Maeda, and Naito in \cite{Kozono_Maeda_Naito_1995} and by Schlatter in \cite{Schlatter_1997}, based in part on predictions of Struwe in \cite{Struwe_1994}. (Some of these estimates may be viewed as local-in-space analogues of the global-in-space \apriori estimate in Lemma \ref{lem:Rade_8-1}.) There are related results due to Chen and Shen in \cite{Chen_Shen_1994}, though they assume that the dimension, $d$, of $X$ is greater than four. In this section, we review the results we shall need from \cite{Kozono_Maeda_Naito_1995, Schlatter_1997} and then develop some additional useful consequences. Many of the local \apriori estimates in \cite[Section 4]{Kozono_Maeda_Naito_1995} which we review here are analogues of earlier results due to Struwe \cite[Section 3]{Struwe_1985} for a harmonic map evolution equation \cite[Equation (1.7)]{Struwe_1985}, a gradient-like flow equation for the energy functional \cite[Equations (1.1) and (1.2)]{Struwe_1985} associated to a map from a Riemann surface into a Riemannian manifold. In the context of the anti-self-dual equation over a long tube, $Y\times [0,\infty)$ where $Y$ is a closed, three-dimensional manifold, results of this kind are due to Morgan, Mrowka, and Ruberman \cite{MMR} and Taubes \cite{TauL2}, keeping in mind that a solution to the anti-self-dual equation over $Y\times [0,\infty)$ with bounded energy may be viewed as the gradient flow for the Chern-Simons functional.

\subsection{Fundamental inequalities}
\label{subsec:Kozono_Maeda_Naito_2}
In this subsection, we review several fundamental \apriori estimates either proved or identified by Kozono, Maeda, and Naito  and which appear inspired in turn by results of Struwe \cite[Section 3]{Struwe_1985}, based on embedding theorems for parabolic Sobolev spaces due to Lady{\v{z}}enskaja, Solonnikov, and Ural$'$ceva \cite{LadyzenskajaSolonnikovUralceva}. Throughout this subsection, $X$ will denote a closed, connected smooth manifold of dimension $d \geq 2$, Riemannian metric, $g$, and covariant derivative $\nabla: C^\infty(TX) \to C^\infty(TX\otimes T^*X)$ defined by the Levi-Civita connection associated to $g$.

\begin{prop}
\label{prop:Kozono_Maeda_Naito_2-1}
\cite[Proposition 2.1]{Kozono_Maeda_Naito_1995}
Let $(X, g)$ be a closed, Riemannian, smooth manifold of dimension $d \geq 2$. Then there exist positive constants $C$ and $R_0$ such that, for any $u, v \in W^{1,2}(X)$ and $r \in (0, R_0]$, we have
$$
\int_X |u| |v|^2 \,d\vol_g
\leq
C\sup_{x\in X}\left( \int_{B_r(x)} |u|^2 \,d\vol_g \right)^{1/2}
\left( \int_X |\nabla v|^2 \,d\vol_g + r^{-2} \int_X |v|^2 \,d\vol_g \right).
$$
\end{prop}

As noted by Kozono, Maeda, and Naito, Proposition \ref{prop:Kozono_Maeda_Naito_2-1} follows in part from the following local version.

\begin{lem}
\label{lem:Kozono_Maeda_Naito_2-2}
\cite[Lemma 2.2]{Kozono_Maeda_Naito_1995}
Let $(X, g)$ be a closed, Riemannian, smooth manifold of dimension $d \geq 2$. Then there exist positive constants $C$ and $R_0$ such that, for any $u, v \in W^{1,2}(X)$ and $r \in (0, R_0]$ and non-increasing, non-negative function $\varphi \in C^\infty(\RR)$ such that $\varphi(s) = 0$ for $s \geq 1$, we have
\begin{multline*}
\int_X |u| |v|^2 \varphi(\dist_g(\cdot, x)/r) \,d\vol_g
\\
\leq C\left( \int_{B_r(x)} |u|^2 \,d\vol_g \right)^{1/2}
\int_X \left( |\nabla v|^2 + r^{-2} |v|^2  \right)\varphi(\dist_g(\cdot, x)/r) \,d\vol_g.
\end{multline*}
\end{lem}

Proposition \ref{prop:Kozono_Maeda_Naito_2-1} is derived from Lemma \ref{lem:Kozono_Maeda_Naito_2-2} via the following lemma proved by Struwe in \cite{Struwe_1985}.

\begin{lem}[Covering lemma]
\label{lem:Kozono_Maeda_Naito_2-3}
\cite[Lemma 2.3]{Kozono_Maeda_Naito_1995}, \cite[Lemma 3.3]{Struwe_1985}
Let $(X, g)$ be a closed, Riemannian, smooth manifold of dimension $d \geq 2$. Then there exist positive constants $K$ and $R_0$, depending only on $(X,g)$, such that for any $r \in (0, R_0]$, there exists a covering of $X$ by geodesic balls $B_{r/2}(x_i)$ such that, for any point $x \in X$, at most $K$ of the balls $B_r(x_i)$ contain $x$.
\end{lem}

We record the following useful pointwise inequalities for the curvature $F_A(t)$ of a family of smooth connections $A(t)$ on a principal $G$-bundle $P$ over $X$ and which depend smoothly on $t \in \RR$.
We begin with the following lemma, whose proof we include since we use a different Bochner-Weitzenb\"ock formula than that of \cite{Kozono_Maeda_Naito_1995} and so our pointwise equalities and inequalities contain extra terms due to our inclusion of the Riemann curvature tensor, $\Riem_g$. The Bianchi identity \cite[Equation (2.16)]{FU},
\begin{equation}
\label{eq:Freed_Uhlenbeck_2-16_Bianchi_identity}
d_A F_A = 0 \quad\hbox{on } X,
\end{equation}
valid for any connection $A$ on a principal $G$ bundle $P$ over a manifold $X$ of any dimension $d \geq 2$, is frequently employed (at least implicitly) in the proof of Lemma \ref{lem:Kozono_Maeda_Naito_2-4} below.

\begin{lem}[Equalities and inequalities for smooth solutions to Yang-Mills gradient flow]
\label{lem:Kozono_Maeda_Naito_2-4}
\cite[Lemma 2.4]{Kozono_Maeda_Naito_1995}
Let $(X, g)$ be a closed, Riemannian, smooth manifold of dimension $d \geq 2$. Then there is a positive constant $C$ and, given a non-negative integer $k$, a constant $C_k$ with the following significance. Let $P$ be a principal $G$-bundle over $X$, let $t_1, t_2 \in \RR$ obey $-\infty \leq t_1 < t_2 \leq \infty$, and let $A(t)$ be a solution to the Yang-Mills gradient flow equation \eqref{eq:Yang-Mills_gradient_flow_equation} over $X$ that is $C^\infty$ in time and space, namely,
$$
\frac{\partial A}{\partial t}(t) = -d_{A(t)}^* F_{A(t)}, \quad \forall\,t \in (t_1, t_2) \subset \RR.
$$
Then the following pointwise equalities and inequalities hold on $(t_1, t_2) \times X$:
\begin{align}
\label{eq:Kozono_Maeda_Naito_2-4}
\frac{\partial F_A}{\partial t} &= -\Delta_A F_A,
\\
\label{eq:Kozono_Maeda_Naito_2-5}
\frac{\partial F_A}{\partial t} &= -\nabla_A^*\nabla_A F_A - \Riem_g\times F_A - F_A\times F_A,
\\
\label{eq:Kozono_Maeda_Naito_2-6}
\frac{\partial |F_A|}{\partial t} &\leq \Delta|F_A| + C(|F_A| + |F_A|^2) \quad\hbox{a.e. on } (t_1, t_2)\times X,
\\
\label{eq:Kozono_Maeda_Naito_2-7}
\frac{\partial |\nabla_A^k F_A|}{\partial t} &\leq \Delta|\nabla_A^k F_A|
+ C_k\sum_{j=0}^k \left( |\nabla_A^j F_A| + |\nabla_g^j \Riem_g| \right) |\nabla_A^{k-j}F_A|
\quad\hbox{a.e. on } (t_1, t_2)\times X,
\end{align}
where $\Delta_A = d_A^*d_A + d_Ad_A^*$ is the Hodge Laplace operator \eqref{eq:Lawson_page_93_Hodge_Laplacian} on $\Omega^2(X; \ad P)$ defined by a connection $A$ on $P$ and Riemannian $g$ on $X$, and $\Delta = d^*d$ is the Laplace operator on $C^\infty(X)$ defined by the Riemannian $g$ on $X$.
\end{lem}

\begin{proof}
Equation \eqref{eq:Kozono_Maeda_Naito_2-4} is well-known \cite{DK, Kozono_Maeda_Naito_1995, Rade_1992, Struwe_1994} and follows from the Equation \eqref{eq:Yang-Mills_gradient_flow_equation} for Yang-Mills gradient flow and the Bianchi identity \eqref{eq:Freed_Uhlenbeck_2-16_Bianchi_identity}. We have modified the expressions on the right-hand side of \cite[Equations (2.5), (2.6), and (2.7)]{Kozono_Maeda_Naito_1995} since the terms in the Bochner-Weitzenb\"ock formula \eqref{eq:Lawson_corollary_II-3} involving Riemann curvature tensor, $\Riem_g$, are omitted in \cite{Kozono_Maeda_Naito_1995}. Formula \eqref{eq:Lawson_corollary_II-3} and Equation \eqref{eq:Kozono_Maeda_Naito_2-4} immediately yield Equation \eqref{eq:Kozono_Maeda_Naito_2-5}, while the second-order Kato Inequality \cite[Equation (6.21)]{FU},
\begin{equation}
\label{eq:FU_6-21_second-order_Kato_inequality}
\langle \nabla^*\nabla v, v \rangle \geq |v|\Delta |v| \quad\hbox{a.e. on } X,
\end{equation}
for any $v \in C^2(X; E)$ and Riemannian or Hermitian vector bundle $E$ over $X$ and compatible covariant derivative, $\nabla$, gives
$$
\langle \nabla_A^*\nabla_A F_A, F_A \rangle \geq |F_A|\Delta|F_A| \quad\hbox{a.e. on } (0,T)\times X,
$$
and thus \eqref{eq:Kozono_Maeda_Naito_2-5} implies (compare \cite[page 21]{DK})
$$
|F_A|\frac{\partial |F_A|}{\partial t}
=
\frac{1}{2}\frac{\partial |F_A|^2}{\partial t}
=
\left\langle\frac{\partial F_A}{\partial t}, F_A\right\rangle
\leq
|F_A||\Delta|F_A|| + C(|F_A|^2 + |F_A|^3) \quad\hbox{a.e. on } (t_1, t_2)\times X.
$$
Therefore, on the complement of the set, $|F_A|=0$ in $(t_1, t_2)\times X$, we obtain the Inequality \eqref{eq:Kozono_Maeda_Naito_2-6}. Next, we have
\begin{align*}
\frac{\partial \nabla_A F_A}{\partial t}&= \nabla_A \frac{\partial F_A}{\partial t}
+ \dot A \times F_A
\\
&= \nabla_A \frac{\partial F_A}{\partial t}
- d_A^* F_A \times F_A \quad\hbox{(by \eqref{eq:Yang-Mills_gradient_flow_equation})}
\\
&= -\nabla_A\nabla_A^*\nabla_A F_A - \nabla_A(\Riem_g\times F_A) - \nabla_A(F_A\times F_A)
- d_A^* F_A \times F_A \quad\hbox{(by \eqref{eq:Kozono_Maeda_Naito_2-5})}
\\
&= -\nabla_A^*\nabla_A \nabla_A F_A - \nabla_g\Riem_g\times F_A - \Riem_g\times \nabla_A F_A
- \nabla_A F_A\times F_A.
\end{align*}
By iterating once again and substituting the preceding expression for $\partial\nabla_A F_A/\partial t$ and the expression for $\partial A/\partial t$ given by Yang-Mills gradient flow \eqref{eq:Yang-Mills_gradient_flow_equation}, we obtain
\begin{align*}
\frac{\partial \nabla_A^2 F_A}{\partial t}&= \nabla_A \frac{\partial \nabla_A F_A}{\partial t}
+ \dot A \times \nabla_A F_A
\\
&= -\nabla_A^*\nabla_A \nabla_A^2 F_A - \nabla_g^2\Riem_g\times F_A - \nabla_g\Riem_g\times \nabla_A F_A - \Riem_g\times \nabla_A^2 F_A
\\
&\qquad - \nabla_A^2 F_A\times F_A - \nabla_A F_A\times \nabla_AF_A.
\end{align*}
By iteration, we see that for any integer $k \geq 0$,
\begin{equation}
\label{eq:Kozono_Maeda_Naito_2-8}
\frac{\partial \nabla_A^k F_A}{\partial t} = - \nabla_A^*\nabla_A \nabla_A^k F_A
+ \sum_{j=0}^k \nabla_A^j F_A \times \nabla_A^{k-j} F_A
+ \sum_{j=0}^k \nabla_g^j \Riem_g \times \nabla_A^{k-j} F_A,
\end{equation}
where we have absorbed all signs (except for the term involving the connection Laplacian, $\nabla_A^*\nabla_A$) into the definition of `$\times$' in the preceding expression. It is now straightforward to verify the preceding formula by induction on $k$. Note that our Equation \eqref{eq:Kozono_Maeda_Naito_2-8} slightly from the corresponding \cite[Equation (2.8)]{Kozono_Maeda_Naito_1995} since Kozono, Maeda, and Naito omit the `$\Riem_g$' term in their version of the Bochner-Weitzenb\"ock formula.

Finally, the inequality \eqref{eq:Kozono_Maeda_Naito_2-7} follows from \eqref{eq:Kozono_Maeda_Naito_2-8} by same argument which gave the inequality \eqref{eq:Kozono_Maeda_Naito_2-6} as a consequence of \eqref{eq:Kozono_Maeda_Naito_2-5}.
\end{proof}

\subsection{Local \apriori estimates for a solution to the Yang-Mills gradient flow equation over a ball}
\label{subsec:Kozono_Maeda_Naito_4}
In this subsection, we review the results due to Kozono, Maeda, and Naito \cite[Section 4]{Kozono_Maeda_Naito_1995} which we shall need in our application. Again, it is worth mentioning that many of these ideas owe their origin in part to earlier work of Struwe on solutions to harmonic-map gradient-like flow of a map from a Riemann surface into a Riemannian manifold  \cite[Section 3]{Struwe_1985} and of Uhlenbeck in the case of solutions to the Yang-Mills equation \cite{UhlLp, UhlRem}.

\begin{lem}[Energy decreasing property for a solution to Yang-Mills gradient flow]
\label{lem:Kozono_Maeda_Naito_4-1}
\cite[Lemma 4.1]{Kozono_Maeda_Naito_1995}
Let $P$ be a principal $G$-bundle over a closed, Riemannian, smooth manifold of dimension $d \geq 2$, let $t_1, t_2 \in \RR$ obey $-\infty \leq t_1 < t_2 \leq \infty$, and let $A(t)$ be a solution to the Yang-Mills gradient flow equation \eqref{eq:Yang-Mills_gradient_flow_equation} over $X$. Then the energy function,
\begin{equation}
\label{eq:Kozono_Maeda_Naito_page_108_energy_of_a_connection}
\sE(A(t)) := \frac{1}{2}\int_X |F_A(t)|^2\,d\vol_g,
\end{equation}
is non-increasing the function of time $t \in (t_1, t_2)$.
\end{lem}

\begin{proof}
This is a general property of a solution to a gradient flow equation --- see the proof of \cite[Proposition 3.2]{Huang_2006}.
\end{proof}

By analogy\footnote{We omit the square root for convenience and also for the sake of prioritizing consistency with the proof of \cite[Lemma 5.2 and Theorem 5.4]{Kozono_Maeda_Naito_1995}, which conflicts with the same authors' usage elsewhere in \cite[Sections 4 and 5]{Kozono_Maeda_Naito_1995}.}
with \cite[Equation (4.1)]{Kozono_Maeda_Naito_1995} and \cite[page 571]{Struwe_1985}, for any point $x \in X$ and positive constant, $r$, less than the injectivity radius of $(X,g)$ at $x$, we define a measure of the energy of family of connections, $A(t)$ on a principal $G$-bundle $P\restriction B_r(x) \subset X$, uniform with respect to $t \in [0, T)$,
\begin{equation}
\label{eq:Kozono_Maeda_Naito_4-1}
\eps(r,x) = \sup_{t \in [0, T)} \int_{B_r(x)} |F_A(t)|^2 \,d\vol_g.
\end{equation}
We note that
$$
\eps(r,x) = \eps(r,x; A, T) = \sup_{t \in [0, T)}\sE_r(A(t);x),
$$
where
$$
\sE_r(A;x) := \int_{B_r(x)} |F_A|^2 \,d\vol_g
$$
is the local energy of a connection, $A$ on $P\restriction B_r(x)$, defined by analogy with the definition of the local energy of a map from a Riemann surface into a Riemannian manifold in \cite[page 561]{Struwe_1985}.

In the sequel, we review the local-in-space \apriori estimates due to Kozono, Maeda, and Naito for norms of a solution, $A(t)$ for $t\in [0,T)$, to Yang-Mills gradient flow on $P$ in terms of the initial energy, $\sE(A_0)$, $T$ and $\eps_1$. Here, $\eps_1$ is a positive constant depending only on $(X,g)$ and which will be determined in Lemmata \ref{lem:Kozono_Maeda_Naito_4-2} through \ref{lem:Kozono_Maeda_Naito_4-9}; we will set $\eps_1$ to be the smallest of the positive constants labeled $\eps_1$ appearing in these lemmata. Lemmata \ref{lem:Kozono_Maeda_Naito_4-3} through \ref{lem:Kozono_Maeda_Naito_4-9} may be compared with their analogues for a solution to harmonic-map gradient-like flow of a map from a Riemann surface into a Riemannian manifold  \cite[Lemmata 3.7$'$ -- 3.10$'$]{Struwe_1985}. We begin with the

\begin{lem}[Local-in-space \apriori estimate for the curvature of a solution to Yang-Mills gradient flow]
\label{lem:Kozono_Maeda_Naito_4-2}
\cite[Lemma 4.2]{Kozono_Maeda_Naito_1995}, \cite[Lemma 2.2]{Schlatter_1997}
Let $P$ be a principal $G$-bundle over a closed, four-dimensional, Riemannian, smooth manifold, $X$. Then there exist positive constants, $C$ and $\eps_1$, with the following significance. Let $R_0$ be as in Proposition \ref{prop:Kozono_Maeda_Naito_2-1}. If $A(t)$ is a smooth solution to the Yang-Mills gradient flow equation \eqref{eq:Yang-Mills_gradient_flow_equation} over $(0, T) \times X$ with initial value $A_0$ of class $H^1$, and $r \in (0, R_0]$ and $x_0 \in X$ obey
$$
\eps(r,x_0) < \eps_1,
$$
then
\begin{equation}
\label{eq:Kozono_Maeda_Naito_4-2}
\int_0^T \int_{B_{r/2}(x_0)} |\nabla_A F_A(t)|^2 \,d\vol_g\,dt \leq C(1 + r^{-2}T)\sE(A_0).
\end{equation}
\end{lem}

Lemma \ref{lem:Kozono_Maeda_Naito_4-2} is a local analogue of Lemma \ref{lem:Rade_8-1} when $X$ has dimension four.

\begin{lem}[Local-in-space \apriori estimate for the curvature of a solution to Yang-Mills gradient flow]
\label{lem:Kozono_Maeda_Naito_4-3}
\cite[Lemma 4.3]{Kozono_Maeda_Naito_1995}
Let $P$ be a principal $G$-bundle over a closed, four-dimensional, Riemannian, smooth manifold, $X$. Then there exists a positive constant, $\eps_1$, with the following significance. Given $\sE_0> 0$, and $R_0$ be as in Proposition \ref{prop:Kozono_Maeda_Naito_2-1}, and $r \in (0, R_0]$, and $T>0$, and $\delta \in (0, T)$, there exists a positive constant, $C$, such that the following holds. If $A(t)$ is a smooth solution to the Yang-Mills gradient flow equation \eqref{eq:Yang-Mills_gradient_flow_equation} over $(0, T) \times X$ with initial value $A_0$ of class $H^1$ and $\sE(A_0) \leq \sE_0$, and $r \in (0, R_0]$ and $x_0 \in X$ obey
$$
\eps(r,x_0) < \eps_1,
$$
then
\begin{equation}
\label{eq:Kozono_Maeda_Naito_4-5}
\sup_{\begin{subarray}{c} \delta < t < T \\ x \in B_{r/2}(x_0) \end{subarray}} |F_A(t,x)| \leq C.
\end{equation}
\end{lem}

\begin{lem}[Local-in-space \apriori estimate for the curvature of a solution to Yang-Mills gradient flow]
\label{lem:Kozono_Maeda_Naito_4-4}
\cite[Lemma 4.4]{Kozono_Maeda_Naito_1995}
Assume the hypotheses of Lemma \ref{lem:Kozono_Maeda_Naito_4-3}. Then
\begin{equation}
\label{eq:Kozono_Maeda_Naito_4-10}
\sup_{\begin{subarray}{c} \delta < t < T \\ x \in B_{r/2}(x_0) \end{subarray}} |\nabla_AF_A(t,x)| \leq C.
\end{equation}
\end{lem}

\begin{lem}[Local-in-space \apriori estimate for the curvature of a solution to Yang-Mills gradient flow]
\label{lem:Kozono_Maeda_Naito_4-5}
\cite[Lemma 4.5]{Kozono_Maeda_Naito_1995}
Assume the hypotheses of Lemma \ref{lem:Kozono_Maeda_Naito_4-3} and, in addition, that $k \geq 2$ is an integer. Then
$$
\int_\delta^T \int_{B_{r/2}(x_0)} |\nabla_A^k F_A(t)|^2 \,d\vol_g\,dt \leq C,
$$
where $C$ has the same dependencies as those in Lemma \ref{lem:Kozono_Maeda_Naito_4-3} but now depends in addition on $k$ and $\|F_A(\delta)\|_{W^{k,2}(B_r(x_0))}$.
\end{lem}

\begin{lem}[Local-in-space \apriori estimate for the curvature of a solution to Yang-Mills gradient flow]
\label{lem:Kozono_Maeda_Naito_4-6}
\cite[Lemma 4.6]{Kozono_Maeda_Naito_1995}
Assume the hypotheses of Lemma \ref{lem:Kozono_Maeda_Naito_4-3} and, in addition, that $k \geq 2$ is an integer. Then
$$
\sup_{\begin{subarray}{c} \delta < t < T \\ x \in B_{r/2}(x_0) \end{subarray}} |\nabla_A^kF_A(t,x)| \leq C,
$$
where $C$ has the same dependencies as those in Lemma \ref{lem:Kozono_Maeda_Naito_4-5}.
\end{lem}

\begin{lem}[Local-in-space \apriori estimate for a solution to Yang-Mills gradient flow]
\label{lem:Kozono_Maeda_Naito_4-7}
\cite[Lemma 4.7]{Kozono_Maeda_Naito_1995}
Assume the hypotheses of Lemma \ref{lem:Kozono_Maeda_Naito_4-3} and, in addition, that $p \geq 2$. Then
$$
\sup_{\delta < t < T} \int_{B_{r/2}(x_0)} |A(t) - \Gamma|^p \,d\vol_g \leq C,
$$
where $C$ has the same dependencies as those in Lemma \ref{lem:Kozono_Maeda_Naito_4-3} but now depends in addition on $p$ and\footnote{We correct a typographical error in the statement of \cite[Lemma 4.7]{Kozono_Maeda_Naito_1995}: as is clear from the proof, the dependency is on $\|F_A(\delta)\|_{W^{1,p}(B_r(x_0))}$ and not $\|F_A(\delta)\|_{W^{2,n}(B_r(x_0))}$.}
$\|F_A(\delta)\|_{W^{1,p}(B_r(x_0))}$ and $\|A(\delta) - \Gamma\|_{L^p(B_{r/2}(x_0))}$, where $\Gamma$ is the product connection on $B_{r/2}(x_0) \times G$.
\end{lem}

\begin{lem}[Local-in-space \apriori estimate for a solution to Yang-Mills gradient flow]
\label{lem:Kozono_Maeda_Naito_4-8}
\cite[Lemma 4.8]{Kozono_Maeda_Naito_1995}
Assume the hypotheses of Lemma \ref{lem:Kozono_Maeda_Naito_4-3} and, in addition, that $k \geq 2$ is an integer and $p \geq 2$. Then
$$
\sup_{\delta < t < T} \int_{B_{r/2}(x_0)} |\nabla_\Gamma^k(A(t) - \Gamma)|^p \,d\vol_g \leq C,
$$
where $C$ has the same dependencies as those in Lemma \ref{lem:Kozono_Maeda_Naito_4-3} but now depends in addition on $k$ and $p$ and $\|F_A(\delta)\|_{W^{k+1,2}(B_r(x_0))}$ and $\|A(\delta) - \Gamma\|_{L^p(B_{r/2}(x_0))}$.
\end{lem}

\begin{lem}[Local-in-space \apriori estimate for a solution to Yang-Mills gradient flow]
\label{lem:Kozono_Maeda_Naito_4-9}
\cite[Lemma 4.9]{Kozono_Maeda_Naito_1995}
Assume the hypotheses of Lemma \ref{lem:Kozono_Maeda_Naito_4-3} and, in addition, that $p \geq 2$. Then
$$
\sup_{\delta < t < T} \int_{B_{r/2}(x_0)} \left|\frac{\partial A}{\partial t}\right|^p \,d\vol_g \leq C,
$$
where $C$ has the same dependencies as those in Lemma \ref{lem:Kozono_Maeda_Naito_4-3} but now depends in addition on $k$ and $p$ and $\|F_A(\delta)\|_{W^{2,2}(B_r(x_0))}$ and $\|A(\delta) - \Gamma\|_{L^p(B_{r/2}(x_0))}$.
\end{lem}


An important hypothesis underpinning all of the local-in-space \apriori estimates reviewed thus far in this section is that the local energy for the Yang-Mills gradient flow,
$$
\int_{B_r(x_0)} |F_A(t)|^2 \,d\vol_g, \quad 0\leq t < T,
$$
can be controlled. With that in mind, we have the following version of the \cite[Lemma 4.10]{Kozono_Maeda_Naito_1995} due to Kozono, Maeda, Naito and the \cite[Lemma 2.2]{Schlatter_1997} due to Struwe \cite[Equation (13)]{Struwe_1985}; those lemmata are analogous to a result of Struwe \cite[Lemma 3.6]{Struwe_1985} for a harmonic-map gradient-like flow equation \cite[Equation (1.7)]{Struwe_1985}.

\begin{lem}[Local-in-space \apriori estimate for a solution to Yang-Mills gradient flow]
\label{lem:Kozono_Maeda_Naito_4-10}
Let $P$ be a principal $G$-bundle over a closed, four-dimensional, Riemannian, smooth manifold, $X$. Then there is a positive constant, $c$, with the following significance. Let $R_0$ be as in Proposition \ref{prop:Kozono_Maeda_Naito_2-1}. If $A(t)$ is a smooth solution to the Yang-Mills gradient flow equation \eqref{eq:Yang-Mills_gradient_flow_equation} over $(0, T) \times X$ with initial value $A_0$ of class $H^1$, and $r \in (0, R_0]$ and $x_0 \in X$, then
\begin{multline*}
\int_0^t\int_{B_{r/2}(x_0)} \left|\frac{\partial A}{\partial s}(s)\right|^2\,d\vol_g\,ds + \int_{B_{r/2}(x_0)} |F_A(t)|^2 \,d\vol_g
\\
\leq \int_{B_r(x_0)} |F_A(0)|^2 \,d\vol_g + cr^{-2}t\sE(A_0), \quad\forall\, t \in [0, T).
\end{multline*}
\end{lem}

\begin{proof}
We add details to explain the origin of the inequality immediately preceding \cite[Equation (4.22)]{Kozono_Maeda_Naito_1995} in the proof of \cite[Lemma 4.10]{Kozono_Maeda_Naito_1995} and also strengthen the inequality via the addition of the $L^2$ time and space integral of $\dot A$ on the left-hand side. (The outline of the proof of \cite[Lemma 2.2]{Schlatter_1997} due to Schlatter appears incomplete.)

We begin by fixing a $C^\infty$ cut-off function $\kappa:\RR\to [0,1]$ such that $\kappa(t) = 0$ for $t \geq 1$ and $\kappa(t) = 1$ for $t \leq 1/2$. Now define a $C^\infty$ cut-off function $\chi:X \to [0,1]$ by setting $\chi(x) := \kappa(\dist_g(x,x_0)/r)$ for $x \in X$ and observe that $\chi = 1$ on $B_{r/2}(x_0)$ and $\chi = 0$ on $X \less B_r(x_0)$. Hence, there is a positive constant, $c$, depending at most on the Riemannian metric $g$ on $X$ and a (universal) pointwise bound for $\nabla\kappa$ on $\RR$ such that $|\nabla\chi| \leq cr^{-1}$ on $X$.

By multiplying the Yang-Mills gradient flow equation, $\partial_t A = -d_A^*F_A$ on $(0,T)\times X$, by $\chi$, squaring both sides, and integrating over $X$, we find that
\begin{align*}
{}&\int_X \chi^2 |\partial_t A|^2\,d\vol_g
\\
&\quad = \int_X \chi^2\langle d_A^*F_A, d_A^*F_A\rangle \,d\vol_g
\\
&\quad = \int_X \chi^2\langle d_Ad_A^*F_A, F_A\rangle \,d\vol_g + \int_X \langle 2\chi d\chi\wedge d_A^*F_A, F_A\rangle \,d\vol_g
\\
&\quad = \int_X \chi^2\langle \Delta_AF_A, F_A\rangle \,d\vol_g - \int_X \langle 2\chi d\chi\wedge \partial_t A, F_A\rangle \,d\vol_g
\quad \hbox{(by \eqref{eq:Lawson_page_93_Hodge_Laplacian} and \eqref{eq:Freed_Uhlenbeck_2-16_Bianchi_identity})}
\\
&\quad = -\int_X \chi^2\langle \partial_t F_A, F_A\rangle \,d\vol_g - \int_X \langle 2\chi d\chi\wedge \partial_t A, F_A\rangle \,d\vol_g
\quad \hbox{(by \eqref{eq:Kozono_Maeda_Naito_2-4})}
\\
&\quad = -\frac{1}{2}\int_X \chi^2 \partial_t |F_A|^2 \,d\vol_g - \int_X \langle 2\chi d\chi\wedge \partial_t A, F_A\rangle \,d\vol_g
\\
&\quad \leq -\frac{1}{2}\int_X \chi^2 \partial_t |F_A|^2 \,d\vol_g
+ 2\sup_X|d\chi|\left(\int_X \chi^2|\partial_t A|^2 \,d\vol_g\right)^{1/2}\left(\int_{\supp d\chi} |F_A|^2 \,d\vol_g\right)^{1/2}
\end{align*}
Therefore,
$$
\int_X \chi^2 |\partial_t A|^2\,d\vol_g + \frac{1}{2}\int_X \chi^2 \partial_t |F_A|^2 \,d\vol_g
\leq \frac{1}{2}\int_X \chi^2|\partial_t A|^2 \,d\vol_g + cr^{-2}\int_X |F_A|^2 \,d\vol_g,
$$
and so,
\begin{align*}
\frac{1}{2}\int_X \chi^2 |\partial_t A|^2\,d\vol_g + \frac{1}{2}\frac{d}{dt}\int_X \chi^2|F_A|^2 \,d\vol_g &\leq cr^{-2}\int_X |F_A|^2 \,d\vol_g
\\
&\leq cr^{-2}\int_X |F_A(0)|^2 \,d\vol_g \quad \hbox{(by Lemma \ref{lem:Kozono_Maeda_Naito_4-1})}.
\end{align*}
Integrating with respect to time over $[0, t]$ yields
\begin{multline*}
\int_0^t\int_X \chi^2 |\partial_s A(s)|^2\,d\vol_g\,ds + \int_X \chi^2 |F_A(t)|^2 \,d\vol_g - \int_X \chi^2 |F_A(0)|^2 \,d\vol_g
\\
\leq 2cr^{-2}t\int_X |F_A(0)|^2 \,d\vol_g.
\end{multline*}
The conclusion now follows.
\end{proof}

Of course, by taking $\chi \equiv 1$ on $X$ in the proof of Lemma \ref{lem:Kozono_Maeda_Naito_4-10}, we obtain a familiar global version and simple analogue of the \apriori estimate \cite[Lemma 3.4]{Struwe_1985}, \cite[Lemma 1.1]{Struwe_1996} for harmonic map gradient-like flow and \cite[Lemma 4.1]{Kozono_Maeda_Naito_1995}, \cite[Equation (12)]{Struwe_1994} for Yang-Mills gradient flow.

\begin{cor}[Energy estimate for a solution to Yang-Mills gradient flow]
\label{cor:Kozono_Maeda_Naito_4-10_global}
\cite[Equation (12)]{Struwe_1994}
Let $P$ be a principal $G$-bundle over a closed, four-dimensional, Riemannian, smooth manifold, $X$. Then there is a positive constant, $c$, with the following significance. If $A(t)$ is a smooth solution to the Yang-Mills gradient flow equation \eqref{eq:Yang-Mills_gradient_flow_equation} over $(0, T) \times X$ with initial value $A_0$ of class $H^1$, then
$$
\int_0^t\int_X \left|\frac{\partial A}{\partial s}(s)\right|^2\,d\vol_g\,ds + \int_X |F_A(t)|^2 \,d\vol_g
\leq \int_X |F_A(0)|^2 \,d\vol_g, \quad\forall\, t \in [0, T).
$$
\end{cor}

This completes our review of the local-in-space \apriori estimates for the curvature of a solution to Yang-Mills gradient flow equation from \cite[Section 4]{Kozono_Maeda_Naito_1995}.

\subsection{Evolution of anti-self-dual curvature in Yang-Mills gradient flow}
\label{subsec:Evolution_ASD_curvature_Yang-Mills_gradient_flow}
The Chern-Weil formula yields the following useful refinement of Lemma \ref{lem:Kozono_Maeda_Naito_4-1}, showing that not only is the energy $\sE(A(t)) = \frac{1}{2}\|F_A(t)\|_{L^2(X)}^2$ a non-increasing function of time, $t\geq 0$, but so also are the of the anti-self-dual and self-dual curvature components, $\sE^\pm(A(t)) := \frac{1}{2}\|F_A^\pm(t)\|_{L^2(X)}^2$.

\begin{lem}[Evolution of anti-self-dual and self-dual curvatures in Yang-Mills gradient flow]
\label{lem:Nonincreasing_ASD_curvature_Yang-Mills_gradient_flow}
Let $G$ be a compact Lie group and $P$ a principal $G$-bundle over a closed, connected, four-dimensional, oriented, smooth manifold, $X$, with Riemannian metric, $g$, and $T \in (0,\infty]$. If $A(t)$, for $t\in [0,T)$, is a solution to Yang-Mills gradient flow \eqref{eq:Yang-Mills_gradient_flow} on $P$, then
$$
\int_X |F_A^\pm(t)|_\fg^2\,d\vol_g
\leq
\int_X |F_A^\pm(0)|_\fg^2\,d\vol_g,
\quad \forall\,t\in [0,T).
$$
\end{lem}

\begin{proof}
We recall from Section \ref{sec:Taubes_1982_Appendix} that, for any connection $A$ on $P$, we have
$$
\int_X |F_A|_\fg^2\,d\vol_g
=
\sum_{i=1}^l \frac{1}{r_{\fg_i}}\int_X \left(|F_A^+|_{\fg_i}^2
+ |F_A^-|_{\fg_i}^2\right)\,d\vol_g,
$$
and
\begin{equation}
\label{eq:L2_norm_anti-self-dual_curvature_bound}
4\pi^2|\bkappa(P)|
=
\sum_{i=1}^l \frac{1}{r_{\fg_i}}\int_X \left(|F_A^+|_{\fg_i}^2
- |F_A^-|_{\fg_i}^2\right)\,d\vol_g,
\end{equation}
where $|\bkappa(P)| := \sum_{i=1}^l\kappa_i(P)$, and thus
$$
4\pi^2|\bkappa(P)| + \sum_{i=1}^l \frac{1}{r_{\fg_i}}\int_X |F_A^-|_{\fg_i}^2 \,d\vol_g
=
\sum_{i=1}^l \frac{1}{r_{\fg_i}}\int_X |F_A^+|_{\fg_i}^2 \,d\vol_g.
$$
To analyze the evolution of the `self-dual energy', we apply the identity \eqref{eq:L2_norm_anti-self-dual_curvature_bound} and resulting equality to $A(t)$ yields
\begin{align*}
{}&\sum_{i=1}^l \frac{1}{r_{\fg_i}}\int_X |F_A^+(t)|_{\fg_i}^2 \,d\vol_g
\\
&\quad =
\sum_{i=1}^l \frac{1}{r_{\fg_i}}\int_X |F_A(t)|_{\fg_i}^2 \,d\vol_g
- \sum_{i=1}^l \frac{1}{r_{\fg_i}}\int_X |F_A^-(t)|_{\fg_i}^2 \,d\vol_g
\\
&\quad \leq \sum_{i=1}^l \frac{1}{r_{\fg_i}}\int_X |F_A(0)|_{\fg_i}^2 \,d\vol_g
- \sum_{i=1}^l \frac{1}{r_{\fg_i}}\int_X |F_A^-(t)|_{\fg_i}^2 \,d\vol_g
\quad\hbox{(by Lemma \ref{lem:Kozono_Maeda_Naito_4-1})}
\\
&\quad = \sum_{i=1}^l \frac{1}{r_{\fg_i}}\int_X |F_A^+(0)|_{\fg_i}^2 \,d\vol_g
+ \sum_{i=1}^l \frac{1}{r_{\fg_i}}\int_X |F_A^-(0)|_{\fg_i}^2 \,d\vol_g
- \sum_{i=1}^l \frac{1}{r_{\fg_i}}\int_X |F_A^-(t)|_{\fg_i}^2 \,d\vol_g
\\
&\quad = 2\sum_{i=1}^l \frac{1}{r_{\fg_i}}\int_X |F_A^+(0)|_{\fg_i}^2 \,d\vol_g
- 4\pi^2|\bkappa(P)|
- \sum_{i=1}^l \frac{1}{r_{\fg_i}}\int_X |F_A^-(t)|_{\fg_i}^2 \,d\vol_g
\quad\hbox{(by \eqref{eq:L2_norm_anti-self-dual_curvature_bound})}
\\
&\quad = 2\sum_{i=1}^l \frac{1}{r_{\fg_i}}\int_X |F_A^+(0)|_{\fg_i}^2 \,d\vol_g
- \sum_{i=1}^l \frac{1}{r_{\fg_i}}\int_X |F_A^+(t)|_{\fg_i}^2 \,d\vol_g
\quad\hbox{(by \eqref{eq:L2_norm_anti-self-dual_curvature_bound})},
\end{align*}
yields the desired inequality by combining terms.

For the evolution of the anti-self-dual curvature component, one argues in the same way. We include the details for completeness. Applying the identity \eqref{eq:L2_norm_anti-self-dual_curvature_bound} and the resulting equality
$$
4\pi^2|\bkappa(P)| - \sum_{i=1}^l \frac{1}{r_{\fg_i}}\int_X |F_A^+|_{\fg_i}^2 \,d\vol_g
=
-\sum_{i=1}^l \frac{1}{r_{\fg_i}}\int_X |F_A^-|_{\fg_i}^2 \,d\vol_g,
$$
to $A(t)$ yields
\begin{align*}
{}&\sum_{i=1}^l \frac{1}{r_{\fg_i}}\int_X |F_A^-(t)|_{\fg_i}^2 \,d\vol_g
\\
&\quad =
\sum_{i=1}^l \frac{1}{r_{\fg_i}}\int_X |F_A(t)|_{\fg_i}^2 \,d\vol_g
- \sum_{i=1}^l \frac{1}{r_{\fg_i}}\int_X |F_A^+(t)|_{\fg_i}^2 \,d\vol_g
\\
&\quad \leq \sum_{i=1}^l \frac{1}{r_{\fg_i}}\int_X |F_A(0)|_{\fg_i}^2 \,d\vol_g
- \sum_{i=1}^l \frac{1}{r_{\fg_i}}\int_X |F_A^+(t)|_{\fg_i}^2 \,d\vol_g
\quad\hbox{(by Lemma \ref{lem:Kozono_Maeda_Naito_4-1})}
\\
&\quad = \sum_{i=1}^l \frac{1}{r_{\fg_i}}\int_X |F_A^+(0)|_{\fg_i}^2 \,d\vol_g
+ \sum_{i=1}^l \frac{1}{r_{\fg_i}}\int_X |F_A^-(0)|_{\fg_i}^2 \,d\vol_g
- \sum_{i=1}^l \frac{1}{r_{\fg_i}}\int_X |F_A^+(t)|_{\fg_i}^2 \,d\vol_g
\\
&\quad = 2\sum_{i=1}^l \frac{1}{r_{\fg_i}}\int_X |F_A^-(0)|_{\fg_i}^2 \,d\vol_g
+ 4\pi^2|\bkappa(P)|
- \sum_{i=1}^l \frac{1}{r_{\fg_i}}\int_X |F_A^+(t)|_{\fg_i}^2 \,d\vol_g
\quad\hbox{(by \eqref{eq:L2_norm_anti-self-dual_curvature_bound})}
\\
&\quad = 2\sum_{i=1}^l \frac{1}{r_{\fg_i}}\int_X |F_A^+(0)|_{\fg_i}^2 \,d\vol_g
- \sum_{i=1}^l \frac{1}{r_{\fg_i}}\int_X |F_A^-(t)|_{\fg_i}^2 \,d\vol_g
\quad\hbox{(by \eqref{eq:L2_norm_anti-self-dual_curvature_bound})},
\end{align*}
again yields the desired inequality by combining terms.
\end{proof}

\subsection{Local \apriori estimates due to Schlatter, Struwe, and Uhlenbeck}
\label{subsec:Schlatter_Struwe_and_Uhlenbeck_apriori_estimates}
We recall some of the key \apriori estimates employed by Schlatter, Struwe, and Uhlenbeck in Schlatter's proof of \cite[Theorem 1.2]{Schlatter_1997}. We sometimes include proofs since some the statements in \cite{Schlatter_1997} lack a degree of precision that we shall later require. For consistency with Section \ref{subsec:Kozono_Maeda_Naito_4}, we let $\eps_1$ denote the smallest of the positive constants labeled $\eps_1$ appearing in this section and Section \ref{subsec:Kozono_Maeda_Naito_4} and for convenience we continue to assume that solutions, $A(t)$, to the Yang-Mills gradient flow equation \eqref{eq:Yang-Mills_gradient_flow} are classical, unless specified otherwise, that is smooth for $t > 0$ and if the initial condition is involved, continuous for $t \geq 0$.

We begin with the following refinement of \cite[Lemma 2.3]{Schlatter_1997}, \cite[Lemma 3.3]{Struwe_1994}.

\begin{lem}
\label{lem:Schlatter_2-3_and_Struwe_3-3}
\cite[Lemma 2.3]{Schlatter_1997}, \cite[Lemma 3.3]{Struwe_1994}
Let $G$ be a compact Lie group and $P$ be a principal $G$-bundle over a closed, four-dimensional, smooth manifold, $X$, with Riemannian metric, $g$. Then there are positive constants, $\eps_1$ and $C$, with the following significance. Let $A$ be a connection of class $H^1$ on $P$ such that
\begin{equation}
\label{eq:Struwe_lemma_3-3_sup_x_L2_ball_x_curvature_condition}
\sup_{x \in X} \int_{B_r(x)} |F_A|^2\,d\vol_g < \eps_1 \quad\hbox{for \emph{some} } r \in (0, 1].
\end{equation}
If $a \in H_A^1(X;\Lambda^p\otimes\ad P)$, for an integer $p \geq 0$, then
$$
\|a\|_{L^4(X)}^2 + \|\nabla_A a\|_{L^2(X)}^2
\leq C\left(\|d_A^*a\|_{L^2(X)}^2 + \|d_Aa\|_{L^2(X)}^2 \right) + Cr^{-2}\|a\|_{L^2(X)}^2.
$$
Moreover, the constants $\eps_1$ and $r$ depend on the metric $g$ through the quantities,
$$
\|\Riem(g)\|_{L^4(X,g)}, \quad \|g-g_0\|_{C(X)}, \quad\hbox{and}\quad \Inj(X,g_0),
$$
where $g_0$ is any fixed reference metric on $X$, and this dependence is continuous.
\end{lem}

\begin{rmk}[Dependence of the injectivity radius on the Riemannian metric]
\label{rmk:Dependence_injectivity_radius_Riemannian_metric}
By comparing the sectional curvatures of a pair of Riemannian metrics, $g_1$ and $g_2$, on a closed, smooth manifold $X$ of dimension $d\geq 2$, it is possible to compare their injectivity radii, $\Inj(X,g_1)$ and $\Inj(X,g_2)$, through the Rauch Comparison Theorems \cite[Chapter 1.10]{Cheeger_Ebin_1975}. Continuity properties of the injectivity radius as a function of the metric, $g$, are developed in \cite{Ehrlich_1974, Sakai_1983}. The dependence of the injectivity radius, $\Inj(X,g)$, on curvature of the Riemannian metric, $g$, which can be traced back to fundamental properties of geodesics \cite[Chapter 1]{Aubin}, suggests continuous dependence on $g \in \Met(X)$, the space of all Riemannian metrics on $X$, when equipped with the $C^2$ topology.
\end{rmk}

\begin{proof}[Proof of Lemma \ref{lem:Schlatter_2-3_and_Struwe_3-3}]
The only difference between the statement of Lemma \ref{lem:Schlatter_2-3_and_Struwe_3-3} and that in the cited references is in our clarification of the dependence of the constants on the Riemannian metric, $g$. The proof is otherwise identical.

The first occurrence of constants depending on $g$ is at the top of \cite[p. 132]{Struwe_1994}, through the use of the Sobolev embedding $H^1(X,g) \hookrightarrow L^4(X,g)$ with estimate \cite[Theorem 4.12]{AdamsFournier},
$$
\|f\|_{L^4(X,g)} \leq c_1(g)\|f\|_{H^1(X,g)}
\quad\hbox{and}\quad
\|f\|_{H^1(X,g)} = \left(\int_X\left(f^2 + |df|_g^2\right)d\vol_g\right)^{1/2},
$$
for all $f \in C^\infty(X)$. Plainly, $c_1(g)$ depends continuously on $g$ with respect to the $C^0(X)$ topology on the space of Riemannian metrics on $X$. As usual, when this Sobolev estimate is converted to one for $a \in \Omega^p(X; \ad P)$, with covariant derivative on $\Lambda^p(T^*X)\otimes\ad P$ induced by the Levi-Civita connection on $TX$ and connection $A$ on $P$ with the aid of the Kato Inequality \eqref{eq:FU_6-20_first-order_Kato_inequality}, namely $|d|a|| \leq |\nabla_A a|$, the constant $c_1(g)$ remains unchanged.

The second occurrence of constants depending on $g$ is near the top of \cite[p. 132]{Struwe_1994}, through the use of the Bochner-Weitzenb\"ock formula \eqref{eq:Lawson_corollary_II-2} (this version is stated for one-forms but the formula holds more generally with $\Ric_g$ replaced by $\Riem_g$). Rather than estimate
$$
(\{\Riem_g,a\}, a)_{L^2(X,g)}
\leq
c_0\|\Riem_g\|_{C(X,g)}\|a\|_{L^2(X,g)}^2,
$$
as in \cite{Struwe_1994} (with universal constant $c_0$ independent of $g$), we may use
\begin{align*}
(\{\Riem_g,a\}, a)_{L^2(X,g)} &\leq c_0\|\Riem_g\|_{L^4(X,g)} \|a\|_{L^4(X,g)} \|a\|_{L^2(X,g)}
\\
&\leq \frac{c_0}{2}\|\Riem_g\|_{L^4(X,g)}
\left( \zeta\|a\|_{L^4(X,g)}^2 + \zeta^{-1}\|a\|_{L^2(X,g)}^2 \right),
\end{align*}
for arbitrary $\zeta > 0$. With a choice of $\zeta$ such that
$$
\frac{c_0\zeta}{2}\|\Riem_g\|_{L^4(X,g)} = \frac{1}{2c_1},
$$
observing that the Sobolev embedding constant, $c_1$, is denoted by $C(\eta)$ in \cite[p. 132]{Struwe_1994}, one can use rearrangement and Struwe's proof of \cite[Lemma 3.3]{Struwe_1994} then continues without change.

The third apparent occurrence of constants depending on $g$ is through Struwe's use of the Sobolev embedding $H^1(B_R(x_0),g) \hookrightarrow L^4(B_R(x_0),g)$ with estimate,
$$
\|f\|_{L^4(B_R(x_0),g)} \leq c_2(g,R)\|f\|_{H^1(B_R(x_0),g)}, \quad\forall\, x_0 \in X,
$$
where $R \in (0, \Inj(X,g)]$ and $c_2(g,R)$ is potentially different from $c_1(g)$. However, one could instead use a covering of $X$ by balls $B_{R/2}(x_0)$ and fix, once and for all, a cut-off function $\chi \in C_0^\infty(B_1(0))$ with $B_1(0)\subset\RR^4$ having $0\leq \chi\leq 1$ on $\RR^4$ and $\chi=1$ on $B_{1/2}(0)$ and $|d\chi|\leq 4$. We then define $\chi_R := \chi(\dist_g(\cdot, x_0)/R) \in C_0^\infty(B_R(x_0))$ and observe that $|d\chi_R|_g \leq c(g)R^{-2}$. The constant, $c(g)$, may  be absorbed into a possibly larger Sobolev constant, $c_1(g)$, and this is how Struwe obtains his stated estimate for $(\{F_A,a\}, a)_{L^2(X,g)}$ (in our notation).

The last occurrence of constants depending on $g$ is through Struwe's choice of a radius $R_0 \in (0,\Inj(X,g)]$ and number $L$ independent of $(X,g)$ such that, if $\{B_R(x_i)\}$ is a cover of $X$ with $R \in (0, R_0]$, then at most $L$ balls intersect at each point $x\in X$.

Finally, as we noted in Remark \ref{rmk:Dependence_injectivity_radius_Riemannian_metric}, the dependence of the injectivity radius, $\Inj(X,g)$, on the Riemannian $g$ can be studied with the aid of $C^0$ curvature bounds, implying $C^2$ control over $g$. However, the term
$$
(\{F_A,a\}, a)_{L^2(X,g)},
$$
depends continuously on $g\in \Met(X)$, the space of all Riemannian metrics on $X$, when equipped with the $C^0$ topology. Thus, if $g_0$ is a reference Riemannian metric on $X$, there is a small enough constant, $\eta=\eta(g_0) \in (0,1]$, such that if $\|g-g_0\|_{C(X)} \leq \eta$, then
$$
\frac{1}{2}(\{F_A,a\}, a)_{L^2(X,g_0)} \leq (\{F_A,a\}, a)_{L^2(X,g)} \leq 2(\{F_A,a\}, a)_{L^2(X,g_0)},
$$
for all $C^\infty$ connections $A$ on $P$ and $a\in \Omega^p(X;\ad P)$. Thus, for the constants $\eps_1$ and $r$ appearing in \eqref{eq:Struwe_lemma_3-3_sup_x_L2_ball_x_curvature_condition}, they may be assumed to depend on the injectivity radius of $\Inj(X,g_0)$ and $\|g-g_0\|_{C(X)}$ rather than $\Inj(X,g)$. This completes the proof.
\end{proof}

Our Corollary \ref{cor:Rade_8-1} can be viewed as a variant of the following result due to Struwe \cite{Struwe_1994}.

\begin{lem}
\label{lem:Struwe_3-4}
\cite[Lemma 3.4]{Struwe_1994}
Let $G$ be a compact Lie group and $P$ be a principal $G$-bundle over a closed, four-dimensional, Riemannian, smooth manifold, $X$. Then there are positive constants, $\eps_1$ and $C$, with the following significance. If $A_1$ is a fixed reference connection of class $C^\infty$ on $P$ and $T > 0$ and $A(t) = A_1 + a(t)$ for $t\in [0, T)$ is a classical solution to the Yang-Mills gradient flow \eqref{eq:Yang-Mills_gradient_flow} on $P$ such that $A(t)$, for $t\in [0, T)$, obeys
\begin{equation}
\label{eq:Struwe_15}
\sup_{(t,x) \in (0,T)\times X} \int_{B_r(x)} |F_A(t)|^2\,d\vol_g < \eps_1 \quad\hbox{for \emph{some} } r \in (0, 1],
\end{equation}
then $F_A \in L^3(0,T; L^3(X;\Lambda^2\otimes\ad P))$ and obeys
$$
\|F_A\|_{L^3(0,T; L^3(X))}^2 \leq C\left(1 + r^{-2}T\right)\sE(A(0)).
$$
\end{lem}

Lemma \ref{lem:Struwe_3-4} has the following useful consequence, which in turn leads to the important Lemma \ref{lem:Schlatter_2-6}.

\begin{lem}
\label{lem:Struwe_3-5}
\cite[Lemma 3.5]{Struwe_1994}
Assume the hypotheses of Lemma \ref{lem:Struwe_3-4}. Then $F_A$ obeys
\begin{align*}
d_A^*F_A &\in L_{\loc}^2((0,T]; L^4(X;\Lambda^1\otimes\ad P)),
\\
\frac{\partial F_A}{\partial t} &\in L_{\loc}^2((0,T]; L^2(X;\Lambda^2\otimes\ad P)).
\end{align*}
\end{lem}

\begin{lem}
\label{lem:Schlatter_2-4_and_Struwe_3-6}
\cite[Lemma 2.4]{Schlatter_1997}, \cite[Lemma 3.6]{Struwe_1994}
Let $G$ be a compact Lie group and $P$ be a principal $G$-bundle over a closed, four-dimensional, Riemannian, smooth manifold, $X$. If $A_1$ is a fixed reference connection of class $C^\infty$ on $P$ and $A(t) = A_1 + a(t)$ for $t\in [0, T)$ and $T > 0$ is a classical solution to the Yang-Mills gradient flow \eqref{eq:Yang-Mills_gradient_flow} on $P$ such that $A(t)$, for $t\in [0, T)$, obeys \eqref{eq:Struwe_15}, then $A(t)$ extends from $[0, T)$ to $[0, T]$ in the sense that
$$
a \in C_{\loc}([0,T]; H^1_{A_1}(X;\Lambda^1\otimes\ad P)).
$$
\end{lem}

In the sequel, we state and prove a local version, Lemma \ref{lem:Schlatter_2-4_and_Struwe_3-6_local}, of Lemma \ref{lem:Schlatter_2-4_and_Struwe_3-6} that we shall find useful in our analysis of Uhlenbeck limits. Lemma \ref{lem:Schlatter_2-6} below may be viewed as a local version of the global estimate Lemma \ref{lem:Schlatter_2-3_and_Struwe_3-3}.

\begin{lem}
\label{lem:Schlatter_2-6}
\cite[Lemma 2.6]{Schlatter_1997}
There is a positive constant, $\eps_1$, and, if $\rho > 0$, there is a positive constant $C = C(\rho)$ with the following significance. Let $G$ be a compact Lie group and $A$ be a connection of class $H^2$ on $B_1\times G$, where $B_r \subset \RR^4$ is the ball of radius $r > 0$ centered at the origin. If
$$
\int_{B_1} |F_A|^2 \, d^4x \leq \eps_1,
$$
then
$$
\|F_A\|_{H_A^1(B_\rho)} \leq C\left(\|d_A^*F_A\|_{L^2(B_1)} + \|F_A\|_{L^2(B_1)} \right).
$$
\end{lem}


\begin{lem}
\label{lem:Schlatter_2-7}
\cite[Lemma 2.7]{Schlatter_1997}
There are positive constants, $\eps \in (0, 1]$ and $K$ and, if $\rho > 0$, there is a positive constant $C = C(K,\rho)$ with the following significance. Let $G$ be a compact Lie group and $A$ be a connection of class $H^2$ on $B_1\times G$, where $B_r \subset \RR^4$ is the ball of radius $r > 0$ centered at the origin. If
\begin{enumerate}
  \item $\|A-\Gamma\|_{H_\Gamma^1(B_1)} \leq \eps$;
  \item $d_\Gamma^* (A-\Gamma) = 0$ over $B_1$;
  \item $\|F_A\|_{H_\Gamma^1(B_1)} \leq K$;
\end{enumerate}
where $\Gamma$ denotes the product connection on $B_1\times G$, then
$$
\|A-\Gamma\|_{H_\Gamma^2(B_\rho)}\leq C.
$$
\end{lem}


The following bound on local gauge transformations may be compared with \cite[Corollary 3.3 and Lemma 3.5]{UhlLp}.

\begin{lem}
\label{lem:Schlatter_2-8}
\cite[Lemma 2.8]{Schlatter_1997}
Let $G$ be a compact Lie group. Then there are positive constants, $C$ and $\eps$, with the following significance. For $i=1,2$, let $B_1(y_i)$ be two open balls in $\RR^4$ with radius one and centers, $y_i$, and non-empty overlap, and $A_i$ be connections on $B_1(y_i)\times G$ such that $A_i - \Gamma \in H_\Gamma^2(B_1(y_i); \Lambda^1\otimes\fg)$, where $\Gamma$ is the product connection $\RR^4\times G$, and $\|A_i - \Gamma\|_{H_\Gamma^1(B_1(y_i))} \leq \eps$. If there is a local gauge transformation,
$\varphi \in H_\Gamma^3(U; G)$, such that $\varphi^*A_1 = A_2$ a.e. on $U := B_1(y_1) \cap B_1(y_2)$, then
$$
\|\varphi\|_{H_\Gamma^3(U; G)} \leq C.
$$
\end{lem}


We next recall some fundamental results due to Uhlenbeck.

\begin{thm}[Construction of a $W^{2,p}$ transformation to Coulomb gauge]
\label{thm:Uhlenbeck_Lp_1-3}
\cite[Theorem 1.3]{UhlLp}
There are positive constants, $c$ and $\eps_1$, with the following significance. Let $G$ be a compact Lie group and $A = \Gamma + a$ be a connection on $B_1\times G$, where $B_1 \subset \RR^4$ is the unit ball centered at the origin, such that $a \in W^{1,p}(B_1; \Lambda^1\otimes\fg)$ and $p \geq 2$. If
$$
\int_{B_1} |F_A|^2 \, d^4x \leq \eps_1,
$$
then there is a $W^{2,p}$ gauge transformation, $\varphi:B_1 \to G$, such that the gauge-transformed connection, $\widetilde A = \varphi^*A = \Gamma + \widetilde a$ given by
$$
\widetilde a := \varphi^{-1}a \varphi + \varphi^{-1}d\varphi,
$$
obeys
\begin{subequations}
\begin{align}
\label{eq:Uhlenbeck_Lp_1-3_Coulomb_gauge}
d_\Gamma^* (\widetilde A-\Gamma) &= 0 \quad\hbox{a.e. on } B_1,
\\
\label{eq:Uhlenbeck_Lp_1-3_apriori_estimate}
\|\widetilde A-\Gamma\|_{W_\Gamma^{1,p}(B_1)} &\leq c\|F_A\|_{L^p(B_1)}.
\end{align}
\end{subequations}
\end{thm}

\begin{rmk}[Construction of a $W^{k+1,p}$ transformation to Coulomb gauge]
\label{rmk:Uhlenbeck_theorem_1-3_Wkp}
We note that if $A$ is of class $W^{k,p}$, for an integer $k \geq 1$ and $p \geq 2$, then the gauge transformation, $\varphi$, in Theorem \ref{thm:Uhlenbeck_Lp_1-3} is of class $W^{k+1,p}$; see \cite[page 32]{UhlLp}, the proof of \cite[Lemma 2.7]{UhlLp} via the Implicit Function Theorem for smooth functions on Banach spaces, and our proof of \cite[Theorem 1.1]{FeehanSlice} --- a global version of Theorem \ref{thm:Uhlenbeck_Lp_1-3}.
\end{rmk}

\begin{rmk}[Equivalent scale invariant version of the \apriori estimate \eqref{eq:Uhlenbeck_Lp_1-3_apriori_estimate}]
\label{rmk:Uhlenbeck_theorem_1-3_apriori_estimate_scale_invariant}
We recall from the Sobolev Embedding Theorem \cite[Theorem 4.12]{AdamsFournier} in dimension four that $W^{1,p}(B_1) \hookrightarrow L^q(B_1)$, for $p \in [2, 4)$ and $q \in [4, \infty)$ defined by $1/p = 1/4 + 1/q$. With that in mind, it is often useful to replace the \apriori estimate \eqref{eq:Uhlenbeck_Lp_1-3_apriori_estimate} in Theorem \ref{thm:Uhlenbeck_Lp_1-3} by an equivalent scale invariant version,
\begin{equation}
\label{eq:Uhlenbeck_Lp_1-3_apriori_estimate_scale_invariant}
\|\widetilde A-\Gamma\|_{L^q(B_1)} + \|\cov_\Gamma(\widetilde A-\Gamma)\|_{L^p(B_1)} \leq c\|F_A\|_{L^p(B_1)}.
\end{equation}
\end{rmk}

Uhlenbeck has two versions of her Removable Singularities Theorem, one for finite-energy Yang-Mills connections of class $W^{1,p}$ (with $p > 2$) over a punctured ball in $\RR^4$ and a more delicate version for finite-energy Sobolev connections of class $H^1$ over a punctured ball in $\RR^d$ with $d \geq 2$. We shall have occasion to use both in dimension four.

\begin{thm}[Removability of singularities for finite-energy Yang-Mills connections]
\label{thm:Uhlenbeck_removable_singularity_Yang-Mills}
\cite[Theorem 4.1]{UhlRem}
Let $G$ be a compact Lie group and $A$ be a Yang-Mills connection, of class $W^{1,p}$ with $p>2$, on a principal $G$-bundle over $B_1\less\{0\}$, where $B_1\subset \RR^4$ is the unit ball centered at the origin, such that
$$
\int_{B_1} |F_A|^2 \, d^4x < \infty,
$$
Then there is a $W^{2,p}$ gauge transformation, $\varphi: P\restriction B_1\less\{0\} \to P\restriction B_1\less\{0\}$, such that the gauge-transformed connection, $\varphi^*A$, extends to a $W^{1,p}$ Yang-Mills connection, $\bar A$, on a principal $G$-bundle $\bar P$ over $B_1$.
\end{thm}

\begin{rmk}[On the regularity of the connections in Uhlenbeck's Removable Singularities Theorem]
\label{rmk:Uhlenbeck_removable_singularity_Yang-Mills}
Uhlenbeck assumes that the given connection, $A$, in \cite[Theorem 4.1]{UhlRem} is $C^\infty$. This assumption implies no loss of generality and reduces notational clutter in the statement and proof of her theorem. Nevertheless, as we wish to keep track of the regularity of all connections and gauge transformations, we include the minimal regularity in both the Yang-Mills and non-Yang-Mills assertions in Theorem \ref{thm:Uhlenbeck_removable_singularity_Yang-Mills}.
\end{rmk}

We shall only need the case $d=4$ of \cite[Theorem 2.1]{UhlChern}, which we record below.

\begin{thm}[Removability of singularities for finite-energy Sobolev connections]
\label{thm:Uhlenbeck_removable_singularity_Sobolev}
\cite[Theorems 2.1 and 4.5 and Corollary 4.6]{UhlChern}
Let $G$ be a compact Lie group. Then there are positive constants, $c$ and $\eps_1$, with the following significance. If $A$ is a connection of class $H^1$ on the product bundle, $G \times B_1\less\{0\}$, where $B_1\subset \RR^4$ is the unit ball centered at the origin, such that
$$
\int_{B_1} |F_A|^2 \, d^4x \leq \eps_1,
$$
then there is a gauge transformation, $\varphi \in H_{\Gamma,\loc}^2 (B_1\less\{0\}; G)$, such that the gauge-transformed connection, $\varphi^*A$, obeys
\begin{subequations}
\begin{align}
\label{eq:Uhlenbeck_removable_singularity_Sobolev_Coulomb_gauge}
d_\Gamma^*(\varphi^*A - \Gamma) &= 0 \quad\hbox{a.e. on } B_1,
\\
\label{eq:Uhlenbeck_removable_singularity_Sobolev_A_is_H1_on_ball}
\varphi^*A - \Gamma &\in H_\Gamma^1(B_1; \Lambda^1\otimes\fg),
\\
\label{eq:Uhlenbeck_removable_singularity_Sobolev_apriori_estimate}
\|\varphi^*A - \Gamma\|_{H_\Gamma^1(B_1)} &\leq c\|F_A\|_{L^2(B_1)}.
\end{align}
\end{subequations}
If the Coulomb gauge condition is omitted, then the gauge transformation, $\varphi$, may be chosen to be $C^\infty$ on $\{x \in \RR^4: \rho \leq |x| \leq 1-\rho$, for $\rho \in (0, 1/2)$.
\end{thm}

\begin{rmk}[Construction of a $W^{k+1,p}$ transformation to Coulomb gauge and removability of singularities]
\label{rmk:Uhlenbeck_removable_singularity_Sobolev_Wkp}
As in the case of Theorem \ref{thm:Uhlenbeck_removable_singularity_Yang-Mills}, the proof of Theorem \ref{thm:Uhlenbeck_removable_singularity_Sobolev} relies on the Implicit Function Theorem for smooth functions on Banach spaces to construct the required gauge transformation. If $A$ is of class $W^{k,p}$, for an integer $k \geq 1$ and $p \geq 2$, one may rework the proof of Theorem \ref{thm:Uhlenbeck_removable_singularity_Sobolev} to show that the gauge transformation, $\varphi$, in Theorem \ref{thm:Uhlenbeck_removable_singularity_Sobolev} is of class $W^{k+1,p}$.
\end{rmk}

Just as we noted in Remark \ref{rmk:Uhlenbeck_theorem_1-3_apriori_estimate_scale_invariant}, we may replace \eqref{eq:Uhlenbeck_removable_singularity_Sobolev_apriori_estimate} by
\begin{equation}
\label{eq:Uhlenbeck_removable_singularity_Sobolev_apriori_estimate_scale_invariant}
\|\widetilde A-\Gamma\|_{L^4(B_1)} + \|\cov_\Gamma(\widetilde A-\Gamma)\|_{L^2(B_1)} \leq c\|F_A\|_{L^2(B_1)},
\end{equation}
an equivalent scale invariant version that is often useful in the sequel.

\subsection{A local version of the continuous extension Lemma \ref{lem:Schlatter_2-4_and_Struwe_3-6} of Struwe}
\label{subsec:Local_version_continuous_extension_lemma_Schlatter_2-4}
We shall need a local version of Lemma \ref{lem:Schlatter_2-4_and_Struwe_3-6}.

\begin{lem}
\label{lem:Schlatter_2-4_and_Struwe_3-6_local}
Let $G$ be a compact Lie group and $P$ be a principal $G$-bundle over a closed, four-dimensional, Riemannian, smooth manifold, $X$. Let $\eps_1$ and $R_0$ be the smaller of the corresponding pair of constants in Section \ref{subsec:Kozono_Maeda_Naito_4}. Suppose that $A_1$ is a fixed reference connection of class $C^\infty$ on $P$ and $U \subset X$ is an open subset and $A(t) = A_1 + a(t)$ for $t\in [0, T)$ with $0 < T < \infty$, is a classical solution to the Yang-Mills gradient flow equation \eqref{eq:Yang-Mills_gradient_flow} on $P \restriction U$. If $U' \Subset U$ is an open subset and\footnote{If $U = X$, then $d_0 =  \infty$ and so Lemma \ref{lem:Schlatter_2-4_and_Struwe_3-6_local} generalizes Lemma \ref{lem:Schlatter_2-4_and_Struwe_3-6}.}
$d_0 := \dist(U', \partial U)$ and $A(t)$ obeys
\begin{equation}
\label{eq:Struwe_15_local}
\sup_{(t,x) \in (0,T)\times U'} \int_{B_r(x)} |F_A(t)|^2\,d\vol_g < \eps_1 \quad\hbox{for \emph{some} } r \in (0, d_0\wedge R_0],
\end{equation}
then the solution $A(t)$ extends continuously from $(0, T)$ to $(0, T]$ over $U'$ in the sense that
$$
a \in C_{\loc}((0,T]; H_{A_1}^1(U';\Lambda^1\otimes\ad P)).
$$
\end{lem}

\begin{rmk}[On the dependence of the constants in Lemma \ref{lem:Schlatter_2-4_and_Struwe_3-6_local} on the Riemannian metric]
\label{rmk:Schlatter_2-4_and_Struwe_3-6_local_equivalent_Riemannian_metric}
If $\bar g$ is another Riemannian metric on $X$ which is $C^2$ equivalent to the given Riemannian metric $g$ on $X$ in the sense that $\|g - \bar g\|_{C^2(X,g)} \leq K$, for some positive constant, $K$, that is uniform with respect to $\|g - \bar g\|_{C^0(X,g)} \in (0,1]$, then $C^0$ norms of curvatures, distances between pairs of points in $X$, injectivity radii, tensor norms, volumes of open subsets, and the positive constants $\bar\eps_1$ and $\bar R_0$ defined by $\bar g$ will be comparable to the corresponding quantities defined by $g$, with comparison constants depending on $g$ and $K$.

In the metric $\bar g$ is not $C^2$ equivalent to $g$, then we note that the weaker equivalence described in Lemma \ref{lem:Schlatter_2-3_and_Struwe_3-3} carries over \mutatis to the statement and proof of Lemma \ref{lem:Schlatter_2-4_and_Struwe_3-6_local}.
\end{rmk}

\begin{proof}[Proof of Lemma \ref{lem:Schlatter_2-4_and_Struwe_3-6_local}]
We shall adapt Struwe's proof of his \cite[Lemma 3.6]{Struwe_1994}. Since $\bar U'$ is a compact subset of $U$, it suffices to consider a ball $B_r(x) \subset U$, for a point $x \in U'$ and radius $r \in (0, R_0]$, where $R_0$ is the positive constant appearing in Section \ref{subsec:Kozono_Maeda_Naito_4}. For any $\delta \in (0, T)$, Lemma \ref{lem:Kozono_Maeda_Naito_4-9} provides a positive constant, $C$, such that
\begin{equation}
\label{eq:Kozono_Maeda_Naito_4-9_apriori_estimate}
\sup_{t\in (\delta,T)}\left\|\frac{\partial a}{\partial t}(t)\right\|_{L^4(B_{r/2}(x))} \leq C.
\end{equation}
Therefore, because $-\infty < \delta < T < \infty$, we obtain $a \in C([\delta,T]; L^4(B_{r/2}(x); \Lambda^1\otimes\ad P))$. (The dependencies of $C$ include $\|F_A(\delta)\|_{W_\Gamma^{2,2}(B_r(x))}$ and $\|a(\delta)\|_{L^4(B_{r/2}(x))}$.) Recalling that $F_{A_1+a} = F_{A_1} + d_{A_1}a + [a, a]$, we obtain
$$
\frac{\partial d_{A_1} a}{\partial t} = \frac{\partial F_A}{\partial t} + \left[\frac{\partial a}{\partial t}, a\right] + \left[a, \frac{\partial a}{\partial t}\right].
$$
Also, Lemma \ref{lem:Kozono_Maeda_Naito_4-7} provides a positive constant, $C$, such that
\begin{equation}
\label{eq:Kozono_Maeda_Naito_4-7_apriori_estimate_L-infinity_time_and_L4_space}
\sup_{t\in (\delta,T)} \|a(t)\|_{L^4(B_{r/2}(x))} \leq C,
\end{equation}
where the dependencies of $C$ include $\|F_A(\delta)\|_{W_\Gamma^{1,4}(B_r(x))}$ and $\|a(\delta)\|_{L^4(B_{r/2}(x))}$. Because
\begin{align*}
\frac{\partial a}{\partial t} &\in L^\infty(\delta, T; L^4(B_{r/2}(x); \Lambda^1\otimes\ad P)) \quad\hbox{(by \eqref{eq:Kozono_Maeda_Naito_4-9_apriori_estimate})} \quad\hbox{and}
\\
a &\in L^\infty(\delta, T; L^4(B_{r/2}(x); \Lambda^1\otimes\ad P)) \quad\hbox{(by \eqref{eq:Kozono_Maeda_Naito_4-7_apriori_estimate_L-infinity_time_and_L4_space})},
\end{align*}
we see that
$$
\left[\frac{\partial a}{\partial t}, a\right] + \left[a, \frac{\partial a}{\partial t}\right] \in L^\infty(\delta, T; L^2(B_{r/2}(x); \Lambda^1\otimes\ad P)).
$$
Moreover, as
\begin{align}
\label{eq:FA_is_L-infinity_time_and_L-infinity_space}
F_A &\in L^\infty(\delta, T; L^\infty(B_{r/2}(x); \Lambda^1\otimes\ad P)) \quad\hbox{(by Lemma \ref{lem:Kozono_Maeda_Naito_4-3})},
\\
\label{eq:nablaA_FA_is_L2_time_and_L2_space}
\nabla_A F_A &\in L^2(\delta, T; L^2(B_{r/2}(x); \Lambda^1\otimes\ad P)) \quad\hbox{(by Lemma \ref{lem:Kozono_Maeda_Naito_4-2})},
\\
\label{eq:nablaA2_FA_is_L2_time_and_L2_space}
\nabla_A^2F_A &\in L^2(\delta, T; L^2(B_{r/2}(x); \Lambda^1\otimes\ad P)) \quad\hbox{(by Lemma \ref{lem:Kozono_Maeda_Naito_4-5})},
\end{align}
and recalling from \eqref{eq:Kozono_Maeda_Naito_2-5} that
$$
\frac{\partial F_A}{\partial t} = -\nabla_A^*\nabla_A F_A - \Riem_g\times F_A - F_A\times F_A,
$$
we obtain
$$
\frac{\partial F_A}{\partial t} \in L^2(\delta, T; L^2(B_{r/2}(x); \Lambda^1\otimes\ad P)).
$$
Therefore, by combining the preceding observations, we find that
$$
\frac{\partial d_{A_1} a}{\partial t} \in L^2(\delta, T; L^2(B_{r/2}(x); \Lambda^1\otimes\ad P)),
$$
and so, because $[\delta,T]$ is a finite-length interval,
$$
d_{A_1} a \in C([\delta, T]; L^2(B_{r/2}(x); \Lambda^1\otimes\ad P)).
$$
Furthermore,
\begin{align*}
\frac{\partial d_{A_1}^* a}{\partial t} &= d_{A_1}^* \frac{\partial A}{\partial t}
\\
&= -d_{A_1}^*d_A^*F_A  \quad\hbox{(by \eqref{eq:Yang-Mills_gradient_flow})}
\\
&= -d_A^*d_A^*F_A - *[a, *d_A^*F_A] \quad\hbox{(because $A = A_1+a$)}
\\
&= (d_A^2)^*F_A  - *[a, *d_A^*F_A]
\\
&=  F_A^*F_A  - *[a, *d_A^*F_A].
\end{align*}
For any $b \in \Omega^1(B_\rho; \fg)$, where $B_\rho \subset \RR^4$ is the ball of radius $\rho > 0$ and center at the origin, define $\tilde b \in \Omega^1(B_1; \fg)$ by $b(\rho y) = \tilde b(y)$, for $y \in B_1$, and for the connection, $A_1$ on $B_\rho\times G$, we similarly define $\widetilde A_1$ on $B_1\times G$. Then, by the scaling behavior of the $L^4$ norm on one-forms and $L^2$ norm on two-forms,
$$
\|b\|_{L^4(B_\rho)}
=
\|\tilde b\|_{L^4(B_1)}
\leq
c\left(\|\nabla_{\widetilde A_1}\tilde b\|_{L^2(B_1)} + \|\tilde b\|_{L^2(B_1)}\right),
$$
by the Sobolev embedding $H^1(X;\RR) \hookrightarrow L^4(X;\RR)$ \cite[Theorem 4.12]{AdamsFournier} and the Kato Inequality \eqref{eq:FU_6-20_first-order_Kato_inequality}, and so
\begin{equation}
\label{eq:L^4_norm_one-form_bounded_by_H1_norm_ball_radius_rho}
\|b\|_{L^4(B_\rho)} \leq c\left(\|\nabla_{A_1}b\|_{L^2(B_\rho)} + \rho^{-1}\|b\|_{L^2(B_\rho)}\right),
\end{equation}
where $c$ is a universal positive constant. Applying the H\"older inequality, we discover that
\begin{align*}
{}&\|[a, *d_A^*F_A]\|_{L^2(B_{r/2}(x))}
\\
&\quad \leq c\|a\|_{L^4(B_{r/2}(x))} \|d_A^*F_A\|_{L^4(B_{r/2}(x))}
\\
&\quad \leq c\|a\|_{L^4(B_{r/2}(x))} \left(\|d_A^*F_A\|_{L^2(B_{r/2}(x))} + r^{-1}\|\nabla_A d_A^*F_A\|_{L^2(B_{r/2}(x))}\right)
\quad\hbox{(by \eqref{eq:L^4_norm_one-form_bounded_by_H1_norm_ball_radius_rho})}
\\
&\quad \leq c\|a\|_{L^4(B_{r/2}(x))} \left(\|\nabla_AF_A\|_{L^2(B_{r/2}(x))} + r^{-1}\|\nabla_A^2F_A\|_{L^2(B_{r/2}(x))}\right),
\end{align*}
where $c$ is a positive constant depending at most on the Riemannian metric on $X$. Thus, by \eqref{eq:Kozono_Maeda_Naito_4-7_apriori_estimate_L-infinity_time_and_L4_space}, \eqref{eq:FA_is_L-infinity_time_and_L-infinity_space}, \eqref{eq:nablaA_FA_is_L2_time_and_L2_space}, and \eqref{eq:nablaA2_FA_is_L2_time_and_L2_space}, we obtain
$$
\frac{\partial d_{A_1}^* a}{\partial t} \in L^2(\delta, T; L^2(B_{r/2}(x); \Lambda^1\otimes\ad P)),
$$
and therefore
$$
d_{A_1}^* a \in C([\delta, T]; L^2(B_{r/2}(x); \Lambda^1\otimes\ad P)).
$$
Let $\chi \in C^\infty(X)$ be a cut-off function such that $0 \leq \chi \leq 1$ on $X$ with $\supp\chi \subset B_{r/2}(x)$ and $\chi = 1$ on $B_{r/4}(x)$. We then observe that, for all $t \in (\delta,T)$ and writing $a = a(t)$,
\begin{align*}
\|\nabla_{A_1}(\chi a)\|_{L^2(X)}^2 &= (\nabla_{A_1}^*\nabla_{A_1}(\chi a), \chi a)_{L^2(X)}
\\
&= ((d_{A_1}^*d_{A_1} + d_{A_1}d_{A_1}^*)(\chi a), \chi a)_{L^2(X)} + (F_{A_1} \times \chi a, \chi a)_{L^2(X)}
\\
&\quad + (\Ric_g \times \chi a, \chi a)_{L^2(X)}  \quad\hbox{(by \eqref{eq:Lawson_corollary_II-2})}
\\
&= \|d_{A_1}(\chi a)\|_{L^2(X)}^2 + \|d_{A_1}^*(\chi a)\|_{L^2(X)}^2 + (F_{A_1} \times \chi a, \chi a)_{L^2(X)}
\\
&\quad + (\Ric_g \times \chi a, \chi a)_{L^2(X)}.
\end{align*}
Hence, there is a positive constant, $c$, depending at most on the Riemannian metric on $X$ such that, for all $t \in (\delta,T)$ and writing $a = a(t)$,
\begin{align*}
\|\nabla_{A_1}(\chi a)\|_{L^2(X)} &\leq c\left(\|d_{A_1}a\|_{L^2(B_{r/2}(x))} + \|d_{A_1}^*a\|_{L^2(B_{r/2}(x))} \right.
\\
&\quad + \|d\chi\|_{L^\infty(X)}\|a\|_{L^2(B_{r/2}(x))} + \|F_{A_1}\|_{L^\infty(X)}^{1/2} \|a\|_{L^2(B_{r/2}(x))}
\\
&\quad + \left. \|\Ric_g\|_{L^\infty(X)}^{1/2} \|a\|_{L^2(B_{r/2}(x))}\right).
\end{align*}
Thus, $\chi a \in C([\delta, T]; H_{A_1}^1(B_{r/2}(x); \Lambda^1\otimes\ad P))$ and so $a \in C([\delta, T]; H_{A_1}^1(B_{r/4}(x); \Lambda^1\otimes\ad P))$, as desired. The conclusion now follows.
\end{proof}

It is clear from the proof of Lemma \ref{lem:Schlatter_2-4_and_Struwe_3-6_local} that under its hypotheses one should expect higher-order continuity results of the form \eqref{eq:Schlatter_2-4_local_H2_or_W1p} for any integer $k\geq 1$ or $p \geq 2$. It will suffice for our application to consider the case of $k=2$ and $p=2$ or $k = 1$ and $p > 2$ since, for that situation, one has $W^{k+1,p}(X;\RR) \hookrightarrow C(X)$ by the Sobolev Embedding Theorem \cite[Theorem 4.12]{AdamsFournier}. We leave further extensions to the reader.

\begin{lem}
\label{lem:Schlatter_2-4_and_Struwe_3-6_local_H2_or_W1p}
Assume the hypotheses of Lemma \ref{lem:Schlatter_2-4_and_Struwe_3-6_local}. Then for every precompact open subset $U' \Subset U$, the solution $A(t)$ extends continuously from $(0, T)$ to $(0, T]$ over $U'$ in the sense that, $k=2$ and $p=2$ or $k = 1$ and $2 \leq p < \infty$,
\begin{equation}
\label{eq:Schlatter_2-4_local_H2_or_W1p}
a \in C_{\loc}((0,T]; W_{A_1}^{k, p}(U';\Lambda^1\otimes\ad P)).
\end{equation}
\end{lem}

\begin{proof}
Suppose first that $k = 1$ and $p \geq 2$. Since $A(t) = A_1 + a(t)$, we observe that
\begin{align*}
\frac{\partial \nabla_{A_1}a}{\partial t} &= \nabla_{A_1}\frac{\partial a}{\partial t}
\\
&= \nabla_A\frac{\partial a}{\partial t} + a\otimes \frac{\partial a}{\partial t}
\\
&= -\nabla_Ad_A^*F_A + a\otimes \frac{\partial a}{\partial t} \quad\hbox{(by \eqref{eq:Yang-Mills_gradient_flow})}.
\end{align*}
But for any integer $l \geq 0$, the curvature $F_A(t)$ (of the classical solution, $A(t)$) satisfies
\begin{multline}
\label{eq:nablaAl_FA_is_L-infinity_time_and_L-infinity_space}
\nabla_A^lF_A \in L^\infty(\delta, T; L^\infty(B_{r/2}(x); (T^*X)^{\otimes (l+1)}\otimes\ad P))
\\
\quad\hbox{(by Lemmata \ref{lem:Kozono_Maeda_Naito_4-3}, \ref{lem:Kozono_Maeda_Naito_4-4}, and \ref{lem:Kozono_Maeda_Naito_4-6})},
\end{multline}
and so \eqref{eq:nablaAl_FA_is_L-infinity_time_and_L-infinity_space} implies that
$$
\nabla_Ad_A^*F_A \in L^\infty(\delta, T; L^\infty(B_{r/2}(x); (T^*X)^{\otimes 2}\otimes\ad P)).
$$
On the other hand, for any $q \in [2, \infty)$, the terms $a(t)$ and $\partial a/\partial t$ satisfy
\begin{align}
\label{eq:Kozono_Maeda_Naito_4-7_regularity_L-infinity_time_and_Lq_space}
a &\in L^\infty(\delta, T; L^q(B_{r/2}(x); \Lambda^1\otimes\ad P)) \quad\hbox{(by Lemma \ref{lem:Kozono_Maeda_Naito_4-7})},
\\
\label{eq:Kozono_Maeda_Naito_4-9_regularity_L-infinity_time_and_Lq_space}
\frac{\partial a}{\partial t} &\in L^\infty(\delta, T; L^q(B_{r/2}(x); \Lambda^1\otimes\ad P)) \quad\hbox{(by Lemma \ref{lem:Kozono_Maeda_Naito_4-9})},
\end{align}
and so, choosing $q = 2p$, we see that
$$
a\otimes \frac{\partial a}{\partial t} \in L^\infty(\delta, T; L^p(B_{r/2}(x); (T^*X)^{\otimes 2}\otimes\ad P)).
$$
By combining the preceding observations, we discover that
\begin{equation}
\label{eq:Kozono_Maeda_Naito_4-9_nablaA1_a_regularity_L-infinity_time_and_Lp_space}
\frac{\partial \nabla_{A_1}a}{\partial t} \in L^\infty(\delta, T; L^p(B_{r/2}(x); (T^*X)^{\otimes 2}\otimes\ad P)).
\end{equation}
Therefore, by \eqref{eq:Kozono_Maeda_Naito_4-9_regularity_L-infinity_time_and_Lq_space} and \eqref{eq:Kozono_Maeda_Naito_4-9_nablaA1_a_regularity_L-infinity_time_and_Lp_space}, we see that
$$
\frac{\partial a}{\partial t} \in L^\infty(\delta, T; W_{A_1}^{1,p}(B_{r/2}(x); (T^*X)^{\otimes 2}\otimes\ad P)),
$$
and because $[\delta, T]$ is a finite-length interval, we have
$$
a \in C([\delta, T]; W_{A_1}^{1,p}(B_{r/2}(x); \Lambda^1\otimes\ad P)),
$$
which gives \eqref{eq:Schlatter_2-4_local_H2_or_W1p} for the case $k = 1$ and $p \geq 2$.

If $k = 1$ and $p = 2$, we similarly observe that
\begin{align*}
\frac{\partial \nabla_{A_1}^2a}{\partial t} &= \nabla_{A_1}^2\frac{\partial a}{\partial t}
\\
&= \nabla_A^2\frac{\partial a}{\partial t} + \nabla_{A_1}a \times \frac{\partial a}{\partial t}
 + a \times \nabla_{A_1}\frac{\partial a}{\partial t}  + a \times a \times \frac{\partial a}{\partial t}
\\
&= -\nabla_A^2d_A^*F_A + \nabla_{A_1}a \times \frac{\partial a}{\partial t}
 + a \times \frac{\partial \nabla_{A_1}a}{\partial t}  + a \times a \times \frac{\partial a}{\partial t}
\quad\hbox{(by \eqref{eq:Yang-Mills_gradient_flow})}.
\end{align*}
By applying \eqref{eq:nablaAl_FA_is_L-infinity_time_and_L-infinity_space} (with $l=3$), \eqref{eq:Kozono_Maeda_Naito_4-7_regularity_L-infinity_time_and_Lq_space}, \eqref{eq:Kozono_Maeda_Naito_4-9_regularity_L-infinity_time_and_Lq_space}, and \eqref{eq:Kozono_Maeda_Naito_4-9_nablaA1_a_regularity_L-infinity_time_and_Lp_space} to the preceding identity, we discover that
\begin{equation}
\label{eq:Kozono_Maeda_Naito_4-9_nablaA12_a_regularity_L-infinity_time_and_L2_space}
\frac{\partial \nabla_{A_1}^2a}{\partial t} \in L^\infty(\delta, T; L^2(B_{r/2}(x); (T^*X)^{\otimes 2}\otimes\ad P)).
\end{equation}
Therefore, by \eqref{eq:Kozono_Maeda_Naito_4-9_regularity_L-infinity_time_and_Lq_space}, \eqref{eq:Kozono_Maeda_Naito_4-9_nablaA1_a_regularity_L-infinity_time_and_Lp_space}, and \eqref{eq:Kozono_Maeda_Naito_4-9_nablaA12_a_regularity_L-infinity_time_and_L2_space}, we see that
$$
\frac{\partial a}{\partial t} \in L^\infty(\delta, T; H_{A_1}^2(B_{r/2}(x); (T^*X)^{\otimes 2}\otimes\ad P)),
$$
and because $[\delta, T]$ is a finite-length interval, we have
$$
a \in C([\delta, T]; H_{A_1}^2(B_{r/2}(x); \Lambda^1\otimes\ad P)),
$$
which gives \eqref{eq:Schlatter_2-4_local_H2_or_W1p} for the case $k = 2$ and $p = 2$.
\end{proof}

\subsection{An extension of the Lemma 4.2 of Kozono, Maeda, and Naito and analogue of the Lemma 8.1 of R\r{a}de}
\label{subsec:Lemma_Kozono_Maeda_Naito_4-3_global_extension}
As we know from the statement of Lemma \ref{lem:Rade_8-1}, a direct attempt to derive an four-dimensional analogue of R\r{a}de's \cite[Lemma 8.1]{Rade_1992} (which holds when $X$ has dimension two or three) only yields a conclusion which is useful when the initial curvature, $F_A(0)$, is $L^2(X)$-small and that condition is usually impossible to satisfy. However, the proof of Lemma \ref{lem:Rade_8-1} did suggest that a local analogue should hold and, indeed, such a result was established by Kozono, Maeda, and Naito in their \cite[Lemma 4.2]{Kozono_Maeda_Naito_1995}, which we recall in our monograph as Lemma \ref{lem:Kozono_Maeda_Naito_4-2}.

In this section, we establish a simple extension of Lemma \ref{lem:Kozono_Maeda_Naito_4-2} which will serve as a replacement for Lemma \ref{lem:Rade_8-1} when $X$ has dimension four but the initial curvature, $F_A(0)$, is not $L^2(X)$-small.

\begin{lem}[\Apriori estimate for the curvature of a solution to Yang-Mills gradient flow]
\label{lem:Kozono_Maeda_Naito_4-2_global}
Let $P$ be a principal $G$-bundle over a closed, four-dimensional, Riemannian, smooth manifold, $X$. Then there exist positive constants, $c$ and $\eps_1$, with the following significance. Let $R_0$ be as in Proposition \ref{prop:Kozono_Maeda_Naito_2-1}. If $A(t)$ is a smooth solution to the Yang-Mills gradient flow equation \eqref{eq:Yang-Mills_gradient_flow_equation} over $(0, T) \times X$ with initial value $A_0$ of class $H^1$, and $r \in (0, R_0]$ and $U \Subset V \Subset X$ are open subsets obeying
\begin{subequations}
\begin{align}
\label{eq:Lemma_Kozono_Maeda_Naito_4-2_global_condition_small_energy_on_small_balls}
\sup_{x\in U}\eps(r,x;A,T) &< \eps_1,
\\
\label{eq:Lemma_Kozono_Maeda_Naito_4-2_global_condition_boundaries_U_and_V_r_separated}
\dist_g(\partial U, \partial V) &\geq r,
\end{align}
\end{subequations}
where $\eps(r,x) = \eps(r,x;A,T)$ is as defined in \eqref{eq:Kozono_Maeda_Naito_4-1}, and
$$
K := \inf \left\{n \in \NN: \exists \{x_1,\ldots,x_n\} \subset U \hbox{ such that } U \Subset \bigcup_{i=1}^n B_{r/2}(x_i) \Subset V\right\},
$$
then
\begin{equation}
\label{eq:Kozono_Maeda_Naito_4-2_global}
\int_0^T \int_U |\nabla_A F_A(t)|^2 \,d\vol_g\,dt \leq cK(1 + r^{-2}T)\sE(A_0).
\end{equation}
\end{lem}

\begin{proof}
The conclusion is an immediate consequence of Lemma \ref{lem:Kozono_Maeda_Naito_4-2} and the given properties of $U$ and $V$ and the definition of the covering number, $K$.
\end{proof}

\subsection{Estimate of the contribution of a non-flat Riemannian metric to the size of the perturbation in Yang-Mills gradient-like flow}
\label{subsec:Estimate_gradient-like_flow_perturbation_due_to_non-flat_Riemannian_metric}
In the preceding subsections, we have reviewed or developed \apriori estimates for Yang-Mills gradient flow, that are local in time and space, which allow us to estimate the size of the error term, $R(t)$, in Yang-Mills gradient-like flow resulting from cutting off Yang-Mills gradient over a small annulus, where the energy of the connections, $A(t)$, is sufficiently small though the annulus may enclose a ball where the energy of $A(t)$ is becoming large. However, when we graft the gradient flow over a small ball in $X$ to a gradient-like flow over $S^4$, we also introduce an error due to the fact the Riemannian metric, $g$, is not necessarily flat near the center of the ball, $x_0 \in X$. When we choose geodesic normal coordinates, the components of the metric obey \cite[Definition 1.24, Proposition 1.25, and Corollary 1.32]{Aubin}
\begin{equation}
\label{eq:Riemannian_metric_components_in_geodesic_normal_coordinates}
\Gamma_{\alpha\beta}^\gamma(x_0) = 0,  \quad \Gamma_{\alpha\beta}^\gamma(x) = O(r), \quad\hbox{and}\quad g_{\mu\nu}(x) = \delta_{\mu\nu} + O(r^2),
\end{equation}
where $r := \dist_g(x, x_0)$. If $g_0$ denotes the flat metric near $x_0 \in X$ and $g$ is the given metric,
then\footnote{To obtain the correct sign, it is helpful to recall the definition of the adjoint of the exterior derivative, $d:\Omega^p(X) \to \Omega^{p+1}(X)$, in terms of the Hodge-star operator, $*_g:\Omega^p(X) \to \Omega^{d-p}(X)$, for integers $0 \leq p \leq d$, giving $d^{*_g}:\Omega^p(X) \to \Omega^{p-1}(X)$ as $d^{*_g} = (-1)^{-d(p+1)+1}*_gd*_g$.}
$$
d_A^{*_g}F_A = -*_g d_A *_g F_A.
$$
Therefore, over an neighborhood of $x_0 \in X$,
\begin{align*}
d_A^{*_g}F_A - d_A^{*_{g_0}}F_A &= -*_g d_A *_g F_A + *_{g_0} d_A *_{g_0} F_A
\\
&= -*_g d_A (*_g - *_{g_0})F_A + (*_{g_0} - *_g)d_A *_{g_0} F_A,
\end{align*}
and so, using $c$ to denote a constant depending at most on the Riemannian metric, $g$, on $X$, the following pointwise inequality holds on a ball $B_\rho(x_0)$, for a positive constant $\rho$ less than the injectivity radius of $g$ on $X$,
\begin{align*}
|d_A^{*_g}F_A - d_A^{*_{g_0}}F_A|
&\leq c\|g - g_0\|_{C^2(\bar B_\rho(x_0))} \left( r|F_A| + r^2 |\nabla_A F_A| + r^2 |d_A *_{g_0} F_A|\right)
\\
&\leq c\|\Riem_g\|_{C(\bar B_\rho(x_0))} \left( r|F_A| + r^2 |\nabla_A F_A| + r^2 |d_A^{*_{g_0}} F_A|\right).
\end{align*}
and thus,
\begin{equation}
\label{eq:Pointwise_norm_dA*curvedmetricFA_minus_dA*flatmetricFA}
|d_A^{*_g}F_A - d_A^{*_{g_0}}F_A|
\leq
c\|\Riem_g\|_{C(\bar B_\rho(x_0))} \left( r|F_A| + r^2 |\nabla_A F_A| + r^2 |d_A^{*_{g_0}} F_A|\right)
\quad\hbox{on } X.
\end{equation}
Therefore, we obtain an estimate for the $L^2$ norm of the difference between the Yang-Mills energy functional gradients for the metrics $g$ and $g_0$ on $B_\rho(x_0)$,
\begin{multline}
\label{eq:L2_norm_ball_dA*curvedmetricFA_minus_dA*flatmetricFA}
\|d_A^{*_g}F_A - d_A^{*_{g_0}}F_A\|_{L^2(B_\rho(x_0))}
\leq
c\|\Riem_g\|_{C(\bar B_\rho(x_0))} \left( \rho \|F_A\|_{L^2(B_\rho(x_0))} \right.
\\
+ \left. \rho^2 \|\nabla_A F_A\|_{L^2(B_\rho(x_0))}
+ \rho^2 \|d_A^{*_{g_0}} F_A\|_{L^2(B_\rho(x_0))} \right).
\end{multline}
Unfortunately, Lemma \ref{lem:Kozono_Maeda_Naito_4-2_global} cannot be used to estimate
$$
\int_0^T \|\nabla_{A(t)} F_A(t)\|_{L^2(B_\rho(x_0))}^2 \,dt,
$$
when $A(t)$ is a solution to Yang-Mills gradient flow, since the hypothesis \eqref{eq:Lemma_Kozono_Maeda_Naito_4-2_global_condition_small_energy_on_small_balls} may fail at $x_0$. Instead, we shall circumvent the problem of estimating \eqref{eq:L2_norm_ball_dA*curvedmetricFA_minus_dA*flatmetricFA}, when $A$ is replaced with Yang-Mills gradient flow, $A(t)$, by assuming that $g$ is locally conformally flat near $x_0$ and proving that flow would bubble at a point $x_0$ for the given Riemannian metric on $X$ at time $T<\infty$ if and only if that is true for the flow defined by a small perturbation of the Riemannian metric.

\section[Uhlenbeck and bubble-tree limits for Yang-Mills gradient flow]{Uhlenbeck and bubble-tree limits for Yang-Mills gradient flow over a four-dimensional manifold}
\label{sec:Uhlenbeck_and_bubbletree_limits_Yang-Mills_gradient_flow_four-manifold}
In Section \ref{subsec:Kozono_Maeda_Naito_theorems_5-1_through_5-6_Schlatter_theorems_1-2_and_1-3_Struwe_theorems_2-3_and_2-4}, we summarize the weak global existence results of Kozono, Maeda, and Naito \cite{Kozono_Maeda_Naito_1995}, Schlatter \cite{Schlatter_1997}, and Struwe \cite{Struwe_1994}, with particular emphasis on the behavior of the solution to Yang-Mills gradient flow approaching a time $T \in (0, \infty]$ where energy bubbling may occur. In Section \ref{subsec:Bubbletree_limits_Yang-Mills_gradient_flow_four-manifold},
we describe the bubble-tree limit process for a solution to Yang-Mills gradient flow over a closed, four-dimensional, Riemannian manifold.

The bubble-tree convergence possibility was first noted by Sacks and Uhlenbeck \cite[p. 3]{Sacks_Uhlenbeck_1981} in the context of sequences of harmonic maps from $S^2$, even before further development by Taubes in the context of sequences of Yang-Mills connections in dimension four \cite[Section 5]{TauFrame}, Parker and Wolfson in the context of sequences of harmonic maps from Riemann surfaces \cite{ParkerHarmonic, ParkerWolfson}, and the author in the context of sequences of anti-self-dual connections \cite{FeehanGeometry}.

\subsection{Condition (C) of Palais and Smale}
\label{subsec:Palais-Smale_condition_C}
For a base manifold, $X$, it is known from research of Daskalopoulos \cite{Daskalopoulos_1992} when $X$ has dimension two and R\r{a}de \cite{Rade_1992} when $X$ has dimension two or three, that the Yang-Mills energy functional obeys a gauge-equivariant version the condition (C) of Palais and Smale \cite{Daskalopoulos_1992, Rade_1992}. This is a consequence of Uhlenbeck's Weak Compactness Theorem \cite[Theorems 1.5 and 3.6]{UhlLp} for a sequence of connections with bounded $L^2$-energy.

\begin{prop}[Palais-Smale Condition C for the Yang-Mills energy functional in dimensions two and three]
\label{prop:Daskalopoulos_4-1}
\cite[Proposition 4.1]{Daskalopoulos_1992}, \cite[Proposition 7.1]{Rade_1992}
Let $G$ be a compact Lie group and $P$ be a principal $G$-bundle with $C^\infty$ reference connection, $A_1$, over a closed, Riemannian manifold of dimension $2$ or $3$. Let $\{A_m\}_{m\in\NN}$ be a sequence of $H^2$-connections on $P$ with
$$
\sup_{m\in\NN}\sE(A_m) < \infty \quad\hbox{and}\quad d_{A_m}^*F_{A_m} \to 0 \quad\hbox{as } m \to \infty.
$$
Then there is a sequence of $H^3$-gauge transformations, $\{\varphi_m\}_{m\in \NN} \subset \Aut P$, and a $C^\infty$ critical point, $A_\infty$ on $P$, such that, after passing to a subsequence,
$$
\varphi_m^*A_m \to A_\infty \quad\hbox{strongly in } H_{A_1}^1(X) \quad\hbox{as } m \to \infty.
$$
\end{prop}

When $X$ has the critical dimension four, of course, Proposition \ref{prop:Daskalopoulos_4-1} no longer holds. However, following Taubes \cite{TauFrame}, one can instead try to exploit the underlying reason for the failure of the Palais-Smale Condition C. We begin by examining the nature of the failure of the Palais-Smale Condition C in the following subsections.

\subsection{Uhlenbeck limits for Yang-Mills gradient flow over a four-dimensional manifold}
\label{subsec:Kozono_Maeda_Naito_theorems_5-1_through_5-6_Schlatter_theorems_1-2_and_1-3_Struwe_theorems_2-3_and_2-4}
In this section, we review the results of Kozono, Maeda, and Naito \cite[Theorems 5.1, 5.3, 5.4, and 5.6]{Kozono_Maeda_Naito_1995}, Schlatter \cite[Theorems 1.2 and 1.3]{Schlatter_1997}, and Struwe \cite[Theorem 2.3]{Struwe_1994}. The \cite[Theorems 1.2 and 1.3]{Schlatter_1997} due to Schlatter were conjectured earlier by Struwe as \cite[Theorem 2.4]{Struwe_1994} and are based in turn on Struwe's results for harmonic-map gradient-like flow for maps of a Riemann surface into a closed Riemannian manifold \cite[Theorems 4.1, 4.2, and 4.3]{Struwe_1985}.

Let $T = T(A_0)$ be the maximal existence time for a strong solution to the Yang-Mills gradient flow equation \eqref{eq:Yang-Mills_gradient_flow_equation} with initial condition \eqref{eq:Yang-Mills_heat_or_gradient_flow_equation_initial_condition} defined by a connection, $A_0$ on $P$, of class $W^{1,p}$ with $p \geq 1$. Kozono, Maeda, and Naito, Schlatter, and Struwe all allow $A_0$ to be of class $H^1$. Theorems \ref{thm:Kozono_Maeda_Naito_5-1} and \ref{thm:Kozono_Maeda_Naito_5-3} below may be viewed as analogues of Sedlacek's \cite[Theorem 3.1]{Sedlacek} and Taubes' \cite[Proposition 4.5]{TauPath}. See \cite[Appendix]{TauPath} for a careful and detailed analysis of the convergence of a Palais-Smale sequence of gauge-equivalence classes connections, taking into account the important distinction between sequences of gauge transformations of class $H^3$ and class $H^2$ which, in dimension four, may be discontinuous.

\begin{thm}[Characterization of the maximal lifetime and formation of bubble singularities for a solution to the Yang-Mills gradient flow over $X$]
\label{thm:Kozono_Maeda_Naito_5-1}
\cite[Theorem 5.1]{Kozono_Maeda_Naito_1995}, \cite[Theorem 2.3]{Struwe_1994}
Let $(X, g)$ be a closed, four-dimensional, Riemannian, smooth manifold. Then there is a positive constant $\eps_1$ with the following significance.
Let $A_0$ be a connection of class $H^1$ on a principal $G$-bundle $P$ over $X$ and let $A$ be a classical solution
to the Yang-Mills gradient flow equation \eqref{eq:Yang-Mills_gradient_flow_equation} on a maximal interval $(0, T)$ with initial condition \eqref{eq:Yang-Mills_heat_or_gradient_flow_equation_initial_condition} defined by $A_0$. Then the following hold:
\begin{enumerate}
\item The maximal lifetime, $T = T(A_0)$, is characterized by\footnote{We omit the square root appearing in the corresponding expression appearing in \cite[Theorem 5.1]{Kozono_Maeda_Naito_1995} since we also omit the square root in our redefinition \eqref{eq:Kozono_Maeda_Naito_4-1} of $\eps(x,r)$ in \cite[Equation (4.1)]{Kozono_Maeda_Naito_1995}.}
$$
\limsup_{t\nearrow T}\sup_{x \in X} \int_{B_r(x)} |F_A(t)|^2 \, d\vol_g \geq \eps_1, \quad\forall\, r \in (0, R_0],
$$
where the positive constants $\eps_1$ and $R_0$ are as in Section \ref{subsec:Kozono_Maeda_Naito_4}.

\item If $T < \infty$, then the solution, $A(t)$, is smooth on $X\times (0, T]$ except for a finite set of points, $\Sigma_T := \{(x_l, T): 1 \leq l \leq L\}$.

\item A singular point, $(x_l , T)$, is characterized by
$$
\limsup_{t\nearrow T} \int_{B_r(x_l)} |F_A(t)|^2 \, d\vol_g \geq \eps_1, \quad\forall\, r \in (0, R_0].
$$

\item The energy, $\sE(A(t))$, is a non-increasing function of $t \in [0,T]$.
\end{enumerate}
\end{thm}

\begin{proof}
We include selected details from the proofs provided by Kozono, Maeda, and Naito \cite[Theorem 5.1]{Kozono_Maeda_Naito_1995} and Struwe \cite[Theorem 2.3]{Struwe_1994} since we shall need to exploit those details in the sequel. As in \cite[Section 7]{Struwe_1994}, we argue by contradiction and suppose that there exists a constant $r \in (0, 1]$ such that \eqref{eq:Struwe_15} holds (see Lemma \ref{lem:Struwe_3-4}). But then, for any fixed $C^\infty$ reference connection $A_1$ on $P$ and writing $A(t) = A_1 + a(t)$ for $t \in [0, T)$, we note that Lemma \ref{lem:Schlatter_2-4_and_Struwe_3-6} implies $a \in C_{\loc}((0, T]; H_{A_1}^1(X; \Lambda^1\otimes \ad P))$ and, in particular, $\lim_{t\to T} a(t) = a(T) \in H_{A_1}^1(X; \Lambda^1\otimes \ad P)$. Consequently, for $t_0 \in (0, T]$ sufficiently large, there exists a solution $\breve A(t)$
for $t \in [t_0, t_1)$ to the Yang-Mills heat equation \eqref{eq:Yang-Mills_heat_equation_as_perturbation_rough_Laplacian_plus_one_heat_equation} with initial data, $\breve A(t_0) = A(T)$ and $t_1 > T$.
Hence, by Donaldson's version of the DeTurck trick (see Section \ref{subsec:Struwe_4-4_Donaldson-DeTurck trick}), there exists a family of gauge transformations, $\Phi(t) \in \Aut P$
for $t \in (t_0, t_1)$, such that and $\lim_{t\to t_0}\Phi(t) = \id_P$ and $\bar A(t) = \Phi(t)^*\breve A(t)$ for $t \in (t_0, t_1)$ is a solution to the Yang-Mills gradient flow equation \eqref{eq:Yang-Mills_gradient_flow_equation} with initial data $\bar A(t_0) = A(t_0)$,
contradicting the maximality of $T$.

The fact that the number, $L$, of bubble points $x_l \in \Sigma_T$ is finite follows from \cite[Lemma 5.2]{Kozono_Maeda_Naito_1995}, whose proof is based in turn on the proof of the corresponding result for the evolution equation for harmonic maps \cite[page 577]{Struwe_1985}.
\end{proof}

It is interesting to compare Theorem \ref{thm:Kozono_Maeda_Naito_5-3} below with the description \cite[Section 4.4]{DK} due to Donaldson and Kronheimer for the Uhlenbeck compactification of the moduli space of anti-self-dual connections on a principal $\SU(2)$-bundle $P$ over $X$ and of Taubes for the limiting behavior of \emph{good sequences} of connections (in the sense of \cite[Definition 4.1]{TauPath}) on a principal $G$-bundle $P$ with bounded energy and $L^2$-small self-dual curvature \cite[Propositions 4.4 and 4.5]{TauPath}, \cite[Proposition 5.1]{TauPath, TauFrame}. Recall that a sequence of connections, $\{A^j\}_{j\in\NN}$, on $P$ is \emph{good} in the sense of \cite[Definition 4.1]{TauPath}, \cite[Proposition 5.1]{TauPath} if
$$
\sup_{j\in\NN} \sE(A_j) < \infty \quad\hbox{and}\quad \|\sE'(A_j)\| \to 0 \quad\hbox{as } j \to \infty,
$$
where we regard the gradient of $\sE$ at $A_j$ as an element $\sE'(A_j) \in H^{-1}_{A_j}(X;\Lambda^1\otimes\ad P)$ and so its norm\footnote{There is a typographical error in the definition \cite[Equation (2.10)]{TauPath} of the dual space norm --- the infimum should be replaced by a supremum.}
is given by \cite[pages 340--341]{TauFrame}
$$
\|\sE'(A_j)\| = \|\sE'(A_j)\|_{H^{-1}_{A_j}(X)}
= \sup_{\begin{subarray}{c} b \in H_{A_j}(X;\Lambda^1\otimes\ad P) \\ \|b\|_{H^1_{A_j}(X)} \leq 1\end{subarray}} \frac{|\sE'(A_j)(b)|}{\|b\|_{H^1_{A_j}(X)}},
$$
where, as usual,
$$
\|b\|_{H^1_{A_j}(X)} = \left( \|\nabla_{A_j}b\|_{L^2(X)}^2 + \|b\|_{L^2(X)}^2 \right)^{1/2}.
$$
See Section \ref{sec:Taubes_1982_Appendix} for the classification of principal $G$-bundles, $P$ over $X$, where $G$ is a compact, connected, semi-simple Lie group, together with related Chern-Weil formulae. If $K$ is a compact Hausdorff space, then $(C(K;\RR))'$, the dual space of $C(K;\RR)$, may be identified with the space of Radon measures on $(K, \sB(K))$, where $\sB(K)$ is sigma algebra of Borel sets in $K$.


\begin{thm}[Uhlenbeck convergence over $X$ for a solution to the Yang-Mills gradient flow approaching a singularity at time $T < \infty$]
\label{thm:Kozono_Maeda_Naito_5-3}
\cite[Theorems A, 5.3, 5.6, and 7.1]{Kozono_Maeda_Naito_1995}, \cite[Theorem 1.2]{Schlatter_1997}, and \cite[Theorem 2.4]{Struwe_1994}
Assume the hypotheses of Theorem \ref{thm:Kozono_Maeda_Naito_5-1} and, in addition, that $T < \infty$. Then there exists a principal $G$-bundle $P^\flat$ over $X$, a connection, $A^\flat$, on $P^\flat$ of class $H^2$, and constants, $\sE_1, \ldots, \sE_L$ satisfying\footnote{To be consistent with our definition \eqref{eq:Yang-Mills_energy_functional} of the energy of a connection, it would be logical to include a factor of $1/2$ in the definition of the multiplicities but we omit that in favor of notational simplicity.}
$$
\eps_1 \leq \sE_l \leq \int_X |F_A(0)|^2 \, d\vol_g, \quad 1\leq l \leq L,
$$
and defined by the limit suprema,
\begin{equation}
\label{eq:Limit_r_to_zero_limsup_t_to_T_energy_over_ball_equals_E}
\sE_l := \lim_{r\to 0} \limsup_{t\nearrow T} \int_{B_r(x_l)} |F_A(t)|^2 \, d\vol_g, \quad 1 \leq l \leq L,
\end{equation}
such that the following hold:
\begin{enumerate}
\item
\label{item:Theorem_Kozono_Maeda_Naito_5-3_T_is_finite_Wkploc_strong_convergence_At}
For any fixed $C^\infty$ reference connection, $A_1$, on $P$ and $k=2$, $p=2$ or $k=1$, $p\in [2,\infty)$, one has $A(t) \to A^\flat$ strongly in $W_{A_1,\loc}^{k,p}(X \less \Sigma_T; \Lambda^1\otimes\ad P)$ as $t\nearrow T$;

\item
\label{item:Yang-Mills_gradient_flow_finite_time_bubbling_obstruction_preserved}
The principal $G$-bundle obstruction is preserved, $\eta(P^\flat) = \eta(P)$;

\item
\label{item:Limit_r_to_zero_t_to_T_energy_over_ball_equals_E}
The positive constants, $\sE_l$, are computed by the limits
\begin{equation}
\label{eq:Limit_r_to_zero_t_to_T_energy_over_ball_equals_E}
\sE_l = \lim_{r\to 0} \lim_{t\nearrow T} \int_{B_r(x_l)} |F_A(t)|^2 \, d\vol_g, \quad 1 \leq l \leq L;
\end{equation}

\item
\label{item:Theorem_Kozono_Maeda_Naito_5-3_T_is_finite_weak_convergence_FAt_measures}
$|F_A(t)|^2 \rightharpoonup |F_{A^\flat}|^2 + \sum_{l=1}^L \sE_l\, \delta_{x_l}$ in $(C(X;\RR))'$ as $t\nearrow T$, where $\delta_{x_l}$ denotes the Dirac delta measure with unit mass centered at the point $x_l \in X$.
\end{enumerate}
\end{thm}

\begin{proof}
The proof is omitted in \cite{Kozono_Maeda_Naito_1995, Struwe_1994} and only briefly described in \cite{Schlatter_1997}, so we include additional details here since we shall need them in the sequel.

Consider Item \eqref{item:Theorem_Kozono_Maeda_Naito_5-3_T_is_finite_Wkploc_strong_convergence_At}. Following Schlatter in his proof of \cite[Theorem 1.2]{Schlatter_1997}, let $U \Subset X\less\Sigma_T$ be a precompact open subset. Theorem \ref{thm:Kozono_Maeda_Naito_5-1} implies that there are constants, $r_1 \in (0, R_0]$ and $T_0 \in (0, T)$ (depending on $U$), such that (see, for example, the proof of Lemma \ref{lem:Kozono_Maeda_Naito_4-2_global} by a covering argument)
\begin{equation}
\label{eq:sup_t_and_x_in_(delta_T)_and_U_energy_Yang-Mills_gradient_flow}
\sup_{(t,x) \in (T_0, T)\times U} \int_{B_r(x)} |F_A(t)|^2\,d\vol_g < \eps_1, \quad\forall\, r \in (0, r_1].
\end{equation}
According to Lemma \ref{lem:Schlatter_2-4_and_Struwe_3-6_local}, we then have
$$
A \in C_{\loc}((0,T]; H_{A_1}^1(U; \Lambda^1\otimes\ad P)).
$$
In particular, this ensures that $A(t) \restriction U \to A(T) \restriction U$ in $H_{A_1}^1(U; \Lambda^1\otimes\ad P)$ as $t \nearrow T$. Since $U \Subset X\less\Sigma_T$ was arbitrary, this yields
\begin{equation}
\label{eq:Kozono_Maeda_Naito_5-3_A(t)_converges_to_A(T)_strongly_in_H1loc_on_X_away_from_bubble_points}
A(t) \to A(T) \quad\hbox{strongly in } H_{A_1,\loc}^1(X\less\Sigma_T; \Lambda^1\otimes\ad P) \quad\hbox{as } t \nearrow T.
\end{equation}
To obtain the indicated stronger convergence of $A(t)$ to $A(T)$ over precompact subsets $U \Subset X\less \Sigma_T$ as $t \nearrow T$, we apply Lemma \ref{lem:Schlatter_2-4_and_Struwe_3-6_local_H2_or_W1p}.

Uhlenbeck's Removable Singularity Theorem for finite-energy Sobolev connections, Theorem \ref{thm:Uhlenbeck_removable_singularity_Sobolev}, and Remark \ref{rmk:Uhlenbeck_removable_singularity_Sobolev_Wkp} ensure the existence of a principal $G$-bundle, $P^\flat$, over $X$ and a principal $G$-bundle isomorphism, $\Phi: P^\flat\restriction X\less \Sigma_T \cong P\restriction X\less \Sigma_T$ of class $H^3$, and a connection, $A^\flat$, of class $H^2$ on $P^\flat$ such that $\Phi^*A(T) = A^\flat$ on $P^\flat\restriction X\less \Sigma_T$. This completes the proof of Item \eqref{item:Theorem_Kozono_Maeda_Naito_5-3_T_is_finite_Wkploc_strong_convergence_At}.

Consider Item \eqref{item:Yang-Mills_gradient_flow_finite_time_bubbling_obstruction_preserved}. It follows from \cite[Theorem 5.5]{Sedlacek}, or more precisely its proof, that $\eta(P^\flat) = \eta(P)$. The only relevant difference between our setting and that of \cite[Section 5]{Sedlacek} is that the limiting connection, $A^\flat$, need not be Yang-Mills (when $T < \infty$), but that fact is only used in Sedlacek's application of Uhlenbeck's Removable Singularity Theorem for Yang-Mills connections, which we replace here by her Removable Singularity Theorem for finite-energy Sobolev connections. This completes the proof of Item \eqref{item:Yang-Mills_gradient_flow_finite_time_bubbling_obstruction_preserved}.

Consider Items \eqref{item:Limit_r_to_zero_t_to_T_energy_over_ball_equals_E} and \eqref{item:Theorem_Kozono_Maeda_Naito_5-3_T_is_finite_weak_convergence_FAt_measures}. As a consequence of the definition \eqref{eq:Limit_r_to_zero_limsup_t_to_T_energy_over_ball_equals_E} of the constants, $\sE_l$, there is a sequence, $\{t_m\}_{m\in\NN} \subset (0, T)$, such that $t_m \to T$ as $m \to \infty$ and
\begin{equation}
\label{eq:Limit_r_to_zero_m_to_infinity_energy_tm_over_ball_equals_E}
\sE_l =  \lim_{r\to 0} \lim_{m\to\infty} \int_{B_r(x_l)} |F_A(t_m)|^2 \, d\vol_g, \quad 1 \leq l \leq L.
\end{equation}
In order to prove the stronger assertion \eqref{eq:Limit_r_to_zero_t_to_T_energy_over_ball_equals_E} we need, of course, to show that the limits on the right-hand side of \eqref{eq:Limit_r_to_zero_m_to_infinity_energy_tm_over_ball_equals_E} are equal to $\sE_l$ for \emph{any} choice of sequence $\{t_m\}_{m\in\NN} \subset (0, T)$ such that $t_m \nearrow T$.

By the expression \eqref{eq:Limit_r_to_zero_m_to_infinity_energy_tm_over_ball_equals_E} for $\sE_l$, we may fix $m_0$ large enough and $\rho_0 \in (0, R_0]$ small enough such that
\begin{equation}
\label{eq:Energy_Yang-Mills_gradient_flow_over_ball_at_tm_near_E_for_m_geq_m0}
\left|\sE_l - \int_{B_\rho(x_l)} |F_A(t_m)|^2 \,d\vol_g\right| < \frac{\eps}{8L}, \quad\forall\, m \geq m_0, \quad \rho \in (0, \rho_0].
\end{equation}
Suppose $f \in C(X)$ and $\eps>0$. Fix a positive constant, $\zeta$, such that
\begin{equation}
\label{eq:zeta_lessthan_epsilon_over_twice_initial_energy}
\zeta \int_X |F_A(0)|^2 \,d\vol_g < \frac{\eps}{4}.
\end{equation}
We may decrease the size of $\rho_0$, if necessary, to also ensure that
\begin{equation}
\label{eq:Max_bubble_balls_C0_f_minus_fxl_lessthan_zeta}
\max_{1\leq l \leq L} \|f - f(x_l)\|_{C(\bar B_\rho(x_l))} < \zeta, \quad \forall\, \rho \in (0, \rho_0].
\end{equation}
Writing
\begin{equation}
\label{eq:Yang-Mills_gradient_flow_energy_f_decomposition_union_bubble_balls_and_complement}
\int_X f |F_A(t)|^2 \,d\vol_g
= \int_{X\less \cup_{l=1}^LB _\rho(x_l)} f |F_A(t)|^2 \,d\vol_g
+ \sum_{l=1}^L \int_{B_\rho(x_l)} f |F_A(t)|^2 \,d\vol_g,
\end{equation}
we wish to compare the right-hand side of the preceding identity with
$$
\int_X f |F_{A^\flat}|^2 \,d\vol_g + \sum_{l=1}^L \sE_l f(x_l),
$$
as $t \nearrow T$.

For the first term on the right-hand side of \eqref{eq:Yang-Mills_gradient_flow_energy_f_decomposition_union_bubble_balls_and_complement}, the convergence in Item \eqref{item:Theorem_Kozono_Maeda_Naito_5-3_T_is_finite_Wkploc_strong_convergence_At} implies that
$$
\int_{X\less \cup_{l=1}^LB _\rho(x_l)} f |F_A(t)|^2 \,d\vol_g \to \int_{X\less \cup_{l=1}^LB _\rho(x_l)} f |F_{A^\flat}|^2 \,d\vol_g,
\quad\hbox{as } t \to T,
$$
and so there is a large enough $T_0 \in (0, T)$ such that
\begin{equation}
\label{eq:Sup_t_in_0T_complement_union_bubble_balls_f_Yang-Mills_gradient_flow_energy_minus_FAflat}
\left|\int_{X\less \cup_{l=1}^LB_\rho(x_l)} f \left( |F_A(t)|^2 - |F_{A^\flat}|^2 \right) \,d\vol_g \right| < \frac{\eps}{4},
\quad \forall\, t \in [T_0, T).
\end{equation}
For the second term on the right-hand side of \eqref{eq:Yang-Mills_gradient_flow_energy_f_decomposition_union_bubble_balls_and_complement}, we have
\begin{align*}
\sup_{t\in(0,T)} \left|\sum_{l=1}^L \int_{B_\rho(x_l)} (f - f(x_l))|F_A(t)|^2 \,d\vol_g \right|
&< \zeta \sup_{t\in(T_0,T)} \sum_{l=1}^L\int_{B_\rho(x_l)} |F_A(t)|^2 \,d\vol_g
\\
&\leq \zeta \sup_{t\in(0,T)} \int_X |F_A(t)|^2 \,d\vol_g
\\
&= \zeta \int_X |F_A(0)|^2 \,d\vol_g \quad\hbox{(by Lemma \ref{lem:Kozono_Maeda_Naito_4-1})}.
\end{align*}
Therefore, by our choice of $\zeta$ in \eqref{eq:zeta_lessthan_epsilon_over_twice_initial_energy}, we obtain
\begin{equation}
\label{eq:Sup_t_in_0T_sum_bubble_balls_f_minus_fxl_Yang-Mills_gradient_flow_energy}
\sup_{t\in(0,T)} \left| \sum_{l=1}^L \int_{B_\rho(x_l)} (f - f(x_l))|F_A(t)|^2 \,d\vol_g \right| < \frac{\eps}{4}.
\end{equation}
We may further decrease the size of $\rho_0 \in (0, R_0]$ (depending at least on $A^\flat$, $f$, $L$, and the points $x_l \in X$), if necessary, to ensure that
\begin{equation}
\label{eq:Sum_bubble_balls_f_Aflat_energy_lessthan_epsilon_over_4}
\sum_{l=1}^L \int_{B_\rho(x_l)} |f| |F_{A^\flat}|^2 \,d\vol_g < \frac{\eps}{4},
\quad\forall\, \rho \in (0, \rho_0], \quad 1 \leq l \leq L.
\end{equation}
In order to prove that the following terms become suitably close as $t \nearrow T$ and then $\rho \searrow 0$,
$$
\int_{B_\rho(x_l)} f(x_l) |F_A(t)|^2 \,d\vol_g \quad\hbox{and}\quad \sE_l f(x_l), \quad\hbox{for } 1 \leq l \leq L,
$$
it obviously suffices to prove that \eqref{eq:Limit_r_to_zero_t_to_T_energy_over_ball_equals_E} holds. We begin by recalling that (just as in the proof of Lemma \ref{lem:Kozono_Maeda_Naito_4-1}),
\begin{align*}
\frac{1}{2}\frac{d}{dt} \int_X |F_A(t)|^2 \,d\vol_g &= \int_X \left\langle \frac{\partial F_A}{\partial t}(t), F_A(t)\right\rangle \,d\vol_g
\\
&= -\int_X \langle \Delta_A F_A(t), F_A(t) \rangle \,d\vol_g
\\
&= -\int_X \langle d_Ad_A^* F_A(t), F_A(t) \rangle \,d\vol_g
\\
&= -\int_X |d_A^* F_A(t)|^2 \,d\vol_g \leq 0,
\end{align*}
and thus $\sE(A(t))$ is a non-increasing function of $t \in [0, T)$. Fix a $C^\infty$ cut-off function $\kappa:\RR\to [0,1]$ such that $\kappa(s) = 0$ for $s \geq 2$ and $\kappa(s) = 1$ for $s \leq 1$. Now define $C^\infty$ cut-off functions, $\chi_{l,\rho}:X\to [0,1]$, depending on the points $x_l \in X$ and the small parameter $\rho \in (0, R_0]$, to be further constrained later, by setting
$$
\chi_{l,\rho}(x) := \kappa(\dist_g(x, x_l)/\rho), \quad\forall\, x \in X, \quad 1 \leq l \leq L.
$$
From the proof of Lemma \ref{lem:L2_and_L4_nabla_chi_and_L2_nabla2_chi_bounds} (see Remark \ref{rmk:Construction_cut-off_function_chi_N_rho}), there is a positive constant, $c$, depending at most on the Riemannian metric, $g$, on $X$ such that
\begin{equation}
\label{eq:L4_nabla_chi_bound_4manifold_basic}
\|\cov\chi_{l,\rho}\|_{L^4(X)} \leq c, \quad\forall\,\rho \in (0, R_0], \quad 1 \leq l \leq L.
\end{equation}
Denote $\Omega(x_l; \rho, 2\rho) := B_{2\rho}(x_l) \less \bar B_\rho(x_l)$. From the convergence given by Item \eqref{item:Theorem_Kozono_Maeda_Naito_5-3_T_is_finite_Wkploc_strong_convergence_At}, we have
\begin{equation}
\label{eq:limit_t_to_T_energy_Yang-Mills_gradient_flow_over_annulus}
\lim_{t \nearrow T} \int_{\Omega(x_l; \rho, 2\rho)} |F_A(t)|^2\,d\vol_g
= \int_{\Omega(x_l; \rho, 2\rho)} |F_{A^\flat}|^2\,d\vol_g, \quad 1 \leq l \leq L.
\end{equation}
In particular, given a small positive constant, $\eta$, to be determined later, there are a small enough constant, $\rho \in (0, \rho_0]$ (where $\rho_0$ is as in \eqref{eq:Energy_Yang-Mills_gradient_flow_over_ball_at_tm_near_E_for_m_geq_m0}) depending at least on $A^\flat$ and $\eta$, such that
$$
\int_{\Omega(x_l; \rho, 2\rho)} |F_{A^\flat}|^2\,d\vol_g < \frac{\eta}{2}, \quad 1 \leq l \leq L,
$$
and a large enough constant, $T_0 \in (0, T)$ depending at least on $\rho \in (0, \rho_0]$, such that
$$
\sup_{t \in (T_0, T)} \int_{\Omega(x_l; \rho, 2\rho)} \left| |F_A(t)|^2 - |F_{A^\flat}|^2\right| \,d\vol_g < \frac{\eta}{2}, \quad 1 \leq l \leq L,
$$
and hence,
\begin{equation}
\label{eq:sup_t_in_(T0_T)_energy_Yang-Mills_gradient_flow_over_annulus}
\sup_{t \in (T_0, T)} \int_{\Omega(x_l; \rho, 2\rho)} |F_A(t)|^2\,d\vol_g < \eta, \quad 1 \leq l \leq L.
\end{equation}
Again from the convergence given by Item \eqref{item:Theorem_Kozono_Maeda_Naito_5-3_T_is_finite_Wkploc_strong_convergence_At}, there are a possibly larger constant, $T_0 \in (0, T)$ depending at least on $\rho \in (0, \rho_0]$, and a positive constant, $C$, depending at least on $A^\flat$ but independent of $\rho \in (0, \rho_0]$ or $T_0 \in (0, T)$, such that
\begin{equation}
\label{eq:sup_t_in_(T0_T)_Linfinity_annulus_FA_and_nablaA_FA}
\sup_{t \in (T_0, T)} \left( \|F_A(t)\|_{L^\infty(\Omega(x_l; \rho, 2\rho))} + \|\nabla_AF_A(t)\|_{L^\infty(\Omega(x_l; \rho, 2\rho))} \right) \leq C,
\quad 1 \leq l \leq L.
\end{equation}
From the expression \eqref{eq:Warner_6-2} for $d_A^*$, we see that
\begin{align*}
\frac{1}{2}\frac{d}{dt} \int_X \chi_\rho^2 |F_A(t)|^2 \,d\vol_g &= -\int_X \langle d_Ad_A^* F_A(t), \chi_\rho^2 F_A(t) \rangle \,d\vol_g
\\
&= -\int_X \chi_\rho^2|d_A^* F_A(t)|^2 \,d\vol_g + 2\int_X \langle d_A^*F_A(t), *(d\chi_\rho \wedge *F_A(t)) \rangle \,d\vol_g.
\end{align*}
To further examine the second term on the right-hand side of the preceding identity, observe that
\begin{align*}
{}& \left|\int_X \langle d_A^*F_A(t), *(d\chi_\rho \wedge *F_A(t)) \rangle \,d\vol_g\right|
\\
&\quad \leq
\|d_A^*F_A(t)\|_{L^4(\Omega(x_l; \rho, 2\rho))} \|d\chi_\rho\|_{L^4(X)} \|F_A(t)\|_{L^2(\Omega(x_l; \rho, 2\rho))}
\\
&\quad \leq C\sqrt{\eta}, \quad\forall\, t\in [T_0, T)
\quad\hbox{(by \eqref{eq:L4_nabla_chi_bound_4manifold_basic}, \eqref{eq:sup_t_in_(T0_T)_energy_Yang-Mills_gradient_flow_over_annulus}, and \eqref{eq:sup_t_in_(T0_T)_Linfinity_annulus_FA_and_nablaA_FA})}.
\end{align*}
Hence, by combining the preceding inequalities,
$$
\frac{d}{dt} \int_X \chi_\rho^2 |F_A(t)|^2 \,d\vol_g \leq 4C\sqrt{\eta}, \quad\forall\, t\in [T_0, T),
$$
and so, for any $t', t'' \in [T_0, T)$ such that $t' \leq t''$, we have
$$
\int_X \chi_\rho^2 |F_A(t'')|^2 \,d\vol_g \leq \int_X \chi_\rho^2 |F_A(t')|^2 \,d\vol_g + 4C\sqrt{\eta} (t'' - t').
$$
Consequently, by definition of $\chi_\rho$, we obtain the approximate monotonicity relationship,
$$
\int_{B_\rho(x_l)} |F_A(t'')|^2 \,d\vol_g
\leq \int_{B_\rho(x_l)} |F_A(t')|^2 \,d\vol_g + \int_{\Omega(x_l; \rho, 2\rho))} |F_A(t'')|^2 \,d\vol_g + 4C\sqrt{\eta} (t'' - t'),
$$
and combining this with the energy bound \eqref{eq:sup_t_in_(T0_T)_energy_Yang-Mills_gradient_flow_over_annulus} yields,
\begin{multline}
\label{eq:Approximate_monotonicity_energy_Yang-Mills_gradient_flow_over_ball}
\int_{B_\rho(x_l)} |F_A(t'')|^2 \,d\vol_g \leq \int_{B_\rho(x_l)} |F_A(t')|^2 \,d\vol_g + \eta + 4C\sqrt{\eta}(t'' - t'),
\\
T_0 \leq t' \leq t'' < T, \quad 1 \leq l \leq L.
\end{multline}
Hence, for any $t \in [t_{m_0}, T)$ (where $m_0$ is as in \eqref{eq:Energy_Yang-Mills_gradient_flow_over_ball_at_tm_near_E_for_m_geq_m0}),
\begin{align*}
\int_{B_\rho(x_l)} |F_A(t)|^2 \,d\vol_g
&\leq \int_{B_\rho(x_l)} |F_A(t_{m_0})|^2 \,d\vol_g + \eta + 4C\sqrt{\eta}(t - t_{m_0})
\quad\hbox{(by \eqref{eq:Approximate_monotonicity_energy_Yang-Mills_gradient_flow_over_ball})}
\\
&< \sE_l + \frac{\eps}{8L} + \eta + 4C\sqrt{\eta}(t - t_{m_0}), \quad 1 \leq l \leq L
\quad\hbox{(by \eqref{eq:Energy_Yang-Mills_gradient_flow_over_ball_at_tm_near_E_for_m_geq_m0}).}
\end{align*}
Let us now choose $\eta = \eta(C, L, T, \eps) >0$ small enough that,
\begin{equation}
\label{eq:half_eta_plus_2C_sqrt_eta_T_lessthan_epsilon_over_8L}
\eta + 4C\sqrt{\eta}\,T < \frac{\eps}{8L},
\end{equation}
and thus, by the preceding inequality,
$$
\int_{B_\rho(x_l)} |F_A(t)|^2 \,d\vol_g < \sE_l + \frac{\eps}{4L}, \quad\forall\, t \in [t_{m_0}, T), \quad 1 \leq l \leq L.
$$
On the other hand, for any $t \in [t_{m_0}, T)$, we may choose $m_1 \geq m_0$ large enough that $t_{m_1} \in [t, T)$ and then we observe that
\begin{align*}
\int_{B_\rho(x_l)} |F_A(t)|^2 \,d\vol_g
&\geq  \int_{B_\rho(x_l)} |F_A(t_{m_1})|^2 \,d\vol_g - \eta - 4C\sqrt{\eta}(t_{m_1} - t)
\quad\hbox{(by \eqref{eq:Approximate_monotonicity_energy_Yang-Mills_gradient_flow_over_ball})}
\\
&> \sE_l - \frac{\eps}{8L} - \eta - 4C\sqrt{\eta}(t_{m_1} - t)
\quad\hbox{(by \eqref{eq:Energy_Yang-Mills_gradient_flow_over_ball_at_tm_near_E_for_m_geq_m0})}
\\
&> \sE_l - \frac{\eps}{4L}, \quad\forall\, t \in [t_{m_0}, T)  \quad\hbox{(by \eqref{eq:half_eta_plus_2C_sqrt_eta_T_lessthan_epsilon_over_8L})}.
\end{align*}
Assembling the preceding upper and lower bounds gives,
\begin{equation}
\label{eq:Max_bubble_balls_E_minus_Yang-Mills_gradient_flow_energy_ball_rho}
\max_{1\leq l\leq L}\left|\sE_l - \int_{B_\rho(x_l)} |F_A(t)|^2 \,d\vol_g \right|
< \frac{\eps}{4L}, \quad\forall\, t \in [t_{m_0}, T).
\end{equation}
This yields the desired strong characterization \eqref{eq:Limit_r_to_zero_t_to_T_energy_over_ball_equals_E} of the positive constants, $\sE_l$ for $1\leq l\leq L$, and completes the proof of Item \eqref{item:Limit_r_to_zero_t_to_T_energy_over_ball_equals_E}.

Finally, by combining \eqref{eq:Sup_t_in_0T_complement_union_bubble_balls_f_Yang-Mills_gradient_flow_energy_minus_FAflat},
\eqref{eq:Sup_t_in_0T_sum_bubble_balls_f_minus_fxl_Yang-Mills_gradient_flow_energy},
\eqref{eq:Sum_bubble_balls_f_Aflat_energy_lessthan_epsilon_over_4}, and \eqref{eq:Max_bubble_balls_E_minus_Yang-Mills_gradient_flow_energy_ball_rho}, we discover that
\begin{align*}
{}& \left| \int_X f |F_A(t)|^2 \,d\vol_g - \int_X f |F_{A^\flat}|^2 \,d\vol_g - \sum_{l=1}^L \sE_l f(x_l) \right|
\\
&\quad \leq \int_{X\less \cup_{l=1}^LB_\rho(x_l)} |f| \left| |F_A(t)|^2 -  |F_{A^\flat}|^2 \right| \,d\vol_g
\\
&\qquad + \sum_{l=1}^L \int_{B_\rho(x_l)} |f - f(x_l)| |F_A(t)|^2 \,d\vol_g
\\
&\qquad + \sum_{l=1}^L |f(x_l)|\left|\int_{B_\rho(x_l)} |F_A(t)|^2 \,d\vol_g - \sE_l\right|
+ \sum_{l=1}^L \int_{B_\rho(x_l)} |f| |F_{A^\flat}|^2 \,d\vol_g
\\
&\quad < \eps, \quad \forall\, t \in [T_0 \vee t_{m_0}, T).
\end{align*}
The preceding inequality completes our verification of Item \eqref{item:Theorem_Kozono_Maeda_Naito_5-3_T_is_finite_weak_convergence_FAt_measures} and the proof of Theorem \ref{thm:Kozono_Maeda_Naito_5-3}.
\end{proof}

When $T = \infty$, then the conclusions of Theorem \ref{thm:Kozono_Maeda_Naito_5-3} take the following weaker form because we can no longer rely on Lemmata \ref{lem:Schlatter_2-4_and_Struwe_3-6_local} or \ref{lem:Schlatter_2-4_and_Struwe_3-6_local_H2_or_W1p} to give continuity of $\|A(t)\|_{H_{A_1}^1(U)}$ or $\|A(t)\|_{W_{A_1}^{1,p}(U)}$ as $t \nearrow T$, for $U \Subset X \less \Sigma_T$, when $T < \infty$.

\begin{thm}[Uhlenbeck convergence over $X$ for a solution to the Yang-Mills gradient flow approaching a singularity at time $T = \infty$]
\label{thm:Kozono_Maeda_Naito_5-3_T_is_infinite}
\cite[Theorems A, 5.6, and 7.1, and Remark following Theorem 5.3]{Kozono_Maeda_Naito_1995}, \cite[Theorems 1.2 and 1.3]{Schlatter_1997}, and \cite[Theorem 2.4]{Struwe_1994}
Assume the hypotheses of Theorem \ref{thm:Kozono_Maeda_Naito_5-1} and, in addition, that $T = \infty$. Then there exists a principal $G$-bundle, $P_\infty^\flat$ over $X$, a sequence of times, $\{t_m\}_{m \in \NN} \subset [0,\infty)$, such that $t_m \to \infty$ as $m \to \infty$, a sequence of bundle isomorphisms, $\Phi_m:P^\flat\restriction X\less \Sigma_\infty \cong P\restriction X\less \Sigma_\infty$, of class $H^3$, a smooth \emph{Yang-Mills} connection, $A_\infty^\flat$ on $P_\infty^\flat$, and constants, $\sE_1, \ldots, \sE_L$, satisfying
$$
\eps_1 \leq \sE_l \leq \int_X |F_A(0)|^2 \, d\vol_g, \quad 1\leq l \leq L,
$$
and defined by the limit suprema,
\begin{equation}
\label{eq:Limit_r_to_zero_limsup_t_to_infinity_energy_over_ball_equals_E}
\sE_l := \lim_{r\to 0} \limsup_{t\to \infty} \int_{B_r(x_l)} |F_A(t)|^2 \, d\vol_g, \quad 1 \leq l \leq L,
\end{equation}
such that the following hold:
\begin{enumerate}
\item
\label{item:Theorem_Kozono_Maeda_Naito_5-3_T_is_infinite_H2loc_weak_and_W1ploc_strong_convergence_Am}
For any fixed $C^\infty$ reference connection, $A_1$, on $P$ and $p\in [2,4)$, one has
\begin{subequations}
\label{eq:Kozono_Maeda_Naito_5-3_Phi(tm)A(tm)_converges_to_Ainfinity_flat}
\begin{align}
\label{eq:Kozono_Maeda_Naito_5-3_Phi(tm)A(tm)_converges_to_Ainfinity_flat_weakly_in_H2loc_on_X_away_from_bubble_points}
\Phi_m^*A(t_m) &\rightharpoonup A_\infty^\flat \quad\hbox{weakly in } H_{A_1,\loc}^2(X\less \Sigma_\infty; \Lambda^1\otimes\ad P),
\\
\label{eq:Kozono_Maeda_Naito_5-3_Phi(tm)A(tm)_converges_to_Ainfinity_flat_strongly_in_W1ploc_on_X_away_from_bubble_points}
\Phi_m^*A(t_m) &\to A_\infty^\flat \quad\hbox{strongly in } W_{A_1,\loc}^{1,p}(X\less \Sigma_\infty; \Lambda^1\otimes\ad P), \quad\hbox{as } m \to \infty.
\end{align}
\end{subequations}

\item
\label{item:Yang-Mills_gradient_flow_infinite_time_bubbling_obstruction_preserved}
The principal $G$-bundle obstruction is preserved, $\eta(P_\infty^\flat) = \eta(P)$;

\item
\label{item:Limit_r_to_zero_m_to_infinity_energy_over_ball_equals_E_T_is_infinite}
The positive constants, $\sE_l$, are computed by the limas,
\begin{equation}
\label{eq:Limit_r_to_zero_m_to_infinity_energy_over_ball_equals_E_T_is_infinite}
\sE_l = \lim_{r\to 0} \lim_{m \to \infty} \int_{B_r(x_l)} |F_A(t_m)|^2 \, d\vol_g, \quad 1 \leq l \leq L;
\end{equation}

\item
\label{item:Theorem_Kozono_Maeda_Naito_5-3_T_is_infinite_weak_convergence_FAm_measures}
$|F_A(t_m)|^2 \rightharpoonup |F_{A_\infty^\flat}|^2 + \sum_{l=1}^L \sE_l\, \delta_{x_l}$ in $(C(X;\RR))'$ as $m \to \infty$.
\end{enumerate}
\end{thm}

\begin{proof}
Item \eqref{item:Limit_r_to_zero_m_to_infinity_energy_over_ball_equals_E_T_is_infinite} is a consequence the definition \eqref{eq:Limit_r_to_zero_limsup_t_to_infinity_energy_over_ball_equals_E} of the constants, $\sE_l$. In our proof of Item \eqref{item:Theorem_Kozono_Maeda_Naito_5-3_T_is_infinite_H2loc_weak_and_W1ploc_strong_convergence_Am}, we shall pass to subsequences, with relabeling, of the choice of sequence $\{t_m\}_{m \in \NN} \subset [0,\infty)$ used to compute the limits \eqref{eq:Limit_r_to_zero_m_to_infinity_energy_over_ball_equals_E_T_is_infinite}.

Consider Item \eqref{item:Theorem_Kozono_Maeda_Naito_5-3_T_is_infinite_H2loc_weak_and_W1ploc_strong_convergence_Am}. Let $\{B_{r_n}(x_n^*)\}_{n\in \NN}$ be a countable open cover of $X\less \Sigma_\infty$ by open balls in $X\less \Sigma_\infty$ such that $\{B_{r_n/4}(x_n^*)\}_{n\in \NN}$ is a countable subcover. For any such ball, $B_r(x_*) \subset X\less \Sigma_\infty$, we may suppose by Theorem \ref{thm:Kozono_Maeda_Naito_5-1} that $r \in (0, R_0]$ is small enough and $t_0 \in (0, \infty)$ is large enough that
\begin{equation}
\label{eq:L-infinity_time_energy_A(t)_ball_is_small_t_geq_t_zero}
\sup_{t \in (t_0, \infty)} \int_{B_r(x_*)} |F_A(t)|^2\,d\vol_g < \eps_1.
\end{equation}
Theorem \ref{thm:Uhlenbeck_Lp_1-3} and Remarks \ref{rmk:Uhlenbeck_theorem_1-3_Wkp} and \ref{rmk:Uhlenbeck_theorem_1-3_apriori_estimate_scale_invariant} imply that there are a positive constant $c$ (depending at most on the Riemannian metric on $X$) and a family of local sections, $\sigma(t):B_r(x_*) \to P \restriction B_r(x_*)$ for $t \in (t_0, \infty)$, such that $\sigma(t)^*A(t) \in H_\Gamma^2(B_r(x_*); \Lambda^1\otimes\ad P)$ and obeys a Coulomb gauge condition and a spatially scale-invariant estimate, for $p \in [2, 4)$ and $q \in [4,\infty)$ defined by $1/p=1/4+1/q$,
\begin{subequations}
\label{eq:a(t)_Coulomb_gauge_and_W_1-q_bound_ball_t_geq_t_zero}
\begin{align}
\label{eq:d_Gamma*a(t)_is_zero_ball_t_geq_t_zero}
d_\Gamma^*\sigma(t)^*A(t) &= 0 \quad\hbox{on } B_r(x_*),
\\
\label{eq:W_Gamma_1-q_norm_a(t)_leq_Lp_norm_FA(t)_ball_t_geq_t_zero}
\|\sigma(t)^*A(t)\|_{L^q(B_r(x_*)} + \|\nabla_\Gamma \sigma(t)^*A(t)\|_{L^p(B_r(x_*)} &\leq c\|F_A(t)\|_{L^p(B_r(x_*)}, \quad\forall\, t \in (t_0, \infty),
\end{align}
\end{subequations}
where $\Gamma$ is the product connection on $B_r(x_*) \times G$.

We would like to apply Lemmata \ref{lem:Kozono_Maeda_Naito_4-3} and \ref{lem:Kozono_Maeda_Naito_4-4}
in order to control the spatial $L^\infty$ norms of $F_A(t)$ and $\nabla_AF_A(t)$, on the balls $B_r(x_*)$, uniformly for all $t\geq t_0$ and therefore use \eqref{eq:a(t)_Coulomb_gauge_and_W_1-q_bound_ball_t_geq_t_zero} to obtain a stronger Sobolev bound on the local connection one-forms, $\sigma(t)^*A(t)$, than that given in \eqref{eq:W_Gamma_1-q_norm_a(t)_leq_Lp_norm_FA(t)_ball_t_geq_t_zero}. However, the estimate constant, $C$, in those lemmata depend explicitly on a \emph{finite} time, $T$, to parlay a small energy bound for $A(t)$ over $B_r(x_*)$, uniform with respect to $t\in(0,T)$, into stronger Sobolev estimates for $F_A(t)$ over $B_{r/2}(x_*)$, uniform with respect to $t\in (\delta,T)$, for $0<\delta<T$. Thus, we instead consider the sequence of intervals, $(t_k, t_k + 2)$ with $t_k := t_0+k$ for $k \in \NN$, over which we have a uniformly small energy bound for $A(t)$ over $B_r(x_*)$ by \eqref{eq:L-infinity_time_energy_A(t)_ball_is_small_t_geq_t_zero} and obtain stronger Sobolev estimates for $F_A(t)$ over $B_{r/2}(x_*)$, uniform with respect to $t\in (t_k+1, t_k + 2)$, for all $k\in\NN$, with estimate constant, $C$, \emph{independent} of $k\in\NN$.

With the preceding remarks understood, Lemmata \ref{lem:Kozono_Maeda_Naito_4-3} and \ref{lem:Kozono_Maeda_Naito_4-4} provide a positive constant, $C$, such that
\begin{equation}
\label{eq:L_infty_time_and_space_FA_bound_good_ball}
\sup_{t \in (t_0+1, \infty)} \left( \|F_A(t)\|_{L^\infty(B_{r/2}(x_*))} + \|\nabla_AF_A(t)\|_{L^\infty(B_{r/2}(x_*))} \right) \leq C.
\end{equation}
Therefore, Lemma \ref{lem:Schlatter_2-7} yields a (possibly larger) positive constant, $C$, and the bound,
$$
\|\sigma(t)^*A(t)\|_{H_\Gamma^2(B_{r/4}(x_*))} \leq C, \quad\forall\, t \in (t_0+1, \infty).
$$
Hence, passing to a subsequence and relabeling, there are a sequence of times, $\{t_m\} \subset (t_0+1, \infty)$, and a local connection one-form, $a_\infty \in H_\Gamma^2(B_{r/4}(x_*); \Lambda^1\otimes\ad P)$, such that $t_m \to \infty$ as $m \to \infty$ and
\begin{subequations}
\label{eq:Kozono_Maeda_Naito_5-3_sigma(tm)A(tm)_converges_to_a_infty}
\begin{align}
\label{eq:Kozono_Maeda_Naito_5-3_sigma(tm)A(tm)_converges_to_a_infty_weakly_in_H2_on_good_ball}
\sigma(t_m)^*A(t_m) &\rightharpoonup a_\infty \quad\hbox{weakly in } H_\Gamma^2(B_{r/4}(x_*); \Lambda^1\otimes\ad P),
\\
\label{eq:Kozono_Maeda_Naito_5-3_sigma(tm)A(tm)_converges_to_a_infty_strongly_in_W1p_on_good_ball}
\sigma(t_m)^*A(t_m) &\to a_\infty \quad\hbox{strongly in } W_\Gamma^{1,p}(B_{r/4}(x_*); \Lambda^1\otimes\ad P), \quad\hbox{as } m \to \infty,
\end{align}
\end{subequations}
where $2 \leq p < 4$. A standard diagonal subsequence and patching argument (see Step \ref{step:Global_ideal_limit_sequence_rescaled_connections_over_4space} in the proof of Theorem \ref{thm:Kozono_Maeda_Naito_5-4}) now applies \mutatis to give a (relabeled) subsequence of times $\{t_m\} \subset (0, \infty)$, a sequence of gauge transformations, $\{\Phi_m\}_{m\in\NN} \subset \Aut (P \restriction X\less \Sigma_\infty)$ of class $H^3$, and a finite-energy connection, $A_\infty$, of class $H^2$ on $P \restriction X\less \Sigma_\infty$ such that $t_m \to \infty$ as $m \to \infty$ and
\begin{subequations}
\label{eq:Kozono_Maeda_Naito_5-3_Phi(tm)A(tm)_converges_to_Ainfinity}
\begin{align}
\label{eq:Kozono_Maeda_Naito_5-3_Phi(tm)A(tm)_converges_to_Ainfinity_weakly_in_H2loc_on_X_away_from_bubble_points}
\Phi_m^*A(t_m) &\rightharpoonup A_\infty \quad\hbox{weakly in } H_{A_1,\loc}^2(X\less \Sigma_\infty; \Lambda^1\otimes\ad P),
\\
\label{eq:Kozono_Maeda_Naito_5-3_Phi(tm)A(tm)_converges_to_Ainfinity_strongly_in_W1ploc_on_X_away_from_bubble_points}
\Phi_m^*A(t_m) &\to A_\infty \quad\hbox{strongly in } W_{A_1,\loc}^{1,p}(X\less \Sigma_\infty; \Lambda^1\otimes\ad P), \quad\hbox{as } m \to \infty.
\end{align}
\end{subequations}
Uhlenbeck's Removable Singularity Theorem for finite-energy Sobolev connections, Theorem \ref{thm:Uhlenbeck_removable_singularity_Sobolev}, and Remark \ref{rmk:Uhlenbeck_removable_singularity_Sobolev_Wkp} ensure the existence of a principal $G$-bundle, $P^\flat$, over $X$ and a principal $G$-bundle isomorphism, $\Phi: P^\flat\restriction X\less \Sigma_\infty \cong P\restriction X\less \Sigma_\infty$ of class $H^3$, and a connection, $A_\infty^\flat$, of class $H^2$ on $P^\flat$ such that $\Phi^*A_\infty = A_\infty^\flat$ on $P^\flat\restriction X\less \Sigma_\infty$. The conclusions of Item \eqref{item:Theorem_Kozono_Maeda_Naito_5-3_T_is_infinite_H2loc_weak_and_W1ploc_strong_convergence_Am} are obtained by composing these bundle isomorphisms.

Item \eqref{item:Yang-Mills_gradient_flow_infinite_time_bubbling_obstruction_preserved} is given by \cite[Theorem 5.5]{Sedlacek}. Item \eqref{item:Theorem_Kozono_Maeda_Naito_5-3_T_is_infinite_weak_convergence_FAm_measures} is a simple consequence of Items \eqref{item:Theorem_Kozono_Maeda_Naito_5-3_T_is_infinite_H2loc_weak_and_W1ploc_strong_convergence_Am} and \eqref{item:Limit_r_to_zero_m_to_infinity_energy_over_ball_equals_E_T_is_infinite}. This completes the proof of Theorem \ref{thm:Kozono_Maeda_Naito_5-3_T_is_infinite}.
\end{proof}

The proof of Theorem \ref{thm:Kozono_Maeda_Naito_5-3_T_is_infinite} yields the following useful

\begin{cor}[\Apriori estimates away from bubble points for a solution to the Yang-Mills gradient flow approaching a singularity at time $T = \infty$]
\label{cor:Kozono_Maeda_Naito_5-3_T_is_infinite}
Assume the hypotheses of Theorem \ref{thm:Kozono_Maeda_Naito_5-3_T_is_infinite} for a solution, $A(t)$, to Yang-Mills gradient flow for $t\in [T_0,\infty)$ and some $T_0 \in (0,\infty)$. If $U \Subset X\less\Sigma_\infty$ is a precompact open subset, then there are constants, $t_0 = t_0(g,T_0,U) \in (T_0, \infty)$ and $C = C(A_1,g,t_0,U)$, such that
\begin{equation}
\label{eq:L_infty_time_and_space_FA_bound_precompact_open_subset_X_less_Sigma_infinity}
\|F_A(t)\|_{L^\infty(U)} + \|\nabla_AF_A(t)\|_{L^\infty(U)} \leq C, \quad\forall\, t \geq t_0+1,
\end{equation}
and a family of gauge transformations, $\Phi(t) \in \Aut (P \restriction U)$ for $t\in (t_0,\infty)$, with the same regularity with respect to $t$ as that of $A(t)$ and of class $H^3_{\loc}$ in the spatial variables, such that
\begin{equation}
\label{eq:L_infty_time_H2_space_A_minus_A1_bound_precompact_open_subset_X_less_Sigma_infinity}
\|A(t)-A_1\|_{H_{A_1}^2(U)} \leq C, \quad\forall\, t \geq t_0+1.
\end{equation}
The set, $\Sigma_\infty$, is uniquely determined by the Yang-Mills gradient flow, $A(t)$, for $t\in [T_0,\infty)$.
\end{cor}

Our next result, Theorem \ref{thm:Kozono_Maeda_Naito_5-4}, is an analogue of the \cite[Theorem 5.4]{Kozono_Maeda_Naito_1995}, due to Kozono, Maeda, and Naito, and \cite[Theorem 1.2]{Schlatter_1997}, due to Schlatter. We provide a considerably more detailed statement and proof than that provided in \cite{Kozono_Maeda_Naito_1995} or \cite{Schlatter_1997} for several reasons. For example, in the application of their \cite[Lemma 4.2]{Kozono_Maeda_Naito_1995} in the proof of \cite[Theorem 5.4]{Kozono_Maeda_Naito_1995}, an explanation of how to verify the hypothesis on the smallness of the energy measure, $\eps(x,r;A,T)$, is omitted. Most significantly, however, we provide further details based on Taubes \cite[Section 4]{TauPath}, \cite[Section 5]{TauFrame} and the author \cite[Section 4.2]{FeehanGeometry} regarding the rescaling process: the proofs and original statements of \cite[Theorem 5.4]{Kozono_Maeda_Naito_1995} and \cite[Theorem 1.2]{Schlatter_1997} appear to suggest that, after one level of rescaling, no additional bubbling can occur on $S^4$ --- that is, there are no non-trivial bubble-tree limits. However, based on results for harmonic map gradient-like flow, one expects that non-trivial bubble-tree limits may occur in Yang-Mills gradient flow and the proofs due to Kozono, Maeda, and Naito for \cite[Theorem 5.4]{Kozono_Maeda_Naito_1995} and Schlatter for \cite[Theorem 1.2]{Schlatter_1997} apparently overlook that possibility. We refer the reader to Section \ref{subsec:Bubbletree_limits_Yang-Mills_gradient_flow_four-manifold} for additional discussion of bubble-tree limits in harmonic heat flow and for sequences of connections with bounded energy.

If $\Omega$ is a locally compact Hausdorff space, then $(C_0(\Omega;\RR))'$, the dual space of $C_0(\Omega;\RR)$, may be identified with the space of Radon measures on $(\Omega, \sB(\Omega))$. The behavior of the Yang-Mills gradient flow, $A(t)$, near the flow singularities, $\{(x_l, T): 1 \leq l \leq L\}$, is elucidated in the following (compare \cite[Section 4.2]{FeehanGeometry})

\begin{thm}[Uhlenbeck convergence over $S^4$ for a rescaled solution to the Yang-Mills gradient flow approaching a singularity at time $T \leq \infty$]
\label{thm:Kozono_Maeda_Naito_5-4}
\cite[Theorem 5.4]{Kozono_Maeda_Naito_1995}, \cite[Theorem 1.2]{Schlatter_1997}, and \cite[Theorem 2.4]{Struwe_1994}
Assume the hypotheses of Theorem \ref{thm:Kozono_Maeda_Naito_5-1} and that $L \geq 1$ with $\Sigma_T := \{x_l: 1 \leq l \leq L\}$ denoting a set of singular points for the Yang-Mills gradient flow on the principal $G$-bundle $P$ over $X$, for a time $T \in (0, \infty]$. Let $g$ denote the smooth Riemannian metric on $X$ and $\rho \in (0, 1]$ and $R \geq 1$ be constants obeying the conditions \eqref{eq:Kozono_Maeda_Naito_Theorem_5-4_fixed_small_singular_point_radial_parameter} and \eqref{eq:Kozono_Maeda_Naito_Theorem_5-4_proof_choice_large_Euclidean_4ball_radius_R} in the sequel and, for $1 \leq l \leq L$, choose $g$-orthonormal frames, $f_l$, for the tangent spaces, $(TX)_{x_l}$, and points $p_l \in P_{x_l}$. For $1 \leq l \leq L$, there exist
\begin{itemize}
\item A geodesic normal coordinate chart, $\varphi_l: \RR^4 \supset B_{8R\rho}(0) \cong B_{8R\rho}(x_l) \subset X$ with $\varphi_l(0) = x_l$, defined by the frame, $f_l$;

\item Times, $\{t_m\}_{m \in \NN} \subset (0, T)$, such that $t_m \nearrow T$ as $m \to \infty$;

\item Mass centers, $\{x_{l,m}\}_{m \in \NN} \subset X$, such that $x_{l,m} \to x_l$ as $m \to \infty$;

\item Scales, $\{\lambda_{l,m}\}_{m \in \NN} \subset (0, R_0]$, such that $\lambda_{l,m} \searrow 0$ as $m \to \infty$;

\item Sequences of geodesic normal coordinate charts, $\varphi_{l,m}: \RR^4 \supset B_{2\rho}(0) \cong B_{2\rho}(x_{l,m}) \subset X$ with $\varphi_{l,m}(0) = x_{l,m}$, defined by the $g$-orthonormal frame $f_{l,m}$ for $(TX)_{x_{l,m}}$ obtained by parallel translation of $f_l$ along the geodesic joining $x_l$ to $x_{l,m}$;

\item Sequences of local sections, $\sigma_{l,m}: B_{2\rho}(x_{l,m}) \to P\restriction B_{2\rho}(x_{l,m})$, constructed via parallel translation, with respect to the connections, $A(t_m)$ on $P$, of the point $p_l \in P_{x_l}$ to a point $p_{l,m} \in P_{x_{l,m}}$ along the geodesic joining $x_l$ to $x_{l,m}$, followed by parallel translation along radial geodesics emanating from $x_{l,m} \in X$;

\item Times, $\{s_{l,m}\}_{m\in\NN} \subset (-\eps,0)$, for some constant $\eps \in (0,1]$;

\item Sequences of $H^3$ isomorphisms of principal $G$-bundles, $\Phi_{l,m}: P \restriction S^4 \less \{\infty\} \cong G \times \RR^4$.
\end{itemize}
For $1 \leq l \leq L$, define a sequence of flows, $\{A_{l,m}\}_{m \in \NN}$, on a sequence of product $G$-bundles, $Q_{l,m}\times G \subset \RR^5\times G$, by setting
\begin{equation}
\label{eq:Schlatter_25_theoremstatement}
\begin{aligned}
A_{l,m} &:= \Gamma + a_{l,m} \quad\hbox{and}
\\
a_{l,m}(s,y) &:= \varphi_{l,m}^*\sigma_{l,m}^*A(t_m + \lambda_{l,m}^2s, \lambda_{l,m} y),
\quad\forall\, (s,y) \in U_{l,m}, \quad m \in \NN,
\end{aligned}
\end{equation}
where $\Gamma$ denotes the product connection on $\RR^4\times G$ and the sequence of increasing open subsets, $Q_{l,m} \subset \RR^5$, is defined by
\begin{equation}
\label{eq:Schlatter_25_s-y_domain_theoremstatement}
\{(s,y)\in\RR^5: t := t_m + \lambda_{l,m}^2s \in (0, T) \hbox{ and } x := \lambda_{l,m} y \in B_\rho(0)\}.
\end{equation}
For $1 \leq l \leq L$, define a sequence of smooth Riemannian metrics, $\{g_{l,m}\}_{m\in\NN}$, on a sequence of increasing open balls in $\RR^4$ with center at the origin, by setting
\begin{equation}
\label{eq:Feehan_3-33_Euclidean_space}
g_{l,m}(y) := \varphi_{l,m}^*g(\lambda_{l,m} y), \quad\forall\, y \in B_{\rho/\lambda_{l,m}}(0).
\end{equation}
Then the following hold for $1 \leq l \leq L$:
\begin{enumerate}
\item
\label{item:Theorem_Kozono_Maeda_Naito_5-4_Cinfinity_loc_convergence_gm}
If $\delta$ denotes the standard Euclidean metric on $\RR^4$, then $g_{l,m} \to \delta$ in $C_{\loc}^\infty(\RR^4)$ as $m \to \infty$;

\item
\label{item:Theorem_Kozono_Maeda_Naito_5-4_Am_obeys_Yang-Mills_gradient_flow_Euclidean_space}
The family of connections, $A_{l,m}(s)$, obeys the Yang-Mills gradient flow equation \eqref{eq:Yang-Mills_gradient_flow_equation} over $Q_{l,m} \subset \RR^5$ with respect to the metric $g_{l,m}$;

\item
\label{item:Theorem_Kozono_Maeda_Naito_5-4_H2loc_weak_and_W1ploc_strong_convergence_Am}
There is a smooth Yang-Mills connection, $A_{l,\infty}$, on a principal $G$-bundle, $P_l$, over $S^4$ with its standard round metric of radius one, an integer $J_l \geq 0$ and, if $J_l \geq 1$, a set of points, $\Sigma_{T;l} := \{y_{l,j}\}_{j=1}^{J_l} \subset \RR^4\less\{\infty\}$ such that for any $p\in [2,4)$,
\begin{subequations}
\label{eq:Kozono_Maeda_Naito_5-4_PhimAlm(sm)_converges_to_Al_infinity_flat}
\begin{align}
\label{eq:Kozono_Maeda_Naito_5-4_PhimAlm(sm)_converges_to_Al_infinity_flat_weakly_in_H2loc_on_X_away_from_bubble_points}
\Phi_{l,m}^*A_{l,m}(s_{l,m}) &\rightharpoonup A_{l,\infty} \quad\hbox{weakly in } H_{\Gamma,\loc}^2(\RR^4\less \Sigma_{T;l}; \Lambda^1\otimes\fg),
\\
\label{eq:Kozono_Maeda_Naito_5-4_PhimAlm(sm)_converges_to_Al_infinity_flat_strongly_in_W1ploc_on_X_away_from_bubble_points}
\Phi_{l,m}^*A_{l,m}(s_{l,m}) &\to A_{l,\infty} \quad\hbox{strongly in } W_{\Gamma,\loc}^{1,p}(\RR^4\less \Sigma_{T;l}; \Lambda^1\otimes\fg),
\quad\hbox{as } m \to \infty.
\end{align}
\end{subequations}

\item
\label{item:Theorem_Kozono_Maeda_Naito_5-4_weak_convergence_FAm_measures}
If $J_l \geq 1$, there is a set of constants, $\sE_{l,j}$, given by the limits,
\begin{equation}
\label{eq:Limit_r_to_zero_m_to_infinity_energy_over_ball_equals_Elj_Euclidean_space}
\sE_{l,j} = \lim_{r\to 0} \lim_{m\to \infty} \int_{B_r(y_{l,j})} |F_{A_{l,m}}(s_{l,m})|^2 \, d^4y, \quad 1 \leq j \leq J_l,
\end{equation}
such that
$$
\eps_1 \leq \sE_{l,j} \leq \sE_l, \quad 1 \leq j \leq J_l,
$$
and $|F_{A_{l,m}}(s_{l,m})|^2 \rightharpoonup |F_{A_{l,\infty}}|^2 + \sum_{j=1}^{J_l} \sE_{l,j}\, \delta_{y_{l,j}}$ in $(C_0(\RR^4;\RR))'$ as $m \to \infty$. The ideal limit, $(A_{l,\infty}, \Sigma_{T;l})$ is non-trivial for in the sense that either $J_l \geq 1$ or, if $J_l = 0$,
$$
\int_{\RR^4} |F_{A_{l,\infty}}|^2 \,d^4y \geq \frac{\eps_1}{2}.
$$

\item
\label{item:Kozono_Maeda_Naito_5-4_limit_m_to_infinity_energy_over_increasing_balls_in_Euclidean_space_equals_E}
When $T< \infty$, the positive constants, $\sE_l$, in Theorem \ref{thm:Kozono_Maeda_Naito_5-3} may be computed by the alternative formula,
\begin{equation}
\label{eq:Limit_m_to_infinity_energy_over_increasing_balls_in_Euclidean_space_equals_E}
\sE_l = \lim_{\zeta \to 0} \lim_{m \to \infty} \int_{B_{\zeta/\lambda_{l,m}}(0)} |F_{A_{l,m}}(s_{l,m},y)|^2 \, d^4y.
\end{equation}
\end{enumerate}
\end{thm}

\begin{proof}
Item \eqref{item:Theorem_Kozono_Maeda_Naito_5-4_Cinfinity_loc_convergence_gm} is a consequence (compare \cite[Lemma 3.12]{FeehanGeometry}) of the definition \eqref{eq:Feehan_3-33_Euclidean_space} and more refined $C^0$ and $C^1$ convergence rate estimates follow from standard properties \eqref{eq:Riemannian_metric_components_in_geodesic_normal_coordinates} of geodesic normal coordinates.

Because the radii of all geodesic balls in $X$ encountered in our proof of Theorem \ref{thm:Kozono_Maeda_Naito_5-4} will be small relative to the injectivity radius of the Riemannian metric, $g$, on $X$ (see \eqref{eq:Kozono_Maeda_Naito_Theorem_5-4_fixed_small_singular_point_radial_parameter} in the sequel), to reduce notational clutter in the remainder of the proof we shall suppose without loss of generality that the Riemannian metric, $g$, is \emph{flat} in small balls around the points $x_l \in \Sigma_T$. Given the choice of an orthonormal frame, $f_l$, for $(TX)_{x_l}$ and corresponding geodesic normal coordinate chart, $\varphi_l$, we have an identification $\varphi_l: \RR^4 \supset B_\zeta(0) \cong B_\zeta(x_l) \subset X$ of geodesic balls, for any $\zeta > 0$ less than the injectivity radius of $(X,g)$. In the sequel, we shall therefore identify the balls, $B_\zeta(x_l) = B_\zeta(0)$, and the volume elements, $d\vol_g = d^4x$, where $x^\mu$ (for $\mu=1,\ldots,4$) denote the standard Euclidean coordinates on $\RR^4$. Of course, this simplifying assumption of flatness for the metric $g$ near the points $x_l$ is not essential and indeed in our construction in \cite[Section 4.2]{FeehanGeometry}, we treat the case of a Riemannian metric that is not necessarily flat near bubble points, $x_l \in X$, in the context of bubbling sequences of $g$-anti-self-dual connections.

\setcounter{step}{0}
\begin{step}[Definitions of the sequences of mass centers, scales, and times]
\label{step:Definition_sequences_mass_centers_scales_times}
We recall from Theorem \ref{thm:Kozono_Maeda_Naito_5-1} that a singular point $x_l \in \Sigma_T$ is characterized by
\begin{equation}
\label{eq:Kozono_Maeda_Naito_5-1}
\limsup_{t\nearrow T}\int_{B_r(x_l)} |F_A(t)|^2 \, d\vol_g \geq \eps_1, \quad\forall\, r \in (0, R_0].
\end{equation}
Of course, for each $l \in \{1,\ldots,L\}$, since
$$
\int_{B_r(x_l)} |F_A(t)|^2 \, d\vol_g \leq \int_X |F_A(t)|^2 \, d\vol_g \leq \int_X |F_A(0)|^2 \, d\vol_g, \quad\forall\, t \in [0, T),
$$
the limit in \eqref{eq:Kozono_Maeda_Naito_5-1} is finite. By virtue of \eqref{eq:Kozono_Maeda_Naito_5-1}, there is a monotone increasing sequences of times, $\{t_m\}_{m \in \NN} \subset [0, T)$, such that (compare \cite[Equation (4.3)]{FeehanGeometry})
\begin{equation}
\label{eq:Feehan_4-3}
\sE_{x_l} :=  \lim_{r\to 0} \lim_{m\to\infty} \int_{B_r(x_l)} |F_A(t_m)|^2 \, d\vol_g \geq \eps_1, \quad 1 \leq l \leq L.
\end{equation}
For a sufficiently large positive constant,
\begin{equation}
\label{eq:Kozono_Maeda_Naito_Theorem_5-4_proof_large_Euclidean_4ball_radius_R_dependencies}
R = R(\eps_1, \sE_{x_1}, \ldots, \sE_{x_L}) \geq 1,
\end{equation}
to be determined in the sequel (see \eqref{eq:Kozono_Maeda_Naito_Theorem_5-4_proof_choice_large_Euclidean_4ball_radius_R}), let $\rho \in (0, 1]$ be small enough that
\begin{equation}
\label{eq:Kozono_Maeda_Naito_Theorem_5-4_fixed_small_singular_point_radial_parameter}
\begin{gathered}
\dist_g(x_i, x_j) \geq 8\rho R \quad\forall\, i, j \in \{1, \ldots, L\} \hbox{ such that } i \neq j,
\\
8\rho R \ll \hbox{Injectivity radius of $(X,g)$},
\\
8\rho R \leq R_0,
\end{gathered}
\end{equation}
where $R_0$ is the (small) positive constant appearing in Section \ref{subsec:Kozono_Maeda_Naito_4} and
\begin{equation}
\label{eq:Kozono_Maeda_Naito_Theorem_5-4_proof_small_background_connection}
\int_{B_{8\rho R}(x_l)} |F_{A^\flat}|^2 \, d\vol_g \leq \frac{\eps_1}{16}.
\end{equation}
To analyze the behavior of the flow, $A(t)$, near the set of points $\Sigma_T = \{x_1,\ldots,x_L\}$, it suffices to focus on each point individually. Therefore, to simplify notation during the remainder of the proof, we denote that singular point by $x_0 \in \Sigma_T$ and suppress the explicit dependence on the singular point label $l \in \{1\ldots,L\}$.

We now define sequences of relative mass centers, scales, and times by adapting our construction in \cite[Section 4.2]{FeehanGeometry}. Theorem \ref{thm:Kozono_Maeda_Naito_5-3} implies that
\begin{equation}
\label{eq:Kozono_Maeda_Naito_Theorem_5-3_2_punctured_ball}
F_{A(t)} \to F_{A^\flat} \quad\hbox{strongly in } L^2_{\loc}(B_\rho(x_0)\less\{x_0\}), \quad\hbox{as } t\nearrow T,
\end{equation}
and thus we also have (compare \cite[Equation (4.5)]{FeehanGeometry})
\begin{multline}
\label{eq:Feehan_4-5}
\sE_{x_0} = \lim_{m\to\infty} \int_{B_\rho(x_0)} \left(|F_A(t_m)|^2 - |F_{A^\flat}|^2\right)\, d\vol_g
\\
= \lim_{r\to 0} \lim_{m\to\infty} \int_{B_r(x_0)} \left(|F_A(t_m)|^2 - |F_{A^\flat}|^2\right) \, d\vol_g.
\end{multline}
By virtue of \eqref{eq:Feehan_4-5} we may also assume, by taking large enough $m$ and then relabeling the sequence if needed, that (compare \cite[Equations (4.9) and (4.10)]{FeehanGeometry})
\begin{equation}
\label{eq:Feehan_4-9_and_4-10_difference_measure}
\frac{7}{8}\sE_{x_0} \leq \int_{B_\rho(x_0)} \left||F_A(t_m)|^2 - |F_{A^\flat}|^2\right|\, d\vol_g \leq \frac{9}{8}\sE_{x_0}, \quad\forall\, m \in \NN.
\end{equation}
We define a sequence of relative mass centers, $\{x_m\}_{m\in\NN} \subset B_{9\rho/8}(x_0)$, with coordinates (compare \cite[Equation (4.6)]{FeehanGeometry})
\begin{equation}
\label{eq:Feehan_4-6}
x_m^\mu := \frac{1}{\sE_{x_0}}\int_{B_\rho(x_0)} x^\mu\left(|F_A(t_m)|^2 - |F_{A^\flat}|^2\right)\, d\vol_g, \quad 1\leq \mu \leq 4, \quad\forall\, m \in \NN,
\end{equation}
and scales, $\{\lambda_m\}_{m \in \NN} \subset (0, 3\rho/(2\sqrt{2}))$ (compare \cite[Equation (4.7)]{FeehanGeometry}),
\begin{equation}
\label{eq:Feehan_4-7}
\lambda_m^2 := \frac{1}{\sE_{x_0}}\int_{B_\rho(x_0)} \dist_g(x, x_m)^2 \left||F_A(t_m)|^2 - |F_{A^\flat}|^2\right|\, d\vol_g, \quad\forall\, m \in \NN,
\end{equation}
noting that $\dist_g(x, x_m) = |x - x_m|$ (Euclidean distance function on $\RR^4$) by our simplifying assumption that the Riemannian metric is flat on $B_{8\rho R}(x_0)$. (The upper bound, $3\rho/(2\sqrt{2})$, on the size of the scales, $\lambda_m$, is due to the upper bound $9\sE_{x_0}/8$ in \eqref{eq:Feehan_4-9_and_4-10_difference_measure}.)

Equation \eqref{eq:Feehan_4-7} leads to the following Chebychev-type inequality (compare \cite[Equation (4.8)]{FeehanGeometry}) for all $R \geq 1$,
\begin{equation}
\label{eq:Feehan_4-8}
\int_{B_\rho(x_0) \less B_{\lambda_m R}(x_m)} \left||F_A(t_m)|^2 - |F_{A^\flat}|^2\right|\, d\vol_g \leq R^{-2}\sE_{x_0}, \quad\forall\, m \in \NN.
\end{equation}
To see this, we recall that for any measure space, $(\Omega,\Sigma,\nu)$, and measurable function $f:\Omega \to [-\infty,\infty]$, then the classical Chebychev inequality \cite[Section 5.7]{Feller2} asserts that
\begin{equation}
\label{eq:Chebychev_inequality}
\int_{\{\omega \in \Omega: |f(\omega)| \geq r^2 \}} d\nu \leq \frac{1}{r^2} \int_\Omega f^2 \,d\nu.
\end{equation}
Applying \eqref{eq:Chebychev_inequality} with $\Omega  = B_\rho(x_0)$ and $r = \lambda_m R$ and $f(x) = \dist_g(x, x_m)$ for $x \in B_\rho(x_0)$ and
$$
d\nu = \left||F_A(t_m)|^2 - |F_{A^\flat}|^2\right|\, d\vol_g,
$$
yields, for all $m \in \NN$,
\begin{align*}
{}&\int_{B_\rho(x_0) \less B_{\lambda_m R}(x_m)} \left||F_A(t_m)|^2 - |F_{A^\flat}|^2\right|\, d\vol_g
\\
&\quad \leq \frac{1}{\lambda_m^2 R^2} \int_{B_\rho(x_0)} \dist_g(x, x_m)^2\left||F_A(t_m)|^2 - |F_{A^\flat}|^2\right|\, d\vol_g \quad\hbox{(by \eqref{eq:Chebychev_inequality})}
\\
&\quad = R^{-2}\sE_{x_0} \quad\hbox{(by \eqref{eq:Feehan_4-7})},
\end{align*}
which is \eqref{eq:Feehan_4-8}, as claimed.

We also assert that our definitions \eqref{eq:Feehan_4-6} and \eqref{eq:Feehan_4-7} of the mass centers and scales ensure that they have the properties described in Claim \ref{claim:Convergence_mass_centers_and_scales_and_no_external_bubbles} below.

\begin{claim}[Convergence of mass centers and scales and absence of bubbling over $B_\rho(x_0) \less B_{\lambda_m R}(x_m)$]
\label{claim:Convergence_mass_centers_and_scales_and_no_external_bubbles}
Continue the notation of Step \ref{step:Definition_sequences_mass_centers_scales_times}. Then,
\begin{align}
\label{eq:Mass_centers_converge_to_bubble_point_in_X}
x_m &\to x_0,
\\
\label{eq:Scales_converge_to_zero}
\lambda_m &\to 0, \quad\hbox{as } m \to \infty.
\end{align}
If in addition the constant $R \geq 1$ is large enough that
\begin{equation}
\label{eq:Kozono_Maeda_Naito_Theorem_5-4_proof_choice_large_Euclidean_4ball_radius_R}
R^{-2}\sE_{x_0} \leq \frac{\eps_1}{2}, \quad\forall\, x_0 \in \Sigma_T,
\end{equation}
then no bubbling can occur for the sequence $A(t_m)$ over $B_\rho(x_0) \less B_{\lambda_m R}(x_m)$ as $m \to\infty$.
\end{claim}

\begin{proof}[Proof of Claim \ref{claim:Convergence_mass_centers_and_scales_and_no_external_bubbles}]
The final assertion that the sequence $A(t_m)$ cannot bubble over $B_\rho(x_0) \less B_{\lambda_m R}(x_m)$ as $m \to\infty$ follows from inequalities \eqref{eq:Feehan_4-8} and \eqref{eq:Kozono_Maeda_Naito_Theorem_5-4_proof_choice_large_Euclidean_4ball_radius_R} and Lemma \ref{lem:Kozono_Maeda_Naito_4-3}.

To verify \eqref{eq:Mass_centers_converge_to_bubble_point_in_X}, let $\delta > 0$ and note that
\begin{align*}
|x_m^\mu| &\leq \frac{1}{\sE_{x_0}}\int_{B_\rho(x_0)} |x^\mu|\left||F_A(t_m)|^2 - |F_{A^\flat}|^2\right|\, d\vol_g
\quad\hbox{(by \eqref{eq:Feehan_4-6})}
\\
&=
\frac{1}{\sE_{x_0}}\int_{B_\rho(x_0)\less B_\delta(x_0)} |x^\mu|\left||F_A(t_m)|^2 - |F_{A^\flat}|^2\right|\, d\vol_g
\\
&\quad + \frac{1}{\sE_{x_0}}\int_{B_\delta(x_0)} |x^\mu|\left||F_A(t_m)|^2 - |F_{A^\flat}|^2\right|\, d\vol_g
\\
&\leq \frac{\rho}{\sE_{x_0}}\int_{B_\rho(x_0)\less B_\delta(x_0)} \left||F_A(t_m)|^2 - |F_{A^\flat}|^2\right|\, d\vol_g
+ \frac{9\delta}{8} \quad\hbox{(by \eqref{eq:Feehan_4-9_and_4-10_difference_measure})}.
\end{align*}
The integral term on the right-hand side of the preceding inequality tends to zero as $m \to \infty$ by \eqref{eq:Kozono_Maeda_Naito_Theorem_5-3_2_punctured_ball} and as $\delta >0$ is arbitrary, we find that $x_m \to x_0$ as $m \to \infty$. The verification of \eqref{eq:Scales_converge_to_zero} follows by the same argument, using \eqref{eq:Feehan_4-7} in place of \eqref{eq:Feehan_4-6}. This completes the proof of Claim \ref{claim:Convergence_mass_centers_and_scales_and_no_external_bubbles}.
\end{proof}

(We remark that in the proofs of \cite[Theorem 5.4]{Kozono_Maeda_Naito_1995} or \cite[Theorem 1.2]{Schlatter_1997}, it is simply asserted that there are sequences of points in $X$ converging to $x_0$ and radii converging to zero, as $m \to \infty$, using which one can perform rescaling; however, the existence alone of such sequences is insufficient to guarantee a bubble-tree limit.) This completes Step \ref{step:Definition_sequences_mass_centers_scales_times} and our definition of the sequences of times, $\{t_m\}_{m\in\NN}$, mass centers, $\{x_m\}_{m\in\NN}$, and scales, $\{\lambda_m\}_{m\in\NN}$, along with a development of their basic properties.
\end{step}


\begin{step}[Good small balls in four-dimensional Euclidean space after rescaling]
\label{step:Good_small_balls_four-dimensional_Euclidean_space_after_rescaling}
Because the spatial scale invariance of the $L^2$ norm on covariant two-tensors in dimension four yields\footnote{We abuse notation slightly by not explicitly denoting the composition with spatial translation and rescaling and writing, for example, $F_{\hat A}(t_m,y) = F_A(t_m,x)$, when $x = x_m + \lambda_m y$.}
$$
\int_{B_{\lambda_m R}(x_m)}|F_A(t_m,x)|^2\,d^4x = \int_{B_R(0)}|F_A(t_m,y)|^2\,d^4y, \quad\forall\, m \in \NN,
$$
and the fact that \eqref{eq:Feehan_4-3} at most gives,
\begin{equation}
\label{eq:Upper_bound_energy_Am_ball_lambda_mR}
\limsup_{m \to \infty} \int_{B_{\lambda_m R}(x_m)}|F_A(t_m,x)|^2\,d^4x \leq \sE_{x_0},
\end{equation}
we cannot \apriori eliminate the possibility of additional bubbling for the sequence $A(t_m)$, after rescaling by $\lambda_m$, over the ball $B_R(0) \subset \RR^4$ if $\sE_{x_0}$ is sufficiently large relative to the constant $\eps_1$ appearing in Theorem \ref{thm:Kozono_Maeda_Naito_5-1}. We shall return to the question of possible additional bubbling in Step \ref{step:Bad_small_balls_four-dimensional_Euclidean_space_after_rescaling}. For now, we consider a point $y_* \in \RR^4$ such that there are an integer $m_* = m_*(y_*)$ with the property that
$$
y_* \in B_{\rho/\lambda_m}(0), \quad\forall\, m \geq m_*,
$$
and then a small enough positive constant, $r_* \in (0, R_0]$, whose dependencies include $y_*$ and $\eps_1$, such that
$$
B_{r_*}(y_*) \subset B_{\rho/\lambda_m}(0), \quad\forall\, m \geq m_*,
$$
and the following energy inequality is obeyed,
$$
\int_{B_{r_*}(y_*)}|F_A(t_m,y)|^2\,d^4y \leq \eps_1, \quad\forall\, m \geq m_*.
$$
By decreasing $r_*$ further, if necessary, and relabeling the sequence, if necessary, we shall require in addition that a \emph{good small ball} in $\RR^4$ obeys the stronger energy inequality\footnote{We include the factor of $7/8$ in the inequality for convenience later in the proof.},
\begin{equation}
\label{eq:Kozono_Maeda_Naito_5-1_nested_small_good_ball_in_Euclidean_4pace}
B_{r_*}(y_*) \subset B_{\rho/\lambda_m}(0) \less \Xi_T
\quad\hbox{and}\quad
\sup_{t \in [0, t_m]} \int_{B_{r_*}(y_*)}|F_A(t,y)|^2\,d^4y \leq \frac{7\eps_1}{8},
\quad\forall\, m \in \NN.
\end{equation}
We recall from Step \ref{step:Definition_sequences_mass_centers_scales_times} that $x_m \in B_{9\rho/8}(x_0)$, for all $m \in \NN$. Setting
\begin{equation}
\label{eq:Relationship_between_small_ball_and_Euclidean_space_good_ball_centers}
x_m^* := x_m + \lambda_m y_*, \quad\forall\, m \in \NN,
\end{equation}
and applying the spatial scale invariance of the $L^2$ norm on covariant two-tensors in dimension four, we see that the characterization \eqref{eq:Kozono_Maeda_Naito_5-1_nested_small_good_ball_in_Euclidean_4pace} of a good small ball in $\RR^4$ is equivalent to the following characterization of a sequence of \emph{good very small balls} in $B_{\rho}(x_m)$,
\begin{equation}
\label{eq:Kozono_Maeda_Naito_5-1_nested_very_small_good_ball_in_X}
B_{\lambda_m r_*}(x_m^*) \subset B_{\rho}(x_m)
\quad\hbox{and}\quad
\sup_{t\in[0,t_m]}\int_{B_{\lambda_m r_*}(x_m^*)}|F_A(t,x)|^2\,d^4x \leq \frac{7\eps_1}{8},
\quad\forall\, m \in \NN.
\end{equation}
This concludes Step \ref{step:Good_small_balls_four-dimensional_Euclidean_space_after_rescaling}.
\end{step}

\begin{step}[Lower bounds for the energy of the Yang-Mills gradient flow over sequences of balls enclosing mass centers]
\label{step:Lower_bounds_energy_Yang-Mills_gradient_flow_balls_enclosing_mass_centers}
We first observe that
\begin{align*}
{}&\int_{B_{\lambda_m R}(x_m)}|F_A(t_m)|^2\,d^4x
\\
&\quad = \int_{B_\rho(x_0)}|F_A(t_m)|^2\,d^4x - \int_{B_\rho(x_0)\less B_{\lambda_m R}(x_m)}|F_A(t_m)|^2 \,d^4x
\\
&\quad = \int_{B_\rho(x_0)} \left(|F_A(t_m)|^2 - |F_{A^\flat}|^2\right) \,d^4x
- \int_{B_\rho(x_0)\less B_{\lambda_m R}(x_m)} \left(|F_A(t_m)|^2 - |F_{A^\flat}|^2\right) \,d^4x
\\
&\qquad - \int_{B_{\lambda_m R}(x_m)}|F_{A^\flat}|^2\,d^4x
\\
&\geq \frac{7\sE_{x_0}}{8} - R^{-2}\sE_{x_0}  - \int_{B_{8\rho R}(x_0)}|F_{A^\flat}|^2\,d^4x
\quad\hbox{(by \eqref{eq:Feehan_4-9_and_4-10_difference_measure} and \eqref{eq:Feehan_4-8})}
\\
&\geq \frac{7\sE_{x_0}}{8} - R^{-2}\sE_{x_0} - \frac{\eps_1}{16}
\quad\hbox{(by \eqref{eq:Kozono_Maeda_Naito_Theorem_5-4_proof_small_background_connection})},
\end{align*}
where to obtain the second last inequality we also used the fact that $B_{\lambda_m R}(x_m) \subset B_{8\rho R}(x_0)$ since $\lambda_m \in (0, 3\rho/2)$ and $x_m \in B_{9\rho/8}(x_0)$ from Step \ref{step:Definition_sequences_mass_centers_scales_times}, for all $m \in \NN$. Choosing $R \geq 4$ and noting that $\sE_{x_0} \geq \eps_1$ by \eqref{eq:Feehan_4-3}, we obtain (compare \cite[Equation (22)]{Schlatter_1997})
\begin{equation}
\label{eq:Schlatter_22_lower_bound}
\int_{B_{\lambda_m R}(x_m)}|F_A(t_m)|^2\,d^4x \geq \frac{3\eps_1}{4}, \quad\forall\, m \in \NN.
\end{equation}
By Lemma \ref{lem:Kozono_Maeda_Naito_4-10}, denoting the constant on the right-hand side by $C_1$ and recalling that $8\rho R \leq R_0$ by \eqref{eq:Kozono_Maeda_Naito_Theorem_5-4_fixed_small_singular_point_radial_parameter} (and thus also $2\lambda_m R \leq R_0$ for all $m \in \NN$), for $\eps \in (0, 1]$ to be determined and any $t \in [t_m - \lambda_m^2\eps, t_m]$ we have
\begin{align*}
\int_{B_{\lambda_m R}(x_m)}|F_A(t_m)|^2\,d^4x &\leq \int_{B_{2\lambda_m R}(x_m)}|F_A(t)|^2\,d^4x + C_1(2\lambda_m R)^{-2}(t_m - t) \sE(A(t))
\\
&\leq \int_{B_{2\lambda_m R}(x_m)}|F_A(t)|^2\,d^4x + C_1(2\lambda_m R)^{-2}(\lambda_m^2\eps) \sE(A(0))
\\
&= \int_{B_{2\lambda_m R}(x_m)}|F_A(t)|^2\,d^4x + 4C_1R^{-2}\eps\sE(A(0)),
\end{align*}
where we use the fact that $\sE(A(t)) \leq \sE(A(0))$ by Lemma \ref{lem:Kozono_Maeda_Naito_4-1} when $t \geq 0$. Thus, recalling that $R \geq 4$ from the choice leading to \eqref{eq:Schlatter_22_lower_bound} and choosing $\eps \in (0, 1]$ such that
$$
4C_1R^{-2}\eps\sE(A(0)) \leq 4C_1\eps\sE(A(0))/16 = C_1\eps\sE(A(0))/4 \leq \eps_1/8,
$$
that is,
\begin{equation}
\label{eq:Kozono_Maeda_Naito_page_120_definition_of_eps_for_lower_bound_energy}
0 < \eps \leq \frac{\eps_1}{2C_1\sE(A(0))},
\end{equation}
and combining the preceding energy inequality with the energy inequality \eqref{eq:Schlatter_22_lower_bound}, we see that
\begin{equation}
\label{eq:Schlatter_lemma_3-1_inequality}
\int_{B_{2\lambda_m R}(x_m)}|F_A(t)|^2\,d^4x \geq \frac{5\eps_1}{8}, \quad\forall\, t \in [t_m - \lambda_m^2\eps, t_m] \hbox{ and } m \in \NN.
\end{equation}
This concludes Step \ref{step:Lower_bounds_energy_Yang-Mills_gradient_flow_balls_enclosing_mass_centers} and our derivation of a lower bound for the energy of the flow, $A(t)$, in a sequence of balls enclosing mass center points.
\end{step}

\begin{step}[Upper bounds for the energy of the Yang-Mills gradient flow over a sequence of good very small balls]
\label{step:Upper_bounds_energy_Yang-Mills_gradient_flow_over_very_small_good_ball}
Suppose now, as in Step \ref{step:Good_small_balls_four-dimensional_Euclidean_space_after_rescaling}, that
$$
B_{\lambda_m r_*}(x_m^*) \subset B_\rho(x_m), \quad \forall\, m\in\NN,
$$
is a sequence of good very small balls in $X$ corresponding to a good small ball in $\RR^4$, namely
$$
B_{r_*}(y_*) \subset B_{\rho/\lambda_m}(0), \quad\forall\, m\in \NN,
$$
where $x_m^* = x_m + \lambda_m y_*$. Again as in Step \ref{step:Good_small_balls_four-dimensional_Euclidean_space_after_rescaling}, we assume without loss of generality that the above inclusions hold for all $m \in \NN$ by considering $m \geq m_* = m_*(y_*)$ for some large enough $m_* \geq 1$ and then relabeling the sequence. We recall from \eqref{eq:Kozono_Maeda_Naito_5-1_nested_small_good_ball_in_Euclidean_4pace} and \eqref{eq:Kozono_Maeda_Naito_5-1_nested_very_small_good_ball_in_X} that
\begin{equation}
\label{eq:Kozono_Maeda_Naito_page_120_eps1_energy_inequality_t_lessthan_tm}
\int_{B_{\lambda_m r_*}(x_m^*)}|F_A(t,x)|^2\,d^4x = \int_{B_{r_*}(y^*)}|F_A(t,y)|^2\,d^4y \leq \frac{7\eps_1}{8},
\quad\forall\, m \in \NN \hbox{ and } t \in [0, t_m].
\end{equation}
By Lemma \ref{lem:Kozono_Maeda_Naito_4-10} and again denoting the constant on the right-hand side by $C_1$, for any $t \in [t_m - \lambda_m^2\eps, t_m]$ and noting that $t_m - \lambda_m^2\eps \geq 0$ (for large enough $m \in \NN$) we have
\begin{align*}
\int_{B_{\lambda_m r_*/2}(x_m^*)}|F_A(t,x)|^2\,d^4x &\leq \int_{B_{\lambda_m r_*}(x_m^*)}|F_A(t_m - \lambda_m^2\eps, x)|^2\,d^4x
\\
&\quad + C_1 (\lambda_m r_*)^{-2}(t - (t_m - \lambda_m^2\eps)) \sE(A(t_m - \lambda_m^2\eps))
\\
&\leq \int_{B_{\lambda_m r_*}(x_m^*)}|F_A(t_m - \lambda_m^2\eps,x)|^2\,d^4x + C_1 r_*^{-2}\eps \sE(A(0))
\quad \hbox{(by Lemma \ref{lem:Kozono_Maeda_Naito_4-1})}
\\
&\leq \frac{7\eps_1}{8} + C_1 r_*^{-2}\eps \sE(A(0)) \quad \hbox{by \eqref{eq:Kozono_Maeda_Naito_page_120_eps1_energy_inequality_t_lessthan_tm}.}
\end{align*}
By further decreasing the size of $\eps \in (0, 1]$ chosen to obey \eqref{eq:Kozono_Maeda_Naito_page_120_definition_of_eps_for_lower_bound_energy}, if necessary, we can assume that
\begin{equation}
\label{eq:Kozono_Maeda_Naito_page_120_dependencies_of_eps_for_upper_bound_energy}
\eps = \eps(C_1,\sE(A(0)),r_*,\eps_1)
\end{equation}
obeys $C_1r_*^{-2}\eps\sE(A(0)) \leq \eps_1/8$, namely
\begin{equation}
\label{eq:Kozono_Maeda_Naito_page_120_definition_of_eps_for_upper_bound_energy}
0 < \eps \leq \frac{\eps_1  r_*^2}{8C_1\sE(A(0))}.
\end{equation}
Hence, the preceding energy inequality gives
\begin{equation}
\label{eq:Energy_inequality_condition_needed_to_apply_lemma_Kozono_Maeda_Naito_4-2}
\int_{B_{\lambda_m r_*/2}(x_m^*)}|F_A(t,x)|^2\,d^4x \leq \eps_1, \quad\forall\, t \in [t_m - \lambda_m^2\eps, t_m].
\end{equation}
The small-energy bound \eqref{eq:Energy_inequality_condition_needed_to_apply_lemma_Kozono_Maeda_Naito_4-2} ensures that the \apriori estimates in Sections \ref{subsec:Kozono_Maeda_Naito_4} and \ref{subsec:Schlatter_Struwe_and_Uhlenbeck_apriori_estimates} are applicable on $B_{\lambda_m r_*/2}(x_m^*)$. This concludes Step
\ref{step:Upper_bounds_energy_Yang-Mills_gradient_flow_over_very_small_good_ball}.
\end{step}

\begin{step}[Definition of the sequence of rescaled Yang-Mills gradient flows over subsets of $\RR^4$ and convergence in $L^2$ of the limit of their time derivatives]
\label{step:Definition_sequence_rescaled_Yang-Mills_gradient_flows_Euclidean_space_and_vanishing_L2_limit_time_derivatives}
Keeping in mind the hypotheses of Theorem \ref{thm:Kozono_Maeda_Naito_5-4}, we now define a sequence of rescaled flows of connections, $\{A_m(s)\}_{m\in\NN}$, over a sequence of increasing nested subsets of $\RR\times\RR^4$, by setting
\begin{equation}
\label{eq:Schlatter_25}
A_m(s,y) := \Gamma + a_m(s,y), \quad\hbox{with } a_m(s,y) := a(t_m+\lambda_m^2s, x_m + \lambda_m y) \quad\forall\, m \in \NN,
\end{equation}
where $\Gamma$ is the product connection on $\RR^4\times G$, and $(s,y) \in \RR \times \RR^4$ satisfies
\begin{equation}
\label{eq:Schlatter_25_s-y_domain}
t := t_m + \lambda_m^2s \in [0, T) \quad\hbox{and}\quad x := x_m + \lambda_m y \in B_\rho(0),
\end{equation}
where we identify $X \supset B_\rho(x_0) \cong B_\rho(0) \subset \RR^4$. The rescaled flows of connections, $A_m(s)$, are strong solutions to the Yang-Mills gradient flow equation \eqref{eq:Yang-Mills_gradient_flow_equation} with respect to the standard Euclidean metric on $\RR^4$ (because of our simplifying assumption that $g$ is flat near $x_0 \in X$), for all $m \in \NN$,
\begin{multline}
\label{eq:Schlatter_29}
\frac{\partial A_m}{\partial s}(s,y) = -d_{A_m}^*F_{A_m}(s,y),
\\
\quad\hbox{for a.e. } (s,y) \in [-t_m/\lambda_m^2, (T-t_m)/\lambda_m^2)\times B_{\rho/\lambda_m}(-\lambda_m^{-1}x_m) \subset \RR\times\RR^4.
\end{multline}
This proves Item \eqref{item:Theorem_Kozono_Maeda_Naito_5-4_Am_obeys_Yang-Mills_gradient_flow_Euclidean_space}. We observe that, for all $m \in \NN$,
\begin{align*}
\int_{-\eps}^0 \int_{|y+\lambda_m^{-1}x_m| < \rho/\lambda_m} \left|\frac{\partial A_m}{\partial s}(s,y)\right|^2\,d^4y\,ds
&=
\int_{t_m-\lambda_m^2\eps}^{t_m} \int_{B_\rho(x_0)} \left|\frac{\partial A}{\partial t}(t,x)\right|^2\,d^4x\,dt
\\
&\leq \int_0^T \int_X \left|\frac{\partial A}{\partial t}\right|^2\,d\vol_g\,dt
\\
&\leq \int_X |F_A(0)|^2\,d\vol_g \quad\hbox{(by Corollary \ref{cor:Kozono_Maeda_Naito_4-10_global})}
\\
&= 2\sE(A(0)) < \infty,
\end{align*}
and therefore,
$$
\lim_{m\to\infty}\int_{-\eps}^0 \int_{|y+\lambda_m^{-1}x_m| < \rho/\lambda_m} |\partial_s A_m(s,y)|^2\,d^4y\,ds = 0.
$$
By the Mean Value Theorem for integrals, there is a sequence $\{s_m\}_{m\in\NN} \subset (-\eps, 0)$ such that
$$
\int_{|y+\lambda_m^{-1}x_m| < \rho/\lambda_m} |\partial_s A_m(s_m,y)|^2\,d^4y
=
\frac{1}{\eps}\int_{-\eps}^0 \int_{|y+\lambda_m^{-1}x_m| < \rho/\lambda_m} |\partial_s A_m(s,y)|^2\,d^4y\,ds, \quad\forall\, m \in \NN,
$$
and hence,
\begin{equation}
\label{eq:Schlatter_28}
\lim_{m\to\infty} \int_{|y+\lambda_m^{-1}x_m| < \rho/\lambda_m} |\partial_s A_m(s_m,y)|^2\,d^4y = 0.
\end{equation}
Therefore,
\begin{equation}
\label{eq:Schlatter_28_L2loc_4space}
\frac{\partial A_m}{\partial s}(s_m) \to 0 \quad\hbox{strongly in } L^2_{\loc}(\RR^4; \Lambda^1\otimes\fg), \quad\hbox{as } m \to \infty.
\end{equation}
It is important to note that each $s_m$ inherits the dependencies of the positive constant, $\eps$, from \eqref{eq:Kozono_Maeda_Naito_page_120_dependencies_of_eps_for_upper_bound_energy} and thus
\begin{equation}
\label{eq:Kozono_Maeda_Naito_page_120_dependencies_of_sm_for_upper_bound_energy}
s_m = s_m(C_1,\sE(A(0)),r_*,\eps_1), \quad\forall\, m \in \NN.
\end{equation}
This concludes Step
\ref{step:Definition_sequence_rescaled_Yang-Mills_gradient_flows_Euclidean_space_and_vanishing_L2_limit_time_derivatives}.
\end{step}

\begin{step}[Local integral-norm bounds on the curvatures of the sequence of rescaled connections and local convergence of a sequence of connections in Coulomb gauge over a good small ball]
\label{step:Local_bounds_curvatures_rescaled_connections_and_convergence_sequence_connections_Coulomb_gauge}
Because the rescaled flows of connections, $A_m(s)$, obey the Yang-Mills gradient flow equation \eqref{eq:Schlatter_29} for all $m \in \NN$, we also obtain from \eqref{eq:Schlatter_28} that
\begin{equation}
\label{eq:Schlatter_31}
\lim_{m\to\infty} \int_{|y+\lambda_m^{-1}x_m| < \rho/\lambda_m} |d_{A_m}^*F_{A_m}(s_m,y)|^2\,d^4y = 0.
\end{equation}
Setting $\tau_m := t_m - \lambda_m^2s_m \in (t_m - \lambda_m^2\eps, t_m)$ for $m \in \NN$, the energy inequality \eqref{eq:Energy_inequality_condition_needed_to_apply_lemma_Kozono_Maeda_Naito_4-2} obeyed over a good very small ball, $B_{\lambda_mr_*}(x_m^*) \subset B_\rho(x_0)$, and scale invariance of the $L^2$-norms of covariant two-tensors on four-dimensional Riemannian manifolds implies that
\begin{equation}
\label{eq:Schlatter_26_time_sm_half_ball}
\int_{B_{r_*/2}(y_*)}|F_{A_m}(s_m,y)|^2\,d^4y
= \int_{B_{\lambda_mr_*/2}(x_m^*)}|F_A(\tau_m,x)|^2\,d^4x \leq \eps_1, \quad\forall\, m \in \NN.
\end{equation}
Therefore, thanks to \eqref{eq:Schlatter_31} and \eqref{eq:Schlatter_26_time_sm_half_ball}, Lemma \ref{lem:Schlatter_2-6} yields the \apriori estimate,
\begin{equation}
\label{eq:Schlatter_32}
\|F_{A_m}(s_m)\|_{L^2(B_{3r_*/8}(y_*))} + \|\nabla_{A_m}F_{A_m}(s_m)\|_{L^2(B_{3r_*/8}(y_*))} \leq C,
\quad\forall\, m \in \NN,
\end{equation}
for a positive constant, $C = C(r_*)$. By the Sobolev embedding $H^1(B_{3r_*/8}(0);\RR) \hookrightarrow L^4(B_{3r_*/8}(0);\RR)$ \cite[Theorem 4.12]{AdamsFournier} with constant $c_2 = c_2(r_*)$ and the Kato Inequality \eqref{eq:FU_6-20_first-order_Kato_inequality}, we obtain the bound
\begin{equation}
\label{eq:Schlatter_32_L4}
\|F_{A_m}(s_m)\|_{L^4(B_{3r_*/8}(y_*))} \leq c_2C, \quad\forall\, m \in \NN.
\end{equation}
To proceed further, however, we need more control over the sequence of connections on the product bundle, $B_{3r_*/8}(y_*)\times G$.

Because $A_m(s_m)$ is at least of class $W^{1,p}$ for $2 \leq p \leq 4$ and as the positive constant $\eps_1$ in \eqref{eq:Schlatter_26_time_sm_half_ball} may be assumed without loss of generality to be sufficiently small, Theorem \ref{thm:Uhlenbeck_Lp_1-3} provides a universal positive constant, $c_1$, and a sequence of $W^{2,4}$ gauge transformations, $\varphi_m: B_{3r_*/8}(y_*) \to G$ (depending of course on the point $y_* \in \RR^4$), such that the gauge-transformed sequence of connections, $\widetilde A_m(s_m) := \varphi_m^*A_m(s_m)$, obeys for all $m \in \NN$,
\begin{equation}
\label{eq:Theorem_Uhlenbeck_Lp_1-3_Coulomb_gauge}
d_\Gamma^* (\widetilde A_m(s_m) - \Gamma) = 0 \quad\hbox{a.e. on } B_{3r_*/8}(y_*),
\end{equation}
and
\begin{multline}
\label{eq:Schlatter_27}
\|\widetilde A_m(s_m) - \Gamma\|_{L^4(B_{3r_*/8}(y_*))} + \|\nabla_\Gamma(\widetilde A_m(s_m) - \Gamma)\|_{L^2(B_{3r_*/8}(y_*))}
\\
\leq c_1\|F_{A_m}(s_m)\|_{L^2(B_{3r_*/8}(y_*))}, \quad\forall\, m \in \NN.
\end{multline}
Since $A_m(s_m)$ is actually at least of class $H^2$, the sequence of gauge transformations, $\varphi_m: B_{3r_*/8}(y_*) \to G$, can be assumed by Remark \ref{rmk:Uhlenbeck_theorem_1-3_Wkp} to be at least of class $H^3$, and thus continuous by the Sobolev embedding, $H^3(B_{3r_*/8}(0);\RR) \hookrightarrow C(\bar B_{3r_*/8}(0);\RR)$ \cite[Theorem 4.12]{AdamsFournier}.

By writing $\tilde a_m(s_m) = \widetilde A_m(s_m) - \Gamma$ and $\nabla_\Gamma$ for the covariant derivative defined by the product connection on $\RR^4\times G$ and standard Euclidean metric on $\RR^4$, we discover that, for a universal positive constant, $c_0$, and remaining constants as above,
\begin{align*}
{}& \|\nabla_\Gamma F_{\widetilde A_m}(s_m)\|_{L^2(B_{3r_*/8}(y_*))}
\\
&\quad =  \|\nabla_\Gamma F_{\widetilde A_m}(s_m) + [\tilde a_m, F_{\widetilde A_m}(s_m)]\|_{L^2(B_{3r_*/8}(y_*))}
\\
&\quad \leq \|\nabla_{\widetilde A_m}F_{\widetilde A_m}(s_m)\|_{L^2(B_{3r_*/8}(y_*))}
+ c_0\|F_{\widetilde A_m}(s_m)\|_{L^4(B_{3r_*/8}(y_*))} \|\tilde a_m(s_m)\|_{L^4(B_{3r_*/8}(y_*))}
\\
&\quad \leq \|\nabla_{\widetilde A_m} F_{\widetilde A_m}(s_m)\|_{L^2(B_{3r_*/8}(y_*))}
+ c_0c_1c_2C \|F_{\widetilde A_m}(s_m)\|_{L^2(B_{r_*/2}(y_*))} \quad\hbox{(by \eqref{eq:Schlatter_32_L4} and \eqref{eq:Schlatter_27})}
\\
&\quad \leq \|\nabla_{\widetilde A_m} F_{\widetilde A_m}(s_m)\|_{L^2(B_{3r_*/8}(y_*))}
+ c_0c_1c_2C\sqrt{\eps_1} \quad\hbox{(by \eqref{eq:Schlatter_26_time_sm_half_ball})}
\\
&\quad \leq C\left(1 + c_0c_1c_2\sqrt{\eps_1}\right) \quad\hbox{(by \eqref{eq:Schlatter_32})}.
\end{align*}
Therefore, setting $K :=  C\left(1 + c_0c_1c_2\sqrt{\eps_1}\right) + \sqrt{\eps_1}$, the preceding inequality and \eqref{eq:Schlatter_26_time_sm_half_ball} yields
\begin{equation}
\label{eq:Schlatter_34}
\|F_{\widetilde A_m}(s_m)\|_{H_\Gamma^1(B_{3r_*/8}(y_*))} \leq K, \quad\forall\, m \in \NN.
\end{equation}
Hence, by virtue of the combination of the estimates \eqref{eq:Schlatter_26_time_sm_half_ball} and \eqref{eq:Schlatter_27} for small enough $\eps_1$, the Coulomb gauge condition \eqref{eq:Theorem_Uhlenbeck_Lp_1-3_Coulomb_gauge}, and the $H_\Gamma^1$ curvature bound \eqref{eq:Schlatter_34}, the \apriori estimate provided by Lemma \ref{lem:Schlatter_2-7} supplies a positive constant, $C_2 = C_2(K,r_*)$, such that
\begin{equation}
\label{eq:Schlatter_35}
\|\widetilde A_m(s_m) - \Gamma\|_{H^2_\Gamma(B_{r_*/4}(y_*))} \leq C_2, \quad\forall\, m \in \NN.
\end{equation}
By the Banach-Alaoglu Theorem, after passing to a subsequence, we obtain a weak limit, $a_\infty \in H^2_\Gamma(B_{r_*/4}(y_*); \Lambda^1\otimes\fg)$ and connection $A_\infty := \Gamma + a_\infty$ on $B_{r_*/4}(y_*) \times G$ of class $H^2$, such that
\begin{equation}
\label{eq:Schlatter_36}
\widetilde A_m(s_m) \rightharpoonup A_\infty \quad\hbox{weakly in } H_\Gamma^2(B_{r_*/4}(y_*); \Lambda^1\otimes\fg), \quad\hbox{as } m \to \infty.
\end{equation}
The Rellich-Kondrachov Embedding Theorem \cite[Theorem 6.3]{AdamsFournier} now ensures that, for any $1 \leq q < 4$,
\begin{equation}
\label{eq:Schlatter_36_strong_Wq1_convergence}
\widetilde A_m(s_m) \to A_\infty \quad\hbox{strongly in } W_\Gamma^{q,1}(B_{r_*/4}(y_*); \Lambda^1\otimes\fg), \quad\hbox{as } m \to \infty.
\end{equation}
This concludes Step \ref{step:Local_bounds_curvatures_rescaled_connections_and_convergence_sequence_connections_Coulomb_gauge}.
\end{step}

\begin{step}[Limit of the sequence of rescaled connections over a good small ball is Yang-Mills with respect to the standard Euclidean metric]
\label{step:Limit_sequence_rescaled_connections_is_Yang-Mills}
For any $b \in L^2_{\loc}(\RR^4; \Lambda^1\otimes\fg)$ and good small ball $B_{r_*/4}(y_*) \subset \RR^4$, we discover that
\begin{align*}
(d_{A_\infty}^*F_{A_\infty}, b)_{L^2(B_{r_*/4}(y_*))}
&=
\lim_{m\to\infty}(d_{A_m}^*F_{A_m}(s_m), b)_{L^2(B_{r_*/4}(y_*))}  \quad\hbox{(by \eqref{eq:Schlatter_36})}
\\
&= \lim_{m\to\infty}(\partial_sA_m(s_m), b)_{L^2(B_{r_*/4}(y_*))} \quad\hbox{(by \eqref{eq:Schlatter_29})}
\\
&= 0 \quad\hbox{(by \eqref{eq:Schlatter_28_L2loc_4space}),}
\end{align*}
and thus, as $b$ was arbitrary,
\begin{equation}
\label{eq:Schlatter_43}
d_{A_\infty}^*F_{A_\infty} = 0 \quad\hbox{a.e. on } B_{r_*/4}(y_*),
\end{equation}
and so $A_\infty$ is a Yang-Mills connection on $B_{r_*/4}(y_*)$ with respect to the standard Euclidean metric on $\RR^4$. This concludes Step \ref{step:Limit_sequence_rescaled_connections_is_Yang-Mills}.
\end{step}

\begin{step}[Bad small balls in four-dimensional Euclidean space after rescaling]
\label{step:Bad_small_balls_four-dimensional_Euclidean_space_after_rescaling}
Theorem \ref{thm:Kozono_Maeda_Naito_5-1}, or more precisely its proof applied to the sequence of connections $A_m(s_m)$ on $\RR^4\times G$, implies that there is an at most finite (and possibly empty) subset $\Xi_T \equiv \{y_1, \ldots, y_J\}\subset B_R(0) \subset \RR^4$ of \emph{descendant bubble points} characterized by the following analogue of \eqref{eq:Kozono_Maeda_Naito_5-1},
\begin{equation}
\label{eq:Kozono_Maeda_Naito_5-1_nested_small_bad_ball_in_Euclidean_4pace}
\limsup_{m\to\infty}\int_{B_r(y_j)} |F_{A_m}(s_m,y)|^2 \, d^4y \geq \eps_1, \quad\forall\, r \in (0, R_0].
\end{equation}
The set of descendant bubble points, $\Xi_T \equiv \{y_1, \ldots, y_J\}\subset B_R(0) \subset \RR^4$, is associated to the \emph{parent bubble point} $x_0 \in \Sigma_T \subset X$. (If we temporarily restored the label $l$ for the points $x_l \in \Sigma_T \equiv \{x_1,\ldots, x_L\} \subset X$, we would write $\Sigma_{T;l} \equiv \{y_{l,1}, \ldots, y_{l,L_l}\}\subset B_R(0) \subset \RR^4$ for the set of descendant bubble points associated with the parent bubble point $x_l \in \Sigma_T \subset X$.)

We again recall from Step \ref{step:Definition_sequences_mass_centers_scales_times} that $x_m \in B_{9\rho/8}(x_0)$, for all $m \in \NN$. It is convenient to define, by analogy  with the scaling and translation rule \eqref{eq:Relationship_between_small_ball_and_Euclidean_space_good_ball_centers},
\begin{equation}
\label{eq:Relationship_between_small_ball_and_Euclidean_space_bubble_points}
x_{j,m} := x_m + \lambda_m y_j, \quad\forall\, j \in \{1,\ldots,J\} \hbox{ and } m \in \NN.
\end{equation}
Consequently, by the spatial scale invariance of the $L^2$ norm on covariant two-tensors in dimension four, the characterization \eqref{eq:Kozono_Maeda_Naito_5-1_nested_small_bad_ball_in_Euclidean_4pace} of the \emph{bad small balls}, $B_r(y_j) \subset B_R(0)$, is equivalent to the following characterization of a sequence of \emph{bad very small balls}, $B_{\lambda_m r}(x_{j,m}) \subset B_{\lambda_mR}(x_m)$,
\begin{equation}
\label{eq:Kozono_Maeda_Naito_5-1_nested_very_small_bad_ball_in_X}
\limsup_{m\to\infty}\int_{B_{\lambda_m r}(x_{j,m})} |F_A(\tau_m,x)|^2 \, d^4x \geq \eps_1, \quad\forall\, r \in (0, R_0] \hbox{ and } j \in \{1,\ldots,J\},
\end{equation}
where we recall that $\tau_m = t_m + \lambda_m^2s_m$.

Returning to the primary characterization \eqref{eq:Kozono_Maeda_Naito_5-1_nested_small_bad_ball_in_Euclidean_4pace}, after passing to a subsequence, we may assume that the following limits exist and define the positive constants (relabelled $\sE_{l,j}$) appearing in \eqref{eq:Limit_r_to_zero_m_to_infinity_energy_over_ball_equals_Elj_Euclidean_space},
\begin{equation}
\label{eq:Limit_r_to_zero_m_to_infinity_energy_over_ball_equals_Eyj_Euclidean_space}
\sE_{y_j} = \lim_{r\to 0} \lim_{m\to \infty} \int_{B_r(y_j)} |F_{A_m}(s_m)|^2 \, d^4y, \quad 1 \leq j \leq J.
\end{equation}
This concludes Step \ref{step:Bad_small_balls_four-dimensional_Euclidean_space_after_rescaling}.
\end{step}

\begin{step}[Global ideal limit of the sequence of rescaled connections over $\RR^4$]
\label{step:Global_ideal_limit_sequence_rescaled_connections_over_4space}
We now seek to globalize the conclusions of Step  \ref{step:Local_bounds_curvatures_rescaled_connections_and_convergence_sequence_connections_Coulomb_gauge} with the aid of Uhlenbeck patching arguments. The patching arguments comprising this step are standard and described in detail, in the case of sequences of anti-self-dual connections for example, by Donaldson and Kronheimer in \cite[Section 4.4.2]{DK}, by Schlatter in the present context in the proof of his \cite[Proposition 3.2]{Schlatter_1997}, and of course by Uhlenbeck in \cite[Section 3]{UhlLp} for sequences of connections with an $L^p$ bound on their curvatures. We shall just indicate the necessary changes required to adapt the version of the patching argument described in \cite[Section 4.4.2]{DK} (for sequences of anti-self-dual connections) to the present situation and supply a few details to augment the proof of \cite[Proposition 3.2]{Schlatter_1997}.

We choose countable dense set of points $\{y_n^* \}_{n\in\NN} \subset \RR^4 \less \Xi_T$, where $\Xi_T$ is the finite set of bubble points identified in Step \ref{step:Bad_small_balls_four-dimensional_Euclidean_space_after_rescaling}, and a sequence of radii, $\{r_n^*\}_{n\in\NN} \subset (0,R_0]$, such that
$$
\RR^4 \less \Xi_T = \bigcup_{m\in\NN} B_{r_n^*/4}(y_n^*) = \bigcup_{m\in\NN} B_{r_n^*}(y_n^*),
$$
\emph{and} each small ball $B_{r_n^*/2}(y_n^*) \subset \RR^4 \less \Xi_T$ is \emph{good} in the sense that the energy inequality \eqref{eq:Schlatter_26_time_sm_half_ball} holds for all $m \geq m_*(n)$, where the possible dependence of $m_*(n)$ on $n$ originates with its dependence on $y_n^*$ and the (mild) dependence of $s_m$ in \eqref{eq:Kozono_Maeda_Naito_page_120_dependencies_of_sm_for_upper_bound_energy} on $r_n^*$.

For each good small ball, $B_{r_n^*/2}(y_n^*) \subset \RR^4 \less \Xi_T$,
Step \ref{step:Local_bounds_curvatures_rescaled_connections_and_convergence_sequence_connections_Coulomb_gauge} provides a sequence of local gauge transformations, $\{\varphi_{m;n}\}_{m\in\NN} \subset H_\Gamma^3(B_{r_n^*/4}(y_n^*); G)$, and a connection, $A_{\infty;n}$, of class $H^2$ on $B_{r_n^*/4}(y_n^*)\times G$ such that
\begin{subequations}
\label{eq:Schlatter_37}
\begin{align}
\label{eq:Schlatter_37_weakly_in_H2}
\varphi_{m;n}^*A_m(s_m) &\rightharpoonup A_{\infty;n} \quad\hbox{weakly in }
H_\Gamma^2(B_{r_n^*/4}(y_n^*); \Lambda^1\otimes\fg),
\\
\label{eq:Schlatter_39_strongly_in_W1q}
\varphi_{m;n}^*A_m(s_m) &\to A_{\infty;n} \quad\hbox{strongly in } W_{\Gamma,\loc}^{1,q}(B_{r_n^*/4}(y_n^*); \Lambda^1\otimes\fg), \quad\hbox{as } m \to \infty.
\end{align}
\end{subequations}
By passing to a diagonal subsequence, we may arrange that the preceding convergence holds \emph{simultaneously} for all $n \in \NN$, as $m \to \infty$.  Lemma \ref{lem:Schlatter_2-8} implies that, after passing to a subsequence again, we obtain convergence on each overlap $U_{n,n'} := B_{r_n^*/4}(y_n^*)\cap B_{r_n^*/4}(y_{n'}^*)$ as $m \to \infty$ for the induced transition functions,
$$
\phi_{m; n,n'} := \varphi_{m,n}\circ\varphi_{m,n'}^{-1} \in H_\Gamma^3(U_{n,n'}; G),
$$
in the sense that,
\begin{subequations}
\label{eq:Schlatter_39}
\begin{align}
\label{eq:Schlatter_39_weakly_in_H3}
\phi_{m; n,n'} &\rightharpoonup \phi_{\infty; n,n'} \quad\hbox{weakly in } H_\Gamma^3(U_{n,n'}; G),
\\
\label{eq:Schlatter_39_strongly_in_W2q}
\phi_{m; n,n'} &\to \phi_{\infty; n,n'} \quad\hbox{strongly in } W_\Gamma^{2,q}(U_{n,n'}; G),
\\
\label{eq:Schlatter_39_strongly_in_C0}
\phi_{m; n,n'} &\to \phi_{\infty; n,n'} \quad\hbox{strongly in } C(\bar U_{n,n'}; G), \quad\hbox{as } m \to \infty,
\end{align}
\end{subequations}
where the strong convergence in $C^0$ follows from the Sobolev embedding $W^{2,q}(U_{n,n'};\RR) \hookrightarrow C(\bar U_{n,n'};\RR)$ \cite[Theorem 4.12]{AdamsFournier}.

With the mild exception of Donaldson and Kronheimer's statement of the Uhlenbeck Removable Singularity Theorem for a finite-energy anti-self-dual connection over a punctured ball (see \cite[Theorem 4.4.12]{DK}, the patching argument clearly described in \cite[Section 4.4.2]{DK} applies to any sequence of $G$ connections of class $H^2$ or $W^{1,p}$ with $p > 2$ (rather than $C^\infty$) and compact Lie group, $G$, rather than the unitary group, $\U(r)$. More precisely, the argument in \cite[Section 4.4.2]{DK} applies as follows:
\begin{enumerate}
\item In \cite[Lemma 4.4.5 (page 159) and Lemma 4.4.6 (page 160)]{DK}, `convergence' of a sequence of connections can be taken to mean strong convergence in $W_{\Gamma,\loc}^{1,p}(\Omega)$, with $p > 2$, rather than $C^\infty(\Omega)$, for any open subset $\Omega \subseteqq \RR^4 \less \Xi_T$.

\item In \cite[Lemma 4.4.6 (page 160)]{DK}, we choose $\{\Omega_i\}_{i\in\NN}$ to be an exhaustion of an open subset $\Omega \subseteqq \RR^4 \less \Xi_T$ by an increasing sequence of precompact open subsets,
$$
\Omega_1 \Subset \cdots \Subset \Omega_i \Subset \Omega_{i+1} \Subset \cdots \Omega
\quad\hbox{and}\quad
\bigcup_{i\in\NN} \Omega_i = \Omega.
$$
\item In \cite[Lemma 4.4.7]{DK}, the sequences of connections and gauge transformations over the union of a pair of overlapping open subsets of $\RR^4$ can be taken to be of class $W^{1,p}$ and $W^{2,p}$, respectively, for $p > 2$, rather than $C^\infty(\Omega)$.

\item In \cite[Corollary 4.4.8]{DK}, the sequence of connections on $P\restriction \Omega$ and global gauge transformations on $P\restriction \Omega$ can be taken to be of class $W^{1,p}$ and $W^{2,p}$, respectively, for $p > 2$, rather than $C^\infty(\Omega)$, for an open subset $\Omega \subseteqq \RR^4 \less \Xi_T$.

\item In place of \cite[Proposition 4.4.9]{DK}, we note that \cite[Corollary 4.4.8]{DK} and Uhlenbeck's local Coulomb gauge-fixing result, namely Theorem \ref{thm:Uhlenbeck_Lp_1-3} and Remark \ref{rmk:Uhlenbeck_theorem_1-3_Wkp} instead of \cite[Theorem 2.3.9]{DK}, and the scale invariance of the $L^2$ norm on covariant two-tensors and $L^4$ norm on one-forms in dimension four yield the following analogue of \cite[Proposition 4.4.9]{DK}: If for each point $y_* \in \Omega \subseteqq \RR^4 \less \Xi_T$, there is a small enough ball $B_{r_*}(y_*) \subset \Omega$ such that,
$$
\int_{B_{r_*/2}(y_*)} |F_{A_m}(s_m)|^2 \,d^4y \leq \eps_1, \quad \forall\, m \geq m_* = m_*(y_*),
$$
where $\eps_1$ is the small positive constant in Theorem \ref{thm:Uhlenbeck_Lp_1-3}. Taking $\Omega = \RR^4 \less \Xi_T$ and passing to a diagonal subsequence, the local convergence \eqref{eq:Schlatter_37} and \eqref{eq:Schlatter_39} over sequences of good small balls yields a sequence of global gauge transformations, $\{\Phi_m\}_{m\in\NN} \subset H_{\Gamma,\loc}^3(\Omega; G)$, and a connection, $A_\infty$, of class $H^2$ on $\Omega\times G$ such that
\begin{subequations}
\label{eq:Schlatter_37_global}
\begin{align}
\label{eq:Schlatter_37_global_weakly_in_H2}
\Phi_m^*A_m(s_m) &\rightharpoonup A_\infty \quad\hbox{weakly in } H_{\Gamma,\loc}^2(\Omega; \Lambda^1\otimes\fg),
\\
\label{eq:Schlatter_37_global_strongly_in_W1q}
\Phi_m^*A_m(s_m) &\to A_\infty \quad\hbox{strongly in } W_{\Gamma,\loc}^{1,q}(\Omega; \Lambda^1\otimes\fg), \quad\hbox{as } m \to \infty.
\end{align}
\end{subequations}
The connection, $A_\infty$ on $\Omega\times G$, is Yang-Mills with respect to the standard Euclidean metric on $\RR^4$ by Step \ref{step:Limit_sequence_rescaled_connections_is_Yang-Mills}. This proves Item \eqref{item:Theorem_Kozono_Maeda_Naito_5-4_H2loc_weak_and_W1ploc_strong_convergence_Am}.

\item In place of \cite[Theorem 4.4.12]{DK}, we apply Uhlenbeck's removable singularity result for finite-energy Yang-Mills connections, namely Theorem \ref{thm:Uhlenbeck_removable_singularity_Yang-Mills}, to provide a gauge transformation, $\Phi:\Omega\times G \to \Omega\times G$, and an extension of the gauge-transformed Yang-Mills connection, $\Phi^*A_\infty$ on the product bundle, $\Omega\times G$, over the punctured sphere, $\Omega = S^4 \less \{y_1,\ldots,y_K,\infty\} \cong \RR^4\less\Xi_T$, to a Yang-Mills connection, $\bar A_\infty$ on a principal $G$-bundle, $\bar P_\infty$, over $S^4$.
\end{enumerate}

In the sequel, we denote $\bar A_\infty$ on $\bar P$ simply by $A_\infty$ on $P$. This concludes Step \ref{step:Global_ideal_limit_sequence_rescaled_connections_over_4space}.
\end{step}

\begin{step}[Convergence of measures on compact subsets of $\RR^4$]
\label{step:Convergence_measures_compact_subsets_Euclidean_space}
Item \eqref{item:Theorem_Kozono_Maeda_Naito_5-4_weak_convergence_FAm_measures} --- aside from non-triviality of the ideal limit, $(A_\infty, \{y_{l,1}, \ldots, y_{l,J_l}\})$, and which is proved in Step \ref{step:Non-triviality_ideal_limit_sequence_rescaled_connections} --- follows \mutatis from the proofs of Theorems \ref{thm:Kozono_Maeda_Naito_5-1} and \ref{thm:Kozono_Maeda_Naito_5-3_T_is_infinite}, recalling that the set of descendant bubble points, $\{y_{l,1}, \ldots, y_{l,J_l}\} \subset B_R(0) \subset \RR^4$, is characterized by \eqref{eq:Kozono_Maeda_Naito_5-1_nested_small_bad_ball_in_Euclidean_4pace} and thus $\sE_{l,j} \geq \eps_1$ for all $l, j$. The inequalities, $\sE_{l,j} \leq \sE_l$ for all $l, j$, are provided by \eqref{eq:Upper_bound_energy_Am_ball_lambda_mR}.
\end{step}

\begin{step}[Non-triviality of the ideal limit of the sequence of rescaled connections]
\label{step:Non-triviality_ideal_limit_sequence_rescaled_connections}
Lastly, we observe that the ideal limit, $(A_\infty, \{y_1,\ldots,y_{L_1}\})$, is non-trivial in the sense that at least \emph{one} of the following must hold:
\begin{inparaenum}[\itshape a\upshape)]
\item $A_\infty$ is not a flat connection on $P_\infty$, or
\item $L_1 \geq 1$.
\end{inparaenum}
The first possibility above implies that $A_\infty$ is gauge equivalent to the product connection on $P_\infty \cong S^4 \times G$. (We emphasize that it is incorrect, as assumed in \cite[Lemma 3.1]{Schlatter_1997}, that $A_\infty$ cannot be a flat connection: neither the proofs of \cite[Theorem 5.4]{Kozono_Maeda_Naito_1995} and \cite[Theorem 1.2]{Schlatter_1997} nor the more detailed argument we have presented here exclude the possibility that $A_\infty$ is a flat connection on a principal $G$-bundle $P_\infty$ over $S^4$.) Denoting the (positive) multiplicities of the points $y_j$ by $\sE_{y_j}$ as in \eqref{eq:Limit_r_to_zero_m_to_infinity_energy_over_ball_equals_Eyj_Euclidean_space} for $1 \leq j \leq L_1$, the inequality \eqref{eq:Schlatter_lemma_3-1_inequality} and fact that $\tau_m = t_m + \lambda_m^2s_m \in [t_m - \lambda_m^2\eps, t_m]$ for all $m \in\NN$ and the scale invariance of $L^2$-norms of covariant two-tensors in dimension four yield,
\begin{align*}
\int_{B_{2R}(0)}|F_{A_\infty}(y)|^2\,d^4y + \sum_{y_j \in B_{2R}(0)}\sE_{y_j}
&= \lim_{m\to\infty}\int_{B_{2R}(0)}|F_{A_m}(s_m,y)|^2\,d^4y
\\
&= \lim_{m\to\infty}\int_{B_{2\lambda_m R}(x_m)}|F_A(\tau_m,x)|^2\,d^4x
\\
&\geq \frac{5\eps_1}{8}.
\end{align*}
Hence, the ideal limit, $(A_\infty, \{y_1,\ldots,y_{L_1}\})$, is non-trivial. This completes the proof of the remaining part of Item \eqref{item:Theorem_Kozono_Maeda_Naito_5-4_weak_convergence_FAm_measures} and concludes Step
\ref{step:Non-triviality_ideal_limit_sequence_rescaled_connections}.
\end{step}

\begin{step}[Computation of the positive constants, $\sE_l$, when $T < \infty$ via the curvatures of the sequence of connections, $A_m$, over an exhaustion of $\RR^4$]
\label{step:Limit_m_to_infinity_energy_over_increasing_balls_in_Euclidean_space_equals_E}
We now verify the limit formula \eqref{eq:Limit_m_to_infinity_energy_over_increasing_balls_in_Euclidean_space_equals_E} in Item \eqref{item:Kozono_Maeda_Naito_5-4_limit_m_to_infinity_energy_over_increasing_balls_in_Euclidean_space_equals_E} which, in our simplified notation, becomes
\begin{equation}
\label{eq:Limit_m_to_infinity_energy_over_increasing_balls_in_Euclidean_space_equals_E_simplified_notation}
\sE_{x_0} = \lim_{\zeta\to 0}\lim_{m \to \infty} \int_{B_{\zeta/\lambda_m}(0)} |F_{A_m}(s_m,y)|^2 \, d^4y.
\end{equation}
By Lemma \ref{lem:Kozono_Maeda_Naito_4-10}, recalling that $\tau_m \equiv t_m + \lambda_m^2 s_m \in [t_m - \lambda_m^2\eps, t_m]$ and $s_m \in (-\eps, 0)$, we have
\begin{align*}
\int_{B_{\zeta/\lambda_m}(0)} |F_{A_m}(s_m,y)|^2 \, d^4y &= \int_{B_\zeta(x_m)} |F_A(\tau_m,x)|^2 \, d^4x
\\
&\geq \int_{B_{\zeta/2}(x_m)} |F_A(t_m,x)|^2 \, d^4x - c(\zeta/2)^{-2}(t_m - \tau_m) \sE(A(\tau_m))
\\
&\geq \int_{B_{\zeta/2}(x_m)} |F_A(t_m,x)|^2 \, d^4x - 4c\zeta^{-2}\lambda_m^2\eps \sE(A(0)),
\end{align*}
where we again use the fact that $\sE(A(\tau_m)) \leq \sE(A(0))$ by Lemma \ref{lem:Kozono_Maeda_Naito_4-1} when $t \geq 0$. On the other hand, for any large enough integer $k = k(m) \geq 0$ such that $t_m \leq \tau_{m+k}  \equiv t_{m+k} + \lambda_{m+k}^2 s_{m+k}$, where again $s_{m+k} \in (-\eps, 0)$, we have
\begin{align*}
\int_{B_\zeta(x_m)} |F_A(\tau_{m+k},x)|^2 \, d^4x &\leq \int_{B_{2\zeta}(x_m)} |F_A(t_m,x)|^2 \, d^4x + c(\zeta/2)^{-2}(\tau_{m+k} - t_m) \sE(A(t_m))
\\
&\leq \int_{B_{2\zeta}(x_m)} |F_A(t_m,x)|^2 \, d^4x + 4c\zeta^{-2}(t_{m+k} - t_m - \lambda_{m+k}^2\eps)\sE(A(0)).
\end{align*}
Hence, for $0 < T < \infty$, we have $t_m - \lambda_{m+k(m)} \to 0$ as $m \to \infty$, since $t_m \nearrow T$ as $m \to \infty$, and because $\lambda_m \searrow 0$ as $m \to \infty$, we obtain, for any fixed $\zeta > 0$ (less than one half of the injectivity radius of $(X,g)$),
\begin{align*}
\limsup_{m \to \infty} \int_{B_{\zeta/\lambda_m}(0)} |F_{A_m}(s_m,y)|^2 \, d^4y &\leq  \limsup_{m \to \infty} \int_{B_{2\zeta}(x_m)} |F_A(t_m,x)|^2 \, d^4x,
\\
\liminf_{m \to \infty} \int_{B_{\zeta/\lambda_m}(0)} |F_{A_m}(s_m,y)|^2 \, d^4y &\geq  \liminf_{m \to \infty} \int_{B_{\zeta/2}(x_m)} |F_A(t_m,x)|^2 \, d^4x.
\end{align*}
The conclusion \eqref{eq:Limit_m_to_infinity_energy_over_increasing_balls_in_Euclidean_space_equals_E_simplified_notation}, for $0 < T < \infty$, now follows from the definition \eqref{eq:Feehan_4-3} of $\sE_{x_0}$.
\end{step}

This completes the proof of Theorem \ref{thm:Kozono_Maeda_Naito_5-4}.
\end{proof}

\subsection{Bubble-tree limits for Yang-Mills gradient flow over a four-dimensional manifold}
\label{subsec:Bubbletree_limits_Yang-Mills_gradient_flow_four-manifold}
In this section, we describe the bubble-tree limit process for a solution to Yang-Mills gradient flow over a closed, four-dimensional, Riemannian manifold, building on our prior results \cite{FeehanGeometry} for the bubble-tree compactification of the moduli space of anti-self-dual connections and which is in turn inspired by work of Parker and Wolfson \cite{ParkerHarmonic, ParkerWolfson} on the bubble-tree limit process for sequences of harmonic maps  of Riemann surfaces into a Riemannian manifold and Taubes \cite{TauPath, TauFrame} on the bubble-tree limit process for Yang-Mills connections over a closed, four-dimensional, Riemannian manifold.

In general, when starting from arbitrary initial energy, we only know that the bubble connections on principal $G$-bundles over $S^4$ converge (after rescaling) to centered Yang-Mills connections with Pontrjagin numbers $k_i \geq 1$ (for suitable $G$). If $k_i\geq 2$, these $S^4$ connections may either bubble further (to two or more or $S^4$ connections with instanton numbers $1\leq k_{ij} < k_i$) or, like the case $k_i=1$, converge (after rescaling) to smooth Yang-Mills connections over $S^4$. Note that if we have a bubble-tree limit, then we may need to consider a nested sequence of gradient-like flows until, at the top of the tree, we have convergence to smooth Yang-Mills limits over $S^4$, show these gradient-like flows converge, and then successively step down the tree and inductively establish gradient-like convergence as we move down each level to the tree base.

See the proof of our \cite[Lemma 4.21]{FL1} and references contained therein, including Donaldson and Kronheimer \cite{DK}, Friedman and Morgan \cite{FrM}, Schlatter \cite{Schlatter_1996, Schlatter_1997}, Struwe \cite{Struwe_1994} (and researchers from Japan with contemporaneous papers), and Taubes \cite{TauPath, TauFrame}.

It is important to note that the results of Kozono, Maeda, Naito \cite{Kozono_Maeda_Naito_1995}, Schlatter \cite{Schlatter_1997}, and Struwe \cite{Struwe_1994} do not eliminate the possibility that, in the case of gradient flow, there may be `bubbles on bubbles'.

Topping \cite[Theorem 3]{Topping_2000} leaves open the possibility of bubble tree occurring at $T < \infty$ for harmonic maps $S^2 \hookrightarrow (M,g)$ but shows can occur when $T=\infty$; Struwe \cite{Struwe_1985} leaves open possibility of bubble tree occurring at $T < \infty$; van der Hout \cite{VanderHout_2003} uses symmetry to rule out bubble trees when $T < \infty$. For further results on bubble trees in harmonic map gradient flow, we refer the reader to Bertsch, van der Hout, and Hulshof \cite{Bertsch_vanderHout_Hulshof_2011} and Topping \cite{Topping_2004am}.

\section[Yang-Mills heat equation and variations in the Riemannian metric]{Continuity and stability of solutions to the Yang-Mills heat equation with respect to variations in the Riemannian metric}
\label{sec:Continuity_solution_Yang-Mills_heat_equation_wrt_Riemannian_metric}
To motivate this section, it is useful to recall Taubes' construction (see, for example, \cite[Section 8]{TauSelfDual}) of a family of almost $g$-anti-self-dual connections, $A'$, on a principal $G$-bundle $P \to X$, where $X$ is assumed to be oriented for the purpose of this introduction. The error, $\|F_{A'}^{+,g}\|_{L^p(X)}$ for $p \in [1,\infty)$, arising in this construction is small enough to allow one to perturb $A'$ to a nearby $g$-anti-self-dual connection, $A'+a$, by solving the first-order anti-self-dual equation for $a \in \Omega^1(X; \ad P)$,
$$
F_{A'+a}^{+,g} = 0  \quad\hbox{over } X.
$$
There are two reasons for the fact that the error, $\|F_{A'}^{+,g}\|_{L^p(X)}$, is non-zero. One is because of the cut-off function, $\chi = \chi_{N,\lambda}$, required to splice a family of anti-self-dual connections, $A_1$, on a principal $G$-bundle $P_1 \to S^4$ with the product connection, $\Gamma$, on $X\times G$, where the splicing takes place over a small
annulus\footnote{Our usage of the scale parameter, $\lambda$, differs from that of Taubes in \cite[Section 8]{TauSelfDual}.},
$$
\Omega(x_0;\lambda/N, \lambda/2) := B_{\lambda/2}(x_0) \less \bar B_{\lambda/N}(x_0) \subset X,
$$
for some large auxiliary parameter, $N \geq 4$, and where the parameter, $\lambda \in (0, \lambda_0]$, is used to rescale the family of anti-self-dual connections, $A_1$, over $S^4$ and $\lambda_0$ is a positive constant that is small relative to the injectivity radius of $(X,g)$. This cut-off function error is supported in $\Omega$. The second source of error, supported in $B_{\lambda/2}(x_0)$, arises whenever the given Riemannian metric, $g$, is not conformally flat near the point $x_0 \in X$ whereas, of course, the standard round metric on $S^4 \less \{s\}$ is conformally flat. However, a choice of geodesic normal coordinates, $\{x^\mu\}_{\mu=1}^4$, near $x_0 \in X$ in the splicing construction ensures that one can regard $g$ as almost flat on $B_{\lambda/2}(x_0)$ and this metric error can also be easily estimated in a useful way via Taubes' analysis \cite[Section 8]{TauSelfDual}.

In the setting of Yang-Mills gradient or heat flow, $A(t)$ on $P \to X$, we shall have occasion to graft $A(t)$ in the reverse direction and create an approximation solution, $\widehat A(t)$, to the Yang-Mills gradient or heat flow equations on $\widehat P \to S^4$ defined by the standard round metric on $S^4$. Equivalently, we shall create a solution to the Yang-Mills gradient-like flow equation on $\widehat P \to S^4$,
$$
\frac{\partial \widehat A(t)}{\partial t} = -d_{\widehat A(t)}^*F_{\widehat A(t)} + R(t) \quad\hbox{over } S^4.
$$
The error, $R(t) = R_\chi(t) + R_g(t)$, is again non-zero because of cut-off function and metric errors, $R_\chi(t)$ and $R_g(t)$, respectively. The cut-off function error in this reverse grafting operation, $R_\chi(t)$, is supported in $\Omega$ and can be easily controlled. The metric error, $R_g(t)$, is supported in the ball, $B_{\lambda/2}(x_0)$, and this unfortunately causes far greater difficulties than in \cite[Section 8]{TauSelfDual}. These difficulties are due, in part, to the fact that the Yang-Mills gradient (and heat) flow equation is a second-order in the spatial directions, whereas the anti-self-dual equation is first-order. A more detailed explanation of the difficulty is explained in Section \ref{sec:Yang-Mills_gradient-like_flow_over_four_sphere}. It is not apparent to the author how to show that the combined error, $R(t) = R_\chi(t) + R_g(t)$, is small enough for $\widehat A(t)$ to be a Yang-Mills gradient-like flow in the sense of Section \ref{sec:Huang_3_gradientlike_system}.

To circumvent this difficulty, we recall that the motivation for our analysis of Yang-Mills gradient-like flow over $S^4$ was to show that the Yang-Mills gradient flow over $X$ does not acquire bubble singularities in finite time. Suppose that bubbling occurs in a solution, $A(t)$, to the Yang-Mills heat equation, for the given Riemannian metric $g$ and initial data, $A_0$, at a \emph{finite} time $T \in (0, \infty)$. Let $\Sigma = \{x_1,\ldots,x_L\} \subset X$ denote the set of bubble points for the solution, $A(t)$, at time $T$. Now consider a Riemannian metric, $\bar g$, obtained from $g$ by `flattening' $g$ on small balls, $B_\rho(x_l)$, around the bubble points, $x_l \in \Sigma$. Let $\bar A(t)$ be the solution to the Yang-Mills heat equation for the Riemannian metric $\bar g$ and initial data $A_0$ and note that $\|g - \bar g\|_{C^1(X)}$ is small, so it is reasonable to expect (as we shall indeed show) that $\bar A(t)$ remains suitably close to $A(t)$ for $t \in [0,T)$. Our goal in this section is to show that $\bar A(t)$ would necessarily bubble in at least one ball, $B_{\lambda/N}(x_l) \subset B_\rho(x_l)$, leading to a contradiction, since we can independently show that a Yang-Mills gradient-like flow, $\widehat A(t)$ over $S^4$, cannot bubble in finite time when the error term, $R_g(t)$, is zero.

\subsection{Overview}
\label{subsec:Overview_continuity_solution_Yang-Mills_heat_equation_wrt_Riemannian_metric}
We begin in Section \ref{subsec:Krylov_2-3-3_covariant_derivative_vector_bundle} by recalling the $L^p$ \apriori global estimate for a solution to the Cauchy problem for the linear heat equation defined by a connection Laplace operator for a smooth reference connection on vector bundle over a closed, Riemannian, smooth manifold, $X$. In our application in Section \ref{sec:Continuity_solution_Yang-Mills_heat_equation_wrt_Riemannian_metric}, we shall only need the case $p=2$ (see Theorem \ref{thm:Krylov_2-5-2_covariant_derivative_vector_bundle}), but we include the statement for $p \in (1, \infty)$ (see Theorem \ref{thm:Krylov_5-2-10_covariant_derivative_vector_bundle}) for the sake of completeness. Section \ref{subsec:Donaldson_Kronheimer_6-3-1} contains a preparatory description of the nonlinearities in the Yang-Mills heat and gradient flow equations which we shall need in our subsequent analysis. Section \ref{subsec:Global_comparison_estimates_two_solutions_Yang-Mills_heat_equations_for_two_metrics} contains our first major result, Theorem \ref{thm:Global_apriori_estimate_difference_solutions_Yang-Mills_heat_equations_pair_metrics}, providing an \apriori estimate for the difference between two solutions, $A(t)$ and $\bar A(t)$, to the Yang-Mills heat equations defined by two nearby Riemannian metrics, $g$ and $\bar g$, on $X$ with common initial data, $A_0$, together with a strict lower bound for the maximal lifetime, $\bar\sT$, of the solution, $\bar A(t)$, in terms of data associated with the solution, $A(t)$, and the pair of Riemannian metrics. In Section \ref{subsec:Local_comparison_estimates_two_solutions_Yang-Mills_heat_equations_for_two_metrics}, we localize the argument in Section \ref{subsec:Global_comparison_estimates_two_solutions_Yang-Mills_heat_equations_for_two_metrics} and establish a local version, Theorem \ref{thm:Local_apriori_estimate_difference_solutions_Yang-Mills_heat_equations_pair_metrics}, of Theorem \ref{thm:Global_apriori_estimate_difference_solutions_Yang-Mills_heat_equations_pair_metrics}, which we will ultimately apply to an open subset, $U \subset X$, given by the complement of a collection of small balls centered at points in $X$ where $A(t)$ develops bubble singularities as $t \nearrow T$. In Section \ref{subsec:Bubbling_for_metric_g_implies_bubbling_for_metric_barg}, we show (Theorem \ref{thm:Bubbling_for_g_implies_bubbling_for_barg}) that bubbling in $A(t)$ as $t \nearrow T$ implies for a given Riemannian metric, $g$, implies bubbling in $\bar A(t)$ as $t \nearrow T$ for a nearby Riemannian metric, $\bar g$. Finally, in Section \ref{subsec:Stability_bubbling_wrt_local_flattening_Riemannian_metric}, we specialize the results of Sections \ref{subsec:Local_comparison_estimates_two_solutions_Yang-Mills_heat_equations_for_two_metrics} and \ref{subsec:Bubbling_for_metric_g_implies_bubbling_for_metric_barg} to the case of a Riemannian metric, $\bar g$, obtained by flattening $g$ near the bubble points for the flow, $A(t)$, as exemplified in Corollaries \ref{cor:Stability_bubbling_wrt_local_flattening_Riemannian_metric} and \ref{cor:Bubbling_for_g_implies_bubbling_for_nearby_barg_locally_flattened_finite_number_small_balls}.

\subsection{The \apriori global $L^2$ and $L^p$ estimates for a solution to the linear heat equation}
\label{subsec:Krylov_2-3-3_covariant_derivative_vector_bundle}
We take $X$ to be a $C^\infty$ closed manifold of dimension $d \geq 2$ with Riemannian metric $g$, and $E$ a complex Hermitian (real Riemannian) vector bundle over $X$, and $A$ a $C^\infty$ Hermitian (Riemannian) connection on $E$ with covariant derivative denoted by $\nabla_A$. We define the corresponding heat operator on $C^\infty(\RR\times X; E)$ by
\begin{equation}
\label{eq:Heat_operator_on_sections_vectorbundle_over_manifold}
L_A := \frac{\partial}{\partial t} + \nabla_A^*\nabla_A \quad\hbox{on } C^\infty(\RR \times X; E),
\end{equation}
The following \apriori global $L^2$ estimate is an analogue of \cite[Theorems 2.3.1, 2.3.2, and 2.5.2]{Krylov_LecturesSobolev} for the scalar parabolic equation \cite[Equation (2.3.1)]{Krylov_LecturesSobolev} on $(-\infty, T) \times \RR^d$, where $T \in (-\infty, \infty]$ (compare \cite[pages 46 and 68]{Krylov_LecturesSobolev}).


\begin{thm}[Global \apriori $L^2$ estimate for the linear heat operator over a closed Riemannian manifold]
\label{thm:Krylov_2-3-1_and_2-3-2_covariant_derivative_vector_bundle}
Let $X$ be a $C^\infty$ closed manifold of dimension $d \geq 2$ with Riemannian metric $g$, and $E$ a complex Hermitian (real Riemannian) vector bundle over $X$, and $A$ a $C^\infty$ Hermitian (Riemannian) connection on $E$. Then there exist a constant $\mu_0 \geq 1$ and, given an integer $k \geq 0$, a positive constant, $C$, with the following significance. If $T \in (-\infty, \infty]$ and $u \in C_0^\infty((-\infty,T]\times X; E)$ and $\mu \geq \mu_0$, then
\begin{multline}
\label{eq:Krylov_2-3-2_covariant_derivative_vector_bundle}
\sum_{j=0}^k\left(\mu\|\nabla_A^j u\|_{L^2((-\infty, T) \times X)} + \sqrt{\mu}\|\nabla_A^{j+1} u\|_{L^2((-\infty, T) \times X)}
\right.
\\
+ \left. \|\nabla_A^{j+2} u\|_{L^2((-\infty, T) \times X)} + \|\nabla_A^j \partial_t u\|_{L^2((-\infty, T) \times X)}\right)
\\
\leq
C\sum_{j=0}^k \|\nabla_A^j(L_A + \mu)u\|_{L^2((-\infty, T) \times X)}.
\end{multline}
\end{thm}

Of course, by approximation, the estimate \eqref{eq:Krylov_2-3-2_covariant_derivative_vector_bundle} continues to hold for $u$ belonging to the Sobolev space of sections of $E$ over $(-\infty,T)\times X$ for which the left-hand side is defined. The \apriori global $L^2$ estimate for a solution to the Cauchy problem is provided by the following analogue of \cite[Theorem 2.5.2]{Krylov_LecturesSobolev} for the scalar parabolic equation \cite[Equation (2.3.1)]{Krylov_LecturesSobolev} on $(0, T) \times \RR^d$, where $T \in (-\infty, \infty]$ in the case $k=0$.

\begin{thm}[Global \apriori $L^2$ estimate for the Cauchy problem for a linear heat operator over a closed Riemannian manifold]
\label{thm:Krylov_2-5-2_covariant_derivative_vector_bundle}
Let $X$ be a $C^\infty$ closed manifold of dimension $d \geq 2$ with Riemannian metric $g$, and $E$ a complex Hermitian (real Riemannian) vector bundle over $X$, and $A$ a $C^\infty$ Hermitian (Riemannian) connection on $E$. Then there exist a constant $\mu_0 \geq 1$ and, given an integer $k \geq 0$, a positive constant, $C$, with the following significance. If $T \in (-\infty, \infty]$ and $u \in C_0^\infty([0,T]\times X; E)$ and $\mu \geq \mu_0$, then
\begin{multline}
\label{eq:Krylov_theorem_2-5-2_covariant_derivative_vector_bundle}
\sum_{j=0}^k\left(\mu\|\nabla_A^j u\|_{L^2((0, T) \times X)} + \sqrt{\mu}\|\nabla_A^{j+1} u\|_{L^2((0, T) \times X)} \right.
\\
+ \left. \|\nabla_A^{j+2} u\|_{L^2((0, T) \times X)} + \|\nabla_A^j \partial_t u\|_{L^2((0, T) \times X)}\right)
\\
\leq
C\sum_{j=0}^k \|\nabla_A^j(L_A + \mu)u\|_{L^2((0, T) \times X)} + C\|u(0)\|_{H_A^{k+2}(X)}.
\end{multline}
\end{thm}

Again, by approximation, the estimate \eqref{eq:Krylov_theorem_2-5-2_covariant_derivative_vector_bundle} continues to hold for $u$ belonging to the Sobolev space of sections of $E$ over $(0,T)\times X$ for which the left-hand side is defined and initial data, $u(0)$, which is acquired in the Sobolev sense.

\begin{rmk}[On the proofs of Theorems \ref{thm:Krylov_2-3-1_and_2-3-2_covariant_derivative_vector_bundle} and \ref{thm:Krylov_2-5-2_covariant_derivative_vector_bundle}]
\label{rmk:Proofs_theorem_Krylov_2-3-1_and_2-3-2_and_2-5-2_covariant_derivative_vector_bundle}
We shall only need the case of $k=0$ in our application of Theorem \ref{thm:Krylov_2-5-2_covariant_derivative_vector_bundle} in our application to the question of continuity of the solution, $A(t)$, to the Yang-Mills heat equation \eqref{eq:Yang-Mills_heat_equation_as_perturbation_rough_Laplacian_heat_equation} with respect to variations in the Riemannian metric, $g$, on $X$ and that \apriori estimate is well known --- see, for example, Corollary \ref{cor:Sell_You_42-12_heat_equation_alpha_is_one_apriori_estimates} or Lemma \ref{lem:Struwe_3-2} when $\mu=1$. The argument for the cases $k \geq 1$ follows in a standard way by mimicking the proof of \cite[Theorem 2.3.2]{Krylov_LecturesSobolev}; alternatively, one can appeal to Corollary \ref{cor:Sell_You_42-12_heat_equation_alpha_is_real_apriori_estimates} when $\alpha=k$ (and $\mu=1$).
\end{rmk}

\begin{rmk}[On the regularity of the initial data in Theorem \ref{thm:Krylov_2-5-2_covariant_derivative_vector_bundle}]
\label{rmk:Regularity_initial_data_theorem_Krylov_2-5-2_covariant_derivative_vector_bundle}
The norm, $\|u(0)\|_{H_A^{k+2}(X)}$, of the initial data, $u(0)$, in the \apriori estimate \eqref{eq:Krylov_theorem_2-5-2_covariant_derivative_vector_bundle} is what one expects from first deriving the \apriori estimate for homogeneous initial data, $u(0)=0$, which is primary case we shall need. However, when $k=0$, one can improve \eqref{eq:Krylov_theorem_2-5-2_covariant_derivative_vector_bundle} and replace $\|u(0)\|_{H_A^2(X)}$ by $\|u(0)\|_{H_A^1(X)}$ by adapting Evan's proof of his \cite[Theorem 7.1.5 (i)]{Evans2} and, when $k \geq 1$, replace $\|u(0)\|_{H_A^{k+2}(X)}$ by $\|u(0)\|_{H_A^{k+1}(X)}$ by adapting his proof of \cite[Theorem 7.1.6]{Evans2}.
\end{rmk}

The following \apriori global $L^p$ estimate is an analogue of \cite[Theorems 5.2.1 and 5.2.2]{Krylov_LecturesSobolev} for the scalar parabolic equation \cite[Equation (5.2.1)]{Krylov_LecturesSobolev} on $\RR^{d+1}$.

\begin{thm}[Global \apriori $L^p$ estimate for a linear heat operator over a closed Riemannian manifold]
\label{thm:Krylov_5-2-1_and_2_covariant_derivative_vector_bundle}
Let $p \in (1, \infty)$, and $X$ be a $C^\infty$ closed manifold of dimension $d \geq 2$ with Riemannian metric $g$, and $E$ a complex Hermitian (real Riemannian) vector bundle over $X$, and $A$ a $C^\infty$ Hermitian (Riemannian) connection on $E$. Then there exist a constant $\mu_0 \geq 1$ and, given an integer $k \geq 0$, a positive constant, $C$, with the following significance. If $T \in (-\infty, \infty]$ and $u \in C_0^\infty((-\infty,T]\times X; E)$ and $\mu \geq \mu_0$, then
\begin{multline}
\label{eq:Krylov_5-2-4_covariant_derivative_vector_bundle}
\sum_{j=0}^k\left(\mu\|\nabla_A^j u\|_{L^p((-\infty,T)\times X)} + \sqrt{\mu}\|\nabla_A^{j+1} u\|_{L^p((-\infty,T)\times X)} \right.
\\
+ \left. \|\nabla_A^{j+2} u\|_{L^p((-\infty,T)\times X)} + \|\nabla_A^j \partial_t u\|_{L^p((-\infty,T)\times X)}\right)
\\
\leq
C\sum_{j=0}^k \|\nabla_A^j(L_A + \mu)u\|_{L^p((-\infty,T)\times X)}.
\end{multline}
\end{thm}

The \apriori global $L^p$ estimate for a solution to the Cauchy problem is provided by the following analogue of \cite[Theorem 5.2.10]{Krylov_LecturesSobolev} for the scalar parabolic equation \cite[Equation (2.3.1)]{Krylov_LecturesSobolev} on $(0, T) \times \RR^d$, where $T \in (-\infty, \infty]$ in the case $k=0$.

\begin{thm}[Global \apriori $L^p$ estimate for the Cauchy problem for a linear heat operator over a closed Riemannian manifold]
\label{thm:Krylov_5-2-10_covariant_derivative_vector_bundle}
Let $p \in (1, \infty)$, and $X$ be a $C^\infty$ closed manifold of dimension $d \geq 2$ with Riemannian metric $g$, and $E$ a complex Hermitian (real Riemannian) vector bundle over $X$, and $A$ a $C^\infty$ Hermitian (Riemannian) connection on $E$. Then there exist a constant $\mu_0 \geq 1$ and, given an integer $k \geq 0$, a positive constant, $C$, with the following significance. If $T \in (-\infty, \infty]$ and $u \in C_0^\infty([0,T]\times X; E)$ and $\mu \geq \mu_0$, then
\begin{multline}
\label{eq:Krylov_5-2-10_theorem_covariant_derivative_vector_bundle}
\sum_{j=0}^k\left(\mu\|\nabla_A^j u\|_{L^p((-\infty,T)\times X))} + \sqrt{\mu}\|\nabla_A^{j+1} u\|_{L^p((-\infty,T)\times X))}
\right.
\\
+ \left. \|\nabla_A^{j+2} u\|_{L^p((-\infty,T)\times X))} + \|\nabla_A^j \partial_t u\|_{L^p((-\infty,T)\times X))}\right)
\\
\leq
C\sum_{j=0}^k \|\nabla_A^j(L_A + \mu)u\|_{L^p((-\infty,T)\times X))} + C\|u(0)\|_{W_A^{k+2,p}((-\infty,T)\times X))}.
\end{multline}
\end{thm}

As usual, by approximation, the estimates \eqref{eq:Krylov_5-2-4_covariant_derivative_vector_bundle} and \eqref{eq:Krylov_5-2-10_theorem_covariant_derivative_vector_bundle} continue to hold for $u$ belonging to the Sobolev space of sections of $E$ over $(0,T)\times X$ for which the left-hand sides are defined or initial data, $u(0)$, which is acquired in the Sobolev sense.

\begin{rmk}[On the proof of Theorem \ref{thm:Krylov_5-2-1_and_2_covariant_derivative_vector_bundle}]
\label{rmk:Proof_theorem_Krylov_5-2-1_and_2_covariant_derivative_vector_bundle}
We include the statement of Theorem \ref{thm:Krylov_5-2-1_and_2_covariant_derivative_vector_bundle} for completeness; the simpler Theorem \ref{thm:Krylov_2-3-1_and_2-3-2_covariant_derivative_vector_bundle} will suffice for our application. Theorem \ref{thm:Krylov_5-2-1_and_2_covariant_derivative_vector_bundle} may be established by adapting Krylov's proofs of his \cite[Theorems 5.2.1 and 5.2.2]{Krylov_LecturesSobolev} with the aid of the vector-valued Calder\'on-Zygmund inequality --- see, for example, Taylor \cite{Taylor_PDO_and_NLPDE}.
\end{rmk}

\subsection{The Yang-Mills heat and gradient flow equations revisited}
\label{subsec:Donaldson_Kronheimer_6-3-1}
Recall that $A(t) = A_1 + a(t)$ solves the nonlinear Yang-Mills heat equation, viewed as a perturbation,
\begin{equation}
\label{eq:Yang-Mills_heat_equation_as_perturbation_rough_Laplacian_heat_equation}
\frac{\partial a}{\partial t} + \nabla_{A_1}^*\nabla_{A_1}a(t) + \sF(a(t)) = 0, \quad\forall\, t > 0,
\end{equation}
of the \emph{linear} heat equation defined by the \emph{connection Laplacian}, $\nabla_{A_1}^*\nabla_{A_1}$, that is,
\begin{equation}
\label{eq:Linear_heat_equation_with_rough_Laplacian_on_Omega_1_adP}
\frac{\partial a}{\partial t} + \nabla_{A_1}^*\nabla_{A_1}a(t) = f(t)
\quad\hbox{in } \Omega^1(X; \Lambda^1\otimes\ad P), \quad \hbox{for } t > 0,
\end{equation}
with $f(t) := -\sF(a(t))$, we obtain the following schematic expression for the \emph{Yang-Mills heat equation non linearity},
\begin{multline}
\label{eq:Yang-Mills_heat_equation_nonlinearity_relative_rough_Laplacian_plus_sign}
\sF(a) := d_{A_1}^*F_{A_1} + F_{A_1}\times a + \Ric(g)\times a + \nabla_{A_1}a\times a + a\times a\times a,
\\
\forall\, a \in \Omega^1(X; \ad P),
\end{multline}
where we recall that $\Ric(g)$ denotes the Ricci curvature tensor of a Riemannian metric, $g$, on $X$. We shall be comparing solutions to Yang-Mills heat equations defined by metrics $g$ and $\mathring{g}$ that are $C^1(X)$ close but that may not be $C^2(X)$ close and for which it may not be true that $\Ric(g)$ and $\Ric(\mathring{g})$ are $C^0(X)$-close, although they may be $L^q(X)$-close for any $q \in [1, \infty)$. We shall need a more explicit expression for the nonlinearity, $\sF(a)$, and one that makes the role of the Hodge star operator clear.

Therefore, we shall first view the Yang-Mills heat equation as a perturbation of the linear heat equation defined by the Hodge Laplace operator \eqref{eq:Lawson_page_93_Hodge_Laplacian} on $\Omega^1(X; \ad P)$. Namely, writing $A(t) = A_1 + a(t)$ and writing the Yang-Mills heat equation as in \cite[Equation (6.3.3)]{DK},
\begin{equation}
\label{eq:Donaldson_Kronheimer_6-3-3}
\frac{\partial A}{\partial t} + d_A^*F_A + d_Ad_A^* a = 0,
\end{equation}
we expand the elliptic part of \eqref{eq:Donaldson_Kronheimer_6-3-3} as
\begin{align*}
d_A^*F_A + d_Ad_A^* a &= d_{A_1+a}^*F_{A_1+a} + d_{A_1 + a}d_{A_1+a}^* a
\\
&= d_{A_1+a}^*(F_{A_1} + d_{A_1}a + [a, a]) + d_{A_1 + a}(d_{A_1}^*a - *[a, *a])
\\
&= d_{A_1}^*F_{A_1} + d_{A_1}^*d_{A_1}a + d_{A_1}^*([a, a]) - *[a, *F_{A_1}] - *[a, *d_{A_1}a] - *[a, *[a, a]]
\\
&\quad + d_{A_1}d_{A_1}^*a - d_{A_1}(*[a, *a]) + [a, d_{A_1}^*a] - *[a, *[a, a]]
\\
&= (d_{A_1}^*d_{A_1} + d_{A_1}d_{A_1}^*)a + \sG(a)
\\
&= \Delta_{A_1}^*a + \sG(a),
\end{align*}
where the induced nonlinearity is defined by the expression,
\begin{align*}
\sG(a) &:= d_{A_1}^*F_{A_1} + d_{A_1}^*([a, a]) - *[a, *F_{A_1}] - *[a, *d_{A_1}a] - 2*[a, *[a, a]]
\\
&\quad - d_{A_1}(*[a, *a]) + [a, d_{A_1}^*a].
\end{align*}
While the signs and values of universal constants (factors of $2$ or $1/2$ and so on) appearing in the expression for $\sG(a)$ are unimportant, we do need to keep track of the presence of the Hodge $*$-operators. We have
\begin{equation}
\label{eq:Donaldson_Kronheimer_6-3-3_Yang-Mills_heat_equation_with_Hodge_Laplacian_nonlinearity}
\frac{\partial a}{\partial t} + \Delta_{A_1}^*a + \sG(a) = 0.
\end{equation}
The Bochner-Weitzenb\"ock formula \eqref{eq:Lawson_corollary_II-2} allows us to further refine our expression for the Yang-Mills heat equation as
\begin{equation}
\label{eq:Donaldson_Kronheimer_6-3-3_Yang-Mills_heat_equation_with_rough_Laplacian_nonlinearity}
\frac{\partial a}{\partial t} + \nabla_{A_1}^*\nabla_{A_1}a + \sF(a) = 0,
\end{equation}
where
\begin{multline}
\label{eq:Donaldson_Kronheimer_6-3-3_rough_Laplacian_nonlinearity}
\sF(a) := d_{A_1}^*F_{A_1} + d_{A_1}^*([a, a]) - *[a, *F_{A_1}] - *[a, *d_{A_1}a] - 2*[a, *[a, a]]
\\
- d_{A_1}(*[a, *a]) + [a, d_{A_1}^*a] + \{\Ric(g), a\} + \{F_{A_1}, a\}
\end{multline}
where we recall that $\{, \}$ denotes a universal bilinear expression (independent of the Riemannian metric on $X$).

\begin{rmk}
See \cite[page 235]{DK} for the identity $d_A^*d_A^*F_A = \{F_A, F_A\} \equiv 0$. See \cite[Section 6.3.1]{DK} for a detailed explanation of the Donaldson-DeTurck trick and the role of irreducibility of the connections, $A(t)$, for $t\in [0, T)$.
\end{rmk}

\subsection{Global comparison estimate for two solutions to the Yang-Mills heat equations defined by two nearby Riemannian metrics}
\label{subsec:Global_comparison_estimates_two_solutions_Yang-Mills_heat_equations_for_two_metrics}
We wish to compare the pair of solutions, $A(t)$ and $\bar A(t)$, to the Yang-Mills heat equation \eqref{eq:Yang-Mills_heat_equation_as_perturbation_rough_Laplacian_heat_equation} defined by two close Riemannian metrics, $g$ and $\bar g$, and common initial data, $A_0$. Our eventual goal is to consider the case where $\bar g$ is obtained by flattening $g$ on small balls around the bubble points in the set $\Sigma = \{x_1, \ldots, x_L\} \subset X$, when $A(t)$ acquires bubble singularities at time $T$. However, it will prove very useful to first consider the simpler case where
$$
A - A_1 \in C([0,T]; H_{A_1}^1(X;\Lambda^1\otimes\ad P)),
$$
and so the solution, $A(t)$, exists for $t \in [0, T+\delta)$ and some positive constant, $\delta$. For this simpler case, our aim is to show that if $\bar g$ is close enough to $g$, then
$$
\bar A - A_1 \in C([0,T]; H_{A_1}^1(X;\Lambda^1\otimes\ad P)),
$$
and so $\bar A(t)$ also exists for $t \in [0, T+\bar\delta)$ and some positive constant, $\bar\delta$, and remains close to $A(t)$ for $t \in [0, T]$. In fact, even the assumption that $A - A_1 \in C([0,T]; H_{A_1}^1(X;\Lambda^1\otimes\ad P))$ can be relaxed and still yield the conclusion, $\bar A - A_1 \in C([0,T]; H_{A_1}^1(X;\Lambda^1\otimes\ad P))$, as the statement and proof of Theorem \ref{thm:Global_apriori_estimate_difference_solutions_Yang-Mills_heat_equations_pair_metrics} below illustrates. Recall that $\Riem(g) \equiv \Riem_g$ denotes the full Riemann curvature tensor of a Riemannian metric, $g$, on $X$. In Theorem \ref{thm:Global_apriori_estimate_difference_solutions_Yang-Mills_heat_equations_pair_metrics} below and in the sequel, we regard $g$ as the \emph{reference} Riemannian metric and use that in our definitions of H\"older and Sobolev norms; in Corollary \ref{cor:Global_apriori_estimate_difference_solutions_YM_heat_eqns_pair_metrics_C0_close}, we draw the same conclusion, with additional work, for a pair of metrics that are close in a coarser topology.

\begin{thm}[Global \apriori estimate for the difference between solutions to the Yang-Mills heat equations defined by two close Riemannian metrics and common initial data]
\label{thm:Global_apriori_estimate_difference_solutions_Yang-Mills_heat_equations_pair_metrics}
Let $G$ be a compact Lie group, $P$ a principal $G$-bundle over a closed, four-dimensional, smooth manifold with Riemannian metric $g$, and $A_1$ a fixed reference connection of class $C^\infty$ on $P$, and $K$ and $T$ positive constants. Then there are a small enough constant $\eta = \eta(A_1,g,K,T) \in (0, 1]$, a large enough constant $\mu_0 = \mu_0(A_1,g,K) \in [1, \infty)$, and a positive constant $z_1 = z_1(g,K)$ with the following significance. Let $\bar g$ be a Riemannian metric on $X$ such that\footnote{The Sobolev Embedding Theorem \cite[Theorem 4.12]{AdamsFournier} with $n=4$, $j=0$, $m=2$, and $p=4$ gives $W^{2,4}(X)\hookrightarrow C^\alpha(X)$, for any $\alpha \in [0,1)$.}
\begin{align}
\label{eq:Riem_barg_C0_norm_lessthan_bound}
\|\Riem(\bar g)\|_{C(X)} &\leq K,
\\
\label{eq:Riemannian_metrics_g_minus_barg_C1_plus_W24_norm_lessthan_small_positive_constant}
\|g - \bar g\|_{C^1(X)} + \|g - \bar g\|_{W^{2,4}(X)} &\leq \eta,
\\
\label{eq:Riemannian_metrics_g_minus_barg_C2_lessthan_bound}
\|g - \bar g\|_{C^2(X)} &\leq K.
\end{align}
Let $A(t)$, for $t\in [0, T)$, and $\bar A(t)$, for $t\in [0, \bar\sT)$ with maximal lifetime $\bar\sT \in (0, \infty]$, be strong solutions to the Yang-Mills heat equation \eqref{eq:Yang-Mills_heat_equation_as_perturbation_rough_Laplacian_heat_equation} on $P$ for the Riemannian metrics, $g$ and $\bar g$, respectively, and common initial data, $A_0$. Suppose that $A(t)$ obeys
\begin{equation}
\label{eq:Yang-Mills_gradient_flow_A_minus_A1_Linfty_time_H2_space_lessthan_bound}
\|A - A_1\|_{L^\infty(0,T;H_{A_1}^2(X))} \leq K,
\end{equation}
and that $A_1$ obeys
\begin{equation}
\label{eq:F_A1_C0_norm_lessthan_bound}
\|F_{A_1}\|_{C(X))} \leq K.
\end{equation}
Then $\bar A - A_1 \in C([0,T]; H_{A_1}^1(X;\Lambda^1\otimes\ad P))$ and $\bar\sT > T$ and $\bar A(t)$ obeys
\begin{multline}
\label{eq:Apriori_estimate_for_A_minus_barA_in_terms_C1_and_W24_norms_g_minus_barg}
\|A - \bar A\|_{L^2(0,T;H_{A_1}^2(X))} + \|\partial_t (A - \bar A)\|_{L^2(0,T;L^2(X))}
+ \|A - \bar A\|_{L^\infty(0,T;H_{A_1}^1(X))}
\\
\leq z_1e^{\mu_0T}\left(\|g - \bar g\|_{C^1(X)} + \|g - \bar g\|_{W^{2,4}(X)}\right).
\end{multline}
\end{thm}


\begin{proof}
We find it convenient to write the pair of solutions as $A(t) = A_1 + a(t)$ and $\bar A(t) = A_1 + \bar a(t)$, with common initial data expressed as $a(0) = a_0 = \bar a(0)$, where $a_0 = A_0 - A_1$. For the metric $\bar g$ (and connection $A_1$), we then write the Hodge star operator, connection Laplace operator, exterior covariant derivative adjoint, and Yang-Mills heat equation nonlinearity as $\bar *$, and $\bar\nabla_{A_1}^*\nabla_{A_1}$, and $\bar d_{A_1}^*$, and $\bar\sF(\bar a)$, respectively.

\setcounter{step}{0}
\begin{step}[Derivation of the quasilinear parabolic equation for $a - \bar a$]
\label{step:Derivation_quasilinear_parabolic_equation_for_a_minus_bara}
Beginning with the pair of Cauchy problems for Yang-Mills heat equations \eqref{eq:Yang-Mills_heat_equation_as_perturbation_rough_Laplacian_heat_equation} defined by the metrics $g$ and $\bar g$,
\begin{align*}
\frac{\partial a}{\partial t} + \nabla_{A_1}^*\nabla_{A_1}a + \sF(a) &= 0
\quad\hbox{a.e. on }(0,T)\times X, \quad a(0) = a_0,
\\
\frac{\partial \bar a}{\partial t} + \bar\nabla_{A_1}^*\nabla_{A_1}\bar a + \bar\sF(\bar a) &= 0
\quad\hbox{a.e. on }(0,\bar\sT)\times X, \quad \bar a(0) = a_0,
\end{align*}
we consider the induced quasilinear Cauchy problem for $a - \bar a$,
\begin{multline}
\label{eq:Quasilinear_Cauchy_problem_a_minus_bara}
\frac{\partial}{\partial t}(a - \bar a) + \nabla_{A_1}^*\nabla_{A_1}a -\bar\nabla_{A_1}^*\nabla_{A_1}\bar a + \sF(a) - \bar\sF(\bar a) = 0 \quad\hbox{a.e. on }(0,T_1)\times X,
\\
\quad (a - \bar a)(0) = 0,
\end{multline}
where $T_1 > 0$ obeys $T_1 < \bar\sT$ and $T_1 \leq T$.

To proceed further, we shall need to express the nonlinearity, $\sF(a) - \bar\sF(\bar a)$, as one in powers of $a-\bar a$, with coefficients involving powers of $a(t)$, and also make its dependence on $g-\bar g$ explicit. We shall do this in several stages, first writing the (linear) difference between the connection Laplace operators as
$$
\nabla_{A_1}^*\nabla_{A_1}a -\bar\nabla_{A_1}^*\nabla_{A_1}\bar a
\\
= \bar\nabla_{A_1}^*\nabla_{A_1}(a - \bar a) + \left(\nabla_{A_1}^*\nabla_{A_1} - \bar\nabla_{A_1}^*\nabla_{A_1}\right)a,
$$
and writing the nonlinearity as
\begin{multline*}
\sF(a) - \bar\sF(\bar a)
=
-*[a, *F_{A_1}] + \bar*[\bar a, \bar*F_{A_1}] + \{F_{A_1}, a - \bar a\} + \{\Ric(g) - \Ric(\bar g), a\} + \{\Ric(\bar g), a - \bar a\}
\\
+ d_{A_1}^*([a, a]) - \bar d_{A_1}^*([\bar a, \bar a]) - *[a, *d_{A_1}a] + \bar*[\bar a, \bar*d_{A_1}\bar a]
\\
- d_{A_1}(*[a, *a]) + d_{A_1}(\bar*[\bar a, \bar*\bar a]) + [a, d_{A_1}^*a] - [\bar a, \bar d_{A_1}^*\bar a]
\\
- 2*[a, *[a, a]] + 2\bar*[\bar a, \bar*[\bar a, \bar a]].
\end{multline*}
We wish to write our quasilinear parabolic equation \eqref{eq:Quasilinear_Cauchy_problem_a_minus_bara} for $a - \bar a$ in the form,
\begin{multline*}
\frac{\partial}{\partial t}(a - \bar a) + \bar\nabla_{A_1}^*\nabla_{A_1}(a - \bar a) = \hbox{Sum of powers of $a-\bar a$}
\\
+ \hbox{Sum of powers of $a$ with small $g-\bar g$ coefficients.}
\end{multline*}
To accomplish this, we expand the difference terms appearing in the preceding expression for $\sF(a) - \bar\sF(\bar a)$. The first difference term is
\begin{align*}
\hbox{Term 1} :=&\ *[a, *F_{A_1}] - \bar*[\bar a, \bar*F_{A_1}]
\\
= &\ (* - \bar*)[a, *F_{A_1}] + \bar*([a, *F_{A_1}] - [\bar a, \bar*F_{A_1}])
\\
= &\ (* - \bar*)[a, *F_{A_1}] + \bar*([a, (* - \bar*)F_{A_1}] + [a-\bar a, \bar*F_{A_1}]).
\end{align*}
The second difference term is
\begin{align*}
\hbox{Term 2} :=&\ d_{A_1}^*([a, a]) - \bar d_{A_1}^*([\bar a, \bar a])
\\
=&\  (d_{A_1}^* - \bar d_{A_1}^*)([a, a]) + \bar d_{A_1}^*([a, a] - [\bar a, \bar a])
\\
=&\  (d_{A_1}^* - \bar d_{A_1}^*)([a, a]) + \bar d_{A_1}^*([a, a - \bar a] + [a - \bar a, \bar a])
\\
=&\  (d_{A_1}^* - \bar d_{A_1}^*)([a, a]) + \bar d_{A_1}^*([a, a - \bar a] + [a - \bar a, \bar a - a] + [a - \bar a, a]).
\end{align*}
The third difference term is
\begin{align*}
\hbox{Term 3} :=&\ *[a, *d_{A_1}a] - \bar*[\bar a, \bar*d_{A_1}\bar a]
\\
=&\  (* - \bar*)[a, *d_{A_1}a] + \bar*([a, *d_{A_1}a] - [\bar a, \bar*d_{A_1}\bar a])
\\
=&\  (* - \bar*)[a, *d_{A_1}a] + \bar*([a - \bar a, *d_{A_1}a] + [\bar a, *d_{A_1}a -\bar*d_{A_1}\bar a])
\\
=&\  (* - \bar*)[a, *d_{A_1}a] + \bar*([a - \bar a, *d_{A_1}a] + [a, *d_{A_1}a -\bar*d_{A_1}\bar a])
+ [\bar a - a, *d_{A_1}a -\bar*d_{A_1}\bar a])
\\
=&\  (* - \bar*)[a, *d_{A_1}a] + \bar*([a - \bar a, *d_{A_1}a]
+ [a, (* - \bar*)d_{A_1}a]) + [a, \bar*d_{A_1}(a - \bar a)])
\\
&\quad + [\bar a - a, (* - \bar*)d_{A_1}a]) + [\bar a - a, \bar*d_{A_1}(a - \bar a)]).
\end{align*}
The fourth difference term is
\begin{align*}
\hbox{Term 4} :=&\ d_{A_1}(*[a, *a]) - d_{A_1}(\bar*[\bar a, \bar*\bar a])
\\
=&\  d_{A_1}((* - \bar*)[a, *a]) + d_{A_1}(\bar*([a, *a] - [\bar a, \bar*\bar a]))
\\
=&\  d_{A_1}((* - \bar*)[a, *a]) + d_{A_1}(\bar*[a, *a - \bar*\bar a])
+ d_{A_1}(\bar*[a - \bar a, \bar*\bar a])
\\
=&\  d_{A_1}((* - \bar*)[a, *a])
+ d_{A_1}(\bar*[a, (* - \bar*)a])+ d_{A_1}(\bar*[a, \bar*(a - \bar a)])
+ d_{A_1}(\bar*[a - \bar a, \bar*\bar a])
\\
=&\  d_{A_1}((* - \bar*)[a, *a])
+ d_{A_1}(\bar*[a, (* - \bar*)a])+ d_{A_1}(\bar*[a, \bar*(a - \bar a)])
+ d_{A_1}(\bar*[a - \bar a, \bar*a])
\\
&\quad + d_{A_1}(\bar*[a - \bar a, \bar*(\bar a - a)]).
\end{align*}
The fifth difference term is
\begin{align*}
\hbox{Term 5} :=&\ [a, d_{A_1}^*a] - [\bar a, \bar d_{A_1}^*\bar a]
\\
=&\  [a, (d_{A_1}^*a - \bar d_{A_1}^*\bar a)] + [a - \bar a, \bar d_{A_1}^*\bar a]
\\
=&\  [a, (d_{A_1}^* - \bar d_{A_1}^*)a] + [a, \bar d_{A_1}^*(a - \bar a)]
+ [a - \bar a, \bar d_{A_1}^*\bar a]
\\
=&\  [a, (d_{A_1}^* - \bar d_{A_1}^*)a] + [a, \bar d_{A_1}^*(a - \bar a)]
+ [a - \bar a, \bar d_{A_1}^*a]+ [a - \bar a, \bar d_{A_1}^*(\bar a - a)].
\end{align*}
The sixth difference term is
\begin{align*}
\hbox{Term 6} :=&\ *[a, *[a, a]] - \bar*[\bar a, \bar*[\bar a, \bar a]]
\\
=&\  (* - \bar*)[a, *[a, a]] + \bar*([a, *[a, a]] - [\bar a, \bar*[\bar a, \bar a]])
\\
=&\  (* - \bar*)[a, *[a, a]] + \bar*[a - \bar a, *[a, a]] + \bar*[\bar a, *[a, a] - \bar*[\bar a, \bar a]]
\\
=&\  (* - \bar*)[a, *[a, a]] + \bar*[a - \bar a, *[a, a]] + \bar*[\bar a, (* - \bar*)[a, a]]
+ \bar*[\bar a, \bar*([a, a] - [\bar a, \bar a])]
\\
=&\  (* - \bar*)[a, *[a, a]] + \bar*[a - \bar a, *[a, a]] + \bar*[\bar a, (* - \bar*)[a, a]]
\\
&\qquad + \bar*[\bar a, \bar*[a, a - \bar a]] + \bar*[\bar a, \bar*[a - \bar a, \bar a]]
\\
=&\  (* - \bar*)[a, *[a, a]] + \bar*[a - \bar a, *[a, a]] + \bar*[a, (* - \bar*)[a, a]] + \bar*[\bar a - a, (* - \bar*)[a, a]]
\\
&\qquad + \bar*[a, \bar*[a, a - \bar a]] + \bar*[\bar a - a, \bar*[a, a - \bar a]]
 + \bar*[a, \bar*[a - \bar a, \bar a]] + \bar*[\bar a - a, \bar*[a - \bar a, \bar a]]
\\
=&\ (* - \bar*)[a, *[a, a]] + \bar*[a - \bar a, *[a, a]] + \bar*[a, (* - \bar*)[a, a]] + \bar*[\bar a - a, (* - \bar*)[a, a]]
\\
&\qquad + \bar*[a, \bar*[a, a - \bar a]] + \bar*[\bar a - a, \bar*[a, a - \bar a]]
 + \bar*[a, \bar*[a - \bar a, a]] + \bar*[a, \bar*[a - \bar a, \bar a - a]]
\\
&\qquad + \bar*[\bar a - a, \bar*[a - \bar a, a]] + \bar*[\bar a - a, \bar*[a - \bar a, \bar a - a]].
\end{align*}
With the aid of the formula \eqref{eq:Donaldson_Kronheimer_6-3-3_rough_Laplacian_nonlinearity} for $\sF(a)$ (and thus also $\bar\sF(\bar a)$), the quasilinear parabolic equation \eqref{eq:Quasilinear_Cauchy_problem_a_minus_bara} for $a-\bar a$ takes the following schematic form,
\begin{multline}
\label{eq:Quasilinear_Cauchy_problem_a_minus_bara_schematic_sum_terms_1_to_6}
\frac{\partial}{\partial t}(a - \bar a) + \bar\nabla_{A_1}^*\nabla_{A_1}(a - \bar a)
+ \{F_{A_1}, a - \bar a\} + \{\Ric(\bar g), a - \bar a\}
\\
= -\left(\nabla_{A_1}^*\nabla_{A_1} - \bar\nabla_{A_1}^*\nabla_{A_1}\right)a - \{\Ric(g) - \Ric(\bar g), a\}
+ \hbox{Sum of Terms 1 through 6}.
\end{multline}
When expressed schematically, the Terms 1 through 6 which we have identified simplify as follows:
$$
\hbox{Term 1} = (* - \bar*)F_{A_1} \times a + F_{A_1} \times (a-\bar a);
$$
and
\begin{align*}
\hbox{Term 2} &= a \times (\nabla_{A_1}^* - \bar\nabla_{A_1}^*)a
+ (a - \bar a) \times \bar\nabla_{A_1}^*a + a \times \bar\nabla_{A_1}^*(a - \bar a)
\\
&\quad + (a - \bar a) \times (a - \bar a) + a \times (a - \bar a);
\end{align*}
and
\begin{align*}
\hbox{Term 3} &= (* - \bar*)a \times \nabla_{A_1}a + (a - \bar a) \times \nabla_{A_1}a + a \times \nabla_{A_1}(a - \bar a)
\\
&\quad + (a - \bar a) \times \nabla_{A_1}(a - \bar a);
\end{align*}
and
\begin{align*}
\hbox{Term 4} &= (\nabla(* - \bar*))a \times a + (* - \bar*)a \times \nabla_{A_1}a
+ a \times \nabla_{A_1}(a - \bar a)
\\
&\quad + (a - \bar a) \times \nabla_{A_1}a + (a - \bar a) \times \nabla_{A_1}(a - \bar a);
\end{align*}
and
\begin{align*}
\hbox{Term 5} &= a \times (\nabla_{A_1}^* - \bar\nabla_{A_1}^*)a + a \times \bar\nabla_{A_1}^*(a - \bar a)
\\
&\quad + (a - \bar a) \times \bar\nabla_{A_1}^*a + (a - \bar a) \times \bar\nabla_{A_1}^*(a - \bar a);
\end{align*}
and
\begin{align*}
\hbox{Term 6} &= (* - \bar*)a \times a \times a + a \times a \times (a - \bar a) + a \times (a - \bar a) \times (a - \bar a)
\\
&\quad + (a - \bar a) \times (a - \bar a) \times (a - \bar a).
\end{align*}
We can now reorganize our expression \eqref{eq:Quasilinear_Cauchy_problem_a_minus_bara_schematic_sum_terms_1_to_6} for the quasilinear parabolic equation for $a-\bar a$ in order to better focus on its essential structure:
\begin{equation}
\label{eq:Quasilinear_parabolic_equation_a_minus_bara}
\begin{aligned}
{}&\frac{\partial}{\partial t}(a - \bar a) + \bar\nabla_{A_1}^*\nabla_{A_1}(a - \bar a)
+ \Ric(\bar g) \times (a - \bar a) + F_{A_1} \times (a-\bar a)
\\
&\quad + (a - \bar a) \times \nabla_{A_1}a +  a \times (a - \bar a)
+ a \times \nabla_{A_1}(a - \bar a) + a \times a \times (a - \bar a)
\\
&= -\left(\nabla_{A_1}^*\nabla_{A_1} - \bar\nabla_{A_1}^*\nabla_{A_1}\right)a
- (\Ric(g) - \Ric(\bar g)) \times a
\\
&\quad + (* - \bar*)F_{A_1} \times a + a \times (\nabla_{A_1}^* - \bar\nabla_{A_1}^*)a
 + (* - \bar*)a \times \nabla_{A_1}a + (\nabla(* - \bar*))a \times a
\\
&\quad + (* - \bar*)a \times a \times a
\\
&\quad + (a - \bar a) \times (a - \bar a) + (a - \bar a) \times \nabla_{A_1}(a - \bar a)
+ a \times (a - \bar a) \times (a - \bar a)
\\
&\quad + (a - \bar a) \times (a - \bar a) \times (a - \bar a) \quad\hbox{a.e. on }(0,T_1)\times X.
\end{aligned}
\end{equation}
This concludes Step \ref{step:Derivation_quasilinear_parabolic_equation_for_a_minus_bara}.
\end{step}

\begin{step}[Passage to an equivalent exponentially shifted quasilinear parabolic equation for $b - \bar b$]
\label{step:Exponentially_shifted_quasilinear_parabolic_equation}
It is useful at this stage to digress briefly and recall the exponential shift trick in the present context. For convenience, define $X_T := (0, T)\times X$. Given $f \in \Omega^1(X_T; \ad P)$, suppose $u \in \Omega^1(X_T; \ad P)$ obeys
$$
\frac{\partial u}{\partial t} + \nabla_{A_1}^*\nabla_{A_1}u = f
\quad\hbox{on } X_T.
$$
If $u =: e^{\mu t}v$ for a constant $\mu \geq 0$, then $\partial_t u = e^{\mu t}(\mu v + \partial_t v)$ and thus $v \in \Omega^1(X_T; \ad P)$ obeys
$$
\frac{\partial v}{\partial t} + \left(\nabla_{A_1}^*\nabla_{A_1} + \mu\right)v = \tilde f
\quad\hbox{on } X_T,
$$
where $\tilde f := e^{-\mu t}f$. For a quasilinear parabolic equation with nonlinearity, $\sK$,
$$
\frac{\partial u}{\partial t} + \nabla_{A_1}^*\nabla_{A_1}u + \sK(t,u) = 0
\quad\hbox{on } X_T,
$$
we obtain
$$
\frac{\partial v}{\partial t} + \left(\nabla_{A_1}^*\nabla_{A_1} + \mu\right)v
+ \tilde\sK(t,v) = 0 \quad\hbox{on } X_T,
$$
where $\tilde\sK(t,v) := e^{-\mu t}\sK(t,e^{\mu t}v)$. Hence, even in the case of quasilinear parabolic equation, the corresponding exponentially-shifted version is entirely equivalent for $t \in [0, T)$ when $T < \infty$. We end this digression and return to our quasilinear parabolic equation.

By applying the exponential shift trick to \eqref{eq:Quasilinear_parabolic_equation_a_minus_bara} by writing $a(t) =: e^{\mu t}b(t)$ and $\bar a(t) =: e^{\mu t}\bar b(t)$, for $\mu \geq \mu_0$ and $\mu_0 \geq 1$ to be determined in the sequel, we see that $b-\bar b$ obeys
\begin{equation}
\label{eq:Quasilinear_parabolic_equation_b_minus_barb}
\begin{aligned}
{}&\frac{\partial}{\partial t}(b - \bar b) + \left(\bar\nabla_{A_1}^*\nabla_{A_1} + \mu\right)(b - \bar b)
+ \Ric(\bar g) \times (b - \bar b) + F_{A_1} \times (b-\bar b)
\\
&\quad + e^{\mu t}\left[(b - \bar b) \times \nabla_{A_1}b +  b \times (b - \bar b)
+ b \times \nabla_{A_1}(b - \bar b)\right] + e^{2\mu t}b \times b \times (b - \bar b)
\\
&= -\left(\nabla_{A_1}^*\nabla_{A_1} - \bar\nabla_{A_1}^*\nabla_{A_1}\right)b
- (\Ric(g) - \Ric(\bar g)) \times b
\\
&\quad + (* - \bar*)F_{A_1} \times b + e^{\mu t}\left[b \times (\nabla_{A_1}^* - \bar\nabla_{A_1}^*)b
 + (* - \bar*)b \times \nabla_{A_1}b + (\nabla(* - \bar*))b \times b\right]
\\
&\quad + (* - \bar*)e^{2\mu t}b \times b \times b
\\
&\quad + e^{\mu t}\left[(b - \bar b) \times (b - \bar b) + (b - \bar b) \times \nabla_{A_1}(b - \bar b)\right]
+ e^{2\mu t}b \times (b - \bar b) \times (b - \bar b)
\\
&\quad + e^{2\mu t}(b - \bar b) \times (b - \bar b) \times (b - \bar b) \quad\hbox{a.e. on }(0, T_1)\times X.
\end{aligned}
\end{equation}
This concludes Step \ref{step:Exponentially_shifted_quasilinear_parabolic_equation}.
\end{step}

\begin{step}[\Apriori estimate for $b-\bar b$]
\label{step:Apriori_estimate_for_b_minus_barb}
From Theorem \ref{thm:Krylov_2-3-1_and_2-3-2_covariant_derivative_vector_bundle}, we have the \apriori estimate
\begin{multline}
\label{eq:Apriori_estimate_for_b_minus_barb_in_terms_of_L2_time_space_heat_operator_on_b_minus_barb}
\mu\|b - \bar b\|_{L^2(X_{T_1})} + \sqrt{\mu}\|\nabla_{A_1}(b - \bar b)\|_{L^2(X_{T_1})} + \|\nabla_{A_1}^2(b - \bar b)\|_{L^2(X_{T_1})} + \|\partial_t(b - \bar b)\|_{L^2(X_{T_1})}
\\
\quad \leq C_1\|(L_{A_1} + \mu)(b - \bar b)\|_{L^2(X_{T_1})},
\end{multline}
where $C_1 = C(A_1,g)$.

On the other hand, from the quasilinear parabolic equation \eqref{eq:Quasilinear_parabolic_equation_b_minus_barb} for $b - \bar b$ we see that\footnote{It is worth taking note of the occurrence of the term $\|e^{\mu t}\nabla_{A_1}b\|_{L^\infty(0, T; L^4(X))}$ on the right-hand side of \eqref{eq:Estimate_for_L2_time_space_heat_operator_on_b_minus_barb_raw}, as this motivates our hypothesis of the bound \eqref{eq:Yang-Mills_gradient_flow_A_minus_A1_Linfty_time_H2_space_lessthan_bound} for $\|a\|_{L^\infty(0, T; H_{A_1}^2(X))}$ and thus $\|\nabla_{A_1}a\|_{L^\infty(0, T; L^4(X))}$.}
\begin{equation}
\label{eq:Estimate_for_L2_time_space_heat_operator_on_b_minus_barb_raw}
\begin{aligned}
{} & \|(L_{A_1} + \mu)(b - \bar b)\|_{L^2(X_{T_1})}
\\
&\quad\leq c\left(\|\Ric(\bar g)\|_{C(X)} + \|F_{A_1}\|_{C(X)} \right)\|b - \bar b\|_{L^2(X_{T_1})}
\\
&\qquad + c\left(\|e^{\mu t}b\|_{L^\infty(0, T; L^4(X))} + c\|e^{\mu t}\nabla_{A_1}b\|_{L^\infty(0, T; L^4(X))} \right)
 \|b - \bar b\|_{L^2(0,{T_1}; L^4(X))}
\\
&\qquad + c\|e^{2\mu t}b\times b\|_{L^\infty(0, T; L^4(X))} \|b - \bar b\|_{L^2(0,T_1; L^4(X))}
\\
&\qquad + c\|e^{\mu t}b\|_{L^\infty(0, T; L^6(X))} \|\nabla_{A_1}(b - \bar b)\|_{L^2(0,{T_1}; L^3(X))}
\\
&\qquad + \|\left(\nabla_{A_1}^*\nabla_{A_1} - \bar\nabla_{A_1}^*\nabla_{A_1}\right)b\|_{L^2(X_T)}
+ c\|(* - \bar*)F_{A_1}\|_{C(X)}\|b\|_{L^2(X_T)}
\\
&\qquad + c\|\Ric(g) - \Ric(\bar g)\|_{L^4(X)} \|b\|_{L^2(0,T; L^4(X))}
\\
&\qquad + c\|e^{\mu t}b \|_{L^\infty(0, T; L^4(X))} \|(\nabla_{A_1}^* - \bar\nabla_{A_1}^*)b\|_{L^2(0,T; L^4(X))}
\\
&\qquad + c\|(* - \bar*)e^{\mu t}b \|_{L^\infty(0, T; L^4(X))} \|\nabla_{A_1}b\|_{L^2(0,T; L^4(X))}
\\
&\qquad + c\|(\nabla(* - \bar*))e^{\mu t}b \|_{L^\infty(0, T; L^4(X))} \|b\|_{L^2(0,T; L^4(X))}
\\
&\qquad + c\|(* - \bar*)e^{\mu t}b \|_{L^\infty(0, T; L^4(X))} \|e^{\mu t}b \|_{L^\infty(0, T; H_{A_1}^1(X))}
\|b \|_{L^2(0, T; H_{A_1}^2(X))}
\quad \hbox{(by \eqref{eq:Struwe_Yang-Mills_heat_equation_estimate_term_VIII_W14})}
\\
&\qquad + e^{\mu T}\|(b - \bar b) \times (b - \bar b)\|_{L^2(X_{T_1})}
+ e^{\mu T}\|(b - \bar b) \times \nabla_{A_1}(b - \bar b)\|_{L^2(X_{T_1})}
\\
&\qquad + ce^{2\mu T}\|b\|_{L^\infty(0, T; L^4(X))} \|(b - \bar b) \times (b - \bar b) \|_{L^2(0, {T_1}; L^4(X))}
\\
&\qquad + e^{2\mu T}\|(b - \bar b) \times (b - \bar b) \times (b - \bar b) \|_{L^2(X_{T_1})}.
\end{aligned}
\end{equation}
In deriving the preceding estimate, we used H\"older inequalities such as,
\begin{align*}
\|fh\|_{L^2(X_T)} &\leq \|f\|_{L^\infty(0, T; L^4(X))} \|h\|_{L^2(0,T; L^4(X))},
\\
\|fh\|_{L^2(X_T)} &\leq \|f\|_{L^\infty(0, T; L^6(X))} \|h\|_{L^2(0,T; L^3(X))}.
\end{align*}
Terms such as $\|e^{\mu t}b(t)\|_{L^6(X)}$ may be estimated using the Sobolev embeddings $W^{1, \frac{12}{5}}(X) \hookrightarrow L^6(X)$ (where $k=1$ and $p=12/5)$) and $H^1(X) \hookrightarrow L^4(X)$ \cite[Theorem 4.12]{AdamsFournier} and the Kato Inequality \eqref{eq:FU_6-20_first-order_Kato_inequality} to give,
$$
\|e^{\mu t}b(t)\|_{L^6(X)} \leq c\|e^{\mu t}b(t)\|_{W_{A_1}^{1, \frac{12}{5}}(X)} \leq c\|e^{\mu t}b(t)\|_{W_{A_1}^{1, 4}(X)}
\leq c\|e^{\mu t}b(t)\|_{H_{A_1}^2(X)}.
$$
For the $g - \bar g$ terms in the bound \eqref{eq:Estimate_for_L2_time_space_heat_operator_on_b_minus_barb_raw} for $\|(L_{A_1} + \mu)(b - \bar b)\|_{L^2(X_{T_1})}$, we see that
\begin{equation}
\label{eq:Norm_g_minus_barg_b_terms}
\begin{aligned}
\|\left(\nabla_{A_1}^*\nabla_{A_1} - \bar\nabla_{A_1}^*\nabla_{A_1}\right)b\|_{L^2(X_T)}
&\leq
c\|g-\bar g\|_{C(X)} \|\nabla_{A_1}^2b\|_{L^2(X_T)}
\\
&\quad + c\|g-\bar g\|_{C^1(X)} \|\nabla_{A_1}b\|_{L^2(X_T)},
\\
\|(* - \bar*)F_{A_1}\|_{C(X)}
&\leq
c\|g-\bar g\|_{C(X)} \|F_{A_1}\|_{C(X)},
\\
\|\Ric(g) - \Ric(\bar g)\|_{L^4(X)}
&\leq
c\|g-\bar g\|_{W^{2,4}(X)},
\\
\|(\nabla_{A_1}^* - \bar\nabla_{A_1}^*)b\|_{L^2(0,T; L^4(X))}
&\leq
c\|g-\bar g\|_{C(X)} \|\nabla_{A_1}b\|_{L^2(0,T; L^4(X))}
\\
&\quad + c\|g-\bar g\|_{C^1(X)} \|b\|_{L^2(0,T; L^4(X))},
\\
\|(* - \bar*)e^{\mu t}b\|_{L^\infty(0,T; L^4(X))}
&\leq
c\|g-\bar g\|_{C(X)} \|e^{\mu t}b\|_{L^\infty(0,T; L^4(X))},
\\
\|(\nabla(* - \bar*))e^{\mu t}b \|_{L^\infty(0,T; L^4(X))}
&\leq
c\|g-\bar g\|_{C^1(X)} \|e^{\mu t}b\|_{L^\infty(0,T; L^4(X))}.
\end{aligned}
\end{equation}
By applying the Sobolev embedding $H^1(X) \hookrightarrow L^4(X)$ \cite[Theorem 4.12]{AdamsFournier} and the Kato Inequality \eqref{eq:FU_6-20_first-order_Kato_inequality}, we also obtain
\begin{align*}
\|e^{2\mu t}b\times b\|_{L^\infty(0, T; L^4(X))}
&\leq
c\|e^{\mu t}b \times e^{\mu t}b\|_{L^\infty(0, T; H_{A_1}^1(X))}
\\
&\leq
c\|e^{\mu t}b \times e^{\mu t}b\|_{L^\infty(0, T; L^2(X))}
+ c\|\nabla_{A_1}(e^{\mu t}b \times e^{\mu t}b)\|_{L^\infty(0, T; L^2(X))}
\\
&\leq  c\|e^{\mu t}b\|_{L^\infty(0, T; L^4(X))}^2
+ c\|e^{\mu t}b\|_{L^\infty(0, T; L^4(X))} \|e^{\mu t}\nabla_{A_1}b\|_{L^\infty(0, T; L^4(X))}
\\
&\leq c\|e^{\mu t}b\|_{L^\infty(0, T; H_{A_1}^2(X))}^2.
\end{align*}
We estimate the quadratic and cubic powers of $b - \bar b$ in the right-hand side of the bound \eqref{eq:Estimate_for_L2_time_space_heat_operator_on_b_minus_barb_raw} for $\|(L_{A_1} + \mu)(b - \bar b)\|_{L^2(X_{T_1})}$,
\begin{align}
\label{eq:L2_time_and_space_norm_quadratic_b_minus_barb}
\|(b - \bar b) \times (b - \bar b)\|_{L^2(X_{T_1})}
&\leq
c\|b - \bar b\|_{L^\infty(0,T_1;L^4(X))} \|b - \bar b\|_{L^2(0,T_1;L^4(X))}
\\
\notag
&\leq c\|b - \bar b\|_{L^\infty(0,T_1;H_{A_1}^1(X))} \|b - \bar b\|_{L^2(0,T_1;H_{A_1}^1(X))},
\\
\label{eq:L2_time_and_space_norm_b_minus_barb_times_nablaA1_b_minus_barb}
\|(b - \bar b) \times \nabla_{A_1}(b - \bar b)\|_{L^2(X_{T_1})}
&\leq
c\|b - \bar b\|_{L^\infty(0,T_1;L^4(X))} \|\nabla_{A_1}(b - \bar b)\|_{L^2(0,T_1;L^4(X))}
\\
\notag
&\leq c\|b - \bar b\|_{L^\infty(0,T_1;H_{A_1}^1(X))} \|b - \bar b\|_{L^2(0,T_1;H_{A_1}^2(X))},
\\
\label{eq:L2_time_and_L4_space_norm_quadratic_b_minus_barb}
\|(b - \bar b) \times (b - \bar b)\|_{L^2(0,T_1;L^4(X))}
&\leq
c\|b - \bar b\|_{L^\infty(0,T_1;H_{A_1}^1(X))} \|b - \bar b\|_{L^2(0,T_1;H_{A_1}^2(X))},
\end{align}
where \eqref{eq:L2_time_and_L4_space_norm_quadratic_b_minus_barb} follows
from \eqref{eq:Struwe_Yang-Mills_heat_equation_estimate_term_VII_H2}, and
\begin{equation}
\label{eq:L2_time_and_space_norm_cubic_b_minus_barb}
\begin{aligned}
{}& \|(b - \bar b) \times (b - \bar b) \times (b - \bar b) \|_{L^2(X_{T_1})}
\\
&\quad \leq c\|b - \bar b\|_{L^\infty(0,T_1;L^4(X))}\|(b - \bar b) \times (b - \bar b)\|_{L^2(0,T_1;L^4(X))}
\\
&\quad \leq c\|b - \bar b\|_{L^\infty(0,T_1;H_{A_1}^1(X))}^2 \|b - \bar b\|_{L^2(0,T_1;H_{A_1}^2(X))}.
\end{aligned}
\end{equation}
Note that $(b - \bar b)(0) = 0$, by definition of the solutions $A(t)$ and $\bar A(t)$, and so the estimate \eqref{eq:Struwe_9} yields
$$
\|b - \bar b\|_{L^\infty(0,T_1;H_{A_1}^1(X))}
\leq \sqrt{2}\left( \|b - \bar b\|_{L^2(0,T_1;H_{A_1}^2(X))} + \|\partial_t(b - \bar b)\|_{L^2(0,T_1;L^2(X))} \right).
$$
We now substitute the preceding inequalities together with the hypotheses \eqref{eq:Riem_barg_C0_norm_lessthan_bound}, \eqref{eq:Yang-Mills_gradient_flow_A_minus_A1_Linfty_time_H2_space_lessthan_bound}, \eqref{eq:F_A1_C0_norm_lessthan_bound}, the bounds \eqref{eq:Norm_g_minus_barg_b_terms},  \eqref{eq:L2_time_and_space_norm_quadratic_b_minus_barb},
\eqref{eq:L2_time_and_space_norm_b_minus_barb_times_nablaA1_b_minus_barb}
\eqref{eq:L2_time_and_L4_space_norm_quadratic_b_minus_barb}, \eqref{eq:L2_time_and_space_norm_cubic_b_minus_barb},
and the Sobolev embedding $H^1(X) \hookrightarrow L^4(X)$ \cite[Theorem 4.12]{AdamsFournier} and Kato Inequality \eqref{eq:FU_6-20_first-order_Kato_inequality} into our bound \eqref{eq:Estimate_for_L2_time_space_heat_operator_on_b_minus_barb_raw} for $\|(L_{A_1} + \mu)(b - \bar b)\|_{L^2(X_{T_1})}$ to give
\begin{equation}
\label{eq:Estimate_for_L2_time_space_heat_operator_on_b_minus_barb_before_interpolation_and_rearrangement}
\begin{aligned}
{}& \|(L_{A_1} + \mu)(b - \bar b)\|_{L^2(X_{T_1})}
\\
&\quad\leq C_2\left(\|b - \bar b\|_{L^2(X_{T_1})} + \|\nabla_{A_1}(b - \bar b)\|_{L^2(X_{T_1})}\right)
+ C_2\|\nabla_{A_1}(b - \bar b)\|_{L^2(0,T_1; L^3(X))}
\\
&\qquad + C_2\|g-\bar g\|_{C^1(X)} + c\|g-\bar g\|_{W^{2,4}(X)}
\\
&\qquad + ce^{\mu T}\|b - \bar b\|_{L^\infty(0,T_1;H_{A_1}^1(X))} \|b - \bar b\|_{L^2(0,T_1;H_{A_1}^1(X))}
\\
&\qquad + ce^{\mu T} \|b - \bar b\|_{L^\infty(0,T_1;H_{A_1}^1(X))} \|b - \bar b\|_{L^2(0,T_1;H_{A_1}^2(X))}
\\
&\qquad + C_2 e^{2\mu T} \|b - \bar b\|_{L^\infty(0,T_1;H_{A_1}^1(X))} \|b - \bar b\|_{L^2(0,T_1;H_{A_1}^2(X))}
\\
&\qquad + ce^{2\mu T}\|b - \bar b\|_{L^\infty(0,T_1;H_{A_1}^1(X))}^2 \|b - \bar b\|_{L^2(0,T_1;H_{A_1}^2(X))},
\end{aligned}
\end{equation}
where $c = c(g)$ and $C_2 = C_2(K, g)$ are positive constants. The term $\|\nabla_{A_1}(b - \bar b)\|_{L^2(0,T_1; L^3(X))}$ may be estimated with the aid of the interpolation inequality \cite[Equation (7.10)]{GilbargTrudinger}, for any $\eps \in (0, 1]$,
\begin{align*}
\|\nabla_{A_1}(b - \bar b)\|_{L^2(0,T_1; L^3(X))}
&\leq
\eps\|\nabla_{A_1}(b - \bar b)\|_{L^2(0,T_1; L^4(X))} + \eps^{-2}\|\nabla_{A_1}(b - \bar b)\|_{L^2(0,T_1; L^2(X))}
\\
&\leq c\eps\|\nabla_{A_1}(b - \bar b)\|_{L^2(0,T_1; H_{A_1}^1(X))} + \eps^{-2}\|\nabla_{A_1}(b - \bar b)\|_{L^2(0,T_1; L^2(X))}
\\
&\leq c\eps\|\nabla_{A_1}^2(b - \bar b)\|_{L^2(0,T_1; L^2(X))}
+ (\eps^{-2} + c\eps)\|\nabla_{A_1}(b - \bar b)\|_{L^2(0,T_1; L^2(X))},
\end{align*}
where the second inequality follows from the Sobolev embedding $H^1(X) \hookrightarrow L^4(X)$ \cite[Theorem 4.12]{AdamsFournier} and the Kato Inequality \eqref{eq:FU_6-20_first-order_Kato_inequality}. We now choose $\eps = \eps(C_1,C_2,g) = \eps(C_1,K,g) \in (0, 1]$ small enough that $cC_2\eps \leq C_1/2$, where we use $C_1$ to label the constant appearing on the right-hand side of the \apriori estimate \eqref{eq:Apriori_estimate_for_b_minus_barb_in_terms_of_L2_time_space_heat_operator_on_b_minus_barb} for $b-\bar b$. We can thus combine our estimate \eqref{eq:Estimate_for_L2_time_space_heat_operator_on_b_minus_barb_before_interpolation_and_rearrangement} for $\|(L_{A_1} + \mu)(b - \bar b)\|_{L^2(X_{T_1})}$ with the \apriori estimate \eqref{eq:Apriori_estimate_for_b_minus_barb_in_terms_of_L2_time_space_heat_operator_on_b_minus_barb} and absorb the resulting right-hand side term, $cC_2\eps\|\nabla_{A_1}^2(b - \bar b)\|_{L^2(0,T_1; L^2(X))}$, into the left-hand side of the inequality to give
\begin{align*}
{} & \mu\|b - \bar b\|_{L^2(X_{T_1})} + \sqrt{\mu}\|\nabla_{A_1}(b - \bar b)\|_{L^2(X_{T_1})} + \|\nabla_{A_1}^2(b - \bar b)\|_{L^2(X_{T_1})} + \|\partial_t(b - \bar b)\|_{L^2(X_{T_1})}
\\
&\quad\leq 2C_2\|b - \bar b\|_{L^2(X_{T_1})} + C_3\|\nabla_{A_1}(b - \bar b)\|_{L^2(X_{T_1})}
\\
&\qquad + 2C_2\|g-\bar g\|_{C^1(X)} + 2c\|g-\bar g\|_{W^{2,4}(X)}
\\
&\qquad + 4ce^{\mu T} \|b - \bar b\|_{L^\infty(0,T_1;H_{A_1}^1(X))} \|b - \bar b\|_{L^2(0,T_1;H_{A_1}^2(X))}
\\
&\qquad + 2C_2 e^{2\mu T} \|b - \bar b\|_{L^\infty(0,T_1;H_{A_1}^1(X))} \|b - \bar b\|_{L^2(0,T_1;H_{A_1}^2(X))}
\\
&\qquad + 2ce^{2\mu T}\|b - \bar b\|_{L^\infty(0,T_1;H_{A_1}^1(X))}^2 \|b - \bar b\|_{L^2(0,T_1;H_{A_1}^2(X))},
\end{align*}
where $C_3 = C_3(C_1,C_2,g) = C_3(C_1,K,g)$. We therefore choose $\mu_0 = \mu(C_2,C_3,g) = \mu(C_1,K,g) \geq 1$ large enough that
$$
2C_2 \leq \frac{\mu_0}{2} \quad\hbox{and}\quad C_3 \leq \frac{\sqrt{\mu_0}}{2},
$$
and obtain
\begin{align*}
{} & \mu\|b - \bar b\|_{L^2(X_{T_1})} + \sqrt{\mu}\|\nabla_{A_1}(b - \bar b)\|_{L^2(X_{T_1})} + \|\nabla_{A_1}^2(b - \bar b)\|_{L^2(X_{T_1})} + \|\partial_t(b - \bar b)\|_{L^2(X_{T_1})}
\\
&\quad\leq 4C_2\|g-\bar g\|_{C^1(X)} + 4c\|g-\bar g\|_{W^{2,4}(X)}
\\
&\qquad + 8ce^{\mu T} \|b - \bar b\|_{L^\infty(0,T_1;H_{A_1}^1(X))} \|b - \bar b\|_{L^2(0,T_1;H_{A_1}^2(X))}
\\
&\qquad + 4C_2 e^{2\mu T} \|b - \bar b\|_{L^\infty(0,T_1;H_{A_1}^1(X))} \|b - \bar b\|_{L^2(0,T_1;H_{A_1}^2(X))}
\\
&\qquad + 4c^2e^{2\mu T}\|b - \bar b\|_{L^\infty(0,T_1;H_{A_1}^1(X))}^2 \|b - \bar b\|_{L^2(0,T_1;H_{A_1}^2(X))},
\quad \forall\, \mu \geq \mu_0,
\end{align*}
where we may allow that $c\geq 1$ without loss of generality when writing the final term in the preceding inequality. We now fix $\mu = \mu_0$ for the remainder of the proof. If we further suppose that $\|b - \bar b\|_{L^\infty(0,T_1;H_{A_1}^1(X))}$ is small enough that
$$
\left(8ce^{\mu_0 T} + 4C_2e^{2\mu_0 T}\right)\|b - \bar b\|_{L^\infty(0,T_1;H_{A_1}^1(X))} \leq 1/4
$$
that is,
\begin{equation}
\label{eq:Linfty_time_H1_space_b_minus_barb_bound_allowing_rearrangement_quadratic_and_cubic_terms}
\|b - \bar b\|_{L^\infty(0,T_1;H_{A_1}^1(X))} \leq \sigma := \frac{1}{16\left(2ce^{\mu_0 T} + C_2e^{2\mu_0 T}\right)},
\end{equation}
so $\sigma = \sigma(C_2,g,\mu_0,T) = \sigma(C_1,g,K,T)$ (recall that $C_1 = C_1(A_1,g)$ was defined in Step \ref{step:Apriori_estimate_for_b_minus_barb}), then we discover the inequality,
\begin{multline*}
\mu_0\|b - \bar b\|_{L^2(X_{T_1})} + \sqrt{\mu_0}\|\nabla_{A_1}(b - \bar b)\|_{L^2(X_{T_1})} + \|\nabla_{A_1}^2(b - \bar b)\|_{L^2(X_{T_1})} + \|\partial_t(b - \bar b)\|_{L^2(X_{T_1})}
\\
\leq 4C_2\|g-\bar g\|_{C^1(X)} + 4c\|g-\bar g\|_{W^{2,4}(X)} + \frac{9}{16}\|b - \bar b\|_{L^2(0,T_1;H_{A_1}^2(X))}.
\end{multline*}
Hence, noting that $\mu_0 \geq 1$, we find that
\begin{multline*}
\|b - \bar b\|_{L^2(X_{T_1})} + \|\nabla_{A_1}(b - \bar b)\|_{L^2(X_{T_1})}
+ \|\nabla_{A_1}^2 (b - \bar b)\|_{L^2(X_{T_1})} + \|\partial_t (b - \bar b)\|_{L^2(X_{T_1})}
\\
+ \|b - \bar b\|_{L^\infty(0,T_1;H_{A_1}^1(X))} \leq 16C_2\|g-\bar g\|_{C^1(X)} + 16c\|g-\bar g\|_{W^{2,4}(X)}.
\end{multline*}
that is, after relabeling the constant $16\max\{c, C_2\}$ as $C_2$,
\begin{multline}
\label{eq:Apriori_estimate_for_b_minus_barb_in_terms_C1_and_W24_norms_g_minus_barg}
\|b - \bar b\|_{L^2(0,T_1;H_{A_1}^2(X))} + \|\partial_t (b - \bar b)\|_{L^2(0,T_1;L^2(X))}
+ \|b - \bar b\|_{L^\infty(0,T_1;H_{A_1}^1(X))}
\\
\leq C_2\left(\|g-\bar g\|_{C^1(X)} + \|g-\bar g\|_{W^{2,4}(X)}\right).
\end{multline}
This concludes Step \ref{step:Apriori_estimate_for_b_minus_barb}.
\end{step}

\begin{step}[Norm criterion for continuous temporal extension of $\bar A(t)$]
\label{step:Norm_criterion_for_no_bubbling_in_barA}
We next establish the following key

\begin{claim}[Norm criterion for continuous temporal extension of $\bar A(t)$]
\label{claim:Norm_criterion_for_no_bubbling_in_barA}
There is a positive constant, $\zeta = \zeta(A_1,g,K,T) \in (0, 1]$ with the following significance. If $T_2 \in (0, T]$ has the property that $\bar A - A_1 \in C_{\loc}([0,T_2); H_{A_1}^1(X;\Lambda^1\otimes\ad P))$ and
\begin{equation}
\label{eq:Linfty_time_H1_space_b_minus_barb_leq_small_constant_ensuring_no_bubbling}
\|b(t) - \bar b(t)\|_{H_{A_1}^1(X)} <  \zeta, \quad \forall\, t \in [0,T_2),
\end{equation}
then $\bar A - A_1 \in C([0,T_2]; H_{A_1}^1(X;\Lambda^1\otimes\ad P))$.
\end{claim}

\begin{proof}[Proof of Claim \ref{claim:Norm_criterion_for_no_bubbling_in_barA}]
Observe that, for any positive constant, $r \leq \Inj(X,g)$, and point $x \in X$,
\begin{align*}
\|F_{\bar A}(t)\|_{L^2(B_r(x))} &\leq \|F_A(t)\|_{L^2(B_r(x))} + \|F_A(t)-F_{\bar A}(t)\|_{L^2(B_r(x))}
\\
&\leq \|F_A(t)\|_{L^2(B_r(x))} + \|F_A(t)-F_{\bar A}(t)\|_{L^2(X)}, \quad\forall\, t \in [0, T_2),
\end{align*}
where $B_r(x_l) := \{x\in X:\dist_g(x,x_l) < r\}$, and thus,
\begin{equation}
\label{eq:Linfty_time_L2_ball_FbarA_leq_Linfty_time_L2_ball_FA_plus_Linfty_time_L2_space_FA_minus_FbarA}
\sup_{x\in X}\|F_{\bar A}(t)\|_{L^2(B_r(x))}
\leq
\sup_{x\in X}\|F_A\|_{L^2(B_r(x))} +  \|F_A - F_{\bar A}\|_{L^2(X)}, \quad\forall\, t \in [0, T_2).
\end{equation}
Writing $F_A(t) = F_{A_1} + d_{A_1}a(t) + [a(t), a(t)]$ and $F_{\bar A}(t) = F_{A_1} + d_{A_1}\bar a(t) + [\bar a(t), \bar a(t)]$, we obtain
\begin{align*}
F_A(t) - F_{\bar A}(t) &= d_{A_1}(a(t) - \bar a(t)) + [a(t) - \bar a(t), a(t)] + [\bar a(t), a(t) - \bar a(t)]
\\
&= d_{A_1}(a(t) - \bar a(t)) + [a(t) - \bar a(t), a(t)] + [a(t), a(t) - \bar a(t)]
\\
&\quad + [\bar a(t) - a(t), a(t) - \bar a(t)],
\end{align*}
and hence,
\begin{align*}
\|F_A(t)-F_{\bar A}(t)\|_{L^2(X)} &\leq c\|\nabla_{A_1}(a(t) - \bar a(t))\|_{L^2(X)}
+ c\|a(t)\|_{L^4(X)} \|a(t) - \bar a(t)\|_{L^4(X)}
\\
&\quad + c\|a(t) - \bar a(t)\|_{L^4(X)}^2, \quad\forall\, t\in [0, T_2),
\end{align*}
where $c$ is a positive constant depending at most on $g$. Therefore, recalling that $a(t) = e^{\mu_0 t}b(t)$ and $\bar a(t) = e^{\mu_0 t}\bar b(t)$,
\begin{align*}
\|F_A(t) - F_{\bar A}(t)\|_{L^2(X)} &\leq e^{\mu_0 T}\|\nabla_{A_1}(b - \bar b)(t)\|_{L^2(X)}
\\
&\quad + ce^{\mu_0 T}\|e^{\mu_0 t}b(t)\|_{L^4(X)} \|(b - \bar b)(t)\|_{L^4(X)}
\\
&\quad + ce^{2\mu_0 T}\|(b - \bar b)(t)\|_{L^4(X)}^2, \quad\forall\, t\in [0, T_2).
\end{align*}
The Sobolev embedding, $H^1(X) \hookrightarrow L^4(X)$ \cite[Theorem 4.12]{AdamsFournier}, and the Kato Inequality \eqref{eq:FU_6-20_first-order_Kato_inequality} yield,
\begin{multline}
\label{eq:Linfty_time_L2_space_FA_minus_FbarA_raw}
\|F_A(t) - F_{\bar A}(t)\|_{L^2(X)}
\\
\leq
c_1(1 + K)e^{\mu_0 T}\|(b - \bar b)(t)\|_{H_{A_1}^1(X)} + c_1e^{2\mu_0 T}\|(b - \bar b)(t)\|_{H_{A_1}^1(X)}^2,
\quad\forall\, t\in [0, T_2),
\end{multline}
where $c_1$ is a positive constant depending at most on $g$. Provided
$$
\|b - \bar b\|_{L^\infty(0,T_2;H_{A_1}^1(X))} \leq 1,
$$
which will be assured by our hypothesis \eqref{eq:Linfty_time_H1_space_b_minus_barb_leq_small_constant_ensuring_no_bubbling} (since $\zeta \in (0, 1]$), we thus obtain,
$$
\|F_A(t) - F_{\bar A}(t)\|_{L^2(X)}
\leq
c_1\left(1 + K + e^{\mu_0 T}\right) e^{\mu_0 T} \|b - \bar b\|_{L^\infty(0,T_2;H_{A_1}^1(X))},
\quad\forall\, t\in [0, T_2).
$$
Let $\bar\eps_1 = \bar\eps_1(\bar g) \in (0, 1]$ and $\bar r_1 = \bar r_1(\bar g) \in (0, 1]$ be the positive constants in Equation \eqref{eq:Struwe_15} and Lemma \ref{lem:Schlatter_2-4_and_Struwe_3-6_local} and which we shall apply to the flow, $\bar A(t)$ for the metric, $\bar g$, and initial data, $A_0$. By Remark \ref{rmk:Schlatter_2-4_and_Struwe_3-6_local_equivalent_Riemannian_metric} and the fact that $g$ and $\bar g$ are related by \eqref{eq:Riemannian_metrics_g_minus_barg_C1_plus_W24_norm_lessthan_small_positive_constant} and \eqref{eq:Riemannian_metrics_g_minus_barg_C2_lessthan_bound}, the dependence of $\bar\eps_1$ and $\bar r_1$ on $\bar g$ have the following equivalents,
\begin{equation}
\label{eq:Struwe_lemma_3-6_eps1_and_r1_dependence_on_metrics_C1_equivalent_to_g}
\bar\eps_1(\bar g) = \bar\eps_1(g,K) \quad\hbox{and}\quad \bar r_1(\bar g) = \bar r_1(g,K).
\end{equation}
Recalling that $\mu_0 = \mu_0(A_1,g,K)$, we define $\zeta \in (0, 1]$ by
\begin{equation}
\zeta \equiv \zeta(A_1,g,K,T) := \frac{\sqrt{\bar\eps_1}}{4c_1\left(1 + K + e^{\mu_0 T}\right)e^{\mu_0 T}}.
\end{equation}
Therefore, provided $b-\bar b$ obeys the inequality \eqref{eq:Linfty_time_H1_space_b_minus_barb_leq_small_constant_ensuring_no_bubbling}, we discover that
\begin{equation}
\label{eq:Linfty_time_L2_space_FA_minus_FbarA}
\|F_A(t) - F_{\bar A}(t)\|_{L^2(X)} < \frac{\sqrt{\bar\eps_1}}{4}, \quad\forall\, t\in [0, T_2).
\end{equation}
Next observe that, for any $x \in X$ and $r \in (0, 1 \vee \Inj(X,g)]$, and writing $F_A(t) = F_{A_1} + d_{A_1}a(t) + [a(t), a(t)]$ over $X$,
\begin{align*}
\|F_A(t)\|_{L^2(B_r(x))} &\leq \|F_{A_1}\|_{L^2(B_r(x))} + \|d_{A_1}a(t)\|_{L^2(B_r(x))} + c\|[a(t), a(t)]\|_{L^2(B_r(x))}
\\
&\leq \left(\Vol_g(B_r(x))\right)^{1/2}\|F_{A_1}\|_{C(X)} + \left(\Vol_g(B_r(x))\right)^{1/4}\|d_{A_1}a(t)\|_{L^4(X)}
\\
&\quad + c\left(\Vol_g(B_r(x))\right)^{1/4} \|a(t)\|_{L^4(X)}^2
\\
&\leq c(Kr^2 + Kr + K^2r) \quad\hbox{(by \eqref{eq:Yang-Mills_gradient_flow_A_minus_A1_Linfty_time_H2_space_lessthan_bound} and \eqref{eq:F_A1_C0_norm_lessthan_bound})},
\end{align*}
where $c$ is a positive constant depending at most on the Riemannian metric, $g$, on $X$ and we apply the Sobolev embedding $H^1(X) \hookrightarrow L^4(X)$ \cite[Theorem 4.12]{AdamsFournier} and the Kato Inequality \eqref{eq:FU_6-20_first-order_Kato_inequality} in the usual way. Thus (for a possibly larger $c$),
$$
\sup_{x\in X}\|F_A(t)\|_{L^2(B_r(x))} \leq c(1 + K)Kr, \quad 0 < r \leq 1 \vee \Inj(X,g).
$$
Therefore, provided $r_0 = r_0(g,K,\bar\eps_1) = r_0(g,K) \in (0, 1 \vee \Inj(X,g)]$ is small enough that
$$
c(1 + K)Kr_0 < \frac{\sqrt{\bar\eps_1}}{4},
$$
where we recall that $\bar\eps_1 = \bar\eps_1(g,K) \in (0, 1]$, then
$$
\sup_{x\in X}\|F_A(t)\|_{L^2(B_r(x))} < \frac{\sqrt{\bar\eps_1}}{4}, \quad\forall\, t\in [0, T) \hbox{ and }
0 < r \leq r_0.
$$
Combining the preceding inequality with the bounds  \eqref{eq:Linfty_time_L2_ball_FbarA_leq_Linfty_time_L2_ball_FA_plus_Linfty_time_L2_space_FA_minus_FbarA} and \eqref{eq:Linfty_time_L2_space_FA_minus_FbarA}, yields
$$
\sup_{x\in X}\|F_{\bar A}\|_{L^2(B_r(x))}
<
\frac{\sqrt{\bar\eps_1}}{2}, \quad\forall\, t\in [0, T_2) \hbox{ and } 0 < r_0\vee \bar r_1.
$$
But Lemma \ref{lem:Schlatter_2-4_and_Struwe_3-6} now implies that $\bar A$ extends continuously from $[0, T_2)$ to $[0, T_2]$, that is, $\bar A - A_1 \in C([0,T_2]; H_{A_1}^1(X;\Lambda^1\otimes\ad P))$. This completes the proof of Claim \ref{claim:Norm_criterion_for_no_bubbling_in_barA}.
\end{proof}

This concludes Step \ref{step:Norm_criterion_for_no_bubbling_in_barA}.
\end{step}

\begin{step}[Conclusion that the maximal lifetime $\bar\sT$ of $\bar A$ is greater than $T$]
\label{step:Maximal_lifetime_barA_geq_T}
We claim that
\begin{equation}
\label{eq:H1_space_norm_b_minus_barb_strictly_lessthan_nobubbling_and_rearrangement_constant}
\|(b - \bar b)(t)\|_{H_{A_1}^1(X)} < \min\{\zeta, \sigma\}, \quad\forall\, t \in [0, T),
\end{equation}
and thus $\bar A - A_1 \in C([0,T]; H_{A_1}^1(X;\Lambda^1\otimes\ad P))$ and $\bar\sT > T$, where $\zeta$ is the positive constant in Claim \ref{claim:Norm_criterion_for_no_bubbling_in_barA} and $\sigma$ is the positive constant in \eqref{eq:Linfty_time_H1_space_b_minus_barb_bound_allowing_rearrangement_quadratic_and_cubic_terms}.

Suppose the claim is false. Then there must be a large enough time, $T_1 \in (0, T)$, such that
\begin{subequations}
\label{eq:maximal_T1_such_that_H1_space_norm_minus_barb_lessthan_nobubbling_and_rearrangement_constant}
\begin{align}
\label{eq:maximal_T1_such_that_H1_space_norm_minus_barb_strictly_lessthan_nobubbling_and_rearrangement_constant}
\|(b - \bar b)(t)\|_{H_{A_1}^1(X)} &< \min\{\zeta, \sigma\}, \quad\forall\, t \in [0, T_1), \quad\hbox{but}
\\
\label{eq:maximal_T1_such_that_H1_space_norm_minus_barb_geq_nobubbling_and_rearrangement_constant}
\|(b - \bar b)(T_1)\|_{H_{A_1}^1(X)} &\geq \min\{\zeta, \sigma\}.
\end{align}
\end{subequations}
For if not, then the strict inequality \eqref{eq:H1_space_norm_b_minus_barb_strictly_lessthan_nobubbling_and_rearrangement_constant} would hold and Claim \ref{claim:Norm_criterion_for_no_bubbling_in_barA} (with $T_2 = T$) would imply that $\bar A - A_1 \in C([0,T]; H_{A_1}^1(X;\Lambda^1\otimes\ad P))$ and $\bar\sT > T$.

We now choose the positive constant, $\eta$, in the hypotheses of Theorem \ref{thm:Global_apriori_estimate_difference_solutions_Yang-Mills_heat_equations_pair_metrics} by setting
\begin{equation}
\label{eq:Definition_smallness_constant_g_minus_barg_metric_difference}
\eta \equiv \eta(A_1,g,K,T) \leq \frac{\min\{\zeta, \sigma\}}{4C_2},
\end{equation}
where $C_2$ is the constant on the right-hand side of the \apriori estimate \eqref{eq:Apriori_estimate_for_b_minus_barb_in_terms_C1_and_W24_norms_g_minus_barg}, and so
$$
C_2\left(\|g-\bar g\|_{C^1(X)} + \|g-\bar g\|_{W^{2,4}(X)}\right)
\leq
C_2\eta \leq \frac{1}{4}\min\{\zeta, \sigma\}.
$$
Since the condition \eqref{eq:Linfty_time_H1_space_b_minus_barb_bound_allowing_rearrangement_quadratic_and_cubic_terms} is satisfied, the \apriori estimate \eqref{eq:Apriori_estimate_for_b_minus_barb_in_terms_C1_and_W24_norms_g_minus_barg} thus yields
$$
\|(b - \bar b)(t)\|_{H_{A_1}^1(X)} \leq \frac{1}{4}\min\{\zeta, \sigma\}, \quad\forall\, t \in [0, T_1),
$$
and ensures that \eqref{eq:Linfty_time_H1_space_b_minus_barb_leq_small_constant_ensuring_no_bubbling} holds with $T_2=T_1$. Consequently, Claim \ref{claim:Norm_criterion_for_no_bubbling_in_barA} (with $T_2 = T_1$) implies that $\bar A - A_1 \in C([0,T_1]; H_{A_1}^1(X;\Lambda^1\otimes\ad P))$ and
$$
\|(b - \bar b)(t)\|_{H_{A_1}^1(X)} \leq \frac{1}{4}\min\{\zeta, \sigma\}, \quad\forall\, t \in [0, T_1],
$$
which contradicts \eqref{eq:maximal_T1_such_that_H1_space_norm_minus_barb_geq_nobubbling_and_rearrangement_constant}. Hence, the strict inequality \eqref{eq:H1_space_norm_b_minus_barb_strictly_lessthan_nobubbling_and_rearrangement_constant} holds and consequently $\bar A - A_1 \in C([0,T]; H_{A_1}^1(X;\Lambda^1\otimes\ad P))$ and $\bar\sT > T$. This concludes Step \ref{step:Maximal_lifetime_barA_geq_T}.
\end{step}


By virtue of Step \ref{step:Maximal_lifetime_barA_geq_T}, the condition \eqref{eq:Linfty_time_H1_space_b_minus_barb_bound_allowing_rearrangement_quadratic_and_cubic_terms} and \apriori estimate \eqref{eq:Apriori_estimate_for_b_minus_barb_in_terms_C1_and_W24_norms_g_minus_barg} hold with $T_1 = T$. The \apriori estimate \eqref{eq:Apriori_estimate_for_A_minus_barA_in_terms_C1_and_W24_norms_g_minus_barg} now follows from \eqref{eq:Apriori_estimate_for_b_minus_barb_in_terms_C1_and_W24_norms_g_minus_barg}, after setting $T_1 = T$, relabeling the constant $C_2$, and recalling that $b(t) - \bar b(t) = e^{-\mu_0 t}(A(t) - \bar A(t))$. This completes the proof of Theorem \ref{thm:Global_apriori_estimate_difference_solutions_Yang-Mills_heat_equations_pair_metrics}.
\end{proof}

When we introduce the cut-off functions, $\chi_l$, we can choose them in such a way that $d\chi_l$ is supported in the interior of where $\bar g$ is flat and we know that both $A(t)$ and $\bar A(t)$ do not bubble there and so have good control over the connection one-forms.

The hypotheses \eqref{eq:Riem_barg_C0_norm_lessthan_bound}, \eqref{eq:Riemannian_metrics_g_minus_barg_C1_plus_W24_norm_lessthan_small_positive_constant}, and \eqref{eq:Riemannian_metrics_g_minus_barg_C2_lessthan_bound} on $\bar g$
in Theorem \ref{thm:Global_apriori_estimate_difference_solutions_Yang-Mills_heat_equations_pair_metrics} can be relaxed and this is useful in situations where we know that $\bar g$ can be chosen $C^0$ close as desired to $g$, but not necessarily $C^1$-close as desired to $g$, while $\Riem(\bar g)$ remains uniformly $L^4$ bounded but not necessarily uniformly $C^0$ bounded.

\begin{cor}[Global \apriori estimate for the difference between solutions to the Yang-Mills heat equations defined by two weakly close Riemannian metrics and common initial data]
\label{cor:Global_apriori_estimate_difference_solutions_YM_heat_eqns_pair_metrics_C0_close}
Assume the hypotheses of Theorem \ref{thm:Global_apriori_estimate_difference_solutions_Yang-Mills_heat_equations_pair_metrics}, except replace the hypotheses \eqref{eq:Riem_barg_C0_norm_lessthan_bound}, \eqref{eq:Riemannian_metrics_g_minus_barg_C1_plus_W24_norm_lessthan_small_positive_constant}, and \eqref{eq:Riemannian_metrics_g_minus_barg_C2_lessthan_bound} on $\bar g$ by the following:
\begin{align}
\label{eq:Riem_barg_L4_norm_lessthan_bound}
\|\Riem(\bar g)\|_{L^4(X)} &\leq K,
\\
\label{eq:Riemannian_metrics_g_minus_barg_C0_plus_W24_norm_lessthan_small_positive_constant}
\|g - \bar g\|_{W^{2,4}(X)} &\leq \eta.
\end{align}
Then the conclusions of Theorem \ref{thm:Global_apriori_estimate_difference_solutions_Yang-Mills_heat_equations_pair_metrics} continue to hold.
\end{cor}

\begin{proof}
We only need to slightly modify the proof of Theorem \ref{thm:Global_apriori_estimate_difference_solutions_Yang-Mills_heat_equations_pair_metrics} by employing different combinations of bounds for the $g-\bar g$ terms in the collection of estimates \eqref{eq:Norm_g_minus_barg_b_terms} arising in Step \ref{step:Apriori_estimate_for_b_minus_barb}. We indicate below the terms that require changes:
\begin{equation}
\label{eq:Norm_g_minus_barg_b_terms_C0-close_metrics}
\begin{aligned}
\|\left(\nabla_{A_1}^*\nabla_{A_1} - \bar\nabla_{A_1}^*\nabla_{A_1}\right)b\|_{L^2(X_T)}
&\leq
c\|g-\bar g\|_{C(X)} \|\nabla_{A_1}^2b\|_{L^2(X_T)}
\\
&\quad + c\|g-\bar g\|_{W^{1,4}(X)} \|\nabla_{A_1}b\|_{L^2(0,T;L^4(X))},
\\
\|(\nabla_{A_1}^* - \bar\nabla_{A_1}^*)b\|_{L^2(0,T; L^4(X))}
&\leq
c\|g-\bar g\|_{C(X)} \|\nabla_{A_1}b\|_{L^2(0,T; L^4(X))}
\\
&\quad + c\|g-\bar g\|_{W^{1,8}(X)} \|b\|_{L^2(0,T; L^8(X))},
\\
\|(\nabla(* - \bar*))e^{\mu t}b \|_{L^\infty(0,T; L^4(X))}
&\leq
c\|g-\bar g\|_{W^{1,8}(X)} \|e^{\mu t}b\|_{L^\infty(0,T; L^8(X))}.
\end{aligned}
\end{equation}
Noting that $a(t) = e^{\mu t}b(t)$ (see Step \ref{step:Exponentially_shifted_quasilinear_parabolic_equation} in the proof of Theorem \ref{thm:Global_apriori_estimate_difference_solutions_Yang-Mills_heat_equations_pair_metrics}) and $W^{1,8/3}(X) \hookrightarrow L^8(X)$ by \cite[Theorem 4.12]{AdamsFournier} and applying the Kato Inequality \eqref{eq:FU_6-20_first-order_Kato_inequality}, we have
\begin{align*}
\|a(t)\|_{L^8(X)} &\leq c(g)\|a(t)\|_{W_{A_1}^{1,8/3}(X)}
\leq c(g)\|a(t)\|_{W_{A_1}^{2,4}(X)},
\\
\|b(t)\|_{L^8(X)} &\leq c(g)\|b(t)\|_{W_{A_1}^{2,4}(X)}, \quad\forall\, t\in [0,T),
\end{align*}
while, noting that $W^{1,q}(X) \hookrightarrow C(X)$ for any $q>4$ by \cite[Theorem 4.12]{AdamsFournier},
\begin{align*}
\|g-\bar g\|_{C(X)} &\leq c(g)\|g-\bar g\|_{W^{1,8}(X)},
\\
\|g-\bar g\|_{W^{1,8}(X)} &\leq c(g)\|g-\bar g\|_{W^{2,8/3}(X)}
 \leq c(g)\|g-\bar g\|_{W^{2,4}(X)}.
\end{align*}
In Step \ref{step:Norm_criterion_for_no_bubbling_in_barA} (see Equation \eqref{eq:Struwe_lemma_3-6_eps1_and_r1_dependence_on_metrics_C1_equivalent_to_g} in the proof of Claim \ref{claim:Norm_criterion_for_no_bubbling_in_barA}), we appealed to Remark \ref{rmk:Schlatter_2-4_and_Struwe_3-6_local_equivalent_Riemannian_metric} to note that $C^0$ norms of the curvature, $\Riem(\bar g)$, distances between pairs of points in $X$, $\dist_{\bar g}(x_0,x_1)$, the injectivity radius $\Inj(X,\bar g)$, tensor norms, volumes of open subsets, and the positive constants $\bar\eps_1$ and $\bar R_0$  defined by $\bar g$ in Lemma \ref{lem:Schlatter_2-4_and_Struwe_3-6_local} will be comparable to the corresponding quantities defined by $g$, with comparison constants depending on $g$ and $K$. However, just as in Lemma \ref{lem:Schlatter_2-3_and_Struwe_3-3}, those quantities depend on the metric $g$ only through
$$
\|\Riem(g)\|_{L^4(X,g)}, \quad \|g-g_0\|_{C(X)}, \quad\hbox{and}\quad \Inj(X,g_0),
$$
where $g_0$ is any fixed reference metric on $X$, and this dependence is continuous. In particular, the positive constants $\bar\eps_1$ and $\bar R_0$  defined by $\bar g$ in Lemma \ref{lem:Schlatter_2-4_and_Struwe_3-6_local} will be comparable to the corresponding quantities defined by $g$, but with comparison constants depending on $g$ and
$$
\|\Riem(\bar g)\|_{L^4(X,g)}, \quad \|\bar g-g\|_{C(X)}, \quad\hbox{and}\quad \Inj(X,g),
$$
and this dependence is continuous. The remainder of the proof of Theorem \ref{thm:Global_apriori_estimate_difference_solutions_Yang-Mills_heat_equations_pair_metrics} is unaffected by the change in the hypotheses, so this completes the proof of Corollary \ref{cor:Global_apriori_estimate_difference_solutions_YM_heat_eqns_pair_metrics_C0_close}.
\end{proof}

\subsection{Local comparison estimate for two solutions to the Yang-Mills heat equations defined by two nearby Riemannian metrics}
\label{subsec:Local_comparison_estimates_two_solutions_Yang-Mills_heat_equations_for_two_metrics}
Ultimately we wish to compare the pair of solutions, $A(t)$ and $\bar A(t)$, to the Yang-Mills heat equation \eqref{eq:Yang-Mills_heat_equation_as_perturbation_rough_Laplacian_heat_equation} defined by two close Riemannian metrics, $g$ and $\bar g$, and common initial data, $A_0$, over a precompact open subset $U \Subset X\less\Sigma$, where $\Sigma = \{x_1,\ldots,x_L\}$ is a finite set of points. In our application, $\Sigma$ will be a set of points where the flow, $A(t)$, acquires bubble singularities as $t \nearrow T$. In Theorem \ref{thm:Local_apriori_estimate_difference_solutions_Yang-Mills_heat_equations_pair_metrics} below, we provide a slightly more general statement that yields the result we need.

\begin{thm}[Local \apriori estimate for the difference between solutions to the Yang-Mills heat equations defined by two close Riemannian metrics and common initial data]
\label{thm:Local_apriori_estimate_difference_solutions_Yang-Mills_heat_equations_pair_metrics}
Let $G$, $P$, $X$, $g$, $A_1$, $K$, $T$ be as in the hypotheses of Theorem \ref{thm:Global_apriori_estimate_difference_solutions_Yang-Mills_heat_equations_pair_metrics} and $\mu_0$ and $z_1$ as in the assertions of Theorem \ref{thm:Global_apriori_estimate_difference_solutions_Yang-Mills_heat_equations_pair_metrics}. Then there are a positive constant, $z_2=z_2(g)$, and a small enough positive constants, $\eps = \eps(A_1,g,K,T) \in (0, 1]$ and $\eta = \eta(A_1,g,K,T) \in (0, 1]$, with the following significance. Let $\bar g$ be a Riemannian metric on $X$ obeying \eqref{eq:Riem_barg_C0_norm_lessthan_bound}, \eqref{eq:Riemannian_metrics_g_minus_barg_C1_plus_W24_norm_lessthan_small_positive_constant}, and \eqref{eq:Riemannian_metrics_g_minus_barg_C2_lessthan_bound}. Furthermore, let $U' \Subset U \subset X$ be open subsets and $\chi \in C_0^\infty(U)$ a cut-off function such that $0 \leq \chi \leq 1$ on $X$ with $\supp\chi \subset U$ and $\chi = 1$ on $U'$ and $\Omega := U \less \bar U'$ obeying
\begin{subequations}
\label{eq:Volume_Omega_and_cut-off_function_H2_norm_leq_epsilon}
\begin{align}
\label{eq:Volume_Omega_leq_epsilon}
\Vol_g(\Omega) &\leq \eps,
\\
\label{eq:Cutoff_function_L4_norm_nabla_chi_plus_L2_norm_nabla2_chi_leq_epsilon}
\|\nabla\chi\|_{L^4(X)} + \|\nabla^2\chi\|_{L^2(X)} &\leq \eps,
\end{align}
\end{subequations}
Let $A(t)$, for $t\in [0, T)$, and $\bar A(t)$, for $t\in [0, \bar\sT_U)$ with maximal lifetime $\bar\sT_U \in (0, \infty]$, be strong solutions to the Yang-Mills heat equation \eqref{eq:Yang-Mills_heat_equation_as_perturbation_rough_Laplacian_heat_equation} on $P \restriction U$ for the Riemannian metrics, $g$ and $\bar g$, respectively, and common initial data, $A_0$. Suppose that $A(t)$ and $\bar A(t)$ obey
\begin{subequations}
\label{eq:Yang-Mills_gradient_flow_A_minus_A1_and_barA_minus_A_bounds_local}
\begin{align}
\label{eq:Yang-Mills_gradient_flow_A_minus_A1_Linfty_time_Linfty_space_plus_L4nabla_space_plus_L2nabla2_space_leq_bound_local}
\|A - A_1\|_{L^\infty(0,T;L^\infty(U))} + \|\nabla_{A_1}(A - A_1)\|_{L^\infty(0,T;L^4(U))} \quad &
\\
\notag
+ \|\nabla_{A_1}^2(A - A_1)\|_{L^\infty(0,T;L^2(U))} &\leq K,
\\
\label{eq:Yang-Mills_gradient_flow_barA_minus_A1_annulus_bound}
\|\bar A - A_1\|_{L^\infty(0,T;L^\infty(\Omega))} + \|\nabla_{A_1}(\bar A - A_1)\|_{L^\infty(0,T;L^4(\Omega))}
&\leq K,
\end{align}
\end{subequations}
and that $A_1$ obeys \eqref{eq:F_A1_C0_norm_lessthan_bound}. Then $\bar A - A_1 \in C([0,T]; H_{A_1}^1(U';\Lambda^1\otimes\ad P))$, so $\bar\sT_{U'} > T$, where $\bar\sT_{U'}$ is the maximal lifetime of $\bar A$ on $P \restriction U'$, and $\bar A(t)$ obeys
\begin{multline}
\label{eq:Apriori_estimate_for_A_minus_barA_in_terms_C1_and_W24_norms_g_minus_barg_local}
\|A - \bar A\|_{L^2(0,T;H_{A_1}^2(U'))} + \|\partial_t (A - \bar A)\|_{L^2(0,T;L^2(U'))}
\\
+ \|A - \bar A\|_{L^\infty(0,T;L^4(U'))} + \|A - \bar A\|_{L^\infty(0,T;H_{A_1}^1(U'))}
\\
\leq z_1e^{\mu_0T}\left(\|g - \bar g\|_{C^1(X)} + \|g - \bar g\|_{W^{2,4}(X)}\right) + z_2\eps^{1/4}K\sqrt{T}.
\end{multline}
\end{thm}

\begin{proof}
We proceed by modifying the proof of Theorem \ref{thm:Global_apriori_estimate_difference_solutions_Yang-Mills_heat_equations_pair_metrics} at the appropriate stages and continue the notation adopted there.

\setcounter{step}{0}
\begin{step}[Derivation of the quasilinear parabolic equation for $b - \bar b$]
\label{step:Derivation_quasilinear_parabolic_equation_for_b_minus_barb_cut-off}
Consider $T_1 > 0$ such that $T_1 \leq T$ and $T_1 < \bar\sT_U$. We begin by multiplying the quasilinear parabolic equation \eqref{eq:Quasilinear_parabolic_equation_b_minus_barb} for $b-\bar b$ by the cut-off function, $\chi$, and commuting it with the differential operators acting on $b-\bar b$ and the quadratic and cubic $b-\bar b$ terms in that equation to give, a.e. on $(0, T_1)\times X$,
\begin{equation}
\label{eq:Quasilinear_parabolic_equation_b_minus_barb_cut-off}
\begin{aligned}
{}&\frac{\partial}{\partial t}\chi(b - \bar b) + \left(\bar\nabla_{A_1}^*\nabla_{A_1} + \mu\right)\chi(b - \bar b)
+ \Ric(\bar g) \times \chi(b - \bar b) + F_{A_1} \times \chi(b-\bar b)
\\
&\quad + e^{\mu t}\left[\chi(b - \bar b) \times \nabla_{A_1}b +  b \times \chi(b - \bar b)
+ b \times \nabla_{A_1}\chi(b - \bar b)\right] + e^{2\mu t}b \times b \times \chi(b - \bar b)
\\
&= -\chi\left(\nabla_{A_1}^*\nabla_{A_1} - \bar\nabla_{A_1}^*\nabla_{A_1}\right)b
- (\Ric(g) - \Ric(\bar g)) \times \chi b
\\
&\quad + (* - \bar*)F_{A_1} \times \chi b + e^{\mu t}\left[\chi b \times (\nabla_{A_1}^* - \bar\nabla_{A_1}^*)b
 + (* - \bar*)\chi b \times \nabla_{A_1}b + \chi(\nabla(* - \bar*))b \times b\right]
\\
&\quad + (* - \bar*)e^{2\mu t}\chi b \times b \times b
\\
&\quad + e^{\mu t}\left[\chi(b - \bar b) \times \chi(b - \bar b)
+ \chi(b - \bar b) \times \nabla_{A_1}\chi(b - \bar b)\right]
+ e^{2\mu t}b \times \chi(b - \bar b) \times \chi(b - \bar b)
\\
&\quad + e^{2\mu t}\chi(b - \bar b) \times \chi(b - \bar b) \times \chi(b - \bar b)
\\
&\quad + \hbox{Sum of commutator terms},
\end{aligned}
\end{equation}
where
\begin{equation}
\label{eq:Quasilinear_parabolic_equation_b_minus_barb_cut-off_sum_commutator_terms}
\begin{aligned}
{}&\hbox{Sum of commutator terms}
\\
&\quad= - \nabla\chi \times \nabla_{A_1}(b-\bar b) - \nabla^2\chi \times (b-\bar b)
- e^{\mu t}b \times \nabla\chi \times (b - \bar b)
\\
&\qquad - e^{\mu t}\chi(b - \bar b) \times \left[ (\chi-1)(b - \bar b)
+ (\chi-1)\nabla_{A_1}(b - \bar b) + \nabla\chi\times(b - \bar b) \right]
\\
&\qquad - e^{2\mu t}b \times \chi(b - \bar b) \times (\chi-1)(b - \bar b)
\\
&\qquad - e^{2\mu t} \chi(b - \bar b) \times \left[(\chi-1)(b - \bar b) \times (b - \bar b)
+ \chi(b - \bar b) \times (\chi-1)(b - \bar b) \right].
\end{aligned}
\end{equation}
This concludes Step \ref{step:Derivation_quasilinear_parabolic_equation_for_b_minus_barb_cut-off}.
\end{step}

Observe that, aside from the `sum of commutator terms' --- which we shall treat as a non-zero source, the quasilinear parabolic equation \eqref{eq:Quasilinear_parabolic_equation_b_minus_barb_cut-off} for $\chi(b-\bar b)$ has exactly the same structure as the equation \eqref{eq:Quasilinear_parabolic_equation_b_minus_barb} for $b-\bar b$.

\begin{step}[Estimates for the terms in \eqref{eq:Quasilinear_parabolic_equation_b_minus_barb_cut-off} excluding the sum of commutators]
\label{step:Estimates_for_non_commutator_terms_in_quasilinear_parabolic_equation_b_minus_barb_cut-off}
We modify our development of our estimates for the corresponding terms in \eqref{eq:Quasilinear_parabolic_equation_b_minus_barb} in order to make use of the stronger hypotheses on $a = A-A_1$ on $U$, noting that $\supp\chi \subset U$ and continuing to write $X_T := (0,T)\times X$ and similarly for $U_T$ or $X_{T_1}$. Thus, in place of the estimate \eqref{eq:Estimate_for_L2_time_space_heat_operator_on_b_minus_barb_raw}, we have the slightly simpler bound,
\begin{equation}
\label{eq:Estimate_for_L2_time_space_heat_operator_on_b_minus_barb_raw_cut-off}
\begin{aligned}
{} & \|(L_{A_1} + \mu)\chi(b - \bar b)\|_{L^2(X_{T_1})}
\\
&\quad\leq c\left(\|\Ric(\bar g)\|_{C(X)} + \|F_{A_1}\|_{C(X)} \right)\|\chi(b - \bar b)\|_{L^2(X_{T_1})}
\\
&\qquad + c\left(\|e^{\mu t}b\|_{L^\infty(0, T; L^4(U))} + c\|e^{\mu t}\nabla_{A_1}b\|_{L^\infty(0, T; L^4(U))} \right)
 \|\chi(b - \bar b)\|_{L^2(0,{T_1}; L^4(X))}
\\
&\qquad + c\|e^{\mu t}b\|_{L^\infty(U_T)}^2 \|\chi(b - \bar b)\|_{L^2(X_{T_1})}
+ c\|e^{\mu t}b\|_{L^\infty(U_T)} \|\nabla_{A_1}\chi(b - \bar b)\|_{L^2(X_{T_1})}
\\
&\qquad + \|\chi\left(\nabla_{A_1}^*\nabla_{A_1} - \bar\nabla_{A_1}^*\nabla_{A_1}\right)b\|_{L^2(X_T)}
+ c\|(* - \bar*)F_{A_1}\|_{C(X)}\|\chi b\|_{L^2(X_T)}
\\
&\qquad + c\|\Ric(g) - \Ric(\bar g)\|_{L^4(X)} \|\chi b\|_{L^2(0,T; L^4(X))}
\\
&\qquad + c\|e^{\mu t}b \|_{L^\infty(0, T; L^4(U))} \|\chi(\nabla_{A_1}^* - \bar\nabla_{A_1}^*)b\|_{L^2(0,T; L^4(X))}
\\
&\qquad + c\|(* - \bar*)e^{\mu t}b \|_{L^\infty(0, T; L^4(U))} \|\chi\nabla_{A_1}b\|_{L^2(0,T; L^4(X))}
\\
&\qquad + c\|(\nabla(* - \bar*))e^{\mu t}b \|_{L^\infty(0, T; L^4(U))} \|\chi b\|_{L^2(0,T; L^4(X))}
\\
&\qquad + c\|(* - \bar*)e^{\mu t}b \|_{L^\infty(U_T)} \|e^{\mu t}b \|_{L^\infty(U_T)} \|\chi b \|_{L^2(X_T)}
\\
&\qquad + e^{\mu T}\|\chi(b - \bar b) \times \chi(b - \bar b)\|_{L^2(X_{T_1})}
+ e^{\mu T}\|\chi(b - \bar b) \times \nabla_{A_1}\chi(b - \bar b)\|_{L^2(X_{T_1})}
\\
&\qquad + ce^{2\mu T}\|b\|_{L^\infty(U_T)} \|\chi(b - \bar b) \times \chi(b - \bar b) \|_{L^2(X_{T_1})}
\\
&\qquad + e^{2\mu T}\|\chi(b - \bar b) \times \chi(b - \bar b) \times \chi(b - \bar b) \|_{L^2(X_{T_1})}
\\
&\qquad + \|\hbox{Sum of commutator terms}\|_{L^2(X_T)}.
\end{aligned}
\end{equation}
By appealing to our hypothesis \eqref{eq:Yang-Mills_gradient_flow_A_minus_A1_Linfty_time_Linfty_space_plus_L4nabla_space_plus_L2nabla2_space_leq_bound_local} for $a(t) = A(t)-A_1$ on $U_T$, and thus $b(t) = e^{-\mu t}a(t)$ on $U_T$, we can now repeat our proof of the \apriori estimate \eqref{eq:Apriori_estimate_for_b_minus_barb_in_terms_C1_and_W24_norms_g_minus_barg} for $b - \bar b$ \mutatis to give the following \apriori estimate for $\chi(b - \bar b)$, valid for $\mu \geq \mu_0$,
\begin{multline}
\label{eq:Apriori_estimate_for_b_minus_barb_in_terms_C1_and_W24_norms_g_minus_barg_cut-off_raw}
\|\chi(b - \bar b)\|_{L^2(0,T_1;H_{A_1}^2(X))} + \|\partial_t \chi(b - \bar b)\|_{L^2(0,T_1;L^2(X))}
\\
+ \|\chi(b - \bar b)\|_{L^\infty(0,T_1;L^4(X))} + \|\chi(b - \bar b)\|_{L^\infty(0,T_1;H_{A_1}^1(X))}
\\
\leq C_2\|g-\bar g\|_{C^1(X)} + c\|g-\bar g\|_{W^{2,4}(X)}
\\
+ \|\hbox{Sum of commutator terms}\|_{L^2(X_T)},
\end{multline}
and \emph{provided} $\chi(b - \bar b)$ obeys the following analogue of the condition \eqref{eq:Linfty_time_H1_space_b_minus_barb_bound_allowing_rearrangement_quadratic_and_cubic_terms}, namely
\begin{equation}
\label{eq:Linfty_time_H1_space_b_minus_barb_bound_allowing_rearrangement_quadratic_and_cubic_terms_cut-off}
\|\chi(b - \bar b)\|_{L^\infty(0,T_1;H_{A_1}^1(X))} \leq \sigma,
\end{equation}
where $\sigma = \sigma(A_1,g,K,T)$ is as defined in \eqref{eq:Linfty_time_H1_space_b_minus_barb_bound_allowing_rearrangement_quadratic_and_cubic_terms}. In writing \eqref{eq:Apriori_estimate_for_b_minus_barb_in_terms_C1_and_W24_norms_g_minus_barg_cut-off_raw}, keeping in mind an application in the sequel, we have explicitly added the term, $\|\chi(b - \bar b)\|_{L^\infty(0,T_1; L^4(X))}$, even though this is of course bounded by $c\|\chi(b - \bar b)\|_{L^\infty(0,T_1;H_{A_1}^1(X))}$, for some positive constant, $c = c(g)$. This concludes Step \ref{step:Estimates_for_non_commutator_terms_in_quasilinear_parabolic_equation_b_minus_barb_cut-off}.
\end{step}

\begin{step}[Estimate for the sum of commutator terms in \eqref{eq:Quasilinear_parabolic_equation_b_minus_barb_cut-off}]
\label{step:Estimates_for_sum_commutator_terms_in_quasilinear_parabolic_equation_b_minus_barb_cut-off}
The additional commutator terms (not present in \eqref{eq:Quasilinear_parabolic_equation_b_minus_barb}) arise because we have used the fact that the commutator of the connection Laplace operator and $\chi$ is given by
\begin{align*}
\bar\nabla_{A_1}^*\nabla_{A_1} \chi(b-\bar b)
&= -\bar*\bar\nabla_{A_1}\bar*\left( \chi\nabla_{A_1}(b-\bar b) + d\chi\wedge(b-\bar b) \right)
\\
&= \chi\bar\nabla_{A_1}^*\nabla_{A_1}(b-\bar b) + \nabla\chi \times \nabla_{A_1}(b-\bar b)
+ \nabla^2\chi \times (b-\bar b),
\end{align*}
and the commutators for the quadratic and cubic $b-\bar b$ terms and $\chi$ are given by
\begin{align*}
b \times \nabla_{A_1}\chi(b - \bar b)
&=
b \times \chi\nabla_{A_1}(b - \bar b) + b \times \nabla\chi \times (b - \bar b),
\\
\chi(b - \bar b) \times \chi(b - \bar b)
&=
\chi(b - \bar b) \times (b - \bar b) + \chi(b - \bar b) \times (\chi-1)(b - \bar b),
\\
\chi(b - \bar b) \times \nabla_{A_1}\chi(b - \bar b)
&=
\chi(b - \bar b) \times \nabla_{A_1}(b - \bar b) + \chi(b - \bar b) \times (\chi-1)\nabla_{A_1}(b - \bar b)
\\
&\quad + \chi(b - \bar b) \times \nabla\chi\times(b - \bar b),
\\
\chi(b - \bar b) \times \chi(b - \bar b) \times \chi(b - \bar b)
&= \chi(b - \bar b) \times (b - \bar b) \times (b - \bar b)
\\
&\quad + \chi(b - \bar b) \times (\chi-1)(b - \bar b) \times (b - \bar b)
\\
&\quad + \chi(b - \bar b) \times \chi(b - \bar b) \times (\chi-1)(b - \bar b).
\end{align*}
We now fix $\mu = \mu_0$, where $\mu_0 \in [1,\infty)$ is as in Theorem \ref{thm:Global_apriori_estimate_difference_solutions_Yang-Mills_heat_equations_pair_metrics}, and make the

\begin{claim}
\label{claim:Quasilinear_parabolic_equation_b_minus_barb_cut-off_sum_norms_commutators_lessthan_sigma_over_4}
There is a positive constant, $z_2 = z_2(g)$, such that the sum of commutator terms \eqref{eq:Quasilinear_parabolic_equation_b_minus_barb_cut-off_sum_commutator_terms} obeys
\begin{equation}
\label{eq:Quasilinear_parabolic_equation_b_minus_barb_cut-off_sum_norms_commutators_lessthan_sigma_over_4}
\|\hbox{\emph{Sum of commutator terms}}\|_{L^2(X_T)} \leq z_2\eps^{1/4} K\sqrt{T},
\end{equation}
where $\eps$ and $K$ are the constants in the hypotheses of Theorem \ref{thm:Local_apriori_estimate_difference_solutions_Yang-Mills_heat_equations_pair_metrics}.
\end{claim}

\begin{proof}[Proof of Claim \ref{claim:Quasilinear_parabolic_equation_b_minus_barb_cut-off_sum_norms_commutators_lessthan_sigma_over_4}]
Using our hypothesis \eqref{eq:Cutoff_function_L4_norm_nabla_chi_plus_L2_norm_nabla2_chi_leq_epsilon} on $\nabla\chi$ and $\nabla^2\chi$ as needed and the hypotheses \eqref{eq:Yang-Mills_gradient_flow_A_minus_A1_Linfty_time_Linfty_space_plus_L4nabla_space_plus_L2nabla2_space_leq_bound_local} and \eqref{eq:Yang-Mills_gradient_flow_barA_minus_A1_annulus_bound} for the norms of $a$ and $\bar a$ on $\Omega$, respectively, we estimate the $L^2(X_T; \Lambda^1\otimes\ad P)$ norms of the commutator terms via
\begin{align*}
\|\nabla\chi \times \nabla_{A_1}(b-\bar b)\|_{L^2(X_T)}
&\leq
c\|\nabla\chi\|_{L^4(X)} \|\nabla_{A_1}(b-\bar b)\|_{L^2(0,T;L^4(\Omega))}
\\
&\leq c\eps \|e^{-\mu t}\nabla_{A_1}(a-\bar a)\|_{L^2(0,T;L^4(\Omega))}
\\
&\leq 2c\eps K,
\end{align*}
and
\begin{align*}
\|\nabla^2\chi \times (b-\bar b)\|_{L^2(X_T)}
&\leq
c\|\nabla^2\chi\|_{L^2(X)} \|b-\bar b\|_{L^2(0,T;L^\infty(\Omega))}
\\
&\leq c\eps \|e^{-\mu t}(a-\bar a)\|_{L^2(0,T;L^\infty(\Omega))}
\\
&\leq 2c\eps K\sqrt{T},
\end{align*}
and
\begin{align*}
\|e^{\mu t}b \times \nabla\chi \times (b - \bar b)\|_{L^2(X_T)}
&\leq c\|\nabla\chi\|_{L^4(X)} \|e^{\mu t}b \times (b - \bar b)\|_{L^2(0,T;L^4(\Omega))}
\\
&\leq c\|\nabla\chi\|_{L^4(X)} \|b - \bar b\|_{L^\infty(0,T;L^4(\Omega))}
\|e^{\mu t}b\|_{L^2(0,T;L^\infty(\Omega))}
\\
&= c\eps \|e^{-\mu t}(a-\bar a)\|_{L^\infty(0,T;L^4(\Omega))} \|a\|_{L^2(0,T;L^\infty(\Omega))}
\\
&\leq 2c\eps K\sqrt{T}(\Vol_g(\Omega))^{1/4}
\\
&\leq 2c\eps^{5/4} K\sqrt{T},
\end{align*}
where the last inequality follows from the hypothesis \eqref{eq:Volume_Omega_leq_epsilon} on $\Vol_g(\Omega)$, and
\begin{align*}
\|e^{\mu t}\chi(b - \bar b) \times \nabla\chi\times(b - \bar b)\|_{L^2(X_T)}
&\leq c\|\nabla\chi\|_{L^4(X)} \|e^{\mu t}(b - \bar b) \times (b - \bar b)\|_{L^2(0,T;L^4(\Omega))}
\\
&\leq c\eps \|e^{\mu t}(b - \bar b)\|_{L^\infty(0,T;L^4(\Omega))}
\|b - \bar b\|_{L^2(0,T;L^\infty(\Omega))}
\\
&= c\eps \|a-\bar a\|_{L^\infty(0,T;L^4(\Omega))} \|e^{-\mu t}(a-\bar a)\|_{L^2(0,T;L^\infty(\Omega))}
\\
&\leq 2c\eps^{5/4} K\sqrt{T}.
\end{align*}
The remaining quadratic commutator terms are estimated by
\begin{align*}
\|e^{\mu t}\chi(b - \bar b) \times (\chi-1)(b - \bar b)\|_{L^2(X_T)}
&\leq c\|e^{\mu t}(b - \bar b)\|_{L^\infty(0,T;L^\infty(\Omega))} \|b - \bar b\|_{L^2(0,T;L^2(\Omega))}
\\
&= c\|a-\bar a\|_{L^\infty(0,T;L^\infty(\Omega))} \|e^{-\mu t}(a-\bar a)\|_{L^2(0,T;L^2(\Omega))}
\\
&\leq 2c\eps^{1/2} K\sqrt{T},
\end{align*}
and
\begin{align*}
{}& \|e^{\mu t}\chi(b - \bar b) \times (\chi-1)\nabla_{A_1}(b - \bar b)\|_{L^2(X_T)}
\\
&\quad \leq c\|e^{\mu t}(b - \bar b)\|_{L^2(0,T;L^4(\Omega))} \|\nabla_{A_1}(b - \bar b)\|_{L^\infty(0,T;L^4(\Omega))}
\\
&\quad = c\|a-\bar a\|_{L^2(0,T;L^4(\Omega))} \|e^{-\mu t}\nabla_{A_1}(a-\bar a)\|_{L^\infty(0,T;L^4(\Omega))}
\\
&\quad \leq 2c\eps^{1/4} K\sqrt{T},
\end{align*}
and the cubic commutator terms are estimated by
\begin{align*}
{}& \|e^{2\mu t}b \times \chi(b - \bar b) \times (\chi-1)(b - \bar b)\|_{L^2(X_T)}
\\
&\quad \leq c\|e^{\mu t}b\|_{L^\infty(0,T;L^\infty(\Omega))}
\|e^{\mu t}(b - \bar b)\|_{L^\infty(0,T;L^\infty(\Omega))} \|b - \bar b\|_{L^2(0,T;L^2(\Omega))}
\\
&\quad = c\|a\|_{L^\infty(0,T;L^\infty(\Omega))}
\|a-\bar a\|_{L^\infty(0,T;L^\infty(\Omega))} \|e^{-\mu t}(a-\bar a)\|_{L^2(0,T;L^2(\Omega))}
\\
&\quad \leq 2c\eps^{1/2} K\sqrt{T},
\end{align*}
and
\begin{align*}
{}& \|e^{2\mu t} \chi(b - \bar b) \times (\chi-1)(b - \bar b) \times (b - \bar b)\|_{L^2(X_T)}
\\
&\quad \leq c\|e^{\mu t}(b - \bar b)\|_{L^\infty(0,T;L^\infty(\Omega))}^2
\|b - \bar b\|_{L^2(0,T;L^2(\Omega))}
\\
&\quad = c\|a-\bar a\|_{L^\infty(0,T;L^\infty(\Omega))}^2 \|e^{-\mu t}(a-\bar a)\|_{L^2(0,T;L^2(\Omega))}
\\
&\quad \leq 2c\eps^{1/2} K\sqrt{T}.
\end{align*}
The conclusion \eqref{eq:Quasilinear_parabolic_equation_b_minus_barb_cut-off_sum_norms_commutators_lessthan_sigma_over_4} now follows by combining the preceding estimates. This completes the proof of Claim \ref{claim:Quasilinear_parabolic_equation_b_minus_barb_cut-off_sum_norms_commutators_lessthan_sigma_over_4}.
\end{proof}
This concludes Step \ref{step:Estimates_for_sum_commutator_terms_in_quasilinear_parabolic_equation_b_minus_barb_cut-off}.
\end{step}

\begin{step}[\Apriori estimate for $\chi(b-\bar b)$]
\label{step:Apriori_estimate_for_b_minus_barb_cut-off}
The \apriori estimate \eqref{eq:Apriori_estimate_for_b_minus_barb_in_terms_C1_and_W24_norms_g_minus_barg_cut-off_raw} and Claim \ref{claim:Quasilinear_parabolic_equation_b_minus_barb_cut-off_sum_norms_commutators_lessthan_sigma_over_4} now imply that if $\chi(b - \bar b)$ obeys the condition \eqref{eq:Linfty_time_H1_space_b_minus_barb_bound_allowing_rearrangement_quadratic_and_cubic_terms_cut-off}, then
\begin{multline}
\label{eq:Apriori_estimate_for_b_minus_barb_in_terms_C1_and_W24_norms_g_minus_barg_cut-off}
\|\chi(b - \bar b)\|_{L^2(0,T_1;H_{A_1}^2(X))} + \|\partial_t \chi(b - \bar b)\|_{L^2(0,T_1;L^2(X))}
\\
+ \|\chi(b - \bar b)\|_{L^\infty(0,T_1;L^4(X))} + \|\chi(b - \bar b)\|_{L^\infty(0,T_1;H_{A_1}^1(X))}
\\
\leq C_2\left(\|g-\bar g\|_{C^1(X)} + \|g-\bar g\|_{W^{2,4}(X)}\right) + z_2\eps^{1/4} K\sqrt{T},
\end{multline}
where $C_2 = C_2(g,K)$ is as in \eqref{eq:Apriori_estimate_for_b_minus_barb_in_terms_C1_and_W24_norms_g_minus_barg} and we have fixed $\mu = \mu_0$. This concludes Step \ref{step:Apriori_estimate_for_b_minus_barb_cut-off}.
\end{step}

\begin{step}[Norm criterion for continuous temporal extension of $\bar A(t)$ over $U'$]
\label{step:Norm_criterion_for_no_bubbling_in_barA_local}
We next establish the following local analogue of Claim \ref{claim:Norm_criterion_for_no_bubbling_in_barA} from the proof of Theorem \ref{thm:Global_apriori_estimate_difference_solutions_Yang-Mills_heat_equations_pair_metrics}.

\begin{claim}[Norm criterion for continuous temporal extension of $\bar A(t)$ over $U'$]
\label{claim:Norm_criterion_for_no_bubbling_in_barA_local}
There is a positive constant, $\zeta = \zeta(A_1,g,K,T) \in (0, 1]$ with the following significance. If $T_2 \in (0, T]$ has the property that $\bar A - A_1 \in C_{\loc}([0,T_2); H_{A_1}^1(U;\Lambda^1\otimes\ad P))$ and
\begin{equation}
\label{eq:Linfty_time_H1_space_b_minus_barb_leq_small_constant_ensuring_no_bubbling_local}
\|b(t) - \bar b(t)\|_{L^4(U)} + \|\nabla_{A_1}(b(t) - \bar b(t))\|_{L^2(U)} <  \zeta, \quad \forall\, t \in [0,T_2),
\end{equation}
then $\bar A - A_1 \in C([0,T_2]; H_{A_1}^1(U';\Lambda^1\otimes\ad P))$ and $\bar\sT_{U'} > T_2$, where $\bar\sT_{U'}$ is the maximal lifetime of $\bar A$ on $P \restriction U'$.
\end{claim}

\begin{proof}[Proof of Claim \ref{claim:Norm_criterion_for_no_bubbling_in_barA_local}]
Although the argument is formally similar to that used to establish Claim \ref{claim:Norm_criterion_for_no_bubbling_in_barA}, there differences. With subsequent applications in mind, it will be very important to be aware of any dependency, if present, of the constant, $\zeta = \zeta(A_1,g,K,T)$, in the hypotheses of Claim \ref{claim:Norm_criterion_for_no_bubbling_in_barA} on the open subsets, $U'$ or $U$. For that reason, we shall include the details of the argument in order to check whether such dependencies arise.

Let $\bar\eps_1 = \bar\eps_1(\bar g)= \bar\eps_1(g,K) \in (0, 1]$ and $\bar r_1 = \bar r_1(\bar g) = \bar r_1(g,K) \in (0, 1]$ be the positive constants described in \eqref{eq:Struwe_lemma_3-6_eps1_and_r1_dependence_on_metrics_C1_equivalent_to_g}. It will be convenient to define the separation constant,
$$
\delta := \dist_g(U', \partial U).
$$
Observe that, for any positive constant, $r \in (0, \delta \vee \Inj(X,g)]$, and any point $x \in U'$,
\begin{align*}
\|F_{\bar A}(t)\|_{L^2(B_r(x))} &\leq \|F_A(t)\|_{L^2(B_r(x))} + \|F_A(t)-F_{\bar A}(t)\|_{L^2(B_r(x))}
\\
&\leq \|F_A(t)\|_{L^2(B_r(x))} + \|F_A(t)-F_{\bar A}(t)\|_{L^2(U)}, \quad\forall\, t \in [0, T_2),
\end{align*}
and thus,
\begin{equation}
\label{eq:Linfty_time_L2_ball_FbarA_leq_Linfty_time_L2_ball_FA_plus_Linfty_time_L2_space_FA_minus_FbarA_local}
\sup_{x\in U'}\|F_{\bar A}(t)\|_{L^2(B_r(x))}
\leq
\sup_{x\in U'}\|F_A\|_{L^2(B_r(x))} +  \|F_A(t) - F_{\bar A}(t)\|_{L^2(U)}, \quad\forall\, t \in [0, T_2).
\end{equation}
Exactly as in the proof of Claim \ref{claim:Norm_criterion_for_no_bubbling_in_barA}, we discover that
\begin{multline}
\label{eq:Linfty_time_L2_space_FA_minus_FbarA_raw_local}
\|F_A(t) - F_{\bar A}(t)\|_{L^2(U)}
\\
\leq c_1e^{\mu_0 T}\|\nabla_{A_1}(b - \bar b)(t)\|_{L^2(U)} + c_1e^{\mu_0 T}\|e^{\mu_0 t}b(t)\|_{L^4(U)}\|(b - \bar b)(t)\|_{L^4(U)}
\\
+ c_1e^{2\mu_0 T}\|(b - \bar b)(t)\|_{L^4(U)}^2, \quad\forall\, t\in [0, T_2),
\end{multline}
where $c_1 \in [1,\infty)$ is a positive constant depending at most on $g$. Provided
$$
\|b - \bar b\|_{L^\infty(0,T_2;L^4(U))} \leq 1,
$$
which will be assured by our hypothesis \eqref{eq:Linfty_time_H1_space_b_minus_barb_leq_small_constant_ensuring_no_bubbling_local} (since $\zeta \in (0, 1])$) and recalling that, by \eqref{eq:Yang-Mills_gradient_flow_A_minus_A1_Linfty_time_Linfty_space_plus_L4nabla_space_plus_L2nabla2_space_leq_bound_local} and the fact that $a(t) = e^{\mu_0 t}b(t)$,
$$
\|e^{\mu_0 t}b(t)\|_{L^4(U)} = \|a(t)\|_{L^4(U)} \leq K, \quad\forall\, t \in [0, T),
$$
we thus obtain from \eqref{eq:Linfty_time_L2_space_FA_minus_FbarA_raw_local} that,
\begin{multline*}
\|F_A(t) - F_{\bar A}(t)\|_{L^2(U)}
\\
\leq
c_1\left(1 + K + e^{\mu_0 T}\right) e^{\mu_0 T}
\left( \|\nabla_{A_1}(b - \bar b)(t)\|_{L^2(U)} + \|(b - \bar b)(t)\|_{L^4(U)}\right),
\quad\forall\, t\in [0, T_2).
\end{multline*}
Recalling that $\mu_0 = \mu_0(A_1,g,K)$, we define $\zeta \in (0, 1]$ by
\begin{equation}
\zeta \equiv \zeta(A_1,g,K,T) := \frac{\sqrt{\bar\eps_1}}{4c_1\left(1 + K + e^{\mu_0 T}\right)e^{\mu_0 T}}.
\end{equation}
Therefore, provided $b-\bar b$ obeys the hypothesis \eqref{eq:Linfty_time_H1_space_b_minus_barb_leq_small_constant_ensuring_no_bubbling_local}, we discover that
\begin{equation}
\label{eq:Linfty_time_L2_space_FA_minus_FbarA_local}
\|F_A(t) - F_{\bar A}(t)\|_{L^2(U)} < \frac{\sqrt{\bar\eps_1}}{4}, \quad\forall\, t\in [0, T_2).
\end{equation}
Next observe that, for any $x \in U'$ and $r \in (0, 1 \vee \delta \vee \Inj(X,g)]$, and writing $F_A(t) = F_{A_1} + d_{A_1}a(t) + [a(t), a(t)]$ over $U$,
\begin{align*}
\|F_A(t)\|_{L^2(B_r(x))} &\leq \|F_{A_1}\|_{L^2(B_r(x))} + \|d_{A_1}a(t)\|_{L^2(B_r(x))} + c\|[a(t), a(t)]\|_{L^2(B_r(x))}
\\
&\leq \left(\Vol_g(B_r(x))\right)^{1/2}\|F_{A_1}\|_{C(\bar B_r(x))} + \left(\Vol_g(B_r(x))\right)^{1/4}\|d_{A_1}a(t)\|_{L^4(B_r(x))}
\\
&\quad + c\left(\Vol_g(B_r(x))\right)^{1/2} \|a(t)\|_{C(\bar B_r(x))}^2
\\
&\leq c(Kr^2 + Kr + K^2r^2), \quad\hbox{(by \eqref{eq:F_A1_C0_norm_lessthan_bound} and \eqref{eq:Yang-Mills_gradient_flow_A_minus_A1_Linfty_time_Linfty_space_plus_L4nabla_space_plus_L2nabla2_space_leq_bound_local})}
\end{align*}
where $c$ is a positive constant depending at most on the Riemannian metric, $g$, on $X$, and thus (for a possibly larger $c$),
$$
\sup_{x\in U'}\|F_A(t)\|_{L^2(B_r(x))} \leq c(1 + K)Kr, \quad 0 < r \leq 1 \vee \delta \vee \Inj(X,g).
$$
Therefore, provided $r_0 = r_0(g,K,\bar\eps_1) = r_0(g,K) \in (0, 1 \vee \Inj(X,g)]$ is small enough that
$$
c(1 + K)Kr_0 < \frac{\sqrt{\bar\eps_1}}{4},
$$
where we recall that $\bar\eps_1 = \bar\eps_1(g,K) \in (0, 1]$, then
$$
\sup_{x\in U'}\|F_A(t)\|_{L^2(B_r(x))} < \frac{\sqrt{\bar\eps_1}}{4}, \quad\forall\, t\in [0, T) \hbox{ and }
0 < r \leq r_0 \vee \delta.
$$
Combining the preceding inequality with the bounds  \eqref{eq:Linfty_time_L2_ball_FbarA_leq_Linfty_time_L2_ball_FA_plus_Linfty_time_L2_space_FA_minus_FbarA_local} and \eqref{eq:Linfty_time_L2_space_FA_minus_FbarA_local} yields,
$$
\sup_{x\in U'}\|F_{\bar A}\|_{L^2(B_r(x))}
<
\frac{\sqrt{\bar\eps_1}}{2}, \quad\forall\, t\in [0, T_2) \hbox{ and }
0 < r \leq r_0 \vee \bar r_1 \vee \delta.
$$
But Lemma \ref{lem:Schlatter_2-4_and_Struwe_3-6_local} now implies that $\bar A$ extends continuously from $[0, T_2)$ to $[0, T_2]$, that is, $\bar A - A_1 \in C([0,T_2]; H_{A_1}^1(U';\Lambda^1\otimes\ad P))$ and $\bar\sT_{U'} > T_2$. This completes the proof of Claim \ref{claim:Norm_criterion_for_no_bubbling_in_barA_local}.
\end{proof}

This concludes Step \ref{step:Norm_criterion_for_no_bubbling_in_barA_local}.
\end{step}

\begin{step}[Conclusion that the maximal lifetime of $\bar A$ on $U'$ is greater than $T$]
\label{step:Maximal_lifetime_barA_geq_T_local}
We claim that
\begin{equation}
\label{eq:H1_space_norm_b_minus_barb_strictly_lessthan_nobubbling_and_rearrangement_constant_local}
\|(b - \bar b)(t)\|_{L^4(U)} + \|\nabla_{A_1}(b - \bar b)(t)\|_{L^2(U)}
< \min\{\zeta, \sigma/N\}, \quad\forall\, t \in [0, T),
\end{equation}
and $\bar\sT_{U'} > T$, where $\zeta$ is the positive constant in Claim \ref{claim:Norm_criterion_for_no_bubbling_in_barA_local} and $\sigma$ is the positive constant in \eqref{eq:Linfty_time_H1_space_b_minus_barb_bound_allowing_rearrangement_quadratic_and_cubic_terms} and $N := \left(1+\Vol_g(X)\right)^{1/4}$.

Suppose the claim is false. Then there must be a large enough time, $T_1 \in (0, T)$, such that
\begin{subequations}
\label{eq:maximal_T1_such_that_H1_space_norm_minus_barb_lessthan_nobubbling_and_rearrangement_constant_local}
\begin{align}
\label{eq:maximal_T1_such_that_H1_space_norm_minus_barb_strictly_lessthan_nobubbling_and_rearrangement_constant_local}
\|(b - \bar b)(t)\|_{L^4(U)} + \|\nabla_{A_1}(b - \bar b)(t)\|_{L^2(U)} &< \min\{\zeta, \sigma/N\},
\quad\forall\, t \in [0, T_1), \quad\hbox{but}
\\
\label{eq:maximal_T1_such_that_H1_space_norm_minus_barb_geq_nobubbling_and_rearrangement_constant_local}
\|(b - \bar b)(t)\|_{L^4(U)} + \|\nabla_{A_1}(b - \bar b)(t)\|_{L^2(U)} &\geq \min\{\zeta, \sigma/N\}.
\end{align}
\end{subequations}
For if not, then the strict inequality \eqref{eq:H1_space_norm_b_minus_barb_strictly_lessthan_nobubbling_and_rearrangement_constant_local} would hold and Claim \ref{claim:Norm_criterion_for_no_bubbling_in_barA_local} (with $T_2 = T$) would imply that $\bar A - A_1 \in C([0,T]; H_{A_1}^1(U';\Lambda^1\otimes\ad P))$ and $\bar\sT_{U'} > T$.

We choose the positive constant, $\eta \in (0, 1]$, in the hypothesis of Theorem  \ref{thm:Local_apriori_estimate_difference_solutions_Yang-Mills_heat_equations_pair_metrics}, small enough that
\begin{equation}
\label{eq:Definition_smallness_constant_g_minus_barg_metric_difference_cut-off}
\eta \equiv \eta(A_1,g,K,T) \leq \frac{\min\{\zeta, \sigma/N\}}{8(1+\kappa)C_2},
\end{equation}
where $C_2$ is the constant on the right-hand side of the \apriori estimate \eqref{eq:Apriori_estimate_for_b_minus_barb_in_terms_C1_and_W24_norms_g_minus_barg_cut-off} and $\kappa = \kappa(g)$ is the Sobolev constant for the embedding $H^1(X) \hookrightarrow L^4(X)$. We also choose the positive constant, $\eps \in (0, 1]$, in the hypotheses of Theorem \ref{thm:Local_apriori_estimate_difference_solutions_Yang-Mills_heat_equations_pair_metrics} by requiring that
\begin{equation}
\label{eq:Definition_smallness_constant_sum_commutator_terms}
\eps^{1/4} \equiv \eps^{1/4}(A_1,g,K,T) \leq \frac{\min\{\zeta, \sigma/N\}}{8(1+\kappa)z_2K\sqrt{T}},
\end{equation}
where $z_2$ is the constant on the right-hand side of the \apriori estimate \eqref{eq:Apriori_estimate_for_b_minus_barb_in_terms_C1_and_W24_norms_g_minus_barg_cut-off}. Therefore, by \eqref{eq:Riemannian_metrics_g_minus_barg_C1_plus_W24_norm_lessthan_small_positive_constant}, \eqref{eq:Definition_smallness_constant_g_minus_barg_metric_difference_cut-off}, and \eqref{eq:Definition_smallness_constant_sum_commutator_terms}, we see that
\begin{multline}
\label{eq:priori_estimate_for_b_minus_barb_in_terms_C1_and_W24_norms_g_minus_barg_cut-off_RHS_small}
C_2\left(\|g-\bar g\|_{C^1(X)} + \|g-\bar g\|_{W^{2,4}(X)}\right) + z_2\eps^{1/4} K\sqrt{T}
\\
\leq
C_2\eta + z_2\eps^{1/4} K\sqrt{T} =  \frac{\min\{\zeta, \sigma/N\}}{4(1+\kappa)}.
\end{multline}
Next, we observe that
\begin{align*}
\|\chi(b - \bar b)(t)\|_{H_{A_1}^1(X)} &\leq \|\chi(b - \bar b)(t)\|_{L^2(X)} + \|\nabla_{A_1}\chi(b - \bar b)(t)\|_{L^2(X)}
\\
&\leq \left(\Vol_g(X)\right)^{1/4} \|(b - \bar b)(t)\|_{L^4(U)} + \|\nabla_{A_1}(b - \bar b)(t)\|_{L^2(U)}
\\
&\quad + \|\nabla\chi\|_{L^4(X)}\|(b - \bar b)(t)\|_{L^4(\Omega)}
\\
&\leq \left(1+\Vol_g(X)\right)^{1/4} \|(b - \bar b)(t)\|_{L^4(U)} + \|\nabla_{A_1}(b - \bar b)(t)\|_{L^2(U)},
\end{align*}
where the last inequality follows from our hypothesis \eqref{eq:Cutoff_function_L4_norm_nabla_chi_plus_L2_norm_nabla2_chi_leq_epsilon} that $\|\nabla\chi\|_{L^4(X)} \leq \eps$, for some $\eps \in (0, 1]$. Consequently, as $N = \left(1+\Vol_g(X)\right)^{1/4}$, the preceding inequality yields
\begin{align*}
\|\chi(b - \bar b)\|_{L^\infty(0,T_1;H_{A_1}^1(X))}
&\leq
N\left( \|(b - \bar b)(t)\|_{L^4(U)} + \|\nabla_{A_1}(b - \bar b)(t)\|_{L^2(U)} \right)
\\
&\leq \sigma \quad\hbox{(by \eqref{eq:maximal_T1_such_that_H1_space_norm_minus_barb_strictly_lessthan_nobubbling_and_rearrangement_constant_local}).}
\end{align*}
Hence, the condition \eqref{eq:Linfty_time_H1_space_b_minus_barb_bound_allowing_rearrangement_quadratic_and_cubic_terms_cut-off} is satisfied and so the \apriori estimate \eqref{eq:Apriori_estimate_for_b_minus_barb_in_terms_C1_and_W24_norms_g_minus_barg_cut-off} and the bound \eqref{eq:priori_estimate_for_b_minus_barb_in_terms_C1_and_W24_norms_g_minus_barg_cut-off_RHS_small} on its right-hand side yield,
$$
\|\chi(b - \bar b)(t)\|_{H_{A_1}^1(X)} \leq \frac{\min\{\zeta, \sigma/N\}}{4(1+\kappa)}, \quad\forall\, t \in [0, T_1),
$$
and ensures that, applying the Sobolev embedding $H^1(X) \hookrightarrow L^4(X)$ \cite[Theorem 4.12]{AdamsFournier} and the Kato Inequality \eqref{eq:FU_6-20_first-order_Kato_inequality} in the usual way,
\begin{align*}
{}& \|(b - \bar b)(t)\|_{L^4(U')} + \|\nabla_{A_1}(b - \bar b)(t)\|_{L^2(U')}
\\
&\quad \leq \|\chi(b - \bar b)(t)\|_{L^4(X)} + \|\nabla_{A_1}\chi(b - \bar b)(t)\|_{L^2(X)}
\\
&\quad \leq \kappa\left(\|\chi(b - \bar b)(t)\|_{L^2(X)} + \|\nabla_{A_1}\chi(b - \bar b)(t)\|_{L^2(X)}\right)
+ \|\nabla_{A_1}\chi(b - \bar b)(t)\|_{L^2(X)}
\\
&\quad \leq 2(1+\kappa)\|\chi(b - \bar b)(t)\|_{H_{A_1}^1(X)}.
\end{align*}
Therefore, by combining the preceding two inequalities, we obtain
\begin{equation}
\label{eq:H1_Uprime_b_minus_barb_leq_half_zeta_sigma_over_N_on_0_T1_open}
\|(b - \bar b)(t)\|_{H_{A_1}^1(U')} \leq \frac{1}{2}\min\{\zeta, \sigma/N\}, \quad\forall\, t \in [0, T_1).
\end{equation}
On the other hand, for the annulus, $\Omega$, we have
\begin{align*}
\|(b - \bar b)(t)\|_{H_{A_1}^1(\Omega)}
&\leq
\|(b - \bar b)(t)\|_{L^2(\Omega)} + \|\nabla_{A_1}(b - \bar b)(t)\|_{L^2(\Omega)}
\\
&\leq \left(\Vol_g(\Omega)\right)^{1/2}\|(b - \bar b)(t)\|_{L^\infty(\Omega)}
+ \left(\Vol_g(\Omega)\right)^{1/4}\|\nabla_{A_1}(b - \bar b)(t)\|_{L^4(\Omega)}
\\
&\leq 2\eps^{1/4}K, \quad\forall\, t \in [0, T_1] \quad\hbox{(by \eqref{eq:Volume_Omega_leq_epsilon} and
\eqref{eq:Yang-Mills_gradient_flow_A_minus_A1_Linfty_time_Linfty_space_plus_L4nabla_space_plus_L2nabla2_space_leq_bound_local}).}
\end{align*}
By decreasing the size constant of the constant $\eps \in (0, 1]$ in \eqref{eq:Definition_smallness_constant_sum_commutator_terms}, if necessary, we may further suppose that $\eps$ is chosen small enough that
\begin{equation}
\label{eq:Definition_smallness_constant_b_minus_barb_annulus_term}
2\eps^{1/4}K \leq \frac{1}{2}\min\{\zeta, \sigma/N\},
\end{equation}
and therefore we also have
\begin{equation}
\label{eq:H1_Omega_b_minus_barb_leq_half_zeta_sigma_over_N_on_0_T1_closed}
\|(b - \bar b)(t)\|_{H_{A_1}^1(\Omega)} \leq \frac{1}{2}\min\{\zeta, \sigma/N\}, \quad\forall\, t \in [0, T_1].
\end{equation}
By combining \eqref{eq:H1_Uprime_b_minus_barb_leq_half_zeta_sigma_over_N_on_0_T1_open} and \eqref{eq:H1_Omega_b_minus_barb_leq_half_zeta_sigma_over_N_on_0_T1_closed}, we obtain
\begin{align*}
\|(b - \bar b)(t)\|_{H_{A_1}^1(U)} &\leq \|(b - \bar b)(t)\|_{H_{A_1}^1(U')} + \|(b - \bar b)(t)\|_{H_{A_1}^1(\Omega)}
\\
&\leq \frac{1}{2}\min\{\zeta, \sigma/N\}, \quad\forall\, t \in [0, T_1).
\end{align*}
Therefore, the condition \eqref{eq:Linfty_time_H1_space_b_minus_barb_leq_small_constant_ensuring_no_bubbling_local} holds with $T_2=T_1$. Consequently, Claim \ref{claim:Norm_criterion_for_no_bubbling_in_barA_local} (with $T_2 = T_1$) implies that $\bar A - A_1 \in C([0,T_1]; H_{A_1}^1(U';\Lambda^1\otimes\ad P))$ and
\begin{equation}
\label{eq:H1_Uprime_b_minus_barb_leq_half_zeta_sigma_over_N_closed}
\|(b - \bar b)(t)\|_{H_{A_1}^1(U')} \leq \frac{1}{2}\min\{\zeta, \sigma/N\}, \quad\forall\, t \in [0, T_1].
\end{equation}
Thus, by combining \eqref{eq:H1_Omega_b_minus_barb_leq_half_zeta_sigma_over_N_on_0_T1_closed} and  \eqref{eq:H1_Uprime_b_minus_barb_leq_half_zeta_sigma_over_N_closed}, we obtain
$$
\|(b - \bar b)(t)\|_{H_{A_1}^1(U)} \leq \frac{1}{2}\min\{\zeta, \sigma/N\}, \quad\forall\, t \in [0, T_1],
$$
which contradicts \eqref{eq:maximal_T1_such_that_H1_space_norm_minus_barb_geq_nobubbling_and_rearrangement_constant_local}. Hence, the strict inequality \eqref{eq:H1_space_norm_b_minus_barb_strictly_lessthan_nobubbling_and_rearrangement_constant_local} holds and $\bar\sT_{U'} > T$. This concludes Step \ref{step:Maximal_lifetime_barA_geq_T_local}.
\end{step}


By virtue of Step \ref{step:Maximal_lifetime_barA_geq_T_local}, the condition \eqref{eq:Linfty_time_H1_space_b_minus_barb_bound_allowing_rearrangement_quadratic_and_cubic_terms_cut-off} and \apriori estimate \eqref{eq:Apriori_estimate_for_b_minus_barb_in_terms_C1_and_W24_norms_g_minus_barg_cut-off} hold with $T_1 = T$. The \apriori estimate \eqref{eq:Apriori_estimate_for_A_minus_barA_in_terms_C1_and_W24_norms_g_minus_barg_local} now follows from \eqref{eq:Apriori_estimate_for_b_minus_barb_in_terms_C1_and_W24_norms_g_minus_barg_cut-off}, after setting $T_1 = T$, relabeling the constant $C_2$, and recalling that $b(t) - \bar b(t) = e^{-\mu_0 t}(A(t) - \bar A(t))$. This completes the proof of Theorem \ref{thm:Local_apriori_estimate_difference_solutions_Yang-Mills_heat_equations_pair_metrics}.
\end{proof}

Just as in Corollary \ref{cor:Global_apriori_estimate_difference_solutions_YM_heat_eqns_pair_metrics_C0_close}, the hypotheses \eqref{eq:Riem_barg_C0_norm_lessthan_bound}, \eqref{eq:Riemannian_metrics_g_minus_barg_C1_plus_W24_norm_lessthan_small_positive_constant}, and \eqref{eq:Riemannian_metrics_g_minus_barg_C2_lessthan_bound} on $\bar g$
in Theorem \ref{thm:Global_apriori_estimate_difference_solutions_Yang-Mills_heat_equations_pair_metrics} can be relaxed to give the following corollary of the proof of Theorem \ref{thm:Local_apriori_estimate_difference_solutions_Yang-Mills_heat_equations_pair_metrics}

\begin{cor}[Local \apriori estimate for the difference between solutions to the Yang-Mills heat equations defined by two weakly close Riemannian metrics and common initial data]
\label{cor:Local_apriori_estimate_difference_solutions_YM_heat_equations_pair_metrics_C0_close}
Assume the hypotheses of Theorem \ref{thm:Local_apriori_estimate_difference_solutions_Yang-Mills_heat_equations_pair_metrics}, except replace the hypotheses \eqref{eq:Riem_barg_C0_norm_lessthan_bound}, \eqref{eq:Riemannian_metrics_g_minus_barg_C1_plus_W24_norm_lessthan_small_positive_constant}, and \eqref{eq:Riemannian_metrics_g_minus_barg_C2_lessthan_bound} on $\bar g$ by \eqref{eq:Riem_barg_L4_norm_lessthan_bound} and \eqref{eq:Riemannian_metrics_g_minus_barg_C0_plus_W24_norm_lessthan_small_positive_constant}.
Then the conclusions of Theorem \ref{thm:Local_apriori_estimate_difference_solutions_Yang-Mills_heat_equations_pair_metrics} continue to hold.
\end{cor}

\subsection{Bubbling for a given Riemannian metric implies bubbling for a nearby Riemannian metric}
\label{subsec:Bubbling_for_metric_g_implies_bubbling_for_metric_barg}
We have seen in Section \ref{subsec:Local_comparison_estimates_two_solutions_Yang-Mills_heat_equations_for_two_metrics} that if a solution, $A(t)$, to the Yang-Mills heat equation for a Riemannian metric, $g$, and initial data, $A_0$, \emph{does not} bubble in an open subset $U \subset X$ for $t\in [0, T]$, then a solution, $\bar A(t)$, to the Yang-Mills heat equation for a nearby Riemannian metric, $\bar g$, and initial data, $A_0$, does \emph{not} bubble in an open subset $U' \Subset U$ for $t\in [0, T]$. In this section, we prove that if $A(t)$ \emph{does} bubble at \emph{every} point in a set $\Sigma = \{x_1,\ldots,x_L\}\subset X$ and $\bar g$ is sufficiently close to $g$, then also $\bar A(t)$ \emph{does} bubble in a small neighborhood of \emph{at least one} point in $\Sigma$. More precisely, we have the

\begin{thm}[Bubble singularities for a given metric imply at least one bubble singularity for a close metric]
\label{thm:Bubbling_for_g_implies_bubbling_for_barg}
Let $G$ be a compact Lie group, $P$ a principal $G$-bundle over a closed, four-dimensional, smooth manifold and $A_1$ a fixed reference connection of class $C^\infty$ on $P$, and $K$ and $T$ and $\delta$ positive constants. Let $g$ be a Riemannian metric on $X$ such that
\begin{equation}
\label{eq:Riem_g_C0_norm_lessthan_bound}
\|\Riem(g)\|_{C(X)} \leq K.
\end{equation}
Then there are a large enough constant $\mu_0 = \mu_0(A_1,g,K) \in [1, \infty)$ and small enough constants $\eta = \eta(A_1,g,K,T) \in (0, 1]$ and $\eps = \eps(A_1,g,K,T) \in (0, 1]$ with the following significance. Given $\eta \in (0, 1]$ and a connection, $A_0$, of class $C^\infty$
on $P$, there is a small enough constant $r = r(A_0,g,T,\delta,\eta) \in (0, \delta \wedge \Inj(X,g)]$ such that the following holds. Let $\Sigma := \{x_1,\ldots,x_L\}\subset X$ with $L \geq 1$ and $\dist_g(x_i,x_j) \geq 2\delta$ for all $i\neq j$ and define $U := X \less \cup_{l=1}^L \bar B_r(x_l)$, where $B_r(x_l) := \{x\in X:\dist_g(x,x_l) < r\}$. Let $U' \Subset U$ and $\Omega := U \less \bar U'$ satisfy the conditions \eqref{eq:Volume_Omega_and_cut-off_function_H2_norm_leq_epsilon} for some cut-off function, $\chi \in C^\infty_0(U)$.

Let $\bar g$ be a Riemannian metric on $X$ obeying \eqref{eq:Riem_barg_C0_norm_lessthan_bound}, \eqref{eq:Riemannian_metrics_g_minus_barg_C1_plus_W24_norm_lessthan_small_positive_constant}, and \eqref{eq:Riemannian_metrics_g_minus_barg_C2_lessthan_bound}. Let $A(t)$ and $\bar A(t)$, for $t\in [0, T)$, be strong solutions to the Yang-Mills heat equation
\eqref{eq:Yang-Mills_heat_equation_as_perturbation_rough_Laplacian_heat_equation} on $P \restriction U$ for the Riemannian metrics, $g$ and $\bar g$, respectively, and common initial data, $A_0$. Assume that $A_1$ obeys \eqref{eq:F_A1_C0_norm_lessthan_bound}, $A(t)$ and $\bar A(t)$ obey \eqref{eq:Yang-Mills_gradient_flow_A_minus_A1_and_barA_minus_A_bounds_local}. Let $\eps_1 = \eps_1(g) \in (0, 1]$ and $R_0 = R_0(g) \in (0, \Inj(X,g)]$ denote the positive constants in Section \ref{subsec:Kozono_Maeda_Naito_4}, Lemma \ref{lem:Schlatter_2-4_and_Struwe_3-6}, and Theorem \ref{thm:Kozono_Maeda_Naito_5-1}.

If $A(t)$ develops a bubble singularity as $t \nearrow T$ at every point in the set $\Sigma$ in the sense that
\begin{equation}
\label{eq:Theorem_Kozono_Maeda_Naito_5-1-1_bubble_characterization}
\limsup_{t \nearrow T} \int_{B_\varrho(x_l)} |F_A(t)|^2 \, d\vol_g \geq \eps_1,
\quad \forall\, \varrho \in (0, R_0] \hbox{ and } 1 \leq l \leq L,
\end{equation}
then $\bar A(t)$ develops a bubble singularity in \emph{at least one} ball $B_r(x_l)$, for \emph{some} $l \in \{1, \ldots, L\}$.
\end{thm}

\begin{proof}
We shall argue by contradiction and suppose that $\bar A(t)$ does not develop a bubble singularity in \emph{any} $g$-ball, $B_r(x_l)$, for $1 \leq l \leq L$ and $t \in [0, T]$. Consequently, there is a small enough positive constant, $\bar r_0 = \bar r_0(A_0,g,\bar g,L,T,\delta) = \bar r_0(A_0,g,L,T,\delta,\eta) \leq \delta\wedge\Inj(X,g)$, such that
\begin{equation}
\label{eq:Energy_barA_over_small_bubble_ball_is_small}
\int_{B_\varrho(x_l)} |F_{\bar A}(t)|^2 \, d\vol_g \leq \frac{\bar\eps_1}{8L},
\quad\forall\, \varrho \in (0, \bar r_0] \hbox{ and } t \in [0, T] \hbox{ and } 1 \leq l \leq L,
\end{equation}
where $\bar\eps_1 = \bar\eps_1(\bar g) = \bar\eps_1(g,K)$ is the positive constant in Lemma \ref{lem:Schlatter_2-4_and_Struwe_3-6} for the flow, $\bar A(t)$, and Riemannian metric, $\bar g$. For our definition of $U = X \less \bigcup_{l=1}^L \bar B_r(x_l)$, we now fix $r \in (0, \bar r_0 \wedge R_0]$ in the hypotheses by setting $r = \bar r_0 \wedge R_0$.

By combining our hypothesis \eqref{eq:Yang-Mills_gradient_flow_barA_minus_A1_annulus_bound} on $\bar A(t)$ over $\Omega$ and \eqref{eq:F_A1_C0_norm_lessthan_bound} for $A_1$ over $X$, we also have
$$
\int_\Omega |F_{\bar A}(t)|^2 \, d\vol_g \leq cK^2\eps, \quad\forall\, t \in [0, T],
$$
for some positive constant, $c$, depending at most on the Riemannian metric, $g$. We may choose $\eps \in (0, 1]$ in the hypotheses small enough that $cK^2\eps \leq \eps_1/4$ and so
\begin{equation}
\label{eq:Energy_barA_over_annulus_is_small}
\int_\Omega |F_{\bar A}(t)|^2 \, d\vol_g \leq \frac{\bar\eps_1}{8}, \quad\forall\, t \in [0, T].
\end{equation}
Theorem \ref{thm:Local_apriori_estimate_difference_solutions_Yang-Mills_heat_equations_pair_metrics} implies that $\bar A - A_1 \in C([0, T]; H_{A_1}^1(U'; \Lambda^1\otimes\ad P))$. More specifically, the \apriori estimate \eqref{eq:Apriori_estimate_for_A_minus_barA_in_terms_C1_and_W24_norms_g_minus_barg_local} for $A-\bar A$ coupled with the hypothesis
\eqref{eq:Yang-Mills_gradient_flow_A_minus_A1_Linfty_time_Linfty_space_plus_L4nabla_space_plus_L2nabla2_space_leq_bound_local} on $A$ yield a positive constant, $C_0 = C_0(A_1,g,K,T)$, such that
$$
\|F_{\bar A}\|_{L^\infty(0,T;L^4(U'))} \leq C_0,
$$
and hence, for positive constant, $C_1 = C_1(A_1,g,K,T)$,
$$
\int_{B_\varrho(x)\cap U'} |F_{\bar A}(t)|^2 \, d\vol_g \leq C_1\varrho^2,
\quad\forall\, x \in X \hbox{ and } \varrho \in (0, \Inj(X,g)].
$$
Hence, choosing $\bar r_1 = \bar r_1(A_1,g,\bar g,K,T) = \bar r_1(A_1,g,K,T)\in (0, \Inj(X,g)]$ small enough that $C_1\bar r_1^2 \leq \bar\eps_1/8$ yields
\begin{equation}
\label{eq:Energy_barA_over_small_ball_intersect_U_is_small}
\int_{B_\varrho(x)\cap U'} |F_{\bar A}(t)|^2 \, d\vol_g \leq \frac{\bar\eps_1}{8},
\quad\forall\, x \in X \hbox{ and } \varrho \in (0, \bar r_1] \hbox{ and } t \in [0, T].
\end{equation}
Therefore, by writing
\begin{multline*}
\int_{B_\varrho(x)} |F_{\bar A}(t)|^2 \, d\vol_g
=
\int_{B_\varrho(x)\cap U'} |F_{\bar A}(t)|^2 \, d\vol_g
+ \int_{B_\varrho(x)\cap \Omega} |F_{\bar A}(t)|^2 \, d\vol_g
\\
+ \sum_{l=1}^L\int_{B_\varrho(x)\cap B_r(x_l)} |F_{\bar A}(t)|^2 \, d\vol_g,
\end{multline*}
the combination of inequalities \eqref{eq:Energy_barA_over_small_bubble_ball_is_small}, \eqref{eq:Energy_barA_over_annulus_is_small}, and \eqref{eq:Energy_barA_over_small_ball_intersect_U_is_small} gives,
$$
\int_{B_\varrho(x)} |F_{\bar A}(t)|^2 \, d\vol_g \leq \frac{3\bar\eps_1}{8},
\quad\forall\, x \in X \hbox{ and } \varrho \in (0, \bar r_0\wedge\bar r_1\wedge R_0] \hbox{ and } t \in [0, T].
$$
But $\bar g$ is $C^1$-close to $g$ by hypothesis \eqref{eq:Riemannian_metrics_g_minus_barg_C1_plus_W24_norm_lessthan_small_positive_constant} and so the preceding inequality gives (for small enough $\eta \in (0, 1]$),
$$
\int_{B_\varrho(x;\bar g)} |F_{\bar A}(t)|^2 \, d\vol_{\bar g} \leq \frac{\bar\eps_1}{2},
\quad\forall\, x \in X \hbox{ and } 0 < \varrho \leq \frac{1}{2}\bar r_0\wedge\bar r_1\wedge R_0 \hbox{ and } t \in [0, T],
$$
where $B_\varrho(x;\bar g) := \{x\in X:\dist_{\bar g}(x,x_l) < \varrho\}$. Thus, Lemma \ref{lem:Schlatter_2-4_and_Struwe_3-6} implies that
\begin{equation}
\label{eq:barA_continuous_map_from_closed_interval_0_to_T_to_H1}
\bar A - A_1 \in C([0, T]; H_{A_1}^1(X; \Lambda^1\otimes\ad P)),
\end{equation}
and $\bar\sT > T$, where $\bar\sT$ is the maximal lifetime of $\bar A$ over $X$. Our regularity theory (Section \ref{subsec:Struwe_page_137_contraction_mapping_regularity_initial_data_in_H1}) for solutions to the Yang-Mills gradient flow ensures that $\bar A$ obeys the following analogue of \eqref{eq:Yang-Mills_gradient_flow_A_minus_A1_Linfty_time_H2_space_lessthan_bound}, when $A(0)-A_1 \in H_{A_1}^2(X;\Lambda^1\otimes\ad P)$, namely\footnote{One can either assume initial data, $A_0$, of class at least $H_{A_1}^2(X)$ or (we assume $A_0$ is $C^\infty$, which is unnecessary), without loss of generality, restrict attention to an interval $(t_0, T)$, for some $t_0 \in (0, T)$, and take advantage of the smoothing property of the Yang-Mills heat equation in order to achieve this bound.}
$$
\|\bar A - A_1\|_{L^\infty(0,T;H_{A_1,\bar g}^2(X,\bar g))} \leq \bar K,
$$
for some finite, positive constant $\bar K$. But then Theorem \ref{thm:Global_apriori_estimate_difference_solutions_Yang-Mills_heat_equations_pair_metrics}, by interchanging the roles of $(A,g)$ and $(\bar A, \bar g)$ --- which is possible because of the symmetry between the roles of $g$ and $\bar g$ by virtue of the pair of conditions \eqref{eq:Riem_barg_C0_norm_lessthan_bound} and \eqref{eq:Riem_g_C0_norm_lessthan_bound} and the conditions \eqref{eq:Riemannian_metrics_g_minus_barg_C1_plus_W24_norm_lessthan_small_positive_constant} and \eqref{eq:Riemannian_metrics_g_minus_barg_C2_lessthan_bound} --- ensures that
$$
A - A_1 \in C([0, T]; H_{A_1}^1(X; \Lambda^1\otimes\ad P)),
$$
and $A$ has maximal lifetime greater then $T$, contradicting our hypothesis that $A(t)$ bubbles at each of the points in $\Sigma$ as $t\nearrow T$.
\end{proof}

Just as in Corollaries \ref{cor:Global_apriori_estimate_difference_solutions_YM_heat_eqns_pair_metrics_C0_close} and \ref{cor:Local_apriori_estimate_difference_solutions_YM_heat_equations_pair_metrics_C0_close}, the hypotheses \eqref{eq:Riem_barg_C0_norm_lessthan_bound}, \eqref{eq:Riemannian_metrics_g_minus_barg_C1_plus_W24_norm_lessthan_small_positive_constant}, and \eqref{eq:Riemannian_metrics_g_minus_barg_C2_lessthan_bound} on $\bar g$
in Theorem \ref{thm:Global_apriori_estimate_difference_solutions_Yang-Mills_heat_equations_pair_metrics}
and \eqref{eq:Riem_g_C0_norm_lessthan_bound} on $g$ in Theorem \ref{thm:Bubbling_for_g_implies_bubbling_for_barg} can be relaxed to give the following corollary of the proof of Theorem \ref{thm:Bubbling_for_g_implies_bubbling_for_barg}.

\begin{cor}[Bubble singularities for a given metric imply at least one bubble singularity for a weakly close metric]
\label{cor:Bubbling_for_g_implies_bubbling_for_barg_C0_close}
Assume the hypotheses of Theorem \ref{thm:Bubbling_for_g_implies_bubbling_for_barg}, except replace the hypotheses \eqref{eq:Riem_barg_C0_norm_lessthan_bound}, \eqref{eq:Riemannian_metrics_g_minus_barg_C1_plus_W24_norm_lessthan_small_positive_constant}, \eqref{eq:Riemannian_metrics_g_minus_barg_C2_lessthan_bound} on $\bar g$ by \eqref{eq:Riem_barg_L4_norm_lessthan_bound} and \eqref{eq:Riemannian_metrics_g_minus_barg_C0_plus_W24_norm_lessthan_small_positive_constant} and replace the hypothesis \eqref{eq:Riem_g_C0_norm_lessthan_bound} on $g$ by
\begin{equation}
\label{eq:Riem_g_L4_norm_lessthan_bound}
\|\Riem(g)\|_{L^4(X)} \leq K.
\end{equation}
Then the conclusions of Theorem \ref{thm:Bubbling_for_g_implies_bubbling_for_barg} continue to hold.
\end{cor}

\subsection{Continuity and stability of Yang-Mills heat flow with respect to flattening a given Riemannian metric over finitely many small balls}
\label{subsec:Stability_bubbling_wrt_local_flattening_Riemannian_metric}
We now turn to applications of Theorems \ref{thm:Local_apriori_estimate_difference_solutions_Yang-Mills_heat_equations_pair_metrics} and \ref{thm:Bubbling_for_g_implies_bubbling_for_barg} to the case where the Riemannian metric, $\bar g$, is constructed explicitly by a `local flattening' procedure.

We begin by defining a Riemannian metric, $\bar g$, which is flat near each one of a finite set of points in $X$.

\begin{defn}[Riemannian metric flattened over a finite collection of small balls]
\label{defn:Locally_flattened_Riemannian_metric}
Let $X$ be a closed, smooth manifold with Riemannian metric $g$ and dimension $d \geq 2$ and $\Sigma := \{x_1,\ldots,x_L\} \subset X$ a finite set of points and $\rho \in (0, \Inj(X,g))$ a parameter, where $\Inj(X,g)$ denotes the injectivity radius of $(X,g)$. Let $f_l$ be an orthonormal frame for $(TX)_{x_l}$, for $1 \leq l \leq L$, and
$$
\varphi_l^{-1}: X \supset B_r(x_l) \cong B_r(0) \subset \RR^d
$$
the corresponding geodesic normal coordinate chart for any $r \in (0, \Inj(X,g))$, and $\delta$ the standard Euclidean metric on $\RR^d$. Let $\chi_0:\RR\to[0,1]$ be a $C^\infty$ cut-off function such that $\chi_0(r) = 1$ for $r\leq 1/2$ and $\chi_0(r) = 0$ for $r\geq 1$. The \emph{locally flattened Riemannian metric} corresponding to $g$ for the data $\rho$, $\{f_l\}_{l=1}^L$, and $\chi_0$ is defined by
\begin{equation}
\label{eq:Locally_flattened_Riemannian_metric}
\bar g
:=
\begin{cases}
g &\hbox{over } X \less \cup_{l=1}^L B_\rho(x_l),
\\
(1 - \chi_0(\dist_g(\cdot, x_l)/\rho))g &{}
\\
\quad + \chi_0(\dist_g(\cdot, x_l)/\rho)(\varphi_l^{-1})^*\delta
&\hbox{over } B_\rho(x_l) \less B_{\rho/2}(x_l), \quad 1 \leq l \leq L,
\\
(\varphi_l^{-1})^*\delta &\hbox{over } B_{\rho/2}(x_l), \quad 1 \leq l \leq L.
\end{cases}
\end{equation}
\end{defn}

\begin{rmk}[On the construction of a Riemannian metric flattened over a finite collection of small balls]
\label{rmk:Construction_locally_flattened_Riemannian_metric}
The set of positive definite, symmetric matrices is a convex cone and therefore the matrices, $(\bar g_{\mu\nu}(x))$ for $x \in B_\rho(x_l)$, defined by the expression \eqref{eq:Locally_flattened_Riemannian_metric} are positive definite and symmetric, so $\bar g$ is a well-defined Riemannian metric for all $\rho \in (0, \Inj(X,g))$.
\end{rmk}

A measure of closeness between $g$ and $\bar g$ is provided by the

\begin{lem}[Bound for the difference between the given Riemannian metric, $g$, and a locally flattened version]
\label{lem:Bound_norm_difference_g_minus_barg}
Let $X$ be a closed, smooth manifold with Riemannian metric $g$ and dimension $d \geq 2$. Then there is a positive constant, $c_0$, with the following significance. If $\rho \in (0,\Inj(X,g))$, and $\Sigma \subset X$, and $\bar g$ are as in Definition \ref{defn:Locally_flattened_Riemannian_metric} and $p \in [1,\infty)$, then
\begin{align}
\label{eq:Riemannian_metrics_g_minus_locally_flattened_g_C0_bound}
\|g - \bar g\|_{C(X)} &\leq c_0\rho^2,
\\
\label{eq:Riemannian_metrics_g_minus_locally_flattened_g_C1_bound}
\|g - \bar g\|_{C^1(X)} &\leq c_0\rho,
\\
\label{eq:Riemannian_metrics_g_minus_locally_flattened_g_W2p_norm_bound}
\|g - \bar g\|_{W^{2,p}(X)} &\leq c_0L\rho^{d/p},
\\
\label{eq:Riemannian_metrics_g_minus_locally_flattened_g_C2_bound}
\|g - \bar g\|_{C^2(X)} &\leq c_0,
\\
\label{eq:Riemannian_curvature_locally_flattened_g_C0_norm_bound}
\|\Riem(\bar g)\|_{C(X)} &\leq c_0.
\end{align}
\end{lem}

\begin{proof}
The estimates follow easily from the pointwise estimates for the components, $g_{\mu\nu}$, of the metric, $g$, with respect to the geodesic normal coordinate system, $\{x^\mu\}_{\mu=1}^d$, described in Section \ref{subsec:Estimate_gradient-like_flow_perturbation_due_to_non-flat_Riemannian_metric}.
\end{proof}

Lemma \ref{lem:Bound_norm_difference_g_minus_barg} provides an example of a Riemannian metric, $\bar g$, fulfilling the hypotheses of Theorems \ref{thm:Global_apriori_estimate_difference_solutions_Yang-Mills_heat_equations_pair_metrics} and \ref{thm:Local_apriori_estimate_difference_solutions_Yang-Mills_heat_equations_pair_metrics}.

\begin{cor}[Continuity and stability of Yang-Mills heat flow with respect to flattening of a given Riemannian metric over finitely many small balls]
\label{cor:Stability_bubbling_wrt_local_flattening_Riemannian_metric}
Let $G$ be a compact Lie group, $P$ a principal $G$-bundle over a closed, four-dimensional, smooth manifold with Riemannian metric $g$, and $A_1$ a fixed reference connection of class $C^\infty$ on $P$, and $K$ and $T$ positive constants and $L$ a positive integer. Then there are small enough constants $r = r(A_1,g,K,L,T) \in (0, \Inj(X,g)]$ and $\rho = \rho(A_1,g,K,L,T) \in (0, \Inj(X,g)]$, large enough constants $N = N(A_1,g,K,L,T) \in [4, \infty)$ and $\mu_0 = \mu_0(A_1,g,K) \in [1, \infty)$, and positive constants $z_1 = z_1(g,K)$ and $z_2=z_2(g)$ with the following significance. Let $\Sigma = \{x_1,\ldots,x_L\}\subset X$ and $\bar g$ be the Riemannian metric on $X$ constructed as in Definition \ref{defn:Locally_flattened_Riemannian_metric} corresponding to $g$, $\rho$, and $\Sigma$. Define
\begin{equation}
\label{eq:Definition_Uprime_U_Omega_in_terms_of_Sigma_points_and_N_rho}
\begin{gathered}
U \equiv U_{r/N} := X \less \bigcup_{l=1}^L \bar B_{r/N}(x_l),
\quad
U' \equiv U_{r/2} := X \less \bigcup_{l=1}^L \bar B_{r/2}(x_l),
\\
\hbox{and}\quad \Omega := U \less \bar U' = \bigcup_{l=1}^L \Omega(x_l;r/N,r/2),
\end{gathered}
\end{equation}
where $\Omega(x_l;r/N,r/2) := B_{r/2}(x_l) \less \bar B_{r/N}(x_l)$ for $1 \leq l \leq L$. Let $A(t)$, for $t\in [0, T)$, and $\bar A(t)$, for $t\in [0, \bar\sT_U)$ with maximal lifetime $\bar\sT_U \in (0, \infty]$, be strong solutions to the Yang-Mills heat equation \eqref{eq:Yang-Mills_heat_equation_as_perturbation_rough_Laplacian_heat_equation} on $P \restriction U$ for the Riemannian metrics, $g$ and $\bar g$, respectively, and common initial data, $A_0$. Assume that $A_1$ obeys \eqref{eq:F_A1_C0_norm_lessthan_bound} and $A(t)$ and $\bar A(t)$ obey \eqref{eq:Yang-Mills_gradient_flow_A_minus_A1_and_barA_minus_A_bounds_local}. Then $\bar A - A_1 \in C([0,T]; H_{A_1}^1(U';\Lambda^1\otimes\ad P))$ and $\bar\sT_{U'} > T$, where $\bar\sT_{U'}$ is the maximal lifetime of $\bar A$ on $P \restriction U'$, and $\bar A(t)$ obeys the \apriori estimate \eqref{eq:Apriori_estimate_for_A_minus_barA_in_terms_C1_and_W24_norms_g_minus_barg_local}.
\end{cor}

\begin{rmk}[On the choice of the radii $r$ and $\rho$ in Corollary \ref{cor:Stability_bubbling_wrt_local_flattening_Riemannian_metric}]
\label{rmk:Stability_bubbling_wrt_local_flattening_Riemannian_metric_choice_r_and_rho}
In our applications of Corollary \ref{cor:Stability_bubbling_wrt_local_flattening_Riemannian_metric}, we shall \emph{usually} choose $r \in (0, \rho]$, but that is not a requirement, as is evident from the proof.
\end{rmk}

\begin{proof}[Proof of Corollary \ref{cor:Stability_bubbling_wrt_local_flattening_Riemannian_metric}]
First, Lemma \ref{lem:Bound_norm_difference_g_minus_barg} assures us that $\bar g$ obeys the conditions \eqref{eq:Riem_barg_C0_norm_lessthan_bound} and \eqref{eq:Riemannian_metrics_g_minus_barg_C2_lessthan_bound} and, for $\rho \in (0, \Inj(X,g)]$ small enough that
$$
c_0(1+L)\rho \leq \eta,
$$
where $\eta \in (0, 1]$ is the constant in the hypotheses of Theorem \ref{thm:Local_apriori_estimate_difference_solutions_Yang-Mills_heat_equations_pair_metrics}, then $\bar g$ also obeys
\eqref{eq:Riemannian_metrics_g_minus_barg_C1_plus_W24_norm_lessthan_small_positive_constant}.

Second, we define a cut-off function, $\chi \equiv \chi_{N,r,\Sigma} \in C^\infty(X)$, following the recipe described in Remark \ref{rmk:Construction_cut-off_function_chi_N_rho} in a neighborhood of each point $x_l\in\Sigma$, so
$$
\chi
:=
\begin{cases}
1 &\hbox{on } U_{r/2},
\\
0 &\hbox{on } \cup_{l=1}^L B_{r/N}(x_l).
\end{cases}
$$
We observe that $\supp\chi \subset \bar U_{r/N}$ and (also for future reference),
$$
\supp\nabla\chi \subset \bigcup_{l=1}^L \bar\Omega(x_l;r/N,r/2).
$$
Then Lemma \ref{lem:L2_and_L4_nabla_chi_and_L2_nabla2_chi_bounds} assures us that $\chi, U', U$ have the properties specified in \eqref{eq:Volume_Omega_and_cut-off_function_H2_norm_leq_epsilon} by choosing $r \in (0, \Inj(X,g)]$ small enough to ensure that
$$
\Vol_g(\Omega) \leq \sum_{l=1}^L \Vol_g(B_{r/2}(x_l)) \leq cLr^4 \leq \eps,
$$
where $c$ is a positive constant depending at most on the Riemannian metric, $g$, and $\eps \in (0, 1]$ is the constant in the hypotheses of Theorem \ref{thm:Local_apriori_estimate_difference_solutions_Yang-Mills_heat_equations_pair_metrics}, and $N \geq 4$ is a large enough constant (independent of $r$) such that
$$
\|\nabla\chi\|_{L^4(X)} + \|\nabla^2\chi\|_{L^2(X)}
\leq cL\left( (\log N)^{-3/4} + (\log N)^{-1/2} \right) \leq \eps.
$$
The conclusions now follow from Theorem \ref{thm:Local_apriori_estimate_difference_solutions_Yang-Mills_heat_equations_pair_metrics}.
\end{proof}

Similarly, combining Corollary \ref{cor:Stability_bubbling_wrt_local_flattening_Riemannian_metric} and Theorem \ref{thm:Bubbling_for_g_implies_bubbling_for_barg} leads to the following useful

\begin{cor}[Bubble singularities for a given metric imply at least one bubble singularity for a nearby metric flattened over finitely many small balls]
\label{cor:Bubbling_for_g_implies_bubbling_for_nearby_barg_locally_flattened_finite_number_small_balls}
Assume the hypotheses of Corollary \ref{cor:Stability_bubbling_wrt_local_flattening_Riemannian_metric} and, in addition, that the constant, $r = (A_0,A_1,g,K,L,T,\delta,\rho) \in (0, \delta\wedge\Inj(X,g)]$, is chosen small enough to also satisfy the hypotheses of Theorem \ref{thm:Bubbling_for_g_implies_bubbling_for_barg}. If $A(t)$ develops a bubble singularity as $t \nearrow T$ at every point in the set $\Sigma$ in the sense of \eqref{eq:Theorem_Kozono_Maeda_Naito_5-1-1_bubble_characterization}, then $\bar A(t)$ develops a bubble singularity in \emph{at least one} ball $B_{r/N}(x_l)$, for \emph{some} $l \in \{1, \ldots, L\}$.
\end{cor}

\begin{proof}
We recall from the hypothesis \eqref{eq:Definition_Uprime_U_Omega_in_terms_of_Sigma_points_and_N_rho} of Corollary \ref{cor:Stability_bubbling_wrt_local_flattening_Riemannian_metric} that
$$
U := X \less \bigcup_{l=1}^L \bar B_{r/N}(x_l), \quad U' := X \less \bigcup_{l=1}^L \bar B_{r/2}(x_l),
\quad\hbox{and}\quad \Omega := \bigcup_{l=1}^L B_{r/2}(x_l) \less \bar B_{r/N}(x_l).
$$
It now suffices to observe that, exactly as in the proof of Corollary \ref{cor:Stability_bubbling_wrt_local_flattening_Riemannian_metric}, we can make the constants $\eps \in (0, 1]$ and $\eta \in (0, 1]$ sufficiently small in the hypotheses of Theorem \ref{thm:Bubbling_for_g_implies_bubbling_for_barg} by choosing $\rho \in (0, \Inj(X,g)]$ sufficiently small and $N \in [4, \infty)$ sufficiently large. The conclusion thus follows from Theorem \ref{thm:Bubbling_for_g_implies_bubbling_for_barg}, for a constant $\rho = \rho(A_1,g,K,L,T) \in (0, \Inj(X,g)]$ chosen small enough and constants $N = N(A_1,g,K,L,T) \in [4, \infty)$ and $\mu_0 = \mu_0(A_1,g,K) \in [1, \infty)$ chosen large enough to ensure that the hypotheses of Theorem \ref{thm:Bubbling_for_g_implies_bubbling_for_barg} are obeyed with small enough constants $\eps \in (0, 1]$ and $\eta \in (0, 1]$, exactly as in the proof of Corollary \ref{cor:Stability_bubbling_wrt_local_flattening_Riemannian_metric}.
\end{proof}


\subsection{Continuity and stability of Yang-Mills heat flow with respect to weakly flattening a given Riemannian metric over finitely many small balls}
\label{subsec:Stability_bubbling_wrt_local_flattening_Hermitian_metric}
In this section, we rederive the results of Section \ref{subsec:Stability_bubbling_wrt_local_flattening_Riemannian_metric} using a weaker notion of `locally flattened' Riemannian metric, one that is appropriate for the case of a Riemannian metric induced by a Hermitian, non-K\"ahler metric on a complex manifold.

Lemma \ref{lem:Bound_norm_difference_g_minus_barg} allows us, given a Riemannian metric, $g$, on $X$, to construct a Riemannian metric, $\bar g$, on $X$ where $\|g-\bar g\|_{C^1(X,g)}$ can be made as small as desired while $\|g-\bar g\|_{C^2(X,g)}$ remains uniformly bounded. However, that splicing construction relies on (real) geodesic normal coordinates, which are only available for real Riemannian manifolds. For complex manifolds with a Hermitian metric, the analogous assumption regarding properties of local coordinates would imply that the Hermitian metric is K\"ahler, the very hypothesis we aim to relax. There, we shall introduce a weaker version of Lemma \ref{lem:Bound_norm_difference_g_minus_barg} that will nonetheless suffice for our application to Hermitian, but non-K\"ahler complex manifolds. Note that Lemma \ref{lem:Bound_norm_difference_g_minus_barg} applies to a smooth manifold of any (real) dimension $d \geq 2$ whereas Lemma \ref{lem:Bound_norm_difference_g_minus_barg_Hermitian} is specific to (real) dimension four.

\begin{lem}[Bound for the difference between the given Riemannian metric, $g$, and a weakly locally flattened version]
\label{lem:Bound_norm_difference_g_minus_barg_Hermitian}
Let $X$ be a closed, four-dimensional, smooth manifold with Riemannian metric $g$. Then there are a positive constant, $c$, and, if $\alpha\in[0,1)$, a positive constant, $c_\alpha$ with the following significance. Let $\rho \in (0,\Inj(X,g))$, and $\Sigma \subset X$, and $\bar g$ be as in Definition \ref{defn:Locally_flattened_Riemannian_metric}, except that
\begin{enumerate}
  \item The local coordinates, $\{x^\mu\}$, around each point $x_i \in \Sigma$ are only assumed to obey $g_{\mu\nu}(x_i) = \delta_{\mu\nu}$, and are not assumed to be geodesic normal coordinates;

  \item The cut-off function $\chi_0$ and annuli $\Omega(x_i;\rho/2,\rho)$ are replaced by the cut-off function $\chi = 1-\chi_{N,\rho}$, where $N \geq 4$ and $\chi_{N,\rho}$ is as in Remark \ref{rmk:Construction_cut-off_function_chi_N_rho}, and $\Omega(x_i;\rho/N,\rho/2)$, respectively, for each $x_i \in \Sigma$.
\end{enumerate}
Then the following hold:
\begin{align}
\label{eq:Riemannian_metrics_g_minus_locally_flattened_g_Calpha_bound_Hermitian}
\|g - \bar g\|_{C^\alpha(X,g)} &\leq c_\alpha\left(\rho + (\log N)^{-3/4}\right),
\\
\label{eq:Riemannian_metrics_g_minus_locally_flattened_g_W24_norm_Hermitian}
\|g - \bar g\|_{W^{2,4}(X,g)} &\leq c\left(\rho + (\log N)^{-3/4}\right),
\\
\label{eq:Riemannian_curvature_locally_flattened_g_L4_norm_bound_Hermitian}
\|\Riem(\bar g)\|_{L^4(X,g)} &\leq c.
\end{align}
\end{lem}

\begin{proof}
It is enough to consider one of the points, $x_0 \in \Sigma$. We observe that, for any real local coordinates, $\{x^\mu\}$, such that $g_{\mu\nu}(x_0) = \delta_{\mu\nu}$, we have
$$
|g_{\mu\nu}(x) - \delta_{\mu\nu}| \leq c_1(g)\dist_g(x,x_0),
$$
by Taylor's formula, where $c_1(g)$ is a positive constant depending at most on
$$
\max_{\mu,\nu,\lambda}
\left( \|g_{\mu\nu} - \delta_{\mu\nu}\|_{C(B_\rho(x_0),g)}
+ \left\|\frac{\partial g_{\mu\nu}}{\partial x^\lambda}\right\|_{C(B_\rho(x_0),g)}\right).
$$
We interpolate between $g$ on $X\less B_\rho(x_0)$ and $\delta$ on $B_{\rho/N}(x_0)$ just as in Definition \ref{defn:Locally_flattened_Riemannian_metric}, but with $\rho/2$ replaced by $\rho/N$ and $\rho$ replaced by $\rho$). We identify $B_\rho(x_0) \subset X$ with $B_\rho(0) \subset \RR^4$ and employ the cut-off function $\chi = 1-\chi_{N,\rho}$ provided by Remark \ref{rmk:Construction_cut-off_function_chi_N_rho}, so $\chi \in C_0^\infty(X)$ with
$$
\chi
=
\begin{cases}
1 &\hbox{on } B_{\rho/N}(x_0),
\\
0 &\hbox{on } \Omega(x_0;\rho/2,\rho).
\end{cases}
$$
Then, on $B_{\rho/2}(x_0)$ we have
\begin{align*}
\bar g_{\mu\nu} &= (1-\chi)g_{\mu\nu} + \chi\delta_{\mu\nu}
\\
&= g_{\mu\nu} + \chi\left(\delta_{\mu\nu} - g_{\mu\nu}\right),
\\
\frac{\partial\bar g_{\mu\nu}}{\partial x^\lambda}
&=
(1-\chi)\frac{\partial g_{\mu\nu}}{\partial x^\lambda}
+
\frac{\partial \chi}{\partial x^\lambda}(\delta_{\mu\nu} - g_{\mu\nu}),
\\
\frac{\partial^2\bar g_{\mu\nu}}{\partial x^\lambda\partial x^\xi}
&=
(1-\chi)\frac{\partial^2 g_{\mu\nu}}{\partial x^\lambda\partial x^\xi}
-
\frac{\partial \chi}{\partial x^\xi}\frac{\partial g_{\mu\nu}}{\partial x^\lambda}
-
\frac{\partial \chi}{\partial x^\lambda}\frac{\partial g_{\mu\nu}}{\partial x^\xi}
+
\frac{\partial^2 \chi}{\partial x^\lambda\partial x^\xi}(\delta_{\mu\nu} - g_{\mu\nu}),
\end{align*}
and although we no longer assume that $(\partial g_{\mu\nu}/\partial x^\lambda)(x_0) = 0$, we still have
$$
\bar g_{\mu\nu} = g_{\mu\nu} \quad\hbox{on } \bar\Omega(x_0;\rho/2,\rho).
$$
We recall from Equation \eqref{eq:Feehan_Leness_5-4} in Remark \ref{rmk:Construction_cut-off_function_chi_N_rho} that, denoting $r = \dist_g(x,x_0)$,
$$
|\nabla\chi_{N,\rho}(x)| \leq \frac{c_0}{r\log N}
\quad\hbox{and}\quad
|\nabla^2\chi_{N,\rho}(x)| \leq \frac{c_0}{r^2\log N}, \quad\forall\, x \in B_{\rho}(x_0),
$$
where $c_0(g)$ is a positive constant depending at most on
$$
\max_{\mu,\nu}\|g_{\mu\nu} - \delta_{\mu\nu}\|_{C(B_\rho(x_0),g)},
$$
and, unless a different Riemannian metric is explicitly noted, we use
$$
B_{\rho}(x_0) = \left\{x\in X: \dist_g(x,x_0) < \rho\right\},
$$
with $\rho \in (0, \Inj(X,g)]$. Consequently, Taylor's formula and Remark \ref{rmk:Construction_cut-off_function_chi_N_rho} yield the following pointwise estimates, for $r = \dist_g(x,x_0)$ and all $x \in B_\rho(x_0)$,
\begin{align*}
|\bar g_{\mu\nu} - g_{\mu\nu}| &\leq c_1(g)r\,1_{\{r\leq\rho/2\}},
\\
\left|\frac{\partial\bar g_{\mu\nu}}{\partial x^\lambda}(x)
- \frac{\partial g_{\mu\nu}}{\partial x^\lambda}(x)\right|
&\leq
\chi\left\|\frac{\partial g_{\mu\nu}}{\partial x^\lambda}\right\|_{C(B_\rho(x_0),g)}
+  \frac{c_0(g)}{\log N}1_{\{\rho/N\leq r\leq\rho/2\}},
\\
\left|\frac{\partial^2\bar g_{\mu\nu}}{\partial x^\lambda\partial x^\xi}(x)
- \frac{\partial^2 g_{\mu\nu}}{\partial x^\lambda\partial x^\xi}(x)\right|
&\leq
\chi\left\|\frac{\partial^2 g_{\mu\nu}}{\partial x^\lambda\partial x^\xi} \right\|_{C(B_\rho(x_0),g)}
\\
&\quad + \frac{c_0(g)}{r\log N}1_{\{\rho/N\leq r\leq\rho/2\}}\sum_{\lambda=1}^4
\left(1+\left\|\frac{\partial g_{\mu\nu}}{\partial x^\lambda}\right\|_{C(B_\rho(x_0),g)}\right),
\end{align*}
where $c_0(g)$ and $c_1(g)$ are as above. Thus,
\begin{align*}
\|\bar g_{\mu\nu} - g_{\mu\nu}\|_{C(B_\rho(x_0),g)}
&\leq c_1(g)\rho,
\\
\|\bar g_{\mu\nu} - g_{\mu\nu}\|_{L^q(B_\rho(x_0),g)}
&\leq
c_1(g)\rho^{1+4/q}, \quad\forall\, q\geq 4,
\\
\left\|\frac{\partial\bar g_{\mu\nu}}{\partial x^\lambda}
- \frac{\partial g_{\mu\nu}}{\partial x^\lambda}
\right\|_{L^q(B_\rho(x_0),g)}
&\leq
c_1(g)\rho^{4/q}\left(1 + (\log N)^{-1}\right), \quad\forall\, q\geq 4,
\\
\left\|\frac{\partial^2\bar g_{\mu\nu}}{\partial x^\lambda\partial x^\xi}
- \frac{\partial^2 g_{\mu\nu}}{\partial x^\lambda\partial x^\xi}
\right\|_{L^4(B_\rho(x_0),g)}
&\leq
c_2(g)\left(\rho + (\log N)^{-3/4}\right),
\end{align*}
where $c_2(g)$ is a positive constant depending at most on
$$
\max_{\mu,\nu,\lambda\xi}
\left( \|g_{\mu\nu} - \delta_{\mu\nu}\|_{C(B_\rho(x_0),g)}
+ \left\|\frac{\partial g_{\mu\nu}}{\partial x^\lambda}\right\|_{C(B_\rho(x_0),g)}
+ \left\|\frac{\partial^2 g_{\mu\nu}}{\partial x^\lambda\partial x^\xi}
\right\|_{C(B_\rho(x_0),g)}\right).
$$
The estimate for the second-order partial derivatives of $g_{\mu\nu} - \bar g_{\mu\nu}$ is obtained using
\begin{align*}
\left(\int_{\Omega(x_0;\rho/N,\rho/2)} \frac{c_1^4}{r^4(\log N)^4} d\vol_g\right)^{1/4}
&\leq
\frac{c_1}{\log N} \left(\int_{S^3}\int_{\rho/N}^{\rho/2} \frac{1}{r}dr\,d\theta\right)^{1/4}
\\
&= c_1\Vol(S^3)\frac{\left(\log N - \log 2\right)^{1/4}}{\log N}.
\end{align*}
Consequently, using $W^{k,p}(X,g)$ to denote the Sobolev space defined by the Riemannian metric, $g$, on $X$ and Levi-Civita connection for $g$ on $TX$ and associated vector bundles, we have
$$
\|g - \bar g\|_{W^{2,4}(X,g)} \leq c(g)\left(\rho + (\log N)^{-3/4}\right),
$$
for all $\rho \in (0, 1\wedge \Inj(X,g)]$ and $N\geq 4$, which gives \eqref{eq:Riemannian_metrics_g_minus_locally_flattened_g_W24_norm_Hermitian}. The Sobolev Embedding Theorem \cite[Theorem 4.12]{AdamsFournier} with $n=4$, $j=0$, $m=2$, and $p=4$ gives $W^{2,4}(X,g)\hookrightarrow C^\alpha(X,g)$, for any $\alpha \in [0,1)$, and thus
$$
\|g - \bar g\|_{C^\alpha(X,g)} \leq c_\alpha(g)\left(\rho + (\log N)^{-3/4}\right),
$$
which gives \eqref{eq:Riemannian_metrics_g_minus_locally_flattened_g_Calpha_bound_Hermitian}. The bound \eqref{eq:Riemannian_curvature_locally_flattened_g_L4_norm_bound_Hermitian} follows from \eqref{eq:Riemannian_metrics_g_minus_locally_flattened_g_W24_norm_Hermitian}, noting that $c(g)$ is allowed to depend on $\|\Riem(g)\|_{C^0(X)}$. This completes the proof of Lemma \ref{lem:Bound_norm_difference_g_minus_barg_Hermitian}.
\end{proof}

Lemma \ref{lem:Bound_norm_difference_g_minus_barg_Hermitian} provides an example of a Riemannian metric, $\bar g$, fulfilling the hypotheses of Corollaries \ref{cor:Global_apriori_estimate_difference_solutions_YM_heat_eqns_pair_metrics_C0_close} and \ref{cor:Local_apriori_estimate_difference_solutions_YM_heat_equations_pair_metrics_C0_close}. Consequently, we obtain the following analogues of Corollaries \ref{cor:Stability_bubbling_wrt_local_flattening_Riemannian_metric} and \ref{cor:Bubbling_for_g_implies_bubbling_for_nearby_barg_locally_flattened_finite_number_small_balls}.

\begin{cor}[Continuity and stability of Yang-Mills heat flow with respect to weakly flattening a given Riemannian metric over finitely many small balls]
\label{cor:Stability_bubbling_wrt_local_flattening_Riemannian_metric_C0_close}
Assume the hypotheses of Corollary \ref{cor:Stability_bubbling_wrt_local_flattening_Riemannian_metric}, except that $\bar g$ is the Riemannian metric on $X$ constructed as in Lemma \ref{lem:Bound_norm_difference_g_minus_barg_Hermitian} corresponding to $g$, $\rho$, $\Sigma$, and a constant $N_1 \geq 4$. With the additional hypothesis of a large enough constant $N_1=N_1(A_1,g,K,L,T) \in [4,\infty)$, the conclusions of Corollary \ref{cor:Stability_bubbling_wrt_local_flattening_Riemannian_metric} continue to hold.
\end{cor}

\begin{cor}[Bubble singularities for a given metric imply at least one bubble singularity for a weakly close metric flattened over finitely many small balls]
\label{cor:Bubbling_for_g_implies_bubbling_for_nearby_barg_weakly_locally_flattened_finite_number_small_balls}
Assume the hypotheses of Corollary \ref{cor:Stability_bubbling_wrt_local_flattening_Riemannian_metric_C0_close} and, in addition, that the constant, $r = (A_0,A_1,g,K,L,T,\delta,\rho) \in (0, \delta\wedge\Inj(X,g)]$, is chosen small enough to also satisfy the hypotheses of Theorem \ref{thm:Bubbling_for_g_implies_bubbling_for_barg} and Corollary \ref{cor:Bubbling_for_g_implies_bubbling_for_barg_C0_close}. Then the conclusions of Corollary \ref{cor:Stability_bubbling_wrt_local_flattening_Riemannian_metric_C0_close} continue to hold.
\end{cor}


\chapter[Yang-Mills gradient-like flow over four-dimensional manifolds]{Yang-Mills gradient-like flow over four-dimensional manifolds and applications}
\label{chapter:Yang-Mills_gradient-like_flow_four_manifolds_applications}
Our analysis thus far of the asymptotic behavior of Yang-Mills gradient flow, $A(t)$, over a four-dimensional manifold, $X$, suggests that it may be profitable to try isolate the balls, $B_\rho(x_i) \subset X$, where the energy density of the flow, $|F_A(t)|^2$, is becoming unbounded as $t \nearrow T$, where $T \in (0,\infty]$. Indeed, Theorems \ref{thm:Kozono_Maeda_Naito_5-1}, \ref{thm:Kozono_Maeda_Naito_5-3} (for $T<\infty)$, and \ref{thm:Kozono_Maeda_Naito_5-3_T_is_infinite} (for $T=\infty)$, Corollary \ref{cor:Kozono_Maeda_Naito_5-3_T_is_infinite}, and Theorem \ref{thm:Kozono_Maeda_Naito_5-4} provide precise information about the flow, both on the balls $B_\rho(x_i)$ and on their complement, $X \less \cup_{i=1}^L B_\rho(x_i)$. Thus, a natural strategy is to cut-off the flow over small annuli near each bubble point, $x_i$, and consider the resulting Yang-Mills gradient-\emph{like} flows on standard, conformally flat manifolds, such as $S^4$ when $X$ is a real four-dimensional manifold or $\CC\PP^1\times\CC\PP^1$ when $X$ is a complex surface, as well as the residual background flow over $X$. Of course, even if the Riemannian metrics match exactly, any such cutting off procedure introduces a cut-off function error over the annuli and the success or otherwise of this strategy depends on whether those errors are sufficiently small or, at the very least, bounded as $t\nearrow T$. In \cite[Theorem 2]{Simon_1983}, Simon assumes that the source term, $f \in C^\infty([0,\infty)\times X)$, in his gradient-like flow equation (namely, \cite[Equation (0.1)]{Simon_1983}) obeys
$$
\|f(t)\|_{H^l(X)} + \|\dot f(t)\|_{H^l(X)} + \|\ddot f(t)\|_{L^2(X)} \leq \delta e^{-\eps t}, \quad \forall\, t \in [0,\infty),
$$
for positive constants $\eps$ and $\delta$ and an integer $l$ chosen large enough to ensure that\footnote{In \cite[Equation (1.27)]{Simon_1983}, the Sobolev embedding should be reversed.}
$H^{l-1}(X) \hookrightarrow C^2(X)$. Unfortunately, it appears quite difficult to achieve this kind of exponential decay through a cutting off procedure.

Given the preceding caveat, we shall describe in this chapter some results that can obtained using the methods developed in this monograph, together with some applications. Section \ref{sec:Yang-Mills_gradient-like_flow_over_four_sphere} explores properties of Yang-Mills gradient-like flow over the four-dimensional sphere while Section \ref{sec:Yang-Mills_gradient-like_flow_Kahler_surface} develops more refined properties of Yang-Mills gradient-like flow over K\"ahler surfaces (such as $\CC\PP^1\times\CC\PP^1$ with its product Fubini-Study metric).

\section{Yang-Mills gradient-like flow over the four-dimensional sphere}
\label{sec:Yang-Mills_gradient-like_flow_over_four_sphere}
A basic construction underlying much of Taubes' work in gauge theory concerns cutting off and splicing connections \cite{TauSelfDual, TauPath, TauIndef, TauFrame, TauStable}. We describe such a construction here, using the notation and conventions of \cite[Section 3.3]{FLKM1}.

\begin{defn}[Construction of a principal $G$-bundle and a spliced connection over the sphere]
\label{defn:Construction_spliced_connection_over_sphere}
Let $P$ be a principal $G$-bundle over a closed, oriented, smooth manifold of dimension $d \geq 2$ with Riemannian metric, $g$. Assume we are given a point $x_0 \in X$; a positive constant, $\rho$, obeying
$$
\rho \leq \max\{1, \rho_0\},
$$
where $\rho_0 < r_0$ and $r_0 = r_0(x_0)$ is the injectivity radius of $(X, g)$ at $x_0$; an orthonormal frame $f_0$ for the tangent space $(TX)_{x_0}$ and corresponding geodesic normal coordinate chart,
$$
\varphi_0^{-1}: X \supset B_\rho(x_0) \cong B_\rho(0) \subset \RR^d \cong (TX)_{x_0},
$$
for the ball $B_\rho(x_0) \subset X$ and defined via the exponential map,
$$
\varphi_0 \equiv \exp_{f_0}: (TX)_{x_0} \supset B_\rho(0) \cong B_\rho(x_0) \subset X;
$$
an orthonormal frame $f_n$ for the tangent space $(TS^d)_n$ over the north pole $n \in S^d \cong \RR^d\cup\{\infty\}$ (identified with the origin in $\RR^d$) and corresponding coordinate chart,
$$
\varphi_n: \RR^d \cong S^d \less \{s\},
$$
that is inverse to a stereographic projection from the south pole $s \in S^d \subset \RR^{d+1}$ (identified with the point at infinity in $\RR^d\cup\{\infty\}$); a $C^\infty$ reference connection $A_1$ on $P$; a point $p_0$ in the fiber $P_{x_0}$ and a $W^{k+1, p}$ section,
$$
\sigma_0 \equiv \sigma_0(A_1, p_0): B_\rho(x_0) \to P \restriction B_\rho(x_0),
$$
defined by parallel transport via the connection, $A_1$, along geodesics emanating from $x_0$; a trivialization,
$$
\tau_0: P \restriction B_\rho(x_0) \cong B_\rho(x_0) \times G,
$$
defined by the section $\sigma_0$ via $\tau_0^{-1}(x, u) = \sigma_0(x)u$ for all $(x,u) \in B_\rho(x_0) \times G$; a constant $N \geq 4$ and a $C^\infty$ cut-off function $\beta = \beta_{N,\rho}:\RR \to [0, 1]$ such that $\beta(r) = 1$ for $r \leq \rho/N$ and $\beta(r) = 0$ for $r \geq \rho/2$. Define a principal $G$-bundle, $\widehat P$, over $S^d$ by setting
\begin{equation}
\label{eq:Construction_spliced_principle_G-bundle_over_sphere}
\widehat P
:=
\begin{cases}
(\varphi_0\circ \varphi_n^{-1})^*P &\hbox{on } B_\rho(n),
\\
\left(S^d \less B_{\rho/2}(n)\right) \times G &\hbox{on } S^d \less B_{\rho/2}(n).
\end{cases}
\end{equation}
If $A$ is a connection on $P$ of class $W^{k, p}$ for an integer $k\geq 1$ and $1 \leq p \leq \infty$, define a connection, $\widehat A$, on $\widehat P$ of class $W^{k, p}$ by setting
\begin{equation}
\label{eq:Construction_spliced_connection_over_sphere}
\widehat A
:=
\begin{cases}
(\varphi_0\circ \varphi_n^{-1})^*A &\hbox{over } \varphi_n\left(B_{\rho/N}(0)\right),
\\
\Gamma + (\varphi_0\circ \varphi_n^{-1})^*\left(\chi\sigma_0^*A\right) &\hbox{over } \varphi_n\left(B_{\rho/2}(0) \less B_{\rho/N}(0)\right),
\\
\Gamma &\hbox{over } \varphi_n\left(B_{\rho/2}(0)\right),
\end{cases}
\end{equation}
where $\chi:X \to [0, 1]$ is the $C^\infty$ cut-off function defined by
\begin{equation}
\label{eq:Definition_chi_in_terms_of_beta}
\chi = \chi_{N,\rho,g, x_0} := \beta_{N,\rho}(\dist_g(\cdot, x_0)),
\end{equation}
and $\Gamma$ is the product connection on $(S^d \less \{n\}) \times G$.
\end{defn}

\begin{rmk}[Topology of the principal $G$-bundle $\widehat P$]
\label{rmk:Sphere_principal_bundle_topology}
The Pontrjagin classes of the bundle $\widehat P$ can be determined from the connection $\widehat A$ in \eqref{eq:Construction_spliced_connection_over_sphere} (see, for example, \cite{UhlChern}) constructed in Definition \ref{defn:Construction_spliced_connection_over_sphere} or by the transition function for $\widehat P$ implied by the construction \eqref{eq:Construction_spliced_principle_G-bundle_over_sphere}. See Section \ref{sec:Taubes_1982_Appendix} for a discussion of the classification of principal $G$-bundles when $G$ is a compact, connected, semi-simple Lie group.
\end{rmk}

The splicing procedure prescribed in Definition \ref{defn:Construction_spliced_connection_over_sphere} can be applied to a solution, $A(t)$, to Yang-Mills gradient flow on a principal $G$-bundle $P$ over a Riemannian, smooth manifold to create a solution, $\widehat A(t)$, to Yang-Mills gradient-like flow on a principal $G$-bundle $\widehat P$ over $S^d$.

\begin{defn}[Splicing construction of a Yang-Mills gradient-like flow over the sphere]
\label{defn:Splicing_construction_Yang-Mills_gradientlike_flow_over_sphere}
Assume the notation and set-up of Definition \ref{defn:Construction_spliced_connection_over_sphere}. Let $A(t)$, for $t \in [0, T)$, be a solution\footnote{When the precise sense of what we mean by `solution' --- such as strong, classical, and so on --- is unimportant for the discussion at hand, we shall leave that sense unspecified.}
to the Yang-Mills gradient flow equation \eqref{eq:Yang-Mills_gradient_flow_equation} on $P$ with initial condition \eqref{eq:Yang-Mills_heat_or_gradient_flow_equation_initial_condition}. Define a family of connections $\widehat A(t)$, for $t \in [0, T)$, on a principal $G$-bundle $\widehat P$ over $S^d$ by the cutting off procedure prescribed in Definition \ref{defn:Construction_spliced_connection_over_sphere}. Let $g_1$ denote the standard round metric of radius one on $S^d$. Define the perturbation, $R(t)$, for $t \in [0, T)$, in \eqref{eq:Simon_0-1} by
\begin{equation}
\label{eq:Yang-Mills_gradientlike_flow_equation_sphere}
\frac{\partial\widehat A}{\partial t} = -d_{\widehat A(t)}^{*, g_1}F_{\widehat A(t)} + R(t), \quad\forall\, t \in (0, T),
\end{equation}
with
\begin{equation}
\label{eq:Yang-Mills_gradientlike_flow_perturbation_sphere}
R(t)
=
\begin{cases}
-d_{A(t)}^{*, g}F_{A(t)} + d_{A(t)}^{*, g_1}F_{A(t)} &\hbox{over } \varphi_n\left(B_{\rho/N}(0)\right),
\\
d_{\widehat A(t)}^{*,g_1}F_{\widehat A}(t) - \chi d_{A(t)}^{*,g}F_A(t) &\hbox{over } \varphi_n\left(B_{\rho/2}(0) \less B_{\rho/N}(0)\right),
\\
0 &\hbox{over } S^d \less \varphi_n\left(B_{\rho/2}(0)\right), \quad\forall\, t \in (0, T),
\end{cases}
\end{equation}
where, abusing notation, we let $g$ and $A(t)$ denote the pull-backs of the given Riemannian metric and Yang-Mills gradient flow connections from $B_\rho(x_0) \subset X$ to $\varphi_n(B_\rho(0)) \subset S^d$. The initial condition for \eqref{eq:Yang-Mills_gradientlike_flow_equation_sphere} is given by
\begin{equation}
\label{eq:Yang-Mills_gradientlike_flow_equation_sphere_initial_condition}
\widehat A(0) = \widehat A_0,
\end{equation}
where $A_0$ is the initial connection for the Yang-Mills gradient flow equation \eqref{eq:Yang-Mills_gradient_flow_equation} on $P$ and $\widehat A_0$ is produced by the cutting off procedure prescribed in Definition \ref{defn:Construction_spliced_connection_over_sphere}.
\end{defn}

To understand the origin of the expression for $R(t)$ in \eqref{eq:Yang-Mills_gradientlike_flow_perturbation_sphere} over the annulus $\varphi_n\left(B_{\rho/2}(0) \less B_{\rho/N}(0)\right)$, observe that the definition of $\widehat A(t)$ over $S^d$ in \eqref{eq:Construction_spliced_connection_over_sphere} yields
\begin{align*}
\frac{\partial\widehat A}{\partial t} &= \chi \frac{\partial A}{\partial t}  \quad\hbox{(by \eqref{eq:Construction_spliced_connection_over_sphere})}
\\
&= -\chi d_{A(t)}^{*,g}F_A(t) \quad\hbox{(by \eqref{eq:Yang-Mills_gradient_flow_equation})}
\\
&= -d_{\widehat A(t)}^{*,g_1}F_{\widehat A}(t) + \left(d_{\widehat A(t)}^{*,g_1}F_{\widehat A}(t) - \chi d_{A(t)}^{*,g}F_A(t)\right).
\end{align*}
Next, we compute an $L^2(S^d)$ estimate for the perturbation term, $R(t)$, in \eqref{eq:Yang-Mills_gradientlike_flow_perturbation_sphere}. Given positive constants $r_1 < r_2$, with $r_2$ less than one half the injectivity radius, $r_0(x_0)$, of the Riemannian metric $g$ at a point $x_0 \in X$, it will be convenient denote the open annulus in $X$ with radii $r_1$ and $r_2$ and center $x_0$ by
\begin{equation}
\label{eq:annulus}
\Omega(x_0, r_1, r_2) := B_{r_2}(x_0) \less \bar B_{r_1}(x_0) \subset X,
\end{equation}
and similarly let $\varphi_n(\Omega(0, r_1, r_2)) \subset S^4$ denote the open annulus in $S^4$ with center at the north pole $n \in S^4$, where $\Omega(0, r_1, r_2) = B_{r_2}(0) \less \bar B_{r_1}(0) \subset \RR^4$ is the open annulus in $\RR^4$ with radii $r_1$ and $r_2$ and center at the origin. We identify $\Omega(x_0, r_1, r_2) \cong \varphi_n(\Omega(0, r_1, r_2))$ via the orientation-preserving diffeomorphism $\varphi_n\circ\varphi_0^{-1}$ in the setting of Definition \ref{defn:Construction_spliced_connection_over_sphere}.

\begin{prop}[$L^2$ \apriori estimate for the perturbation term in Yang-Mills gradient-like flow over the four-dimensional sphere]
\label{prop:Yang-Mills_gradientlike_flow_L2_spatial_norm_4sphere_R}
Let $K$ be a positive constant, $P$ be a principal $G$-bundle over a closed, smooth four-dimensional manifold, $X$, with Riemannian metric, $g$, and $T > 0$, and $A_1$ be a $C^\infty$ reference connection on $P$, and $A(t) = A_1 + a(t)$ be a strong solution to the Yang-Mills gradient flow equation \eqref{eq:Yang-Mills_gradient_flow_equation} over $(0, T)\times X$ with
$$
a(t) \in L^\infty(X; \Lambda^1\otimes \ad P) \cap H^2_{A_1}(X; \Lambda^1\otimes \ad P), \quad\hbox{a.e. } t \in (0, T).
$$
Given a point $x_0 \in X$ and constant $\rho \in (0, \rho_0]$, suppose that $g$ is \emph{locally conformally flat} near $x_0$ in the sense that
\begin{equation}
\label{eq:Four_manifold_locally_conformally_flat_near_point}
(\varphi_0\circ\varphi_n^{-1})^*g = g_1 \quad\hbox{on } B_{\rho/2}(x_0),
\end{equation}
where $g_1$ is the standard round metric of radius one on $S^4$ and $\rho_0$, $\varphi_0$, and $\varphi_n$ are as in Definition \ref{defn:Construction_spliced_connection_over_sphere}. Suppose further that $a(t)$ obeys, for a constant $N \geq 4$,
\begin{subequations}
\label{eq:W1infinity_bounds_Yang-Mills_gradient_flow_and_A1_minus_Gamma_annulus}
\begin{align}
\label{eq:W1infinity_bound_Yang-Mills_gradient_flow_annulus}
\|a(t)\|_{W^{1,\infty}_{A_1}(\Omega)} &\leq K, \quad\forall\, t \in [0, T),
\\
\label{eq:W1infinity_bound_A1_minus_Gamma_annulus}
\|A_1 - \Gamma\|_{W^{1,\infty}_\Gamma(\Omega)} &\leq K,
\end{align}
\end{subequations}
where $\Omega = \Omega(x_0, \rho/N, \rho/2)$ and $\Gamma$ is the product connection on $(B_\rho(x_0) \less \{x_0\})\times G \cong P \restriction (B_\rho(x_0) \less \{x_0\})$ arising in Definition \ref{defn:Construction_spliced_connection_over_sphere}. Define a flow $\widehat A(t)$ for $t \in [0, T)$ on a principal $G$-bundle $\widehat P$ over $S^4$ by the cutting off procedure prescribed in Definitions \ref{defn:Construction_spliced_connection_over_sphere} and \ref{defn:Splicing_construction_Yang-Mills_gradientlike_flow_over_sphere}. Then the resulting perturbation term, $R(t)$, in \eqref{eq:Yang-Mills_gradientlike_flow_perturbation_sphere} obeys, for all $t \in [0, T)$,
\begin{equation}
\label{eq:L2_spatial_norm_sphere_R}
\|R(t)\|_{L^2(S^4)}
\leq
c\left(\rho^2 + (\log N)^{-1/2}\right)K(1 + K^2),
\end{equation}
where $c$ depends at most on the Riemannian metric, $g$, on $X$.
\end{prop}

Before we commence the proof of Proposition \ref{prop:Yang-Mills_gradientlike_flow_L2_spatial_norm_4sphere_R}, it is useful to recall the following elementary lemma from \cite{FLKM1}, which is turn a simple extension of \cite[Lemma 7.2.10]{DK}.

\begin{lem}[Integral bounds on the derivatives of a radial cut-off function]
\label{lem:L2_and_L4_nabla_chi_and_L2_nabla2_chi_bounds}
\cite[Lemma 5.8]{FLKM1}
There is a positive constant, $c$, with the following significance. For any constants $N\geq 4$ and $\rho>0$, there is a $C^\infty$ cut-off function $\chi \equiv \chi_{N,\rho}$ on $\RR^4$
such that
$$
\chi(x) =
\begin{cases}
1 &\text{if }|x|\geq \rho/2,
\\
0 &\text{if }|x|\leq \rho/N,
\end{cases}
$$
and satisfying the following estimates,
\begin{subequations}
\label{eq:L2_and_L4_nabla_chi_and_L2_nabla2_chi_bounds}
\begin{align}
\label{eq:L2_nabla_chi_bound}
\|\cov\chi\|_{L^2(\RR^4)} &\leq c\rho(\log N)^{-1},
\\
\label{eq:L4_nabla_chi_bound}
\|\cov\chi\|_{L^4(\RR^4)} &\leq c(\log N)^{-3/4},
\\
\label{eq:L2_nabla2_chi_bound}
\|\cov^2\chi\|_{L^2(\RR^4)} &\leq c(\log N)^{-1/2}.
\end{align}
\end{subequations}
\end{lem}

\begin{rmk}[On the construction of the cut-off function, $\chi_{N,\rho}$]
\label{rmk:Construction_cut-off_function_chi_N_rho}
It will be useful to recall from \cite{FLKM1} the construction of the cut-off function $\chi_{N,\rho}$ in Lemma \ref{lem:L2_and_L4_nabla_chi_and_L2_nabla2_chi_bounds}. We begin by fixing a $C^\infty$ cut-off function $\kappa:\RR\to [0,1]$ such that $\kappa(t) = 1$ for $t \geq 1$ and $\kappa(t) = 0$ for $t \leq 0$. Now define a $C^\infty$ cut-off function $\alpha_N:\RR\to [0,1]$, depending on the parameter $N$, by setting
$$
\alpha_N(t) := \kappa\left(\frac{\log N + t}{\log N - \log 2}\right).
$$
Finally, define $C^\infty$ cut-off functions $\beta = \beta_{N,\rho}:\RR\to [0,1]$ and $\chi = \chi_{N,\rho}:\RR^4\to [0,1]$, depending on the parameters $N$ and $\rho$, by setting
\begin{align}
\label{eq:Definition_beta_N_rho}
\beta_{N,\rho}(s) &:= \alpha_N(\log s - \log \rho), \quad s \in \RR,
\\
\label{eq:Definition_chi_in_terms_of_beta_Euclidean_space}
\chi_{N,\rho}(x) &:= \beta_{N,\rho}(|x|), \quad x \in \RR^4.
\end{align}
From the proof of \cite[Lemma 5.8]{FLKM1}, we recall that the derivatives of $\chi_{N,\rho}$ obey the following pointwise bounds on $\RR^4$,
\begin{equation}
\label{eq:Feehan_Leness_5-4}
|\nabla\chi_{N,\rho}(x)| \leq \frac{c}{|x|\log N}
\quad\hbox{and}\quad
|\nabla^2\chi_{N,\rho}(x)| \leq \frac{c}{|x|^2\log N}, \quad\forall\, x \in \RR^4\less\{0\},
\end{equation}
where $c$ is a universal positive constant (that is, depending only on our initial choice of cut-off function, $\kappa$).

Of course, the construction of $\chi = \chi_{N, \rho}$ in \eqref{eq:Definition_chi_in_terms_of_beta_Euclidean_space} transfers easily to any Riemannian manifold, $(X, g)$, given a point $x_0 \in X$ and requirement that $\rho$ is less than the injectivity radius of $(X,g)$, just as in \eqref{eq:Definition_chi_in_terms_of_beta}. The resulting $C^\infty$ cut-off function, $\chi = \chi_{N,\rho,g,x_0}: X \to [0,1]$, again satisfies the integral-norm bounds in Lemma \ref{lem:L2_and_L4_nabla_chi_and_L2_nabla2_chi_bounds} and the pointwise bounds \eqref{eq:Feehan_Leness_5-4} from which they follow, albeit now with a positive constant, $c$, which depends on the Riemannian metric on $g$ (but not the point $x_0$, since we always assume that $X$ is closed).
\end{rmk}

Given Lemma \ref{lem:L2_and_L4_nabla_chi_and_L2_nabla2_chi_bounds}, we can proceed with the

\begin{proof}[Proof of Proposition \ref{prop:Yang-Mills_gradientlike_flow_L2_spatial_norm_4sphere_R}]
The fact that perturbation term, $R(t)$, in the Yang-Mills gradient-like flow equation \eqref{eq:Yang-Mills_gradientlike_flow_equation_sphere} is non-zero is due, in general, to two sources of error, namely the effect of
\begin{inparaenum}[\itshape a\upshape)]
\item the difference between the standard round metric of radius one, $g_1$, over $S^4$ and the pull back, via the orientation-preserving diffeomorphism, $\varphi_0\circ \varphi_n^{-1}$, of the given Riemannian metric, $g$, on the ball, $B_{\rho/2}(x_0) \subset X$; and
\item the cut-off function, $\chi(\dist_g(\cdot, x_0)/\rho)$, over the annulus $\varphi_n(B_{\rho/2}(0) \less B_{\rho/N}(0))$ in the construction of $\widehat A(t)$ in Definition \ref{defn:Splicing_construction_Yang-Mills_gradientlike_flow_over_sphere}, together with the difference between the metrics $g$ and $g_1$.
\end{inparaenum}
Thus, it is natural to write the perturbation term, $R(t)$, in equation \eqref{eq:Yang-Mills_gradientlike_flow_perturbation_sphere} as
$$
R(t)
=
\begin{cases}
R_g(t) &\hbox{over } \varphi_n(B_{\rho/N}(0)),
\\
R_\chi(t) &\hbox{over } \varphi_n(B_{\rho/2}(0) \less B_{\rho/N}(0)),
\\
0 &\hbox{over } S^4\less \varphi_n(B_{\rho/2}(0), \quad\forall\, t \in [0, T),
\end{cases}
$$
where
\begin{align}
\label{eq:Yang-Mills_gradientlike_flow_perturbation_sphere_g}
R_g(t) &:= -d_{A(t)}^{*, g}F_A(t) + d_{A(t)}^{*, g_1}F_A(t) \quad\hbox{over } \varphi_n(B_{\rho/N}(0)),
\\
\label{eq:Yang-Mills_gradientlike_flow_perturbation_sphere_chi}
R_\chi(t) &:= d_{\widehat A(t)}^{*,g_1}F_{\widehat A}(t) - \chi d_{A(t)}^{*,g}F_A(t)
\quad\hbox{over } \varphi_n(B_{\rho/2}(0) \less B_{\rho/N}(0)), \quad\forall\, t \in [0, T).
\end{align}
Because $\widehat A(t) \equiv \Gamma$ and $F_{\widehat A}(t) \equiv 0$ over $S^d \less \varphi_n(B_{\rho/2}(0))$ for all $t \in [0, T)$, the cut-off flow, $\widehat A(t)$, is a (trivial) solution to the Yang-Mills gradient flow equation over $S^4 \less \varphi_n(B_{\rho/2}(0))$, with respect to the Riemannian metric, $g_1$.

By hypothesis \eqref{eq:Four_manifold_locally_conformally_flat_near_point} on $g$, the term $R_g$ is identically zero on the ball $\varphi_n(B_{\rho/2}(0))$ and the term $R_\chi$ is non-zero solely due to the presence of the cut-off function, $\chi$.

Thus, we need only estimate $R_\chi(t)$. To minimize clutter, we suppress the explicit notation for the diffeomorphism $\varphi_n\circ\varphi_0^{-1}$ defined by the coordinate charts and write $A(t) = \Gamma + a_1 + a(t)$, and $A_1 = \Gamma + a_1$, and $\widehat A(t) = \Gamma + \hat a(t)$ with $\hat a(t) := \chi(a_1 + a(t))$ as in equation \eqref{eq:Construction_spliced_connection_over_sphere}, where $\Gamma$ denotes the product connection on $(B_\rho(x_0) \less \{x_0\}) \times G$. Consequently, writing $b(t) := a_1 + a(t)$
and $\Omega = \varphi_n(B_{\rho/2}(0) \less \bar B_{\rho/N}(0)) \cong \Omega(x_0,\rho/N,\rho/2)$ for brevity, we have
$$
\widehat A(t) = \Gamma + \chi b(t) \quad\hbox{over } \Omega,
$$
and so our expression \eqref{eq:Yang-Mills_gradientlike_flow_perturbation_sphere_chi} for $R_\chi(t)$ becomes
\begin{align*}
R_\chi(t) &= d_{\widehat A(t)}^{*,g_1}F_{\widehat A}(t) - \chi d_{A(t)}^{*,g}F_A(t)
\\
&= d_{\widehat A(t)}^{*,g}F_{\widehat A}(t) - \chi d_{A(t)}^{*,g}F_A(t)
\quad\hbox{(by \eqref{eq:Four_manifold_locally_conformally_flat_near_point})}
\\
&= d_{\Gamma + \chi b(t)}^{*,g}F_{\Gamma + \chi b}(t)
- \chi d_{\Gamma + b(t)}^{*,g}F_{\Gamma + b}(t)
\\
&= d_\Gamma^{*,g}\left(d_\Gamma(\chi b(t)) + [\chi b(t), \chi b(t)] \right)
-*_g \left[\chi b(t), *_gd_\Gamma(\chi b(t)) + *_g[\chi b(t), \chi b(t)] \right]
\\
&\quad - \chi d_\Gamma^{*,g}\left(d_\Gamma b(t) + [b(t), b(t)] \right)
- *_g\left[b(t), *_gd_\Gamma b(t) + *_g[b(t), b(t)] \right].
\end{align*}
Formally expanding the preceding expression for $R_\chi(t)$, combining like terms, and dropping the distinction between $\chi$ and $1-\chi$ yields
\begin{multline}
\label{eq:R_chi_expansion}
R_\chi(t) = \nabla\chi \times \nabla_\Gamma b(t)
+ \nabla^2\chi \times b(t)
\\
+ \chi\nabla_\Gamma b(t) \times \chi b(t) + \nabla\chi \times b(t) \times \chi b(t)
+ \chi b(t) \times \chi b(t) \times \chi b(t).
\end{multline}
Therefore, applying the H\"older inequality we obtain
\begin{align*}
\|R_\chi(t)\|_{L^2(\Omega)} &\leq c\|\nabla\chi\|_{L^4(\Omega)} \|\nabla_\Gamma b(t)\|_{L^4(\Omega)}
+ c\|\nabla^2\chi\|_{L^2(\Omega)} \|b(t)\|_{L^\infty(\Omega)}
\\
&\quad + c\|\nabla_\Gamma b(t)\|_{L^4(\Omega)} \|b(t)\|_{L^4(\Omega)}
+ c\|\nabla\chi\|_{L^4(\Omega)} \|b(t)\|_{L^4(\Omega)} \|b(t)\|_{L^\infty(\Omega)}
\\
&\quad + c\|b(t)\|_{L^4(\Omega)}^2 \|b(t)\|_{L^\infty(\Omega)},
\end{align*}
where $c$ depends at most on the Riemannian metric, $g$, on $X$.

For any integer $k \geq 1$, the construction of $\chi$ given in Remark \ref{rmk:Construction_cut-off_function_chi_N_rho} yields the following pointwise bounds,
\begin{equation}
\label{eq:nabla_k_chi_rho_pointwise_bound}
|\nabla^k\chi(x)|
\leq
\begin{cases}
c_{k,N}/\rho^k, &x \in \Omega,
\\
0, &x \in X \less \Omega,
\end{cases}
\end{equation}
where $c_{k,N}$ depends at most on the Riemannian metric, $g$, on $X$ and $k\geq 1$ and $N\geq 4$. However, the more precise pointwise bounds given in \eqref{eq:Feehan_Leness_5-4} result in the following useful integral-norm bounds provided by Lemma \ref{lem:L2_and_L4_nabla_chi_and_L2_nabla2_chi_bounds} (when $X$ has dimension four),
$$
\|\nabla\chi\|_{L^4(\Omega)} \leq \frac{c}{(\log N)^{3/4}}
\quad\hbox{and}\quad
\|\nabla^2\chi\|_{L^2(\Omega)} \leq \frac{c}{(\log N)^{1/2}},
$$
where $c$ depends at most on the Riemannian metric, $g$, on $X$. Hence, applying the preceding estimates and combining terms,
\begin{align*}
\|R_\chi(t)\|_{L^2(\Omega)} &\leq c\left((\log N)^{-3/4} + \|b(t)\|_{L^4(\Omega)}\right)\|\nabla_\Gamma b(t)\|_{L^4(\Omega)}
\\
&\quad + c\left((\log N)^{-1/2} + (\log N)^{-3/4}\|b(t)\|_{L^4(\Omega)} + \|b(t)\|_{L^4(\Omega)}^2 \right) \|b(t)\|_{L^\infty(\Omega)},
\end{align*}
where $c$ depends at most on the Riemannian metric, $g$, on $X$. By our hypothesis \eqref{eq:W1infinity_bound_Yang-Mills_gradient_flow_annulus} on $A(t) = A_1 + a(t)$ over the annulus $\Omega \cong \Omega(x_0, \rho/N, \rho/2)$ for $t \in [0, T)$,
$$
\|a(t)\|_{W^{1,\infty}_{A_1}(\Omega)} \leq K, \quad\forall\, t \in [0, T),
$$
and as $b(t) = a_1 + a(t)$ on $B_\rho(x_0) \less \{x_0\}$, we have
$$
\|b(t)\|_{W^{1,\infty}_{A_1}(\Omega)} \leq K + \|a_1\|_{W^{1,\infty}_{A_1}(\Omega)},
\quad\forall\, t \in [0, T),
$$
where also our hypothesis \eqref{eq:W1infinity_bound_A1_minus_Gamma_annulus} on $A_1 = \Gamma + a_1$ gives
$$
\|a_1\|_{W^{1,\infty}_\Gamma(\Omega)} \leq K.
$$
Thus, noting that $0 < \rho \leq 1$ and $N \geq 4$ and $A_1 = \Gamma + a_1$ on $B_{4\rho}(x_0) \less \{x_0\}$, and
$$
\nabla_{A_1} b(t) = \nabla_\Gamma b(t) + [a_1, b(t)],
$$
we obtain
$$
\|R_\chi(t)\|_{L^2(\Omega)}
\leq
c\left((\log N)^{-3/4} + K\rho \right)K\rho
+ c\left((\log N)^{-1/2} + (\log N)^{-3/4}K\rho + K^2\rho^2 \right) K
$$
and thus, simplifying slightly to give a more manageable bound,
\begin{equation}
\label{eq:L2_norm_sphere_Rchi}
\|R_\chi(t)\|_{L^2(\Omega)}
\leq
c\left(\rho^2 + (\log N)^{-1/2}\right)K(1 + K^2),
\quad\forall\, t \in [0, T),
\end{equation}
where $c$ depends at most on the Riemannian metric, $g$, on $X$ and this yields the desired inequality \eqref{eq:L2_spatial_norm_sphere_R}.
\end{proof}

We now verify that the Yang-Mills gradient-like flow, $\widehat A(t)$, on the principal $G$-bundle $\widehat P$ over $S^d$ satisfies \apriori estimates of the form \eqref{eq:Abstract_apriori_interior_estimate_trajectory_with_small_constant_error} given in Hypothesis \ref{hyp:Abstract_apriori_interior_estimate_trajectory_with_perturbation}. It is useful at this point to define a metric on $S^d$ which agrees with the metric $g$ on $X$, after identifying the ball $B_\rho(x_0) \subset X$ with the ball $\varphi_n(B_\rho(0)) \subset S^d$ via the coordinate charts.

\begin{defn}[Splicing construction of an almost round Riemannian metric on the sphere]
\label{defn:Splicing_construction_Riemannian_metric_on_sphere}
Assume the notation and set-up of Definition \ref{defn:Construction_spliced_connection_over_sphere}. Let $g_1$ denote the standard round metric of radius one on $S^d$ and let $g_{1,\rho}$ denote the almost round metric with radius approximately equal to one obtained by pulling back the given Riemannian metric, $g$, on $B_\rho(x_0) \subset X$ via the specified coordinate charts,
\begin{equation}
\label{eq:Construction_spliced_metric_over_sphere}
g_{1,\rho}
:=
\begin{cases}
(\varphi_0\circ \varphi_n^{-1})^*g &\hbox{over } \varphi_n\left(B_{\rho/2}(0)\right),
\\
\left(1-\chi(\dist_g(\cdot, x_0)/(2\rho)\right)g_1 &{}
\\
\quad + (\varphi_0\circ \varphi_n^{-1})^*(\chi g)
&\hbox{over } \varphi_n\left(B_\rho(0) \less B_{\rho/2}(0)\right),
\\
g_1 &\hbox{over } \varphi_n\left(B_\rho(0)\right).
\end{cases}
\end{equation}
where $\chi = \chi_{N,\rho,g,x_0}$ is as in \eqref{eq:Definition_chi_in_terms_of_beta} and may be regarded as a $C^\infty$ cut-off function on $X$ or $S^d$ via the identifications implied by the local coordinate charts in Definition \ref{defn:Splicing_construction_Riemannian_metric_on_sphere}.
\end{defn}

With the aid of the Riemannian metric on $S^d$ in Definition \eqref{defn:Splicing_construction_Riemannian_metric_on_sphere}, we can derive an \apriori estimate for a Yang-Mills gradient-like flow over $S^d$ where the only perturbation source term that appears is due to the effect of the cut-off function, $\chi$.

\begin{lem}[Higher-order \apriori estimates for the cut-off function perturbation term in Yang-Mills gradient-like flow over the four-dimensional sphere]
\label{lem:Yang-Mills_gradientlike_flow_higher-order_apriori_estimates_Rchi}
Assume the hypotheses of Proposition \ref{prop:Yang-Mills_gradientlike_flow_L2_spatial_norm_4sphere_R}, but strengthen the conditions \eqref{eq:W1infinity_bounds_Yang-Mills_gradient_flow_and_A1_minus_Gamma_annulus} to
\begin{subequations}
\label{eq:W4infinity_bounds_Yang-Mills_gradient_flow_and_A1_minus_Gamma_and_W2infinity_dot_flow_annulus}
\begin{align}
\label{eq:W4infinity_bound_Yang-Mills_gradient_flow_annulus}
\|a(t)\|_{W^{4,\infty}_\Gamma(\Omega)} &\leq K,
\\
\label{eq:W2infinity_bound_time_derivative_Yang-Mills_gradient_flow_annulus}
\|\dot a(t)\|_{W^{2,\infty}_\Gamma(\Omega)} &\leq K, \quad\forall\, t \in [0, T),
\\
\label{eq:W4infinity_bound_A1_minus_Gamma_annulus}
\|A_1 - \Gamma\|_{W^{4,\infty}_\Gamma(\Omega)} &\leq K,
\end{align}
\end{subequations}
where $\Omega = \Omega(x_0, \rho/N, \rho/2)$ and $\Gamma$ is the product connection on $(B_\rho(x_0) \less \{x_0\})\times G \cong P \restriction (B_\rho(x_0) \less \{x_0\})$ arising in Definition \ref{defn:Construction_spliced_connection_over_sphere}. Then the perturbation term, $R_\chi(t)$, in \eqref{eq:Yang-Mills_gradientlike_flow_perturbation_sphere_chi} obeys
\begin{align}
\label{eq:H2_Gamma_spatial_norm_sphere_Rchi}
\|R_\chi(t)\|_{H^2_\Gamma(\Omega)} &\leq  C(1 + \rho^{-2})K(1 + K^2),
\\
\label{eq:Linfinity_spatial_norm_sphere_Rchi}
\|R_\chi\|_{L^\infty(\Omega)} &\leq C(1 + \rho^{-2}) K(1 + K^2),
\\
\label{eq:H2_A1_spatial_norm_sphere_Rchi}
\|R_\chi(t)\|_{H^2_{A_1}(\Omega)} &\leq  C(1 + \rho^{-2})K(1 + K^4),
\\
\label{eq:L2_spatial_norm_sphere_time_derivative_Rchi}
\|\dot R_\chi(t)\|_{L^2(\Omega)} &\leq cK(1 + K^2),  \quad\forall\, t \in [0, T),
\end{align}
where the positive constant, $c$, depends at most on the Riemannian metric, $g$, on $X$ and the positive constant, $C$, depends at most on $g$ and $N$.
\end{lem}

\begin{proof}
We keep the notation of the proof of Proposition \ref{prop:Yang-Mills_gradientlike_flow_L2_spatial_norm_4sphere_R} and note that the pointwise estimates \eqref{eq:nabla_k_chi_rho_pointwise_bound} for $\nabla^k \chi$ on $\Omega = \Omega(x_0, \rho/N, \rho/2) \cong \varphi_n(\Omega(0, \rho/N, \rho/2))$ imply the $L^2$ bounds,
\begin{equation}
\label{eq:nabla_k_chi_rho_L2_bound}
\|\nabla^k\chi\|_{L^2(\Omega)} \leq c_{k,N}\rho^{2-k},
\end{equation}
where $c_{k,N}$ depends at most on the integer $k \geq 1$, parameter $N$, and the Riemannian metric, $g$, on $X$. By repeatedly applying the H\"older inequality to the expression \eqref{eq:R_chi_expansion} for $R_\chi$ and its first and second covariant derivatives with respect to $\Gamma$, we obtain
\begin{multline*}
\|R_\chi\|_{L^\infty(\Omega)}  \leq \|\nabla_\Gamma^2 b\|_{L^\infty(\Omega)}
+ \frac{C}{\rho}\|\nabla_\Gamma b\|_{L^\infty(\Omega)}
+ \frac{C}{\rho^2}\|b\|_{L^\infty(\Omega)}
+ C\|b\|_{L^\infty(\Omega)} \|\nabla_\Gamma b\|_{L^\infty(\Omega)}
\\
+ \frac{C}{\rho}\|b\|_{L^\infty(\Omega)}^2
+ C\|b\|_{L^\infty(\Omega)}^3
+ \|\chi d_{A(t)}^{*,g}F_A(t)\|_{L^\infty(\Omega)},
\end{multline*}
and
\begin{align*}
\|\nabla_\Gamma R_\chi\|_{L^2(\Omega)}  &\leq \|\nabla_\Gamma^3 b\|_{L^2(\Omega)}
+ \frac{C}{\rho} \|\nabla_\Gamma^2 b\|_{L^2(\Omega)}
+ \frac{C}{\rho^2} \|\nabla_\Gamma b\|_{L^2(\Omega)}
+ \frac{C}{\rho^3} \|b\|_{L^2(\Omega)}
\\
&\quad + \frac{C}{\rho} \|\nabla_\Gamma b\|_{L^4(\Omega)} \|b\|_{L^4(\Omega)}
+ C\|\nabla_\Gamma^2 b\|_{L^4(\Omega)} \|b\|_{L^4(\Omega)}
+ C\|\nabla_\Gamma b\|_{L^4(\Omega)}^2
\\
&\quad + \frac{C}{\rho^2} \|b\|_{L^4(\Omega)}^2
+ \frac{C}{\rho} \|b\|_{L^4(\Omega)}^2\|b\|_{L^\infty(\Omega)}
+ C\|\nabla_\Gamma b\|_{L^4(\Omega)} \|b\|_{L^4(\Omega)} \|b\|_{L^\infty(\Omega)}
\\
&\quad + \|\nabla_\Gamma (\chi d_{A(t)}^{*,g}F_A(t))\|_{L^2(\Omega)},
\end{align*}
and
\begin{align*}
\|\nabla_\Gamma^2 R_\chi\|_{L^2(\Omega)}  &\leq \|\nabla_\Gamma^4 b\|_{L^2(\Omega)}
+ \frac{C}{\rho} \|\nabla_\Gamma^3 b\|_{L^2(\Omega)}
+ \frac{C}{\rho^2} \|\nabla_\Gamma^2 b\|_{L^2(\Omega)}
+ \frac{C}{\rho^3} \|\nabla_\Gamma b\|_{L^2(\Omega)}
\\
&\quad + \frac{C}{\rho^4} \|b\|_{L^2(\Omega)}
+ \frac{C}{\rho^2} \|\nabla_\Gamma b\|_{L^4(\Omega)} \|b\|_{L^4(\Omega)}
+ \frac{C}{\rho} \|\nabla_\Gamma^2 b\|_{L^4(\Omega)} \|b\|_{L^4(\Omega)}
\\
&\quad + \frac{C}{\rho} \|\nabla_\Gamma b\|_{L^4(\Omega)}^2
+ C\|\nabla_\Gamma^3 b\|_{L^4(\Omega)} \|b\|_{L^4(\Omega)}
+ C\|\nabla_\Gamma^2 b\|_{L^4(\Omega)} \|\nabla_\Gamma b\|_{L^4(\Omega)}
\\
&\quad + \frac{C}{\rho^3} \|b\|_{L^4(\Omega)}^2
+ \frac{C}{\rho^2} \|b\|_{L^4(\Omega)}^2\|b\|_{L^\infty(\Omega)}
+ \frac{C}{\rho} \|\nabla_\Gamma b\|_{L^4(\Omega)}  \|b\|_{L^4(\Omega)} \|b\|_{L^\infty(\Omega)}
\\
&\quad +  C\|\nabla_\Gamma b\|_{L^4(\Omega)}^2 \|b\|_{L^\infty(\Omega)}
+  C\|\nabla_\Gamma^2 b\|_{L^4(\Omega)}  \|b\|_{L^4(\Omega)} \|b\|_{L^\infty(\Omega)}
\\
&\quad + \|\nabla_\Gamma^2 (\chi d_{A(t)}^{*,g}F_A(t))\|_{L^2(\Omega)},
\end{align*}
where $C$ depends at most on the Riemannian metric, $g$, on $X$, and parameter $N$. Because $b(t) = a_1 + a(t)$, the bound \eqref{eq:W4infinity_bound_Yang-Mills_gradient_flow_annulus} for $a(t)$ and bound \eqref{eq:W4infinity_bound_A1_minus_Gamma_annulus} for $a_1$ allow us to simplify the preceding estimates and give
$$
\|R_\chi\|_{L^\infty(\Omega)} \leq C\left(K + \frac{K}{\rho} + \frac{K}{\rho^2} + K^2 + \frac{K^2}{\rho} + K^3\right) + \|\chi d_{A(t)}^{*,g}F_A(t)\|_{L^\infty(\Omega)},
$$
and
\begin{align*}
\|\nabla_\Gamma R_\chi\|_{L^2(\Omega)}
&\leq C\left(K\rho^2 + K\rho + K + \frac{K}{\rho} + K^2\rho+ K^2\rho^2 + K^2 + K^3\rho + K^3\rho^2\right)
\\
&\quad + \|\nabla_\Gamma (\chi d_{A(t)}^{*,g}F_A(t))\|_{L^2(\Omega)},
\end{align*}
and
\begin{align*}
{}&\|\nabla_\Gamma^2 R_\chi\|_{L^2(\Omega)}
\\
&\leq C\left(K\rho^2 + K\rho + K + \frac{K}{\rho} + \frac{K}{\rho^2}
+ K^2 + K^2\rho+ K^2\rho^2 + \frac{K^2}{\rho} + K^3 + K^3\rho + K^3\rho^2\right)
\\
&\quad + \|\nabla_\Gamma^2 (\chi d_{A(t)}^{*,g}F_A(t))\|_{L^2(\Omega)}
\end{align*}
where $C$ depends at most on the Riemannian metric, $g$, on $X$, and parameter $N$. Noting that $0<\rho\leq 1$ by hypothesis, the preceding bounds simplify to give
\begin{align*}
\|R_\chi\|_{L^\infty(\Omega)} &\leq CK(1 + K^2)(1 + \rho^{-2}) + \|\chi d_{A(t)}^{*,g}F_A(t)\|_{L^\infty(\Omega)},
\\
\|\nabla_\Gamma R_\chi\|_{L^2(\Omega)}
&\leq CK(1 + K^2)(1 + \rho^{-1}) + \|\nabla_\Gamma (\chi d_{A(t)}^{*,g}F_A(t))\|_{L^2(\Omega)},
\\
\|\nabla_\Gamma^2 R_\chi\|_{L^2(\Omega)}
&\leq CK(1 + K^2)(1 + \rho^{-2}) + \|\nabla_\Gamma^2 (\chi d_{A(t)}^{*,g}F_A(t))\|_{L^2(\Omega)},
\end{align*}
where $C$ depends at most on the Riemannian metric, $g$, on $X$, and parameter $N$.

Next, we estimate the indicated norms of $\chi d_{A(t)}^{*,g}F_A(t)$. Recall that $A(t) = A_1 + a(t) = \Gamma + a_1 + a(t) = \Gamma + b(t)$ over $\Omega$ and so, writing $F_A(t) = d_\Gamma b(t) + [b(t), b(t)]$, we obtain the formal expressions,
$$
d_A^{*,g}F_A = \nabla_\Gamma^2 b + b \times \nabla_\Gamma b + b\times b \times b,
$$
and
\begin{align*}
\nabla_\Gamma (\chi d_A^{*,g}F_A) &= \chi \left(\nabla_\Gamma^3 b + \nabla_\Gamma b \times \nabla_\Gamma b
+ b \times \nabla_\Gamma^2 b + \nabla_\Gamma b\times b \times b\right)
\\
&\quad + \nabla\chi\times \left(\nabla_\Gamma^2 b + b \times \nabla_\Gamma b + b\times b \times b\right),
\end{align*}
and
\begin{align*}
\nabla_\Gamma^2 (\chi d_A^{*,g}F_A) &= \chi \left(\nabla_\Gamma^4 b + \nabla_\Gamma^2 b \times \nabla_\Gamma b
+ b \times \nabla_\Gamma^3 b + \nabla_\Gamma^2 b\times b \times b + \nabla_\Gamma b\times \nabla_\Gamma b \times b\right)
\\
&\quad + \nabla\chi\times \left(\nabla_\Gamma^3 b + \nabla_\Gamma b \times \nabla_\Gamma b
+ b \times \nabla_\Gamma^2 b + \nabla_\Gamma b\times b \times b\right)
\\
&\quad + \nabla^2\chi\times \left(\nabla_\Gamma^2 b + b \times \nabla_\Gamma b + b\times b \times b\right).
\end{align*}
Consequently, by repeatedly applying the H\"older inequality, the bounds in Lemma \ref{lem:L2_and_L4_nabla_chi_and_L2_nabla2_chi_bounds} for the $L^4$ norm of $\nabla\chi$ and $L^2$ norm of $\nabla^2\chi$, and the bounds \eqref{eq:W4infinity_bounds_Yang-Mills_gradient_flow_and_A1_minus_Gamma_and_W2infinity_dot_flow_annulus} for the $W^{4,\infty}_\Gamma(\Omega)$ norms of $a(t)$ and $a_1$, we obtain the estimates,
\begin{align*}
\|\chi d_A^{*,g}F_A\|_{L^2(\Omega)} &\leq c\rho^2 K(1 + K^2),
\\
\|\chi d_A^{*,g}F_A\|_{L^\infty(\Omega)} &\leq cK(1 + K^2),
\\
\|\nabla_\Gamma (\chi d_A^{*,g}F_A)\|_{L^2(\Omega)} &\leq c\rho(\rho + (\log N)^{-3/4}) K(1 + K^2),
\\
\|\nabla_\Gamma^2 (\chi d_A^{*,g}F_A)\|_{L^2(\Omega)} &\leq c\rho(\rho + (\log N)^{-3/4} + (\log N)^{-1/2}) K(1 + K^2),
\end{align*}
where the positive constant, $c$, depends at most on the Riemannian metric, $g$, on $X$. In particular, noting that $0 < \rho \leq 1$ and $N \geq 4$ by hypothesis, we see that
$$
\|\chi d_{A(t)}^{*,g}F_A(t)\|_{H^2_\Gamma(\Omega)} \leq c\rho(\rho + (\log N)^{-1/2}) K(1 + K^2), \quad\forall\, t \in [0, T).
$$
Combining the preceding inequalities with the $L^2(\Omega)$ estimate \eqref{eq:L2_norm_sphere_Rchi} for $R_\chi(t)$ obtained in the proof of Proposition \ref{prop:Yang-Mills_gradientlike_flow_L2_spatial_norm_4sphere_R} yields the following $H_\Gamma^2(\Omega)$ estimate for $R_\chi(t)$,
$$
\|R_\chi(t)\|_{H^2_\Gamma(\Omega)} \leq C(1 + \rho^{-2}) K(1 + K^2), \quad\forall\, t \in [0, T),
$$
where the positive constant, $C$, depends at most on the Riemannian metric, $g$, on $X$, and parameter $N$; this is the desired $H^2_\Gamma(\Omega)$ estimate \eqref{eq:H2_Gamma_spatial_norm_sphere_Rchi} for $R_\chi(t)$.

Similarly, gathering the $L^\infty(\Omega)$ estimates in the preceding inequalities yields
$$
\|R_\chi\|_{L^\infty(\Omega)} \leq C(1 + \rho^{-2}) K(1 + K^2), \quad\forall\, t \in [0, T),
$$
where the positive constant, $C$, depends at most on the Riemannian metric, $g$, on $X$, and parameter $N$; this is the desired $L^\infty(\Omega)$ estimate \eqref{eq:Linfinity_spatial_norm_sphere_Rchi} for $R_\chi(t)$.

The bound \eqref{eq:W4infinity_bounds_Yang-Mills_gradient_flow_and_A1_minus_Gamma_and_W2infinity_dot_flow_annulus} for the $W^{4,\infty}_\Gamma(\Omega)$ norm of $a_1$ ensures that the derivation of the preceding estimate can be repeated, \emph{mutatis mutandis}, with $\nabla_\Gamma$ replaced by $\nabla_{A_1}$ to give a bound with respect to the $H^2_{A_1}(\Omega)$ norm,
$$
\|R_\chi(t)\|_{H^2_{A_1}(\Omega)} \leq C(1 + \rho^{-2}) K(1 + K^4), \quad\forall\, t \in [0, T),
$$
where $C$ depends at most on the Riemannian metric, $g$, on $X$, and parameter $N$; this is the desired $H_{A_1}^2(\Omega)$ estimate \eqref{eq:H2_A1_spatial_norm_sphere_Rchi} for $R_\chi(t)$.

By taking the time derivative of the expression \eqref{eq:R_chi_expansion} for $R_\chi(t)$, the derivation leading to the $L^2(\Omega)$ bound \eqref{eq:L2_norm_sphere_Rchi} for $R_\chi$ now gives, again keeping in mind the bounds \eqref{eq:W4infinity_bounds_Yang-Mills_gradient_flow_and_A1_minus_Gamma_and_W2infinity_dot_flow_annulus},
$$
\|\dot R_\chi(t)\|_{L^2(\Omega)} \leq cK(1 + K^2) + \|\partial_t(d_{A(t)}^{*,g}F_A(t))\|_{L^2(\Omega)}, \quad\forall\, t \in [0, T),
$$
where $c$ depends at most on the Riemannian metric, $g$, on $X$. Expanding $\partial_t(d_{A(t)}^{*,g}F_A(t))$ schematically yields,
$$
\frac{\partial}{\partial t}(d_{A(t)}^{*,g}F_A(t)) = \nabla_\Gamma^2 \dot b(t) + b(t)\times \nabla_\Gamma \dot b(t) + \dot b(t)\times \nabla_\Gamma b(t) + b(t)\times b(t)\times \dot b(t).
$$
By applying the bounds \eqref{eq:W4infinity_bounds_Yang-Mills_gradient_flow_and_A1_minus_Gamma_and_W2infinity_dot_flow_annulus} once again, we see that
$$
\|\partial_t(d_{A(t)}^{*,g}F_A(t))\|_{L^2(\Omega)} \leq cK(1 + K^2), \quad\forall\, t \in [0, T),
$$
where $c$ depends at most on the Riemannian metric, $g$, on $X$. Combining the preceding inequalities yields the desired $L^2(\Omega)$ estimate \eqref{eq:L2_spatial_norm_sphere_time_derivative_Rchi} for $\dot R_\chi(t)$.
\end{proof}




We then have the

\begin{prop}[\Apriori estimates for Yang-Mills gradient-like flow over the four-dimensional sphere]
\label{prop:Apriori_estimate_Yang-Mills_gradient-like_flow_S4}
Assume the hypotheses of Proposition \ref{prop:Yang-Mills_gradientlike_flow_L2_spatial_norm_4sphere_R} and, in addition, that $\widehat A_1$ is a $C^\infty$ reference connection on $\widehat P$ such that\footnote{The connection, $\widehat A_1$, can be constructed explicitly using the splicing procedure in Definition \ref{defn:Construction_spliced_connection_over_sphere}.}
$$
\widehat A_1
=
\begin{cases}
\Gamma &\hbox{on } \widehat P \restriction S^4 \less \varphi_n(B_\rho(0)),
\\
A_1 &\hbox{on } \widehat P \restriction \varphi_n(B_{\rho/2}(0)),
\end{cases}
$$
where $\Gamma$ is the product connection on $\widehat P \restriction S^4 \less \{n\} \cong (S^4 \less \{n\})\times G$. Then there are positive constants, $C = C(A_1, g)$ and $\eps_1 = \eps_1(A_1, g) \in (0, 1]$, and, given $\beta \in [3/4,1)$, positive constants, $C = C(A_1, g, \beta)$ and $\eps_1 = \eps_1(A_1, g, \beta) \in (0, 1]$, with the following significance.
\begin{enumerate}
\item If the Yang-Mills gradient-like flow, $\widehat A(t)$, obeys
\begin{equation}
\label{eq:Linfinity_in_time_H1_in_space_small_norm_At_minus_A_1_condition_lemma_7.3_S4}
\|\widehat A(t) - \widehat A_1\|_{H_{\widehat A_1}^1(S^4)} \leq \eps_1 \quad\forall\, t \in (S, T),
\end{equation}
then
\begin{multline}
\label{eq:L1_interior_time_H1_space_apriori_estimate_Yang-Mills_gradient-like_flow_S4}
\int_{S+\delta}^T \|\partial_t\widehat A \|_{H_{\widehat A_1}^1(S^4)}\,dt
\leq C(1 + \delta^{-1/2})\int_S^T \|\partial_t\widehat A \|_{L^2(S^4)}\,dt
\\
+ C\sqrt{\delta}(T-S)(1 + \rho^{-2})K(1 + K^2).
\end{multline}

\item If $\beta \in [3/4,1)$ and the Yang-Mills gradient-like flow, $\widehat A(t)$, obeys
\begin{equation}
\label{eq:Linfinity_in_time_H2beta_in_space_small_norm_At_minus_A_1_condition_lemma_7.3_S4}
\|\widehat A(t) - \widehat A_1\|_{H_{\widehat A_1}^{2\beta}(S^4)}
\leq \eps_1 \quad\forall\, t \in (S, T),
\end{equation}
then
\begin{multline}
\label{eq:L1_interior_time_H2beta_space_apriori_estimate_Yang-Mills_gradient-like_flow_S4}
\int_{S+\delta}^T \|\partial_t\widehat A \|_{H_{\widehat A_1}^{2\beta}(S^4)}\,dt
\leq C(1 + \delta^{-1})\int_S^T \|\partial_t\widehat A \|_{L^2(S^4)}\,dt
\\
+ C(T-S)(1 + \rho^{-2})K(1 + K^2).
\end{multline}
\end{enumerate}
\end{prop}

\begin{rmk}[On the role of the connection $\widehat A_1$ in Proposition \ref{prop:Apriori_estimate_Yang-Mills_gradient-like_flow_S4}]
The proof of Proposition relies on Lemmata \ref{rmk:Rade_7-3} and \ref{lem:Rade_7-3_L1_in_time_H2beta_in_space_apriori_estimate_by_L1_in_time_L2_in_space} and those lemmata in turn hypothesize that the conditions \eqref{eq:Linfinity_in_time_H1_in_space_small_norm_At_minus_A_1_condition_lemma_7.3_S4} and \eqref{eq:Linfinity_in_time_H2beta_in_space_small_norm_At_minus_A_1_condition_lemma_7.3_S4}, respectively, hold when $\widehat A_1$ is replaced by a Yang-Mills connection. See Remark \ref{rmk:Rade_7-3} for additional discussion.
\end{rmk}


\begin{proof}[Proof of Proposition \ref{prop:Apriori_estimate_Yang-Mills_gradient-like_flow_S4}]
By appealing to the construction in Definition \ref{defn:Splicing_construction_Riemannian_metric_on_sphere} of the Riemannian metric, $g_{1,\rho}$, on $S^4$, we can consider $\widehat A$ to be a solution to the Yang-Mills gradient-like flow equation where the full perturbation, $R(t)$, is replaced by the \emph{local-in-space} perturbation term, $R_\chi(t)$, which is identically zero outside the annulus, $\Omega = \varphi_n(B_{\rho/2}(0)\less \bar B_{\rho/N}(0)) \subset S^4$,
\begin{equation}
\label{eq:Yang-Mills_gradientlike_flow_equation_sphere_Rchi}
\frac{\partial\widehat A}{\partial t} = -d_{\widehat A}^{*,g_{1,\rho}}F_{\widehat A} + R_\chi \quad\hbox{a.e. over } (S, T)\times S^4,
\end{equation}
and $\widehat A(t) = \Gamma$ is a (trivial solution) to the Yang-Mills gradient flow equation over $S^4 \less \varphi_n(B_{\rho/2}(0))$.

We first apply Lemma \ref{lem:Rade_7-3} to the Yang-Mills gradient-like flow equation \eqref{eq:Yang-Mills_gradientlike_flow_equation_sphere_Rchi} and observe that the \apriori estimate \eqref{eq:Rade_apriori_interior_estimate_lemma_7.3} gives
\begin{multline*}
\int_{S+\delta}^T \|\partial_t\widehat A \|_{H_{\widehat A_1}^1(S^4)}\,dt
\leq C(1 + \delta^{-1/2})\int_S^T \|\partial_t\widehat A \|_{L^2(S^4)}\,dt
\\
+ C\sqrt{\delta}\int_{S+\delta/2}^T \left(\|R_\chi(t)\|_{H^2_{\widehat A_1}(S^4)} + \|R_\chi(t)\|_{L^\infty(S^4)}
+ \|\partial_t R_\chi\|_{L^2(S^4)}\right)\,dt.
\end{multline*}
We now apply Lemma \ref{lem:Yang-Mills_gradientlike_flow_higher-order_apriori_estimates_Rchi}, substituting our $H^2_{A_1}(S^4)\cap L^\infty(S^4)$ bounds \eqref{eq:H2_Gamma_spatial_norm_sphere_Rchi} and \eqref{eq:Linfinity_spatial_norm_sphere_Rchi} for $R_\chi$ and $L^2(S^4)$ bound \eqref{eq:L2_spatial_norm_sphere_time_derivative_Rchi} for $\partial_t R_\chi$ into the preceding inequality to give the estimate \eqref{eq:L1_interior_time_H1_space_apriori_estimate_Yang-Mills_gradient-like_flow_S4}.

Next, we apply Lemma \ref{lem:Rade_7-3_L1_in_time_H2beta_in_space_apriori_estimate_by_L1_in_time_L2_in_space} to the Yang-Mills gradient-like flow equation \eqref{eq:Yang-Mills_gradientlike_flow_equation_sphere_Rchi} and observe that the \apriori estimate \eqref{eq:Rade_7-3_L1_in_time_H2beta_in_space_apriori_estimate_by_L1_in_time_L2_in_space} gives
\begin{multline*}
\int_{S+\delta}^T \|\partial_t\widehat A \|_{H_{\widehat A_1}^{2\beta}(X)}\,dt
\leq C(1 + \delta^{-1})\int_S^T \|\partial_t\widehat A \|_{L^2(S^4)}\,dt
\\
+ C\int_{S+\delta/2}^T\left(\|R_\chi(t)\|_{H^2_{\widehat A_1}(S^4)} + \|R_\chi(t)\|_{L^\infty(S^4)} + \|\partial_t R_\chi(t)\|_{L^2(S^4)} \right)\,dt.
\end{multline*}
Again applying Lemma \ref{lem:Yang-Mills_gradientlike_flow_higher-order_apriori_estimates_Rchi} yields the estimate \eqref{eq:L1_interior_time_H2beta_space_apriori_estimate_Yang-Mills_gradient-like_flow_S4}. This completes the proof.
\end{proof}

\begin{rmk}[Application of Propositions \ref{prop:Yang-Mills_gradientlike_flow_L2_spatial_norm_4sphere_R} and \ref{prop:Apriori_estimate_Yang-Mills_gradient-like_flow_S4}]
\label{rmk:Application_Yang-Mills_gradient-like_flow_S4}
If the bounds provided by Propositions \ref{prop:Yang-Mills_gradientlike_flow_L2_spatial_norm_4sphere_R} and \ref{prop:Apriori_estimate_Yang-Mills_gradient-like_flow_S4} could be further sharpened, they might be used to verify the hypotheses of Theorem \ref{thm:Huang_5-1-1_gradientlike_uniform_continuity} for abstract gradient-like flows and hence provide continuity in $t \in [0, T]$ for the Yang-Mills gradient-like flow, $\widehat A(t)$, and thus a `no bubbling in finite time' result for Yang-Mills gradient flow.
A measure of the sharpness of the estimates required is evident from Section \ref{subsec:Rade_proposition_7-2_dimension_four}, where we discuss the dependence on the Riemannian metric of the constants in the {\L}ojasiewicz-Simon gradient inequality for the Yang-Mills energy functional in dimension four and the effect of rescaling.
\end{rmk}

\begin{rmk}[On eliminating the possibility of bubbling at $T = \infty$]
\label{rmk:No_bubbling_at_T_equalto_infinity}
In order to further refine the potential application described in Remark \ref{rmk:Application_Yang-Mills_gradient-like_flow_S4} to give no bubbling at $T = \infty$, one would need to show that the perturbation, $R(t)$, converges to zero sufficiently fast as $t \to \infty$ with respect to the $L^2(S^4)$,  $H^2_{A_1}(S^4)$, and $L^\infty(S^4)$ spatial norms and that the time derivative, $\dot R(t)$, also converges to zero sufficiently fast as $t \to \infty$ with respect to the $L^2(S^4)$ spatial norms. For example, in \cite[Theorem 2]{Simon_1983}, Simon assumes exponential-in-time decay for $R(t)$ to zero as $t \to \infty$, although with different choices of spatial norms. When $\theta = 1/2$ in the {\L}ojasiewicz-Simon gradient inequality \eqref{eq:Simon_2-2} with positive constants $c$, $\sigma$, and $\theta \in [1/2, 1)$, then our abstract convergence-rate result, Theorem \ref{thm:Huang_3-4-8}, yields exponential-in-time convergence (with rate $e^{-c^2t/2}$) for the abstract gradient flow, $u(t)$, to a critical point, $u_\infty$, with respect to a norm for a Banach space, $\sX$. Also when $\theta = 1/2$, our generalization, Theorem \ref{thm:Huang_3-4-8_Yang-Mills}, of R\r{a}de's \cite[Proposition 7.4]{Rade_1992} yields exponential-in-time convergence (with rate $e^{-c^2t/2}$) of the Yang-Mills gradient flow, $A(t)$ on a principal $G$-bundle $P$ over $X$, to a Yang-Mills connection, $A_\infty$, on $P$, with respect to a spatial norm, $H_{A_1}^1(X)$. While one might be able to localize such exponential convergence results on fixed small balls in the complement of potential bubble points in a four-dimensional manifold, $X$, one in general only has power-law convergence in time, with rate $t^{-(1-\theta)/(2\theta -1)}$, for  some $\theta \in (1/2, 1)$, and in general such a convergence rate would be too slow for our application. When $X$ is a Riemann surface, $G = U(n)$, and $A_\infty$ is irreducible, then R\r{a}de has shown that $\theta = 1/2$ in \eqref{eq:Simon_2-2} and thus obtains exponential convergence, but it is unclear that one can draw a similar conclusion when $X$ is a four-dimensional manifold. When $G = \SU(2)$ or $\SO(3)$ and $b^+(X) > 0$ and $A_\infty$ is a $g$-anti-self-dual Yang-Mills connection, then \cite[Corollary 4.3.15]{DK} implies that $A_\infty$ is necessarily irreducible when $g$ is a generic Riemannian metric.
\end{rmk}

\subsection{Refined {\L}ojasiewicz-Simon gradient inequality for the Yang-Mills energy functional in dimension four}
\label{subsec:Rade_proposition_7-2_dimension_four}
In order to examine the dependence on the Riemannian metric of the {\L}ojasiewicz-Simon triple of constants, $(c,\theta,\sigma)$, we shall specialize the {\L}ojasiewicz-Simon gradient inequality for the Yang-Mills energy functional (Theorem \ref{thm:Rade_proposition_7-2}) to dimension four. As we shall most clearly see in Corollary \ref{cor:Rade_proposition_7-2_S4} below, the {\L}ojasiewicz-Simon constant, $c$, has an \emph{unfavorable} behavior with respect to conformal rescaling.

For this purpose, it is convenient to define a norm for $\Omega^1(X;\ad P)$ which only depends on the Riemannian metric on $X$ through its conformal equivalence class, $[g]$, namely,
\begin{equation}
\label{eq:H1_norm_one-forms_conformally_invariant}
\|a\|_{H_A^1(X,[g])}
:=
\|\nabla_A a\|_{L^2(X,g)} + \|a\|_{L^4(X,g)}, \quad\forall\, a \in \Omega^1(X;\ad P).
\end{equation}
For simplicity in the statement of the {\L}ojasiewicz-Simon gradient inequality, we take $A_1 = A_\ym$ in Theorem \ref{thm:Rade_proposition_7-2} (as well as $d=4$).

\begin{cor}[Refined {\L}ojasiewicz-Simon gradient inequality for the Yang-Mills energy functional in dimension four]
\label{cor:Rade_proposition_7-2_dimension_four}
Let $(X,g)$ be a closed, four-dimensional, Riemannian manifold, $G$ a compact Lie group, and $A_\ym$ a Yang-Mills connection of class $C^\infty$, on a principal $G$-bundle $P$ over $X$. Then there are positive constants $c$, $\sigma$, and $\theta \in [1/2,1)$, depending on $A_\ym$, $g$, $G$, $P$, and $X$ with the following significance.  If $A$ is a connection of class $H^1$ on $P$ and
\begin{equation}
\label{eq:Rade_7-1_neighborhood_dimension_four}
\|A - A_\ym\|_{H_{A_\ym}^1(X,[g])} < \sigma,
\end{equation}
then
\begin{equation}
\label{eq:Rade_7-1_dimension_four}
\|d_A^*F_A\|_{L^2(X,g)} \geq c|\sE(A;[g]) - \sE(A_\ym;[g])|^\theta,
\end{equation}
where $\sE(A;[g]) = \sE(A)$ is given by \eqref{eq:Potential_yang_mills}. The constants $\sigma$ and $\theta$ depend on the Riemannian metric on $X$ only through its conformal equivalence class, $[g]$.
\end{cor}

We further specialize Corollary \ref{cor:Rade_proposition_7-2_dimension_four} to the case of the four-dimensional sphere, $S^4$, with its standard round metric of radius one. Following the conventions of Definition \ref{defn:Construction_spliced_connection_over_sphere} and \cite[Section 3.2]{FeehanGeometry}, for any $\lambda>0$, let $f_\lambda:S^4\to S^4$ denote the conformal diffeomorphism of $S^4$ defined by
\begin{equation}
\label{eq:Feehan_1995_page_478_rescaling_diffeomorphism_S4}
f_\lambda := \varphi_n \circ \bar f_\lambda \circ \varphi_n^{-1},
\end{equation}
where $\bar f_\lambda:\RR^4\to\RR^4$ is defined by $x\mapsto y = x/\lambda$ and $\varphi_n:\RR^4\to S^4\less\{s\}$ is inverse to a stereographic projection from the south pole $s \in S^4 \subset \RR^5$ (identified with the point at infinity in $\RR^4\cup\{\infty\}$). If $a \in \Omega^1(S^4,\ad P)$ and $\omega \in \Omega^2(S^4,\ad P)$, it follows by direct calculation and the conformal invariance of the $L^2$ norm on two-forms and $L^4$ norm on one-forms that
\begin{subequations}
\label{eq:L2_and_L4_norms_diff_forms_S4_and_conformal_group_action}
\begin{align}
\label{eq:Rescaling_L2_norm_one-form_S^4}
\|f_\lambda^*a\|_{L^2(S^4)} &= \lambda\|a\|_{L^2(S^4)},
\\
\label{eq:Conformal_invariance_L4_norm_one-form_S^4}
\|f^*a\|_{L^4(S^4)} &= \|a\|_{L^4(S^4)},
\\
\label{eq:Conformal_invariance_L2_norm_two-form_S^4}
\|f^*\omega\|_{L^2(S^4)} &= \|\omega\|_{L^2(S^4)},
\end{align}
\end{subequations}
where $f_\lambda$ is as in \eqref{eq:Feehan_1995_page_478_rescaling_diffeomorphism_S4}, while $f$ is \emph{any} conformal diffeomorphism of $S^4$.

We recall from \cite[Equation (3.10)]{FeehanGeometry}, which is based in turn on Taubes \cite[p. 343]{TauFrame}, the definition of the \emph{center}, $q = \mathrm{Center}[A] \in \RR^4$, and \emph{scale}, $\lambda = \mathrm{Scale}[A] \in (0,\infty)$, of a connection, $A$, on a principal $G$-bundle $P$ over $S^4 \cong \RR^4\cup\{\infty\}$:
\begin{subequations}
\label{eq:Feehan_1995_3-10}
\begin{align}
\label{eq:Feehan_1995_3-10_center}
q^\mu &:= \frac{1}{\sE(A)} \int_{\RR^4} x^\mu |F_A|^2\, d^4x, \quad 1\leq \mu \leq 4,
\\
\label{eq:Feehan_1995_3-10_scale}
\lambda^2 &:= \frac{1}{\sE(A)} \int_{\RR^4} |x-q|^2 |F_A|^2\, dx.
\end{align}
\end{subequations}
Alternative definitions of center and scale are possible --- for example, motivated by Section \ref{sec:Taubes_1982_Appendix}, one might replace $\sE(A)$ by $4\pi^2|\bkappa(P)|$, the topological invariant representing the minimum value of the energy attainable by a connection $A$ on $P$.

The map $f_\lambda$ lifts to an action on the affine space of connections on $P\to S^4$, as discussed in \cite[Section 3.2]{FeehanGeometry}. If a connection, $A$, on $P$ has scale one in the sense of \eqref{eq:Feehan_1995_3-10_scale}, then $f_\lambda^*A$ will have scale $\lambda$. It is useful to recall from \cite[Lemma 3.1]{TauFrame} that there is a \emph{universal} constant, $z \in [1,\infty)$, such that
$$
z^{-1}\|a\|_{H_A^1(S^4)} \leq \|f^*a\|_{H_A^1(S^4)} \leq \|a\|_{H_A^1(S^4)}, \quad\forall\, a\in \Omega^1(S^4,\ad P),
$$
where $f:S^4\to S^4$ is any conformal diffeomorphism. Thus, it is convenient to define a norm for $\Omega^1(S^4;\ad P)$ which is exactly conformally invariant and also equivalent to the standard definition of the norm on $H_A^1(S^4;\ad P)$ via universal constants, as suggested by \eqref{eq:H1_norm_one-forms_conformally_invariant}, that is,
\begin{equation}
\label{eq:H1_norm_one-forms_conformally_invariant_S4}
\|a\|_{H_A^1[S^4]}
:=
\|\nabla_A a\|_{L^2(S^4)} + \|a\|_{L^4(S^4)}, \quad\forall\, a \in \Omega^1(S^4;\ad P).
\end{equation}
With the preceding comments in mind, we can now provide the promised specialization of Corollary \ref{cor:Rade_proposition_7-2_dimension_four} to the four-sphere:

\begin{cor}[Refined {\L}ojasiewicz-Simon gradient inequality for the Yang-Mills energy functional over the four-dimensional sphere]
\label{cor:Rade_proposition_7-2_S4}
Let $G$ a compact Lie group and $A_\ym$ a Yang-Mills connection of class $C^\infty$ on a principal $G$-bundle $P$ over $S^4$ with its standard round metric of radius one. Then there are positive constants $c$, $\sigma$, and $\theta \in [1/2,1)$, depending on $A_\ym$, with the following significance.  If $A$ is a connection of class $H^1$ on $P$ and
\begin{equation}
\label{eq:Rade_7-1_neighborhood_S4}
\|A - A_\ym\|_{H_{A_\ym}^1[S^4]} < \sigma,
\end{equation}
then
\begin{equation}
\label{eq:Rade_7-1_dimension_S4}
\|d_A^*F_A\|_{L^2(S^4)} \geq c|\sE(A) - \sE(A_\ym)|^\theta,
\end{equation}
where $\sE(A)$ is given by \eqref{eq:Potential_yang_mills}. Moreover, if $\lambda>0$ and $f_\lambda:S^4\to S^4$ is the conformal diffeomorphism \eqref{eq:Feehan_1995_page_478_rescaling_diffeomorphism_S4}, then
\begin{equation}
\label{eq:Rade_7-1_dimension_S4_rescaled}
\|d_{f_\lambda^*A}^*F_{f_\lambda^*A}\|_{L^2(S^4)} \geq c\lambda |\sE(f_\lambda^*A) - \sE(f_\lambda^*A_\ym)|^\theta,
\end{equation}
In particular, if $(c,\theta,\sigma)$ are the {\L}ojasiewicz-Simon constants for $A_\ym$, then $(c\lambda,\theta,\sigma)$ are the {\L}ojasiewicz-Simon constants for $f_\lambda^*A_\ym$, for any $\lambda>0$.
\end{cor}

\begin{proof}
It suffices to observe that
\begin{align*}
\|d_{f_\lambda^*A}^*F_{f_\lambda^*A}\|_{L^2(S^4)}
&=
\lambda\|d_A^*F_A\|_{L^2(S^4)}
\quad\hbox{(by \eqref{eq:Rescaling_L2_norm_one-form_S^4})}
\\
&\geq c\lambda|\sE(A) - \sE(A_\ym)|^\theta
\quad\hbox{(by \eqref{eq:Rade_7-1_dimension_S4})}
\\
&= c\lambda|\sE(f_\lambda^*A) - \sE(f_\lambda^*A_\ym)|^\theta
\quad\hbox{(by \eqref{eq:Conformal_invariance_L2_norm_two-form_S^4})}.
\end{align*}
This completes the proof.
\end{proof}

\begin{rmk}[Identifying the {\L}ojasiewicz-Simon exponent for a Yang-Mills connection over $S^4$]
\label{rmk:Rade_section_10_explicit_Lojasiewicz-Simon exponents}
If $X$ is closed two-dimensional Riemannian manifold and $G=\U(n)$, R\r{a}de has shown that any Yang-Mills connection, $A$, on a principal $G$-bundle $P$ over $X$ has {\L}ojasiewicz-Simon exponent in the range $1/2 \leq \theta \leq 3/4$ and, if $A$ is also irreducible, then one has the optimal exponent, $\theta=1/2$ \cite[Section 10]{Rade_1992}. Thus, it seems reasonable that one may be able to generalize R\r{a}de's argument and identify more precise bounds on the {\L}ojasiewicz-Simon exponent when $A$ is a Yang-Mills connection on a principal $G$-bundle over $S^4$, at least for certain compact Lie groups, $G$.

It is interesting to note that if $u:[0,\infty)\to\sX$ is a trajectory in the sense of Definition \ref{defn:Huang_3-1-1} and $\int_0^\infty \|\dot u(t)\|_\sX\,dt < \infty$, then a gradient map, $\sE':\sU\subset\sX \to \sH$, necessarily satisfies a gradient inequality in the orbit $O(u) = \{u(t):t\geq 0\}$,
$$
\phi(\sE(u(t)) \leq \|\sE'(u(t)\|_\sH, \quad\forall\, t \in [0,\infty),
$$
where $\phi:\RR\to[0,\infty)$ is a measurable function such that $1/\phi \in L^1(\RR)$ (see  \cite[Proposition 3.3.1]{Huang_2006}).
\end{rmk}

\section[Yang-Mills gradient-like flow over K{\"a}hler surfaces and applications]{Yang-Mills gradient-like flow over compact K{\"a}hler surfaces and applications}
\label{sec:Yang-Mills_gradient-like_flow_Kahler_surface}
Our goal in this section is to adapt certain results developed by Donaldson and Kronheimer in \cite[Sections 6.2 and 6.3]{DK} from the context of (pure) Yang-Mills gradient flow on a Hermitian vector bundle over a compact, K\"ahler surface to the case of Yang-Mills gradient-like (also called pseudo-gradient) flow. Specifically, we shall establish a continuous extension property (Theorem \ref{thm:Donaldson_Kronheimer_6-3-11_complete}) for smooth solutions, $A(t)$, to Yang-Mills gradient-like given an initial Hermitian connection $A_0$ on $E$ which has curvature of type $(1,1)$ over a given open subset $U \subset X$, but not necessarily over all of $X$. Lastly, we apply the continuous extension result for Yang-Mills gradient-like flow over K\"ahler surfaces to prove global existence for pure Yang-Mills gradient flow with arbitrary initial energy over complex surfaces with non-K\"ahler metrics (Theorem \ref{mainthm:Yang-Mills_gradient_flow_global_existence_and_convergence_started_with arbitrary_initial_energy}).

\subsection{Preliminaries on connections with curvature of type $(1,1)$, the K\"ahler identities, and the Bochner-Weitzenb\"ock formula}
\label{subsec:Preliminaries_curvature_1-1_Kaehler_identities_Bochner}
We recall some basic facts concerning connections on Hermitian vector bundles over complex surfaces from \cite[Sections 2.1.5, 6.1.1, and 6.1.3]{DK}, for ease of reference and to confirm sign conventions. We closely follow the development in \cite{DK}.

For a compact Lie group $G$ and principal $G$-bundle $P$ over a compact, complex surface, $X$, we let $E := P\times_\varrho \CC^n$, where we assume $G \subset \U(n)$ and let $\varrho:G\to\GL(n,\CC)$ denote the standard representation. As in \cite[p. 33]{DK}, we write $\fg_E = P\times_{\ad}\fg = \ad P$, for the bundle of Lie algebras associated to the adjoint representation, so $\fg_E$ is a real subbundle of $\End E = E\otimes E^*$. Thus, if $G=\U(n)$, then $\fg_E$ consists of skew-adjoint endomorphisms of $E$ and if $G=\SU(n)$, then those endomorphisms are also trace-free.

As in \cite[Equation (2.1.47)]{DK} we write, for any connection $A$ on $E$,
\begin{equation}
\label{eq:Donaldson_Kronheimer_2.1.47}
\nabla_A = d_A := \partial_A \oplus \bar\partial_A:
\Omega^0(X;\fg_E) \to \Omega^{1,0}(X;\fg_E) \oplus \Omega^{0,1}(X;\fg_E),
\end{equation}
with respect to the decomposition
$$
\Omega^1(X;\fg_E) = \Omega^{1,0}(X;\fg_E) \oplus \Omega^{0,1}(X;\fg_E).
$$
If $A$ is any connection on $E$, then its curvature may be expressed in terms of its components
$$
F_A = F_A^{2,0} + F_A^{1,1} + F_A^{0,2},
$$
with respect to the decomposition \cite[p. 45]{DK},
$$
\Omega^2(X;\fg_E) = \Omega^{0,2}(X;\fg_E) \oplus \Omega^{1,1}(X;\fg_E) \oplus \Omega^{0,2}(X;\fg_E),
$$
where \cite[p. 45]{DK}
\begin{equation}
\label{eq:FA02_is_barpartialA_squared}
F_A^{0,2} = \bar\partial_A^2.
\end{equation}
Thus, \cite[Theorem 2.1.53]{DK} asserts that $\bar\partial_A$ defines a holomorphic structure on $E$ if and only if $F_A^{0,2}=0$.

If $A$ is a unitary connection, then one has $F_A^{0,2} = - (F_A^{2,0})^*$ by \cite[p. 46]{DK}. One says that the curvature of $A$ is type $(1,1)$ if $F_A^{0,2} = 0$ (and thus $F_A^{2,0} = 0$), so $F_A = F_A^{1,1}$. Consequently, a unitary connection $A$ on $E$ is compatible with a holomorphic structure on $E$ if and only if its curvature is of type $(1,1)$ and, in this case, the connection is uniquely determined by the Hermitian metric and holomorphic structure on $E$ (see \cite[Lemma 2.1.54]{DK} and \cite[Proposition 2.1.56]{DK}).

Recall also from \cite[Lemma 2.1.57]{DK} that if $X$ is a complex surface with its natural orientation and Hermitian metric, $h$, then
\begin{equation}
\label{eq:Donaldson_Kronheimer_lemma_2-1-57}
\begin{aligned}
\Omega^2(X;\CC) &= \Omega^+(X;\CC) \oplus \Omega^-(X;\CC),
\\
\Omega^+(X;\CC) &= \Omega^{2,0}(X;\CC) \oplus \Omega^0(X;\CC)\omega \oplus \Omega^{0,2}(X;\CC),
\\
\Omega^-(X;\CC) &= \Omega_0^{1,1}(X;\CC),
\end{aligned}
\end{equation}
where $\Omega_0^{1,1}(X;\CC)$ consists of forms in $\Omega^{1,1}(X;\CC)$ that are pointwise orthogonal to $\omega$, the $(1,1)$ form defined by the complex structure and the Hermitian metric, $h$, on $X$.

If $A$ is any connection on $E$ then, following \cite[p. 47]{DK}, we write
\begin{equation}
\label{eq:Donaldson_Kronheimer_2.1.58}
\widehat F_A := \langle F_A,\omega\rangle \in \Omega^0(X;\fg_E).
\end{equation}
In particular, one has the following key relationship (\cite[Proposition 2.1.59]{DK}): If $A$ is an anti-self-dual connection on a complex vector bundle $E$ over a complex, Hermitian surface, $X$, then $\bar\partial_A$ defines a holomorphic structure on $E$. Conversely if $\partial_\sE$ is a holomorphic structure on $E$, and $A$ is a unitary connection compatible with that holomorphic structure (that is, $\bar\partial_A=\bar\partial_\sE$), then $A$ is anti-self-dual if and only
if $\widehat F_A = 0$.

As in \cite[p. 210]{DK}, we let $\sG_E^\CC$ denote the group of all general linear
automorphisms of the complex vector bundle, $E$, covering the identity map on $X$. This contains as a subgroup the group $\sG$ of automorphisms preserving the Hermitian metric on $E$, and so $\sG_E^\CC$ can be thought of as the complexification of $\sG_E$. The action of $\sG_E^\CC$ on $E$ preserves the subspace $\sA^{1,1}(E) \subset \sA(E)$ of unitary connections on $E$ whose curvature is type (1,1).

When $A(t)$, for $t\in [0,T)$, is a solution to pure Yang-Mills gradient flow on a Hermitian vector bundle $E$ over a complex Hermitian surface $X$,
\begin{equation}
\label{eq:Yang-Mills_gradient_flow_complex_surface}
\frac{\partial A}{\partial t} = -d_{A(t)}^*F_{A(t)}, \quad\forall\,t\in (0,T),
\end{equation}
where
\begin{equation}
\label{Yang-Mills_gradient_flow_initial_data_type_1-1}
A(0) = A_0 \in \sA^{1,1}(E),
\end{equation}
we recall from \cite[p. 218]{DK} that the gradient vectors $d_A^*F_A$ lie in the tangent spaces to the $\sG_E^\CC$ orbits. The Lie algebra of $\sG_E^\CC$ is $\Omega^0(X;\End_0(E)) = \Omega^0(\fg_E \otimes \CC)$ and the derivative of the $\sG_E^\CC$ action at $A$ is \cite[Equation (6.2.6)]{DK},
\begin{equation}
\label{eq:Donaldson_Kronheimer_6-2-6}
\phi \mapsto -\bar\partial_A\phi + \partial_A(\phi^*),
\end{equation}
for $\phi \in \Omega^0(X;\End_0(E))$. The tangent vector $-d_A^*F_A$ represents the
infinitesimal action of the element $i\widehat F_A$ of the Lie algebra of $\sG_E^\CC$.

For any connection $A$ on $E$, the operators in \eqref{eq:Donaldson_Kronheimer_2.1.47} extend to give
\begin{equation}
\label{eq:Donaldson_Kronheimer_2.1.47_pq_forms}
d_A := \partial_A \oplus \bar\partial_A:
\Omega^{p,q}(X;\fg_E) \to \Omega^{p+1,q}(X;\fg_E) \oplus \Omega^{p,q+1}(X;\fg_E).
\end{equation}
Therefore, the Bianchi identity $d_AF_A = 0$ on a complex surface is equivalent to
\begin{equation}
\label{eq:Complex_Bianchi_identities}
\partial_AF_A = 0 \quad\hbox{and}\quad \bar\partial_AF_A = 0.
\end{equation}
As in \cite[p. 212]{DK}, we write
$$
\Lambda:\Omega^{p,q}(X;\CC) \to \Omega^{p-1,q-1}(X;\CC)
$$
for the adjoint of the operator,
$$
\Omega^{p,q}(X;\CC) \ni \alpha \mapsto \alpha\wedge\omega \in \Omega^{p+1,q+1}(X;\CC).
$$
For any connection, $A$, on $E$ we have
\begin{equation}
\label{eq:hatFA_is_LambdaFA}
\widehat F_A = \Lambda F_A.
\end{equation}
When the Hermitian metric, $h$, on $X$ is \emph{K\"ahler}, with $(1,1)$ form $\omega$, one has the \emph{K\"ahler identities} on $\Omega^{p,q}(X; E)$ \cite[Equation (6.2.6)]{DK}\footnote{We correct a missing $*$ in the expression for $\partial_A^*$},
\begin{equation}
\label{eq:Donaldson_Kronheimer_6-1-8}
\bar\partial_A^* = i[\partial_A,\Lambda], \quad \partial_A^* = -i[\bar\partial_A,\Lambda].
\end{equation}
On $\Omega^{1,0}(X; E)$ and $\Omega^{0,1}(X; E)$, respectively, these reduce to\footnote{We correct the signs.}
\begin{equation}
\label{eq:Donaldson_Kronheimer_6-1-8_one-zero_and_zero-one_forms}
\bar\partial_A^* = -i\Lambda\partial_A, \quad \partial_A^* = i\Lambda\bar\partial_A.
\end{equation}
Since $d_A^* = \partial_A^* + \bar\partial_A^*$ by \eqref{eq:Donaldson_Kronheimer_2.1.47_pq_forms}, we have
\begin{align*}
d_A^*F_A &= \partial_A^*F_A + \bar\partial_A^*F_A
\\
&= -i[\bar\partial_A,\Lambda]F_A + i[\partial_A,\Lambda]F_A
\quad\hbox{(by \eqref{eq:Donaldson_Kronheimer_6-1-8})}
\\
&= -i\bar\partial_A\Lambda F_A + i\partial_A\Lambda F_A,
\quad\hbox{(by \eqref{eq:Complex_Bianchi_identities})}
\end{align*}
and so \cite[Equation (6.2.4)]{DK}\footnote{We correct a sign.}
\begin{equation}
\label{eq:Donaldson_Kronheimer_6-2-4}
d_A^*F_A = i(\partial_A - \bar\partial_A)\widehat F_A,
\end{equation}
by \eqref{eq:hatFA_is_LambdaFA}. We recall from \cite[Section 6.1.3]{DK} that

\begin{lem}
\label{lem:Donaldson_Kronheimer_6-1-7}
\cite[Lemma 6.1.7]{DK}
For any connection $A$ on a Hermitian vector bundle $E$ over a K\"ahler surface $X$, one has
\begin{align*}
\bar\partial_A^*\bar\partial_A &= \frac{1}{2}\nabla_A^*\nabla_A + i\widehat F_A
\quad\hbox{on }\Omega^0(X;E),
\\
\partial_A^*\partial_A &= \frac{1}{2}\nabla_A^*\nabla_A - i\widehat F_A
\quad\hbox{on }\Omega^0(X;E).
\end{align*}
\end{lem}

In Lemma \ref{lem:Donaldson_Kronheimer_6-1-7}, the term $i\widehat F_A$ acts on $\fg_E$ by the standard algebraic action of $\fg_E$ on $E$. Our only application of Lemma \ref{lem:Donaldson_Kronheimer_6-1-7} will be through the simpler,
\begin{equation}
\label{eq:Donaldson_Kronheimer_lemma_6-1-7_corollary}
\bar\partial_A^*\bar\partial_A + \partial_A^*\partial_A = \nabla_A^*\nabla_A
\quad\hbox{on } \Omega^0(X;E).
\end{equation}
This completes our discussion of preliminaries.

\subsection{Global existence for Yang-Mills gradient-like flow}
We shall adapt certain results from \cite[Sections 6.2 and 6.3]{DK} for pure Yang-Mills gradient flow to the case of Yang-Mills gradient-like flow, for some $T>0$,
\begin{equation}
\label{eq:Yang-Mills_gradient-like_flow_Kaehler_surface}
\frac{\partial A}{\partial t} = -d_{A(t)}^*F_{A(t)} + R(t), \quad\forall\,t\in (0,T),
\end{equation}
with initial data,
\begin{equation}
\label{Yang-Mills_gradient-like_flow_initial_data}
A(0) = A_0 \in \sA(E),
\end{equation}
where $A_0$ is a unitary connection on a Hermitian vector bundle, $E$, over a compact surface, $X$, with K\"ahler metric, $h$, and
\begin{equation}
\label{Yang-Mills_gradient-like_flow_error}
R \in C^\infty([0,\infty)\times X; \Lambda^1\otimes\fg_E).
\end{equation}
While $F_{A_0}$ is not assumed to have type $(1,1)$ over all of $X$, we shall assume that
\begin{equation}
\label{Yang-Mills_gradient-like_flow_initial_data_product_or_type_1-1}
A_0 \in \sA^{1,1}(E\restriction \sU'),
\end{equation}
where $\sU'' \Subset \sU' \Subset \sU \subset X$ are given open subsets (typical choices and their roles are explained in Remark \ref{rmk:Applications_Yang-Mills_gradient-like_flow_Kaehler_surface}). The solution, $A(t)$, to \eqref{eq:Yang-Mills_gradient-like_flow_Kaehler_surface} and \eqref{Yang-Mills_gradient-like_flow_initial_data} will not necessarily have curvature, $F_A(t)$, of type $(1,1)$ on the subset
$X \less \sU'$, as distinct from the situation considered in \cite[Sections 6.2 and 6.3]{DK}, but in our applications we will be able to impose auxiliary conditions which ensure that flow, $A(t)$, remains well behaved over $X \less \sU''$ for $t\in [0,T)$.

\begin{rmk}[Applications of Yang-Mills gradient-like flow]
\label{rmk:Applications_Yang-Mills_gradient-like_flow_Kaehler_surface}
A term such as $R(t)$ in \eqref{eq:Yang-Mills_gradient-like_flow_Kaehler_surface} arises in our application by considering a \emph{pure} Yang-Mills gradient flow, $A(t)$ on a Hermitian bundle $E$ over a complex, possibly \emph{non-K\"ahler} surface, $X$. Given a point $x_0\in X$ (for example, one where the curvature of $A(t)$ is assumed to concentrate as $t\nearrow T$), we take
$$
\sU'' := B_{r_2}(x_0),
\quad \sU' := B_{r_1}(x_0),
\quad\hbox{and}\quad \sU = B_{r_0}(x_0), \quad\hbox{for } 0 < r_2<r_1<r_0<\Inj(X,h),
$$
and cut off $A(t)$ over the annulus,
$$
\Omega(x_0;r_1,r_0)
:=
B_{r_0}(x_0) \less \bar B_{r_1}(x_0)
=
\sU \less \bar \sU' \subset X,
$$
with the aid of a cut-off function, $\chi \in C_0^\infty(B_{r_0}(x_0))$ such that $0\leq \chi\leq 1$ on $X$, and a trivialization for $E \restriction B_{r_0}(x_0)$. The resulting flow, $A_\chi(t)$ on a Hermitian bundle $E_\chi$ over $X$, coincides with $A(t)$ on $E$ over $B_{r_1}(x_0)$ and with the product connection, $\Gamma$, on the product bundle, $\CC^n\times X\less B_{r_0}(x_0)$. We will consider an application where the Hermitian metric $h$ on $X$ is \emph{K\"ahler} on a ball $B_{r_0}(x_0)\subset X$, identified by an isometric biholomorphism with a ball $B_{r_0}(x_0)\subset Z$, where $Z$ is a compact K\"ahler surface. In particular, $A_\chi(t)$ may be regarded as a Yang-Mills gradient-\emph{like} flow on the compact K\"ahler surface, $Z$. The flow $A_\chi(t)$ obeys \eqref{eq:Yang-Mills_gradient-like_flow_Kaehler_surface} with error term,
$$
R(t) := d_{A_\chi(t)}^*F_{A_\chi(t)} - \chi d_{A(t)}^*F_{A(t)}
\in \Omega^1(Z;\fg_{E_\chi}), \quad\forall\, t\in[0,T),
$$
due to the effect of the cut-off function and thus supported in $\Omega(x_0;r_1,r_0) \subset Z$. Of course, even if the flow, $A(t)$, has curvature of type $(1,1)$ for $t\in[0,T)$, that will no longer necessarily be true for the cut-off flow, $A_\chi(t)$, over $\Omega(x_0;r_1,r_0)$ but in our application, we shall be able to assume good control on $A(t)$ over $B_{r_0}(x_0)\less B_{r_2}(x_0)$ for $t\in[0,T)$.
\end{rmk}

Given a one-parameter family, $A(t)$ for $t\in [0,T)$, of connections on $E$, we define\footnote{We use $u(t)$ rather than $e(t)$ as in \cite[Sections 6.2 and 6.3]{DK}.}
\begin{equation}
\label{eq:Donaldson_Kronheimer_6-2-10}
u(t) := |\widehat F_A(t)|^2, \quad\forall\, t \in [0,T).
\end{equation}
We have the following analogue of \cite[Lemma 6.2.11]{DK}.

\begin{lem}
\label{lem:Donaldson_Kronheimer_6-2-11}
Let $E$ be a Hermitian vector bundle over a compact, K\"ahler surface, $X$. If $A(t)$, for $t\in (0,T)$, is a smooth solution to the Yang-Mills gradient-like flow equation \eqref{eq:Yang-Mills_gradient-like_flow_Kaehler_surface} and $u$ is as in \eqref{eq:Donaldson_Kronheimer_6-2-10}, then
\begin{equation}
\label{eq:Donaldson_Kronheimer_6-2-11}
\frac{\partial u}{\partial t} + \Delta u = - |d_A^*F_A|^2 + f
\quad\hbox{on } (0,T)\times X,
\end{equation}
where
\begin{equation}
\label{eq:Definition_f_in_terms_A_and_R}
f := 2\langle\widehat F_A,\Lambda d_AR\rangle \in C^\infty((0,T)\times X).
\end{equation}
\end{lem}

\begin{proof}
When $R$ is identically zero, the identity \eqref{eq:Donaldson_Kronheimer_6-2-11} is given by \cite[Lemma 6.2.11]{DK} for pure Yang-Mills gradient flow. When $R$ is non-zero, we have
\begin{align*}
\frac{\partial \widehat F_A}{\partial t} &= \Lambda d_A\left(\frac{\partial A}{\partial t}\right)
\\
&= -\Lambda d_Ad_A^*F_A + \Lambda d_AR
\quad\hbox{(by \eqref{eq:Yang-Mills_gradient-like_flow_Kaehler_surface})}
\\
&= -i\Lambda d_A (\partial_A - \bar\partial_A) \widehat F_A + \Lambda d_AR
\quad\hbox{(by \eqref{eq:Donaldson_Kronheimer_6-2-4})}
\\
&= -i\Lambda (\bar\partial_A\partial_A - \partial_A\bar\partial_A) \widehat F_A + \Lambda d_AR
\quad\hbox{(as $\Lambda\partial_A^2=0=\Lambda\bar\partial_A^2$ on $\Omega^0(X;\fg_E)$)}
\\
&= -(\partial_A^*\partial_A + \bar\partial_A^*\bar\partial_A) \widehat F_A + \Lambda d_AR
\quad\hbox{(by \eqref{eq:Donaldson_Kronheimer_6-1-8})},
\end{align*}
and so by \eqref{eq:Donaldson_Kronheimer_lemma_6-1-7_corollary},
\begin{equation}
\label{eq:ddt_hatFA_Yang-Mills_gradient-like_flow}
\frac{\partial \widehat F_A}{\partial t} = -\nabla_A^*\nabla_A \widehat F_A + \Lambda d_AR.
\end{equation}
As in the proof of \cite[Lemma 6.2.11]{DK}, we have
\begin{align}
\label{eq:Freed_Uhlenbeck_6-18_hatFA}
\Delta(|\widehat F_A|^2) &= 2\langle\widehat F_A, \nabla_A^*\nabla_A\widehat F_A\rangle
- |\nabla_A\widehat F_A|^2
\quad\hbox{(by \eqref{eq:Freed_Uhlenbeck_6-18})},
\\
\label{eq:Length_nablaA_hatA_is_length_dA*FA}
|\nabla_A\widehat F_A|^2 &= |d_A^*F_A|^2,
\end{align}
where the second identity follows from the fact that $d_A^*F_A = i(\partial_A - \bar\partial_A)\widehat F_A$ by \eqref{eq:Donaldson_Kronheimer_6-2-4} while $\nabla_A\widehat F_A = (\partial_A + \bar\partial_A)\widehat F_A$ by \eqref{eq:Donaldson_Kronheimer_2.1.47}. Thus,
\begin{align*}
\frac{\partial |\widehat F_A|^2}{\partial t}
&=
2\left\langle \widehat F_A, \frac{\partial \widehat F_A}{\partial t}\right\rangle
\\
&= -2\langle\widehat F_A, \nabla_A^*\nabla_A\widehat F_A\rangle + 2\langle\widehat F_A,\Lambda d_AR\rangle
\quad\hbox{(by \eqref{eq:ddt_hatFA_Yang-Mills_gradient-like_flow})}
\\
&= -\Delta|\widehat F_A|^2 - |d_A^*F_A|^2 + 2\langle\widehat F_A,\Lambda d_AR\rangle
\quad\hbox{(by \eqref{eq:Freed_Uhlenbeck_6-18_hatFA} and \eqref{eq:Length_nablaA_hatA_is_length_dA*FA})},
\end{align*}
which yields \eqref{eq:Donaldson_Kronheimer_6-2-11} and $f = 2\langle\widehat F_A,\Lambda d_AR\rangle$, as claimed.
\end{proof}

We next give an analogue of \cite[Corollary 6.2.12]{DK}. We find it convenient to denote
$$
X_T := (0,T)\times X \quad\hbox{and}\quad s^+ = s \vee 0, \quad\forall\, s \in \RR.
$$
While \cite[Corollary 6.2.12]{DK} asserts that $\sup_X u(t)$ is a decreasing function of $t\in[0,T)$ in the case of pure Yang-Mills gradient flow, the weaker Corollary \ref{cor:Donaldson_Kronheimer_6-2-12} will suffice for our applications.

\begin{cor}
\label{cor:Donaldson_Kronheimer_6-2-12}
Assume the hypotheses of Lemma \ref{lem:Donaldson_Kronheimer_6-2-11} and, in addition, that $A(t)$ is continuous up to $t=0$. If $u$ is as in \eqref{eq:Donaldson_Kronheimer_6-2-10} and $f$ is as in \eqref{eq:Definition_f_in_terms_A_and_R}, then
\begin{equation}
\label{eq:Donaldson_Kronheimer_6-2-12}
u(t,x) \leq e^{2(T-t)}\sup_{X_T} f^+ \vee \sup_X u(0), \quad\forall\, (t,x)\in X_T.
\end{equation}
\end{cor}

\begin{rmk}[Application of Corollary \ref{cor:Donaldson_Kronheimer_6-2-12}]
\label{rmk:Donaldson_Kronheimer_6-2-12}
We note that, unlike \cite[Corollary 6.2.12]{DK} which gives the simpler bound $u \leq \sup_X u(0)$ on $X_T$, independent of $T \in (0,\infty]$, the bound in Corollary \ref{cor:Donaldson_Kronheimer_6-2-12} is clearly only useful when $T<\infty$.
\end{rmk}

\begin{proof}[Proof of Corollary \ref{cor:Donaldson_Kronheimer_6-2-12}]
The conclusion follows from Lemma \ref{lem:Donaldson_Kronheimer_6-2-11} and the maximum principle for the heat operator, $L := \partial_t + \Delta$, on $C^\infty(X_T)$. Indeed, we have
$$
Lu = - |d_A^*F_A|^2 + f \leq f \quad\hbox{on }X_T.
$$
If $u\geq c_0$ on $X_T$ for a positive constant $c_0$, then
\cite[Proposition 2.6 (Item 2)]{Feehan_parabolicmaximumprinciple} provides the \apriori maximum principle estimate,
$$
u \leq 0 \vee \frac{1}{c_0}\sup_{X_T} Lu \vee \sup_{\mydirac X_T}u
\quad\hbox{on } X_T,
$$
where $\mydirac X_T = \{0\}\times X$, where $\mydirac X_T$ is the parabolic boundary of $X_T$ in the notation of \cite{Feehan_parabolicmaximumprinciple}. In general, we will not have $u\geq c_0$ on $X_T$ for some $c_0>0$, but we can instead appeal to \cite[Lemma 2.7]{Feehan_parabolicmaximumprinciple}, noting that $u \geq 0$ on $X_T$ by definition and so (taking $K_0 = 1$ in \cite[Lemma 2.7]{Feehan_parabolicmaximumprinciple}) we can apply the preceding estimate with $1/c_0$ replaced by $e^{2(T-t)}$, that is
$$
u(t,x) \leq 0 \vee e^{2(T-t)}\sup_{X_T} f \vee \sup_X u(0), \quad\forall\, (t,x) \in X_T,
$$
or more simply \eqref{eq:Donaldson_Kronheimer_6-2-12}, as claimed.
\end{proof}

\begin{rmk}[Temporal sign conventions]
We note that in \cite{Feehan_parabolicmaximumprinciple}, because of the motivating applications to optimal control theory, we considered terminal-value problems and parabolic problems with terminal data at time $T$ rather than initial data at time $0$, as we do here and for that reason, the heat operator is written as $-\partial_t + \Delta$ in \cite{Feehan_parabolicmaximumprinciple} and not $\partial_t + \Delta$ as in the proof of Corollary \ref{cor:Donaldson_Kronheimer_6-2-12}.
\end{rmk}

While Corollary \ref{cor:Donaldson_Kronheimer_6-2-12} provides a useful \emph{global} upper bound for $u$ over all of $X_T$, albeit with right-hand side depending unfavorably on $T$, we can also apply Lemma \ref{lem:Donaldson_Kronheimer_6-2-11} to give a bound for $u$ on $\sU'_T := (0,T)\times \sU'$ which is \emph{local} but with right-hand side potentially \emph{independent of $T$}.

\begin{cor}
\label{cor:Donaldson_Kronheimer_6-2-12_local}
Assume the hypotheses of Lemma \ref{lem:Donaldson_Kronheimer_6-2-11} and, in addition, that $A(t)$ is continuous up to $t=0$. If $u$ is as in \eqref{eq:Donaldson_Kronheimer_6-2-10}, then
\begin{equation}
\label{eq:Donaldson_Kronheimer_6-2-12_local}
u(t,x) \leq  \sup_{(0,T)\times\partial \sU'} u \vee \sup_{\sU'} u(0),
\quad\forall\, (t,x)\in X_T.
\end{equation}
\end{cor}

\begin{rmk}[On applications of Corollary \ref{cor:Donaldson_Kronheimer_6-2-12_local}]
\label{rmk:Application_corollary_Donaldson_Kronheimer_6-2-12_local}
The importance of the \apriori estimate \eqref{eq:Donaldson_Kronheimer_6-2-12_local} in our applications is that, when $T=\infty$, we will have uniform bound for $u$ on $(0,\infty)\times\partial \sU'$.
\end{rmk}

\begin{proof}
The conclusion follows from the maximum principle in the shape of \cite[Proposition 2.6 (Item 1)]{Feehan_parabolicmaximumprinciple}, noting that Lemma \ref{lem:Donaldson_Kronheimer_6-2-11} now yields
$$
Lu \leq 0 \quad\hbox{on }\sU'_T,
$$
since $R \equiv 0$ and thus $f \equiv 0$ on $\sU'_T$, and noting that $\mydirac \sU'_T = \{0\}\times \bar \sU' \cup (0,T)\times \partial \sU'$ is the parabolic boundary in this case.
\end{proof}

We define, for a one-parameter family of connections, $A(t)$ on a Hermitian vector bundle $E$ over $X$ (which can be any Riemannian manifold in this definition) for $0\leq t <T$,
\begin{equation}
\label{eq:Donaldson_Kronheimer_6-3-10}
\delta(r)
:=
\sup_{\begin{subarray}{c} x\in X,\\ t\in [0,T)\end{subarray}}
\int_{B_r(x)} |F_A(t)|^2 \, d\vol_g, \quad\forall\, r \in (0,\Inj(X)].
\end{equation}
We have the following analogue of \cite[Proposition 6.3.11]{DK}.

\begin{prop}[Continuous extension for a solution to Yang-Mills gradient-like flow over a compact, K\"ahler surface]
\label{prop:Donaldson_Kronheimer_6-3-11}
Let $E$ be a Hermitian vector bundle over a compact, K\"ahler surface, $X$. Let $A(t)$, for $t\in [0,T)$, be a solution to the Yang-Mills gradient-like flow equation \eqref{eq:Yang-Mills_gradient-like_flow_Kaehler_surface} that is smooth over $(0,T)\times X$ and continuous over $[0,T)\times X$. If $\delta(r) \to 0$ as $r\to 0$, then $A(t)$ extends to a smooth solution to \eqref{eq:Yang-Mills_gradient-like_flow_Kaehler_surface} on an interval $0\leq t<T+\eta$, for some $\eta > 0$.
\end{prop}

\begin{proof}
The proof is outlined for the case of pure Yang-Mills gradient flow in \cite[p. 236]{DK} and carries over \mutatis to the case of Yang-Mills gradient-like flow.
\end{proof}

Before proceeding to verify the hypothesis in Proposition \ref{prop:Donaldson_Kronheimer_6-3-11}, we recall from \cite[p. 235]{DK} that if $A(t)$ is a solution to \emph{pure} Yang-Mills gradient flow \eqref{eq:Yang-Mills_gradient_flow_complex_surface} over $X$ with initial data $A_0 \in \sA^{1,1}(E)$ as in \eqref{Yang-Mills_gradient_flow_initial_data_type_1-1}, so $F_{A_0}$ has type $(1,1)$, we can use \eqref{eq:Yang-Mills_gradient_flow_complex_surface} to compute $\partial F_A/\partial t$ and see that
\begin{equation}
\label{eq:Yang-Mills_gradient_flow_implies_FA_type_1-1_all_time}
F_A(t) \in \Omega^{1,1}(X;\fg_E), \quad\forall\, t\in [0,T).
\end{equation}
One can then define a one-parameter family of complex gauge transformations\footnote{In \cite[Chapter 6]{DK}, the authors use $g$ to denote a complex gauge transformation, but this conflicts with our extensive use of $g$ for a Riemannian metric.}, $\varphi_t=\varphi(t)$ for $t\in[0,T)$, by solving the ordinary differential equation,
\begin{equation}
\label{eq:Donaldson_Kronheimer_6-3-7}
\frac{\partial \varphi}{\partial t}
=
i\widehat F_A \varphi\quad\hbox{on } X, \quad \varphi(0) = \id_E,
\end{equation}
and obtain
\begin{equation}
\label{eq:Yang-Mills_gradient_flow_A_evolves_by_complex_gauge_transformations}
A(t) = \varphi_t(A_0) \quad\hbox{on } E, \quad\forall\, t\in [0,T).
\end{equation}
Returning now to the case where $A(t)$ is a solution to Yang-Mills gradient-\emph{like} flow \eqref{eq:Yang-Mills_gradient-like_flow_Kaehler_surface} over $X$ with initial data satisfying the conditions \eqref{Yang-Mills_gradient-like_flow_initial_data_product_or_type_1-1}
we observe that \eqref{eq:Yang-Mills_gradient_flow_implies_FA_type_1-1_all_time} and \eqref{eq:Yang-Mills_gradient_flow_A_evolves_by_complex_gauge_transformations} will be replaced by their local counterparts,
\begin{align}
\label{eq:Yang-Mills_gradient_flow_implies_FA_type_1-1_all_time_local}
F_A(t) &\in \Omega^{1,1}(\sU';\fg_E), \quad\forall\, t\in [0,T),
\\
\label{eq:Yang-Mills_gradient_flow_A_evolves_by_complex_gauge_transformations_local}
A(t) &= \varphi_t(A_0) \quad\hbox{on } E \restriction \sU', \quad\forall\, t\in [0,T).
\end{align}
Just as in \cite[pp. 236--7]{DK}, the proof of Lemma \ref{lem:Donaldson_Kronheimer_6-3-11_hypothesis_holds} below relies on the

\begin{lem}
\label{lem:Donaldson_Kronheimer_6-3-13}
\cite[Lemma 6.3.13]{DK}
Suppose that $A$ is a finite-energy anti-self-dual connection on the product bundle, $\CC^2\times\CC^n$, over $\CC^2$ which can be written as $\varphi(\Theta)$, where $\Theta$ is the product connection on $\CC^2\times\CC^n$, for a complex gauge transformation, $\varphi:\CC^2 \to \GL(n,\CC)$, with $\varphi$ and $\varphi^{-1}$ bounded. Then $A$ is a flat connection.
\end{lem}

In \cite[Lemma 6.3.13]{DK}, it is assumed that $n=2$, but this is only done for consistency with the applications to differential topology in \cite{DK} which require $E$ to have complex rank two and is not required in the proof of \cite[Lemma 6.3.13]{DK}. We have the following analogue of the conclusion asserted in \cite[pp. 236--7]{DK}.

\begin{lem}
\label{lem:Donaldson_Kronheimer_6-3-11_hypothesis_holds}
Let $E$ be a Hermitian vector bundle over a compact, K\"ahler surface, $X$. Let $A(t)$, for $t\in [0,T)$, be a solution to the Yang-Mills gradient-like flow equation \eqref{eq:Yang-Mills_gradient-like_flow_Kaehler_surface} that is smooth over $(0,T)\times X$ and continuous over $[0,T)\times X$. If $A(t)$ obeys
\begin{subequations}
\label{eq:Yang-Mills_gradient-like_flow_A_minus_Gamma_and_FA_uniform_local_bounds}
\begin{align}
\label{eq:Yang-Mills_gradient-like_flow_A_minus_Gamma_L-infinity_bounded_annulus}
\sup_{t\in [0,T)}\|A - A_1\|_{L^\infty(\sU\less \sU')} < \infty,
\\
\label{eq:Yang-Mills_gradient-like_flow_FA_L-infinity_bounded_local}
\sup_{t\in [0,T)}\|F_A(t)\|_{L^\infty(X\less \sU'')} < \infty,
\end{align}
\end{subequations}
where $A_1$ is a fixed $C^\infty$ reference connection on $E$, then $\delta(r) \to 0$ as $r\to 0$.
\end{lem}

\begin{proof}
Since $R$ obeys \eqref{Yang-Mills_gradient-like_flow_error}, we know in particular that
$$
R \in C^1([0,T]\times X; \Lambda^1\otimes\fg_E).
$$
Consequently, the definition \eqref{eq:Definition_f_in_terms_A_and_R} of $f$ and the conditions \eqref{eq:Yang-Mills_gradient-like_flow_A_minus_Gamma_and_FA_uniform_local_bounds} ensure that
$$
f \in C_b([0,T)\times X).
$$
Therefore, the role of \cite[Corollary 6.2.12]{DK} in the argument in \cite[pp. 236--7]{DK} can be replaced by that of Corollary \ref{cor:Donaldson_Kronheimer_6-2-12}.

Because $A(t)$ obeys \eqref{eq:Yang-Mills_gradient_flow_implies_FA_type_1-1_all_time_local} and \eqref{eq:Yang-Mills_gradient_flow_A_evolves_by_complex_gauge_transformations_local} over $\sU'$, the argument in the last two paragraphs of \cite[p. 236]{DK} and first paragraph of \cite[p. 237]{DK} carries over without change for balls $B_r(x) \subset \sU'$. On the other hand, our hypothesis \eqref{eq:Yang-Mills_gradient-like_flow_FA_L-infinity_bounded_local} immediately takes care of balls $B_r(x) \subset X\less \sU''$ in the expression \eqref{eq:Donaldson_Kronheimer_6-3-10} for $\delta(r)$. This concludes the proof of Lemma \ref{lem:Donaldson_Kronheimer_6-3-11_hypothesis_holds}.
\end{proof}

We can now combine Proposition \ref{prop:Donaldson_Kronheimer_6-3-11} and Lemma \ref{lem:Donaldson_Kronheimer_6-3-11_hypothesis_holds} to give the desired continuous extension property for Yang-Mills gradient-like flow on a Hermitian vector bundle over a compact, K\"ahler surface.

\begin{thm}[Continuous extension for a solution to Yang-Mills gradient-like flow over a compact, K\"ahler surface]
\label{thm:Donaldson_Kronheimer_6-3-11_complete}
Let $E$ be a Hermitian vector bundle over a compact, K\"ahler surface, $X$. Let $A(t)$, for $t\in [0,T)$, be a solution to the Yang-Mills gradient-like flow equation \eqref{eq:Yang-Mills_gradient-like_flow_Kaehler_surface} that is smooth over $[0,T)\times X$. If $A(t)$ obeys \eqref{eq:Yang-Mills_gradient-like_flow_A_minus_Gamma_and_FA_uniform_local_bounds}, then $A(t)$ extends to a smooth solution to \eqref{eq:Yang-Mills_gradient-like_flow_Kaehler_surface} on an interval $0\leq t<T+\eta$, for some $\eta > 0$.
\end{thm}

\subsection{Application to the proof of Theorem \ref{mainthm:Yang-Mills_gradient_flow_global_existence_and_convergence_started_with arbitrary_initial_energy}}
Given Theorem \ref{thm:Donaldson_Kronheimer_6-3-11_complete}, we are at last ready to provide the

\begin{proof}[Proof of Theorem \ref{mainthm:Yang-Mills_gradient_flow_global_existence_and_convergence_started_with arbitrary_initial_energy}]
We shall argue by contradiction, so we suppose that the maximal lifetime, $\sT$, of a solution $A(t)$ to Yang-Mills gradient flow \eqref{eq:Yang-Mills_gradient_flow_complex_surface} is finite. Theorem \ref{thm:Kozono_Maeda_Naito_5-1} then provides a set of points, $\Sigma = \{x_1,\ldots,x_L\} \subset X$, where the curvature, $F_A(t)$, concentrates as $t \nearrow \sT$.

Suppose first that the given complex Hermitian metric $h$ on $X$ is \emph{K\"ahler} on small open neighborhoods $\sU_i\subset X$ of the points $x_i$ and that there are biholomorphisms,
\begin{equation}
\label{eq:Local_biholomorphism_preserving_Kaehler_metrics}
\psi_i: X \supset \sU_i \cong \phi_i(\sU_i) \subset Z_i, \quad\hbox{for } 1 \leq i \leq L,
\end{equation}
from $\sU_i$ onto open neighborhoods of the points $y_i = \psi_i(x_i) \in Z_i$, where each $Z_i$ is a compact, complex surface with K\"ahler metric $h_i$ such that $h=\psi_i^*h_i$ on $\sU_i$. For small enough $\varrho \in (0,\Inj(X,h)]$, the open neighborhoods, $\sU_i$, contain geodesic balls $B_\varrho(x_i) \subset X$ for $1\leq i\leq L$. We now cut off $A(t)$ over each the annuli, $\Omega(x_i;\varrho/2,\varrho/4)$, as described in Remark \ref{rmk:Applications_Yang-Mills_gradient-like_flow_Kaehler_surface}, for $1\leq i\leq L$ and choose $\varrho \in (0,\Inj(X,g)]$ small enough that $\varrho \leq \min_{i\neq j}\dist_g(x_i,x_j)/4$, where $g$ is the real Riemannian metric on $X$ corresponding to $h$.

It is enough to consider one of the points in $\Sigma$, which we relabel as $x_0$, and consider the Yang-Mills gradient-like flow, $A_\chi(t)$ on $E_\chi$ over $X$, which is given by $A(t)$ on $E_\chi := E \restriction B_{\varrho/2}(x_0)$ and the product connection, $\Gamma$, on
$E_\chi := X \less B_\varrho(x_0) \times \CC^n$. (This construction also depends on a choice of trivialization
$E\restriction \Omega(x_0;\varrho/2,\varrho) \cong \Omega(x_0;\varrho/2,\varrho)\times \CC^n$, which we suppress from our notation.) Because $A_\chi(t) = \Gamma$ over $X \less B_\varrho(x_0)$ for all $t\in [0,\sT)$, we can equally well regard $A_\chi(t)$ as a flow on a Hermitian vector bundle, $E_\chi$, over a compact, K\"ahler surface, $Z$, containing a copy of the ball $B_\varrho(x_0)$ via the local identification \eqref{eq:Local_biholomorphism_preserving_Kaehler_metrics}.

Over the annulus $\Omega(x_0;\varrho/2,\varrho)\subset X$, the Yang-Mills gradient flow, $A(t)$, on $E$ and hence the Yang-Mills gradient-like flow, $A_\chi(t)$, on $E_\chi$ over $\Omega(x_0;\varrho/2,\varrho)\subset Z$ obey the conditions \eqref{eq:Yang-Mills_gradient-like_flow_A_minus_Gamma_and_FA_uniform_local_bounds} for
$t\in [0,\sT)$ by Theorem \ref{thm:Kozono_Maeda_Naito_5-3}. Moreover, the error term $R(t)$ in \eqref{eq:Yang-Mills_gradient-like_flow_Kaehler_surface} obeys the condition
$$
R \in C^\infty([0,\sT+\eta)\times Z; \Lambda^1\otimes\fg_{E_\chi})),
$$
for some positive constant $\eta$, which, although weaker than \eqref{Yang-Mills_gradient-like_flow_error}, is all that is required for the proof of Theorem \ref{thm:Donaldson_Kronheimer_6-3-11_complete}.

Consequently, Theorem \ref{thm:Donaldson_Kronheimer_6-3-11_complete} implies that $A_\chi(t)$ to extends to a solution to the Yang-Mills gradient-like flow equation \eqref{eq:Yang-Mills_gradient-like_flow_Kaehler_surface} on $E_\chi$ that is smooth over $[0,\sT+\eta)\times X$ for some (possibly smaller) positive constant $\eta$. But this means that $A(t)$ extends to a solution to the Yang-Mills gradient flow equation \eqref{eq:Yang-Mills_gradient_flow_complex_surface} on $E$ that is smooth over $[0,\sT+\eta)\times B_{\varrho/2}(x_0)$ for some positive constant $\eta$ and, since $x_0\in\Sigma$ was arbitrary, this holds for any point $x_i\in\Sigma$. Thus, $A(t)$ extends to a solution to the Yang-Mills gradient flow equation \eqref{eq:Yang-Mills_gradient_flow_complex_surface} on $E$ that is smooth over $[0,\sT+\eta)\times X$ for some positive constant $\eta$, contradicting the maximality of $\sT$.

If the given Hermitian metric $h$ on $X$ is \emph{not} K\"ahler near each of the points $x_i \in \Sigma$, we can apply Remark \ref{rmk:Stability_bubbling_wrt_local_flattening_Riemannian_metric_choice_r_and_rho} and Corollary \ref{cor:Bubbling_for_g_implies_bubbling_for_nearby_barg_weakly_locally_flattened_finite_number_small_balls}
(with $r=\varrho$) to assure us that if $A(t)$ has maximal lifetime $\sT$ with bubble singularity set $\Sigma=\{x_1,\ldots,x_L\}$, for the given Hermitian metric $h$ and initial data $A_0$, then $\widetilde A(t)$ has the same maximal lifetime, $\sT$, with at least \emph{one} bubble point in a ball $B_{\varrho/\sN}(x_i)$, for some $x_i \in \Sigma$ and $\sN\geq 4$, and initial data $A_0$ and a perturbed Hermitian metric $\tilde h$ on $X$ such that $\tilde h$ is \emph{K\"ahler} on $\cup_{i=1}^L B_\varrho(x_i)$. Writing $\varrho = \rho/N$, for $\rho\in (0,\Inj(X,h)]$ and suitable $N\geq 4$ (the size of the constants $\rho$ and $N$ are determined by Lemma \ref{lem:Bound_norm_difference_g_minus_barg_Hermitian} and Corollary \ref{cor:Bubbling_for_g_implies_bubbling_for_nearby_barg_weakly_locally_flattened_finite_number_small_balls}
the perturbed Hermitian metric $\tilde h$ on $X$ is chosen so that $\tilde h=h$ on $X\less \cup_{i=1}^L B_{\rho/2}(x_i)$ and $\tilde h(x_i) = h(x_i)$ for $i \in \{1,\ldots,L\}$.

To actually construct $\tilde h$, consider one of the points $x_0 \in \Sigma$. We choose local holomorphic coordinates, $\{z^\alpha\}$, on an open neighborhood, $\sU$, of the point $x_0$ and thus obtain a biholomorphism,
$$
\psi: X \supset \sU \cong \psi(\sU) \supset \CC\PP^1\times\CC\PP^1,
$$
where $\CC\PP^1\times\CC\PP^1$ has its standard complex structure, product Fubini-Study metric, $h^0$, and $\phi(x_0) = (0,0) \in \CC\PP^1\times\CC\PP^1$. The local holomorphic coordinates, $\{z^\alpha\}$, may be chosen such that $h_{\alpha\bar\beta}(x_0) = \delta_{\alpha\bar\beta} = h_{\alpha\bar\beta}^0(0,0)$, where $\delta_{\alpha\bar\beta}$ are the components of the standard flat K\"ahler metric on $\CC^2$ with respect to the coordinates $\{z^\alpha\}$. We now define a Hermitian metric, $\tilde h$, on $B_\rho(x_0)$ by the interpolation construction in Definition \ref{defn:Locally_flattened_Riemannian_metric}, so that
$$
\tilde h
=
\begin{cases}
\psi^*h^0 &\hbox{on } B_{\rho/N}(x_0),
\\
h &\hbox{on } \Omega(x_0;\rho/2,\rho).
\end{cases}
$$
Repeating this construction for each point $x_0\in \Sigma$, we obtain a Hermitian metric on $X$ that on each ball $B_\varrho(x_i) = B_{\rho/N}(x_i)\subset X$ is the biholomorphic pull-back of the standard K\"ahler metric on $Z=\CC\PP^1\times\CC\PP^1$ via $\phi_i$ as in \eqref{eq:Local_biholomorphism_preserving_Kaehler_metrics}. Lemma \ref{lem:Bound_norm_difference_g_minus_barg_Hermitian} ensures that $\tilde h$ is close enough to $h$ that the hypotheses of Corollary \ref{cor:Bubbling_for_g_implies_bubbling_for_nearby_barg_weakly_locally_flattened_finite_number_small_balls}
are satisfied.

We can now appeal to the argument for the case when the Hermitian metric, $h$, on $X$ was assumed to be the biholomorphic pull-back on the balls $B_\varrho(x_i)$ for $1\leq i\leq L$ of a K\"ahler metric on a compact, complex surface, $Z$, to establish global existence for Yang-Mills gradient flow $A(t)$. Uniqueness of the set, $\Sigma$, follows from Corollary \ref{cor:Kozono_Maeda_Naito_5-3_T_is_infinite}.

Finally, the assertions regarding convergence of $A(t_m)$, modulo gauge transformations, for a subsequence $\{t_m\}_{m\in\NN} \subset [0,\infty)$ follow immediately from Theorem \ref{thm:Kozono_Maeda_Naito_5-3}. This completes the proof of Theorem \ref{mainthm:Yang-Mills_gradient_flow_global_existence_and_convergence_started_with arbitrary_initial_energy}.
\end{proof}

\begin{rmk}[On the approximation of Hermitian metrics]
Our first version of Theorem \ref{thm:Global_apriori_estimate_difference_solutions_Yang-Mills_heat_equations_pair_metrics} assumed a convenient but rather strong condition on the pair of nearby Riemannian metrics, $g$ and $\tilde g$, in particular that $\|\tilde g-g\|_{C^1(X)}$ could be as small as desired while retaining a uniform bound on $\|\tilde g-g\|_{C^2(X)}$. For a given Riemannian metric, $g$, on $X$, we constructed $\tilde g$ using (real) normal geodesic coordinates, $\{x^\mu\}$, on neighborhoods of each of the points $x_i \in \Sigma$ to replace $g$ by the standard Euclidean metric on a finite collection of balls, $B_\rho(x_i) \subset X$, where $g_{\mu\nu}(x_i) = \delta_{\mu\nu}$ and $(\partial g_{\mu\nu}/\partial x^\eta)(x_i) = 0$, for all $\mu,\nu,\eta \in \{1,\ldots,4\}$. However, in the case of a complex manifold $(X,J)$ with Hermitian $h$, the analogous approximation condition with respect to local holomorphic coordinates, $\{z^\alpha\}$, is equivalent to the assumption that $h$ is actually K\"ahler. See, for example, \cite[Lemma, p. 107]{GriffithsHarris}, \cite[Theorem 11.6]{Moroianu_2007}, or \cite[Remark, p. 196]{Wells3}: \emph{A Hermitian metric, $h$, on a complex manifold, $(X, J)$, is K\"ahler if and only if around each point of $X$ there exist holomorphic coordinates
in which $h$ osculates to the standard Hermitian metric to order $2$}.
\end{rmk}

\begin{rmk}[On non-K\"ahler complex surfaces]
Interest in non-K\"ahler, complex Hermitian manifolds has grown in recent years, in part due to motivations from theoretical physics. Basic results on Hermitian manifolds can be found in the text by Yano \cite{Yano_1965}. The class of compact, locally conformal K\"ahler manifolds lies between the class of K\"ahler and complex Hermitian manifolds and their properties are described in the monograph by Dragomir and Ornea \cite{Dragomir_Ornea_1998}. Just a small selection of many recent results is provided by Liu and Yang \cite{Liu_Yang_2014arxiv, Liu_Yang_2012} and Streets and Tian  \cite{Streets_Tian_2010} and references contained therein.
\end{rmk}

The following result is proved using an argument similar to that of \cite[Theorem 1]{Daskalopoulos_Wentworth_2007} and, with indicated convergence hypothesis, gives a partial characterization of certain analytic multiplicities associated with points, $x_i$, in the set of bubble points, $\Sigma \subset X$, in Theorem \ref{mainthm:Yang-Mills_gradient_flow_global_existence_and_convergence_started_with arbitrary_initial_energy}. The hypotheses on the regularity of the initial data, $A_0$, and the convergence of $\Phi(t)^*A(t)$ over $X\less\Sigma$ may be relaxed, albeit at the expense of an increase in technicalities.

\begin{lem}[Characterization of analytic multiplicities]
\label{lem:Daskalopoulos_Wentworth_2007_lemma_5_analogue}
Let $E$ and $E_\infty$ be Hermitian vector bundles over a closed, connected, four-dimensional, oriented, smooth manifold, $X$, with Riemannian metric, $g$. Let $A_0$ and $A_\infty$ be $C^\infty$ unitary connections on $E$ and $E_\infty$, respectively, and $\Sigma = \{x_1,\ldots,x_L\} \subset X$ a finite set of points such that $E\restriction X\less\Sigma \cong E_\infty\restriction X\less\Sigma$. Let $A(t)$, for $t\in [0,\infty)$, be a solution to the Yang-Mills gradient flow \eqref{eq:Yang-Mills_gradient_flow} with initial data, $A(0) = A_0$, such that
$$
A - A_0 \in C^\infty([0,\infty)\times X; \Lambda^1\otimes\fg_E),
$$
and $\Phi(t) \in \Aut E$, for $t \geq 0$, a $C^\infty$ path of $C^\infty$ automorphisms of the Hermitian vector bundle $E$ such that,
as $t \to \infty$,
$$
\Phi(t)^*A(t) \to A_\infty
\quad\hbox{in } H_{A_\infty,\loc}^2(X \less \Sigma; \Lambda^1\otimes\fg_E).
$$
Then the \emph{analytic multiplicity},
$$
\mu_{\mathrm{an}}(x_i)
:=
\lim_{t\to\infty}
\int_{B_\rho(x_i)} \tr(F_A(t)\wedge F_A(t)) - \tr(F_{A_\infty}\wedge F_{A_\infty}) \in \ZZ,
$$
associated with each point, $x_i \in \Sigma$, is well-defined and independent of
$$
\rho \in \left(0,\frac{1}{4}\min_{k\neq l}\dist_g(x_k,x_l)\wedge\Inj(X,g)\right].
$$
\end{lem}

\begin{rmk}[Extension to compact Lie groups]
\label{rmk:Daskalopoulos_Wentworth_2007_lemma_5_analogue_compact_Lie_group}
A similar result holds, more generally, for any compact Lie group, as is clear from the discussion in Section \ref{sec:Taubes_1982_Appendix}. We restrict our attention to the case where $G=\U(n)$ for notational convenience in the proof and for the sake of consistency with \cite[Section 2.1.4]{DK}.
\end{rmk}

\begin{proof}[Proof of Lemma \ref{lem:Daskalopoulos_Wentworth_2007_lemma_5_analogue}]
Consider one of the points, $x_i \in \Sigma$, and choose $\rho \in (0, \Inj(X,g)]$ small enough that $\dist_g(x_i,x_j) \geq 4\rho$ for all $j\neq i$. For each $r \in (0,\rho)$, define a $C^\infty$ unitary connection, $A_\infty^r$, on $E$ such that $A_\infty^r = A_\infty$ on $E \restriction X\less \cup_{j=1}^L B_r(x_j)$ and $A_\infty^r = A_0$ on $E \restriction \cup_{j=1}^L B_{r/2}(x_j)$. This may be accomplished by choosing orthonormal frames, $p_j$, for the fibers, $E_{x_j}$, defining trivializations, $E\restriction B_\rho(x_j) \cong B_\rho(x_j)\times \CC^n$ for $1\leq j\leq L$ using $p_j$ and $A_0$, choosing smooth cut-off functions over the annuli, $\Omega(x_j;r/2,r)$, and using the cut-off functions and trivializations to splice the connections $A_\infty$ and $A_0$ over those annuli.

Denoting the trace on complex $n\times n$ matrices by $\tr(\cdot)$, we recall from \cite[Equation (2.1.28]{DK} that, for any $C^\infty$ unitary connection $A$ on $E$,
$$
c_2(E) - \frac{1}{2}c_1(E)^2 = \frac{1}{8\pi^2} [\tr(F_A\wedge F_A)] \in H^4(X;\ZZ),
$$
and
$$
\left\langle c_2(E) - c_1(E)^2/2, [X] \right\rangle
=
\frac{1}{8\pi^2} \int_X \tr(F_A\wedge F_A) \in \ZZ.
$$
The four-form,
$$
\tr(F_A\wedge F_A) \in \Omega^4(X),
$$
is \emph{closed} \cite[p. 39]{DK}. In particular, writing $\Phi^*A(t) = A_\infty^r + a(t)$ over $X$, then \cite[Equation (2.1.27]{DK} gives, for all $t\geq 0$,
\begin{align*}
\tr(F_A(t)\wedge F_A(t)) - \tr(F_{A_\infty^r}\wedge F_{A_\infty^r})
&=
\tr(F_{\Phi^*A}(t)\wedge F_{\Phi^*A}(t)) - \tr(F_{A_\infty^r}\wedge F_{A_\infty^r})
\\
&= d\left(\tr(a(t)\wedge d_{A_\infty^r}a(t)) + \frac{2}{3}a(t)\wedge a(t) \wedge a(t) \right).
\end{align*}
Therefore, noting that $A_\infty^r = A_\infty$ on $\partial B_\rho(x_i)$, we obtain
\begin{multline*}
\int_{B_\rho(x_i)} \tr(F_A(t)\wedge F_A(t)) - \tr(F_{A_\infty^r}\wedge F_{A_\infty^r})
\\
= \int_{\partial B_\rho(x_i)}
\left(\tr(a(t)\wedge d_{A_\infty}a(t)) + \frac{2}{3}a(t)\wedge a(t) \wedge a(t) \right).
\end{multline*}
If $Y \subset X$ is a $C^\infty$, closed, three-dimensional submanifold, it follows from \cite[Theorem 4.12]{AdamsFournier} that there are continuous embeddings,
$$
H^2(X) \hookrightarrow W^{1,3}(Y) \quad\hbox{and}\quad W^{1,3}(Y) \hookrightarrow L^q(Y),
\quad 1 \leq q < \infty,
$$
and similarly for the Sobolev spaces, $H_{A_0}^2(X; \Lambda^1\otimes\fg_E)$, and $W_{A_\infty}^{1,3}(X; \Lambda^1\otimes\fg_E)$, and $L^q(X; \Lambda^1\otimes\fg_E)$. Moreover, similar embeddings hold if $X$ is replaced by a tubular neighborhood, $U$, of $Y$. Consequently, we see that
\begin{align*}
{} & \left|\int_{\partial B_\rho(x_i)}
\left(\tr(a(t)\wedge d_{A_\infty}a(t)) + \frac{2}{3}a(t)\wedge a(t) \wedge a(t) \right) \right|
\\
&\leq c_\rho \left(\|a(t)\|_{L^{3/2}(\partial B_\rho(x_i))}
\|d_{A_\infty}a(t)\|_{L^3(\partial B_\rho(x_i))}
+ \|a(t)\|_{L^3(\partial B_\rho(x_i))}^3 \right)
\\
&\leq c_\rho \left(\|a(t)\|_{W_{A_\infty}^{1,3}(\partial B_\rho(x_i))}^2
+ \|a(t)\|_{L^3(\partial B_\rho(x_i))}^3 \right)
\\
&\leq C \left(\|a(t)\|_{H_{A_\infty}^2(\Omega(x_i;\rho/2,2\rho))}^2
+ \|a(t)\|_{H_{A_\infty}^2(\Omega(x_i;\rho/2,2\rho))}^3 \right),
\end{align*}
where $c_\rho$ is a positive constant depending at most on $g$, $n$, and $\rho$, while $C = C(A_\infty,g,n,\rho)$. Our convergence hypothesis on $A(t)$ thus ensures that
$$
\lim_{t\to\infty}
\int_{B_\rho(x_i)} \tr(F_A(t)\wedge F_A(t)) - \tr(F_{A_\infty^r}\wedge F_{A_\infty^r})
= 0, \quad\forall\, r \in (0,\rho).
$$
The resulting identity,
\begin{multline}
\label{eq:Daskalopoulos_Wentworth_2007_lemma_6_analogue}
\lim_{t\to\infty}
\int_{B_\rho(x_i)}
\left( \tr(F_A(t)\wedge F_A(t)) - \tr(F_{A_\infty}\wedge F_{A_\infty}) \right)
\\
=
\int_{B_\rho(x_i)}
\left( \tr(F_{A_\infty^r}\wedge F_{A_\infty^r}) - \tr(F_{A_\infty}\wedge F_{A_\infty}) \right),
\quad\forall\, r \in (0,\rho),
\end{multline}
interprets $\mu_{\mathrm{an}}(x_i)$ in terms of $A_0$ and $A_\infty$.
The convergence property of $A(t)$ ensures that the left-hand side of the preceding identity (and hence the right-hand side) is independent of $\rho$ in the stated range, while the right-hand side is independent of $r$ in the stated range because that is manifestly true for the left-hand side.
\end{proof}

\subsection{Proof of Corollary \ref{cor:Donaldson_Kronheimer_6-2-7_and_6-2-14_convergence}}
\label{subsec:Proof_corollary_Donaldson_Kronheimer_6-2-7_and_6-2-14_convergence}
Lastly, we can give the

\begin{proof}[Proof of Corollary \ref{cor:Donaldson_Kronheimer_6-2-7_and_6-2-14_convergence}]
The conclusion follows from Theorem \ref{mainthm:Donaldson_Kronheimer_6-2-7_and_6-2-14_plus} and our Theorem \ref{thm:Huang_3-3-6_Yang-Mills}. Indeed, Theorem \ref{mainthm:Donaldson_Kronheimer_6-2-7_and_6-2-14_plus} implies that, for a sufficiently large integer, $m_0$, we have
$$
\|A(t_{m_0}) - A_\infty'\|_{H_{A_\infty'}^1(X)} < \eps,
$$
where $\eps$ is the positive constant in Theorem \ref{thm:Huang_3-3-6_Yang-Mills} and $A_\infty' = (\Phi_{m_0}^{-1})^*A_\infty$. Moreover,
$$
\sE(A(t)) \geq \lim_{t\to\infty}\sE(A(t)) = \sE(A_\infty'), \quad\forall\, t \geq 0,
$$
and so the second alternative in Theorem \ref{thm:Huang_3-3-6_Yang-Mills} necessarily holds. Therefore,
$$
\|A(t) - A_\infty'\|_{H_{A_\infty'}^1(X)} \to 0 \quad\hbox{as } t \to \infty,
$$
and the conclusion follows by relabeling the limiting connection, $A_\infty$, in Theorem \ref{mainthm:Donaldson_Kronheimer_6-2-7_and_6-2-14_plus}.
\end{proof}

As we noted in Remark \ref{rmk:Donaldson_Kronheimer_6-2-7_and_6-2-14_convergence_mod_gauge}, it is possible to prove a weaker version of Corollary \ref{cor:Donaldson_Kronheimer_6-2-7_and_6-2-14_convergence}, namely convergence modulo a $C^\infty$ path of $C^\infty$ gauge transformations without appeal Theorem \ref{thm:Huang_3-3-6_Yang-Mills}. To see this, we first claim that for \emph{every} subsequence, $\{t_m\}_{m\in\NN} \subset (0, \infty)$ with $t_m \to \infty$ as $m \to \infty$, there is a sequence of gauge transformations, $\{\Phi_m\}_{m\in\NN} \subset \sG_E$, such that the sequence, $\{\Phi_m^*A(t_m)\}_{m\in\NN}$, converges to $A_\infty$ over $X$ as $m \to \infty$ in the sense of $H_{A_0}^1(X;\Lambda^1\otimes\fg_E)$. If not, then there is a sequence, $\{t_m\}_{m\in\NN} \subset (0, \infty)$, such that for \emph{every} sequence, $\{\Phi_m\}_{m\in\NN} \subset \sG_E$,
\begin{equation}
\label{eq:Wentworth_explanation_theorem_1}
\|\Phi_m^*A(t_m)-A_\infty\|_{H_{A_0}^1(X)} \geq \eps, \quad\forall\, m\in\NN,
\end{equation}
for \emph{some} constant, $\eps \in (0,1]$. However, by \cite[Section 6.2.4]{DK} and the fact that $\Sigma = \emptyset$,
$$
\|\Psi_m^*A(t_m)-A_\infty\|_{H_{A_0}^2(X)} \leq K, \quad\forall\, m\in\NN,
$$
for \emph{some} constant, $K \in [1,\infty)$, and \emph{some} sequence, $\{\Psi_m\}_{m\in\NN} \subset \sG_E$. (The preceding pair of inequalities can be expressed more elegantly using distance functions on $\sB(E) = \sA(E)/\sG_E$, as in \cite{FeehanSlice}.) According to \cite{UhlLp},
we can pass to a subsequence, relabeled as $\{m\}$, such that $\{\Psi_m^*A(t_m)\}_{m\in\NN}$, converges to a limit $A_\infty'$ on the bundle $E$ over $X$ as $m \to \infty$ strongly in the sense of $W_{A_0}^{1,p}(X;\Lambda^1\otimes\fg_E)$ for any $p \in [1,4)$. The connection, $A_\infty'$, is anti-self-dual and, by the uniqueness assertion of Theorem \ref{mainthm:Donaldson_Kronheimer_6-2-7_and_6-2-14_plus}, is equivalent to $A_\infty$ by an element of $\sG_E$, contradicting \eqref{eq:Wentworth_explanation_theorem_1}.

Consequently, we have
$$
\dist_{H_{[A_0]}^2(X)}\left([A(t)], [A_\infty]\right) \to 0, \quad\hbox{as } t \to \infty,
$$
and similarly if $A_0$ is replaced by $A_\infty$. Therefore, by the Slice Theorem for $\sB_E$ (see \cite[Theorem 1.1]{FeehanSlice}), for large enough $T \in (0,\infty)$ and each $t\in [T,\infty)$, there is a gauge transformation $\Phi(t) \in \sG_E$, unique up to $\{\pm\id_E\}$ (since $A_\infty$ is irreducible) such that
$$
d_{A_\infty}^*\left(\Phi(t)^*A(t) - A_\infty\right) = 0, \quad\forall\, t \geq T,
$$
and
$$
\|\Phi(t)^*A(t)-A_\infty\|_{H_{A_\infty}^1(X)} \leq C\dist_{H_{[A_\infty]}^2(X)}\left([A(t)], [A_\infty]\right), \quad\forall\, t \geq T,
$$
for a constant $C = C(A_\infty,h) \in [1,\infty)$, where $h$ is the K\"ahler metric on $X$. Moreover, by the proof of \cite[Theorem 1.1]{FeehanSlice}, the family of gauge transformations, $\Phi(t) \in \sG_E$, can be shown to depend smoothly on $t\in [T,\infty)$. Thus, we obtain a version of Corollary \ref{cor:Donaldson_Kronheimer_6-2-7_and_6-2-14_convergence} with $A(t)$ replaced by $\Phi(t)^*A(t)$, for a $C^\infty$ path of $C^\infty$ gauge transformations, $\Phi(t)$.

\chapter{Solution to the anti-self-dual equation and applications}
\label{chapter:Solution_to_anti-self-dual_equation_applications}

\section{Solution to the anti-self-dual equation}
\label{sec:Solution_to_anti-self-dual_equation}
In this section, we review a particular case of the gluing theorem for anti-self-dual and self-dual connections \cite{DK, FLKM1, MorganMrowkaTube, TauSelfDual, TauIndef, TauFrame, TauGluing}. If in Theorems \ref{thm:Kozono_Maeda_Naito_5-3} and \ref{thm:Kozono_Maeda_Naito_5-3_T_is_infinite} and Theorem \ref{thm:Kozono_Maeda_Naito_5-4}, after sufficiently many iterations, the limiting connections at a singularity time, $T\leq \infty$, are all smooth anti-self-dual connections on principal $G$ bundles over $X$ or $S^4$, then the limiting process can be reversed, under suitable conditions, with the aid of the gluing theorem to produce an anti-self-dual connection, $\tilde A$ on $P$, which is $H_A^1(X)$-close to $A(t_1)$ for a time $t_1 \in [0, T)$ close enough to $T$. The analogous statement holds for self-dual connections and so it is enough to restrict our attention to the case of anti-self-dual connections.

A proof of a `gluing theorem', at least in the context of anti-self-dual connections with no small-eigenvalue obstructions to gluing, has the following principal steps:
\begin{enumerate}
  \item \emph{Perturbation of an approximate anti-self-dual connection to a solution of the anti-self-dual equation.} Given an approximate anti-self-dual connection, $A$, with sufficiently small $\|F_A^+\|$ (for some suitable norm $\|\cdot\|$) and a positive lower bound for the small eigenvalues of the self-adjoint operator, $d_A^+d_A^{+,*}$, on $L^2(X;\Lambda^+\otimes\ad P)$, solve $F_{A+a}^+ = 0$ for $a = d_A^{+,*}v \in \Omega^1(X;\ad P)$, namely,
  \begin{equation}
  \label{eq:Feehan_Leness_7-1}
  d_A^+a + (a \wedge a)^+ = - F^+_A \quad\hbox{on } X,
  \end{equation}
  that is, solve the following elliptic, quasilinear, second-order equation for $v \in \Omega^+(X;\ad P) \subset \Omega^2(X; \ad P)$,
  \begin{equation}
  \label{eq:Feehan_Leness_7-2}
  d_A^+d_A^{+,*}v + (d_A^{+,*}v \wedge d_A^{+,*}v)^+ = - F^+_A \quad\hbox{on } X;
  \end{equation}

  \item \emph{Splicing construction of an approximate anti-self-dual connection.} Construct an approximate anti-self-dual connection, $A$, by cut-and-paste splicing of an (approximate) anti-self-dual connection on a principal $G$-bundle (with different topology) over $X$ with a set of anti-self-dual connections on trees of principal $G$-bundles over $S^4$;

  \item \emph{Local \apriori estimates for anti-self-dual connections.} The estimates are applied over annuli where connections are cut off in the splicing construction and are used to give the required small bound for $\|F_A^+\|$;

  \item \emph{Positive lower bound for small eigenvalues.} The positive lower bound is for the eigenvalues of the Laplace operator $d_A^+d_A^{+,*}$ for a spliced connection, $A$, which is approximately anti-self-dual.
\end{enumerate}
A method of solving the anti-self-dual equation in the case where there is a positive lower bound for the small eigenvalues of $d_A^+d_A^{+,*}$ on $\Omega^+(X; \ad P)$ --- the only case we shall consider in this monograph --- was established by Taubes and comprises his \cite[Theorem 3.2]{TauSelfDual}, proved in \cite[Sections 4 and 5]{TauSelfDual}. Later refinements of Taubes' fundamental result hinge on
\begin{inparaenum}[\itshape a\upshape)]
\item relaxing the constraint that the small eigenvalues of $d_A^+d_A^{+,*}$ have a positive lower bound \cite{DonConn, DK, FLKM1, MorganMrowkaTube, TauIndef, TauFrame, TauStable, TauGluing, Taylor_2002};
\item weakening the measure, $\|F_A^+\|$, of deviation from anti-self-duality
\cite{FLKM1, TauStable, TauGluing}.
\end{inparaenum}

With the exception of \cite{TauIndef, Taylor_2002}, all results which relax the constraint that the small eigenvalues of $d_A^+d_A^{+,*}$ have a positive lower bound have achieved their goal by appealing to the concept underlying the Kuranishi method \cite{Kuranishi} and replacing the anti-self-dual equation with its projection onto the complement of the eigenspace defined by the small eigenvalues, so they are no longer solving the true anti-self-dual equation. The resulting solutions, `generalized' or `extended' anti-self-dual connections, are not necessarily even Yang-Mills connections, so this approach does not have an obvious application in our monograph. In \cite{TauIndef}, Taubes extended his existence results \cite{TauSelfDual} from the case where $b^+(X)=0$ to $b^+(X)>0$ (though now assuming $G=\SU(2)$ or $\SO(3)$) but at cost of splicing finitely many additional anti-self-dual $\SU(2)$-connections on the four-dimensional sphere, $S^4$, with second Chern class one (or `one-instantons'). Taubes' idea was extended by Taylor \cite{Taylor_2002} (a Ph.D. student of Mrowka) to the case where $G = \SO(n)$ and $n \geq 4$. Unfortunately, neither Taubes' nor Taylor's results for the case $b^+(X)>0$ appear to be applicable in our monograph since they involve splicing in additional instantons on $S^4$ in order to achieve existence of a solution to the anti-self-dual connection on a new principal $G$-bundle with different topology than the one supporting the given or initially constructed approximate anti-self-dual connection.

Since our monograph is concerned, more generally, with Yang-Mills connections which need not be anti-self-dual, it is natural to ask whether it is possible to solve the Yang-Mills equation (even on $S^4$) by mimicking Taubes' method for solving the anti-self-dual equation. Certainly, one can use the splicing construction to build approximate Yang-Mills connections over a closed, Riemannian, smooth manifold of arbitrary dimension using local \apriori estimates for Yang-Mills connections due to Taubes (for example, \cite[Section 9]{TauFrame} when $X$ has dimension four) and Uhlenbeck \cite{UhlRem} to estimate the splicing error. However, the perturbation step is complicated by the fact that the Hodge Laplace operator, $\Delta_A = d_Ad_A^* + d_A^*d_A$ on $\Omega^1(X; \ad P)$, will have a non-zero kernel, even in the simplest case of a Yang-Mills connection, $A$, over $S^4$. (A similar problem arises when gluing $\SO(3)$ monopoles \cite{FL3}.) Therefore, one is forced (apparently in all cases) to apply the Kuranishi method and restrict attention to the projection of the Yang-Mills equation onto the complement of the eigenspace defined by the small eigenvalues of $\Delta_A$ on $\Omega^1(X; \ad P)$. When $X$ has dimension four, this approach was implemented by Taubes in \cite[Lemma 7.3]{TauFrame} and extended to higher-dimensional manifolds by Brendle in \cite[Theorem 1.1]{Brendle_2003arxiv}. Unfortunately, neither result appears to be applicable in this monograph since there is no obvious condition (for example, on the Riemannian metric) that would provide a positive lower bound for the small eigenvalues of $\Delta_A$ on $\Omega^1(X; \ad P)$.

\subsection{Perturbation of an approximate anti-self-dual connection to a solution of the anti-self-dual equation}
\label{subsec:Perturbation_and_solution_anti-self-dual_equation}
In this subsection, we review a special case of the perturbation theorem for the anti-self-dual equation pioneered by Taubes \cite{TauSelfDual, TauIndef, TauFrame, TauGluing}. The actual version we shall recall and use in our monograph is closest to one we established with Leness in \cite{FLKM1}.



\begin{thm}[Perturbation of an approximate anti-self-dual connection to a unique solution of the anti-self-dual equation]
\label{thm:Proposition_Feehan_Leness_7-6}
\cite[Proposition 7.6]{FLKM1}
Let $X$ be a closed, four-dimensional, oriented, smooth manifold with Riemannian metric, $g$, and let $M \geq 1$ and $\Lambda \geq 1$ be constants.  Then there are a constant, $\eps = \eps(g, M, \Lambda) \in (0, 1]$, and a positive constant, $C = C(g, M, \Lambda)$, with the following significance. Let $G$ be a compact Lie group, $P$ a principal $G$-bundle over $X$, and $A$ an $H^4$ connection on $P$ with $\|F_A^+\|_{L^{\sharp,2}(X)}\leq \eps$, and $\|F_A\|_{L^2(X)}\leq M$, and $\mu$ a positive constant such that $\Lambda^{-1} \leq \mu \leq \Lambda$. Then there is a unique solution, $v\in H_A^2(X;\Lambda^+\otimes \ad P)\cap C(X; \Lambda^+\otimes\ad P)$, to equation \eqref{eq:Feehan_Leness_7-2} such that
\begin{subequations}
\label{eq:Proposition_Feehan_Leness_7-6}
\begin{align}
\label{eq:Proposition_Feehan_Leness_7-6_1}
\|v\|_{H_A^1(X)} &\leq C\|F_A^+\|_{L^{4/3}(X)},
\\
\label{eq:Proposition_Feehan_Leness_7-6_2}
\|v\|_{H_A^2(X)} + \|v\|_{C(X)} &\leq C\|F_A^+\|_{L^{\sharp,2}(X)}.
\end{align}
\end{subequations}
Moreover, if $A$ is an $H^k$ connection, for $k\geq 4$, then
$v\in H_A^{k+1}(X;\Lambda^+\otimes \ad P)$.
\end{thm}

Theorem \ref{thm:Proposition_Feehan_Leness_7-6} requires that $F_A^+$ is small in the sense that $\|F_A^+\|_{L^{\sharp,2}(X)} \leq \eps$, for a suitable $\eps \in (0, 1]$. The reason for the presence of the stronger (but still scale-invariant) measure of $F_A^+$ in Theorem \ref{thm:Proposition_Feehan_Leness_7-6} is that its use allows a statement with constants, $C,\eps$, that depend only on the energy $A$, through a bound on $\|F_A\|_{L^2(X)}$, rather than $\|F_A\|_{L^p(X)}$ as in \cite[Theorem 3.2]{TauSelfDual}, together, of course, with a positive lower bound for the small eigenvalues of $d_A^+d_A^{+,*}$. A version of Theorem \ref{thm:Proposition_Feehan_Leness_7-6} that only required $\|F_A^+\|_{L^2(X)} \leq \eps$ would be more convenient but it is appears difficult to couple that requirement with a simultaneous demand that $C,\eps$ depend only on $\|F_A\|_{L^p(X)}$ when $p=2$.

We recall the following result of Taubes, which we have translated from the self-dual setting considered by Taubes to the anti-self-dual setting. Note that $d_A^+d_A^{+,*}$ has discrete spectrum and that its eigenvalues are non-negative and real. We omit Taubes' hypothesis that $G$ is semi-simple since that plays no role in his proof of Theorem \ref{thm:Taubes_1982_3-2} but rather was included for the sake of consistency with results of \cite{ADHM, AHS}.

\begin{defn}[Least eigenvalue of $d_A^+d_A^{+,*}$]
\label{defn:Taubes_1982_3-1}
\cite[Definition 3.1]{TauSelfDual}
Let $G$ be a compact Lie group, $P$ a principal $G$-bundle a closed, four-dimensional, oriented, Riemannian, smooth manifold, and $A$ an $H^1$ connection on $P$. The least eigenvalue of $d_A^+d_A^{+,*}$ on $L^2(X; \Lambda^+\otimes\ad P)$ is
\begin{equation}
\label{eq:Least_eigenvalue_dA+dA+*}
\mu(A) := \inf_{v \in \Omega^+(X;\ad P)\less\{0\}}
\frac{\|d_A^{+,*}v\|_{L^2(X)}^2}{\|v\|_{L^2(X)}^2}.
\end{equation}
\end{defn}

We then have the

\begin{thm}[Perturbation of an approximate anti-self-dual connection to a unique solution of the anti-self-dual equation]
\label{thm:Taubes_1982_3-2}
\cite[Theorem 3.2]{TauSelfDual}
Let $X$ be a closed, four-dimensional, oriented, smooth manifold with Riemannian metric, $g$.  Then there are a constant, $\eps = \eps(g) \in (0, 1]$, and a positive constant, $c = c(g)$, with the following significance. Let $G$ be a compact Lie group, $P$ a principal $G$-bundle over $X$, and $A$ an $H^2$ connection on $P$ obeying
\begin{equation}
\label{eq:Taubes_1982_3-6}
\mu(A) > 0 \quad\hbox{and}\quad \delta(A) < \eps,
\end{equation}
where
\begin{subequations}
\label{eq:Taubes_1982_definition_3-1}
\begin{align}
\label{eq:Taubes_1982_3-5b}
\delta(A) &:= \|F_A\|_{L^2(X)} + \zeta(A)\|F_A^+\|_{L^{4/3}(X)}\left(1 + \|F_A\|_{L^4(X)}\right),
\\
\label{eq:Taubes_1982_3-5a}
\zeta(A) &:= \mu(A)^{-1/2}\left(1 + \mu(A) + \|F_A^+\|_{L^3(X)}^3\right)^{-1/2},
\end{align}
\end{subequations}
where $\mu(A)$ is as in \eqref{eq:Least_eigenvalue_dA+dA+*}. Then there is a unique solution, $v\in H_A^2(X;\Lambda^+\otimes \ad P)$, to equation \eqref{eq:Feehan_Leness_7-2} such that
\begin{equation}
\label{eq:Taubes_1982_3-7}
\|d_A^{+,*}v\|_{H_A^1(X)} \leq c\|F_A^+\|_{L^2(X)}.
\end{equation}
Moreover, if $A$ is an $H^k$ connection, for $k\geq 3$, then
$v\in H_A^{k+1}(X;\Lambda^+\otimes \ad P)$.
\end{thm}

A stronger version of Theorem \ref{thm:Taubes_1982_3-2} in its essential ingredients (with a requirement that $\|F_A^+\|_{L^2(X)}\leq \eps$ and $\|F_A^+\|_{L^p(X)} \leq M$ for some $p>2$, where $\eps\in(0,1]$ and $M\geq 1$ can be specified independently), can be recovered from Theorem \ref{thm:Proposition_Feehan_Leness_7-6} using the Sobolev embedding for $L^p(X)\hookrightarrow L^\sharp(X)$ with $p>2$ given by Lemma \ref{lem:Feehan_4-1} and a Sobolev norm interpolation inequality \cite[Equation (7.9)]{GilbargTrudinger} for $\|\cdot\|_{L^p(X)}$. Indeed, Lemma \ref{lem:Feehan_4-1} provides a positive constant, $c = c(g)$, such that
$$
\|F_A^+\|_{L^\sharp(X)} \leq c\|F_A^+\|_{L^{8/3}(X)},
$$
while the interpolation inequality \cite[Equation (7.9)]{GilbargTrudinger} with $p=2$, $q = 8/3$, $r=4$, and $\lambda \in (0,1)$ defined by $3/8 = \lambda/2 + (1-\lambda)/4$, namely $\lambda = 1/2$, gives
$$
\|F_A^+\|_{L^\sharp(X)} \leq c\|F_A^+\|_{L^2(X)}^{1/2}\|F_A^+\|_{L^4(X)}^{1/2}.
$$
Thus, if $\|F_A^+\|_{L^2(X)} \leq \eps$ and $\|F_A^+\|_{L^4(X)} \leq M$, we obtain
$$
\|F_A^+\|_{L^\sharp(X)} \leq cM\eps^{1/2}.
$$
Obviously, the preceding bound is only useful if we can assume that $\eps \in (0,1]$ can be made arbitrarily small while $M$ is fixed. As \cite[Proposition 8.6]{TauSelfDual} illustrates, this is sometimes possible, where the Taubes' family of approximately anti-self-dual connections, $A_\lambda$ for $\lambda \in (0, \frac{1}{2}\Inj(X,g)]$, obeys the bounds
\begin{subequations}
\begin{align}
\label{eq:Taubes_1982_8-19a}
\|F_{A_\lambda}^+\|_{L^p(X)} &\leq c\lambda^{2/p},
\\
\label{eq:Taubes_1982_8-19b}
\|F_{A_\lambda}\|_{L^p(X)} &\leq c\lambda^{(4/p)-2}, \quad \forall\, p \in [1,\infty),
\end{align}
\end{subequations}
for a positive constant $c$ depending at most on the Riemannian metric, $g$.

It is also possible to perturb an approximately anti-self-dual connection, $A$, to one that is exactly anti-self-dual in the setting of Theorem \ref{thm:Proposition_Feehan_Leness_7-6}, but with a hypothesis that $\|F_A^+\|_{L^2(X)}\leq \eps$ rather than $\|F_A^+\|_{L^{2,\sharp}(X)}\leq \eps$. However, in this case, the perturbation need not be small or obey the \apriori estimates \eqref{eq:Proposition_Feehan_Leness_7-6} or \eqref{eq:Taubes_1982_3-7}. To construct this perturbation, Taubes' considers the Cauchy problem for the \emph{anti-self-dual curvature flow} \cite[Equations (5.3) and (5.4)]{TauStable} (which we have translated from Taubes' self-dual flow) for a path of connections, $A(t)$, on $P$ for $t\geq 0$,
\begin{equation}
\label{eq:Taubes_1989_5-3_and_4}
\frac{\partial A}{\partial t} = d_A^{+,*}(d_A^+d_A^{+,*} + 1)^{-1}F_A^+, \quad A(0) = A_0.
\end{equation}
In \cite[Lemma 5.4]{TauStable}, Taubes establishes the existence and uniqueness of a global solution, $A$ on $[0,\infty)\times P$, to \eqref{eq:Taubes_1989_5-3_and_4} provided $\|F_{A_0}^+\|_{L^2(X)}\leq \eps$ and $\eps = \eps(g,G) \in (0,1]$ is sufficiently small. From
\cite[Equations (5.5), (5.6), and (5.7)]{TauStable}, one knows that both $\|F_A^+(t)\|_{L^2(X)}$ and $\|F_A^+(t)\|_{L^\sharp(X)}$ are non-increasing functions of $t \in [0,\infty)$ when $A(t)$ evolves according to \eqref{eq:Taubes_1989_5-3_and_4}. It is difficult, without further hypotheses, to prove that the solution, $A(t)$, converges as $t\to\infty$ to a limit, $A_\infty$, on $P$ (which would then necessarily obey $F_{A_\infty}^+=0$).

In an earlier article, where $(X,g)$ is the four-dimensional sphere, $S^4$, with its standard round metric of radius one, Taubes considers the following, essentially equivalent, version \cite[Equations (3.2) and (3.3)]{TauPath} of \eqref{eq:Taubes_1989_5-3_and_4}, namely
\begin{equation}
\label{eq:Taubes_1984b_3-2_and_3}
\frac{\partial A}{\partial t} = d_A^{+,*}(d_A^+d_A^{+,*})^{-1}F_A^+, \quad A(0) = A_0.
\end{equation}
In writing $(d_A^+d_A^{+,*})^{-1}$ in \eqref{eq:Taubes_1984b_3-2_and_3}, Taubes' exploits the well-known Bochner-Weitzenb\"ock formula \cite[Appendix C]{FU} for $d_A^+d_A^{+,*}$, which we recall here. For a Riemannian metric $g$ on a four-dimensional, oriented manifold, $X$, let $R_g(x)$ denote its scalar curvature at a point $x \in X$ and $\sW_g^\pm(x) \in \End(\Lambda_x^\pm)$, denote its self-dual and anti-self-dual Weyl curvature tensors at $x$, where $\Lambda_x^2 = \Lambda_x^+\oplus \Lambda_x^-$. Set
$$
w_g^\pm(x) := \text{Largest eigenvalue of } \sW_g^\pm(x), \quad\forall\, x \in X.
$$
One then has the following Bochner-Weitzenb\"ock formula \cite[Equation (6.26) and Appendix C, p. 174]{FU}, \cite[Equation (5.2)]{GroisserParkerSphere},
\begin{equation}
\label{eq:Freed_Uhlenbeck_6-26}
2d_A^+d_A^{+,*}v = \nabla_A^*\nabla_Av + \left(\frac{1}{3}R - 2w_+\right)v + \{F_A^+, v\},
\quad\forall\, v \in \Omega^+(X; \ad P).
\end{equation}
We shall call a Riemannian metric, $g$, \emph{positive} if there is a metric $\tilde g$, conformally equivalent to $g$ (so $\tilde g  = e^{2f}g$ for some $f \in C^\infty(X)$), such that
\begin{equation}
\label{eq:Freed_Uhlenbeck_page_174_positive_metric}
\frac{1}{3}R_{\tilde g} - 2w_{\tilde g}^+ > 0 \quad\hbox{on } X,
\end{equation}
that is, the operator $R_{\tilde g}/3 - 2\sW_{\tilde g}^+ \in \End(\Lambda^+)$ is pointwise positive definite. Of course, the simplest example is the standard round metric of radius one on $S^4$, where $R=1$ and $w^+=0$.

When $X$ has a positive Riemannian metric, $g$, Taubes establishes existence of a global solution, $A(t)$ on $[0,\infty)\times P$, together with convergence of $A(t)$ as $t\to\infty$ to a limit, $A_\infty$, on $P$ in his \cite[proof of Proposition 3.1, p. 352]{TauPath} and \cite[Lemma 3.3]{TauPath}, respectively, when $\|F_{A_0}^+\|_{L^2(X)}\leq \eps$. While his \cite[Proposition 3.1]{TauPath} is only stated for the case of $S^4$ with its standard round metric of radius one, all of Taubes' analysis in \cite[Section 3]{TauPath} is applied in this more general context.

A related version of Taubes' anti-self-dual flow was investigated by Hong and Zheng in \cite{Hong_Zheng_2008, Zheng_2007}, who may have been unaware of Taubes' earlier work on this topic. They consider a path $v(t) \in \Omega^+(X;\ad P)$ for $t\geq 0$ defined by
\begin{equation}
\label{eq:Hong_Zheng_1-2}
\frac{\partial v}{\partial t} = F^+(A_0+d_{A_0}^{+,*}v), \quad v(0) = v_0.
\end{equation}
However, while they obtain some interesting partial results for the flow \eqref{eq:Hong_Zheng_1-2}, they do not establish global existence or convergence.

\subsection{\Apriori estimates for the operators $d_A^{+,*}$ and $d_A^+d_A^{+,*}$}
\label{subsec:Apriori_estimates_dA+*_and_dA+d_A+*}
We begin by recalling the following useful \apriori estimates from \cite[Lemma 6.6]{FLKM1}, based in turn on estimates due to Taubes in \cite[Lemma 5.2]{TauSelfDual} and in
\cite[Appendix A]{TauIndef}. For the sake of consistency, we shall assume that the section, $v$, of $\Lambda^+\otimes\ad P$ and connection, $A$, on $P$ are both $C^\infty$; \apriori estimates with weaker regularities follow by approximation.

\begin{lem}[An \apriori $L^2$ estimate for $d_A^{+,*}$ and $L^{4/3}$ estimate for $d_A^+d_A^{+,*}$]
\label{lem:Feehan_Leness_6-6}
\cite[Lemma 6.6]{FLKM1}
Let $X$ be a closed, four-dimensional, oriented, smooth manifold with Riemannian metric, $g$. Then there are positive constants, $c = c(g)$ and $\eps = \eps(g) \in (0,1]$, with the following significance. If $G$ is a compact Lie group, $A$ is a connection of class $C^\infty$ on a principal $G$-bundle $P$ over $X$ with
\begin{equation}
\label{eq:L2norm_FA+_leq_small}
\|F_A^+\|_{L^2(X)} \leq \eps,
\end{equation}
and $v \in \Omega^+(X;\ad P)$, then\footnote{We correct a typographical error in the statement of inequality (2) in \cite[Lemma 6.6]{FLKM1}, where the term $\|v\|_{L^2(X)}$ was omitted on the right-hand side.}
\begin{align}
\label{eq:Feehan_Leness_6-6-1_H1Av_L2normdA+*v}
\|v\|_{H_A^1(X)}
&\leq
c(\|d_A^{+,*}v\|_{L^2(X)}+\|v\|_{L^2(X)}),
\\
\label{eq:Feehan_Leness_6-6-1_L4v_L2dA+*v}
\|v\|_{L^4(X)}
&\leq
c(\|d_A^{+,*}v\|_{L^2(X)}+\|v\|_{L^2(X)}),
\\
\label{eq:Feehan_Leness_6-6-2_L2dA+*v_L4over3dA+dA+*v}
\|d_A^{+,*}v\|_{L^2(X)}
&\leq
c\|d_A^+d_A^{+,*}v\|_{L^{4/3}(X)}+\|v\|_{L^2(X)}),
\\
\label{eq:Feehan_Leness_6-6-2_H1Av_L4over3dA+dA+*v}
\|v\|_{H_A^1(X)}
&\leq
c(\|d_A^+d_A^{+,*}v\|_{L^{4/3}(X)}+\|v\|_{L^2(X)}).
\end{align}
\end{lem}

\begin{proof}
The \apriori estimates \eqref{eq:Feehan_Leness_6-6-1_H1Av_L2normdA+*v} and \eqref{eq:Feehan_Leness_6-6-2_L2dA+*v_L4over3dA+dA+*v} are given by \cite[Lemma 6.6]{FLKM1} and \eqref{eq:Feehan_Leness_6-6-2_H1Av_L4over3dA+dA+*v} is a trivial consequence of those. The \apriori estimate \eqref{eq:Feehan_Leness_6-6-1_L4v_L2dA+*v} is obtained by combining \eqref{eq:Feehan_Leness_6-6-1_H1Av_L2normdA+*v} with the Kato Inequality \eqref{eq:FU_6-20_first-order_Kato_inequality} and the Sobolev embedding $H^1(X) \hookrightarrow L^4(X)$.
\end{proof}

We next establish the following $L^p$ analogue of the \apriori $L^\infty$ estimate \cite[Lemma 5.3, Item (1)]{FeehanSlice}.

\begin{lem}[An \apriori $L^p$ estimate for the connection Laplace operator]
\label{lem:Feehan_5-3-1_Lp}
Let $X$ be a closed, smooth manifold of dimension $d\geq 4$ and Riemannian metric, $g$, and $q \in (d,\infty)$. Then there is a positive constant, $c = c(g,q)$, with the following significance. Let $r \in (d/3, d/2)$ be defined by $1/r = 2/d + 1/q$. Let $A$ be a Riemannian connection of class $C^\infty$ on a Riemannian vector bundle $E$ over $X$ with covariant derivative $\nabla_A$ and curvature $F_A$. If $v \in C^\infty(X;E)$, then
\begin{equation}
\label{eq:Feehan_5-3-1_dimension_Lp}
\|v\|_{L^q(X)}
\leq
c\left(\|\cov_A^*\cov_Av\|_{L^r(X)} + \|v\|_{L^r(X)}\right).
\end{equation}
\end{lem}

\begin{proof}
We adapt the proof of the estimate \cite[Lemma 5.3, Item (1)]{FeehanSlice}. For any $v\in C^\infty(X;E)$, we have the pointwise identity \eqref{eq:Freed_Uhlenbeck_6-18}, namely
$$
|\cov_Av|^2 + \frac{1}{2} d^*d|v|^2 = \langle\cov_A^*\cov_A v,v\rangle \quad\hbox{on } X,
$$
and thus,
$$
|\cov_Av|^2 + \frac{1}{2}(1+ d^*d)|v|^2 = \langle\cov_A^*\cov_A v,v\rangle + \frac{1}{2}|v|^2
\quad\hbox{on }X.
$$
As in \cite[Section 5.1]{FeehanSlice}, we let $\cG \in C^\infty(X\times X\less\Delta;\mathbb{R})$ denote the Green kernel for the augmented Laplace operator, $d^*d+1$, on $C^\infty(X;\mathbb{R})$, where $\Delta$ denotes the diagonal of $X\times X$. Using the preceding identity and the fact that
$$
\int_X \cG(x,\cdot)(d^*d+1)|v|^2\,dV = |v|^2(x), \quad\forall\, x \in X,
$$
we obtain
\begin{multline*}
\int_X \cG(x,\cdot)|\cov_Av|^2\,d\vol + \frac{1}{2} |v|^2(x)
\\
\leq
\int_X G(x,\cdot)|\langle \cov_A^*\cov_Av,v\rangle|\,d\vol
+ \frac{1}{2}\int_X G(x,\cdot)|v|^2\,d\vol, \quad\forall\, x \in X.
\end{multline*}
Writing the Green operator, $\cG := (d^*d + 1)^{-1}$, as
$$
(\cG v)(x) := \int_X \cG(x,\cdot)v\,d\vol, \quad\forall\, x \in X,
$$
we observe that $\cG$ extends to define a bounded operator,
$$
\cG:L^s(X) \to L^t(X),
$$
when $s \in (1, d/2)$ and $t\in (d/2,\infty)$ satisfy $1/s = 2/d + 1/t$ since \afortiori $\cG$ extends to a bounded operator,
$$
\cG:L^s(X) \to W^{2,s}(X),
$$
and $W^{2,s}(X) \hookrightarrow L^t(X)$ is a continuous embedding by \cite[Theorem 4.12]{AdamsFournier} when $t$ is defined as above. In particular, there is a positive constant, $c=c(g,s)$, such that
$$
\|\cG f\|_{L^t(X)} \leq c\|f\|_{L^s(X)}, \quad \forall\, f \in L^s(X).
$$
Therefore, expressing the preceding inequality for $v$ more compactly as
$$
\cG(|\cov_Av|^2) + \frac{1}{2} |v|^2
\leq
\cG|\langle \cov_A^*\cov_Av,v\rangle| + \frac{1}{2}\cG(|v|^2) \quad\hbox{on } X,
$$
and dropping the first term on the left, we find that, for $c = c(g,s)$,
\begin{align*}
\||v|^2\|_{L^t(X)}
&\leq
2\|\cG|\langle \cov_A^*\cov_Av,v\rangle|\|_{L^t(X)}
+ \|\cG(|v|^2)\|_{L^t(X)}
\\
&\leq c\|\langle \cov_A^*\cov_Av,v\rangle\|_{L^s(X)} + c\||v|^2\|_{L^s(X)}.
\end{align*}
Using
$$
\frac{1}{s} = \frac{2}{d} + \frac{1}{t} = \frac{2}{d} + \frac{1}{2t} + \frac{1}{2t}
= \frac{d+4t}{2dt} + \frac{1}{2t},
$$
we see that
\begin{align*}
\|\langle \cov_A^*\cov_Av,v\rangle\|_{L^s(X)}
&\leq
\|\cov_A^*\cov_Av\|_{L^{2dt/(d+4t)}(X)} \|v\|_{L^{2t}(X)},
\\
\||v|^2\|_{L^s(X)} &\leq \|v\|_{L^{2dt/(d+4t)}(X)} \|v\|_{L^{2t}(X)}.
\end{align*}
Therefore,
$$
\|v\|_{L^{2t}(X)}^2 \leq c\|\cov_A^*\cov_Av\|_{L^{2dt/(d+4t)}(X)} \|v\|_{L^{2t}(X)}
+ c\|v\|_{L^{2dt/(d+4t)}(X)} \|v\|_{L^{2t}(X)},
$$
and thus, for $v$ not identically zero,
$$
\|v\|_{L^{2t}(X)}
\leq
c\left(\|\cov_A^*\cov_Av\|_{L^{2dt/(d+4t)}(X)} + \|v\|_{L^{2dt/(d+4t)}(X)}\right),
$$
for any $t \in (d/2,\infty)$. But $(d+4t)/(2dt)=  1/(2t) + 2/d$, so writing $q = 2t \in (d,\infty)$ and $r = 2dt/(d+4t) \in (d/3, d/2)$ (for $t \in (d/2,\infty)$) yields \eqref{eq:Feehan_5-3-1_dimension_Lp}, where $1/r = 2/d + 1/q$.
\end{proof}

We now apply Lemma \ref{lem:Feehan_5-3-1_Lp} to $E = \Lambda^+\otimes\ad P$ with Riemannian connection induced by a connection $A$ on $P$ and the Levi-Civita connection $TX$.

\begin{lem}[An \apriori $L^p$ estimate for $d_A^+d_A^{+,*}$]
\label{lem:Feehan_5-3-1_Lp_dA+dA+*}
Let $X$ be a closed, four-dimensional, oriented, smooth manifold, $X$, with Riemannian metric, $g$, and $q \in [4,\infty)$. Then there are positive constants, $c = c(g,q)\in [1,\infty)$ and $\eps = \eps(g,q) \in (0,1]$, with the following significance. Let $r \in [4/3, 2)$ be defined by $1/r = 1/2 + 1/q$. Let $G$ be a compact Lie group and $A$ a connection of class $C^\infty$ on a principal bundle $P$ over $X$ that obeys the curvature bound \eqref{eq:L2norm_FA+_leq_small}. If $v \in \Omega^+(X;\ad P)$, then
\begin{equation}
\label{eq:Feehan_5-3-1_Lp_dA+dA+*}
\|v\|_{L^q(X)}
\leq
c\left(\|d_A^+d_A^{+,*}v\|_{L^r(X)} + \|v\|_{L^r(X)}\right).
\end{equation}
\end{lem}

\begin{proof}
We first dispose of the simplest case, $q=4$ and $r=4/3$. We combine the Kato Inequality \eqref{eq:FU_6-20_first-order_Kato_inequality}, Sobolev embedding $H^1(X) \hookrightarrow L^4(X)$, and \apriori estimate \eqref{eq:Feehan_Leness_6-6-2_H1Av_L4over3dA+dA+*v} to give
$$
\|v\|_{L^4(X)}
\leq
c(\|d_A^+d_A^{+,*}v\|_{L^{4/3}(X)}+\|v\|_{L^2(X)}).
$$
Substituting the interpolation inequality, $\|v\|_{L^2(X)} \leq \|v\|_{L^{4/3}(X)}^{1/2}\|v\|_{L^4(X)}^{1/2}$, in the preceding estimate yields, for any $\zeta > 0$,
\begin{align*}
\|v\|_{L^4(X)}
&\leq
c\|d_A^+d_A^{+,*}v\|_{L^{4/3}(X)} + c\|v\|_{L^{4/3}(X)}^{1/2}\|v\|_{L^4(X)}^{1/2}
\\
&\leq c\|d_A^+d_A^{+,*}v\|_{L^{4/3}(X)} + \frac{c}{2\zeta}\|v\|_{L^{4/3}(X)} + \frac{c\zeta}{2}\|v\|_{L^4(X)}.
\end{align*}
We obtain \eqref{eq:Feehan_5-3-1_Lp_dA+dA+*} when $q=4$ and $r=4/3$ from the preceding inequality by choosing $\zeta = 1/c$.

For the remainder of the proof, we assume $q \in (4,\infty)$ and $r \in (4/3,2)$. The Bochner-Weitzenb\"ock formula \eqref{eq:Freed_Uhlenbeck_6-26}, namely,
$$
2d_A^+d_A^{+,*} = \cov_A^*\cov_A + \left(\frac{R}{3} - 2w^+\right) + \{F_A^+, \cdot\},
$$
yields, for $v \in \Omega^+(X;\ad P)$,
$$
\|\cov_A^*\cov_Av\|_{L^r(X)}
\leq
2\|d_A^+d_A^{+,*}v\|_{L^r(X)}
+ c\|v\|_{L^r(X)} + \|\{F_A^+, v\}\|_{L^r(X)},
$$
and some $c = c(g)$. Since $1/r = 1/2 + 1/q$ by hypothesis, we see that
$$
\|\{F_A^+, v\}\|_{L^r(X)} \leq c\|F_A^+\|_{L^2(X)} \|v\|_{L^q(X)},
$$
for some $c = c(g)$. Combining the preceding inequalities with the estimate \eqref{eq:Feehan_5-3-1_dimension_Lp} yields
$$
\|v\|_{L^q(X)}
\leq
c\|d_A^+d_A^{+,*}v\|_{L^r(X)} + c\|v\|_{L^r(X)} + c\|F_A^+\|_{L^2(X)} \|v\|_{L^q(X)},
$$
for some $c = c(g,q)$. Provided $c\|F_A^+\|_{L^2(X)} \leq 1/2$, rearrangement gives \eqref{eq:Feehan_5-3-1_Lp_dA+dA+*} when $r \in (4/3,2)$.
\end{proof}

Finally, we recall the following analogues of Lemmata \ref{lem:Feehan_5-3-1_Lp} and \ref{lem:Feehan_5-3-1_Lp_dA+dA+*} when $q=\infty$ and $r=2$.

\begin{lem}
\label{lem:Feehan_5-8}
\cite[Lemma 5.8]{FeehanSlice}
Let $X$ be a closed, four-dimensional, smooth manifold with Riemannian metric, $g$. Then there is a positive constant, $c=c(g)$, with the following significance. Let $A$ be a Riemannian connection of class $C^\infty$ on a Riemannian vector bundle $E$ over $X$ with covariant derivative $\nabla_A$ and curvature $F_A$. If $v \in C^\infty(X;E)$, then
$$
\|v\|_{H_A^2(X)} + \|v\|_{C(X)}
\leq c\left(1+\|F_A\|_{L^2(X)}\right)
\left(\|\cov_A^*\cov_Av\|_{L^{\sharp,2}(X)} + \|v\|_{L^2(X)}\right).
$$
\end{lem}

We apply Lemma \ref{lem:Feehan_5-8} to $E = \Lambda^+\otimes\ad P$ with Riemannian connection induced by a connection $A$ on $P$ and the Levi-Civita connection $TX$.

\begin{lem}
\label{lem:Feehan_5-9}
\cite[Lemma 5.9]{FeehanSlice}
Let $X$ be a closed, four-dimensional, oriented, smooth manifold with Riemannian metric, $g$. Then there is a positive constant, $c=c(g)$, with the following significance. Let $G$ be a compact Lie group and $A$ a connection of class $C^\infty$ on a principal bundle $P$ over $X$ that obeys the curvature bound,
\begin{equation}
\label{eq:L2norm_FA+sharp_leq_small}
\|F_A^+\|_{L^{\sharp,2}(X)} <\eps.
\end{equation}
If $v \in \Omega^+(X;\ad P)$, then
\begin{equation}
\label{eq:Lsharp2_dA+dA+*}
\|v\|_{H_A^2(X)} + \|v\|_{C(X)}
\leq
c\left(1+\|F_A\|_{L^2(X)}\right)
\left(\|d_A^+d_A^{+,*}v\|_{L^{\sharp,2}(X)} + \|v\|_{L^2(X)}\right).
\end{equation}
\end{lem}

If one seeks an $L^p$ analogue of the \apriori elliptic estimate in Lemma \ref{lem:Feehan_5-9} where the norm on the right-hand-side of \eqref{eq:Lsharp2_dA+dA+*} is replaced by $W^{2,p}$ for $p>2$, replacement of the condition \eqref{eq:L2norm_FA+sharp_leq_small} by one on $\|F_A^+\|_{L^p(X)}$ is to be expected but unfortunately it no longer appears to be possible for the estimate constant to depend at most on the connection, $A$, through its energy, $\sE(A) = \frac{1}{2}\|F_A\|_{L^2(X)}^2$.

Lemma \ref{lem:Apriori_Lp_estimate_solution_ASD_equation} in the sequel appears to provide an $L^p$ extension of Lemma \ref{lem:Feehan_5-9} which is optimal with regard to dependence on $A$ when $p > 2$. Before stating that result, we shall need an elementary lemma allowing us to compare Sobolev norms defined by two nearby connections, $A_0$ and $A_1$ on $P$.

\begin{lem}[Comparison estimates for $W^{k,p}$ Sobolev norms defined by connections that are $W^{k-1,4}$-close]
\label{lem:Equivalence_W1p_or_W2p_norms_connections_A0_and_A1_L4_or_W14_close}
Let $G$ be a compact Lie group and $P$ a principal $G$-bundle over a closed, connected, four-dimensional, smooth manifold, $X$, with Riemannian metric, $g$, and $p \in [2,4)$. Then there is a positive constant, $\delta = \delta(g,p) \in (0,1]$, with the following significance.
\begin{enumerate}
\item If $A_0, A_1$ are $L^4$ connections on $P$ such that
$$
\|A_0 - A_1\|_{L^4(X)} \leq \delta,
$$
then
\begin{equation}
\label{eq:Equivalence_W1p_norms_connections_A0_and_A1_L4_close}
\frac{1}{4}\|b\|_{W_{A_1}^{1,p}(X)}
\leq
\|b\|_{W_{A_0}^{1,p}(X)}
\leq
4\|b\|_{W_{A_1}^{1,p}(X)}, \quad\forall\, b \in W_{A_0}^{1,p}(X;\Lambda^1\otimes\ad P).
\end{equation}
\item If $A_0, A_1$ are $W^{1,4}$ connections on $P$ such that
$$
\|A_0 - A_1\|_{W_{A_0}^{1,4}(X)} \leq \delta
\quad\hbox{or}\quad
\|A_0 - A_1\|_{W_{A_1}^{1,4}(X)} \leq \delta,
$$
then
\begin{equation}
\label{eq:Equivalence_W2p_norms_connections_A0_and_A1_W14_close}
\frac{1}{8}\|b\|_{W_{A_1}^{2,p}(X)}
\leq
\|b\|_{W_{A_0}^{2,p}(X)}
\leq
8\|b\|_{W_{A_1}^{2,p}(X)}, \quad\forall\, b \in W_{A_0}^{2,p}(X;\Lambda^1\otimes\ad P).
\end{equation}
\end{enumerate}
\end{lem}

\begin{proof}
Define $q\in [4,\infty)$ by $1/p=1/4+1/q$. For any $b \in W_{A_0}^{1,p}(X; \Lambda^1\otimes\ad P)$, we have
$$
\nabla_{A_1}b = \nabla_{A_0}b + [A_1-A_0, b].
$$
The H\"older Inequality, Sobolev embedding $W^{1,p}(X) \hookrightarrow L^q(X)$ \cite[Theorem 4.12]{AdamsFournier}, with embedding constant denoted by $\kappa_p = \kappa_p(g,p)$, and Kato Inequality \eqref{eq:FU_6-20_first-order_Kato_inequality} give
\begin{align*}
\|\nabla_{A_1}b\|_{L^p(X)}
&\leq
\|\nabla_{A_0}b\|_{L^p(X)} + 2\|A_0 - A_1\|_{L^4(X)}\|b\|_{L^q(X)}
\\
&\leq
\|\nabla_{A_0}b\|_{L^p(X)}
+ 2\kappa_p\|A_0 - A_1\|_{L^4(X)} \|b\|_{W_{A_1}^{1,p}(X)}
\\
&\leq
\|\nabla_{A_0}b\|_{L^p(X)}
+ 2\kappa_p\|A_0 - A_1\|_{L^4(X)}
\left(\|b\|_{L^p(X)} + \|\nabla_{A_1}b\|_{L^p(X)}\right).
\end{align*}
We now choose $\delta = \delta(g,p) \in (0,1]$ so that $2\kappa_p\delta \leq 1/2$ and hence
$$
\|\nabla_{A_1}b\|_{L^p(X)}
\leq \|\nabla_{A_0}b\|_{L^p(X)}
+ \frac{1}{2}\left(\|b\|_{L^p(X)} + \|\nabla_{A_1}b\|_{L^p(X)}\right),
$$
so that, after rearranging,
\begin{equation}
\label{eq:Equivalence_W1p_norms_connections_A0_and_A1_L4_close_prelim}
\|\nabla_{A_1}b\|_{L^p(X)}
\leq
2\|\nabla_{A_0}b\|_{L^p(X)} + \|b\|_{L^p(X)}.
\end{equation}
In particular, we find that
\begin{align*}
\|b\|_{W_{A_1}^{1,p}(X)}
&=
\left( \|\nabla_{A_1}b\|_{L^p(X)}^p + \|b\|_{L^p(X)}^p \right)^{1/p}
\\
&\leq \|\nabla_{A_1}b\|_{L^p(X)} + \|b\|_{L^p(X)}
\\
&\leq 2\left(\|\nabla_{A_0}b\|_{L^p(X)} + \|b\|_{L^p(X)}\right)
\\
&\leq 4\|b\|_{W_{A_0}^{1,p}(X)},
\end{align*}
which gives \eqref{eq:Equivalence_W1p_norms_connections_A0_and_A1_L4_close}. Similarly,
\begin{align*}
\|b\|_{W_{A_1}^{2,p}(X)}
&=
\left( \|\nabla_{A_1}^2b\|_{L^p(X)}^p + \|\nabla_{A_1}b\|_{L^p(X)}^p + \|b\|_{L^p(X)}^p \right)^{1/p}
\\
&\leq \|\nabla_{A_1}^2b\|_{L^p(X)} + \|\nabla_{A_1}b\|_{L^p(X)} + \|b\|_{L^p(X)}
\\
&\leq \|\nabla_{A_1}^2b\|_{L^p(X)} + 2\left(\|\nabla_{A_0}b\|_{L^p(X)} + \|b\|_{L^p(X)}\right)
\quad\hbox{(by \eqref{eq:Equivalence_W1p_norms_connections_A0_and_A1_L4_close_prelim}).}
\end{align*}
Writing $A_1=A_0+a$ for brevity in the calculations below, we have
\begin{align*}
\nabla_{A_1}^2b
&=
\nabla_{A_0+a}(\nabla_{A_0}b + [a,b])
\\
&= \nabla_{A_0}^2b + \nabla_{A_0}[a,b] + [a,[a,b]]
\\
&= \nabla_{A_0}^2b + [\nabla_{A_0}a,b] - [a,\nabla_{A_0}b] + [a,[a,b]].
\end{align*}
Therefore, the H\"older Inequality, Sobolev embedding $W^{1,p}(X) \hookrightarrow L^q(X)$ \cite[Theorem 4.12]{AdamsFournier}, and Kato Inequality \eqref{eq:FU_6-20_first-order_Kato_inequality} yield
\begin{align*}
\|\nabla_{A_1}^2b\|_{L^p(X)}
&\leq
\|\nabla_{A_0}^2b\|_{L^p(X)} + 2\|\nabla_{A_0}a\|_{L^4(X)} \|b\|_{L^q(X)}
+ 2\|a\|_{L^4(X)} \|\nabla_{A_0}b\|_{L^q(X)}
\\
&\quad + 2\|a\|_{L^4(X)} \|[a,b]\|_{L^q(X)}
\\
&\leq
\|\nabla_{A_0}^2b\|_{L^p(X)} + 2\kappa_p\|\nabla_{A_0}a\|_{L^4(X)} \|b\|_{W_{A_0}^{1,p}(X)}
+ 2\kappa_p\|a\|_{L^4(X)} \|\nabla_{A_0}b\|_{W_{A_0}^{1,p}(X)}
\\
&\quad + 2\kappa_p\|a\|_{L^4(X)} \|[a,b]\|_{W_{A_0}^{1,p}(X)}.
\end{align*}
Since $\nabla_{A_0}[a,b] = [\nabla_{A_0}a,b] - [a,\nabla_{A_0}b]$, we similarly find that
\begin{align*}
\|[a,b]\|_{W_{A_0}^{1,p}(X)}
&\leq
\|\nabla_{A_0}[a,b]\|_{L^p(X)} + \|[a,b]\|_{L^p(X)}
\\
&\leq \|[\nabla_{A_0}a,b]\|_{L^p(X)} + \|[a,\nabla_{A_0}b]\|_{L^p(X)}
+ 2\|a\|_{L^4(X)} \|b\|_{L^q(X)}
\\
&\leq 2\|\nabla_{A_0}a\|_{L^4(X)} \|b\|_{L^q(X)}
+ 2\|a\|_{L^4(X)} \|\nabla_{A_0}b\|_{L^q(X)}
+ 2\|a\|_{L^4(X)} \|b\|_{L^q(X)}
\\
&\leq 2\kappa_p\|a\|_{W_{A_0}^{1,4}(X)}
\left( \|\nabla_{A_0}b\|_{W_{A_0}^{1,p}(X)} + 2\|b\|_{W_{A_0}^{1,p}(X)}\right)
\\
&\leq 6\kappa_p\|a\|_{W_{A_0}^{1,4}(X)}\|b\|_{W_{A_0}^{2,p}(X)}.
\end{align*}
Combining the preceding inequalities yields,
\begin{align*}
\|b\|_{W_{A_1}^{2,p}(X)}
&\leq
\|\nabla_{A_1}^2b\|_{L^p(X)} + 2\left(\|\nabla_{A_0}b\|_{L^p(X)} + \|b\|_{L^p(X)}\right)
\\
&\leq \|\nabla_{A_0}^2b\|_{L^p(X)}
+ 2\kappa_p\|\nabla_{A_0}a\|_{L^4(X)} \|b\|_{W_{A_0}^{1,p}(X)}
+ 2\kappa_p\|a\|_{L^4(X)} \|b\|_{W_{A_0}^{2,p}(X)}
\\
&\quad + 12\kappa_p^2\|a\|_{L^4(X)} \|a\|_{W_{A_0}^{1,4}(X)} \|b\|_{W_{A_0}^{2,p}(X)}
+ 4\|b\|_{W_{A_0}^{1,p}(X)}.
\end{align*}
Hence, for small enough $\delta = \delta(g,p) \in (0,1]$ and $\|a\|_{W_{A_0}^{1,4}(X)} \leq \delta$, we obtain
$$
\|b\|_{W_{A_1}^{2,p}(X)} \leq 8\|b\|_{W_{A_0}^{2,p}(X)},
$$
which gives \eqref{eq:Equivalence_W2p_norms_connections_A0_and_A1_W14_close}. By symmetry, we also have \eqref{eq:Equivalence_W2p_norms_connections_A0_and_A1_W14_close} if $\|a\|_{W_{A_1}^{1,4}(X)} \leq \delta$. This completes the proof of Lemma \ref{lem:Equivalence_W1p_or_W2p_norms_connections_A0_and_A1_L4_or_W14_close}.
\end{proof}

Of course, the norms on $W_{A_0}^{k,p}(X;\Lambda^1\otimes\ad P)$ and $W_{A_1}^{k,p}(X;\Lambda^1\otimes\ad P)$ remain equivalent regardless of the size of $\|A_0 - A_1\|_{W_{A_0}^{k-1,p}(X)}$ or $\|A_0 - A_1\|_{W_{A_1}^{k-1,p}(X)}$, but then the equivalences will no longer take the form given in Lemma \ref{lem:Equivalence_W1p_or_W2p_norms_connections_A0_and_A1_L4_or_W14_close} with universal numerical constants.

\begin{lem}[\Apriori $L^2$ and $L^p$ estimates for solutions to the anti-self-dual equation]
\label{lem:Apriori_Lp_estimate_solution_ASD_equation}
Let $G$ be a compact Lie group, $P$ a principal $G$-bundle over a closed, connected, four-dimensional, smooth manifold, $X$, with Riemannian metric, $g$, and $\mu_0 \in (0,\infty)$ a constant.
\begin{enumerate}
\item If $E_0 \in (0,\infty)$ is a constant, then there are constants, $\delta = \delta(g) \in (0,1]$ and $C_0 = C_0(E_0,g,\mu_0) \in [1,\infty)$ with the following significance. If $A_0$ is a $W^{1,4}$ connection on $P$ such that
\begin{subequations}
\label{eq:FA0+_L2sharp_leq_small_and_FA0_L2_leq_bound_and_muA0_geq_mu0}
\begin{align}
\label{eq:FA0+_L2sharp_leq_small}
\|F_{A_0}^+\|_{L^{\sharp,2}(X)} &\leq \eps,
\\
\label{eq:FA0_L2_leq_bound}
\|F_{A_0}\|_{L^2(X)} &\leq E_0,
\\
\label{eq:muA0_geq_mu0}
\mu(A_0) &\geq \mu_0,
\end{align}
\end{subequations}
where $\mu(A_0)$ is as in \eqref{eq:Least_eigenvalue_dA+dA+*}, and
\begin{equation}
\label{eq:dA0+*v_L2sharp_and_L4_norm_leq_small}
\|d_{A_0}^{+,*}v\|_{L^{2\sharp,4}(X)} \leq \delta,
\end{equation}
where $L^{2\sharp,4}(X) = L^{2\sharp}(X) \cap L^4(X)$ is as in \eqref{eq:WholeFamilyOfSharpNorms}, then
\begin{subequations}
\begin{align}
\label{eq:Apriori_estimate_dA0+*v_L2sharp_and_4}
\|d_{A_0}^{+,*}v\|_{L^{2\sharp,4}(X)} &\leq C_0\|F_{A_0}^+\|_{L^{\sharp,2}(X)},
\\
\label{eq:Apriori_estimate_dA0+*v_HA1}
\|d_{A_0}^{+,*}v\|_{H_{A_1}^1(X)} &\leq C_0\|F_{A_0}^+\|_{L^{\sharp,2}(X)},
\\
\label{eq:Apriori_estimate_v_HA2}
\|v\|_{H_{A_0}^2(X)} + \|v\|_{C(X)} &\leq C_0\|F_{A_0}^+\|_{L^{\sharp,2}(X)}.
\end{align}
\end{subequations}
\item If $A_1$ is a $W^{1,4}$ connection on $P$ and $p \in [2,4)$ and $q \in [4,\infty)$ is defined by $1/p=1/4+1/q$, then there are constants, $\delta = \delta(g,p) \in (0,1]$ and $C = C([A_1],g,p,\mu_0) \in [1,\infty)$ with the following significance. If $A_0$ is a $W^{1,4}$ connection on $P$ that obeys \eqref{eq:muA0_geq_mu0} and
\begin{equation}
\label{eq:dA0+*v_L4_norm_leq_small}
\|d_{A_0}^{+,*}v\|_{L^4(X)} \leq \delta,
\end{equation}
then
\begin{subequations}
\begin{align}
\label{eq:Apriori_estimate_dA0+*v_Lq}
\|d_{A_0}^{+,*}v\|_{L^q(X)} &\leq C\|F_{A_0}^+\|_{L^p(X)},
\\
\label{eq:Apriori_estimate_dA0+*v_WA1_p}
\|d_{A_0}^{+,*}v\|_{W_{A_1}^{1,p}(X)} &\leq C\|F_{A_0}^+\|_{L^p(X)},
\\
\label{eq:Apriori_estimate_v_WA2_p}
\|v\|_{W_{A_1}^{2,p}(X)} &\leq C\|F_{A_0}^+\|_{L^p(X)},
\\
\label{eq:Apriori_estimate_FAasd_Lp}
\|F(A_0 + d_{A_0}^{+,*}v)\|_{L^p(X)} &\leq \|F_{A_0}\|_{L^p(X)} + C\|F_{A_0}^+\|_{L^p(X)}.
\end{align}
\end{subequations}
\end{enumerate}
\end{lem}

\begin{proof}
By expanding the $g$-anti-self-dual equation, $F^+(A_0 + d_{A_0}^{+,*}v) = 0$, in the usual way, one obtains,
\begin{equation}
\label{eq:Second_order_ASD_equation}
d_{A_0}^+d_{A_0}^{+,*}v + (d_{A_0}^{+,*}v\wedge d_{A_0}^{+,*}v)^+ = -F_{A_0}^+,
\end{equation}
In the sequel, we let $c=c(g)$ denote a positive constant depending at most on the Riemannian metric, $g$. We first consider the case $p=2$:
\begin{align*}
&\|v\|_{H_{A_0}^2(X)} + \|v\|_{C(X)}
\\
&\quad\leq c\left(1+\|F_{A_0}\|_{L^2(X)}\right)
\left(\|d_{A_0}^+d_{A_0}^{+,*}v\|_{L^{\sharp,2}(X)} + \|v\|_{L^2(X)}\right)
\quad\hbox{(by Lemma \ref{lem:Feehan_5-9})}
\\
&\quad\leq c\left(1+\|F_{A_0}\|_{L^2(X)}\right)\left(1 + \mu(A_0)^{-1}\right)
\|d_{A_0}^+d_{A_0}^{+,*}v\|_{L^{\sharp,2}(X)}
\quad\hbox{(by \eqref{eq:Least_eigenvalue_dA+dA+*})}
\\
&\quad\leq c\left(1+\|F_{A_0}\|_{L^2(X)}\right)\left(1 + \mu(A_0)^{-1}\right)
\left( \|F_{A_0}^+\|_{L^{\sharp,2}(X)}
+ \|d_{A_0}^{+,*}v\wedge d_{A_0}^{+,*}v\|_{L^{\sharp,2}(X)} \right)
\\
&\qquad\hbox{(by \eqref{eq:Second_order_ASD_equation})}
\\
&\quad\leq c\left(1+\|F_{A_0}\|_{L^2(X)}\right)\left(1 + \mu(A_0)^{-1}\right)
\left( \|F_{A_0}^+\|_{L^{\sharp,2}(X)}
+ \|d_{A_0}^{+,*}v\|_{L^{2\sharp,4}(X)}^2 \right)
\\
&\qquad\hbox{(by Lemma \ref{lem:Feehan_4-3})}
\\
&\quad\leq c\left(1+\|F_{A_0}\|_{L^2(X)}\right)\left(1 + \mu(A_0)^{-1}\right)
\\
&\qquad \times \left( \|F_{A_0}^+\|_{L^{\sharp,2}(X)}
+ (\kappa_\sharp+\kappa_2)\|d_{A_0}^{+,*}v\|_{L^{2\sharp,4}(X)}\|d_{A_0}^{+,*}v\|_{H_{A_0}^1(X)} \right),
\end{align*}
where we apply the Sobolev embedding $H^1(X) \hookrightarrow L^{2\sharp}(X)$ from Lemma \ref{lem:Feehan_4-1} with constant $\kappa_\sharp=\kappa_\sharp(2)$, the usual embedding $H^1(X) \hookrightarrow L^4(X)$ \cite[Theorem 4.12]{AdamsFournier} with constant $\kappa_2=\kappa_2(g)$, and the Kato Inequality \eqref{eq:FU_6-20_first-order_Kato_inequality} to obtain the last inequality. Hence, for $\|d_{A_0}^{+,*}v\|_{L^{2\sharp,4}(X)} \leq \delta$ as in \eqref{eq:dA0+*v_L2sharp_and_L4_norm_leq_small} and small enough $\delta = \delta(E_0,g,\mu(A_0)) \in (0,1]$, rearrangement gives
$$
\|v\|_{H_{A_0}^2(X)} + \|v\|_{C(X)}
\leq
C_0\|F_{A_0}^+\|_{L^{\sharp,2}(X)},
$$
with constant $C_0=C_0(E_0,g,\mu(A_0))$, which is the desired estimate \eqref{eq:Apriori_estimate_v_HA2}. Since
\begin{align*}
\|d_{A_0}^{+,*}v\|_{L^{2\sharp,4}(X)}
&\leq
(\kappa_\sharp+\kappa_2)\|d_{A_0}^{+,*}v\|_{H_{A_0}^1(X)},
\\
\|d_{A_0}^{+,*}v\|_{H_{A_0}^1(X)} &\leq \|v\|_{H_{A_0}^2(X)},
\end{align*}
we also obtain the estimates \eqref{eq:Apriori_estimate_dA0+*v_L2sharp_and_4} and \eqref{eq:Apriori_estimate_dA0+*v_HA1}.

For $p \in [2,4)$ and $q\in [4,\infty)$, equation \eqref{eq:Second_order_ASD_equation} gives
\begin{align*}
\|d_{A_0}^+d_{A_0}^{+,*}v\|_{L^p(X)}
&\leq
\|F_{A_0}^+\|_{L^p(X)} + 2\|d_{A_0}^{+,*}v\|_{L^4(X)} \|d_{A_0}^{+,*}v\|_{L^q(X)}
\\
&\leq \|F_{A_0}^+\|_{L^p(X)} + 2\delta \|d_{A_0}^{+,*}v\|_{L^q(X)}
\quad\hbox{(by \eqref{eq:dA0+*v_L4_norm_leq_small})}
\\
&\leq \|F_{A_0}^+\|_{L^p(X)} + 2\delta \|\nabla_{A_0}v\|_{L^q(X)}
\\
&\leq \|F_{A_0}^+\|_{L^p(X)}
+ 2\delta\kappa_p\|\nabla_{A_0}v\|_{W_{A_0}^{1,p}(X)}
\\
&\leq \|F_{A_0}^+\|_{L^p(X)}
+ 2\delta\kappa_p\|v\|_{W_{A_0}^{2,p}(X)},
\end{align*}
where we have applied the Sobolev embedding $W^{1,p}(X) \hookrightarrow L^q(X)$ with embedding constant $\kappa_p$ \cite[Theorem 4.12]{AdamsFournier} and Kato Inequality \eqref{eq:FU_6-20_first-order_Kato_inequality} above. Because the norms on $W_{A_0}^{2,p}(X;\Lambda^1\otimes\ad P)$ and $W_{A_1}^{2,p}(X;\Lambda^1\otimes\ad P)$ are equivalent via universal constants by Lemma \ref{lem:Equivalence_W1p_or_W2p_norms_connections_A0_and_A1_L4_or_W14_close}, the preceding inequality yields
\begin{equation}
\label{eq:Lp_estimate_dA0+dA0+*v}
\|d_{A_0}^+d_{A_0}^{+,*}v\|_{L^p(X)}
\leq
\|F_{A_0}^+\|_{L^p(X)} + 16\delta\kappa_p\|v\|_{W_{A_1}^{2,p}(X)}.
\end{equation}
By Theorem \ref{thm:Krylov_Sobolev_lectures_8-5-3}, we have an \apriori $L^p$ estimate for the elliptic operator, $d_{A_1}^+d_{A_1}^{+,*}$, namely
\begin{equation}
\label{eq:Apriori_Lp_estimate_elliptic_operator_dA1+dA1+*}
\|v\|_{W_{A_1}^{2,p}(X)}
\leq
C_1\left(\|d_{A_1}^+d_{A_1}^{+,*}v\|_{L^p(X)} + \|v\|_{L^p(X)}\right),
\quad\forall\, v \in W_{A_1}^{2,p}(X;\Lambda^1\otimes\ad P),
\end{equation}
with constant $C_1 = C_1([A_1],g,p)$. A simple modification of the calculations in the proof of Lemma \ref{lem:Equivalence_W1p_or_W2p_norms_connections_A0_and_A1_L4_or_W14_close} that lead to the pair of inequalities \eqref{eq:Equivalence_W2p_norms_connections_A0_and_A1_W14_close} also give
\begin{equation}
\label{eq:Equivalence_dA+dA+*_connections_A0_and_A1_W14_close}
\frac{1}{8}\|d_{A_1}^+d_{A_1}^{+,*}v\|_{L^p(X)}
\leq
\|d_{A_0}^+d_{A_0}^{+,*}v\|_{L^p(X)}
\leq
8\|d_{A_1}^+d_{A_1}^{+,*}v\|_{L^p(X)}, \quad\forall\, v \in W_{A_0}^{2,p}(X;\Lambda^+\otimes\ad P).
\end{equation}
Therefore, using $\|v\|_{L^{4/3}(X)} \leq c\|v\|_{L^2(X)} \leq c\mu(A_0)^{-1}\|d_{A_0}^+d_{A_0}^{+,*}v\|_{L^2(X)}$ in the calculations below,
\begin{align*}
\|v\|_{W_{A_1}^{2,p}(X)} &\leq 8C_1\left(\|d_{A_1}^+d_{A_1}^{+,*}v\|_{L^p(X)} + \|v\|_{L^p(X)}\right)
\quad\hbox{(by \eqref{eq:Apriori_Lp_estimate_elliptic_operator_dA1+dA1+*})}
\\
&\leq 8C_1\left(8\|d_{A_0}^+d_{A_0}^{+,*}v\|_{L^p(X)} + \|v\|_{L^p(X)}\right)
\quad\hbox{(by \eqref{eq:Equivalence_dA+dA+*_connections_A0_and_A1_W14_close})}.
\end{align*}
Since $p \in [2,4)$ (and in particular, $p \leq 4$), we can apply Lemma \ref{lem:Feehan_5-3-1_Lp_dA+dA+*} to deduce that
\begin{align*}
\|v\|_{L^p(X)} &\leq c_p\left(\|d_{A_0}^+d_{A_0}^{+,*}v\|_{L^{4/3}(X)}
+ \|v\|_{L^{4/3}(X)}\right)
\\
&\leq c_p\|d_{A_0}^+d_{A_0}^{+,*}v\|_{L^p(X)} + \|v\|_{L^2(X)}
\\
&\leq c_p\|d_{A_0}^+d_{A_0}^{+,*}v\|_{L^p(X)}
+ \mu(A_0)^{-1}\|d_{A_0}^+d_{A_0}^{+,*}v\|_{L^2(X)}
\\
&\leq c_p\left(1 + \mu(A_0)^{-1}\right)\|d_{A_0}^+d_{A_0}^{+,*}v\|_{L^p(X)}.
\end{align*}
Combining the preceding inequalities gives
\begin{align*}
\|v\|_{W_{A_1}^{2,p}(X)}
&\leq
64C_1\|d_{A_0}^+d_{A_0}^{+,*}v\|_{L^p(X)}
+ 8cc_pC_1\left(1+\mu(A_0)^{-1}\right)\|d_{A_0}^+d_{A_0}^{+,*}v\|_{L^2(X)}
\\
&\leq 8C_1(8+cc_p)\left(1+\mu(A_0)^{-1}\right)\|d_{A_0}^+d_{A_0}^{+,*}v\|_{L^p(X)}
\\
&\leq 8C_1(8+cc_p)\left(1+\mu(A_0)^{-1}\right)
\left(\|F_{A_0}^+\|_{L^p(X)} + 16\delta\kappa_p\|v\|_{W_{A_1}^{2,p}(X)}\right)
\quad\hbox{(by \eqref{eq:Lp_estimate_dA0+dA0+*v})}
\end{align*}
Thus, for small enough $\delta = \delta([A_1],g,p,\mu(A_0)) \in (0,1]$, rearrangement yields
$$
\|v\|_{W_{A_1}^{2,p}(X)} \leq 16C_1(8+cc_p)\left(1+\mu(A_0)^{-1}\right)\|F_{A_0}^+\|_{L^p(X)},
$$
and this is the desired estimate \eqref{eq:Apriori_estimate_v_WA2_p}. The inequalities \eqref{eq:Apriori_estimate_dA0+*v_Lq} and \eqref{eq:Apriori_estimate_dA0+*v_WA1_p} follow from
\begin{align*}
\|d_{A_0}^{+,*}v\|_{L^q(X)}
&\leq
\kappa_p\|d_{A_0}^{+,*}v\|_{W_{A_1}^{1,p}(X)},
\\
\|d_{A_0}^{+,*}v\|_{W_{A_1}^{1,p}(X)} &\leq \|v\|_{W_{A_1}^{2,p}(X)}.
\end{align*}
It remains to estimate the $L^p$ norm of the curvature of $A_\asd := A_0 + d_{A_0}^{+,*}v$. From
$$
F(A_\asd) = F_{A_0} + d_{A_0}d_{A_0}^{+,*}v + d_{A_0}^{+,*}v \wedge d_{A_0}^{+,*}v,
$$
we obtain
\begin{align*}
\|F(A_\asd)\|_{L^p(X)} &\leq \|F_{A_0}\|_{L^p(X)} + \|d_{A_0}d_{A_0}^{+,*}v\|_{L^p(X)}
+ 2\|d_{A_0}^{+,*}v\|_{L^4(X)} \|d_{A_0}^{+,*}v\|_{L^q(X)}
\\
&\leq \|F_{A_0}\|_{L^p(X)} + \|F_{A_0}^+\|_{L^p(X)} + 16\delta\kappa_p\|v\|_{W_{A_1}^{2,p}(X)}
+ 2\delta \|d_{A_0}^{+,*}v\|_{L^q(X)}
\\
&\qquad\hbox{(by \eqref{eq:Lp_estimate_dA0+dA0+*v} and \eqref{eq:dA0+*v_L4_norm_leq_small})}
\\
&\leq \|F_{A_0}\|_{L^p(X)} + \|F_{A_0}^+\|_{L^p(X)} + 16\delta\kappa_p\|v\|_{W_{A_1}^{2,p}(X)}
+ 2\delta\kappa_p \|d_{A_0}^{+,*}v\|_{W_{A_0}^{1,p}(X)}
\\
&\leq \|F_{A_0}\|_{L^p(X)} + \|F_{A_0}^+\|_{L^p(X)} + 16\delta\kappa_p\|v\|_{W_{A_1}^{2,p}(X)}
+ 2\delta\kappa_p \|v\|_{W_{A_0}^{2,p}(X)}
\\
&\leq \|F_{A_0}\|_{L^p(X)} + \|F_{A_0}^+\|_{L^p(X)} + 32\delta\kappa_p\|v\|_{W_{A_1}^{2,p}(X)}
\quad\hbox{(by Lemma \ref{lem:Equivalence_W1p_or_W2p_norms_connections_A0_and_A1_L4_or_W14_close})}.
\end{align*}
Combining the preceding estimate with our bound \eqref{eq:Apriori_estimate_v_WA2_p} for $\|v\|_{W_{A_1}^{2,p}(X)}$ yields
\begin{multline*}
\|F(A_\asd)\|_{L^p(X)} \leq \|F_{A_0}\|_{L^p(X)} + \|F_{A_0}^+\|_{L^p(X)}
\\
+ 32\delta\kappa_p\left(128C_1(8+cc_p)\left(1+\mu(A_0)^{-1}\right)\|F_{A_0}^+\|_{L^p(X)}\right),
\end{multline*}
and thus we obtain \eqref{eq:Apriori_estimate_FAasd_Lp}, with constant $C = C([A_1],g,p,\mu(A_0))$.
\end{proof}

\subsection{Continuity of the least eigenvalue of $d_A^+d_A^{+,*}$ with respect to the connection}
\label{subsec:Continuity_least_eigenvalue_d_A+d_A+*_wrt_connection}
In order to extend Lemma \ref{lem:Positive_metric_implies_positive_lower_bound_small_eigenvalues}
in the sequel from the simple case of a \emph{positive} Riemannian metric, $g$, though arbitrary compact Lie group, $G$, to the more difficult case of a \emph{generic} Riemannian metric, $g$, and Lie groups $G=\SU(2)$ or $\SO(3)$, we shall need to closely examine continuity properties of the least eigenvalue of the elliptic operator $d_A^+d_A^{+,*}$ with respect to the connection, $A$. We begin with the simplest results and then develop continuity properties and bounds of increasing generality.

\begin{lem}[$L^4$ lower semicontinuity of the least eigenvalue of $d_A^+d_A^{+,*}$ with respect to the connection]
\label{lem:L4-continuity_least_eigenvalue_wrt_connection}
Let $X$ be a closed, four-dimensional, oriented, smooth manifold with Riemannian metric, $g$. Then there are positive constants, $c = c(g)$ and $\eps = \eps(g) \in (0,1]$, with the following significance. Let $G$ be a compact Lie group and $P$ a principal $G$-bundle over $X$. If $A_0$ is an $H^1$ connection on $P$ obeying the curvature bound \eqref{eq:L2norm_FA+_leq_small}
and $A$ is an $H^1$ connection on $P$ such that
\begin{equation}
\label{eq:L4_norm_AminusA0_leq_small}
\|A-A_0\|_{L^4(X)} \leq \eps,
\end{equation}
then
$$
\mu(A) \geq \left(1 - c\|A-A_0\|_{L^4(X)}\right)\mu(A_0) - c\|A-A_0\|_{L^4(X)}.
$$
where $\mu(A)$ is as in \eqref{eq:Least_eigenvalue_dA+dA+*}.
\end{lem}

\begin{proof}
For convenience, write $a := A-A_0 \in L^4(X;\Lambda^1\otimes\ad P)$. For $v \in H_{A_0}^1(X;\Lambda^+\otimes\ad P)$, we have $d_A^{+,*}v = d_{A_0}^{+,*}v - *[a, v]$ and
\begin{align*}
\|d_A^{+,*}v\|_{L^2(X)} &= \|d_{A_0}^{+,*}v - *[a, v]\|_{L^2(X)}
\\
&\geq \|d_{A_0}^{+,*}v\|_{L^2(X)} - \|[a, v]\|_{L^2(X)}
\\
&\geq \|d_{A_0}^{+,*}v\|_{L^2(X)} - 2\|a\|_{L^4(X)}\|v\|_{L^4(X)}
\\
&\geq \|d_{A_0}^{+,*}v\|_{L^2(X)} - 2c_1\|a\|_{L^4(X)}\|v\|_{H_{A_0}^1(X)},
\end{align*}
where $c_1=c_1(g)$ is the Sobolev embedding constant for $H^1(X) \hookrightarrow L^4(X)$ \cite[Theorem 4.12]{AdamsFournier} and we apply that and the Kato Inequality \eqref{eq:FU_6-20_first-order_Kato_inequality} to achieve the last inequality. Applying the \apriori estimate \eqref{eq:Feehan_Leness_6-6-1_H1Av_L2normdA+*v} for $\|v\|_{H_{A_0}^1(X)}$ from Lemma \ref{lem:Feehan_Leness_6-6}, with $c=c(g)$ and small enough $\eps=\eps(g) \in (0,1]$, yields
$$
\|v\|_{H_{A_0}^1(X)} \leq c\left(\|d_{A_0}^{+,*}v\|_{L^2(X)} + \|v\|_{L^2(X)}\right).
$$
Combining the preceding two inequalities gives
$$
\|d_A^{+,*}v\|_{L^2(X)}
\geq
\|d_{A_0}^{+,*}v\|_{L^2(X)} - 2c_1\|a\|_{L^4(X)}\|d_{A_0}^{+,*}v\|_{L^2(X)}
- cc_1\|a\|_{L^4(X)}\|v\|_{L^2(X)}.
$$
Now take $v$ to be an eigenvector of $d_A^+d_A^{+,*}$ with eigenvalue $\mu(A)$ and $\|v\|_{L^2(X)} = 1$ and also suppose that $\|A-A_0\|_{L^4(X)}$ is small enough that $2c_1\|a\|_{L^4(X)} \leq 1$. The preceding inequality then gives
$$
\mu(A)
\geq
\left(1 - 2c_1\|a\|_{L^4(X)}\right)\|d_{A_0}^{+,*}v\|_{L^2(X)} - cc_1\|a\|_{L^4(X)}.
$$
Since $\|v\|_{L^2(X)} = 1$, we have $\|d_{A_0}^{+,*}v\|_{L^2(X)} \geq \mu(A_0)$ by \eqref{eq:Least_eigenvalue_dA+dA+*} and the conclusion follows, for a possibly larger, relabelled constant.
\end{proof}

Lemma \ref{lem:L4-continuity_least_eigenvalue_wrt_connection} is not symmetric with respect to the connections $A$ and $A_0$. However, by replacing the $L^4(X)$ norm with the $H_{A_0}^1(X)$ norm, we obtain a continuity result that is symmetric with respect to $A$ and $A_0$.

\begin{lem}[$H^1$ continuity of the least eigenvalue of $d_A^+d_A^{+,*}$ with respect to the connection]
\label{lem:H^1-continuity_least_eigenvalue_wrt_connection}
Let $X$ be a closed, four-dimensional, oriented, smooth manifold with Riemannian metric, $g$. Then there are positive constants, $c = c(g)$ and $\eps = \eps(g) \in (0,1]$, with the following significance. Let $G$ be a compact Lie group and $P$ a principal $G$-bundle over $X$. If $A_0$ is an $H^1$ connection on $P$ that obeys the curvature bound \eqref{eq:L2norm_FA+_leq_small} and $A$ is an $H^1$ connection on $P$ such that
\begin{equation}
\label{eq:HA01_norm_AminusA_0_leq_small}
\|A-A_0\|_{H_{A_0}^1(X)} \leq \eps,
\end{equation}
then
\begin{multline*}
\left(1 - c\|A-A_0\|_{L^4(X)}\right)\mu(A_0) - c\|A-A_0\|_{L^4(X)}
\\
\leq
\mu(A)
\leq
\left(1 + c\|A-A_0\|_{L^4(X)}\right)^{-1}\left(\mu(A_0) + c\|A-A_0\|_{L^4(X)}\right).
\end{multline*}
\end{lem}

\begin{proof}
For convenience, write $a := A-A_0 \in H_{A_0}^1(X;\Lambda^1\otimes\ad P)$. The lower bound for $\mu(A)$ follows from Lemma \ref{lem:L4-continuity_least_eigenvalue_wrt_connection}. To obtain the upper bound for $\mu(A)$, observe that $F_A^+ = F_{A_0}^+ + d_{A_0}^+a + (a\wedge a)^+$ and thus
\begin{align*}
\|F_A^+\|_{L^2(X)}
&\leq
\|F_{A_0}^+\|_{L^2(X)} + \|d_{A_0}^+a\|_{L^2(X)} + 2\|a\|_{L^4(X)}^2
\\
&\leq \|F_{A_0}^+\|_{L^2(X)} + c'\|a\|_{H_{A_0}^1(X)}\left(1 + \|a\|_{H_{A_0}^1(X)}\right)
\\
&\leq (1+c')\eps_0 \quad\hbox{(by \eqref{eq:L2norm_FA+_leq_small} and \eqref{eq:HA01_norm_AminusA_0_leq_small})},
\end{align*}
where $c' = c'(g)$ and we applied, as usual, the Kato Inequality \eqref{eq:FU_6-20_first-order_Kato_inequality} and Sobolev embedding $H^1(X)\hookrightarrow L^4(X)$ \cite[Theorem 4.12]{AdamsFournier}. Hence, $A$ obeys the condition \eqref{eq:L2norm_FA+_leq_small}, for a constant $\eps' := (1+c')\eps \in (0,1]$ (for small enough $\eps$). Therefore, applying Lemma \ref{lem:L4-continuity_least_eigenvalue_wrt_connection} with the roles of $A$ and $A_0$ interchanged yields the inequality,
$$
\mu(A_0) \geq \left(1 - c\|a\|_{L^4(X)}\right)\mu(A) - c\|a\|_{L^4(X)},
$$
and hence we obtain the desired upper bound for $\mu(A)$.
\end{proof}

While Lemma \ref{lem:L4-continuity_least_eigenvalue_wrt_connection} provides a basic continuity result for $\mu(A)$, it is nonetheless very useful to further weaken the topology in which $A$ is close to $A_0$. To that end, we provide the

\begin{lem}[$L^p$ continuity of the least eigenvalue of $d_A^+d_A^{+,*}$ with respect to the connection for $2<p\leq\infty$]
\label{lem:Lp-continuity_least_eigenvalue_wrt_connection}
Let $X$ be a closed, four-dimensional, oriented, smooth manifold with Riemannian metric, $g$, and $p \in (2,\infty]$. Then there are positive constants, $c = c(g,p)$ and $\eps = \eps(g,p) \in (0,1]$, with the following significance. Let $G$ be a compact Lie group and $P$ a principal $G$-bundle over $X$. If $A_0$ and $A$ are connections on $P$ of class $H^1\cap L^p$ that each obey the curvature bound \eqref{eq:L2norm_FA+_leq_small}, then
\begin{multline}
\label{eq:Lp-continuity_least_eigenvalue_wrt_connection}
\sqrt{\mu(A_0)} - c\left(1 + \mu(A_0)
+ (1+\mu(A_0))^2\|A-A_0\|_{L^p(X)}^2\right)\|A-A_0\|_{L^p(X)}
\\
\leq
\sqrt{\mu(A)}
\leq
\sqrt{\mu(A_0)} + c(1+\mu(A_0))\|A-A_0\|_{L^p(X)},
\end{multline}
where $\mu(A)$ is as in \eqref{eq:Least_eigenvalue_dA+dA+*}.
\end{lem}

\begin{proof}
For convenience, write $a := A-A_0 \in L^p(X;\Lambda^1\otimes\ad P)$. Define $q \in [2,\infty)$ by $1/2 = 1/p+1/q$ and consider $v \in H_A^1(X;\Lambda^+\otimes\ad P)$. We use $d_A^{+,*}v = d_{A_0}^{+,*}v - *[a, v]$ and the triangle and H\"older inequalities to give
\begin{align*}
\|d_A^{+,*}v\|_{L^2(X)} &= \|d_{A_0}^{+,*}v - *[a, v]\|_{L^2(X)}
\\
&\geq \|d_{A_0}^{+,*}v\|_{L^2(X)} - \|[a, v]\|_{L^2(X)}
\\
&\geq \|d_{A_0}^{+,*}v\|_{L^2(X)} - 2\|a\|_{L^p(X)}\|v\|_{L^q(X)}.
\end{align*}
For $p > 4$, we have $2\leq q < 4$ and $\|v\|_{L^q(X)} \leq (\Vol_g(X))^{1/q-1/4}\|v\|_{L^4(X)}$, while for $2<p\leq 4$, we have $4\leq q<\infty$. Therefore, it suffices to consider the case $4\leq q<\infty$. Applying the \apriori estimate \eqref{eq:Feehan_5-3-1_Lp_dA+dA+*} for $\|v\|_{L^q(X)}$ in terms of $\|d_A^+d_A^{+,*}v\|_{L^r(X)}$ from Lemma \ref{lem:Feehan_5-3-1_Lp_dA+dA+*}, with small enough $\eps=\eps(g,q) \in (0,1]$ and $r \in [4/3,2)$ defined by $1/r = 1/2 + 1/q$, yields
$$
\|v\|_{L^q(X)} \leq c_1\left(\|d_A^+d_A^{+,*}v\|_{L^r(X)} + \|v\|_{L^r(X)}\right),
$$
for $c_1 = c_1(g,q)$. By combining the preceding two inequalities we find that
\begin{align*}
\|d_A^{+,*}v\|_{L^2(X)}
&\geq
\|d_{A_0}^{+,*}v\|_{L^2(X)} - 2c_1\|a\|_{L^p(X)}\|d_A^+d_A^{+,*}v\|_{L^r(X)}
- 2c_1\|a\|_{L^p(X)}\|v\|_{L^r(X)}
\\
&\geq
\|d_{A_0}^{+,*}v\|_{L^2(X)}
- c_0\|a\|_{L^p(X)}\left(\|d_A^+d_A^{+,*}v\|_{L^2(X)} + \|v\|_{L^2(X)}\right),
\end{align*}
for $c_0 = c_0(g,q) = 2c_1(\Vol_g(X))^{1/q}$, using the fact that $\|v\|_{L^r(X)} \leq (\Vol_g(X))^{1/q}\|v\|_{L^2(X)}$ for $r\in [4/3,2)$ and $1/r=1/2+1/q$. By taking $v \in H_A^1(X;\Lambda^+\otimes\ad P)$ to be an eigenvector of $d_A^+d_A^{+,*}$ with eigenvalue $\mu(A)$ such that $\|v\|_{L^2(X)} = 1$ and noting that
$$
\|d_A^+d_A^{+,*}v\|_{L^2(X)} = \mu(A)
\quad\hbox{and}\quad
\|d_A^{+,*}v\|_{L^2(X)} = \sqrt{\mu(A)},
$$
by \eqref{eq:Least_eigenvalue_dA+dA+*}, we obtain the bound,
$$
\sqrt{\mu(A)}
\geq
\|d_{A_0}^{+,*}v\|_{L^2(X)} - c_0(1+\mu(A))\|a\|_{L^p(X)}.
$$
But $\|d_{A_0}^{+,*}v\|_{L^2(X)} \geq \sqrt{\mu(A_0)}$ by \eqref{eq:Least_eigenvalue_dA+dA+*} and thus we have the inequality,
$$
\sqrt{\mu(A)}
\geq
\sqrt{\mu(A_0)} - c_0(1+\mu(A))\|a\|_{L^p(X)}.
$$
Interchanging the roles of $A$ and $A_0$ in the preceding derivation yields,
$$
\sqrt{\mu(A)}
\leq
\sqrt{\mu(A_0)} + c_0(1+\mu(A_0))\|a\|_{L^p(X)},
$$
the desired upper bound for $\sqrt{\mu(A)}$ in \eqref{eq:Lp-continuity_least_eigenvalue_wrt_connection}, after relabelling the constant. By substituting the resulting upper bound for $\mu(A)$ in the preceding lower bound for $\sqrt{\mu(A)}$, we discover that
$$
\sqrt{\mu(A)}
\geq
\sqrt{\mu(A_0)} - c_0\left(1 + 2\mu(A_0)
+ 2c_0^2(1+\mu(A_0))^2\|a\|_{L^p(X)}^2\right)\|a\|_{L^p(X)}.
$$
This gives the desired lower bound for $\sqrt{\mu(A)}$ in \eqref{eq:Lp-continuity_least_eigenvalue_wrt_connection} after relabelling the constants.
\end{proof}

We now wish to extend Lemma \ref{lem:Lp-continuity_least_eigenvalue_wrt_connection} so we can accommodate the weak notion of convergence described by Sedlacek in his \cite[Theorem 3.1]{Sedlacek} which, in contrast to the Uhlenbeck convergence as defined in \cite[Condition 4.4.2]{DK}, essentially replaces the usual strong convergence in $W_{\loc}^{k,p}(X\less\Sigma)$ of connections with $k\geq 1$ and $p\geq 2$ obeying $kp>4$ by weak convergence in $H_{\loc}^1(X\less\Sigma)$ and strong convergence in $L_{\loc}^p(X\less\Sigma)$ for $p\in [2,4)$, where $\Sigma = \{x_1,\ldots,x_l\} \subset X$ is a finite set of points where the curvature densities, $|F_{A_m}|^2$, concentrate as $m\to\infty$.

\begin{prop}[$L_{\loc}^p$ continuity of the least eigenvalue of $d_A^+d_A^{+,*}$ with respect to the connection for $2\leq p < 4$]
\label{prop:Lp_loc_continuity_least_eigenvalue_wrt_connection}
Let $X$ be a closed, connected, four-dimensional, oriented, smooth manifold with Riemannian metric, $g$. Then there are a positive constant $c = c(g) \in [1, \infty)$ and a constant $\eps = \eps(g) \in (0,1]$ such that the following holds. Let $G$ be a compact Lie group, $A_0$ a connection of class $H^1$ on a principal $G$-bundle $P_0$ over $X$ obeying the curvature bound \eqref{eq:L2norm_FA+_leq_small} with constant $\eps$, and $L\geq 1$ an integer, and $p \in [2,4)$. Then there are constants $c_p = c_p(g,p) \in [1, \infty)$ and $\delta = \delta(\mu(A_0),g,L,p) \in (0,1]$ and $\rho_0 = \rho_0(\mu(A_0),g,L)\in (0, 1\wedge \Inj(X,g)]$ with the following significance. Let $\rho \in (0,\rho_0]$ and $\Sigma = \{x_1,\ldots,x_L\} \subset X$ be such that
$$
\dist_g(x_l,x_k) \geq \rho \quad\hbox{for all } k\neq l,
$$
and let $U \subset X$ be the open subset given by
$$
U := X \less \bigcup_{l=1}^L \bar B_{\rho/2}(x_l).
$$
Let $P$ be a principal $G$-bundle over $X$ such that there is an isomorphism of principal $G$-bundles, $u:P\restriction X\less\Sigma \cong P_0\restriction X\less\Sigma$, and identify $P\restriction X\less\Sigma$ with $P_0\restriction X\less\Sigma$ using this isomorphism. Let $A$ be a connection of class $H^1$ on $P$ obeying the curvature bound \eqref{eq:L2norm_FA+_leq_small} with constant $\eps$ such that
\begin{equation}
\label{eq:Lp_norm_AminusA0_U_leq_small}
\|A-A_0\|_{L^p(U)} \leq \delta.
\end{equation}
Then $\mu(A)$ in \eqref{eq:Least_eigenvalue_dA+dA+*} satisfies the \emph{lower} bound,
\begin{multline}
\label{eq:Lower_bound_sqrt_muA_for_A_Lp_loc_near_A0}
\sqrt{\mu(A)}
\geq
\sqrt{\mu(A_0)} - c\sqrt{L}\,\rho^{1/6}(\mu(A)+1)
\\
- cL\rho\left(\sqrt{\mu(A)}+1\right) - c_p\|A-A_0\|_{L^p(U)}(\mu(A) + 1),
\end{multline}
and \emph{upper} bound,
\begin{multline}
\label{eq:Upper_bound_sqrt_muA_for_A_Lp_loc_near_A0}
\sqrt{\mu(A)}
\leq
\sqrt{\mu(A_0)} + c\sqrt{L}\,\rho^{1/6}(\mu(A_0)+1)
\\
+ cL\rho\left(\sqrt{\mu(A_0)}+1\right) + c_p\|A-A_0\|_{L^p(U)}(\mu(A_0) + 1).
\end{multline}
\end{prop}

\begin{proof}
Since the argument is lengthy, we divide it into several steps.

\setcounter{step}{0}
\begin{step}[The eigenvalue identity]
By hypothesis, $X = U \cup (\cup_{l=1}^L \bar B_{\rho/2}(x_l))$ and we may choose a $C^\infty$ partition of unity, $\{\chi_l\}_{l=0}^L$, for $X$ subordinate to the open cover of $X$ given by $U$ and the open balls, $B_\rho(x_l)$ for $1\leq l\leq L$, such that $\chi_l = 1$ on $B_{\rho/2}(x_l)$ and $\supp\chi_l \subset B_\rho(x_l)$ for $1\leq l\leq L$ and $\supp\chi_0 \subset U$, while $\sum_{l=0}^L\chi_l = 1$ on $X$ and $0\leq \chi_l\leq 1$ for $0\leq l\leq L$. In addition, we may suppose that there is a constant, $c=c(g)$, such that
\begin{equation}
\label{eq:Pointwise_bound_dchi_l}
|d\chi_l| \leq \frac{c}{\rho} \quad\hbox{on } X, \quad\hbox{for } 0 \leq l \leq L.
\end{equation}
We consider $v \in \Omega^+(X;\ad P)$ and write $v = \sum_{l=0}^L\chi_l v$. Because $\supp\chi_k \cap \supp\chi_l = \emptyset$ for all $k\neq l$ with $1 \leq k,l \leq L$ and $\chi_0 = 1-\chi_l$ on $B_\rho(x_l)$ for $1\leq l \leq l$, we see that
\begin{align*}
\|d_A^{+,*}v\|_{L^2(X)}^2
&=
\sum_{l=0}^L\|d_A^{+,*}(\chi_lv)\|_{L^2(X)}^2
+ 2\sum_{k<l} \left(d_A^{+,*}(\chi_kv), d_A^{+,*}(\chi_lv)\right)_{L^2(X)}
\\
&=
\sum_{l=0}^L\|d_A^{+,*}(\chi_lv)\|_{L^2(X)}^2
+ 2\sum_{l=1}^L \left(d_A^{+,*}(\chi_0v), d_A^{+,*}(\chi_lv)\right)_{L^2(X)}
\\
&=
\sum_{l=0}^L\|d_A^{+,*}(\chi_lv)\|_{L^2(X)}^2
+ 2\sum_{l=1}^L \left(d_A^{+,*}((1-\chi_l)v), d_A^{+,*}(\chi_lv)\right)_{L^2(X)},
\end{align*}
and hence,
$$
\|d_A^{+,*}v\|_{L^2(X)}^2
= \|d_A^{+,*}(\chi_0v)\|_{L^2(X)}^2 - \sum_{l=1}^L\|d_A^{+,*}(\chi_lv)\|_{L^2(X)}^2
+ 2\sum_{l=1}^L \left(d_A^{+,*}v, d_A^{+,*}(\chi_lv)\right)_{L^2(X)}.
$$
We now choose $v \in H_A^1(X;\Lambda^+\otimes\ad P)$ with $\|v\|_{L^2(X)}=1$ to be an eigenvector for the least eigenvalue $\mu(A)$ of $d_A^+d_A^{+,*}$. Hence, the preceding identity and \eqref{eq:Least_eigenvalue_dA+dA+*} yield
\begin{equation}
\label{eq:Eigenvalue_identity_muA_for_A_Lp_loc_near_A0}
\mu(A) = \|d_A^{+,*}(\chi_0v)\|_{L^2(X)}^2 - \sum_{l=1}^L\|d_A^{+,*}(\chi_lv)\|_{L^2(X)}^2
+ 2\sum_{l=1}^L \left(d_A^{+,*}v, d_A^{+,*}(\chi_lv)\right)_{L^2(X)},
\end{equation}
the basic eigenvalue identity for $\mu(A)$.
\end{step}

To proceed further, we need a lower bound for the expression $\|d_A^{+,*}(\chi_0v)\|_{L^2(X)}^2$ in \eqref{eq:Eigenvalue_identity_muA_for_A_Lp_loc_near_A0} in terms of $\mu(A_0)$ and small upper bounds for the remaining terms on the right-hand side of \eqref{eq:Eigenvalue_identity_muA_for_A_Lp_loc_near_A0}.

\begin{step}[Upper bound for the terms $(d_A^{+,*}v, d_A^{+,*}(\chi_lv))_{L^2(X)}$ when $1\leq l\leq L$]
Observe that
$$
\left(d_A^{+,*}v, d_A^{+,*}(\chi_lv)\right)_{L^2(X)}
=
\left(d_A^+d_A^{+,*}v, \chi_lv\right)_{L^2(X)}
=
\mu(A)\left(v, \chi_lv\right)_{L^2(X)},
$$
and thus, for $1\leq l\leq L$,
\begin{align*}
\left|\left(d_A^{+,*}v, d_A^{+,*}(\chi_lv)\right)_{L^2(X)}\right|
&\leq
\mu(A) \|\sqrt{\chi_l} v\|_{L^2(X)}^2
\\
&\leq
\mu(A) \left(\Vol_g(\supp\chi_l)\right)^{1/2}\|v\|_{L^4(X)}^2
\\
&\leq c\rho^2\mu(A)\|v\|_{L^4(X)}^2
\\
&\leq c\rho^2\mu(A)\left(\|d_A^{+,*}v\|_{L^2(X)}+\|v\|_{L^2(X)}\right)^2
\quad\hbox{(by \eqref{eq:Feehan_Leness_6-6-1_L4v_L2dA+*v})}.
\end{align*}
Thus, noting that $\|d_A^{+,*}v\|_{L^2(X)} = \sqrt{\mu(A)}\|v\|_{L^2(X)}$ by \eqref{eq:Least_eigenvalue_dA+dA+*} and $\|v\|_{L^2(X)}=1$, we have
\begin{equation}
\label{eq:Upper_bound_dA*v_dA*chiv_L2innerproduct_for_A_Lp_loc_near_A0}
\left|\left(d_A^{+,*}v, d_A^{+,*}(\chi_lv)\right)_{L^2(X)}\right|
\leq
c\rho^2\mu(A)(1+\mu(A)), \quad\hbox{for } 1\leq l\leq L,
\end{equation}
where $c=c(g)$.
\end{step}

\begin{step}[Upper bound for the terms $\|d_A^{+,*}(\chi_lv)\|_{L^2(X)}^2$ when $1\leq l\leq L$]
For $0\leq l\leq L$, we have
$$
d_{A}^{+,*}(\chi_{l}v)
=
-*d_{A}*(\chi_{l}v)
=
-*(d\chi_{l}\wedge v + \chi_{l}d_{A}v)
=
-*(d\chi_{l}\wedge v) + \chi_{l}d_{A}^{+,*}v.
$$
Hence, for $s \in (2,4)$ and $t\in (4,\infty)$ obeying $1/2=1/s+1/t$, we have
$$
\left| \|d_A^{+,*}(\chi_lv)\|_{L^2(X)} - \|\chi_ld_A^{+,*}v\|_{L^2(X)} \right|
\leq
\|d\chi_l\|_{L^s(X)} \|v\|_{L^t(X)}.
$$
The pointwise bound \eqref{eq:Pointwise_bound_dchi_l} for $d\chi_l$ implies that there is a positive constant, $c = c(g)$ when $1\leq l\leq L$ and $c = Lc_0(g)$ when $l=0$, such that, for any $u \in [1,\infty]$,
\begin{equation}
\label{eq:Lp_norm_dchi_l}
\|d\chi_l\|_{L^u(X)} \leq c\rho^{(4/u)-1}, \quad 0\leq l\leq L.
\end{equation}
We choose $s=u=3$ and $t=6$ and combine the two preceding inequalities to give
\begin{align*}
{}&\left| \|d_A^{+,*}(\chi_lv)\|_{L^2(X)} - \|\chi_ld_A^{+,*}v\|_{L^2(X)} \right|
\\
&\quad \leq
c\rho^{1/3} \|v\|_{L^6(X)}
\\
&\quad \leq c\rho^{1/3} \left(\|d_A^+d_A^{+,*}v\|_{L^{3/2}(X)} + \|v\|_{L^{3/2}(X)}\right)
\quad\hbox{(by \eqref{eq:Feehan_5-3-1_Lp_dA+dA+*} with $r=3/2$)}
\\
&\quad = c\rho^{1/3}(\mu(A) + 1)\|v\|_{L^{3/2}(X)},
\end{align*}
where we used the fact that $d_A^+d_A^{+,*}v = \mu(A)v$ in the last equality. Thus, for a positive constant, $c = c(g)$ when $1\leq l\leq L$ and $c = Lc_0(g)$ when $l=0$,
\begin{equation}
\label{eq:L2norm_dA+*chi_ellv_minus_L2norm_chi_elldA+*v}
\left| \|d_A^{+,*}(\chi_lv)\|_{L^2(X)} - \|\chi_ld_A^{+,*}v\|_{L^2(X)} \right|
\leq c\rho^{1/3}(\mu(A) + 1), \quad\hbox{for } 0\leq l\leq L.
\end{equation}
Restricting now to $1\leq l\leq L$, we have
\begin{align*}
\|\chi_ld_A^{+,*}v\|_{L^2(X)}^2
&=
\left(\chi_l^2 d_A^{+,*}v, d_A^{+,*}v\right)_{L^2(X)}
\\
&=
\left(\chi_l^2 d_A^+d_A^{+,*}v
+ 2\chi_l\left(d\chi_l\wedge d_A^{+,*}v\right)^+, v\right)_{L^2(X)}
\\
&=
\mu(A)\|\chi_lv\|_{L^2(X)}^2 + 2\left(\chi_ld\chi_l\wedge d_A^{+,*}v, v\right)_{L^2(X)}
\\
&\leq \mu(A)(\Vol_g(\supp\chi_l))^{1/2}\|v\|_{L^4(X)}^2
+ 2\left(\chi_ld\chi_l\wedge d_A^{+,*}v, v\right)_{L^2(X)}.
\end{align*}
Therefore, applying \eqref{eq:Feehan_Leness_6-6-1_L4v_L2dA+*v} in the preceding inequality together with the facts that $\|d_A^{+,*}v\|_{L^2(X)} = \sqrt{\mu(A)}\|v\|_{L^2(X)}$ by \eqref{eq:Least_eigenvalue_dA+dA+*} and $\|v\|_{L^2(X)} = 1$,
\begin{equation}
\label{eq:Upper_bound_chidA*v_L2norm_for_A_Lp_loc_near_A0}
\|\chi_ld_A^{+,*}v\|_{L^2(X)}^2
\leq
c\rho^2\mu(A)\left(\sqrt{\mu(A)} + 1\right)^2
+ 2\left(\chi_ld\chi_l\wedge d_A^{+,*}v, v\right)_{L^2(X)},
\quad\hbox{for } 1\leq l\leq L.
\end{equation}
The inner product term in \eqref{eq:Upper_bound_chidA*v_L2norm_for_A_Lp_loc_near_A0} is bounded via
\begin{align*}
\left|\left(\chi_ld\chi_l\wedge d_A^{+,*}v, v\right)_{L^2(X)}\right|
&\leq
\|d\chi_l\|_{L^3(X)} \|d_A^{+,*}v\|_{L^2(X)} \|v\|_{L^6(X)}
\\
&=
\sqrt{\mu(A)}\, \|d\chi_l\|_{L^3(X)} \|v\|_{L^6(X)}
\quad\hbox{(by \eqref{eq:Least_eigenvalue_dA+dA+*})}
\\
&\leq c\rho^{1/3}\sqrt{\mu(A)}(\mu(A)+1)
\quad\hbox{(by \eqref{eq:Feehan_5-3-1_Lp_dA+dA+*} and \eqref{eq:Lp_norm_dchi_l}),}
\end{align*}
with $q=6$ and $r=3/2$ in \eqref{eq:Feehan_5-3-1_Lp_dA+dA+*}.
Hence, substituting the preceding inequality in \eqref{eq:Upper_bound_chidA*v_L2norm_for_A_Lp_loc_near_A0} yields
$$
\|\chi_ld_A^{+,*}v\|_{L^2(X)}^2
\leq
c\left(\rho^2\mu(A) + \rho^{1/3}\sqrt{\mu(A)}\right)(\mu(A)+1),
\quad\hbox{for } 1\leq l\leq L.
$$
By combining the preceding estimate with \eqref{eq:L2norm_dA+*chi_ellv_minus_L2norm_chi_elldA+*v} (and the elementary inequality, $x^2 \leq 2(x-y)^2 + 2y^2$ for $x,y\in\RR$) we obtain
\begin{align*}
\|d_A^{+,*}(\chi_lv)\|_{L^2(X)}^2
&\leq
2\left| \|d_A^{+,*}(\chi_lv)\|_{L^2(X)} - \|\chi_ld_A^{+,*}v\|_{L^2(X)} \right|^2
+ 2\|\chi_ld_A^{+,*}v\|_{L^2(X)}^2
\\
&\leq c\rho^{2/3}(\mu(A) + 1)^2
+ c\left(\rho^2\mu(A) + \rho^{1/3}\sqrt{\mu(A)}\right)(\mu(A)+1),
\end{align*}
and thus, noting that $\rho \in (0,1]$,
\begin{equation}
\label{eq:Upper_bound_dA*chiv_L2norm_for_A_Lp_loc_near_A0}
\|d_A^{+,*}(\chi_lv)\|_{L^2(X)}^2
\leq
c\rho^{1/3}(\mu(A)+1)^2, \quad\hbox{for } 1\leq l\leq L,
\end{equation}
for $c=c(g)$. This completes our analysis of all terms on the right-hand side of \eqref{eq:Eigenvalue_identity_muA_for_A_Lp_loc_near_A0} with $l\neq 0$.
\end{step}

\begin{step}[Lower bound for the term $\|d_A^{+,*}(\chi_0v)\|_{L^2(X)}$ and preliminary lower bound for $\mu(A)$]
Without loss of generality in the remainder of the proof, we may restrict attention to  $p \in (2,4]$. For convenience, we write  $a := A-A_0 \in H_{A_0}^1(X;\Lambda^1\otimes\ad P)$. For the term in \eqref{eq:Eigenvalue_identity_muA_for_A_Lp_loc_near_A0} with $l=0$, we note that $d_A^{+,*}v = d_{A_0}^{+,*}v - *(a\wedge v)$ on $X\less\Sigma$ and thus, for $p\in (2,4]$ and $q\in [4,\infty)$ defined by $1/2=1/p+1/q$ and $r\in [4/3,2)$ defined by $1/r = 1/2+1/q$,
\begin{align*}
{}&\left|\|d_A^{+,*}(\chi_0v)\|_{L^2(X)} - \|d_{A_0}^{+,*}(\chi_0v)\|_{L^2(X)}\right|
\\
&\quad \leq
\|*(a\wedge \chi_0v)\|_{L^2(X)}
\\
&\quad \leq 2\|a\|_{L^p(\supp\chi_0)} \|v\|_{L^q(X)}
\\
&\quad \leq c_p\|a\|_{L^p(U)} \left(\|d_A^+d_A^{+,*}v\|_{L^r(X)} + \|v\|_{L^r(X)}\right)
\quad\hbox{(by \eqref{eq:Feehan_5-3-1_Lp_dA+dA+*})}
\\
&\quad = c_p\|a\|_{L^p(U)} (\mu(A) + 1)\|v\|_{L^r(X)}
\quad\hbox{(by \eqref{eq:Least_eigenvalue_dA+dA+*})},
\end{align*}
where we used the fact that $\supp\chi_0 \subset U$ by construction and $c_p=c_p(g,p)$. Thus, noting that
$$
\|v\|_{L^r(X)} \leq \left(\Vol_g(X)\right)^{1/q}\|v\|_{L^2(X)},
$$
and $\|v\|_{L^2(X)}=1$, we see that
\begin{equation}
\label{eq:Upper_bound_dA*chi0v_L2norm_minus_dA0*chi0v_L2norm_for_A_Lp_loc_near_A0}
\left|\|d_A^{+,*}(\chi_0v)\|_{L^2(X)} - \|d_{A_0}^{+,*}(\chi_0v)\|_{L^2(X)}\right|
\leq
c_p\|a\|_{L^p(U)} (\mu(A) + 1),
\end{equation}
for a positive constant, $c_p=c_p(g,p)$. Next, we observe that
\begin{align*}
\|d_{A_0}^{+,*}(\chi_0v)\|_{L^2(X)}
&\geq
\sqrt{\mu(A_0)}\|\chi_0v\|_{L^2(X)} \quad\hbox{(by \eqref{eq:Least_eigenvalue_dA+dA+*})}
\\
&\geq
\sqrt{\mu(A_0)}\left( \|v\|_{L^2(X)} - \sum_{l=1}^L \|\chi_lv\|_{L^2(X)}\right)
\\
&\geq
\sqrt{\mu(A_0)}\left( \|v\|_{L^2(X)}
- \|v\|_{L^4(X)}\sum_{l=1}^L \left(\Vol_g(\supp\chi_l)\right)^{1/4} \right)
\\
&\geq \sqrt{\mu(A_0)}\left( \|v\|_{L^2(X)}
- cL\rho\|v\|_{L^4(X)}\right)
\\
&\geq \sqrt{\mu(A_0)}\|v\|_{L^2(X)}
- cL\rho\left(\|d_A^{+,*}v\|_{L^2(X)}+\|v\|_{L^2(X)}\right),
\end{align*}
where $c=c(g)$ and we used \eqref{eq:Feehan_Leness_6-6-1_L4v_L2dA+*v} to obtain the preceding inequality. Therefore, because $\|d_A^{+,*}v\|_{L^2(X)} = \sqrt{\mu(A)}\|v\|_{L^2(X)}$ by \eqref{eq:Least_eigenvalue_dA+dA+*} and $\|v\|_{L^2(X)}=1$,
\begin{equation}
\label{eq:Lower_bound_dA*chi0v_L2norm_for_A_Lp_loc_near_A0}
\|d_{A_0}^{+,*}(\chi_0v)\|_{L^2(X)}
\geq \sqrt{\mu(A_0)} - cL\rho\left(\sqrt{\mu(A)}+1\right).
\end{equation}
Observe that
\[
\|d_{A_0}^{+,*}(\chi_0v)\|_{L^2(X)}
\leq
\|d_A^{+,*}(\chi_0v)\|_{L^2(X)}
+ \left| \|d_{A_0}^{+,*}(\chi_0v)\|_{L^2(X)} - \|d_A^{+,*}(\chi_0v)\|_{L^2(X)} \right|.
\]
We rewrite the preceding inequality and combine with \eqref{eq:Upper_bound_dA*chi0v_L2norm_minus_dA0*chi0v_L2norm_for_A_Lp_loc_near_A0} and
\eqref{eq:Lower_bound_dA*chi0v_L2norm_for_A_Lp_loc_near_A0} to give
\begin{align*}
\|d_A^{+,*}(\chi_0v)\|_{L^2(X)}
&\geq
\|d_{A_0}^{+,*}(\chi_0v)\|_{L^2(X)}
- \left|\|d_A^{+,*}(\chi_0v)\|_{L^2(X)} - \|d_{A_0}^{+,*}(\chi_0v)\|_{L^2(X)}\right|
\\
&\geq \sqrt{\mu(A_0)} - cL\rho\left(\sqrt{\mu(A)}+1\right)
- c_p\|a\|_{L^p(U)}(\mu(A) + 1),
\end{align*}
for positive constants, $c=c(g)$ and $c_p=c_p(g,p)$. We substitute the preceding inequality, together with \eqref{eq:Upper_bound_dA*v_dA*chiv_L2innerproduct_for_A_Lp_loc_near_A0} and \eqref{eq:Upper_bound_dA*chiv_L2norm_for_A_Lp_loc_near_A0}, in \eqref{eq:Eigenvalue_identity_muA_for_A_Lp_loc_near_A0} to discover that $\mu(A)$ obeys
\begin{multline}
\label{eq:Lower_bound_sqrt_muA_for_A_Lp_loc_near_A0_preliminary}
\mu(A)
\geq
\left(\sqrt{\mu(A_0)} - cL\rho\left(\sqrt{\mu(A)}+1\right)
- c_p\|a\|_{L^p(U)}(\mu(A) + 1)\right)^2
\\
- cL\rho^2\mu(A)(1+\mu(A)) - cL\rho^{1/3}(\mu(A)+1)^2,
\end{multline}
a preliminary lower bound for $\mu(A)$, where $c=c(g)$ and $c_p=c_p(g,p)$.
\end{step}

\begin{step}[Upper and lower bounds for $\mu(A)$]
The inequality \eqref{eq:Lower_bound_sqrt_muA_for_A_Lp_loc_near_A0_preliminary} implies an upper bound for $\mu(A_0)$ in terms of $\mu(A)$ and hence, by interchanging the roles of $A$ and $A_0$, an upper bound for $\mu(A)$ in terms of $\mu(A_0)$. Therefore, regarding $A_0$ as fixed, for small enough $\delta = \delta(\mu(A_0),g,L,p) \in (0,1]$ and $\rho_0 = \rho_0(\mu(A_0),g,L)\in (0, 1\wedge \Inj(X,g)]$, recalling that $\rho \in (0,\rho_0]$ and $\|a\|_{L^p(U)} \leq \delta$ by hypothesis, we may suppose that
$$
\sqrt{\mu(A_0)} - cL\rho\left(\sqrt{\mu(A)}+1\right)
- c_p\|a\|_{L^p(U)}(\mu(A) + 1) \geq 0.
$$
Thus, using the elementary inequality, $(x+y)^{1/2} \leq x^{1/2} + y^{1/2}$ for $x,y \geq 0$, we obtain from \eqref{eq:Lower_bound_sqrt_muA_for_A_Lp_loc_near_A0_preliminary} that
\begin{multline*}
\sqrt{\mu(A)} + \left(cL\rho^2\mu(A)(1+\mu(A)) + cL\rho^{1/3}(\mu(A)+1)^2\right)^{1/2}
\\
\geq
\sqrt{\mu(A_0)} - cL\rho\left(\sqrt{\mu(A)}+1\right)
- c_p\|a\|_{L^p(U)}(\mu(A) + 1).
\end{multline*}
The preceding inequality yields the desired lower bound \eqref{eq:Lower_bound_sqrt_muA_for_A_Lp_loc_near_A0} for $\mu(A)$, after an another application of the elementary inequality, $(x+y)^{1/2} \leq x^{1/2} + y^{1/2}$ for $x,y \geq 0$ to give
\begin{align*}
{}&\left(cL\rho^2\mu(A)(1+\mu(A)) + cL\rho^{1/3}(\mu(A)+1)^2\right)^{1/2}
\\
&\quad \leq c\sqrt{L}\,\rho\sqrt{\mu(A)}\left(1+\sqrt{\mu(A)}\right) + c\sqrt{L}\,\rho^{1/6}(\mu(A)+1)
\\
&\quad \leq c\sqrt{L}\,\rho^{1/6}(\mu(A)+1).
\end{align*}
Interchanging the roles of $A$ and $A_0$ in the preceding inequality yields the desired upper bound \eqref{eq:Upper_bound_sqrt_muA_for_A_Lp_loc_near_A0} for $\mu(A)$.
\end{step}

This completes the proof of Proposition \ref{prop:Lp_loc_continuity_least_eigenvalue_wrt_connection}.
\end{proof}

\subsection{Minimization of the Yang-Mills energy functional and application of the weak compactness result of Sedlacek}
\label{subsec:Minimization_Yang-Mills_energy_weak_compactness_Sedlacek}
We next recall two important results due to Sedlacek. In his \cite[Theorem 3.1]{Sedlacek}, Sedlacek assumes that each $A_m$ is $C^\infty$ but that assumption can easily be weakened.

\begin{thm}
\label{thm:Sedlacek_3-1}
\cite[Theorem 3.1 and Lemma 3.4]{Sedlacek}
Let $G$ be a compact Lie group and $P$ a principal $G$-bundle over a closed, connected, Riemannian, smooth four-dimensional manifold, $X$. If $\{A_m\}_{m\geq 1}$ is a sequence of connections on $P$ of class $W^{k,p}$, with $k\geq 1$ and $p\geq 2$ obeying $kp>4$, such that $\sup_{m\in\NN}\sE(A_m) < \infty$, then there are an integer $L\geq 0$ and, for $L\geq 1$, a finite set of points, $\Sigma := \{x_1,\ldots,x_L\} \subset X$, characterized by
\begin{equation}
\label{eq:Sedlacek_proof_of_theorem_3-1_characterization_bad_points}
\lim_{r\searrow 0}\limsup_{m\to\infty} \|F_{A_m}\|_{L^2(B_r(x_l))}^2 \geq \eps_1,
\quad 1\leq l \leq L,
\end{equation}
for a constant $\eps_1 \in (0,1]$ as in \cite[Theorem 3.2]{Sedlacek} (or \cite[Theorem 1.3]{UhlLp}), and sequences of
{\alphenumi
\begin{enumerate}
  \item geodesic balls, $\{B_\alpha\}_{\alpha\in\NN} \subset X\less\Sigma$, such that $\cup_{\alpha\in\NN}B_\alpha = X\less\Sigma$,
  \item sections, $\sigma_{\alpha,m}:B_\alpha \to P$;
  \item local connection one-forms,
  $a_{\alpha,m} := \sigma_{\alpha,m}^*A_m \in H_\Gamma^1(B_\alpha; \Lambda^1\otimes\fg)$; and
  \item transition functions, $u_{\alpha\beta,m} \in W_\Gamma^{1,4}(B_\alpha\cap B_\beta; G)$,
\end{enumerate}
}
such that the following hold for each $\alpha, \beta \in \NN$:
\begin{enumerate}
  \item $d_\Gamma^* a_{\alpha,m} = 0$, for all $m$ sufficiently large (depending on $\alpha$);
  \item $d_\Gamma^* a_\alpha = 0$;
  \item $u_{\alpha\beta,m} \rightharpoonup u_{\alpha\beta}$ weakly in $W_\Gamma^{1,4}(B_\alpha\cap B_\beta; G)$;
  \item $F_{\alpha,m} \rightharpoonup F_\alpha$ weakly in $L^2(B_\alpha; \Lambda^2\otimes\fg)$; and
  \item the sequence $\{a_{\alpha,m}\}_{m\in\NN}$ obeys
  \begin{enumerate}
  \item $\{a_{\alpha,m}\}_{m\in\NN} \subset H_\Gamma^1(B_\alpha; \Lambda^1\otimes\fg)$ is bounded;

  \item $a_{\alpha,m} \rightharpoonup a_\alpha$ weakly in $H^1(B_\alpha; \Lambda^1\otimes\fg)$; and

  \item $a_{\alpha,m} \to a_\alpha$ strongly in $L^p(B_\alpha; \Lambda^1\otimes\fg)$ for $1\leq p<4$,
  \end{enumerate}
\end{enumerate}
where $F_{\alpha,m} := da_{\alpha,m} + [a_{\alpha,m}, a_{\alpha,m}] = \sigma_{\alpha,m}^*F_{A_m}$ and $F_\alpha := da_\alpha + [a_\alpha, a_\alpha]$, and $d^*$ is the formal adjoint of $d$ with respect to the flat metric defined by a choice of geodesic normal coordinates on $B_\alpha$.
\end{thm}

In our statement of Theorem \ref{thm:Sedlacek_3-1}, we have strengthened Sedlacek's conclusions concerning the boundedness and convergence properties of the local connection one-forms as follows. Boundedness of the sequence $\{a_{\alpha,m}\}_{m\in\NN} \subset H_\Gamma^1(B_\alpha; \Lambda^1\otimes\fg)$ is given by \cite[Lemma 3.4]{Sedlacek} and, after passing to a subsequence, the strong convergence in $L^p(B_\alpha)$ for $1\leq p < 4$ (and not just $p=2$ as in \cite[Theorem 3.1]{Sedlacek}) is a consequence of compactness of the embedding $H^1(B) \to L^p(B)$ in dimension four by the Rellich-Kondrachov Embedding Theorem \cite[Theorem 6.3]{AdamsFournier}.

Because $W^{1,4}(X)$ does not embed in $C(X)$, it is not all obvious that the collection of local connection one-forms, $\{a_\alpha\}_{\alpha\in\NN}$, produced by Theorem \ref{thm:Sedlacek_3-1} can be patched together to form a connection on $P\restriction X\less\Sigma$. Nevertheless, Sedlacek also deduces this and more when $\{A_m\}_{m\in\NN}$ is a minimizing sequence. Recall that $\sE(A) := \frac{1}{2}\|F_A\|_{L^2(X)}^2$ denotes the energy of a connection of class $H^1$ on $P$.

\begin{defn}
\label{defn:Sedlacek_page_521}
\cite[p. 521]{Sedlacek}
Let $G$ be a compact Lie group and $P$ a principal $G$-bundle with obstruction $\eta=\eta(P)$ \cite[Section 2]{Sedlacek} over a closed, connected, Riemannian, smooth four-dimensional manifold, $X$. Let $\fm(\eta) := \inf\{\sE(A): A$ is a connection on a principal $G$-bundle with obstruction $\eta$ over $X\}$.
\end{defn}

Sedlacek applies Theorem \ref{thm:Sedlacek_3-1} to a minimizing sequence of connections to obtain

\begin{thm}
\label{thm:Sedlacek_4-3}
\cite[Proposition 4.2 and Theorem 4.3]{Sedlacek}
Let $G$ be a compact Lie group and $P$ a principal $G$-bundle over a closed, connected, smooth four-dimensional manifold, $X$, with Riemannian metric, $g$. If $\{A_m\}_{m\geq 1}$ is a \emph{minimizing sequence} of connections on $P$ of class $W^{k,p}$, with $k\geq 1$ and $p\geq 2$ obeying $kp>4$, in the sense that
$$
\sE(A_m) \searrow \fm(\eta), \quad\hbox{as } m \to \infty,
$$
then the following hold. For each $\alpha,\beta\in\NN$,
\begin{enumerate}
  \item $a_\alpha \in C^\infty(B_\alpha; \Lambda^1\otimes\fg))$ and is a solution to the Yang-Mills equations with respect to the metric $g$ over $B_\alpha$;

  \item $u_{\alpha\beta} \in C^\infty(B_\alpha\cap B_\beta; G)$;

  \item The sequences, $\{a_\alpha\}_{\alpha\in\NN}$ and $\{u_{\alpha,\beta}\}_{\alpha,\beta\in\NN}$, define a $C^\infty$ connection, $\tilde A_\infty$, on a principal $G$-bundle, $\tilde P_\infty$, over $X\less\Sigma$, where $\Sigma$ is as in Theorem \ref{thm:Sedlacek_3-1};

  \item The connection, $\tilde A_\infty$, is a solution to the Yang-Mills equation with respect to $g$; and

  \item There is a $C^\infty$ bundle automorphism, $u_\infty \in \Aut(\tilde P_\infty\restriction X\less\Sigma)$, such that $u_\infty^*\tilde A_\infty$ extends to a $C^\infty$ Yang-Mills connection, $A_\infty$, on a principal $G$-bundle, $P_\infty$, over $X$.
\end{enumerate}
\end{thm}

We now have the useful

\begin{cor}[Convergence of the least eigenvalue of $d_{A_m}^+d_{A_m}^{+,*}$ for a minimizing sequence of connections, $\{A_m\}_{m\in\NN}$]
\label{cor:Weak_H_loc1_X_less_Sigma_continuity_least_eigenvalue_wrt_connection}
Assume the hypotheses of Theorem \ref{thm:Sedlacek_4-3} and that $X$ is oriented. Then
\begin{equation}
\label{eq:Weak_H_loc1_X_less_Sigma_continuity_least_eigenvalue_wrt_connection_bounds}
\lim_{m\to\infty}\mu(A_m) = \mu(A_\infty),
\end{equation}
where $\mu(A)$ is as in \eqref{eq:Least_eigenvalue_dA+dA+*}.
\end{cor}

\begin{proof}
Proposition \ref{prop:Lp_loc_continuity_least_eigenvalue_wrt_connection} implies that, for each $\rho \in (0,\rho_0]$, where $\rho_0 = \rho_0(\mu(A_\infty),g,L)\in (0, 1\wedge \Inj(X,g)]$, we have
\begin{align*}
\liminf_{m\to\infty}\sqrt{\mu(A_m)}
&\geq
\sqrt{\mu(A_\infty)} - c\sqrt{L}\,\rho^{1/6}\left(\limsup_{m\to\infty}\mu(A_m)+1\right)
- cL\rho\left(\limsup_{m\to\infty}\sqrt{\mu(A_m)}+1\right),
\\
\limsup_{m\to\infty}\sqrt{\mu(A_m)}
&\leq
\sqrt{\mu(A_\infty)} + c\sqrt{L}\,\rho^{1/6}(\mu(A_\infty)+1)
+ cL\rho\left(\sqrt{\mu(A_\infty)}+1\right).
\end{align*}
The preceding inequalities follow from Proposition \ref{prop:Lp_loc_continuity_least_eigenvalue_wrt_connection} since Theorem \ref{thm:Sedlacek_3-1} implies that, for $p\in [2,4)$,
$$
\|\sigma_{\alpha,m}^*A_m-\sigma_\alpha A_\infty\|_{L^p(B_\alpha)} \to 0
\quad\hbox{strongly in } L^p(B_\alpha;\Lambda^1\otimes\ad P)
\quad \hbox{as } m\to\infty,
$$
for sequences of local sections, $\{\sigma_{\alpha,m}\}_{m\in\NN}$ of $P\restriction B_\alpha$, and a local section, $\sigma_\alpha$ of $P_\infty\restriction B_\alpha$, and
$$
\|A_m-A_\infty\|_{L^p(U)}
\leq
\sum_{\alpha:B_\alpha\cap U\neq\emptyset}
\|\sigma_{\alpha,m}^*A_m - \sigma_\alpha^*A_\infty\|_{L^p(B_\alpha)}
$$
where $U := X - \cup_{l=1}^L \bar B_{\rho/2}(x_l)$ is as in Proposition \ref{prop:Lp_loc_continuity_least_eigenvalue_wrt_connection} and so the index set in the preceding sum is finite. Because the lower bound for $\liminf_{m\to\infty}\sqrt{\mu(A_m)}$ and upper bound for $\limsup_{m\to\infty}\sqrt{\mu(A_m)}$ hold for every $\rho \in (0,\rho_0]$, the conclusion follows.
\end{proof}

\subsection{Minimization of the Yang-Mills energy functional and application of the weak compactness result of Taubes}
\label{subsec:Minimization_Yang-Mills_energy_weak_compactness_Taubes}
The convergence result, Corollary \ref{cor:Weak_H_loc1_X_less_Sigma_continuity_least_eigenvalue_wrt_connection}, can also be proved by replacing the role of Sedlacek's Theorem \ref{thm:Sedlacek_3-1} by that of Taubes' \cite[Proposition 4.5]{TauPath}. The advantage of \cite[Proposition 4.5]{TauPath} over Theorem \ref{thm:Sedlacek_3-1} is that one has strong convergence in $H_{\loc}^1(X\less\Sigma;\Lambda^1\otimes)$ of $u_m^*A_m$ to $A_\infty$ on $P_\infty\restriction X\less\Sigma$ as $m\to\infty$ for a sequence of $\{u_m\}_{m\in\NN}$, of bundle isomorphisms, $u_m: P_\infty\restriction X\less\Sigma \cong P\restriction X\less\Sigma$ of class $H_{\loc}^3(X\less\Sigma)$ and some (possibly empty) subset $\Sigma = \{x_1,\ldots,x_L\} \subset X$, provided $\{A_m\}_{m\in\NN}$ is a \emph{convergent Palais-Smale} sequence in the sense that
\begin{equation}
\label{eq:Convergent_Palais-Smale_sequence}
\sE(A_m) \to \sE(A_\infty) \quad\hbox{and}\quad \sE'(A_m) \to 0
\quad\hbox{as } m \to \infty,
\end{equation}
where $A_\infty$ is a Yang-Mills connection on $P_\infty$, as in Theorem \ref{thm:Sedlacek_4-3} and, as usual, $\sE(A) = \frac{1}{2}\|F_A\|_{L^2(X)}^2$ denotes the energy of a connection, $A$, on $P$. The condition, $\sE'(A_m) \to 0$ as $m\to\infty$, comprises part of Taubes' \cite[Definition 4.1]{TauPath} of a \emph{good} sequence, where our assumption, $\sE(A_m) \to \sE(A_\infty)$, in \eqref{eq:Convergent_Palais-Smale_sequence} is relaxed to a condition that $\{\sE(A_m)\}_{m\in\NN}$ is bounded. The disadvantage of this approach is that one has to first establish the existence of a convergent Palais-Smale sequence as in \eqref{eq:Convergent_Palais-Smale_sequence}, not just a minimizing sequence as in Theorem \ref{thm:Sedlacek_4-3} (and whose existence is an immediate consequence of the Definition \ref{defn:Sedlacek_page_521} of $\fm(\eta)$).

Jeanjean provides an existence theorem \cite[Theorem 1.1]{Jeanjean_1999} for bounded Palais-Smale sequences for suitable abstract functionals on a Banach space, but it would be cumbersome at best to try to appeal to his result here. Fortunately, Taubes provides a similar result in the context of the Yang-Mills energy functional
\cite[Proposition 4.2]{TauPath}. However, the desired result can also be extracted from the results of Kozono, Maeda, and Naito \cite{Kozono_Maeda_Naito_1995} or Schlatter \cite{Schlatter_1997} and Struwe \cite{Struwe_1994} on Yang-Mills gradient flow, as we indicate below in Proposition \ref{prop:Convergent_Palais-Smale_sequence_Yang-Mills_functional}.

\begin{prop}[Existence of a convergent Palais-Smale sequence with non-increasing energies for the Yang-Mills functional]
\label{prop:Convergent_Palais-Smale_sequence_Yang-Mills_functional}
Let $G$ be a compact Lie group and $P$ a principal $G$-bundle over a closed, connected, smooth four-dimensional manifold, $X$, with Riemannian metric, $g$. Then there exists a finite subset, $\Sigma \subset X$, a sequence, $\{A_m\}_{m\in\NN}$, of connections  of class $W^{k,p}$, with $k\geq 1$ and $p\geq 2$ obeying $kp>4$
on a finite sequence of principal $G$-bundles, $P_0=P,P_1,\ldots,P_K$, with obstructions $\eta(P_k)=\eta(P)$ for $k=0,\ldots,K$, where $A_m$ is defined on $P_k$ for $m_k\leq m < m_{k+1}$, with $k=0,\ldots,K$, and $m_{K+1}=\infty$,
that obeys \eqref{eq:Convergent_Palais-Smale_sequence}, for some Yang-Mills connection, $A_\infty$, on a principal $G$-bundle, $P_\infty$, over $X$ with obstruction $\eta(P_\infty)=\eta(P)$, and the sequence, $\{\sE(A_m)\}_{m\in\NN}$, is non-increasing.
\end{prop}

\begin{proof}
Given any initial connection, $A_0$, on $P$ of class $W^{k,p}$, we consider the Yang-Mills gradient flow, $A(t)$ for $t \geq 0$, with initial data. Modulo singularities occurring at finitely many times, $\{T_1,\ldots,T_K\} \subset (0,\infty)$, the flow, $A(t)$, may be extended to a global weak solution in the sense of Kozono, Maeda, and Naito
\cite[Theorem A]{Kozono_Maeda_Naito_1995}, Schlatter \cite[Theorems 1.2 and 1.3]{Schlatter_1997}, and Struwe \cite[Theorem 2.4]{Struwe_1994} on a finite sequence of principal $G$-bundles, $P_0=P,P_1,\ldots,P_K$, with $A(t)$ defined on $[T_k,T_{k+1})$ for $k=0,\ldots,K$ and $T_{K+1}=\infty$. In particular,
\cite[Theorems 1.2 and 1.3 and Equations (56) or (58)]{Schlatter_1997} provide a sequence of times, $\{t_m\}_{m\in\NN} \subset [0,\infty)$, such that $\sE(A(t_m)) \to \sE(A_\infty)$ and  $\sE'(A(t_m)) = d_{A(t_m)}^*F_{A(t_m)}\to 0$ as $m \to \infty$ and $\{\sE(A(t_m))\}_{m\in\NN}$ is non-increasing, which completes the proof by taking $A_m := A(t_m)$ for $m \in \NN$.
\end{proof}

The sequence, $\{A_m\}_{m\in\NN}$, asserted by Schlatter's
\cite[Theorem 1.3]{Schlatter_1997} only converges weakly in
$$
H_{A_\infty,\loc}^1(X\less\Sigma;\Lambda^1\otimes\ad P_\infty),
$$
modulo local gauge transformations to $(A_\infty,\Sigma)$, essentially as in Sedlacek's Theorem \ref{thm:Sedlacek_3-1}. We emphasize that the subset, $\Sigma = \{x_1,\ldots,x_L\} \subset X$ in Proposition \ref{prop:Convergent_Palais-Smale_sequence_Yang-Mills_functional}, is the union of the sets of bubble points, $\Sigma_1,\ldots,\Sigma_K,\Sigma_\infty$, arising in the application of \cite[Theorems 1.2 and 1.3]{Schlatter_1997} in the proof of Proposition \ref{prop:Convergent_Palais-Smale_sequence_Yang-Mills_functional}, and that the connections in the sequence, $\{A_m\}_{m\in\NN}$, while all defined on $P\restriction\Sigma$, are not all defined on $P$.

Following Taubes' \cite[Definition 4.1]{TauPath}, we call a sequence of connections, $\{A_m\}_{m\in\NN}$, on a principal $G$-bundle, $P$, over $X$ \emph{good} if
\begin{equation}
\label{eq:Taubes_definition_4-1}
\sup_{m\in\NN}\sE(A_m) < \infty \quad\hbox{and}\quad \sE'(A_m) \to 0
\quad\hbox{as } m \to \infty,
\end{equation}
Clearly, a convergent Palais-Smale sequence of connections is good. Thanks to Taubes' \cite[Proposition 4.5]{TauPath}, the mode of convergence in Proposition \ref{prop:Convergent_Palais-Smale_sequence_Yang-Mills_functional} may be improved to give

\begin{prop}[Strong convergence in $H_{\loc}^1(X\less\Sigma)$ for a good sequence for the Yang-Mills functional]
\label{prop:Taubes_1984_4-5_for_convergent_Palais-Smale_sequence_Yang-Mills_functional}
Let $G$ be a compact Lie group and $P$ a principal $G$-bundle over a closed, connected, smooth four-dimensional manifold, $X$, with Riemannian metric, $g$, and $\{A_m\}_{m\in\NN}$ a \emph{good} sequence of connections of class $H^1$ on a finite sequence of principal $G$-bundles, $P_0=P,P_1,\ldots,P_K$, with obstructions $\eta(P_k)=\eta(P)$ for $k=0,\ldots,K$, where $A_m$ is defined on $P_k$ for $m_k\leq m < m_{k+1}$, with $k=0,\ldots,K$, and $m_{K+1}=\infty$. After passing to a subsequence of $\{A_m\}_{m\in\NN}$ and relabelling, there is a sequence, $\{u_m\}_{m\in\NN}$, of bundle isomorphisms, $u_m: P_\infty\restriction X\less\Sigma \cong P\restriction X\less\Sigma$ of class $H_{\loc}^3(X\less\Sigma)$ such that
$$
\|u_m^*A_m - A_\infty\|_{H_{A_\infty}^1(U)} \to 0
\quad\hbox{strongly in } H_{A_\infty}^1(U;\Lambda^1\otimes\ad P_\infty)
\quad\hbox{as } m \to \infty,
$$
for every $U \Subset X\less \Sigma$, for some Yang-Mills connection, $A_\infty$, on a principal $G$-bundle, $P_\infty$, over $X$ with obstruction $\eta(P_\infty)=\eta(P)$.
\end{prop}

\begin{proof}
In \cite[Proposition 4.5]{TauPath}, it is assumed that the sequence of connections, $\{A_m\}_{m\in\NN}$, are all defined on the same principal $G$-bundle, $P$, over $X$. However, the proof of \cite[Proposition 4.5]{TauPath} extends \mutatis to the case where $\{A_m\}_{m\in\NN}$ is instead defined on a finite sequence principal $G$-bundles, $P_0=P,P_1,\ldots,P_K$.
\end{proof}

In particular, by the Sobolev embedding $H^1(X) \hookrightarrow L^4(X)$ \cite[Theorem 4.12]{AdamsFournier} and Kato Inequality \eqref{eq:FU_6-20_first-order_Kato_inequality}, Proposition \ref{prop:Taubes_1984_4-5_for_convergent_Palais-Smale_sequence_Yang-Mills_functional} implies that
$$
\|u_m^*A_m - A_\infty\|_{L^4(U)} \to 0
\quad\hbox{strongly in } L^4(U;\Lambda^1\otimes\ad P)
\quad\hbox{as } m \to \infty,
$$
for every $U \Subset X\less \Sigma$. We then have the

\begin{cor}[Convergence of the least eigenvalue of $d_{A_m}^+d_{A_m}^{+,*}$ for a convergent Palais-Smale sequence of connections, $\{A_m\}_{m\in\NN}$, converging strongly in $H_{\loc}^1(X\less\Sigma)$]
\label{cor:Strong_H_loc1_X_less_Sigma_continuity_least_eigenvalue_wrt_connection}
Let $G$ be a compact Lie group and $P$ a principal $G$-bundle over a closed, connected, smooth four-dimensional manifold, $X$, with Riemannian metric, $g$, and $\{A_m\}_{m\in\NN}$ a convergent Palais-Smale sequence of connections of class $H^1$ on $P\restriction\Sigma$ that converges strongly in $H_{\loc}^1(X\less\Sigma)$, modulo a sequence, $\{u_m\}_{m\in\NN}$, of bundle isomorphisms, $u_m: P_\infty\restriction X\less\Sigma \cong P\restriction X\less\Sigma$ of class $H_{\loc}^3(X\less\Sigma)$, to a Yang-Mills connection, $A_\infty$, on a principal $G$-bundle, $P_\infty$, over $X$. Then \eqref{eq:Weak_H_loc1_X_less_Sigma_continuity_least_eigenvalue_wrt_connection_bounds} holds, that is,
$$
\lim_{m\to\infty}\mu(A_m) = \mu(A_\infty),
$$
where $\mu(A)$ is as in \eqref{eq:Least_eigenvalue_dA+dA+*}.
\end{cor}

\subsection{The generic metric theorems of Freed and Uhlenbeck and uniform positive lower bounds for the least eigenvalue of $d_A^+d_A^{+,*}$}
\label{subsec:Generic_metric_theorems_Freed_Uhlenbeck}
Theorem \ref{thm:Proposition_Feehan_Leness_7-6} makes no assumption on the topology of $P$ or $X$ or the Riemannian metric, $g$. However, application of Theorem \ref{thm:Proposition_Feehan_Leness_7-6} requires a positive lower bound, $\mu_0$, for the smallest eigenvalue, $\mu(A)$, of $d_A^+d_A^{+,*}$ on $L^2(X; \Lambda^+\otimes\ad P)$, while $\eps \in (0,1]$ can be simultaneously assumed to be as small as desired in order to appeal to the Contraction Mapping Principle. The analogous remarks apply to Theorem \ref{thm:Taubes_1982_3-2}.

By way of a first step in this direction, we shall describe two approaches to establishing a uniform positive lower bound for $\mu(A)$ as $[A]$ varies over the moduli space of $g$-anti-self-dual connections on $P$,
\begin{equation}
\label{eq:Moduli_space_anti-self-dual_connections}
M(P,g) := \{[A] \in \sB(P,g): F_A^{+,g} = 0\},
\end{equation}
where $\sB(P,g) := \sA(P)/\Aut(P)$ and $\sA(P)$ is the affine (Banach) space of $W^{k,p}$ connections on $P$ and $\Aut(P)$ the Banach Lie group of gauge transformations of $P$ of class $W^{k+1,p}$, where $k\geq 1$ and $p\geq 2$ and $(k+1)p > 4$. We refer the reader to \cite{DK, FU} for discussions of the Banach manifold structures on $\sB^*(P,g)$ (the open subset of irreducible connections on $P$) and $\Aut(P)$.

If the Riemannian metric, $g$, on $X$ is positive in the sense of \eqref{eq:Freed_Uhlenbeck_page_174_positive_metric}, then the Bochner-Weitenb\"ock formula \eqref{eq:Freed_Uhlenbeck_6-26} ensures that the least eigenvalue function,
\begin{equation}
\label{eq:Least_eigenvalue_dA+dA+*_function_on_MPg}
\mu_g[\,\cdot\,]: M(P,g) \to [0, \infty),
\end{equation}
defined by $\mu_g(A)=\mu(A)$ in \eqref{eq:Least_eigenvalue_dA+dA+*}, admits a uniform positive lower bound, $\mu_0 = \mu(g)$,
$$
\mu_g(A) \geq \mu_0,  \quad\forall\, [A] \in M(P,g),
$$
as the following well-known elementary lemma illustrates and which underlies Taubes' proof of his \cite[Theorem 1.4]{TauSelfDual}.

\begin{lem}[Positive lower bound for the least eigenvalue of $d_A^+d_A^{+,*}$ on a four-manifold with a positive Riemannian metric and $L^2$-small $F_A^+$]
\label{lem:Positive_metric_implies_positive_lower_bound_small_eigenvalues}
Let $X$ be a closed, four-dimensional, oriented, smooth manifold with Riemannian metric, $g$, that is \emph{positive} in the sense of \eqref{eq:Freed_Uhlenbeck_page_174_positive_metric}. Then there is a positive constant, $\eps = \eps(g) \in (0,1]$, with the following significance. Let $G$ be a compact Lie group and $P$ a principal $G$-bundle over $X$. If $A$ is an $H^1$ connection on $P$ such that
$$
\|F_A^+\|_{L^2(X)} \leq \eps,
$$
and $\mu(A)$ is as in \eqref{eq:Least_eigenvalue_dA+dA+*}, then
\begin{equation}
\label{eq:Positive_metric_implies_positive_lower_bound_small_eigenvalues_L2small_FA+}
\mu(A) \geq \inf_{x\in X}\left(\frac{1}{3}R(x) - 2w_+(x)\right) > 0.
\end{equation}
\end{lem}

\begin{proof}
The result is well-known, but we include the proof for the sake of completeness. Let $v \in H_A^1(X;\Lambda^+\otimes\ad P)$ be an eigenvector of $d_A^+d_A^{+,*}$ for the eigenvalue $\mu(A)$ with $\|v\|_{L^2(X)}=1$. By applying the Bochner-Weitenb\"ock formula and integration by parts, we see that
\begin{align*}
\mu(A) &= (d_A^+d_A^{+,*}v,v)_{L^2(X)}
\\
&= (\nabla_A^*\nabla_Av,v)_{L^2(X)}
+ \left(\left(\frac{1}{3}R - 2w_+\right)v,v\right)_{L^2(X)}
+ (\{F_A^{+,g}, v\}, v)_{L^2(X)} \quad\hbox{(by \eqref{eq:Freed_Uhlenbeck_6-26})}
\\
&\geq \inf_{x\in X}\left(\frac{1}{3}R(x) - 2w_+(x)\right)\|v\|_{L^2(X)}^2
- \|F_A^{+,g}\|_{L^2(X)}\|v\|_{L^4(X)}^2
\\
&\geq \inf_{x\in X}\left(\frac{1}{3}R(x) - 2w_+(x)\right)
- c\eps\left(\|d_A^{+,*}v\|_{L^2(X)}+\|v\|_{L^2(X)}\right)
\quad\hbox{(by hypothesis and \eqref{eq:Feehan_Leness_6-6-1_L4v_L2dA+*v})}
\\
&= \inf_{x\in X}\left(\frac{1}{3}R(x) - 2w_+(x)\right)
- c\eps\left(\sqrt{\mu(A)} + 1\right) \quad\hbox{(by \eqref{eq:Least_eigenvalue_dA+dA+*})}.
\end{align*}
Using $c\eps\sqrt{\mu(A)} \leq \frac{1}{2}(c^2\eps^2 + \mu(A))$, we obtain
$$
\frac{1}{2}\mu(A) \geq \inf_{x\in X}\left(\frac{1}{3}R(x) - 2w_+(x)\right)
- c\eps - \frac{c^2\eps^2}{2},
$$
Choosing $\eps \in (0,1]$ small enough that
$$
c\eps + \frac{c^2\eps^2}{2} \leq \frac{1}{2}\inf_{x\in X}\left(\frac{1}{3}R(x) - 2w_+(x)\right),
$$
yields the desired result.
\end{proof}

While Lemma \ref{lem:Positive_metric_implies_positive_lower_bound_small_eigenvalues} has the benefit that it only requires $F_A^{+,g}$ to be $L^2$-small, not zero, and so has direct application to the existence results, Theorem \ref{thm:Proposition_Feehan_Leness_7-6} and \ref{thm:Taubes_1982_3-2}, the positivity hypothesis on $g$ imposes a strong constraint on the topology of $X$ since, when applied to the product connection on $X\times G$ and Levi-Civita connection on $TX$, it implies that $b^+(X) = \Ker d^+d^{+,*} = 0$ and thus $X$ is necessarily a four-dimensional manifold with negative definite intersection form, $Q_X$, on $H_2(X;\ZZ)$. Indeed, we recall from \cite[Section 1.1.6]{DK} that, given \emph{any} Riemannian metric $g$ on $X$, we have an isomorphism of real vector spaces,
$$
H^2(X;\RR) \cong \sH^{+,g}(X) \oplus \sH^{-,g}(X;\RR)
$$
where $\sH^{\pm,g}(X) := \ker\{d^{+,g}d^{+,g,*}:\Omega^2(X) \to C^\infty(X)\}$, the real vector spaces of harmonic self-dual and anti-self-dual two-forms defined by the metric, $g$, and
$$
b^\pm = \dim\sH^{\pm,g}(X).
$$
Another approach to ensuring the existence of a uniform positive lower bound for the least eigenvalue function \eqref{eq:Least_eigenvalue_dA+dA+*_function_on_MPg} is more subtle and relies on the generic metric theorems of Freed and Uhlenbeck \cite[pp. 69--73]{FU}, together with certain extensions due to Donaldson and Kronheimer \cite[Sections 4.3.3]{DK}. Under suitable hypotheses on $P$ and a generic Riemannian metric, $g$, on $X$, their results collectively ensure that $\mu(A) > 0$ for all $[A]$ in both $M(P,g)$ and every moduli space, $M(P_l,g)$, appearing in its \emph{Uhlenbeck compactification} (see \cite[Theorem 4.4.3]{DK}),
\begin{equation}
\label{eq:Uhlenbeck_compactification}
\bar M(P,g) \subset \bigcup_{l=0}^L \left(M(P_l,g)\times\Sym^l(X)\right).
\end{equation}
While the statement of \cite[Theorem 4.4.3]{DK} assumes that $G=\SU(2)$ or $\SO(3)$ --- see \cite[pages 157 and 158]{DK} --- the proof applies to any compact Lie group via the underlying analytic results due Uhlenbeck \cite{UhlRem, UhlLp}; alternatively, one may appeal directly to the general compactness result due to Taubes \cite[Proposition 4.4]{TauPath}, \cite[Proposition 5.1]{TauFrame}. Every principal $G$-bundle, $P_l$, over $X$ appearing in \eqref{eq:Uhlenbeck_compactification} has the property that $\eta(P_l) = \eta(P)$ by \cite[Theorem 5.5]{Sedlacek}.

The generic metric theorems of Freed and Uhlenbeck \cite[pp. 69--73]{FU} and Donaldson and Kronheimer \cite[Sections 4.3.3]{DK} are normally phrased in terms of existence of a Riemannian metric, $g$, on $X$ such that $\Coker d_A^{+,g}=0$ for all $[A] \in M(P,g)$, a property of $g$ which is equivalent to $\mu_g(A) > 0$ for all $[A] \in M(P,g)$, as we use in the following restatement of their results.

\begin{thm}[Generic metrics theorem for simply-connected four-manifolds]
\label{thm:Donaldson-Kronheimer_Corollary_4-3-15_and_18_and_Proposition_4-3-20}
Let $G$ be a compact Lie group and $P$ a principal $G$-bundle over a closed, connected, Riemannian, smooth, four-dimensional manifold, $X$. Then there is an open dense subset, $\sC(P)$, of the Banach space, $\sC(X)$, of conformal equivalence classes, $[g]$, of $C^r$ Riemannian metrics on $X$ (for some $r\geq 3$) with the following significance. Assume that $[g] \in \sC(P)$ and $\pi_1(X)$ is trivial and at least \emph{one} of the following holds:
\begin{enumerate}
  \item $b^+(X) = 0$ and $G = \SU(2)$ or $G = \SO(3)$; or

  \item $b^+(X) > 0$, $G = \SO(3)$, and the second Stiefel-Whitney class, $w_2(P) \in H^2(X;\ZZ/2\ZZ)$, is non-trivial;
\end{enumerate}
Then every point $[A] \in M(P,g)$ has the property that $\mu_g(A) > 0$.
\end{thm}

\begin{proof}
Let us first consider points, $[A]$, in the open subset, $M^*(P,g) \subset M(P,g)$, of gauge-equivalence classes of irreducible $g$-anti-self-dual connections. For $b^+(X) \geq 0$, every irreducible $g$-anti-self-dual connection, $A$, on a principal $\SU(2)$ or $\SO(3)$-bundle over $X$ has $\Coker d_A^{+,g} = 0$ by \cite[Corollary 4.3.18]{DK} and thus $\mu(A)>0$ for all $[A] \in M^*(P,g)$ when $G = \SU(2)$ or $G = \SO(3)$.

It remains to consider points $[A] \in M^{\red}(P,g) := M(P,g) \less M^*(P,g)$. If $b^+(X) = 0$ and $[A] \in M^{\red}(P,g)$, for $G = \SU(2)$ or $\SO(3)$, then $\Coker d_A^{+,g} = 0$ and thus $\mu(A)>0$ by \cite[Proposition 4.3.20]{DK}. If $b^+(X) > 0$, the only reducible $g$-anti-self-dual connection on a principal $\SU(2)$ or $\SO(3)$-bundle over $X$ is the product connection on the product bundle $P = X\times G$ by \cite[Corollary 4.3.15]{DK} and the latter possibility is excluded by our hypothesis in this case that $G = \SO(3)$ and the second Stiefel-Whitney class, $w_2(P) \in H^2(X;\ZZ/2\ZZ)$, is non-trivial.
\end{proof}

The hypothesis in Theorem \ref{thm:Donaldson-Kronheimer_Corollary_4-3-15_and_18_and_Proposition_4-3-20} that the four-manifold $X$ is simply-connected may be relaxed.

\begin{cor}[Generic metrics theorem for four-manifolds with no non-trivial representations of $\pi_1(X)$ in $G$]
\label{cor:Donaldson-Kronheimer_Corollary_4-3-15_and_18_and_Proposition_4-3-20}
Assume the hypotheses of Theorem \ref{thm:Donaldson-Kronheimer_Corollary_4-3-15_and_18_and_Proposition_4-3-20}, except replace the hypothesis that $X$ is simply-connected by the requirement that the fundamental group, $\pi_1(X)$, has no non-trivial representations in $G$. Then the conclusions of Theorem \ref{thm:Donaldson-Kronheimer_Corollary_4-3-15_and_18_and_Proposition_4-3-20} continue to hold.
\end{cor}

\begin{proof}
If $X$ has the property that there are no non-trivial representations, $\pi_1(X) \to G$, and a principal $G$-bundle $P$ admits a flat connection, $A$, then $P$ is isomorphic to the product bundle, $X \times G$, and that isomorphism identifies $A$ with the product connection by \cite[Proposition 2.2.3]{DK}. The conclusions now follow from the proofs of Theorem \ref{thm:Donaldson-Kronheimer_Corollary_4-3-15_and_18_and_Proposition_4-3-20} and the underlying results \cite[Corollaries 4.3.15 and 4.3.18 and Proposition 4.3.20]{DK} developed by Donaldson and Kronheimer.
\end{proof}

Our results in Section \ref{subsec:Continuity_least_eigenvalue_d_A+d_A+*_wrt_connection} assure the continuity of $\mu_g[\,\cdot\,]$ with respect to the Uhlenbeck topology and they may be profitably applied here. Before doing this, since there are many potential combinations of conditions on $G$, $P$, $X$, and $g$ which imply that $\mu_g(A) > 0$, when $A$ is $g$-anti-self-dual, it is convenient to introduce the

\begin{defn}[Good Riemannian metric]
\label{defn:Good_Riemannian_metric}
Let $G$ be a compact Lie group and $P$ a principal $G$-bundle over a closed, four-dimensional, oriented, smooth manifold, $X$. We say that a Riemannian metric, $g$, on $X$ is \emph{good} if $\mu_g(A_0)>0$ whenever $A_0$ is a connection with $F_{A_0}^{+,g}=0$ on a principal $G$-bundle $P_0$ over $X$ with $\eta(P_0)=\eta(P)$, where $\mu_g(A_0)$ is as in \eqref{eq:Least_eigenvalue_dA+dA+*} and $\eta(P) \in H^2(X; \pi_1(G))$ is the obstruction \cite[Section 2]{Sedlacek}.
\end{defn}

We then have the

\begin{thm}[Positive lower bound for the least eigenvalue of $d_A^+d_A^{+,*}$ on a four-manifold with a good Riemannian metric and anti-self-dual connection $A$]
\label{thm:Good_metric_implies_positive_lower_bound_small_eigenvalues}
Let $G$ be a compact Lie group and $P$ a principal $G$-bundle over a closed, four-dimensional, oriented, smooth manifold, $X$, with Riemannian metric, $g$. Assume that $g$ is \emph{good} in the sense of Definition \ref{defn:Good_Riemannian_metric}. Then there is a positive constant, $\mu_0 = \mu_0(P,g)$ with the following significance. If $A$ is an $H^1$ connection on $P$ such that
$$
F_A^{+,g} = 0,
$$
and $\mu_g(A)$ is as in \eqref{eq:Least_eigenvalue_dA+dA+*}, then
\begin{equation}
\label{eq:Positive_metric_implies_positive_lower_bound_small_eigenvalues_good_Riem_metric_ASD_conn}
\mu_g(A) \geq \mu_0.
\end{equation}
\end{thm}

\begin{proof}
The conclusion is a consequence of the fact that $\bar M(P,g)$ is compact, the extension,
\begin{equation}
\label{eq:Least_eigenvalue_function_on_Uhlenbeck_compactification}
\bar\mu_g[\,\cdot\,]: \bar M(P,g) \ni ([A],\bx) \to \mu_g[A] \in [0, \infty),
\end{equation}
to $\bar M(P,g)$ of the function \eqref{eq:Least_eigenvalue_dA+dA+*_function_on_MPg} defined by \eqref{eq:Least_eigenvalue_dA+dA+*} is continuous by Proposition \ref{prop:Lp_loc_continuity_least_eigenvalue_wrt_connection}, the fact that $\mu_g(A) > 0$ for $[A] \in M(P_l,g)$ and $P_l$ is a principal $G$-bundle over $X$ appearing in the Uhlenbeck compactification \eqref{eq:Uhlenbeck_compactification}, and $g$ is good by hypothesis.
\end{proof}

We remark that, unlike Lemma \ref{lem:Positive_metric_implies_positive_lower_bound_small_eigenvalues}, Theorem \ref{thm:Good_metric_implies_positive_lower_bound_small_eigenvalues} requires that $F_A^+=0$, not merely that $\|F_A^+\|_{L^2(X)} \leq \eps$, for a small enough $\eps(g) \in (0,1]$. However, by appealing to the results of Sedlacek, we can obtain an extension of Theorem \ref{thm:Good_metric_implies_positive_lower_bound_small_eigenvalues} that replaces the condition $F_A^+=0$ by $\|F_A^+\|_{L^2(X)} \leq \eps$, for a small enough $\eps(P,g) \in (0,1]$.

Consider the open subset of the space $\sB(P,g)$ defined by
\begin{equation}
\label{eq:Open_neighborhood_asd_moduli_space_FA+_L2_small}
\sB_\eps(P,g) := \{[A] \in \sB(P,g): \fA_g(A) < \eps\},
\end{equation}
where (compare \cite[Equation (2.2) and p. 349]{TauPath}),
\begin{equation}
\label{eq:Taubes_1984b_2-2_and_page_349}
\fA_g(A) := \frac{1}{2}\int_X |F_A^{+,g}|^2\,d\vol_g.
\end{equation}
We then have the

\begin{thm}[Positive lower bound for the least eigenvalue of $d_A^+d_A^{+,*}$ on a four-manifold with a good Riemannian metric and almost anti-self-dual connection $A$]
\label{thm:Good_metric_implies_positive_lower_bound_small_eigenvalues_almost_ASD}
Let $G$ be a compact Lie group and $P$ a principal $G$-bundle over a closed, four-dimensional, oriented, smooth manifold, $X$, with Riemannian metric, $g$. Assume that $g$ is \emph{good} in the sense of Definition \ref{defn:Good_Riemannian_metric}. Then there is a positive constant, $\eps = \eps(P,g) \in (0,1]$, with the following significance. If $A$ is an $H^1$ connection on $P$ such that
$$
\|F_A^{+,g}\|_{L^2(X)} \leq \eps,
$$
and $\mu_g(A)$ is as in \eqref{eq:Least_eigenvalue_dA+dA+*}, then
\begin{equation}
\label{eq:Positive_metric_implies_positive_lower_bound_small_eigenvalues_good_Riem_metric_almost_ASD_conn}
\mu_g(A) \geq \frac{\mu_0}{2},
\end{equation}
where $\mu_0 = \mu_0(P,g)$ is the positive constant in Theorem \ref{thm:Good_metric_implies_positive_lower_bound_small_eigenvalues}.
\end{thm}

\begin{proof}
Suppose that the constant, $\eps = \eps(P,g) \in (0,1]$, does not exist. We may then choose a minimizing sequence, $\{A_m\}_{m\in\NN}$, of connections on $P$ such that $\|F_{A_m}^+\|_{L^2(X)} \to 0$ and $\mu(A_m) \to 0$ as $m\to\infty$. According to Sedlacek's Theorem \ref{thm:Sedlacek_4-3}, there is a $g$-anti-self-dual connection, $A_0$, on a principal $G$-bundle, $P_0$, over $X$ such that $\eta(P)=\eta(P_0)$ and a finite set of points, $\Sigma = \{x_1,\ldots,x_L\} \subset X$, such that $A_m$ converges to $A_0$ modulo local gauge transformations weakly in $H_{\loc}^1(X\less\Sigma)$ in the sense of Theorem \ref{thm:Sedlacek_3-1}. But then $\mu(A_m) \to \mu(A_0) \geq \mu_0$ by Corollary \ref{cor:Weak_H_loc1_X_less_Sigma_continuity_least_eigenvalue_wrt_connection} and $\mu(A_0)>0$ because $g$ is good by hypothesis, contradicting our initial assumption regarding the sequence $\{A_m\}_{m\in\NN}$. In particular, the preceding argument shows that the desired constant, $\eps$, exists.
\end{proof}

\begin{cor}[Uniform positive lower bound for the smallest eigenvalue function when $g$ is generic, $G$ is $\SU(2)$ or $\SO(3)$, and $\pi_1(X)$ has no non-trivial representations in $G$]
\label{cor:Least_eigenvalue_function_positive_on_Uhlenbeck_compactification_pi1_X}
Assume the hypotheses of Corollary \ref{cor:Donaldson-Kronheimer_Corollary_4-3-15_and_18_and_Proposition_4-3-20} and that $g$ is generic. Then there are constants, $\eps = \eps(P,g) \in (0,1]$ and $\mu_0 = \mu_0(P,g) > 0$, such that
\begin{align*}
\mu_g(A) &\geq \mu_0, \quad\forall\, [A] \in M(P,g),
\\
\mu_g(A) &\geq \frac{\mu_0}{2}, \quad\forall\, [A] \in \sB_\eps(P,g).
\end{align*}
\end{cor}

\begin{proof}
It suffices to observe that $g \in \cap_{l=1}^L \sC(P_l)$, where $\sC(P_l)$ is as in Theorem \ref{thm:Donaldson-Kronheimer_Corollary_4-3-15_and_18_and_Proposition_4-3-20}, and thus $g$ is good in the sense of Definition \ref{defn:Good_Riemannian_metric}. The conclusions now follow from Theorem \ref{thm:Good_metric_implies_positive_lower_bound_small_eigenvalues} and \ref{thm:Good_metric_implies_positive_lower_bound_small_eigenvalues_almost_ASD}.
\end{proof}

The proof of Corollary \ref{cor:Least_eigenvalue_function_positive_on_Uhlenbeck_compactification_pi1_X} extends without change to give

\begin{cor}[Uniform positive lower bound for the smallest eigenvalue function when $g$ is generic and $G$ is $\SO(3)$]
\label{cor:Least_eigenvalue_function_positive_on_Uhlenbeck_compactification}
Assume the hypotheses of Corollary \ref{cor:Least_eigenvalue_function_positive_on_Uhlenbeck_compactification_pi1_X}, but replace the hypothesis on $\pi_1(X)$ by the requirement that $G = \SO(3)$ and no principal $\SO(3)$-bundle $P_l$ over $X$ appearing in the Uhlenbeck compactification \eqref{eq:Uhlenbeck_compactification} admits a flat connection. Then the conclusions of Corollary \ref{cor:Least_eigenvalue_function_positive_on_Uhlenbeck_compactification_pi1_X} continue to hold.
\end{cor}

Despite its apparently technical nature, the alternative hypothesis in Corollary \ref{cor:Least_eigenvalue_function_positive_on_Uhlenbeck_compactification} is easy to achieve in practice, albeit at the cost of blowing-up the given four-manifold, $X$, and modifying the given principal $G$-bundle, $P$. Indeed, we recall the following facts discussed in \cite{MorganMrowkaPoly}:
\begin{enumerate}
\item If $H_1(X;\ZZ)$ has no $2$-torsion, then every principal $\SO(3)$-bundle $\bar P$ over $X$ lifts to a principal $\U(2)$-bundle $P$, every $\SO(3)$-connection on $\bar P$ lifts to a $\U(2)$-connection on $P$, and every $\SO(3)$-gauge transformation lifts to an $\SU(2)$-gauge transformation \cite[Remark (ii), p. 225]{MorganMrowkaPoly};

\item If $\widetilde X$ is the connected sum, $X \# \CC\PP^2$, and $\widetilde P$ is the connected sum of a principal $\U(2)$-bundle $P$ over $X$ with the principal $\U(2)$-bundle $Q$ over $\CC\PP^2$ with $c_2(Q)=0$ and $c_1(Q)$ equal to the Poincar\'e dual of $e = [\CC\PP^1] \in H^2(\CC\PP^2; \ZZ)$, and $\bar P$ is the principal $\SO(3)$-bundle associated to $\widetilde P$, then the following holds: No principal $\SO(3)$-bundle $P'$ over $X$ with $w_2(P') = w_2(\bar P)$ admits a flat connection
    \cite[Paragraph prior to Corollary 2.2]{MorganMrowkaPoly}.
\end{enumerate}

An enhancement due to Kronheimer and Mrowka \cite[Corollary 2.5]{KMStructure} of the generic metrics theorem of Freed and Uhlenbeck \cite[pp. 69--73]{FU} and Donaldson and Kronheimer \cite[Sections 4.3.3, 4.3.4, and 4.3.5]{DK}, for a generic Riemannian metric, $g$, on $X$ ensures the following holds even when $\pi_1(X)$ is non-trivial: If $b^+(X) > 0$, then\footnote{This result may also be true when $b^+(X) = 0$, but a verification would require a close examination of their proof, as this possibility is not explicitly allowed by their statement of \cite[Corollary 2.5]{KMStructure}.}
every $g$-anti-self-dual connection, $A$, on a principal $\SO(3)$ bundle over $X$ has $\Coker d_A^{+,g} = 0$, unless $A$ is reducible or flat.



\subsection{Tubular neighborhood of the moduli space of anti-self-dual connections}
\label{subsec:Tubular_neighborhood_moduli_space_asd_connections}
Suppose that $G$ is a compact Lie group and $P$ is a principal $G$-bundle over a closed, four-dimensional, oriented, smooth manifold, $X$, with Riemannian metric, $g$, such that $M(P,g)$ is non-empty. In \cite[Proposition 3.1]{TauPath}, Taubes applies his construction of anti-self-dual curvature flow to show that $M(P,g)$ is a strong deformation retract of the open subset $\sB_\eps(P,g)$ in \eqref{eq:Open_neighborhood_asd_moduli_space_FA+_L2_small} of the space $\sB(P,g)$, for a small enough $\eps \in (0,1]$, provided there is a positive constant, $\mu_0$, such that $\mu(A) \geq \mu_0$ for all $[A] \in \sB_\eps(P,g)$.

Taubes develops a more general version of his \cite[Proposition 3.1]{TauPath} in \cite[Lemma 5.4]{TauStable}, where the constraint $\mu(A) \geq \mu_0$ is omitted at the expense of replacing the constant $\eps$ in \eqref{eq:Open_neighborhood_asd_moduli_space_FA+_L2_small} by $z(\mu(A))$, where $z:[0,\infty)\to [0,\infty)$ is a function such that $z(0)=0$. However, it is worth noting that the generic metrics result (due to Freed and Uhlenbeck \cite{FU}) he appeals to in \cite[page 194]{TauStable} only applies when $G$ is $\SU(2)$ or $\SO(3)$ and not the arbitrary compact Lie group, $G$, he allows in \cite[Lemma 5.4]{TauStable}.

By slightly strengthening the norm employed in the definition \eqref{eq:Open_neighborhood_asd_moduli_space_FA+_L2_small} to give
\begin{equation}
\label{eq:Open_neighborhood_asd_moduli_space_FA+_Lsharp2_small}
\sB_{\sharp;\eps}(P,g)
:= \{[A] \in \sB(P,g): \|F_A^{+,g}\|_{L^{\sharp,2}(X)} < 2\eps^2\},
\end{equation}
we obtain

\begin{prop}[Existence of a tubular neighborhood of the moduli space of anti-self-dual connections]
\label{prop:Existence_tubular_neighborhood_moduli_space_asd_connections}
Let $G$ be a compact Lie group and $P$ a principal $G$-bundle over a closed, four-dimensional, oriented, smooth manifold, $X$, with Riemannian metric, $g$. Then there are a positive constant $c = c(P,g)$ and a constant $\eps = \eps(P,g) \in (0, 1]$ with the following significance. If $g$ is \emph{good} in the sense of Theorem \ref{thm:Good_metric_implies_positive_lower_bound_small_eigenvalues} and $M(P,g)$ is non-empty, then $\sB_{\sharp;\eps}(P,g)$ is a \emph{tubular neighborhood} of $M(P,g)$ in the sense that there is a continuous map, smooth upon restriction to $\sB_{\sharp;\eps}^*(P,g)$,
\begin{equation}
\label{eq:Tubular_neighborhood_moduli_space_asd_connections_projection}
\pi:\sB_{\sharp;\eps}(P,g) \to M(P,g), \quad [A] \mapsto [A+a],
\end{equation}
such that $\pi$ is the identity map on $M(P,g)$ and
\begin{equation}
\label{eq:Tubular_neighborhood_moduli_space_asd_connections_projection_bound}
\|a\|_{H_A^1(X)} \leq c\|F_A^{+,g}\|_{L^{\sharp,2}(X)}.
\end{equation}
\end{prop}

\begin{proof}
Theorem \ref{thm:Good_metric_implies_positive_lower_bound_small_eigenvalues_almost_ASD} provides a constant, $\eps = \eps(P,g) \in (0,1]$, such that if $[A] \in \sB_{\sharp;\eps}(P,g)$, then $\mu(A) \geq \mu_0/2$, where $\mu_0 = \mu(P,g)$ is the positive constant produced by Theorem \ref{thm:Good_metric_implies_positive_lower_bound_small_eigenvalues}.

We may further suppose that $\eps \in (0,1]$ is chosen small enough to also satisfy the hypotheses of Theorem \ref{thm:Proposition_Feehan_Leness_7-6}. Consequently, if $[A] \in \sB_{\sharp;\eps}(P,g)$, there exists a unique $a \in H_A^1(X;\Lambda^1\otimes\ad P)$ such that $F_{A+a}^{+,g} = 0$ and \eqref{eq:Tubular_neighborhood_moduli_space_asd_connections_projection_bound} holds,
so the map \eqref{eq:Tubular_neighborhood_moduli_space_asd_connections_projection}
is well-defined. The methods of \cite{FLKM1} may be applied to show that the map $\pi$ in \eqref{eq:Tubular_neighborhood_moduli_space_asd_connections_projection} is continuous and smooth upon restriction to $\sB_{\sharp;\eps}^*(P,g)$. Clearly, $\pi$ restricts to the identity map on $M(P,g)$.
\end{proof}

It is natural to ask what conclusions of Proposition \ref{prop:Existence_tubular_neighborhood_moduli_space_asd_connections} can be salvaged if $\sB_{\sharp;\eps}(P,g)$ is replaced by $\sB_\eps(P,g)$. Of course, one knows from Taubes' \cite[Proposition 3.1]{TauPath} that $M(P,g)$ is a strong deformation retract of $\sB_\eps(P,g)$ but his proof does not provide any analogue of the bound \eqref{eq:Tubular_neighborhood_moduli_space_asd_connections_projection_bound}. To that end, we have the

\begin{prop}[Existence of a connection with $L^2$-small self-dual curvature and a nearby anti-self-dual connection]
\label{prop:Existence_asd_connection_near_A_with_L2-small_FA+}
Let $G$ be a compact Lie group and $P$ a principal $G$-bundle over a closed, four-dimensional, oriented, smooth manifold, $X$, with Riemannian metric, $g$. Then there are a positive constant $c = c(P,g)$ and a constant $\delta = \delta(P,g) \in (0, 1]$ with the following significance. If $g$ is \emph{good} in the sense of Theorem \ref{thm:Good_metric_implies_positive_lower_bound_small_eigenvalues} and $M(P,g)$ is non-empty, there is a connection $A$ on $P$ of class $H^1$ such that $\|F_A^{+,g}\|_{L^2(X)} < \delta$ and a perturbation $a \in H_A^1(X;\Lambda^1\otimes\ad P)$ such that $F_{A+a}^{+,g} = 0$ and
$$
\|a\|_{H_A^1(X)} < c\|F_A^{+,g}\|_{L^{\sharp,2}(X)}.
$$
\end{prop}

\begin{proof}
Let $\{A_m\}_{m\geq 1}$ be a minimizing sequence of connections on $P$ such that $\|F_{A_m}^+\|_{L^2(X)} \leq 1/m$. Since the energies, $\sE(A_m) = \frac{1}{2}\|F_{A_m}\|_{L^2(X)}^2$, are necessarily bounded (see Section \ref{sec:Taubes_1982_Appendix}) then, after passing to a subsequence, we may suppose that the sequence, $\{A_m\}_{m\geq 1}$, converges in the sense of Sedlacek (Theorem \ref{thm:Sedlacek_3-1}), that is weakly in $H_{\loc}^1(X\less\Sigma)$ modulo local gauge transformations, to a limit $(A_0, \Sigma_0)$, for some (possibly empty) subset $\Sigma_0 = \{x_1,\ldots,x_L\} \subset X$ and smooth $g$-anti-self-dual connection, $A_0$, on a principal $G$-bundle, $P_0$, over $X$ with obstruction $\eta(P_0)=\eta(P)$.

Unfortunately, the mode of convergence in Theorem \ref{thm:Sedlacek_3-1} does not allow us to conclude that, given $\eps\in (0,1]$, we have $\|u_m^*A_m-A_0\|_{H_{A_0}^1(U)}<\eps$, for $U\Subset X\less\Sigma$ and all $m \geq m_0(U,\eps)$, with $m_0(U,\eps)$ sufficiently large, for some sequence of $G$-bundle isomorphisms, $u_m:P_0\restriction X\less\Sigma \cong P\restriction X\less \Sigma$, of class $H_{\loc}^3(X\less\Sigma)$.

Therefore, we apply Proposition \ref{prop:Convergent_Palais-Smale_sequence_Yang-Mills_functional} to produce sequences, $\{A_{mn}\}_{n\geq 1}$ on $P\restriction\Sigma_m$ for a finite subset $\Sigma_m\subset X$, for each $m \geq 1$, that are convergent Palais-Smale in the sense of \eqref{eq:Convergent_Palais-Smale_sequence} and hence, by passing to the diagonal subsequence and relabelling, a sequence, $\{A_m\}_{n\geq 1}$ on $P\restriction\Sigma$ for a (necessarily) finite subset of points, $\Sigma = \cup_{m\geq 0}\Sigma_m \subset X$, that obeys, as $m \to \infty$,
\begin{gather*}
\|u_m^*A_m-A_\infty\|_{H_{A_\infty}^1(U)} \to 0, \quad\forall\, U\Subset X\less\Sigma,
\\
\|F_{A_m}^{+,g}\|_{L^2(X)} \to 0
\\
\sE(A_{m+1}) \leq \sE(A_m) \quad\forall\, m \geq 1,
\\
\sE(A_m) \searrow \sE(A_\infty),
\\
\sE'(A_m) \to 0,
\end{gather*}
where $u_m:P_\infty\restriction X\less\Sigma \cong P\restriction X\less \Sigma$ is a sequence of $G$-bundle isomorphisms of class $H_{\loc}^3(X\less\Sigma)$ and $A_\infty$ is a $g$-anti-self-dual connection on a principal $G$-bundle, $P_\infty$, over $X$ with obstruction $\eta(P_\infty)=\eta(P)$..

According to \cite[Theorem 4.14]{FeehanGeometry} and the development of the bubble-tree compactification for $M(P,g)$ in \cite[Section 4.2]{FeehanGeometry}, there are
\begin{enumerate}
  \item a tree of centered connections (in the sense of \cite[p. 479]{FeehanGeometry}), $A_{l_1}$, $A_{l_1l_2}$, $A_{l_1l_2l_3}$, \ldots, on a tree of principal $G$-bundles $P_{l_1}$, $P_{l_1l_2}$, $P_{l_1l_2l_3}$, \ldots, over copies of $S^4$ which are anti-self-dual with respect to the standard round metric of radius one on $S^4$;

  \item sequences of scales, $\lambda_{m;l_1}$, $\lambda_{m;l_1l_2}$, $\lambda_{m;l_1l_2l3}$, \ldots, belonging to $(0,1]$ and monotonically decreasing to zero;

  \item sequences of mass centers, $\{x_{m;l_1}\} \subset X$, $\{x_{m;l_1l_2}\} \subset S^4$, $\{x_{m;l_1l_2l3}\} \subset S^4$, \ldots, convergent to points in $X$ or $S^4\less\{s\}$, respectively;

  \item sequences of $g$-orthonormal frames, $\{f_{m;l_1}\}$, for $(TX)_{x_{m;l_1}}$;

  \item sequences of fiber points, $p_{m;l_1} \in P|_{x_{m;l_1}}$, $p_{m;l_1s} \in P_{l_1}|_s$, $p_{m;l_1l_2s} \in P_{l_1l_2}|_s$, \ldots;
\end{enumerate}
such that, by repeatedly applying the splicing construction described \cite[Section 3]{FLKM1} with the preceding set of gluing data, one obtains a sequence of connections, $A_m'$, on $P$ with the properties that, for sufficiently large $m = m(\delta)$,
$$
\|F_{A_m'}^{+,g}\|_{L^{\sharp,2}(X)} < \delta,
\quad
\|A_m - A_m'\|_{H_{A_m'}^1(X)} < \delta,
\quad\hbox{and}\quad
\mu(A_m') \geq \frac{\mu(A_m)}{2} \geq \frac{\mu(A_\infty)}{4} \geq \frac{\mu_0}{4},
$$
where $\mu_0 = \mu_0(P,g)$ is the positive constant produced by Theorem \ref{thm:Good_metric_implies_positive_lower_bound_small_eigenvalues}. A similar result is provided by Taubes in \cite[Propositions 5.3 and 5.4]{TauFrame}. The fact $\mu(A_m) \geq \mu(A_\infty)/2$ follows from Corollary \ref{cor:Weak_H_loc1_X_less_Sigma_continuity_least_eigenvalue_wrt_connection}, while the fact that $\mu(A_m') \geq \mu(A_m)/2$ follows from Proposition \ref{prop:Lp_loc_continuity_least_eigenvalue_wrt_connection}.

We may suppose that $\delta \in (0,1]$ is chosen small enough (and hence $m(\delta)$ large enough) to satisfy the hypotheses of Theorem \ref{thm:Proposition_Feehan_Leness_7-6}. Consequently, there exists a unique $a_m \in H_{A_m'}^1(X;\Lambda^1\otimes\ad P)$ such that $F_{A_m'+a_m}^{+,g} = 0$ and
$$
\|a_m\|_{H_{A_m'}^1(X)}
\leq
c\|F_{A_m'}^{+,g}\|_{L^{\sharp,2}(X)} < c\delta,
$$
where $c = c(P,g)$. Choosing $A = A_m'$ and $a = a_m$ concludes the proof.
\end{proof}

Clearly, Proposition \ref{prop:Existence_asd_connection_near_A_with_L2-small_FA+} is a considerably weaker result than one would ideally like, namely, that given a connection $A$ on $P$ with sufficiently small $\|F_A^{+,g}\|_{L^2(X)}$, there is a perturbation, $a \in H_A^1(X;\Lambda^1\otimes\ad P)$ such that $F_{A+a}^{+,g} = 0$.

\subsection{Strong deformation retraction of the moduli space of anti-self-dual connections}
\label{subsec:Strong_deformation_retraction_moduli_space_asd_connections}
As an alternative to our tubular neighborhood result, Proposition \ref{prop:Existence_tubular_neighborhood_moduli_space_asd_connections}, we recall the following result of Taubes \cite{TauPath}.

\begin{prop}[Existence of a strong deformation retraction of the moduli space of anti-self-dual connections]
\label{prop:Existence_strong_deformation_retraction_moduli_space_asd_connections}
\cite[Proposition 3.1]{TauPath}
Let $G$ be a compact Lie group and $P$ a principal $G$-bundle over a closed, four-dimensional, oriented, smooth manifold, $X$, with Riemannian metric, $g$. Then there are a positive constant $c = c(P,g)$ and a constant $\eps = \eps(P,g) \in (0, 1]$ with the following significance. If $g$ is \emph{good} in the sense of Theorem \ref{thm:Good_metric_implies_positive_lower_bound_small_eigenvalues} and $M(P,g)$ is non-empty, then $\sB_\eps(P,g)$ in \eqref{eq:Open_neighborhood_asd_moduli_space_FA+_L2_small} is a \emph{strong deformation retraction} of $M(P,g)$: There is a continuous map,
\begin{equation}
\label{eq:Strong_deformation_retraction_moduli_space_asd_connections_projection}
H:\sB_\eps(P,g)\times[0,1] \to M(P,g),
\end{equation}
such that $H(\cdot, 0)$ is the identity map on $\sB_\eps(P,g)$ and $H(\cdot,1): \sB_\eps(P,g) \to M(P,g)$ is continuous. Moreover, $H(\cdot,1)$ restricts to the identity map on $M(P,g)$.
\end{prop}

Unfortunately, Taubes' proof of Proposition \ref{prop:Existence_strong_deformation_retraction_moduli_space_asd_connections} does not yield an analogue of the \apriori estimate \eqref{eq:Tubular_neighborhood_moduli_space_asd_connections_projection_bound} provided by Proposition \ref{prop:Existence_tubular_neighborhood_moduli_space_asd_connections}. What one can extract from his proof is summarized below in the


\begin{prop}[Global existence, convergence, and stability of the anti-self-dual curvature flow]
\label{prop:Taubes_1984b_section3}
\cite[Section 3]{TauPath}
Let $G$ be a compact Lie group and $P$ a principal $G$-bundle over a closed, four-dimensional, oriented, smooth manifold, $X$, with Riemannian metric, $g$, and let $\mu_0$ be a positive constant. Then there is a positive constant, $\eps = \eps(g,P,\mu_0) \in (0,1]$, with the following significance. Suppose that $M(P,g)$ is non-empty and let $\sA_\eps(P)$ denote the open subset of the affine space, $\sA(P)$, of connections, $A$ on $P$, of class $W^{1,3}$ obeying
$\fA_g(A) < \eps$, so that $\sB_\eps(P) = \sA_\eps/\Aut P$ in \eqref{eq:Open_neighborhood_asd_moduli_space_FA+_L2_small}. Let $A_1$ be a fixed $C^\infty$ reference connection on $P$. If
$$
\mu(A) \geq \mu_0 \quad\forall\, A \in \sA_\eps(P),
$$
where $\mu(A)$ is as in \eqref{eq:Least_eigenvalue_dA+dA+*}, then the following hold for each initial connection, $A_0 \in \sA_\eps(P)$.

\begin{enumerate}
  \item
  \label{item:Taubes_1984b_anti-self-dual_curvature_flow_is_gradient_flow_for_fA}
  The anti-self-dual curvature flow \eqref{eq:Taubes_1984b_3-2_and_3} is the gradient flow for the functional $\fA$ in \eqref{eq:Taubes_1984b_2-2_and_page_349}.

  \item
  \label{item:Taubes_1984b_step_in_proof_proposition_3-1}
  There is a unique solution, $A(t) = A_0+\alpha(t)$ for $t\in[0,\infty]$, to \eqref{eq:Taubes_1984b_3-2_and_3}, with
  $$
  \alpha \in C([0,\infty]; W^{1,3}_{A_1}(X; \Lambda^1\otimes\ad P);
  $$

  \item
  \label{item:Taubes_1984b_lemma_3-3}
  The solution, $A(t)$, obeys
  \begin{equation}
  \label{eq:Taubes_1984b_lemma_3-3}
  \|A(t) - A_\infty\|_{W_{A_1}^{1,3}(X)} \to 0, \quad\hbox{as } t \to \infty,
  \end{equation}
  where $A_\infty$ is a $g$-anti-self-dual connection on $P$ of class $W^{1,3}$;

  \item
  \label{item:Taubes_1984b_lemma_3-4}
  The convergence in \eqref{eq:Taubes_1984b_lemma_3-3} is uniform with respect to initial connections, $A_0 \in \sA_\eps(P)$, in closed $W^{1,3}_{A_1}(X; \Lambda^1\otimes\ad P)$ balls;

  \item
  \label{item:Taubes_1984b_lemma_3-11}
  There is an open $W^{1,3}_{A_1}(X; \Lambda^1\otimes\ad P)$-neighborhood, $\fN(A_0) \subset \sA(P)$, of $A_0$ and a positive constant, $c = c(\fN(A_0))$, such that, for all $A_0' \in \fN(A_0)$, the solution $A(t)$ to \eqref{eq:Taubes_1984b_3-2_and_3} with $A(0)=A_0'$ obeys
  \begin{equation}
  \label{eq:Taubes_1984b_equation_3-22_lemma_3-11}
  \|A(t) - A_\infty\|_{L^4(X)} + \|\nabla_{A_0'}(A(t) - A_\infty)\|_{L^3(X)}
  \leq ce^{-t}, \quad\forall\, t \in [0,\infty);
  \end{equation}

  \item
  \label{item:Taubes_1984b_equations_3-7_and_8}
  One has
  \begin{equation}
  \label{eq:Taubes_1984b_3-7}
  \frac{\partial}{\partial t}F_{A(t)}^+ = -F_{A(t)}^+, \quad\hbox{a.e. } t \in (0,\infty),
  \end{equation}
  and, for $p \in (1,3]$,
  \begin{equation}
  \label{eq:Taubes_1984b_3-8}
  \|F_{A(t)}^+\|_{L^p(X)} \leq e^{-t}\|F_{A(0)}^+\|_{L^p(X)}, \quad\forall\, t \in [0,\infty).
  \end{equation}
\end{enumerate}
\end{prop}

\begin{proof}
Proposition \ref{prop:Taubes_1984b_section3} summarizes results obtained by Taubes in the course of the proof of his \cite[Proposition 3.1]{TauPath}, thus our task here is just to translate from his notation and verify a few elementary consequences. In \cite{TauPath}, Taubes considers `based connections', $(A,h)$, consisting of pairs of connections $A$ on $P$ and points $h$ in a fiber $P_{x_0}$, for some point $x_0 \in X$). We need only consider the simpler case of connections, $A$, on $P$.

See \cite[p. 349]{TauPath} for Item \eqref{item:Taubes_1984b_anti-self-dual_curvature_flow_is_gradient_flow_for_fA}.
We observe that Taubes' proof of \cite[Proposition 3.1]{TauPath} (see \cite[p. 352]{TauPath}) establishes Item \eqref{item:Taubes_1984b_step_in_proof_proposition_3-1}. His \cite[Lemma 3.3]{TauPath} asserts that for a fixed $A_0$, the flow $A(t)$ converges strongly in $W^{1,3}_{A_1}(X; \Lambda^1\otimes\ad P)$ as $t \to \tau \in (0, \infty]$ and, as we may take $\tau = \infty$, this gives Item \eqref{item:Taubes_1984b_lemma_3-3}. His \cite[Lemma 3.4]{TauPath} asserts that the convergence of $A(t)$, as $t \to \tau \in (0, \infty]$, described in \cite[Lemma 3.4]{TauPath} is uniform on closed $W^{1,3}_{A_1}(X; \Lambda^1\otimes\ad P)$-balls in $\sA(P)$ for $\tau \in (0,\infty]$. As we may take $\tau = \infty$, this gives Item \eqref{item:Taubes_1984b_lemma_3-4}.

The convergence rate estimate \eqref{eq:Taubes_1984b_equation_3-22_lemma_3-11} follows from
\cite[Equation (3.22) and Lemma 3.11]{TauPath} which together assert that
$$
\|A(t) - A(t')\|_{L^4(X)} + \|\nabla_{A_0'}(A(t) - A(t')\|_{L^3(X)}
\leq
ce^{-t}, \quad\forall\, t, t' \in [0,\tau),
$$
where $c$ is independent of $\tau$ and $t,t' \in [0,\tau)$. We obtain \eqref{eq:Taubes_1984b_equation_3-22_lemma_3-11} from the preceding estimate by noting that we can take $\tau = \infty$. Item \eqref{item:Taubes_1984b_equations_3-7_and_8} follows from \cite[Equations (3.7) and (3.8)]{TauPath}. This concludes the proof of Proposition \ref{prop:Taubes_1984b_section3}.
\end{proof}

\begin{rmk}[An application of Proposition \ref{prop:Taubes_1984b_section3} to reduce the $L^p$ norm of the self-dual curvature]
\label{rmk:Reducing_Lp_norm_FA+_via_anti_self_curvature_flow}
Our Theorem \ref{thm:Proposition_Feehan_Leness_7-6} requires the hypothesis,
$$
\|F_{A_0}^+\|_{L^{\sharp,2}(X)} \leq \eps,
$$
in order to be able to solve the $g$-anti-self-dual equation, $F_{A_0+a}^{+,g}=0$ for a perturbation, $a \in H_{A_0}^1(X;\Lambda^1\otimes \ad P)$. However, if $A_0$ only obeys
$$
\|F_{A_0}^+\|_{L^2(X)} \leq \eps,
$$
then we can apply Proposition \ref{prop:Taubes_1984b_section3} to find a large enough $T \in [0,\infty)$ such that $A(T)$ obeys
$$
\|F_{A(T)}^+\|_{L^3(X)} \leq \eps,
$$
where $A(t)$ is the anti-self-dual curvature flow with initial data, $A(0) = A_0$. But there is a continuous Sobolev embedding, $L^3(X) \hookrightarrow L^\sharp(X)$ by Lemma \ref{lem:Feehan_4-1}, and thus
$$
\|F_{A(T)}^+\|_{L^3(X)} \leq \eps,
$$
for a positive constant, $c = c(g)$. Theorem \ref{thm:Proposition_Feehan_Leness_7-6} will then apply to the connection, $A(T)$.
\end{rmk}

\section[Proofs of main corollaries]{Proofs of main corollaries of Theorem \ref{mainthm:Yang-Mills_gradient_flow_global_existence_and_convergence_started_near_minimum}}
\label{sec:Yang-Mills_gradient_flow_initial_connection_almost_minimal_energy_proofs}
In this section, we provide the proofs of the main corollaries of Theorem \ref{mainthm:Yang-Mills_gradient_flow_global_existence_and_convergence_started_near_minimum} described in Section \ref{subsec:Yang-Mills_gradient_flow_initial_connection_almost_minimal_energy}. Those results are stated, for the sake simplicity, for the case where the initial data, $A_0$, is a $C^\infty$ connection and we continue to make this assumption here. More general results for $A_0$ belonging to a suitable Sobolev class may be easily extracted \mutatis from the proofs of this section by referring to the proof of Theorem \ref{mainthm:Yang-Mills_gradient_flow_global_existence_and_convergence_started_near_minimum} in Section \ref{sec:Application_abstract_gradient_system_results_Yang-Mills_energy_functional}.

\begin{proof}[Proof of Corollary \ref{maincor:Yang-Mills_gradient_flow_global_existence_and_convergence_started_near_connection with almost minimal energy_and_vanishing_H2+}]
Theorem \ref{thm:Proposition_Feehan_Leness_7-6} provides, for small enough $\eta = \eta(E_0,g,\mu_0) \in (0,1]$, a constant, $C_0 = C_0(E_0,g,\mu_0) \in (0,\infty)$, and a unique $v \in \Omega^{+,g}(X;\ad P)$ such that $A_{\asd} := A_0 + d_{A_0}^{+,g}v$ is a $g$-anti-self-dual connection on $P$ that obeys
$$
\|A_0 - A_{\asd}\|_{H_{A_0}^1(X)} \leq C_0\|F_{A_0}^{+,g}\|_{L^{\sharp,2}(X)}.
$$
The remaining assertions now follow immediately from Theorem \ref{mainthm:Yang-Mills_gradient_flow_global_existence_and_convergence_started_near_minimum}.
\end{proof}

\begin{proof}[Proof of Corollary \ref{maincor:Yang-Mills_gradient_flow_global_existence_and_convergence_started_near_connection_with_almost_minimal_energy_and_g_is_good_Lp_version}]
We have $\|F_{A_0}^{+,g}\|_{L^{\sharp,2}(X)} \leq c_p\|F_{A_0}^{+,g}\|_{L^p(X)}$, for a positive constant $c_p = c_p(g,p)$, by Lemma \ref{lem:Feehan_4-1}, and $\|F_{A_0}^{+,g}\|_{L^p(X)} \leq \eta$ by hypothesis \eqref{eq:FA0+_Lp_small_and_FA0_Lp_bounded}, and therefore $\|F_{A_0}^{+,g}\|_{L^{\sharp,2}(X)} \leq c_p\eta$. Since $g$ is good by hypothesis, for a small enough $\eta = \eta(g,P,p) \in (0,1]$, Theorem \ref{thm:Good_metric_implies_positive_lower_bound_small_eigenvalues_almost_ASD} provides a positive constant, $\mu_0 = \mu_0(g,P,p)$, such that
$$
\mu_g(A_0) \geq \frac{\mu_0}{2}.
$$
For a small enough $\eta = \eta(\sE(A_0),g,P,p) \in (0,1]$, Theorem \ref{thm:Proposition_Feehan_Leness_7-6} thus provides a constant, $C_0 = C_0(\sE(A_0),g,P) \in (0,\infty)$, and a unique $v \in \Omega^{+,g}(X;\ad P)$ such that $A_{\asd} := A_0 + d_{A_0}^{+,g}v$ is a $g$-anti-self-dual connection on $P$ that obeys
\begin{align}
\label{eq:A0_minus_Aasd_L4_norm_leq_C1_FA0+_Lsharp2_norm}
\|A_0 - A_{\asd}\|_{L^4(X)} &\leq C_0\|F_{A_0}^{+,g}\|_{L^{\sharp,2}(X)} \leq C_0\eta,
\\
\label{eq:A0_minus_Aasd_H1_norm_leq_C1_FA0+_Lsharp2_norm}
\|A_0 - A_{\asd}\|_{H_{A_0}^1(X)} &\leq C_0\|F_{A_0}^{+,g}\|_{L^{\sharp,2}(X)} \leq C_0\eta,
\end{align}
where, as usual, $\sE(A_0) = \frac{1}{2}|F_{A_0}\|_{L^2(X)}^2$.

Moreover, Lemma \ref{lem:Apriori_Lp_estimate_solution_ASD_equation} implies that, for small enough $\eta = \eta([A_1],g,P,p) \in (0,1]$, there is a positive constant, $C_1 = C_1([A_1],g,P,p)$, such that
\begin{align}
\label{eq:A0_minus_Aasd_Lq_norm_leq_C1_FA0+_Lp_norm}
\|A_0 - A_{\asd}\|_{L^q(X)} &\leq C_1\|F_{A_0}^{+,g}\|_{L^p(X)} \leq C_1\eta,
\\
\label{eq:A0_minus_Aasd_W1p_norm_leq_C1_FA0+_Lp_norm}
\|A_0 - A_{\asd}\|_{W_{A_1}^{1,p}(X)} &\leq C_1\|F_{A_0}^{+,g}\|_{L^p(X)} \leq C_1\eta,
\end{align}
where $q\in (4,\infty)$ is defined by $1/p = 1/4+1/q$. Using the curvature expansion,
$$
F(A_\asd) = F_{A_0} + d_{A_0}(A_0 - A_{\asd}) + (A_0 - A_{\asd})\wedge (A_0 - A_{\asd}),
$$
and the inequality $\|d_{A_0}(A_0 - A_{\asd})\|_{L^p(X)} \leq 4\|d_{A_1}(A_0 - A_{\asd})\|_{L^p(X)}$ by Lemma \ref{lem:Equivalence_W1p_or_W2p_norms_connections_A0_and_A1_L4_or_W14_close},
\begin{align*}
\|F(A_\asd)\|_{L^p(X)}
&\leq \|F_{A_0}\|_{L^p(X)} + \|d_{A_0}(A_0 - A_{\asd})\|_{L^p(X)}
+ 2\|A_0 - A_{\asd}\|_{L^4(X)} \|A_0 - A_{\asd}\|_{L^q(X)}
\\
&\leq \|F_{A_0}\|_{L^p(X)} + 4C_1\|F_{A_0}^{+,g}\|_{L^p(X)}
+ 2C_0C_1\|F_{A_0}^{+,g}\|_{L^{\sharp,2}(X)}\|F_{A_0}^{+,g}\|_{L^p(X)}
\\
&\qquad \hbox{(by \eqref{eq:A0_minus_Aasd_L4_norm_leq_C1_FA0+_Lsharp2_norm}, \eqref{eq:A0_minus_Aasd_Lq_norm_leq_C1_FA0+_Lp_norm}, and \eqref{eq:A0_minus_Aasd_W1p_norm_leq_C1_FA0+_Lp_norm})}
\\
&\leq \|F_{A_0}\|_{L^p(X)} + 4C_1\eta + 2c_pC_0C_1\eta^2
\\
&\leq K_1,
\end{align*}
for $K_1 := K + 4C_1 + 2c_pC_0C_1$. Thus, $[A_\asd]$ belongs to the compact subset $\bar\sU(g,K_1,P,p) \subset M(P,g)$ given by the closure (in the Uhlenbeck topology \cite[Definition 4.4.1]{DK} on $M(P,g)$) of the precompact open subset,
\begin{equation}
\label{eq:Precompact_subset_MPg_FA_Lp-norm_leq_bound}
\sU(K_1,P,p,g) := \{[A] \in M(P,g): \|F_A\|_{L^p(X)} < K_1\} \Subset M(P,g).
\end{equation}
For each $[A] \in M(P,g)$, there is a {\L}ojasiewicz-Simon constant, $\sigma([A_1],[A],g) \in (0,1]$, defined by Theorem \ref{mainthm:Yang-Mills_gradient_flow_global_existence_and_convergence_started_near_minimum}. Because the set $\bar\sU(K_1,P,p,g)$ is compact, it has a finite cover by open balls,
$$
B([A_\alpha],\sigma_\alpha)
:=
\left\{[A]\in M(P,g): \dist_{H_{A_1}^1}([A],[A_\alpha]) < \sigma_\alpha\right\},
$$
with radii $\sigma_\alpha = \sigma_\alpha([A_1],[A_\alpha],g) \in (0,1]$, where
$$
\dist_{H_{A_1}^1}([A],[A_\alpha])
:=
\inf_{u \in \Aut P}\|u^*A - A_\alpha\|_{H_{A_1}^1(X)},
$$
and the infimum is taken over all $C^\infty$ gauge transformations of $P$. We let
$$
\sigma([A_1],g,K,P,p) := \min_\alpha \sigma([A_1],[A_\alpha],g),
$$
where we take the minimum over the index set corresponding to this finite open cover. We have
\begin{align*}
\|A_0 - A_{\asd}\|_{H_{A_1}^1(X)}
&\leq
4\|A_0 - A_{\asd}\|_{H_{A_0}^1(X)}
\quad\hbox{(by Lemma \ref{lem:Equivalence_W1p_or_W2p_norms_connections_A0_and_A1_L4_or_W14_close})}
\\
&\leq 4C_0\|F_{A_0}^{+,g}\|_{L^{\sharp,2}(X)}
\quad\hbox{(by \eqref{eq:A0_minus_Aasd_H1_norm_leq_C1_FA0+_Lsharp2_norm})}
\\
&\leq 4c_pC_0\|F_{A_0}^{+,g}\|_{L^p(X)}
\quad\hbox{(by Lemma \ref{lem:Feehan_4-1})}
\\
&\leq 4c_pC_0\eta
\quad\hbox{(by \eqref{eq:FA0+_Lp_small})}.
\end{align*}
Thus, for a small enough $\eta([A_1],g,K,P,p) \in (0,1]$, we may suppose that $4c_pC_0\eta < \sigma$ and therefore,
$$
\|A_0 - A_{\asd}\|_{H_{A_1}^1(X)} < \sigma.
$$
The conclusions now follow from Theorem \ref{mainthm:Yang-Mills_gradient_flow_global_existence_and_convergence_started_near_minimum}.
This completes the proof of Corollary \ref{maincor:Yang-Mills_gradient_flow_global_existence_and_convergence_started_near_connection_with_almost_minimal_energy_and_g_is_good_Lp_version}.
\end{proof}

\begin{proof}[Proof of Corollary \ref{maincor:Yang-Mills_gradient_flow_global_existence_and_convergence_started_near_connection_with_almost_minimal_energy_and_g_is_good_L2sharp_version}]
Observe that $E_0 := \|F_{A_0}\|_{L^2(X)} \leq c_p\|F_{A_0}\|_{L^p(X)} \leq c_pK$ by hypothesis \eqref{eq:FA0_Lp_bounded}. Just as in the proof of Corollary \ref{maincor:Yang-Mills_gradient_flow_global_existence_and_convergence_started_near_connection_with_almost_minimal_energy_and_g_is_good_Lp_version}, Theorem \ref{thm:Proposition_Feehan_Leness_7-6} provides, for $\mu_0 = \mu(A_0)$ and small enough $\eta = \eta(E_0,g,\mu_0) \in (0,1]$, a constant, $C_0 = C_0(E_0,g,\mu_0) \in (0,\infty)$, and a unique $v \in \Omega^{+,g}(X;\ad P)$ such that $A_{\asd} = A_0 + d_{A_0}^{+,*}v$ is a $g$-anti-self-dual connection on $P$ and
$$
\|v\|_{H_{A_0}^2(X)} + \|v\|_{C(X)} \leq C_0\|F_{A_0}^+\|_{L^{\sharp,2}(X)}.
$$
The conclusions now follow from the proof of Corollary \ref{maincor:Yang-Mills_gradient_flow_global_existence_and_convergence_started_near_connection_with_almost_minimal_energy_and_g_is_good_Lp_version}, noting that $[A_\asd]$ belongs to the compact subset $\bar\sU(g,K,P,p) \subset M(P,g)$ defined in \eqref{eq:Precompact_subset_MPg_FA_Lp-norm_leq_bound} by hypothesis \eqref{eq:FAasd_Lp_bounded}.
\end{proof}

\setcounter{sectioncontinuation}{\thesection}
\appendix

\appendix

\chapter[Embeddings for fractional-order Sobolev spaces]{Continuous and compact embeddings for fractional-order Sobolev spaces}
\label{chap:Embedding_multiplication_theorems_fractional_order_Sobolev_spaces}
While comprehensive references for continuous and compact embedding theorems, together with multiplication theorems, are easily found for standard integer order Sobolev spaces --- for example Adams and Fournier \cite{AdamsFournier} in the case of embedding theorems and Freed and Uhlenbeck \cite{FU} and Palais \cite{PalaisFoundationGlobal} for multiplication theorems --- the analogous references for fractional order Sobolev spaces are less accessible. A further complication is that there are numerous definitions of fractional order Sobolev spaces and they are not necessarily equivalent in all cases. While the desired results can usually be guessed by extrapolation of known results for integer order Sobolev spaces, their justification can be difficult and dependent on the particular definition for a fractional order Sobolev space. For these reasons, we gather the most useful results in this Appendix and provide references to proofs in the literature for cases of
\begin{inparaenum}[\itshape i\upshape)]
\item $H^s$ spaces, following Lions and Magenes \cite{Lions_Magenes_v1_1972} and Taylor \cite{Taylor_PDE1}
\item Bessel potential spaces $H^{s,p}$, following Taylor \cite{Taylor_PDE3}, and
\item Slobodeckij spaces $W^{s,p}$, following Di Nezza, Palatucci, and Valdinoci \cite{DiNezza_Palatucci_Valdinoci_2012}, Haroske and Triebel \cite{Haroske_Triebel_2008}, and Maz'ya \cite{Mazya_2011}.
\end{inparaenum}

\setcounter{section}{\thesectioncontinuation}
\section{Fractional order Sobolev spaces}
\label{subsec:Fractional_order_Sobolev_spaces}
We begin with a review of some of the different approaches to defining fractional order Sobolev spaces.

\subsection{Standard Sobolev spaces}
\label{sec:Standard_Sobolev_spaces}
We first recall the standard definition of Sobolev spaces of integer order in terms of derivatives \cite[Sections 3.2 and 3.7]{AdamsFournier}. Given a domain $\Omega \subset \RR^d$ with $d \geq 2$, for $1 \leq p \leq \infty$ and integer $k \geq 0$, one sets
\[
W^{k,p}(\Omega)
:=
\left\{u \in L^p(\Omega): D^\alpha u \in L^p(\Omega),\quad\text{for any } \alpha \in \NN^k \text{ with } |\alpha| \leq k\right\},
\]
where the derivatives $D^\alpha u$ are interpreted in the weak (distributional) sense. For $1 < p \leq \infty$ and integer $k < 0$, one may define $W^{k,p}(\Omega)$ by duality (see \cite[Theorem 3.9 and Section 3.14]{AdamsFournier} and \cite[Theorem 3.12 and Section 3.13]{AdamsFournier} for the dual of $W_0^{k,p}(\Omega)$)
\[
W^{k,p}(\Omega) := (W_0^{-k,p'}(\Omega))',
\]
where $1 \leq p' < \infty$ is the dual exponent defined by $1/p + 1/p' = 1$ and $W_0^{k,p}(\Omega)$ is the closure of $C_0^\infty(\Omega)$ in $W^{k,p}(\Omega)$ \cite[Section 3.2]{AdamsFournier}. When $p = 2$, one denotes $W^{k,2}(\Omega) = H^k(\Omega)$ for $k \in \ZZ$.

One approach to development of embedding results for Sobolev spaces on domains involves first establishing results on Euclidean space then appealing to extension results to generalize to the case of domains \cite[Section 5.17]{AdamsFournier}. See \cite[Theorems 5.21, 5.22, and Remark 5.33.1]{AdamsFournier} for existence of a total extension operator when $\Omega$ is uniformly $C^k$-regular for all integers $k \geq 0$ \cite[Section 4.10]{AdamsFournier} and has a bounded boundary and see \cite[Theorem 5.24]{AdamsFournier} and \cite[Chapter 6]{Stein} for the Stein Extension Theorem for existence of a total extension operator when $\Omega$ obeys the strong local Lipschitz condition \cite[Section 4.9]{AdamsFournier}.

\subsection{Complex interpolation}
\label{subsec:Complex_interpolation}
Following \cite[Section 7.57]{AdamsFournier}, one can define a scale of fractional order Sobolev spaces by complex interpolation \cite[Section 7.51]{AdamsFournier}, \cite[Section 1.5]{Lions_Magenes_v1_1972}, \cite[Section 4.2]{Taylor_PDE1} between $L^p$ and (integer-order) Sobolev spaces. Specifically, if $s > 0$ and $k = \lceil s \rceil$ is the smallest integer greater than $s$ and $\Omega \subset \RR^d$ is a domain, one defines
\[
W^{s,p}(\Omega) := [L^p(\Omega), W^{k,p}(\Omega)]_{s/k}.
\]
If $s$ is a positive integer and $\Omega$ has a suitable extension property, then $W^{s,p}(\Omega)$ coincides with the usual Sobolev space with the same notation \cite[Sections 3.2 and 3.7]{AdamsFournier}. When it comes to characterizing Sobolev spaces of fractional order more concretely, there are several approaches in the literature \cite[Section 1]{DeNapoli_Drelichman_2014arxiv} as we recall next.

\subsection{$H^s$ spaces}
\label{subsec:Hs_spaces}
A definition of fractional order Sobolev spaces via Fourier transforms is straightforward when $p=2$ and $\Omega=\RR^d$; see \cite[Section 7.62]{AdamsFournier}, \cite[Sections 1.7 and 1.9]{Lions_Magenes_v1_1972}, \cite[Section 4.1]{Taylor_PDE1}.

For real $s \in \RR$, let \cite[Section 1.7]{Lions_Magenes_v1_1972}
\[
H^s(\RR^d) := \left\{u \in L^2(\RR^d) : \int_{\RR^d} |\hat u(\omega)|^2 \left(1 + |\omega|^2\right)^s \, d\omega < \infty\right\},
\]
where $\hat u$ denotes the Fourier transform \cite[Section 7.53]{AdamsFournier} of $u \in L^1(\RR^d)$,
\[
\hat u(\omega) := \sF(u)(\omega) = \frac{1}{(2\pi)^{d/2}}\int_{\RR^d} e^{ix\cdot\omega} u(x)\,dx, \quad\forall\, \omega \in \RR^d.
\]
For integer $s$, one has $H^s(\RR^d) = W^{k,2}(\RR^d)$ \cite[Section 7.62]{AdamsFournier}, \cite[Theorem 1.9.7]{Lions_Magenes_v1_1972}.

For real $s \geq 0$ and a domain $\Omega \subset \RR^d$, one can define $H^s(\Omega)$ by complex interpolation \cite[Section 9.1]{Lions_Magenes_v1_1972},
\begin{equation}
\label{eq:Lions_Magenes_1-9-1}
H^s(\Omega) := [L^2(\Omega), H^k(\Omega)]_\theta,
\end{equation}
where $k\geq 0$ is an integer and $\theta \in (0,1)$ obeys $s = k\theta$. The properties of $H^s(\Omega)$ are elucidated in the

\begin{thm}[Properties of $H^s$ spaces on domains]
\label{thm:Lions_Magenes_1-9-1}
\cite[Theorem 1.9.1]{Lions_Magenes_v1_1972}
Let $d \geq 2$ and $\Omega \subset \RR^d$ be a bounded domain with a $C^\infty$ boundary \cite[Equations (1.7.10) and (1.7.11)]{Lions_Magenes_v1_1972}. Then $H^s(\Omega)$ coincides (algebraically) with the space of restrictions to $\Omega$ of functions in $H^s(\RR^d)$, the definition \eqref{eq:Lions_Magenes_1-9-1} of $H^s(\Omega)$ is independent of $k$, and when $s=k$ is an integer, then $H^s(\Omega) = H^k(\Omega)$ (standard Sobolev space).
\end{thm}

To define $H^s(\Omega)$ for $s < 0$, one sets \cite[Section 12.1]{Lions_Magenes_v1_1972}
\[
H^s(\Omega) := \left(H_0^{-s}(\Omega)\right)'.
\]

\subsection{Bessel potential spaces}
\label{subsec:Bessel_potential_spaces}
When $1<p<2$ or $2<p<\infty$, the generalization of the preceding method introduces Fourier transforms of certain Bessel functions and the resulting spaces $H^{s,p}(\RR^d)$ are called \emph{spaces of Bessel potentials} \cite[Section 7.6.3]{AdamsFournier}, \cite[Section 1]{DeNapoli_Drelichman_2014arxiv}, \cite[Section 13.6]{Taylor_PDE3}. For $1<p<\infty$, one defines the spaces of Bessel potentials (or classical potential spaces) by
\[
H^{s,p}(\RR^d) := \left\{u : u = (1 + \Delta)^{−s/2}f \text{ for some } f \in L^p(\RR^d) \right\},
\]
where we use the sign convention $\Delta = -\sum_{i=1}^d \partial_{x_i}^2$ and the fractional power $(1 + \Delta)^{−s/2}$ can be defined by means of the Fourier transform (for functions in the Schwartz class of rapidly decreasing functions \cite[Section 7.59]{AdamsFournier}),
\[
(1 + \Delta)^{−s/2}f := \sF^{-1}\left(\left(1 + |\omega|^2\right)^{-s/2}\sF(f)\right) = G_s \star f,
\]
where
\begin{align*}
G_s(x)
&:=
\sF^{-1}\left(1 + |\omega|^2\right)^{-s/2}(x)
\\
&\,=
\frac{1}{(4\pi)^{d/2}\Gamma(s/2)} \int_0^\infty e^{-t} e^{-|x|^2/(4t)} t^{(s - d)/2}\, \frac{dt}{t}, \quad\forall\, x \in \RR^d,
\end{align*}
is called the \emph{Bessel potential}. Classical references on Bessel potential spaces include Aronszajn, Mulla, Szeptycki, and Smith \cite{Aronszajn_Mulla_Szeptycki_1963, Aronszajn_Smith_1961}, Calder{\'o}n \cite{Calderon_1961}, and Stein \cite{Stein}. For integer $k \geq 0$, one has  \cite[Section 7.62]{AdamsFournier}
\[
H^{k,p}(\RR^d) = W^{k,p}(\RR^d).
\]
Indeed, $H^{k,p}(\RR^d) = W^{k,p}(\RR^d)$ for $1 < p < \infty$ by \cite[Theorem 7]{Calderon_1961} when $k=1$ and by \cite[Proposition 13.6.1]{Taylor_PDE3} any integer $k \geq 1$.

Also $H^{s,2}(\RR^d) = H^s(\RR^d)$ for any real $s \geq 0$ by Plancherel's Theorem \cite[Theorem 7.61]{AdamsFournier}, \cite[Section 1]{DeNapoli_Drelichman_2014arxiv}.

\subsection{Slobodeckij spaces}
\label{subsec:Slobodeckij_spaces}
Another (non-equivalent) definition of fractional order Sobolev spaces is given by the \emph{Aronszajn}, \emph{Gagliardo}, or \emph{Slobodeckij} spaces (the nomenclature varies) \cite[Section 1]{DeNapoli_Drelichman_2014arxiv}. Following Di Nezza, Palatucci, and Valdinoci, for $0 < s < 1$ and $1 \leq p < \infty$, given a domain $\Omega \subset \RR^d$, one sets \cite[Equation (2.1)]{DiNezza_Palatucci_Valdinoci_2012}
\[
W^{s,p}(\Omega) := \left\{u \in L^p(\Omega) : \frac{|u(x) - u(y)|}{|x - y|^{(d/p) + s}} \in L^p(\Omega \times \Omega) \right\},
\]
that is, an intermediary Banach space between $L^p(\Omega)$ and $W^{1,p}(\Omega)$, endowed with the natural norm defined by \cite[Equation (2.2)]{DiNezza_Palatucci_Valdinoci_2012}
\[
\|u\|_{W^{s,p}(\Omega)}^p := \|u\|_{W^{s,p}(\Omega)}^p + [u]_{W^{s,p}(\Omega)}^p,
\]
where the term
\[
[u]_{W^{s,p}(\Omega)}
:=
\left( \int_\Omega\int_\Omega \frac{|u(x)-u(y)|^p}{|x-y|^{d+sp}} \,dx\,dy \right)^{1/p}
\]
is called the \emph{Gagliardo (semi-)norm} of $u$. Maz'ya defines the space $W_0^{s,p}(\RR^d)$ as the completion of $C_0^\infty(\RR^d)$ with respect to the Gagliardo norm $[\cdot]_{W^{s,p}(\Omega)}$.

It is known that $H^{s,2}(\RR^d) = W^{s,2}(\RR^d)$ for any $0 < s < 1$ (for example, see \cite[Proposition 3.4]{DiNezza_Palatucci_Valdinoci_2012} for the fact that $H^s(\RR^d) = W^{s,2}(\RR^d)$), but for $p \neq 2$, the spaces $H^{s,p}(\RR^d)$ and $W^{s,p}(\RR^d)$ are different.

When $s \geq 1$ and is not an integer, one writes $s = k + \sigma$ , where $k\geq 0$ is an integer and $\sigma \in (0, 1)$. In this case the space $W^{s,p}(\Omega)$ consists of those equivalence classes of functions $u \in W^{k,p}(\Omega)$ whose distributional derivatives $D^\alpha u$, with $|\alpha| = k$, belong to $W^{\sigma,p}(\Omega)$, namely \cite[Equation (2.10)]{DiNezza_Palatucci_Valdinoci_2012}
\[
W^{s,p}(\Omega) := \left\{ u \in W^{k,p}(\Omega): D^\alpha u \in W^{\sigma,p}(\Omega) \text{ for any } \alpha
\text{ such that } |\alpha| = k\right\},
\]
and this is a Banach space with respect to the norm \cite[Equation (2.11)]{DiNezza_Palatucci_Valdinoci_2012}
\[
\|u\|_{W^{s,p}(\Omega)} := 	\left( \|u\|_{W^{k,p}(\Omega)}^p + \sum_{|\alpha|=k} \|D^\alpha u\|_{W^{\sigma,p}(\Omega)}^p \right)^{1/p}.
\]
If $s = k$ is an integer, then the space $W^{s,p}(\Omega)$ coincides with the Sobolev space $W^{k,p}(\Omega)$.

According to \cite[Theorem 7.38]{Adams}, \cite[Theorem 2.4]{DiNezza_Palatucci_Valdinoci_2012}, for any $s > 0$, the space $C_0^\infty(\RR^d)$ of smooth functions with compact support is dense in $W^{s,p}(\Omega)$. One defines the space $W_0^{s,p}(\Omega)$ as the completion of $C_0^\infty(\Omega)$ with respect to the Slobodeckij norm $\|\cdot\|_{W^{s,p}(\Omega)}$. If $s < 0$ and $p \in (1,\infty)$, one can define \cite[Remark 2.5]{DiNezza_Palatucci_Valdinoci_2012}
\[
W^{s,p}(\Omega) := \left(W_0^{-s,p'}(\Omega)\right)'
\]
where $1/p + 1/p' = 1$. In this case, the space $W^{s,p}(\Omega)$ is a space of distributions on $\Omega$, since it is the dual of a space having $C_0^\infty(\Omega)$ as dense subset.

\section{Continuous embeddings}
\label{sec:Continuous_embeddings}
The standard Sobolev Embedding \cite[Theorem 4.12]{AdamsFournier} (for integer $k$) may be extended to non-integral $s\in\RR$.

\subsection{Continuous embeddings for $H^s$ spaces}
\label{subsec:Continuous_embeddings_Hs_spaces}
For the fractional order Sobolev spaces, $H^s(\Omega)$, we recall the

\begin{thm}[Continuous embedding for $H^s$ spaces on domains]
\label{thm:Lions_Magenes_1-9-8}
\cite[Theorem 1.9.8 and Corollary 1.9.1]{Lions_Magenes_v1_1972}
\cite[Proposition 4.3]{Taylor_PDE1}
Let $d \geq 1$ and $\Omega \subset \RR^d$ be a bounded domain with a $C^\infty$ boundary \cite[Equations (1.7.10) and (1.7.11)]{Lions_Magenes_v1_1972}. If $s \in \RR$ obeys $s > d/2$, then the following embedding is continuous:
\[
H^s(\Omega) \subset C(\bar\Omega).
\]
Moreover, if $s = k + \sigma$ for an integer $k \geq 0$ and $\sigma \in \RR$ obeys $\sigma > d/2$, then the following embedding is continuous:
\[
H^s(\Omega) \subset C^k(\bar\Omega).
\]
Finally, if $s = d/2 + \alpha$, for $\alpha \in (0,1)$, then the following embedding is continuous:
\[
H^s(\Omega) \subset C^\alpha(\bar\Omega).
\]
\end{thm}

Here, $C^k(\Omega)$ denotes the set of functions having all derivatives of order at most $k$ continuous in $\Omega$, and $C^k(\bar\Omega)$ denotes the set of functions in $C^k(\Omega)$ all of whose derivatives of order at most $k$ have continuous extensions to the closure $\bar\Omega$ \cite[Chapter 1]{GilbargTrudinger}, while $C^\alpha(\Omega)$ and $C^\alpha(\bar\Omega)$ denote the H\"older spaces \cite[Section 4.1]{GilbargTrudinger}.

\subsection{Continuous embeddings for Bessel potential spaces}
\label{subsec:Continuous_embeddings_Bessel_potential_spaces}
For the Bessel potential spaces, $H^{s,p}(\Omega)$, we recall the

\begin{prop}[Continuous embedding for Bessel potential spaces on Euclidean space]
\label{prop:Taylor_PDE_proposition_13-6-3and4}
\cite[Propositions 13.6.3 and 13.6.4]{Taylor_PDE3}
If $d \geq 1$ is an integer, $p \in (1,\infty)$, and $s \in \RR$, then the following embeddings are continuous:
\begin{equation}
\label{eq:Taylor_PDE_proposition_13-6-3and4}
H^{s,p}(\RR^d)
\subset
\begin{cases}
L^q(\RR^d), &\hbox{for } 0 \leq sp < d \text{ and } q \in [p, p^*],
\\
C_b(\RR^d), &\hbox{for } sp > d,
\end{cases}
\end{equation}
where $p^* := dp/(d-sp)$.
\end{prop}

\subsection{Continuous embeddings for Bessel potential spaces via complex interpolation theory}
\label{subsec:Continuous_embeddings_Bessel_potential_spaces_complex_interpolation}
The Calder{\'o}n method of complex interpolation theory \cite{Bergh_Lofstrom_1976, Triebel_interpolation_theory_1978} provides a convenient way to derive continuous embedding results for fractional order Sobolev spaces from related results for integer order Sobolev spaces, whether for functions on Euclidean space, domains, or closed manifolds, or sections of vector bundles over closed manifolds. Following Taylor \cite[Section 13.6]{Taylor_PDE3} and recalling that \cite[Equation (13.6.1)]{Taylor_PDE3}, for $p \in (1,\infty)$ and $s \in \RR$,
\[
H^{s,p}(\RR^d) := (\Delta + 1)^{-s/2}L^p(\RR^d)
\]
one has the following interpolation result.

\begin{prop}
\label{prop:Taylor_PDE_13-6-2}
If $d \geq 1$ and $s \in \RR$, and $\theta \in (0,1)$, and $p \in (1,\infty)$, then
\[
[L^p(\RR^d), H^{s,p}(\RR^d)]_\theta = H^{\theta s,p}(\RR^d).
\]
\end{prop}

If $X$ is a closed, smooth manifold of dimension $d \geq 2$, one can use the spaces $H^{s,p}(\RR^d)$ to $H^{s,p}(X)$ via local coordinate charts and a partition of unity subordinate to a cover of $X$ by those charts as in \cite[Section 4.3]{Taylor_PDE1}. One can show that the spaces $H^{s,p}(X)$ also obey the interpolation inequality \cite[Proposition 4.3.1]{Taylor_PDE1}, \cite[Equation (13.6.6)]{Taylor_PDE3},
\begin{equation}
\label{eq:Taylor_PDE_13-6-6}
H^{s,p}(X) = [L^p(X), H^{k,p}(X)]_\theta,
\quad \theta \in (0,1), \quad s = k\theta.
\end{equation}
If $\Omega$ is a compact subdomain with smooth boundary of a closed manifold $X$, one defines $H^{k,p}(\Omega) := W^{k,p}(\Omega)$ (the standard Sobolev space). There exists a (total) extension operator \cite[Equation (4.4.12)]{Taylor_PDE1}, \cite[Equation (13.6.6)]{Taylor_PDE3}
\[
E:H^{k,p}(\Omega) \to H^{k,p}(X), \quad 0 \leq k \leq N,
\]
for every finite integer $N$, since $\partial\Omega$ is a closed, smooth manifold. One defines $H^{s,p}(\Omega)$ for $s > 0$ by \cite[Equation (13.6.7)]{Taylor_PDE3}
\begin{equation}
\label{eq:Taylor_PDE_13-6-7}
H^{s,p}(\Omega) := [L^p(\Omega), H^{k,p}(\Omega)]_\theta,
\quad \theta \in (0,1), \quad s = k\theta,
\end{equation}
and obtains an extension operator, $E:H^{s,p}(\Omega) \to H^{s,p}(X)$. One can show that the characterization \eqref{eq:Taylor_PDE_13-6-7} of $H^{s,p}(\Omega)$ is independent of $k\geq s$ (see \cite[Equation (4.4.16)]{Taylor_PDE1} and \cite[Equation (13.6.8)]{Taylor_PDE3}). The space $H^{s,p}(\Omega)$ agrees with the standard characterization when $s = k$ is a positive integer.

For the spaces $H^s(\RR^d)$ and $H^s(X)$, there are natural isomorphisms $(H^s(\RR^d))' \cong H^{-s}(\RR^d)$ and $(H^s(X))' \cong H^{-s}(X)$ for $s \in \RR$ by \cite[Proposition 3.3.2]{Taylor_PDE1}. These isomorphisms extend to the case $p \in (1,\infty)$, giving $(H^{s,p}(\RR^d))' \cong H^{-s,p'}(\RR^d)$ and $(H^{s,p}(X))' \cong H^{-s,p'}(X)$ for $s \in \RR$, where $p'$ is the dual exponent for $p$ defined by $1/p+1/p'=1$, using $(L^p(\RR^d))' \cong L^{p'}(\RR^d)$ and $(L^p(X))' \cong L^{p'}(X)$.

Complex interpolation can also be used to extend the Sobolev Embedding \cite[Theorem 4.12]{AdamsFournier} for the standard Sobolev spaces $W^{k,p}(\Omega)$ to the case of Bessel potential spaces $H^{s,p}(\Omega)$. (This is the approach of Palais to the derivation of embedding and multiplication theorems for fractional order Sobolev spaces \cite[Section 9]{PalaisFoundationGlobal}.) For this purpose, we recall the

\begin{thm}[Exact interpolation]
\label{thm:Bergh_Lofstrom_4-1-2}
\cite[Section 7.52]{AdamsFournier},
\cite[Theorem 4.1.2]{Bergh_Lofstrom_1976}
Given two \emph{compatible pairs} of Banach spaces, $(E_0, E_1)$ and $(F_0, F_1)$, the pair $([E_0, E_1]_\theta, [F_0, F_1]_\theta)$ is an \emph{exact interpolation pair} of exponent $\theta$, that is, if $T : E_0 + E_1 \to F_0 + F_1$ is a linear operator that is bounded from $E_j$ to $F_j$, for $j = 0, 1$, then for each $\theta \in [0,1]$, the operator $T$ is bounded from $[E_0, E_1]_\theta$ to $[F_0, F_1]_\theta$ and
\[
\|T\|_\theta \leq \|T\|_0^{1 - \theta} \|T\|_1^\theta.
\]
\end{thm}

In Theorem \ref{thm:Bergh_Lofstrom_4-1-2}, we denote the operator norms for $T:E_0\to F_0$, $T:E_0\to F_0$, and $T:[E_0, E_1]_\theta \to [F_0, F_1]_\theta$ by $\|T\|_0$, $\|T\|_1$, and $\|T\|_\theta$, respectively. Note that $[E_0, E_1]_0 = E_0$ and $[E_0, E_1]_1 = E_1$.

As an application of Proposition \ref{prop:Taylor_PDE_13-6-2}, Theorem \ref{thm:Bergh_Lofstrom_4-1-2}, and the standard Sobolev Embedding \cite[Theorem 4.12]{AdamsFournier}, we note the following extension of Proposition \ref{prop:Taylor_PDE_proposition_13-6-3and4}. The analogous results hold for domains $\Omega \subset \RR^d$ and closed manifolds $X$ in place of $\RR^d$.

\begin{cor}[Continuous embedding for Bessel potential spaces on Euclidean space]
\label{cor:Taylor_PDE_proposition_13-6-3and4}
If $d \geq 1$ and $k \geq 0$ are integers, $p \in (1,\infty)$, and $\sigma > 0$, then the following embeddings are continuous:
\begin{equation}
\label{eq:Taylor_PDE_proposition_13-6-3and4_k_geq_zero}
H^{k+\sigma,p}(\RR^d)
\subset
\begin{cases}
H^{k,q}(\RR^d), &\hbox{for } 0 \leq \sigma p < d \text{ and } q \in [p, p^*],
\\
H^{k,q}(\RR^d), &\hbox{for } \sigma p = d \text{ and } q \in [p, \infty),
\\
C_b^k(\RR^d), &\hbox{for } \sigma p > d,
\end{cases}
\end{equation}
where $p^* := dp/(d-\sigma p)$.
\end{cor}

\subsection{Continuous embeddings for Slobodeckij spaces}
\label{subsec:Continuous_embeddings_Slobodeckij_spaces}
For the Slobodeckij spaces on Euclidean space, Maz'ya and Di Nezza, Palatucci, and Valdinoci provide the following embedding theorem.

\begin{thm}[Continuous embeddings for Slobodeckij spaces on Euclidean space for subcritical exponents]
\label{thm:DiNezza_Palatucci_Valdinoci_6-5}
\cite[Theorem 6.5]{DiNezza_Palatucci_Valdinoci_2012}
\cite[Theorem 10.2.1]{Mazya_2011}
Let $d \geq 1$ be an integer and $p\geq 1$. Then there is a constant $C = C(d,p) \in [1,\infty)$ with the following significance. Let $0 < s < 1$ be such that $sp < d$. If $u \in W_0^{s,p}(\RR^d)$, then
\begin{equation}
\|u\|_{L^{p^*}(\RR^d)}^p \leq c(d,p)\frac{s(1-s)}{(d-sp)^{p-1}} [u]_{W^{s,p}(\RR^d)}^p,
\end{equation}
where $p^* := dp/(d-sp)$ and $c(d,p)$ is a function of $d$ and $p$. In particular, there are continuous embeddings,
\[
W^{s,p}(\RR^d) \subset L^q(\RR^d), \quad\forall\, q \in [p, p^*].
\]
\end{thm}

The preceding embedding results extend to the case of suitably regular domains, as discussed by Di Nezza, Palatucci, and Valdinoci \cite{DiNezza_Palatucci_Valdinoci_2012}. For any $s \in (0, 1)$ and any $p \in [1,\infty)$, one says that an open set $\Omega \subset \RR^d$ is an \emph{extension domain} for $W^{s,p}$ if there exists a positive constant $C = C(d, p, s, \Omega)$ with the following significance: for every function $u \in W^{s,p}(\Omega)$ there exists $\tilde u \in W^{s,p}(\RR^d)$ with $\tilde u(x) = u(x)$ for all $x \in \Omega$ and
$\|\tilde u\|_{W^{s,p}(\RR^d)} \leq C\|u\|_{W^{s,p}(\Omega)}$. We now recall the

\begin{thm}[Existence of extension domains for Slobodeckij spaces]
\label{thm:DiNezza_Palatucci_Valdinoci_5-4}
Let $d \geq 1$ be an integer, $p \in [1,\infty)$, and $s \in (0, 1)$, and $\Omega \subset \RR^d$ be an open set of class $C^{0,1}$ with bounded boundary. Then there is a constant $C = C(d,p,s,\Omega) \in [1,\infty)$ with the following significance. For any $u \in W^{s,p}(\Omega)$, there exists $\tilde u \in W^{s,p}(\RR^d)$ such that $\tilde u \restriction \Omega = u$ and
\[
\|\tilde u\|_{W^{s,p}(\RR^d)} \leq C\|u\|_{W^{s,p}(\Omega)},
\]
and there is a continuous embedding,
\[
W^{s,p}(\Omega) \subset W^{s,p}(\RR^d).
\]
\end{thm}

\begin{thm}[Continuous embeddings for Slobodeckij spaces on domains for subcritical exponents]
\label{thm:DiNezza_Palatucci_Valdinoci_6-7}
\cite[Theorem 6.7]{DiNezza_Palatucci_Valdinoci_2012}
Let $d \geq 1$ be an integer, $p \in [1,\infty)$, and $s \in (0, 1)$ be such that $sp < d$. Let $\Omega \subset \RR^d$ be an extension domain for $W^{s,p}$. Then there exists a positive constant $C = C(d,p,s,\Omega) \in [1,\infty)$ such that, for any $u \in W^{s,p}(\Omega)$ and $q \in [p, p^*]$, we have
\begin{equation}
\label{eq:DiNezza_Palatucci_Valdinoci_6-27}
\|u\|_{L^q(\Omega)} \leq C\|u\|_{W^{s,p}(\Omega)},
\end{equation}
and the following embeddings are continuous:
\[
W^{s,p}(\Omega) \subset L^q(\Omega), \quad\forall\, q \in [p, p^*].
\]
If, in addition, $\Omega$ is bounded, then the preceding embeddings are continuous for any $q \in [1, p^*]$.
\end{thm}

Analogues of Theorems \ref{thm:DiNezza_Palatucci_Valdinoci_6-5} and \ref{thm:DiNezza_Palatucci_Valdinoci_6-7} hold when $sp=d$:

\begin{thm}[Continuous embeddings for Slobodeckij spaces on Euclidean space for the critical exponent]
\label{thm:DiNezza_Palatucci_Valdinoci_6-9}
\cite[Theorem 6.9]{DiNezza_Palatucci_Valdinoci_2012}
Let $d \geq 1$ be an integer, and $p \in [1,\infty)$, and $0 < s < 1$ be such that $sp = d$. Then there is a constant $C = C(d,p,s) \in [1,\infty)$ with the following significance. If $q \in [p, \infty)$ and $u:\RR^d \to \RR$ is measurable and compactly supported, then
\begin{equation}
\label{eq:DiNezza_Palatucci_Valdinoci_6-30}
\|u\|_{L^q(\RR^d)} \leq C\|u\|_{W^{s,p}(\RR^d)},
\end{equation}
and there is a continuous embedding,
\[
W^{s,p}(\RR^d) \subset L^q(\RR^d).
\]
\end{thm}

\begin{thm}[Continuous embeddings for Slobodeckij spaces on domains for the critical exponent]
\label{thm:DiNezza_Palatucci_Valdinoci_6-10}
\cite[Theorem 6.10]{DiNezza_Palatucci_Valdinoci_2012}
Let $d \geq 1$ be an integer, $p \in [1,\infty)$, and $0 < s < 1$ be such that $sp = d$. Let $\Omega \subset \RR^d$ be an extension domain for $W^{s,p}$. Then there is a constant $C = C(d,p,s,\Omega) \in [1,\infty)$ with the following significance. If $q \in [p, \infty)$ and $u \in W^{s,p}(\Omega)$, then
\begin{equation}
\label{eq:DiNezza_Palatucci_Valdinoci_6-31}
\|u\|_{L^q(\Omega)} \leq C\|u\|_{W^{s,p}(\Omega)},
\end{equation}
and there is a continuous embedding,
\[
W^{s,p}(\Omega) \subset L^q(\Omega).
\]
If, in addition, $\Omega$ is bounded, then the preceding embedding is continuous for any $q \in [1, \infty)$.
\end{thm}

Finally, in the case of supercritical exponents, one has the

\begin{thm}[Continuous embeddings for Slobodeckij spaces on domains for supercritical exponents]
\label{thm:DiNezza_Palatucci_Valdinoci_8-2}
\cite[Theorem 8.2]{DiNezza_Palatucci_Valdinoci_2012}
Let $d \geq 1$ be an integer, $p \in [1,\infty)$, and $s \in (0, 1)$ be such that $sp > d$. Let $\Omega \subset \RR^d$ be an extension domain for $W^{s,p}$ with no external cusps. Then there is a constant $C = C(d,p,s,\Omega) \in [1,\infty)$ with the following significance. If $u \in L^p(\Omega)$, then
\begin{equation}
\label{eq:DiNezza_Palatucci_Valdinoci_8-2}
\|u\|_{C^\alpha(\bar\Omega)} \leq C\|u\|_{W^{s,p}(\Omega)}.
\end{equation}
with $\alpha := (sp - d)/p \in (0,1)$.
\end{thm}

\section{Compact embeddings}
\label{sec:Compact_embeddings}
Analogues of the classical Rellich-Kondrachov Embedding \cite[Theorem 6.3]{AdamsFournier} hold for fractional order Sobolev spaces. We note that Theorems \ref{thm:Lions-Magenes_1-16-1}, \ref{thm:Taylor_PDE_equation_13-6-9}, and \ref{thm:DiNezza_Palatucci_Valdinoci_corollary_7-2} are proved directly, without appeal to interpolation theory for Hilbert or Banach spaces \cite{Bergh_Lofstrom_1976, Triebel_interpolation_theory_1978}, where the question of whether one can infer compactness of embeddings is known to be delicate \cite{Cwikel_Janson_2006}, unlike the simpler question of continuity.

\subsection{Compact embeddings for $H^s$ spaces}
\label{subsec:Compact_embeddings_Hs_spaces}
Lions and Magenes provide the

\begin{thm}[Compact embedding for $H^s$ spaces on domains]
\label{thm:Lions-Magenes_1-16-1}
\cite[Theorem 1.16.1]{Lions_Magenes_v1_1972}
Let $d \geq 1$ be an integer and $\Omega \subset \RR^d$ be a bounded domain with a $C^\infty$ boundary \cite[Equations (1.7.10) and (1.7.11)]{Lions_Magenes_v1_1972}. If $s, t \in \RR$ obey $-\infty < s < t < \infty$, then the following embedding is compact:
\[
H^t(\Omega) \Subset H^s(\Omega).
\]
\end{thm}

Haroske and Triebel \cite[Theorems 4.17 and 7.8]{Haroske_Triebel_2008} and Taylor \cite[Proposition 4.4]{Taylor_PDE1} establish Theorem \ref{thm:Lions-Magenes_1-16-1} for the case $0 \leq s < t < \infty$. The case $-\infty < s < t < 0$ then follows by duality and \cite[Theorem 6.4]{Brezis}, which asserts that if $T:E\to F$ is a compact linear map of Banach spaces, then the adjoint operator $T^*:F'\to E'$ is also compact. The case $-\infty < s < 0 < t < \infty$ is an immediate consequence of either one of the preceding compact embedding results.

\subsection{Compact embeddings for Bessel potential spaces}
\label{subsec:Compact_embeddings_Bessel_potential_spaces}
Taylor provides the

\begin{thm}[Compact embedding for Bessel potential spaces on domains]
\label{thm:Taylor_PDE_equation_13-6-9}
\cite[Equation (13.6.9)]{Taylor_PDE3}
Let $d \geq 1$ be an integer and $\Omega \subset \RR^d$ be a bounded domain with a $C^\infty$ boundary. If $p \in (1,\infty)$ and $s, t \in \RR$ obey $0 \leq s < t < \infty$, then the following embedding is compact:
\[
H^{p,t}(\Omega) \Subset H^{p,s}(\Omega).
\]
\end{thm}

Just as in Section \ref{subsec:Compact_embeddings_Hs_spaces}, Theorem \ref{thm:Taylor_PDE_equation_13-6-9} extends to the cases of $-\infty < s < t < 0$ and $-\infty < s < 0 < t < \infty$ by duality.

\subsection{Compact embeddings for Slobodeckij spaces}
\label{subsec:Compact_embeddings_Slobodeckij_spaces}
For the Slobodeckij spaces, $W^{s,p}(\Omega)$ with  $p \in [1, \infty)$ and $s \in (0,1)$, Di Nezza, Palatucci, Valdinoci provide the

\begin{thm}[Compact embedding for Slobodeckij spaces]
\label{thm:DiNezza_Palatucci_Valdinoci_corollary_7-2}
\cite[Corollary 7.2]{DiNezza_Palatucci_Valdinoci_2012}
Let $d \geq 1$ be an integer and $p \in [1, \infty)$ and $s \in (0,1)$ obey $sp < d$. Let $p^* := dp/(d-sp)$ denote the critical exponent and $q \in [1, p^*)$. If $\Omega \subset \RR^d$ is a bounded extension domain for $W^{s,p}$, then the following embedding is compact:
\[
W^{s,p}(\Omega) \Subset L^q(\Omega).
\]
\end{thm}

Theorem \ref{thm:DiNezza_Palatucci_Valdinoci_corollary_7-2} immediately yields the

\begin{cor}[Compact embedding for Slobodeckij spaces on domains]
\label{cor:DiNezza_Palatucci_Valdinoci_corollary_7-2}
Let $d \geq 1$ be an integer and $p \in [1, \infty)$ and $\sigma > 0$ obey $\sigma p < d$. Let $p^* := dp/(d-\sigma p)$ denote the critical exponent and $q \in [1, p^*)$. Let $k \geq 0$ be an integer. If $\Omega \subset \RR^d$ is a bounded extension domain for $W^{k + \sigma,p}$, then the following embedding is compact:
\[
W^{k + \sigma,p}(\Omega) \Subset W^{k,q}(\Omega).
\]
\end{cor}

\section[Fractional order Sobolev spaces of sections of vector bundles]{Fractional order Sobolev spaces of sections of vector bundles over closed manifolds}
\label{subsec:Fractional_order_Sobolev_spaces_sections_vector_bundles_closed_manifolds}
The definitions of Bessel potential spaces in Section \ref{subsec:Bessel_potential_spaces} or Slobodeckij spaces in Section \ref{subsec:Slobodeckij_spaces}, over Euclidean space or domains, together with the associated continuous and compact embedding results, extend in a standard way using a partition of unity to the case of spaces of functions over closed, Riemannian, smooth manifolds, and hence to spaces of sections of vector bundles over such manifolds. This is described in Section \ref{subsec:Continuous_embeddings_Bessel_potential_spaces_complex_interpolation}, following Taylor \cite[Section 4.3]{Taylor_PDE1}, for the case of Bessel potential spaces.

In the case of functions over closed Riemannian smooth manifolds, one also has the option of defining Bessel potential spaces directly in terms of the fractional powers, $(\Delta+1)^{s/2}$ for $s\in \RR$, where $\Delta$ is the Laplace operator on $C^\infty(X)$ defined by the Riemannian metric on $TX$. More generally, if $E$ is a Riemannian vector bundle over $X$ with a metric connection $A$, one may define Bessel potential spaces directly in terms of the fractional powers, $(\nabla_A^*\nabla_A+1)^{s/2}$ for $s\in \RR$, where $\nabla_A^*\nabla_A$ is the connection Laplace operator defined by the connection $A$ on $E$ and Levi-Civita connection on $TX$. When $s$ is an integer, these fractional order Sobolev spaces agree with their standard definitions in terms of covariant derivatives, as in Aubin \cite[Section 2.1]{Aubin_1998}.

\chapter[Fredholm properties of elliptic partial differential operators]{Fredholm properties of elliptic partial differential operators on Sobolev and H\"older spaces}
\label{chap:Fredholm_properties_elliptic_PDEs_Sobolev_Holder_spaces}
We prove two corollaries of \cite[Lemma 1.4.5]{Gilkey2} for the Fredholm property of elliptic pseudo-differential operators, $P:H^{s+m}(X;V) \to H^s(X;W)$, where $V$ and $W$ are vector bundles over a closed, smooth manifold, $X$, and $P:C^\infty(X;V)\to C^\infty(X;W)$ is an elliptic pseudo-differential operator of order $m\in \RR$ \cite[Sections 1.2 and 1.3]{Gilkey2} and $H^s(X;V) = W^{s,2}(X;V)$ is the Sobolev space of order $s\in\RR$ (see Section \ref{subsec:Fractional_order_Sobolev_spaces} for a survey of definitions of fractional-order Sobolev spaces). The conclusions of Lemmata \ref{lem:Gilkey_1-4-5_Sobolev} and \ref{lem:Gilkey_1-4-5_Holder} are widely assumed without comment (for example, based on the results for $H^s$ spaces in Gilkey \cite{Gilkey2} or H\"ormander \cite{Hormander_v3}), but they are not immediate consequences and require justification.

\section[Fredholm properties of elliptic operators on Sobolev spaces]{Fredholm properties of elliptic partial differential operators on Sobolev spaces}
\label{sec:Fredholm_properties_elliptic_PDEs_Sobolev_spaces}
We have the following corollary of a well-known result \cite[Lemma 1.4.5]{Gilkey2} for the Fredholm property of elliptic pseudo-differential operators on $H^s$ spaces. Compare \cite[Theorem 5.7.2]{Volpert1}.

\begin{lem}[Fredholm property for elliptic partial differential operators on Sobolev spaces]
\label{lem:Gilkey_1-4-5_Sobolev}
Let $V$ and $W$ be finite-rank vector bundles over a closed, smooth manifold, $X$, and $k\geq 0$ an integer, and $p \in (1,\infty)$. If $P:C^\infty(X;V)\to C^\infty(X;W)$ is an elliptic partial differential operator of order $m\geq 0$, then $P:W^{k+m,p}(X;V) \to W^{k,p}(X;W)$ is a Fredholm operator with index
\begin{align*}
\Ind(P) &= \dim\Ker\left(P:C^\infty(X;V)\to C^\infty(X;W)\right)
\\
&\quad - \dim\Ker\left(P^*:C^\infty(X;W)\to C^\infty(X;V)\right),
\end{align*}
where $P^*$ is the formal ($L^2$) adjoint of $P$.
\end{lem}

\begin{proof}
We use \cite[Lemma 1.3.6]{Gilkey2} to find a pseudo-differential operator $S:C^\infty(X;W) \to C^\infty(X;V)$ of order $-m$ so that
\[
SP - I \in \Psi^{-\infty}(X;V) \quad\text{and}\quad PS - I \in \Psi^{-\infty}(X;W),
\]
where $\Psi^{-\infty}(X;V)$ is the vector space of infinitely smoothing pseudo-differential operators \cite[Sections 1.2 and 1.3]{Gilkey2}.

The operator $P:W^{k+m,p}(X;V) \to W^{k,p}(X;W)$ is continuous since $P$ is an elliptic partial differential operator of order $m \geq 0$ and by definition
of the Sobolev space $W^{k+m,p}(X;V)$. Combining the preceding observation with the \apriori elliptic estimate (see Theorem \ref{thm:Krylov_Sobolev_lectures_8-5-3}, for example),
\[
\|v\|_{W^{k+m,p}(X)} \leq C\left(\|v\|_{W^{k,p}(X)} + \|Pv\|_{W^{k,p}(X)}\right),
\]
implies that the expression $\|v\|_{W^{k,p}(X)} + \|Pv\|_{W^{k,p}(X)}$ defines a norm for $v \in C^\infty(X;V)$ which is equivalent to the standard norm on $W^{k+m,p}(X;V)$. Furthermore, the operators
\begin{align*}
SP - I:W^{k,p}(X;V) &\to W^{k+1,p}(X;V),
\\
PS - I:W^{k,p}(X;W) &\to W^{k+1,p}(X;W),
\end{align*}
are continuous because the operators
\begin{align*}
SP - I: H^s(X;V) &\to H^{s+t}(X;V),
\\
PS - I: H^s(X;W) &\to H^{s+t}(X;W),
\end{align*}
are continuous by \cite[Lemma 1.3.5]{Gilkey2}, for any $s\in\RR$ and $t\geq 0$, and the Sobolev Embedding Theorem for Sobolev spaces of fractional order (for example, see \cite[Theorem 4.12]{AdamsFournier} for Sobolev spaces of non-negative integer order and Section \ref{sec:Continuous_embeddings} for Sobolev spaces of fractional order), which provides continuous embeddings, $W^{k,p}(X;V) \subset H^s(X;V)$ for sufficiently small $s < s_0(d,k,p)$ and $H^{s+t}(X;W) \subset W^{k+1,p}(X;W)$ for sufficiently large $t > t_0(d,k,p)$, where $d$ is the dimension of $X$. (The finite constants $s_0(d,k,p)$ and $t_0(d,k,p)$ can of course be determined explicitly from the full statement of the Sobolev Embedding Theorem, but their precise values are unimportant here.) Thus,
\begin{align*}
\|Sw\|_{W^{k+m,p}(X;W)} &\leq C\left(\|PSw\|_{W^{k,p}(X;W)} + \|w\|_{W^{k,p}(X;W)}\right)
\\
&\leq C\left(\|(PS-I)w\|_{W^{k,p}(X;W)} + \|w\|_{W^{k,p}(X;W)}\right)
\\
&\leq C\|w\|_{W^{k,p}(X;W)} \quad\text{(by continuity of $PS - I$ on $W^{k,p}(X;W)$),}
\end{align*}
and so the operator $S:W^{k,p}(X;V) \to W^{k+m,p}(X;W)$ is continuous.

The embeddings $W^{k+1,p}(X;V) \subset W^{k,p}(X;V)$ and $W^{k+1,p}(X;W) \subset W^{k,p}(X;W)$ are compact by the Rellich-Kondrachov Theorem (see \cite[Theorem 6.3]{AdamsFournier}) and so the operators (now viewed as compositions with compact embeddings)
\begin{align*}
SP - I:W^{k,p}(X;V) &\to W^{k,p}(X;V),
\\
PS - I:W^{k,p}(X;W) &\to W^{k,p}(X;W),
\end{align*}
are compact by \cite[Proposition 6.3]{Brezis}. Hence, $P:W^{k+m,p}(X;V) \to W^{k,p}(X;W)$ is Fredholm by \cite[Section 1.4.2, Definition]{Gilkey2}.

Because $P:W^{k+m,p}(X;V) \to W^{k,p}(X;W)$ is Fredholm, its index may be computed by \cite[Lemma 4.38]{Abramovich_Aliprantis_2002}
\begin{align*}
\Ind(P) &= \dim\Ker\left(P:W^{k+m,p}(X;V) \to W^{k,p}(X;W)\right)
\\
&\quad - \dim\Ker\left(P^*:W^{-k,p'}(X;W) \to W^{-k-m,p'}(X;V)\right),
\end{align*}
where $p' \in (1, \infty)$ is the dual exponent defined by $1/p+1/p'=1$ and we have appealed to the characterization of Banach duals of Sobolev spaces \cite[Sections 3.7 to 3.14]{AdamsFournier}. If $w \in (\Ker P^*)\cap W^{-k,p'}(X;W)$, then $w \in (\Ker P^*)\cap H^s(X;W)$ for all $s < s_0(d,k,p)$, where $s_0$ is given by the Sobolev Embedding Theorem, and consequently (since the Banach space dual $P^*$ is defined by the realization of the formal adjoint and $P^*$ is thus an elliptic partial differential operator of order $m$) we see that $w \in (\Ker P^*)\cap C^\infty(X;W)$ by elliptic regularity \cite[Lemma 1.3.2 and Section 1.3.5]{Gilkey2}. Of course, if $v \in (\Ker P)\cap W^{k+m,p}(X;V)$, then $v \in (\Ker P)\cap C^\infty(X;V)$ by the same argument. This yields the stated formula for the index of $P:W^{k+m,p}(X;V) \to W^{k,p}(X;W)$.
\end{proof}

\section[Fredholm properties of elliptic operators on H\"older spaces]{Fredholm properties of elliptic partial differential operators on H\"older spaces}
\label{sec:Fredholm_properties_elliptic_PDEs_Holder_spaces}
We now turn to the corresponding assertions for elliptic partial differential operators on H\"older spaces, but we specialize to the case of second-order operators since there fewer suitable references for \apriori estimates or regularity theory in the case of elliptic partial differential operators of arbitrary order on H\"older spaces. Compare \cite[Theorem 5.8.3]{Volpert1}.

\begin{lem}[Fredholm property for second-order elliptic partial differential operators on H\"older spaces]
\label{lem:Gilkey_1-4-5_Holder}
Let $V$ and $W$ be finite-rank vector bundles over a closed, smooth manifold, $X$, and $k\geq 0$ an integer and $\alpha\in (0,1)$. If $P:C^\infty(X;V)\to C^\infty(X;W)$ is a second-order, elliptic partial differential operator, then $P:C^{k+2,\alpha}(X;V) \to C^{k,\alpha}(X;W)$ is a Fredholm operator with index
\begin{align*}
\Ind(P) &= \dim\Ker\left(P:C^\infty(X;V)\to C^\infty(X;W)\right)
\\
&\quad - \dim\Ker\left(P^*:C^\infty(X;W)\to C^\infty(X;V)\right),
\end{align*}
where $P^*$ is the formal ($L^2$) adjoint of $P$.
\end{lem}

\begin{proof}
We again use \cite[Lemma 1.3.6]{Gilkey2} to find a pseudo-differential operator $S:C^\infty(X;W) \to C^\infty(X;V)$ of order $-m$ so that
\[
SP - I \in \Psi^{-\infty}(X;V) \quad\text{and}\quad PS - I \in \Psi^{-\infty}(X;W),
\]
where $\Psi^{-\infty}(X;V)$ is the vector space of infinitely smoothing pseudo-differential operators \cite[Sections 1.2 and 1.3]{Gilkey2}.

The operator $P:C^{k+2,\alpha}(X;V) \to C^{k,\alpha}(X;W)$ is continuous since $P$ is a second-order elliptic partial differential operator and by definition of the H\"older space $C^{k+2,\alpha}(X;V)$. Combining this continuity property with the \apriori elliptic estimate (see \cite[Theorem 6.2, Corollary 6.3, or Theorem 6.6]{GilbargTrudinger}, for example, in the case of scalar operators),
\[
\|v\|_{C^{k+2,\alpha}(X)} \leq C\left(\|v\|_{C^{k,\alpha}(X)} + \|Pv\|_{C^{k,\alpha}(X)}\right),
\]
implies that the expression $\|v\|_{C^{k,\alpha}(X)} + \|Pv\|_{C^{k,\alpha}(X)}$ defines a norm for $v \in C^\infty(X;V)$ which is equivalent to the standard norm on $C^{k+2,\alpha}(X;V)$. Furthermore, the operators
\begin{align*}
SP - I:C^{k,\alpha}(X;V) &\to C^{k,\beta}(X;V),
\\
PS - I:C^{k,\alpha}(X;W) &\to C^{k,\beta}(X;W),
\end{align*}
are continuous for any $\beta$ obeying $\alpha\leq\beta\leq 1$ because the operators
\begin{align*}
SP - I: H^s(X;V) &\to H^{s+t}(X;V),
\\
PS - I: H^s(X;W) &\to H^{s+t}(X;W),
\end{align*}
are continuous by \cite[Lemma 1.3.5]{Gilkey2}, for any $s\in\RR$ and $t\geq 0$, and the Sobolev Embedding Theorem for Sobolev spaces of fractional order (for example, see \cite[Theorem 4.12]{AdamsFournier} for Sobolev spaces of non-negative integer order and Section \ref{sec:Continuous_embeddings}
for Sobolev spaces of fractional order), which provides continuous embeddings, $C^{k,\alpha}(X;V) \subset H^s(X;V)$ for sufficiently small $s < s_0(d,k,\alpha)$ and $H^{s+t}(X;W) \subset C^{k,\beta}(X;W)$ for sufficiently large $t > t_0(d,k,\beta)$, where $d$ is the dimension of $X$. Thus,
\begin{align*}
\|Sw\|_{C^{k+2,\alpha}(X;W)} &\leq C\left(\|PSw\|_{C^{k,\alpha}(X;W)} + \|w\|_{C^{k,\alpha}(X;W)}\right)
\\
&\leq C\left(\|(PS-I)w\|_{C^{k,\alpha}(X;W)} + \|w\|_{C^{k,\alpha}(X;W)}\right)
\\
&\leq C\|w\|_{C^{k,\alpha}(X;W)} \quad\text{(by continuity of $PS - I$ on $C^{k,\alpha}(X;W)$),}
\end{align*}
and so the operator $S:C^{k,\alpha}(X;V) \to C^{k+2,\alpha}(X;W)$ is continuous.

When $\beta>\alpha$, the embeddings $C^{k,\beta}(X;V) \subset C^{k,\alpha}(X;V)$ and $C^{k,\beta}(X;W) \subset C^{k,\alpha}(X;W)$ are compact as a consequence of the Arzel\`a-Ascoli Theorem (see \cite[Theorem 1.3.4]{AdamsFournier}) and so the operators (now viewed as compositions with compact embeddings)
\begin{align*}
SP - I:C^{k,\alpha}(X;V) &\to C^{k,\alpha}(X;V),
\\
PS - I:C^{k,\alpha}(X;W) &\to C^{k,\alpha}(X;W),
\end{align*}
are compact by \cite[Proposition 6.3]{Brezis}. Hence, $P:C^{k+2,\alpha}(X;V) \to C^{k,\alpha}(X;W)$ is Fredholm by \cite[Section 1.4.2, Definition]{Gilkey2}.

Because $P:C^{k+2,\alpha}(X;V) \to C^{k,\alpha}(X;W)$ is Fredholm, its index may be computed by \cite[Lemma 4.38]{Abramovich_Aliprantis_2002}
\begin{align*}
\Ind(P) &= \dim\Ker\left(P:C^{k+2,\alpha}(X;V) \to C^{k,\alpha}(X;W)\right)
\\
&\quad - \dim\Ker\left(P^*:\left(C^{k,\alpha}(X;W)\right)^* \to \left(C^{k+2,\alpha}(X;V)\right)^*\right).
\end{align*}
If $v \in (\Ker P)\cap C^{k+2,\alpha}(X;V)$, then $v \in (\Ker P)\cap H^s(X;V)$ for any $s < s_0(d,k,\alpha)$ small enough (and determined by the Sobolev Embedding Theorem) that $C^{k+2,\alpha}(X;V) \subset H^s(X;V)$ is a continuous embedding and consequently $v \in (\Ker P)\cap C^\infty(X;V)$ by elliptic regularity \cite[Lemma 1.3.2 and Section 1.3.5]{Gilkey2}. Similarly, if $w \in (\Ker P^*)\cap (C^{k,\alpha}(X;W))^*$, then we observe that $H^t(X;W) \subset C^{k,\alpha}(X;W)$ is a continuous (and dense) embedding by the Sobolev Embedding Theorem for any $t > t_0(d,k,\alpha)$ large enough and thus $(C^{k,\alpha}(X;W))^* \subset (H^t(X;W))^* = H^{-t}(X;W)$ is a continuous embedding and so $w \in (\Ker P^*)\cap H^{-t}(X;W)$. Therefore, $w \in (\Ker P^*)\cap C^\infty(X;V)$ by elliptic regularity. This yields the stated formula for the index of $P:C^{k+2,\alpha}(X;V) \to C^{k,\alpha}(X;W)$.
\end{proof}

\backmatter

%
%

\bibliography{master,mfpde}

\def\cprime{$'$} \def\cprime{$'$}
  \def\ocirc#1{\ifmmode\setbox0=\hbox{$#1$}\dimen0=\ht0 \advance\dimen0
  by1pt\rlap{\hbox to\wd0{\hss\raise\dimen0
  \hbox{\hskip.2em$\scriptscriptstyle\circ$}\hss}}#1\else {\accent"17 #1}\fi}
  \def\cprime{$'$} \def\cprime{$'$} \def\cprime{$'$} \def\cprime{$'$}
  \def\polhk#1{\setbox0=\hbox{#1}{\ooalign{\hidewidth
  \lower1.5ex\hbox{`}\hidewidth\crcr\unhbox0}}} \def\cprime{$'$}
  \def\cprime{$'$} \def\cprime{$'$}
  \def\lfhook#1{\setbox0=\hbox{#1}{\ooalign{\hidewidth
  \lower1.5ex\hbox{'}\hidewidth\crcr\unhbox0}}} \def\cprime{$'$}
  \def\cprime{$'$} \def\cprime{$'$} \def\cprime{$'$} \def\cprime{$'$}
\providecommand{\bysame}{\leavevmode\hbox to3em{\hrulefill}\thinspace}
\providecommand{\MR}{\relax\ifhmode\unskip\space\fi MR }
\providecommand{\MRhref}[2]{%
  \href{http://www.ams.org/mathscinet-getitem?mr=#1}{#2}
}
\providecommand{\href}[2]{#2}
\begin{thebibliography}{100}

\bibitem{AMR}
R.~Abraham, J.~E. Marsden, and T.~Ratiu, \emph{Manifolds, tensor analysis, and
  applications}, second ed., Springer, New York, 1988. \MR{960687 (89f:58001)}

\bibitem{Abramovich_Aliprantis_2002}
Y.~A. Abramovich and C.~D. Aliprantis, \emph{An invitation to operator theory},
  Graduate Studies in Mathematics, vol.~50, American Mathematical Society,
  Providence, RI, 2002. \MR{1921782 (2003h:47072)}

\bibitem{Ache_2011arxiv}
A.~G. Ache, \emph{On the uniqueness of asymptotic limits of the {R}icci flow},
  arXiv:1211.3387.

\bibitem{Adams}
R.~A. Adams, \emph{Sobolev spaces}, Academic Press [A subsidiary of Harcourt
  Brace Jovanovich, Publishers], New York-London, 1975, Pure and Applied
  Mathematics, Vol. 65. \MR{0450957 (56 \#9247)}

\bibitem{AdamsFournier}
R.~A. Adams and J.~J.~F. Fournier, \emph{Sobolev spaces}, second ed.,
  Elsevier/Academic Press, Amsterdam, 2003. \MR{2424078 (2009e:46025)}

\bibitem{Agmon_1962}
S.~Agmon, \emph{On the eigenfunctions and on the eigenvalues of general
  elliptic boundary value problems}, Comm. Pure Appl. Math. \textbf{15} (1962),
  119--147. \MR{0147774 (26 \#5288)}

\bibitem{AgmonLecturesEllipticBVP}
\bysame, \emph{Lectures on elliptic boundary value problems}, AMS Chelsea
  Publishing, Providence, RI, 2010, Revised edition of the 1965 original.
  \MR{2589244 (2011c:35004)}

\bibitem{AgmonDouglisNirenberg1}
S.~Agmon, A.~Douglis, and L.~Nirenberg, \emph{Estimates near the boundary for
  solutions of elliptic partial differential equations satisfying general
  boundary conditions. {I}}, Comm. Pure Appl. Math. \textbf{12} (1959),
  623--727.

\bibitem{AgmonDouglisNirenberg2}
\bysame, \emph{Estimates near the boundary for solutions of elliptic partial
  differential equations satisfying general boundary conditions. {II}}, Comm.
  Pure Appl. Math. \textbf{17} (1964), 35--92.

\bibitem{Amann_1983}
H.~Amann, \emph{Dual semigroups and second order linear elliptic boundary value
  problems}, Israel J. Math. \textbf{45} (1983), 225--254. \MR{719122
  (85i:35043)}

\bibitem{Amann_1995}
\bysame, \emph{Linear and quasilinear parabolic problems. {V}ol. {I}},
  Monographs in Mathematics, vol.~89, Birkh\"auser Boston Inc., Boston, MA,
  1995, Abstract linear theory. \MR{1345385 (96g:34088)}

\bibitem{Ambrosetti_Rabinowitz_1973}
A.~Ambrosetti and P.~H. Rabinowitz, \emph{Dual variational methods in critical
  point theory and applications}, J. Functional Analysis \textbf{14} (1973),
  349--381. \MR{0370183 (51 \#6412)}

\bibitem{Andrews_Baker_2010}
B.~Andrews and C.~Baker, \emph{Mean curvature flow of pinched submanifolds to
  spheres}, J. Differential Geom. \textbf{85} (2010), no.~3, 357--395.
  \MR{2739807 (2012a:53122)}

\bibitem{Andrews_Hopper_2011}
B.~Andrews and C.~Hopper, \emph{The {R}icci flow in {R}iemannian geometry},
  Lecture Notes in Mathematics, vol. 2011, Springer, Heidelberg, 2011, A
  complete proof of the differentiable 1/4-pinching sphere theorem. \MR{2760593
  (2012d:53208)}

\bibitem{Angiuli_Pallara_Paronetto_2010}
L.~Angiuli, D.~Pallara, and F.~Paronetto, \emph{Analytic semigroups generated
  in {$L^1(\Omega)$} by second order elliptic operators via duality methods},
  Semigroup Forum \textbf{80} (2010), 255--271. \MR{2601763 (2011i:35037)}

\bibitem{Aronszajn_Mulla_Szeptycki_1963}
N.~Aronszajn, F.~Mulla, and P.~Szeptycki, \emph{On spaces of potentials
  connected with {$L^{p}$} classes}, Ann. Inst. Fourier (Grenoble) \textbf{13}
  (1963), 211--306. \MR{0180846 (31 \#5076)}

\bibitem{Aronszajn_Smith_1961}
N.~Aronszajn and K.~T. Smith, \emph{Theory of {B}essel potentials. {I}}, Ann.
  Inst. Fourier (Grenoble) \textbf{11} (1961), 385--475. \MR{0143935 (26
  \#1485)}

\bibitem{Atiyah_Bott_1983}
M.~F. Atiyah and R.~Bott, \emph{The {Y}ang-{M}ills equations over {R}iemann
  surfaces}, Philos. Trans. Roy. Soc. London Ser. A \textbf{308} (1983),
  523--615. \MR{702806 (85k:14006)}

\bibitem{ADHM}
M.~F. Atiyah, N.~J. Hitchin, V.~G. Drinfel{\cprime}d, and Yu.~I. Manin,
  \emph{Construction of instantons}, Phys. Lett. A \textbf{65} (1978),
  185--187. \MR{598562 (82g:81049)}

\bibitem{AHS}
M.~F. Atiyah, N.~J. Hitchin, and I.~M. Singer, \emph{Self-duality in
  four-dimensional {R}iemannian geometry}, Proc. Roy. Soc. London Ser. A
  \textbf{362} (1978), no.~1711, 425--461. \MR{506229 (80d:53023)}

\bibitem{Aubin_1976}
T.~Aubin, \emph{\'{E}quations diff\'erentielles non lin\'eaires et probl\`eme
  de {Y}amabe concernant la courbure scalaire}, J. Math. Pures Appl. (9)
  \textbf{55} (1976), 269--296. \MR{0431287 (55 \#4288)}

\bibitem{Aubin}
\bysame, \emph{Nonlinear analysis on manifolds. {M}onge-{A}mp\`ere equations},
  Springer, New York, 1982. \MR{681859 (85j:58002)}

\bibitem{Aubin_1998}
\bysame, \emph{Some nonlinear problems in {R}iemannian geometry}, Springer,
  Berlin, 1998. \MR{1636569 (99i:58001)}

\bibitem{Bahri_1993na}
A.~Bahri, \emph{Another proof of the {Y}amabe conjecture for locally
  conformally flat manifolds}, Nonlinear Anal. \textbf{20} (1993), 1261--1278.
  \MR{1219242 (94e:53033)}

\bibitem{Bahri_1993lnpam}
\bysame, \emph{Proof of the {Y}amabe conjecture, without the positive mass
  theorem, for locally conformally flat manifolds}, Einstein metrics and
  {Y}ang-{M}ills connections ({S}anda, 1990), Lecture Notes in Pure and Appl.
  Math., vol. 145, Dekker, New York, 1993, pp.~1--26. \MR{1215275 (94m:53050)}

\bibitem{BakerThesis}
C.~Baker, \emph{The mean curvature flow of submanifolds of high codimension},
  {Ph.D}. thesis, Australian National University, Canberra, November 2010,
  arXiv:1104.4409.

\bibitem{Banasiak_Arlotti_2006}
J.~Banasiak and L.~Arlotti, \emph{Perturbations of positive semigroups with
  applications}, Springer Monographs in Mathematics, Springer-Verlag London
  Ltd., London, 2006. \MR{2178970 (2006i:47076)}

\bibitem{Barbu_Precupanu_2012}
V.~Barbu and T.~Precupanu, \emph{Convexity and optimization in {B}anach
  spaces}, fourth ed., Springer Monographs in Mathematics, Springer, Dordrecht,
  2012. \MR{3025420}

\bibitem{Bellettini_lecture_notes_on_mean_curvature_flow}
G.~Bellettini, \emph{Lecture notes on mean curvature flow, barriers and
  singular perturbations}, Appunti. Scuola Normale Superiore di Pisa (Nuova
  Serie) [Lecture Notes. Scuola Normale Superiore di Pisa (New Series)],
  vol.~12, Edizioni della Normale, Pisa, 2013. \MR{3155251}

\bibitem{Benjamini_Chavel_Feldman_1996}
I.~Benjamini, I.~Chavel, and E.~A. Feldman, \emph{Heat kernel lower bounds on
  {R}iemannian manifolds using the old ideas of {N}ash}, Proc. London Math.
  Soc. (3) \textbf{72} (1996), 215--240. \MR{1357093 (97c:58150)}

\bibitem{Berger_1977}
M.~Berger, \emph{Nonlinearity and functional analysis}, Academic Press, New
  York, 1977. \MR{0488101 (58 \#7671)}

\bibitem{Berger_Gauduchon_Mazet_1971}
M.~Berger, P.~Gauduchon, and E.~Mazet, \emph{Le spectre d'une vari\'et\'e
  riemannienne}, Lecture Notes in Mathematics, Vol. 194, Springer-Verlag,
  Berlin-New York, 1971. \MR{0282313 (43 \#8025)}

\bibitem{Bergh_Lofstrom_1976}
J.~Bergh and J.~L{\"o}fstr{\"o}m, \emph{Interpolation spaces. {A}n
  introduction}, Springer-Verlag, Berlin-New York, 1976, Grundlehren der
  Mathematischen Wissenschaften, No. 223. \MR{0482275 (58 \#2349)}

\bibitem{BerlineGetzlerVergne}
N.~Berline, E.~Getzler, and M.~Vergne, \emph{Heat kernels and {D}irac
  operators}, Grundlehren Text Editions, Springer-Verlag, Berlin, 2004,
  Corrected reprint of the 1992 original. \MR{2273508 (2007m:58033)}

\bibitem{Bertsch_vanderHout_Hulshof_2011}
M.~Bertsch, R.~van~der Hout, and J.~Hulshof, \emph{Energy concentration for
  2-dimensional radially symmetric equivariant harmonic map heat flows},
  Commun. Contemp. Math. \textbf{13} (2011), 675--695. \MR{2826442
  (2012j:53085)}

\bibitem{Besse_1987}
A.~L. Besse, \emph{Einstein manifolds}, Springer, Berlin, 1987. \MR{867684
  (88f:53087)}

\bibitem{Bethuel_1991}
F.~Bethuel, \emph{The approximation problem for {S}obolev maps between two
  manifolds}, Acta Math. \textbf{167} (1991), 153--206. \MR{1120602
  (92f:58023)}

\bibitem{Biernat_2015}
P.~Biernat, \emph{Non-self-similar blow-up in the heat flow for harmonic maps
  in higher dimensions}, Nonlinearity \textbf{28} (2015), no.~1, 167--185,
  arXiv:1404.2209. \MR{3297131}

\bibitem{BierstoneMilman}
E.~Bierstone and P.~D. Milman, \emph{Semianalytic and subanalytic sets}, Inst.
  Hautes \'Etudes Sci. Publ. Math. (1988), no.~67, 5--42. \MR{972342
  (89k:32011)}

\bibitem{Biswas_Wilkin_2010}
I.~Biswas and G.~Wilkin, \emph{Morse theory for the space of {H}iggs
  {$G$}-bundles}, Geom. Dedicata \textbf{149} (2010), 189--203. \MR{2737688
  (2012c:58023)}

\bibitem{Bizon_Wasserman_2015}
P.~Bizo{\'n} and A.~Wasserman, \emph{Nonexistence of shrinkers for the harmonic
  map flow in higher dimensions}, Int. Math. Res. Not. IMRN (2015), no.~17,
  7757--7762, arXiv:1404.7381. \MR{3403999}

\bibitem{Blatt_2016arxiv}
S.~Blatt, \emph{The gradient flow of {O}'{H}ara's knot energies},
  arXiv:1601.02840.

\bibitem{Blatt_2012jntr}
\bysame, \emph{Boundedness and regularizing effects of {O}'{H}ara's knot
  energies}, J. Knot Theory Ramifications \textbf{21} (2012), no.~1, 1250010,
  9. \MR{2887901}

\bibitem{Blatt_2012cvpde}
\bysame, \emph{The gradient flow of the {M}\"obius energy near local
  minimizers}, Calc. Var. Partial Differential Equations \textbf{43} (2012),
  no.~3-4, 403--439. \MR{2875646}

\bibitem{Boling_Kelleher_Streets_2015arxiv}
J.~Boling, C.~Kelleher, and J.~Streets, \emph{Entropy, stability, and harmonic
  map flow}, arXiv:1506.07567.

\bibitem{Bor_1992}
G.~Bor, \emph{Yang-{M}ills fields which are not self-dual}, Comm. Math. Phys.
  \textbf{145} (1992), 393--410. \MR{1162805 (93e:58030)}

\bibitem{Bor_Montgomery_1990}
G.~Bor and R.~Montgomery, \emph{{${\rm SO}(3)$} invariant {Y}ang-{M}ills fields
  which are not self-dual}, Hamiltonian systems, transformation groups and
  spectral transform methods ({M}ontreal, {PQ}, 1989), Univ. Montr\'eal,
  Montreal, QC, 1990, pp.~191--198. \MR{1110384 (92f:58041)}

\bibitem{Bourguignon_1981}
J.-P. Bourguignon, \emph{Formules de {W}eitzenb\"ock en dimension {$4$}},
  Riemannian geometry in dimension 4 ({P}aris, 1978/1979), Textes Math.,
  vol.~3, CEDIC, Paris, 1981, pp.~308--333. \MR{769143}

\bibitem{Bourguignon_1990}
\bysame, \emph{The ``magic'' of {W}eitzenb\"ock formulas}, Variational methods,
  Proceedings of the Conference on Variational Problems, Paris, 1988
  (H.~Berestycki, J-M. Coron, and I.~Ekeland, eds.), Progress in nonlinear
  differential equations and their applications, vol.~4, Birkh\"auser, Boston,
  MA, 1990, pp.~251--271.

\bibitem{Bourguignon_Lawson_1981}
J-P. Bourguignon and H.~B. Lawson, Jr., \emph{Stability and isolation phenomena
  for {Y}ang-{M}ills fields}, Comm. Math. Phys. \textbf{79} (1981), 189--230.
  \MR{612248 (82g:58026)}

\bibitem{Bourguignon_Lawson_Simons_1979}
J-P. Bourguignon, H.~B. Lawson, Jr., and J.~Simons, \emph{Stability and gap
  phenomena for {Y}ang-{M}ills fields}, Proc. Nat. Acad. Sci. U.S.A.
  \textbf{76} (1979), 1550--1553. \MR{526178 (80h:53028)}

\bibitem{Brendle_2003arxiv}
S.~Brendle, \emph{On the construction of solutions to the {Y}ang-{M}ills
  equations in higher dimensions}, arXiv:math/0302093.

\bibitem{Brendle_2005}
\bysame, \emph{Convergence of the {Y}amabe flow for arbitrary initial energy},
  J. Differential Geom. \textbf{69} (2005), 217--278. \MR{2168505
  (2006e:53119)}

\bibitem{Brendle_2007invent}
\bysame, \emph{Convergence of the {Y}amabe flow in dimension 6 and higher},
  Invent. Math. \textbf{170} (2007), 541--576. \MR{2357502 (2008k:53136)}

\bibitem{Brendle_2008jams}
\bysame, \emph{Blow-up phenomena for the {Y}amabe equation}, J. Amer. Math.
  Soc. \textbf{21} (2008), 951--979. \MR{2425176 (2009m:53084)}

\bibitem{Brendle_2008sdg}
\bysame, \emph{On the conformal scalar curvature equation and related
  problems}, Surveys in differential geometry. {V}ol. {XII}. {G}eometric flows,
  Surv. Differ. Geom., vol.~12, Int. Press, Somerville, MA, 2008, pp.~1--19.
  \MR{2488950 (2009k:53077)}

\bibitem{Brendle_2011jjm}
\bysame, \emph{Evolution equations in {R}iemannian geometry}, Jpn. J. Math.
  \textbf{6} (2011), 45--61. \MR{2835361 (2012j:53086)}

\bibitem{Brendle_Marques_2009}
S.~Brendle and F.~C. Marques, \emph{Blow-up phenomena for the {Y}amabe
  equation. {II}}, J. Differential Geom. \textbf{81} (2009), 225--250.
  \MR{2472174 (2010k:53050)}

\bibitem{Brezis}
H.~Br{\'e}zis, \emph{Functional analysis, {S}obolev spaces and partial
  differential equations}, Universitext, Springer, New York, 2011. \MR{2759829
  (2012a:35002)}

\bibitem{Brezis_Nirenberg_1983}
H.~Br{\'e}zis and L.~Nirenberg, \emph{Positive solutions of nonlinear elliptic
  equations involving critical {S}obolev exponents}, Comm. Pure Appl. Math.
  \textbf{36} (1983), 437--477. \MR{709644 (84h:35059)}

\bibitem{Brezis_Strauss_1973}
H.~Br{\'e}zis and W.~A. Strauss, \emph{Semi-linear second-order elliptic
  equations in {$L^{1}$}}, J. Math. Soc. Japan \textbf{25} (1973), 565--590.
  \MR{0336050 (49 \#826)}

\bibitem{Browder_1961}
F.~E. Browder, \emph{On the spectral theory of elliptic differential operators.
  {I}}, Math. Ann. \textbf{142} (1960/1961), 22--130. \MR{0209909 (35 \#804)}

\bibitem{Buchdahl_1988}
N.~P. Buchdahl, \emph{Hermitian-{E}instein connections and stable vector
  bundles over compact complex surfaces}, Math. Ann. \textbf{280} (1988),
  625--648. \MR{939923 (89d:53092)}

\bibitem{Buck_Orloff_1995}
G.~Buck and J.~Orloff, \emph{A simple energy function for knots}, Topology
  Appl. \textbf{61} (1995), no.~3, 205--214. \MR{1317077 (95k:58037)}

\bibitem{Calderon_1961}
A.-P. Calder{\'o}n, \emph{Lebesgue spaces of differentiable functions and
  distributions}, Partial differential equations (C.~B. Morrey, ed.), Proc.
  {S}ympos. {P}ure {M}ath., vol.~IV, American Mathematical Society, Providence,
  R.I., 1961, pp.~33--49. \MR{0143037 (26 \#603)}

\bibitem{Cannarsa_2014_private}
P.~Cannarsa, \emph{Private communcation}, February 14, 2014.

\bibitem{Cannarsa_Terreni_Vespri_1985}
P.~Cannarsa, B.~Terreni, and V.~Vespri, \emph{Analytic semigroups generated by
  nonvariational elliptic systems of second order under {D}irichlet boundary
  conditions}, J. Math. Anal. Appl. \textbf{112} (1985), 56--103. \MR{812793
  (87e:35034)}

\bibitem{Cannarsa_Vespri_1988}
P.~Cannarsa and V.~Vespri, \emph{Generation of analytic semigroups in the
  {$L^p$} topology by elliptic operators in {${\bf R}^n$}}, Israel J. Math.
  \textbf{61} (1988), 235--255. \MR{941240 (89h:35087)}

\bibitem{Cao_Hamilton_Ilmanen_2004arxiv}
H.-D. Cao, R.S. Hamilton, and T.~Ilmanen, \emph{Gaussian densities and
  stability for some {R}icci solitons}, arXiv:math/0404165.

\bibitem{Cao_Zhu_2006arxiv}
H.-D. Cao and X.-P. Zhu, \emph{Hamilton-{P}erelman's proof of the {P}oincar\'e
  conjecture and the geometrization conjecture}, arXiv:0612069.

\bibitem{Cao_Zhu_2006}
\bysame, \emph{A complete proof of the {P}oincar\'e and geometrization
  conjectures---application of the {H}amilton-{P}erelman theory of the {R}icci
  flow}, Asian J. Math. \textbf{10} (2006), 165--492. \MR{2233789
  (2008d:53090)}

\bibitem{Cao_Zhu_2006erratum}
\bysame, \emph{Erratum to: ``{A} complete proof of the {P}oincar\'e and
  geometrization conjectures---application of the {H}amilton-{P}erelman theory
  of the {R}icci flow'' [{A}sian {J}. {M}ath. {\bf 10} (2006), no. 2,
  165--492]}, Asian J. Math. \textbf{10} (2006), no.~4, 663. \MR{2282358
  (2008d:53091)}

\bibitem{Caps_evolution_equations_scales_banach_spaces}
O.~Caps, \emph{Evolution equations in scales of {B}anach spaces}, Teubner-Texte
  zur Mathematik [Teubner Texts in Mathematics], vol. 140, B. G. Teubner,
  Stuttgart, 2002. \MR{2352809 (2008m:34129)}

\bibitem{Carlotto_Chodosh_Rubinstein_2015}
A.~Carlotto, O.~Chodosh, and Y.~A. Rubinstein, \emph{Slowly converging {Y}amabe
  flows}, Geom. Topol. \textbf{19} (2015), no.~3, 1523--1568, arXiv:1401.3738.
  \MR{3352243}

\bibitem{Chang_Ding_1991}
K-C. Chang and W.~Y. Ding, \emph{A result on the global existence for heat
  flows of harmonic maps from {$D^2$} into {$S^2$}}, Nematics ({O}rsay, 1990)
  (J-M. Coron, J-M. Ghidaglia, and F.~H{\'e}lein, eds.), NATO Adv. Sci. Inst.
  Ser. C Math. Phys. Sci., vol. 332, Kluwer Acad. Publ., Dordrecht, 1991,
  pp.~37--47. \MR{1178085 (94c:58048)}

\bibitem{Chang_Ding_Ye_1992}
K-C. Chang, W.~Y. Ding, and R.~Ye, \emph{Finite-time blow-up of the heat flow
  of harmonic maps from surfaces}, J. Differential Geom. \textbf{36} (1992),
  507--515. \MR{1180392 (93h:58043)}

\bibitem{Chavel}
I.~Chavel, \emph{Eigenvalues in {R}iemannian geometry}, Pure and Applied
  Mathematics, vol. 115, Academic Press Inc., Orlando, FL, 1984. \MR{768584
  (86g:58140)}

\bibitem{Cheeger_Ebin_1975}
J.~Cheeger and D.~G. Ebin, \emph{Comparison theorems in {R}iemannian geometry},
  North-Holland Publishing Co., Amsterdam-Oxford; American Elsevier Publishing
  Co., Inc., New York, 1975, North-Holland Mathematical Library, Vol. 9.
  \MR{0458335 (56 \#16538)}

\bibitem{Chen_Shen_1993}
Y-M. Chen and C-L. Shen, \emph{Evolution of {Y}ang-{M}ills connections},
  Differential geometry ({S}hanghai, 1991), World Sci. Publ., River Edge, NJ,
  1993, pp.~33--41. \MR{1341596 (96f:58039)}

\bibitem{Chen_Shen_1994}
\bysame, \emph{Monotonicity formula and small action regularity for
  {Y}ang-{M}ills flows in higher dimensions}, Calc. Var. Partial Differential
  Equations \textbf{2} (1994), 389--403. \MR{1383915 (98c:58030)}

\bibitem{Chen_Shen_1995}
\bysame, \emph{Evolution problem of {Y}ang-{M}ills flow over {$4$}-dimensional
  manifold}, Variational methods in nonlinear analysis ({E}rice, 1992), Gordon
  and Breach, Basel, 1995, pp.~63--66. \MR{1451148 (99g:58030)}

\bibitem{Chen_Shen_Zhou_2002}
Y-M. Chen, C-L. Shen, and Q.~Zhou, \emph{Asymptotic behavior of {Y}ang-{M}ills
  flow in higher dimensions}, Differential geometry and related topics, World
  Sci. Publ., River Edge, NJ, 2002, pp.~16--38. \MR{2066798 (2005c:53082)}

\bibitem{Chen_Wu_1998}
Y-Z. Chen and L-C. Wu, \emph{Second order elliptic equations and elliptic
  systems}, Translations of Mathematical Monographs, vol. 174, American
  Mathematical Society, Providence, RI, 1998. \MR{1616087 (99i:35016)}

\bibitem{Chen_Zhang_2015}
Z.~Chen and Y.~Zhang, \emph{Stabilities of homothetically shrinking
  {Y}ang-{M}ills solitons}, Trans. Amer. Math. Soc. \textbf{367} (2015), no.~7,
  5015--5041, arXiv:1410.5150. \MR{3335408}

\bibitem{Chill_2003}
R.~Chill, \emph{On the {{\L}}ojasiewicz-{S}imon gradient inequality}, J. Funct.
  Anal. \textbf{201} (2003), 572--601. \MR{1986700 (2005c:26019)}

\bibitem{Chill_2006}
\bysame, \emph{The {{\L}}ojasiewicz-{S}imon gradient inequality in {H}ilbert
  spaces}, Proceedings of the 5th European-Maghrebian Workshop on Semigroup
  Theory, Evolution Equations, and Applications (M.~A. Jendoubi, ed.), 2006,
  pp.~25--36.

\bibitem{Chill_Fiorenza_2006}
R.~Chill and A.~Fiorenza, \emph{Convergence and decay rate to equilibrium of
  bounded solutions of quasilinear parabolic equations}, J. Differential
  Equations \textbf{228} (2006), 611--632. \MR{2289546 (2007k:35226)}

\bibitem{Chill_Haraux_Jendoubi_2009}
R.~Chill, A.~Haraux, and M.~A. Jendoubi, \emph{Applications of the
  {{\L}}ojasiewicz-{S}imon gradient inequality to gradient-like evolution
  equations}, Anal. Appl. (Singap.) \textbf{7} (2009), 351--372. \MR{2572850
  (2011a:35557)}

\bibitem{Chill_Jendoubi_2003}
R.~Chill and M.~A. Jendoubi, \emph{Convergence to steady states in
  asymptotically autonomous semilinear evolution equations}, Nonlinear Anal.
  \textbf{53} (2003), 1017--1039. \MR{1978032 (2004d:34103)}

\bibitem{Chill_Jendoubi_2007}
\bysame, \emph{Convergence to steady states of solutions of non-autonomous heat
  equations in {$\Bbb R^N$}}, J. Dynam. Differential Equations \textbf{19}
  (2007), 777--788. \MR{2350247 (2009h:35208)}

\bibitem{Chow_1992}
B.~Chow, \emph{The {Y}amabe flow on locally conformally flat manifolds with
  positive {R}icci curvature}, Comm. Pure Appl. Math. \textbf{45} (1992),
  1003--1014. \MR{1168117 (93d:53045)}

\bibitem{Chow_etal_Ricci_Flow_I}
B.~Chow, S.-C. Chu, D.~Glickenstein, C.~Guenther, J.~Isenberg, .~Ivey,
  D.~Knopf, P.~Lu, F.~Luo, and L.~Ni, \emph{The {R}icci flow: techniques and
  applications. {P}art {I}}, Mathematical Surveys and Monographs, vol. 135,
  American Mathematical Society, Providence, RI, 2007, Geometric aspects.
  \MR{2302600 (2008f:53088)}

\bibitem{Chow_etal_Ricci_Flow_II}
\bysame, \emph{The {R}icci flow: techniques and applications. {P}art {II}},
  Mathematical Surveys and Monographs, vol. 144, American Mathematical Society,
  Providence, RI, 2008, Analytic aspects. \MR{2365237 (2008j:53114)}

\bibitem{Chow_Knopf_2004}
B.~Chow and D.~Knopf, \emph{The {R}icci flow: an introduction}, Mathematical
  Surveys and Monographs, vol. 110, American Mathematical Society, Providence,
  RI, 2004. \MR{2061425 (2005e:53101)}

\bibitem{Chow_Lu_Ni_2006}
B.~Chow, P.~Lu, and L.~Ni, \emph{Hamilton's {R}icci flow}, Graduate Studies in
  Mathematics, vol.~77, American Mathematical Society, Providence, RI, 2006.
  \MR{2274812 (2008a:53068)}

\bibitem{Colding_Minicozzi_2014sdg}
T.~H. Colding and W.~P. Minicozzi, II, \emph{{\L}ojasiewicz inequalities and
  applications}, Surveys in Differential Geometry \textbf{XIX} (2014), 63--82,
  arXiv:1402.5087.

\bibitem{Cwikel_Janson_2006}
M.~Cwikel and S.~Janson, \emph{Complex interpolation of compact operators: an
  update}, Proc. Estonian Acad. Sci. Phys. Math. \textbf{55} (2006), no.~3,
  164--169. \MR{2257738 (2007c:46072)}

\bibitem{Daskalopoulos_1992}
G.~D. Daskalopoulos, \emph{The topology of the space of stable bundles on a
  compact {R}iemann surface}, J. Differential Geom. \textbf{36} (1992),
  699--746. \MR{1189501 (93i:58026)}

\bibitem{Daskalopoulos_Wentworth_2004}
G.~D. Daskalopoulos and R.~A. Wentworth, \emph{Convergence properties of the
  {Y}ang-{M}ills flow on {K}\"ahler surfaces}, J. Reine Angew. Math.
  \textbf{575} (2004), 69--99. \MR{2097548 (2005m:53116)}

\bibitem{Daskalopoulos_Wentworth_2004correction}
\bysame, \emph{Correction to the paper: ``{C}onvergence properties of the
  {Y}ang-{M}ills flow on {K}\"ahler surfaces'' [{J}. {R}eine {A}ngew. {M}ath.
  {\bf 575} (2004), 69--99; mr2097548]}, J. Reine Angew. Math. \textbf{606}
  (2007), 235--236. \MR{2337650 (2008f:53089)}

\bibitem{Daskalopoulos_Wentworth_2007}
\bysame, \emph{On the blow-up set of the {Y}ang-{M}ills flow on {K}\"ahler
  surfaces}, Math. Z. \textbf{256} (2007), 301--310. \MR{2289876 (2008h:53122)}

\bibitem{Davies_1989}
E.~B. Davies, \emph{Heat kernels and spectral theory}, Cambridge Tracts in
  Mathematics, vol.~92, Cambridge University Press, Cambridge, 1989.

\bibitem{DeNapoli_Drelichman_2014arxiv}
P.~L. De~Napoli and I.~Drelichman, \emph{Elementary proofs of embedding
  theorems for potential spaces of radial functions}, arXiv:1404.7468.

\bibitem{Deimling_1985}
K.~Deimling, \emph{Nonlinear functional analysis}, Springer-Verlag, Berlin,
  1985. \MR{787404 (86j:47001)}

\bibitem{Denk_Dreher_2011}
R.~Denk and M.~Dreher, \emph{Resolvent estimates for elliptic systems in
  function spaces of higher regularity}, Electron. J. Differential Equations
  (2011), No. 109, 12. \MR{2832284 (2012e:35049)}

\bibitem{Denk_Hieber_Pruss_2003}
R.~Denk, M.~Hieber, and J.~Pr{\"u}ss, \emph{{$\mathscr{R}$}-boundedness,
  {F}ourier multipliers and problems of elliptic and parabolic type}, Mem.
  Amer. Math. Soc. \textbf{166} (2003), no.~788. \MR{2006641 (2004i:35002)}

\bibitem{DeTurck_1983}
D.~M. DeTurck, \emph{Deforming metrics in the direction of their {R}icci
  tensors}, J. Differential Geom. \textbf{18} (1983), 157--162. \MR{697987
  (85j:53050)}

\bibitem{DiBlasio_1991}
G.~Di~Blasio, \emph{Analytic semigroups generated by elliptic operators in
  {$L^1$} and parabolic equations}, Osaka J. Math. \textbf{28} (1991),
  367--384. \MR{1132171 (93a:46055)}

\bibitem{DiNezza_Palatucci_Valdinoci_2012}
E.~Di~Nezza, G.~Palatucci, and E.~Valdinoci, \emph{Hitchhiker's guide to the
  fractional {S}obolev spaces}, Bull. Sci. Math. \textbf{136} (2012), no.~5,
  521--573. \MR{2944369}

\bibitem{DonNS}
S.~K. Donaldson, \emph{A new proof of a theorem of {N}arasimhan and
  {S}eshadri}, J. Differential Geom. \textbf{18} (1983), 269--277. \MR{710055
  (85a:32036)}

\bibitem{DonASD}
\bysame, \emph{Anti self-dual {Y}ang-{M}ills connections over complex algebraic
  surfaces and stable vector bundles}, Proc. London Math. Soc. (3) \textbf{50}
  (1985), 1--26. \MR{765366 (86h:58038)}

\bibitem{DonConn}
\bysame, \emph{Connections, cohomology and the intersection forms of
  four-manifolds}, J. Differential Geom. \textbf{24} (1986), 275--341.

\bibitem{DonInfDet}
\bysame, \emph{Infinite determinants, stable bundles and curvature}, Duke Math.
  J. \textbf{54} (1987), 231--247.

\bibitem{DonApprox}
\bysame, \emph{The approximation of instantons}, Geom. Funct. Anal. \textbf{3}
  (1993), 179--200.

\bibitem{DK}
S.~K. Donaldson and P.~B. Kronheimer, \emph{The geometry of four-manifolds},
  Oxford University Press, New York, 1990. \MR{1079726 (92a:57036)}

\bibitem{Dong_Kim_2011cvpde}
H.~Dong and D.~Kim, \emph{{$L_p$} solvability of divergence type parabolic and
  elliptic systems with partially {BMO} coefficients}, Calc. Var. Partial
  Differential Equations \textbf{40} (2011), 357--389. \MR{2764911
  (2011m:35393)}

\bibitem{Dong_Kim_2011arma}
\bysame, \emph{On the {$L_p$}-solvability of higher order parabolic and
  elliptic systems with {BMO} coefficients}, Arch. Ration. Mech. Anal.
  \textbf{199} (2011), 889--941. \MR{2771670 (2012h:35152)}

\bibitem{Dragomir_Ornea_1998}
S.~Dragomir and L.~Ornea, \emph{Locally conformal {K}\"ahler geometry},
  Progress in Mathematics, vol. 155, Birkh\"auser Boston, Inc., Boston, MA,
  1998. \MR{1481969 (99a:53081)}

\bibitem{Duncan_2015preprint}
D.~Duncan, \emph{The {Y}ang-{M}ills heat flow for cylindrical end
  $4$-manifolds}, December 18, 2015 preprint,
  \url{http://ms.mcmaster.ca/~duncand}.

\bibitem{Ebin_1968}
D.~G. Ebin, \emph{On the space of {R}iemannian metrics}, Bull. Amer. Math. Soc.
  \textbf{74} (1968), 1001--1003. \MR{0231410 (37 \#6965)}

\bibitem{Ebin_1970}
\bysame, \emph{The manifold of {R}iemannian metrics}, Global {A}nalysis
  ({P}roc. {S}ympos. {P}ure {M}ath., {V}ol. {XV}, {B}erkeley, {C}alif., 1968),
  Amer. Math. Soc., Providence, R.I., 1970, pp.~11--40. \MR{0267604 (42
  \#2506)}

\bibitem{Eells_Lemaire_1978}
J.~Eells and L.~Lemaire, \emph{A report on harmonic maps}, Bull. London Math.
  Soc. \textbf{10} (1978), 1--68. \MR{495450 (82b:58033)}

\bibitem{Eells_Sampson_1964ajm}
J.~Eells, Jr. and J.~H. Sampson, \emph{Harmonic mappings of {R}iemannian
  manifolds}, Amer. J. Math. \textbf{86} (1964), 109--160. \MR{0164306 (29
  \#1603)}

\bibitem{Ehrlich_1974}
P.~E. Ehrlich, \emph{Continuity properties of the injectivity radius function},
  Compositio Math. \textbf{29} (1974), 151--178. \MR{0417977 (54 \#6022)}

\bibitem{Engel_Nagel_2000}
K.-J. Engel and R.~Nagel, \emph{One-parameter semigroups for linear evolution
  equations}, Graduate Texts in Mathematics, vol. 194, Springer-Verlag, New
  York, 2000, With contributions by S. Brendle, M. Campiti, T. Hahn, G.
  Metafune, G. Nickel, D. Pallara, C. Perazzoli, A. Rhandi, S. Romanelli and R.
  Schnaubelt. \MR{1721989 (2000i:47075)}

\bibitem{Engel_Nagel_2006}
\bysame, \emph{A short course on operator semigroups}, Universitext, Springer,
  New York, 2006. \MR{2229872 (2007e:47001)}

\bibitem{Evans}
L.~C. Evans, \emph{Partial differential equations}, American Mathematical
  Society, Providence, RI, 1998.

\bibitem{Evans2}
\bysame, \emph{Partial differential equations}, second ed., Graduate Studies in
  Mathematics, vol.~19, American Mathematical Society, Providence, RI, 2010.
  \MR{2597943 (2011c:35002)}

\bibitem{Fan_1999}
H.~Fan, \emph{Existence of the self-similar solutions in the heat flow of
  harmonic maps}, Sci. China Ser. A \textbf{42} (1999), no.~2, 113--132.
  \MR{1694169 (2001f:58033)}

\bibitem{Feehan_yangmillsenergygapflat}
P.~M.~N. Feehan, \emph{Energy gap for {Y}ang-{M}ills connections, {II}:
  Arbitrary closed {R}iemannian manifolds}, arXiv:1502.00668.

\bibitem{Feehan_parabolicmaximumprinciple}
\bysame, \emph{Maximum principles for boundary-degenerate linear parabolic
  differential operators}, arXiv:1306.5197.

\bibitem{FeehanGeometry}
\bysame, \emph{Geometry of the ends of the moduli space of anti-self-dual
  connections}, J. Differential Geom. \textbf{42} (1995), 465--553,
  arXiv:1504.05741. \MR{1367401 (97d:58034)}

\bibitem{FeehanSlice}
\bysame, \emph{Critical-exponent {S}obolev norms and the slice theorem for the
  quotient space of connections}, Pacific J. Math. \textbf{200} (2001),
  71--118, arXiv:dg-ga/9711004.

\bibitem{FLKM1}
P.~M.~N. Feehan and T.~G. Leness, \emph{{D}onaldson invariants and
  wall-crossing formulas. {I}: Continuity of gluing and obstruction maps},
  arXiv:math/9812060.

\bibitem{FL3}
\bysame, \emph{{$\rm PU(2)$} monopoles. {III}: {E}xistence of gluing and
  obstruction maps}, arXiv:math/9907107.

\bibitem{FL1}
\bysame, \emph{{$\rm PU(2)$} monopoles. {I}. {R}egularity, {U}hlenbeck
  compactness, and transversality}, J. Differential Geom. \textbf{49} (1998),
  265--410. \MR{1664908 (2000e:57052)}

\bibitem{Feehan_Maridakis_Lojasiewicz-Simon_harmonic_maps}
P.~M.~N. Feehan and M.~Maridakis, \emph{{\L}ojasiewicz-{S}imon gradient
  inequalities for analytic and {M}orse-{B}ott functionals on {B}anach spaces
  and applications to harmonic maps}, arXiv:1510.03817.

\bibitem{Feehan_Maridakis_Lojasiewicz-Simon_coupled_Yang-Mills}
\bysame, \emph{{\L}ojasiewicz-{S}imon gradient inequalities for coupled
  {Y}ang-{M}ills energy functionals}, arXiv:1510.03815.

\bibitem{Feireisl_Laurencot_Petzeltova_2007}
E.~Feireisl, P.~Lauren{\c{c}}ot, and H.~Petzeltov{\'a}, \emph{On convergence to
  equilibria for the {K}eller-{S}egel chemotaxis model}, J. Differential
  Equations \textbf{236} (2007), 551--569. \MR{2322024 (2008c:35121)}

\bibitem{Feireisl_Simondon_2000}
E.~Feireisl and F.~Simondon, \emph{Convergence for semilinear degenerate
  parabolic equations in several space dimensions}, J. Dynam. Differential
  Equations \textbf{12} (2000), 647--673. \MR{1800136 (2002g:35116)}

\bibitem{Feireisl_Takac_2001}
E.~Feireisl and P.~Tak{\'a}{\v{c}}, \emph{Long-time stabilization of solutions
  to the {G}inzburg-{L}andau equations of superconductivity}, Monatsh. Math.
  \textbf{133} (2001), no.~3, 197--221. \MR{1861137 (2003a:35022)}

\bibitem{Feller2}
W.~Feller, \emph{An introduction to probability theory and its applications.
  {V}ol. {II}.}, Second edition, John Wiley \& Sons, Inc., New
  York-London-Sydney, 1971. \MR{0270403 (42 \#5292)}

\bibitem{Folland_realanalysis}
G.~B. Folland, \emph{Real analysis}, second ed., Wiley, New York, 1999.

\bibitem{FU}
D.~S. Freed and K.~K. Uhlenbeck, \emph{Instantons and four-manifolds}, second
  ed., Mathematical Sciences Research Institute Publications, vol.~1, Springer,
  New York, 1991. \MR{1081321 (91i:57019)}

\bibitem{FriedmanBundleBook}
R.~Friedman, \emph{Algebraic surfaces and holomorphic vector bundles},
  Springer-Verlag, New York, 1998.

\bibitem{FrM}
R.~Friedman and J.~W. Morgan, \emph{Smooth four-manifolds and complex
  surfaces}, Springer, Berlin, 1994.

\bibitem{Friedman_Morgan_1998}
R.~Friedman and J.~W. Morgan (eds.), \emph{Gauge theory and the topology of
  four-manifolds}, IAS/Park City Mathematics Series, vol.~4, American
  Mathematical Society, Providence, RI, 1998, Papers from the Summer School
  held in Park City, UT, July 10--30, 1994. \MR{1611444 (98j:57001)}

\bibitem{Frigeri_Grasselli_Krejcic_2013}
S.~Frigeri, M.~Grasselli, and P.~Krej{\v{c}}{\'{\i}}, \emph{Strong solutions
  for two-dimensional nonlocal {C}ahn-{H}illiard-{N}avier-{S}tokes systems}, J.
  Differential Equations \textbf{255} (2013), no.~9, 2587--2614. \MR{3090070}

\bibitem{GilbargTrudinger}
D.~Gilbarg and N.~Trudinger, \emph{Elliptic partial differential equations of
  second order}, second ed., Springer, New York, 1983.

\bibitem{Gilkey2}
P.~B. Gilkey, \emph{Invariance theory, the heat equation, and the
  {A}tiyah-{S}inger index theorem}, second ed., Studies in Advanced
  Mathematics, CRC Press, Boca Raton, FL, 1995. \MR{1396308 (98b:58156)}

\bibitem{Grasselli_Wu_2013}
M.~Grasselli and H.~Wu, \emph{Long-time behavior for a hydrodynamic model on
  nematic liquid crystal flows with asymptotic stabilizing boundary condition
  and external force}, SIAM J. Math. Anal. \textbf{45} (2013), no.~3,
  965--1002. \MR{3048212}

\bibitem{Grasselli_Wu_Zheng_2009}
M.~Grasselli, H.~Wu, and S.~Zheng, \emph{Convergence to equilibrium for
  parabolic-hyperbolic time-dependent {G}inzburg-{L}andau-{M}axwell equations},
  SIAM J. Math. Anal. \textbf{40} (2008/09), no.~5, 2007--2033.

\bibitem{Greene_Jacobowitz_1971}
R.~E. Greene and H.~Jacobowitz, \emph{Analytic isometric embeddings}, Ann. of
  Math. (2) \textbf{93} (1971), 189--204. \MR{0283728 (44 \#958)}

\bibitem{GriffithsHarris}
P.~Griffiths and J.~Harris, \emph{Principles of algebraic geometry},
  Wiley-Interscience [John Wiley \& Sons], New York, 1978, Pure and Applied
  Mathematics. \MR{507725 (80b:14001)}

\bibitem{Grigoryan_2009}
Alexander Grigor'yan, \emph{Heat kernel and analysis on manifolds}, AMS/IP
  Studies in Advanced Mathematics, vol.~47, American Mathematical Society,
  Providence, RI, 2009. \MR{2569498 (2011e:58041)}

\bibitem{Grillakis_Shatah_Strauss_1987}
M.~Grillakis, J.~Shatah, and W.~Strauss, \emph{Stability theory of solitary
  waves in the presence of symmetry. {I}}, J. Funct. Anal. \textbf{74} (1987),
  160--197. \MR{901236 (88g:35169)}

\bibitem{Grillakis_Shatah_Strauss_1990}
\bysame, \emph{Stability theory of solitary waves in the presence of symmetry.
  {II}}, J. Funct. Anal. \textbf{94} (1990), 308--348. \MR{1081647 (92a:35135)}

\bibitem{GroisserParkerSphere}
D.~Groisser and T.~H. Parker, \emph{The {R}iemannian geometry of the
  {Y}ang-{M}ills moduli space}, Comm. Math. Phys. \textbf{112} (1987),
  663--689. \MR{910586 (89b:58024)}

\bibitem{Grotowski_2001}
J.~F. Grotowski, \emph{Finite time blow-up for the {Y}ang-{M}ills heat flow in
  higher dimensions}, Math. Z. \textbf{237} (2001), 321--333. \MR{1838314
  (2002e:53102)}

\bibitem{Grotowski_Shatah_2007}
J.~F. Grotowski and J.~Shatah, \emph{Geometric evolution equations in critical
  dimensions}, Calc. Var. Partial Differential Equations \textbf{30} (2007),
  499--512. \MR{2332425 (2008e:53122)}

\bibitem{Grubb_1996}
G.~Grubb, \emph{Functional calculus of pseudodifferential boundary problems},
  second ed., Progress in Mathematics, vol.~65, Birkh\"auser Boston Inc.,
  Boston, MA, 1996. \MR{1385196 (96m:35001)}

\bibitem{Guidetti_1993}
D.~Guidetti, \emph{On elliptic systems in {$L^1$}}, Osaka J. Math. \textbf{30}
  (1993), 397--429. \MR{1240004 (94j:35042)}

\bibitem{Haase_2006}
M.~Haase, \emph{The functional calculus for sectorial operators}, Operator
  Theory: Advances and Applications, vol. 169, Birkh\"auser Verlag, Basel,
  2006. \MR{2244037 (2007j:47030)}

\bibitem{Haller-Dintelmann_Heck_Hieber_2006}
R.~Haller-Dintelmann, H.~Heck, and M.~Hieber, \emph{{$L^p$}-{$L^q$} estimates
  for parabolic systems in non-divergence form with {VMO} coefficients}, J.
  London Math. Soc. (2) \textbf{74} (2006), 717--736. \MR{2286441
  (2008b:35119)}

\bibitem{Hamilton_1975}
R.~S. Hamilton, \emph{Harmonic maps of manifolds with boundary}, Lecture Notes
  in Mathematics, vol. 471, Springer, New York, 1975. \MR{0482822 (58 \#2872)}

\bibitem{Hamilton_1982bams}
\bysame, \emph{The inverse function theorem of {N}ash and {M}oser}, Bull. Amer.
  Math. Soc. (N.S.) \textbf{7} (1982), no.~1, 65--222. \MR{656198 (83j:58014)}

\bibitem{Hamilton_1982}
\bysame, \emph{Three-manifolds with positive {R}icci curvature}, J.
  Differential Geom. \textbf{17} (1982), no.~2, 255--306. \MR{664497
  (84a:53050)}

\bibitem{Hamilton_1995}
\bysame, \emph{The formation of singularities in the {R}icci flow}, Surveys in
  differential geometry, {V}ol.\ {II} ({C}ambridge, {MA}, 1993), Int. Press,
  Cambridge, MA, 1995, pp.~7--136. \MR{1375255 (97e:53075)}

\bibitem{Haraux_2012}
A.~Haraux, \emph{Some applications of the {{\L}}ojasiewicz gradient
  inequality}, Commun. Pure Appl. Anal. \textbf{11} (2012), 2417--2427.
  \MR{2912754}

\bibitem{Haraux_Jendoubi_1998}
A.~Haraux and M.~A. Jendoubi, \emph{Convergence of solutions of second-order
  gradient-like systems with analytic nonlinearities}, J. Differential
  Equations \textbf{144} (1998), 313--320. \MR{1616968 (99a:35182)}

\bibitem{Haraux_Jendoubi_2007}
\bysame, \emph{On the convergence of global and bounded solutions of some
  evolution equations}, J. Evol. Equ. \textbf{7} (2007), 449--470. \MR{2328934
  (2008k:35480)}

\bibitem{Haraux_Jendoubi_2011}
\bysame, \emph{The {{\L}}ojasiewicz gradient inequality in the
  infinite-dimensional {H}ilbert space framework}, J. Funct. Anal. \textbf{260}
  (2011), 2826--2842. \MR{2772353 (2012c:47168)}

\bibitem{Haraux_Jendoubi_Kavian_2003}
A.~Haraux, M.~A. Jendoubi, and O.~Kavian, \emph{Rate of decay to equilibrium in
  some semilinear parabolic equations}, J. Evol. Equ. \textbf{3} (2003),
  463--484. \MR{2019030 (2004k:35187)}

\bibitem{HardyLittlewoodPolya}
G.~H. Hardy, J.~E. Littlewood, and G.~P{\'o}lya, \emph{Inequalities}, Cambridge
  University Press, Cambridge, 1952, 2nd ed.

\bibitem{Haroske_Triebel_2008}
D.~D. Haroske and H.~Triebel, \emph{Distributions, {S}obolev spaces, elliptic
  equations}, EMS Textbooks in Mathematics, European Mathematical Society
  (EMS), Z\"urich, 2008. \MR{2375667 (2009a:46003)}

\bibitem{Hartman_2002}
P.~Hartman, \emph{Ordinary differential equations}, Classics in Applied
  Mathematics, vol.~38, Society for Industrial and Applied Mathematics (SIAM),
  Philadelphia, PA, 2002, Corrected reprint of the second (1982) edition.
  \MR{1929104 (2003h:34001)}

\bibitem{Haslhofer_2012cvpde}
R.~Haslhofer, \emph{Perelman's lambda-functional and the stability of
  {R}icci-flat metrics}, Calc. Var. Partial Differential Equations \textbf{45}
  (2012), 481--504. \MR{2984143}

\bibitem{Haslhofer_Muller_2014}
R.~Haslhofer and R.~M{\"u}ller, \emph{Dynamical stability and instability of
  {R}icci-flat metrics}, Math. Ann. \textbf{360} (2014), no.~1-2, 547--553,
  arXiv:1301.3219. \MR{3263173}

\bibitem{Heck_Hieber_2003}
H.~Heck and M.~Hieber, \emph{Maximal {$L^p$}-regularity for elliptic operators
  with {VMO}-coefficients}, J. Evol. Equ. \textbf{3} (2003), 332--359.
  \MR{1980981 (2005e:35095)}

\bibitem{Heck_Hieber_Stavrakidis_2010}
H.~Heck, M.~Hieber, and K.~Stavrakidis, \emph{{$L^\infty$}-estimates for
  parabolic systems with {VMO}-coefficients}, Discrete Contin. Dyn. Syst. Ser.
  S \textbf{3} (2010), 299--309. \MR{2610566 (2011c:35230)}

\bibitem{Helein_harmonic_maps}
F.~H{\'e}lein, \emph{Harmonic maps, conservation laws and moving frames},
  second ed., Cambridge Tracts in Mathematics, vol. 150, Cambridge University
  Press, 2002. \MR{1913803 (2003g:58024)}

\bibitem{Henry_geometric_theory_semilinear_parabolic_equations}
D.~Henry, \emph{Geometric theory of semilinear parabolic equations}, Lecture
  Notes in Mathematics, vol. 840, Springer-Verlag, Berlin-New York, 1981.
  \MR{610244 (83j:35084)}

\bibitem{HewittStromberg}
E.~Hewitt and K.~Stromberg, \emph{Real and abstract analysis}, Springer-Verlag,
  New York-Heidelberg, 1975, A modern treatment of the theory of functions of a
  real variable, Third printing, Graduate Texts in Mathematics, No. 25.
  \MR{0367121 (51 \#3363)}

\bibitem{Hille_Phillips}
E.~Hille and R.~S. Phillips, \emph{Functional analysis and semi-groups},
  American Mathematical Society, Providence, R. I., 1974, Third printing of the
  revised edition of 1957, American Mathematical Society Colloquium
  Publications, Vol. XXXI. \MR{0423094 (54 \#11077)}

\bibitem{Hong_Schabrun_2010}
M-C. Hong and L.~Schabrun, \emph{Global existence for the {S}eiberg-{W}itten
  flow}, Comm. Anal. Geom. \textbf{18} (2010), 433--473. \MR{2747435
  (2012b:53139)}

\bibitem{Hong_Tian_2004}
M-C. Hong and G.~Tian, \emph{Asymptotical behaviour of the {Y}ang-{M}ills flow
  and singular {Y}ang-{M}ills connections}, Math. Ann. \textbf{330} (2004),
  441--472. \MR{2099188 (2006h:53063)}

\bibitem{Hong_Zheng_2008}
M-C. Hong and Y.~Zheng, \emph{Anti-self-dual connections and their related flow
  on 4-manifolds}, Calc. Var. Partial Differential Equations \textbf{31}
  (2008), 325--349. \MR{2366128 (2009c:53090)}

\bibitem{Hormander_v1}
L.~H{\"o}rmander, \emph{The analysis of linear partial differential operators,
  {I}. distribution theory and fourier analysis}, Springer-Verlag, Berlin,
  2003. \MR{1996773}

\bibitem{Hormander_v3}
\bysame, \emph{The analysis of linear partial differential operators, {III}.
  pseudo-differential operators}, Springer, Berlin, 2007. \MR{2304165
  (2007k:35006)}

\bibitem{Huang_2006}
S.-Z. Huang, \emph{Gradient inequalities}, Mathematical Surveys and Monographs,
  vol. 126, American Mathematical Society, Providence, RI, 2006. \MR{2226672
  (2007b:35035)}

\bibitem{Huang_Takac_2001}
S.-Z. Huang and P.~Tak{\'a}{\v{c}}, \emph{Convergence in gradient-like systems
  which are asymptotically autonomous and analytic}, Nonlinear Anal.
  \textbf{46} (2001), 675--698. \MR{1857152 (2002f:35125)}

\bibitem{Huybrechts_2005}
D.~Huybrechts, \emph{Complex geometry}, Universitext, Springer-Verlag, Berlin,
  2005, An introduction. \MR{2093043 (2005h:32052)}

\bibitem{Ilmanen_1994}
T.~Ilmanen, \emph{Elliptic regularization and partial regularity for motion by
  mean curvature}, Mem. Amer. Math. Soc. \textbf{108} (1994), no.~520.
  \MR{1196160 (95d:49060)}

\bibitem{IrwinThesis}
C.~A. Irwin, \emph{Bubbling in the harmonic map heat flow}, {Ph.D.} thesis,
  Stanford University, Palo Alto, CA, 1998. \MR{2698290}

\bibitem{Jacob_2015conm}
A.~Jacob, \emph{Stable {H}iggs bundles and {H}ermitian-{E}instein metrics on
  non-{K}\"ahler manifolds}, Analysis, complex geometry, and mathematical
  physics: in honor of {D}uong {H}. {P}hong, Contemp. Math., vol. 644, Amer.
  Math. Soc., Providence, RI, 2015, pp.~117--140. \MR{3372463}

\bibitem{Jacob_v1}
N.~Jacob, \emph{Pseudo-differential operators and {M}arkov processes}, vol. 1:
  Fourier Analysis and Semigroups, Imperial College Press, London, 2001.

\bibitem{Jacob_v2}
\bysame, \emph{Pseudo-differential operators and {M}arkov processes}, vol. 2:
  Generators and Their Potential Theory, Imperial College Press, London, 2001.

\bibitem{Jacob_v3}
\bysame, \emph{Pseudo-differential operators and {M}arkov processes}, vol. 3:
  Markov Processes And Applications, Imperial College Press, London, 2001.

\bibitem{Jeanjean_1999}
L.~Jeanjean, \emph{On the existence of bounded {P}alais-{S}male sequences and
  application to a {L}andesman-{L}azer-type problem set on {${\bf R}^N$}},
  Proc. Roy. Soc. Edinburgh Sect. A \textbf{129} (1999), 787--809. \MR{1718530
  (2001c:35034)}

\bibitem{Jendoubi_1998jde}
M.~A. Jendoubi, \emph{Convergence of global and bounded solutions of the wave
  equation with linear dissipation and analytic nonlinearity}, J. Differential
  Equations \textbf{144} (1998), 302--312. \MR{1616964 (99e:35149)}

\bibitem{Jendoubi_1998jfa}
\bysame, \emph{A simple unified approach to some convergence theorems of {L}.
  {S}imon}, J. Funct. Anal. \textbf{153} (1998), 187--202. \MR{1609269
  (99c:35101)}

\bibitem{Jost_1991}
J.~Jost, \emph{Nonlinear methods in {R}iemannian and {K}\"ahlerian geometry},
  second ed., DMV Seminar, vol.~10, Birkh\"auser Verlag, Basel, 1991.
  \MR{1138210 (93k:58065)}

\bibitem{Jost_two_dim_geom_var_probs}
\bysame, \emph{Two-dimensional geometric variational problems}, Pure and
  Applied Mathematics (New York), John Wiley \& Sons, Ltd., Chichester, 1991, A
  Wiley-Interscience Publication. \MR{1100926 (92h:58045)}

\bibitem{Jost_riemannian_geometry_geometric_analysis}
\bysame, \emph{Riemannian geometry and geometric analysis}, sixth ed.,
  Universitext, Springer, Heidelberg, 2011. \MR{2829653}

\bibitem{Kato}
T.~Kato, \emph{Perturbation theory for linear operators}, second ed., Springer,
  New York, 1984.

\bibitem{Kelleher_Streets_2014arxiv}
C.~Kelleher and J.~Streets, \emph{Entropy, stability, and {Y}ang-{M}ills flow},
  Communications in Contemporary Mathematics, arXiv:1410.4547.

\bibitem{Kelleher_Streets_2016preprint}
\bysame, \emph{Singularity formation of the {Y}ang-{M}ills flow}, preprint,
  February 1, 2016, \url{math.uci.edu/~clkelleh/about.html}.

\bibitem{Kleiner_Lott_2008}
B.~Kleiner and J.~Lott, \emph{Notes on {P}erelman's papers}, Geometry and
  Topology \textbf{12} (2008), 2587--2855, arXiv:math/0605667.

\bibitem{Knapp_2002}
A.~W. Knapp, \emph{Representation theory of semisimple groups}, Princeton
  Mathematical Series, vol.~36, Princeton University Press, Princeton, NJ,
  1986. \MR{855239 (87j:22022)}

\bibitem{Kobayashi}
S.~Kobayashi, \emph{Differential geometry of complex vector bundles},
  Publications of the Mathematical Society of Japan, vol.~15, Princeton
  University Press, Princeton, NJ, 1987, Kan{\^o} Memorial Lectures, 5.
  \MR{909698 (89e:53100)}

\bibitem{Kobayashi_Nomizu_v1}
S.~Kobayashi and K.~Nomizu, \emph{Foundations of differential geometry. {V}ol
  {I}}, Interscience Publishers, a division of John Wiley \& Sons, New
  York-London, 1963. \MR{0152974 (27 \#2945)}

\bibitem{Kobayashi_Nomizu_v2}
\bysame, \emph{Foundations of differential geometry. {V}ol. {II}}, Interscience
  Tracts in Pure and Applied Mathematics, No. 15 Vol. II, Interscience
  Publishers John Wiley \& Sons, Inc., New York-London-Sydney, 1969.
  \MR{0238225 (38 \#6501)}

\bibitem{Kono_Nagasawa_1988}
K.~Kono and T.~Nagasawa, \emph{Weak asymptotical stability of {Y}ang-{M}ills
  gradient flow}, Tokyo J. Math. \textbf{11} (1988), 339--357. \MR{976570
  (90a:58034)}

\bibitem{Koshelev_1995}
A.~Koshelev, \emph{Regularity problem for quasilinear elliptic and parabolic
  systems}, Lecture Notes in Mathematics, vol. 1614, Springer-Verlag, Berlin,
  1995. \MR{1442954 (98j:35004)}

\bibitem{Kozono_Maeda_Naito_1995}
H.~Kozono, Y.~Maeda, and H.~Naito, \emph{Global solution for the {Y}ang-{M}ills
  gradient flow on {$4$}-manifolds}, Nagoya Math. J. \textbf{139} (1995),
  93--128. \MR{1355271 (97a:58038)}

\bibitem{Krein_Petunin_Semenov_interpolation_linear_operators}
S.~G. Kre{\u\i}n, Yu.~{\=I}. Petun{\={\i}}n, and E.~M. Sem{\"e}nov,
  \emph{Interpolation of linear operators}, Translations of Mathematical
  Monographs, vol.~54, American Mathematical Society, Providence, R.I., 1982,
  Translated from the Russian by J. Sz{\H{u}}cs. \MR{649411 (84j:46103)}

\bibitem{Kroncke_2013arxiv}
K.~Kr{\"o}ncke, \emph{Ricci flow, {E}instein metrics and the {Y}amabe
  invariant}, arXiv:1312.2224.

\bibitem{Kroncke_2015cvpde}
\bysame, \emph{Stability and instability of {R}icci solitons}, Calc. Var.
  Partial Differential Equations \textbf{53} (2015), no.~1-2, 265--287,
  arXiv:1403.3721. \MR{3336320}

\bibitem{KMStructure}
P.~B. Kronheimer and T.~S. Mrowka, \emph{Embedded surfaces and the structure of
  {D}onaldson's polynomial invariants}, J. Differential Geom. \textbf{41}
  (1995), 573--734. \MR{1338483 (96e:57019)}

\bibitem{Krylov_LecturesHolder}
N.~V. Krylov, \emph{Lectures on elliptic and parabolic equations in {H}\"older
  spaces}, American Mathematical Society, Providence, RI, 1996.

\bibitem{Krylov_LecturesSobolev}
\bysame, \emph{Lectures on elliptic and parabolic equations in {S}obolev
  spaces}, American Mathematical Society, Providence, RI, 2008.

\bibitem{Kuranishi}
M.~Kuranishi, \emph{New proof for the existence of locally complete families of
  complex structures}, Proc. {C}onf. {C}omplex {A}nalysis ({M}inneapolis, 1964)
  (A.~Aeppli, E.~Calabi, and H.~R{\"o}hrl, eds.), Springer, Berlin, 1965,
  pp.~142--154. \MR{0176496 (31 \#768)}

\bibitem{KwonThesis}
H.~Kwon, \emph{Asymptotic convergence of harmonic map heat flow}, {Ph.D}.
  thesis, Stanford University, Palo Alto, CA, 2002. \MR{2703296}

\bibitem{LadyzenskajaSolonnikovUralceva}
O.~A. Lady{\v{z}}enskaja, V.~A. Solonnikov, and N.~N. Ural{\cprime}ceva,
  \emph{Linear and quasilinear equations of parabolic type}, Translated from
  the Russian by S. Smith. Translations of Mathematical Monographs, Vol. 23,
  American Mathematical Society, Providence, R.I., 1968. \MR{0241822 (39
  \#3159b)}

\bibitem{Ladyzhenskaya_Uraltseva_1968}
O.~A. Ladyzhenskaya and N.~N. Ural{\cprime}tseva, \emph{Linear and quasilinear
  elliptic equations}, Translated from the Russian by Scripta Technica, Inc.
  Translation editor: Leon Ehrenpreis, Academic Press, New York, 1968.
  \MR{0244627 (39 \#5941)}

\bibitem{Lawson}
H.~B. Lawson, \emph{The theory of gauge fields in four dimensions}, Conf. Board
  Math. Sci., vol.~58, Amer. Math. Soc., Providence, RI, 1985.

\bibitem{LM}
H.~B. Lawson, Jr. and M.-L. Michelsohn, \emph{Spin geometry}, Princeton
  Mathematical Series, vol.~38, Princeton University Press, Princeton, NJ,
  1989. \MR{1031992 (91g:53001)}

\bibitem{Leng_Zhao_Zhao_2014}
Y.~Leng, E.~Zhao, and H.~Zhao, \emph{Notes on the extension of the mean
  curvature flow}, Pacific J. Math. \textbf{269} (2014), no.~2, 385--392.
  \MR{3238481}

\bibitem{Li_Zhang_2011}
J.~Li and X.~Zhang, \emph{The gradient flow of {H}iggs pairs}, J. Eur. Math.
  Soc. (JEMS) \textbf{13} (2011), 1373--1422. \MR{2825168 (2012m:53043)}

\bibitem{Li_Zhu_2012}
J.~Li and X.~Zhu, \emph{Energy identity for the maps from a surface with
  tension field bounded in {$L^p$}}, Pacific J. Math. \textbf{260} (2012),
  181--195. \MR{3001790}

\bibitem{Lieberman}
G.~M. Lieberman, \emph{Second order parabolic differential equations}, World
  Scientific Publishing Co. Inc., River Edge, NJ, 1996.

\bibitem{Lin_Wang_2008}
F.-H. Lin and C.~Wang, \emph{The analysis of harmonic maps and their heat
  flows}, World Scientific, Hackensack, NJ, 2008. \MR{2431658 (2011a:58030)}

\bibitem{Lions_Magenes_v1_1972}
J.-L. Lions and E.~Magenes, \emph{Non-homogeneous boundary value problems and
  applications. {V}ol. {I}}, Die Grundlehren der mathematischen Wissenschaften,
  Band 181, Springer-Verlag, New York, 1972. \MR{0350177 (50 \#2670)}

\bibitem{Liu_Yang_2014arxiv}
K-F. Liu and X-K. Yang, \emph{Ricci curvatures on hermitian manifolds},
  arXiv:1404.2481.

\bibitem{Liu_Yang_2012}
\bysame, \emph{Geometry of {H}ermitian manifolds}, Internat. J. Math.
  \textbf{23} (2012), 1250055, 40. \MR{2925476}

\bibitem{Liu_Yang_2010}
Q.~Liu and Y.~Yang, \emph{Rigidity of the harmonic map heat flow from the
  sphere to compact {K}\"ahler manifolds}, Ark. Mat. \textbf{48} (2010),
  121--130. \MR{2594589 (2011a:53066)}

\bibitem{Lockhart_McOwen_1985}
R.~B. Lockhart and R.~C. McOwen, \emph{Elliptic differential operators on
  noncompact manifolds}, Ann. Scuola Norm. Sup. Pisa Cl. Sci. (4) \textbf{12}
  (1985), 409--447. \MR{837256 (87k:58266)}

\bibitem{Lojasiewicz_1963}
S.~{\L}ojasiewicz, \emph{Une propri\'et\'e topologique des sous-ensembles
  analytiques r\'eels}, Les \'{E}quations aux {D}\'eriv\'ees {P}artielles
  ({P}aris, 1962), \'Editions du Centre National de la Recherche Scientifique,
  Paris, 1963, pp.~87--89. \MR{0160856 (28 \#4066)}

\bibitem{Lojasiewicz_1965}
\bysame, \emph{Ensembles semi-analytiques},  (1965), Publ. Inst. Hautes Etudes
  Sci., Bures-sur-Yvette, preprint, 112 pages,
  \url{perso.univ-rennes1.fr/michel.coste/Lojasiewicz.pdf}.

\bibitem{Lojasiewicz_1984}
\bysame, \emph{Sur les trajectoires du gradient d'une fonction analytique},
  Geometry seminars, 1982--1983 ({B}ologna, 1982/1983), Univ. Stud. Bologna,
  Bologna, 1984, pp.~115--117. \MR{771152 (86m:58023)}

\bibitem{Lojasiewicz_1993}
\bysame, \emph{Sur la g\'eom\'etrie semi- et sous-analytique}, Ann. Inst.
  Fourier (Grenoble) \textbf{43} (1993), 1575--1595. \MR{1275210 (96c:32007)}

\bibitem{Lorenzi_Bertoldi_2007}
L.~Lorenzi and M.~Bertoldi, \emph{Analytical methods for {M}arkov semigroups},
  Chapman \& Hall/CRC, Boca Raton, FL, 2007.

\bibitem{Lubke_Teleman_1995}
M.~L{\"u}bke and A.~Teleman, \emph{The {K}obayashi-{H}itchin correspondence},
  World Scientific Publishing Co., Inc., River Edge, NJ, 1995. \MR{1370660
  (97h:32043)}

\bibitem{Lunardi_1995}
A.~Lunardi, \emph{Analytic semigroups and optimal regularity in parabolic
  problems}, Progress in Nonlinear Differential Equations and their
  Applications, 16, Birkh\"auser Verlag, Basel, 1995. \MR{1329547 (96e:47039)}

\bibitem{Lunardi_Metafune_2004}
A.~Lunardi and G.~Metafune, \emph{On the domains of elliptic operators in
  {$L^1$}}, Differential Integral Equations \textbf{17} (2004), 73--97.
  \MR{2035496 (2004m:35104)}

\bibitem{Luo_2012}
Y.~Luo, \emph{Energy identity and removable singularities of maps from a
  {R}iemann surface with tension field unbounded in {$L^2$}}, Pacific J. Math.
  \textbf{256} (2012), 365--380. \MR{2944981}

\bibitem{Mantegazza_2011_lectures_mcf}
C.~Mantegazza, \emph{Lecture notes on mean curvature flow}, Progress in
  Mathematics, vol. 290, Birkh\"auser/Springer Basel AG, Basel, 2011.
  \MR{2815949}

\bibitem{Maugeri_Palagachev_Softova_2000}
A.~Maugeri, D.~K. Palagachev, and L.~G. Softova, \emph{Elliptic and parabolic
  equations with discontinuous coefficients}, Mathematical Research, vol. 109,
  Wiley-VCH Verlag Berlin GmbH, Berlin, 2000. \MR{2260015 (2007f:35001)}

\bibitem{Mawhin_Willem_2010}
J.~Mawhin and M.~Willem, \emph{Origin and evolution of the {P}alais-{S}male
  condition in critical point theory}, J. Fixed Point Theory Appl. \textbf{7}
  (2010), 265--290. \MR{2729392 (2011j:58016)}

\bibitem{Mazya_2011}
V.~Maz'ya, \emph{Sobolev spaces with applications to elliptic partial
  differential equations}, augmented ed., Grundlehren der Mathematischen
  Wissenschaften [Fundamental Principles of Mathematical Sciences], vol. 342,
  Springer, Heidelberg, 2011. \MR{2777530 (2012a:46056)}

\bibitem{McNamara_Zhao_2014arxiv}
J.~McNamara and Y.~Zhao, \emph{Limiting behavior of donaldson’s heat flow on
  non-k{\"a}hler}, arXiv:1403.8037.

\bibitem{Megginson_1998}
R.~E. Megginson, \emph{An introduction to {B}anach space theory}, Graduate
  Texts in Mathematics, vol. 183, Springer-Verlag, New York, 1998. \MR{1650235
  (99k:46002)}

\bibitem{Milgram_Rosenbloom_1951a}
A.~N. Milgram and P.~C. Rosenbloom, \emph{Harmonic forms and heat conduction.
  {I}. {C}losed {R}iemannian manifolds}, Proc. Nat. Acad. Sci. U. S. A.
  \textbf{37} (1951), 180--184. \MR{0042769 (13,160a)}

\bibitem{Min-Oo_1982}
M.~Min-Oo, \emph{An {$L_{2}$}-isolation theorem for {Y}ang-{M}ills fields},
  Compositio Math. \textbf{47} (1982), 153--163. \MR{0677017 (84b:53033a)}

\bibitem{Morgan_Fong_2010}
J.~W. Morgan and F.~T.-H. Fong, \emph{Ricci flow and geometrization of
  3-manifolds}, University Lecture Series, vol.~53, American Mathematical
  Society, Providence, RI, 2010. \MR{2597148 (2011d:53160)}

\bibitem{MorganMrowkaTube}
J.~W. Morgan and T.~S. Mrowka, \emph{The gluing construction for anti-self-dual
  connections over manifolds with long tubes}, preprint, October 13, 1994, 115
  pages.

\bibitem{MorganMrowkaPoly}
\bysame, \emph{A note on {D}onaldson's polynomial invariants}, Internat. Math.
  Res. Notices (1992), 223--230. \MR{1191573 (93m:57032)}

\bibitem{MMR}
J.~W. Morgan, T.~S. Mrowka, and D.~Ruberman, \emph{The {$L^2$}-moduli space and
  a vanishing theorem for {D}onaldson polynomial invariants}, Monographs in
  Geometry and Topology, vol.~2, International Press, Cambridge, MA, 1994.
  \MR{1287851 (95h:57039)}

\bibitem{Morgan_Tian_2007}
J.~W. Morgan and G.~Tian, \emph{Ricci flow and the {P}oincar\'e conjecture},
  Clay Mathematics Monographs, vol.~3, American Mathematical Society,
  Providence, RI, 2007. \MR{2334563 (2008d:57020)}

\bibitem{Moroianu_2007}
A.~Moroianu, \emph{Lectures on {K}\"ahler geometry}, London Mathematical
  Society Student Texts, vol.~69, Cambridge University Press, Cambridge, 2007.
  \MR{2325093 (2008f:32028)}

\bibitem{Morrey}
C.~B. Morrey, Jr., \emph{Multiple integrals in the calculus of variations},
  Classics in Mathematics, Springer-Verlag, Berlin, 2008, Reprint of the 1966
  edition. \MR{2492985}

\bibitem{Moser_1966b}
J.~Moser, \emph{A rapidly convergent iteration method and non-linear
  differential equations. {II}}, Ann. Scuola Norm. Sup. Pisa (3) \textbf{20}
  (1966), 499--535. \MR{0206461 (34 \#6280)}

\bibitem{Moser_1966a}
\bysame, \emph{A rapidly convergent iteration method and non-linear partial
  differential equations. {I}}, Ann. Scuola Norm. Sup. Pisa (3) \textbf{20}
  (1966), 265--315. \MR{0199523 (33 \#7667)}

\bibitem{MrowkaThesis}
T.~S. Mrowka, \emph{Local {M}ayer-{V}ietoris principle for {Y}ang-{M}ills
  moduli spaces}, {Ph.D}. thesis, Harvard University, Cambridge, MA, 1988.

\bibitem{Muller_2006}
R.~M{\"u}ller, \emph{Differential {H}arnack inequalities and the {R}icci flow},
  EMS Series of Lectures in Mathematics, European Mathematical Society (EMS),
  Z\"urich, 2006. \MR{2251315 (2007d:53112)}

\bibitem{Munkres2}
J.~R. Munkres, \emph{Topology: a first course}, second ed., Prentice-Hall,
  Englewood Cliffs, N.J., 2000.

\bibitem{Nagasawa_1989}
T.~Nagasawa, \emph{Asymptotical stability of {Y}ang-{M}ills' gradient flow},
  Geometry of manifolds (K.~Shiohama, ed.), Perspectives in Mathematics,
  vol.~8, Academic Press, Boston, MA, 1989, Papers presented at the
  Thirty-fifth Symposium on Differential Geometry, July 25--30, 1988 at Shinshu
  University, Matsumoto, Japan., pp.~499--517. \MR{1040545 (91c:58029)}

\bibitem{Naito_1994}
H.~Naito, \emph{Finite time blowing-up for the {Y}ang-{M}ills gradient flow in
  higher dimensions}, Hokkaido Math. J. \textbf{23} (1994), 451--464.
  \MR{1299637 (95i:58054)}

\bibitem{Nakajima_1987}
H.~Nakajima, \emph{Removable singularities for {Y}ang-{M}ills connections in
  higher dimensions}, J. Fac. Sci. Univ. Tokyo Sect. IA Math. \textbf{34}
  (1987), 299--307. \MR{914024 (89d:58029)}

\bibitem{Nash_1956}
J.~Nash, \emph{The imbedding problem for {R}iemannian manifolds}, Ann. of Math.
  (2) \textbf{63} (1956), 20--63. \MR{0075639 (17,782b)}

\bibitem{Nash_1966}
\bysame, \emph{Analyticity of the solutions of implicit function problems with
  analytic data}, Ann. of Math. (2) \textbf{84} (1966), 345--355. \MR{0205266
  (34 \#5099)}

\bibitem{Neuberger_2010}
J.~W. Neuberger, \emph{Sobolev gradients and differential equations}, second
  ed., Lecture Notes in Mathematics, vol. 1670, Springer-Verlag, Berlin, 2010.
  \MR{2573187 (2011m:35005)}

\bibitem{Neves_2013}
A.~Neves, \emph{Finite time singularities for {L}agrangian mean curvature
  flow}, Ann. of Math. (2) \textbf{177} (2013), 1029--1076. \MR{3034293}

\bibitem{OHara_1991}
J.~O'Hara, \emph{Energy of a knot}, Topology \textbf{30} (1991), no.~2,
  241--247. \MR{1098918 (92c:58017)}

\bibitem{OHara_1992}
\bysame, \emph{Family of energy functionals of knots}, Topology Appl.
  \textbf{48} (1992), no.~2, 147--161. \MR{1195506 (94h:58064)}

\bibitem{OHara_1994}
\bysame, \emph{Energy functionals of knots. {II}}, Topology Appl. \textbf{56}
  (1994), no.~1, 45--61. \MR{1261169 (94m:58028)}

\bibitem{Olver_Lozier_Boisvert_Clark}
F.~W.~J. Olver, D.~W. Lozier, R.~F. Boisvert, and C.~W. Clark, \emph{{NIST}
  handbook of mathematical functions}, National Institute of Standards and
  Technology and Cambridge University Press, New York, 2010.

\bibitem{Palais_1961}
R.~S. Palais, \emph{On the existence of slices for actions of non-compact {L}ie
  groups}, Ann. of Math. (2) \textbf{73} (1961), 295--323. \MR{0126506 (23
  \#A3802)}

\bibitem{Palais_1966LS}
\bysame, \emph{Lusternik-{S}chnirelman theory on {B}anach manifolds}, Topology
  \textbf{5} (1966), 115--132. \MR{0259955 (41 \#4584)}

\bibitem{PalaisFoundationGlobal}
\bysame, \emph{Foundations of global non-linear analysis}, Benjamin, New York,
  1968. \MR{0248880 (40 \#2130)}

\bibitem{ParkerGauge}
T.~H. Parker, \emph{Gauge theories on four-dimensional {R}iemannian manifolds},
  Comm. Math. Phys. \textbf{85} (1982), 563--602. \MR{677998 (84b:58036)}

\bibitem{ParkerHarmonic}
\bysame, \emph{Bubble tree convergence for harmonic maps}, J. Differential
  Geom. \textbf{44} (1996), 595--633. \MR{1431008 (98k:58069)}

\bibitem{ParkerTaubes}
T.~H. Parker and C.~H. Taubes, \emph{On {W}itten's proof of the positive energy
  theorem}, Comm. Math. Phys. \textbf{84} (1982), 223--238. \MR{661134
  (83m:83020)}

\bibitem{ParkerWolfson}
T.~H. Parker and J.~G. Wolfson, \emph{Pseudo-holomorphic maps and bubble
  trees}, J. Geom. Anal. \textbf{3} (1993), 63--98. \MR{1197017 (95c:58032)}

\bibitem{Pazy_1983}
A.~Pazy, \emph{Semigroups of linear operators and applications to partial
  differential equations}, Applied Mathematical Sciences, vol.~44,
  Springer-Verlag, New York, 1983. \MR{710486 (85g:47061)}

\bibitem{Perelman_2002}
G.~Perelman, \emph{The entropy formula for the {R}icci flow and its geometric
  applications}, arXiv:math/0211159.

\bibitem{Perelman_2003b}
\bysame, \emph{Finite extinction time for the solutions to the {R}icci flow on
  certain three-manifolds}, arXiv:math/0307245.

\bibitem{Perelman_2003a}
\bysame, \emph{Ricci flow with surgery on three-manifolds}, arXiv:math/0303109.

\bibitem{Pulemotov_2008}
A.~Pulemotov, \emph{The {L}i-{Y}au-{H}amilton estimate and the {Y}ang-{M}ills
  heat equation on manifolds with boundary}, J. Funct. Anal. \textbf{255}
  (2008), 2933--2965. \MR{2464197 (2009k:58049)}

\bibitem{Ramanathan_Subramanian_1988}
A.~Ramanathan and S.~Subramanian, \emph{Einstein-{H}ermitian connections on
  principal bundles and stability}, J. Reine Angew. Math. \textbf{390} (1988),
  21--31. \MR{953674 (90b:32057)}

\bibitem{Renardy_Rogers_2004}
M.~Renardy and R.~C. Rogers, \emph{An introduction to partial differential
  equations}, second ed., Texts in Applied Mathematics, vol.~13,
  Springer-Verlag, New York, 2004. \MR{2028503 (2004j:35001)}

\bibitem{RadeThesis}
J.~R\r{a}de, \emph{On the {Y}ang-{M}ills heat equation in two and three
  dimensions}, {Ph.D}. thesis, The University of Texas, Austin, TX, 1991.
  \MR{2686594}

\bibitem{Rade_1992}
\bysame, \emph{On the {Y}ang-{M}ills heat equation in two and three
  dimensions}, J. Reine Angew. Math. \textbf{431} (1992), 123--163. \MR{1179335
  (94a:58041)}

\bibitem{Rudin}
W.~Rudin, \emph{Functional analysis}, Mc{G}raw-{H}ill, New York, NY, 1973.

\bibitem{Rybka_Hoffmann_1998}
P.~Rybka and K.-H. Hoffmann, \emph{Convergence of solutions to the equation of
  quasi-static approximation of viscoelasticity with capillarity}, J. Math.
  Anal. Appl. \textbf{226} (1998), 61--81. \MR{1646449 (99h:35146)}

\bibitem{Rybka_Hoffmann_1999}
\bysame, \emph{Convergence of solutions to {C}ahn-{H}illiard equation}, Comm.
  Partial Differential Equations \textbf{24} (1999), 1055--1077. \MR{1680877
  (2001a:35028)}

\bibitem{Sacks_Uhlenbeck_1981}
J.~Sacks and K.~K. Uhlenbeck, \emph{The existence of minimal immersions of
  {$2$}-spheres}, Ann. of Math. (2) \textbf{113} (1981), 1--24. \MR{604040
  (82f:58035)}

\bibitem{Sacks_Uhlenbeck_1982}
\bysame, \emph{Minimal immersions of closed {R}iemann surfaces}, Trans. Amer.
  Math. Soc. \textbf{271} (1982), 639--652. \MR{654854 (83i:58030)}

\bibitem{Sadun_1987thesis}
L.~Sadun, \emph{Continuum regularized {Y}ang-{M}ills theory}, {Ph.D}. thesis,
  The University of California, Berkeley, CA, 1987.

\bibitem{Sadun_1994}
\bysame, \emph{A symmetric family of {Y}ang-{M}ills fields}, Comm. Math. Phys.
  \textbf{163} (1994), no.~2, 257--291. \MR{1284785 (95f:53059)}

\bibitem{Sadun_Segert_1991}
L.~Sadun and J.~Segert, \emph{Non-self-dual {Y}ang-{M}ills connections with
  nonzero {C}hern number}, Bull. Amer. Math. Soc. (N.S.) \textbf{24} (1991),
  163--170. \MR{1067574 (91m:58038)}

\bibitem{Sadun_Segert_1992cmp}
\bysame, \emph{Non-self-dual {Y}ang-{M}ills connections with quadrupole
  symmetry}, Comm. Math. Phys. \textbf{145} (1992), 363--391. \MR{1162804
  (94b:53056)}

\bibitem{Sadun_Segert_1992cpam}
\bysame, \emph{Stationary points of the {Y}ang-{M}ills action}, Comm. Pure
  Appl. Math. \textbf{45} (1992), 461--484. \MR{1161540 (93d:58034)}

\bibitem{Sadun_Segert_1993}
\bysame, \emph{Constructing non-self-dual {Y}ang-{M}ills connections on {$S^4$}
  with arbitrary {C}hern number}, Differential geometry: geometry in
  mathematical physics and related topics ({L}os {A}ngeles, {CA}, 1990), Proc.
  Sympos. Pure Math., vol.~54, Amer. Math. Soc., Providence, RI, 1993,
  pp.~529--537. \MR{1216561}

\bibitem{Sakai_1983}
T.~Sakai, \emph{On continuity of injectivity radius function}, Math. J. Okayama
  Univ. \textbf{25} (1983), 91--97. \MR{701970 (84g:53070)}

\bibitem{Schlatter_1996}
A.~E. Schlatter, \emph{Global existence of the {Y}ang-{M}ills flow in four
  dimensions}, J. Reine Angew. Math. \textbf{479} (1996), 133--148. \MR{1414392
  (97k:58050)}

\bibitem{Schlatter_1997}
\bysame, \emph{Long-time behaviour of the {Y}ang-{M}ills flow in four
  dimensions}, Ann. Global Anal. Geom. \textbf{15} (1997), 1--25. \MR{1443269
  (98i:58066)}

\bibitem{Schlatter_Struwe_Tahvildar-Zadeh_1998}
A.~E. Schlatter, M.~Struwe, and A.~S. Tahvildar-Zadeh, \emph{Global existence
  of the equivariant {Y}ang-{M}ills heat flow in four space dimensions}, Amer.
  J. Math. \textbf{120} (1998), 117--128. \MR{1600272 (98k:58063)}

\bibitem{Schoen_1984}
R.~M. Schoen, \emph{Conformal deformation of a {R}iemannian metric to constant
  scalar curvature}, J. Differential Geom. \textbf{20} (1984), 479--495.
  \MR{788292 (86i:58137)}

\bibitem{Schoen_1991}
\bysame, \emph{On the number of constant scalar curvature metrics in a
  conformal class}, Differential geometry, Pitman Monogr. Surveys Pure Appl.
  Math., vol.~52, Longman Sci. Tech., Harlow, 1991, pp.~311--320. \MR{1173050
  (94e:53035)}

\bibitem{Schoen_Yau_1981}
R.~M. Schoen and S.-T. Yau, \emph{Proof of the positive mass theorem. {II}},
  Comm. Math. Phys. \textbf{79} (1981), 231--260. \MR{612249 (83i:83045)}

\bibitem{Schoen_Yau_1988}
\bysame, \emph{Conformally flat manifolds, {K}leinian groups and scalar
  curvature}, Invent. Math. \textbf{92} (1988), 47--71. \MR{931204 (89c:58139)}

\bibitem{Schulze_2014}
F.~Schulze, \emph{Uniqueness of compact tangent flows in mean curvature flow},
  J. Reine Angew. Math. \textbf{690} (2014), 163--172, arXiv:1107.4643.
  \MR{3200339}

\bibitem{Schwetlick_Struwe_2003}
H.~Schwetlick and M.~Struwe, \emph{Convergence of the {Y}amabe flow for
  ``large'' energies}, J. Reine Angew. Math. \textbf{562} (2003), 59--100.
  \MR{2011332 (2004h:53097)}

\bibitem{Sedlacek}
S.~Sedlacek, \emph{A direct method for minimizing the {Y}ang-{M}ills functional
  over {$4$}-manifolds}, Comm. Math. Phys. \textbf{86} (1982), 515--527.
  \MR{679200 (84e:81049)}

\bibitem{Sell_You_2002}
G.~R. Sell and Y.~You, \emph{Dynamics of evolutionary equations}, Applied
  Mathematical Sciences, vol. 143, Springer, New York, 2002. \MR{MR1873467
  (2003f:37001b)}

\bibitem{Sesum_2006}
N.~Sesum, \emph{Linear and dynamical stability of {R}icci-flat metrics}, Duke
  Math. J. \textbf{133} (2006), 1--26. \MR{2219268 (2007c:53089)}

\bibitem{Showalter}
R.~E. Showalter, \emph{Monotone operators in {B}anach space and nonlinear
  partial differential equations}, American Mathematical Society, Providence,
  RI, 1996.

\bibitem{SibnerSibnerUhlenbeck}
L.~M. Sibner, R.~J. Sibner, and K.~K. Uhlenbeck, \emph{Solutions to
  {Y}ang-{M}ills equations that are not self-dual}, Proc. Nat. Acad. Sci.
  U.S.A. \textbf{86} (1989), 8610--8613. \MR{1023811 (90j:58032)}

\bibitem{Simon_1983}
L.~Simon, \emph{Asymptotics for a class of nonlinear evolution equations, with
  applications to geometric problems}, Ann. of Math. (2) \textbf{118} (1983),
  525--571. \MR{727703 (85b:58121)}

\bibitem{Simon_1985}
\bysame, \emph{Isolated singularities of extrema of geometric variational
  problems}, Lecture Notes in Math., vol. 1161, Springer, Berlin, 1985.
  \MR{821971 (87d:58045)}

\bibitem{Simon_1996}
\bysame, \emph{Theorems on regularity and singularity of energy minimizing
  maps}, Lectures in Mathematics ETH Z\"urich, Birkh{\"a}user, Basel, 1996.
  \MR{1399562 (98c:58042)}

\bibitem{Smoczyk_2012}
K.~Smoczyk, \emph{Mean curvature flow in higher codimension: introduction and
  survey}, Global differential geometry, Springer Proc. Math., vol.~17,
  Springer, Heidelberg, 2012, pp.~231--274. \MR{3289845}

\bibitem{Stein}
E.~M. Stein, \emph{Singular integrals and differentiability properties of
  functions}, Princeton Mathematical Series, No. 30, Princeton University
  Press, Princeton, N.J., 1970. \MR{0290095}

\bibitem{Stern_2010}
M.~Stern, \emph{Geometry of minimal energy {Y}ang-{M}ills connections}, J.
  Differential Geom. \textbf{86} (2010), 163--188. \MR{2772548 (2012h:53055)}

\bibitem{Stewart_1974}
H.~B. Stewart, \emph{Generation of analytic semigroups by strongly elliptic
  operators}, Trans. Amer. Math. Soc. \textbf{199} (1974), 141--162.
  \MR{0358067 (50 \#10532)}

\bibitem{Stewart_1980}
\bysame, \emph{Generation of analytic semigroups by strongly elliptic operators
  under general boundary conditions}, Trans. Amer. Math. Soc. \textbf{259}
  (1980), 299--310. \MR{561838 (82h:35048)}

\bibitem{Streets_Tian_2010}
J.~Streets and G.~Tian, \emph{A parabolic flow of pluriclosed metrics}, Int.
  Math. Res. Not. IMRN (2010), 3101--3133. \MR{2673720 (2011h:53091)}

\bibitem{Struwe_1984}
M.~Struwe, \emph{A global compactness result for elliptic boundary value
  problems involving limiting nonlinearities}, Math. Z. \textbf{187} (1984),
  511--517. \MR{760051 (86k:35046)}

\bibitem{Struwe_1985}
\bysame, \emph{On the evolution of harmonic mappings of {R}iemannian surfaces},
  Comment. Math. Helv. \textbf{60} (1985), 558--581. \MR{826871 (87e:58056)}

\bibitem{Struwe_1994}
\bysame, \emph{The {Y}ang-{M}ills flow in four dimensions}, Calc. Var. Partial
  Differential Equations \textbf{2} (1994), 123--150. \MR{1385523 (97i:58034)}

\bibitem{Struwe_1996}
\bysame, \emph{Geometric evolution problems}, Nonlinear partial differential
  equations in differential geometry ({P}ark {C}ity, {UT}, 1992), IAS/Park City
  Math. Ser., vol.~2, Amer. Math. Soc., Providence, RI, 1996, pp.~257--339.
  \MR{1369591 (97e:58057)}

\bibitem{Struwe_variational_methods}
\bysame, \emph{Variational methods}, fourth ed., Springer, Berlin, 2008.
  \MR{2431434 (2009g:49002)}

\bibitem{Sun_Wang_2015}
S.~Sun and Y.~Wang, \emph{On the {K}\"ahler-{R}icci flow near a
  {K}\"ahler-{E}instein metric}, J. Reine Angew. Math. \textbf{699} (2015),
  143--158. \MR{3305923}

\bibitem{Takac_2000}
P.~Tak{\'a}{\v{c}}, \emph{Stabilization of positive solutions for analytic
  gradient-like systems}, Discrete Contin. Dynam. Systems \textbf{6} (2000),
  947--973. \MR{1788263 (2001i:35162)}

\bibitem{Tanabe_1979}
H.~Tanabe, \emph{Equations of evolution}, Monographs and Studies in
  Mathematics, vol.~6, Pitman (Advanced Publishing Program), Boston, MA, 1979,
  Translated from the Japanese by N. Mugibayashi and H. Haneda. \MR{533824
  (82g:47032)}

\bibitem{Tanabe_1997}
\bysame, \emph{Functional analytic methods for partial differential equations},
  Monographs and Textbooks in Pure and Applied Mathematics, vol. 204, Marcel
  Dekker Inc., New York, 1997. \MR{1413304 (97i:35002)}

\bibitem{Tao_2008lecture1}
T.~Tao, \emph{Math 285{G} - {L}ecture 1: {F}lows on {R}iemannian manifolds},
  Internet, April 2008, \url{http://terrytao.wordpress.com}.

\bibitem{Tao_2008lecture8}
\bysame, \emph{Math 285{G} - {L}ecture 8: {R}icci flow as a gradient flow,
  log-{S}obolev inequalities, and {P}erelman entropy}, Internet, April 2008,
  \url{http://terrytao.wordpress.com}.

\bibitem{TauSelfDual}
C.~H. Taubes, \emph{Self-dual {Y}ang-{M}ills connections on non-self-dual
  {$4$}-manifolds}, J. Differential Geom. \textbf{17} (1982), 139--170.
  \MR{658473 (83i:53055)}

\bibitem{TauPath}
\bysame, \emph{Path-connected {Y}ang-{M}ills moduli spaces}, J. Differential
  Geom. \textbf{19} (1984), 337--392. \MR{755230 (85m:58049)}

\bibitem{TauIndef}
\bysame, \emph{Self-dual connections on {$4$}-manifolds with indefinite
  intersection matrix}, J. Differential Geom. \textbf{19} (1984), 517--560.
  \MR{755237 (86b:53025)}

\bibitem{TauFrame}
\bysame, \emph{A framework for {M}orse theory for the {Y}ang-{M}ills
  functional}, Invent. Math. \textbf{94} (1988), 327--402. \MR{958836
  (90a:58035)}

\bibitem{TauStable}
\bysame, \emph{The stable topology of self-dual moduli spaces}, J. Differential
  Geom. \textbf{29} (1989), 163--230. \MR{978084 (90f:58023)}

\bibitem{TauCasson}
\bysame, \emph{Casson's invariant and gauge theory}, J. Differential Geom.
  \textbf{31} (1990), 547--599. \MR{1037415 (91m:57025)}

\bibitem{TauConf}
\bysame, \emph{The existence of anti-self-dual conformal structures}, J.
  Differential Geom. \textbf{36} (1992), 163--253. \MR{1168984 (93j:53063)}

\bibitem{TauL2}
\bysame, \emph{{$L^2$} moduli spaces on 4-manifolds with cylindrical ends},
  Monographs in Geometry and Topology, I, International Press, Cambridge, MA,
  1993. \MR{1287854 (96b:58018)}

\bibitem{TauGluing}
\bysame, \emph{Metrics, connections, and gluing theorems}, Conf. Board Math.
  Sci., vol.~89, Amer. Math. Soc., Providence, RI, 1996.

\bibitem{Taylor_2002}
G.~Taylor, \emph{A topological existence proof for {${\rm
  SO}(n)$}-anti-self-dual connections}, Topology Appl. \textbf{118} (2002),
  275--291. \MR{1874551 (2002i:57043)}

\bibitem{Taylor_PDO_and_NLPDE}
M.~E. Taylor, \emph{Pseudodifferential operators and nonlinear {PDE}}, Progress
  in Mathematics, vol. 100, Birkh\"auser Boston Inc., Boston, MA, 1991.

\bibitem{Taylor_PDE1}
\bysame, \emph{Partial differential equations {I}. {B}asic theory}, second ed.,
  Applied Mathematical Sciences, vol. 115, Springer, New York, 2011.
  \MR{2744150 (2011m:35001)}

\bibitem{Taylor_PDE3}
\bysame, \emph{Partial differential equations {III}. {N}onlinear equations},
  second ed., Applied Mathematical Sciences, vol. 117, Springer, New York,
  2011. \MR{2744149 (2011m:35003)}

\bibitem{Temam_1982}
R.~Temam, \emph{Behaviour at time {$t=0$} of the solutions of semilinear
  evolution equations}, J. Differential Equations \textbf{43} (1982), 73--92.
  \MR{645638 (83c:35058)}

\bibitem{ToppingThesis}
P.~Topping, \emph{The harmonic map heat flow from surfaces}, {Ph.D.} thesis,
  University of Warwick, United Kingdom, April 1996.

\bibitem{Topping_1997}
\bysame, \emph{Rigidity in the harmonic map heat flow}, J. Differential Geom.
  \textbf{45} (1997), 593--610. \MR{1472890 (99d:58050)}

\bibitem{Topping_2000}
\bysame, \emph{An example of a nontrivial bubble tree in the harmonic map heat
  flow}, Harmonic morphisms, harmonic maps, and related topics ({B}rest, 1997),
  Chapman \& Hall/CRC Res. Notes Math., vol. 413, Chapman \& Hall/CRC, Boca
  Raton, FL, 2000, pp.~185--191. \MR{1735698 (2001a:58019)}

\bibitem{Topping_2004am}
\bysame, \emph{Repulsion and quantization in almost-harmonic maps, and
  asymptotics of the harmonic map flow}, Ann. of Math. (2) \textbf{159} (2004),
  465--534. \MR{2081434 (2005g:58029)}

\bibitem{Topping_2004mz}
\bysame, \emph{Winding behaviour of finite-time singularities of the harmonic
  map heat flow}, Math. Z. \textbf{247} (2004), 279--302. \MR{2064053
  (2004m:53120)}

\bibitem{Topping_Lectures_ricci_flow}
\bysame, \emph{Lectures on the {R}icci flow}, London Mathematical Society
  Lecture Note Series, vol. 325, Cambridge University Press, Cambridge, 2006.
  \MR{2265040 (2007h:53105)}

\bibitem{Triebel_interpolation_theory_1978}
H.~Triebel, \emph{Interpolation theory, function spaces, differential
  operators}, North-Holland Mathematical Library, vol.~18, North-Holland
  Publishing Co., Amsterdam-New York, 1978. \MR{503903 (80i:46032b)}

\bibitem{Troianiello}
G.~M. Troianiello, \emph{Elliptic differential equations and obstacle
  problems}, Plenum Press, New York, 1987.

\bibitem{Trudinger_1968}
N.~S. Trudinger, \emph{Remarks concerning the conformal deformation of
  {R}iemannian structures on compact manifolds}, Ann. Scuola Norm. Sup. Pisa
  (3) \textbf{22} (1968), 265--274. \MR{0240748 (39 \#2093)}

\bibitem{UhlLp}
K.~K. Uhlenbeck, \emph{Connections with {$L^{p}$} bounds on curvature}, Comm.
  Math. Phys. \textbf{83} (1982), 31--42. \MR{648356 (83e:53035)}

\bibitem{UhlRem}
\bysame, \emph{Removable singularities in {Y}ang-{M}ills fields}, Comm. Math.
  Phys. \textbf{83} (1982), 11--29. \MR{648355 (83e:53034)}

\bibitem{UhlChern}
\bysame, \emph{The {C}hern classes of {S}obolev connections}, Comm. Math. Phys.
  \textbf{101} (1985), 449--457. \MR{815194 (87f:58028)}

\bibitem{Uhlenbeck_Yau_1986}
K.~K. Uhlenbeck and S.-T. Yau, \emph{On the existence of
  {H}ermitian-{Y}ang-{M}ills connections in stable vector bundles}, Comm. Pure
  Appl. Math. \textbf{39} (1986), no.~S, suppl., S257--S293, Frontiers of the
  mathematical sciences: 1985 (New York, 1985). \MR{861491 (88i:58154)}

\bibitem{Uhlenbeck_Yau_1989}
\bysame, \emph{A note on our previous paper: ``{O}n the existence of
  {H}ermitian-{Y}ang-{M}ills connections in stable vector bundles'' [{C}omm.\
  {P}ure {A}ppl.\ {M}ath.\ {\bf 39} (1986), {S}257--{S}293; {MR}0861491
  (88i:58154)]}, Comm. Pure Appl. Math. \textbf{42} (1989), 703--707.
  \MR{997570 (90i:58029)}

\bibitem{VanderHout_2003}
R.~van~der Hout, \emph{On the nonexistence of finite time bubble trees in
  symmetric harmonic map heat flows from the disk to the 2-sphere}, J.
  Differential Equations \textbf{192} (2003), 188--201. \MR{1987090
  (2004d:53084)}

\bibitem{Vespri_1991}
V.~Vespri, \emph{Analytic semigroups generated by ultraweak operators}, Proc.
  Roy. Soc. Edinburgh Sect. A \textbf{119} (1991), 87--105. \MR{1130598
  (92j:47074)}

\bibitem{Volpert1}
V.~Volpert, \emph{Elliptic partial differential equations. {V}olume 1:
  {F}redholm theory of elliptic problems in unbounded domains}, Monographs in
  Mathematics, vol. 101, Birkh\"auser/Springer Basel AG, Basel, 2011.
  \MR{2778694 (2012g:35003)}

\bibitem{Vrabie_2003}
I.~I. Vrabie, \emph{{$C_0$}-semigroups and applications}, North-Holland
  Mathematics Studies, vol. 191, North-Holland Publishing Co., Amsterdam, 2003.
  \MR{1972224 (2004c:47088)}

\bibitem{Waldron_2014arxiv}
A.~Waldron, \emph{Instantons and singularities in the {Y}ang-{M}ills flow},
  arXiv:1402.3224.

\bibitem{WaldronThesis}
\bysame, \emph{Self-duality and singularities in the {Y}ang-{M}ills flow},
  {Ph.D}. thesis, Columbia University, New York, NY, 2014. \MR{3232377}

\bibitem{Wang_2008sdg}
M.-T. Wang, \emph{Some recent developments in {L}agrangian mean curvature
  flows}, Surveys in differential geometry. {V}ol. {XII}. {G}eometric flows,
  Surv. Differ. Geom., vol.~12, Int. Press, Somerville, MA, 2008,
  arXiv:1104.3355, pp.~333--347. \MR{2488942 (2009m:53181)}

\bibitem{Wang_Zhang_2008}
Y.~Wang and X.~Zhang, \emph{The coupled {Y}ang-{M}ills-{H}iggs flow}, J. Math.
  Anal. Appl. \textbf{339} (2008), 153--174. \MR{2370640 (2008k:53142)}

\bibitem{Warner}
F.~W. Warner, \emph{Foundations of differentiable manifolds and {L}ie groups},
  Graduate Texts in Mathematics, vol.~94, Springer, New York, 1983. \MR{722297
  (84k:58001)}

\bibitem{Weinstein_1987}
M.~I. Weinstein, \emph{Solitary waves of nonlinear dispersive evolution
  equations with critical power nonlinearities}, J. Differential Equations
  \textbf{69} (1987), 192--203. \MR{899159 (89d:35172)}

\bibitem{Wells3}
R.~O. Wells, Jr., \emph{Differential analysis on complex manifolds}, third ed.,
  Graduate Texts in Mathematics, vol.~65, Springer, New York, 2008, with a new
  appendix by Oscar Garcia-Prada. \MR{2359489 (2008g:32001)}

\bibitem{Wilkin_2008}
G.~Wilkin, \emph{Morse theory for the space of {H}iggs bundles}, Comm. Anal.
  Geom. \textbf{16} (2008), 283--332. \MR{2425469 (2010f:53115)}

\bibitem{Wu_Xu_2013}
H.~Wu and X.~Xu, \emph{Strong solutions, global regularity, and stability of a
  hydrodynamic system modeling vesicle and fluid interactions}, SIAM J. Math.
  Anal. \textbf{45} (2013), no.~1, 181--214. \MR{3032974}

\bibitem{Wu_1988}
H.~H. Wu, \emph{The {B}ochner technique in differential geometry}, Math. Rep.
  \textbf{3} (1988), no.~2, i--xii and 289--538. \MR{1079031 (91h:58031)}

\bibitem{Yagi_abstract_parabolic_evolution_equations_applications}
A.~Yagi, \emph{Abstract parabolic evolution equations and their applications},
  Springer Monographs in Mathematics, Springer-Verlag, Berlin, 2010.
  \MR{2573296 (2011c:35008)}

\bibitem{Yamabe_1960}
H.~Yamabe, \emph{On a deformation of {R}iemannian structures on compact
  manifolds}, Osaka Math. J. \textbf{12} (1960), 21--37. \MR{0125546 (23
  \#A2847)}

\bibitem{Yang_2003aim}
B.~Yang, \emph{The uniqueness of tangent cones for {Y}ang-{M}ills connections
  with isolated singularities}, Adv. Math. \textbf{180} (2003), 648--691.
  \MR{2020554 (2004m:58026)}

\bibitem{Yano_1965}
K.~Yano, \emph{Differential geometry on complex and almost complex spaces},
  International Series of Monographs in Pure and Applied Mathematics, Vol. 49,
  A Pergamon Press Book. The Macmillan Co., New York, 1965. \MR{0187181 (32
  \#4635)}

\bibitem{Ye_1994jdg}
R.~Ye, \emph{Global existence and convergence of {Y}amabe flow}, J.
  Differential Geom. \textbf{39} (1994), 35--50. \MR{1258912 (95d:53044)}

\bibitem{Yosida}
K.~Yosida, \emph{Functional analysis}, sixth ed., Springer, New York, 1980.

\bibitem{Zaal_2015}
M.~M. Zaal, \emph{The gradient structure of the mean curvature flow}, Adv.
  Calc. Var. \textbf{8} (2015), no.~3, 183--202. \MR{3365740}

\bibitem{Zeidler_nfaa_v3}
E.~Zeidler, \emph{Nonlinear functional analysis and its applications. {III}},
  Springer-Verlag, New York, 1985, Variational methods and optimization.
  \MR{768749 (90b:49005)}

\bibitem{Zeidler_nfaa_v1}
\bysame, \emph{Nonlinear functional analysis and its applications, {I}.
  {F}ixed-point theorems}, Springer, New York, 1986. \MR{816732 (87f:47083)}

\bibitem{Zhang_2012}
Y.~Zhang, \emph{{$\mathscr{F}$}-stability of self-similar solutions to harmonic
  map heat flow}, Calc. Var. Partial Differential Equations \textbf{45} (2012),
  347--366. \MR{2984136}

\bibitem{Zheng_2007}
Y.~Zheng, \emph{On the study of one flow for {ASD} connection}, Commun.
  Contemp. Math. \textbf{9} (2007), 545--569. \MR{2348843 (2008j:53119)}

\end{thebibliography}
\bibliographystyle{amsplain}


\end{document}